\documentclass[a4paper, 11pt, notitlepage]{report}

\usepackage[utf8]{inputenc}
\usepackage[T1]{fontenc}

\usepackage[xindy]{imakeidx}

\usepackage[
	backend = biber,
	hyperref = auto,
	backref = true,
	doi = false,
	url = false
]{biblatex}

\usepackage{amsthm}
\usepackage{amssymb}
\usepackage{mathtools}
\usepackage{stmaryrd}	\SetSymbolFont{stmry}{bold}{U}{stmry}{m}{n}	% To avoid warnings
\usepackage{bm}
\usepackage{tikz-cd}
\usetikzlibrary{positioning, shapes.geometric, patterns, graphs, decorations.pathmorphing}

\usepackage{paralist}

\usepackage[colorlinks]{hyperref}
\hypersetup{
	linkcolor=blue,
	citecolor=teal!70!lime,
	urlcolor=cyan,
	pdftitle = {Stabilization of the trace formula for metaplectic groups},
	pdfauthor = {Wen-Wei Li}
}

\usepackage{nccrules}
\usepackage{mathrsfs}

\newcommand{\Z}{\ensuremath{\mathbb{Z}}}
\newcommand{\Q}{\ensuremath{\mathbb{Q}}}
\newcommand{\R}{\ensuremath{\mathbb{R}}}
\newcommand{\CC}{\ensuremath{\mathbb{C}}}
\newcommand{\A}{\ensuremath{\mathbb{A}}}

\newcommand{\gr}{\operatorname{gr}}
\newcommand{\Sym}{\operatorname{Sym}}
\newcommand{\Aut}{\operatorname{Aut}}
\newcommand{\Tr}{\operatorname{tr}}

\newcommand{\Ind}{\operatorname{Ind}}

\newcommand{\Gal}{\operatorname{Gal}}

\newcommand{\Weil}[1]{\ensuremath{\mathrm{W}_{#1}}}	% The Weil group
\newcommand{\WD}[1]{\ensuremath{\mathrm{WD}_{#1}}}	% The Weil-Deligne group
\newcommand{\Resprod}{\ensuremath{{\prod}'}}	% Restricted product

\newcommand{\dd}{\mathop{}\!\mathrm{d}}
\newcommand{\lrangle}[1]{\ensuremath{\langle #1 \rangle}}
\newcommand{\sgn}{\operatorname{sgn}}
\newcommand{\Stab}{\operatorname{Stab}}
\newcommand{\mes}{\operatorname{mes}}

\newcommand{\identity}{\ensuremath{\mathrm{id}}}

\newcommand{\Hom}{\operatorname{Hom}}

\newcommand{\rightiso}{\ensuremath{\stackrel{\sim}{\rightarrow}}}

\newcommand{\Ker}{\operatorname{ker}}
\newcommand{\Coker}{\operatorname{coker}}
\newcommand{\Image}{\operatorname{im}}
\newcommand{\Hm}{\operatorname{H}}
\newcommand{\dotimes}[1]{\ensuremath{\underset{#1}{\otimes}}}

\newcommand{\Lie}{\operatorname{Lie}}
\newcommand{\Ad}{\operatorname{Ad}}
\newcommand{\ad}{\operatorname{ad}}
\newcommand{\Spec}{\operatorname{Spec}}
\newcommand{\Gm}{\ensuremath{\mathbb{G}_\mathrm{m}}}
\newcommand{\Ga}{\ensuremath{\mathbb{G}_\mathrm{a}}}
\newcommand{\Res}{\operatorname{Res}}
\newcommand{\utimes}[1]{\ensuremath{\overset{#1}{\times}}}
\newcommand{\dtimes}[1]{\ensuremath{\underset{#1}{\times}}}
\newcommand{\Supp}{\operatorname{Supp}}

\newcommand{\GL}{\operatorname{GL}}
\newcommand{\SO}{\operatorname{SO}}
\newcommand{\so}{\ensuremath{\mathfrak{so}}}
\newcommand{\Spin}{\operatorname{Spin}}

\newcommand{\SL}{\operatorname{SL}}
\newcommand{\Sp}{\operatorname{Sp}}
\newcommand{\GSp}{\operatorname{GSp}}
\newcommand{\syp}{\ensuremath{\mathfrak{sp}}}
\newcommand{\Mp}{\ensuremath{\widetilde{\mathrm{Sp}}}}
\newcommand{\MMp}{\ensuremath{\overline{\mathrm{Sp}}}}
\newcommand{\bmu}{\ensuremath{\bm\mu}}
\newcommand{\Lgrp}[1]{\ensuremath{{}^{\mathrm{L}} #1}}	% The L-group

\newcommand{\Endo}{\ensuremath{\mathcal{E}}}
\newcommand{\orbI}{\ensuremath{\mathcal{I}}}
\newcommand{\cusp}{\operatorname{cusp}}
\newcommand{\elli}{\operatorname{ell}}
\newcommand{\desc}{\operatorname{desc}}	% Descente de Harish-Chandra
\newcommand{\asp}{\ensuremath{\dashrule[.7ex]{2 2 2 2}{.4}}} % Le symbole pour anti-spécifiques
\newcommand{\rev}{\ensuremath{\mathbf{p}}} % Le symbole pour revêtements
\newcommand{\Trans}{\ensuremath{\mathcal{T}}}	% Transfert géométrique
\newcommand{\trans}{\ensuremath{\check{\mathcal{T}}}}	% Transfert de distributions

\theoremstyle{plain}
\newtheorem{proposition}{Proposition}
\newtheorem{lemma}[proposition]{Lemma}
\newtheorem{theorem}[proposition]{Theorem}
\newtheorem{corollary}[proposition]{Corollary}
\theoremstyle{definition}
\newtheorem{definition}[proposition]{Definition}
\newtheorem{definition-theorem}[proposition]{Definition--Theorem}
\newtheorem{definition-proposition}[proposition]{Definition--Proposition}
\newtheorem{remark}[proposition]{Remark}
\newtheorem{hypothesis}[proposition]{Hypothesis}
\newtheorem{notation}[proposition]{Notation}

\newtheorem{example}[proposition]{Example}

\theoremstyle{definition}

\theoremstyle{plain}

\numberwithin{equation}{section}
\numberwithin{proposition}{section}
\numberwithin{conj}{section}	% In the Introduction only

\renewcommand{\Re}{\operatorname{Re}}	% Symbol for the real and imaginary parts
\renewcommand{\Im}{\operatorname{Im}}
\renewcommand{\emptyset}{\ensuremath{\varnothing}}	% Symbol for the emptyset

\usepackage{amsmath}
\usepackage[nottoc]{tocbibind}

\usepackage{geometry}
\usepackage{lmodern}

\title{Stabilization of the trace formula for metaplectic groups}
\author{Wen-Wei Li}
\date{}
\DeclareFieldFormat{postnote}{#1}	% Disable pagination in postnote for biblatex
\makeindex

\begin{document}

\maketitle

\begin{abstract}
	We stabilize the full Arthur--Selberg trace formula for the metaplectic covering of symplectic groups over a number field. This provides a decomposition of the invariant trace formula for metaplectic groups, which encodes information about the genuine $L^2$-automorphic spectrum, into a linear combination of stable trace formulas of products of split odd orthogonal groups via endoscopic transfer. By adapting the strategies of Arthur and Mo\!eglin--Waldspurger from the linear case, the proof is built on a long induction process that mixes up local and global, geometric and spectral data. As a by-product, we also stabilize the local trace formula for metaplectic groups over any local field of characteristic zero.
\end{abstract}

\vspace{1em}
{\scriptsize
\begin{tabular}{ll}
	\textbf{MSC (2010)} & \textbf{11F72}, 11F27, 11F70 \\
	\textbf{Keywords} & Arthur-Selberg trace formula, stable trace formula, metaplectic group
\end{tabular}}

\vfill

{\footnotesize\begin{tabular}{l}
	Wen-Wei Li \\
	School of Mathematical Sciences / BICMR, Peking University \\
	5 Yiheyuan Road, Beijing 100871, People's Republic of China \\
	E-mail: \texttt{wwli@bicmr.pku.edu.cn}
\end{tabular}}

% This work is supported by NSFC-11922101.

\tableofcontents

\chapter{Introduction}
Let $F$ be a field of characteristic zero, and let $\Sp(W)$ be the symplectic group associated to a symplectic $F$-vector space $(W, \lrangle{\cdot|\cdot})$. When a symplectic basis is chosen, it is common to write $\Sp(2n)$ instead of $\Sp(W)$, with $2n = \dim_F W$.

Suppose that $F$ is a number field. Upon choosing a non-trivial additive character $\psi$ of $\A_F / F$ and denoting by $\bmu_m$ the group of $m$-th roots of unity in $\CC^\times$, there is a canonical central extension of locally compact groups
\[ 1 \to \bmu_2 \to \Mp^{(2)}(W, \A_F) \to \Sp(W, \A_F) \to 1 \]
which splits canonically over the $F$-points of $\Sp(W)$. This is the well-known metaplectic covering introduced by A.\ Weil \cite{Weil64}. When $F$ is a local field and $\psi$ is a chosen additive character of $F$, we also have the corresponding covering $\Mp^{(2)}(W)$ of $\Sp(W)$, where we identify $\Sp(W)$ with its group of $F$-points. These coverings are nontrivial except when $F = \CC$.

We may consider the automorphic forms and automorphic representations of $\Mp^{(2)}(W, \A_F)$ which do not factor through $\Sp(W, \A_F)$. They are naturally tied up with the study of $\theta$-series and $\theta$-lifting. These coverings are instances of the Brylinski--Deligne coverings of connected reductive groups, hence provide a nontrivial testing ground for the Langlands--Weissman program \cite{GG18} for such covering groups, in both local and global aspects.

To be more precise, what we study in this work are the \emph{eightfold coverings}
\[ 1 \to \bmu_8 \to \Mp(W, \A_F) \xrightarrow{\rev} \Sp(W, \A_F) \to 1 \]
in the global case, and
\[ 1 \to \bmu_8 \to \Mp(W) \xrightarrow{\rev} \Sp(W) \to 1 \]
in the local case, defined as push-outs via the inclusion $\bmu_2 \hookrightarrow \bmu_8$; they will be abbreviated as \emph{metaplectic groups} in what follows. The advantages of the eightfold way will be clarified soon.

In this setting, we are led to study automorphic forms (resp.\ representations) on $\Mp(W, \A_F)$ (resp.\ of $\Mp(W, \A_F)$ or its local avatars) which are \emph{genuine} in the sense that $\bmu_8$ acts as $z \mapsto z \cdot \identity$. As for various distributions that arise in this context, such as the character-distribution $\Theta_\pi$ attached to a genuine admissible representation $\pi$, what matter are the \emph{anti-genuine} test functions $f$ in the sense that $f(z\tilde{x}) = z^{-1} f(\tilde{x})$ for all $z \in \bmu_8$.

Let us denote $G := \Sp(W)$, and $\tilde{G} := \Mp(W, \A_F)$ (resp.\ $\tilde{G} := \Mp(W)$) when $F$ is global (resp.\ local). The goal of this work is to stabilize the \emph{Arthur--Selberg trace formula} for $\tilde{G}$ in its entirety (Theorems \ref{prop:STF-intro}, \ref{prop:STF-disc-intro}), akin to the case of connected reductive groups established in the monumental works \cite{Ar02, Ar01, Ar03-3} of Arthur. There are at least two motivations.
\begin{itemize}
	\item It is the first fully stabilized trace formula in the Langlands--Weissman program which involves a nontrivial theory of endoscopy.
	\item It furnishes a vehicle capable of accessing the genuine automorphic spectrum $L^2_-(G(F) \backslash \tilde{G})$, as exemplified by \cite{Ar13} in the case of classical groups. The stabilization of the elliptic terms of the geometric side in \cite{Li15} does not suffice for this purpose.

	For instance, one might hope to strengthen the multiplicity formula of Gan--Ichino \cite{GI18, GI21} using the stabilized full trace formula, but this is surely beyond the scope of the present work.
\end{itemize}

The stabilization is based on the earlier works \cite{Li11, Li12a, Li15, Li19, Luo18} on the endoscopy for metaplectic groups. Surely, the techniques from Arthur and Mo\!eglin--Waldspurger \cite{MW16-1, MW16-2} are also indispensable in this work. A detailed introduction goes as follows.

\section{Invariant trace formula for metaplectic groups}
Consider a global metaplectic covering $\rev: \tilde{G} \to G(\A_F)$. Fix a symplectic basis for $(W, \lrangle{\cdot|\cdot})$ and the corresponding minimal Levi subgroup $M_0$ of $G$. Let $L$ be the lattice in $W$ generated by that basis over the ring of integers $\mathfrak{o}_F$. Let $V_{\mathrm{ram}} = V_{\mathrm{ram}}(\tilde{G})$ be the finite set consisting of all Archimedean places, the dyadic places and the non-Archimedean places $v$ such that $L_v$ is not self-dual with respect to $\psi_v \circ \lrangle{\cdot|\cdot}$. 

For every subset $E \subset G(\A_F)$, set $\tilde{E} := \rev^{-1}(E) \subset \tilde{G}$.

Let $V$ be a finite set of places such that $V \supset V_{\mathrm{ram}}$. Fix maximal compact subgroups $K_v \subset G(F_v)$ at each $v \in V$ in good positions relative to $M_0$. For $v \notin V$, let $K_v$ be the hyperspecial subgroup associated with $L_v$. Let $K_V := \prod_{v \in V} K_v$. Define $\tilde{G}_V := \rev^{-1}(G(F_V))$ and $\tilde{K}_V$ by the earlier convention, and define $\tilde{G}_v$, $\tilde{K}_v$ similarly. It is known that when $v \notin V$,
\begin{itemize}
	\item the covering $\rev_v: \tilde{G}_v \to G(F_v)$ splits canonically over $K_v$, and
	\item there is a Satake isomorphism for the anti-genuine spherical Hecke algebra for $K_v \backslash \tilde{G}_v / K_v$.
\end{itemize}

Fix Haar measures. Following Arthur's paradigm, the \emph{invariant trace formula} for coverings \cite{Li14a, Li12b, Li13, Li14b} is an equality between two genuine invariant distributions on $\tilde{G}_V$, called the geometric and spectral sides of the trace formula:
\[ I^{\tilde{G}}_{\mathrm{geom}}(f) = I^{\tilde{G}}_{\mathrm{spec}}(f) \]
where $f$ is an anti-genuine, $C^\infty_c$ and $\tilde{K}_V \times \tilde{K}_V$-finite function on $\tilde{G}_V$. Oftentimes the superscripts $\tilde{G}$ will be omitted. We also write $\dot{f} := f \prod_{v \notin V} f_{K_v}$ where $f_{K_v}$ is the unit in the anti-genuine spherical Hecke algebra of $\tilde{G}_v$.

Let $\mathcal{L}(M_0)$ denote the set of semi-standard Levi subgroups with respect to $M_0$. For every $M \in \mathcal{L}(M_0)$, let $\mathcal{P}(M)$ be the set of semi-standard parabolic subgroups with Levi factor $M$.

\begin{description}
	\item[Core of the spectral side] We have a decomposition $I_{\mathrm{spec}} = \sum_{t \geq 0} I_t$; in turn, each $I_t$ decomposes into $\sum_\nu I_\nu$ where $\nu$ stands for the infinitesimal characters at Archimedean places, subject to $\|\Im(\nu)\| = t$ where $\|\cdot\|$ is some chosen height function. The right hand side of
	\[ I_{\mathrm{spec}} = \sum_{t \geq 0} \sum_{\|\Im(\nu)\|=t} I_\nu \]
	converges as an iterated sum. The absolute convergence would follow from a metaplectic analogue of the result of Finis--Lapid--Müller \cite{FLM11}, but we do not make this assumption.
	
	Each $I_t$ admits a further decomposition, and there is a core part encoding automorphic information, namely the $t$-discrete part
	\begin{align*}
		I_{\mathrm{disc}, t}(\dot{f}) = I^{\tilde{G}}_{\mathrm{disc}, t}(\dot{f}) & := \sum_{M \in \mathcal{L}(M_0)} \sum_{s \in W^G(M)_{\mathrm{reg}}} \frac{|W^M_0|}{|W^G_0|} \\
		& \quad \left| \det(1-s | \mathfrak{a}^G_M) \right|^{-1} \Tr\left( M_{P|P}(s, 0) \mathcal{I}_{\tilde{P}, \mathrm{disc}, t}(0, \dot{f}) \right),
	\end{align*}
	where
	\begin{itemize}
		\item $W^G_0$ (resp.\ $W^M_0$) is the Weyl group of $G$ (resp.\ $M$) relative to $M_0$;
		\item $W^G(M) := \Stab_{W^G_0}(M) / W^M_0$, and $W^G(M)_{\mathrm{reg}}$ is its subset consisting of elements $s$ such that $s-1$ is invertible on the $\R$-vector space $\mathfrak{a}^G_M$ omnipresent in Arthur's theory (see \S\ref{sec:combinatorics-reductive}, for example);
		\item $P \in \mathcal{P}(M)$ is arbitrary and $\mathcal{I}_{\tilde{P}, \mathrm{disc}, t}(0, \cdot)$ is the normalized parabolic induction to $\tilde{G}$ of $L^2_{-, t}(M(F) \backslash \tilde{M}^1)$, the $t$-discrete part of $L^2(M(F) \backslash \tilde{M}^1)$, with $\tilde{M}^1$ being the kernel of Harish-Chandra's homomorphism for $\tilde{M}$ (see below);
		\item $M_{P|P}(s, 0)$ is the standard intertwining operator.
	\end{itemize}
	We refer to \cite[(5.10)]{Li14b} for a complete explanation of all terms. We shall view $I_{\mathrm{disc}, t}$ as a formal linear combination of irreducible genuine characters of $\tilde{G}$.
	
	The terms with $M=G$ in $I_{\mathrm{disc}, t}$ are simply the contribution of the $t$-part of $L^2_{-, t}(G(F) \backslash \tilde{G})$, but there are also some ``shadows'' from the discrete spectrum of $\tilde{M} \neq \tilde{G}$.
	
	The discussions above can be refined to $I_{\mathrm{disc}, \nu}$, when $\nu$ is a specified infinitesimal character. They give rise to distributions on $\tilde{G}_V$ through $f \mapsto \dot{f}$. Furthermore, $I_{\mathrm{disc}, \nu} = \sum_{c^V} I_{\mathrm{disc}, \nu, c^V}$ as distributions on $\tilde{G}_V$, where $c^V$ ranges over Satake parameters off $V$.
		
	\item[Core of the geometric side] Assume temporarily that $f$ has \emph{admissible support}; given $f$, this can be achieved by replacing $V$ by a sufficiently large finite set $S \supset V$, and replacing $f$ by $f_S := f \prod_{v \in S \smallsetminus V} f_{K_v}$ accordingly. Then $I_{\mathrm{geom}}(f)$ contains the elliptic part
	\begin{equation*}
		I_{\elli}(f) = I^{\tilde{G}}_{\elli}(f) := \sum_{\mathcal{O}} I^{\tilde{G}_V}\left(A(V, \mathcal{O})_{\elli}, f \right)
	\end{equation*}
	where
	\begin{itemize}
		\item $\mathcal{O}$ ranges over finitely many elliptic semisimple conjugacy classes in $G(F)$ (see Definition \ref{def:elliptic-element}), in a way that depends only on $\Supp(f)$;
		\item $A(V, \mathcal{O})_{\elli}$ is a linear combination of conjugacy classes in $\tilde{G}_V$ extracted from $G(F)$ in a suitable manner, whose semisimple part comes from $\mathcal{O}$;
		\item $I^{\tilde{G}_V}(\cdot, f)$ denotes the orbital integral of $f$, evaluated on $A(V, \mathcal{O})_{\elli}$ by linearity.
	\end{itemize}
	The classes $\mathcal{O}$ which are regular semisimple yield the sum of elliptic regular orbital integrals weighted by Tamagawa numbers, if the orbital integrals are defined by Tamagawa measures. However, $I_{\elli}$ also contains contributions from other types of conjugacy classes.
\end{description}

The ``cores'' above are just pieces of the part with $M=G$ in the general expansions of $I_{\mathrm{spec}}$ and $I_{\mathrm{geom}}$, explained briefly as follows.

\begin{description}
	\item[Full spectral side] For all $t \geq 0$ and test functions $f$, we have
	\begin{equation*}
		I_t \left(f \right) = \sum_{M \in \mathcal{L}(M_0)} \frac{|W^M_0|}{|W^G_0|} I^{\mathrm{glob}}_{\tilde{M}}\left( I^{\tilde{M}}_{\mathrm{disc}, t}, f \right),
	\end{equation*}
	where $I^{\mathrm{glob}}_{\tilde{M}}(\cdot, f)$ stands for the (invariant) global weighted character, to be introduced in \S\ref{sec:spec-nr}. The global weighted character is expressed in terms of
	\begin{itemize}
		\item certain factors $r^{\tilde{L}}_{\tilde{M}}(c^V)$ involving Satake parameters $c^V$ off $V$ which are automorphic, and unramified normalizing factors for standard intertwining operators, where $M \subset L \subset G$,
		\item semi-local invariant weighted characters $I^{\tilde{G}_V}_{\tilde{L}_V}(\pi_V, \lambda, X, f)$ --- ``semi-local'' as it involves a finite set $V$ of places of $F$ --- where $\lambda \in \mathfrak{a}^*_L$ governs a shift of contours and $X \in \mathfrak{a}_L$ (see \S\ref{sec:weighted-characters}).
	\end{itemize}
	The evaluation of $I^{\mathrm{glob}}_{\tilde{M}}(\cdot, f)$ on the formal linear combination $I^{\tilde{M}}_{\mathrm{disc}, t}$ is well-defined. Note that the formulation above is slightly different from that in \cite{Li14b}.
	
	Similarly, when an infinitesimal character $\nu$ is specified, there is a refined version expressing $I_\nu(f)$ in terms of $I^{\tilde{M}}_{\mathrm{disc}, \nu}$ for various $M$.
	
	\item[Full geometric side] For all test functions $f$, we have
	\begin{equation*}
		I_{\mathrm{geom}}(f) = \sum_{M \in \mathcal{L}(M_0)} \frac{|W^M_0|}{|W^G_0|} \sum_{\substack{\mathcal{O}^M \in M(F) /\text{conj} \\ \text{semisimple}}} I_{\tilde{M}_V}\left(A^{\tilde{M}}(V, \mathcal{O}^M), f\right),
	\end{equation*}
	where $I_{\tilde{M}_V}(\cdot, f)$ stands for the semi-local invariant weighted orbital integral.
\end{description}

The terms with $M=G$ are defined from the aforementioned cores $I^{\tilde{M}}_{\mathrm{disc}, t}$ and $A^{\tilde{M}}(S, \mathcal{O}^M)_{\elli}$, for various $M \in \mathcal{L}(M_0)$ and $S \supset V$, by a procedure called \emph{compression of coefficients}; we refer to \cite[\S 5]{Li14b} for an overview.

Without stepping into the details, we emphasize some key features of this formalism.
\begin{enumerate}
	\item In either case, we obtain a sum indexed by $M \in \mathcal{L}(M_0)$ of
	\begin{itemize}
		\item semi-local distributions,
		\item global coefficients, namely $I^{\tilde{M}}_{\mathrm{disc}, t}$, $I^{\tilde{M}}_{\mathrm{disc}, \nu}$ and $\sum_{\mathcal{O}^M} A^{\tilde{M}}(V, \mathcal{O}^M)$ --- the term ``coefficient'' is used in a loose sense, since they will be defined to be certain distributions on $\tilde{M}_V$.
	\end{itemize}
	The semi-local distributions are expressible in terms of their local avatars through splitting formulas.
	
	Typically, the global coefficients depend on the group $\tilde{G}$ or its Levi subgroups, recognizable as the superscript. The semi-local or local distributions depend on $\tilde{G}$ (superscript) as well as a Levi subgroup $\tilde{M}$ (subscript). Superscripts will often be omitted.
	
	\item It is clearly necessary to consider all the Levi subgroups, not just the metaplectic group $\tilde{G}$ itself. Up to conjugacy, the Levi subgroups $M$ of $G$ take the form
	\[ M = \prod_{i \in I} \GL(n_i) \times \Sp(W^\flat) \]
	where $I$ is a finite set, $n_i \in \Z_{\geq 1}$ for all $i$ and $W^\flat$ is a symplectic subspace of $W$. Thanks to the use of eightfold coverings and the Schrödinger model for Weil representations (which depends on the choice of $\psi$), the preimage $\tilde{M}$ of $M(\A_F)$ decomposes canonically into
	\[ \tilde{M} = \prod_{i \in I} \GL(n_i, \A_F) \times \Mp(W^\flat, \A_F). \]
	Ditto over local fields. Coverings of this form are said to be of \emph{metaplectic type}. It is justified to put $V_{\mathrm{ram}}(\tilde{M}) := V_{\mathrm{ram}}(\Mp(W^\flat, \A_F))$.
	
	This brings a recursive structure into our study of trace formulas, reducing everything to smaller metaplectic groups.
	
	\item According to the previous item, we can and must consider the invariant trace formula attached to groups of metaplectic type. This causes a minor issue on the choice of test functions. Let $\tilde{G}$ be of metaplectic type, and let $\tilde{G}^1$ be the kernel of Harish-Chandra's homomorphism $H_G: G(\A_F) \to \mathfrak{a}_G$ composed with $\rev: \tilde{G} \to G(\A_F)$. There is a central subgroup $A_{G, \infty}$ in $G(F_\infty)$ such that $H_G: A_{G, \infty} \rightiso \mathfrak{a}_G$ as Lie groups, and with $\tilde{G}_V^1 := \tilde{G}_V \cap \tilde{G}^1$ there are canonical decompositions
	\[ \tilde{G} = \tilde{G}^1 \times A_{G, \infty} , \quad \tilde{G}_V = \tilde{G}_V^1 \times A_{G, \infty}. \]
	The trace formula for $\tilde{G}$ involves only the representations (resp.\ conjugacy classes) that are trivial on $A_{G, \infty}$ (resp.\ lying in $\tilde{G}^1$). In \cite{Li14b}, we considered test functions $f^1$ on $\tilde{G}_V/A_{G, \infty}$. In this work, we switch to the viewpoint of \cite{Ar88-2, MW16-2} by using $C^\infty_c$-test functions $f$ on $\tilde{G}_V$, but the resulting distributions will depend only on $f|_{\tilde{G}^1}$.
	
	The transition from $f^1$ to $f$ has been sketched in \cite[Théorème 6.4]{Li14b}. It will be reviewed in \S\ref{sec:geom-side}.
\end{enumerate}

We also indicate a feature of groups of metaplectic type that distinguishes them from general Brylinski--Deligne coverings. Let $\tilde{G}$ be of metaplectic type, say over a local field $F$. Then all elements $\tilde{x}$ of $\tilde{G}$ are \emph{good} in the sense that $\tilde{x}\tilde{y} = \tilde{y}\tilde{x}$ if and only if their images in $G(F)$ satisfies $xy = yx$. In general, only good conjugacy classes can support invariant genuine distributions.

\section{Stabilization: desiderata}
Let $F$ be a number field. For a connected reductive $F$-group $L$, the goal of stabilization is to express both sides of the invariant trace formula $I^L_{\mathrm{geom}} = I^L_{\mathrm{spec}}$ in terms of the stable trace formulas $S^{L'}_{\mathrm{geom}} = S^{L'}_{\mathrm{spec}}$ for suitable quasisplit groups $L'$. The precise recipe hinges on the theory of endoscopy, and the groups $L'$ in question are called the elliptic endoscopic groups; they underlie the elliptic endoscopic data $\mathbf{L}'$ of $L$, taken up to equivalence. The test functions on $L$ and $L'$ are related by the transfer of orbital integrals. This program is accomplished in Arthur's works \cite{Ar02, Ar01, Ar03-3}.

We wish to do the same for $\tilde{G}$, in a manner such that the endoscopic groups will be quasisplit $F$-groups, not coverings.\footnote{Some authors prefer to take endoscopic groups to be coverings. We will say a few words about this in \S\ref{sec:endo-covering}.} The endoscopic groups (resp.\ data) will be denoted as $G^!$ (resp.\ $\mathbf{G}^!$).

Based on ideas of J.\ Adams \cite{Ad98} and D.\ Renard \cite{Re99} over $\R$, such a formalism of endoscopy is proposed in \cite{Li11}. The rough idea goes as follows. Let $\tilde{G} = \Mp(W, \A_F)$ (global case) or $\tilde{G} = \Mp(W)$ (local case); assume $\dim_F W = 2n$.
\begin{itemize}
	\item The dual group of $\tilde{G}$ is defined as the $\CC$-group $\tilde{G}^\vee := \Sp(2n, \CC)$ with trivial Galois action, and this is compatible with Satake isomorphisms. Note that $\Sp(2n, \CC)$ is also the dual of the split $\SO(2n+1)$.
	\item Using $\tilde{G}^\vee$, the elliptic endoscopic data are defined in the same manner as in Langlands--Shelstad \cite{LS87}, with one exception: we disregard the effect of $Z_{\tilde{G}^\vee} = \{\pm 1\}$ on endoscopic data.
	
	As a consequence, the subtle issue of automorphisms of endoscopic data disappears, and the elliptic endoscopic data turn out to be described by pairs $(n', n'') \in \Z_{\geq 0}^2$ such that $n' + n'' = n$. The pair $(n', n'')$ corresponds to an elliptic semisimple element $s \in \tilde{G}^\vee$ up to conjugacy: $2n'$ (resp.\ $2n''$) is the multiplicity of the eigenvalue $+1$ (resp.\ $-1$) of $s$.
	
	In contrast, the equivalence classes of elliptic endoscopic data of $\SO(2n+1)$ are in bijection with pairs $(n', n'')$ up to $(n', n'') \leftrightarrow (n'', n')$, the symmetry stemming from the effect of $Z_{\Sp(2n, \CC)} = \{\pm 1\}$.
	
	\item The endoscopic group $G^!$ corresponding to the datum $\mathbf{G}^! \leftrightarrow (n', n'')$ is defined as the split group
	\[ G^! = \SO(2n' + 1) \times \SO(2n'' + 1). \]
	The elliptic endoscopic groups of $\SO(2n+1)$ take the same form, but the ``endoscopic outer automorphisms'' in the theory for $\SO(2n+1)$ (i.e.\ swapping two factors when $n' = n''$) play no role here.
	
	In the local setting, the local Langlands correspondence is known for $G^!$; see \cite[Theorem 1.5.1]{Ar13}. In the global setting, the stable trace formula for $G^!$ is relatively well-understood, for example the decomposition of $S^{G^!}_{\mathrm{disc}}$ is completely described by the \emph{stable multiplicity formula} in \cite[Section 4.1]{Ar13}.
	
	\item Denote by $\Endo_{\elli}(\tilde{G})$ the set of elliptic endoscopic data of $\tilde{G}$, without worrying about automorphisms. For every $\mathbf{G}^! \in \Endo_{\elli}(\tilde{G})$, the correspondence between semisimple conjugacy classes of $G$ and $G^!$ and the transfer factors are defined in \cite{Li11}. When $F$ is local, the existence of local transfer $f \mapsto f^!$ and the fundamental lemma for units are established in \textit{loc.\ cit.} The fundamental lemma for the whole spherical anti-genuine Hecke algebra is obtained by C.\ Luo \cite{Luo18}. Let us identify $f^!_1, f^!_2 \in C^\infty_c(G^!(F))$ if they have the same stable orbital integrals along every regular semisimple orbit; modulo this equivalence relation, the transfer is a well-defined linear map, denoted as
	\[ f \mapsto f^! = \Trans_{\mathbf{G}^!, \tilde{G}}(f), \quad \mathbf{G}^! \in \Endo_{\elli}(\tilde{G}). \]
\end{itemize}

It is routine to extend the formalism to groups of metaplectic type, since its counterparts on each $\GL$-factor are trivial. For every $\mathbf{G}^! \in \Endo_{\elli}(\tilde{G})$, put
\[ \iota(\tilde{G}, G^!) := \left( Z_{(G^!)^\vee} : Z_{\tilde{G}^\vee}^\circ \right)^{-1}
= 2^{- \#(\SO\text{-factors in}\; G^!)}. \]
For a number field $F$, we also have $\iota(G, G^!) = \tau(G) \tau(G^!)^{-1}$ where $\tau(\cdot)$ is the normalized Tamagawa number.

\begin{definition}\label{def:I-Endo}
	Consider a number field $F$, a covering of metaplectic type $\rev: \tilde{G} \to G(\A_F)$ and $V \supset V_{\mathrm{ram}}$. For each invariant genuine distribution $I = I^{\tilde{G}}$ on $\tilde{G}_V$ that intervenes in the invariant trace formula for $\tilde{G}$, we may consider its stable counterparts $S^{G^!}$ on $G^!(F_V)$, for various $\mathbf{G}^! \in \Endo_{\elli}(\tilde{G})$. By taking $\otimes$-products, an anti-genuine test function $f$ on $\tilde{G}_V$ may be transferred to $G^!(F_V)$, still denoted as $f \mapsto f^!$. Define the \emph{endoscopic counterpart} of $I$ as
	\[ I^{\Endo}(f) := \sum_{\mathbf{G}^! \in \Endo_{\elli}(\tilde{G})} \iota(G, G^!) S^{G^!}(f^!). \]
\end{definition}

In particular, we obtain $I_{\mathrm{geom}}^{\Endo}$, $I_t^{\Endo}$, $I^{\Endo}_\nu$ and so on. For the definition of $I^{\Endo}_\nu$, note that there is a well-defined map of infinitesimal characters $\nu^! \mapsto \nu$ for each $\mathbf{G}^! \in \Endo_{\elli}(\tilde{G})$, which allows us to define $S^{G^!}_\nu := \sum_{\nu^! \mapsto \nu} S^{G^!}_{\nu^!}$. Similarly, we have a map of Satake parameters $c^{V, !} \mapsto c^V$, which leads to $I^{\Endo}_{\mathrm{disc}, \nu, c^V}$, etc.

The main results of this works are stated as follows.

\begin{theorem}[\textit{infra.} Theorem \ref{prop:spectral-stabilization}]\label{prop:STF-intro}
	For all infinitesimal character $\nu$, we have
	\[ I_\nu = I^{\Endo}_\nu \]
	as genuine invariant distributions on $\tilde{G}_V$. Consequently,
	\[ I_t = I^{\Endo}_t, \quad I_{\mathrm{spec}} = I^{\Endo}_{\mathrm{spec}} \]
	for all $t \geq 0$.
\end{theorem}

Moreover, one also stabilizes the ``core'' or the global coefficients in $I_{\mathrm{spec}}$, namely the discrete part.

\begin{theorem}[\textit{infra.} Theorem \ref{prop:spectral-stabilization-disc}]\label{prop:STF-disc-intro}
	For all infinitesimal character $\nu$ and Satake parameter $c^V$ outside $V$, we have
	\[ I_{\mathrm{disc}, \nu, c^V} = I^{\Endo}_{\mathrm{disc}, \nu, c^V} \]
	as genuine invariant distributions on $\tilde{G}_V$. Consequently, we also have
	\[ I_{\mathrm{disc}, \nu} = I^{\Endo}_{\mathrm{disc}, \nu}, \quad I_{\mathrm{disc}, t} = I^{\Endo}_{\mathrm{disc}, t} \]
	for all $t \geq 0$.
\end{theorem}

More generally, we may consider the stabilization of a distribution, a map or a function, in either the global, local or semi-local settings, with the understanding that
\begin{center}
	Stabilization := matching some term with its endoscopic counterpart. 
\end{center}
This point of view will be applied to all the terms in the invariant trace formula.

\section{Refined matching}
The terms in the invariant trace formula $I^{\tilde{G}}_{\mathrm{geom}} = I^{\tilde{G}}_{\mathrm{spec}}$ can be catalogued as follows.

\begin{center}\begin{tabular}{|c|c|c|c|} \hline
	Type & Global distributions & Global ``coefficients'' & Semi-local distributions \\ \hline
	Notation & $I^{\tilde{G}}$ & $A^{\tilde{G}}$ & $I^{\tilde{G}_V}_{\tilde{M}_V}$ \\ \hline
	Examples & $I^{\tilde{G}}_{\mathrm{geom}}$ & $A^{\tilde{G}}(V, \mathcal{O})$ & $I^{\tilde{G}_V}_{\tilde{M}_V}(\tilde{\gamma}, f)$ \\
	& $I^{\tilde{G}}_{\mathrm{spec}}$ & $I^{\tilde{G}}_{\mathrm{disc}, t}$ & $I^{\tilde{G}_V}_{\tilde{M}_V}(\pi_V, \nu, X, f)$ \\ \hline
\end{tabular}\end{center}
Here $V$ is a finite set of places containing $V_{\mathrm{ram}}$. The semi-local distributions are of a relative nature: they depend on a Levi subgroup $M \subset G$, and each semi-distribution admits a splitting formula into local avatars $I^{\tilde{G}_v}_{\tilde{M}_v}$ where $v \in V$. As usual, the superscripts $\tilde{G}$ are often omitted.

Thanks to Arthur's works, we have a similar pattern for the stable trace formula $S^{G^!}_{\mathrm{geom}} = S^{G^!}_{\mathrm{spec}}$ for every $\mathbf{G}^! \in \Endo_{\elli}(\tilde{G})$, namely
\begin{center}\begin{tabular}{|c|c|c|c|} \hline
	Type & Global distributions & Global ``coefficients'' & Semi-local distributions \\ \hline
	Notation & $S^{G^!}$ & $SA^{G^!}$ & $S^{G^!_V}_{M^!_V}$ \\ \hline
	Examples & $S^{G^!}_{\mathrm{geom}}$ & $SA^{G^!}(V, \mathcal{O}^!)$ & $S^{G^!_V}_{M^!_V}(\delta, f^!)$ \\
	& $S^{G^!}_{\mathrm{spec}}$ & $S^{G^!}_{\mathrm{disc}, t}$ & $S^{G^!_V}_{M^!_V}(\pi^!_V, \nu, X, f)$ \\ \hline
\end{tabular}\end{center}
Again, stable semi-local distributions can be expressed by local ones via stable splitting formulas.

Following Arthur, the strategy towards Theorems \ref{prop:STF-intro}, \ref{prop:STF-disc-intro} is to define an endoscopic analogue for each term in the invariant trace formula from the stable counterparts, affected with a superscript $\Endo$, and try to prove a \textsc{term by term matching}. This means:
\[ A^{\tilde{G}, \Endo} = A^{\tilde{G}}, \quad I^{\tilde{G}_V, \Endo}_{\tilde{M_V}} = I^{\tilde{G}_V}_{\tilde{M}_V}, \quad I^{\tilde{G}, \Endo} = I^{\tilde{G}}. \]
Furthermore, we wish to reduce the matching for semi-local distributions to the purely local case, by showing that $I^{\tilde{G}_V, \Endo}_{\tilde{M}_V}$ and $I^{\tilde{G}_V}_{\tilde{M}_V}$ satisfy splitting formulas of the same form.

The formation of the global objects $I^{\tilde{G}, \Endo}$ and $A^{\tilde{G}, \Endo}$ follows the pattern of Definition \ref{def:I-Endo}. The formation of $I^{\tilde{G}_V, \Endo}_{\tilde{M}_V}$ and their local avatars are more roundabout: one has to consider all $\mathbf{M}^! \in \Endo_{\elli}(\tilde{M})$, together with all possible ways to fit $(\mathbf{M}^!, \tilde{M}, \tilde{G})$ into the situation
\begin{equation}\label{eqn:s-situation-pre}\begin{tikzcd}
	G^! \arrow[dashed, leftrightarrow, r, "\text{ell.}", "\text{endo.}"'] & \tilde{G} \\
	M^! \arrow[dashed, leftrightarrow, r, "\text{ell.}", "\text{endo.}"'] \arrow[hookrightarrow, u, "\text{Levi}"] & \tilde{M} \arrow[hookrightarrow, u, "\text{Levi}"']
\end{tikzcd}. \end{equation}
The possible $\mathbf{G}^!$'s are parameterized by a set $\Endo_{\mathbf{M}^!}(\tilde{G})$ that can either be defined in terms of dual groups as Arthur did, or by a simple combinatorial recipe (Definition \ref{def:Endo-s}). For each $s \in \Endo_{\mathbf{M}^!}(\tilde{G})$, we obtain the corresponding $\mathbf{G}^![s] \in \Endo_{\elli}(\tilde{G})$. Set
\[ i_{M^!}(\tilde{G}, G^![s]) := \dfrac{\left( Z_{\check{M}^!} : Z_{\tilde{M}^\vee}^\circ \right)}{\left( Z_{G^![s]^\vee} : Z_{\tilde{G}^\vee}^\circ \right)}. \]
To appreciate this definition, note that these factors appear in the key combinatorial summation formula discussed in \S\ref{sec:combinatorial-summation}.

Take the local invariant weighted orbital integral $I^{\tilde{G}}_{\tilde{M}}(\tilde{\gamma}, f)$ for example, where $\rev: \tilde{G} \to G(F)$ is a local metaplectic covering and $\tilde{\gamma}$ is a conjugacy class in $\tilde{M}$. Given $\mathbf{M}^!$, a stable conjugacy class $\delta$ in $M^!(F)$ together with various choices of Haar measures (for the sake of simplicity), we form
\[ I^{\tilde{G}, \Endo}_{\tilde{M}}(\mathbf{M}^!, \delta, f) = \sum_{s \in \Endo_{\mathbf{M}^!}(\tilde{G})} i_{M^!}(\tilde{G}, G^![s]) S^{G^![s]}_{M^!}(\delta[s], B^{\tilde{G}}, f^!) \]
where
\begin{itemize}
	\item $f^!$ is the transfer to $G^!$ of $f$;
	\item $B^{\tilde{G}}$ is a \emph{system of $B$-functions} on $G^![s]$ determined by $\tilde{G}$ --- this is a notion borrowed from twisted endoscopy, see \S\ref{sec:local-geometric} for a discussion;
	\item $\delta[s] := \delta \cdot z[s]$ where $z[s]$ is a central element in $M^!(F)$ determined by $s$.
\end{itemize}

The system $B^{\tilde{G}}$ governs a rescaling of roots of $G^![s]$. It fades away in the global setting.

The \emph{central twist} $\delta \mapsto \delta[s]$ is another subtlety in the metaplectic theory. Roughly speaking, it comes from the fact that if the square \eqref{eqn:s-situation-pre} is chased in two ways, the resulting correspondence of semisimple classes will differ by $z[s]$.

Furthermore, we have to transform various $I^{\tilde{G}, \Endo}_{\tilde{M}}(\mathbf{M}^!, \delta, \cdot)$ into an expression $I^{\tilde{G}, \Endo}_{\tilde{M}}(\tilde{\gamma}, \cdot)$ that no longer refers to $\mathbf{M}^!$. This requires a complete understanding of the image of transfer. The same pattern applies to other semi-local and local distributions in the trace formula, and nontrivial arguments will be needed in each case.

In the study of $I^{\tilde{G}, \Endo}_{\tilde{M}}(\mathbf{M}^!, \delta, \cdot)$, the starting point will be the $\tilde{G}$-regular case: we say $\delta$ is $\tilde{G}$-regular, or simply $G$-regular, if it corresponds to a stable semisimple class $\gamma \in M(F)$ that is regular semisimple in $G$.

\begin{remark}\label{rem:G-equal-M}
	We conclude this section by two quick observations.
	\begin{itemize}
		\item All the matching statements above reduce to the case $\tilde{G} = \Mp(W, \A_F)$ (resp. $\tilde{G} = \Mp(W)$) in the global (resp.\ local) case, by treating the $\Mp$ and $\GL$ factors in $\tilde{G}$ separately.
		\item When $G = M$, the matching of local distributions $I^{\tilde{G}, \Endo}_{\tilde{M}} = I^{\tilde{G}}_{\tilde{M}}$ reduces to the definition of transfer.
	\end{itemize}
\end{remark}

\section{Strategy of the proof}
In order to prove the term by term matching, we may make the following inductive assumptions.
\begin{description}
	\item[Global assumption] Let $F$ be a number field. For each proper Levi subgroup $L$ of $G$, we have
	\begin{itemize}
		\item the matching of global distributions $I^{\tilde{L}} = I^{\tilde{L}, \Endo}$,
		\item the matching of global coefficients $A^{\tilde{L}} = A^{\tilde{L}, \Endo}$.
	\end{itemize}
	\item[Local assumption] Let $F$ be a local field of characteristic zero. Given a Levi subgroup $M$ of $G$, we have the matching of various local distributions
	\[ I^{\tilde{L}}_{\tilde{S}} = I^{\tilde{L}, \Endo}_{\tilde{S}} \]
	where $L \supset S$ are Levi subgroups of $G$ with $S \supset M$, such that either $L \neq G$ or $S \neq M$.
\end{description}

The local assumption is justified by Remark \ref{rem:G-equal-M}: the case $G=M$ reduces to the definition of transfer. Also note that once the splitting formulas are proved, the local assumption implies a similar matching $I^{\tilde{L}_V}_{\tilde{S}_V} = I^{\tilde{L}_V, \Endo}_{\tilde{S}_V}$ in the semi-local setting.

Ideally, one would hope to start with the geometric side, trying to stabilize all the local distributions and global coefficients therein; the spectral stabilization would then follow smoothly. The actual scenario is more complicated. The difficulty is that we know little beyond the local matching in the case $M=G$. How can one bootstrap on such a ground?

Arthur accomplishes this task for connected reductive groups by a long, local-global argument passing back and forth between spectral and geometric sides. In particular, the stabilization of local invariant weighted orbital integrals makes substantial use of the spectral side of global trace formula. Mo\!eglin--Waldspurger \cite{MW16-1, MW16-2} followed a similar strategy in the twisted case. During this process, one obtains a plethora of new distributions and maps, each of them is to be stabilized.

The route is relatively shorter in the metaplectic case. First, over a number field $F$, the stabilization of global coefficients $A^{\tilde{G}}(V, \mathcal{O})$ (Theorem \ref{prop:matching-coeff-A}) can be done directly. This is achieved in \S\ref{sec:gdesc} by first reducing to a matching for elliptic coefficients $A^{\tilde{G}}(S, \mathcal{O})_{\elli}$, and then descending to unipotent coefficients $A_{\mathrm{unip}}(S)$ on centralizers. With due care on Galois cohomologies and the behavior of transfer factors, we will be able to apply the results of Arthur (resp.\ Mo\!eglin--Waldspurger) for the standard (resp.\ nonstandard) endoscopic matching of $A_{\mathrm{unip}}(S)$. This procedure is called global descent in \cite{Ar01}.

A crucial ingredient here is the \emph{weighted fundamental lemma} (Theorem \ref{prop:LFP-general}), asserting that in the local unramified setting
\[ r^{\tilde{G}, \Endo}_{\tilde{M}}(\mathbf{M}^!, \delta) = r^{\tilde{G}}_{\tilde{M}}(\tilde{\gamma}, K) \]
whenever $\delta$ transfers to $\tilde{\gamma}$. The $\tilde{G}$-regular case is in \cite{Li12a}
\footnote{The proof of weighted fundamental lemma in \cite{Li12a} proceeds by reducing to scenario of \cite{CL10, CL12} on Lie algebras via Harish-Chandra descent. Strictly speaking, the authors of \cite{CL10, CL12} only treated the standard case for split groups, but they asserted that similar methods will apply in general.},
and the general case follows via Shalika germs. It implies that $A^{\tilde{G}}(S, \mathcal{O})_{\elli}$ and its endoscopic counterpart $A^{\tilde{G}, \Endo}(S, \mathcal{O})_{\elli}$ are ``compressed'' into $A^{\tilde{G}}(V, \mathcal{O})$ and $A^{\tilde{G}, \Endo}(V, \mathcal{O})$ in a compatible way.

The stabilization of local invariant weighted orbital integrals $I^{\tilde{G}}_{\tilde{M}}(\tilde{\gamma}, \cdot)$ is our local geometric Theorem \ref{prop:local-geometric}; here $F$ is a local field of characteristic zero. The argument is much more roundabout.
\begin{enumerate}
	\item First, one reduces to the case where $\tilde{\gamma}$ is $\tilde{G}$-regular, by using a family of germs $\rho^{\tilde{G}}_J$ which arise when one perturbs $\tilde{\gamma}$. The index $J$ here ranges over lattices in $\mathfrak{a}^{G, *}_M$ spanned by linearly independent restricted roots $\alpha_1, \ldots, \alpha_d \in \Sigma(A_M)$ with $d := \dim \mathfrak{a}^G_M$. The stabilization of $\rho^{\tilde{G}}_J$ is also a prerequisite, and can be reduced to the $\tilde{G}$-regular matching as well.
	
	The precise arguments are given in \S\ref{sec:germs}. The stabilization of $\rho^{\tilde{G}}_J$ will involve new coverings of the form
	\[ \widetilde{G_J} = \Mp(2n_1) \utimes{\bmu_8} \cdots \utimes{\bmu_8} \Mp(2n_r), \quad n_1 + \cdots + n_r = n, \]
	a quotient group of $\prod_{i=1}^r \Mp(2n_i)$. Fortunately, they are not too different from metaplectic groups, and we know how to define endoscopy and transfer in this context.

	For Archimedean $F$, the account above is actually over-simplifying: we have to consider $I^{\tilde{G}}_{\tilde{M}}(\tilde{\gamma}, f)$ for $\tilde{\gamma}$ which is not a conjugacy class, but contains normal derivatives along a given conjugacy class; a priori, such a distribution can arise as endoscopic transfer. A program is proposed in \S\ref{sec:ext-Arch} in order to tackle this issue, which also assumes the $\tilde{G}$-regular case. The same program applies and simplifies for non-Archimedean $F$, yielding the required reduction.
		
	\item Secondly, the obstruction to the matching is encapsulated into a cuspidal anti-genuine test function $\epsilon_{\tilde{M}}(f) = \epsilon^{\tilde{G}}_{\tilde{M}}(f)$ on $\tilde{M}$, of almost compact support (a notion defined near the end of \S\ref{sec:orbital-integrals}). It is characterized by
	\[ I^{\tilde{G}, \Endo}_{\tilde{M}}(\tilde{\gamma}, f) - I^{\tilde{G}}_{\tilde{M}}(\tilde{\gamma}, f) = I^{\tilde{M}}(\tilde{\gamma}, \epsilon_{\tilde{M}}(f)) \]
	for all $f$ and $\tilde{G}$-regular $\tilde{\gamma}$, where $I^{\tilde{M}} = I^{\tilde{M}}_{\tilde{M}}$ is the orbital integral on $\tilde{M}$. In order to match the left hand side with an orbital integral, one has to stabilize the weighted Shalika germs (resp.\ jump relations and differential equations) for $\tilde{\gamma} \mapsto I^{\tilde{G}}_{\tilde{M}}(\tilde{\gamma}, f)$ when $F$ is non-Archimedean (resp.\ Archimedean). One also has to pass to compactly-supported analogues ${}^c I^{\tilde{G}}_{\tilde{M}}$ of $I^{\tilde{G}}_{\tilde{M}}$. They are related to $I^{\tilde{G}}_{\tilde{M}}$ via a family of maps ${}^c \theta^{\tilde{G}}_{\tilde{M}}$, and both ${}^c I^{\tilde{G}}_{\tilde{M}}$ and ${}^c \theta^{\tilde{G}}_{\tilde{M}}$ require their own stabilization (Theorems \ref{prop:cpt-supported-equalities}, \ref{prop:theta-matching-arch}); these can all be reduced to the postulated $\tilde{G}$-regular matching.
	
	For non-Archimedean $F$, the precise arguments are given in \S\ref{sec:cpt-nonarch}. The Archimedean case brings the extra difficulties below.
	\begin{itemize}
		\item The formation of ${}^c I^{\tilde{G}}_{\tilde{M}}(\tilde{\gamma}, f)$ requires explicit normalization of intertwining operators for $\tilde{M} \subset \tilde{G}$; the normalizing factors give rise to factors $r^{\tilde{G}}_{\tilde{M}}(\pi)$, which are stabilized by a direct attack (Theorem \ref{prop:easy-stabilization}).
		\item We also have to address the $\widetilde{K^M} \times \widetilde{K^M}$-finiteness of $\epsilon_{\tilde{M}}(f)$ in order to plug them into global trace formulas, where $K^M := K \cap M(F)$. We refer to \S\ref{sec:epsilon-arch} for full details.
	\end{itemize}
\end{enumerate}

The local spectral Theorem \ref{prop:local-spectral}, or the stabilization of local invariant weighted characters $I^{\tilde{G}}_{\tilde{M}}(\pi, \lambda, X, f)$ in general will be a consequence of $I^{\tilde{G}, \Endo}_{\tilde{M}}(\tilde{\gamma}, f) = I^{\tilde{G}}_{\tilde{M}}(\tilde{\gamma}, f)$. However, it is possible to prove it directly in the special case that $\pi^!$ is a stable virtual representation of $M^!(F)$ with unitary central character (Corollary \ref{prop:spectral-trick} (ii)), using a spectral trick due to Arthur.

Now consider the global setting. Given the local result above, the stabilization of the global weighted characters appearing in $I_{\mathrm{spec}}$ reduces to the stabilization of certain factors $r^{\tilde{G}}_{\tilde{M}}(c^V)$ associated with quasi-automorphic Satake parameters of $\tilde{M}$ (Definition \ref{def:quasi-automorphic-Satake}). For every $\mathbf{M}^! \in \Endo_{\elli}(\tilde{M})$ and quasi-automorphic Satake parameter $c^{V, !}$ for $M^!$, with transfer $c^V$, the required matching takes the form
\[ r^{\tilde{G}, \Endo}_{\tilde{M}}(\mathbf{M}^!, c^{V, !}) = r^{\tilde{G}}_{\tilde{M}}(c^V). \]
This is of a combinatorial nature, and is directly provable (Theorem \ref{prop:spectral-fundamental-lemma}).

This reduces the spectral stabilization $I^{\tilde{G}, \Endo}_{\mathrm{spec}} = I^{\tilde{G}}_{\mathrm{spec}}$ to the stabilization of ``coefficients'' $I^{\tilde{G}, \Endo}_{\mathrm{disc}, \star} = I^{\tilde{G}}_{\mathrm{disc}, \star}$, where $\star$ is a placeholder for $t$, $\nu$, $c^V$, etc. In turn, the latter equality will follow from the geometric stabilization $I^{\tilde{G}, \Endo}_{\mathrm{geom}} = I^{\tilde{G}}_{\mathrm{geom}}$. The arguments are given in \S\ref{sec:end-stabilization}.

Let us omit the superscripts $\tilde{G}$, and recapitulate the foregoing reductions by the diagram:

% The codes for "on top" are taken from https://tex.stackexchange.com/questions/497037/cross-over-of-arrows-in-a-complex-diagram
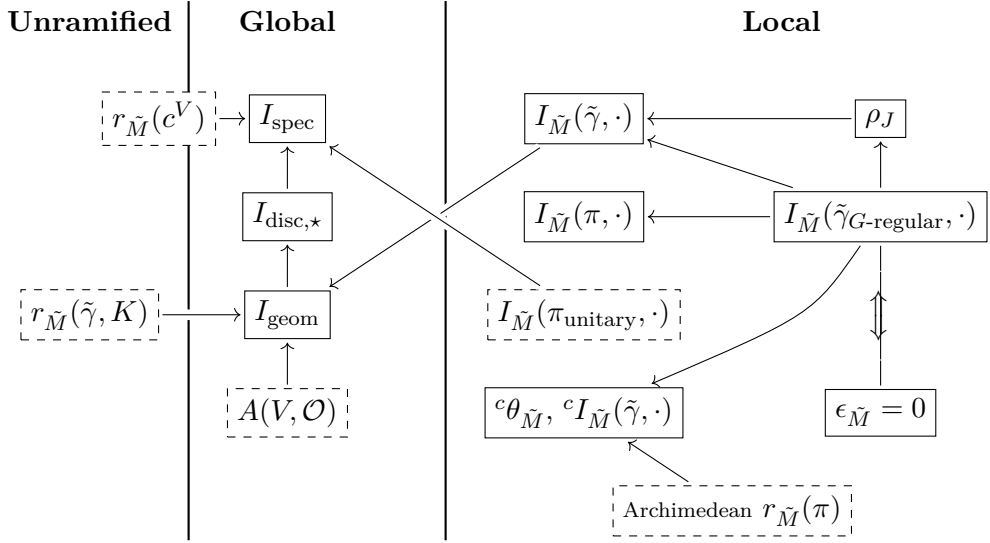
\begin{center}\begin{tikzpicture}[
		every node/.style={draw, outer sep=1.5pt},
		on top/.style={preaction={draw=white,-,line width=#1}},
		on top/.default=4pt,
		xscale=1.3, yscale=1.3
	]
	\draw[thick] (-1.4, 4.2) -- (-1.4, -1.3);
	\draw[thick] (-4, 4.2) -- (-4, -1.3);
		
	\node[draw=white] at (2, 4) {\textbf{Local}};
	\node (RHO) at (3, 3) {$\rho_J$};
	\node (GAMMA) at (0, 3) {$I_{\tilde{M}}(\tilde{\gamma}, \cdot)$};
	\node (GAMMAREG) at (3, 2) {$I_{\tilde{M}}(\tilde{\gamma}_{G\text{-regular}}, \cdot)$};
	\node (PIGEN) at (0, 2) {$I_{\tilde{M}}(\pi, \cdot)$};
	\node[dashed] (PI) at (0, 1) {$I_{\tilde{M}}(\pi_{\mathrm{unitary}}, \cdot)$};
	\node (THETA) at (0, 0) {${}^c \theta_{\tilde{M}}$, ${}^c I_{\tilde{M}}(\tilde{\gamma}, \cdot)$};
	\node (EPSILON) at (3, 0) {$\epsilon_{\tilde{M}} = 0$};
	\node[dashed] (ARCHNOR) at (1.5, -1) {{\scriptsize Archimedean} $r_{\tilde{M}}(\pi)$};
		
	\node[draw=white] at (-3, 4) {\textbf{Global}};
	\node (ISPEC) at (-3, 3) {$I_{\mathrm{spec}}$};
	\node (IDISC) at (-3, 2) {$I_{\mathrm{disc}, \star}$};
	\node (IGEOM) at (-3, 1) {$I_{\mathrm{geom}}$};
	\node[dashed] (A) at (-3, 0) {$A(V, \mathcal{O})$};
		
	\node[draw=white] at (-5, 4) {\textbf{Unramified}};
	\node[dashed, fill=white] (NRNOR) at (-4.3, 3) {$r_{\tilde{M}}(c^V)$};
	\node[dashed] (WFL) at (-5, 1) {$r_{\tilde{M}}(\tilde{\gamma}, K)$};
		
	\path[draw=white] (GAMMAREG) edge node[draw=white, midway, sloped] {$\iff$} (EPSILON);
		
	\draw[->] (GAMMAREG) -- (RHO);
	\draw[->] (GAMMAREG) -- (PIGEN);
	\draw[->] (RHO) -- (GAMMA);
	\draw[->] (GAMMAREG) -- (GAMMA);
	\draw[->] (GAMMAREG) .. controls ([xshift=6em] PI) .. (THETA);
	\draw[->] (ARCHNOR) -- (THETA);
	
	\draw[->] (IGEOM) -- (IDISC);
	\draw[->] (IDISC) -- (ISPEC);
	\draw[->] (A) -- (IGEOM);
	\draw[on top, ->] (GAMMA) -- (IGEOM);
	\draw[->] (NRNOR) -- (ISPEC);
	\draw[on top, ->] (PI) -- (ISPEC);
	\draw[on top, ->] (WFL) -- (IGEOM);
\end{tikzpicture}\end{center}
\begin{itemize}
	\item Each node (except $\epsilon_{\tilde{M}} = 0$) represents a term to be stabilized.
	\item Arrows $A_1, \ldots, A_n \to B$ means that the stabilization of $A_1, \ldots, A_n$ implies that of $B$.
	\item The nodes with dashed borders represent terms that will be stabilized directly in this work.
	\item The inductive assumptions remain in force for each node in this diagram, when we shrink $G$ or enlarge $M$.
\end{itemize}

This reduces the stabilization of trace formula for $\tilde{G}$ to $\epsilon_{\tilde{M}} = 0$, for various Levi subgroups $\tilde{M}$ up to conjugacy. To finalize the proof, we formulate the key geometric Hypothesis \ref{hyp:key-geometric}. For each test function $f$ on $\tilde{G}$, denote by $f_{\tilde{M}}$ its parabolic descent to $\tilde{M}$. The hypothesis asserts that for any given $\mathbf{M}^! \in \Endo_{\elli}(\tilde{M})$, the transfers satisfy
\[ \Trans_{\mathbf{M}^!, \tilde{M}}\left(\epsilon_{\tilde{M}}(f)\right)(\delta) = \epsilon(\mathbf{M}^!, \delta) \Trans_{\mathbf{M}^!, \tilde{M}}\left(f_{\tilde{M}}\right)(\delta) \]
for some $\epsilon(\mathbf{M}^!, \delta) \in \CC$, where
\begin{itemize}
	\item $\delta$ ranges over $\tilde{G}$-regular stable conjugacy classes in $M^!(F)$, and $g(\delta)$ means the stable orbital integral of a test function $g$ along $\delta$, relative to reasonably chosen Haar measures;
	\item $\epsilon(\mathbf{M}^!, \delta)$ is smooth in $\delta$, and independent of $f$.
\end{itemize}

Now enters the final ingredient, the invariant local trace formula in the form of \cite{Li12b}:
\[ I_{\mathrm{geom}}(\overline{f_1}, f_2) = I_{\mathrm{disc}}(\overline{f_1}, f_2). \]
Here $f_1$, $f_2$ are anti-genuine test functions on $\tilde{G}$ and $\overline{f_1}$ denotes the complex conjugate of $f_1$; they can actually be taken from the Schwartz--Harish-Chandra space (Remark \ref{rem:SHC-LTF}). The local trace formula is actually a statement about two different coverings, namely $\tilde{G}$ and its antipodal $\tilde{G}^\dagger$.

There is an endoscopic counterpart $I^{\Endo}_{\mathrm{geom}}(\overline{f_1}, f_2)$ of $I_{\mathrm{geom}}(\overline{f_1}, f_2)$; ditto for $I^{\Endo}_{\mathrm{disc}}(\overline{f_1}, f_2)$. The stabilization of $I_{\mathrm{geom}}(\overline{f_1}, f_2)$ will be a consequence of the local geometric Theorem \ref{prop:local-geometric}.

Using the local trace formula, we will show in Lemma \ref{prop:epsilon-vanishing-1} that
\[ \epsilon(\mathbf{M}^!, \delta) + \overline{\epsilon(\mathbf{M}^!, \delta)} = 0. \]
On the other hand, the MVW involution \cite[p.36]{MVW87} for metaplectic groups relates $\tilde{G}^\dagger$ to $\tilde{G}$. This observation entails that endoscopic transfer is ``isomorphic to its complex conjugate'' in a precise sense, which leads to Lemma \ref{prop:epsilon-vanishing-2}:
\[ \epsilon(\mathbf{M}^!, \delta) = \overline{\epsilon(\mathbf{M}^!, \delta)}. \]
Their conjunction implies $\epsilon(\mathbf{M}^!, \delta) = 0$, hence $\epsilon_{\tilde{M}} = 0$ by varying $(\mathbf{M}^!, \delta)$, as desired.

Ultimately, the key geometric hypothesis will be established in \S\ref{sec:proof-key-geometric} via a local-global argument, using the global spectral stabilization in a special case together with the stabilization of global coefficients $A^{\tilde{M}}(V, \mathcal{O})$ for regular $\mathcal{O}$. As a by-product, we also stabilize the local trace formula, which ought to have applications elsewhere.

\begin{theorem}[\textit{infra.}\ Theorems \ref{prop:stabilization-geom-LTF}, \ref{prop:stabilization-disc-LTF}]
	In the local circumstance, we have
	\[ I^{\Endo}_{\mathrm{geom}}(\overline{f_1}, f_2) = I_{\mathrm{geom}}(\overline{f_1}, f_2), \quad I^{\Endo}_{\mathrm{disc}}(\overline{f_1}, f_2) = I_{\mathrm{disc}}(\overline{f_1}, f_2) \]
	for all anti-genuine test functions $f_1$ and $f_2$ on $\tilde{G}$.
\end{theorem}

To recap:
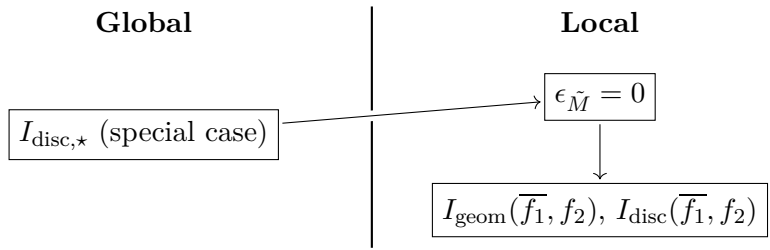
\begin{center}\begin{tikzpicture}[
	every node/.style={draw, outer sep=1.5pt},
	on top/.style={preaction={draw=white,-,line width=#1}},
	on top/.default=4pt
	]
	\node[draw=white] at (3, 2.5) {\textbf{Local}};
	\node (EPSILON) at (3, 1.5) {$\epsilon_{\tilde{M}} = 0$};
	\node (LTF) at (3, 0) {$I_{\mathrm{geom}}(\overline{f_1}, f_2)$, $I_{\mathrm{disc}}(\overline{f_1}, f_2)$};
		
	\node[draw=white] at (-3, 2.5) {\textbf{Global}};
	\node (IDISC) at (-3, 1) {$I_{\mathrm{disc}, \star}$ (special case)};
		
	\draw[thick] (0, 2.7) -- (0, -0.5);
		
	\draw[on top, ->] (IDISC) -- (EPSILON);
	\draw[->] (EPSILON) -- (LTF);
\end{tikzpicture}\end{center}

We close this overview with a few words on the spectral stabilization in a special case (Proposition \ref{prop:Idisc-stabilization-special}) that plays a key role in the bootstrapping process. For a global metaplectic covering $\tilde{G}$ and a given $M \neq G$, the assumptions are that the test function $f = \prod_{v \in V} f_v$ on $\tilde{G}$ is $\tilde{M}$-cuspidal at two places (Definition \ref{def:M-cuspidal}), and has vanishing singular orbital integrals at some non-Archimedean place. The core of the proof involves a famous principle from Jacquet--Langlands: a distribution must vanish if it has both continuous and discrete spectral expansions. This is exactly the method employed by Arthur \cite[\S 7]{Ar03-3}. The same principle is also used in the proof of Lemma \ref{prop:epsilon-vanishing-1}, in the local setting.

\section{Metaplectic features}
Thus far, we have followed the blueprint laid out by Arthur and Mo\!eglin--Waldspurger \cite{MW16-1, MW16-2}, simplified in various aspects, with the existence of transfer and fundamental lemma for $\tilde{G}$ being the main inputs. Let us indicate some new, metaplectic features below.

Recall that there are three kinds of objects in question. The first kind (often denoted by letter $I$) are objects in the invariant trace formula for $\tilde{G}$. They are mostly of an analytic nature, and the properties are analogous to the ones for connected reductive groups. Even though the adaptation to coverings can be nontrivial in some cases, the required works have been done in \cite{Li14a, Li12b, Li13, Li14b} and \cite{Luo17, Luo20}.

The second kind (often denoted by letter $S$) of objects live on the stable side, that is, on the endoscopic groups of $\tilde{G}$. Their definitions and required properties are furnished by Arthur's works, but we also need supplements from Mo\!eglin--Waldspurger since rescaling of roots may also occur on the stable side.

The considerations above for objects of type $I$ and $S$ account for the comparatively short length of this work.

The third kind (with superscript $\Endo$) of objects are of endoscopic nature. Their construction and properties involve combinatorial arguments and/or Galois cohomologies. This is what we have to address seriously in this work.

Although the same trichotomy also exists for connected reductive groups, one cannot simply copy the existing arguments to $\tilde{G}$, for at least two reasons. Assume for simplicity that $\tilde{G} = \Mp(W)$ with $\dim W = 2n$, and denote by $\Gamma_F$ the Galois group of $F$.

\begin{enumerate}
	\item Although $\tilde{G}^\vee$ is defined as $\Sp(2n, \CC)$, the action of $Z_{\tilde{G}^\vee} = \{\pm 1\}$ is to be disregarded; in general, this means that all occurrence of $Z_{M^\vee}^{\Gamma_F}$ in Arthur's works must be replaced by $Z_{\tilde{M}^\vee}^\circ$.
	
	Once phrased in the correct way, this will actually make the metaplectic theory simpler than $\SO(2n+1)$. For example, the elliptic endoscopic groups of $\tilde{G}$ no longer carry endoscopic outer automorphism.
	
	\item Another metaplectic feature is the central twist by $z[s]$ in various endoscopic constructions, for example in $I^{\Endo}_{\tilde{M}}(\mathbf{M}^!, \delta, f)$. Together with the ``centerless'' philosophy above, these force us to redo the proofs of various endoscopic descent and splitting formulas.
	
	\item As mentioned above, the root system on the endoscopic side is to be rescaled when we perform Harish-Chandra descent around a semisimple element. Although this phenomenon is absent in Arthur's case, it does play a critical role in the stabilization twisted trace formula, and brings nonstandard endoscopy into the picture.

	Our results will be formulated in a way such that the nonstandard results in \cite{MW16-1, MW16-2} can help. In fact, we only need the nonstandard endoscopy in type $\mathrm{B}_n \leftrightarrow \mathrm{C}_n$.  The relevance of nonstandard endoscopy to metaplectic groups is first observed in \cite{Li11}.
	
	\item In the study of weighted characters over an Archimedean local field $F$, we have to use normalized intertwining operators instead of the canonically normalized ones in \cite{Li14b}, in order to have a meromorphic family with finitely many polar hyperplanes. In \cite{MW16-2}, this is done via a complicated ``rational normalization'' and comparison of germs. In our case, \cite{Li12b} already provides us with the normalizing factors predicted by Langlands. They match the $\SO(2n+1)$-counterparts up to a simple power of $\sqrt{2}$, as proved by explicit computations in \S\ref{sec:normalized-R-arch}.

	\item Despite the subtle difference between $\Mp(W)$ and $\SO(2n+1)$, several matching statements can actually be reduced to the $\SO(2n+1)$. This is the case for
	\begin{itemize}
		\item the factors $r_{\tilde{M}}(c^V)$ associated with a quasi-automorphic Satake parameter $c^V$ for $\tilde{M}$ in the global setting;
		\item the factors $r_{\tilde{M}}(\pi)$ associated with normalizing factors in the Archimedean local setting, where $\pi$ is a genuine tempered representation of $\tilde{M}$.
	\end{itemize}
	In fact, these factors are related to global and local $L$-functions of Langlands--Shahidi type, respectively. Once the meromorphy of such $L$-functions is known, the corresponding statements are of a combinatorial nature, and reduce to the cases settled in \cite{Ar99b}.
	
	\item Another outstanding metaplectic feature is the use of MVW involution in \S\ref{sec:key-geometric}. Its rationale is deeply rooted in the character formula of the Weil representation in terms of Maslov indices, which enters into the transfer factors for $\Mp(W)$.
\end{enumerate}

As an aside, note that certain endoscopic coverings appear in the Langlands--Shahidi method \cite{Gao18} for $\Mp(W)$, which is used as a black box in the study of $r_{\tilde{M}}(c^V)$. In \S\ref{sec:endo-covering}, we will digress and explain how these coverings are related to our endoscopic groups.

\section{Structure of this work}
Below is only a brief sketch. Detailed overviews will be given in the beginning of each chapter.

In \S\ref{sec:prelim}, we collect the basic definitions about covering groups, their genuine representations, and the distributions of geometric / spectral origin on coverings. Most of them are known results, but have to be formulated in the fashion of \cite{MW16-1} for later use. We also review various constructions about $R$-groups and the trace Paley--Wiener spaces. We conclude this part with a review of Galois cohomologies and their abelianization, in order to fix notation.

In \S\ref{sec:endoscopy-Mp}, we introduce the groups of metaplectic type and relevant constructions in linear algebra, mainly in the local case. We also review the formalism of endoscopy for such covering groups. This part also contains some new results on the Archimedean transfer, including the determination of its image and the spectral transfer factors; these are mild enhancements of the results in \cite{Li19}. Another key ingredient introduced here is the combinatorial summation formula in \S\ref{sec:combinatorial-summation}.

In \S\ref{sec:HC-descent}, we begin with a review of Harish-Chandra descent of endoscopy. We then define the ``diagrams'', which serve to rigidify the endoscopic correspondence of semisimple classes. Diagrams already played a prominent role in \cite{MW16-1, MW16-2}, but they have not been formally introduced in the metaplectic case. Finally, we study Harish-Chandra descent under the presence of a Levi subgroup, and introduce various constants related to the descent of transfer factors in this setting.

The first goal of \S\ref{sec:geom-matching} is to state the local geometric Theorem \ref{prop:local-geometric}, that is, the stabilization of invariant weighted orbital integrals; the Archimedean case requires special care. To this end, we also need to define the system $B^{\tilde{G}}$ of $B$-functions on the endoscopic groups. The second goal is to prove an endoscopic descent formula (Proposition \ref{prop:descent-orbint-Endo}). The combinatorial pattern in its proof will reoccur frequently throughout this work. Finally, we state the weighted fundamental lemma for general conjugacy classes in the unramified case.

Note that we follow \cite[II.1.5]{MW16-1} to modify original the definition of weighted orbital integrals, in order to obtain the desired behavior under induction. See \S\ref{sec:orbint-weighted-nonArch}.

The local geometric Theorem \ref{prop:local-geometric} is reduced to the case of $G$-regular elements in \S\ref{sec:germs}. To this end, we first develop a theory of germs as in \cite[III, V]{MW16-1}. In order to extend the definition of invariant weighted orbital integrals to certain non-orbital distributions in $D_{\mathrm{geom}, -}$ in the Archimedean case, we propose a program in \S\ref{sec:ext-Arch}, which also accomplishes the reduction to $G$-regular case. In \S\ref{sec:rho-main-summary} we explain how the program can be adapted to the non-Archimedean case in a much simpler manner. In the remaining parts, we prove the stabilization of weighted Shalika germs by reducing to known results on centralizers by descent. Based on these ingredients, the weighted fundamental lemma is reduced to the known case of $G$-regular classes.

In \S\ref{sec:cpt-nonarch}, we work over a non-Archimedean $F$, review the compactly-supported version ${}^c I_{\tilde{M}}(\tilde{\gamma}, f)$ of $I_{\tilde{M}}(\tilde{\gamma}, f)$, the associated maps ${}^c \phi_{\tilde{M}}$ and ${}^c \theta_{\tilde{M}}$, and formulate the corresponding matching theorems. We then construct the map $\epsilon_{\tilde{M}}$ under inductive assumptions, the main ingredients being ${}^c I_{\tilde{M}}(\tilde{\gamma}, f)$ and the stabilization of weighted Shalika germs.

The construction of $\epsilon_{\tilde{M}}$ for Archimedean $F$ is more complicated. In \S\ref{sec:weighted-char-arch}, we begin with a study of the normalizing factors for intertwining operators, and match them with their $\SO(2n+1)$-counterparts. We then study ${}^c I_{\tilde{M}}(\tilde{\gamma}, f)$ and the associated maps as in the non-Archimedean case, but we also have to pass between the objects defined by canonically normalization and those defined by normalizing factors. Finally, we state the matching theorems, together with all the preliminary reductions.

The Archimedean construction of $\epsilon_{\tilde{M}}$ continues in \S\ref{sec:epsilon-arch}. The hard part is to stabilize the differential equations and jump relations for invariant weighted orbital integrals; they replace the role of Shalika germs in the non-Archimedean case. After constructing $\epsilon_{\tilde{M}}$, we show that it lands in the space of $\widetilde{K^M} \times \widetilde{K^M}$-finite test functions on $\tilde{M}$ by using a symmetric form of local trace formula.

In \S\ref{sec:TF-geom}, we switch to the global trace formula. After a summary on adélic coverings, the expansion of $I_{\mathrm{geom}}$ is reviewed and reformulated, with various supplements to the previous works \cite{Li14a, Li14b}. In particular, we express the global geometric coefficients in the form of $A^{\tilde{G}}(V, \mathcal{O})$, indexed by semisimple conjugacy classes $\mathcal{O}$ in $M(F)$ (or classes in $M(F_V)$, or stable classes, etc.) We also review the semi-local distributions in $I_{\mathrm{geom}}$, their endoscopic counterparts, as well as the corresponding endoscopic descent formulas. With these preparations, we state the global matching theorems and make the preliminary reductions. Finally, we state and prove the semi-local weighted fundamental lemma, by reducing it to the local case.

The matching of global coefficients $A^{\tilde{G}}(V, \mathcal{O})$ is settled in \S\ref{sec:gdesc}. The procedure is called global descent, requiring an in-depth study of transfer factors and Galois cohomologies. In the end, we will be able to exploit the matching of unipotent coefficients in standard and nonstandard endoscopy, due to Arthur and Mo\!eglin--Waldspurger respectively.

In \S\ref{sec:LTF}, we begin by defining the antipodal covering $\tilde{G}^\dagger$ and relate it to $\tilde{G}$ via MVW involution. Next, we consider the invariant local trace formula in the form of \cite{Li12b}, and write down its endoscopic counterpart. The stabilization of the local trace formula is then reduced to showing that $\epsilon_{\tilde{M}} = 0$ under inductive assumptions. In turn, the vanishing of $\epsilon_{\tilde{M}}$ is reduced to the key geometric Hypothesis \ref{hyp:key-geometric} by using the local trace formula itself.

In the first part of \S\ref{sec:spectral-side}, we introduce the local and semi-local invariant weighted characters, state their stabilization (the local spectral Theorem \ref{prop:local-spectral}) and prove the case of unitary central characters. We then move to the spectral side of the global trace formula, define $I_{\mathrm{disc}, t}$, $I_{\mathrm{disc}, \nu}$, $I_{\mathrm{disc}, \nu, c^V}$ and state the global spectral stabilization. In Proposition \ref{prop:Idisc-stabilization-special}, $I_{\mathrm{disc}, \nu}$ is stabilized for special choices of test functions. Using these results and a local-global argument, we prove in \S\ref{sec:proof-key-geometric} the key geometric Hypothesis \ref{hyp:key-geometric}.

Finally, we will be able to prove all main theorems and resolve all inductive assumptions in \S\ref{sec:end-stabilization}.

\section{Acknowledgements}
This work is an outgrowth of the author's doctoral dissertation defended in 2011. The author is deeply indebted to Prof.\ Jean-Loup Waldspurger for his patience and guidance throughout these years. Besides Arthur's works, many of the proofs here are inspired by, or even modeled on the opus of Mo\!eglin--Waldspurger \cite{MW16-1, MW16-2}. The author is also grateful to Prof.\ Waldspurger for numerous corrections and suggestions on earlier drafts of this work.

Thanks also go to Prof.\ Wee Teck Gan, for his long-running concern for the status of this project despite multiple delays, as well as to Prof.\ Fan Gao for his comments that lead up to \S\ref{sec:endo-covering}.

Last but not least, the author is grateful to the referees for their incredibly careful reading of the draft, which leads to further improvements and corrections.

This work is supported by NSFC-11922101.

\chapter{Preliminaries}\label{sec:prelim}
Below is a presentation of the general notation and definitions that will be used throughout this work. We basically follow Arthur's paradigm, with due care about coverings. The notation for abelianized Galois cohomologies follows Labesse \cite{Lab99}.

\section{General notations}
\subsection{Basics}
The cardinality of a set $E$ is written as $|E|$ or $\# E$. The measure of $E$ is denoted by $\mes(E)$. The index of a subgroup $H$ of $G$ is denoted by $(G:H)$.
\index{mes@$\mes$}

If $E$ is a subset of $\Omega$, let $\mathbf{1}_E: \Omega \to \{0,1\}$ denote its characteristic function.

Representations of a locally compact group $H$ are always realized on $\CC$-vector spaces. Such a representation is customarily written as $(\pi, V)$ or simply $\pi$. The central character of $\pi$, if exists, is denoted by $\omega_\pi$. The contragredient of $\pi$ will always be denoted as $\check{\pi}$, taken in a suitable category that depends on the context (eg.\ smooth admissible representations, etc.)
\index{omega-pi@$\omega_\pi$}

For any commutative ring $R$ and $m \in \Z_{\geq 1}$, let $\mu_m(R) := \left\{ r \in R^\times: r^m = 1 \right\}$. We set
\index{mu-m@${\bmu_m}$}
\[ \bmu_m := \mu_m(\CC). \]

When $\Omega$ is a real manifold (resp.\ totally disconnected space), $C^\infty(\Omega)$ denotes the space of smooth (resp.\ locally constant) functions on $\Omega$. If $\Omega_1$ is a real manifold and $\Omega_2$ is totally disconnected, then $C^\infty(\Omega_1 \times \Omega_2) := C^\infty(\Omega_1) \otimes C^\infty(\Omega_2)$. Ditto for $C^\infty_c$ instead of $C^\infty$, by adding the requirement of having compact support.

\subsection{Local and global fields}
For a field $F$ with chosen separable closure $\overline{F}$, we denote by $\Gamma_F = \Gal(\overline{F}|F)$ its absolute Galois group. If $F$ is a non-Archimedean local field, we denote by $I_F$ the inertia subgroup of $\Gamma_F$. For a finite-dimensional étale $F$-algebra $E$, we denote by $N_{E|F}$ and $\Tr_{E|F}$ the corresponding norm and trace maps, respectively.
\index{GammaF@$\Gamma_F$}

The normalized absolute value of a local field $F$ is $|\cdot| = |\cdot|_F$. For a local or global field $F$, we denote by $\Weil{F}$ (resp.\ $\WD{F}$) its Weil group (resp.\ Weil--Deligne group), and the ring of integers in $F$ (a number field or non-Archimedean local field) is denoted by $\mathfrak{o}_F$. The ring of adèles of a global field $F$ is the restricted product $\A_F := \Resprod_v F_v$ relative to $\mathfrak{o}_v := \mathfrak{o}_{F_v}$ for almost all $v \nmid \infty$. More generally, for any set of places $S$ (resp.\ sets of places $S \supset V$) of $F$, we put
\[ F_S := \Resprod_{v \in S} F_v, \quad (\text{resp. } F^V_S := \Resprod_{v \in S \smallsetminus V} F_v ) \]
\index{WeilF@$\Weil{F}, \WD{F}$}
\index{FVS@$F^V_S$}

By an additive character of a local (resp.\ global) field $F$, we mean a non-trivial continuous homomorphism $\psi: F \to \mathbb{S}^1$ (resp.\ $\psi: \A_F/F \to \mathbb{S}^1$), where $\mathbb{S}^1$ denotes the unit circle in $\CC^\times$.

When $F$ is a non-Archimedean local field, we denote by $\mathfrak{o} = \mathfrak{o}_F$ its ring of integers and $\mathfrak{p} = \mathfrak{p}_F$ the maximal ideal thereof. We say an additive character $\psi$ of $F$ is of conductor $\mathfrak{p}_F^d$ if $\psi|_{\mathfrak{p}_F^d} = 1$ but $\psi|_{\mathfrak{p}_F^{d-1}} \neq 1$.

If $F$ is a global field, we still denote by $\mathfrak{o}$ its ring of integers. More generally, if $S$ is a finite set of places of the global field $F$ containing all the Archimedean ones, $\mathfrak{o}_S$ stands for the ring of $S$-integers in $F$.

\subsection{Varieties and groups}
Let $A$ be a commutative ring and $B$ be an $A$-algebra. The set of $B$-points of an $A$-scheme $X$ is denoted by $X(B)$, and we have the $B$-scheme
\[ X_B := X \dtimes{\Spec A} \Spec B . \]
By an $A$-group we mean an affine group scheme over $A$. Now fix a field $F$. For any $F$-group $H$, its subgroups are always meant to be closed $F$-subgroups, and we denote:
\begin{center}\begin{tabular}{c|c}
	$H^\circ$ & the identity connected component \\
	$Z_H$ & the center \\
	$Z_H(\cdot)$ & centralizers in $H$ \\
	$N_H(\cdot)$ & normalizers in $H$ \\
	$\mathfrak{h} = \Lie H$ & the Lie algebra
\end{tabular}\end{center}

Denote the adjoint action $x \mapsto hxh^{-1}$ as $\Ad(h)$.

The multiplicative group scheme is denoted by $\Gm$. For any affine $F$-group $G$, we set $X^*(G) := \Hom_{F\text{-group}}(G, \Gm)$, whose group structure is written additively. For an $F$-torus $T$, we define $X_*(T) := \Hom_{F\mathrm{-group}}(\Gm, T)$. The canonical pairing $X^*(T) \times X_*(T) \to \Z$ will be denoted by $\lrangle{\cdot, \cdot}$.

When $F$ is algebraically closed, for instance $F = \CC$, we often identify $G$ with $G(F)$.

For a connected reductive $F$-group $G$, we let $G_{\mathrm{SC}} \to G_{\mathrm{der}}$ stand for the simply connected covering of its derived subgroup. We set $G_{\mathrm{ab}} := G/G_{\mathrm{der}}$ and $G_{\mathrm{AD}} := G/Z_G$. Therefore $X^*(G) \simeq X^*(G/G_{\mathrm{der}})$. Furthermore, we define the finite-dimensional $\R$-vector spaces in duality
\begin{gather*}
	\mathfrak{a}^*_G := X^*(G) \otimes \R , \quad \mathfrak{a}_G := \Hom(X^*(G), \R), \\
	\lrangle{\cdot, \cdot}: \mathfrak{a}^*_G \times \mathfrak{a}_G \to \R \;\text{is induced by evaluation.}
\end{gather*}

The automorphism group of $G$ sits in an exact sequence
\[ 1 \to G_{\mathrm{AD}} \to \Aut(G) \to \mathrm{Out}(G) \to 1, \]
which splits after choosing a pinning of $G$.

Let $G$ be as above. For a topological field $F$, we endow $G(F)$ with the topology induced from $F$. For a global field $F$, we topologize $G(\A_F)$ accordingly. When $F$ is local, we have the Harish-Chandra homomorphism
\begin{gather*}
	H_G: G(F) \to \mathfrak{a}_G, \\
	\lrangle{H_G(g), \chi} = |\chi(g)|, \quad \chi \in X^*(G).
\end{gather*}
We also set $G(F)^1 := \Ker(H_G)$.
\index{HG@$H_G$}

When $F$ is global, we can also define $H_G: G(\A_F) \to \mathfrak{a}_G$ by using the absolute value $|\cdot| := \prod_v |\cdot|_v$ on $\A_F$; in this case $H_G|_{G(F)} = 0$.

\subsection{Linear algebra}\label{sec:linear-algebra-prelim}
Let $F$ be a field of characteristic $\neq 2$. By a quadratic (resp.\ symplectic) vector space, we mean a pair $(V, q)$ where $V$ is a finite-dimensional $F$-vector space and $q: V \times V \to F$ is a symmetric (resp.\ alternating) non-degenerate bilinear form. We write $q(v) := q(v|v)$ for all $v \in V$.
\index{symplectic vector space}
\index{quadratic vector space}

Assume $(V, q)$ is a quadratic (resp.\ symplectic) $F$-vector space. A subspace $V_0 \subset V$ is called totally isotropic if $V_0 \subset V_0^\perp$. When $V_0$ is totally isotropic, the subquotient $V_0^\perp / V_0$ is still a quadratic (resp.\ symplectic) vector space. 

Maximal totally isotropic subspaces are called Lagrangian subspaces. In the symplectic case, $V_0 \subset V$ is Lagrangian if and only if $V_0 = V_0^\perp$.

Hereafter, the notation $(V,q)$ (resp.\ $(W, \lrangle{\cdot|\cdot})$) will be reserved for quadratic (resp.\ symplectic) $F$-vector spaces. The orthogonal (resp.\ symplectic) group associated to a quadratic (resp.\ symmetric) $F$-vector space $(V, q)$ (resp.\ $(W, \lrangle{\cdot|\cdot})$) is denoted by $\mathrm{O}(V, q)$ (resp.\ $\Sp(W) = \Sp(W, \lrangle{\cdot|\cdot})$). We set $\SO(V, q) := \mathrm{O}(V, q) \cap \SL(V)$. We view $\SO(V,q)$ and $\Sp(W)$ as connected reductive $F$-groups.

Likewise, for a quadratic extension $E|F$ of field and $\epsilon \in \{\pm 1\}$, one can consider $\epsilon$-hermitian $E$-vector spaces $(W, h)$ where $W$ is a finite-dimensional $E$-vector space, $h: W \times W \to E$ is non-degenerate, sesquilinear and satisfies $h(w_1|w_2) = \epsilon\tau(h(w_2|w_1))$, where $\tau$ is the nontrivial automorphism in $\Gal(E|F)$. The corresponding unitary groups is denoted by $\mathrm{U}(W, h)$.

The formalism of $\epsilon$-hermitian forms can be extended to include the split case $E = F \times F$, with $\tau \in \Aut_F(E)$ acting by $(x,y) \mapsto (y,x)$. In this case, one can show that $W = V \times V$ for some $F$-vector space $V$ and $U(W, h) \simeq \GL(V)$.

Now let $F$ be local and fix an additive character $\psi$; we always assume $\psi$ non-trivial. Let $(V,q)$ be a quadratic $F$-vector space. The \emph{Weil constant} $\gamma_\psi(q)$ in \cite[\S 14]{Weil64} is characterized by the identity
\[ \int_V \left( \;\int_V \phi(x-y) \psi\left(\frac{q(x)}{2}\right) \dd x \right) \dd y = \gamma_\psi(q) \int_V \phi(x) \dd x \]
for all Schwartz--Bruhat function $\phi$ on $V$, where we use the self-dual Haar measure on $V$ with respect to the bi-character $\psi \circ q(\cdot|\cdot): V \times V \to \CC^\times$. It only depends on $\psi_F \circ q: V \to \CC^\times$.

For all $a \in F^\times$, let $q_a$ be the quadratic form on $V=F$ with $q(x) = ax^2$, and set $\gamma_\psi(a) := \gamma_\psi(q_a)$.
\index{gamma-psi@$\gamma_\psi$}

\subsection{Combinatorics}\label{sec:combinatorics-reductive}
Let $G$ be a connected reductive $F$-group where $F$ is a given field. The parabolic subgroups of $G$ are assumed to be defined over $F$, and the Levi subgroups are their Levi components. The Levi decompositions are written in the form $P = M_P U_P$ where $U_P := R_u(P)$, and so forth, where $R_u(\cdot)$ stands for the unipotent radical. We will often fix a minimal Levi subgroup $M_0$ of $G$. The parabolic (resp.\ Levi) subgroups containing $M_0$ are called semi-standard. If a minimal parabolic subgroup $P_0 = M_0 U_0$ is chosen, the parabolic subgroups containing $P_0$ are called standard, and their Levi factors are called standard Levi subgroups.

For each connected reductive group $G$, denote by $A_G$ be the maximal split torus in $Z_G^\circ$. When $M$ is a Levi subgroup of $G$, we have $M = Z_G(A_M)$, and $A_0 := A_{M_0}$ is a maximal split torus in $G$.
\index{AM@$A_M$, $A_G$}

Letting $M \subset G$ be a Levi subgroup, we set
\begin{align*}
	W^G(M) & := \Stab_{W^G_0}(M) / W^M_0, \\
	\mathcal{L}^G(M) & := \left\{ L: \text{Levi subgroup}, \; L \supset M \right\}, \\
	\mathcal{P}^G(M) & := \left\{ P: \text{parabolic subgroup with $M$ as Levi factor} \right\}, \\
	\mathcal{F}^G(M) & := \left\{ P: \text{parabolic subgroup}, \; P \supset M \right\}.
\end{align*}
\index{LM@$\mathcal{L}^G(M)$, $\mathcal{P}^G(M)$, $\mathcal{F}^G(M)$}

Here $W^G_0 := N_G(M_0)/M_0 = W^G(M_0)$. We will often omit the superscript $G$ and write $\mathcal{L}(M)$, $\mathcal{P}(M)$ and so on. Every $P \in \mathcal{F}(M_0)$ has a canonical Levi decomposition $P = M_P U_P$ with $M_P \in \mathcal{L}(M_0)$.
\index{W(M)@$W^G(M), W^G_0$}

Let $P$ be semi-standard with Levi factor $M$. Set $\mathfrak{a}_P = \mathfrak{a}_M$. For all semi-standard Levi subgroups $L \supset M$, one defines the $\R$-vector space $\mathfrak{a}^L_M$ fitting into a canonically split short exact sequence:
\[ 0 \to \mathfrak{a}_L \to \mathfrak{a}_M \leftrightarrows \mathfrak{a}^L_M \to 0. \]
\index{aGM@$\mathfrak{a}^G_M$, $\mathfrak{a}^{G, *}_M$}

Putting $\mathfrak{a}_0 := \mathfrak{a}_{M_0}$, one can view $\mathfrak{a}^L_M$ as a subspace of $\mathfrak{a}_0$. Dualization yields the $\R$-vector spaces $\mathfrak{a}^{L,*}_M$, $\mathfrak{a}_M^*$, etc.\ and the corresponding canonically split short exact sequences. Denote their complexifications as
\[ \mathfrak{a}^L_{M, \CC} := \mathfrak{a}^L_M \dotimes{\R} \CC = \mathfrak{a}^L_M \oplus i\mathfrak{a}^L_M, \]
and similarly for $\mathfrak{a}_{M, \CC}$, $\mathfrak{a}^{L, *}_{M, \CC}$, etc.

For any $M \in \mathcal{L}(M_0)$, let $\Sigma(A_M)$ be the set of nonzero restrictions of roots for $(G, A_0)$ to $A_M$, and let $\Sigma_{\mathrm{red}}(A_M) \subset \Sigma(A_M)$ be the subset of indivisible elements. It should be noted that $\Sigma_{\mathrm{red}}(A_M)$ does not necessarily give rise to a root system. Let $P = M_P U_P$ be a parabolic subgroup. Define the subset $\Sigma_P \subset X^*(A_P) \subset \mathfrak{a}_P^*$ which indexes the decomposition
\[ \mathfrak{u}_P := \Lie(U_P) = \bigoplus_{\alpha \in \Sigma_P} \mathfrak{u}_\alpha \]
into $A_P$-eigenspaces. Denote by $\Sigma_P^{\mathrm{red}} \subset \Sigma_P$ the subset of ``reduced roots'', i.e.\ the indivisible elements. In $\mathfrak{a}_P^*$ we have the following element
\[ \rho_P := \frac{1}{2} \sum_{\alpha \in \Sigma_P} (\dim \mathfrak{u}_\alpha) \alpha. \]
\index{SigmaAM@$\Sigma(A_M), \Sigma_{\mathrm{red}}(A_M), \Sigma_P$}
\index{rhoP@$\rho_P$}

Fix a minimal parabolic $P_0 \in \mathcal{L}(M_0)$. The set $\Delta_0 = \Delta^G_0$ of simple roots with respect to $(A_0, P_0)$ forms a basis of $\mathfrak{a}^{G, *}_0$. The Weyl group attached to this datum coincides with $W^G_0$.

The standard parabolic subgroups are in order-preserving bijection $P \leftrightarrow \Delta_0^P $ with the subsets of $\Delta_0$. Recall the characterization of $\Delta_0^P$: if $P=MU$ is a standard parabolic and $\Delta_P := \Delta_0 \smallsetminus \Delta^P_0$, then we may identify $\Delta_P$ with a subset of $\Sigma_P^{\mathrm{red}}$ by restricting to $A_P$, and every element of $\Sigma_P$ can be expressed uniquely as a $\Z_{\geq 0}$-linear combination from $\Delta_P$.

From the based root datum associated with $(A_0, P_0)$ and $G$, in addition to $\Delta_0$ we define
\begin{align*}
	\Delta_0^\vee & \subset \mathfrak{a}_0^G  : \quad \text{simple coroots},\\
	\widehat{\Delta_0} & \subset \mathfrak{a}^{G, *}_0  : \quad \text{the dual basis of } \Delta_0^\vee, \\
	\widehat{\Delta_0^\vee} & \subset \mathfrak{a}^G_0  : \quad \text{the dual basis of } \Delta_0 .
\end{align*}
\index{DeltaPQ@$\Delta^P_Q, {\Delta^P_Q}^\vee, \widehat{\Delta^P_Q}, \widehat{{\Delta^P_Q}^\vee}$}

As before, there are relative versions for standard parabolic subgroups $P \supset Q$. Namely, we have the following bases
\begin{gather*}
	\Delta^P_Q \subset \mathfrak{a}^{P,*}_Q, \quad {\Delta^P_Q}^\vee \subset \mathfrak{a}^P_Q, \\
	\widehat{\Delta^P_Q} \subset \mathfrak{a}^{P,*}_Q, \quad \widehat{{\Delta^P_Q}^\vee} \subset \mathfrak{a}^P_Q.
\end{gather*}

In practice, we will fix a $W^G_0$-invariant positive definite quadratic form on $\mathfrak{a}^G_0$. This determines Haar measures on $\mathfrak{a}^L_M$ for all $L, M \in \mathcal{L}(M_0)$ such that $L \supset M$. Let $M \in \mathcal{L}(M_0)$ and $L_1, L_2 \in \mathcal{L}(M)$. Equip $\mathfrak{a}^{L_1}_M \oplus \mathfrak{a}^{L_2}_M$ with the product measure. Consider the canonical linear map $\Sigma: \mathfrak{a}^{L_1}_M \oplus \mathfrak{a}^{L_2}_M \to \mathfrak{a}^G_M$. Define the non-negative real number (see eg.\ \cite[\S 4]{Ar99}):
\begin{equation}\label{eqn:d}
	d^G_M(L_1, L_2) := \begin{cases}
		\dfrac{\text{the Haar measure on }\; \mathfrak{a}^G_M}{\Sigma_* \left( \text{that on }\; \mathfrak{a}^{L_1}_M \oplus \mathfrak{a}^{L_2}_M \right) }, & \text{if} \; \Sigma \;\text{is isomorphism} \\
		0, & \text{otherwise.}
	\end{cases}
\end{equation}
\index{dGM@$d^G_M(L_1, L_2)$}
An equivalent way is to take the push-forward of the Haar measure on $\mathfrak{a}^G_{L_1} \oplus \mathfrak{a}^G_{L_2}$ to be the numerator, and the Haar measure on $\mathfrak{a}^G_M$ as the denominator.

\subsection{\texorpdfstring{$(G, M)$}{(G, M)}-families}
Consider a connected reductive $F$-group $G$ with a chosen minimal Levi subgroup $M_0$. Fix an invariant positive-definite quadratic form on $\mathfrak{a}_0$ and endow all vector subspaces of $\mathfrak{a}_0$ with the corresponding Haar measures.

Let $M \in \mathcal{L}(M_0)$ and $P \in \mathcal{P}(M)$. Define
\begin{equation}\label{eqn:GMtheta}
	\theta_P(\Lambda) = \theta^G_P(\Lambda) := \mes\left( \mathfrak{a}^G_M / \lrangle{\Delta_P^\vee } \right)^{-1} \prod_{\check{\alpha} \in \Delta_P^\vee} \lrangle{\Lambda, \check{\alpha}},
\end{equation}
\index{thetaP@$\theta_P$}
\index{GMfamily@$(G, M)$-family}
and let $(\mathfrak{a}^*_M)_P^+$ be the chamber in $\mathfrak{a}_M^*$ cut out by $\lrangle{\cdot, \check{\alpha}} > 0$ for all $\check{\alpha} \in \Delta_P^\vee$. Here we are assuming that a minimal parabolic subgroup $P_0 \in \mathcal{P}(M_0)$ with $P \supset P_0$ is chosen, so that $\Delta_P$ and $\Delta_P^\vee$ makes sense. One readily checks that $\theta_P$ and $(\mathfrak{a}^*_M)_P^+$ are independent of this choice.

A $(G, M)$-family with values in a Banach space $E$ over $\CC$ is a family of smooth functions
\[ c_P: i\mathfrak{a}^*_M \to E, \quad P \in \mathcal{P}(M), \]
such that $c_P(\lambda) = c_{P'}(\lambda)$ for all adjacent $P, P' \in \mathcal{P}(M)$ and all $\lambda$ lying on the wall separating the chambers $i(\mathfrak{a}_M^*)_P^+$ and $i(\mathfrak{a}_M^*)_{P'}^+$.

A standard fact is that the function
\begin{equation*}
	c_M(\lambda) := \sum_{P \in \mathcal{P}(M)} c_P(\lambda) \theta_P(\lambda)^{-1}
\end{equation*}
extends to a smooth function on $i\mathfrak{a}_M^*$. We put
\begin{equation}\label{eqn:GMlim}
	c_M := c_M(0) = \lim_{\lambda \to 0} \sum_{P \in \mathcal{P}(M)} c_P(\lambda) \theta_P(\lambda)^{-1}.
\end{equation}

For the proofs of these facts as well as further operations on $(G, M)$-families, such as the descent and splitting formulas, we refer to \cite[\S 17]{ArIntro}.

\subsection{Roots and \texorpdfstring{$L$}{L}-groups}
Let $G$ be a connected reductive $F$-group with maximal $F$-torus $T \subset G$.
\begin{itemize}
	\item We say $T$ is \emph{elliptic} in $G$ if $T/A_G$ is anisotropic.
	\item We say $T$ is \emph{fundamental} if $A_T$ is as small as possible.
\end{itemize}
Elliptic maximal tori are fundamental, and the converse is true when $F$ is global or local non-Archimedean.
\index{elliptic}

We denote by $\Sigma(G, T) \subset X^*(T_{\overline{F}})$ (resp.\ $\check{\Sigma}(G, T) \subset X_*(T_{\overline{F}})$) the corresponding set of roots (resp.\ coroots). In contrast with the conventions in \S\ref{sec:combinatorics-reductive}, the roots and coroots here are absolute, i.e.\ over $\overline{F}$. In order to define the Galois action, we should embed $T$ in a Borel subgroup $B$ over $\overline{F}$, and define the \emph{quasisplit} Galois action on $T_{\overline{F}}$, $X^*(T_{\overline{F}})$, etc.: namely, $\sigma \in \Gamma_F$ acts by $\Ad(g_\sigma) \sigma$ where $g_\sigma \in G(\overline{F})$ maps $\sigma(B, T_{\overline{F}})$ to $(B, T_{\overline{F}})$. Equivalently, one can replace $T$ by \emph{the abstract Cartan subgroup}, which is the inverse limit of all $B/[B, B]$, with transition maps being inner automorphisms.
\index{SigmaGT@$\Sigma(G, T)$, $\check{\Sigma}(G, T)$}

The dual group $\check{G}$ (sometimes $G^\vee$) is the connected reductive $\CC$-group defined from the root datum dual to
\[ (X^*(T_{\overline{F}}), \Sigma(G, T), X_*(T_{\overline{F}}), \check{\Sigma}(G, T)). \]
It comes with a pinning $\left( \check{B}, \check{T}, (E_{\check{\alpha}})_{\check{\alpha}} \right)$, and $\Gamma_F$ acts on $\check{G}$ by pinned automorphisms; this is the contragredient of the quasisplit action alluded to above.
\index{Gcheck@$\check{G}$, $G^\vee$, $\Lgrp{G}$}

The $L$-group is $\Lgrp{G} := \check{G} \rtimes \Gamma_F$. Occasionally we also consider its finite Galois form $\check{G} \rtimes \Gal(E|F)$ where $E|F$ is a finite Galois extension splitting $G$.

Recall that $\check{G}$ and $\Lgrp{G}$ depend only on the quasisplit inner form $G^*$ of $G$. Also, every Levi subgroup $M \subset G$ corresponds to a Levi subgroup $\Lgrp{M} = \check{M} \rtimes \Gamma_F \subset \Lgrp{G}$, well-defined up to conjugacy, together with an identification $W^{\check{G}}(\check{M}) \simeq W^G(M)$.

Next, assume that $F$ is either a local field, and recall the notion of $L$-parameters from \cite[\S 8.2]{Bo79}. Let $\Phi(G)$ be the set of $L$-parameters $\phi: \WD{F} \to \Lgrp{G}$ (resp.\ $\phi: \Weil{F} \to \Lgrp{G}$) when $F$ is non-Archimedean (resp.\ Archimedean), subject to the conditions in \textit{loc.\ cit.}, taken up to $\check{G}$-conjugacy. Given an $L$-parameter $\phi$, we define
\[ S_\phi := Z_{\check{G}}(\Image(\phi)). \]
Then $S_\phi^\circ$ is a connected reductive $\CC$-group.
\index{L-parameter@$L$-parameter}

We say $\phi$ is bounded if the projection of $\phi(\Weil{F})$ to $\check{G}$ is relatively compact. This notion depends only on the equivalence class of $\phi$. Define
\begin{align*}
	\Phi_{\mathrm{bdd}}(G) & := \left\{\phi \in \Phi(G) : \text{bounded} \right\}, \\
	\Phi_{2, \mathrm{bdd}}(G) & := \left\{\phi \in \Phi_\text{bdd}(G) : \Image(\phi) \text{ not contained in any proper Levi of } \Lgrp{G} \right\}.
\end{align*}
\index{Phi@$\Phi$, $\Phi_{\mathrm{bdd}}$, $\Phi_{2, \mathrm{bdd}}$}

When $G$ is quasisplit, the factor $d^G_M(L_1, L_2)$ in \eqref{eqn:d} has the following stable avatar (see \cite[(6.1)]{Ar99}): denote by $\check{L}_i \subset \check{G}$ the standard Levi subgroup corresponding to $L_i \subset G$, then
\begin{equation}\label{eqn:e}\begin{aligned}
	k^G_M(L_1, L_2) & := \left( Z_{\check{L}_1}^{\Gamma_F} \cap Z_{\check{L}_2}^{\Gamma_F} : Z_{\check{G}}^{\Gamma_F} \right) \\
	& = \# \Ker\left[ Z_{\check{L}_1}^{\Gamma_F} \big/ Z_{\check{G}}^{\Gamma_F} \to Z_{\check{M}}^{\Gamma_F} \big/ Z_{\check{L}_2}^{\Gamma_F} \right], \\
	e^G_M(L_1, L_2) & := \begin{cases}
		k^G_M(L_1, L_2)^{-1} d^G_M(L_1, L_2), & \text{if}\; \Sigma: \mathfrak{a}^{L_1}_M \oplus \mathfrak{a}^{L_2}_M \rightiso \mathfrak{a}^G_M \\
		0, & \text{otherwise.}
	\end{cases}
\end{aligned}\end{equation}
Observe that $\Sigma$ is an isomorphism implies that $\mathfrak{a}_{L_1}^* \cap \mathfrak{a}_{L_2}^* = \mathfrak{a}_G^*$.
\index{kGM@$k^G_M(L_1, L_2)$}
\index{eGM@$e^G_M(L_1, L_2)$}

\subsection{Conjugacy classes}
What follows will apply to any connected reductive $F$-group $G$.

\begin{definition}\index{strongly regular semisimple}
	We say an element $\gamma \in G(F)$ is \emph{strongly regular semisimple} if $\gamma$ is regular semisimple and $Z_G(\gamma)$ is connected.
\end{definition}

Denote by $G_{\mathrm{reg}} \subset G$ the open subset of strongly regular semisimple elements. Let $G_{\mathrm{unip}} \subset G$ be the closed subvariety of unipotent elements. Let $G(F)_{\mathrm{ss}} \subset G(F)$ be the subset of semisimple elements.

\index{GregGssGunip@$G_{\mathrm{reg}}$, $G(F)_{\mathrm{ss}}$, $G_{\mathrm{unip}}$}

For every $\gamma \in G(F)$, let $G_\gamma := Z_G(\gamma)^\circ$. We always have $A_G \subset A_{G_\gamma}$.
\index{Ggamma@$G_\gamma$}

\begin{definition}\label{def:elliptic-element}\index{elliptic}
	We say that $\gamma \in G_{\mathrm{ss}}(F)$ is \emph{elliptic}\footnote{In some other references, $\gamma$ is said to be elliptic in $G$ if it belongs to an elliptic maximal torus. If $F$ is global or non-Archimedean, or if $\gamma \in G_{\mathrm{reg}}$, these two notions are equivalent. If $F = \R$, our definition of ellipticity is weaker in general. } if $A_{G_\gamma} = A_G$.
\end{definition}

Let $\Gamma(G)$ be the set of conjugacy classes in $G(F)$; it has the subsets
\[ \Gamma_{\mathrm{ss}}(G) \supset \Gamma_{\mathrm{reg}}(G) \supset \Gamma_{\mathrm{reg, ell}}(G) \]
of semisimple, strongly regular semisimple, and elliptic strongly regular semisimple classes, respectively.
\index{GammassGammareg@$\Gamma_{\mathrm{ss}}$, $\Gamma_{\mathrm{reg}}$, $\Gamma_{\mathrm{reg, ell}}$}

\begin{definition}\label{def:stable-conjugacy}
	\index{stable conjugacy}
	Assume that $G$ is quasisplit. Following \cite[\S 3]{Ko82}, we say that $\delta_1, \delta_2 \in G(F)$ are \emph{stably conjugate} if there exists $y \in G(\overline{F})$ such that $\delta_2 = y^{-1} \delta_1 y$ and $y\sigma(y)^{-1} \in G_{\delta_1}(\overline{F})$ for all $\sigma \in \Gamma_F$.
	
	Note that $\delta_2 = y^{-1} \delta_1 y \in G(F)$ only implies that $y\sigma(y)^{-1} \in Z_G(\delta_1)(\overline{F})$; if there exists such a $y$, we say that $\delta_1$ and $\delta_2$ are \emph{geometrically conjugate}. For elements in $G_{\mathrm{reg}}(F)$, all these notions coincide. Define the sets
	\[  \Sigma_{\mathrm{ss}}(G) \supset  \Sigma_{\mathrm{reg}}(G) \supset  \Sigma_{\mathrm{reg, ell}}(G) \]
	of stable conjugacy classes in $G(F)$ as above, by noting that these notions depend only stable (in fact, geometric) classes.
\end{definition}
\index{SigmassSigmareg@$\Sigma_{\mathrm{ss}}$, $\Sigma_{\mathrm{reg}}$, $\Sigma_{\mathrm{reg, ell}}$}

Let $M \subset G$ be a Levi subgroup. Let $M_{G\text{-reg}} := M \cap G_{\mathrm{reg}}$. Define the subsets $\Gamma_{G\text{-reg}}(M) \subset \Gamma_{\mathrm{reg}}(M)$ and $\Sigma_{G\text{-reg}}(M) \subset  \Sigma_{\mathrm{reg}}(M)$ (when $G$ is quasisplit) accordingly.
\index{MGreg@$M_{G\text{-reg}}$}
\index{GammaGreg@$\Gamma_{G\text{-reg}}$}

\section{Local covering groups in general}\label{sec:covering}
Let $F$ be a local field of characteristic zero and consider a connected reductive $F$-group $G$. For $m \in \Z_{\geq 1}$, an $m$-fold covering of $G(F)$ is a central extension of locally compact groups
\[ 1 \to \bmu_m \to \tilde{G} \xrightarrow{\rev} G(F) \to 1. \]

We will adopt the following convention: whenever $\tilde{\delta} \in \tilde{G}$ we write $\delta := \rev(\tilde{\delta})$, and for any subset $\Omega \subset G(F)$ we write $\widetilde{\Omega} := \rev^{-1}(\Omega) \subset \tilde{G}$.

For any $g \in G(F)$, the automorphism $\tilde{\delta} \mapsto \tilde{g}\tilde{\delta}\tilde{g}^{-1}$ of $\tilde{G}$ depends only on $g = \rev(\tilde{g})$; we shall write it as $\Ad(g)(\tilde{\delta}) = g\tilde{\delta}g^{-1}$.

Note that for every connected reductive subgroup $L \subset G$, the pull-back $\tilde{L} := \rev^{-1}(L(F))$ is an $m$-fold covering of $L(F)$.

We put $H_{\tilde{G}} := H_G \circ \rev: \tilde{G} \to \mathfrak{a}_G$. More generally, choose a minimal Levi $M_0 \subset G$ and a maximal compact subgroup $K \subset G(F)$ in good position relative to $M_0$. For every $P = MU \in \mathcal{F}(M_0)$ we have the map $H_P: G(F) \to \mathfrak{a}_M$ defined by $H_P(muk) = H_M(m)$ for $m \in M(F)$, $u \in U(F)$ and $k \in K$. Put
\[ H_{\tilde{P}} := H_P \circ \rev: \tilde{G} \to \mathfrak{a}_M. \]
\index{HGtilde@$H_{\tilde{G}}$}
\index{HP@$H_P$}
The modulus function $\delta_{\tilde{P}}$ equals $\exp\left(\lrangle{ 2\rho_P, H_{\tilde{P}}(\cdot)} \right)$, which also equals $\delta_P \circ \rev$.

We say $\tilde{\delta} \in \tilde{G}$ is semisimple, strongly regular semisimple, elliptic, etc.\ if $\delta$ is. Denote by $\tilde{G}_{\mathrm{ss}}$ (resp.\ $\tilde{G}_{\mathrm{reg}}$) the inverse images of $G(F)_{\mathrm{ss}}$ (resp.\ $G_{\mathrm{reg}}(F)$) in $\tilde{G}$.

\begin{definition}\label{def:good}
	\index{good}
	An element $\gamma \in G(F)$ is said to be \emph{good} if for some (equivalently, any) $\tilde{\gamma} \in \rev^{-1}(\gamma)$, we have
	\[ \forall x \in G(F), \quad x\tilde{\gamma}x^{-1} = \tilde{\gamma} \iff x\gamma x^{-1} = \gamma. \]
	This notion depends only on the conjugacy class of $\gamma$. We say $\tilde{\gamma}$ is good if $\gamma$ is.
\end{definition}

Let $\Gamma(\tilde{G})$ be the set of good $G(F)$-conjugacy classes in $\tilde{G}$. Define its subsets
\[ \Gamma_{\mathrm{ss}}(\tilde{G}) \supset \Gamma_{\mathrm{reg}}(\tilde{G}) \supset \Gamma_{\mathrm{reg, ell}}(\tilde{G}) \]
consisting of classes which are semisimple, etc.
\index{GammassGammareg}

The relevance of good conjugacy classes is that they support orbital integrals; see Remark \ref{rem:why-good}. This property will automatically hold in the coverings in question (Proposition \ref{prop:auto-good}).

If $M \subset G$ is a Levi subgroup, then $N_{G(F)}(M(F))$ acts on $\tilde{M}$ by conjugation. A key property \cite[Appendice I]{MW94} is that for any parabolic subgroup $P = MU$ of $G$, the covering $\rev$ has a unique $M(F)$-equivariant splitting over $U(F)$. This yields a canonical decomposition
\[ \tilde{P} = \tilde{M} \cdot U(F). \]
We say that $\tilde{P}$ (resp.\ $\tilde{M}$) is a parabolic (resp.\ Levi) subgroup of $\tilde{G}$. This makes it possible to perform parabolic induction and take Jacquet functors in the usual manner. Moreover, Harish-Chandra's theory of intertwining operators, Plancherel formula, etc.\ carry over to the case of coverings. See \S\ref{sec:rep-coverings}

We can also define the subsets $\Gamma_{G\text{-reg}}(\tilde{M}) \subset \Gamma_{\mathrm{reg}}(\tilde{M})$, etc.

For a function $f: \widetilde{\Omega} \to \CC$ and $z \in \bmu_m$, set
\begin{equation}\label{eqn:fz}
	f^z := \left[ \tilde{M} \ni \tilde{x} \mapsto f(z\tilde{x}) \right].
\end{equation}
We say that $f$ is \emph{genuine} (resp.\ \emph{anti-genuine}) if $f^z = zf$ (resp.\ $f^z = z^{-1} f$) for all $z \in \bmu_m$. In particular, when $\Omega$ is open in $G(F)$, we set
\begin{align*}
	C^\infty_{c, \asp}(\widetilde{\Omega}) & := \left\{ f \in C^\infty_c(\widetilde{\Omega}): \text{anti-genuine} \right\}, \\
	C^\infty_{c, -}(\widetilde{\Omega}) & := \left\{ f \in C^\infty_c(\widetilde{\Omega}): \text{genuine} \right\}.
\end{align*}
\index{anti-genuine function}
\index{Cinftyasp@$C^\infty_{c, \asp}$, $C^\infty_{c, -}$}

Fix a Haar measure on $G(F)$. The corresponding Haar measure on $\tilde{G}$ is taken to be:
\begin{equation}\label{eqn:measure-cover}
	\mes(E) = \mes(\rev^{-1}(E)), \quad E \subset G(F): \text{any compact subset}.
\end{equation}

\begin{remark}
	The same convention \eqref{eqn:measure-cover} is adopted in \cite{Li14a, Li12b}. The point is that for all $\bmu_m$-invariant $L^1$-function $\phi$ on $\rev^{-1}(E)$, say $\phi = \phi_0 \circ \rev$, we have $\int_{\rev^{-1} E} \phi = \int_E \phi_0$. In particular, the convolution of anti-genuine (resp.\ genuine) functions on $\tilde{G}$ is expressed by integrals over $G(F)$ with the original Haar measure. For example, if $F$ is non-Archimedean, $K \subset G(F)$ is a compact open subgroup equipped with a section $K \hookrightarrow \tilde{G}$ of $\rev$, then the function $f_K \in C_{c, \asp}(K \backslash \tilde{G} / K)$ with $\Supp(f_K) = \rev^{-1} K$ and $f_K(1) = 1$ is the unit of $C_{c, \asp}(K \backslash \tilde{G} / K)$ under convolution, provided that $K$ has total mass $1$ (resp.\ $\frac{1}{m}$) as a subset of $G(F)$ (resp.\ $\tilde{G}$).
\end{remark}

Below is a review of the Jordan decomposition \cite[Proposition 2.2.2]{Li14a} for $\tilde{G}$. Since $\mathrm{char}(F) = 0$, the covering splits canonically over $G_{\mathrm{unip}}(F)$ by \cite[Proposition 2.2.1]{Li14a}. Let $\tilde{\gamma} \in \tilde{G}$ with $\gamma := \rev(\tilde{\gamma})$. Take the Jordan decomposition $\gamma = \gamma_{\mathrm{ss}} \gamma_{\mathrm{u}}$ in $G(F)$. View $\gamma_{\mathrm{u}}$ as an element of $\tilde{G}$ by the splitting above. Then there exists a unique $\tilde{\gamma}_{\mathrm{ss}} \in \rev^{-1}(\gamma_{\mathrm{ss}})$ such that $\tilde{\gamma} = \tilde{\gamma}_{\mathrm{ss}} \gamma_{\mathrm{u}}$, and $\tilde{\gamma}_{\mathrm{ss}} \gamma_{\mathrm{u}} = \gamma_{\mathrm{u}} \tilde{\gamma}_{\mathrm{ss}}$ since the splitting over $G_{\mathrm{unip}}(F)$ is $\Ad(\gamma_{\mathrm{ss}})$-invariant.

Let $M \subset G$ be a Levi subgroup. Consider $\gamma \in M(F)$ and write $\eta := \gamma_{\mathrm{ss}}$. By \cite[II.1.2 (1)]{MW16-1}, we have
\[ M_\gamma = G_\gamma \iff M_\eta = G_\eta. \]

\begin{definition}[Equisingularity; see \textit{loc.\ cit.}]\label{def:equisingular}
	\index{equisingular}
	Let $\gamma \in M(F)$ be as above. When $M_\gamma = G_\gamma$, we say that $\gamma$ is $G$-\emph{equisingular} (or: $\tilde{\gamma}$ is $\tilde{G}$-equisingular if $\tilde{\gamma} \in \rev^{-1}(\gamma)$). This depends only on the conjugacy class of $\gamma_{\mathrm{ss}}$ in $M(F)$.
\end{definition}

As a trivial example, every $\gamma \in M_{G\text{-reg}}(F)$ is $G$-equisingular.

\section{Orbital integrals}\label{sec:orbital-integrals}
Let $F$ be a local field of characteristic zero. First off, we introduce a notation about measures.

\begin{definition}
	\index{mes}
	For every unimodular locally compact group $\Gamma$, let $\mes(\Gamma)$ denote the $\CC$-line generated by Haar measures on $\Gamma$, and let $\mes(\Gamma)^\vee$ denote its dual. If $M$ is a connected reductive $F$-group, we use the shorthand
	\[ \mes(M) := \mes(M(F)), \quad \mes(M)^\vee := \mes(M(F))^\vee . \]
\end{definition}

\index{DM@$D^M$}
Let $M$ be a connected reductive $F$-group. The \emph{Weyl discriminant} is denoted by $D^M$, which is a regular function on the variety $M$. The same symbol also denotes the version for $\mathfrak{m}$.

For this moment, fix a Haar measure on $M(F)$, or equivalently a trivialization over $\R$ of the line $\mes(M)$. Assume hereafter that $\gamma \in M_{\mathrm{reg}}(F)$. Fix a Haar measure on $M_\gamma(F)$ and define the normalized orbital integral
\begin{align*}
	f_M(\gamma) & = I^M(\gamma, f) = I^M_M(\gamma, f) \\
	& := |D^M(\gamma)|^{\frac{1}{2}} \int_{M_\gamma(F) \backslash M(F)} f(x^{-1} \gamma x) \dd x
\end{align*}
where we take the quotient measure on $M_\gamma(F) \backslash M(F)$. We will use both notations throughout this work, for the following reasons.

\begin{itemize}
	\item Choose the Haar measures so that for all $\gamma \in M_{\mathrm{reg}}(F)$ and $x \in M(F)$, the measures on $M_\gamma(F)$ and $M_{x\gamma x^{-1}}(F)$ match under $\Ad(x)$. Then $\gamma \mapsto f_M(\gamma)$ may be regarded as a function $\Gamma_{\mathrm{reg}}(M) \to \CC$.
	\item On the other hand, the notation $I^M(\gamma, f)$ is meant to emphasize that $f \mapsto I^M(\gamma, f)$ is an invariant linear functional in $f$.
\end{itemize}

For a quasisplit $F$-group $M^!$, we may even assume that the Haar measures on various $M^!_\gamma(F)$ are compatible under stable conjugacy, where $\gamma \in M^!_{\mathrm{reg}}(F)$. We have the stable normalized orbital integrals
\begin{align*}
	f^{M^!}(\delta) & = S^{M^!}(\delta, f) = S^{M^!}_{M^!}(\delta, f) \\
	& := \sum_{\substack{\gamma \in \Gamma_{\mathrm{reg}}(M^!) \\ \gamma \mapsto \delta}} f_{M^!}(\gamma)
\end{align*}
where $\delta \in \Sigma_{\mathrm{reg}}(M^!)$. One may regard $f^{M^!}$ as a function $\Sigma_{\mathrm{reg}}(M^!) \to \CC$.
\index{fMdelta@$f^{M^{"!}}(\delta)$}
\index{SMdelta@$S^{M^{"!}}(\delta, f)$}

The same definitions can be adapted to the Lie algebra $\mathfrak{m}$. Define the Weyl discriminant of $X \in \mathfrak{m}_{\mathrm{reg}}(F)$ as
\[ D^M(X) := \det(\ad(X)| \mathfrak{m}/\mathfrak{m}_X) \; \in F^\times. \]
By choosing Haar measures as before, we still have $f_M = I^M_M(\cdot, f): \Gamma_{\mathrm{reg}}(\mathfrak{m}) \to \CC$. When $M^!$ is quasisplit, we similarly have $f^{M^!} = S^{M^!}_{M^!}(\cdot, f): \Sigma_{\mathrm{reg}}(\mathfrak{m}^!) \to \CC$.

Consider now a covering as in \S\ref{sec:covering}
\[ 1 \to \bmu_m \to \tilde{M} \xrightarrow{\rev} M(F) \to 1. \]
Keep the same conventions on Haar measures and recall \eqref{eqn:measure-cover}. Define for every $f \in C^\infty_{\asp}(\tilde{M})$ and $\tilde{\gamma} \in \Gamma_{\mathrm{reg}}(\tilde{M})$,
\begin{align*}
	f_{\tilde{M}}(\tilde{\gamma}) & = I^{\tilde{M}}(\tilde{\gamma}, f) = I^{\tilde{M}}_{\tilde{M}}(\tilde{\gamma}, f) \\
	& := |D^M(\gamma)|^{\frac{1}{2}} \int_{M_\gamma(F) \backslash M(F)} f(x^{-1} \tilde{\gamma} x) \dd x .
\end{align*}
\index{fMgamma@$f_{\tilde{M}}(\tilde{\gamma})$}
\index{IMgamma@$I^{\tilde{M}}(\tilde{\gamma}, f)$}

For every $z \in \bmu_m$, we have $f_{\tilde{M}}(z\tilde{\gamma}) = z^{-1} f_{\tilde{M}}(\tilde{\gamma})$. The same definition also applies to $f \in C^\infty_{-}(\tilde{M})$, in which case $f_{\tilde{M}}(z\tilde{\gamma}) = z f_{\tilde{M}}(\tilde{\gamma})$.

\begin{remark}\label{rem:why-good}
	By definition, conjugacy classes from $\Gamma_{\mathrm{reg}}(\tilde{G})$ are good. The reason for restricting to good classes in $\tilde{G}_{\mathrm{reg}}$ is that for general $\tilde{\gamma}$, the ratio between various $x\tilde{\gamma} x^{-1}$, when $x \in Z_M(\gamma)(F)$, is a $\bmu_m$-valued character. Therefore, among all conjugacy classes in $M(F)$, only the inverse images of the good ones can support nonzero invariant genuine distributions. A comprehensive treatment of the orbital integrals on coverings can be found in \cite[\S 4]{Li12b}.
\end{remark}

\begin{remark}\label{rem:measure-tori}
	In defining $f_{\tilde{M}}(\tilde{\gamma})$, etc., the Haar measures on $M_\gamma(F)$ play only a minor role. In fact, there are even canonical Haar measures on tori; see \cite[\S 3.3]{Li15}. On the other hand, we will have to keep track of Haar measures of $M(F)$.
\end{remark}

\begin{definition}\label{def:orbI}
	\index{Iasp@$\orbI_{\asp}$, $\orbI_-$}
	Denote the $\CC$-vector space of normalized anti-genuine orbital integrals on $\tilde{M}$ by
	\[ \orbI_{\asp}(\tilde{M}) := C^\infty_{c, \asp}(\tilde{M}) \bigg/ \bigcap_{\tilde{\gamma} \in \Gamma_{\mathrm{reg}}(\tilde{M})} \Ker\left[ f \mapsto f_{\tilde{M}}(\tilde{\gamma}) \right]. \]
	The genuine version $\orbI_-(\tilde{M})$ is defined analogously.
\end{definition}

\begin{definition}\label{def:SorbI}
	\index{SI@$S\orbI$}
	Let $M^!$ be a quasisplit connected reductive $F$-group. Denote the $\CC$-vector space of normalized stable orbital integrals on $M^!(F)$ by
	\[ S\orbI(M^!) := C^\infty_c(M^!(F)) \bigg/ \bigcap_{\delta \in \Sigma_{\mathrm{reg}}(M^!)} \Ker\left[ f \mapsto f^{M^!}(\delta) \right]. \]
\end{definition}

The same definitions apply to orbital integrals on Lie algebras.

Assume that the Haar measures on $M_\gamma(F)$ are chosen coherently for all good $\gamma \in M_{\mathrm{reg}}(F)$. Then $\tilde{\delta} \mapsto f_{\tilde{M}}(\tilde{\delta})$ is well-defined for all $f \in \orbI_{\asp}(\tilde{M}) \otimes \mes(M)$, and $f \mapsto f_{\tilde{M}}$ induces an isomorphism
\[
	\orbI_{\asp}(\tilde{M}) \otimes \mes(M) \rightiso \left\{ f_{\tilde{M}}: \Gamma_{\mathrm{reg}}(\tilde{M}) \to \CC , \; f \in \orbI_{\asp}(\tilde{M}) \otimes \mes(M) \right\}.
\]
Thus the elements of $\orbI_{\asp}(\tilde{M}) \otimes \mes(M)$ may be viewed as normalized orbital integrals. Ditto for $\orbI_-(\tilde{M}) \otimes \mes(M)$ and $S\orbI(M^!) \otimes \mes(M^!)$ (for quasisplit $M^!$).

When $F$ is Archimedean, $\orbI_{\asp}(\tilde{M})$, $\orbI_-(\tilde{M})$ and $S\orbI(M)$ are topologized so that they are nuclear LF-spaces, see \cite[pp. 580--581]{Bo94b}. The maps $C^\infty_{c, \asp}(\tilde{M}) \twoheadrightarrow \orbI_{\asp}(\tilde{M})$ are open surjections by the remarks in \cite[p.581]{Bo94b}.

Consider now a connected reductive $F$-group $G$ and a covering $\rev: \tilde{G} \to G(F)$. Let $P = MU$ be a parabolic subgroup of $G$ and choose a maximal compact subgroup $K$ of $G(F)$ in good position relative to $M$ (see \cite[Définition 2.4.1]{Li14a}). Given Haar measures on $G(F)$ and $M(F)$, one can choose Haar measures on $K$ and $U(F)$ so that
\begin{equation}\label{eqn:Iwasawa-integration}
	\int_{\tilde{G}} f(\tilde{x}) \dd \tilde{x} = \iiint_{U(F) \times \tilde{M} \times \tilde{K}} \delta_P(m)^{-1} f(u\tilde{m}\tilde{k}) \dd u \dd \tilde{m} \dd \tilde{k}
\end{equation}
for all $f \in C^\infty_c(\tilde{G})$. Define $f_{\tilde{P}} \in C^\infty_c(\tilde{M})$ as
\begin{equation}\label{eqn:f_P}
	f_{\tilde{P}}: \tilde{x} \longmapsto \delta_P(x)^{\frac{1}{2}} \iint_{K \times U(F)} f(k^{-1} \tilde{x}u k) \dd u \dd k, \quad x \in \tilde{M}.
\end{equation}
\index{fP@$f_{\tilde{P}}$}

The element of $\mes(K) \otimes \mes(U)$ determined by \eqref{eqn:Iwasawa-integration} is proportional to that of $\mes(G) \otimes \mes(M)^\vee$ determined by the chosen measures. Hence $f \mapsto f_{\tilde{P}}$ induces the \emph{parabolic descent} map
\begin{align*}
	\orbI_{\asp}(\tilde{G}) \otimes \mes(G) & \longrightarrow \orbI_{\asp}(\tilde{M})^{W^G(M)} \otimes \mes(M) \\
	f_{\tilde{G}} & \mapsto f_{\tilde{M}}
\end{align*}
which is surjective, and similarly for $\orbI_-(\cdots)$. Here the actions of $W^G(M)$ on $\orbI_{\asp}(\tilde{M})$, etc.\ are induced by conjugation.

Although $f_{\tilde{P}}$ depends on the choice of $P$ and $K$, the elements $f_{\tilde{M}}$, etc.\ do not. In fact we have the

\begin{proposition}\label{prop:descent-orbint}
	For all $\tilde{\gamma} \in \Gamma_{G\mathrm{-reg}}(\tilde{M})$ and the Haar measures chosen above,
	\[ f_{\tilde{G}}(\tilde{\gamma}) = f_{\tilde{M}}(\tilde{\gamma}). \]
	Moreover, if $M \subset L$ are Levi subgroups of $G$, then the transitivity holds: $(f_{\tilde{L}})_{\tilde{M}} = f_{\tilde{M}}$.
\end{proposition}

There is also a stable version for a quasisplit connected reductive $F$-group $G^!$ and its Levi subgroup $M^!$: we have the surjection (see for instance \cite[p.57]{MW16-1})
\begin{align*}
	S\orbI(G^!) \otimes \mes(G^!) & \longrightarrow S\orbI(M^!)^{W^{G^!}(M^!)} \otimes \mes(M^!) \\
	f^{G^!} & \longmapsto f^{M^!}.
\end{align*}
\index{fMshrek@$f^{M^{"!}}$}

\begin{definition}
	\index{Iaspcusp@$\orbI_{\cusp, \asp}$, $\orbI_{\cusp, -}$}
	\index{SIcusp@$S\orbI_{\cusp}$}
	Denote by $\orbI_{\cusp, \asp}(\tilde{G}) \subset \orbI_{\asp}(\tilde{G})$ the subspace consisting of $f_{\tilde{G}}$ such that $f_{\tilde{M}} = 0$ for all proper Levi subgroup $M$.
	
	Likewise, we have $\orbI_{\cusp, -}(\tilde{G}) \subset \orbI_-(\tilde{G})$ and $S\orbI_{\cusp}(G^!)$, when $G^!$ is a quasisplit connected reductive $F$-group.
\end{definition}

All these constructions readily generalize to Lie algebras, and we have the spaces $\orbI(\mathfrak{m})$ and $S\orbI(\mathfrak{m})$ (when $M$ is quasisplit), as well as the parabolic descent maps.

\begin{definition}[The $\tilde{K} \times \tilde{K}$-finite variants]\label{def:K-finite}
	\index{Iasp}
	Assume $F$ is Archimedean. Let $K$ be a maximal compact subgroup of $G(F)$ and $\tilde{K} := \rev^{-1}(K)$. We say a function $f: \tilde{G} \to \CC$ is $\tilde{K} \times \tilde{K}$-finite if its translates under left and right $\tilde{K}$-actions span a finite-dimensional vector space. Let
	\begin{align*}
		C^\infty_{c, \asp}(\tilde{G}, \tilde{K}) & := \left\{ f \in C^\infty_{c, \asp}(\tilde{G}) : \tilde{K} \times \tilde{K}\text{-finite} \right\}, \\
		\orbI_{\asp}(\tilde{G}, \tilde{K}) & := \text{the image of $C^\infty_{c, \asp}(\tilde{G}, \tilde{K})$ in $\orbI_{\asp}(\tilde{G})$}.
	\end{align*}
	If $K$ is in good position relative to a Levi subgroup $M$ (see the discussions preceding \eqref{eqn:Iwasawa-integration}), then $K^M := K \cap M(F)$ is also maximal compact in $M(F)$, and $f_{\tilde{G}} \mapsto f_{\tilde{M}}$ induces
	\[ \orbI_{\asp}(\tilde{G}, \tilde{K}) \otimes \mes(G) \twoheadrightarrow \orbI_{\asp}(\tilde{M}, \widetilde{K^M})^{W^G(M)} \otimes \mes(M). \]
\end{definition}

Define $C^\infty_{c,-}(\tilde{G}, \tilde{K})$ and $\orbI_{-}(\tilde{G}, \tilde{K})$ in a similar vein. For a maximal compact subgroup $K^!$ in a quasisplit $G^!$, we also have $S\orbI(G^!, K^!)$ and the corresponding parabolic descent maps in the $K^! \times K^!$-finite context.

\index{fb@$f^b$}
Let us return to a general covering $\tilde{G}$ over a local field $F$. For all function $f$ on $\tilde{G}$ and $b \in C^\infty_c(\mathfrak{a}_G)$, set $f^b := b(H_{\tilde{G}}(\cdot)) f$. Set
\[ C^\infty_{\mathrm{ac}}(\tilde{G}) := \left\{ f \in C^\infty(\tilde{G}): \forall b \in C^\infty_c(\mathfrak{a}_G),\; f^b \in C^\infty_c(\tilde{G}) \right\}. \]
and one defines $C^\infty_{\mathrm{ac}, \asp}(\tilde{G})$, etc.\ accordingly. Such functions are said to be of \emph{almost compact support}. The regular semisimple orbital integrals can be evaluated on elements of $C^\infty_{\mathrm{ac}, \asp}(\tilde{G})$ once the Haar measures are chosen. Set
\[ \orbI_{\mathrm{ac}, \asp}(\tilde{G}) := C^\infty_{\mathrm{ac}, \asp}(\tilde{G}) \bigg/ \bigcap_{\tilde{\gamma} \in \Gamma_{\mathrm{reg}}(\tilde{G})} \Ker\left[ f \mapsto f_{\tilde{G}}(\tilde{\gamma}) \right]. \]
Hence $\orbI_{\mathrm{ac}, \asp}(\tilde{G}) \supset \orbI_{\asp}(\tilde{G})$.
\index{Iacasp@$\orbI_{\mathrm{ac}, \asp}$}

Note that $\mathfrak{a}_G$ is a direct factor of the group $\tilde{G}$ when $F$ is Archimedean, and the definitions above are simplified.

When $F$ is Archimedean, there are $\tilde{K} \times \tilde{K}$-finite versions $C^\infty_{\mathrm{ac}, \asp}(\tilde{G}, \tilde{K})$ and $\orbI_{\mathrm{ac}, \asp}(\tilde{G}, \tilde{K})$ (note that $H_G|_K = 0$). Since the maximal compact subgroups are mutually $G(F)$-conjugate, $\orbI_{\mathrm{ac}, \asp}(\tilde{G}, \tilde{K})$ is independent of the choice of $K$. Moreover, these spaces carry obvious topologies.

\begin{definition}\label{def:concentrated}
	\index{concentration}
	Let $Z \in \mathfrak{a}_G$. We say a linear functional $T: C^\infty_{c, \asp}(\tilde{G}) \to \CC$ is \emph{concentrated} at $Z$ if $T(f) = T(f^b)$ for all $f$ and $b$ such that $b(Z) = 1$; in this case, $T$ extends naturally to $C^\infty_{\mathrm{ac}, \asp}(\tilde{G}) \to \CC$. Ditto for linear functionals on $\orbI_{\asp}(\tilde{G})$, $\orbI_{\asp}(\tilde{G}, \tilde{K})$, etc.
\end{definition}

\section{Distributions}\label{sec:distributions}
Consider a local field $F$ of characteristic zero and an $m$-fold covering as in \S\ref{sec:covering}
\[ 1 \to \bmu_m \to \tilde{M} \xrightarrow{\rev} M(F) \to 1. \]

\begin{definition}
	\index{distribution}
	Let $\Omega \subset M(F)$ be open. A \emph{distribution} on $\widetilde{\Omega}$ means a linear functional $C^\infty_c(\widetilde{\Omega}) \to \CC$, continuous if $F$ is Archimedean. Recall that $C^\infty_c(\widetilde{\Omega})$ for $F \in \{\R, \CC\}$ is topologized as the direct limit of $\{f \in C^\infty(\widetilde{\Omega}): \Supp(f) \subset C \}$ with norm $f \mapsto \sup_C |f|$, where $C$ ranges over compact subsets.
	
	Given a distribution $\Lambda$ and $z \in \bmu_m$, define
	\[ {}^z \Lambda: f \mapsto \lrangle{\Lambda, f^z}, \]
	where $f^z$ is as in \eqref{eqn:fz}. We say that $\Lambda$ is genuine (resp.\ anti-genuine) if ${}^z \Lambda = z^{-1} \Lambda$ (resp.\ ${}^z \Lambda = z\Lambda$) for all $z \in \bmu_m$.
\end{definition}

\begin{compactitem}
	\item A genuine (resp.\ anti-genuine) distribution is determined by its values on $C^\infty_{\asp}(\widetilde{\Omega})$ (resp.\ $C^\infty_-(\widetilde{\Omega})$).
	\item A genuine (resp.\ anti-genuine) function defines a genuine (resp.\ anti-genuine) distribution provided that it is locally integrable, and that a Haar measure is chosen.
\end{compactitem}

\begin{definition}[Invariant distributions]\label{def:invariant-distributions}
	Let $D_-(\tilde{M})$ denote the linear dual (resp.\ continuous dual) of $\orbI_{\asp}(\tilde{M})$ when $F$ is non-Archimedean (resp.\ Archimedean). We define $D_{\asp}(\tilde{M})$ by dualizing $\orbI_-(\tilde{M})$ instead.
	
	Similarly, for a quasisplit connected reductive group $M^!$, let $SD(M^!)$ denote the linear dual (resp.\ continuous dual) of $S\orbI(M^!)$ when $F$ is non-Archimedean (resp.\ Archimedean).
\end{definition}

The space $D_-(\tilde{M})$ (resp.\ $D_{\asp}(\tilde{M})$) can be identified with the space of genuine (resp.\ anti-genuine) invariant distributions on $\tilde{M}$ if a Haar measure is chosen; see \cite[Théorème 5.8.10]{Li12b} for the non-Archimedean case, and \cite[\S 3]{Bo94b} for the Archimedean case.
\index{Dist@$D_-$, $D_{\asp}$}

Similarly, $SD(M^!)$ can be identified with the space of stable distributions on $M^!(F)$ if a Haar measure is chosen. This is a definition if $F$ is non-Archimedean, and for $F = \R$ this is \cite[Théorème 6.2.1 (ii)]{Bo94b}.
\index{SDist@$SD$}

Consider a covering $\tilde{G}$ together with a Levi subgroup $M \subset G$. The dual of $f_{\tilde{G}} \mapsto f_{\tilde{M}}$ is the \emph{parabolic induction} map
\begin{equation}\label{eqn:parabolic-ind-dist}
	\begin{tikzcd}[row sep=tiny]
		\Ind^{\tilde{G}}_{\tilde{M}}: D_-(\tilde{M}) \otimes \mes(M)^\vee \arrow[r] & D_-(\tilde{G}) \otimes \mes(G)^\vee \\
		\Lambda \arrow[mapsto, r] & \Lambda^{\tilde{G}} := \Ind^{\tilde{G}}_{\tilde{M}}(\Lambda)
	\end{tikzcd}
\end{equation}
\index{Ind@$\Ind$}

Similarly, if $M^!$ is a Levi subgroup of $G^!$, where $G^!$ is quasisplit, the dual of $f_{G^!} \mapsto f_{M^!}$ and its stable avatar $f^{G^!} \mapsto f^{M^!}$ yield
\begin{equation*}\begin{tikzcd}[row sep=tiny]
	\Ind^{G^!}_{M^!}: D(M^!) \otimes \mes(M^!)^\vee \arrow[r] & D(G^!) \otimes \mes(G^!)^\vee \\
	SD(M^!) \otimes \mes(M^!)^\vee \arrow[r] \arrow[phantom, u, "\subset" description, sloped] & SD(G^!) \otimes \mes(G^!)^\vee \arrow[phantom, u, "\subset" description, sloped] \\
	\Lambda \arrow[mapsto, r] & \Lambda^{G^!} := \Ind^{G^!}_{M^!}(\Lambda) .
\end{tikzcd}\end{equation*}

\begin{notation}\label{nota:dist-as-I}
	\index{IGtau@$I^{\tilde{G}}(\tau, f)$}
	\index{SFsigma@$S^{G^{"!}}(\sigma, f)$}
	For all $\tau \in D_-(\tilde{G}) \otimes \mes(G)^\vee$ and $f \in \orbI_{\asp}(\tilde{G}) \otimes \mes(G)$, we will often write
	\[ I^{\tilde{G}}(\tau, f) = I^{\tilde{G}}_{\tilde{G}}(\tau, f) := \lrangle{\tau, f}. \]
	Similarly, for $\sigma \in SD(G^!) \otimes \mes(G^!)^\vee$ and $f \in S\orbI(G^!) \otimes \mes(G^!)$, where $G^!$ is a quasisplit $F$-group, we often write
	\[ S^{G^!}(\sigma, f) = S^{G^!}_{G^!}(\sigma, f) := \lrangle{\sigma, f}. \]
	The purpose is to keep compatibility with the standard notation for orbital integrals and characters in the trace formula.
\end{notation}

\section{Representations of coverings}\label{sec:rep-coverings}
Let $F$ be a local field and consider a covering
\[ 1 \to \bmu_m \to \tilde{G} \to G(F) \to 1 \]
where $G$ is a connected reductive $F$-group and $m \in \Z_{\geq 1}$. We are interested in the representation theory of $\tilde{G}$. Observe that $\tilde{G}$ is totally disconnected locally compact (resp.\ in the Harish-Chandra class) when $F$ is non-Archimedean (resp.\ Archimedean), so it makes sense to consider the category of admissible representations of $\tilde{G}$. This means a Harish-Chandra module when $F$ is Archimedean, for a chosen maximal compact subgroup $K \subset G(F)$ and $\tilde{K} := \rev^{-1}(K)$.

We say a representation $\rho$ of $\tilde{G}$ is \emph{genuine} (resp.\ \emph{anti-genuine}) if $\bmu_m$ acts by $z \mapsto z \cdot \identity$ (resp.\ $z \mapsto z^{-1} \cdot \identity$). We will mainly be interested in the study of genuine representations of coverings.
\index{genuine representation}

In fact, Harish-Chandra's theory carries over to $\tilde{G}$, including the Plancherel formula: see \cite{Li12b, Li14a} for detailed explanations. Below are several instances.
\begin{itemize}
	\item Given $\lambda \in \mathfrak{a}^*_{G, \CC}$ and an admissible representation $\pi$ of $\tilde{G}$, we can form the twist
	\[ \pi_\lambda := \pi \otimes \exp\left( \lrangle{\lambda, H_{\tilde{G}}(\cdot)} \right). \]
	\item Let $P = MU$ be a parabolic subgroup of $G$. Then the normalized parabolic induction of an admissible representation $\sigma$ of $\tilde{M}$ is
	\[ I_{\tilde{P}}(\sigma) = I^{\tilde{G}}_{\tilde{P}}(\sigma) := \Ind^{\tilde{G}}_{\tilde{P}}\left(\sigma \otimes \delta_P^{1/2}\right), \]
	where $\sigma$ is inflated to $\tilde{P}$ via $\tilde{M} = \tilde{P}/U(F)$, and $\delta_P$ is the modulus character pulled back to $\tilde{P}$. This operation preserved genuine (resp.\ anti-genuine) representations.
	\index{IP@$I_{\tilde{P}}$}
	
	\item The standard intertwining operators between parabolic inductions are denoted by
	\[ J_{\tilde{Q}|\tilde{P}}(\sigma_\lambda): I_{\tilde{P}}(\sigma_\lambda) \to I_{\tilde{Q}}(\sigma_\lambda), \quad P, Q \in \mathcal{P}(M), \]
	which are meromorphic families of intertwining operators in $\lambda \in \mathfrak{a}_{M, \CC}^*$, since $I_{\tilde{P}}(\sigma_\lambda)$ and $I_{\tilde{Q}}(\sigma_\lambda)$ can be realized on spaces which are independent of $\lambda$ (namely by restricting the elements to $\tilde{K}$).
	\index{JQP@$J_{\tilde{Q}{"|}\tilde{P}}$}

	\item The contragredient $\check{\pi}$ of an admissible representation $\pi$ of $\tilde{G}$ is defined; it switches genuine and anti-genuine representations.
	\item Assume that $F$ is of characteristic zero, then Harish-Chandra's regularity theorem holds: the character $\Theta_\pi: f \mapsto \Tr\left( \pi(f) \right)$ of an admissible genuine (resp.\ anti-genuine) representation $\pi$ is a genuine (resp.\ anti-genuine) distribution on $\tilde{G}$, such that
	\begin{compactitem}
		\item $\Theta_\pi$ is smooth on $\tilde{G}_{\mathrm{reg}}$,
		\item $\Theta_\pi$ is representable by a locally integrable function on $\tilde{G}$,
		\item $\left| D^G \right|^{1/2} \Theta_\pi$ is locally bounded on $\tilde{G}$,
		\item $\Theta_\pi$ has a local character expansion as in the uncovered case.
	\end{compactitem}
	For Archimedean $F$, this is already included in Harish-Chandra's setting. The case of non-Archimedean $F$ is established in \cite[Corollaire 4.3.3]{Li12b} by descent.
	\index{Theta-pi@$\Theta_\pi$}
\end{itemize}

Since Harish-Chandra's theory extends to $\tilde{G}$, it is legitimate to define the sets
\[ \Pi_-(\tilde{G}) \supset \Pi_{\mathrm{unit}, -}(\tilde{G}) \supset \Pi_{\mathrm{temp}, -}(\tilde{G}) \supset \Pi_{2, -}(\tilde{G}) \]
of isomorphism classes of genuine representations of $\tilde{G}$, the unitary, tempered, and essentially square-integrable ones, respectively. Note that $i\mathfrak{a}^*_G$ preserves each of these subsets.
\index{PiG@$\Pi_-(\tilde{G})$, $\Pi_{\mathrm{unit}, -}(\tilde{G})$, $\Pi_{\mathrm{temp}, -}(\tilde{G})$, $\Pi_{2, -}(\tilde{G})$}

Likewise, we denote by $\Pi_{\asp}(\tilde{G})$, etc.\ their anti-genuine avatars.

\section{Spectral distributions: characters}\label{sec:spectral-distributions}
Fix a local field $F$ of characteristic zero and consider a covering
\[ 1 \to \bmu_m \to \tilde{G} \to G(F) \to 1. \]
Fix a minimal Levi subgroup $M_0$ of $G$. We are going to review certain elements in $D_-(\tilde{G})$, à la Arthur.

We refer to \cite[\S 5.4]{Li12b} and \cite{Luo20} for the Knapp--Stein theory of $R$-groups for coverings. First, for each $M \in \mathcal{L}(M_0)$ we have a space $\tilde{T}_{\elli, -}(\tilde{M})$ formed by triples $(L, \pi, r)$, where
\begin{itemize}
	\item $L \in \mathcal{L}^M(M_0)$, $\pi \in \Pi_{2, -}(\tilde{L})$ and $r$ is an element in $\tilde{R}_\pi$, where
	\[ 1 \to \mathbb{S}^1 \to \tilde{R}_\pi \to R_\pi \to 1 \]
	\index{Rpi@$\tilde{R}_\pi$, $R_\pi$}
	is the central extension of the $R$-group $R_\pi$ attached to $\tilde{L}$ and $\pi$;
	\item the triple $(L, \pi, r)$ must be essential in the sense that for all $z \in \mathbb{S}^1$, we have $zr$ conjugate to $r$ in $\tilde{R}_\pi$ if and only if $z = 1$;
	\item we choose normalizing factors in order to associate a normalized intertwining operator to $r$ --- the general theory is settled in \cite[\S 3]{Li12b}, whereas an explicit recipe for metaplectic groups in the Archimedean case will be given in \S\ref{sec:normalized-R-arch}.
\end{itemize}

The precise conditions for a triplet $(L, \pi, r)$ to belong to $\tilde{T}_{\elli, -}(\tilde{G})$ will be given in \eqref{eqn:Tell-def}: it is a regularity condition on $r$. We also have:
\begin{compactitem}
	\item $\mathbb{S}^1$ acts on $\tilde{T}_{\elli,-}(\tilde{M})$ by scaling $r$, making $\tilde{T}_{\elli,-}(\tilde{M})$ an $\mathbb{S}^1$-torsor;
	\item $W^M_0$ acts on $\tilde{T}_{\elli,-}(\tilde{M})$, and we put
	\[ T_{\elli, -}(\tilde{M}) := \tilde{T}_{\elli, -}(M) / W^M_0 ; \]
	\item $i\mathfrak{a}^*_M$ acts on $\tilde{T}_{\elli,-}(\tilde{M})$ by mapping $\tau := (L, \pi, r)$ to $\tau_\lambda := (L, \pi_\lambda, r)$ where $\lambda \in i\mathfrak{a}^*_M$, which makes $\tilde{T}_{\elli,-}(\tilde{M})$ into a disjoint union of compact tori. Here we used \cite[Lemme 5.6.9]{Li12b} to ensure $R_\pi = R_{\pi_\lambda}$.
\end{compactitem}
\index{TG@$T_-(\tilde{G})$, $T_{\elli, -}(\tilde{G})$, $T_-(\tilde{G})_{\CC}$, $T_{\elli, -}(\tilde{G})_{\CC}$}

The action of $i\mathfrak{a}^*_M$ factors through the quotient group $i\mathfrak{a}^*_{M, F}$, which is defined to be the Pontryagin dual of $H_M(M(F))$. It equals $i\mathfrak{a}^*_M$ (resp.\ is a compact torus) when $F$ is Archimedean (resp.\ non-Archimedean).

This picture can be complexified: let $\tilde{T}_{\elli, -}(\tilde{M})_{\CC}$ be the space of pairs $(L, \pi_\lambda, r)$ with $(L, \pi, r)$ as before and $\lambda \in \mathfrak{a}^*_{M, \CC}$, so that $\mathfrak{a}^*_{M, \CC}$ acts on $\tilde{T}_{\elli, -}(\tilde{M})_{\CC}$. This makes $\tilde{T}_{\elli, -}(\tilde{M})_{\CC}$ into a disjoint union of complex tori. Put $T_{\elli, -}(\tilde{M})_{\CC} := \tilde{T}_{\elli, -}(\tilde{M})_{\CC} / W^M_0$.

Next, we define
\begin{align*}
	\tilde{T}_-(\tilde{G}) & := \bigsqcup_{M \in \mathcal{L}(M_0)} \tilde{T}_{\elli, -}(\tilde{M}), \\
	\tilde{T}_-(\tilde{G})_{\CC} & := \bigsqcup_{M \in \mathcal{L}(M_0)} \tilde{T}_{\elli, -}(\tilde{M})_{\CC}, \\
	T_-(\tilde{G}) & := \tilde{T}_-(\tilde{G}) / W^G_0 = \bigsqcup_{M \in \mathcal{L}(M_0)/W^G_0} T_{\elli, -}(\tilde{M})/W^G(M), \\
	T_-(\tilde{G})_{\CC} & := \tilde{T}_-(\tilde{G})_{\CC} / W^G_0 = \bigsqcup_{M \in \mathcal{L}(M_0)/W^G_0} T_{\elli, -}(\tilde{M})_{\CC} / W^G(M).
\end{align*}
Here we fix normalizing factors for each $M \in \mathcal{L}(M_0)$, in a $W^G_0$-invariant manner. Note that $\mathbb{S}^1$ acts on each of spaces above.

To each $\tau \in (L, \pi, r) \in \tilde{T}_-(\tilde{G})$, one defines the virtual tempered genuine character of $\tilde{G}$ in terms of the normalized intertwining operator $R_{\tilde{P}}(r, \pi)$ (see \cite[p.78]{Li19}):
\begin{align*}
	\Theta_\tau & = \Theta^{\tilde{G}}_\tau : \orbI_{\asp}(\tilde{G}) \otimes \mes(G) \to \CC \\
	& = \Tr\left( R_{\tilde{P}}(r, \pi) I_{\tilde{P}}(\pi, \cdot) \right).
\end{align*}
\index{Theta-tau@$\Theta_\tau$}

The factors in the decomposition of $\Theta_\tau$ into irreducibles share the same central character and infinitesimal character (in the Archimedean case), determined from those of $\pi$. The special case $L = G$ yields the characters $\Theta_\pi$ of $\pi \in \Pi_{2, -}(\tilde{G})$.

For $\tau \in \tilde{T}_-(\tilde{G})$, one views $\Theta_\tau$ as an element of $D_-(\tilde{G}) \otimes \mes(G)^\vee$ since we do not fix the Haar measure on $G(F)$. For $f_{\tilde{G}} \in \orbI_{\asp}(\tilde{G}) \otimes \mes(G)$, we write
\[ f_{\tilde{G}}(\tau) := \Theta_\tau(f_{\tilde{G}}), \quad \tau \in \tilde{T}_-(\tilde{G}). \]
\index{fGtau@$f_{\tilde{G}}(\tau)$}

Following the convention of \cite[\S 6.1]{Li19}, the $\mathbb{S}^1$-action satisfies $R_{\tilde{P}}(\pi, zr) = z^{-1} R_{\tilde{P}}(\pi, r)$, thus
\[ f_{\tilde{G}}(z\tau) = z^{-1} f_{\tilde{G}}(\tau), \quad z \in \mathbb{S}^1, \quad \tau \in \tilde{T}_-(\tilde{G}) . \]

This construction can also be complexified: if $\tau \in \tilde{T}_{\elli, -}(\tilde{M})$, then $\lambda \mapsto f_{\tilde{G}}(\tau_\lambda)$ extends meromorphically to all $\lambda \in \mathfrak{a}^*_{M, \CC}$. This follows from the meromorphic continuation of normalized intertwining operators..

For $\tau \in \tilde{T}_-(\tilde{M})$, in the notation of \eqref{eqn:parabolic-ind-dist} we have
\[ \Ind^{\tilde{G}}_{\tilde{M}} \left( \Theta^{\tilde{M}}_\tau \right) = \Theta^{\tilde{G}}_\tau \; \in D_-(\tilde{G}) \otimes \mes(G)^\vee . \]

Hence we obtain a linear map
\begin{equation}\label{eqn:PW-map}\begin{tikzcd}[row sep=tiny]
		\orbI_{\asp}(\tilde{G}) \otimes \mes(G) \arrow[r] & \left\{ \text{functions}\; \alpha: \tilde{T}_-(\tilde{G}) \to \CC \right\} \\
		f_{\tilde{G}} \arrow[mapsto, r] & {\left[ \tau \mapsto f_{\tilde{G}}(\tau) \right]}.
\end{tikzcd}\end{equation}

\begin{remark}\label{rem:archimedean-lambda-twist}
	\index{TG0@$\tilde{T}_-(\tilde{G})_0$}
	In the Archimedean case, each $i\mathfrak{a}^*_M$-orbit (resp.\ $\mathfrak{a}^*_{M, \CC}$-orbit) in $\tilde{T}_{\elli, -}(\tilde{M})$ (resp.\ $\tilde{T}_{\elli, -}(\tilde{M})_{\CC}$) contains a unique member $\tau$ whose central character is trivial on $A_M(F)^\circ \simeq \mathfrak{a}_M$. Denote the subset of such base-points as $\tilde{T}_{\elli, -}(\tilde{M})_0$ and put $\tilde{T}_-(\tilde{G})_0 := \bigsqcup_{M \in \mathcal{L}(M_0)} \tilde{T}_{\elli, -}(\tilde{M})_0$, and so forth.
\end{remark}

Now comes the stable avatars. Let $G^!$ be a direct product of groups of the form $\GL(m)$ or $\SO(2m+1)$. Fix a minimal Levi subgroup $M^!_0$ of $G^!$. The local Langlands correspondence for the special case of $G^!$ is known after Arthur \cite{Ar13}. In particular, $\Pi_{\mathrm{temp}, -}(G^!)$ (resp.\ $\Pi_{2, -}(G^!)$) is partitioned into $L$-packets indexed by $\Phi_{\mathrm{bdd}}(G^!)$ (resp.\ $\Phi_{2, \mathrm{bdd}}(G^!)$), and we have
\[ \Phi_{\mathrm{bdd}}(G^!) = \bigsqcup_{M^! \in \mathcal{L}(M^!_0)/W^{G^!}_0} \Phi_{2, \mathrm{bdd}}(M^!) / W^{G^!}(M^!). \]
As before, $i\mathfrak{a}^*_{M^!}$ acts on each $\Phi_{2, \mathrm{bdd}}(M^!)$, making $\Phi_{\mathrm{bdd}}(G^!)$ into a disjoint union of quotients of compact tori by some finite groups. Again, there are complexified version with $\mathfrak{a}^*_{M^!, \CC}$-action, etc.

By abusing notation, we choose representatives $\phi$ for the elements of $\Phi_{\mathrm{bdd}}(G^!)$. To each $\phi$, we have the stable tempered character $S\Theta^{G^!}_\phi = S\Theta_\phi$ on $G^!(F)$; it belongs to $SD(G^!) \otimes \mes(G^!)^\vee$. We write $f^{G^!}(\phi) := S\Theta_\phi(f^{G^!})$ for all $f^{G^!} \in S\orbI(G^!) \otimes \mes(G^!)$. This yields a linear map
\[\begin{tikzcd}[row sep=tiny]
	S\orbI(G^!) \otimes \mes(G^!) \arrow[r] & \left\{ \text{functions}\; \beta: \Phi_{\mathrm{bdd}}(G^!) \to \CC \right\} \\
	f^{G^!} \arrow[mapsto, r] & {\left[ \phi \mapsto f^{G^!}(\phi) \right] } .
\end{tikzcd}\]
\index{STheta-phi@$S\Theta_\phi$}
\index{fGphi@$f^{G^{"!}}(\phi)$}

\begin{definition}\label{def:Dspec}
	\index{Dspec@$D_{\mathrm{spec}, -}$, $D_{\mathrm{temp}, -}$, $D_{\elli, -}$}
	Let $F$ be a local field of characteristic zero. Let $D_{\mathrm{spec}, -}(\tilde{G}) \subset D_-(\tilde{G})$ be the linear subspace generated by genuine irreducible characters. Define its subspaces
	\[ D_{\mathrm{temp}, -}(\tilde{G}) \supset D_{\elli, -}(\tilde{G}) \]
	generated by the $\Theta_\tau$ for $\tau$ in $\tilde{T}_-(\tilde{G})$ and $\tilde{T}_{\elli, -}(\tilde{G})$, respectively; the first one is also the subspace generated by all tempered genuine irreducible characters.
	
	When $F$ is Archimedean, we define $D_{\mathrm{spec}, \mu, -}(\tilde{G})$ as the subspace of $D_{\mathrm{spec}, -}(\tilde{G})$ generated by genuine irreducible characters of infinitesimal character $\mu$. Similarly we have $D_{\elli, \mu, -}(\tilde{G})$. We also denote by $D_{\elli, -}(\tilde{G})_0$ the subspace generated by $\Theta_\tau$ where $\tau \in \tilde{T}_{\elli, -}(\tilde{G})_0$.
	
	By taking $\tau$ from the complexified version $\tilde{T}_{\elli, -}(\tilde{G})_{\CC}$, etc., one can define the subspaces $D_{\elli, -}(\tilde{G})_{\CC}$, etc.\ of $D_{\mathrm{spec}, -}(\tilde{G})$.
\end{definition}

For example, the classification of admissible genuine representations translates into
\begin{equation}\label{eqn:Dspec-Ind}\begin{aligned}
	D_{\mathrm{temp}, -}(\tilde{G}) \otimes \mes(G)^\vee & = \left(\bigoplus_{M \in \mathcal{L}(M_0)} \Ind^{\tilde{G}}_{\tilde{M}} D_{\elli, -}(\tilde{M}) \otimes \mes(M)^\vee \right)^{W^G_0}, \\
	D_{\mathrm{spec}, -}(\tilde{G}) \otimes \mes(G)^\vee & = \left(\bigoplus_{M \in \mathcal{L}(M_0)} \Ind^{\tilde{G}}_{\tilde{M}} D_{\elli, -}(\tilde{M})_{\CC} \otimes \mes(M)^\vee \right)^{W^G_0} .
\end{aligned}\end{equation}
Indeed, this is just the variant of \cite[IV.1.2 (2)]{MW16-1} for coverings.

For a quasisplit group $G^!$, we follow \cite[IV.2.2]{MW16-1} to define the stable avatars
\begin{align*}
	SD_{\mathrm{spec}}(G^!) & := D_{\mathrm{spec}}(G^!) \cap SD(G^!), \\
	SD_{\mathrm{temp}}(G^!) & := D_{\mathrm{temp}}(G^!) \cap SD(G^!), \\
	SD_{\elli}(G^!) & := D_{\elli}(G^!) \cap SD(G^!),
\end{align*}
as well as the complexified version $SD_{\elli}(G^!)_{\CC}$, etc. When $F$ is Archimedean, we also have $SD_{\elli}(G^!)_0$, $SD_{\mathrm{spec}, \mu}(G^!)$, etc.\ by specifying an infinitesimal character $\mu$.
\index{SDspec@$SD_{\mathrm{spec}}$, $SD_{\mathrm{temp}}$, $SD_{\elli}$}
\index{SD0@$SD(G^{"!})_0$}

According to \cite[IV.2.8 Corollaire]{MW16-1} in the Archimedean case, or as an immediate consequence of \cite[XI.3.1 Corollaire]{MW16-2} in the non-Archimedean case, we have
\begin{equation}\label{eqn:SDspec-Ind}\begin{aligned}
	SD_{\mathrm{temp}}(G^!) \otimes \mes(G^!)^\vee & = \left( \bigoplus_{M^! \in \mathcal{L}(M^!_0)} \Ind^{G^!}_{M^!} SD_{\elli}(M^!) \otimes \mes(M^!)^\vee \right)^{W^{G^!}_0}, \\
	SD_{\mathrm{spec}}(G^!) \otimes \mes(G^!)^\vee & = \left( \bigoplus_{M^! \in \mathcal{L}(M^!_0)} \Ind^{G^!}_{M^!} SD_{\elli}(M^!)_{\CC} \otimes \mes(M^!)^\vee \right)^{W^{G^!}_0} .
\end{aligned}\end{equation}

\section{Discrete spectral parameters}\label{sec:disc-parameter}
Take $\rev: \tilde{G} \to G(F)$ and $M_0 \subset G$ be as in \S\ref{sec:spectral-distributions}. We will have to consider a subset $\tilde{T}_{\mathrm{disc}, -}(\tilde{G})$ of $\tilde{T}_-(\tilde{G})$, consisting of parameters which occur discretely in the local trace formula; see \S\ref{sec:LTF}.

Recall from \cite[p.834]{Li12b} or \cite[\S 6.1]{Li19} that
\begin{align*}
	W^G(L)_{\mathrm{reg}} & := \left\{ w \in W^G(L): \det(w-1 | \mathfrak{a}^G_L) \neq 0 \right\}, \quad L \in \mathcal{L}(M_0), \\
	\tilde{T}_{\mathrm{disc}, -}(\tilde{G}) & := \left\{ (L, \pi, r) \in \tilde{T}_-(\tilde{G}) : W_\pi^0 r \cap W^G(L)_{\mathrm{reg}} \neq \emptyset \right\},
\end{align*}
\index{Wreg@$W^G(L)_{\mathrm{reg}}$}
\index{Tdisc@$T_{\mathrm{disc}, -}(\tilde{G})$}
through which we can express
\begin{equation}\label{eqn:Tell-def}\begin{aligned}
	\tilde{T}_{\elli, -}(\tilde{G}) & = \left\{ (L, \pi, r) \in \tilde{T}_-(\tilde{G}) : r \in \tilde{R}_{\pi, \mathrm{reg}} \right\}, \\
	& = \left\{ (L, \pi, r) \in \tilde{T}_{\mathrm{disc}, -}(\tilde{G}) : W^0_\pi = \{1\} \right\}
\end{aligned}\end{equation}
where
\begin{itemize}
	\item $P \in \mathcal{P}(M)$ is chosen in order to define the normalized intertwining operators $R_{\tilde{P}}(w, \pi)$,
	\item $W_\pi := \Stab_{W^G(L)}(\pi)$ and $W_\pi^0 \subset W_\pi$ is the subgroup consisting of elements $w$ satisfying $R_{\tilde{P}}(w, \pi) \in \CC^\times \identity$, so that $R_\pi \simeq W_\pi / W_\pi^0$ and there is a decomposition $W_\pi \simeq W_\pi^0 \rtimes R_\pi$,
	\index{Rpi}
	\item $\tilde{R}_{\pi, \mathrm{reg}} \subset \tilde{R}_\pi$ is the preimage of $R_\pi \cap W^G(L)_{\mathrm{reg}}$.
\end{itemize}
The second equality in \eqref{eqn:Tell-def} follows by combining \cite[2.11 Lemme]{MW18} and \cite[p.90]{Ar93}.

The group $W^G_0 \times i\mathfrak{a}^*_G \times \mathbb{S}^1$ acts on $\tilde{T}_{\mathrm{disc}, -}(\tilde{G})$. Define $T_{\mathrm{disc}, -}(\tilde{G}) := \tilde{T}_{\mathrm{disc}, -}(\tilde{G}) / W^G_0$.

Equip $\mathbb{S}^1 \backslash T_{\mathrm{disc}, -}(\tilde{G})$ with the Radon measure such that
\begin{equation}\label{eqn:Tdisc-measure}
	\int_{\mathbb{S}^1 \backslash T_{\mathrm{disc}, -}(\tilde{G})} \alpha(\tau) \dd\tau = \sum_{\substack{\tau \in \mathbb{S}^1 \backslash T_{\mathrm{disc},-}(\tilde{G}) / i\mathfrak{a}_G^* \\ \tau = (L, \pi, r) }} |Z_{R_\pi}(r)|^{-1} \int_{i\mathfrak{a}_G^* / \Stab(\tau)} \alpha(\tau_\lambda) \dd\lambda
\end{equation}
for suitable test functions $\alpha$ on $\mathbb{S}^1 \backslash T_{\mathrm{disc, -}(\tilde{G})}$, where we endow $\Stab(\tau)$ (resp.\ $i\mathfrak{a}^*_G / \Stab(\tau)$) with the counting (resp.\ quotient) measure. This is legitimate since $R_\pi = R_{\pi_\lambda}$ for all $\lambda \in i\mathfrak{a}_G^*$, by \cite[Lemme 5.6.9]{Li12b}.

\begin{definition}
	When $F$ is Archimedean, we denote by $\tilde{T}_{\mathrm{disc}, -}(\tilde{G})_0$ the subset consisting of triplets such that $\omega_\pi|_{A_G(\R)^\circ} = 1$, and so forth.
\end{definition}

\begin{lemma}[Cf.\ {\cite[p.1118]{MW16-2}}]\label{prop:elli-disc-finiteness}
	For every $\tau \in T_{\mathrm{disc}, -}(\tilde{G})$, there exist $R \in \mathcal{L}(M_0)$ and $\sigma \in T_{\elli, -}(\tilde{R})$ such that $\tau$ is the parabolic induction of $\sigma$ to $\tilde{G}$. Moreover, the pair $(R, \sigma)$ is unique up to $W^G_0$.
	
	In the opposite direction, assume $F$ is Archimedean. Let $R \in \mathcal{L}(M_0)$ and $\sigma \in T_{\elli, -}(\tilde{R})_0$. There exist at most finitely many pairs $(L, \tau)$ where $L \in \mathcal{L}(R)$ and $\tau \in T_{\mathrm{disc}, -}(\tilde{L})_0$, such that $\tau$ is the induction of $\sigma_\nu$ for some $\nu \in i\mathfrak{a}_R^*$.
\end{lemma}
\begin{proof}
	This is combinatorial in principle. The first assertion is a consequence of \cite[2.11 Lemme]{MW18}. Uniqueness is built into the definition of $T_-(\tilde{G})$ as the quotient of $\bigsqcup_R \tilde{T}_{\elli, -}(\tilde{R})$ by $W^G_0$.

	As to the second assertion, it suffices to bound the possibilities of $\nu$ such that $\sigma_\nu$ induces to $T_{\mathrm{disc}, -}(\tilde{L})_0$, for each given $L \in \mathcal{L}(R)$. Represent $\sigma$ by a triplet $(S, \pi, r)$ where $S \in \mathcal{L}^R(M_0)$. It is necessary that $\exists w \in W^0_\pi r \subset W^L(S)$ such that
	\[ w(\pi_\nu) \simeq \pi_\nu, \quad w \in W^L(S)_{\mathrm{reg}}, \quad \nu \in i\mathfrak{a}^{L, *}_R . \]
	For such a $w$, the first condition determines $(w - 1)\nu$, which lies in $i\mathfrak{a}^{L, *}_S$ by the third condition. In turn, this determines $\nu$ by the second condition.
\end{proof}

\section{Paley--Wiener spaces}
Retain the notations from \S\ref{sec:spectral-distributions}.

\begin{definition}\label{def:PW-spaces}
	\index{PWasp@$\mathrm{PW}_{\asp}(\tilde{G})$}
	Assume $F$ is non-Archimedean. Let $\mathrm{PW}_{\asp}(\tilde{G})$ be the space of functions $\alpha: \tilde{T}_-(\tilde{G}) \to \CC$ satisfying
	\begin{itemize}
		\item $W^G_0$-invariance,
		\item $\alpha(z\tau) = z^{-1}\alpha(\tau)$ for all $\tau \in \tilde{T}_-(\tilde{G})$ and $z \in \mathbb{S}^1$,
		\item $\alpha$ is supported on finitely many connected components of $\tilde{T}_-(\tilde{G})$,
		\item on the connected component containing a given $\tau \in \tilde{T}_{\elli, -}(\tilde{M})$, which is a compact torus, the function $\lambda \mapsto \alpha(\tau_\lambda)$ on $i\mathfrak{a}^*_M / \Stab(\tau)$ is Paley--Wiener.
	\end{itemize}
\end{definition}

For each $M \in \mathcal{L}(M_0)$, let $\mathrm{PW}_{\asp}(\tilde{G})_{\tilde{M}}$ be the subspace of functions supported on $W^G_0 \cdot \tilde{T}_{\elli,-}(\tilde{M})$. Then $\mathrm{PW}_{\asp}(\tilde{G}) = \bigoplus_{M \in \mathcal{L}(M_0)/W^G_0} \mathrm{PW}_{\asp}(\tilde{G})_{\tilde{M}}$. Restriction to $\tilde{T}_-(\tilde{M})$ induces a surjection $\mathrm{PW}_{\asp}(\tilde{G}) \twoheadrightarrow \mathrm{PW}_{\asp}(\tilde{M})^{W^G(M)}$, and an isomorphism
\[ \mathrm{PW}_{\asp}(\tilde{G})_{\tilde{M}} \rightiso \mathrm{PW}_{\asp}(\tilde{M})_{\tilde{M}}^{W^G(M)}. \]

\begin{theorem}[Trace Paley--Wiener theorem]\label{prop:trace-PW}
	\index{trace Paley--Wiener theorem}
	Assuming that $F$ is non-Archimedean, the map \eqref{eqn:PW-map} induces an isomorphism $\orbI_{\asp}(\tilde{G}) \otimes \mes(G) \rightiso \mathrm{PW}_{\asp}(\tilde{G})$. Moreover, the diagram
	\[\begin{tikzcd}
		\orbI_{\asp}(\tilde{G}) \otimes \mes(G) \arrow[r, "\sim"] \arrow[d, "{f_{\tilde{G}} \mapsto f_{\tilde{M}}}"'] & \mathrm{PW}_{\asp}(\tilde{G}) \arrow[d, "\text{restriction}"] \\
		\orbI_{\asp}(\tilde{M})^{W^G(M)} \otimes \mes(M) \arrow[r, "\sim"] & \mathrm{PW}_{\asp}(\tilde{M})^{W^G(M)}
	\end{tikzcd}\]
	commutes for each $M \in \mathcal{L}(M_0)$.
\end{theorem}

In particular, the commutative diagram implies that $\orbI_{\asp, \cusp}(\tilde{G}) \otimes \mes(G) \rightiso \mathrm{PW}_{\asp}(\tilde{G})_{\tilde{G}}$. We refer to \cite[pp.78--79]{Li19} for more discussions and references on the trace Paley--Wiener theorem for coverings.

By the discussion above and the decomposition of $\mathrm{PW}_{\asp}(\tilde{G})$ into various $\mathrm{PW}_{\asp}(\tilde{G})_{\tilde{M}}$, we obtain an isomorphism denoted as follows:
\begin{equation}\label{eqn:orbI-gr}\begin{aligned}
		\orbI_{\asp}(\tilde{G}) \otimes \mes(G) & \rightiso \bigoplus_{M \in \mathcal{L}(M_0)/W^G_0} \orbI_{\asp, \cusp}(\tilde{M})^{W^G(M)} \otimes \mes(M) \\
		f_{\tilde{G}} & \mapsto f_{\tilde{G}, \gr}.
\end{aligned}\end{equation}

Let us move to the stable version as in \S\ref{sec:spectral-distributions}. Namely, we take $G^!$ to be a direct product of groups of the form $\GL(m)$ or $\SO(2m+1)$, and fix a minimal Levi subgroup $M^!_0$ of $G^!$.

Define the stable avatar $S\mathrm{PW}(G^!) \subset \left\{ \text{functions}\; \beta: \Phi_{\mathrm{bdd}}(G^!) \to \CC \right\}$ of the Paley--Wiener spaces of Definition \ref{def:PW-spaces}, namely by replacing $\tilde{T}_-(\cdot)$ (resp.\ $\tilde{T}_{\elli, -}(\cdot)$) by $\Phi_{\mathrm{bdd}}(\cdot)$ (resp.\ $\Phi_{2, \mathrm{bdd}}(\cdot)$) and removing the condition about $\mathbb{S}^1$-action.

\begin{theorem}[Stable trace Paley--Wiener theorem, cf.\ {\cite[IV.2.3]{MW16-1}} for the Archimedean counterpart]\label{prop:stable-PW}
	Assuming $F$ is non-Archimedean, the map above induces $S\orbI(G^!) \otimes \mes(G^!) \rightiso S\mathrm{PW}(G^!)$, and it fits into commutative diagrams similar to those in Theorem \ref{prop:trace-PW}.
\end{theorem}
\begin{proof}
	In view of the local Langlands correspondence for $G^!$ established in \cite[Chapter 6]{Ar13}, this is an instance of the general formalism of \cite[\S\S 4---5]{Ar96}.
\end{proof}

\begin{remark}\label{rem:real-PW}
	Consider now the case of Archimedean $F$. The spaces $\tilde{T}_-(\tilde{G})$, etc.\ and the virtual tempered characters $f_{\tilde{G}} \mapsto f_{\tilde{G}}(\tau)$ are defined as before. However, the connected components of $\tilde{T}_-(\tilde{G})$ are now Euclidean spaces by Remark \ref{rem:archimedean-lambda-twist}, and there are two trace Paley--Wiener theorems, one for $\orbI_{\asp}(\tilde{G})$ and the other for $\orbI_{\asp}(\tilde{G}, \tilde{K})$, where $\tilde{K} = \rev^{-1}(K)$ and $K \subset G(F)$ is a maximal compact subgroup in good position relative to $M_0$; their Paley--Wiener spaces are different. Let us sketch briefly.
	\begin{itemize}
		\item 	The $\tilde{K} \times \tilde{K}$-finite version, denoted by $\mathrm{PW}_{\asp}(\tilde{G})$, is defined in the same way as in the non-Archimedean case, namely any $\alpha \in \mathrm{PW}_{\asp}(\tilde{G})$ is supported on finitely many connected components, and is Paley--Wiener over each of them.
		
		Note that $\mathrm{PW}_{\asp}(\tilde{G})$ does not depend on the choice of $\tilde{K}$.
		
		For the $\tilde{K} \times \tilde{K}$-finite trace Paley--Wiener theorem, assume that $\tilde{G}$ is a direct product of $\GL(m, \R)$ and the metaplectic groups $\Mp(2m)$. In other words, $\tilde{G}$ is of \emph{metaplectic type} to be introduced in Definition \ref{def:metaplectic-type}. This is all what we need in this work. This theorem asserts that $f_{\tilde{G}} \mapsto [\tau \mapsto f_{\tilde{G}}(\tau)]$ induces an isomorphism
		\[ \orbI_{\asp}(\tilde{G}, \tilde{K}) \otimes \mes(G) \rightiso \mathrm{PW}_{\asp}(\tilde{G}) \]
		compatibly with parabolic descent.
		
		\item The $C^\infty_c$-version $\mathrm{PW}_{\asp}(\tilde{G})^\infty$ is again a space of functions $\alpha: \tilde{T}_-(\tilde{G}) \to \CC$. The difference is that instead of being supported on finitely many connected components, we only require $\alpha(\tau)$ to have rapid decay in terms of the infinitesimal character of $\tau$; the precise definition can be found in \cite[IV.1.3]{MW16-1}. The map $f_{\tilde{G}} \mapsto [\tau \mapsto f_{\tilde{G}}(\tau)]$ induces an isomorphism
		\[ \orbI_{\asp}(\tilde{G}) \otimes \mes(G) \rightiso \mathrm{PW}_{\asp}(\tilde{G})^\infty. \]
		This holds for general coverings.
		\index{PWaspinfty@$\mathrm{PW}_{\asp}(\tilde{G})^\infty$}
		\item In either case, one realizes the Paley--Wiener space as a limit of Fréchet spaces with explicitly given semi-norms. The isomorphisms above are then homeomorphisms.
	\end{itemize}
	
	As for the proofs, the $\tilde{K} \times \tilde{K}$-finite case is done in \cite[\S 3.4]{Li14b}, which includes the groups of metaplectic type and some other examples. The $C^\infty_c$-case is due to \cite{Bo94b}, which applies to general coverings; see also \cite[IV.1]{MW16-1}.
\end{remark}

\begin{remark}
	When $G^!$ is a quasisplit connected reductive $F$-group which is a direct product of groups of the form $\GL(m)$ or $\SO(2m+1)$, we also have trace Paley--Wiener theorems for $S\orbI(G^!, K^!)$ and $S\orbI(G^!)$, with target spaces $S\mathrm{PW}(G^!)$ and $S\mathrm{PW}(G^!)^\infty$ defined in the similar way. For the proofs in the Archimedean case, we refer to \cite[IV.2.3 and IV.3.4]{MW16-1}.
\end{remark}

Finally, all the trace Paley--Wiener theorems above extend to spaces with subscript ``ac'' reviewed in \S\ref{sec:orbital-integrals}; see \cite[\S 3]{Li14b} for details about the metaplectic case.

\section{Galois cohomology}\label{sec:Galois-cohomology}
Let $F$ be a field with separable closure $\overline{F}$. For an algebraic $F$-group $I$, we denote by $\Hm^i(F, I) = \Hm^i(\Gamma_F, I)$ the Galois cohomologies of $I$, where $i \in \{0, 1\}$. They are pointed sets, functorial in $I$.

Assume $I$ to be connected reductive. For any diagonalizable $F$-subgroup $T \subset I$, denote by $T_{\mathrm{sc}}$ its preimage under $I_{\mathrm{SC}} \to I$. Following Borovoi and Labesse \cite{Lab99}, we define the complex of diagonalizable groups
\begin{equation}\label{eqn:Iab}
	I_{\mathrm{ab}} = [I_{\mathrm{SC}} \to I] := [Z_{G, \mathrm{sc}} \to Z_G]
\end{equation}
in degrees $-1$, $0$; it is also quasi-isomorphic to $T_{\mathrm{sc}} \to T$ for any maximal $F$-subtorus $T \subset I$. We then define
\begin{align*}
	\Hm^i_{\mathrm{ab}}(F, I) & := \Hm^i(F, I_{\mathrm{ab}}) = \Hm^i(F, I_{\mathrm{SC}} \to I ) \\
	& = \Hm^i(\Gamma_F, \underbracket{T_{\mathrm{sc}} \to T}_{\deg = -1, 0} ) ,
\end{align*}
the abelianized Galois cohomology groups, for all $i \in \Z_{\geq 0}$. It is functorial in $I$. There are also variants with $\Gamma_F$ replaced by $\Weil{F}$, or with $I$ replaced by some dual objects.
\index{HmiF@$\Hm^i(F, I)$, $\Hm^i_{\mathrm{ab}}(F, I)$}

The canonical \emph{abelianization maps} are $\mathrm{ab}_I^i: \Hm^i(F, I) \to \Hm^i_{\mathrm{ab}}(F, I)$, where $i \in \{0, 1\}$.

When $F$ is a global field and $I$ is connected reductive, one defines
\begin{align*}
	\Ker^i(F, I) & := \Ker\left[ \Hm^i(F, I) \to \prod_v \Hm^i(F_v, I) \right], \quad i \in \{0, 1\}, \\
	\Ker^i_{\mathrm{ab}}(F, I) & := \Ker\left[ \Hm_{\mathrm{ab}}^i(F, I) \to \prod_v \Hm_{\mathrm{ab}}^i(F_v, I) \right], \quad i \in \Z.
\end{align*}
\index{KeriF@$\Ker^i(F, I)$, $\Ker^i_{\mathrm{ab}}(F, I)$}

One also defines the abelian groups
\[ \Hm^i_{\mathrm{ab}}(\A_F, \cdot) = \prod'_v \Hm^i_{\mathrm{ab}}(F_v, \cdot), \quad \Hm^i_{\mathrm{ab}}(\A_F/F, \cdot), \quad \Hm^i_{\mathrm{ab}}(\mathfrak{o}_{F_v}, \cdot), \]
where
\begin{compactitem}
	\item in the last definition, the groups are assumed to be unramified at the place $v \nmid \infty$,
	\item the restricted product $\prod'_v$ is taken with respect to $\Hm^i_{\mathrm{ab}}(\mathfrak{o}_{F_v}, \cdot)$ for almost all $v$
\end{compactitem}
For the precise definitions, we refer to \cite{Lab99} or \cite[\S 3.1]{Li15}.
\index{HmiA@$\Hm^i_{\mathrm{ab}}(\A_F, \cdot)$, $\Hm^i_{\mathrm{ab}}(\A_F/F, \cdot)$}

\begin{example}
	The following special instance will be needed later on. Take a quadratic extension $E|F$ of fields with characteristic $\neq 2$, and the unitary group $I$ for a non-trivial hermitian or anti-hermitian form relative to $E|F$, then it is well-known (see \cite[Exemple 3.1.2]{Li15}) that
	\begin{align*}
		\Hm^1_{\mathrm{ab}}(F, I) & = \{\pm 1\}, \quad F:\; \text{local}, \\
		\Hm^1_{\mathrm{ab}}(\A_F/F, I) & = \{\pm 1\}, \quad F:\; \text{global}.
	\end{align*}
	On the other hand, these $\Hm^1_{\mathrm{ab}}$ are trivial when $I = \GL(m)$ for some $m$. The same holds when if we replace $I$ by a restriction of scalars thereof.
\end{example}

For affine $F$-groups $I \subset H$, we define the pointed set
\[ \mathfrak{D}(I, H; F) := \Ker\left[ \Hm^1(F, I) \to \Hm^1(F, H) \right]. \]
It is finite when $F$ is local.
\index{DIH@$\mathfrak{D}(I, H; \cdot)$}

Assuming that $I$ and $H$ are both connected reductive over any field $F$, we define
\[ \mathfrak{E}(I, H; F) := \Ker\left[ \Hm^1_{\mathrm{ab}}(F, I) \to \Hm^1_{\mathrm{ab}}(F, H) \right], \]
and $\mathrm{ab}^1$ induces a map $\mathfrak{D}(I, H; F) \to \mathfrak{E}(I, H; F)$.
\index{EIH@$\mathfrak{E}(I, H; \cdot)$}

Now assume $F$ is global. We have the adélic variants $\mathfrak{D}(I, H; \A_F)$ and $\mathfrak{E}(I, H; \A_F)$. Moreover, set
\begin{align*}
	\Hm_{\mathrm{ab}}^0(\A_F/ F, I \to H) & := \Hm^0(\A_F/F, \underbracket{I_{\mathrm{ab}} \to H_{\mathrm{ab}}}_{\deg = -1, 0} ), \\
	\mathfrak{E}(I, H; \A_F / F ) & := \Coker\left[ \Hm_{\mathrm{ab}}^0(\A_F, H) \to \Hm_{\mathrm{ab}}^0(\A_F/F, I \to H) \right] ,
\end{align*}
where $I_{\mathrm{ab}}$ and $H_{\mathrm{ab}}$ are viewed as complexes of diagonalizable groups, as in \eqref{eqn:Iab}, and $I_{\mathrm{ab}} \to H_{\mathrm{ab}}$ means the mapping cone complex.

By \cite[Proposition 1.8.4]{Lab99}, there is a natural exact sequence
\begin{equation}\label{eqn:E-local-global}
	\mathfrak{E}(I, H; F) \to \mathfrak{E}(I, H; \A_F) \to \mathfrak{E}(I, H; \A_F/F).
\end{equation}

\begin{theorem}[{\cite[Corollaire 1.8.6]{Lab99}}]\label{prop:Langlands-obstruction}
	Let $I \subset H$ be connected reductive groups over a global field $F$. Let $B(I, H)$ be the kernel of $\Ker^1(F, I) \to \Ker^1(F, H)$. It fits into the exact sequence
	\[ 1 \to B(I, H) \to \mathfrak{D}(I, H; F) \to \mathfrak{D}(I, H; \A_F) \to \mathfrak{E}(I, H; \A_F/F). \]
\end{theorem}

Define $\mathfrak{R}(I, H; F)$ (resp.\ $\mathfrak{R}(I, H; F)_1$) to be the Pontryagin dual of $\Hm_{\mathrm{ab}}^0(\A_F/F, I \to H)$ (resp.\ of $\mathfrak{E}(I, H; \A_F/F)$). Hence $\mathfrak{R}(I, H; F)_1 \subset \mathfrak{R}(I, H; F)$.
\index{RIH@$\mathfrak{R}(I, H; F)$, $\mathfrak{R}(I, H; F)_1$}

\begin{corollary}\label{prop:Langlands-obs}
	A class in $\mathfrak{D}(I, H; \A_F)$ comes from $\mathfrak{D}(I, H; F)$ if and only if it is annihilated by $\mathfrak{R}(I, H; F)_1$.
\end{corollary}

The pointed sets above are related to stable conjugacy (Definition \ref{def:stable-conjugacy}) in the following manner. For every $F$ and $\eta \in G(F)_{\mathrm{ss}}$, there is a canonical bijection
\[ \left\{ \eta' \in G(F): \text{stably conjugate to}\; \eta \right\} \big/ \text{conj} \xleftrightarrow{1:1} \mathfrak{D}(G_\eta, G; F) \]
mapping the $G(F)$-conjugacy class of $\eta[y] := y^{-1} \eta y$ to the class of the $1$-cocycle $c_y: \Gamma_F \to G_\eta(\overline{F})$ given by $c_y(\sigma) = y \sigma(y)^{-1}$.

We proceed to refine this parametrization as follows.

\begin{definition}\label{def:Y-set}
	\index{Yeta@$\mathcal{Y}(\eta)$, $\dot{\mathcal{Y}}(\eta)$}
	Let $G$ be a connected reductive $F$-group. Given $\eta \in G(F)_{\mathrm{ss}}$, define the pointed set
	\[ \mathcal{Y}(\eta) = \mathcal{Y}^G(\eta) := \left\{ y \in G(\overline{F}) : \forall \sigma \in \Gamma_F, \; y\sigma(y)^{-1} \in G_\eta(\overline{F}) \right\}, \]
	the base point being $y=1$. It is acted by $G_\eta(\overline{F})$ (resp.\ $G(F)$) on the left (resp.\ right), so that
	\[\begin{tikzcd}[row sep=tiny]
		G_\eta(\overline{F}) \backslash \mathcal{Y}(\eta) \arrow[leftrightarrow, r, "1:1"] & \left\{ \eta' \in G(F)_{\mathrm{ss}}: \text{stably conjugate to}\; \eta \right\} \\
		y \arrow[mapsto, r] & {\eta[y]} := y^{-1}\eta y .
	\end{tikzcd}\]

	This gives rises to
	\begin{align*}
		\dot{\mathcal{Y}}(\eta) & = \dot{\mathcal{Y}}^G(\eta) := G_\eta(\overline{F}) \backslash \mathcal{Y}(\eta) / G(F) \\
		& \xrightarrow{1:1} \left\{ \text{stable conjugates of}\; \eta \right\} \big/ G(F)\text{-conj.}
	\end{align*}

	Note that $G_{\eta[y]}$ and $G_\eta$ are related by the pure inner twist $\Ad(y)$, for all $y \in \mathcal{Y}(\eta)$.
\end{definition}

\begin{lemma}[See {\cite[I.5.11]{MW16-1}}]
	When $M \subset G$ is a Levi subgroup and $\eta \in M(F)_{\mathrm{ss}}$, the evident map $\dot{\mathcal{Y}}^M(\eta) \to \dot{\mathcal{Y}}(\eta)$ is injective.
\end{lemma}

In fact, each $y \in \mathcal{Y}(\eta)$ defines a $1$-cocycle $c_y: \Gamma_F \to G_\eta(\overline{F})$ alluded to above. In this way, we have the following commutative diagram of pointed sets
\begin{equation}\label{eqn:Y-inclusion}\begin{tikzcd}
		\dot{\mathcal{Y}}(\eta) \arrow[leftrightarrow, r, "1:1"] & \mathfrak{D}(G_\eta, G; F) \arrow[phantom, r, "\subset" description] & \Hm^1(F, G_\eta) \\
		\dot{\mathcal{Y}}^M(\eta) \arrow[leftrightarrow, r, "1:1"] \arrow[u] & \mathfrak{D}(M_\eta, M; F) \arrow[u] \arrow[phantom, r, "\subset" description] & \Hm^1(F, M_\eta) \arrow[u]
\end{tikzcd}\end{equation}
whose vertical arrows are induced by $M \to G$. It remains to observe that $\Hm^1(F, M_\eta) \to \Hm^1(F, G_\eta)$ is injective since $M_\eta$ is a Levi subgroup of $G_\eta$.

\begin{remark}
	In later applications, $G$ will take the form $\prod_{i \in I} \GL(n_i) \times \Sp(W)$ for some finite set $I$, $(n_i)_{i \in I} \in \Z_{\geq 1}^I$ and a symplectic $F$-vector space $(W, \lrangle{\cdot|\cdot})$. For such groups, we have
	\begin{itemize}
		\item $G_\eta = Z_G(\eta)$ for all $\eta \in G(F)_{\mathrm{ss}}$, hence stable conjugacy reduces to geometric conjugacy;
		\item $\Hm^1(F, G)$ is trivial.
	\end{itemize}
	The Levi subgroups of $G$ also take a similar form, thus satisfy the properties above. Hence the inclusions in \eqref{eqn:Y-inclusion} are equalities.
\end{remark}
\chapter{Endoscopy for the metaplectic group}\label{sec:endoscopy-Mp}
For a local field $F$ of characteristic zero endowed with a non-trivial additive character $\psi$, we will review the metaplectic coverings
\[ 1 \to \bmu_8 \to \tilde{G} \to G(F) \to 1, \quad G := \Sp(W) \]
where $(W, \lrangle{\cdot|\cdot})$ is a symplectic $F$-vector space, and more generally the class of covering groups of metaplectic type, which includes all Levi subgroups of $\Mp(W)$. These will mostly be phrased in terms of abstract linear algebra. Afterwards, we will review the endoscopy, the notion of transfer $\Trans_{\mathbf{G}^!, \tilde{G}}: \orbI_{\asp}(\tilde{G}) \otimes \mes(G) \to S\orbI(G^!) \otimes \mes(G^!)$, and the metaplectic fundamental lemma for the anti-genuine spherical Hecke algebra due to \cite{Luo18}.

Besides these summaries, we also supply the following new materials.
\begin{itemize}
	\item A precise notion of ``diagrams'' in \S\ref{sec:diagram}, describing the correspondence between semisimple conjugacy classes.
	\item A combinatorial summation formula in \S\ref{sec:combinatorial-summation}, which will play a key role throughout the stabilization process.
	\item For Archimedean $F$, a characterization of the image of transfer in \S\ref{sec:geom-transfer} and the tempered local character relations in \S\ref{sec:spectral-transfer}.
	\item A definition of the spaces $D_{\mathrm{geom}, -}(\tilde{G})$, $D_{\mathrm{geom}, -}(\tilde{G}, \mathcal{O})$ of geometric distributions in \S\ref{sec:geom-dist-transfer}, where $\mathcal{O}$ is a finite union of semisimple conjugacy classes in $G(F)$, as well as their behavior under transfer.
\end{itemize}

Although the formalism pertains to the global case as well, the review of adélic coverings of metaplectic type is deferred to \S\S\ref{sec:adelic-coverings}--\ref{sec:metaplectic-type-adelic}.

\section{More linear algebra}\label{sec:linear-algebra}
Our aim here is to fix notations for symplectic groups, following \S\ref{sec:linear-algebra-prelim}. Consider a field $F$ with $\mathrm{char}(F) \neq 2$. Fix a symplectic vector space $(W, \lrangle{\cdot|\cdot})$ and consider the symplectic group $\Sp(W)$.

A decomposition $W = \ell \oplus \ell'$ with $\ell, \ell'$ being Lagrangians is called a \emph{polarization} of $W$. Using $\lrangle{\cdot|\cdot}$, we identify $\ell'$ with the dual of $\ell$, and $\GL(\ell)$ (or $\GL(\ell')$) can be identified as a subgroup of $\Sp(W)$.

By a \emph{symplectic basis} for $(W, \lrangle{\cdot|\cdot})$, we mean an ordered basis $e_{-n}, \ldots, e_{-1}, e_1, \ldots, e_n$ of $W$ such that $\bigoplus_{i=1}^n F e_i$ and $\bigoplus_{i=1}^n Fe_{-i}$ are Lagrangians, and $\lrangle{e_i | e_{-j}} = \delta_{i,j}$ for all $1 \leq i, j \leq n$. Symplectic bases exist, and by fixing such a basis, we may write $\Sp(2n)$ instead of $\Sp(W)$, if necessary. 
\index{symplectic basis}
\index{Sp(W)@$\Sp(W)$, $\Sp(2n)$}

The parabolic subgroups of $\Sp(W)$ are in bijection with isotropic flags in $W$, say
\[ \{0\} = V_0 \subsetneq V_1 \subsetneq \cdots \subsetneq V_r , \]
which prolongs into the flag $\cdots \subsetneq V_r^\perp \subsetneq \cdots \subsetneq V_0^\perp = W$; the corresponding parabolic subgroup is simply its stabilizer in $\Sp(W)$. Its Levi quotient is $\prod_{i=1}^r \GL(\ell_i) \times \Sp(W^\flat)$ where $\ell_i := V_i/V_{i-1}$ and $W^\flat := V_r^\perp/V_r$.

The Levi subgroups are in bijection with orthogonal decompositions $W = W^\flat \oplus \bigoplus_{i=1}^r (\ell_i \oplus \ell'_i)$ where $W^\flat$ and $\ell_i \oplus \ell'_i$ (unordered) are symplectic subspaces and $\ell_i$, $\ell'_i$ give a polarization for $\ell_i \oplus \ell'_i$ for each $i$. To each decomposition is attached the subgroup $\prod_{i=1}^r \GL(\ell_i) \times \Sp(W^\flat)$. These bijections are $\Sp(W)$-equivariant in the obvious sense.

In particular, by choosing a symplectic basis, we obtain the standard Borel pair $(B, T)$ where $T$ corresponds to the decomposition $W = \bigoplus_{i=1}^n (Fe_i \oplus Fe_{-i})$, and $B$ corresponds to the isotropic flag $\{0\} \subsetneq Fe_1 \subsetneq \cdots \subsetneq \bigoplus_{i=1}^n Fe_i$. The Weyl group of $T$ is identified with $\mathfrak{S}_n \ltimes \{\pm 1\}^n$, where
\begin{itemize}
	\item $\mathfrak{S}_n$ acts on $\{\pm 1\}^n$ (resp.\ on $T \simeq \Gm^n$) by permutation;
	\item $(z_i)_{i=1}^n \in \{\pm 1\}^n$ acts on $T \simeq \Gm^n$ by mapping $(t_i)_{i=1}^n$ to $\left( t_i^{z_i} \right)_{i=1}^n$.
\end{itemize}

Using the isomorphism $T \simeq \Gm^n$ above, we obtain the basis $e^*_1, \ldots, e^*_n$ for $X^*(T)$. The dual basis of $X_*(T)$ is denoted as $e_{1, *}, \ldots, e_{n, *}$. The action of $\mathfrak{S_n} \ltimes \{\pm 1\}^n$ on them has a similar description.
\index{ei@$e_i^*$, $e_{i, *}$}

Now consider the orthogonal counterparts. Consider a quadratic $F$-vector space $(V, q)$. A decomposition $V = \ell \oplus \ell'$ with $\ell, \ell'$ being Lagrangians is called a polarization. As before, parabolic subgroups of $\SO(V, q)$ are in equivariant bijection with isotropic flags in $V$. The Levi subgroups are in equivariant bijection with orthogonal decompositions $V = V^\flat \oplus \bigoplus_{i=1}^r (\ell_i \oplus \ell'_i)$ where $V^\flat$ and $\ell_i \oplus \ell'_i$ (unordered) are quadratic subspaces, polarized by $\ell_i$, $\ell'_i$ for each $i$. The corresponding Levi subgroup is then $\prod_{i=1}^r \GL(\ell_i) \times \SO(V^\flat, q^\flat)$ with $q^\flat := q|_{V^\flat}$, which is minimal exactly when $(V^\flat, q^\flat)$ is anisotropic and $\dim \ell_i = 1$ for all $i$.

\index{SO-odd@$\SO(2n+1)$}
Let $n \in \Z_{\geq 0}$. The split group $\SO(2n+1)$ is associated with the quadratic vector space $(V, q)$ endowed with the basis
\[ v_1, \ldots, v_n, v_0, v_{-n}, \ldots, v_{-1} \]
such that $q$ is specified by
\[ q(v_i|v_{-j}) = \delta_{i, j}, \quad -n \leq i, j \leq n. \]
The data $(V^\flat, q^\flat)$ in the description of Levi subgroups also take the same form. In particular, our choice of basis yields the standard Borel pair $(B_{\SO}, T_{\SO})$, such that
\begin{itemize}
	\item the Weyl group of $T_{\SO}$ is identified with $\mathfrak{S}_n \ltimes \{\pm 1\}^n$, and $T_{\SO} \simeq \Gm^n$ canonically;
	\item we also have the standard bases
	\[ e^*_{1, \SO}, \ldots, e^*_{n, \SO} \quad \text{for}\; X^*(T_{\SO}), \]
	and its dual basis $e_{1, \SO, *}, \ldots, e_{n, \SO, *}$ for $X_*(T_{\SO})$.
\end{itemize}

The following result is thus evident.

\begin{proposition}
	Fix $n \in \Z_{\geq 1}$. Let $W$ (resp.\ $V$) be a $2n$-dimensional symplectic $F$-vector space (resp.\ the quadratic $F$-vector space giving rise to the split $\SO(2n+1)$). There is a natural bijection $P \leftrightarrow P_{\SO}$ (resp.\ $M \leftrightarrow M_{\SO}$) between conjugacy classes of parabolic subgroups (resp.\ Levi subgroups) of $\Sp(W)$ and $\SO(V, q)$, such that if $M \simeq \prod_{i \in I} \GL(n_i) \times \Sp(W^\flat)$, then $M_{\SO} \simeq \prod_{i \in I} \GL(n_i) \times \SO(V^\flat, q^\flat)$, where $(n_i)_{i \in I} \in \Z_{\geq 1}^I$ satisfies
	\[ \frac{1}{2} \dim W^\flat = n - \sum_{i \in I} n_i = \frac{1}{2} \left( \dim V^\flat - 1 \right) . \]
	The groups $W(M)$ and $W(M_{\SO})$ are also identified under this bijection: as groups of outer automorphisms of $\prod_{i \in I} \GL(n_i)$, both are generated by the transpose-inverse of $\GL(n_i)$ for various $i \in I$, together with permutations of factors of the same size.

	The maximal tori $T$ and $T_{\SO}$ are canonically isomorphic with respect to the given bases, compatibly with the identification of their Weyl groups $\mathfrak{S}_n \ltimes \{\pm 1\}^n$ .
\end{proposition}

The subgroups $B$, $T$ (resp.\ $B_{\SO}$, $T_{\SO}$) arising in this way are \emph{standard}. Following the general practice, we also denote these maximal tori as $M_0$; i.e.\ as the chosen minimal Levi subgroup.

The comparison of half-sum of positive roots is standard. In $X^*(T) \otimes \Q$ and $X^*(T_{\SO}) \otimes \Q$, respectively, we have
\begin{equation}\label{eqn:rho-comparison}\begin{aligned}
	\rho_B & = ne^*_1 + (n-1)e^*_2 + \cdots + e^*_n , \\
	\rho_{B_{\SO}} & = \left( n - \frac{1}{2} \right) e^*_{1, \SO} + \left(n - \frac{3}{2} \right) e^*_{2, \SO} + \cdots + \frac{1}{2} e^*_{n, \SO}.
\end{aligned}\end{equation}

The study of weighted orbital integrals will require positive-definite invariant quadratic forms on the spaces $\mathfrak{a}_M$. This is straightforward in our concrete setting.

\begin{definition}\label{def:invariant-quadratic-form}
	With the data above, we fix the integral quadratic form on $X^*(T)$ given by $X_1^2 + \cdots + X_n^2$, where $X_i$ is the $i$-th coordinate on $X^*(T) \simeq \Z^n$; it is Weyl-invariant and positive-definite. This induces quadratic forms on $\mathfrak{a}_M$, for each Levi subgroup $M$ of $\Sp(W)$, in a compatible manner.

	Ditto for $X^*(T_{\SO})$ and the Levi subgroups of $\SO(V)$.
\end{definition}

The standard choice of invariant quadratic forms is also available for groups of the form
\[ \prod_{i \in I} \GL(n_i) \times \Sp(W), \quad \prod_{i \in I} \GL(n_i) \times \SO(V), \]
by embedding them as Levi subgroups of a larger symplectic (resp.\ odd special orthogonal) group.

\section{Description of semisimple classes}\label{sec:ss-classes}
Let $F$ be a field of characteristic $\neq 2$, and let $(W, \lrangle{\cdot|\cdot})$ be a symplectic $F$-vector space. The following reproduces the parametrization of semisimple conjugacy classes in $\Sp(W)$ in \cite[\S 3.1]{Li11}. They are described by data $(K, K^\sharp, x, W_K, h_K, W_\pm, \lrangle{\cdot|\cdot}_\pm)$, where
\begin{itemize}
	\item $K$ is an étale finite-dimensional $F$-algebra endowed with an involution $\tau$;
	\item $K^\sharp$ is the $\tau$-fixed subalgebra of $K$;
	\item $x \in K^\times$, $\tau(x) = x^{-1}$, $K=F[x]$, and we also assume $x^2 - 1 \in K^\times$;
	\item $W_K$ is a faithful $K$-module and $h_K$ is an anti-hermitian form on $W_K$ with respect to $\tau$;
	\item $(W_\pm, \lrangle{\cdot|\cdot}_\pm)$ are symplectic vector spaces.
\end{itemize}
They are subject to $\dim_F W_K + \dim_F W_+ + \dim_F W_- = \dim_F W$. The invertibility of $x^2 - 1$ was missing in \cite{Li11}; if this condition is omitted, one has to absorb $(W_\pm, \lrangle{\cdot|\cdot}_\pm)$ into $(W_K, h_K)$.

In this case $K^\sharp$ is also an étale $F$-algebra, thus is a product $\prod_{i \in I} K_i^\sharp$ of fields. Accordingly, $K = \prod_{i \in I} K_i$ such that for each $i$, either
\begin{compactenum}[(a)]
	\item $K_i$ is a quadratic extension of $K_i^\sharp$ with $\Gal(K_i/K_i^\sharp) = \{\identity, \tau|_{K_i} \}$, or
	\item $K_i \simeq K_i^\sharp \times K_i^\sharp$ with $\tau(u,v)=(v,u)$.
\end{compactenum}
In both cases $K_i/K_i^\sharp$ determines $\tau|_{K_i}$, hence $(K, K^\sharp)$ determines $\tau$. By definition, $(W_K, h_K)$ also decomposes into $\prod_{i \in I} (W_i, h_i)$ such that each $h_i$ is anti-hermitian with respect to $(K_i, \tau)$.

Isomorphisms between these data are defined using the equivalences between $(K, K^\sharp)$, $x$, and the corresponding forms. To each datum we attach a symplectic form
\[ \left( \Tr_{K|F} \circ h_K, \lrangle{\cdot|\cdot}_+ , \lrangle{\cdot|\cdot}_- \right): \left(W_K \oplus W_+ \oplus W_-\right)^2 \to F. \]
By reasons of dimension, there exists an isomorphism
\[ \iota: (W, \lrangle{\cdot|\cdot}) \rightiso \left(W_K \oplus W_+ \oplus W_- , \left( \Tr_{K|F} \circ h_K, \lrangle{\cdot|\cdot}_+ , \lrangle{\cdot|\cdot}_- \right) \right). \]

For every $x \in K^{\times}$ with $x\tau(x) = 1$, the automorphism $m_x: t \mapsto xt$ of $W_K$ preserves $h_K$. Therefore
\[ \delta := \iota^{-1} (m_x, \identity, -\identity) \iota \;\in \Sp(W). \]
The resulting conjugacy of $\delta$ in $\Sp(W)$ is semisimple and independent of $\iota$.

\begin{proposition}\label{prop:ss-parameters}
	The procedure above sets up a bijection between the set of data
	\[ (K, K^\sharp, x, W_K, h_K, W_\pm, \lrangle{\cdot|\cdot}_\pm) \quad \text{modulo}\; \simeq \]
	and $\Gamma_{\mathrm{ss}}(\Sp(W))$. Furthermore, by choosing $\iota$ as above, we have an isomorphism of $F$-groups
	\[\begin{tikzcd}[row sep=tiny]
		\mathrm{U}(W_K, h_K) \times \Sp(W_+) \times \Sp(W_-) \arrow[r, "\sim"] & Z_{\Sp(W)}(\delta) \\
		 (y, g_+, g_-) \arrow[mapsto, r] & \iota^{-1} (m_y, g_+, g_-) \iota .
	\end{tikzcd}\]
	Here we regard $\mathrm{U}(W_K, h_K)$ as a connected reductive $F$-group by Weil restriction from various $K_i^\sharp$.
\end{proposition}

In particular $Z_{\Sp(W)}(\delta)$ is connected for all semisimple $\delta$, which is well-known. The datum $(K, K^\sharp, \ldots)$ is called the parameter of $\delta$.

\begin{notation}\label{nota:hyperbolic-i}
	In the decomposition $K = \prod_{i \in I} K_i$, the indices with $K_i \simeq K_i^\sharp \times K_i^\sharp$ (resp.\ $K_i$ is a field) are called hyperbolic (resp.\ elliptic). The norm map is
	\begin{align*}
		N_{K/K^\sharp}: K & \to K^\sharp \\
		t & \mapsto t\tau(t).
	\end{align*}
	The norm-one $F$-torus is, modulo Weil restrictions,
	\[ K^1 := \left\{ y \in K^\times: N_{K/K^\sharp}(y)=1 \right\}. \]
\end{notation}

\begin{corollary}\label{prop:parameter-elliptic}
	Let $\delta \in \Gamma_{\mathrm{ss}}(\Sp(W))$. Then $\delta$ is elliptic if and only if its parameter involves no hyperbolic indices in $K$.
\end{corollary}
\begin{proof}
	It suffices to look at $Z_{\Sp(W)}(\delta)$. It is a product of symplectic groups, Weil restrictions of unitary groups (taken over elliptic $i \in I$) and those of general linear groups (taken over hyperbolic $i \in I$).
\end{proof}

Consider now the regular semisimple case, i.e.\ when $Z_{\Sp(W)}(\delta)$ is a torus. The corresponding parameter is characterized by $W_\pm = \{0\}$, $W_K = K$, so $h_K$ is represented by $c \in K^\times$ such that $\tau(c) = -c$. Two data $(K, K^\sharp, x, c)$ and $(L, L^\sharp, y, d)$ are isomorphic if and only if there is an isomorphism $\varphi: K \rightiso L$ such that the involutions are preserved, $\varphi(x) = y$, and $\varphi(c)d^{-1} \in N_{L/L^\sharp}(L^\times)$. The following is thus an immediate consequence of Proposition \ref{prop:ss-parameters}.

\begin{proposition}\label{prop:reg-ss-parameters}
	The recipe in Proposition \ref{prop:ss-parameters} gives a bijection between isomorphism classes of data $(K, K^\sharp, x, c)$ and $\Gamma_{\mathrm{reg}}(\Sp(W))$. If $(K, K^\sharp, x, c)$ is the parameter of $\delta$, then $Z_{\Sp(W)}(\delta) \simeq K^1$.
\end{proposition}

There is a similar notion of isomorphisms among data $(K, K^\sharp, c)$. Given such a datum, $j: y \mapsto \iota^{-1} m_y \iota$ embeds $K^1$ as a maximal $F$-torus $j(K^1) \subset \Sp(W)$; its conjugacy class is independent of the choice of $\iota$.

\begin{proposition}\label{prop:parameter-tori}
	The mapping $(K, K^\sharp, c) \mapsto j$ yields a bijection between the $(K, K^\sharp, c)$ modulo equivalence and the embeddings of maximal $F$-tori $j: T \hookrightarrow \Sp(W)$ modulo conjugacy.
\end{proposition}

\begin{remark}\label{rem:parameter-classical-group}
	Semisimple classes in classical groups can be parameterized by the same paradigm. For example, let $(V, q)$ be a quadratic vector space of odd dimension. The semisimple conjugacy classes in $\SO(V, q)$ are parameterized by data
	\[ \left( K, K^\sharp, x, V_K, h_K, V_\pm, q_\pm \right), \]
	the only differences being that $(V_K, h_K)$ is hermitian instead of anti-hermitian, $(V_\pm, q_\pm)$ are quadratic vector spaces, $\dim V_+$ (resp.\ $\dim V_-$) is odd (resp.\ even), and we suppose that there exists
	\[ \iota: (V, q) \rightiso \left(V_K \oplus V_+ \oplus V_- , \left( \Tr_{K|F} \circ h_K, q_+ , q_- \right) \right). \]
	Then $\delta := \iota^{-1} (m_x, \identity, -\identity) \iota \in \SO(V, q)$. We also have
	\[ Z_{\mathrm{O}(V, q)}(\delta) \simeq \mathrm{U}(V_K, h_K) \times \mathrm{O}(V_+, q_+) \times \mathrm{O}(V_-, q_-) . \]
	
	Regular semisimple classes correspond to the case $V_K \simeq K$, $V_- = \{0\}$ and $\dim V_+ = 1$, so that their parameters reduce to $(K, K^\sharp, x, c)$ with $\tau(c) = c$.
\end{remark}

\begin{remark}\label{rem:parameter-elliptic-classical-group}
	For the parameterization of semisimple classes in classical groups, if $\delta$ is elliptic then $K$ involves no hyperbolic indices. The converse holds for symplectic groups (Corollary \ref{prop:parameter-elliptic}), but not true in general, as one sees in the case of orthogonal groups.
\end{remark}

\section{Groups of metaplectic type: local case}\label{sec:metaplectic-type}
Hereafter, $F$ is a local field of characteristic zero and $\psi$ is an additive character of $F$. The constructions reviewed below can be found in \cite[\S 2]{Li11}.

Given a symplectic $F$-vector space $(W, \lrangle{\cdot|\cdot})$, we construct the central extension of locally compact groups
\[ 1 \to \CC^\times \to \MMp_\psi(W) \to \Sp(W) \to 1 \]
using the Schrödinger models of irreducible representations of the Heisenberg group $\mathrm{H}(W)$ with central character $\psi$. This is not the group that we work with. Instead, we reduce $\MMp_\psi(W)$ to
\[ 1 \to \bmu_8 \to \Mp(W) \xrightarrow{\rev} \Sp(W) \to 1, \]
using the fact that Weil's constants $\gamma_\psi$, which figure in the Maslov cocycle describing $\MMp_\psi(W)$, take value in $\bmu_8$. The eightfold coverings $\rev: \Mp(W) \twoheadrightarrow \Sp(W)$ in the sense of \S\ref{sec:covering} are the \emph{metaplectic groups} considered in this work.
\index{SpWtilde@$\Mp(W)$}
\index{metaplectic group}

\begin{itemize}
	\item The group $\Mp(W)$ carries the \emph{Weil representation} $\omega_\psi$, which is genuine, unitary and admits a canonical orthogonal decomposition into irreducibles, called the even ($+$) and odd ($-$) pieces of $\omega_\psi$:
	\[ \omega_\psi = \omega_\psi^+ \oplus \omega_\psi^- . \]
	\item There is a canonical central element of order two in $\rev^{-1}(-1)$, hereafter denoted as $-1$; see \cite[Définition 2.8]{Li11}. We shall write
	\begin{equation}\label{eqn:minus-1-lifting}
		-\tilde{\delta} := (-1) \cdot \tilde{\delta}, \quad \tilde{\delta} \in \Mp(W).
	\end{equation}
	\item If $F \neq \CC$ and $W \neq \{0\}$, then $\Mp(W)$ can be reduced to a twofold covering
	\[ 1 \to \bmu_2 \to \Mp(W)^{(2)} \to \Sp(W) \to 1. \]
	In fact, $\Mp(W)^{(2)}$ is the derived subgroup of $\Mp(W)$, and this is the non-trivial twofold covering of $\Sp(W)$, unique up to a unique isomorphism.
	\item If $F = \CC$ or $W = \{0\}$, then the central extension splits canonically: $\Mp(W) \simeq \Sp(W) \times \bmu_8^\times$. Here we define $\Sp(\{0\})$ to be the trivial group.
\end{itemize}

Note that $\Mp(W)$ depends on $\psi \circ \lrangle{\cdot|\cdot}$; same for $\omega_\psi$ and the $-1\in \Mp(W)$. On the other hand, $\Mp(W)^{(2)}$ is insensitive to $\psi$ and is the one preferred by many authors.

\begin{definition}[{\cite[Définition 3.1.1]{Li12a}}]\label{def:metaplectic-type}
	\index{group of metaplectic type}
	By a \emph{group of metaplectic type}, we mean a covering of the form
	\[ 1 \to \bmu_8 \to \prod_{i \in I} \GL(n_i, F) \times \Mp(W^\flat) \xrightarrow{\identity \times \rev} \prod_{i \in I} \GL(n_i, F) \times \Sp(W^\flat) \to 1 \]
	where $\rev: \Mp(W^\flat) \twoheadrightarrow \Sp(W^\flat)$ is the metaplectic covering associated with a symplectic $F$-vector space $(W^\flat, \lrangle{\cdot|\cdot})$ and $\psi$.
\end{definition}

The relevance of groups of metaplectic type is explained by the following
\begin{proposition}[See {\cite[\S 5.4]{Li11}}]
	Consider a Levi subgroup $M = \prod_{i \in I} \GL(n_i) \times \Sp(W^\flat)$ of $\Sp(W)$. Then $\tilde{M}$ is canonically isomorphic to $\prod_{i \in I} \GL(n_i, F) \times \Mp(W^\flat)$ as coverings. The isomorphism can be characterized in terms of the values of the character of $\omega_\psi$ at regular semisimple elements.
\end{proposition}

The formalism of orbital integrals and invariant distributions from \S\S\ref{sec:orbital-integrals}---\ref{sec:distributions} simplifies enormously for groups of metaplectic type, because of the following fact.

\begin{proposition}\label{prop:auto-good}
	Let $\tilde{G}$ be a group of metaplectic type. Then every element of $G(F)$ is good in the sense of Definition \ref{def:good}.
\end{proposition}
\begin{proof}
	This reduces immediately to the well-known case $\tilde{G} = \Mp(W)$ in \cite[II.5 Lemme]{MVW87}.
\end{proof}

Now comes the unramified setting. Suppose that $F$ is non-Archimedean with residual characteristic $p > 2$. Let $L \subset W$ be an $\mathfrak{o}_F$-lattice. Define its dual lattice as
\[ L^* := \left\{ x \in W: \forall y \in L, \;\psi(\lrangle{x|y}) = 1 \right\}. \]
We say that $L$ is self-dual with respect to $\psi\lrangle{\cdot|\cdot}$ if $L^* = L$. The choice of $L$ affords an $\mathfrak{o}_F$-model for $\Sp(W)$. When $L$ is self-dual, we know that
\begin{itemize}
	\item $K := \Sp(W, \mathfrak{o}_F)$ is hyperspecial,
	\item the lattice model for $\omega_\psi$ yields a splitting $s: K \to \Mp(W)$ of $\rev$,
	\item the image of $-1 \in K$ coincides with the $-1$ defined before \eqref{eqn:minus-1-lifting}.
\end{itemize}

The datum $(\rev, K, s)$ verifies the ``unramified conditions'' in the sense of \cite[Définition 3.1.1]{Li14a}. We will return to the unramified situation in \S\ref{sec:LF}.

\section{Endoscopic data}\label{sec:endoscopic-data}
Consider a local or global field $F$ of characteristic zero and fix an additive character $\psi$. Given a symplectic $F$-vector space $(W, \lrangle{\cdot|\cdot})$ of dimension $2n$, we put $G := \Sp(W)$. We can form
\begin{itemize}
	\item the metaplectic covering $1 \to \bmu_8 \to \tilde{G} \xrightarrow{\rev} G(F) \to 1$ in the local case, or
	\item the covering $\tilde{G} \xrightarrow{\rev} G(\A_F)$ in the global case, which is to be reviewed in \S\ref{sec:metaplectic-type-adelic}.
\end{itemize}
However, we will not deal with adélic points in this section.

In order to execute inductive arguments, the definitions below will also apply to groups of metaplectic type.

\begin{definition}
	\index{endoscopic data}
	\index{Endoell@$\Endo_{\elli}(\tilde{G})$}
	The \emph{elliptic endoscopic data} of $\tilde{G}$ are triplets $(G^!, n', n'')$ where
	\begin{gather*}
		(n', n'') \in \Z_{\geq 0}^2, \quad n' + n'' = n, \\
		G^! := \SO(2n'+1) \times \SO(2n''+1).
	\end{gather*}
	We call $G^!$ the corresponding \emph{endoscopic group}. Define $\Endo_{\elli}(\tilde{G})$ to be the set of elliptic endoscopic data of $\tilde{G}$. It is sometimes convenient to designate the elements of $\Endo_{\elli}(\tilde{G})$ simply as $(n', n'')$.
	
	More generally, for $\tilde{M} = \prod_{i \in I} \GL(n_i, F) \times \Mp(W^\flat)$ with $\dim W^\flat = 2n^\flat$, we set
	\[ \Endo_{\elli}(\tilde{M}) := \left\{ (M^!, n', n'') : n' + n'' = n^\flat \right\}, \]
	where $M^! := \prod_{i \in I} \GL(n_i) \times \SO(2n'+1) \times \SO(2n''+1)$ is the corresponding endoscopic group.
	
	\index{Gshrekbold@$\mathbf{G}^{"!}$}
	The endoscopic data will also be typeset in boldface, such as $\mathbf{G}^! \in \Endo_{\elli}(\tilde{G})$, and use the normal font $G^!$ to denote the corresponding endoscopic group.
\end{definition}

The definition above is compatible with the fact that $\GL(n_i)$ has no non-tautological elliptic endoscopic data. Interpretations on the dual side will be given in \S\ref{sec:dual}.

The general endoscopic data are defined to be elliptic endoscopic data of Levi subgroups, taken up to conjugacy.

\begin{definition}
	\index{EndoG@$\Endo(\tilde{G})$}
	We set
	\[ \Endo(\tilde{G}) := \left( \bigsqcup_{M \in \mathcal{L}(M_0)} \Endo_{\elli}(\tilde{M}) \right) \big/ W^G_0. \]
	Likewise, we can define $\Endo(\tilde{L})$ whenever $\tilde{L}$ is a group of metaplectic type.
\end{definition}

Next, we recall the correspondence between semisimple elements as follows. Fix an elliptic endoscopic datum $(n', n'')$ for $\tilde{G} = \Mp(W)$ and consider $G^! = \SO(2n' + 1) \times \SO(2n''+1)$.
\begin{itemize}
	\item Take the standard (as in \S\ref{sec:linear-algebra}) maximal tori $T' \subset \SO(2n'+1)$ (resp.\ $T'' \subset \SO(2n'' + 1)$) and $T'_{\Sp} \subset \Sp(2n')$ (resp.\ $T''_{\Sp} \subset \Sp(2n'')$). Take the standard maximal torus $T \subset G$.
	\item Denote by $W'_{\SO}$ (resp.\ $W''_{\SO}$), $W'_{\Sp}$ (resp.\ $W''_{\Sp}$) and $W_{\Sp}$ their Weyl groups, respectively.
\end{itemize}
Then we have $T' \simeq T'_{\Sp}$ (resp.\ $T'' \simeq T''_{\Sp}$), equivariantly with respect to $W'_{\SO} \simeq W'_{\Sp}$, as well as $\iota: T'_{\Sp} \times T''_{\Sp} \rightiso T$ equivariantly with respect to $W'_{\Sp} \times W''_{\Sp} \hookrightarrow W_{\Sp}$.

The standard maximal tori depend on symplectic bases, but different choices of bases yield conjugate data.

\begin{definition}\label{def:corr-orbits}
	For $(n', n'')$ as above, the composite map
	\[\begin{tikzcd}
		\Psi: T' \times T'' \arrow[r, "\sim"] & T'_{\Sp} \times T''_{\Sp} \arrow[rr, "{(g',g'') \mapsto \iota(g', -g'')}" inner sep=0.8em] & & T
	\end{tikzcd}\]
	which is equivariant with respect to $W'_{\SO} \times W''_{\SO} \rightiso W'_{\Sp} \times W''_{\Sp} \hookrightarrow W_{\Sp}$, induces a finite morphism between adjoint quotients of Chevalley, whence a canonical map
	\[ \Sigma_{\mathrm{ss}}(G^!) \to \Sigma_{\mathrm{ss}}(G). \]
	We shall write $\gamma^! \leftrightarrow \gamma$ to indicate that the stable class of $\gamma^!$ maps to that of $\gamma$.
\end{definition}

\begin{remark}
	The correspondence can be made explicit by eigenvalues. If $\gamma^! = (\gamma', \gamma'')$ where $\gamma'$ (resp.\ $\gamma''$) has eigenvalues $(a'_1)^{\pm 1}, \ldots, (a'_n)^{\pm 1}, 1$ (resp.\ $(a''_1)^{\pm 1}, \ldots, (a''_n)^{\pm 1}, 1$) over $\overline{F}$, then $\gamma^! \leftrightarrow \gamma$ if and only if $\gamma$ has eigenvalues $(a'_1)^{\pm 1}, \ldots, (a'_n)^{\pm 1}, -(a''_1)^{\pm 1}, \ldots, -(a''_n)^{-1}$ over $\overline{F}$.
\end{remark}

The correspondence extends by the following routine procedure.
\begin{compactitem}
	\item One gets a correspondence for elliptic endoscopic data when $\tilde{G}$ is replaced by a group of metaplectic type $\tilde{M} = \prod_{i \in I} \GL(n_i, F) \times \Mp(W^\flat)$, by taking the tautological correspondence on each $\GL$ factor.
	\item We deduce a correspondence for all endoscopic data for $\tilde{G}$ as follows. Suppose that the endoscopic datum comes from $\mathbf{M}^! \in \Endo_{\elli}(\tilde{M})$ where $M$ is a Levi subgroup of $G$. Then we take the map
	\[ \Sigma_{\mathrm{ss}}(M^!) \xrightarrow{\text{endoscopy}} \Sigma_{\mathrm{ss}}(\tilde{M}) \xrightarrow{Levi} \Sigma_{\mathrm{ss}}(G). \]
\end{compactitem}

Next, consider a Levi subgroup $M \subset G$. As reviewed in \S\ref{sec:linear-algebra}, to $M$ is attached a symplectic $F$-vector subspace $W^\flat \subset W$. Consider $\mathbf{M}^! \in \Endo_{\elli}(\tilde{M})$; the situation will also be summarized as
$\begin{tikzcd}
	M^! \arrow[dashed, leftrightarrow, r, "\text{ell.}", "\text{endo.}"'] & \tilde{M}
\end{tikzcd}$.
We are interested in completing the diagram
\[\begin{tikzcd}
	& \tilde{G} \\
	M^! \arrow[dashed, leftrightarrow, r, "\text{ell.}", "\text{endo.}"'] & \tilde{M} \arrow[hookrightarrow, u, "\text{Levi}"']
\end{tikzcd} \quad \text{into} \quad \begin{tikzcd}
	G^! \arrow[dashed, leftrightarrow, r, "\text{ell.}", "\text{endo.}"'] & \tilde{G} \\
	M^! \arrow[dashed, leftrightarrow, r, "\text{ell.}", "\text{endo.}"'] \arrow[hookrightarrow, u, "\text{Levi}"] & \tilde{M} \arrow[hookrightarrow, u, "\text{Levi}"']
\end{tikzcd} . \]

Suppose that $\mathbf{M}^!$ arises from the pair $(m', m'')$. Writing
\[ M = \prod_{i \in I} \GL(n_i) \times \Sp(W^\flat), \quad M^! = \prod_{i \in I} \GL(n_i) \times \SO(2m' + 1) \times \SO(2m'' + 1), \]
those $\mathbf{G}^! \in \Endo_{\elli}(\tilde{G})$, or more precisely the pair $(n', n'')$, together with the embedding $M^! \hookrightarrow G^!$ can thus be obtained from ordered partitions $I = I' \sqcup I''$, where
\[ I' := \left\{ i \in I: \GL(n_i) \hookrightarrow \SO(2n'+1) \right\}, \quad I'' := \left\{ i \in I: \GL(n_i) \hookrightarrow \SO(2n''+1) \right\}. \]

\begin{remark}
	In the situation above, the identification of $\GL$-factors gives an embedding
	\[ W^{G^!}(M^!) \hookrightarrow W^G(M), \]
	and an equivariant isomorphism
	\[ \mathfrak{a}_{M^!} \simeq \mathfrak{a}_M . \]
	The isomorphism is compatible with the invariant quadratic forms induced from $G^!$ and $G$, respectively; see Definition \ref{def:invariant-quadratic-form}.
\end{remark}

\begin{definition}\label{def:Endo-s}
	\index{EndoMshrekG@$\Endo_{\mathbf{M}^{"!}}(\tilde{G})$}
	Given a Levi subgroup $M \subset G$ and $\mathbf{M}^! \in \Endo_{\elli}(\tilde{M})$ as above, write $M = \prod_{i \in I} \GL(n_i) \times \Sp(W^\flat)$ and define
	\[ \Endo_{\mathbf{M}^!}(\tilde{G}) := \left\{ \text{ordered partitions}\; I = I' \sqcup I'' \right\}. \]
	To each $s \in \Endo_{\mathbf{M}^!}(\tilde{G})$, we have the corresponding element of $\Endo_{\elli}(\tilde{G})$ is denoted by $\mathbf{G}^![s]$. The reason for denoting an ordered partition by $s$ will be made clear in Proposition \ref{prop:Endo-s-interpretation}.
\end{definition}

Therefore, each $s \in \Endo_{\mathbf{M}^!}(\tilde{G})$ gives rise to the situation
\begin{equation}\label{eqn:s-situation}\begin{tikzcd}
	G^![s] \arrow[dashed, leftrightarrow, r, "\text{ell.}", "\text{endo.}"'] & \tilde{G} \\
	M^! \arrow[dashed, leftrightarrow, r, "\text{ell.}", "\text{endo.}"'] \arrow[hookrightarrow, u, "\text{Levi}"] & \tilde{M} \arrow[hookrightarrow, u, "\text{Levi}"']
\end{tikzcd}. \end{equation}
Recall that the correspondence of semisimple elements is obtained by going from $M^!$ to $M$ and then to $G$ in the diagram.
%along the way
%\begin{tikzpicture}[baseline=(O), scale=0.5]
%	\draw[->] (0,0) -- (1,0) -- (1, 0.7);
%	\coordinate (O) at (0, 0);
%\end{tikzpicture}.

\begin{definition}[The central twist]\label{def:central-twist}
	\index{central twist}
	\index{gamma-s@$\gamma[s]$}
	\index{z-s@$z[s]$}
	Suppose $s \in \Endo_{\mathbf{M}^!}(\tilde{G})$ corresponds to $I = I' \sqcup I''$. Define
	\[ z[s] := \left( (z_i)_{i \in I}, 1, 1 \right) \in Z_{M^!}(F), \quad
	z_i := \begin{cases} 1, & i \in I' \\ -1, & i \in I'' . \end{cases} \]
	For each $\gamma \in M^!(F)$, we write $\gamma[s] := \gamma \cdot z[s]$.
\end{definition}

\begin{remark}
	More generally, when $\tilde{G}$ is of metaplectic type and $\mathbf{M}^! \in \Endo_{\elli}(\tilde{M})$, the set $\Endo_{\mathbf{M}^!}(\tilde{G})$ and the central elements $z[s] \in Z_{M^!}(F)$ are also defined as before. The recipe is routine: decompose $\tilde{G}$ into $\prod_{j \in J} \GL(h_j, F) \times \Mp(W)$ and then reduce to the case of $\Mp(W)$, by leaving the factors $\GL(h_j)$ intact.
\end{remark}

The main relevance of the twist $\gamma \mapsto \gamma[s]$ in this work is explained as follows.
\begin{proposition}[See {\cite[Proposition 3.3.4]{Li12a}}]\label{prop:central-twist-corr}
	Consider a group of metaplectic type $\tilde{G}$, its Levi subgroup $\tilde{M}$ and $\mathbf{M}^! \in \Endo_{\elli}(\tilde{M})$. For all $s \in \Endo_{\mathbf{M}^!}(\tilde{G})$, the assignments
	\begin{equation*}
		\begin{tikzcd}
			\Sigma_{\mathrm{ss}}(G^![s]) \arrow[r, "{\text{via}\; \mathbf{G}^![s]}" inner sep=0.6em] & \Sigma_{\mathrm{ss}}(G) \\
			\Sigma_{\mathrm{ss}}(M^!) \arrow[u, "\text{Levi}"] & \\
			\gamma \arrow[phantom, u, "\in" description, sloped] &
		\end{tikzcd} \quad \text{and} \quad \begin{tikzcd}
			& \Sigma_{\mathrm{ss}}(G) \\
			\Sigma_{\mathrm{ss}}(M^!) \arrow[r, "{\text{via}\; \mathbf{M}^!}"' inner sep=0.6em] & \Sigma_{\mathrm{ss}}(M) \arrow[u, "\text{Levi}"'] \\
			{\gamma[s]} \arrow[phantom, u, "\in" description, sloped] &
	\end{tikzcd}\end{equation*}
	have the same image.
\end{proposition}

\section{Dual groups}\label{sec:dual}
Consider a metaplectic covering $\tilde{G}$ as in \S\ref{sec:endoscopic-data}, associated with $(W, \lrangle{\cdot|\cdot})$ where $\dim W = 2n$.

\begin{definition}\label{def:dual-group}
	\index{Gtildevee@$\tilde{G}^\vee$, $\Lgrp{\tilde{G}}$}
	The \emph{dual group} of $\tilde{G}$ is defined to be $\tilde{G}^\vee := \Sp(2n, \CC)$, equipped with the trivial $\Gamma_F$-action. Define $\Lgrp{\tilde{G}} := \tilde{G}^\vee \times \Gamma_F$.
	
	More generally, for $\tilde{M} = \prod_{i \in I} \GL(n_i, F) \times \Mp(W^\flat)$ with $\dim W^\flat = 2n^\flat$, we set
	\[ \tilde{M}^\vee := \prod_{i \in I} \GL(n_i, \CC) \times \Sp(2n^\flat, \CC), \quad \Lgrp{\tilde{M}} := \tilde{M}^\vee \times \Gamma_F . \]
\end{definition}

For any Levi subgroup $M \subset G$, we obtain a corresponding embedding $\tilde{M}^\vee \hookrightarrow \tilde{G}^\vee$ and $\Lgrp{\tilde{M}} \hookrightarrow \Lgrp{\tilde{G}}$ of standard Levi subgroups.

Elliptic endoscopic data are interpreted on the dual side via the bijection
\begin{equation}\label{eqn:endoscopy-dual-bijection-0}\begin{tikzcd}
	\Endo_{\elli}(\tilde{G}) \arrow[r, "1:1"] & \left\{ s \in \tilde{G}^\vee_{\mathrm{ss}} : s^2 = 1 \right\} \big/ \text{conj} .
\end{tikzcd}\end{equation}
Specifically, to each $(n' ,n'') \in \Endo_{\elli}(\tilde{G})$ we attach an element $s = s_{n' , n''} \in \tilde{G}^\vee$ with eigenvalues $1$ and $-1$, with multiplicities $2n'$ and $2n''$ respectively. We can actually take $s$ in $T^\vee$, the standard maximal torus in $\Sp(2n, \CC)$. If $\mathbf{G}^!$ corresponds to $s$, then $\check{G}^! = Z_{\tilde{G}^\vee}(s) \simeq \Sp(2n', \CC) \times \Sp(2n'', \CC)$ and we obtain $\check{G}^! \hookrightarrow \tilde{G}^\vee$, canonically up to $\tilde{G}^\vee$-conjugacy.

In general, for a group of metaplectic type $\tilde{M}$, a semisimple element $s \in \tilde{M}^\vee$ is called \emph{elliptic} if $C := Z_{\tilde{M}^\vee}(s)$ is not contained in any proper Levi of $\tilde{M}^\vee$. This notion depends only on $s Z_{\tilde{M}^\vee}$ modulo conjugation. In the setting of \eqref{eqn:endoscopy-dual-bijection-0}, a semisimple $s \in \tilde{G}^\vee$ is elliptic if and only if $s^2 = 1$. Hence the bijection generalizes as follows.

\begin{lemma}
	Let $\tilde{M} := \prod_{i \in I} \GL(n_i, F) \times \tilde{M}^\flat$, where $M^\flat = \Mp(W^\flat)$ for some symplectic vector space $W^\flat$. Then $(n' ,n'') \mapsto s_{n', n''} \in (\widetilde{M^\flat})^\vee \subset \tilde{M}^\vee$ yields a bijection
	\[\begin{tikzcd}
		\Endo_{\elli}(\tilde{M}) \arrow[r, "1:1"] & Z_{\tilde{M}^\vee}^\circ \big\backslash \left\{ s \in \tilde{M}^\vee_{\mathrm{ss}} : \mathrm{elliptic} \right\} \big/ \text{conj} .
	\end{tikzcd}\]
	If $\mathbf{M}^! \in \Endo_{\elli}(\tilde{M})$ arises from $s$, then we obtain $\check{M}^! \rightiso Z_{\tilde{M}^\vee}(s) \subset \tilde{M}^\vee$ which is canonical up to conjugacy.
\end{lemma}

In view of the paradigm for linear reductive groups, it would be more reasonable to replace $Z_{\tilde{M}^\vee}^\circ$ by $Z_{\tilde{M}^\vee}^{\Gamma_F, \circ}$ in the formulas above. Since the $\Gamma_F$-action are trivial here, we omit these superscripts.

For every $s \in \tilde{G}^\vee_{\mathrm{ss}}$, take $\tilde{M}^\vee$ to be the centralizer of the maximal central torus in $Z_{\tilde{G}^\vee}(s)$. Then $s$ is elliptic in $\tilde{M}^\vee$. Hence we obtain the following parametrization.

\begin{proposition}\label{prop:endoscopic-data-vs-dual}
	The bijection in \eqref{eqn:endoscopy-dual-bijection-0} extends to a surjection
	\[ \tilde{G}^\vee_{\mathrm{ss}} \twoheadrightarrow \Endo(\tilde{G}). \]
	The same assertion also extends to the case when $\tilde{G}$ is of metaplectic type.
\end{proposition}

\begin{remark}\label{rem:centerless}
	Compared to the endoscopy of reductive groups, the metaplectic feature here is that $Z_{\tilde{G}^\vee}$ does not impose any symmetries on endoscopic data. Consider the split $\SO(2n+1)$ for example: it has the same dual group as $\tilde{G}^\vee$, with trivial $\Gamma_F$-action, but its elliptic endoscopic data are in bijection with \emph{unordered} pairs $(n', n'')$ with $n' + n'' = n$, i.e.\ with elliptic semisimple classes in $\Sp(2n, \CC)$ modulo $\{\pm 1\}$. The slogan here is that one should replace every occurrence of $Z_{M^\vee}^{\Gamma_F}$ in Arthur's theory by $Z_{\tilde{M}^\vee}^{\circ}$, when one wants to deal with a group $\tilde{M}$ of metaplectic type.
	
	Following this heuristic, the metaplectic version $\mathrm{Out}(\mathbf{G}^!) = \mathrm{Out}^{\tilde{G}}(\mathbf{G}^!)$ of endoscopic outer automorphisms of $G^!$ should be trivial when $\mathbf{G}^! \in \Endo_{\elli}(\tilde{G})$. In general, we should define $\mathrm{Out}^{\tilde{G}}(\mathbf{M}^!)$ to be $W^G(M)$ if $\mathbf{M}^! \in \Endo_{\elli}(\tilde{M})$ where $M$ is a Levi subgroup of $G$.
\end{remark}

Using Proposition \ref{prop:endoscopic-data-vs-dual}, the Definition \ref{def:Endo-s} of $\Endo_{\mathbf{M}^!}(\tilde{G})$ can also be understood from the dual picture, which is Arthur's original recipe: see \cite[\S 4]{Ar98} or \S\ref{sec:endoscopy-linear}. Suppose that $\mathbf{M}^! \in \Endo_{\elli}(\tilde{M})$ arises from an elliptic semisimple $s^\flat \in \tilde{M}^\vee$. Consider the map
\begin{equation}\label{eqn:Endo-s-dual}\begin{tikzcd}[row sep=tiny]
	s^\flat Z_{\tilde{M}^\vee}^\circ \big/ Z_{\tilde{G}^\vee}^\circ \arrow[r] & \Endo(\tilde{G}) \\
	s = s^\flat t \arrow[mapsto, r] & \text{the $\mathbf{G}^![s]$ associated with } s^\flat t.
\end{tikzcd}\end{equation}
The quotient by $Z_{\tilde{G}^\vee}^\circ$ might appear superfluous; we add it merely in order to accommodate the groups of metaplectic type.

Given $\mathbf{M}^!$, we may and do assume $s^\flat \in (\widetilde{M^\flat})^\vee$ and $(s^\flat)^2 = 1$. Write any $t \in Z_{\tilde{M}^\vee}^\circ$ as
\begin{gather*}
	t = t_+ t_- t_u, \\
	t_\pm \;\text{has eigenvalue}\; \pm 1, \quad t_u \;\text{has eigenvalues}\; \neq 1.
\end{gather*}
Then $s$ is contained in the Levi subgroup $\tilde{L}^\vee := Z_{\tilde{G}^\vee}(t_u) \supset \tilde{M}^\vee$. Suppose that $L \in \mathcal{L}^G(M)$ corresponds to $\tilde{L}^\vee \in \mathcal{L}^{\tilde{G}^\vee}(\tilde{M}^\vee)$. By construction, we are thus in the situation
\[\begin{tikzcd}
	& \tilde{G} \\
	{G^![s]} \arrow[dashed, r, "\text{endo.}", "\text{ell.}"'] \arrow[ru, "\text{endo.}", sloped] & \tilde{L} \arrow[hookrightarrow, u, "\text{Levi}"'] \\
	M^! \arrow[hookrightarrow, u, "\text{Levi}"] \arrow[dashed, r, "\text{endo.}", "\text{ell.}"'] & \tilde{M} \arrow[hookrightarrow, u, "\text{Levi}"'] .
\end{tikzcd}\]

\begin{proposition}\label{prop:Endo-s-interpretation}
	With the notation above, there is a canonical bijection
	\[\begin{tikzcd}[row sep=tiny]
		\left\{s \in s^\flat Z_{\tilde{M}^\vee}^\circ \big/ Z_{\tilde{G}^\vee}^\circ : \mathbf{G}^![s] \in \Endo_{\elli}(\tilde{G}) \right\} \arrow[r, "1:1"] & \Endo_{\mathbf{M}^!}(\tilde{G}) \\
		s^\flat t \arrow[mapsto, r] & (I', I'')
	\end{tikzcd}\]
	where $I' := \left\{ i \in I: t_i = 1 \right\}$ and $I'' := \left\{ i \in I: t_i = -1 \right\}$, where we write $t = (t_i)_{i \in I}$ according to $Z_{\tilde{M}^\vee}^\circ = (\CC^\times)^I$.
	
	The description above carries over verbatim to the case when $\tilde{G}$ is of metaplectic type.
\end{proposition}
\begin{proof}
	Indeed, $\mathbf{G}^![s]$ is elliptic if and only if $t_u = 1$, i.e.\ $t_i \in \{\pm 1\}$ for all $i \in I$. Suppose that $\mathbf{G}^![s]$ corresponds to the pair $(n', n'')$. Then the signs $t_i$ determine whether $\GL(n_i) \subset M^!$ embeds into the $\SO(2n' + 1)$ or $\SO(2n''+1)$ factor of $G^![s]$. This amounts to a partition $I = I' \sqcup I''$ and matches the definition of $\Endo_{\mathbf{M}^!}(\tilde{G})$.
	
	The generalization to groups of metaplectic type is routine.
\end{proof}

Next, we define certain constants of a combinatorial nature: the symplectic $F$-vector spaces and additive characters behind these covering will be immaterial. Specifically, given a group of metaplectic type $\tilde{G}$, what matters below is just the dual group $\tilde{G}^\vee$, with trivial Galois action; moreover, the formalism applies to any field $F$ with $\mathrm{char}(F) \neq 2$.

Given a group of metaplectic type $\tilde{G}$ and $G^! \in \Endo_{\elli}(\tilde{G})$, there is a natural embedding
\[ Z_{\tilde{G}^\vee}^\circ \hookrightarrow Z_{\check{G}^!}, \]
or equivalently $Z_{\tilde{G}^\vee}^{\Gamma_F, \circ} \hookrightarrow Z_{\check{G}^!}^{\Gamma_F}$. In contrast with the endoscopy of linear groups, one has to limit to the identity connected component of $Z_{\tilde{G}^\vee} = Z_{\tilde{G}^\vee}^{\Gamma_F}$. This is consistent with Remark \ref{rem:centerless}.

\begin{definition}[Cf.\ {\cite[\S 4.2]{Li12a}}]\label{def:i-const}
	\index{iMG@$i_{M^{"!}}(\tilde{G}, G^{"!})$}
	Let $\tilde{G}$ be a group of metaplectic type. In the situation
	\[\begin{tikzcd}
		G^! \arrow[dashed, dash, r, "\text{ell.}", "\text{endo.}"'] & \tilde{G} \\
		M^! \arrow[hookrightarrow, u, "\text{Levi}"] \arrow[dashed, dash, r, "\text{ell.}", "\text{endo.}"'] & \tilde{M} \arrow[hookrightarrow, u, "\text{Levi}"']
	\end{tikzcd}\]
	as in \eqref{eqn:s-situation}, we set
	\[ i_{M^!}(\tilde{G}, G^!) := \dfrac{\left( Z_{\check{M}^!} : Z_{\tilde{M}^\vee}^\circ \right)}{\left( Z_{\check{G}^!} : Z_{\tilde{G}^\vee}^\circ \right)}. \]
\end{definition}

\begin{lemma}\label{prop:i-transitivity}
	Let $\tilde{G}$ be a group of metaplectic type and $L \supset M$ be Levi subgroups of $G$. Let $s \in \Endo_{M^!}(\tilde{L})$ and $t \in \Endo_{L^![s]}(\tilde{G})$, so that we are in the situation
	\[\begin{tikzcd}
		G^![t] \arrow[dashed, dash, r, "\text{ell.}", "\text{endo.}"'] & \tilde{G} \\
		L^![s] \arrow[hookrightarrow, u, "\text{Levi}"] \arrow[dashed, dash, r, "\text{ell.}", "\text{endo.}"'] & \tilde{L} \arrow[hookrightarrow, u, "\text{Levi}"'] \\
		M^! \arrow[hookrightarrow, u, "\text{Levi}"] \arrow[dashed, dash, r, "\text{ell.}", "\text{endo.}"'] & \tilde{M} \arrow[hookrightarrow, u, "\text{Levi}"'] .
	\end{tikzcd}\]
	Then $i_{M^!}(\tilde{G}, G^![t]) = i_{M^!}(\tilde{L}, L^![s]) i_{L^![s]}(\tilde{G}, G^![t])$.
\end{lemma}
\begin{proof}
	Immediate from Definition \ref{def:i-const}.
\end{proof}

\section{Diagrams}\label{sec:diagram}
The notion of \emph{diagrams} enhances the correspondence of semisimple classes under endoscopy. It is used throughout \cite{Wal08, MW16-1, MW16-2}. Since the metaplectic version has not appeared in earlier works, we give a short account here.

Let $F$ be any field of characteristic zero. Let $\tilde{G}$ be a covering of metaplectic type of $G(F)$ (resp.\ $G(\A_F)$) for local (resp.\ global) $F$. Let $\mathbf{G}^! \in \Endo_{\elli}(\tilde{G})$. We fix the standard Borel pairs
\[ (B_G, T_G) \;\text{for}\; G, \quad (B_{G^!}, T_{G^!}) \;\text{for}\; G^! \]
by the recipe of \S\ref{sec:linear-algebra}. They come from the standard choices of bases for the symplectic and odd orthogonal $F$-vector spaces in question. In particular, they are defined over $F$. Denote by $W^G$ and $W^{G^!}$ the corresponding Weyl groups. There is then a natural inclusion $W^{G^!} \hookrightarrow W^G$.

Moreover, we have a standard isomorphism $\Phi: T_{G^!} \rightiso T_G$ between split $F$-tori, equivariant with respect to the actions of Weyl groups.

\begin{definition}[cf.\ {\cite[I.1.10]{MW16-1}}]\label{def:diagram}
	\index{diagrams}
	\index{xiTT@$\xi_{T^{"!}, T}$, $\tilde{\xi}_{T^{"!}, T}$}
	Let $\epsilon \in G^!(F)_{\mathrm{ss}}$ and $\eta \in G(F)_{\mathrm{ss}}$. A \emph{diagram} joining $\epsilon$ and $\eta$ is a sextuplet $(\epsilon, B^!, T^!, B, T, \eta)$ such that
	\begin{itemize}
		\item $T^!$ (resp.\ $T$) is a maximal $F$-torus of $G^!_\epsilon$ (resp.\ of $G_\eta$);
		\item $(B^!, T^!)$ and $(B, T)$ are Borel pairs over $\overline{F}$ for $G^!$ and $G$, respectively;
		\item take $h \in G^!(\overline{F})$ and $g \in G(\overline{F})$ such that
		\[ (h B^! h^{-1}, h T^! h^{-1}) = (B_{G^!}, T_{G^!}), \quad (gBg^{-1}, gTg^{-1}) = (B_G, T_G), \]
		then we have an isomorphism between $F$-tori
		\[ \xi_{T^!, T}: \Ad(g)^{-1} \Phi \Ad(h): T^! \to T; \]
		\item $\tilde{\xi}_{T^!, T} := \Ad(g)^{-1} \Psi \Ad(h)$ sends $\epsilon$ to $\eta$, where $\Psi: T_{G^!} \rightiso T_G$ is the isomorphism between $F$-varieties to be explained below.
	\end{itemize}
	Observe that $\xi_{T^!, T}$ and $\tilde{\xi}_{T^!, T}$ are independent of the choice of $(h, g)$.
\end{definition}

Recall that we have decompositions $T_G = T_{\GL} \times T_{\Sp}$ and $T_{G^!} = T_{\GL} \times T' \times T''$, according to
\begin{align*}
	G & = \prod_{i \in I} \GL(n_i) \times \Sp(W), \\
	G^! & = \prod_{i \in I} \GL(n_i) \times \SO(2n'+1) \times \SO(2n''+1).
\end{align*}

By the discussions preceding Definition \ref{def:corr-orbits}, extended easily to groups of metaplectic type, $\Phi$ decomposes as $\identity_{T_{\GL}}$ times
\[ T' \times T'' \rightiso T'_{\Sp} \times T''_{\Sp} = T_{\Sp} \]
where the maximal $F$-tori $T'_{\Sp} \subset \Sp(2n')$, $T''_{\Sp} \subset \Sp(2n'')$, the embedding $\Sp(2n') \times \Sp(2n'') \hookrightarrow \Sp(W)$ and the isomorphism above are induced by the chosen symplectic basis. The map $g'' \mapsto -g''$ makes sense in $\Sp(2n'', F)$. Define the twisted version of $\Phi$ as follows
\begin{equation*}
	\Psi: T_G = T_{\GL} \times T_{\Sp} = T_{\GL} \times T'_{\Sp} \times T''_{\Sp} \xrightarrow[\sim]{(t, g', g'') \mapsto (t, g', -g'')} T_{\GL} \times T'_{\Sp} \times T''_{\Sp} = T_{G^!}.
\end{equation*}
Then $\Psi$ is still equivariant with respect to $W^{G^!} \hookrightarrow W^G$.

Observe that $\Psi$ is the map used in the Definition \ref{def:corr-orbits} of the correspondence between semisimple classes, thus justifies the last item in Definition \ref{def:diagram}. It is not a group homomorphism.

\begin{remark}\label{rem:xi-center-diagram}
	The maps $\xi_{T^!, T}$ and $\Psi$ depend only on $B$ and $B^!$ in the diagram, not on the choices of $g$ and $h$. If $\delta \in T^!(F)$, then $\tilde{\xi}_{T^!, T}(\delta)$ corresponds to $\delta$ via endoscopy. In particular, $\epsilon$ corresponds to $\delta$ under endoscopy.
	
	Note that $\Ad(g)$, $\Ad(h)$ and $\Phi$ do not affect the center of the common factor $\prod_{i \in I} \GL(n_i)$ of $G^!$ and $G$. Hence $\xi_{T, T^!} := \xi_{T^!, T}^{-1}: T \rightiso T^!$ restricts to the isomorphism $\xi: A_G \rightiso A_{G^!}$ afforded by endoscopy.
\end{remark}

\begin{remark}\label{rem:diagram-dual}
	Another usage of a diagram $(\epsilon, B^!, T^!, B, T, \eta)$ is that it determines isomorphisms
	\[\begin{tikzcd}[row sep=tiny]
		X^*(T^!_{\overline{F}}) & X^*(T_{G^!}) \arrow[l, "\sim", "{\Ad(h)^*}"'] & X^*(T_G) \arrow[l, "\sim", "{\Phi^*}"'] \arrow[r, "\sim"', "{\Ad(g)^*}"] & X^*(T_{\overline{F}}) \\
		& X_*\left( T_{G^!}^\vee \right) \arrow[r, "\sim"'] \arrow[equal, u] & X_*(T^\vee) \arrow[equal, u] &
	\end{tikzcd}\]
	where $T_{G^!}^\vee \subset (G^!)^\vee$ and $T^\vee \subset \tilde{G}^\vee$ are parts of the root data of the dual groups.
\end{remark}

We will come back to diagrams in \S\ref{sec:diagram-2}. In particular, it will be shown in Lemma \ref{prop:diagram-reg} that $\delta \in G^!_{\mathrm{reg}}(F)$ and $\gamma \in G_{\mathrm{reg}}(F)$ correspond to each other via endoscopy, then they are joined by a diagram.

\section{A summation formula}\label{sec:combinatorial-summation}
Let $F$ be either a local or global field of characteristic zero. The conventions are the same as \S\ref{sec:endoscopic-data}. What follows will be of a combinatorial nature.

\begin{definition}\label{def:iotaGG-general}
	\index{iotaGG@$\iota(\tilde{G}, G^{"!})$}
	Let $\tilde{G}$ be of metaplectic type and $\mathbf{G}^! \in \Endo(\tilde{G})$. Define
	\[ \iota(\tilde{G}, G^!) := \left( Z_{(G^!)^\vee} : Z_{\tilde{G}^\vee}^\circ \right)^{-1}
	= \begin{cases}
		2^{- \#(\SO\text{-factors in}\; G^!)}, & \text{if}\; \mathbf{G}^! \in \Endo_{\elli}(\tilde{G}), \\
		0, & \text{otherwise}.
	\end{cases}\]
	This constant will reappear in Definition \ref{def:iota-const}. Note that $\SO(1) = \{1\}$ does not count as an $\SO$-factor here.
\end{definition}

In what follows, we consider conjugacy classes of Levi subgroups $M$ of $G$, denoted as $M/\text{conj}$, and similarly for $M^! / \text{conj}$ in $G^!$. This ambiguity can be avoided by considering only standard Levi subgroups. When parameterizing them, the $\GL$-part is described by a finite set $I$ and $(n_i)_{i \in I} \in \Z_{\geq 1}^I$, modulo bijections. These data can be rigidified by requiring $I = \{1, \ldots, h\}$ and $n_1 \geq \cdots \geq n_h$. In other words, we consider partitions $\mathbf{n} = [n_1 \geq \cdots \geq n_h]$ of various length $h$, and set $|\mathbf{n}| := \sum_i n_i$.

\begin{proposition}\label{prop:situation-easier}
	Assume that $G = \Sp(W)$ with $\dim W = 2n$. Define
	\begin{itemize}
		\item the set $\mathcal{A}$ consisting of data
		\[ (n', n''), \; (m', m'') \in \Z_{\geq 0}^2, \quad \mathbf{n}', \; \mathbf{n}'' : \text{partitions}, \]
		subject to
		\[ n' + n'' = n, \quad |\mathbf{n}'| + m' = n', \quad |\mathbf{n}''| + m'' = n'' ; \]
		\item the set $\mathcal{B}$ consisting of data
		\[ (m', m'') \in \Z_{\geq 0}^2, \quad \mathbf{n}: \text{partition}, \quad s \in \Endo_{\mathbf{M}^!}(\tilde{G}) \]
		where $\mathbf{M}^!$ is the elliptic endoscopic datum of $\tilde{M}$ corresponding to $(m', m'')$, where
		\[ M := \prod_i \GL(n_i, F) \times \Sp(2(m' + m'')) \hookrightarrow \Sp(2n) \simeq G, \]
		and we require that $|\mathbf{n}| + m' + m'' = n$.
	\end{itemize}
	\begin{enumerate}[(i)]
		\item There are natural bijections
		\begin{align*}
			\mathcal{A} & \xrightarrow{1:1} \left\{ (\mathbf{G}^!, M^! / \text{conj}) \right\}, \\
			\mathcal{B} & \xrightarrow{1:1} \left\{ (M/\text{conj}, \mathbf{M}^!, s) \right\}
		\end{align*}
		where $\mathbf{G}^! \in \Endo_{\elli}(\tilde{G})$, $\mathbf{M}^! \in \Endo_{\elli}(\tilde{M})$ and $s \in \Endo_{\mathbf{M}^!}(\tilde{G})$.
		\item There is a natural map $\mathcal{B} \to \mathcal{A}$; under this map and (i), $M^!$ is the endoscopic group underlying $\mathbf{M}^!$, the elliptic endoscopic datum $\mathbf{G}^!$ is $\mathbf{G}^![s]$, and they fit into the situation \eqref{eqn:s-situation}.
		\item The map $\mathcal{B} \to \mathcal{A}$ is surjective, and two data $(m', m'', \mathbf{n}, s)$ and $(\underline{m}', \underline{m}'', \underline{\mathbf{n}}, \underline{s})$ have the same image if and only if
		\[ m' = \underline{m}' , \quad m'' = \underline{m}'' , \quad \mathbf{n} = \underline{\mathbf{n}} \]
		and $s$, $\underline{s}$ are related by permuting those indices $i$ with the same $n_i$; in particular, they give rise to the same $M/\text{conj}$ in this case.
	\end{enumerate}
\end{proposition}
\begin{proof}
	The bijections in (i) are evident if we interpret the data in $\mathcal{A}$ and $\mathcal{B}$ (except $s$) as parameters classifying Levi subgroups up to conjugacy. The map $\mathcal{B} \to \mathcal{A}$ in (ii) is given by
	\begin{itemize}
		\item splitting the partition $\mathbf{n}$ according to $s$, viewing $s$ as an ordered partition $\{1, \ldots, h\} = I' \sqcup I''$ where $h$ is the length of $\mathbf{n}$;
		\item setting $n' := m' + |\mathbf{n}'|$ and $n'' := m'' + |\mathbf{n}''|$.
	\end{itemize}
	The remaining assertions in (ii) are all clear.
	
	To show the surjectivity of $\mathcal{B} \to \mathcal{A}$, note that $\mathbf{n}'$ and $\mathbf{n}''$ can merge into $\mathbf{n}$ so that their lengths add up. This is not yet an element of $\mathcal{B}$ since the $n_i$ for $i = 1, \ldots, h$ in $\mathbf{n}$ have to be labeled: it must come from either $\mathbf{n'}$ or $\mathbf{n}''$, accordingly $i \in I'$ or $i \in I''$. Such a labeling surely exists (eg.\ by requiring that $I'$ precedes $I''$), and leads to the description of fibers in (iii).
\end{proof}

\begin{remark}\label{rem:s-situation-variant}
	The surjectivity in (iii) implies that a situation like
	\begin{equation*}\begin{tikzcd}
			G^! \arrow[dashed, leftrightarrow, r, "\text{endo.}", "\text{ell.}"'] & \tilde{G} \\
			M^! \arrow[hookrightarrow, u, "\text{Levi}"] &
	\end{tikzcd}\end{equation*}
	can always be completed into \eqref{eqn:s-situation}; this result also generalizes to $\tilde{G}$ of metaplectic type. The non-uniqueness of $(M/\text{conj}, \mathbf{M}^!, s)$ already appears when $n = 2$, $\mathbf{G}^!$ corresponds to $(n', n'') = (1, 1)$ and $M^! = \Gm \times \Gm$ (standard Levi with $I = \{1, 2\}$, $n_1 = n_2 = 1$), in which case $M = \Gm \times \Gm$, but there are exactly two choices of $s$: either $I' = \{1\}$, $I'' = \{2\}$, or oppositely.

	A less precise form\footnote{In \textit{loc.\ cit.}, the datum $(M/\text{conj}, \mathbf{M}^!, s)$ is claimed to be unique, which is true only if the uniqueness is understood modulo $W^G(M)$. This suffices for the purposes of that work.} of this result on extensions appeared in \cite[Lemma 3.3.14]{Li19}.
\end{remark}

Next, we consider a family of complex numbers $S(\mathbf{G}^!, M^!)$ attached to all $\mathbf{G}^! \in \Endo_{\elli}(\tilde{G})$ and all conjugacy classes of Levi subgroups $M^! \subset G^!$.

\begin{proposition}[Cf.\ {\cite[Lemma 9.2]{Ar99}, \cite[Lemma 10.2]{Ar02}}]\label{prop:combinatorial-summation}
	Suppose that a complex number $S(\mathbf{G}^!, M^!)$ is prescribed for each $\mathbf{G}^! \in \Endo_{\elli}(\tilde{G})$ and each conjugacy class of Levi subgroup $M^! \subset G^!$. Then
	\begin{multline*}
		\sum_{\mathbf{G}^! \in \Endo_{\elli}(\tilde{G})} \iota(\tilde{G}, G^!) \sum_{M^! / \text{conj}} |W^{G^!}(M^!)|^{-1} S(\mathbf{G}^!, M^!) \\
		= \sum_{M / \text{conj}} |W^G(M)|^{-1} \sum_{\mathbf{M}^! \in \Endo_{\elli}(\tilde{M})} \iota(\tilde{M}, M^!) \sum_{s \in \Endo_{\mathbf{M}^!}(\tilde{G})} i_{M^!}(\tilde{G}, G^![s]) S(\mathbf{G}^![s], M^!).
	\end{multline*}
	Equivalently, one can sum over semi-standard Levi subgroups, namely
	\begin{multline*}
		\sum_{\mathbf{G}^! \in \Endo_{\elli}(\tilde{G})} \iota(\tilde{G}, G^!) \sum_{M^! \in \mathcal{L}^{G^!}(M_0^!)} \frac{|W^{M^!}_0|}{|W^{G^!}_0|} S(\mathbf{G}^!, M^!) \\
		= \sum_{M \in \mathcal{L}(M_0)} \frac{|W^M_0|}{|W^G_0|} \sum_{\mathbf{M}^! \in \Endo_{\elli}(\tilde{M})} \iota(\tilde{M}, M^!) \sum_{s \in \Endo_{\mathbf{M}^!}(\tilde{G})} i_{M^!}(\tilde{G}, G^![s]) S(\mathbf{G}^![s], M^!),
	\end{multline*}
	where we fix minimal Levi $M_0 \subset G$ and $M_0^! \subset G^!$ for each $\mathbf{G}^!$.
\end{proposition}
\begin{proof}
	The equivalence between these two assertions is well-known. It suffices to assume $G = \Sp(W)$ and treat the first version. First off, we claim that
	\begin{equation*}
		\iota(\tilde{M}, M^!) i_{M^!}(\tilde{G}, G^!) = \iota(\tilde{G}, G^!)
	\end{equation*}
	whenever $\mathbf{M}^! \in \Endo_{\elli}(\tilde{M})$ and $\mathbf{G}^! := \mathbf{G}^![s]$ for an $s \in \Endo_{\mathbf{M}^!}(\tilde{G})$. Indeed,
	\begin{align*}
		\iota(\tilde{G}, G^!) \iota(\tilde{M}, M^!)^{-1} & = \left( Z_{\check{M}^!} : Z_{\tilde{M}^\vee}^\circ \right) \left( Z_{\check{G}^!} : Z_{\tilde{G}^\vee}^\circ \right)^{-1} \\
		& = i_{M^!}(\tilde{G}, G^!).
	\end{align*}
	
	We can now use Proposition \ref{prop:situation-easier}. The sum over $(M / \text{conj}, \mathbf{M}^!, s)$ is indexed by $\mathcal{B}$, namely
	\[ \sum_{\substack{(m', m'', \mathbf{n}, s) \\ \in \mathcal{B}}} |W^G(M)|^{-1} \iota(\tilde{G}, G^![s]) S(\mathbf{G}^![s], M^!). \]
	The term indexed by $(m', m'', \mathbf{n}, s)$ depends only on its image in $\mathcal{A}$. The sum can thus be descended to a sum over $(\mathbf{G}^!, M^! / \mathrm{conj.})$, upon multiplying it by the cardinality of fibers.

	Recall that $W^G(M)$ is generated by the permutations of the $\GL$-factors in $M$ (i.e.\ in $M^!$) of the same size, together with transpose-inverses on each $\GL(n_i)$; the same holds for $W^{G^!}(M^!)$, except that only the $\GL$-factors in the same $\SO$-factor of $G^!$ are permuted. Note that $w \in W^G(M)$ acts on $\Endo_{\mathbf{M}^!}(\tilde{G})$ by permutations. Suppose that $s$ corresponds to $(I', I'')$; if the indices $i$ with the same $n_i$ are permuted within $I'$ and $I''$, then $ws = s$; if there is any mixing, then $ws \neq s$.

	Combined with Proposition \ref{prop:situation-easier} (iii), it follows that the fiber over $(\mathbf{G}^!, M^! / \mathrm{conj.})$ has cardinality $|W^G(M)| / |W^{G^!}(M^!)|$. This yields
	\[ \sum_{\mathbf{G}^!} \sum_{M^! / \text{conj.}} |W^G(M)|^{-1} \cdot \frac{|W^G(M)|}{|W^{G^!}(M^!)|} \iota(\tilde{G}, G^!) S(\mathbf{G}^!, M^!) \]
	which is the expected expression.
\end{proof}

It is probably more natural to prove Proposition \ref{prop:combinatorial-summation} by working in dual groups as in \cite[Lemma 9.2]{Ar99}, provided that we replace $Z_{G^\vee}^{\Gamma_F}$ (resp.\ $Z_{M^\vee}^{\Gamma_F}$)by $Z_{\tilde{G}^\vee}^\circ$ (resp.\ $Z_{\tilde{M}^\vee}^\circ$). A difference with the cited result is that we do not require $S(\mathbf{G}^!, M^!)$ to be invariant under endoscopic outer automorphisms. All these are in line with Remark \ref{rem:centerless}.

\section{Geometric transfer}\label{sec:geom-transfer}
Let $\tilde{G} = \Mp(W)$ where $(W, \lrangle{\cdot | \cdot})$ and $\psi$ are chosen. Choose a minimal Levi subgroup $M_0$.

Fix $\mathbf{G}^! \in \Endo_{\elli}(\tilde{G})$. A semisimple element $\delta$ of $G^!$ is called $G$-regular or $\tilde{G}$-regular if it corresponds to an element in $G_{\mathrm{reg}}$; this implies $\delta$ is also regular in $G^!$. The $G$-regular elements form an open dense subset of $G^!_{\mathrm{reg}}$. Let $\Sigma_{G\text{-reg}}(G^!) \subset \Sigma_{\mathrm{reg}}(G^!)$ be the corresponding subset.

Following \cite[\S 5.3]{Li11} (see also \cite[\S 3.2]{Li12a}), the \emph{transfer factor} is a function
\[ \Delta = \Delta_{\mathbf{G}^!, \tilde{G}}: \Sigma_{G\text{-reg}}(G^!) \times \Gamma_{\mathrm{reg}}(\tilde{G}) \to \CC \]
\index{transfer factor}
\index{Delta@$\Delta_{\mathbf{G}^{"!}, \tilde{G}}(\delta, \tilde{\gamma})$}
satisfying:
\begin{itemize}
	\item $\Delta(\delta, z\tilde{\gamma}) = z \Delta(\delta, \tilde{\gamma})$ when $z \in \bmu_8$;
	\item $\Delta(\delta, \tilde{\gamma}) \neq 0$ unless $\delta \leftrightarrow \gamma$, in which case $\Delta(\delta, \tilde{\gamma}) \in \bmu_8$;
	\item write $\Delta = \Delta_{(n', n'')}$ when $\mathbf{G}^!$ corresponds to $(n', n'')$, then
	\[ \Delta_{(n', n'')}((\delta', \delta''), \tilde{\gamma}) = \Delta_{(n'', n')}((\delta'', \delta'), -\tilde{\gamma}) \]
	where $\delta' \in \SO(2n' + 1, F)$, $\delta'' \in \SO(2n'' + 1, F)$, and $-\tilde{\gamma}$ is defined as in \eqref{eqn:minus-1-lifting};
	\item Let $M \subset G$ be a Levi subgroup and $\mathbf{M}^! \in \Endo_{\elli}(\tilde{M})$, expressed as
	\[ M = \prod_{i \in I} \GL(n_i) \times M^\flat, \quad M^! = \prod_{i \in I} \GL(n_i) \times (M^!)^\flat. \]
	Then
	\[ \Delta_{\mathbf{M}^!, \tilde{G}}\left( ((\delta_i)_i, \delta^\flat), ((\delta_i)_i, \tilde{\gamma}^\flat ) \right) =
	\Delta_{\mathbf{M}^!, \tilde{M}}\left( ((\delta_i)_i, \delta^\flat), ((\delta_i)_i, \tilde{\gamma}^\flat ) \right) =
	\Delta_{(\mathbf{M}^!)^\flat, \tilde{M}^\flat} \left( \delta^\flat, \tilde{\gamma}^\flat \right) \]
	where $(\delta_i)_i \in \prod_{i \in I} \GL(n_i, F)$, $\delta^\flat \in M^\flat(F)$ and $\tilde{\gamma}^\flat \in \tilde{M}^\flat$.
\end{itemize}

The last property can be used to define $\Delta_{\mathbf{G}^!, \tilde{G}}$ for general $\mathbf{G}^! \in \Endo(\tilde{G})$. Note that $\Delta$ depends on $\psi \circ \lrangle{\cdot | \cdot}$; in fact, the Weil constants $\gamma_\psi(\cdot)$ appear in $\Delta$ through the character of the Weil representation $\omega_\psi$.

Before stating the \emph{geometric transfer}, i.e.\ the transfer of regular semisimple orbital integrals, we recall from \cite[\S 5.5]{Li11} that given $\mathbf{G}^!$, if $(\delta, \tilde{\gamma}) \in \Sigma_{G\text{-reg}}(G^!) \times \Gamma_{\mathrm{reg}}(\tilde{G})$ satisfies $\delta \leftrightarrow \gamma$, then $G^!_\delta \simeq G_\gamma$. Therefore we can take matching Haar measures on their groups of $F$-points in the comparison of orbital integrals.

\begin{theorem}[Geometric transfer]\label{prop:geom-transfer}
	\index{transfer}
	\index{TGG@$\Trans_{\mathbf{G}^{"!}, \tilde{G}}$}
	Given $\mathbf{G}^! \in \Endo(\tilde{G})$, there exists a linear map
	\[\begin{tikzcd}[row sep=tiny]
		\Trans = \Trans_{\mathbf{G}^!, \tilde{G}}: \orbI_{\asp}(\tilde{G}) \otimes \mes(G) \arrow[r] & S\orbI(G^!) \otimes \mes(G^!) \\
		f \arrow[mapsto, r]& f^{G^!}
	\end{tikzcd}\]
	such that for all $\delta \in \Sigma_{G\text{-reg}}(G^!)$, we have
	\[ \sum_{\gamma \in \Gamma_{\mathrm{reg}}(G)} \Delta_{\mathbf{G}^!, \tilde{G}}(\delta, \tilde{\gamma}) f_{\tilde{G}}(\tilde{\gamma}) = f^{G^!}(\delta) \]
	where $\tilde{\gamma} \in \rev^{-1}(\gamma)$ is arbitrary, with the aforementioned convention on Haar measures.
	
	When $F$ is Archimedean, $\Trans$ is continuous and it restricts to
	\[ \orbI_{\asp}(\tilde{G}, \tilde{K}) \otimes \mes(G) \to S\orbI(G^!, K^!) \otimes \mes(G^!), \]
	where $K \subset G(F)$ and $K^! \subset G^!(F)$ are maximal compact subgroups.
\end{theorem}
\begin{proof}
	This is \cite[Théorème 5.20]{Li11}. First, it reduces to the case $\mathbf{G}^! \in \Endo_{\elli}(\tilde{G})$. The continuity in the Archimedean case is addressed in \cite[\S 7.1]{Li19}, which is based on the works of Adams and Renard. The $\tilde{K} \times \tilde{K}$-finite transfer in the Archimedean case is \cite[Theorem 7.4.5]{Li19}.
\end{proof}

We shall also write $\Trans = \Trans_{n', n''}$ if the endoscopic datum in question is $(n', n'') \in \Endo_{\elli}(\tilde{G})$.

\begin{definition}[Transfer of distributions]\label{def:transfer-dist}
	\index{TGGdual@$\trans_{\mathbf{G}^{"!}, \tilde{G}}$}
	Given $\mathbf{G}^!$ of $\tilde{G}$, denote the transpose of $\Trans_{\mathbf{G}^!, \tilde{G}}$ as
	\[ \trans = \trans_{\mathbf{G}^!, \tilde{G}} : SD(G^!) \otimes \mes(G^!)^\vee \to D_-(\tilde{G}) \otimes \mes(G)^\vee . \]
\end{definition}

All the results and definitions above extend to the case when $\tilde{G}$ is a group of metaplectic type. It suffices to treat the $\Mp$ and the $\GL$-factors of $\tilde{G}$ separately.

The following phenomenon (see \cite[Theorem 3.4.6]{Li19}) is specific to the metaplectic covers.

\begin{definition}\label{def:g-s}
	In the situation of \eqref{eqn:s-situation}, for any $g^{M^!} \in S\orbI(M^!)$ we define its $s$-twist $g^{M^!}[s] \in S\orbI(M^!)$ by
	\[ g^{M^!}[s](\delta) = g^{M^!}(\delta[s]), \quad \delta \in \Sigma_{\mathrm{reg}}(M^!). \]
\end{definition}

As $z[s] \in M^!(F)$ is central, so $g^{M^!} \mapsto g^{M^!}[s]$ is a well-defined involution of $S\orbI(M^!)$.

\begin{proposition}[{\cite[Theorem 3.4.6]{Li19}}]\label{prop:Levi-central-twist}
	Consider the situation \eqref{eqn:s-situation}. For all $f \in \orbI_{\asp}(\tilde{G}) \otimes \mes(G)$ we have
	\[ \Trans_{\mathbf{M}^!, \tilde{M}}(f_{\tilde{M}}) = \left( \Trans_{\mathbf{G}^![s], \tilde{G}}(f)\right)^{M^!}[s] \]
	in $S\orbI(M^!) \otimes \mes(M^!)$.
\end{proposition}

It is convenient to introduce endoscopic analogues of various objects attached to $\tilde{G}$.

\begin{align*}
	\Gamma^{\Endo}_{\mathrm{reg, ell}}(\tilde{G}) & := \bigsqcup_{\mathbf{G}^! \in \Endo_{\elli}(\tilde{G})} \Sigma_{\mathrm{reg, ell}}(G^!), \\
	\Gamma^{\Endo}_{\mathrm{reg}}(\tilde{G}) & := \bigsqcup_{M \in \mathcal{L}(M_0) / W^G_0} \Gamma^{\Endo}_{\text{$G$-reg, ell}}(\tilde{M}) \big/ W^G(M).
\end{align*}
\index{GammaEndo@$\Gamma^{\Endo}_{\mathrm{reg}}$, $\Gamma^{\Endo}_{\mathrm{reg, ell}}$}
The second definition is parallel to
\[ \Gamma_{\mathrm{reg}}(G) = \bigsqcup_{M \in \mathcal{L}(M_0) /W^G_0} \Gamma_{\text{$G$-reg, ell}}(M) \big/ W^G(M). \]

\begin{definition}\label{def:orbI-Endo}
	\index{IEndo@$\orbI^{\Endo}$, $\orbI^{\Endo}_{\cusp}$}
	Let $\orbI^{\Endo}(\tilde{G}) \subset \bigoplus_{\mathbf{G}^! \in \Endo_{\elli}(\tilde{G})} S\orbI(G^!) \otimes \mes(G^!)$ be the subspace consisting of data $f^{\Endo} = \left( f^{G^!} \right)_{\mathbf{G}^! \in \Endo_{\elli}(\tilde{G})}$ satisfying:
	\begin{itemize}
		\item for every Levi subgroup $M \subset G$ and $\mathbf{M}^! \in \Endo_{\elli}(\tilde{M})$, the element of $S\orbI(M^!) \otimes \mes(M^!)$
		\[ \left( f^{G^![s]} \right)^{M^!}[s], \quad s \in \Endo_{\mathbf{M}^!}(\tilde{G}) \]
		is independent of $s$;
		\item denote the element above as $f^{M^!}$, then $f^{M^!} \in S\orbI(M^!)^{W^G(M)} \otimes \mes(M^!)$.
	\end{itemize}
	It contains the subspace
	\[ \orbI^{\Endo}_{\cusp}(\tilde{G}) := \bigoplus_{\mathbf{G}^! \in \Endo_{\elli}(\tilde{G})} S\orbI_{\cusp}(G^!) \otimes \mes(G^!). \]
	
	When $F$ is Archimedean, we also take maximal compact subgroups to define
	\begin{align*}
		\orbI^{\Endo}(\tilde{G}, \tilde{K}) & := \orbI^{\Endo}(\tilde{G}) \cap \bigoplus_{\mathbf{G}^! \in \Endo_{\elli}(\tilde{G})} S\orbI(G^!, K^!) \otimes \mes(G^!), \\
		\orbI^{\Endo}(\tilde{G}, \tilde{K}) & := \bigoplus_{\mathbf{G}^! \in \Endo_{\elli}(\tilde{G})} S\orbI_{\cusp}(G^!, K^!) \otimes \mes(G^!)
	\end{align*}
	and all these spaces are all LF spaces.
\end{definition}

\begin{definition}\label{def:collective-transfer}
	\index{TransEndo@$\Trans^{\Endo}$, $\trans^{\Endo}$}
	In view of Proposition \ref{prop:Levi-central-twist}, we may define the \emph{collective geometric transfer} $\Trans^{\Endo}$ as
	\begin{equation*}\begin{tikzcd}[row sep=tiny]
		\orbI_{\asp}(\tilde{G}) \otimes \mes(G) \arrow[r, "{\Trans^{\Endo}}" inner sep=0.8em] & \orbI^{\Endo}(\tilde{G}) \\
		\orbI_{\asp, \cusp}(\tilde{G}) \otimes \mes(G) \arrow[phantom, u, "\subset" description, sloped] \arrow[r, "{\Trans^{\Endo}_{\cusp}}"'] & \orbI^{\Endo}_{\cusp}(\tilde{G}) \arrow[phantom, u, "\subset" description, sloped]
	\end{tikzcd}\end{equation*}
	mapping $f$ to $\left( \Trans_{\mathbf{G}^!, \tilde{G}}(f)\right)_{\mathbf{G}^! \in \Endo_{\elli}(\tilde{G})}$. When $F$ is Archimedean, it is continuous and restricts to
	\[ \Trans_{\mathbf{G}^!, \tilde{G}}: \orbI_{\asp}(\tilde{G}, \tilde{K}) \otimes \mes(G) \to \orbI^{\Endo}(\tilde{G}, \tilde{K}); \]
	ditto for the case with subscripts ``$\cusp$'' (see Theorem \ref{prop:geom-transfer}).

	Taking transpose yields the collective transfer of distributions
	\[ \trans^{\Endo}: \bigoplus_{\mathbf{G}^! \in \Endo_{\elli}(\tilde{G})} SD(G^!) \otimes \mes(G^!)^\vee \to D_-(\tilde{G}) \otimes \mes(G)^\vee . \]
	These notions extend immediately to groups of metaplectic type.
\end{definition}

Note that the elements $f^{\Endo} \in \orbI^{\Endo}(\tilde{G})$ may be seen as functions on $\Gamma^{\Endo}_{\mathrm{reg}}(\tilde{G})$: to all $\mathbf{M}^! \in \Endo_{\elli}(\tilde{M})$ and $\delta \in \Sigma_{\text{$G$-reg, ell}}(M^!) / W^G(M)$, take any $s \in \Endo_{\mathbf{M}^!}(\tilde{G})$ to define
\begin{equation}\label{eqn:fEndo-descent}
	f^{M^!} := \left( f^{G^![s]} \right)^{M^!}[s], \quad f^\Endo(\delta) := f^{M^!}(\delta).
\end{equation}

All in all, $\Trans^\Endo$ is the restriction to $\orbI_{\asp}(\tilde{G}) \otimes \mes(G)$ of the map
\begin{equation}\label{eqn:collective-transfer-fcn}\begin{tikzcd}[row sep=tiny]
	\left\{ \text{anti-genuine functions}\; \Gamma_{\mathrm{reg}}(\tilde{G}) \to \CC \right\} \arrow[r] & \left\{ \text{functions}\; \Gamma^{\Endo}_{\mathrm{reg}}(\tilde{G}) \to \CC \right\} \\
	f_{\tilde{G}} \arrow[mapsto, r] & {\left[ \delta \mapsto \displaystyle\sum_\gamma \Delta_{\mathbf{M}^!, \tilde{M}}(\delta, \tilde{\gamma}) f_{\tilde{M}}(\tilde{\gamma}) \right]},
\end{tikzcd}\end{equation}
still denoted as $\Trans^{\Endo}$, where $M \in \mathcal{L}(M_0)$, $\mathbf{M}^! \in \Endo_{\elli}(\tilde{M})$, $\delta \in \Gamma_{\text{$G$-reg, ell}}(M^!)$. In other words, $\mathcal{T}^{\Endo}$ is the operator with kernel given by transfer factors.

Note that we do not need the ``glued'' version $S\orbI(\mathbf{G}^!)$ of $S\orbI(G^!)$ in \cite[I.2.5]{MW16-1}.

It is routine to deduce from the stable Paley--Wiener Theorem \ref{prop:stable-PW} the isomorphism
\begin{equation}\label{eqn:orbIEndo-gr}\begin{aligned}
		\orbI^{\Endo}(\tilde{G}) & \rightiso \bigoplus_{M \in \mathcal{L}(M_0)/W^G_0} \orbI^{\Endo}_{\cusp}(\tilde{M})^{W^G(M)} \\
		f^{\Endo} & \mapsto \left( f^{M^!} \big|_{\Phi_{2, \mathrm{bdd}}(M^!)} \right)_{\substack{M \in \mathcal{L}(M_0)/W^G_0 \\ \mathbf{M}^! \in \Endo_{\elli}(\tilde{M})}}
\end{aligned}\end{equation}
where $f^{\Endo} \mapsto f^{M^!}$ is as in \eqref{eqn:fEndo-descent} and $f^{M^!} \big|_{\Phi_{2, \mathrm{bdd}}(M^!)}$ is interpreted via Theorem \ref{prop:stable-PW}. See \cite[Lemma 6.2.1]{Li19} for details.

We are going to define a map in the other direction.

\begin{definition}[Adjoint transfer factor]\label{def:adjoint-transfer-factor}
	\index{transfer factor!adjoint}
	Let $\tilde{\gamma} \in \Gamma_{\mathrm{reg}}(\tilde{G})$ and $\delta \in \Gamma^{\Endo}_{\mathrm{reg}}(\tilde{G})$.
	\begin{enumerate}[(i)]
		\item First, assuming $\tilde{\gamma}$ and $\delta$ are both elliptic, we define
		\[ \Delta(\tilde{\gamma}, \delta) := \left| \mathfrak{D}(G_\delta, G; F) \right|^{-1} \overline{\Delta(\delta, \tilde{\gamma})} \]
		\index{Delta-gamma-delta@$\Delta(\tilde{\gamma}, \delta)$}
		with $\Delta(\delta, \tilde{\gamma}) := \Delta_{\mathbf{G}^!, \tilde{G}}(\delta, \tilde{\gamma})$, where $\mathbf{G}^! \in \Endo_{\elli}(\tilde{G})$ is the endoscopic datum coming with $\delta$. This definition extends to all groups of metaplectic type.
		\item In general, when $\delta \leftrightarrow \gamma$, we take $M \in \mathcal{L}(M_0)$ such that $\tilde{\gamma} \in \Gamma_{\text{$G$-reg, ell}}(\tilde{M}) \big/ W^G(M)$ and $\delta \in \Gamma^{\Endo}_{\text{$G$-reg, ell}}(\tilde{M}) \big/ W^G(M)$ and set
		\[ \Delta(\tilde{\gamma}, \delta) := \sum_{w \in W^G(M)} \Delta^{\tilde{M}}\left( \tilde{\gamma}, \Ad(w)\delta \right), \]
		where $\Delta^{\tilde{M}}$ denotes the adjoint transfer factor defined for $\tilde{M}$.
		\item If $\delta \not\leftrightarrow \gamma$ we set $\Delta(\tilde{\gamma}, \delta) = 0$.
	\end{enumerate}
\end{definition}

The \emph{adjoint transfer}
\begin{equation}\label{eqn:adjoint-transfer}\begin{tikzcd}[row sep=small]
	\left\{ \text{functions}\; \Gamma^{\Endo}_{\mathrm{reg}}(\tilde{G}) \to \CC \right\} \arrow[r, "{\Trans_{\Endo}}" inner sep=0.8em] & \left\{ \text{anti-genuine functions}\; \Gamma_{\mathrm{reg}}(\tilde{G}) \to \CC \right\}
\end{tikzcd}\end{equation}
is then characterized as in \eqref{eqn:collective-transfer-fcn}, with the new kernel function being the adjoint transfer factor.

\begin{proposition}[See {\cite[Proposition 5.2.4]{Li19}}]\label{prop:inverse-transfer-0}
	We have $\Trans_{\Endo} \left( \Trans^{\Endo} (f_{\tilde{G}}) \right) = f_{\tilde{G}}$ for all $f \in \orbI_{\asp}(\tilde{G}) \otimes \mes(G)$.
\end{proposition}

For the next statement, recall from \cite[\S 5.1, \S 6.2]{Li19}, which is in turn based on \cite{Ar96}, that $\orbI_{\asp, \cusp}(\tilde{G}) \otimes \mes(G)$ and $\orbI^{\Endo}_{\cusp}(\tilde{G})$ are actually pre-Hilbert spaces; both carry natural hermitian inner products $(\cdot | \cdot)$.

\begin{proposition}[See {\cite[Corollary 5.2.5]{Li19}}]\label{prop:transfer-cusp-isometry}
	The map $\Trans^{\Endo}_{\cusp}$ is an isometry from $\orbI_{\asp, \cusp}(\tilde{G}) \otimes \mes(G)$ into $\orbI^{\Endo}_{\cusp}(\tilde{G})$.
\end{proposition}

We are going to prove that $\Trans^{\Endo}_{\cusp}$ is actually surjective. The proofs below hinge upon some results and terminologies from \S\ref{sec:spectral-transfer} and \S\ref{sec:HC-descent}, but there will be no worry of no circularity.

\begin{theorem}\label{prop:image-transfer}
	The linear maps $\Trans^{\Endo}$ and $\Trans^{\Endo}_{\cusp}$ in Definition \ref{def:collective-transfer} are isomorphisms; their inverses are given by the $\mathcal{T}_{\Endo}$ in \eqref{eqn:adjoint-transfer}.
	
	When $F$ is Archimedean, they are isomorphisms between topological vector spaces and induce
	\[ \orbI_{\asp}(\tilde{G}, \tilde{K}) \otimes \mes(G) \rightiso \orbI^{\Endo}(\tilde{G}, \tilde{K}), \quad \orbI_{\asp, \cusp}(\tilde{G}, \tilde{K}) \otimes \mes(G) \rightiso \orbI^{\Endo}_{\cusp}(\tilde{G}, \tilde{K}). \]
\end{theorem}
\begin{proof}
	For non-Archimedean $F$, this is \cite[Theorem 5.3.1 and Corollary 6.3.3]{Li19}. It remains to establish the Archimedean case. Note that the injectivity of $\Trans^{\Endo}$ and the assertion about inverses follow from Proposition \ref{prop:inverse-transfer-0}.
	
	Assume that $F$ is Archimedean. The first step is to reduce the assertion about $\Trans^{\Endo}$ to that about $\Trans^{\tilde{M}, \Endo}_{\cusp}$ for all Levi subgroup $M \subset G$, using \eqref{eqn:gr-Trans}. It suffices to show the surjectivity of $\Trans^{\tilde{M}, \Endo}_{\cusp}$, since Proposition \ref{prop:transfer-cusp-isometry} will then imply that  $\Trans^{\tilde{M}, \Endo}_{\cusp}$ is an isometry between Hilbert spaces. By induction, we reduce further to the case $M = G$.
	
	Secondly, we may assume $F = \R$ since $\orbI_{\cusp, \asp}(\tilde{G}) = \{0\} = \orbI^{\Endo}_{\cusp}(\tilde{G})$ when $F = \CC$.
	
	Fix $\mathbf{G}^! \in \Endo_{\elli}(\tilde{G})$ and let $\varphi^! \in S\orbI_{\cusp}(G^!) \otimes \mes(G^!)$. Define $\varphi := \Trans_{\Endo}(\varphi): \Gamma_{\mathrm{reg}, -}(\tilde{G}) \to \CC$ (Definition \ref{def:adjoint-transfer-factor} and \eqref{eqn:adjoint-transfer}), namely
	\[ \varphi(\tilde{\gamma}) := |\mathfrak{D}(G_\gamma, G; F)|^{-1} \sum_{\delta \leftrightarrow \gamma} \Delta(\delta, \tilde{\gamma})^{-1} \varphi^!(\delta). \]
	Note that $\varphi(z\tilde{\gamma}) = z^{-1} \varphi(\tilde{\gamma})$, and $\varphi(\tilde{\gamma}) = 0$ for non-elliptic $\tilde{\gamma}$. In view of the inversion formula \cite[Lemma 5.2.2]{Li19} for transfer and its adjoint, we are reduced to show that $\varphi \in \orbI_{\asp}(\tilde{G}) \otimes \mes(G)$, since it will then follow that $\Trans^{\Endo}(\varphi) = \varphi^!$ (zero on the components $\neq \mathbf{G}^!$).

	The subsequent arguments follow \cite[I.4.13]{MW16-1} and we give only a brief sketch. Suppose that a diagram $(\epsilon, B^!, T^!, B, T, \eta)$ (Definition \ref{def:diagram}) is given, such that $T^!$ and $T$ are elliptic maximal tori in $G^!$ and $G$, respectively. The data of diagram yield an isomorphism $T \simeq T^!$ of $F$-tori. If $X \in \mathfrak{t}_{\mathrm{reg}}(F)$, we write $Y \in \mathfrak{t}^!_{\mathrm{reg}}(F)$ for its image, and we have $\exp(Y)\epsilon \leftrightarrow \exp(X)\eta$ under endoscopy.
	
	For $X \in \mathfrak{t}_{\mathrm{reg}}(F)$, define
	\begin{equation*}
		\Delta^{G_\eta}_T(X) := \prod_{\substack{\alpha \in \Sigma(G_\eta, T) \\ \alpha > 0 }} \sgn(i \alpha(X)) \; \in \{\pm 1\}.
	\end{equation*}
	The positive roots above are relative to $B$; they are all imaginary since $T$ is elliptic. Likewise, we define $\Delta^{G^!_\epsilon}_{T^!}(Y) \in \{\pm 1\}$
	
	Choose $\tilde{\gamma} \in \rev^{-1}(\gamma)$. As in \cite[p.104]{MW16-1}, using the characterization \cite{Bo94b} of $\orbI_{\asp, \cusp}(\tilde{G})$, it boils down to showing that
	\begin{equation}\label{eqn:image-transfer-aux}
		X \mapsto \frac{\Delta^{G_\eta}_T(X)}{\Delta^{G^!_\epsilon}_{T^!}(Y)} \cdot \Delta(\exp(Y)\epsilon, \exp(X)\tilde{\eta})^{-1}
	\end{equation}
	extends to a $C^\infty$-function in a neighborhood of $0$ in $\mathfrak{t}(F)$.
	
	The next step involves the semisimple descent for $\Delta$ in \cite[\S 7]{Li11}, to be reviewed in \S\S\ref{sec:descent-endoscopy}--\ref{sec:diagram-2}. The upshot is that one has a relevant endoscopic datum $\overline{\mathbf{G}^!_\epsilon}$ of $G_\eta$ together with the relations
	\[\begin{tikzcd}[column sep=large]
		G^!_\epsilon \arrow[dash, dashed, r, "\text{nonstandard}" inner sep=1em, "\text{endoscopy}"' inner sep=1em] & \overline{G^!_\epsilon} \arrow[dash, dashed, r, "\text{endoscopy}" inner sep=1em] & G_\eta .
	\end{tikzcd}\]
	We refer to \S\ref{sec:nonstandard-endoscopy} for the nonstandard endoscopic datum in question. As a consequence of this construction, we may factor $T^! \simeq T$ into $T^! \simeq \overline{T^!} \simeq T$ where $\overline{T^!}$ is an elliptic maximal torus of $\overline{G^!_\epsilon}$ related to $T^!$ by nonstandard endoscopy.
	
	Denote by $\overline{Y} \in \overline{\mathfrak{t}^!}(F)$ the image of $Y$. The descent of transfer factors (see \eqref{eqn:dy}) asserts
	\[ \Delta(\exp(Y)\epsilon, \exp(X)\tilde{\eta}) = \Delta_{\overline{\mathbf{G}^!_\epsilon}, G_\eta}(\overline{Y}, X) \]
	where $\Delta_{\overline{\mathbf{G}^!_\epsilon}, G_\eta}$ is a transfer factor on Lie algebras for the endoscopic datum $\overline{\mathbf{G}^!_\epsilon}$ of $G_\eta$.
	
	Observe that $\Delta^{G^!_\epsilon}_{T^!}(Y) = c \Delta^{\overline{G^!_\epsilon}}_{\overline{T^!}}(\overline{Y})$, where $c$ is a nonzero constant depending solely on the choice of Borel subgroups, since passing from $G^!_\epsilon \supset T^!$ to $\overline{G^!_\epsilon} \supset \overline{T^!}$ only rescales the roots. All in all, \eqref{eqn:image-transfer-aux} reduces to the corresponding equality for the endoscopic datum $\overline{\mathbf{G}^!_\epsilon}$ of $G_\eta$ established in \cite[p.104]{MW16-1}.
	
	Finally, to show the isomorphisms for $\tilde{K} \times \tilde{K}$-finite versions in the Archimedean case, it suffices to treat the case for $\orbI_{\asp, \cusp}(\tilde{G}) \otimes \mes(G)$. In view of Remark \ref{rem:real-PW}, the latter case follows from the fact \cite[Theorem 7.4.3]{Li19} that $\Delta(\phi, \tau)$ (Definition \ref{def:spectral-transfer-factor}) has finite support in each variable.
\end{proof}

Of course, Theorem \ref{prop:image-transfer} extends to all groups of metaplectic type.

\section{Spectral transfer}\label{sec:spectral-transfer}
The conventions about $\tilde{G}$, etc.\ are the same as \S\ref{sec:geom-transfer}.

Following \cite[\S 6.2]{Li19} and the formalism from \S\ref{sec:spectral-distributions}, we define the endoscopic counterparts of $\mathbb{S}^1 \backslash T_{\elli, -}(\tilde{G})$ and $\mathbb{S}^1 \backslash T_-(\tilde{G})$ as
\begin{align*}
	T^{\Endo}_{\elli}(\tilde{G}) & := \bigsqcup_{\mathbf{G}^! \in \Endo_{\elli}(\tilde{G})} \Phi_{2, \mathrm{bdd}}(G^!), \\
	T^{\Endo}(\tilde{G}) & := \bigsqcup_{M \in \mathcal{L}(M_0)/W^G_0} T^{\Endo}_{\elli}(\tilde{M}) / W^G(M).
\end{align*}
\index{TGEndo@$T^{\Endo}(\tilde{G})$, $T^{\Endo}_{\elli}(\tilde{G})$}

Let $f^{\Endo} \in \orbI^{\Endo}(\tilde{G})$ and $\phi \in T^{\Endo}(\tilde{G})$. Take $M \in \mathcal{L}(M_0)$ and $\mathbf{M}^! \in \Endo_{\elli}(\tilde{M})$ such that $\phi$ comes from $\phi_M \in \Phi_{2, \mathrm{bdd}}(M^!)$. We define
\[ f^{\Endo}(\phi) := f^{M^!}(\phi_M) \]
where $f^{\Endo} \mapsto f^{M^!}$ is as in \eqref{eqn:fEndo-descent}.

\begin{definition}[Spectral transfer factors {\cite[Definition 6.2.4]{Li19}}]\label{def:spectral-transfer-factor}
	\index{transfer factor!spectral}
	\index{Delta-spec@$\Delta(\phi, \tau)$}
	For every $(\phi, \tau) \in T^{\Endo}_{\elli}(\tilde{G}) \times T_{\elli, -}(\tilde{G})$, the spectral transfer factor $\Delta(\phi, \tau)$ is characterized by
	\begin{align*}
		\Delta(\phi, z\tau) & = z \Delta(\phi, \tau), \quad z \in \mathbb{S}^1, \\
		\Trans^{\Endo}(f)(\phi) & = \sum_{\tau \in T_{\elli, -}(\tilde{G})/\mathbb{S}^1} \Delta(\phi, \tau) f_{\tilde{G}}(\tau), \quad \phi \in T^{\Endo}_{\elli}(\tilde{G}),
	\end{align*}
	for all $f \in \orbI_{\asp, \cusp}(\tilde{G}) \otimes \mes(G)$.
	
	For general $(\phi, \tau) \in T^{\Endo}(\tilde{G}) \times T_-(\tilde{G})$, suppose that they come from $\phi_M \in T^{\Endo}_{\elli}(\tilde{M})$ and $\tau_M \in T_{\elli, -}(\tilde{M})$ for the same $M \in \mathcal{L}(M_0)/W^G_0$. We put
	\[ \Delta(\phi, \tau) := \sum_{\tau_M^\dagger \in W^G(M) \cdot \tau_M } \Delta^{\tilde{M}}(\phi_M, \tau_M^\dagger) \]
	where $\Delta^{\tilde{M}}$ denotes the spectral transfer factor defined for $\tilde{M}$ in the obvious way. If there is no such $M$, we put $\Delta(\phi, \tau) := 0$.
\end{definition}

\begin{theorem}[Local character relation in the tempered case]\label{prop:local-character-relation}
	For all $\phi \in T^{\Endo}(\tilde{G})$ and $f \in \orbI_{\asp}(\tilde{G}) \otimes \mes(G)$, we have
	\[ \Trans^{\Endo}(f)(\phi) = \sum_{\tau \in T_-(\tilde{G})/\mathbb{S}^1} \Delta(\phi, \tau) f_{\tilde{G}}(\tau) . \]
\end{theorem}
\begin{proof}
	The non-Archimedean case is the main result of \cite{Li19}. Consider the case $F = \R$ next.
	
	Write $f_{\tilde{G}}(\pi) := \Theta_\pi(f_{\tilde{G}})$ for each $\pi \in \Pi_{\mathrm{temp}, -}(\tilde{G})$. The local character relation of \cite[Theorem 7.4.3]{Li19} yields a function $\Delta_{\mathrm{spec}}: T^{\Endo}(\tilde{G}) \times \Pi_{\mathrm{temp}, -}(\tilde{G}) \to \{\pm 1\}$ satisfying
	\begin{equation}\label{eqn:local-character-relation-aux-0}
		\Trans^{\Endo}(f)(\phi) = \sum_{\pi \in \Pi_{\mathrm{temp}, -}(\tilde{G})} \Delta_{\mathrm{spec}}(\phi, \pi) f_{\tilde{G}}(\pi)
	\end{equation}
	for all $\phi$ and $f \in \orbI_{\asp}(\tilde{G}) \otimes \mes(G)$, and $\Delta_{\mathrm{spec}}(\cdot, \pi)$ (resp.\ $\Delta_{\mathrm{spec}}(\phi, \cdot)$) has finite support for each $\pi$ (resp.\ for each $\phi$). These properties characterize $\Delta_{\mathrm{spec}}$.

	Choose a representative in $T_-(\tilde{G})$ for every class in $T_-(\tilde{G})/\mathbb{S}^1$. By the theory of $R$-groups, $T_-(\tilde{G})/\mathbb{S}^1$ gives a basis of $D_{\mathrm{temp}, -}(\tilde{G}) \otimes \mes(G)^\vee$: specifically, we may write
	\[ \Theta_\tau = \sum_\pi \mathrm{mult}(\tau : \pi) \Theta_\pi, \quad \Theta_\pi = \sum_\tau \mathrm{mult}(\pi : \tau) \Theta_\tau \]
	for all $\tau \in T_-(\tilde{G})/\mathbb{S}^1$ with its representative and $\pi \in \Pi_{\mathrm{temp}, -}(\tilde{G})$, for uniquely determined coefficients $\mathrm{mult}(\cdots)$. Switching between bases, \eqref{eqn:local-character-relation-aux-0} uniquely determines
	\[ \Delta^\circ: T^{\Endo}(\tilde{G}) \times T_-(\tilde{G}) \to \mathbb{S}^1 \]
	such that
	\begin{align*}
		\Delta^\circ(\phi, z\tau) & = z\Delta^\circ(\phi, \tau), \quad z \in \mathbb{S}^1 , \\
		\Trans^{\Endo}(f)(\phi) & = \sum_{\tau \in T_-(\tilde{G})/\mathbb{S}^1} \Delta^\circ(\phi, \tau) f_{\tilde{G}}(\tau)
	\end{align*}
	for all $f$. Specifically, $\Delta^\circ(\phi, \tau) = \sum_{\pi \in \Pi_{\mathrm{temp}, -}(\tilde{G})} \Delta_{\mathrm{spec}}(\phi, \pi) \mathrm{mult}(\pi : \tau)$ for all $\tau \in T_-(\tilde{G})/\mathbb{S}^1$.
	
	Our goal is thus to show
	\begin{equation}\label{eqn:local-character-relation-aux-1}
		\Delta^\circ(\phi, \tau) = \Delta(\phi, \tau), \quad (\phi, \tau) \in T^{\Endo}(\tilde{G}) \times T_-(\tilde{G}).
	\end{equation}

	The first step is to reduce to the elliptic setting. We say $\pi \in \Pi_{\mathrm{temp}, -}(\tilde{G})$ is \emph{elliptic} if $\Theta_\pi$ is not identically zero on $\Gamma_{\mathrm{reg, ell}}(\tilde{G})$. In \cite[Definition 7.4.1]{Li19} one defined a subset $\Pi_{2\uparrow, -}(\tilde{G})$ of $\Pi_{\mathrm{temp}, -}(\tilde{G})$. All $\pi \in \Pi_{2\uparrow, -}(\tilde{G})$ are elliptic. Indeed, by \cite[Remark 7.5.1]{Li19} $\pi$ is a non-degenerate limit of discrete series in the sense of Knapp--Zuckerman, and such representations are known to be elliptic; see \textit{loc.\ cit.} for the relevant references.
	
	By \cite[Proposition 5.4.4]{Li12b}, $T_{\elli, -}(\tilde{G})/\mathbb{S}^1$ gives a basis for the space spanned by the characters of all elliptic $\pi$.
	
	Let $\phi \in T^{\Endo}(\tilde{G})$. Take $M \in \mathcal{L}(M_0)$ and $\mathbf{M}^! \in \Endo_{\elli}(\tilde{M})$ such that $\phi$ comes from $\phi_{M^!} \in \Phi_{\mathrm{bdd}, 2}(M^!)$ up to $W^G(M)$. Denote the factors relative to $\tilde{M}$ as $\Delta^{\tilde{M}}$, etc. By \cite[Theorem 7.4.3]{Li19},
	\[ \Delta_{\mathrm{spec}}(\phi, \pi) = \sum_{\pi_M \in \Pi_{2\uparrow, -}(\tilde{M})} \Delta^{\tilde{M}}_{\mathrm{spec}}(\phi_{M^!}, \pi_M) \mathrm{mult}(I_{\tilde{P}}(\pi_M) : \pi) \]
	for all $\pi \in \Pi_{\mathrm{temp}, -}(\tilde{G})$, where $P \in \mathcal{P}(M)$ and $\mathrm{mult}(I_{\tilde{P}}(\pi_M) : \pi)$ denotes the multiplicity of $\pi$ in $I_{\tilde{P}}(\pi_M)$.

	We claim that
	\begin{equation}\label{eqn:local-character-relation-aux-2}
		\Delta^\circ(\phi, \tau) = \sum_{\substack{\tau_M \in T_{\elli, -}(\tilde{M}) \\ \tau_M \mapsto \tau }} \Delta^{\tilde{M}, \circ}(\phi_{M^!}, \tau_M).
	\end{equation}
	Note that the sum is actually over an orbit $W^G(M) \tau_M$ if such a $\tau_M$ exists. We may assume $\tau$ is the representative of some element from $T_-(\tilde{G})/\mathbb{S}^1$; we also choose representatives in $T_-(\tilde{M})$ for $T_-(\tilde{M})/\mathbb{S}^1$ compatibly with induction. We have
	\begin{multline*}
		\Delta^\circ(\phi, \tau) = \sum_{\pi \in \Pi_{\mathrm{temp}, -}(\tilde{G})} \Delta_{\mathrm{spec}}(\phi, \pi) \mathrm{mult}(\pi : \tau) \\
		= \sum_\pi \sum_{\pi_M \in \Pi_{2\uparrow, -}(\tilde{M})} \sum_{\tau_M \in T_{\elli, -}(\tilde{M})/\mathbb{S}^1} \\
		\cdot \mathrm{mult}(\tau_M : \pi_{\tilde{M}}) \mathrm{mult}(I_{\tilde{P}}(\pi_{\tilde{M}}) : \pi ) \mathrm{mult}(\pi : \tau) \Delta^{\tilde{M}, \circ}(\phi_{M^!}, \tau_M).
	\end{multline*}
	Given $\tau_M$, the sum over $(\pi, \pi_M)$ of the triple products of $\mathrm{mult}(\cdots)$ is readily seen to be $1$ if $\tau_M \mapsto \tau$, otherwise it is zero. This proves \eqref{eqn:local-character-relation-aux-2}.
	
	Observe that \eqref{eqn:local-character-relation-aux-2} take the same form as the induction formula in Definition \ref{def:spectral-transfer-factor}. In order to prove \eqref{eqn:local-character-relation-aux-1}, we may assume $\phi \in T^{\Endo}_{\elli}(\tilde{G})$, in which case $\Delta^\circ(\phi, \tau) = 0$ unless $\tau \in T_{\elli, -}(\tilde{G})$, and ditto for $\Delta(\phi, \tau)$. It is thus legitimate to take $f \in \orbI_{\asp, \cusp}(\tilde{G}) \otimes \mes(G)$ in the characterization of $\Delta^\circ(\phi, \cdot)$. All in all, $\Delta^\circ(\phi, \cdot)$ and $\Delta(\phi, \cdot)$ have the same characterization, whence \eqref{eqn:local-character-relation-aux-1}.
	
	The case $F = \CC$ is even simpler because it reduces to the case of split maximal tori via parabolic induction: see \cite[\S 7.6]{Li19}.
\end{proof}

We remark that by the discussions in \cite[\S 6.2]{Li19}, Theorem \ref{prop:local-character-relation} is equivalent to saying that $\Trans^{\Endo}$ equals the composition of
\begin{multline}\label{eqn:gr-Trans}
	\orbI_{\asp}(\tilde{G}) \otimes \mes(G) \xrightarrow[\sim]{\eqref{eqn:orbI-gr}} \bigoplus_{M \in \mathcal{L}(M_0)/W^G_0} \orbI_{\asp, \cusp}(\tilde{M})^{W^G(M)} \otimes \mes(M) \\
	\xrightarrow{\bigoplus_M \Trans^{\tilde{M}, \Endo}_{\cusp}} \bigoplus_{M \in \mathcal{L}(M_0)/W^G_0} \orbI^{\Endo}_{\cusp}(\tilde{M})^{W^G(M)}
	\xrightarrow[\sim]{\eqref{eqn:orbIEndo-gr}} \orbI^{\Endo}(\tilde{G})
\end{multline}
where $\Trans^{\tilde{M}, \Endo}_{\cusp}$ stands for the collective geometric transfer for $\tilde{M}$ restricted to the cuspidal subspace.

For the next result, recall from \cite[\S 7.3]{Li19} that when $F$ is Archimedean, there is a canonical map $\lambda^! \mapsto \lambda$ from the space of infinitesimal characters for $G^!$ to that for $\tilde{G}$, where $\mathbf{G}^! \in \Endo_{\elli}(\tilde{G})$. The map is surjective with finite fibers. Hence we have a correspondence between infinitesimal characters of members of $\Phi_{\mathrm{bdd}}(G^!)$ and those of $T_-(\tilde{G})/\mathbb{S}^1$.

\begin{proposition}
	When $F$ is Archimedean, we have $\Delta(\phi, \tau) \neq 0$ only when $\phi$ and $\tau$ have matching infinitesimal characters.
\end{proposition}
\begin{proof}
	Upon reviewing the proof of Theorem \ref{prop:local-character-relation}, this reduces immediately to the case of the factors $\Delta(\phi, \pi)$. In turn, the latter case reduces to $(\phi, \pi) \in \Phi_{\mathrm{bdd}, 2}(G^!) \times \Pi_{2\uparrow, -}(\tilde{G})$ which is established in \cite[\S 7.3]{Li19}.
\end{proof}

\section{Geometric distributions and their transfer}\label{sec:geom-dist-transfer}
Consider a local field $F$ of characteristic zero and a group of metaplectic type $\rev: \tilde{M} \to M(F)$. The notion of genuine and anti-genuine distributions has been defined in \S\ref{sec:metaplectic-type}. In this section, we will focus on genuine distributions of geometric origin, defined as follows.

For any $\tilde{\gamma} \in \tilde{M}$, we denote by $\gamma = \gamma_{\text{ss}} \gamma_{\text{u}}$ the Jordan decomposition of $\gamma := \rev(\tilde{\gamma})$ in $M(F)$, where $\gamma_{\text{ss}}$ and $\gamma_{\text{u}}$ are the semisimple and unipotent parts, respectively.

It makes sense to define the support of $D \in D_-(\tilde{M})$ as a $M(F)$-invariant closed subset $\Supp(D) \subset \tilde{M}$. Similarly, for a quasisplit connected reductive $F$-group $M^!$ and $D \in SD(M^!)$, one can define $\Supp(D) \subset M^!(F)$.

\begin{definition}\label{def:geom-O-dist}
	\index{Dgeom@$D_{\mathrm{geom}, -}$, $D_{\mathrm{unip}, -}$}
	Let $\mathcal{O}$ be a finite union of semisimple conjugacy classes in $M(F)$. Define
	\begin{align*}
		D_{\mathrm{geom},-}(\tilde{M}, \mathcal{O}) & := \left\{ D \in D_-(\tilde{M}) : \tilde{\gamma} \in \Supp(D) \implies \gamma_{\text{ss}} \in \mathcal{O} \right\}, \\
		D_{\mathrm{geom}, -}(\tilde{M}) & := \bigoplus_{\substack{\mathcal{O} \subset M(F) \\ \text{ss.\ conj.\ class} }} D_{\mathrm{geom},-}(\tilde{M}, \mathcal{O}) \; \subset D_-(\tilde{M}), \\
		D_{\mathrm{unip}, -}(\tilde{M}) & := D_{\mathrm{geom}, -}(\tilde{M}, \{1\}).
	\end{align*}
\end{definition}

The orbital integrals $I^{\tilde{M}}(\tilde{\gamma}, \cdot)$ along conjugacy classes whose closure contains $\rev^{-1}(\mathcal{O})$ yield elements of $D_{\mathrm{geom}, -}(\tilde{M}, \mathcal{O}) \otimes \mes(M)^\vee$.

\begin{definition}
	Denote by $D_{\mathrm{orb}, -}(\tilde{M})$ (resp.\ $D_{\mathrm{orb}, -}(\tilde{M}, \mathcal{O})$) the subspace of $D_{\mathrm{geom}, -}(\tilde{M})$ (resp.\ $D_{\mathrm{geom}, -}(\tilde{M}, \mathcal{O})$) spanned by the orbital integrals above.
\end{definition}

Note that when $F$ is non-Archimedean, $D_{\mathrm{geom},-} = D_{\mathrm{orb}, -}$, but the inclusion is proper for Archimedean $F$.

\begin{definition}
	Let $\mathcal{O}$ be a finite union of semisimple conjugacy classes in $M(F)$. Denote by $\orbI_{\asp}(\tilde{M})_{\mathcal{O}, 0}$ the subspace of elements $f_{\tilde{M}} \in \orbI_{\asp}(\tilde{M})$ such that there exists an invariant neighborhood $\mathcal{V}$ of $\mathcal{O}$ in $M(F)$ satisfying $f_{\tilde{M}}(\tilde{\gamma}) = 0$ for all $\gamma \in \mathcal{V} \cap M_{\mathrm{reg}}(F)$. Define
	\[ \orbI_{\asp}(\tilde{M})_{\mathcal{O}, \mathrm{loc}} := \begin{cases}
		\orbI_{\asp}(\tilde{M}) \big/ \orbI_{\asp}(\tilde{M})_{\mathcal{O}, 0}, & F: \text{non-Archimedean}, \\
		\orbI_{\asp}(\tilde{M}) \big/ \overline{\orbI_{\asp}(\tilde{M})_{\mathcal{O}, 0}}, & F: \text{Archimedean}.
	\end{cases}\]
\end{definition}

\begin{lemma}\label{prop:Dgeom-as-dual}
	For every $\mathcal{O}$ as above, $D_-(\tilde{M}, \mathcal{O})$ is identified with the dual of $\orbI_{\asp}(\tilde{M})_{\mathcal{O}, \mathrm{loc}}$ (continuous if $F$ is Archimedean).
\end{lemma}
\begin{proof}
	Same as \cite[I.5.1, I.5.2]{MW16-1}, based on standard techniques of harmonic analysis.
\end{proof}

Suppose that $\tilde{M}$ is a Levi subgroup of $\tilde{G}$. If $\mathcal{O}$ is a finite union of semisimple conjugacy classes in $M(F)$, we denote by $\mathcal{O}^G$ the finite union of the semisimple conjugacy classes in $G(F)$ intersecting $\mathcal{O}$. Then
\begin{equation}\label{eqn:parabolic-ind-dist-geom}\begin{tikzcd}[row sep=tiny]
	D_-(\tilde{M}) \otimes \mes(M)^\vee \arrow[r, "{\Ind^{\tilde{G}}_{\tilde{M}}}"] & D_-(\tilde{G}) \otimes \mes(G)^\vee \\
	D_{\mathrm{geom}, -}(\tilde{M}) \otimes \mes(M)^\vee \arrow[r] \arrow[phantom, u, "\subset" description, sloped] & D_{\mathrm{geom}, -}(\tilde{G}) \otimes \mes(G)^\vee \arrow[phantom, u, "\subset" description, sloped] \\
	D_{\mathrm{geom}, -}(\tilde{M}, \mathcal{O}) \otimes \mes(M)^\vee \arrow[r] \arrow[phantom, u, "\subset" description, sloped] & D_{\mathrm{geom}, -}(\tilde{G}, \mathcal{O}^G) \otimes \mes(G)^\vee \arrow[phantom, u, "\subset" description, sloped] .
\end{tikzcd}\end{equation}

We recall the definition of the stable avatars below. See also \cite[I.5.4 --- I.5.5]{MW16-1}.

\begin{definition}
	For a quasisplit connected reductive $F$-group $M^!$, we define the subspaces $SD_{\mathrm{geom}}(M^!, \mathcal{O}) \subset SD_{\mathrm{geom}}(M^!)$ of $SD(M^!)$ accordingly; here $\mathcal{O}$ is a finite union of stable semisimple conjugacy classes in $M^!(F)$. We set $SD_{\mathrm{unip}}(M^!) := SD_{\mathrm{geom}}(M^!, \{1\})$.
	
	We also define $SD_{\mathrm{orb}}(M^!) := SD_{\mathrm{geom}}(M^!) \cap D_{\mathrm{orb}}(M^!)$.

	As before, we define $S\orbI(M^!)_{\mathcal{O}, 0} \subset S\orbI(M^!)$ to be the subspace of elements vanishing in an invariant neighborhood of $\mathcal{O}$ in $M^!(F)$, and then define
	\[ S\orbI(M^!)_{\mathcal{O}, \mathrm{loc}} := \begin{cases}
		S\orbI(M^!) \big/ S\orbI(M^!)_{\mathcal{O}, 0}, & F: \text{non-Archimedean}, \\
		S\orbI(M^!) \big/ \overline{S\orbI(M^!)_{\mathcal{O}, 0}}, & F: \text{Archimedean}.
	\end{cases}\]
\end{definition}

The analogue of Lemma \ref{prop:Dgeom-as-dual} holds here, namely $SD_{\mathrm{geom}}(M^!, \mathcal{O})$ is the dual of $S\orbI(M^!)_{\mathcal{O}, \mathrm{loc}}$. See \cite[I.5.4--I.5.5]{MW16-1}. When $M^!$ is a Levi subgroup of $G^!$, we also have the map $\Ind^{G^!}_{M^!}: SD_{\mathrm{geom}}(M^!) \otimes \mes(M^!)^\vee \to SD_{\mathrm{geom}}(G^!) \otimes \mes(G^!)^\vee$, also denoted by $\delta \mapsto \delta^{G^!}$.

\begin{theorem}[Cf.\ {\cite[I.5.7 Proposition]{MW16-1}}]\label{prop:Dgeom-preservation}
	For any stable semisimple conjugacy class $\mathcal{O}$ in $G(F)$ and $\mathbf{G}^! \in \Endo(\tilde{G})$, write
	\[ \mathcal{O}_{G^!} := \left\{ \delta \in G^!(F)_{\mathrm{ss}} : \exists \gamma \in \mathcal{O}, \; \delta \leftrightarrow \gamma \right\}. \]
	Then $\mathcal{O}_{G^!}$ is a finite union of stable conjugacy classes in $G^!(F)$. The map $\trans^{\Endo}$ in Definition \ref{def:collective-transfer} restricts to a surjection
	\[ \bigoplus_{\mathbf{G}^! \in \Endo_{\elli}(\tilde{G})} SD_{\mathrm{geom}}(G^!, \mathcal{O}_{G^!}) \otimes \mes(G^!)^\vee \to D_{\mathrm{geom}, -}(\tilde{G}, \mathcal{O}) \otimes \mes(G)^\vee, \]
	whose kernel is spanned by elements of the form
	\[ \delta[s]^{G^![s]} - \delta[t]^{G^![t]}, \quad \delta \in SD_{\mathrm{geom}}(M^!, \mathcal{O}_{M^!}) \otimes \mes(M^!)^\vee, \]
	where $M \subset G$ is a Levi subgroup, $\mathbf{M}^! \in \Endo_{\elli}(\tilde{M})$ and $s, t \in \Endo_{\mathbf{M}^!}(\tilde{G})$.

	The same results holds if $\tilde{G}$ is replaced by any group of metaplectic type.
\end{theorem}
\begin{proof}
	Since the morphism $\Sigma_{\mathrm{ss}}(G^!) \to \Sigma_{\mathrm{ss}}(G)$ in Definition \ref{def:corr-orbits} and the subsequent discussions has finite fibers, $\mathcal{O}_{G^!}$ is a finite union of stable classes. The remaining arguments are the same as in \cite[I.5.8, I.5.9]{MW16-1}. We shall omit the details and give only the basic ingredients below.
	\begin{itemize}
		\item Theorem \ref{prop:image-transfer}, which describes the image of $\Trans^{\Endo}$.
		\item The adjoint transfer discussed in \S\ref{sec:geom-transfer}, which expresses $f_{\tilde{G}}(\tilde{\gamma})$ in terms of the values of $\Trans^{\Endo}(f)$.
	\end{itemize}
	Our case is actually simpler than \textit{loc.\ cit.}, since the relevance of endoscopic data and the endoscopic outer automorphisms on $G^!$ do not arise here (see Remark \ref{rem:centerless}). In contrast, the twists $\delta[s]$, $\delta[t]$ of $\delta$ are new: the same phenomenon of twists appeared in Proposition \ref{prop:central-twist-corr}.

	Finally, the extension to groups of metaplectic type is straightforward.
\end{proof}

\section{Spherical fundamental lemma}\label{sec:LF}
Consider the local metaplectic covering $\tilde{G} = \Mp(W)$ with $(W, \lrangle{\cdot|\cdot})$ and $\psi$ chosen. In this section, we assume that $F$ is non-Archimedean with residual characteristic $p > 2$, and fix a self-dual $\mathfrak{o}_F$-lattice $L \subset W$ with respect to $\psi \circ \lrangle{\cdot|\cdot}$.

Regard the resulting hyperspecial subgroup $K := G(\mathfrak{o}_F)$ as an open compact subgroup of $\tilde{G}$ as in \S\ref{sec:metaplectic-type}. Fix the Haar measure on $G(F)$ with $\mes(K) = 1$, and take the corresponding measure on $\tilde{G}$ by \eqref{eqn:measure-cover}. Similarly, for each $\mathbf{G}^! \in \Endo(\tilde{G})$, we use the Haar measure with $\mes(K^!) = 1$ for some (equivalently, any) hyperspecial subgroup $K^! \subset G^!(F)$. The lines $\mes(G)$ and $\mes(G^!)$ are thus trivialized.

We call this the \emph{unramified situation}.
\index{unramified situation}
\index{fK@$f_K$}
\index{HKGK@$\mathcal{H}_{\asp}(K \backslash \tilde{G} / K)$}

Consider the anti-genuine spherical Hecke algebra $\mathcal{H}_{\asp}(K \backslash \tilde{G} / K)$ with unit $f_K$. For each $\mathbf{G}^! \in \Endo(\tilde{G})$, we have the spherical Hecke algebra $\mathcal{H}(K^! \backslash G^!(F) / K^!)$ with unit $\mathbf{1}_{K^!}$. The metaplectic Satake isomorphism (see eg.\ \cite[Theorem 3.1.3]{Luo18}) reads
\[ \mathcal{H}_{\asp}(K \backslash \tilde{G} / K) \simeq \CC\left[ \tilde{G}^\vee /\!/ \tilde{G}^\vee \right] \]
between $\CC$-algebras, where $\tilde{G}^\vee /\!/ \tilde{G}^\vee$ is the categorical quotient of $\tilde{G}^\vee$ under adjoint action, and $\CC[\cdots]$ means the algebra of regular functions. Similarly, the classical Satake isomorphism reads
\[ \mathcal{H}(K^! \backslash G^!(F) / K^!) \simeq \CC\left[ \check{G}^! /\!/ \check{G}^! \right]. \]
The inclusion $\check{G}^! \hookrightarrow \tilde{G}^\vee$ being canonical up to $\tilde{G}^\vee$-conjugacy, it induces a canonical homomorphism
\[ b_{\mathbf{G}^!, \tilde{G}}: \mathcal{H}_{\asp}(K \backslash \tilde{G} / K) \to \mathcal{H}(K^! \backslash G^!(F) / K^!) \]
between $\CC$-algebras, In particular, $b_{\mathbf{G}^!, \tilde{G}}(f_K) = \mathbf{1}_{K^!}$.
\index{bGG@$b_{\mathbf{G}^{"!}, \tilde{G}}$}

\begin{theorem}[Spherical fundamental lemma]
	In the unramified situation, the diagram
	\[\begin{tikzcd}
		\orbI_{\asp}(\tilde{G}) \arrow[r, "{\Trans_{\mathbf{G}^!, \tilde{G}}}"] & S\orbI(G^!) \\
		\mathcal{H}_{\asp}(K \backslash \tilde{G} / K) \arrow[u] \arrow[r, "{b_{\mathbf{G}^!, \tilde{G}}}"'] & \mathcal{H}(K^! \backslash G^!(F) / K^!) \arrow[u]
	\end{tikzcd}\]
	commutes for each $\mathbf{G}^! \in \Endo(\tilde{G})$.
\end{theorem}
\begin{proof}
	The case when $\mathbf{G}^! \in \Endo_{\elli}(\tilde{G})$ is the Main Theorem of \cite{Luo18}, and the general case follows by parabolic descent as usual. Note that in \textit{loc.\ cit.}, one assumes that $\psi$ has conductor $\mathfrak{o}_F$ and $L$ is a self-dual lattice with respect to $\lrangle{\cdot|\cdot}$. We can also reduce to that case by replacing $\psi$ (resp.\ $\lrangle{\cdot|\cdot}$) by $\psi_a: x \mapsto ax$ (resp.\ $a^{-1}\lrangle{\cdot|\cdot}$) for some $a \in F^\times$, so that $\psi_a$ has conductor $\mathfrak{o}_F$. Then $L$ is self-dual with respect to $a^{-1}\lrangle{\cdot|\cdot}$. This is legitimate since the following depend only on $\psi \circ \lrangle{\cdot|\cdot}$:
	\begin{compactitem}
		\item the construction of $\tilde{G}$,
		\item the transfer factors, and
		\item the splitting of $\tilde{G}$ over $K$.
	\end{compactitem}
	Details can be found in \cite{Li11}; for instance, the last item can be checked using the character formula for $\omega_\psi$ in \cite[Lemme 4.20]{Li11}.
\end{proof}

All the constructions and results above extend to groups of metaplectic type, in the routine manner.

\chapter{Harish-Chandra descent}\label{sec:HC-descent}
The Harish-Chandra descent, also known as semisimple descent, is one of the most important tools in this work. We will review the descent of test functions, endoscopic data and transfer factors around a semisimple element of the metaplectic group, and set up the relevant notations. As a natural consequence, the nonstandard endoscopy on Lie algebras enters into the picture. These results are extracted from \cite{Li11, Li12a}.

These materials are supplemented with a discussion on diagrams in \S\ref{sec:diagram-2}. The setting with a Levi subgroup, addressed in \S\ref{sec:descent-endoscopy-Levi}, will play a vital role when we study various ``weighted'' objects.

\section{Harmonic analysis}\label{sec:descent-HA}
Let $F$ be a local field of characteristic zero. The following theory works for general coverings $\rev: \tilde{G} \to G(F)$ with $\Ker(\rev) = \bmu_m$; see \cite[\S 4]{Li12b}.

Let $\eta \in G(F)_{\mathrm{ss}}$ and $\tilde{\eta} \in \rev^{-1}(\eta)$. Consider an open neighborhood $\mathcal{U}^\flat$ of $0 \in \mathfrak{g}_\eta(F)$ subject to the following conditions:
\begin{itemize}
	\item $X \in \mathfrak{g}_\eta(F)$ lies in $\mathcal{U}^\flat$ if and only if its semisimple part does;
	\item $\mathcal{U}^\flat$ is invariant under conjugation by $Z_G(\eta)(F)$;
	\item the exponential map is defined on $\mathcal{U}^\flat$ and induces a homeomorphism onto its image.
\end{itemize}
Specifically, there are two exponential maps fitting into a commutative diagram
\begin{equation}\label{eqn:two-exp}\begin{tikzcd}
	\mathcal{U}^\flat \arrow[r, "\exp"] \arrow[rd, "\exp"'] & \widetilde{G_\eta} \arrow[d, "\rev"] \\
	& G_\eta(F)
\end{tikzcd}\end{equation}
and both will be denoted by $\exp$. Let $\tilde{\mathcal{U}} \subset \tilde{G}$ be the open subset consisting of elements which are conjugate to an element in $\rev^{-1}(\eta) \cdot \exp(\mathcal{U}^\flat)$.

Define $\omega : Z_G(\eta)(F) \to \bmu_m$ by $\omega(g) = \tilde{\eta} \tilde{g} \tilde{\eta}^{-1} \tilde{g}^{-1}$ by choosing liftings. We have $\orbI_{\asp}(\tilde{\mathcal{U}}) \subset \orbI_{\asp}(\tilde{G})$ and $\orbI(\mathcal{U}^\flat)_\omega \subset \orbI(\mathfrak{g}_\eta)_\omega$, where $(\cdots)_\omega$ is the space of $\omega$-coinvariants under $Z_G(\eta)(F)$-action; see \cite[\S 4.1]{Li12b}.

Each $f^\flat \in \orbI(\mathfrak{g}_\eta)_\omega$ can be evaluated on $X \in \mathfrak{g}_{\eta, \mathrm{reg}}(F)$ once the Haar measures are chosen, namely by taking orbital integrals twisted by the character $\omega|_{G_\eta(F)}$. The value will be denoted by $f^\flat_{G_\eta, \omega}(X)$.

Following the terminology of \cite[I.4.1]{MW16-1}, the \emph{Harish-Chandra descent} of test functions is the canonical linear isomorphism
\[ \desc^{\tilde{G}}_{\tilde{\eta}}: \orbI_{\asp}(\tilde{\mathcal{U}}) \otimes \mes(G) \rightiso \orbI(\mathcal{U}^\flat)_\omega \otimes \mes(G_\eta) \]
characterized by
\[ f \mapsto f^\flat \iff \forall X \in \mathcal{U}^\flat \cap \mathfrak{g}_{\eta, \mathrm{reg}}(F), \; f_{\tilde{G}}(\tilde{\eta}\exp X) = f^\flat_{G_\eta, \omega}(X). \]
\index{desc@$\desc^{\tilde{G}}_{\tilde{\eta}}$}

Moreover, if $\eta$ is elliptic in $G(F)$, it restricts to an isomorphism
\[ \orbI_{\asp, \cusp}(\tilde{\mathcal{U}}) \otimes \mes(G) \rightiso \orbI_{\cusp}(\mathcal{U}^\flat)_\omega \otimes \mes(G_\eta) \]
where $\orbI_{\cusp}(\cdots)_\omega$ is defined via parabolic descent, as usual.

In the upcoming applications, $\rev: \tilde{G} \to G(F)$ will be of metaplectic type, therefore $\omega$ is trivial and $Z_G(\eta) = G_\eta$.

Consider now the stable version. Let $G^!$ be a quasisplit $F$-group. Let $\epsilon \in G^!(F)_{\mathrm{ss}}$. Up to stable conjugacy, we may assume $G^!_\epsilon$ is quasisplit by \cite[Lemma 3.3]{Ko82}. Define the finite $F$-group scheme
\[ \Xi_\epsilon = \Xi^{G^!}_\epsilon := Z_{G^!}(\epsilon)/G^!_\epsilon. \]
\index{Xi-epsilon@$\Xi_\epsilon$}

Consider an open neighborhood $\mathcal{U}^\flat$ of $0 \in \mathfrak{g}^!_\epsilon(F)$ subject to the conditions below:
\begin{itemize}
	\item $X \in \mathfrak{g}^!_\epsilon(F)$ lies in $\mathcal{U}^\flat$ if and only if its semisimple part does;
	\item $\mathcal{U}^\flat$ is invariant under stable conjugation and $Z_{G^!}(\epsilon)(F)$-action;
	\item the exponential map is defined on $\mathcal{U}^\flat$ and induces a homeomorphism onto its image.
\end{itemize}

Define $\mathcal{U} \subset G^!(F)$ to be the open subset consisting of elements which are conjugate to an element in $\epsilon \cdot \exp(\mathcal{U}^\flat)$. The stable Harish-Chandra descent is a canonical linear isomorphism
\[ S\desc^{G^!}_{\epsilon}: S\orbI(\mathcal{U}) \otimes \mes(G^!) \rightiso S\orbI(\mathcal{U}^\flat)_{\Xi_\epsilon(F)} \otimes \mes(G^!_\epsilon) \]
characterized in terms of stable orbital integrals on Lie algebras, namely
\[ f \mapsto f^\flat \iff \forall Y \in \mathcal{U}^\flat \cap \mathfrak{g}^!_{\epsilon, \mathrm{reg}}(F), \; f^{G^!}(\epsilon \exp(Y)) = f^{\flat, G^!_\epsilon}(Y). \]
Here $(\cdots)_{\Xi_\epsilon(F)}$ stands for the space of coinvariants under $\Xi_\epsilon(F)$-action.
\index{Sdesc@$S\desc^{G^{"!}}_{\epsilon}$}

When $\epsilon$ is elliptic in $G^!(F)$, the map restricts to $S\orbI_{\cusp}(\mathcal{U}) \otimes \mes(G^!) \rightiso S\orbI_{\cusp}(\mathcal{U}^\flat)_{\Xi_\epsilon(F)} \otimes \mes(G^!_\epsilon)$.

These results are established in \cite[I.4.8]{MW16-1}. By taking transposes, we obtain
\begin{gather*}
	\desc^{\tilde{G}, *}_{\tilde{\eta}}: D(\mathcal{U}^\flat)^\omega \otimes \mes(G_\eta)^\vee \to D_-(\tilde{\mathcal{U}}) \otimes \mes(G)^\vee, \\
	S\desc^{G^!, *}_\epsilon: SD(\mathcal{U}^\flat)^{\Xi_\epsilon(F)} \otimes \mes(G^!_\epsilon)^\vee \to D(\mathcal{U}) \otimes \mes(G^!)^\vee,
\end{gather*}
where $(\cdots)^\omega$ (resp.\ $(\cdots)^{\Xi_\epsilon(F)}$) stands for the $\omega$-eigenspace under $Z_G(\eta)(F)$ (resp.\ the space of $\Xi_\epsilon(F)$-invariants).

\section{Endoscopy for linear groups}\label{sec:endoscopy-linear}
The formalism below will be applied to classical groups exclusively.

\subsection{Endoscopic data}
To begin with, let $F$ be a local or global field, and let $L$ be a connected reductive $F$-group. Fix a quasisplit inner twist $L_{\overline{F}} \rightiso L^*_{\overline{F}}$ to define its $L$-group $\Lgrp{L} = \Lgrp{L^*}$, which is endowed with a $\Gamma_F$-stable pinning $\check{T}$, $\check{B}$ and so on. We shall work with Weil forms of $L$-groups.

\index{endoscopic data}
An \emph{endoscopic datum} of $L$ is a quadruplet $\mathbf{H} = (H, \mathcal{H}, s, \hat{\xi})$, or simply $(H, \mathcal{H}, s)$ if no confusion arises, such that
\begin{itemize}
	\item $H$ is a quasisplit $F$-group;
	\item $\mathcal{H}$ sits in a split topological extension of groups
	\[ 1 \to \check{H} \to \mathcal{H} \to \Weil{F} \to 1 \]
	such that the corresponding homomorphism $\Weil{F} \to \mathrm{Out}(\check{H})$ coincides with the one arising from $H$;
	\item $s \in \check{L}$ is semi-simple;
	\item $\hat{\xi}: \mathcal{H} \to \Lgrp{L}$ is an $L$-homomorphism such that
	\begin{itemize}
		\item there exists a $1$-cocycle $a: \Weil{F} \to Z_{\check{L}}$ satisfying $\Ad(s) \circ \hat{\xi}(x) = a(w(x)) \hat{\xi}(x)$ for all $x \in \mathcal{H}$, where $w(x)$ denotes its projection in $\Weil{F}$;
		\item $\hat{\xi}$ induces an isomorphism $\check{H} \rightiso Z_{\check{L}}(s)^\circ$;
		\item the cohomology class of $a$ is trivial (resp.\ locally trivial) if $F$ is local (resp.\ global).
	\end{itemize}
\end{itemize}

The group $H$ is called an endoscopic group of $L$. We say $\mathbf{H}$ is elliptic if $\hat{\xi}\left(Z_{\check{H}}^{\Gamma_F, \circ}\right) \subset Z_{\check{L}}$.

An isomorphism between endoscopic data $\mathbf{H}_1$ and $\mathbf{H}_2$ is an element $\hat{r} \in \check{L}$ such that
\[ \Ad(\hat{r})\hat{\xi}_1\left(\mathcal{H}_1\right) = \hat{\xi}_2 \left( \mathcal{H}_2 \right), \quad \Ad(\hat{r})(s_1) \in s_2 Z_{\check{L}}. \]

Isomorphic endoscopic data are often called equivalent. Denote the set of equivalence classes of endoscopic data of $L$ by $\Endo(L)$, and the subset of elliptic ones as $\Endo_{\elli}(L)$.

Every auto-equivalence of $\mathbf{H}$ induces an outer automorphism of $H$. We denote the resulting subgroup of $\mathrm{Out}(H)$ as $\mathrm{Out}(\mathbf{H})$.

Given $\mathbf{H} \in \Endo(L)$, adjusting by isomorphisms, we may assume that $s \in \check{T}$ and
\[ (\check{T}_H, \check{B}_H) := (\hat{\xi}^{-1}(\check{T}), \hat{\xi}^{-1}(\check{B})) \]
is the Borel pair coming with $\check{H}$. Hence $\hat{\xi}$ induces $\check{T}_H \rightiso \check{T}$. Denote by $W^L$ the absolute Weyl group of $L$. It turns out that there is a $1$-cocycle $\omega_H: \Gamma_F \to W^L$ such that $\hat{\xi} \circ \sigma(\hat{t}) = \omega_H(\sigma) \sigma \hat{\xi}(\hat{t})$ for all $\sigma \in \Gamma_F$ and $\hat{t} \in \check{T}_H$.

The obstruction $\omega_H$ for the $\Gamma_F$-equivariance of $\hat{\xi}|_{\check{T}_H}$ and its dual $\xi: T_{\overline{F}} \rightiso T_{H, \overline{F}}$ is hence $W^L$-valued. Here are some consequences:
\begin{enumerate}
	\item $\hat{\xi}$ restricts to a $\Gamma_F$-equivariant inclusion $Z_{\check{L}} \hookrightarrow Z_{\check{H}}$;
	\item dually, $\xi$ induces an inclusion $\mathfrak{a}_L \hookrightarrow \mathfrak{a}_H$.
\end{enumerate}
In fact, $\mathbf{H}$ is elliptic if and only if $Z_{\check{L}}^{\Gamma_F, \circ} \simeq Z_{\check{H}}^{\Gamma_F, \circ}$, or equivalently $\mathfrak{a}_L \rightiso \mathfrak{a}_H$.

When $F$ is a non-Archimedean local field, $\mathbf{H}$ is said to be \emph{unramified} if
\begin{itemize}
	\item $L$ is unramified and the chosen inner twist $L_{\overline{F}} \rightiso L^*_{\overline{F}}$ is $\identity$;
	\item $H$ is unramified and $\hat{\xi}(\mathcal{H}) \supset Z_{\check{L}}(s)^\circ \rtimes I_F$.
\end{itemize}

\subsection{Transfer}
\index{transfer}
Given $\mathbf{H}$, we have a correspondence of stable semisimple conjugacy classes between $L$ and $H$, or between stable semisimple classes in their Lie algebras.

\begin{definition}\label{def:relevance}
	\index{endoscopic data!relevant}
	An endoscopic datum $\mathbf{H}$ of $L$ is called \emph{relevant} if there exists some $\delta \in \Sigma_{\mathrm{reg}}(H)$ that corresponds to some $\gamma \in \Gamma_{\mathrm{reg}}(L)$.
\end{definition}

When $\mathbf{H}$ is relevant, the $L$-regular locus $H_{L\text{-reg}}$ is open and dense in $H$.

\begin{remark}\label{rem:no-z-extension}
	In all the upcoming applications, $L$ and $L^*$ are classical groups and the endoscopic datum $\mathbf{H}$ is endowed with an isomorphism $\mathcal{H} \simeq \Lgrp{H}$, which will not be altered. As a consequence, we do not need $z$-extensions of $H$. See \cite{Wal10}.
\end{remark}

Assume $F$ to be local of characteristic in what follows. There is a notion of \emph{absolute transfer factors}
\[ \Delta = \Delta_{\mathbf{H}, L}: \Sigma_{\mathrm{reg}}(H) \times \Gamma_{\mathrm{reg}}(L) \to \CC \]
\index{DeltaHL@$\Delta_{\mathbf{H}, L}$}
such that $\Delta(\delta, \gamma) \neq 0$ only when $\delta \leftrightarrow \gamma$. It suffices to consider transfer factors for relevant $\mathbf{H}$.
\index{transfer factor}

In the original definition of Langlands--Shelstad \cite{LS87}, $\Delta_{\mathbf{H}, L}$ is only canonical up to $\CC^\times$. Let us take a relevant $\mathbf{H}$. Once a choice of $\Delta_{\mathbf{H}, L}$ is made, there is the transfer of orbital integrals
\[ \Trans_{\mathbf{H}, L}: \orbI(L) \otimes \mes(L) \to S\orbI(H) \otimes \mes(H) \]
characterized by matching orbital integrals
\[ \Trans_{\mathbf{H}, L}(f_L)(\delta) = \sum_{\gamma \leftrightarrow \delta} \Delta_{\mathbf{H}, L}(\delta, \gamma) f_L(\gamma), \quad \delta \in \Sigma_{\mathrm{reg}}(H). \]
The transfer is continuous for Archimedean $F$.

Dualization yields the transfer of distributions
\[ \trans_{\mathbf{H}, L}: SD(H) \otimes \mes(H)^\vee \to D(L) \otimes \mes(L)^\vee. \]

There are also Lie algebra versions, namely
\begin{align*}
	\Trans_{\mathbf{H}, L}: \orbI(\mathfrak{l}) \otimes \mes(L) & \to S\orbI(\mathfrak{h}) \otimes \mes(H), \\
	\trans_{\mathbf{H}, L}: SD(\mathfrak{h}) \otimes \mes(H)^\vee & \to D(\mathfrak{l}) \otimes \mes(L)^\vee.
\end{align*}

\subsection{Levi subgroups}
Suppose that $R$ is a Levi subgroup of $L$. The following construction is due to Arthur \cite{Ar98}. Consider an endoscopic datum $\mathbf{R}^! = (R^!, \mathcal{R}^!, s^\flat, \hat{\xi})$ of $R$; our notational scheme here is changed in order to conform to \S\ref{sec:endoscopic-data}. Let $s \in s^\flat Z_{\check{R}}^{\Gamma_F}/Z_{\check{L}}^{\Gamma_F}$. Put
\begin{align*}
	L^![s]^\vee & := Z_{\check{L}}(s)^\circ, \\
	\mathcal{L}^![s] & := L^![s]^\vee \cdot \hat{\xi}(\mathcal{R}^!), \\
	\hat{\xi}[s] & : \mathcal{L}^![s] \hookrightarrow \Lgrp{L} \quad \text{is the inclusion}.
\end{align*}
This gives an endoscopic datum $\mathbf{L}^![s] = (L^![s], \mathcal{L}^![s], s, \hat{\xi}[s])$ for $L$. Moreover, it is unramified whenever $\mathbf{R}^!$ is. Define
\[ \Endo_{\mathbf{R}^!}(L) := \{ s \in s^\flat Z_{\check{R}}^{\Gamma_F} / Z_{\check{L}}^{\Gamma_F} : \mathbf{L}^![s] \;\text{is elliptic} \}. \]
This is a finite set by \cite[\S 4]{Ar98}.

Given $s \in \Endo_{\mathbf{R}^!}(L)$, one can view $R^!$ as a Levi subgroup of $L^![s]$. There is a canonical homomorphism $W^{L^![s]}(R^!) \hookrightarrow W^L(R)$, for which $\mathfrak{a}_R \rightiso \mathfrak{a}_{R^!}$ is equivariant. In this case, one can also define the constant
\begin{equation}
	i_{R^!}(L, L^![s]) := \begin{cases}
		\frac{\left( Z_{\check{R}^!}^{\Gamma_F} : Z_{\check{R}}^{\Gamma_F}\right)}{\left( Z_{L^![s]^\vee}^{\Gamma_F} : Z_{\check{L}}^{\Gamma_F} \right)}, & s \in \Endo_{\mathbf{R}^!}(L), \\
		0, & \text{otherwise},
	\end{cases}
\end{equation}
for all $s \in s^\flat Z_{\check{R}}^{\Gamma_F}/Z_{\check{L}}^{\Gamma_F}$.
\index{iRLL@$i_{R^{"!}}(L, L^{"!}[s])$}

When $\mathbf{R}^!$ is relevant, one has the open dense subvariety $R^!_{L-\text{reg}}$ of $L$-regular elements in $R^!$, and similarly on Lie algebras. Every choice of $\Delta_{\mathbf{R}^!, R}$ determines one for $\Delta_{\mathbf{L}^![s], L}$, characterized by
\[ \Delta_{\mathbf{L}^![s], L}(\delta, \gamma) = \Delta_{\mathbf{R}^!, R}(\delta, \gamma) \]
if $\delta \in R^!_{L\text{-reg}}(F)$, $\gamma \in R_{\text{reg}}(F)$ and $\delta \leftrightarrow \gamma$.

\subsection{Complements on absolute transfer factors}\label{sec:endoscopy-rigid}
Let $F$ be a local or global field. Fix the following data:
\begin{itemize}
	\item an $F$-pinning $\mathcal{E}$ of $L^*$,
	\item a pure inner twist $\varphi: L_{\overline{F}} \rightiso L^*_{\overline{F}}$,
	\item the $1$-cocycle $u: \Gamma_K \to L^*(\overline{F})$ associated with $\varphi$.
\end{itemize}

Consider an endoscopic datum $\mathbf{H}$ of $L$. In view of future applications, we still assume $\mathcal{H} \simeq \Lgrp{H}$ and avoid the use of $z$-extensions; see Remark \ref{rem:no-z-extension}.

Assume that $F$ is local hereafter. It has been observed by Kottwitz that $\Delta_{\mathbf{H}, L}$ can be canonically defined using the data $(\mathcal{E}, \varphi, u)$ fixed above. In the case where $L$ is classical, explicit descriptions of the elliptic endoscopic data can be found in \cite[\S 1.5]{Wal10}; the data above can be specified by linear-algebraic data, and there are also explicit formulas for those $\Delta_{\mathbf{H}, L}$. The non-elliptic case follows by passing to Levi subgroups. Whittaker data or additive characters are not used here.

Kottwitz's idea is vastly amplified by Kaletha in \cite[Section 5.3]{Kal16}, but we will only use the simplest case of pure inner twists. Note that Kaletha works with Whittaker data of $L^*$ instead of pinnings.

Given a relevant endoscopic datum $\mathbf{H} = (H, \mathcal{H}, s)$ of $L$, one considers its translate $z\mathbf{H} := (H, \mathcal{H}, zs)$ where $z \in Z_{L^\vee}^{\Gamma_F}$. Although $z\mathbf{H}$ is equivalent to $\mathbf{H}$, Arthur \cite{Ar06} and Hiraga--Saito \cite{HS12} observed that unless $\varphi = \identity_L$, the transfer factors may differ. In order to state the version for pure inner twists we need, denote by
\[ \Delta_{\mathbf{H}} = \Delta_{\mathbf{H}, L^*, u} : H_{\mathrm{reg}}(F) \times L_{\mathrm{reg}}(F) \to \CC \]
the absolute transfer factor determined by $(\mathcal{E}, \varphi, u)$ and $\mathbf{H}$. Denote by $[u] \in \Hm_{\mathrm{ab}}^1(F, L^*)$ the abelianized cohomology class of $u$. We have the perfect pairing
\[ \lrangle{\cdot, \cdot}_{\mathrm{Kott}}: \Hm_{\mathrm{ab}}^1(F, L^*) \times \pi_0\left( Z_{L^\vee}^{\Gamma_F}\right) \to \CC^{\times, \mathrm{tors}} \]
due to Kottwitz where $\CC^{\times, \mathrm{tors}}$ is the group of roots of unity in $\CC$, see \cite[Proposition 1.7.3]{Lab99}. It is normalized to make the following hold.

\begin{theorem}\label{prop:Delta-under-translation}
	For every $z \in Z_{L^\vee}^{\Gamma_F}$ and $\mathbf{H}$, we have
	\[ \Delta_{z\mathbf{H}} = \lrangle{[u], z}_{\mathrm{Kott}} \Delta_{\mathbf{H}}. \]
\end{theorem}
\begin{proof}
	For the general case, see \cite[Section 5.3]{Kal16}, especially (5.1). For elliptic endoscopic data for classical groups over non-Archimedean $F$, this is also documented in \cite[\S 1.11 (1)]{Wal10}.
\end{proof}

On the other hand, once $\mathcal{E}$ and $L^*$ are chosen, $\Delta_{\mathbf{H}, L^*, u}$ depends only on $[u]$; see the last assertion in \cite[Proposition 5.6]{Kal16}.

\begin{remark}\label{rem:principal-endoscopy}
	Even for principal endoscopic data $\mathbf{H}$, namely those with $s \in Z_{L^\vee}^{\Gamma_F}$, translation by $Z_{L^\vee}^{\Gamma_F}$ or change of pure inner twists can still alter $\Delta_{\mathbf{H}}$; see \textit{loc.\ cit.} Thus we cannot naively assume $\Delta = 1$ for all principal endoscopic data of $L$.
\end{remark}

\index{endoscopic data!locally relevant}
Finally, let $F$ be a number field. There is a well-defined absolute adélic or \emph{global transfer factor} $\Delta_{\mathbf{H}, \A_F}$ for every endoscopic data $\mathbf{H}$ of $L$ that is \emph{locally relevant}; see \cite[\S 6]{LS87} and \cite[\S 4]{Ar02}. By fixing a pinned quasisplit $F$-group $L^*$, a pure inner twist $\varphi: L_{\overline{F}} \rightiso L^*_{\overline{F}}$ and the resulting $1$-cocycle $u$ globally (i.e.\ over $F$), $\Delta_{\mathbf{H}, \A_F}$ becomes the product over the local absolute transfer factors. See \cite[\S 5.7]{Kal16}.

\begin{remark}
	Note that $\Delta_{\mathbf{H}, \A_F}$ is invariant under $Z_{L^\vee}^{\Gamma_F}$-twists. It is also independent of the choice of global pure inner twist.
\end{remark}

\section{An instance of nonstandard endoscopy}\label{sec:nonstandard-endoscopy}
\index{nonstandard endoscopy}
To begin with, consider any field $F$ of characteristic $\neq 2$.

\begin{notation}\label{nota:Lbar}
	Given a quasisplit $F$-group $L$ endowed with a decomposition
	\[ L = \SO(2a + 1) \times \SO(2b+1) \times L_0 \]
	where $a, b \in \Z_{\geq 0}$ and $L_0$ is another quasisplit $F$-group. We will write
	\[ \overline{L} := \Sp(2a) \times \Sp(2b) \times L_0. \]
\end{notation}

In order to compare $L$ and $\overline{L}$, one is reduced to the comparison between $\SO(2n+1)$ and $\Sp(2n)$, for various $n \in \Z_{\geq 1}$.

Now recall the \emph{nonstandard endoscopy} between simply connected groups of type $\mathrm{B}_n$ and $\mathrm{C}_n$, which is discussed in \cite[III.6.1]{MW16-1}, \cite[\S 8.2]{Li11} or \cite[\S 5.3]{Li12a}. The situation is
\[\begin{tikzcd}[column sep=large]
	\Spin(2n+1) \arrow[twoheadrightarrow, d, "\text{isogeny}"'] \arrow[dash, dashed, r, "\text{nonstd.\ endo.}" inner sep=0.8em] & \Sp(2n) \\
	\SO(2n+1) &
\end{tikzcd}\]
where the $\SO$ and $\Spin$ groups are always assumed split.

Put $G_1 := \Spin(2n+1)$ and $G_2 := \Sp(2n)$. There is a bijection between $\Sigma_{\mathrm{reg}}(\mathfrak{g}_1)$ and $\Sigma_{\mathrm{reg}}(\mathfrak{g}_2)$, written as $X_1 \leftrightarrow X_2$, namely by matching their eigenvalues. Matching elements have isomorphic centralizers.

We now assume $F$ is a local field of characteristic zero. Consider $f_i \in S\orbI(\mathfrak{g}_i) \otimes \mes(G_i)$ for $i=1,2$. We say $f_1$ and $f_2$ have matching stable orbital integrals if
\[ \forall X_1 \leftrightarrow X_2, \quad f_1^{\Spin(2n+1)}(X_1) = f_2^{\Sp(2n)}(X_2). \]
Therefore, one can say that the transfer factor for nonstandard endoscopy is the constant $1$.

\begin{theorem}[Nonstandard transfer]
	\index{transfer!nonstandard}
	Let $F$ be a local field of characteristic zero. The matching of stable orbital integrals determines an isomorphism
	\[ \Trans_{G_2, G_1}: S\orbI(\mathfrak{g}_1) \otimes \mes(G_1) \rightiso S\orbI(\mathfrak{g}_2) \otimes \mes(G_2) . \]
	It is also a homeomorphism when $F$ is Archimedean.
\end{theorem}
\begin{proof}
	For non-Archimedean $F$, the isomorphism $\Trans_{\overline{L}, L}$ is established in \cite[Proposition 1.8 (ii)]{Wal08}. The case $F=\R$ is included in \cite[V.5.1]{MW16-1}, and the case $F = \CC$ follows by restriction of scalars; the proof is based on Shelstad's description of $S\orbI(\mathfrak{g}_i)$ in terms of jump conditions.
\end{proof}

In studying $SD_{\mathrm{unip}}$ at the level of groups, one can restrict to arbitrarily small and stably invariant neighborhoods of $1$, on which the exponential map is defined. Hence we obtain the dual of the transfer map
\[ \trans_{G_1, G_2}: SD_{\mathrm{unip}}(G_1) \otimes \mes(G_1)^\vee \rightiso SD_{\mathrm{unip}}(G_2) \otimes \mes(G_2)^\vee . \]

The isogeny $\Spin(2n+1, F) \to \SO(2n+1, F)$ being a local homeomorphism, the respective spaces $SD_{\mathrm{unip}}(\cdots)$ and lines $\mes(\cdots)$ can be identified. Hence we obtain
\[ SD_{\mathrm{unip}}(\SO(2n+1)) \otimes \mes(\SO(2n+1))^\vee \rightiso SD_{\mathrm{unip}}(\Sp(2n)) \otimes \mes(\Sp(2n))^\vee . \]

All in all, in the situation of Notation \ref{nota:Lbar}, we arrive at canonical isomorphisms
\begin{gather*}
	\Trans_{\overline{L}, L}: S\orbI(\mathfrak{l}) \otimes \mes(L) \rightiso S\orbI(\overline{\mathfrak{l}}) \otimes \mes(\overline{L}) , \\
	\Trans_{L, \overline{L}} := \Trans_{\overline{L}, L}^{-1}, \\
	\trans_{L, \overline{L}}: SD_{\mathrm{unip}}(L) \otimes \mes(L)^\vee \rightiso SD_{\mathrm{unip}}(\overline{L}) \otimes \mes(\overline{L})^\vee , \\
	\trans_{\overline{L}, L} := \trans_{L, \overline{L}}^{-1}.
\end{gather*}
They will still be called the nonstandard transfer of orbital integrals and distributions, respectively. In this setting, we still have a bijection between $\Sigma_{\mathrm{reg}}(\overline{\mathfrak{l}})$ and $\Sigma_{\mathrm{reg}}(\mathfrak{l})$, and $\Trans_{\overline{L}, L}$ is characterized by matching orbital integrals as before. Everything boils down to the matching between $\so(2n+1)$ and $\syp(2n)$.

\section{Descent of endoscopy}\label{sec:descent-endoscopy}
Fix a local (resp.\ global) field $F$ of characteristic zero. Let $\rev: \tilde{G} \to G(F)$ (resp.\ $\rev: \tilde{G} \to G(\A_F)$) be a group of metaplectic type, and $\mathbf{G}^! \in \Endo_{\elli}(\tilde{G})$. Write $G = \prod_{i \in I} \GL(n_i) \times \Sp(W)$. Suppose that $\epsilon \in G^!(F)_{\mathrm{ss}}$ corresponds to $\eta \in G(F)_{\mathrm{ss}}$ under endoscopy. Upon replacing $\epsilon$ by a stable conjugate, we may and do assume that $G^!_\epsilon$ is quasisplit.

\index{Gepsilon-bar@$\overline{G^{"!}_\epsilon}$}
The group $G^!_\epsilon$ is a product of classical groups by Remark \ref{rem:parameter-classical-group}: it comes endowed with the data as in Notation \ref{nota:Lbar}, hence $\overline{G^!_\epsilon}$ makes sense. Let us write
\begin{align*}
	G & = \prod_{i \in I} \GL(n_i) \times \Sp(W), \\
	G^! & = \prod_{i \in I} \GL(n_i) \times \SO(2n' + 1) \times \SO(2n'' + 1).
\end{align*}
Accordingly, $\epsilon = (\epsilon_{\GL}, \epsilon', \epsilon'')$. The odd orthogonal factors in $G^!_\epsilon$ correspond to the $(+1)$-eigenvalues of $\epsilon'$ and $\epsilon''$.

One of the main results in \cite[\S 7]{Li11} asserts that there is an endoscopic datum $\mathbf{G}_\eta^!$ for $G_\eta$, canonical up to equivalences, such that $G^!_\eta = \overline{G^!_\epsilon}$ and
\[\begin{tikzcd}[column sep=huge]
	G^!_\epsilon \arrow[dashed, dash, "\text{nonstd.\ endo.}", r] & \overline{G^!_\epsilon} \arrow[dashed, dash, "\text{endo.}", r] & G_\eta .
\end{tikzcd}\]

Henceforth, $F$ will be a local field.

Consider elements that are close to zero:
\[ Y \in \mathfrak{g}^!_{\epsilon, \mathrm{reg}}(F), \quad X \in \mathfrak{g}_{\eta, \mathrm{reg}}(F). \]
Denote by $\overline{Y} \in (\overline{\mathfrak{g}^!_\epsilon})_{\mathrm{reg}}(F)$ any element corresponding to $Y$ via nonstandard endoscopy. Then $\exp(Y)\epsilon$ corresponds to $\exp(X)\eta$ via $\mathbf{G}^! \in \Endo_{\elli}(\tilde{G})$ if and only if $\overline{Y}$ corresponds to $X$ via $\mathbf{G}^!_\eta \in \Endo(G_\eta)$.

We may also replace $\eta$ by its stable conjugate $\eta[y]$, where $y \in \mathcal{Y}(\eta)$. Given $\tilde{\eta}$ and $y$, one defines a canonical element $\tilde{\eta}[y] \in \rev^{-1}(\eta[y])$ according to \cite[Définition 4.1.3]{Li15}; $\tilde{\eta}[y]$ and $\eta[y]$ depend only on the class of $y$ in $\dot{\mathcal{Y}}(\eta)$. For each $y \in \dot{\mathcal{Y}}[y]$, we choose a transfer factor $\Delta(y)$ for the endoscopic datum $\mathbf{G}^!_{\eta[y]}$ of $G_{\eta[y]}$, if it is relevant (Definition \ref{def:relevance}).

In the relevant case, we may choose $X[y] \in \mathfrak{g}_{\eta[y], \mathrm{reg}}(F)$ and $\overline{Y} \in \overline{\mathfrak{g}^!_\epsilon}_{,\mathrm{reg}}(F)$ which are sufficiently close to zero and $\overline{Y} \leftrightarrow X[y]$. We also take $Y \in \mathfrak{g}^!_{\epsilon, \mathrm{reg}}(F)$ such that $\overline{Y} \leftrightarrow Y$; in particular $Y$ is also close to zero. Define
\begin{equation}\label{eqn:dy}
	d(y) := \begin{cases}
		\dfrac{ \Delta_{\mathbf{G}^!, \tilde{G}}(\exp(Y) \epsilon, \exp(X[y])\tilde{\eta}[y]) }{ \Delta(y)(\exp \overline{Y}, \exp X[y] ) }, & \mathbf{G}^!_{\eta[y]} \;\text{is relevant}, \\
		0, & \text{otherwise}.
	\end{cases}
\end{equation}
\index{dy@$d(y)$}

The descent of transfer factors \cite[Théorème 7.23]{Li11} asserts that $d(y)$ is constant for $\overline{Y}$ and $X$ sufficiently close to zero, when $\mathbf{G}^!_{\eta[y]}$ is relevant. This holds in both the non-Archimedean and Archimedean cases. We repeat it as follows.

\begin{theorem}
	The value of $d(y)$ is independent of $\overline{Y}$ and $X$. It depends only on the choice of $\Delta(y)$.
\end{theorem}

\begin{lemma}[Cf.\ {\cite[III.5.2]{MW16-1}}]\label{prop:transfer-descent-y}
	Assume $F$ to be non-Archimedean. In the situation above, suppose that $\delta \in SD_{\mathrm{geom}}(\mathfrak{g}^!_\epsilon) \otimes \mes(G^!_\epsilon)$ is given and $\Supp(\delta)$ is close to zero. For each $y \in \dot{\mathcal{Y}}(\eta)$, define
	\begin{align*}
		\delta^{G^!} & := S\desc^{G^!, *}_{\epsilon}(\delta) \; \in SD_{\mathrm{geom}}(G^!) \otimes \mes(G^!)^\vee , \\
		\delta[y] & := \trans_{\mathbf{G}^!_{\eta[y]}, G_{\eta[y]}} \trans_{G^!_\epsilon, \overline{G^!_\epsilon}} (\delta) \; \in D_{\mathrm{geom}}(\mathfrak{g}_{\eta[y]}) \otimes \mes(G_{\eta[y]})^\vee , \\
		\delta[y]^{\tilde{G}} & := \desc^{\tilde{G}, *}_{\tilde{\eta}[y]}\left( \delta[y] \right) \; \in D_{\mathrm{geom}, -}(\tilde{G}) \otimes \mes(G)^\vee .
	\end{align*}
	When $\mathbf{G}^!_{\eta[y]}$ is relevant, $\trans_{\mathbf{G}^!_{\eta[y]}, G_{\eta[y]}}$ (on the level of Lie algebras) is defined relative to the $\Delta(y)$. Then we have
	\[ \trans_{\mathbf{G}^!, \tilde{G}}\left( \delta^{G^!} \right) = \sum_{y \in \dot{\mathcal{Y}}(\eta)} d(y) \delta[y]^{\tilde{G}}. \]
\end{lemma}
\begin{proof}
	To begin with, let us assume $\tilde{G} = \Mp(W)$. Identify $\delta$ with some $Y \in \Sigma_{\mathrm{reg}}(\mathfrak{g}^!_\epsilon)$ via exponential map by fixing Haar measures. For all $f \in \orbI_{\asp}(\tilde{G}) \otimes \mes(G)$, we obtain
	\begin{align*}
		\lrangle{ \trans_{\mathbf{G}^!, \tilde{G}}(\delta^{G^!}), f } & = \lrangle{ \exp(Y)\epsilon, \Trans_{\mathbf{G}^!, \tilde{G}}(f) } \\
		& = \sum_{y \in \dot{\mathcal{Y}}(\eta)} d(y) \lrangle{ \delta[y], \desc^{\tilde{G}}_{\tilde{\eta}[y]} f } \\
		& = \sum_{y \in \dot{\mathcal{Y}}(\eta)} d(y) \lrangle{ \delta[y]^{\tilde{G}}, f }
	\end{align*}
	where the second equality follows from the proof of \cite[Lemme 8.5]{Li11}. It is routine to deduce the general case when $\tilde{G}$ is only of metaplectic type.
\end{proof}

\section{Complements on diagrams}\label{sec:diagram-2}
We resume the discussion on ``diagrams'' in \S\ref{sec:diagram}. In particular, $F$ is allowed to be a field of characteristic zero in the first part of this section.

We begin with the regular semisimple case.

\begin{lemma}\label{prop:diagram-reg}
	Let $\eta \in G(F)_{\mathrm{reg}}$ and $\epsilon \in G^!(F)_{\mathrm{reg}}$. If $\epsilon$ corresponds to $\eta$ via endoscopy, then there exists a diagram joining $\epsilon$ and $\eta$.
\end{lemma}
\begin{proof}
	Let $T^! := G^!_\epsilon$. Choose $h \in G^!(\overline{F})$ such that $\Ad(h)$ induces $T^!_{\overline{F}} \rightiso T_{G^!, \overline{F}}$. We use the following terminology from Galois descent: for any $F$-varieties $A$, $B$ and a morphism $f: A_{\overline{F}} \to B_{\overline{F}}$, we write ${}^\sigma f = \sigma f \sigma^{-1} : A_{\overline{F}} \to B_{\overline{F}}$, for every $\sigma \in \Gamma_F$. Since $T^!$ and $T_{G^!}$ are $F$-tori,
	\[ {}^\sigma \Ad(h)|_{T^!_{\overline{F}}} = w_\sigma \circ \Ad(h)|_{T^!_{\overline{F}}}, \quad w_\sigma \in W^{G^!} \]
	and $\sigma \mapsto w_\sigma$ is a $1$-cocycle. We also have ${}^\sigma \Phi = \Phi$. Since $\Phi$ is equivariant with respect to $W^{G^!} \to W^G$, as embeddings $T^!_{\overline{F}} \to G_{\overline{F}}$ we have
	\[ {}^\sigma\left( \Phi \circ \Ad(h)|_{T^!_{\overline{F}}} \right) = w_\sigma \circ \Phi \circ \Ad(h)|_{T^!_{\overline{F}}}, \quad \sigma \in \Gamma_F. \]
	Using Kottwitz's result on embeddings of tori \cite[Corollary 2.2]{Ko82}, we conclude that there exists $g \in G(\overline{F})$ such that $\Ad(g)^{-1} \Phi \Ad(h)|_{T^!_{\overline{F}}}$ is defined over $F$, written as $\xi_{T^!, T}: T^! \rightiso T$ for some maximal $F$-torus $T \subset G$.
	
	Note that $\Ad(g): T_{\overline{F}} \rightiso T_{G, \overline{F}}$ also satisfies ${}^\sigma \Ad(g)|_{T_{\overline{F}}} = v_\sigma \circ \Ad(g)|_{T_{\overline{F}}}$ for some $v_\sigma \in W^G$. By comparing with the cocycle $w_\sigma$, we see that $v_\sigma = w_\sigma$ since $\Ad(g)^{-1} \Phi \Ad(h)|_{T^!_{\overline{F}}}$ is defined over $F$.
	
	Next, note that $\Psi = \tau \circ \Phi$ where $\tau$ is the translation by the $W^{G^!}$-invariant element $(1, 1, -1) \in T(F)$. Hence
	\begin{align*}
		{}^\sigma \left(\Ad(g)^{-1} \Psi \Ad(h)|_{T^!_{\overline{F}}} \right) & = \Ad(g)^{-1} v_\sigma^{-1} \tau \Phi w_\sigma \Ad(h)|_{T^!_{\overline{F}}} \\
		& = \Ad(g)^{-1} v_\sigma^{-1} w_\sigma \tau \Phi \Ad(h)|_{T^!_{\overline{F}}} \\
		& = \Ad(g)^{-1} \Psi \Ad(h)|_{T^!_{\overline{F}}}.
	\end{align*}
	This yields a morphism $\tilde{\xi}_{T^!, T}: T^! \rightiso T$ defined over $F$. By construction, $\eta' := \tilde{\xi}_{T^!, T}(\epsilon)$ corresponds to $\epsilon$ under endoscopy. Hence there exists $g_1 \in G(\overline{F})$ with
	\[ \Ad(g_1)(\eta) = \eta', \quad \forall \sigma \in \Gamma_F, \; {}^\sigma g_1 \in T(\overline{F}) g_1. \]
	Therefore $\Ad(g_1): T_1 := G_\eta \rightiso T$ is a stable conjugacy between maximal $F$-tori in $G$, sending $\eta$ to $\eta'$ and defined over $F$.
	
	Set $B^! := \Ad(h)^{-1} B_{G^!}$ and $B := \Ad(gg_1)^{-1} B_G$. We conclude that $(\epsilon, B^!, T^!, B, T_1, \eta)$ is a diagram; the corresponding $\xi_{T^!, T_1}$ (resp.\ $\tilde{\xi}_{T^!, T}$) being $\Ad(g_1)^{-1}\xi_{T^!, T}$ (resp.\ $\Ad(g_1)^{-1}\tilde{\xi}_{T^!, T}$).
\end{proof}

\begin{corollary}\label{prop:diagram-from-epsilon}
	For any $\epsilon \in G^!(F)_{\mathrm{ss}}$ and a maximal $F$-torus $T^! \subset G^!$ containing $\epsilon$, there exists a diagram joining $\epsilon$ to some $\eta \in G(F)_{\mathrm{ss}}$, with $T^!$ being part of the data.
\end{corollary}
\begin{proof}
	Pick $\delta \in T^!_{G\text{-reg}}(F)$. It corresponds to some $\gamma \in G_{\mathrm{reg}}(F)$. Lemma \ref{prop:diagram-reg} furnishes a diagram $(\delta, B^!, T^!, B, T, \gamma)$. Now take $\eta := \tilde{\xi}_{T^!, T}(\epsilon)$.
\end{proof}

We give a supplementary result on the transfer of ellipticity.
\begin{lemma}\label{prop:ellipticity-transfer}
	Suppose that $\epsilon \in G^!(F)_{\mathrm{ss}}$ corresponds to $\eta \in G(F)_{\mathrm{ss}}$ under endoscopy. If $\epsilon$ is elliptic in $G^!$, then $\eta$ is elliptic in $G$.
\end{lemma}
\begin{proof}
	Without loss of generality, we may assume $G = \Sp(W)$. If $\epsilon$ is elliptic, its parameter introduced in \S\ref{sec:ss-classes} has no hyperbolic indices by Remark \ref{rem:parameter-elliptic-classical-group}. The same holds for the parameter of $\eta$, hence $\eta$ is elliptic by Corollary \ref{prop:parameter-elliptic}.
\end{proof}

The converse of Lemma \ref{prop:ellipticity-transfer} is not true. For example, take $G = \Sp(2)$ and $G^! = \SO(3)$, arising from $(1, 0) \in \Endo_{\elli}(\tilde{G})$. The elliptic element $\eta = -1 \in G(F)$ corresponds to some $\epsilon \in G^!(F)$ with $G^!_\epsilon$ isomorphic to the split $\SO(2) \simeq \Gm$, so $\epsilon$ is not elliptic.

\begin{proposition}\label{prop:diagram-relevant}
	Assume that $F$ is a local field. Let $\eta \in G(F)_{\mathrm{ss}}$ and $\epsilon \in G^!(F)_{\mathrm{ss}}$. Suppose that $\epsilon$ corresponds to $\eta$ via endoscopy, $G^!_\epsilon$ is quasisplit, and the endoscopic datum $\mathbf{G}^!_\eta$ of $G_\eta$ is relevant (see \S\ref{sec:descent-endoscopy}), then there exists a diagram $(\epsilon, B^!, T^!, B, T, \eta)$ joining $\epsilon$ and $\eta$ for $\mathbf{G}^! \in \Endo_{\elli}(\tilde{G})$, compatible with descent in the sense that
	\begin{itemize}
		\item $T^! \subset G^!_\epsilon$ and $T \subset G_\eta$,
		\item $T^!$ transfer to a maximal torus $\overline{T^!} \subset \overline{G^!_\epsilon}$ under nonstandard endoscopy,
		\item $\overline{T^!}$ transfers to $T$ under $\mathbf{G}^!_\eta \in \Endo(G_\eta)$.
	\end{itemize}
	In fact, we can take any $T^!$, $\overline{T^!}$ and $T$ subject to the requirements above.
\end{proposition}
\begin{proof}
	By assumption and the descent of endoscopy, there exist $Y \in \mathfrak{g}^!_{\epsilon, \mathrm{reg}}(F)$, $\overline{Y} \in \overline{\mathfrak{g}^!_\epsilon}_{\mathrm{reg}}(F)$ and $X \in \mathfrak{g}_{\eta, \mathrm{reg}}(F)$, arbitrarily close to $0$, such that
	\begin{itemize}
		\item $Y \leftrightarrow \overline{Y} \leftrightarrow X$, the first (resp.\ second) correspondence being relative to nonstandard endoscopy (resp.\ $\mathbf{G}^!_\eta \in \Endo(G_\eta)$),
		\item $\exp(Y)\epsilon \leftrightarrow \exp(X)\eta$ relative to $\mathbf{G}^! \in \Endo_{\elli}(\tilde{G})$.
	\end{itemize}

	Let $T^! \subset G^!_\epsilon$, $\overline{T^!} \subset \overline{G^!_\epsilon}$ and $T \subset G_\eta$ be the centralizers of $Y$, $\overline{Y}$ and $X$ respectively. It follows that $T^! = G^!_{\exp(Y)\epsilon}$ and $T = G_{\exp(X)\eta}$.
	
	By Lemma \ref{prop:diagram-reg}, we obtain a diagram with Borel subgroups $B^! \supset T^!$ and $B \supset T$, for each $(Y, \overline{Y}, X)$, such that the $\tilde{\xi}_{T^!, T}$ maps $\exp(Y)\epsilon$ to $\exp(X)\eta$. By the finiteness of Borel subgroups containing a given maximal torus, we may now take sequences $Y_n, \overline{Y}_n, X_n \to 0$, such that there is a single choice of $(B^!, B)$ giving rise to the $\tilde{\xi}_{T^!, T}: \exp(Y_n) \epsilon \mapsto \exp(X_n)\eta$, for all $n$.
	
	Pass to the limit to conclude $\tilde{\xi}_{T^!, T}(\epsilon) = \eta$. This is the required $(\epsilon, B^!, T^!, B, T, \eta)$.
\end{proof}

The relevance of $\mathbf{G}^!_\eta \in \Endo(G_\eta)$ is automatic if we assume $G_\eta$ is quasisplit. This is a direct consequence of \cite[\S 2]{Ko82}. We also make the following observation.

\begin{proposition}\label{prop:relevance-elliptic}
	Let $F$ be a local field, $\eta \in G(F)_{\mathrm{ss}}$ and $\epsilon \in G^!(F)_{\mathrm{ss}}$. Suppose that $\epsilon$ corresponds to $\eta$ via endoscopy, and $\epsilon$ in elliptic in $G^!$. Then the endoscopic datum $\mathbf{G}^!_\eta$ obtained by descent is relevant.
\end{proposition}
\begin{proof}
	In this case, $\eta$ is elliptic by Lemma \ref{prop:ellipticity-transfer}, so $\mathbf{G}^!_\eta$ is elliptic. Without loss of generality, we assume $G = \Sp(W)$ in what follows.

	There exists $y \in \mathcal{Y}(\eta)$ such that $G_{\eta[y]}$ is quasisplit. Since the parameter of $\epsilon$ (see \S\ref{sec:ss-classes}) contains no hyperbolic indices, the concrete description of $G^!_\epsilon$ implies that there is an elliptic maximal torus $T^! \subset G^!_\epsilon$; its nonstandard transfer $\overline{T^!} \subset \overline{G^!_\epsilon}$ is still elliptic. In turn, $\overline{T^!}$ transfers to $T' \subset G_{\eta[y]}$ via standard endoscopy; $T'$ is still elliptic in $G_{\eta[y]}$.
	
%	By \cite[10.2 Lemma]{Ko86}, $T'_{\mathrm{ad}} \subset G_{\eta[y], \mathrm{AD}}$ (the image of $T'$, which is elliptic) transfers to $T_{\mathrm{ad}} \subset G_{\eta, \mathrm{AD}}$ under the inner twist $\Ad(y)$. Denote by $T \subset G_\eta$ the preimage of $T_{\mathrm{ad}}$, then $T'$ transfers to $T$ under inner twist, as desired.

	By \cite[10.2 Lemma]{Ko86}, $T' \subset G_{\eta[y]}$ transfers to some $T \subset G_\eta$ under a pure inner twist which lies in the same class as $\Ad(y)$ in $\Hm^1(F, G_\eta)$. This shows that $\overline{T^!} \subset \overline{G^!_\epsilon}$ transfers to $T \subset G_\eta$, as desired.
\end{proof}

\section{Descent with Levi subgroups}\label{sec:descent-endoscopy-Levi}
Return to the setting in \S\ref{sec:descent-endoscopy} and consider $\rev: \tilde{G} \twoheadrightarrow G(F)$ together with a Levi subgroup $M \subset G$. Let $\mathbf{M}^! \in \Endo_{\elli}(\tilde{M})$. We also fix $\epsilon \in M^!(F)_{\mathrm{ss}}$ and $\eta \in M(F)_{\mathrm{ss}}$ in endoscopic correspondence. Replacing $\epsilon$ by a stable conjugate if necessary, we may assume $M^!_\epsilon$ is quasisplit.

Using $\eta$, we define the pointed sets $\mathcal{Y}^M$, $\mathcal{Y}$ and $\dot{\mathcal{Y}}^M$, $\dot{\mathcal{Y}}$.

For each $s \in \Endo_{\mathbf{M}^!}(\tilde{G})$, we have the corresponding $\mathbf{G}^![s] \in \Endo_{\elli}(\tilde{G})$. Thus $\epsilon[s]$ corresponds to $\eta$ with respect to $\mathbf{G}^![s] \in \Endo_{\elli}(\tilde{G})$.

\begin{lemma}
	For each $s$, the group $G^![s]_{\epsilon[s]}$ is quasisplit.
\end{lemma}
\begin{proof}
	It contains $M^!_{\epsilon[s]} = M^!_\epsilon$ as a Levi subgroup. The quasisplit-ness of $M^!_\epsilon$ amounts to the existence of a maximal $F$-torus $T^!$ which is a simultaneously a minimal Levi. Then $T^!$ has the same properties relative to $G^!_{\epsilon[s]}$, hence $G^!_{\epsilon[s]}$ is quasisplit as well.
\end{proof}

For all $y \in \mathcal{Y}$ and $s \in \Endo_{\mathbf{M}^!}(\tilde{G})$, we obtain an endoscopic datum $\mathbf{G}^!_{\eta[y]}$ of $G_{\eta[y]}$ obtained from $\mathbf{G}^![s]$ by descent at $\epsilon[s]$ and $\eta[y]$, canonical up to equivalence, as reviewed in \S\ref{sec:descent-endoscopy}. Given $y \in \mathcal{Y}$, let us pick a representative in each equivalence class of relevant endoscopic data of $G_{\eta[y]}$ arising from some $s \in \Endo_{\mathbf{M}^!}(\tilde{G})$, and fix a transfer factor for it, denoted hereafter by $\Delta(\overline{s}, y)$.

One defines the following constants of proportionality.

\begin{enumerate}
	\item For each $y \in \mathcal{Y}^M$, we apply \eqref{eqn:dy} to $\mathbf{M}^! \in \Endo_{\elli}(\tilde{M})$, $\epsilon \in M^!(F)_{\mathrm{ss}}$ and $\eta[y] \in M(F)_{\mathrm{ss}}$ to define $d^M(y)$. This depends on the choice of transfer factor $\Delta(y)$ for the $\mathbf{M}^!_{\eta[y]} \in \Endo(M_{\eta[y]})$ obtained by descent, whenever it is relevant.
	\item For each $s \in \Endo_{\mathbf{M}^!}(\tilde{G})$ and $y \in \mathcal{Y}$, we apply \eqref{eqn:dy} to $\mathbf{M}^! \in \Endo_{\elli}(\tilde{M})$ to define
	\begin{equation}\label{eqn:dsy}
		d(s, y) := \begin{cases}
			\dfrac{ \Delta_{\mathbf{G}^![s], \tilde{G}}(\exp(Y) \epsilon[s], \exp(X[y])\tilde{\eta}[y]) }{ \Delta(\overline{s}, y)(\exp(\overline{Y}), \exp(X[y]) ) }, & \mathbf{G}^!_{\eta[y]} \;\text{is relevant}, \\
			0, & \text{otherwise},
		\end{cases}
	\end{equation}
	where $Y \leftrightarrow \overline{Y} \leftrightarrow X[y]$ are all close to zero.
	\index{dsy@$d(s, y)$}
\end{enumerate}
This is similar to \cite[III.5.4 (2)]{MW16-1}; our case is simpler since there are no $z$-extensions.

These factors depend on various choices. The factor $d(s, y)$ is more delicate since different choices of $s$ (which determines the numerator in \eqref{eqn:dsy}) may give rise to the same equivalence class of endoscopic data of $G_{\eta[y]}$ (which determines the denominator). Ergo, there is no straightforward normalization for $d(s, y)$.

Nevertheless, the problem simplifies when $y \in \mathcal{Y}^M$. In this case, $d^M(y) \neq 0$ if the endoscopic datum $\mathbf{M}^!_{\eta[y]}$ of $M_{\eta[y]}$ obtained by descent is relevant. We explicate as follows.

\begin{proposition}\label{prop:dsy}
	Let $y \in \mathcal{Y}^M$. The transfer factors $\Delta(\overline{s}, y)$ may be chosen so that $d(s, y) = d^M(y)$ when $\mathbf{M}^!_{\eta[y]}$ is relevant, for all $s \in \Endo_{\mathbf{M}^!}(\tilde{G})$.
\end{proposition}
\begin{proof}
	Suppose $\mathbf{M}^!_{\eta[y]}$ is relevant. We may choose $Y \in \mathfrak{m}^!_{\epsilon, \mathrm{reg}}(F)$ corresponding to $\overline{Y} \in \overline{\mathfrak{m}^!_\epsilon}_{, \mathrm{reg}}(F)$, such that $\overline{Y}$ transfers to $X[y] \in \mathfrak{m}_{\eta[y], \mathrm{reg}}(F)$. By \cite[Proposition 3.3.6]{Li12a},
	\[ \Delta_{\mathbf{G}^!, \tilde{G}}\left( \exp(Y) \epsilon[s], \exp(X[y]) \tilde{\eta}[y] \right) = \Delta_{\mathbf{M}^!, \tilde{M}}\left( \exp(Y) \epsilon, \exp(X[y]) \tilde{\eta}[y] \right). \]
	
	Equation \eqref{eqn:dsy} simplifies into
	\[ d(s, y) = \frac{\Delta_{\mathbf{M}^!, \tilde{M}}\left( \exp(Y) \epsilon, \exp(X[y]) \tilde{\eta}[y] \right)}{\Delta(\overline{s}, y)(\exp(\overline{Y}), \exp(X[y]))}. \]

	Consider the denominator: we are in the situation
	\[\begin{tikzcd}
		\overline{G^![s]_{\epsilon[s]}} \arrow[dashed, dash, r, "\text{endo.}"] & G_{\eta[y]} \\
		\overline{M^!_\epsilon} = \overline{M^!_{\epsilon[s]}} \arrow[dashed, dash, r, "\text{endo.}"'] \arrow[hookrightarrow, u, "\text{Levi}"] & M_{\eta[y]} \arrow[hookrightarrow, u, "\text{Levi}"'] ;
	\end{tikzcd}\]
	see Remark \ref{rem:descent-homomorphism} below for a more precise description.
	
	By the parabolic descent of transfer factors in standard endoscopy, we may normalize $\Delta(\overline{s}, y)$, for all $s$, so that
	\[ \Delta(\overline{s}, y)(\exp(\overline{Y}), \exp(X[y])) = \Delta(y)\left( \exp(\overline{Y}), \exp(X[y]) \right) \]
	where $\Delta(y)$ is the chosen transfer factor for $\mathbf{M}^!_{\eta[y]}$. A comparison with \eqref{eqn:dy} yields $d(s, y) = d^M(y)$.
\end{proof}

For future reference, we record another definition pertaining to the method of descent on the dual side.

\begin{definition}\label{def:eLRmu}
	\index{eLRmu@$e^L_R(\mu)$}
	Recall from \cite[p.362]{MW16-1} that for all quasisplit $F$-group $L$, a Levi subgroup $R \subset L$ and an element $\mu \in R(F)_{\mathrm{ss}}$, we have a canonical homomorphism
	\[ \mathrm{descent}: Z_{\check{R}}^{\Gamma_F} / Z_{\check{L}}^{\Gamma_F} \to Z_{R_\mu^\vee}^{\Gamma_F} / Z_{L_\mu^\vee}^{\Gamma_F}. \]
	
	When $\mu$ is elliptic in $L$ (hence in $R$), one defines
	\[ e^L_R(\mu) := \left| \Ker\left( Z_{\check{R}}^{\Gamma_F} / Z_{\check{L}}^{\Gamma_F} \xrightarrow{\text{descent}} Z_{R_\mu^\vee}^{\Gamma_F} / Z_{L_\mu^\vee}^{\Gamma_F} \right) \right|^{-1}. \]
	By \textit{loc.\ cit.}, the homomorphism is surjective with finite kernel, thus $e^L_R(\mu)$ is well-defined.
\end{definition}

\begin{remark}\label{rem:descent-homomorphism}
	We will have to consider the descent of metaplectic endoscopic data from the dual side. Given $\mathbf{M}^! \in \Endo_{\elli}(\tilde{M})$, $\epsilon \in M^!(F)_{\mathrm{ss}}$ and $\eta \in M(F)_{\mathrm{ss}}$ as before, such that $\epsilon \leftrightarrow \eta$ and $M^!_\epsilon$ is quasisplit, we take
	\begin{itemize}
		\item an elliptic element $s^\flat \in \tilde{M}^\vee$ that determines $\mathbf{M}^!$, trivial on the $\GL$-factors;
		\item $\overline{s^\flat} \in M_\eta^\vee$: part of the $\mathbf{M}^!_\eta \in \Endo(M_\eta)$ obtained by descent at $\epsilon$ and $\eta$ (see \cite[p.524]{Li12a} for a more explicit description).
	\end{itemize}

	The metaplectic analogue of Definition \ref{def:eLRmu} is the composite
	\[ \mathrm{descent}: Z_{\tilde{M}^\vee}^\circ \big/ Z_{\tilde{G}^\vee}^\circ \simeq Z_{\check{M}}^\circ \big/ Z_{\check{G}}^\circ \to Z_{\check{M}} \big/ Z_{\check{G}} \to Z_{M_\eta^\vee}^{\Gamma_F} \big/ Z_{G_\eta^\vee}^{\Gamma_F} \]
	that was denoted by $\tau$ and explained in detail in \cite[pp.525--526]{Li12a}. Now interpret $\Endo_{\mathbf{M}^!}(\tilde{G})$ as a subset of $s^\flat Z_{\tilde{M}^\vee}^\circ / Z_{\tilde{G}^\vee}^\circ$ by Proposition \ref{prop:Endo-s-interpretation}. For any $s = s^\flat t Z_{\tilde{G}^\vee}^\circ$ in that subset, let
	\begin{equation}\label{eqn:descent-s}
		\overline{s} := \overline{s^\flat} \cdot \mathrm{descent}(t) \; \in \overline{s^\flat} \cdot  Z_{M_\eta^\vee}^{\Gamma_F} \big/ Z_{G_\eta^\vee}^{\Gamma_F}.
	\end{equation}
	
	It is shown in \cite[Lemme 6.3.2]{Li12a} that the $\mathbf{G}_\eta \in \Endo(G_\eta)$ obtained by descent at $\epsilon[s]$ and $\eta$ is equivalent to the one associated with $\overline{s} \in \Endo_{\mathbf{M}^!_\eta}(G_\eta)$ by Arthur's recipe reviewed in \S\ref{sec:endoscopy-linear}.
\end{remark}
\chapter{The local geometric theorem}\label{sec:geom-matching}
Fix a local field $F$ of characteristic zero, a non-trivial additive character $\psi$ of $F$, and a symplectic $F$-vector space $(W, \lrangle{\cdot|\cdot})$ to form the covering
\[ 1 \to \bmu_8 \to \tilde{G} \xrightarrow{\rev} G(F) \to 1, \quad G := \Sp(W). \]
As in \S\ref{sec:linear-algebra}, we also fix a symplectic basis for $W$ to obtain the corresponding minimal Levi subgroup $M_0$. This gives rise to an invariant positive-definite quadratic form on $\mathfrak{a}_0$, thus on various $\mathfrak{a}^G_M$ (Definition \ref{def:invariant-quadratic-form}).

In fact, we will work with covering groups $\tilde{G}$ of metaplectic type in order to perform inductive arguments.

Given $M \in \mathcal{L}(M_0)$, we will introduce the invariant weighted orbital integrals $I_{\tilde{M}}(\tilde{\gamma}, f)$, where $f \in \orbI_{\asp}(\tilde{G}) \otimes \mes(G)$.
\begin{itemize}
	\item When $F$ is non-Archimedean, $\tilde{\gamma}$ can be any element in $D_{\mathrm{geom}, -}(\tilde{M}) \otimes \mes(M)^\vee$.
	\item The case of Archimedean $F$ is more subtle. Although $I_{\tilde{M}}(\tilde{\gamma}, f)$ can be defined for ``orbital'' distributions $\tilde{\gamma} \in D_{\mathrm{orb}, -}(\tilde{M}) \otimes \mes(M)^\vee$, that subspace is not large enough to accommodate all endoscopic transfers. On the other hand, general $\tilde{\gamma} \in D_{\mathrm{geom}, -}(\tilde{M}) \otimes \mes(M)^\vee$ might involve derivatives in the normal direction of orbits. Nonetheless, it is still possible to define $I_{\tilde{M}}(\tilde{\gamma}, f)$ when $\tilde{\gamma}$ is $\tilde{G}$-equisingular in the sense of Definition \ref{def:equisingular}.
\end{itemize}

Given $\mathbf{M}^! \in \Endo_{\elli}(\tilde{M})$, we set (Definition \ref{def:IEndo-1})
\[ I^{\Endo}_{\tilde{M}}(\mathbf{M}^!, \delta, f) := \sum_{s \in \Endo_{\mathbf{M}^!}(\tilde{G})} i_{M^!}(\tilde{G}, G^![s]) S^{G^![s]}_{M^!}\left(\delta[s], B^{\tilde{G}}, f^{G^![s]} \right) \]
where $f^{G^![s]} := \Trans_{\mathbf{G}^![s], \tilde{G}}(f)$ is the transfer, $\delta \in SD_{\mathrm{geom}}(M^!) \otimes \mes(M^!)^\vee$ is assumed to be $\tilde{G}$-equisingular in the Archimedean case, and $\delta[s] := \delta \cdot z[s]$ is the central twist. Here appears a new ingredient $B^{\tilde{G}}$ in the stable weighted orbital integral $S^{G^![s]}_{M^!}(\cdots)$, to be defined in \S\ref{sec:local-geometric}. In brief, $B^{\tilde{G}}$ prescribes the rescaling of roots in $G^![s]_{\delta[s]}$; it first appeared in the study of twisted endoscopy in \cite[II]{MW16-1}.

From $I^{\Endo}_{\tilde{M}}(\mathbf{M}^!, \delta, f)$, one builds the endoscopic counterpart $I^{\Endo}_{\tilde{M}}(\tilde{\gamma}, f)$ of $I_{\tilde{M}}(\tilde{\gamma}, f)$ by expressing $\tilde{\gamma}$ via transfers from various $(\mathbf{M}^!, \delta)$ (Definition--Proposition \ref{def:I-geom-Endo}). The local geometric Theorem \ref{prop:local-geometric} asserts that
\[ I^{\Endo}_{\tilde{M}}(\tilde{\gamma}, f) = I_{\tilde{M}}(\tilde{\gamma}, f) \]
for all $f$ and $\tilde{\gamma}$, the latter assumed to be $\tilde{G}$-equisingular in the Archimedean case.

The proof of this theorem will only be achieved at the end of this work. At this stage, we prove in \S\ref{sec:descent-orbint-Endo} that $I^{\Endo}_{\tilde{M}}(\tilde{\gamma}, f)$  satisfies the same descent formula as $I_{\tilde{M}}(\tilde{\gamma}, f)$; then we reduce the $\tilde{G}$-equisingular case (for both Archimedean and non-Archimedean $F$) to the $\tilde{G}$-regular case in \S\ref{sec:local-geom-reduction-G-reg}. The combinatorial technique for proving the endoscopic descent formula is a recurring theme in this work. It will be invoked repeatedly later on.

In the unramified situation, the weighted fundamental lemma for general $\tilde{\gamma}$ is stated as Theorem \ref{prop:LFP-general}. Its proof as well as further reductions for Theorem \ref{prop:local-geometric} demand the theory of germs.

The results in \S\S\ref{sec:orbint-weighted-nonArch}--\ref{sec:orbint-weighted-Arch} will only be stated for coverings of metaplectic type, although they can also be adapted to general coverings; see Remark \ref{rem:weighted-general-covering}.

\section{Weighted orbital integrals: non-Archimedean case}\label{sec:orbint-weighted-nonArch}
Assume $F$ to be non-Archimedean.

\subsection{Non-invariant version}
Choose a maximal compact subgroup $K \subset G(F)$ in good position relative to $M_0$. Let $M \in \mathcal{L}(M_0)$ and $\tilde{\gamma} \in \tilde{M}$. We write the Jordan decomposition of $\tilde{\gamma}$ reviewed in \S\ref{sec:covering} as
\[ \tilde{\gamma} = \tilde{\eta}u = u\tilde{\eta}, \quad \tilde{\eta} := \tilde{\gamma}_{\mathrm{ss}}. \]

Let
\[ v_M: M(F) \backslash G(F) / K \to \R_{\geq 0} \]
be Arthur's weight factor attached to the $(G, M)$-family $v_P(\lambda, x) := e^{-\lrangle{\lambda, H_P(x)}}$, where $P \in \mathcal{P}(M)$ and $\lambda \in i\mathfrak{a}^*_M$.

Consider $\tilde{\gamma} \in \tilde{M}$ with Jordan decomposition $\tilde{\gamma} = \tilde{\eta}u = u\tilde{\eta}$. Suppose that $\tilde{\gamma}$ is $\tilde{G}$-equisingular (Definition \ref{def:equisingular}), we define
\[ J_{\tilde{M}}(\tilde{\gamma}, f) = J^{\tilde{G}}_{\tilde{M}}(\tilde{\gamma}, f) := |D^G(\gamma)|^{1/2} \int_{M_\gamma(F) \backslash G(F)} f(g^{-1} \tilde{\gamma} g) v_M(g) \dd g \]
where $f \in C^\infty_{c, \asp}(\tilde{G})$; here the Haar measures on $M_\gamma(F)$, $G(F)$ are chosen. The weighted orbital integrals for general coverings have been defined in \cite[\S 6.3]{Li14a}; the case for $\tilde{G}$ simplifies since all elements are ``good'', see Definition \ref{def:good} and Proposition \ref{prop:auto-good}.

Define $A_{M, \mathrm{reg}} \subset A_M$ be the open subset of elements such that
\[ \prod_{\beta \in \Sigma_P^{\mathrm{red}}} \left( \beta(a) - \beta(a)^{-1} \right), \quad P \in \mathcal{P}(M) \]
are all invertible. For general $\tilde{\gamma}$, in \cite{Li14a} the weighted orbital integral was defined as
\[ J_{\tilde{M}}^{\mathrm{Art}}(\tilde{\gamma}, f) := \lim_{\substack{a \in A_{M, \mathrm{reg}}(F) \\ a \to 1}} \sum_{L \in \mathcal{L}(M)} r^{L, \mathrm{Art}}_M(\gamma, a) J_{\tilde{L}}( a\tilde{\gamma}, f), \quad a\gamma: G\text{-equisingular} \]
by noting that $M_{a\gamma} = M_\gamma$, and that $\rev: \widetilde{A_M} \to A_M(F)$ splits canonically in some neighborhood of $1$, namely by comparing two exponential maps (see \eqref{eqn:two-exp}).

Assume $a \to 1$ is in general position as above. Here $r^{\mathrm{Art}}_M(\gamma, a) = r^{\mathrm{Art}}_M(\gamma, a; 0)$ is attached to the $(G, M)$-family
\begin{align*}
	r^{\mathrm{Art}}_P(\gamma, a; \lambda) & = \prod_{\alpha \in \Sigma_P} r_\alpha^{\mathrm{Art}}\left(\gamma, a; \frac{\lambda}{2} \right), \\
	r_\alpha^{\mathrm{Art}}\left(\gamma, a; \lambda\right) & = \prod_{\substack{\beta \in \Sigma^{G_\eta}_{\mathrm{red}}(A_{M_\eta}) \\ \beta_M = \alpha }} \left| \alpha(a) - \alpha(a)^{-1} \right|^{\lrangle{\lambda, \rho^{\mathrm{Art}}(\beta, u) \check{\beta}_M }}
\end{align*}
where $\check{\beta} \in \mathfrak{a}_{M_\eta}$ is the ``coroot'' of $\beta$ defined by Arthur \cite[p.229]{Ar88LB}, $\check{\beta}_M \in \mathfrak{a}_M$ is its natural projection, and $\beta_M := \beta|_{A_M}$. The crucial factor $\rho^{\mathrm{Art}}(\beta, u) \in \R$ is also defined by Arthur via the Grothendieck--Springer resolution. It depends only on $G_\eta \supset M_\eta$, $\beta$, $u$ and has nothing to do with the covering.

Our definition below is based on the corrected version $\rho(\beta, u) \in \mathfrak{a}_{M_\eta}$ proposed in \cite[II.1.4]{MW16-1}, for all $u \in M_{\eta, \mathrm{unip}}(F)$ and $\beta \in \Sigma^{G_\eta}(A_{M_\eta})$. It is designed to make \cite[Corollary 6.3]{Ar88LB} and its variant \cite[Proposition 5.4.2]{Li14a} (a descent property for weighted orbital integrals --- there are minor inaccuracies in \textit{loc.\ cit.}) true; it satisfies
\[ \rho^{\mathrm{Art}}(\beta, u) \check{\beta} = \sum_{n \geq 1} \rho(n\beta, u) \]
for all $\beta \in \Sigma^{G_\eta}_{\mathrm{red}}(A_{M_\eta})$, where $\rho(n\beta, u) := 0$ if $n\beta \notin \Sigma^{G_\eta}(A_{M_\eta})$. Next, for $\alpha \in \Sigma(A_M)$,
\begin{align*}
	\rho(\alpha, \gamma) & := \sum_{\substack{\beta \in \Sigma^{G_\eta}(A_{M_\eta}) \\ \beta_M = \alpha }} \underbracket{\rho(\beta, u)}_{\in \mathfrak{a}_{M_\eta}} \quad \text{projected to } \mathfrak{a}_M, \\
	r_\alpha\left(\gamma, a; \lambda\right) & := \left| \alpha(a) - \alpha(a)^{-1} \right|^{\lrangle{\lambda, \rho(\alpha, \gamma) }}.
\end{align*}
Then we define $r_P(\gamma, a; \lambda)$, $r_M(\gamma, a)$ in the same way. The weighted orbital integral considered in this work is defined as follows.

\begin{definition}\label{def:weighted-J-general}
	For all $\tilde{\gamma} \in \tilde{M}$ and $f$, set
	\[ J_{\tilde{M}}(\tilde{\gamma}, f) := \lim_{\substack{a \in A_{M, \mathrm{reg}}(F) \\ a \to 1}} \sum_{L \in \mathcal{L}(M)} r^L_M(\gamma, a) J_{\tilde{L}}( a\tilde{\gamma}, f), \quad a\gamma: G\text{-equisingular}. \]
\end{definition}

Here are some basic properties of $J_{\tilde{M}}(\tilde{\gamma}, f)$.
\begin{compactitem}
	\item The limit above exists.
	\item It satisfies the usual properties used in the coarse trace formula, such as the descent formulas and splitting formulas.
	\item Moreover, the estimate in \cite[II.1.5 (3)]{MW16-1} also holds in our context: there exists $r > 0$ such that
	\begin{equation}
		\left|J_{\tilde{M}}(\tilde{\gamma}, f) - \sum_{L \in \mathcal{L}(M)} r^L_M(\gamma, a) J_{\tilde{L}}( a\tilde{\gamma}, f) \right| \ll d(a)^r, \quad \text{as}\; a \in A_{M, \mathrm{reg}}(F), \; a \to 1;
	\end{equation}
	Here $d(\exp H) = \|H\|$ for a chosen norm $\|\cdot\|$ on the Lie algebra of $A_M$, when $H$ is close to zero.
\end{compactitem}

\begin{remark}\label{rem:weighted-orbital-integral-measures}
	To prescribe $\tilde{\gamma} \in \tilde{M}$ and a Haar measure on $M_\gamma(F)$ is the same as to prescribe an element of $D_{\mathrm{orb}, -}(\tilde{M}) \otimes \mes(M)^\vee$. In this way, we shall view $J_{\tilde{M}}(\tilde{\gamma}, f)$ as a function of
	\[ \tilde{\gamma} \in D_{\mathrm{geom}, -}(\tilde{M}) \otimes \mes(M)^\vee, \quad f \in C^\infty_{c, \asp}(\tilde{G}) \otimes \mes(G) \]
	since $D_{\mathrm{geom}, -} = D_{\mathrm{orb}, -}$ for non-Archimedean $F$.
\end{remark}

\subsection{Invariant version}
Let $M$ be a Levi subgroup of $G$. Since $J^{\tilde{G}}_{\tilde{M}}(\tilde{\gamma}, \cdot)$ is clearly concentrated at $H_{\tilde{G}}(\tilde{\gamma}) \in \mathfrak{a}_G$ when $\Supp(\tilde{\gamma})$ is a single conjugacy class, it extends to $C^\infty_{\mathrm{ac}, \asp}(\tilde{G}) \otimes \mes(G)$.

Now let $L$ be a Levi subgroup of $G$. We have the linear map from \cite[Théorème 4.13]{Li14b}
\[ \phi_{\tilde{L}}: \orbI_{\mathrm{ac}, \asp}(\tilde{G}) \otimes \mes(G) \to \orbI_{\mathrm{ac}, \asp}(\tilde{L}) \otimes \mes(L), \]
defined through tempered weighted characters. For a summary of these weighted characters, we refer to \S\ref{sec:weighted-characters-nonarch}.

\begin{definition}
	\index{IMgamma-weighted@$I_{\tilde{M}}(\tilde{\gamma}, f)$}
	The invariant version of weighted orbital integrals is defined recursively by
	\[ I_{\tilde{M}}(\tilde{\gamma}, f) = I^{\tilde{G}}_{\tilde{M}}(\tilde{\gamma}, f) := J_{\tilde{M}}(\tilde{\gamma}, f) - \sum_{L \in \mathcal{L}(M)} I^{\tilde{L}}_{\tilde{M}}(\tilde{\gamma}, \phi_{\tilde{L}}(f)). \]
	Following Remark \ref{rem:weighted-orbital-integral-measures}, here we take
	\begin{compactitem}
		\item $\tilde{\gamma} \in D_{\mathrm{geom}, -}(\tilde{M}) \otimes \mes(M)^\vee$,
		\item $f \in \orbI_{\mathrm{ac}, \asp}(\tilde{G}) \otimes \mes(G)$,
	\end{compactitem}
	to get rid of any choice of Haar measures. By induction, $I_{\tilde{M}}(\tilde{\gamma}, \cdot)$ is concentrated at $H_{\tilde{G}}(\tilde{\gamma}) \in \mathfrak{a}_G$ when $\Supp(\tilde{\gamma})$ is a single conjugacy class.
\end{definition}

As in \cite[II.1.7]{MW16-1},
\begin{gather*}
	I_{\tilde{M}}(\tilde{\gamma}, f) = \lim_{\substack{a \in A_{M, \mathrm{reg}}(F) \\ a \to 1}} \sum_{L \in \mathcal{L}(M)} r^L_M(\gamma, a) I_{\tilde{L}}( a\tilde{\gamma}, f), \quad a\gamma: G\text{-equisingular}, \\
	\left| I_{\tilde{M}}(\tilde{\gamma}, f) - \sum_{L \in \mathcal{L}(M)} r^L_M(\gamma, a) I_{\tilde{L}}( a\tilde{\gamma}, f) \right| \ll d(a)^r, \quad \text{as}\; a \in A_{M, \mathrm{reg}}(F), \; a \to 1.
\end{gather*}
The proofs are formal: see \textit{loc.\ cit.}

\begin{remark}
	In fact, $I_{\tilde{M}}(\tilde{\gamma}, f)$ is independent on the choice of $K$. To see this, it suffices to copy the arguments from \cite[Lemma 3.4]{Ar98}.
\end{remark}

Now we can state the descent formula for invariant weighted orbital integrals for $\tilde{G}$.

\begin{proposition}[Cf.\ {\cite[II.1.7 Lemme]{MW16-1}}]\label{prop:orbint-weighted-descent-nonArch}
	Let $L \in \mathcal{L}(M)$, $f \in \orbI_{\asp}(\tilde{G}) \otimes \mes(G)$ and $\tilde{\gamma} \in D_{\mathrm{geom}, -}(\tilde{M}) \otimes \mes(M)^\vee$. Write $\tilde{\gamma}^{\tilde{L}} := \Ind^{\tilde{L}}_{\tilde{M}}(\tilde{\gamma})$; see \eqref{eqn:parabolic-ind-dist-geom}. Then
	\[ I_{\tilde{L}}\left( \tilde{\gamma}^{\tilde{L}}, f \right) = \sum_{L^\dagger \in \mathcal{L}(M)} d^G_M(L, L^\dagger) I^{\tilde{L}^\dagger}_{\tilde{M}}\left( \tilde{\gamma}, f_{\tilde{L}^\dagger} \right). \]
\end{proposition}
\begin{proof}
	Identical as the proof in \textit{loc.\ cit.}, which is based on harmonic analysis for $\tilde{G}$ and general properties of the constants $d^G_M(L, L^\dagger)$.
\end{proof}

\section{Weighted orbital integrals: Archimedean case}\label{sec:orbint-weighted-Arch}
Consider an Archimedean local field $F$. The same recipe yields
\[ J_{\tilde{M}}(\tilde{\gamma}, f) \quad \text{for all}\quad \tilde{\gamma} \in D_{\mathrm{orb}, -}(\tilde{M}) \otimes \mes(M)^\vee, \quad f \in C^\infty_{\mathrm{ac}, \asp}(\tilde{G}) \otimes \mes(G). \]
Moreover, it is continuous in $f$. The construction works when $\tilde{G}$ is replaced by any group of metaplectic type. As before, for these $\tilde{\gamma}$ we can define
\[ I_{\tilde{M}}(\tilde{\gamma}, f) = I^{\tilde{G}}_{\tilde{M}}(\tilde{\gamma}, f) := J_{\tilde{M}}(\tilde{\gamma}, f) - \sum_{L \in \mathcal{L}(M)} I^{\tilde{L}}_{\tilde{M}}(\tilde{\gamma}, \phi_{\tilde{L}}(f)). \]
The $\phi_{\tilde{L}}$ here is a continuous linear map defined through tempered weighted characters (see \S\ref{sec:weighted-characters-arch}). One can proceed in one of the following ways.
\begin{enumerate}
	\item One can choose a maximal compact subgroup $K \subset G(F)$ in good position relative to $M_0$, take $f \in \orbI_{\asp}(\tilde{G}, \tilde{K}) \otimes \mes(G)$ (see Definition \ref{def:K-finite}), and use
	\[ \phi_{\tilde{L}}: C^\infty_{\mathrm{ac}, \asp}(\tilde{G}, \tilde{K}) \otimes \mes(G) \to \orbI_{\mathrm{ac}, \asp}(\tilde{L}, \widetilde{K^L}) \otimes \mes(L) \]
	for all $L \in \mathcal{L}(M_0)$ with $K^L := K \cap L(F)$. This is the map used in \cite[\S 4.3]{Li14b}, based on the $\tilde{K} \times \tilde{K}$-finite trace Paley--Wiener theorem (Remark \ref{rem:real-PW}) for groups of metaplectic type.

	Since $\orbI_{\mathrm{ac}, \asp}(\tilde{G}, \tilde{K})$, etc.\ are independent of $\tilde{K}$, so is $I_{\tilde{M}}(\tilde{\gamma}, f)$.
	
	\item One can also consider the broader setting of $f \in C^\infty_{\mathrm{ac}, \asp}(\tilde{G}) \otimes \mes(G)$ and the corresponding
	\[ \phi_{\tilde{L}}: C^\infty_{\mathrm{ac}, \asp}(\tilde{G}) \otimes \mes(G) \to \orbI_{\mathrm{ac}, \asp}(\tilde{L}) \otimes \mes(L). \]
	This is based on the $C^\infty_c$ trace Paley--Wiener theorem \cite{Bo94b}, valid for general coverings; see \cite[IV.1]{MW16-1}.
\end{enumerate}

Note that $I_{\tilde{M}}(\tilde{\gamma}, \cdot)$ is concentrated at $H_{\tilde{G}}(\tilde{\gamma}) \in \mathfrak{a}_G$ when $\Supp(\tilde{\gamma})$ is a single conjugacy class.

So far we defined $I_{\tilde{M}}(\tilde{\gamma}, \cdot)$ only for $\tilde{\gamma} \in D_{\mathrm{orb}, -}(\tilde{M}) \otimes \mes(M)^\vee$. Following the recipe in \cite[V.1.3]{MW16-1}, we propose a partial extension to some other $\tilde{\gamma} \in D_{\mathrm{geom}, -}(\tilde{M}) \otimes \mes(M)^\vee$. It is based on the notion of equisingularity (Definition \ref{def:equisingular}).

\begin{lemma}\label{prop:weighted-arch-aux}
	Fix a $G$-equisingular semisimple class $\mathcal{O} \subset M(F)$, $\tilde{\gamma} \in \rev^{-1}(\mathcal{O})$ and $f \in \orbI_{\asp}(\tilde{G}) \otimes \mes(G)$. Then there exist $f_1 \in \orbI_{\asp}(\tilde{M}) \otimes \mes(M)$ and an invariant neighborhood $\mathcal{U}$ of $\mathcal{O}$ such that $I^{\tilde{M}}(\tilde{\gamma}', f_1) = I^{\tilde{G}}_{\tilde{M}}(\tilde{\gamma}', f)$ for all $\tilde{\gamma}' \in D_{\mathrm{orb},-}(\tilde{M}) \otimes \mes(M)^\vee$ supported in $\rev^{-1}(\mathcal{U})$.
\end{lemma}
\begin{proof}
	Assume inductively that the property holds when $G$ is replaced by any $L \in \mathcal{L}(M)$ such that $L \neq G$; note that the case $L=M$ is trivial.

	Choose Haar measures to trivialize the lines $\mes(\cdots)$. By \cite[Proposition 6.4.1]{Li14a}, there exist $f_0 \in C^\infty_{c, \asp}(\tilde{G})$ and an invariant neighborhood $\mathcal{V}$ of $\mathcal{O}$ such that $I^{\tilde{M}}(\tilde{\gamma}', f_0) = J^{\tilde{G}}_{\tilde{M}}(\tilde{\gamma}', f)$ for all $\tilde{\gamma}' \in \Gamma_{\mathrm{orb}, -}(\tilde{M})$ supported in $\rev^{-1}(\mathcal{V})$. The assertion then follows from the definition of $I^{\tilde{G}}_{\tilde{M}}(\tilde{\gamma}', f)$ and induction hypotheses.
\end{proof}

\begin{definition}[Extension to $G$-equisingular classes]\label{def:weighted-I-integral-arch}
	Let $\mathcal{O} \subset M(F)$ be a $G$-equisingular semisimple class. For every $\tilde{\gamma} \in D_{\mathrm{geom}, -}(\tilde{M}, \mathcal{O}) \otimes \mes(M)^\vee$ and $f \in \orbI_{\asp}(\tilde{G}) \otimes \mes(G)$, set
	\[ I^{\tilde{G}}_{\tilde{M}}(\tilde{\gamma}, f) := I^{\tilde{M}}(\tilde{\gamma}, f_1) \]
	where $f_1 \in \orbI_{\asp}(\tilde{M}) \otimes \mes(M)$ is given as in Lemma \ref{prop:weighted-arch-aux}.
\end{definition}

Let us show that this is well-defined. For simplicity, assume $\Supp(\tilde{\gamma})$ is a single conjugacy class. Definition \ref{def:weighted-I-integral-arch} is independent of the choice of $f_1$ since $I^{\tilde{M}}(\tilde{\gamma}, f_1)$ can be expressed as a linear combination of derivatives at $X=0$ of
\[ X \mapsto I^{\tilde{M}}(\exp(X)\tilde{\gamma}_{\mathrm{orb}} , f_1) = I^{\tilde{G}}_{\tilde{M}}( \exp(X)\tilde{\gamma}_{\mathrm{orb}}, f ), \quad X \in \mathfrak{t}(F) \]
where
\begin{compactitem}
	\item $\tilde{\gamma}_{\mathrm{orb}} \in D_{\mathrm{orb}, -}(\tilde{M}) \otimes \mes(M)^\vee$ is the ``orbital part'' of $\tilde{\gamma}$;
	\item $T$ is any maximal $F$-torus in $M$ with $\gamma \in T(F)$;
	\item we assume that $\exp(X)\tilde{\gamma}_{\mathrm{orb}} \in \tilde{G}_{\mathrm{reg}}$.
\end{compactitem}

It is clear that $I^{\tilde{G}}_{\tilde{M}}(\tilde{\gamma}, f)$ gives back the original definition when $\tilde{\gamma}$ is already ``orbital''.

The description above implies that $I^{\tilde{G}}_{\tilde{M}}(\tilde{\gamma}, \cdot)$ is concentrated at $H_{\tilde{G}}(\tilde{\gamma})$. Therefore $I_{\tilde{M}}(\tilde{\gamma}, f)$ is defined for all $f \in \orbI_{\mathrm{ac}, \asp}(\tilde{G}) \otimes \mes(G)$ and $\tilde{\gamma} \in D_{\mathrm{geom}, -}(\mathcal{O}) \otimes \mes(M)^\vee$.

\begin{definition}
	\index{Dgeom-equi@$D_{\mathrm{geom}, G\text{-equi}, -}(\tilde{M})$}
	Let $D_{\mathrm{geom}, G\text{-equi}, -}(\tilde{M}) := \bigoplus_{\mathcal{O}} D_{\mathrm{geom}, -}(\tilde{M}, \mathcal{O})$ where $\mathcal{O}$ ranges over the $G$-equisingular semisimple classes in $M(F)$ (Definition \ref{def:equisingular}).
\end{definition}

Varying $\mathcal{O}$, we arrive at the definition of
\[ I^{\tilde{G}}_{\tilde{M}}(\tilde{\gamma}, f), \quad \tilde{\gamma} \in D_{\mathrm{geom}, G\text{-equi}, -}(\tilde{M}) \otimes \mes(M)^\vee, \quad f \in \orbI_{\mathrm{ac}, \asp}(\tilde{G}) \otimes \mes(G). \]

In the group-theoretic result below, $F$ can be any field.

\begin{lemma}\label{prop:d-equisingular}
	Let $L \in \mathcal{L}(M)$ and $\mathcal{O}$ be a $G$-equisingular semisimple conjugacy class in $M(F)$. Denote by $\mathcal{O}^L$ the conjugacy class in $L(F)$ containing $\mathcal{O}$. Assume that $\mathcal{O}^L$ is $G$-equisingular. For all $L^\dagger \in \mathcal{L}(M)$, we have
	\[ d^G_M(L, L^\dagger) \neq 0 \implies \mathcal{O} \;\text{is $L^\dagger$-equisingular}. \]
\end{lemma}
\begin{proof}
	This is part (i) of \cite[V.1.3 Lemme]{MW16-1}.
\end{proof}

Now comes the Archimedean counterpart of Proposition \ref{prop:orbint-weighted-descent-nonArch}.

\begin{proposition}[Cf.\ part (ii) of {\cite[V.1.3 Lemme]{MW16-1}}]\label{prop:orbint-weighted-descent-Arch}
	Let $L \in \mathcal{L}(M)$, $f \in \orbI_{\asp}(\tilde{G}) \otimes \mes(G)$ and $\tilde{\gamma} \in D_{\mathrm{geom}, -}(\tilde{M}, \mathcal{O}) \otimes \mes(M)^\vee$, where $\mathcal{O}$ is a semisimple conjugacy class in $M(F)$ such that $\mathcal{O}^L$ is $G$-equisingular. Write
	\[ \tilde{\gamma}^{\tilde{L}} := \Ind^{\tilde{L}}_{\tilde{M}}(\tilde{\gamma}) \in D_{\mathrm{geom}, -}(\tilde{L}, \mathcal{O}^L) \otimes \mes(L)^\vee \]
	(see \eqref{eqn:parabolic-ind-dist}). Then
	\[ I_{\tilde{L}}\left( \tilde{\gamma}^{\tilde{L}}, f \right) = \sum_{L^\dagger \in \mathcal{L}(M)} d^G_M(L, L^\dagger) I^{\tilde{L}^\dagger}_{\tilde{M}}\left( \tilde{\gamma}, f_{\tilde{L}^\dagger} \right). \]
\end{proposition}
\begin{proof}
	Note that $d^G_M(L, L^\dagger) I^{\tilde{L}^\dagger}_{\tilde{M}}\left( \tilde{\gamma}, f_{\tilde{L}} \right)$ makes sense by Lemma \ref{prop:d-equisingular}. The remaining arguments are identical to \textit{loc.\ cit.}.
\end{proof}

\begin{remark}\label{rem:weighted-general-covering}
	The constructions in the non-Archimedean case \S\ref{sec:orbint-weighted-nonArch} work for all covering group $\tilde{G}$, its Levi subgroup $\tilde{M}$ and all good conjugacy classes $\tilde{\gamma} \in \Gamma(\tilde{M})$. See for example \cite[\S 6.3]{Li14a}.
	
	The same also holds in the Archimedean case. One can avoid assuming the $\tilde{K} \times \tilde{K}$-finite trace Paley--Wiener theorem by using the $C^\infty_c$-version in \cite[\S 7]{Bo94b}, which works for general coverings. 
\end{remark}

\section{Systems of \texorpdfstring{$B$}{B}-functions and the local geometric theorem}\label{sec:local-geometric}
Rescaling of root lengths will be needed after descent around semisimple elements, and the $B$-functions are bookkeeping devices for this purpose.

To begin with, consider a connected reductive $F$-group $L$. Its $\overline{F}$-pinnings $(B_L, T, (E_\beta)_\beta)$ form a torsor under $L_{\mathrm{AD}}(\overline{F})$-action. One defines universal versions of the root data and the Weyl group as compatible systems indexed by pinnings, with transition maps given by $L_{\mathrm{AD}}$. Given $(B_L, T, (E_\beta)_\beta)$, we obtain a Galois action on the universal $\Sigma(L, T)$ and $\check{\Sigma}(L, T)$, denoted hereafter as $\sigma \mapsto \sigma_L^*$ where $\sigma \in \Gamma_F$; it factors through $\Gal(E|F)$ where $E|F$ is any Galois extension splitting $G$. This is also the viewpoint adopted in \cite[I.1.2]{MW16-1}.

These ``universal'' objects are insensitive to inner twists. They can be manipulated concretely by choosing an $\overline{F}$-pinning $(B_L, T, (E_\beta)_\beta)$. If $(B_L, T, (E_\beta)_\beta)$ is an $F$-pinning, then $\sigma \mapsto \sigma_{L^*}$ is induced by the $F$-structure of $T$.

Let us fix a positive definite quadratic form on the $X_*(T) \otimes \R$ that is invariant under $\Gamma_F$ and the Weyl group. Here $X_*(T)$, the Weyl group and $\Gamma_F$-actions are understood in the universal sense above. By duality, we deduce a similar form on $X^*(T) \otimes \R$, denoted as $(\cdot | \cdot)$. The definitions below are from \cite[II.1.8, II.1.9]{MW16-1}.

\begin{definition}\label{def:B-function}
	\index{system of $B$-functions}
	Consider the universal set of absolute roots $\Sigma(L, T)$ in $L$ defined above.
	\begin{enumerate}
		\item A $B$-function on $\Sigma(L, T)$ is a map $\Sigma(L, T) \to \Q_{> 0}$ satisfying
		\begin{itemize}
			\item $B$ is invariant under $\Gamma_F$ and the Weyl group,
			\item $B(-\beta) = B(\beta)$ for all $\beta \in \Sigma(L, T)$,
			\item over any irreducible factor of the root system associated attached to $\Sigma(L, T)$, either $B$ or $\beta \mapsto (\beta|\beta)^{-1} B(\beta)$ is constant.
		\end{itemize}
		\item Given a $B$-function as above, set
		\begin{align*}
			\Sigma(L, T; B) & := \left\{ B(\beta)^{-1} \beta : \beta \in \Sigma(L, T) \right\} \subset X^*(T) \otimes \Q , \\
			\check{\Sigma}(L, T; B) & := \left\{ B(\beta) \check{\beta}: \beta \in \Sigma(L, T) \right\} \subset X_*(T) \otimes \Q.
		\end{align*}
	\end{enumerate}
\end{definition}

According to \cite[pp.204--205]{MW16-1}, $\Sigma(L, T; B)$ form a root system such that $\check{\Sigma}(L, T; B)$ is the set of coroots, and the Weyl group is the same as the original one.

For any Levi subgroup $M \subset L$, denote by $\Sigma(A_M; B)$ the set of nonzero restrictions to $\mathfrak{a}_M$ of the elements of $\Sigma(L, T; B)$. Note that it still depends on $L$.

Now let $G^!$ be a quasisplit connected reductive $F$-group, and take any $F$-pinning
\[ (B^!, T^!, (E_\alpha)_\alpha) \quad \text{of}\; G^!. \]
Fix a positive definite quadratic form $(\cdot|\cdot)$ on $X_*(T^!) \otimes \R$ with the same invariance properties as before. For every $\epsilon \in G^!(F)_{\mathrm{ss}}$, there exists $g \in G^!_{\mathrm{AD}}(\overline{F})$ such that $\epsilon \in gT^!(\overline{F})g^{-1}$; in this manner, we obtain a similar quadratic form for $L := G^!_\epsilon$ by forming the $\overline{F}$-pinning
\[ \left( gB^! g^{-1} \cap L, \; g T^! g^{-1}, \; (g E_\alpha g^{-1})_{\alpha: \alpha(g^{-1}\epsilon g)=1} \right) \]
of $L$. One readily checks that the form is independent of all choices.

\begin{definition}\label{def:system-B-function}
	Let $G^!$ be a quasisplit connected reductive $F$-group, with a chosen quadratic form $(\cdot|\cdot)$ as above. A \emph{system of $B$-functions} is an assignment $B: \epsilon \mapsto B_\epsilon$ for each $\epsilon \in G^!(F)_{\mathrm{ss}}$ such that
	\begin{itemize}
		\item $B_\epsilon$ is a $B$-function for $G^!_\epsilon$;
		\item if $\epsilon, \epsilon' \in G^!(F)_{\mathrm{ss}}$ and $g \in G^!(\overline{F})$ satisfy $g \epsilon g^{-1} = \epsilon'$, then the inner twist $\Ad_g: (G^!_\epsilon)_{\overline{F}} \to (G^!_{\epsilon'})_{\overline{F}}$ identifies $B_\epsilon$ with $B_{\epsilon'}$.
	\end{itemize}
	The assignment $B_\epsilon = 1$ for all $\epsilon \in G^!(F)_{\mathrm{ss}}$ is called the \emph{trivial system of $B$-functions}.
\end{definition}

The system of $B$-functions prescribe a coherent way to rescale the roots. Given a system of $B$-functions for $G^!$, a Levi subgroup $M^!$ and $\gamma = \epsilon u \in M^!(F)$ (Jordan decomposition), there are variants
\[ \rho(\beta, u; B), \quad \rho(\alpha, \gamma; B), \quad r_P(u, a; B, \lambda) \]
of the objects in \S\ref{sec:orbint-weighted-nonArch}, where $\beta \in \Sigma(A_{M^!_\epsilon}; B_\epsilon)$ and $\alpha \in \Sigma(A_{M^!}; B_\epsilon)$. Here $\Sigma(A_{M^!}; B_\epsilon)$ is the set of nonzero restrictions of elements of $\Sigma(A_{M^!_\epsilon}; B_\epsilon)$ to $\mathfrak{a}_{M^!}$. See \cite[pp.206--209]{MW16-1} for details. Take $\rho(\beta, u; B)$ for example, in \cite[II.1.8 (6)]{MW16-1} it is inductively defined in terms of various $\rho(\beta', u)$.

Since $\Sigma(A_{M^!}; B_\epsilon)$ depends only on the semisimple stable conjugacy class $\mathcal{O}^!$ in $M(F)$ generated by $\epsilon$, it may be denoted as $\Sigma(A_{M^!}, B_{\mathcal{O}^!})$.

\begin{remark}\label{rem:B-restriction}
	Given a system of $B$-functions for $G^!$ and a Levi subgroup $L^! \subset G^!$, we obtain a system $B|_{L^!}$ since $L^!_\epsilon \subset G^!_\epsilon$ is also a Levi subgroup for every $\epsilon \in L^!(F)_{\mathrm{ss}}$, and one can restrict $B_\epsilon$ to the roots which are trivial on $Z_{L^!}^\circ$.
\end{remark}

Using the factors $r^{L^!}_{M^!}(u, a; B)$ for various $L^! \in \mathcal{L}(M^!)$, we are then able to define
\[ J^{G^!}_{M^!}(\gamma, B, f^!), \quad I^{G^!}_{M^!}(\gamma, B, f^!) \]
as before, with similar properties.

With the given system of $B$-functions, the \emph{stable weighted orbital integrals} are defined in \cite[II.1.10, V.1.4]{MW16-1}:
\[ S^{G^!}_{M^!}(\delta, B, f^!), \quad f^! \in S\orbI(G^!) \otimes \mes(G^!), \]
\index{SGM-delta-weighted@$S^{G^{"!}}_{M^{"!}}(\delta, B, f^{"!})$}
and $\delta$ ranges over $SD_{\mathrm{geom}}(M^!) \otimes \mes(M^!)^\vee$ for non-Archimedean $F$ (resp.\ $SD_{\mathrm{geom}, G^!\text{-equi}}(M^!) \otimes \mes(M^!)^\vee$ for Archimedean $F$, as in \S\ref{sec:orbint-weighted-Arch}). The distribution is concentrated at $H_{G^!}(\delta)$ when $\Supp(\delta)$ is a single stable class. It turns out that $S^{G^!}_{M^!}(\delta, B, f^!)$ is continuous in $f^!$ when $F$ is Archimedean. The stability of these distributions is the content of \cite[Theorem II.1.10]{MW16-1}.

When $B$ is trivial, all these objects reduce to the usual versions without reference $B$.

\begin{remark}
	When $\gamma$ is $G$-equisingular, we have $S^{G^!}_{M^!}(\delta, B, f^!) = S^{G^!}_{M^!}(\delta, f^!)$ by \cite[II.1.10 (23)]{MW16-1}. This is the case when $\gamma \in M^!_{G^! \text{-reg}}(F)$, for example.
\end{remark}

Return to the case of a group of metaplectic type $\tilde{G} = \prod_{i \in I} \GL(n_i) \times \Mp(W)$. The hardcore is of course the case $\tilde{G} = \Mp(W)$. Let $\mathbf{G}^! \in \Endo_{\elli}(\tilde{G})$. We are going to define a system of $B$-functions for $G^!$, denoted as $B^{\tilde{G}}$. Recall that invariant quadratic forms for $G^!$ have been chosen in Definition \ref{def:invariant-quadratic-form}. By Definition \ref{def:system-B-function} and \cite[Lemma 3.3]{Ko82}, it suffices to prescribe $B^{\tilde{G}}_\epsilon$ for each $\epsilon \in G^!(F)_{\mathrm{ss}}$ such that $G^!_\epsilon$ is quasisplit.

\begin{definition}\label{def:B-metaplectic}
	\index{BtildeG@$B^{\tilde{G}}$}
	Suppose $\tilde{G} = \Mp(W)$ and let $\mathbf{G}^! \in \Endo_{\elli}(\tilde{G})$ be associated with the pair $(n', n'') \in \Z_{\geq 0}^2$. The system $B^{\tilde{G}}$ for $G^! = \SO(2n'+1) \times \SO(2n''+1)$ is defined as follows. Let $\epsilon = (\epsilon', \epsilon'') \in G^!(F)_{\mathrm{ss}}$. According to \S\ref{sec:ss-classes}, there are canonical decompositions
	\[ \SO(2n'+1)_{\epsilon'} = U' \times \SO(V'_+) \times \SO(V'_-), \quad \SO(2n''+1)_{\epsilon''} = U'' \times \SO(V''_+) \times \SO(V''_-), \]
	where $U'$ is a direct product of unitary groups and general linear groups up to restriction of scalars, $V'_+$ (resp.\ $V'_-$) is a quadratic $F$-vector space of odd (resp.\ even) dimension; ditto for $U''$, $V''_+$ and $V''_-$. Since $G^!_\epsilon$ is quasisplit, $\SO(V'_+)$ and $\SO(V''_+)$ are split.
	
	Let $\beta$ be an (absolute) root of $G^!_\epsilon$. If $\beta$ comes from the factor $\SO(V'_+)$ or $\SO(V''_+)$, set
	\[ B^{\tilde{G}}_\epsilon(\beta) := \begin{cases}
		1, & \beta: \;\text{long}, \\
		\frac{1}{2}, & \beta: \;\text{short},
	\end{cases}\]
	and in the other cases we set $B^{\tilde{G}}(\beta) = 1$; note that when $\SO(V'_+)$ or $\SO(V''_+)$ has rank $1$, every root is short by convention.
	
	The recipe extends to groups of metaplectic type and their elliptic endoscopic data: simply set $B^{\tilde{G}}$ to be trivial on the $\GL$-factors. Observe that $B^{\tilde{G}}_{z\epsilon} = B^{\tilde{G}}_\epsilon$ for all $z \in Z_{G^!}(F)$.
\end{definition}

It is straightforward to verify the conditions in Definition \ref{def:system-B-function} for $B^{\tilde{G}}$.

\begin{lemma}\label{prop:B-Levi}
	Suppose $M$ is a Levi subgroup of $G$ and $\mathbf{M}^! \in \Endo_{\elli}(\tilde{M})$. For every $s \in \Endo_{\mathbf{M}^!}(\tilde{G})$, consider the system of $B$-functions $B^{\tilde{G}}$ on $G^![s]$. Then $B^{\tilde{G}}|_{M^!} = B^{\tilde{M}}$ (see Remark \ref{rem:B-restriction}).
\end{lemma}
\begin{proof}
	It is clear that $B^{\tilde{G}}_\epsilon$ restricted to the roots in $M^!_\epsilon(F) = M^!_{\epsilon[s]}(F)$ equals $B^{\tilde{M}}_\epsilon = B^{\tilde{M}}_{\epsilon[s]}$: all are described in terms of root-lengths in $M^!_\epsilon$.
\end{proof}

\begin{definition}\label{def:equisingular-endo}
	\index{equisingular}
	Suppose $M$ is a Levi subgroup of $G$ and $\mathbf{M}^! \in \Endo_{\elli}(\tilde{M})$. Let $\mathcal{O}^!$ be a stable semisimple conjugacy class in $M^!(F)$.
	\begin{itemize}
		\item We say that $\mathcal{O}^!$ is $\tilde{G}$-equisingular if $\gamma$ is $G$-equisingular (Definition \ref{def:equisingular}) for all element $\gamma \in M(F)_{\mathrm{ss}}$ corresponding to $\mathcal{O}^!$ under $\Sigma_{\mathrm{ss}}(M^!) \xrightarrow{\text{endoscopy}} \Sigma_{\mathrm{ss}}(M)$.
		\item Denote by $SD_{\mathrm{geom}, \tilde{G}\text{-equi}}(M^!)$ the subspace $\bigoplus_{\mathcal{O}^!} SD_{\mathrm{geom}}(M^!, \mathcal{O}^!)$ of $SD_{\mathrm{geom}}(M^!)$, where $\mathcal{O}^!$ ranges over the $\tilde{G}$-equisingular classes.
	\end{itemize}
\end{definition}

The next result involves the Definitions \ref{def:Endo-s} and \ref{def:central-twist}.

\begin{lemma}[Cf.\ {\cite[p.483]{MW16-1}}]\label{prop:equisingular-endo}
	Let $M \subset G$, $\mathbf{M}^! \in \Endo_{\elli}(\tilde{M})$ and $\mathcal{O}^! \subset M^!(F)_{\mathrm{ss}}$ be as above. For all $\delta \in \mathcal{O}^!$ and $s \in \Endo_{\mathbf{M}^!}(\tilde{G})$, if $\delta$ is $\tilde{G}$-equisingular then $\delta[s]$ is $G^![s]$-equisingular.
\end{lemma}
\begin{proof}
	There exist $\gamma \in M(F)_{\mathrm{ss}}$ and a diagram $(\delta, B^!, T^!, B, T, \gamma)$, by Corollary \ref{prop:diagram-from-epsilon}. Since $A_M \subset M_\gamma$, we have $\gamma$ is $G$-equisingular if and only if $A_M$ is central in $G_\gamma$, or equivalently that all $\alpha \in \Sigma(G_\gamma, T)$ are trivial on $A_M$. There is a similar criterion for the $G^![s]$-equisingularity of $\delta[s]$.
	
	By Proposition \ref{prop:central-twist-corr}, $\delta[s]$ corresponds to $\gamma$ via $\mathbf{G}^![s] \in \Endo_{\elli}(\tilde{G})$. The descent of endoscopy in \S\ref{sec:descent-endoscopy} implies that up to rescaling by a factor of two, $\Sigma\left( G^![s]_{\delta[s]}, T^! \right)$ embeds into $\Sigma\left(G_\gamma, T \right)$, compatibly with the isomorphisms $T^! \simeq T$ and $A_{M^!} \simeq A_M$ from the chosen diagram. As $A_M$ is connected, the scaling is immaterial, so all $\beta \in \Sigma\left( G^![s]_{\delta[s]}, T^! \right)$ are trivial on $A_{M^!}$.
\end{proof}

\begin{definition}\label{def:IEndo-1}
	\index{IMEndo-delta-weighted@$I^{\Endo}_{\tilde{M}}(\mathbf{M}^{"!}, \delta, f)$, $I^{\Endo}_{\tilde{M}}(\tilde{\gamma}, f)$}
	Let $M \subset G$ be a Levi subgroup and $\mathbf{M}^! \in \Endo_{\elli}(\tilde{M})$. Given $f \in \orbI_{\asp}(\tilde{G}) \otimes \mes(G)$ and
	\[ \delta \in \begin{cases}
		SD_{\mathrm{geom}}(M^!) \otimes \mes(M^!), & \text{when $F$ is non-Archimedean}, \\
		SD_{\mathrm{geom}, \tilde{G}\text{-equi}}(M^!) \otimes \mes(M^!), & \text{when $F$ is Archimedean.}
	\end{cases}\]
	We define
	\begin{align*}
		I^{\Endo}_{\tilde{M}}(\mathbf{M}^!, \delta, f) & = I^{\tilde{G}, \Endo}_{\tilde{M}}(\mathbf{M}^!, \delta, f) \\
		& := \sum_{s \in \Endo_{\mathbf{M}^!}(\tilde{G})} i_{M^!}(\tilde{G}, G^![s]) S^{G^![s]}_{M^!}\left( \delta[s], B^{\tilde{G}}, f^{G^![s]} \right)
	\end{align*}
	where $f^{G^![s]} := \Trans_{\mathbf{G}^![s], \tilde{G}}(f) \in S\orbI(G^![s]) \otimes \mes(G^![s])$.
\end{definition}

Note that $S^{G^![s]}_{M^!}\left( \delta[s], B^{\tilde{G}}, f^{G^![s]} \right)$ is well-defined in the Archimedean case by Lemma \ref{prop:equisingular-endo}.

\begin{definition-proposition}\label{def:I-geom-Endo}
	Let $M$ be a Levi subgroup of $G$. Suppose that
	\[ \delta_{\mathbf{M}^!} \in SD_{\mathrm{geom}}(M^!) \otimes \mes(M^!)^\vee, \quad \text{(resp.\ } SD_{\mathrm{geom}, \tilde{G}\text{-equi}}(M^!) \otimes \mes(M^!)^\vee \text{ )} \]
	is given for each $\mathbf{M}^! \in \Endo_{\elli}(\tilde{M})$ when $F$ is non-Archimedean (resp.\ Archimedean), and set
	\[ \tilde{\gamma} := \sum_{\mathbf{M}^!} \trans_{\mathbf{M}^!, \tilde{M}}\left( \delta_{\mathbf{M}^!} \right) \; \in
	\begin{cases}
		D_{\mathrm{geom}, -}(\tilde{M}) \otimes \mes(M)^\vee, & F:\; \text{non-Archimedean} \\
		D_{\mathrm{geom}, G\text{-equi}, -}(\tilde{M}) \otimes \mes(M)^\vee, & F:\; \text{Archimedean}.
	\end{cases}\]
	Then $I^{\Endo}_{\tilde{M}}(\tilde{\gamma}, f) := \sum_{\mathbf{M}^! \in \Endo_{\elli}(\tilde{G})} I^{\Endo}_{\tilde{M}}\left( \mathbf{M}^!, \delta_{\mathbf{M}^!}, f \right)$ depends only on $f$ and $\tilde{\gamma}$. This defines
	\[ I^{\Endo}_{\tilde{M}}\left( \tilde{\gamma}, f \right) = I^{\tilde{G}, \Endo}_{\tilde{M}}\left( \tilde{\gamma}, f \right) \]
	for all
	\[ f \in \orbI_{\asp}(\tilde{G}) \otimes \mes(G), \quad \tilde{\gamma} \in
		\begin{cases}
			D_{\mathrm{geom}, -}(\tilde{M}) \otimes \mes(M)^\vee, & F:\; \text{non-Archimedean} \\
			D_{\mathrm{geom}, \text{equi}, -}(\tilde{M}) \otimes \mes(M)^\vee, & F:\; \text{Archimedean}.
	\end{cases}\]
	Moreover, it is continuous in $f$ for each $\tilde{\gamma}$ when $F$ is Archimedean.
\end{definition-proposition}
\begin{proof}
	The following arguments require the endoscopic descent formula in Proposition \ref{prop:descent-orbint-Endo}, to be proved later on. It is deferred to \S\ref{sec:descent-orbint-Endo} since some combinatorial setup will be needed.
	
	In view of Theorem \ref{prop:Dgeom-preservation} and the notation thereof, it suffices to show that given a Levi subgroup $R \subset M$ and $\mathbf{R}^! \in \Endo_{\elli}(\tilde{R})$, we have
	\[ I^{\Endo}_{\tilde{M}}\left( \mathbf{M}^![s], \delta[s]^{M^![s]}, f \right) = I^{\Endo}_{\tilde{M}}\left( \mathbf{M}^![t], \delta[t]^{M^![t]}, f \right) \]
	for all $\delta \in SD_{\mathrm{geom}}(R^!, \mathcal{O}^!) \otimes \mes(R)^!$ where $\mathcal{O}^!$ is a stable semisimple class in $R^!(F)$ which is $\tilde{G}$-equisingular when $F$ is Archimedean, $s, t \in \Endo_{\mathbf{R}^!}(\tilde{M})$ and all $f$. Indeed, this is immediate from the upcoming Proposition \ref{prop:descent-orbint-Endo}.
	
	When $F$ is Archimedean, the continuity in $f$ results from that of $I^{\Endo}_{\tilde{M}}(\mathbf{M}^!, \delta, f)$ for various $\delta$. In turn, this reduces to the continuity of $S^{G^![s]}_{M^!}(\delta[s], B^{\tilde{G}}, \cdot)$ for all $s \in \Endo_{\mathbf{M}^!}(\tilde{G})$.
\end{proof}

\begin{theorem}[Local geometric theorem]\label{prop:local-geometric}
	\index{local geometric theorem}
	For all Levi subgroups $M \subset G$ and all
	\[ f \in \orbI_{\asp}(\tilde{G}) \otimes \mes(G), \quad \tilde{\gamma} \in
	\begin{cases}
		D_{\mathrm{geom}, -}(\tilde{M}) \otimes \mes(M)^\vee, & F:\; \text{non-Archimedean} \\
		D_{\mathrm{geom}, G-\mathrm{equi}, -}(\tilde{M}) \otimes \mes(M)^\vee, & F:\; \text{Archimedean,}
	\end{cases}\]
	 we have $I_{\tilde{M}}(\tilde{\gamma}, f) = I^{\Endo}_{\tilde{M}}(\tilde{\gamma}, f)$.
\end{theorem}

The proof of Theorem \ref{prop:local-geometric} will be completed at the end of this work; see \S\ref{sec:end-stabilization}. As the first step, we will reduce it to the $\tilde{G}$-regular case in Proposition \ref{prop:local-geometric-regular-reduction}.

\section{An endoscopic descent formula}\label{sec:descent-orbint-Endo}
Let $M$ be a Levi subgroup of $G$, and consider the preimage $\tilde{M} \subset \tilde{G}$ of $M(F)$ as in \S\ref{sec:local-geometric}. We begin with some combinatorics.

Assuming temporarily that $\tilde{G}$ is a metaplectic group $\Mp(W)$, we may write
\[ M = \prod_{i \in I} \GL(n_i) \times \Sp(W^\flat) \]
where $I$ is a finite set, $(n_i)_i \in \Z_{\geq 1}^I$ and $W^\flat$ is a symplectic subspace of $W$. To prescribe $L \in \mathcal{L}(M)$ is the same as to prescribe
\begin{itemize}
	\item a finite set $K$ and $(l_k)_{k \in K} \in \Z_{\geq 1}^K$,
	\item a decomposition $I = I_0 \sqcup I_1$ into disjoint subsets, together with a surjection $\eta: I_0 \twoheadrightarrow K$,
	\item a symplectic subspace $W^\Diamond$ of $W$ such that $W^\Diamond \supset W^\flat$,
\end{itemize}
satisfying $\dim W = \dim W^\Diamond + 2\sum_{k \in K} l_k$; specifically, we have then
\[ M \subset L = \prod_{k \in K} \GL(l_k) \times \Sp(W^\Diamond) \subset G, \]
where $M \hookrightarrow L$ is given by
\[ \Sp(W^\flat) \hookrightarrow \Sp(W^\Diamond), \quad \GL(n_i) \hookrightarrow \begin{cases}
	\GL(l_{\eta(i)}), & i \in I_0 \\
	\Sp(W^\Diamond), & i \in I_1.
\end{cases}\]

Suppose $M \subset L \subset G$ are parameterized as above, and fix $\mathbf{M}^! \in \Endo_{\elli}(\tilde{M})$ with
\[ M^! = \prod_{i \in I} \GL(n_i) \times \SO(2a'+1) \times \SO(2a''+1), \quad a' + a'' = \frac{1}{2} \dim W^\flat. \]
Given $s^L \in \Endo_{\mathbf{M}^!}(\tilde{L})$ and $s_L \in \Endo_{\mathbf{L}^![s^L]}(\tilde{G})$, they assemble into $s \in \Endo_{\mathbf{M}^!}(\tilde{G})$ with $\mathbf{G}^![s] = \mathbf{G}^![s_L]$. Specifically, $s^L$ divides all $i \in I_1$ (i.e.\ the ``metaplectic'' indices relative to $\tilde{L}$) into two parts (ordered), whilst $s_L$ divides all $i \in I_0$ into two parts in a way that depends only on $\eta(i) \in K$; this divides $I$ into two parts, and the identification of endoscopic data for $\tilde{G}$ is then clear.

Consider now the reciprocal question: given $M \subset L \subset G$ and $\mathbf{M}^!$, when does a given $s \in \Endo_{\mathbf{M}^!}(\tilde{G})$ arise from a pair $(s^L, s_L)$? Recall that $s$ corresponds to a decomposition $I = I' \sqcup I''$. It turns out that following condition on $s$ is necessary and sufficient:
\begin{equation}\label{eqn:descent-ellipticity}
	\text{The fibers of $\eta: I_0 \twoheadrightarrow K$ are either contained in $I'$ or in $I''$.}
\end{equation}
When \eqref{eqn:descent-ellipticity} holds for $(L, s)$, we may define $s^L \in \Endo_{\mathbf{M}^!}(\tilde{L})$ corresponding to
\begin{equation}\label{eqn:descent-ellipticity-I}
	I_1 = I'_1 \sqcup I''_1, \quad I'_1 := I_1 \cap I', \quad I''_1 := I_1 \cap I'',
\end{equation}
as well as $s_L \in \Endo_{\mathbf{L}^![s^L]}(\tilde{G})$ corresponding to
\begin{equation}\label{eqn:descent-ellipticity-K}
	K = K' \sqcup K'', \quad K' := \left\{ k : \eta^{-1}(k) \subset I' \right\}, \quad K' := \left\{ k : \eta^{-1}(k) \subset I'' \right\}.
\end{equation}
Evidently, $\mathbf{G}^![s]$ (relative to $\mathbf{M}^! \in \Endo_{\elli}(\tilde{M})$) equals $\mathbf{G}^![s_L]$ (relative to $\mathbf{L}^![s^L] \in \Endo_{\elli}(\tilde{L})$). The choice of $(s^L, s_L)$ is unique.

\begin{remark}\label{rem:Gs-nonelliptic}
	It is often more flexible to work inside dual groups as in \eqref{eqn:Endo-s-dual}. Represent $\mathbf{M}^!$ by an elliptic semisimple element $s^\flat \in \tilde{M}^\vee$. Consider an element $s$ of $s^\flat Z_{\tilde{M}^\vee}^\circ / Z_{\tilde{G}^\vee}^\circ$. For every $L \in \mathcal{L}(M)$, these data can always be decomposed into $s^L \in s^\flat Z_{\tilde{M}^\vee}^\circ / Z_{\tilde{L}^\vee}^\circ$ and $s_L \in s^L Z_{\tilde{L}^\vee}^\circ / Z_{\tilde{G}^\vee}^\circ$. However, even if the endoscopic datum $\mathbf{G}^![s] = \mathbf{G}^![s_L]$ of $\tilde{G}$ is elliptic (i.e.\ if $s$ lands in $\Endo_{\mathbf{M}^!}(\tilde{G})$), it is not necessarily the case for $\mathbf{L}^![s^L]$ for $\tilde{L}$. This is why we impose conditions like \eqref{eqn:descent-ellipticity}.
\end{remark}

\begin{lemma}\label{prop:sL-Ls}
	Suppose that we are in the situation (cf.\ \eqref{eqn:s-situation}):
	\[\begin{tikzcd}
		& \tilde{G} \\
		M^! \arrow[dashed, leftrightarrow, r, "\text{ell.}", "\text{endo.}"'] & \tilde{M} \arrow[hookrightarrow, u, "\text{Levi}"'] \\
		R^! \arrow[dashed, leftrightarrow, r, "\text{ell.}", "\text{endo.}"'] \arrow[hookrightarrow, u, "\text{Levi}"] & \tilde{R} \arrow[hookrightarrow, u, "\text{Levi}"']
	\end{tikzcd}\]
	for some given $\mathbf{R}^! \in \Endo_{\elli}(\tilde{R})$ and $\mathbf{M}^! = \mathbf{M}^![t]$ with $t \in \Endo_{\mathbf{R}^!}(\tilde{M})$. Make the identifications $\mathfrak{a}_{R^!} = \mathfrak{a}_R$ and $\mathfrak{a}_{M^!} = \mathfrak{a}_M$ accordingly. There is a bijection
	\begin{align*}
		\left\{\begin{array}{r|l}
			(s, L^!) & s \in \Endo_{\mathbf{M}^!}(\tilde{G}) \\
			& L^! \in \mathcal{L}^{G^![s]}(R^!)
		\end{array}\right\} & \xrightarrow{1:1}
		\left\{\begin{array}{r|l}
			(L, s) & L \in \mathcal{L}^G(R), \; s \in \Endo_{\mathbf{M}^!}(\tilde{G}), \\
			& r = r(s, t) \text{ satisfies \eqref{eqn:descent-ellipticity}} \\
			& \text{w.r.t.}\; R \subset L \subset G
		\end{array}\right\}
	\end{align*}
	where $r(s,t) \in \Endo_{\mathbf{R}^!}(\tilde{G})$ is assembled from $t$ and $s$ as above. Specifically,
	\begin{itemize}
		\item it is characterized by $\mathfrak{a}_L = \mathfrak{a}_{L^!}$ (as subspaces of $\mathfrak{a}_R = \mathfrak{a}_{R^!}$);
		\item if $(s, L^!)$ is mapped to $(L, s)$, then $L^! = L^![r^L]$.
	\end{itemize}
\end{lemma}
\begin{proof}
	This is a metaplectic version of a similar assertion in \cite[p.236]{MW16-1} proved by working on the dual side. One can either adapt the cited proof, or use the direct arguments below.

	First, taking $r := r(s, t)$ yields embeddings
	\begin{gather*}
		\left\{\begin{array}{r|l}
			(s, L^!) & s \in \Endo_{\mathbf{M}^!}(\tilde{G}) \\
			& L^! \in \mathcal{L}^{G^![s]}(R^!)
		\end{array}\right\} \hookrightarrow
		\left\{\begin{array}{r|l}
			(s, L^!) & r \in \Endo_{\mathbf{R}^!}(\tilde{G}) \\
			& L^! \in \mathcal{L}^{G^![r]}(R^!)
		\end{array}\right\}, \\
		\left\{\begin{array}{r|l}
			(L, s) & L \in \mathcal{L}^G(R), \; s \in \Endo_{\mathbf{M}^!}(\tilde{G}), \\
			& r = r(s, t) \text{ satisfies \eqref{eqn:descent-ellipticity}} \\
			& \text{w.r.t.} \; R \subset L \subset G
		\end{array}\right\} \hookrightarrow
		\left\{\begin{array}{r|l}
			(L, r) & L \in \mathcal{L}^G(R), \; r \in \Endo_{\mathbf{R}^!}(\tilde{G}), \\
			& r = r(s, t) \text{ satisfies \eqref{eqn:descent-ellipticity}} \\
			& \text{w.r.t.} \; R \subset L \subset G
		\end{array}\right\},
	\end{gather*}
	whose images are characterized by the condition \eqref{eqn:descent-ellipticity} with respect to $R \subset M \subset G$, and the same reason gives injectivity. We are thus reduced to the case $M = R$, and $r = s$ in this case.
	
	Given $(s, L^!)$, write $G^![s] = \SO(2n'+1) \times \SO(2n''+1)$ and
	\[ L^! = \left( \prod_{k \in K'} \GL(l_k) \times \SO(2b'+1) \right) \times \left( \prod_{k \in K''} \GL(l_k) \times \SO(2b''+1) \right) \subset G^![s]. \]
	The inclusions $M^! \hookrightarrow L^! \hookrightarrow G^![s]$ of Levi subgroups yield the data
	\[ I' = I'_0 \sqcup I'_1, \; \eta': I'_0 \twoheadrightarrow K', \quad I'' = I''_0 \sqcup I''_1, \; \eta'': I''_0 \twoheadrightarrow K''. \]
	Take $I_0 := I'_0 \sqcup I''_0$, $I_1 := I'_1 \sqcup I''_1$ (thus $I = I_0 \sqcup I_1$), $K := K' \sqcup K''$ and
	\[ \eta := \eta' \sqcup \eta'': I_0 \twoheadrightarrow K. \]
	Moreover, take a symplectic subspace $W^\Diamond$ of $W$ containing $W^\flat$ with $\frac{1}{2} \dim W^\Diamond = b' + b''$. Then $(K, (l_k)_k, I_0, I_1, \eta)$ determines $L \in \mathcal{L}(M)$. By construction, $\mathfrak{a}_L = \mathfrak{a}_{L^!}$ and \eqref{eqn:descent-ellipticity} is satisfied by $(L, s)$. It is also clear that $L^! = L^![s^L]$ because \eqref{eqn:descent-ellipticity-I} holds automatically.
	
	Reciprocally, given $(L, s)$ satisfying \eqref{eqn:descent-ellipticity} with respect to $M \subset L \subset G$, we construct $s^L \in \Endo_{\mathbf{M}^!}(\tilde{L})$ and $L^! := L^![s^L]$ as performed in \eqref{eqn:descent-ellipticity-I} and \eqref{eqn:descent-ellipticity-K}, so that $\mathfrak{a}_{L^!} = \mathfrak{a}_L$. Since $L$ (resp.\ $L^!$) is uniquely determined by $\mathfrak{a}_L \subset \mathfrak{a}_M$ (resp.\ $\mathfrak{a}_{L^!} \subset \mathfrak{a}_{M^!}$), these two maps are indeed mutually inverse.
\end{proof}

\begin{lemma}\label{prop:GM-LR}
	In the setting of Lemma \ref{prop:sL-Ls}, fix $L \in \mathcal{L}^G(R)$ and assume $d^G_R(M, L) \neq 0$. Then the map
	\begin{align*}
		\left\{\begin{array}{r|l}
			s \in \Endo_{\mathbf{M}^!}(\tilde{G}) & r = r(s, t) \text{ satisfies \eqref{eqn:descent-ellipticity}} \\
			& \text{w.r.t.}\; R \subset L \subset G
		\end{array}\right\} & \to \Endo_{\mathbf{R}^!}(\tilde{L}) \\
		s & \mapsto r^L
	\end{align*}
	is surjective, each fiber having cardinality equal to $\# \Ker(p_2)$ where $p_2: Z_{\tilde{M}^\vee}^\circ / Z_{\tilde{G}^\vee}^\circ \to Z_{\tilde{R}^\vee}^\circ / Z_{\tilde{L}^\vee}^\circ$ is the natural homomorphism. Moreover, $p_2$ is surjective with finite kernel.
\end{lemma}
\begin{proof}
	In what follows, we interpret endoscopic data on dual groups; see \S\ref{sec:dual}. Take dual Levi subgroups $\tilde{M}^\vee \supset \tilde{R}^\vee \subset \tilde{L}^\vee$ inside $\tilde{G}^\vee$. Fix a maximal torus $T^\vee \subset \tilde{R}^\vee$. Suppose that $s^\flat \in T^\vee \subset \tilde{R}^\vee$ determines $\mathbf{R}^!$; we may identify $t$ (resp.\ $s$) as an element of $Z_{\tilde{M}^\vee}^\circ / Z_{\tilde{R}^\vee}^\circ$ (resp.\ $Z_{\tilde{M}^\vee}^\circ / Z_{\tilde{G}^\vee}^\circ$) such that $s^\flat t \in \tilde{M}^\vee$ (resp.\ $s^\flat t s$) determines the elliptic endoscopic datum $\mathbf{M}^![t]$ (resp.\ $\mathbf{G}^![s]$). Ditto for $R \subset L$ and $r^L$, etc.

	We contend that $s \mapsto r^L$ can be interpreted in terms of $p_2$; in particular, $\{s: s \mapsto r^L \}$ has the same size as $p_2^{-1}(r^L)$ if we view $r^L$ as an element of $Z_{\tilde{R}^\vee}^\circ / Z_{\tilde{G}^\vee}^\circ$ giving rise to an elliptic endoscopic datum of $\tilde{L}$.

	Indeed, if one removes the exponent $\circ$, i.e.\ working with endoscopy of $\SO(2n+1)$ instead of $\Mp(2n)$, this is exactly the arguments in \cite[pp.236--237]{MW16-1}. Adding $(\cdots)^\circ$ passes to the metaplectic setting, cf.\ Remark \ref{rem:centerless}. The assertions about kernel and surjectivity of $p_2$ also reduce to the version in \textit{loc.\ cit.}
\end{proof}

\begin{proposition}\label{prop:descent-orbint-Endo}
	Let $R$ be a Levi subgroup of $M$ and $\mathbf{R}^! \in \Endo_{\elli}(\tilde{R})$. Fix a stable semisimple class $\mathcal{O}^!$ in $R^!(F)$ which is $\tilde{G}$-equisingular when $F$ is Archimedean. Then
	\[ I^{\Endo}_{\tilde{M}}\left( \mathbf{M}^![t], \delta[t]^{M^![t]}, f \right) = \sum_{L \in \mathcal{L}(R)} d^G_R(M, L) I^{\tilde{L}, \Endo}_{\tilde{R}}\left( \mathbf{R}^!, \delta, f_{\tilde{L}} \right) \]
	for all $t \in \Endo_{\mathbf{R}^!}(\tilde{M})$, $f \in \orbI_{\asp}(\tilde{G}) \otimes \mes(G)^\vee$ and $\delta \in SD_{\mathrm{geom}}(R^!, \mathcal{O}^!) \otimes \mes(R)^!$.
\end{proposition}
\begin{proof}
	We may and do assume that $\tilde{G} = \Mp(W)$ without loss of generality. To see this, note that if $\tilde{G} = \prod_i \GL(n_i) \times \Mp(W)$, we have $\orbI_{\asp}(\tilde{G}) = \bigotimes_i \orbI(\GL(n_i)) \otimes \orbI_{\asp}(\Mp(W))$ when $F$ is non-Archimedean; in this case, $I_{\tilde{M}}^{\Endo}$, $\mathcal{L}(R)$, $d^G_R(M, L)$ and $I^{\tilde{L}, \Endo}_{\tilde{R}}$ all decompose accordingly so that we are reduced to the case of $\GL(n_i)$ (which is trivial) and the case of $\Mp(W)$. When $F$ is Archimedean, elements of $\bigotimes_i \orbI(\GL(n_i)) \otimes \orbI_{\asp}(\Mp(W))$ are dense in $\orbI_{\asp}(\tilde{G})$, as one sees on the level of $C^\infty_{c, \asp}$; again, we are reduced to the cases of $\GL(n_i)$ and $\Mp(W)$ by continuity.

	For the ease of notations, write $\mathbf{M}^! = \mathbf{M}^![t]$. For every $s \in \Endo_{\mathbf{M}^!}(\tilde{G})$, define $r = r (s, t) \in \Endo_{\mathbf{R}^!}(\tilde{G})$ and $r^L \in \Endo_{\mathbf{R}^!}(\tilde{L})$ as in Lemma \ref{prop:sL-Ls}. By construction, $\mathbf{G}^![s] = \mathbf{G}^![r]$ and we have
	\begin{equation*}
		\delta[t]^{M^!}[s] = (\delta \cdot z[t])^{M^!} \cdot \underbracket{z[s]}_{\in Z_{M^!}(F)} = (\delta \cdot z[s] \cdot z[t])^{M^!} = \delta[r]^{M^!}.
	\end{equation*}
	
	By the stable descent formula \cite[II.1.14 Proposition (ii)]{MW16-1} and Remark \ref{rem:B-restriction}, we have
	\begin{multline*}
		I^{\Endo}_{\tilde{M}}\left( \mathbf{M}^!, \delta[t]^{M^!}, f \right) = \sum_{s \in \Endo_{\mathbf{M}^!}(\tilde{G})} i_{M^!}(\tilde{G}, G^![s]) S^{G^![s]}_{M^!}\left( \delta[t]^{M^!}[s], B^{\tilde{G}}, f^{G^![s]} \right) \\
		= \sum_{s \in \Endo_{\mathbf{M}^!}(\tilde{G})} i_{M^!}(\tilde{G}, G^![s]) S^{G^![r]}_{M^!}\left( \delta[r]^{M^!}, B^{\tilde{G}}, f^{G^![r]} \right) \\
		= \sum_{s \in \Endo_{\mathbf{M}^!}(\tilde{G})} i_{M^!}(\tilde{G}, G^![s]) \sum_{L^! \in \mathcal{L}^{G^![r]}(R^!)} e^{G^![r]}_{R^!}(M^!, L^!) S^{L^!}_{R^!}\left( \delta[r], B^{\tilde{G}}|_{L^!}, (f^{G^![r]})_{L^!} \right)
	\end{multline*}

	Next, use the bijection $(s, L^!) \leftrightarrow (L, s)$ in Lemma \ref{prop:sL-Ls} to rewrite the above as the sum over $(L, s)$ of
	\begin{gather*}
		i_{M^!}(\tilde{G}, G^![s]) e^{G^![r]}_{R^!}(M^!, L^!) S^{L^!}_{R^!} \left( \delta[r^L], B^{\tilde{L}}, (f_{\tilde{L}})^{L^!} \right).
	\end{gather*}
	This is based on the following observations:
	\begin{compactitem}
		\item we have $L^! = L^![r^L]$ and the situation can be summarized as
		\[\begin{tikzcd}[column sep=small]
			& G^![s] & \\
			M^! \arrow[dash, ru, "s"] & & L^! \arrow[dash, lu, "{r_L}"'] \\
			& R^! \arrow[dash, lu, "t"] \arrow[dash, ru, "{r^L}"'] \arrow[dash, uu, "r" description] &
		\end{tikzcd} \]
		\item $B^{\tilde{G}}|_{L^!} = B^{\tilde{L}}$ by Lemma \ref{prop:B-Levi};
		\item $\delta[r] = \delta[r^L] \cdot z[r_L]$ and Proposition \ref{prop:Levi-central-twist} entails
		\[ S^{L^!}_{R^!}\left( \delta[r], B^{\tilde{L}}, (f^{G^![r]})_{L^!} \right) = S^{L^!}_{R^!}\left( \delta[r] \cdot z[r_L]^{-1} , B^{\tilde{L}}, (f_{\tilde{L}})^{L^!} \right). \]
	\end{compactitem}

	Therefore
	\begin{equation}\label{eqn:descent-orbint-endo-aux-0}
		I^{\Endo}_{\tilde{M}}\left( \mathbf{M}^!, \delta[t]^{M^!}, f \right) = \sum_{(L, s)} d^G_R(M, L) X(s) S^{L^!}_{R^!} \left( \delta[r^L], B^{\tilde{L}}, (f_{\tilde{L}})^{L^!} \right)
	\end{equation}
	where $(s, L^!) \leftrightarrow (L, s)$ as before and
	\begin{align*}
		X(s) & := i_{M^!}(\tilde{G}, G^![s]) \left( Z_{{M^!}^\vee} \cap Z_{{L^!}^\vee} : Z_{{G^![s]}^\vee} \right)^{-1} \\
		& = i_{M^!}(\tilde{G}, G^![s]) \cdot \# \Ker\left( \frac{Z_{{M^!}^\vee}}{Z_{{G^![s]}^\vee}} \to \frac{Z_{{R^!}^\vee}}{Z_{{L^!}^\vee}} \right)^{-1}.
	\end{align*}
	We claim that
	\[ i_{M^!}(\tilde{G}, G^![s]) = \# \Ker\left( \phi^G_M: \frac{Z_{\tilde{M}^\vee}^\circ}{Z_{\tilde{G}^\vee}^\circ} \twoheadrightarrow \frac{Z_{{R^!}^\vee}}{Z_{{L^!}^\vee}} \right)^{-1} \]
	where $\phi^G_M$ is the evident map. To see this, apply Snake Lemma to
	\[\begin{tikzcd}
		1 \arrow[r] & Z_{\tilde{G}^\vee}^\circ \arrow[r] \arrow[hookrightarrow, d] & Z_{\tilde{M}^\vee}^\circ \arrow[hookrightarrow, d] \arrow[r] & Z_{\tilde{M}^\vee}^\circ / Z_{\tilde{G}^\vee}^\circ \arrow[d, "\phi^G_M"] \arrow[r] & 1 \\
		1 \arrow[r] & Z_{{G^![s]}^\vee} \arrow[r] & Z_{{M^!}^\vee} \arrow[r] & Z_{{M^!}^\vee}/Z_{{G^![s]}^\vee} \arrow[r] & 0
	\end{tikzcd}\]
	and observe that $\phi^G_M$ is surjective since \cite[Lemma 1.1]{Ar99} implies
	\[ Z_{{M^!}^\vee} = Z_{{M^!}^\vee}^\circ Z_{{G^![s]}^\vee} = Z_{\tilde{M}^\vee}^\circ Z_{{G^![s]}^\vee}. \]
	
	For the same reason, when $\mathfrak{a}^{M^!}_{R^!} \oplus \mathfrak{a}^{L^!}_{R^!} \rightiso \mathfrak{a}^{G^![s]}_{R^!}$ (equivalently, $d^G_R(M, L) \neq 0$) we have
	\[ Z_{{R^!}^\vee} = Z_{{M^!}^\vee} Z_{{R^!}^\vee}^\circ = Z_{{M^!}^\vee} Z_{{M^!}^\vee}^\circ Z_{{L^!}^\vee}^\circ = Z_{{M^!}^\vee} Z_{{L^!}^\vee}.  \]
	Therefore $Z_{{M^!}^\vee} / Z_{{G^![s]}^\vee} \twoheadrightarrow Z_{{R^!}^\vee} / Z_{{L^!}^\vee}$.

	Assuming $d^G_R(M, L) \neq 0$, we arrive at the commutative diagram of abelian groups
	\begin{equation}\label{eqn:descent-orbint-endo-aux-1}\begin{tikzcd}
		Z_{\tilde{M}^\vee}^\circ / Z_{\tilde{G}^\vee}^\circ \arrow[twoheadrightarrow, r, "\phi^G_M"] \arrow[d, "p_2"'] & Z_{{M^!}^\vee} / Z_{{G^![s]}^\vee} \arrow[twoheadrightarrow, d] \\
		Z_{\tilde{R}^\vee}^\circ / Z_{\tilde{L}^\vee}^\circ \arrow[r, "\phi^L_R"'] & Z_{{R^!}^\vee}/Z_{{L^!}^\vee}
	\end{tikzcd}\end{equation}
	with evident arrows. Denote the diagonal composition as $p_1(s)$. The discussions above lead to
	\begin{equation}\label{eqn:descent-orbint-endo-aux-2}
		X(s)^{-1} = \# \Ker\left( p_1(s) \right).
	\end{equation}
	
	On the other hand, Lemma \ref{prop:GM-LR} together with \eqref{eqn:descent-orbint-endo-aux-1} imply that when $d^G_R(M, L) \neq 0$, we have
	\[ \frac{\# \Ker(p_2)}{\# \Ker(p_1(s))} = \left(\# \Ker(\phi^L_R)\right)^{-1} = i_{R^!}(\tilde{L}, L^!). \]
	
	To conclude the proof, we return to \eqref{eqn:descent-orbint-endo-aux-0}. From $L^! = L^![r^L]$, \eqref{eqn:descent-orbint-endo-aux-1} and \eqref{eqn:descent-orbint-endo-aux-2} we see that $X(s)$ depends only on the image $r^L \in \Endo_{\mathbf{R}^!}(\tilde{L})$ of $s$. In view of Lemma \ref{prop:GM-LR}, we conclude that
	\begin{multline*}
		I^{\Endo}_{\tilde{M}}\left( \mathbf{M}^!, \delta[t]^{M^!}, f \right) = \sum_{(L, s)} d^G_R(M, L) X(s) S^{L^!}_{R^!} \left( \delta[r^L], B^{\tilde{L}}, (f_{\tilde{L}})^{L^!} \right) \\
		= \sum_{L \in \mathcal{L}^G(R)} d^G_R(M, L) \sum_{r^L \in \Endo_{\mathbf{R}^!}(\tilde{L})} \frac{\# \Ker(p_2)}{\# \Ker (p_1(s))} S^{L^!}_{R^!} \left( \delta[r^L], B^{\tilde{L}}, (f_{\tilde{L}})^{L^!} \right) \\
		= \sum_{L \in \mathcal{L}^G(R)} d^G_R(M, L) \sum_{r^L \in \Endo_{\mathbf{R}^!}(\tilde{L})} i_{R^!}(\tilde{L}, L^!) S^{L^!}_{R^!} \left( \delta[r^L], B^{\tilde{L}}, (f_{\tilde{L}})^{L^!} \right) \\
		= \sum_{L \in \mathcal{L}^G(R)} d^G_R(M, L)  I^{\tilde{L}, \Endo}_{\tilde{R}}\left( \mathbf{R}^!, \delta, f_{\tilde{L}} \right),
	\end{multline*}
	as asserted.
\end{proof}

\begin{corollary}\label{prop:descent-orbint-Endo-2}
	In the circumstance of Proposition \ref{prop:descent-orbint-Endo}, we have
	\[ I^{\Endo}_{\tilde{M}}\left( \tilde{\gamma}^{\tilde{M}}, f \right) = \sum_{L \in \mathcal{L}(R)} d^G_R(M, L) I^{\tilde{L}, \Endo}_{\tilde{R}}\left( \tilde{\gamma}, f_{\tilde{L}} \right). \]
\end{corollary}
\begin{proof}
	Combine Proposition \ref{prop:descent-orbint-Endo}, Definition--Proposition \ref{def:I-geom-Endo} with Proposition \ref{prop:Levi-central-twist}.
\end{proof}

We will encounter many other endoscopic descent formulas of the same form. Because of the purely combinatorial nature of the arguments presented above, the proof will not be repeated anymore.

\section{Reduction to the \texorpdfstring{$\tilde{G}$}{G}-regular case}\label{sec:local-geom-reduction-G-reg}
Let $F$ be a local field of characteristic zero, and let $\rev: \tilde{G} \twoheadrightarrow G(F)$ be of metaplectic type. We begin by studying the local behavior of $I^{\Endo}_{\tilde{M}}(\tilde{\gamma}, f)$, where $M \subset G$ is a Levi subgroup.

\begin{lemma}\label{prop:weighted-arch-aux-Endo}
	Let $\mathcal{O}$ be a $G$-equisingular stable semisimple conjugacy class in $M(F)$. For every $f \in \orbI_{\asp}(\tilde{G}) \otimes \mes(G)$, there exists $g \in \orbI_{\asp}(\tilde{M}) \otimes \mes(M)$ such that
	\[ I^{\Endo}_{\tilde{M}}(\tilde{\gamma}, f) = I^{\tilde{M}}(\tilde{\gamma}, g) \]
	holds if one of the following conditions holds:
	\begin{enumerate}[(i)]
		\item $\tilde{\gamma} \in D_{\mathrm{orb}, -}(\tilde{M}) \otimes \mes(M)^\vee$, such that $\Supp(\tilde{\gamma})$ is formed by $\tilde{G}$-regular elements which are sufficiently close to $\rev^{-1}(\mathcal{O})$;
		\item $\tilde{\gamma} \in D_{\mathrm{geom}, -}(\tilde{M}, \mathcal{O}) \otimes \mes(M)^\vee$.
	\end{enumerate}
\end{lemma}
\begin{proof}
	Our arguments are modeled on \cite[V.1.9 Lemme]{MW16-1}. As in the proof of Proposition \ref{prop:descent-orbint-Endo}, without loss of generality we may assume $\tilde{G} = \Mp(W)$.

	Consider a Levi subgroup $R \subset M$ and $\mathbf{R}^! \in \Endo_{\elli}(\tilde{R})$. Let $\mathcal{O}_{R^!}$ be the finite union of stable semisimple classes in $R^!(F)$ corresponding to $\mathcal{O}$ via $\Sigma_{\mathrm{ss}}(R^!) \to \Sigma_{\mathrm{ss}}(R) \to \Sigma_{\mathrm{ss}}(M)$.

	For $L \in \mathcal{L}(R)$ with $d^G_R(M, L) \neq 0$ and $s \in \Endo_{\mathbf{R}^!}(\tilde{L})$, we claim that $\mathcal{O}_{R^!}[s]$ is $L^![s]$-equisingular. Indeed, let $\delta \in \mathcal{O}_{R^!}$ and $\gamma \in R(F)_{\mathrm{ss}}$ such that $\delta \leftrightarrow \gamma$ via endoscopy, then $\gamma$ is equisingular relative to $M \subset G$, hence Lemma \ref{prop:d-equisingular} implies that $\gamma$ is equisingular relative to $R \subset L$. In view of Definition \ref{def:equisingular-endo} and Lemma \ref{prop:equisingular-endo}, the $L^![s]$-equisingularity of $\delta[s]$ follows.
	
	By \cite[V.1.4 (2), (3)]{MW16-1} and its non-Archimedean analogues in \cite[III]{MW16-1}, for all quadruplets $(\tilde{R}, \mathbf{R}^!, \tilde{L}, s)$ above, there exists $g_{\tilde{R}, \mathbf{R}^!, \tilde{L}, s} \in S\orbI(R^!) \otimes \mes(R^!)$ such that
	\begin{equation}\label{eqn:weighted-arch-aux-Endo-0}
		S^{L^![s]}_{R^!}\left( \delta[s], B^{\tilde{L}}, \Trans_{\mathbf{L}^![s], \tilde{L}} (f_{\tilde{L}}) \right) = S^{R^!}\left( \delta, g_{\tilde{R}, \mathbf{R}^!, \tilde{L}, s} \right)
	\end{equation}
	holds in each of the following circumstances:
	\begin{enumerate}[(i)$'$]
		\item $\delta \in SD_{\mathrm{orb}}(R^!) \otimes \mes(R)^\vee$, such that $\Supp(\delta)$ is formed by elements which are sufficiently regular and close to $\mathcal{O}_{R^!}$;
		\item $\delta \in SD_{\mathrm{geom}}(R^!, \mathcal{O}_{R^!}) \otimes \mes(R^!)^\vee$.
	\end{enumerate}
	Furthermore, when $R \subset M$ and $\mathbf{R}^!$ are kept fixed, $W^M(R)$ (as a subgroup of $W^G_0$) acts on $R^!$ and transports the data $(\tilde{L}, s)$. By transport of structure, the value of the left hand side of \eqref{eqn:weighted-arch-aux-Endo-0} is unaltered under $W^M(R)$-action. Taking averages, we may assume that each $g_{\tilde{R}, \mathbf{R}^!, \tilde{L}, s}$ is $W^M(R)$-invariant.

	Given $R \subset M$ and $\mathbf{R}^! \in \Endo_{\elli}(\tilde{R})$, put
	\[ g_{\tilde{R}, \mathbf{R}^!} := \sum_{L \in \mathcal{L}(R)} d^G_R(M, L) \sum_{s \in \Endo_{\mathbf{R}^!}(\tilde{L})} i_{R^!}(\tilde{L}, L^![s]) g_{\tilde{R}, \mathbf{R}^!, \tilde{L}, s}. \]

	We claim that, upon multiplying by a cut-off function (see below), $\left( g_{\tilde{M}, \mathbf{M}^!} \right)_{\mathbf{M}^! \in \Endo_{\elli}(\tilde{M})}$ belongs to the space $\orbI^{\Endo}(\tilde{M})$ in Definition \ref{def:orbI-Endo}.
	
	Consider $R \subset M$, $\mathbf{R}^! \in \Endo_{\elli}(\tilde{R})$ and $t \in \Endo_{\mathbf{R}^!}(\tilde{M})$; put $\mathbf{M}^! := \mathbf{M}^![t]$. Take $\delta \in SD_{\mathrm{geom}}(R^!) \otimes \mes(R)^\vee$ subject to (i)$'$ or (ii)$'$. Since $d^G_M(M, L) \neq 0$ only when $L = G$, we have
	\begin{equation}\label{eqn:weighted-arch-aux-Endo}\begin{aligned}
		S^{R^!} \left( \delta, (g_{\tilde{M}, \mathbf{M}^!})^{R^!}[t] \right) & = S^{M^!}\left( \delta[t]^{M^!}, g_{\tilde{M}, \mathbf{M}^!} \right) \\
		& = \sum_{s \in \Endo_{\mathbf{M}^!}(\tilde{G})} i_{M^!}(\tilde{G}, G^![s]) S^{M^!}\left( \delta[t]^{M^!}, g_{\tilde{M}, \mathbf{M}^!, \tilde{G}, s} \right) \\
		& = I^{\tilde{G}, \Endo}_{\tilde{M}} \left( \mathbf{M}^!, \delta[t]^{M^!}, f \right) \\
		& = I^{\tilde{G}, \Endo}_{\tilde{M}}\left( \trans_{\mathbf{M}^!, \tilde{M}}\left(\delta[t]^{M^!}\right), f \right) = I^{\tilde{G}, \Endo}_{\tilde{M}}\left( \trans_{\mathbf{R}^!, \tilde{R}}(\delta)^{\tilde{M}}, f \right)
	\end{aligned}\end{equation}
	where we invoked Definition--Proposition \ref{def:I-geom-Endo} and Proposition \ref{prop:Levi-central-twist}. The last expression in \eqref{eqn:weighted-arch-aux-Endo} is hence independent of $t$ and $W^M(R)$-invariant.
	
	So far we take $\delta$ with $\Supp(\delta)$ close to $\mathcal{O}_{R^!}$, which implies $\Supp\left(\delta[t]^{M^!} \right)$ is close to $\mathcal{O}_{M^!}$. Consider the diagram of $F$-varieties
	\[\begin{tikzcd}[row sep=small]
		M^! \arrow[twoheadrightarrow, d] & M \arrow[twoheadrightarrow, d] \\
		T^! /\!/ W(M^!, T^!) \arrow[twoheadrightarrow, r] & T /\!/ W(M, T)
	\end{tikzcd}\]
	where $T^!$, $T$ are split maximal $F$-tori, the vertical arrows are those of Chevalley, and the horizontal arrow comes from endoscopy. Fix a $C^\infty$ cut-off function on $(T /\!/ W(M,T))(F)$ centered at the point classifying $\mathcal{O}$, of sufficiently small support, and pull it back to $(T^! /\!/ W(M^!, T^!))(F)$. As $g_{\tilde{M}, \mathbf{M}^!} \in S\orbI(M^!) \otimes \mes(M)$ can be viewed as functions on the strongly regular locus of $(T^! /\!/ W^{M^!}(T^!))(F)$, they may be multiplied by the cut-off function to remove the constraint on supports.

	Hence Theorem \ref{prop:image-transfer} affords $g \in \orbI_{\asp}(\tilde{M}) \otimes \mes(M)$ with $\Trans^{\Endo}(g) = \left( g_{\tilde{M}, \mathbf{M}^!} \right)_{\mathbf{M}^! \in \Endo_{\elli}(\tilde{M})}$.

	Finally, let $\tilde{\gamma}$ be as in the assertions. By Theorem \ref{prop:Dgeom-preservation}, we may write
	\[ \tilde{\gamma} = \sum_{(\mathbf{M}^!, \delta)} \trans_{\mathbf{M}^!, \tilde{M}}(\delta) \]
	where $\mathbf{M}^! \in \Endo_{\elli}(\tilde{M})$ and $\delta \in SD_{\mathrm{geom}}(M^!) \otimes \mes(M^!)^\vee$ fulfills either (i)$'$ or (ii)$'$ with $R = M$. Hence \eqref{eqn:weighted-arch-aux-Endo} with $R=M$ implies
	\[ I^{\tilde{M}}(\tilde{\gamma}, g) = \sum_{(\mathbf{M}^!, \delta)} S^{M^!}\left( \delta, g_{\tilde{M}, \mathbf{M}^!} \right) = \sum_{(\mathbf{M}^!, \delta)} I^{\Endo}_{\tilde{M}} \left( \mathbf{M}^!, \delta, f \right). \]
	This yields $I^{\Endo}_{\tilde{M}}(\tilde{\gamma}, f)$ as desired.
\end{proof}

\begin{proposition}[Cf.\ {\cite[V.1.11 Lemme]{MW16-1}}]\label{prop:local-geometric-regular-reduction}
	Let $F$ be a local field of characteristic zero. If the assertion in Theorem \ref{prop:local-geometric} holds whenever $\tilde{\gamma} \in D_{\mathrm{orb}, -}(\tilde{M}) \otimes \mes(M)^\vee$ satisfies $\Supp(\tilde{\gamma}) \subset \tilde{G}_{\mathrm{reg}}$, then it holds for all $\tilde{\gamma} \in D_{\mathrm{geom}, G\mathrm{-equi}, -}(\tilde{M}) \otimes \mes(M)^\vee$.
\end{proposition}
\begin{proof}
	Let $\mathcal{O} \subset M(F)$ be a $G$-equisingular stable semisimple conjugacy class. Fix $f \in \orbI_{\asp}(\tilde{G}) \otimes \mes(G)$. Definition \ref{def:weighted-I-integral-arch} and Lemma \ref{prop:weighted-arch-aux-Endo} afford $g_1, g_2 \in \orbI_{\asp}(\tilde{M}) \otimes \mes(M)$ such that
	\[
		I_{\tilde{M}}\left(\tilde{\gamma}, f \right) = I^{\tilde{M}}\left( \tilde{\gamma}, g_1 \right), \quad
		I^{\Endo}_{\tilde{M}}\left(\tilde{\gamma}, f \right) = I^{\tilde{M}}\left( \tilde{\gamma}, g_2 \right)
	\]
	in each of the following cases:
	\begin{enumerate}[(i)]
		\item $\tilde{\gamma} \in D_{\mathrm{orb}, -}(\tilde{M}) \otimes \mes(M)^\vee$ whose support is sufficiently regular and close to $\rev^{-1}(\mathcal{O})$;
		\item $\tilde{\gamma} \in D_{\mathrm{geom}, -}(\tilde{M}, \mathcal{O}) \otimes \mes(M)^\vee$.
	\end{enumerate}

	Case (i) and the assumptions imply $g_1 = g_2$ in some neighborhood of $\rev^{-1}(\mathcal{O})$. By plugging this into case (ii) and varying $\mathcal{O}$, we obtain $I_{\tilde{M}}\left(\tilde{\gamma}, f \right) = I^{\Endo}_{\tilde{M}}\left(\tilde{\gamma}, f \right)$ for all $\tilde{\gamma} \in D_{\text{geom, $G$-equi}, -}(\tilde{M}) \otimes \mes(M)^\vee$.
\end{proof}

\section{Statement of the weighted fundamental lemma}\label{sec:LFP}
Consider a group of metaplectic type $\tilde{G} \xrightarrow{\rev} G(F)$ in the unramified situation as in \S\ref{sec:LF}. In particular, we have a hyperspecial subgroup $K = G(\mathfrak{o}_F)$ of $G(F)$, over which $\rev$ splits canonically. It is harmless to assume in what follows that $\tilde{G} = \Mp(W)$, for the sake of simplicity.

Take the Haar measure on $G(F)$ with $\mes(K) = 1$. Let $f_K$ be the unit of $\mathcal{H}_{\asp}(K \backslash \tilde{G} / K)$. If $M \subset G$ is a Levi subgroup in good position relative to $K$, we take the Haar measure on $M(F)$ with $\mes(K \cap M(F)) = 1$. This trivializes the lines $\mes(\cdots)$ in various formulas, and we will omit them from the notation.

\begin{definition}\label{def:r-unramified}
	\index{rGMgammaK@$r^{\tilde{G}}_{\tilde{M}}(\tilde{\gamma}, K)$}
	Let $M \subset G$ be a Levi subgroup in good relative position to $K$. For every $\tilde{\gamma} \in D_{\mathrm{geom}, -}(\tilde{M})$, set
	\[ r^{\tilde{G}}_{\tilde{M}}(\tilde{\gamma}, K) := J^{\tilde{G}}_{\tilde{M}}(\tilde{\gamma}, f_K). \]
\end{definition}

The descent formula for weighted orbital integrals \cite[Corollaire 5.4.3]{Li14a} leads to
\[ r^{\tilde{G}}_{\tilde{L}}\left(\tilde{\gamma}^L, K \right) = \sum_{L^\dagger \in \mathcal{L}^G(M)} d^G_M(L, L^\dagger) r^{\tilde{L}}_{\tilde{M}}(\tilde{\gamma}, K \cap L(F)) \]
where $L \in \mathcal{L}^G(M)$ and $\tilde{\gamma}^L$ is as in Proposition \ref{prop:orbint-weighted-descent-nonArch}.

On the other hand, given a quasisplit group $G^!$ and its Levi subgroup $M^!$, in \cite[II.4.2]{MW16-1} are defined the stable avatars
\[ s^{G^!}_{M^!}(\delta) = s^{G^!}_{M^!}(\delta, K^!), \quad \delta \in SD_{\mathrm{geom}}(M^!). \]
Here $K^! \subset G^!(F)$ is a hyperspecial subgroup and $M^!$ is in good position relative to $K^!$. The Haar measures on $G^!(F)$ and $M^!(F)$ are determined by $K^!$ as before. The choice of $K^!$ is immaterial here because they are conjugate under $G^!_{\mathrm{AD}}(F)$, see \cite[II.4.2 (1)]{MW16-1}. On the other hand, these quantities depend on the choice of an invariant positive-definite quadratic form on $\mathfrak{a}^{G^!}_{M^!}$. Such a form is available when we are in the situation \eqref{eqn:s-situation}, namely
\[\begin{tikzcd}
	G^! \arrow[dashed, leftrightarrow, r, "\text{ell.}", "\text{endo.}"'] & \tilde{G} \\
	M^! \arrow[dashed, leftrightarrow, r, "\text{ell.}", "\text{endo.}"'] \arrow[hookrightarrow, u, "\text{Levi}"] & \tilde{M} \arrow[hookrightarrow, u, "\text{Levi}"'] .
\end{tikzcd}\]

\begin{definition}\label{def:rEndo-unramified}
	\index{rGEndoMdelta@$r^{\tilde{G}, \Endo}_{\tilde{M}}(\mathbf{M}^{"!}, \delta)$}
	Given $\tilde{M} \subset \tilde{G}$ and $\mathbf{M}^! \in \Endo_{\elli}(\tilde{M})$ as above, we set
	\[ r^{\tilde{G}, \Endo}_{\tilde{M}}(\mathbf{M}^!, \delta) := \sum_{s \in \Endo_{\mathbf{M}^!}(\tilde{G})} i_{M^!}(\tilde{G}, G^![s]) s^{G^![s]}_{M^!}(\delta[s]) \]
	for all $\delta \in SD_{\mathrm{geom}}(M^!)$.
\end{definition}

Now we can state the \emph{weighted fundamental lemma}.

\begin{theorem}\label{prop:LFP-general}
	\index{weighted fundamental lemma}
	In the unramified situation, suppose moreover that the residual characteristic $p$ of $F$ satisfies $p > 5$. Given $\tilde{M} \subset \tilde{G}$ as above, we have
	\[ r^{\tilde{G}}_{\tilde{M}}\left( \trans_{\mathbf{M}^!, \tilde{M}}(\delta), K \right) = r^{\tilde{G}, \Endo}_{\tilde{M}}\left( \mathbf{M}^!, \delta \right) \]
	for all $\mathbf{M}^! \in \Endo_{\elli}(\tilde{M})$ and $\delta \in SD_{\mathrm{geom}}(M^!)$.
\end{theorem}

When $\delta \in SD_{\mathrm{geom}, \tilde{G}\text{-reg}}(M^!)$, this reduces to \cite[Théorème 4.2.1]{Li12a}. The proof of the general case will be given in \S\ref{sec:proof-LFP}.

\chapter{Theory of germs}\label{sec:germs}
Consider a covering of metaplectic type $\tilde{G}$ over a local field $F$ of characteristic zero. In the first part, we will follow \cite[II, III]{MW16-1} to develop a theory of germs. Let $M \in \mathcal{L}(M_0)$. The germs in question are associated with certain subsets $\underline{J}$ of $\Sigma(A_M)$ with $|\underline{J}| = \dim \mathfrak{a}^G_M$, taken up to the equivalence that $\underline{J} \sim \underline{J}'$ if they generate the same $\Z$-module (in fact, lattice) in $\mathfrak{a}^{G,*}_M$. Fix a semisimple conjugacy class $\mathcal{O}$ in $M(F)$. The germs are represented by linear maps
\[ \rho^{\tilde{G}}_J: D_{\mathrm{orb}, -}(\tilde{M}, \mathcal{O}) \otimes \mes(M)^\vee \to U_J \otimes \left( D_{\mathrm{orb}, -}(\tilde{M}, \mathcal{O}) \otimes \mes(M)^\vee \right) \big/ \mathrm{Ann}^{\tilde{G}}_{\mathcal{O}} \]
where
\begin{itemize}
	\item $J$ denotes an equivalence class as above and $U_J$ is a certain space of germs of functions on $A_M(F)$ near $a=1$,
	\item $\mathrm{Ann}^{\tilde{G}}_{\mathcal{O}}$ denotes the kernel of parabolic induction to $\tilde{G}$,
	\item we only consider the restriction of $\rho^{\tilde{G}}_J$ to distributions of support sufficiently close to $\rev^{-1}(\mathcal{O})$.
\end{itemize}
They are related to invariant weighted orbital integrals through the equivalence
\[ I_{\tilde{M}}\left(a\tilde{\gamma}, f \right) \sim \sum_{L \in \mathcal{L}(M)} \sum_{J \in \mathcal{J}^L_M} I_{\tilde{L}}\left( \rho^{\tilde{L}}_J(\tilde{\gamma}, a)^{\tilde{L}}, f \right)  \]
between germs of functions in $a \in A_M(F)$, $a \to 1$; see Definition \ref{def:germ-equiv} for the meaning of $\sim$. Note that $a\tilde{\gamma}$ is $\tilde{G}$-equisingular for general $a$, which explains the usefulness of germs towards the local geometric Theorem \ref{prop:local-geometric}.

Note that $D_{\mathrm{geom}, -} = D_{\mathrm{orb}, -}$ for non-Archimedean $F$. As for the Archimedean case, in \S\ref{sec:ext-Arch} we will lay down a program to extend the definitions of $\rho^{\tilde{G}}_J(\tilde{\gamma})$, $I_{\tilde{M}}(\tilde{\gamma}, f)$ and $I^{\Endo}_{\tilde{M}}(\tilde{\gamma}, f)$ to a larger subspace $D_{\text{tr-orb}, -}(\tilde{M})$ of $D_{\mathrm{geom}, -}(\tilde{M})$ that accommodates the transfers of similar distributions from various $\mathbf{M}^! \in \Endo_{\elli}(\tilde{M})$. Furthermore, we expect that the $\rho^{\tilde{G}}_J$ and $I_{\tilde{M}}$ so extended match their endoscopic counterparts $\rho^{\tilde{G}, \Endo}_J$ and $I^{\Endo}_{\tilde{M}}$.

This program will be carried out in \S\ref{sec:ext-Arch-proof}, under the Hypothesis \ref{hyp:ext-Arch} that the local geometric Theorem \ref{prop:local-geometric} holds for $\tilde{G}$-regular $\tilde{\gamma}$. In this process, subgroups $G_J \subset G$ and the corresponding coverings $\widetilde{G_J} \subset \tilde{G}$ will appear naturally: in our case, $G_J$ are simply products of smaller symplectic groups.

It turns out that the same argument can also handle the non-Archimedean case, if one replaces $D_{\mathrm{tr-orb}, -}$ by $D_{\mathrm{geom}, -}$ everywhere. In particular, assume that the local geometric Theorem \ref{prop:local-geometric} holds in the $\tilde{G}$-regular case, then it holds in general, and $\rho^{\tilde{G}}_J = \rho^{\tilde{G}, \Endo}_J$. This is the content of \S\ref{sec:rho-main-summary}.

In the second part, we introduce the weighted Shalika germs on $\tilde{G}$ in the non-Archimedean case, and define their endoscopic counterparts. The matching Theorem \ref{prop:germ-matching} for Shalika germs can be proved directly via Harish-Chandra descent: it is ultimately reduced to the unipotent matching settled in \cite[III]{MW16-1}, for both standard and nonstandard endoscopy. The matching of Shalika germs can be employed in the following ways.
\begin{itemize}
	\item It furnishes another way to reduce the local geometric Theorem \ref{prop:local-geometric} to the $\tilde{G}$-regular case (Corollary \ref{prop:local-geometric-regular-reduction-nonArch}).
	\item It will be used in the construction of the map $\epsilon_{\tilde{M}}$ in \S\ref{sec:epsilon-nonarch}.
	\item Certain arguments in the proof will be re-used in the Archimedean case in \S\ref{sec:pf-jump-Endo-lemmas}.
\end{itemize}

Finally, we apply the technique of germs to prove the weighted fundamental lemma (Theorem \ref{prop:LFP-general}) in \S\ref{sec:proof-LFP}.

\section{Algebraic constructions}
This section is largely a review of \cite[II.3.1]{MW16-1}. We begin with a general connected reductive $F$-group $G$ and its Levi subgroup $M$. Here $F$ can be any field of characteristic zero.

\subsection{Definition of \texorpdfstring{$G_J$}{GJ}}
View $\Sigma(A_M) = \Sigma^G(A_M)$ as a subset of $\mathfrak{a}^{G, *}_M$ and define
\[ \mathcal{J}^G_M := \left\{ \underline{J} \subset \Sigma(A_M): |\underline{J}| = \dim \mathfrak{a}^G_M, \; \text{linearly independent} \right\} \big / \sim \]
where we write $\underline{J} \sim \underline{J}'$ if they generate the same $\Z$-submodule in $\mathfrak{a}^{G, *}_M$.
\index{JGM@$\mathcal{J}^G_M$}

For each $J \in \mathcal{J}^G_M$, let $R_J$ be the $\Z$-submodule of $\Sigma(A_M)$ it generates. Make $\mathcal{J}^G_M$ into a partially ordered set by setting $J \leq J'$ if and only if $R_J \subset R_{J'}$. There is then a unique maximal element $J_{\max}$ of $(\mathcal{J}^G_M, \leq)$, characterized by $R_{J_{\max}} = \Sigma(A_M)$.
\index{Jmax@$J_{\max}$}

Fix a Borel pair $(B_G, T)$ for $G_{\overline{F}}$ such that $T$ is defined over $F$, and $M_{\overline{F}}$ is a standard Levi subgroup relative to $(B_G, T)$; in particular $A_M \subset T$. Given $J \in \mathcal{J}^G_M$, we have the root subsystem
\[ \left\{ \beta \in \Sigma(G, T): \beta|_{A_M} \in R_J \right\}. \]
It is the absolute root system (relative to $T$) of a connected reductive $F$-group $G_J$ such that $M \subset G_J \subset G$. One readily checks that $G_J$ is independent of $(B_G, T)$. Moreover, $M$ is still a Levi subgroup of $G_J$.
\index{GJ@$G_J$}

Note that $G_{J_{\max}} = G$. When $J \neq J_{\max}$, we have $\dim G_J < \dim G$.

For all $J \in \mathcal{J}^G_M$, we have $\Sigma^{G_J}(A_M) \subset \Sigma^G(A_M)$, and
\[ \mathcal{J}^{G_J}_M = \left\{ J' \in \mathcal{J}^G_M: R_{J'} \subset \Sigma^{G_J}(A_M) \right\}. \]
Moreover, $J \in \mathcal{J}^{G_J}_M$ and it is the unique maximal element of $\mathcal{J}^{G_J}_M$.

\begin{proposition}\label{prop:GJ-ss-rank}
	For every $J \in \mathcal{J}^G_M$, the groups $G_{\mathrm{der}}$ and $G_{J, \mathrm{der}}$ have the same rank over $\overline{F}$.
\end{proposition}
\begin{proof}
	They share the Levi subgroup $M$ and $\mathfrak{a}^G_M = |J| = \mathfrak{a}^{G_J}_M$.
\end{proof}

\subsection{The case with \texorpdfstring{$B$}{B}-functions}\label{sec:GJ-B}
The following is a review of \cite[II.3.3]{MW16-1}. Let $L$ be a connected reductive $F$-group with Levi subgroup $M$.

Fix a $B$-function $\Sigma(L, T) \to \Q_{> 0}$ as in Definition \ref{def:B-function}. We define the analogue $\mathcal{J}^L_M(B)$ of $\mathcal{J}^L_M$ by the same recipe, replacing $\Sigma(A_M)$, $\Sigma(L, T)$ by $\Sigma(A_M; B)$, $\Sigma(L, T; B)$ systematically.
\index{JLM-B@$\mathcal{J}^L_M(B)$}

Given $J \in \mathcal{J}^L_M(B)$, let
\[ \Sigma_J := \left\{ \alpha \in \Sigma(L, T) : B(\alpha)^{-1}\alpha|_{\mathfrak{a}_M} \in R_J \right\}, \quad
\check{\Sigma}_J := \left\{ \check{\alpha}: \alpha \in \Sigma_J \right\}. \] 
Then we construct a connected reductive $F$-group $L_J$ containing $T$, whose roots and coroots are given by $\Sigma_J$ and $\check{\Sigma}_J$ respectively. It reverts to the prior construction when $B$ is trivial. It still satisfies:
\begin{compactitem}
	\item $M$ is canonically embedded as a Levi subgroup of $L_J$,
	\item $Z_L \subset Z_{L_J}$,
	\item $Z_{\check{L}} \hookrightarrow Z_{\check{L}_J}$, compatibly with $\Gamma_F$-actions.
\end{compactitem}
However, $L_J$ is not necessarily a subgroup of $L$ for non-constant $B$.

Now, suppose that $G^!$ is a quasisplit group endowed with a system of $B$-functions (Definition \ref{def:system-B-function}), and $M^! \subset G^!$ is a Levi subgroup. In the discussions after that definition, we defined a subset $\Sigma(A_{M^!}, B_{\mathcal{O}^!})$ of $\mathfrak{a}_{M^!}^*$ for every semisimple conjugacy class $\mathcal{O}^!$ in $M^!(F)$. It consists of nonzero restrictions to $A_{M^!}$ of elements from $\Sigma(A_{M^!_\epsilon}, B_\epsilon)$ (taken inside $G^!_\epsilon$) where $\epsilon \in \mathcal{O}^!$ is arbitrary. In particular, we can define the analogue $\mathcal{J}^{G^!}_{M^!}(B_{\mathcal{O}^!})$ of $\mathcal{J}^{G^!}_{M^!}$ by considering linearly independent subsets of $\Sigma(A_{M^!}, B_{\mathcal{O}^!})$ of size $\dim \mathfrak{a}^G_M$.

It follows immediately from Definition \ref{def:system-B-function} that $\mathcal{J}^{G^!}_{M^!}(B_{\mathcal{O}^!})$ depends only on the stable semisimple conjugacy class containing $\mathcal{O}^!$. It can be empty.

In contrast with the previous case, in general $G^!_J$ does not make sense for $J \in \mathcal{J}^{G^!}_{M^!}(B_{\mathcal{O}^!})$, except when $\mathcal{O}^!$ is central.

Let us specialize this formalism to the situation
\[\begin{tikzcd}
	G^![s] \arrow[dashed, leftrightarrow, r, "\text{ell.}", "\text{endo.}"'] & \tilde{G} \\
	M^! \arrow[dashed, leftrightarrow, r, "\text{ell.}", "\text{endo.}"'] \arrow[hookrightarrow, u, "\text{Levi}"] & \tilde{M} \arrow[hookrightarrow, u, "\text{Levi}"']
\end{tikzcd}\]
where $\tilde{G}$ is a group of metaplectic type, $\mathbf{M}^! \in \Endo_{\elli}(\tilde{M})$ and $s \in \Endo_{\mathbf{M}^!}(\tilde{G})$. Consider the family of $B$-functions $B^{\tilde{G}}$ for $G^!$ from Definition \ref{def:B-metaplectic}. Let $\mathcal{O}^! \subset M^!(F)$ be a stable semisimple conjugacy class and put
\[ \mathcal{O}^![s] := \mathcal{O}^! \cdot z[s]. \]

Also recall that endoscopy affords an isomorphism $\xi: A_M \rightiso A_{M^!}$ between $F$-tori. In turn, it induces $\xi^*: \mathfrak{a}_{M^!}^* \rightiso \mathfrak{a}_M$. The following result is somewhat ``built-in'' to $B^{\tilde{G}}$.

\begin{lemma}\label{prop:J-transfer}
	In the situation above, $\xi^*: \mathfrak{a}_{M^!}^* \rightiso \mathfrak{a}_M^*$ induces an injection
	\[ \Sigma^{G^![s]}\left( A_{M^!}, B^{\tilde{G}}_{\mathcal{O}^![s]} \right) \hookrightarrow \Sigma^G(A_M), \]
	in turn, this induces an injection $\mathcal{J}^{G^![s]}_{M^!}\left( B^{\tilde{G}}_{\mathcal{O}^![s]} \right) \hookrightarrow \mathcal{J}^G_M$.
\end{lemma}
\begin{proof}
	Take $\epsilon \in \mathcal{O}^!$. Without loss of generality, we may assume $\tilde{G} = \Mp(W)$ and $M^!_\epsilon$ is quasisplit. Choose by Corollary \ref{prop:diagram-from-epsilon} a diagram $(\epsilon, B^!, T^!, B, T, \eta)$ joining $\epsilon$ and some $\eta \in M(F)_{\mathrm{ss}}$; the data include an isomorphism $\xi_{T^!, T}: T^! \rightiso T$ compatible with $\xi^{-1}$; see Definition \ref{def:diagram}. We have $\epsilon \leftrightarrow \eta$ (resp.\ $\epsilon[s] \leftrightarrow \eta$) with respect to $\mathbf{M}^! \in \Endo_{\elli}(\tilde{M})$ (resp.\ $\mathbf{G}^![s] \in \Endo_{\elli}(\tilde{G})$).

	The descent of endoscopic data in \S\ref{sec:descent-endoscopy} shows that the roots of $G^![s]_{\epsilon[s]} \supset T^!$ injects into those of $G_\eta \supset T$ in this way, up to explicit scaling. By definition, the scaling is exactly prescribed by $B^{\tilde{G}}_{\epsilon[s]}$. Restricting the roots to $A_M$ and $A_{M^!}$ gives the first injection. The second follows immediately.
\end{proof}

\subsection{Determination of \texorpdfstring{$G_J$}{GJ} in the symplectic case}
In this subsection, we take $G = \Sp(W)$ and let $M = \prod_{i \in I} \GL(n_i) \times \Sp(W^\flat)$ be a Levi subgroup of $G$, where $W^\flat \subset W$ are symplectic $F$-vector spaces and $I$ is a finite set.

\begin{proposition}\label{prop:GJ-Sp}
	\index{GJ}
	Up to conjugation, the subgroups $G_J \subset G$ attached to various $J \in \mathcal{J}^G_M$ are exactly of the form
	\[ \prod_{i \in \mathbf{I}} \Sp(W^i) \times \Sp(W^\natural) \]
	where
	\begin{compactitem}
		\item $\mathbf{I}$ is a finite set (possibly empty) equipped with a map $f: I \to \mathbf{I} \sqcup \{\natural\}$;
		\item $W^i$ and $W^\natural$ are symplectic subspaces of $W$ such that $W = \bigoplus_{i \in \mathbf{I}} W^i \oplus W^\natural$ (orthogonal sum) and $W^\flat \subset W^\natural$;
		\item $\dim W^i = 2 \sum_{f(k) = i} n_k$ for all $i \in \mathbf{I}$;
		\item $\dim W^\natural = 2 \sum_{f(k) = \natural} n_k + \dim W^\flat$.
	\end{compactitem}
	The embedding $M \hookrightarrow G_J$ is given by $\GL(n_k) \hookrightarrow \Sp(W^{f(k)})$ for $k \in I$ and $\Sp(W^\flat) \hookrightarrow \Sp(W^\natural)$.
\end{proposition}
\begin{proof}
	Since $\mathrm{char}(F) = 0$ and the descriptions of $R_J$, $G_J$, $\mathcal{J}^G_M$ are insensitive to field extensions, it suffices to describe the corresponding Lie subalgebras $\mathfrak{g}_J$ over $\overline{F}$. Following Dynkin, a subalgebra $\mathfrak{f} \subset \mathfrak{g}$ is called \emph{regular} if $[\mathfrak{h}, \mathfrak{f}] \subset \mathfrak{f}$ for some Cartan subalgebra $\mathfrak{h} \subset \mathfrak{g}$. We claim that $\mathfrak{g}_J$ must be a semisimple regular subalgebra of $\mathfrak{g}$. Indeed, $\mathfrak{g}_J$ is regular since it contains $\mathfrak{m}$, and semisimplicity follows from Proposition \ref{prop:GJ-ss-rank}.
	
	Now we can appeal to Dynkin's classification of semisimple regular subalgebras of $\mathfrak{g}$; see \cite[p.383]{Dyn52} or \cite[Chapter 6]{Vi94}: $\mathfrak{g}_J$ must be a product of simple Lie algebras of type $\mathrm{C}_m$, for various $m$, and they arise from an orthogonal decomposition of $W$ as indicated. The other conditions are imposed in order to make $M$ a Levi subgroup of $G_J$.
	
	To show that all subgroups in the assertion	arise as some $G_J$, assume that the data $\mathbf{I}$, $f$, $W^i$ and $W^\natural$ are given. Take the standard $\Z$-basis $\{\epsilon_i : i \in I \}$ for $X^*(A_M)$, put a total order on $I$, and let
	\[ \underline{J} := \left\{ \epsilon_i - \epsilon_j : i < j, \; f(i) = f(j) \right\} \sqcup \left\{ 2\epsilon_i : i \in I \right\}. \]
	This set is linearly independent in $\mathfrak{a}^{G, *}_M$, and the lattice generated by $\underline{J}$ in $X^*(A_M)$ equals
	\[ \left\{ \sum_{k \in I} x_k \epsilon_k : \forall i \in \mathbf{I} \sqcup \{\natural\},\; \sum_{f(k) = i} x_k \in 2\Z \right\}. \]
	Take $J \in \mathcal{J}^G_M$ to be the equivalence class generated by $\underline{J}$. One verifies that $G_J \simeq \prod_{i \in \mathbf{I}} \Sp(W^i) \times \Sp(W^\natural)$, compatibly with the given embedding of $M$.
\end{proof}

Next, consider $G^! := \SO(2n+1)$ with Levi subgroup $M^! := \prod_{i \in I} \GL(n_i) \times \SO(2n^\flat + 1)$, recalling that they are all split. Equip $G^!$ with the $B$-function specified in Definition \ref{def:B-metaplectic} (with $n' = n$, $n'' = 0$ and $\epsilon = 1$).

\begin{proposition}\label{prop:GJ-SO}
	Up to isomorphisms, the groups $G^!_{J^!}$ attached to various $J^! \in \mathcal{J}^{G^!}_{M^!}(B)$ are exactly of the form
	\[ \prod_{i \in \mathbf{I}} \SO(2n^i + 1) \times \SO(2n^\natural + 1) \]
	where
	\begin{compactitem}
		\item $\mathbf{I}$ is a finite set (possibly empty) equipped with a map $f: I \to \mathbf{I} \sqcup \{\natural\}$;
		\item the $n^i$ and $n^\natural$ are in $\Z_{\geq 0}$ such that $n = \sum_{i \in \mathbf{I}} n^i + n^\natural$;
		\item $n^i = 2 \sum_{f(k) = i} n_k$ for all $i \in \mathbf{I}$;
		\item $n^\natural = 2 \sum_{f(k) = \natural} n_k + n^\flat$.
	\end{compactitem}
	The embedding $M^! \hookrightarrow G^!_{J^!}$ is given by $\GL(n_i) \hookrightarrow \SO\left(2n^{f(i)} + 1\right)$ for $i \in I$ and $\SO(2n^\flat + 1) \hookrightarrow \SO(2n^\natural + 1)$.
\end{proposition}
\begin{proof}
	All these groups are split. In view of the remarks in \cite[p.283]{MW16-1}, the problem can be translated to the dual side
	\[ {M^!}^\vee = \prod_{i \in I} \GL(n_i, \CC) \times \Sp(2n^\flat, \CC) \hookrightarrow {G^!}^\vee = \Sp(2n, \CC) \]
	endowed with the trivial $B$-function. Proposition \ref{prop:GJ-Sp} can thus be applied over $\CC$ to give a complete description of $(G^!_{J^!})^\vee$. By duality, the case for $G^!_{J^!}$ follows at once.
\end{proof}

The evident duality between $G_J$ (type $\mathrm{C}$) and $G^!_{J^!}$ (type $\mathrm{B}$) is compatible with the identification between $\mathcal{J}^G_M$ and $\mathcal{J}^{G^!}_{M^!}(B)$ which intervenes in nonstandard endoscopy. See \S\ref{sec:nonstandard-endoscopy}.

\section{The maps \texorpdfstring{$\rho$}{rho}, \texorpdfstring{$\sigma$}{sigma} and \texorpdfstring{$\rho^{\Endo}$}{rhoE}}\label{sec:rho-sigma}
Let $F$ be a local field of characteristic zero.

\subsection{The spaces \texorpdfstring{$U_J$}{UJ}}
In what follows, $a$ will denote points of $A_M(F)$ which are close to $1$. Such elements can be canonically viewed as elements of $\widetilde{A_M}$, by using the two exponential maps in \eqref{eqn:two-exp}.

We will consider germs at $a=1$ of $\CC$-valued functions defined almost everywhere on $A_M(F)$ (or on $\widetilde{A_M}$ as explained above), i.e.\ off the $F$-points of some proper closed subvariety of $A_M$. Hereafter they will be abbreviated as ``germs''.

For each class $J \in \mathcal{J}^G_M$, define the $\CC$-vector space of germs of functions
\begin{equation}\label{eqn:UJ}\begin{aligned}
		U_J & = U_J^{G \supset M} \\
		& := \text{the span of}\; \left\{ a \mapsto \prod_{\alpha \in \underline{J}} \left| \alpha(a) - \alpha(a)^{-1} \right| : \underline{J} \in J \right\}.
\end{aligned}\end{equation}
\index{UJ@$U_J$}

When $M=G$ we have $\mathcal{J}^G_M = \{ \emptyset \}$ and $U_{\emptyset} = \CC$.

\begin{definition}\label{def:germ-equiv}
	Fix a function $a \mapsto d(a)$ for $a \in A_M(F)$ close to $1$ as in the discussions after Definition \ref{def:weighted-J-general}. Let $u$, $u'$ be germs. We say that $u \sim u'$ if there exists $r > 0$ such that over all domains in $A_M(F)$ of the form
	\[ \forall \alpha \in \Sigma(A_M), \; |\alpha(a) - 1| > c \cdot d(a), \]
	there exists $C > 0$ such that for almost all $a$ (in the sense explained above \eqref{eqn:UJ}) in that domain which are close to $1$,
	\[ \left| u'(a) - u(a) \right| < C \cdot d(a)^r . \]
\end{definition}

The following result is established in \cite[II.3.1 and V.2.3]{MW16-1}.
\begin{proposition}\label{prop:U-germ}
	Let $u$ be a germ.
	\begin{enumerate}[(i)]
		\item If $u \sim 0$ and $\alpha_1, \ldots, \alpha_n$ are elements in $\Sigma(A_M)$, then $u \cdot \prod_{i=1}^n \log\left| \alpha_i - \alpha_i^{-1} \right| \sim 0$.
		\item Suppose that $F$ is non-Archimedean. If $u \in \sum_{L \in \mathcal{L}(M)} \sum_{J \in \mathcal{J}^L_M} U_J$, then $u \sim 0$ if and only if $u = 0$.
		\item Suppose that $M \neq G$, $J \in \mathcal{J}^G_M$ and $u \in U_J$. If $u \sim c$ where $c \in \CC$, then $u=0$ and $c=0$. Consequently, if $u_1 + c_1 \sim u_2 + c_2$, where $c_i \in \CC$ and $u_i \in U_J$ for $i=1,2$, then $u_1 = u_2$ and $c_1 = c_2$.
	\end{enumerate}
\end{proposition}

Define the spaces of germs on $A_M(F)$:
\begin{align*}
	U^G_M & := \sum_{L \in \mathcal{L}(M)} \sum_{J \in \mathcal{J}^L_M} U_J^{L \supset M} , \\
	U^{G,+}_M & := \sum_{\substack{L \in \mathcal{L}(M) \\ L \neq M}} \sum_{J \in \mathcal{J}^L_M} U_J^{L \supset M}.
\end{align*}
Note the the summand with $L=M$ is $U_\emptyset = \CC$ (constant functions).
\index{UGM@$U^G_M$, $U^{G, +}_M$}

\begin{lemma}[{\cite[II.4.6 Lemme]{MW16-1}}]\label{prop:U-germ-5}
	Assume that $F$ is non-Archimedean with residual characteristic $p > 5$. Then the sum $U^G_M = \CC + U^{G, +}_M$ is direct.
\end{lemma}

This lemma allows us to define the constant term of any element of $U^G_M$ in the non-Archimedean case.

Now let $G^!$ be a quasisplit $F$-group endowed with a system of $B$-functions, and fix a Levi subgroup $M^! \subset G^!$. Let $\mathcal{O}^!$ be a stable semisimple conjugacy class in $M^!(F)$. For every $J \in \mathcal{J}^{G^!}_{M^!}(B_{\mathcal{O}^!})$, we iterate the same construction as \eqref{eqn:UJ} to define the linear subspace $U_J$ of the space of germs. Similarly, we define $U^{G^!}_{M^!}$ in this setting. It depends on the $B$-function and $\mathcal{O}^!$.

\subsection{The coverings \texorpdfstring{$\widetilde{G_J}$}{GJ} and annihilators}
Consider a covering $\rev: \tilde{G} \twoheadrightarrow G(F)$ in general. Let $M \subset G$ be a Levi subgroup.

\begin{definition}\label{def:GJ-Ann}
	\index{GJ-tilde@$\widetilde{G_J}$}
	Given $J \in \mathcal{J}^G_M$, let
	\begin{align*}
		\mathrm{Ann}^{\tilde{G}}_{\tilde{M}} & := \Ker\left[ D_{\mathrm{geom}, -}(\tilde{M}) \otimes \mes(M)^\vee \xrightarrow{\Ind^{\tilde{G}}_{\tilde{M}}} D_{\mathrm{geom}, -}(\tilde{G}) \otimes \mes(G)^\vee \right], \\
		\widetilde{G_J} & := \rev^{-1}\left( G_J(F) \right).
	\end{align*}
\end{definition}

\begin{lemma}\label{prop:Ann-J}
	We have $\mathrm{Ann}^{\widetilde{G_J}}_{\tilde{M}} \subset \mathrm{Ann}^{\tilde{G}}_{\tilde{M}}$.
\end{lemma}
\begin{proof}
	The arguments in \cite[p.279]{MW16-1} carry over. It boils down to the easy fact that $W^{G_J}(M) \subset W^G(M)$ operates on the trace Paley--Wiener spaces of $\tilde{M}$, $\tilde{G}$ and $\widetilde{G_J}$. Note that in the Archimedean case, we shall work with the $C^\infty_c$-version $\mathrm{PW}_{\asp}(\tilde{G})^\infty$ (Remark \ref{rem:real-PW}).
\end{proof}

Let $\mathcal{O}$ be a finite union of semisimple conjugacy classes in $M(F)$. Define
\[ \mathrm{Ann}^{\tilde{G}}_{\mathcal{O}} := \mathrm{Ann}^{\tilde{G}}_{\tilde{M}} \cap \left( D_{\mathrm{geom}, -}(\tilde{M}, \mathcal{O}) \otimes \mes(M)^\vee \right). \]
\index{AnnGO@$\mathrm{Ann}^{\tilde{G}}_{\mathcal{O}}$, $S\mathrm{Ann}^{G^{"!}}_{\mathcal{O}^{"!}}$}

The uncovered case $\tilde{G} = G(F)$ has been discussed in \textit{loc.\ cit.} In that case, we also write
\[ \mathrm{Ann}^G_{\mathrm{unip}} := \mathrm{Ann}^G_{\{1\}}. \]

Similarly, for a quasisplit group $G^!$ with Levi subgroup $M^!$, we define the stable avatars
\[ S\mathrm{Ann}^{G^!}_{M^!}, \quad S\mathrm{Ann}^{G^!}_{\mathcal{O}^!}, \quad S\mathrm{Ann}^{G^!}_{\mathrm{unip}} := S\mathrm{Ann}^{G^!}_{\{1\}}, \]
where $\mathcal{O}^!$ is any finite union of stable semisimple conjugacy classes in $M^!(F)$. It suffices to replace $D_{\mathrm{geom}, -}(\tilde{M})$ by $SD_{\mathrm{geom}}(M^!)$, and so forth. Clearly, $S\mathrm{Ann}^{G^!}_{\mathcal{O}^!} \subset \mathrm{Ann}^{G^!}_{\mathcal{O}^!}$.

\subsection{Non-Archimedean case}
Assume furthermore that $F$ is non-Archimedean. Let $\rev: \tilde{G} \twoheadrightarrow G(F)$ be a covering, $M \subset G$ a Levi subgroup. Fix a semisimple conjugacy class $\mathcal{O}$ in $M(F)$.

In fact, the coverings we need will either be of metaplectic type, or of the form $\widetilde{G_J}$ (Definition \ref{def:GJ-Ann}) where $\tilde{G}$ is of metaplectic type, hence the situation is not too different from the essential case $\tilde{G} = \Mp(W)$.

\begin{definition-proposition}\label{def:rho-germ}
	\index{rhoGJ@$\rho^{\tilde{G}}_J$}
	For all $J \in \mathcal{J}^G_M$, there exists a unique linear map
	\[ \rho^{\tilde{G}}_J: D_{\mathrm{geom}, -}(\tilde{M}, \mathcal{O}) \otimes \mes(M)^\vee \to U_J \otimes \left( D_{\mathrm{geom}, -}(\tilde{M}, \mathcal{O}) \otimes \mes(M)^\vee \right) \big/ \mathrm{Ann}^{\tilde{G}}_{\mathcal{O}} \]
	such that
	\begin{enumerate}[(i)]
		\item it factorizes into
		\[ D_{\mathrm{geom}, -}(\tilde{M}, \mathcal{O}) \otimes \mes(M)^\vee \xrightarrow{\rho^{\widetilde{G_J}}_J} U_J \otimes (\cdots) \big/ \mathrm{Ann}^{\widetilde{G_J}}_{\mathcal{O}} \twoheadrightarrow U_J \otimes (\cdots) \big/ \mathrm{Ann}^{\tilde{G}}_{\mathcal{O}} \]
		where we used Lemma \ref{prop:Ann-J};
		\item let $\tilde{\gamma} \in D_{\mathrm{geom}, -}(\tilde{M}, \mathcal{O}) \otimes \mes(M)^\vee$ and $f \in \orbI_{\asp}(\tilde{G}) \otimes \mes(G)$, then
		\[ I_{\tilde{M}}\left(a\tilde{\gamma}, f \right) \sim \sum_{L \in \mathcal{L}(M)} \sum_{J \in \mathcal{J}^L_M} I_{\tilde{L}}\left( \rho^{\tilde{L}}_J(\tilde{\gamma}, a)^{\tilde{L}}, f \right) \]
		as germs in $a \in A_M(F)$.
	\end{enumerate}
	Here, $a\tilde{\gamma}$ makes sense since $a$ is naturally lifted to $\rev^{-1}\left( A_M(F) \right)$ for $a \to 1$, and $(\cdots)^{\tilde{L}}$ stands for the map
	\[ \identity \otimes \Ind^{\tilde{L}}_{\tilde{M}}: U_J \otimes \left( D_{\mathrm{geom}, -}(\tilde{M}) \otimes \mes(M)^\vee \right) \big/ \mathrm{Ann}^{\tilde{L}}_{\tilde{M}} \to U_J \otimes D_{\mathrm{geom}, -}(\tilde{L}) \otimes \mes(L)^\vee \]
	where $U_J = U_J^{L \supset M}$.
\end{definition-proposition}
\begin{proof}
	Identical to \cite[II.3.2 Proposition]{MW16-1}. It boils down to the following ingredients:
	\begin{compactitem}
		\item generalities on the $(G, M)$-family defining weighted orbital integrals,
		\item standard properties of $I_{\tilde{M}}(\tilde{\gamma}, \cdots)$, see Proposition \ref{prop:orbint-weighted-descent-nonArch};
		\item standard properties of germs, see Proposition \ref{prop:U-germ}.
	\end{compactitem}
	In particular, the coverings do not really intervene in the arguments. For a more explicit formula for $\rho^{\tilde{G}}_J$, see \cite[II.3.2 (5)]{MW16-1}.
\end{proof}

Note that our conventions imply that $\rho^{\tilde{G}}_{\emptyset} = \identity$ when $M=G$.

\begin{lemma}\label{prop:rho-equisingular-nonarch}
	Let $\tilde{\gamma} \in D_{\mathrm{geom}, -}(\tilde{M}, \mathcal{O}) \otimes \mes(M)^\vee$ where $\mathcal{O}$ is $G$-equisingular. Then $\rho^{\tilde{G}}_J(\tilde{\gamma}) = 0$ unless $M=G$.
\end{lemma}
\begin{proof}
	One uses the explicit construction \cite[p.282, (5)]{MW16-1} of $\rho^{\tilde{G}}_J(\tilde{\gamma})$, which applies to coverings as well. Suppose $M \neq G$. The terms $m(\underline{\alpha}, \tilde{\gamma})$ therein are all zero because so are the terms  $\rho(\alpha, \tilde{\gamma})$, which is in turn a consequence of the definition in \cite[p.194]{MW16-1} of $\rho(\alpha, \tilde{\gamma})$ since $G_{\gamma_{\mathrm{ss}}} = M_{\gamma_{\mathrm{ss}}}$.
\end{proof}

Before stating the next descent formula, we collect some facts from \cite[p.297]{MW16-1}. Let $L, L^\dagger \in \mathcal{L}^G(M)$ be such that $d^G_M(L, L^\dagger) \neq 0$. Then:
\begin{compactitem}
	\item the inclusion $A_L \subset A_M$ induces an injection $\Sigma^{L^\dagger}(A_M) \hookrightarrow \Sigma^G(A_L)$;
	\item in turn, this induces an injection $\mathcal{J}^{L^\dagger}_M \hookrightarrow \mathcal{J}^G_L$;
	\item for each $J \in \mathcal{J}^{L^\dagger}_M$, also viewed as an element in $\mathcal{J}^G_L$, we have
	\[ U_J^{G \supset L} = U_J^{L^\dagger \supset M} \big|_{A_L(F)} \]
	as spaces of germs over $A_L(F)$.
\end{compactitem}

On the other hand, the parabolic induction $\Ind^{\tilde{L}}_{\tilde{M}}$ of \eqref{eqn:parabolic-ind-dist} always maps $\mathrm{Ann}^{\tilde{L}^\dagger}_{\tilde{M}}$ to $\mathrm{Ann}^{\tilde{G}}_{\tilde{L}}$. Therefore, when $d^G_M(L, L^\dagger) \neq 0$, we obtain the linear map
\[ U_J^{L^\dagger \supset M} \otimes \left( D_{\mathrm{geom}, -}(\tilde{M}) \otimes \mes(M)^\vee \right) \big/ \mathrm{Ann}^{\tilde{L}^\dagger}_{\tilde{M}} \to U_J^{G \supset L} \otimes \left( D_{\mathrm{geom}, -}(\tilde{L}) \otimes \mes(L)^\vee \right) \big/ \mathrm{Ann}^{\tilde{G}}_{\tilde{L}}, \]
which will again be denoted by $\rho \mapsto \rho^{\tilde{L}}$.

\begin{proposition}\label{prop:rho-descent}
	Let $L \in \mathcal{L}(M)$, $J \in \mathcal{J}^G_L$ and $\tilde{\gamma} \in D_{\mathrm{geom}, -}(\tilde{M}) \otimes \mes(M)^\vee$. Then
	\[ \rho^{\tilde{G}}_J \left( \tilde{\gamma}^{\tilde{L}}\right) = \sum_{\substack{L^\dagger \in \mathcal{L}^G(M) \\ \text{s.t.}\; J \in \mathcal{J}^{L^\dagger}_M}} d^G_M(L, L^\dagger) \rho^{\tilde{L}^\dagger}_J \left( \tilde{\gamma} \right)^{\tilde{L}} \]
	as germs on $A_M(F)$.
\end{proposition}
\begin{proof}
	This is the analogue of \cite[II.3.10 Lemme]{MW16-1} for coverings, and the proof in \textit{loc.\ cit.} carries over, cf.\ the explanations for Definition--Proposition \ref{def:rho-germ}.
\end{proof}

Consider now a quasisplit $F$-group $G^!$ endowed with a system of $B$-functions and a Levi subgroup $M^! \subset G^!$. Let $\mathcal{O}^!$ be a stable semisimple conjugacy class in $M^!(F)$. The stable avatar of Definition--Proposition \ref{def:rho-germ} is given in \cite[II.3.5 Proposition]{MW16-1}: for every $J \in \mathcal{J}^{G^!}_{M^!}(B_{\mathcal{O}^!})$ we have a canonical linear map
\[ \sigma^{G^!}_J : SD_{\mathrm{geom}}(M^!, \mathcal{O}^!) \otimes \mes(M^!)^\vee \to U_J \otimes \left( SD_{\mathrm{geom}}(M^!, \mathcal{O}^!) \otimes \mes(M)^\vee \right) \big/ S\mathrm{Ann}^{G^!}_{\mathcal{O}^!}. \]
Note that it involves the system of $B$-functions.
\index{sigmaGJ@$\sigma^{G^{"!}}_J$}

We are now ready to define the endoscopic avatar of $\rho^{\tilde{G}}_J$. Assume henceforth that $\tilde{G}$ is a group of metaplectic type, $M \subset G$ is a Levi subgroup, and $\mathbf{M}^! \in \Endo_{\elli}(\tilde{M})$. Let $\mathcal{O}^!$ be a stable semisimple conjugacy class in $M^!(F)$. Denote the identification $A_M \rightiso A_{M^!}$ as $\xi$. Finally, let $\mathcal{O} \subset M(F)$ be the union of semisimple conjugacy classes corresponding to $\mathcal{O}^!$.

\begin{definition}\label{def:trans-s}
	\index{Trans-s@$\Trans_{\mathbf{M}^{"!}, \tilde{M}}^s$}
	For every $s \in \Endo_{\mathbf{M}^!}(\tilde{G})$, set
	\begin{align*}
		\Trans_{\mathbf{M}^!, \tilde{M}}^s: \orbI_{\asp}(\tilde{M}) \otimes \mes(M) & \longrightarrow S\orbI(M^!) \otimes \mes(M^!) \\
		f & \longmapsto \left( \Trans_{\mathbf{M}^!, \tilde{M}}(f) \right)[s]^{-1}
	\end{align*}
	where $g^{M^!} \mapsto g^{M^!}[s]^{-1}$ is the inverse of the automorphism $g^{M^!} \mapsto g^{M^!}[s]$ in Definition \ref{def:g-s}. Denote its transpose by
	\[ \trans_{\mathbf{M}^!, \tilde{M}}^s: SD_{\mathrm{geom}}(M^!) \otimes \mes(M^!)^\vee \longrightarrow D_{\mathrm{geom}, -}(\tilde{M}) \otimes \mes(M)^\vee. \]
\end{definition}

Set $\mathcal{O}^![s] := \mathcal{O}^! z[s]$. It is routine to verify that $\trans_{\mathbf{M}^!, \tilde{M}}^s$ restricts to
\[ SD_{\mathrm{geom}}(M^!, \mathcal{O}^![s]) \otimes \mes(M)^\vee \to D_{\mathrm{geom}, -}(\tilde{M}, \mathcal{O}) \otimes \mes(M)^\vee . \]

\begin{lemma}
	For every $s \in \Endo_{\mathbf{M}^!}(\tilde{G})$, the map $\trans_{\mathbf{M}^!, \tilde{M}}^s$ maps $S\mathrm{Ann}^{G^![s]}_{M^!}$ (resp.\ $S\mathrm{Ann}^{G^![s]}_{\mathcal{O}^![s]}$) to $\mathrm{Ann}^{\tilde{G}}_{\tilde{M}}$ (resp.\ $\mathrm{Ann}^{\tilde{G}}_{\mathcal{O}}$).
\end{lemma}
\begin{proof}
	Proposition \ref{prop:Levi-central-twist} asserts the commutativity of
	\[\begin{tikzcd}
		\orbI_{\asp}(\tilde{M}) \otimes \mes(M) \arrow[r, "{\Trans_{\mathbf{M}^!, \tilde{M}}^s}"] & S\orbI(M^!) \otimes \mes(M^!) \\
		\orbI_{\asp}(\tilde{G}) \otimes \mes(G) \arrow[r, "{\Trans_{\mathbf{G}^![s], \tilde{G}}}"'] \arrow[u] & S\orbI(G^![s]) \otimes \mes(G^![s]) \arrow[u]
	\end{tikzcd}\]
	the vertical arrows being parabolic descent. Taking transpose yields $\trans_{\mathbf{M}^!, \tilde{M}}^s\left( S\mathrm{Ann}^{G^![s]}_{M^!} \right) \subset \mathrm{Ann}^{\tilde{G}}_{\tilde{M}}$. Restriction to $SD_{\mathrm{geom}}(M^!, \mathcal{O}^![s]) \otimes \mes(M^!)$ yields the case with $\mathcal{O}^![s]$ and $\mathcal{O}$.
\end{proof}

Observe that if $J^! \mapsto J$ as in Lemma \ref{prop:J-transfer}, then $\xi^*: u \mapsto u\xi$ induces an injective linear map $U_{J^!} \to U_J$.

Hence for each $J^! \in \mathcal{J}^{G^![s]}_{M^!}(B^{\tilde{G}}_{\mathcal{O}^![s]})$ and $\delta \in SD_{\mathrm{geom}}(M^!, \mathcal{O}^!) \otimes \mes(M^!)^\vee$, we can form
\[ \trans_{\mathbf{M}^!, \tilde{M}}^s \left( \sigma^{G^![s]}_{J^!}(\delta[s])\right) \in  U_J \otimes \left( D_{\mathrm{geom},-}(\tilde{M}, \mathcal{O}) \otimes \mes(M)^\vee \right) \big/ \mathrm{Ann}^{\tilde{G}}_{\mathcal{O}}. \]
Unraveling definitions, it involves two twists inside $SD_{\mathrm{geom}}(M^!)$: one by $z[s]$ (for $\delta \mapsto \delta[s]$) and one by $z[s]^{-1}$ (for $\trans_{\mathbf{M}^!, \tilde{M}}^s$). Due to the presence of $S\mathrm{Ann}^{G^![s]}_{M^!}$ and $\mathrm{Ann}^{\tilde{G}}_{\tilde{M}}$, they do not cancel.

Observe that the constructions above also work for Archimedean $F$.

\begin{definition}
	\index{rhoGEndoM@$\rho^{\tilde{G}, \Endo}_J$}
	For each $\delta \in SD_{\mathrm{geom}}(M^!, \mathcal{O}^!) \otimes \mes(M^!)^\vee$ and $J \in \mathcal{J}^G_M$, set
	\begin{equation*}
		\rho^{\tilde{G}, \Endo}_J(\mathbf{M}^!, \delta, a) := \sum_{s \in \Endo_{\mathbf{M}^!}(\tilde{G})} i_{M^!}(\tilde{G}, G^![s])
		\sum_{\substack{J^! \in \mathcal{J}^{G^![s]}_{M^!}(B^{\tilde{G}}_{\mathcal{O}^![s]}) \\ J^! \mapsto J }} \trans_{\mathbf{M}^!, \tilde{M}}^s \left( \sigma^{G^![s]}_{J^!} \left( \delta[s], \xi(a) \right)\right)
	\end{equation*}
	where $a \in A_M(F)$ is in general position and close to $1$. Viewing $a$ as a variable, we obtain an element
	\[ \rho^{\tilde{G}, \Endo}_J(\mathbf{M}^!, \delta) \in U_J \otimes \left( D_{\mathrm{geom},-}(\tilde{M}, \mathcal{O}) \otimes \mes(M)^\vee \right) \big/ \mathrm{Ann}^{\tilde{G}}_{\mathcal{O}}. \]
\end{definition}

\begin{proposition}\label{prop:rhoEndo-expansion}
	For all $(\mathbf{M}^!, \delta)$ above, we have
	\[ I^{\tilde{G}, \Endo}_{\tilde{M}}\left( \mathbf{M}^!, \xi(a)\delta, f \right) \sim \sum_{L \in \mathcal{L}(M)} \sum_{J \in \mathcal{J}^L_M} I^{\tilde{G}, \Endo}_{\tilde{L}}\left( \rho^{\tilde{L}, \Endo}_J(\mathbf{M}^!, \delta, a)^{\tilde{L}}, f \right) \]
	as germs in $a \in A_M(F)$, $a \to 1$.
\end{proposition}
\begin{proof}
	Using \cite[II.3.7 Proposition]{MW16-1}, we have
	\begin{align*}
		I^{\tilde{G}, \Endo}_{\tilde{M}}\left( \mathbf{M}^!, \xi(a)\delta, f \right) & = \sum_{s \in \Endo_{\mathbf{M}^!}(\tilde{G})} i_{M^!}(\tilde{G}, G^![s]) S^{G^![s]}_{M^!}\left( \xi(a)\delta[s], B^{\tilde{G}}, f^{G^![s]} \right) \\
		& \sim \sum_{s \in \Endo_{\mathbf{M}^!}(\tilde{G})} i_{M^!}(\tilde{G}, G^![s]) \sum_{L^! \in \mathcal{L}^{G^![s]}(M^!)} \sum_{J^!} S^{G^![s]}_{L^!}\left( \sigma^{L^!}_{J^!}(\delta[s], \xi a)^{L^!}, B^{\tilde{G}}, f^{G^![s]} \right).
	\end{align*}
	Apply Lemma \ref{prop:sL-Ls} with $R = M$ to transform $\sum_{s, L^!}$ into $\sum_{L, s}$ where $L \in \mathcal{L}^G(M)$ and $s \in \Endo_{\mathbf{M}^!}(\tilde{G})$ factors uniquely through a pair $(s^L, s_L)$ with
	\[ s^L \in \Endo_{\mathbf{M}^!}(\tilde{L}), \quad s_L \in \Endo_{\mathbf{L}^![s^L]}(\tilde{G}), \]
	so that $L^!$ (resp.\ $\mathbf{G}^![s]$) corresponds to $L^![s^L]$ (resp.\ $\mathbf{G}^![s_L]$).

	Lemma \ref{prop:i-transitivity} asserts $i_{M^!}(\tilde{G}, G^![s]) = i_{M^!}(\tilde{L}, L^!) i_{L^!}(\tilde{G}, G^![s_L])$. We also have $\delta[s] = \delta[s^L][s_L]$ and $z[s_L]$ is central in $L^!(F)$. The last expression becomes
	\begin{multline*}
		\sum_{\substack{L \in \mathcal{L}^G(M) \\ s^L \in \Endo_{\mathbf{M}^!}(\tilde{L}) \\ L^! := L^![s^L] }} i_{M^!}(\tilde{L}, L^!) \sum_{\substack{s_L \in \Endo_{\mathbf{L}^!}(\tilde{G}) \\ G^! := G^![s_L] }} i_{L^!}(\tilde{G}, G^!) \sum_{J^!}
		S^{G^!}_{L^!} \left( \sigma^{L^!}_{J^!}(\delta[s^L], \xi a)^{L^!} [s_L], B^{\tilde{G}}, f^{G^!} \right) \\
		= \sum_{L \in \mathcal{L}^G(M)} \sum_{J \in \mathcal{J}^L_M} \sum_{\substack{s^L \in \Endo_{\mathbf{M^!}(\tilde{L})} \\ L^! := L^![s^L] }} i_{M^!}(\tilde{L}, L^!) I^{\tilde{G}, \Endo}_{\tilde{L}} \left( \mathbf{L}^!, \sum_{J^! \mapsto J} \sigma^{L^!}_{J^!}(\delta[s^L], \xi a)^{L^!}, f \right).
	\end{multline*}
	It remains to recall that $\trans_{\mathbf{L}^!, \tilde{L}} \left(\sigma^{L^!}_{J^!}(\delta[s^L], \xi a)^{L^!} \right)$ equals $\trans_{\mathbf{M}^!, \tilde{M}}^{s^L}\left(\sigma^{L^!}_{J^!}(\delta[s^L], \xi a) \right)^{\tilde{L}}$.
\end{proof}

\begin{theorem}\label{prop:rho-main-nonArch}
	Recall that $F$ is non-Archimedean here. For every $\delta \in SD_{\mathrm{geom}}(M^!, \mathcal{O}^!)$ and $J \in \mathcal{J}^G_M$, we have
	\[ \rho^{\tilde{G}}_J\left( \trans_{\mathbf{M}^!, \tilde{M}}(\delta) \right) = \rho^{\tilde{G}, \Endo}_J\left( \mathbf{M}^!, \delta \right) \]
	as germs on $A_M(F)$.
\end{theorem}

A conditional proof of Theorem \ref{prop:rho-main-nonArch} will be given in \S\ref{sec:rho-main-summary}, under the following
\begin{hypothesis}\label{hyp:local-geometric-nonArch}
	For all $\tilde{\gamma} \in D_{\mathrm{geom}, -}(\tilde{M}) \otimes \mes(M)^\vee$ with $\Supp(\tilde{\gamma}) \subset \tilde{M}_{G\text{-reg}}$, and all $f \in \orbI_{\asp}(\tilde{G}) \otimes \mes(G)$, we have
	\[ I^{\tilde{G}}_{\tilde{M}}(\tilde{\gamma}, f) = I^{\tilde{G}, \Endo}_{\tilde{M}}(\tilde{\gamma}, f). \]
\end{hypothesis}

The hypothesis will be resolved in \S\ref{sec:end-stabilization}.

\subsection{Archimedean case}\label{sec:rho-sigma-Arch}
Assume that $F$ is an Archimedean local field. In this case $D_{\mathrm{orb}, -}(\tilde{M})$ and $D_{\mathrm{geom}, -}(\tilde{M})$ are no longer equal. In \S\ref{sec:orbint-weighted-Arch} we have defined $I^{\tilde{G}}_{\tilde{M}}(\tilde{\gamma}, \cdot)$ for all $\tilde{\gamma} \in D_{\mathrm{orb}, -}(\tilde{M}) \otimes \mes(M)^\vee$. Let $J \in \mathcal{J}^G_M$. The same construction as in the non-Archimedean case yields a linear map
\[ \rho^{\tilde{G}}_J: D_{\mathrm{orb}, -}(\tilde{M}) \otimes \mes(M)^\vee \to U_J \otimes \left( D_{\mathrm{orb}, -}(\tilde{M}, \mathcal{O}) \otimes \mes(M)^\vee \right) \big/ \mathrm{Ann}^{\tilde{G}}_{\mathcal{O}} \]
with the same properties as in Definition--Proposition \ref{def:rho-germ}. There is also a counterpart of Lemma \ref{prop:rho-equisingular-nonarch}.
\index{rhoGJ}

\begin{lemma}\label{prop:rho-equisingular-arch}
	Let $\tilde{\gamma} \in D_{\mathrm{orb}, -}(\tilde{M}, \mathcal{O}) \otimes \mes(M)^\vee$ where $\mathcal{O}$ is $G$-equisingular. Then $\rho^{\tilde{G}}_J(\tilde{\gamma}) = 0$ unless $M=G$.
\end{lemma}

Nevertheless, the constructions in the stable and endoscopic contexts are more tricky. The reason is that we have to allow invariant distributions which are not necessarily orbital. We will extend these constructions to a broader class of $\tilde{\gamma}$ and formulate the Archimedean analogue of Theorem \ref{prop:rho-main-nonArch}; besides, we need to postulate an Archimedean analogue of Hypothesis \ref{hyp:local-geometric-nonArch}. This is the topic of \S\ref{sec:ext-Arch}.

\section{Extension of definitions in the Archimedean case}\label{sec:ext-Arch}
In this section, $F$ is an Archimedean local field.

\subsection{The space \texorpdfstring{$D_{\text{tr-orb}}$}{Dtrorb}}
Let $\rev: \tilde{G} \to G(F)$ be of metaplectic type. Following \cite[V.2]{MW16-1}, we are going to define a subspace $D_{\text{tr-orb}, -}(\tilde{G})$ of $D_{\text{geom}, -}(\tilde{G})$ in the Archimedean setting.

Let us fix Haar measures temporarily to trivialize the lines $\mes(\cdots)$; they will not affect the following definitions.

First off, recall from \cite[V.2.1]{MW16-1} that for every quasisplit $F$-group $L$, one defines inductively a vector subspace
\[ D_{\text{tr-orb}}(L) \subset D_{\mathrm{geom}}(L) \]
generated by $D_{\mathrm{orb}}(L)$ together with the images of
\[ SD_{\text{tr-orb}}(L^!) := D_{\text{tr-orb}}(L^!) \cap SD_{\mathrm{geom}}(L^!) \]
under endoscopic transfer, where $\mathbf{L}^!$ ranges over $\Endo_{\elli}(L)$ with $\dim L^! < \dim L$, modulo some technical details from $z$-extensions (because the transfer in general will involve $z$-extensions on the endoscopic side, see \textit{loc.\ cit.}). If $\mathcal{O}^!$ is a finite union of stable semisimple conjugacy classes in $L(F)$, we set
\[ SD_{\text{tr-orb}}(L, \mathcal{O}^!) := SD_{\text{tr-orb}}(L) \cap SD_{\text{geom}}(L, \mathcal{O}^!). \]
By taking $\mathcal{O}^! = \{1\}$, we obtain the space $SD_{\text{tr-unip}}(L)$.
\index{SDtr-orb@$SD_{\text{tr-orb}}$}

\begin{definition}\label{def:tr-orb}
	\index{Dtr-orb@$D_{\text{tr-orb}, -}(\tilde{G})$}
	Let
	\[ D_{\text{tr-orb}, -}(\tilde{G}) := D_{\mathrm{orb}, -}(\tilde{G}) + \sum_{\mathbf{G}^! \in \Endo_{\elli}(\tilde{G})} \trans_{\mathbf{G}^!, \tilde{G}} \left( SD_{\text{tr-orb}}(G^!) \right). \]
	It is a vector subspace of $D_{\mathrm{geom}, -}(\tilde{G})$. For $\trans_{\mathbf{G}^!, \tilde{G}}$, recall Definition \ref{def:transfer-dist}.
\end{definition}

If $\mathcal{O}$ is a finite union of semisimple conjugacy classes in $G(F)$, we write
\begin{align*}
	D_{\text{tr-orb}, -}(\tilde{G}, \mathcal{O}) & := D_{\text{tr-orb}, -}(\tilde{G}) \cap D_{\text{geom}, -}(\tilde{G}, \mathcal{O}), \\
	D_{\text{tr-unip}, -}(\tilde{G}) & := D_{\text{tr-orb}, -}(\tilde{G}, \{1\}).
\end{align*}

Hereafter, we relax the choice of Haar measures.

\begin{proposition}\label{prop:tr-orb-induction}
	Let $M \subset G$ be a Levi subgroup. Then $\Ind^{\tilde{G}}_{\tilde{M}}$ maps $D_{\text{tr-orb}, -}(\tilde{M}) \otimes \mes(M)^\vee$ to $D_{\text{tr-orb}, -}(\tilde{G}) \otimes \mes(G)^\vee$.
\end{proposition}
\begin{proof}
	We may assume $\tilde{G} = \Mp(W)$. The spaces $SD_{\text{tr-orb}}(L) \otimes \mes(L)^\vee$ are known to be stable under parabolic induction. By considering the situation \eqref{eqn:s-situation} for various $\mathbf{M}^!$ and $s \in \Endo_{\mathbf{M}^!}(\tilde{G})$, then applying Proposition \ref{prop:Levi-central-twist}, it suffices to show that $\delta \mapsto \delta \cdot z[s]$ induces an automorphism of $SD_{\text{tr-orb}}(M^!) \otimes \mes(M^!)^\vee$. This is \cite[V.2.1 (7)]{MW16-1}.
\end{proof}

The shorthand $\tilde{\gamma}^{\tilde{G}} := \Ind^{\tilde{G}}_{\tilde{M}}(\tilde{\gamma})$ will be used.

\subsection{A program}\label{sec:ext-Arch-prog}
This subsection is a variant of the situation (A) in \cite[V.2.4]{MW16-1}. Consider a group of metaplectic type $\rev: \tilde{G} \to G(F)$ and a Levi subgroup $M \subset G$. Our aim is to define the following objects.

\begin{enumerate}[(O1)]
	\item For each $J \in \mathcal{J}^G_M$ and a semisimple conjugacy class $\mathcal{O} \subset M(F)$, a linear map
	\[ \rho^{\tilde{G}}_J: D_{\text{tr-orb}, -}(\tilde{M}, \mathcal{O}) \otimes \mes(M)^\vee \to U_J \otimes \left( D_{\mathrm{geom}, -}(\tilde{M}, \mathcal{O}) \otimes \mes(M)^\vee \right) / \mathrm{Ann}^{\tilde{G}}_{\mathcal{O}}. \]
	It extends by linearity to the case when $\mathcal{O}$ is a finite union of semisimple conjugacy classes.
	
	\item For each $\tilde{\gamma} \in D_{\text{tr-orb}, -}(\tilde{M}) \otimes \mes(M)^\vee$, a continuous linear functional
	\[ f \mapsto I^{\tilde{G}}_{\tilde{M}}(\tilde{\gamma}, f) \]
	where $f \in \orbI_{\asp}(\tilde{G}) \otimes \mes(G)$.

	\item Given
	\begin{compactitem}
		\item $\mathbf{M}^! \in \Endo_{\elli}(\tilde{M})$,
		\item a stable semisimple conjugacy class $\mathcal{O}^! \subset M^!(F)$,
		\item the image $\mathcal{O}$ of $\mathcal{O}^!$, which is a stable semisimple conjugacy class in $M(F)$,
		\item $\delta \in SD_{\text{tr-orb}}(M^!, \mathcal{O}^!) \otimes \mes(M^!)^\vee$,
	\end{compactitem}
	we associate the germ $\rho^{\tilde{G}, \Endo}_J(\mathbf{M}^!, \delta)$ in $a \in A_M(F)$, $a \to 1$ in general position:
	\begin{equation*}
		\rho^{\tilde{G}, \Endo}_J (\mathbf{M}^!, \delta, a) := \sum_{s \in \Endo_{\mathbf{M}^!}(\tilde{G})} i_{M^!}(\tilde{G}, G^![s])
		\sum_{\substack{J^! \in \mathcal{J}^{G^![s]}_{M^!}\left(B^{\tilde{G}}_{\mathcal{O}^![s]}\right) \\ J^! \mapsto J }} \trans_{\mathbf{M}^!, \tilde{M}}^s \left( \sigma^{G^![s]}_{J^!} \left( \delta[s], \xi(a) \right)\right).
	\end{equation*}
	It is an element of $U_J \otimes \left( D_{\mathrm{geom}, -}(\tilde{M}, \mathcal{O}) \otimes \mes(M)^\vee\right) / \mathrm{Ann}^{\tilde{G}}_{\mathcal{O}}$. Note that $\sigma^{G^![s]}_{J^!}\left( \cdots \right)$ is defined in the situation (B) of \cite[V.2.4]{MW16-1}, as $\delta[s] \in SD_{\text{tr-orb}}(M^!, \mathcal{O}^![s]) \otimes \mes(M^!)^\vee$ (see the proof of Proposition \ref{prop:tr-orb-induction}).

	\item For $\mathbf{M}^!$, $\delta$ as above,
	\[ I^{\tilde{G}, \Endo}_{\tilde{M}}(\mathbf{M}^!, \delta, f) = \sum_{s \in \Endo_{\mathbf{M}^!}(\tilde{G})} i_{M^!}(\tilde{G}, G^![s]) S^{G^![s]}_{M^!}\left( \delta[s], B^{\tilde{G}}, f^{G^![s]} \right) \]
	where $f^{G^![s]} := \Trans_{\mathbf{G}^![s], \tilde{G}}(f)$. Again, $S^{G^![s]}_{M^!}(\delta[s], B^{\tilde{G}}, \cdot)$ is defined in the situation (B) of \cite[V.2.4]{MW16-1}.
\end{enumerate}

Observe that the definitions in (O3) and (O4) are unconditional.

The objects postulated above must be subject to the following conditions.
\begin{enumerate}[(C1)]
	\item When $\tilde{\gamma} \in D_{\text{orb}, -}(\tilde{M}, \mathcal{O}) \otimes \mes(M)^\vee$, the $\rho^{\tilde{G}}_J(\tilde{\gamma})$ in (O1) (resp.\ $I^{\tilde{G}}_{\tilde{M}}(\tilde{\gamma}, f)$ in (O2)) coincides with the one defined in \S\ref{sec:rho-sigma-Arch} (resp.\ in \S\ref{sec:orbint-weighted-Arch}).
	
	\item For all $(\mathbf{M}^!, \delta)$ as in (O4) and all $f$, we have
	\[ \underbracket{I^{\tilde{G}}_{\tilde{M}}\left( \trans_{\mathbf{M}^!, \tilde{M}}(\delta), f \right)}_{\text{in (O2)}} = \underbracket{I^{\tilde{G}, \Endo}_{\tilde{M}}(\mathbf{M}^!, \delta, f)}_{\text{in (O4)}}. \]
	
	\item For all $(\mathbf{M}^!, \delta)$ as in (O3), we have the equality of germs
	\[ \underbracket{\rho^{\tilde{G}}_J\left( \trans_{\mathbf{M}^!, \tilde{M}}(\delta)\right)}_{\text{in (O1)}} = \underbracket{\rho^{\tilde{G}, \Endo}_J(\mathbf{M}^!, \delta)}_{\text{in (O3)}} . \]
	
	\item For each $\tilde{\gamma} \in D_{\text{tr-orb}, -}(\tilde{M}, \mathcal{O}) \otimes \mes(M)^\vee$, $J \in \mathcal{J}^G_M$ and $a \in A_M(F)$ in general position, $a \to 1$, we have
	\[ \rho^{\tilde{G}}_J(\tilde{\gamma}, a)^{\tilde{G}} \in D_{\text{tr-orb}, -}(\tilde{G}) \otimes \mes(G)^\vee. \]
	
	\item Let $\tilde{\gamma}$ be as above. The germ at $1$ of $a \mapsto I^{\tilde{G}}_{\tilde{M}}(a\tilde{\gamma}, f)$ is equivalent in the sense of \S\ref{sec:rho-sigma} to
	\[ a \mapsto \sum_{L \in \mathcal{L}(M)} \sum_{J \in \mathcal{J}^L_M} I^{\tilde{G}}_{\tilde{L}}\left( \rho^{\tilde{L}}_J(\tilde{\gamma}, a)^{\tilde{L}}, f \right), \]
	which is well-defined by (C4). Observe that for $a$ in general position, $a\tilde{\gamma}$ is $\tilde{G}$-equisingular and $I^{\tilde{G}}_{\tilde{M}}(a\tilde{\gamma}, f)$ should be understood by Definition \ref{def:weighted-I-integral-arch}; see (C8) below.

	\item For each $\tilde{\gamma} \in D_{\text{tr-orb}, -}(\tilde{M}, \mathcal{O}) \otimes \mes(M)^\vee$ and $L \in \mathcal{L}^G(M)$, we have
	\[ I^{\tilde{G}}_{\tilde{L}}\left( \tilde{\gamma}^{\tilde{L}}, f \right) = \sum_{L^\dagger \in \mathcal{L}^G(M)} d^G_M(L, L^\dagger) I^{\tilde{L}^\dagger}_{\tilde{M}}\left( \tilde{\gamma}, f_{\tilde{L}^\dagger} \right). \]
	Cf.\ Proposition \ref{prop:orbint-weighted-descent-Arch}.
	
	\item For each $\tilde{\gamma} \in D_{\text{tr-orb}, -}(\tilde{M}, \mathcal{O}) \otimes \mes(M)^\vee$, $L \in \mathcal{L}^G(M)$ and $J \in \mathcal{J}^L_M$, we have
	\[ \rho^{\tilde{G}}_J\left( \tilde{\gamma}^{\tilde{L}} \right) = \sum_{\substack{L^\dagger \in \mathcal{L}^G(M) \\ \text{s.t.}\; J \in \mathcal{J}^{L^\dagger}_M }} d^G_M(L, L^\dagger) \rho^{\tilde{L}^\dagger}_J \left( \tilde{\gamma} \right)^{\tilde{L}}. \]
	Cf.\ Proposition \ref{prop:rho-descent}.
	
	\item When $\tilde{\gamma} \in D_{\text{tr-orb}}(\tilde{M}, \mathcal{O}) \otimes \mes(M)^\vee$ and $\mathcal{O}$ is $G$-equisingular, then the $I^{\tilde{G}}_{\tilde{M}}(\tilde{\gamma}, f)$ in (O2) coincides with the $I^{\tilde{G}}_{\tilde{M}}(\tilde{\gamma}, f)$ of Definition \ref{def:weighted-I-integral-arch}.
\end{enumerate}

The condition (C8) is not explicitly stated in \cite{MW16-1}. Nonetheless, the arguments below for (C8) will also apply to the scenario in \textit{loc.\ cit.}. In particular, the stable analogues for (C8) hold on the endoscopic side.

Write $\tilde{G} = \prod_{i \in I} \GL(n_i, F) \times \Mp(W)$. We shall establish these objectives by assuming
\begin{itemize}
	\item the case when $\tilde{G}$ is replaced by a group of metaplectic type of smaller dimension;
	\item the case when $M$ is replaced by some strictly larger Levi subgroup of $G$, which is legitimate since the case $M=G$ is already known (namely $\rho^{\tilde{G}}_{\emptyset} = \identity$, and the distributions $I^{\tilde{G}}_{\tilde{G}}$ are just invariant orbital integrals);
	\item the following hypothesis.
\end{itemize}

\begin{hypothesis}\label{hyp:ext-Arch}
	For all $\tilde{\gamma} \in D_{\mathrm{orb}, -}(\tilde{M}) \otimes \mes(M)^\vee$ with $\Supp(\tilde{\gamma}) \subset \tilde{M}_{G\text{-reg}}$, and all $f \in \orbI_{\asp}(\tilde{G}) \otimes \mes(G)$, we have
	\[ I^{\tilde{G}}_{\tilde{M}}(\tilde{\gamma}, f) = I^{\tilde{G}, \Endo}_{\tilde{M}}(\tilde{\gamma}, f). \]
\end{hypothesis}

Note that the equality above extends automatically to $\tilde{\gamma} \in D_{\text{geom, $G$-equi}, -}(\tilde{M}) \otimes \mes(M)^\vee$ by Proposition \ref{prop:local-geometric-regular-reduction}.

The Hypothesis \ref{hyp:ext-Arch} will be verified near the end of this work, where we establish the local geometric Theorem \ref{prop:local-geometric}; see \S\ref{sec:end-stabilization}.

\subsection{Preliminary reductions}
This part is modeled on \cite[V.2.5]{MW16-1}. Let $\tilde{G} \supset \tilde{M}$ be as in \S\ref{sec:ext-Arch-prog}. First off, recall that $I^{\tilde{G}, \Endo}_{\tilde{M}}(\mathbf{M}^!, \delta, f)$, etc.\ (resp.\ $\rho^{\tilde{G}, \Endo}_J(\mathbf{M}^!, \delta)$, etc.) are well-defined in (O4) (resp.\ (O3)).

\begin{lemma}\label{prop:germ-ext-Arch}
	Let $\mathbf{M}^! \in \Endo_{\elli}(\tilde{M})$, $\mathcal{O}^! \subset M^!(F)$ be a stable semisimple conjugacy class, and $\delta \in SD_{\text{tr-orb}}(M^!, \mathcal{O}^!) \otimes \mes(M^!)^\vee$. For all $f \in \orbI_{\asp}(\tilde{G}) \otimes \mes(G)$, the germ near $a = 1$ of
	\[ a \mapsto I^{\tilde{G}, \Endo}_{\tilde{M}}\left( \mathbf{M}^!, \xi(a)\delta, f \right) \]
	is equivalent to that of
	\begin{multline*}
		a \mapsto I^{\tilde{G}, \Endo}_{\tilde{M}}\left(\mathbf{M}^!, \delta, f\right) + \sum_{J \in \mathcal{J}^G_M} I^{\tilde{G}}_{\tilde{G}}\left( \rho^{\tilde{G}, \Endo}_J(\mathbf{M}^!, \delta, a)^{\tilde{G}}, f \right) \\
		+ \sum_{L \in \mathcal{L}^G(M) \smallsetminus \{M, G\}} \sum_{J \in \mathcal{J}^L_M} I^{\tilde{G}}_{\tilde{L}}\left( \rho^{\tilde{L}}_J\left( \trans_{\mathbf{M}^!, \tilde{M}}(\delta) , a\right)^{\tilde{L}} , f \right),
	\end{multline*}
	where $a \in A_M(F)$ is in general position so that $\xi(a)\delta$ is $\tilde{G}$-equisingular.
\end{lemma}
\begin{proof}
	Arguments from \cite[II.3.9]{MW16-1} show that the germ of $I^{\tilde{G}, \Endo}_{\mathbf{M}^!}\left( \mathbf{M}^!, \xi(a)\delta, f \right)$ is equivalent to
	\begin{multline*}
		\sum_{L \in \mathcal{L}^G(M)} \sum_{J \in \mathcal{J}^L_M} \sum_{s \in \Endo_{\mathbf{M}^!}(\tilde{L})} i_{M^!}(\tilde{L}, L^![s])
		\sum_{\substack{J^! \in \mathcal{J}^{L^![s]}_{M^!}( B^{\tilde{G}}_{\mathcal{O}^![s]} ) \\ J^! \mapsto J}} I^{\tilde{G}, \Endo}_{\tilde{L}}\left(\mathbf{L}^![s], \sigma^{L^![s]}_{J^!}(\delta[s], \xi(a))^{L^![s]} , f \right).
	\end{multline*}
	This step is purely combinatorial, the only inputs being Lemma \ref{prop:i-transitivity} and Lemma \ref{prop:sL-Ls}.

	Note that $B^{\tilde{G}}_{\mathcal{O}^![s]} \big|_{L^![s]} = B^{\tilde{L}}_{\mathcal{O}^![s]}$ by Lemma \ref{prop:B-Levi}. We shall show that the asserted decomposition corresponds to the contributions of $L = M$, $L = G$ and $L \notin \{M, G\}$, respectively. This step requires some care because of the $z[s]$-twists.

	For all $(L, s)$ with $L \neq M$, Proposition \ref{prop:Levi-central-twist} and the induction hypotheses applied to (C2) yield
	\begin{align*}
		I^{\tilde{G}, \Endo}_{\tilde{L}}\left(\mathbf{L}^![s], \sigma^{L^![s]}_{J^!}(\delta[s], \xi(a))^{L^![s]} , f \right) & = I^{\tilde{G}}_{\tilde{L}}\left( \trans_{\mathbf{L}^![s], \tilde{L}}\left(\left( \sigma^{L^![s]}_{J^!}(\delta[s], \xi(a)) \right)^{L^![s]}\right) , f \right) \\
		& = I^{\tilde{G}}_{\tilde{L}}\left( \trans^s_{\mathbf{M}^!, \tilde{M}}\left( \sigma^{L^![s]}_{J^!}(\delta[s], \xi(a)) \right)^{\tilde{L}} , f \right).
	\end{align*}
	For fixed $L \neq M$ and $J$, the corresponding contribution is thus equal to
	\[ I^{\tilde{G}}_{\tilde{L}}\left( \rho^{\tilde{L}, \Endo}_J(\mathbf{M}^!, \delta, a)^{\tilde{L}}, f \right). \]
	\begin{compactitem}
		\item When $L \neq M, G$, the induction hypotheses applied to (C3) show that the contribution is $I^{\tilde{G}}_{\tilde{L}}\left( \rho^{\tilde{L}}_J\left(\trans_{\mathbf{M}^!, \tilde{M}}(\delta), a\right)^{\tilde{L}}, f \right)$.
		\item When $L = G$, the contribution is $I^{\tilde{G}}_{\tilde{G}}\left( \rho^{\tilde{G}, \Endo}(\mathbf{M}^!, \delta, a) \right)$.
	\end{compactitem}

	Finally, when $L=M$ we have: $s$ is trivial, $\mathcal{J}^L_M = \{\emptyset\} = \mathcal{J}^{L^![s]}_{M^!}\left(B^{\tilde{L}}_{\mathcal{O}^![s]}\right)$ and $i_{M^!}(\tilde{L}, L^![s]) = 1$. The contribution is thus $I^{\tilde{G}, \Endo}_{\tilde{M}}\left(\mathbf{M}^!, \delta, f \right)$.
\end{proof}

We need the endoscopic version below of the conditions (C6) and (C7).

\begin{lemma}\label{prop:descent-ext-Arch}
	Let $(\mathbf{M}^!, \delta)$ and $\mathcal{O}^!$ be as before. Suppose that we are in the situation
	\[\begin{tikzcd}
		& \tilde{G} \\
		L^! \arrow[dashed, leftrightarrow, r, "\text{ell.}", "\text{endo.}"'] & \tilde{L} \arrow[hookrightarrow, u, "\text{Levi}"] \\
		M^! \arrow[dashed, leftrightarrow, r, "\text{ell.}", "\text{endo.}"'] \arrow[hookrightarrow, u, "\text{Levi}"] & \tilde{M} \arrow[hookrightarrow, u, "\text{Levi}"']
	\end{tikzcd} \quad \text{with} \quad \mathbf{L}^! := \mathbf{L}^![t] \in \Endo_{\elli}(\tilde{L}), \; t \in \Endo_{\mathbf{M}^!}(\tilde{L}). \]
	Then for all $\delta \in SD_{\text{tr-orb}}(M^!, \mathcal{O}^!) \otimes \mes(M)^!$, $J \in \mathcal{J}^G_L$ and $f \in \orbI_{\asp}(\tilde{G}) \otimes \mes(G)$, we have
	\begin{align*}
		I^{\tilde{G}, \Endo}_{\tilde{L}}\left( \mathbf{L}^!, \delta[t]^{L^!}, f \right) & = \sum_{L^\dagger \in \mathcal{L}^G(M)} d^G_M(L, L^\dagger) I^{\tilde{L}^\dagger, \Endo}_{\tilde{M}}\left( \mathbf{M}^!, \delta, f_{\tilde{L}^\dagger} \right), \\
		\rho^{\tilde{G}, \Endo}_J \left( \mathbf{L}^!, \delta[t]^{L^!} \right) & = \sum_{\substack{L^\dagger \in \mathcal{L}^G(M) \\ \text{s.t.}\; J \in \mathcal{J}^{L^\dagger}_M }} d^G_M(L, L^\dagger) \rho^{\tilde{L}^\dagger, \Endo}_J \left( \mathbf{M}^!, \delta \right)^{\tilde{L}} .
	\end{align*}
\end{lemma}
\begin{proof}
	The arguments are combinatorial. The first formula follows from exactly the same proof of Proposition \ref{prop:descent-orbint-Endo}. The second formula requires extra care, as explained below.

	Let us assume $\tilde{G} = \Mp(W)$ without loss of generality. We claim that
	\begin{multline*}
		\rho^{\tilde{G}, \Endo}_J\left( \mathbf{L}^!, \delta[t]^{L^!} \right) = \sum_{s \in \Endo_{\mathbf{L}^!}(\tilde{G})} i_{L^!}(\tilde{G}, G^![s]) \sum_{\substack{J^! \in \mathcal{J}^{G^![s]}_{L^!}\left( B^{\tilde{G}}_{\mathcal{O}^![t][s]} \right) \\ J^! \mapsto J }} \trans_{\mathbf{L}^!, \tilde{L}}^s \left( \sigma^{G^![s]}_{J^!}( \delta[t]^{L^!} [s] ) \right) \\
		= \sum_{s \in \Endo_{\mathbf{L}^!}(\tilde{G})} i_{L^!}(\tilde{G}, G^![s]) \sum_{\substack{J^! \in \mathcal{J}^{G^![r]}_{L^!}\left( B^{\tilde{G}}_{\mathcal{O}^![r]}\right) \\ J^! \mapsto J }} \sum_{\substack{L^{\dagger !} \in \mathcal{L}^{G^![r]}(M^!) \\ \text{s.t.}\; J^! \in \mathcal{J}^{L^{\dagger !}}_{M^!}\left( B^{\tilde{G}}_{\mathcal{O}^![r]} \big|_{L^{\dagger !}} \right) }} e^{G^![r]}_{M^!}(L^!, L^{\dagger !}) \trans_{\mathbf{L}^!, \tilde{L}}^s \left( \sigma^{L^{\dagger !}}_{J^!}( \delta[r] )^{L^!} \right),
	\end{multline*}
	where $\mathcal{O}^![t][s]$ also stands for the stable class it generates in $L^!(F)$, by abusing notation slightly. Indeed, the first equality follows from the definition. Next, as in the proof of Proposition \ref{prop:descent-orbint-Endo},
	\begin{compactitem}
		\item we write $r = r(s, t) \in \Endo_{\mathbf{M}^!}(\tilde{G})$;
		\item to each $(s, L^{\dagger !})$ we attach a canonical $L^\dagger \in \mathcal{L}^G(M)$ and decompose $r$ into $r^L \in \Endo_{\mathbf{M}^!}(\tilde{L}^\dagger)$, $r_L \in \Endo_{\mathbf{L}^{\dagger !}}(\tilde{G})$ so that $L^{\dagger !} = L^{\dagger !}[r^L]$ is an endoscopic group of $\tilde{L}^\dagger$;
		\item the situation can be summarized as
		\[\begin{tikzcd}[column sep=small]
			& G^![s] & \\
			L^! \arrow[dash, ru, "s"] & & L^{\dagger !} \arrow[dash, lu, "{r_L}"'] \\
			& M^! \arrow[dash, lu, "t"] \arrow[dash, ru, "{r^L}"'] \arrow[dash, uu, "r" description]; &
		\end{tikzcd} \]
	\end{compactitem}
	hence the second equality follows from $\delta[t]^{L^!}[s] = \delta[r]^{L^!}$ and the descent formula \cite[II.3.11 Lemme]{MW16-1} for $\sigma^{G^![s]}_{J^!}$.

	Recall that $\trans_{\mathbf{L}^!, \tilde{L}}^s(\cdots) = \trans_{\mathbf{L}^!, \tilde{L}}(\cdots [s]^{-1})$ and factorize $[s]^{-1}$ into $[r_L]^{-1} [r^L]^{-1} [t]$. We have
	\[ \sigma^{L^{\dagger !}}_{J^!}\left( \delta[r] \right)^{L^!} [s]^{-1} = \left( \sigma^{L^{\dagger !}}_{J^!}\left( \delta[r^L] [r_L] \right) [s]^{-1} \right)^{L^!} = \left( \sigma^{L^{\dagger !}}_{J^!}\left( \delta[r^L] \right) [r^L]^{-1} [t] \right)^{L^!} \]

	Since $\mathbf{L}^! = \mathbf{L}^![t]$, Proposition \ref{prop:Levi-central-twist} thus yields
	\begin{align*}
		\trans_{\mathbf{L}^!, \tilde{L}}^s \left( \sigma^{L^{\dagger !}}_{J^!}( \delta[r] )^{L^!} \right) & = \trans_{\mathbf{L}^!, \tilde{L}} \left( \left( \sigma^{L^{\dagger !}}_{J^!}\left( \delta[r^L] \right) [r^L]^{-1} [t] \right)^{L^!} \right) \\
		& = \trans_{\mathbf{M}^!, \tilde{M}}\left( \sigma^{L^{\dagger !}}_{J^!}\left( \delta[r^L] \right) [r^L]^{-1} \right)^{\tilde{L}} \\
		& = \trans_{\mathbf{M}^!, \tilde{M}}^{r^L}\left( \sigma^{L^{\dagger !}}_{J^!}\left( \delta[r^L]\right) \right)^{\tilde{L}}.
	\end{align*}

	Next, we have $\mathcal{J}^{L^{\dagger !}}_{M^!}\left( B^{\tilde{G}}_{\mathcal{O}^![r]} \big|_{L^{\dagger !}} \right) = \mathcal{J}^{L^{\dagger !}}_{M^!}\left( B^{\tilde{L}}_{\mathcal{O}^![r]} \right) = \mathcal{J}^{L^{\dagger !}}_{M^!} \left( B^{\tilde{L}}_{\mathcal{O}^![r^L]} \right)$ as the central twist by $z[r_L]$ has no effect here. Furthermore, when $e^{G^![r]}_{M^!}(L^!, L^{\dagger !}) \neq 0$ (thus $d^G_M(L, L^\dagger) \neq 0$), Lemma \ref{prop:J-transfer} and the discussions preceding Proposition \ref{prop:rho-descent} lead to the commutative diagram
	\[\begin{tikzcd}[arrows=hookrightarrow]
		\mathcal{J}^{L^{\dagger !}}_{M^!}\left( B^{\tilde{L}}_{\mathcal{O}^![r^L]} \right) \arrow[r] \arrow[d] & \mathcal{J}^{G^![r]}_{L^!}\left( B^{\tilde{G}}_{\mathcal{O}^![r]}\right) \arrow[d] \\
		\mathcal{J}^{L^\dagger}_M \arrow[r] & \mathcal{J}^G_L .
	\end{tikzcd}\]

	With the notations above, we arrive at
	\begin{multline*}
		\rho^{\tilde{G}, \Endo}_J\left( \mathbf{L}^!, \delta[t]^{L^!} \right) = \sum_{(L^{\dagger !}, s) } i_{L^!}(\tilde{G}, G^![s]) e^{G^![r]}_{M^!}\left( L^!, L^{\dagger !} \right) \\
		\xi_J(L^\dagger) \sum_{\substack{J^! \in \mathcal{J}^{\tilde{L}^\dagger}_{M^!}\left( B^{\tilde{L}}_{\mathcal{O}^![r^L]} \right) \\ J^! \mapsto J }} \trans_{\mathbf{M}^!, \tilde{M}}^{r^L}\left( \sigma^{L^{\dagger !}}_{J^!}\left( \delta[r^L]\right) \right)^{\tilde{L}}
	\end{multline*}
	where $\xi_J(L^\dagger) = 1$ when $d^G_M(L, L^\dagger) \neq 0$ and $J \in \mathcal{J}^{L^\dagger}_M$, otherwise $\xi_J(L^\dagger) = 0$.

	The combinatorial proof of Proposition \ref{prop:descent-orbint-Endo} can now be applied to transform this expression into
	\begin{multline*}
		\sum_{L^\dagger \in \mathcal{L}^G(M)} d^G_M(L^\dagger, L) \xi_J(L^\dagger) \sum_{r^L \in \Endo_{\mathbf{M}^!}(\tilde{L}^\dagger)} i_{\mathbf{M}^!}(\tilde{L}^\dagger, L^{\dagger !}) \\
		 \sum_{\substack{J^! \in \mathcal{J}^{\tilde{L}^\dagger}_{M^!}\left( B^{\tilde{L}}_{\mathcal{O}^![r^L]} \right) \\ J^! \mapsto J }} \trans_{\mathbf{M}^!, \tilde{M}}^{r^L}\left( \sigma^{L^{\dagger !}}_{J^!}\left( \delta[r^L]\right) \right)^{\tilde{L}}.
	\end{multline*}
	It remains to recall the definition of $\rho^{\tilde{L}^\dagger, \Endo}_J\left( \mathbf{M}^!, \delta \right)$.
\end{proof}

\begin{proposition}[Cf.\ {\cite[V.2.5]{MW16-1}}]\label{prop:ext-Arch-Jmax}
	Assume that Hypothesis \ref{hyp:ext-Arch} holds. Suppose that for all $J \in \mathcal{J}^G_M \smallsetminus \{J_{\max}\}$ and semisimple conjugacy classes $\mathcal{O} \subset M(F)$, we have $\rho^{\tilde{G}}_J$ as in (O1) meeting the conditions (C1), (C3). Then the program in \S\ref{sec:ext-Arch-prog} can be realized.
\end{proposition}
\begin{proof}
	We may assume $M \neq G$. Let $J \in \mathcal{J}^G_M$. Our conditions force that $\rho^{\tilde{G}}_J$ and $I^{\tilde{G}}_{\tilde{M}}$ must be defined as follows. Express $\tilde{\gamma} \in D_{\text{tr-orb}, -}(\tilde{M}, \mathcal{O}) \otimes \mes(M)^\vee$ as
	\begin{equation}\label{eqn:Jmax-0}
		\tilde{\gamma} = \tilde{\gamma}_{\text{orb}} + \sum_{i=1}^n \trans_{\mathbf{M}^!_i, \tilde{M}} \left( \delta_i \right)
	\end{equation}
	with $\tilde{\gamma} \in D_{\text{orb}, -}(\tilde{M}) \otimes \mes(M)^\vee$ and $\delta_i \in SD_{\text{tr-orb}}(M_i^!, \mathcal{O}^!_i) \otimes \mes(M_i^!)^\vee$, where $\mathbf{M}_i^! \in \Endo_{\elli}(\tilde{M})$ and $\mathcal{O}^!_i \mapsto \mathcal{O}$ for all $i$. Then
	\begin{align*}
		\rho^{\tilde{G}}_J\left( \tilde{\gamma} \right) & = \rho^{\tilde{G}}_J\left( \tilde{\gamma}_{\text{orb}} \right) + \sum_i \rho^{\tilde{G}, \Endo}_J \left( \mathbf{M}_i^!, \delta_i \right), \\
		I^{\tilde{G}}_{\tilde{M}}\left( \tilde{\gamma}, f \right) & = I^{\tilde{G}}_{\tilde{M}}\left( \tilde{\gamma}_{\text{orb}}, f \right) + \sum_i I^{\tilde{G}, \Endo}_{\tilde{M}}\left( \mathbf{M}_i^!, \delta_i, f \right).
	\end{align*}
	What remains is to show that they are independent of \eqref{eqn:Jmax-0} and satisfy the required conditions.

	Consider the germ at $a=1$ of $I^{\tilde{G}}_{\tilde{M}}(a\tilde{\gamma}_{\mathrm{orb}}, f) + \sum_i I^{\tilde{G}, \Endo}_{\tilde{M}}(\mathbf{M}^!_i, \xi_i(a)\delta_i, f)$, where $\xi_i: A_M \rightiso A_{M^!_i}$ comes from endoscopy. The germ of the first term is equivalent to
	\begin{equation}\label{eqn:Jmax-1}
		\sum_{L \in \mathcal{L}(M)} \sum_{J \in \mathcal{J}^L_M} I^{\tilde{G}}_{\tilde{L}}\left( \rho^{\tilde{L}}_J(\tilde{\gamma}_{\mathrm{orb}}, a)^{\tilde{L}}, f \right),
	\end{equation}
	whilst by Lemma \ref{prop:germ-ext-Arch}, the second term is equivalent to
	\begin{multline}\label{eqn:Jmax-2}
		\sum_i I^{\tilde{G}, \Endo}_{\tilde{M}}\left(\mathbf{M}^!_i, \delta_i, f\right) + \sum_{J \in \mathcal{J}^G_M} I^{\tilde{G}}_{\tilde{G}}\left( \sum_i \rho^{\tilde{G}, \Endo}_J(\mathbf{M}^!_i, \delta_i, a)^{\tilde{G}}, f \right) \\
		+ \sum_{\substack{L \in \mathcal{L}^G(M) \\ L \neq M, G}} \sum_{J \in \mathcal{J}^L_M} I^{\tilde{G}}_{\tilde{L}}\left( \rho^{\tilde{L}}_J\left( \sum_i \trans_{\mathbf{M}^!_i, \tilde{M}}(\delta_i) , a\right)^{\tilde{L}} , f \right).
	\end{multline}

	To simplify the expression, put
	\[ \underline{\rho}_J(\tilde{\gamma}) := \rho^{\tilde{G}}_J(\tilde{\gamma}_{\mathrm{orb}}) + \sum_{i=1}^n \rho^{\tilde{G}, \Endo}_J(\mathbf{M}^!_i, \delta_i) . \]
	When $J \in \mathcal{J}^G_M \smallsetminus \{J_{\max}\}$, we have $\underline{\rho}_J(\tilde{\gamma}) = \rho^{\tilde{G}}_J(\tilde{\gamma})$. The sum of \eqref{eqn:Jmax-1} and \eqref{eqn:Jmax-2} becomes
	\begin{equation}\label{eqn:Jmax-sum}\begin{gathered}
		\underbracket{I^{\tilde{G}}_{\tilde{M}}\left( \tilde{\gamma}_{\text{orb}}, f \right)}_{\text{$L=M$ in \eqref{eqn:Jmax-1}}} + \underbracket{\sum_i I^{\tilde{G}, \Endo}_{\tilde{M}}\left( \mathbf{M}_i^!, \delta_i, f \right)}_{\text{the first term in \eqref{eqn:Jmax-2}}} +
		\underbracket{I^{\tilde{G}}_{\tilde{G}}\left( \underline{\rho}_{J_{\max}}(\tilde{\gamma}, a)^{\tilde{G}} \right)}_{\substack{L=G, \; J=J_{\max} \;\text{in \eqref{eqn:Jmax-1}} \\ J = J_{\max} \in \mathcal{J}^G_M \; \text{in \eqref{eqn:Jmax-2}}}} \\
		+ \underbracket{\sum_{\substack{J \in \mathcal{J}^G_M \\ J \neq J_{\max}}} I^{\tilde{G}}_{\tilde{G}}\left( \rho^{\tilde{G}}_J(\tilde{\gamma}, a)^{\tilde{G}}, f \right)}_{\substack{L=G, J \neq J_{\max} \;\text{in \eqref{eqn:Jmax-1}} \\ J \in \mathcal{J}^G_M \smallsetminus \{J_{\max}\} \;\text{in \eqref{eqn:Jmax-2}} }}
		+ \underbracket{\sum_{\substack{L \in \mathcal{L}^G(M) \\ L \neq M, G}} \sum_{J \in \mathcal{J}^L_M} I^{\tilde{G}}_{\tilde{L}}\left( \rho^{\tilde{L}}_J(\tilde{\gamma}, a)^{\tilde{L}}, f \right)}_{\substack{ L \neq M,G \;\text{in \eqref{eqn:Jmax-1}} \\ \text{the final sum in \eqref{eqn:Jmax-2}} }}.
	\end{gathered}\end{equation}

	Up to equivalence of germs, \eqref{eqn:Jmax-sum} is determined by $I^{\tilde{G}}_{\tilde{M}}(a\tilde{\gamma}, f)$ for $a \to 1$ in general position so that $a\tilde{\gamma}$ is $\tilde{G}$-equisingular, by Hypothesis \ref{hyp:ext-Arch}. Hence
	\[ \left( I^{\tilde{G}}_{\tilde{M}}\left( \tilde{\gamma}_{\text{orb}}, f \right) + \sum_i I^{\tilde{G}, \Endo}_{\tilde{M}}\left( \mathbf{M}_i^!, \delta_i, f \right) \right) + I^{\tilde{G}}_{\tilde{G}}\left( \underline{\rho}_{J_{\max}}(\tilde{\gamma})^{\tilde{G}}, f \right) \]
	is independent of \eqref{eqn:Jmax-0} up to equivalence of germs.
	
	Consider the formula displayed above. Since the first term in parentheses is constant in $a$ whereas the second term is in $U_{J_{\max}}$, Proposition \ref{prop:U-germ} (iii) implies that both are independent of \eqref{eqn:Jmax-0}. Hence
	\[ I^{\tilde{G}}_{\tilde{M}}\left( \tilde{\gamma}, f \right) :=I^{\tilde{G}}_{\tilde{M}}\left( \tilde{\gamma}_{\text{orb}}, f \right) + \sum_i I^{\tilde{G}, \Endo}_{\tilde{M}}\left( \mathbf{M}_i^!, \delta_i, f \right) \]
	is well-defined. By varying $f$, we see that $\underline{\rho}_{J_{\max}}(\tilde{\gamma})$ is also independent of \eqref{eqn:Jmax-0} (since it is a coset under $\mathrm{Ann}^{\tilde{G}}_{\tilde{M}}$). One obtains the well-defined germ
	\[ \rho^{\tilde{G}}_{J_{\max}}(\tilde{\gamma}) := \underline{\rho}_{J_{\max}}(\tilde{\gamma}). \]

	The condition (C5) is verified by \eqref{eqn:Jmax-sum}, and (C1), (C2), (C3) all follow by construction. Now consider the condition (C4). The case $\tilde{\gamma} \in D_{\text{orb}, -}(\tilde{M}, \mathcal{O}) \otimes \mes(M)^\vee$ is known. In view of (C3), we may assume $\tilde{\gamma} = \trans_{\mathbf{M}^!, \tilde{M}}(\delta)$ for some $\mathbf{M}^! \in \Endo_{\elli}(\tilde{M})$ and $\delta$. Then by Proposition \ref{prop:Levi-central-twist},
	\begin{align*}
		\rho^{\tilde{G}}_J\left( \tilde{\gamma}, a \right)^{\tilde{G}} & = \rho^{\tilde{G}, \Endo}_J\left(\mathbf{M}^!, \delta, a \right)^{\tilde{G}} \\
		& = \sum_{s \in \Endo_{\mathbf{M}^!}(\tilde{G})} i_{M^!}(\tilde{G}, G^![s]) \sum_{J^! \mapsto J} \left(\trans^s_{\mathbf{M^!}, \tilde{M}}\left( \sigma^{G^![s]}_{J^!}(\delta[s], \xi(a)) \right)\right)^{\tilde{G}} \\
		& = \sum_s i_{M^!}(\tilde{G}, G^![s]) \sum_{J^! \mapsto J} \trans_{\mathbf{G}^![s], \tilde{G}}\left( \sigma^{G^![s]}_{J^!}(\delta[s], \xi(a))^{G^![s]} \right) \\
		& \in D_{\text{tr-orb}, -}(\tilde{M}, \mathcal{O}) \otimes \mes(M)^\vee .
	\end{align*}
	
	The descent formulas (C6) and (C7) are known when $\tilde{\gamma} \in D_{\text{orb}, -}(\tilde{M}, \mathcal{O}) \otimes \mes(M)^\vee$. The case of $\tilde{\gamma} = \trans_{\mathbf{M}^!, \tilde{M}}(\delta)$ for some $(\mathbf{M}^!, \delta)$ is covered by Lemma \ref{prop:descent-ext-Arch}.
	
	Finally, let $\mathcal{O}$ and $\tilde{\gamma}$ be as in (C8). Fix $f$ and denote by $I_1(\tilde{\gamma})$ (resp.\ $I_2(\tilde{\gamma})$) the $I^{\tilde{G}}_{\tilde{M}}(\tilde{\gamma}, f)$ of Definition \ref{def:weighted-I-integral-arch} (resp.\ of (O2)). By (C5) and Lemma \ref{prop:rho-equisingular-arch}, the germ of $a \mapsto I_1(a\tilde{\gamma})$ is equivalent to $I_2(\tilde{\gamma})$.
	
	However, it is clear from Definition \ref{def:weighted-I-integral-arch} and basic properties of $I^{\tilde{M}}_{\tilde{M}}(a\tilde{\gamma}, \cdot)$ that $a \mapsto I_1(a\tilde{\gamma})$ is $C^\infty$ near $a=1$. We conclude by Proposition \ref{prop:U-germ} (iii) that $I_1(\tilde{\gamma}) = I_2(\tilde{\gamma})$, which is (C8).
\end{proof}

Next, let $J \in \mathcal{J}^G_M$. By Proposition \ref{prop:GJ-Sp}, we have
\begin{align*}
	G_J & = \prod_{k \in K} \GL(n_k) \times \Sp(W_1) \times \cdots \times \Sp(W_n), \\
	\widetilde{G_J} & = \prod_{k \in K} \GL(n_k, F) \times \left( \Mp(W_1) \utimes{\bmu_8} \cdots \utimes{\bmu_8} \Mp(W_n) \right),
\end{align*}
where
\begin{compactitem}
	\item $K$ is a finite set and $n_k \in \Z_{\geq 1}$ for all $k \in K$;
	\item $W_1, \ldots, W_n$ are symplectic $F$-vector spaces, $n \in \Z_{\geq 0}$;
	\item $\Mp(W_1) \utimes{\bmu_8} \cdots \utimes{\bmu_8} \Mp(W_n)$ is the quotient of $\prod_{i=1}^n \Mp(W_i)$ by $\left\{ (z_i)_{i=1}^n \in \bmu_8^n : z_1 \cdots z_n = 1 \right\}$.
\end{compactitem}
Whence the evident central extension
\[ 1 \to \bmu_8 \to \widetilde{G_J} \xrightarrow{\rev_J} G_J(F) \to 1 . \]

Although it is not necessarily a group of metaplectic type, there is still a theory of elliptic endoscopy for $\widetilde{G_J}$. Obviously, it is the ``product'' of the theories for $\GL(n_k)$ and those for $\Mp(W_i)$, for all $k, i$; eg.\ the elliptic endoscopic data are indexed by sequences of pairs $(m'_i, m''_i)_{i=1}^n$ with $m'_i + m''_i = \frac{1}{2}\dim W_i$. In particular, it makes sense to talk about the transfers $\Trans_{\mathbf{G_J}^!, \widetilde{G_J}}$ and $\trans_{\mathbf{G_J}^!, \widetilde{G_J}}$ for each $\mathbf{G}_J^! \in \Endo_{\elli}(\widetilde{G_J})$; they are $\otimes$-products of transfers for various $\Mp(W_i)$. The dual group of $\widetilde{G_J}$, with trivial Galois action, is simply
\[ \widetilde{G_J}^\vee = \prod_{k \in K} \GL(n_k, \CC) \times \prod_{i=1}^n \Sp(2m_i, \CC), \quad 2m_i := \dim_F W_i. \]

In the same vein, one defines $\rho^{\widetilde{G_J}, \Endo}_J(\mathbf{M}^!, \delta)$ as in (O3). It satisfies the analogues of Lemmas \ref{prop:germ-ext-Arch} and \ref{prop:descent-ext-Arch} for $\widetilde{G_J}$. If the program in \S\ref{sec:ext-Arch-prog} can be achieved for each $\Mp(W_i)$, the same can be said for $\widetilde{G_J}$.

\begin{proposition}\label{prop:ext-Arch-Jmax-red}
	Assume that Hypothesis \ref{hyp:ext-Arch} holds. Suppose that for all $J \in \mathcal{J}^G_M$, $\mathbf{M}^! \in \Endo_{\elli}(\tilde{M})$ and $\delta \in SD_{\text{tr-orb}}(M^!, \mathcal{O}^!) \otimes \mes(M^!)^\vee$, where $\mathcal{O}^! \subset M^!(F)$ and $\mathcal{O} \subset M(F)$ are stable semisimple conjugacy classes with $\mathcal{O}^! \mapsto \mathcal{O}$, the map induced by Lemma \ref{prop:Ann-J}
	\[ D_{\mathrm{geom}, -}(\tilde{M}, \mathcal{O}) \otimes \mes(M)^\vee \big/ \mathrm{Ann}^{\widetilde{G_J}}_{\mathcal{O}} \twoheadrightarrow D_{\mathrm{geom}, -}(\tilde{M}, \mathcal{O}) \otimes \mes(M)^\vee \big/ \mathrm{Ann}^{\tilde{G}}_{\mathcal{O}} \]
	sends $\rho^{\widetilde{G_J}, \Endo}_J(\mathbf{M}^!, \delta)$ to $\rho^{\tilde{G}, \Endo}_J(\mathbf{M}^!, \delta)$. Then the premises of Proposition \ref{prop:ext-Arch-Jmax} hold.
\end{proposition}
\begin{proof}
	Given $\tilde{\gamma} \in D_{\text{tr-orb}, -}(\tilde{M}, \mathcal{O}) \otimes \mes(M)^\vee$, choose an expression
	\begin{equation*}
		\tilde{\gamma} = \tilde{\gamma}_{\text{orb}} + \sum_{i=1}^n \trans_{\mathbf{M}^!_i, \tilde{M}} \left( \delta_i \right).
	\end{equation*}
	In order to verify the premises of Proposition \ref{prop:ext-Arch-Jmax}, it suffices to show that for all $J \in \mathcal{J}^G_M \smallsetminus \{J_{\max}\}$, the germ
	\[ \rho^{\tilde{G}}_J(\tilde{\gamma}) := \rho^{\tilde{G}}_J(\tilde{\gamma}_{\mathrm{orb}}) + \sum_i \rho^{\tilde{G}, \Endo}_J(\mathbf{M}^!_i, \delta_i) \]
	is independent of the expression above of $\tilde{\gamma}$.

	The Archimedean version of Definition--Proposition \ref{def:rho-germ} implies that $\rho^{\widetilde{G_J}}_J(\tilde{\gamma}_{\mathrm{orb}}) \mapsto \rho^{\tilde{G}}_J(\tilde{\gamma}_{\mathrm{orb}})$ under the natural surjection. By assumption, $\rho^{\widetilde{G_J}, \Endo}_J(\mathbf{M}_i^!, \delta_i) \mapsto \rho^{\tilde{G}, \Endo}_J(\mathbf{M}^!_i, \delta_i)$ for all $i$. Since $\dim G_J < \dim G$, by induction we see that
	\[ \rho^{\widetilde{G_J}}_J(\tilde{\gamma}) := \rho^{\widetilde{G_J}}_J(\tilde{\gamma}_{\mathrm{orb}}) + \sum_i \rho^{\widetilde{G_J}, \Endo}_J(\mathbf{M}_i^!, \delta_i) \]
	is independent of the expression of $\tilde{\gamma}$. Hence so is its image $\rho^{\tilde{G}}_J(\tilde{\gamma})$.
\end{proof}

Note that $\rho^{\widetilde{G_J}, \Endo}_J(\mathbf{M}^!, \delta) \mapsto \rho^{\tilde{G}, \Endo}_J(\mathbf{M}^!, \delta)$ is tautologically true when $J = J_{\max}$, since $G_{J_{\max}} = G$.

\section{Conditional proof of the local theorems: Archimedean case}\label{sec:ext-Arch-proof}
Let $F$ be an Archimedean local field. Consider a group of metaplectic type $\rev: \tilde{G} \to G(F)$ and a Levi subgroup $M \subset G$. The goal here is to realize the program in \S\ref{sec:ext-Arch-prog} under the Hypothesis \ref{hyp:ext-Arch}. We will also take the induction hypotheses for granted. See the discussions preceding Hypothesis \ref{hyp:ext-Arch}.

In view of Propositions \ref{prop:ext-Arch-Jmax} and \ref{prop:ext-Arch-Jmax-red}, it suffices to show that $\rho^{\widetilde{G_J}, \Endo}_J(\mathbf{M}^!, \delta) \mapsto \rho^{\tilde{G}, \Endo}_J(\mathbf{M}^!, \delta)$ under the natural map, for each $J \in \mathcal{J}^G_M$. Here $(\mathbf{M}^!, \delta)$ and $\mathcal{O}^!$ are as in Proposition \ref{prop:ext-Arch-Jmax-red}.

\subsection{Reduction to the elliptic case}
Fix $\epsilon \in \mathcal{O}^!$ and $\eta \in \mathcal{O}$. By \cite[Lemma 3.3]{Ko82}, the $\epsilon$ can be chosen in $\mathcal{O}^!$ such that $M^!_\epsilon$ is quasisplit, which we assume hereafter.

\begin{lemma}\label{prop:ext-Arch-red-ell}
	Consider the statement that $\rho^{\widetilde{G_J}, \Endo}_J(\mathbf{M}^!, \delta) \mapsto \rho^{\tilde{G}, \Endo}_J(\mathbf{M}^!, \delta)$ under the natural map, for various $\tilde{M} \subset \tilde{G}$, $J \in \mathcal{J}^G_M$, $\mathcal{O}^! \subset M^!(F)_{\mathrm{ss}}$ and $(\mathbf{M}^!, \delta)$. If it holds whenever $\mathcal{O}^!$ is elliptic in $M^!$, then it holds in general.
\end{lemma}
\begin{proof}
	Define the Levi subgroup
	\[ R^! := Z_{M^!}(A_{M^!_\epsilon}) \subset M^! \]
	and let $\mathcal{O}_{R^!}$ be the stable semisimple class in $R^!(F)$ containing $\epsilon$. By Proposition \ref{prop:situation-easier}, to $\mathbf{M}^! \in \Endo_{\elli}(\tilde{M})$ and $R^! \hookrightarrow M^!$ is associated a Levi subgroup $R \subset M$, unique up to conjugacy, together with a unique $t \in \Endo_{\mathbf{R}^!}(\tilde{M})$ such that
	\[\begin{tikzcd}
		M^! \arrow[dashed, leftrightarrow, r, "\text{ell.}", "\text{endo.}"'] & \tilde{M} \\
		R^! \arrow[dashed, leftrightarrow, r, "\text{ell.}", "\text{endo.}"'] \arrow[hookrightarrow, u, "\text{Levi}"] & \tilde{R} \arrow[hookrightarrow, u, "\text{Levi}"']
	\end{tikzcd} \quad \text{and} \quad \mathbf{M}^! = \mathbf{M}^![t]. \]

	By \cite[V.4.4 Lemme]{MW16-1}, there exists $\delta_{R^!} \in SD_{\text{tr-orb}}(R^!, \mathcal{O}_{R^!}[t]^{-1}) \otimes \mes(R^!)^\vee$ such that $\delta = \left( \delta_{R^!}[t] \right)^{M^!}$, where $\mathcal{O}_{R^!}[t^{-1}] := \mathcal{O}_{R^!} \cdot z[t]^{-1}$.

	Lemma \ref{prop:descent-ext-Arch} implies
	\[ \rho^{\tilde{G}, \Endo}_J\left( \mathbf{M}^!, \delta \right) = \sum_{\substack{L \in \mathcal{L}^G(R) \\ \text{s.t.}\; J \in \mathcal{J}^L_R }} d^G_R(M, L) \rho^{\tilde{L}, \Endo}_J \left( \mathbf{R}^!, \delta_{R^!} \right)^{\tilde{M}}. \]
	
	The analogue of Lemma \ref{prop:descent-ext-Arch} for $\widetilde{G_J}$ also yields
	\[ \rho^{\widetilde{G_J}, \Endo}_J\left( \mathbf{M}^!, \delta \right) = \sum_{\substack{L^\dagger \in \mathcal{L}^{G_J}(R) \\ \text{s.t.}\; J \in \mathcal{J}^{L^\dagger}_R }} d^{G_J}_R(M, L^\dagger) \rho^{\widetilde{L^\dagger}, \Endo}_J \left( \mathbf{R}^!, \delta_{R^!} \right)^{\tilde{M}}. \]
	We contend that $L \mapsto L_J$ gives a bijection
	\[\begin{tikzcd}[row sep=tiny]
		\left\{ L \in \mathcal{L}^G(R) : J \in \mathcal{J}^L_R \right\} \arrow[r, "1:1"] & \left\{ L^\dagger \in \mathcal{L}^{G_J}(R) : J \in \mathcal{J}^{L^\dagger}_R \right\}
	\end{tikzcd}\]
	and we have $d^G_R(M, L) = d^{G_J}_R(M, L_J)$. Indeed, the desired bijection is in \cite[p.546]{MW16-1}, and so is the equality for $d_R(\cdot, \cdot)$ which is based on $\mathfrak{a}_{G_J} = \mathfrak{a}_G$ as subspaces of $\mathfrak{a}_M$.

	We may assume that $\epsilon$ is not elliptic in $M^!$. Then $R^! \neq M^!$, and those $L$ with $d^G_R(M, L) \neq 0$ must satisfy $L \neq G$. By induction, $\rho^{\widetilde{L_J}, \Endo}_J \left( \mathbf{R}^!, \delta_{R^!} \right)$ maps $\rho^{\tilde{L}, \Endo}_J \left( \mathbf{R}^!, \delta_{R^!} \right)$. It follows that $\rho^{\widetilde{G_J}, \Endo}_J(\mathbf{M}^!, \delta)$ maps to $\rho^{\tilde{G}, \Endo}_J(\mathbf{M}^!, \delta)$.
\end{proof}

\subsection{The elliptic case}\label{sec:ext-Arch-elliptic}
In what follows, we assume $\mathcal{O}^!$ is elliptic. Then $\mathcal{O}$ is also elliptic by Lemma \ref{prop:ellipticity-transfer}.

The endoscopic datum of $M_\eta$ obtained by descent is relevant (Proposition \ref{prop:relevance-elliptic}), hence we can take a diagram $(\epsilon, B^!, T^!, B, T, \eta)$ for the endoscopic datum $\mathbf{M}^!$ of $\tilde{M}$, compatibly with descent in the sense of Proposition \ref{prop:diagram-relevant}. It affords isomorphisms between tori
\[\begin{tikzcd}[row sep=small]
	T^! \arrow[r, "{\xi_{T^!, T}}"] & T \\
	A_{M^!} \arrow[phantom, u, "\subset" description, sloped] \arrow[r, "{\xi^{-1}}"'] & A_M \arrow[phantom, u, "\subset" description, sloped] .
\end{tikzcd}\]
In turn, we obtain $\xi^*: \mathfrak{a}_{M^!}^* \rightiso \mathfrak{a}_M^*$; see Remark \ref{rem:xi-center-diagram}.

\begin{notation}
	Fix $J \in \mathcal{J}^G_M \smallsetminus \{J_{\max}\}$ and $\delta \in SD_{\text{tr-orb}}(M^!) \otimes \mes(M^!)^\vee$. For the sake of simplicity, $B$ (resp.\ $B_J$) will denote the system of $B$-functions $B_{\mathcal{O}^!}$ (resp.\ $B^J_{\mathcal{O}^!}$, defined relative to $\tilde{M} \subset \widetilde{G_J}$).
\end{notation}

Let $s \in \Endo_{\mathbf{M}^!}(\tilde{G})$. Definition \ref{def:B-function} furnishes the following finite subsets of $X^*(T^!_{\overline{F}}) \otimes \Q$:
\begin{align*}
	\Sigma\left(G^![s]_{\epsilon[s]}, T^!, B \right) & := \left\{ \beta' = B(\beta)^{-1}\beta : \beta \in \Sigma\left( G^![s]_{\epsilon[s]}, T^! \right) \right\}, \\
	\Sigma^J\left(G^![s]_{\epsilon[s]}, T^!, B \right) & := \left\{ \beta' \in \Sigma\left(G^![s]_{\epsilon[s]}, T^!, B \right) : \xi^*( \beta' |_{\mathfrak{a}_{M^!}}) \in R_J \right\}, \\
	\Sigma^J \left( G^![s]_{\epsilon[s]}, T^! \right) & := \left\{ \beta \in \Sigma\left( G^![s]_{\epsilon[s]}, T^! \right) : B(\beta)^{-1} \beta \in \Sigma^J\left(G^![s]_{\epsilon[s]}, T^!, B \right) \right\}.
\end{align*}
Note that $\Sigma\left(G^![s]_{\epsilon[s]}, T^!, B \right)$ (resp.\ $\Sigma^J\left(G^![s]_{\epsilon[s]}, T^!, B \right)$) is obtained from $\Sigma\left( G^![s]_{\epsilon[s]}, T^! \right)$ (resp.\ $\Sigma^J \left( G^![s]_{\epsilon[s]}, T^! \right)$) by multiplying by $B(\cdot)^{-1}$.

Similarly, for $t \in \Endo_{\mathbf{M}^!}(\widetilde{G_J})$ we define the subsets of $X^*(T^!_{\overline{F}}) \otimes \Q$:
\[ \Sigma\left( G_J^![t]_{\epsilon[t]}, T^! \right), \quad \Sigma\left( G_J^![t]_{\epsilon[t]}, T^! , B_J\right), \]
which are related by dilation by $B_J(\cdot)^{-1}$.

By Lemma \ref{prop:J-transfer}, there is at most one $J^! \in \mathcal{J}^{G^![s]}_{M^!}(B)$ with $J^! \mapsto J$. If it exists, denote it as $J_s$. Similarly, for each $t \in \Endo_{\mathbf{M}^!}(\widetilde{G_J})$, denote by $J_t \in \mathcal{J}^{G_J^![t]}_{M^!}(B_J)$ the unique element mapping to $J$, if it exists.

Given $s$ and $J_s$, we construct the connected reductive group $G^![s]_{\epsilon[s], J_s}$ using the recipe of \S\ref{sec:GJ-B} and the $A_{M^!_{\epsilon[s]}} = A_{M^!}$. It contains $T^!$ as a maximal $F$-torus, and satisfies
\begin{equation}\label{eqn:J-G1}\begin{aligned}
	\Sigma\left( G^![s]_{\epsilon[s], J_s}, T^! \right) & = \left\{ \beta \in \Sigma\left( G^![s]_{\epsilon[s]}, T^! \right) : B(\beta)^{-1} \beta|_{\mathfrak{a}_{M^!}} \in R_{J_s} \right\} \\
	& = \Sigma^J\left( G^![s]_{\epsilon[s]}, T^! \right),
\end{aligned}\end{equation}
the last equality being a consequence of the (rather roundabout) definitions above.

There is a natural surjection
\begin{align*}
	\Endo_{\mathbf{M}^!}(\tilde{G}) & \twoheadrightarrow \Endo_{\mathbf{M}^!}(\widetilde{G_J}) \\
	s & \mapsto t
\end{align*}
which can be interpreted in two ways.
\begin{enumerate}
	\item Let $s^\flat \in \tilde{M}^\vee$ be an elliptic semisimple element determining $\mathbf{M}^!$. The map above is induced by
	\[ s^\flat Z_{\tilde{M}^\vee}^{\circ} / Z_{\tilde{G}^\vee}^{\circ} \to s^\flat Z_{\tilde{M}^\vee}^{\circ} / Z_{\widetilde{G_J}^\vee}^{\circ} \]
	using $Z_{\widetilde{G_J}^\vee} \supset Z_{\tilde{G}^\vee}$.
	\item Alternatively, let us assume for simplicity that $\tilde{G} = \Mp(W)$ and $\tilde{M} = \prod_{i \in I} \GL(n_i, F) \times \Mp(W^\flat)$. With the notation of Proposition \ref{prop:GJ-Sp}, write
	\[ G_J = \prod_{i \in \mathbf{I}} \Sp(W^i) \times \Sp(W^\natural) \]
	and assume that $M \hookrightarrow G_J$ is given by a map $f: I \to \mathbf{I} \sqcup \{\natural\}$, so that
	\begin{compactitem}
		\item $\Sp(W^\flat) \hookrightarrow \Sp(W^\natural)$,
		\item $\GL(n_k) \hookrightarrow \Sp(W^{f(k)})$ for all $k \in I$.
	\end{compactitem}
	Then we obtain
	\[ I = \bigsqcup_{i \in \mathbf{I} \sqcup \{\natural\}} I_i, \quad I_i := f^{-1}(i). \]
	Given $s \in \Endo_{\mathbf{M}^!}(\tilde{G})$, corresponding to $I = I' \sqcup I''$, define $I'_i := I_i \cap I'$ and $I''_i := I_i \cap I''$ for all $i \in \mathbf{I} \sqcup \{\natural\}$. Then the map $\Endo_{\mathbf{M}^!}(\tilde{G}) \twoheadrightarrow \Endo_{\mathbf{M}^!}(\widetilde{G_J})$ maps $s$ to the element $t$ corresponding to the family of partitions $\left( I_i = I'_i \sqcup I''_i\right)_{i \in \mathbf{I} \sqcup \{\natural\}}$.
\end{enumerate}

The second interpretation leads immediately to the surjectivity of the map, as well as
\begin{equation}\label{eqn:zs-zt}
	z[s] = z[t] \; \in Z_{M^!}(F).
\end{equation}

\begin{definition}\label{def:Z-ZJ}
	For $s \in \Endo_{\mathbf{M}^!}(\tilde{G})$ (resp.\ $t \in \Endo_{\mathbf{M}^!}(\widetilde{G_J})$), write $\exists J_s$ (resp.\ $\exists J_t$) to signify that there exists $J_s$ (resp.\ $J_t$) mapping to $J$ via endoscopy. Let
	\[\begin{tikzcd}[column sep=tiny]
		\mathcal{Z} \arrow[phantom, r, ":=" description] & \left\{ s \in \Endo_{\mathbf{M}^!}(\tilde{G}) : \exists J_s \right\} \arrow[phantom, r, "\subset" description] & \Endo_{\mathbf{M}^!}(\tilde{G}) \arrow[twoheadrightarrow, d] \\
		\mathcal{Z}_J \arrow[phantom, r, ":=" description] & \left\{ t \in \Endo_{\mathbf{M}^!}(\widetilde{G_J}) : \exists J_t \right\} \arrow[phantom, r, "\subset" description] & \Endo_{\mathbf{M}^!}(\widetilde{G_J}) .
	\end{tikzcd}\]
	Moreover, we fix an elliptic semisimple element $s^\flat \in \tilde{M}^\vee$ to identify $\mathcal{Z}$ (resp.\ $\mathcal{Z}_J$) as a subset of $s^\flat Z_{\tilde{M}^\vee}^{\circ} / Z_{\tilde{G}^\vee}^{\circ}$ (resp.\ $s^\flat Z_{\tilde{M}^\vee}^{\circ} / Z_{\widetilde{G_J}^\vee}^{\circ}$).
\end{definition}

\begin{lemma}\label{prop:J-roots}
	If $s \in \Endo_{\mathbf{M}^!}(\tilde{G})$ maps to $t \in \Endo_{\mathbf{M}^!}(\widetilde{G_J})$, then
	\[ \Sigma^J \left( G^![s]_{\epsilon[s]}, T^!, B \right) = \Sigma\left( G_J^![t]_{\epsilon[t]}, T^!, B_J \right). \]
\end{lemma}
The proof is deferred to \S\ref{sec:proof-germ-lemmas}.

\begin{lemma}\label{prop:Z-Z}
	In the circumstance of Definition \ref{def:Z-ZJ}, $\mathcal{Z}$ equals the preimage of $\mathcal{Z}_J$.
\end{lemma}
\begin{proof}
	Let $s \in \Endo_{\mathbf{M}^!}(\tilde{G})$. The condition $s \in \mathcal{Z}$ is equivalent to that there exist linearly independent $\beta_1, \ldots, \beta_r \in \Sigma\left( G^![s]_{\epsilon[s]}, T^!, B \right)$ such that $r = \dim \mathfrak{a}^G_M$ and $\left\{\xi^*(\beta_i|_{\mathfrak{a}_{M^!}})\right\}_{i=1}^r$ generate $R_J$. These conditions force $\beta_i \in \Sigma^J\left( G^![s]_{\epsilon[s]}, T^!, B \right)$ for all $i$.
	
	Likewise, for $t \in \Endo_{\mathbf{M}^!}(\widetilde{G_J})$, the condition $t \in \mathcal{Z}_J$ is equivalent to that there exist linearly independent $\beta_1, \ldots, \beta_r \in \Sigma\left( G^!_J[t]_{\epsilon[t]}, T^!, B_J \right)$ such that $r = \dim \mathfrak{a}^{G_J}_M$ and $\left\{\xi^*(\beta_i|_{\mathfrak{a}_{M^!}})\right\}_{i=1}^r$ generate $R_J$.
	
	Since $\dim \mathfrak{a}^{G_J}_M = \dim \mathfrak{a}^G_M$, we conclude by Lemma \ref{prop:J-roots}.
\end{proof}

Suppose that $s \in \mathcal{Z}$ maps to $t \in \mathcal{Z}_J$. View both $Z_{{G^!_J[t]}^\vee}$ and $Z_{{G^![s]}^\vee}$ as subgroups of $Z_{{M^!}^\vee}$. Note that $Z_{\tilde{G}^\vee}^\circ$ lies in their intersection. Define
\begin{equation}\label{eqn:GJ-Gs}\begin{aligned}
	\left( Z_{{G^!_J[t]}^\vee} : Z_{{G^![s]}^\vee} \right) & := \frac{\left(Z_{{G^!_J[t]}^\vee} : Z_{\tilde{G}^\vee}^\circ\right)}{\left( Z_{{G^![s]}^\vee} : Z_{\tilde{G}^\vee}^\circ \right)} \; \in \Q_{> 0}, \\
	\left( Z_{{G^![s]}^\vee} : Z_{{G^!_J[t]}^\vee} \right) & := \left( Z_{{G^!_J[t]}^\vee} : Z_{{G^![s]}^\vee} \right)^{-1} .
\end{aligned}\end{equation}

\begin{lemma}[Cf.\ {\cite[p.552 (7)]{MW16-1}}]\label{prop:Z-Z-i}
	For $s \in \mathcal{Z}$ mapping to $t \in \mathcal{Z}_J$, we have
	\[ i_{M^!}\left( \widetilde{G_J}, G^!_J[t] \right) \left( Z_{{G^!_J[t]}^\vee} : Z_{{G^![s]}^\vee} \right) = \underbracket{\left( Z_{\widetilde{G_J}^\vee}^\circ : Z_{\tilde{G}^\vee}^\circ \right) }_{=1\; \text{in our case}} i_{M^!}\left( \tilde{G}, G^![s] \right) . \]
\end{lemma}
\begin{proof}
	In the commutative diagram below of diagonalizable $\CC$-groups, all arrows are surjective with finite kernels:
	\[\begin{tikzcd}
		Z_{\tilde{M}^\vee}^\circ / Z_{\tilde{G}^\vee}^\circ \arrow[r, "\mathbf{A}"] \arrow[d, "\mathbf{B}"'] & Z_{\tilde{M}^\vee}^\circ / Z_{\widetilde{G_J}^\vee}^\circ \arrow[d, "\mathbf{C}"] \\
		Z_{{M^!}^\vee} / Z_{\tilde{G}^\vee}^\circ \arrow[r, "\mathbf{D}"'] \arrow[d, "\mathbf{E}"'] & Z_{{M^!}^\vee} / Z_{{G^!_J[t]}^\vee} \\
		Z_{{M^!}^\vee} / Z_{{G^![s]^\vee}} . &
	\end{tikzcd}\]
	Hence
	\[ | \Ker(\mathbf{A}) | \cdot |\Ker(\mathbf{C})| = |\Ker(\mathbf{B})| \cdot |\Ker(\mathbf{D})| = \frac{|\Ker(\mathbf{E}\mathbf{B})|}{|\Ker(\mathbf{E})|} \cdot |\Ker(\mathbf{D})|. \]
	On the other hand, one readily checks that
	\begin{gather*}
		i_{M^!}\left( \widetilde{G_J} : G^!_J[t] \right)  = |\Ker(\mathbf{C})|^{-1}, \quad \left( Z_{{G^!_J[t]}^\vee} : Z_{{G^![s]}^\vee} \right) = \frac{|\Ker(\mathbf{D})|}{|\Ker(\mathbf{E})|}, \\
		\left( Z_{\widetilde{G_J}^\vee}^\circ : Z_{\tilde{G}^\vee}^\circ \right) = \underbracket{|\Ker(\mathbf{A})|}_{=1\; \text{in our case}}, \quad i_{M^!}\left( \tilde{G}, G^![s] \right) = |\Ker(\mathbf{E}\mathbf{B})|^{-1}.
	\end{gather*}
	The assertion follows at once.
\end{proof}

We keep the constants which equal to $1$ in order to maintain the analogy with \cite[V.5]{MW16-1}.

\begin{lemma}\label{prop:sigma-sigma}
	For $s \in \mathcal{Z}$ mapping to $t \in \mathcal{Z}_J$, we have the equality
	\[ \left( Z_{{G^!_J[t]}^\vee} : Z_{{G^![s]}^\vee} \right) \cdot \trans^s_{\mathbf{M}^!, \tilde{M}} \left( \sigma^{G^![s]}_{J_s}(\delta[s]) \right) = \trans^t_{\mathbf{M}^!, \tilde{M}}\left( \sigma^{G^!_J[t]}_{J_t}(\delta[t]) \right) \]
	of germs, where $\delta \in SD_{\text{tr-orb}}(M^!, \mathcal{O}^!) \otimes \mes(M^!)^\vee$.
\end{lemma}

The proof is deferred to \S\ref{sec:proof-germ-lemmas}. Let us show how it completes the program in \S\ref{sec:ext-Arch-prog}.
\begin{theorem}\label{prop:ext-Arch-prog-realize}
	Under the Hypothesis \ref{hyp:ext-Arch}, the program in \S\ref{sec:ext-Arch-prog} can be realized. Specifically, $\rho^{\widetilde{G_J}, \Endo}_J(\mathbf{M}^!, \delta) \mapsto \rho^{\tilde{G}, \Endo}_J(\mathbf{M}^!, \delta)$ under the natural map for all $\tilde{M} \subset \tilde{G}$, $(\mathbf{M}^!, \delta)$, $\mathcal{O}^!$ and $J \in \mathcal{J}^G_M$.
\end{theorem}
\begin{proof}
	We may assume $J \neq J_{\max}$. Lemma \ref{prop:ext-Arch-red-ell} reduces the assertion to the case of elliptic $\mathcal{O}^!$. By (O3) and Lemma \ref{prop:Z-Z-i}, we have the equality of germs
	\begin{multline*}
		\rho^{\tilde{G}, \Endo}_J\left( \mathbf{M}^!, \delta \right) = \sum_{s \in \mathcal{Z}} i_{M^!}(\tilde{G}, G^![s]) \trans^s_{\mathbf{M}^!, \tilde{M}} \left( \sigma^{G^![s]}_{J_s}(\delta[s]) \right) \\
		= \sum_{s \in \mathcal{Z}} \left( Z_{{G^!_J[t]}^\vee} : Z_{{G^![s]}^\vee} \right) i_{M^!}\left( \widetilde{G_J} : G^![t] \right) \trans^s_{\mathbf{M}^!, \tilde{M}} \left( \sigma^{G^![s]}_{J_s}(\delta[s]) \right)
	\end{multline*}
	where $s \mapsto t \in \mathcal{Z}_J$. Lemma \ref{prop:sigma-sigma} transforms the last expression into
	\[ \sum_{s \in \mathcal{Z}} i_{M^!}\left( \widetilde{G_J} : G^![t] \right) \trans^t_{\mathbf{M}^!, \tilde{M}} \left( \sigma^{G_J^![t]}_{J_t}(\delta[t]) \right). \]
	
	In view of Lemma \ref{prop:Z-Z}, $\left( Z_{\widetilde{G_J}^\vee}^\circ : Z_{\tilde{G}^\vee}^\circ \right) = 1$ is the common cardinality of fibers of $\mathcal{Z} \twoheadrightarrow \mathcal{Z}_J$. All in all, we obtain the germ $\rho^{\widetilde{G_J}, \Endo}_J(\mathbf{M}^!, \delta)$.
\end{proof}

\subsection{Proof of crucial lemmas}\label{sec:proof-germ-lemmas}
Retain the notations of $J$, $s \mapsto t$, $J_s$, $J_t$, $\delta$ etc., as in Lemma \ref{prop:sigma-sigma}. By \eqref{eqn:zs-zt}, we have
\[ \epsilon[s] = \epsilon[t], \quad \delta[s] = \delta[t], \quad M^!_{\epsilon[s]} = M^!_\epsilon = M^!_{\epsilon[t]}. \]

Define $\Xi_\xi$ to be the finite $F$-group scheme $Z_{M^!}(\epsilon) / M^!_\epsilon$. Then $\Xi_\epsilon(F)$ acts on $SD_{\mathrm{unip}}(M^!_\epsilon) \otimes \mes(M^!_\epsilon)^\vee$, leaving $SD_{\text{tr-unip}}(M^!_\epsilon) \otimes \mes(M^!_\epsilon)^\vee$ invariant. By the general theory of descent, the dual of stable descent maps
\[\begin{tikzcd}
	SD_{\mathrm{unip}}(M^!_\epsilon) \otimes \mes(M^!_\epsilon) \arrow[r, "{S\desc^{M^!, *}_{\epsilon[s]}}" inner sep=0.8em] \arrow[rd, "{S\desc^{M^!, *}_\epsilon}"'] & SD_{\mathrm{geom}}(M^!, \mathcal{O}^![s]) \otimes \mes(M^!)^\vee & \delta[s] \arrow[phantom, l, "\ni" description, sloped] \\
	& SD_{\mathrm{geom}}(M^!, \mathcal{O}^!) \otimes \mes(M^!)^\vee \arrow[u, "{\cdot z[s]}"', "\simeq"] & \delta \arrow[phantom, l, "\ni" description, sloped] \arrow[mapsto, u]
\end{tikzcd}\]
become injective after restriction to $\Xi_\epsilon(F)$-invariants. Note that the diagram commutes.

This injectivity together with the ``localization'' of tr-orb spaces in \cite[V.4.3 Lemme (ii)]{MW16-1} furnish a unique $\delta_\epsilon \in SD_{\text{tr-unip}}(M^!_\epsilon) \otimes \mes(M^!_\epsilon)^\vee$ such that $\Supp(\delta_\epsilon)$ is close to $1$ and
\[ S\desc^{M^!, *}_\epsilon(\delta_\epsilon) = \delta, \quad \text{thus} \quad S\desc^{M^!, *}_{\epsilon[s]}(\delta_\epsilon) = \delta[s] \]

We assumed $A_{M^!_\epsilon} = A_{M^!}$. In proving Lemma \ref{prop:sigma-sigma}, we may also assume $A_{G^![s]_{\epsilon[s]}} = A_{G^![s]}$, since otherwise $\mathcal{J}^{G^![s]}_{M^!}(B) = \emptyset$, contradicting that $s \in \mathcal{Z}$.

Under these assumptions, we have the bijection $\mathcal{J}^{G^![s]_{\epsilon[s]}}_{M^!_\epsilon}(B) \xrightarrow{1:1} \mathcal{J}^{G^![s]}_{M^!}(B)$ such that $J_s \mapsto J$, and \cite[V.4.6 (9)]{MW16-1} entails
\[ \sigma^{G^![s]}_{J_s}\left( \delta[s] \right) = S\desc^{M^!, *}_{\epsilon[s]} \left( e^{G^![s]}_{M^!}(\epsilon[s]) \sigma^{G^![s]_{\epsilon[s]}}_{J_s}(\delta_\epsilon) \right), \]
where $e^{G^![s]}_{M^!}(\epsilon[s])$ is as in Definition \ref{def:eLRmu}. Analogous constructions for $G_J$ lead to
\[ \sigma^{G_J^![t]}_{J_t}\left( \delta[t] \right) = S\desc^{M^!, *}_{\epsilon[t]} \left( e^{G_J^![t]}_{M^!}(\epsilon[t]) \sigma^{G_J^![t]_{\epsilon[t]}}_{J_t}(\delta_\epsilon) \right) \]
where $J_t$ is identified with an element of $\mathcal{J}^{G_J^![t]_{\epsilon[t]}}_{M^!_\epsilon}(B_J)$.

Note that the situation here is simpler than \cite[V.5.6]{MW16-1} since there are neither simply connected covers nor $z$-extensions. The transition factors $\lambda$ and the constants $c$ in \cite[p.556]{MW16-1} are thus trivial.

Define $\left( Z_{{G^!_J[t]}^\vee} : Z_{{G^![s]}^\vee} \right)$ as in \eqref{eqn:GJ-Gs}. All in all, Lemma \ref{prop:sigma-sigma} reduces to the assertion:
\begin{equation}\label{eqn:sigma-sigma-mid}\begin{split}
	\text{There exists}\quad & \tau \in SD_{\mathrm{unip}}(M^!_\epsilon) \otimes \mes(M^!_\epsilon)^\vee \quad \text{that maps to both} \\
	& \sigma^{G^![s]_{\epsilon[s]}}_{J_s}(\delta_\epsilon) \quad \text{and} \\
	& e^{G^!_J[t]}_{M^!}(\epsilon[t]) e^{G^![s]}_{M^!}(\epsilon[s])^{-1} \left( Z_{{G^!_J[t]}^\vee} : Z_{{G^![s]}^\vee} \right)^{-1} \sigma^{G^!_J[t]_{\epsilon[t]}}_{J_t}(\delta_\epsilon)
\end{split}\end{equation}
under the obvious quotient maps. This is a variant of \cite[V.5.6 (5)]{MW16-1}.

\begin{definition}
	Put
	\begin{align*}
		G^1 & := G^![s]_{\epsilon[s], J_s}, \\
		G^2 & := G^!_J[t]_{\epsilon[t]}.
	\end{align*}
\end{definition}

It has been shown in \cite[II.3.3]{MW16-1} that $S\mathrm{Ann}^{G^![s]_{\epsilon[s], J_s}}_{\mathrm{unip}} \subset S\mathrm{Ann}^{G^![s]_{\epsilon[s]}}_{\mathrm{unip}}$. By \cite[V.3.2 (3) and p.315]{MW16-1} we obtain
\[ \left( Z_{\check{G}^1}^{\Gamma_F} : Z_{{G^![s]_{\epsilon[s]}}^\vee}^{\Gamma_F} \right)^{-1} \sigma^{G^1}_{J_s}(\delta_\epsilon) \mapsto \sigma^{G^![s]_{\epsilon[s]}}_{J_s}(\delta_\epsilon). \]
In view of \eqref{eqn:sigma-sigma-mid}, Lemma \ref{prop:sigma-sigma} is reduced to proving:
\begin{equation}\label{eqn:sigma-sigma-final}\begin{split}
	\text{There exists}\quad & \tau \in SD_{\mathrm{unip}}(M^!_\epsilon) \otimes \mes(M^!_\epsilon)^\vee \quad \text{that maps to both} \\
	& \sigma^{G^1}_{J_s}(\delta_\epsilon) \quad \text{and} \\
	& \left( Z_{\check{G}^1}^{\Gamma_F} : Z_{{G^![s]_{\epsilon[s]}}^\vee}^{\Gamma_F} \right) e^{G^!_J[t]}_{M^!}(\epsilon[t]) e^{G^![s]}_{M^!}(\epsilon[s])^{-1} \left( Z_{{G^!_J[t]}^\vee} : Z_{{G^![s]}^\vee} \right)^{-1}
	\sigma^{G^2}_{J_t}(\delta_\epsilon).
\end{split}\end{equation}
This is a variant of \cite[V.5.6 (6)]{MW16-1}. Also note that $J_s$ and $J_t$ are the same lattice in $\mathfrak{a}_{M^!}^*$.

Before proving \eqref{eqn:sigma-sigma-final}, let us prove Lemma \ref{prop:J-roots} first.

\begin{proof}[Proof of Lemma \ref{prop:J-roots}]
	Suppose that $J \in \mathcal{J}^G_M \smallsetminus \{J_{\max}\}$ and $s \in \Endo_{\mathbf{M}^!}(\tilde{G})$ maps to $t \in \Endo_{\mathbf{M}^!}(\widetilde{G_J})$.	For each $\alpha \in \Sigma(G^![s]_{\epsilon[s]}, T^!)$, we write
	\[ \beta := B(\alpha)^{-1} \alpha \; \in X^*(T^!_{\overline{F}}) \otimes \Q. \]
	By our choice of $B$-functions in \S\ref{sec:local-geometric}, $\alpha \mapsto \beta$ transforms a root system of type $\mathrm{B}_n$ into type $\mathrm{C}_n$, and $\alpha$ is short if and only if $\beta$ is long.
	
	Hereafter, we assume $\tilde{G} = \Mp(W)$ for the ease of notations. In this case $Z_{\tilde{G}^\vee}^\circ = Z_{\widetilde{G_J}^\vee}^\circ = \{1\}$. Viewing $s$ and $t$ as elements in $\tilde{M}^\vee$, we infer from $s \mapsto t$ that $s = t$.
	
	By Remark \ref{rem:diagram-dual}, the chosen diagram $(\epsilon, B^!, T^!, B, T, \eta)$ induces $X^*(T^!_{\overline{F}}) \simeq X^*(T_{\overline{F}}) \simeq X_*(T^\vee)$ where $T^\vee \subset \tilde{G}^\vee$ is part of the root datum for $\tilde{G}^\vee$. In what follows, $\alpha \in X^*(T^!_{\overline{F}}) \simeq X_*(T^\vee)$ denotes a priori a coroot for $\tilde{G}^\vee \supset T^\vee$,. Denote by $\check{\alpha} \in X^*(T^\vee)$ the corresponding root. Then the image of $\beta$ by $X^*(T_{\overline{F}}) \simeq X^*(T^!_{\overline{F}})$ is in $\Sigma(G, T)$. The sets $\Sigma(G^i, T^!)$ arise from these $\alpha$ as follows ($i=1,2$).
	\begin{itemize}
		\item The roots of $G^1 \supset T^!$ are characterized by
		\begin{enumerate}
			\item $\check{\alpha}(s) = 1$ (so that $\alpha \in \Sigma(G^![s], T^!)$),
			\item $\alpha(\epsilon[s]) = 1$ (so that $\alpha \in \Sigma(G^![s]_{\epsilon[s]}, T^!)$),
			\item $\beta|_{\mathfrak{a}_{M^!}} \in R_{J_s}$ (so that $\alpha \in \Sigma(G^1, T^!)$); this is tantamount to $\xi^*(\beta|_{\mathfrak{a}_{M^!}}) \in R_J$.
		\end{enumerate}
		Moreover, let $\Psi: X^*(T^!_{\overline{F}}) \rightiso X_*(T^!_{\overline{F}})$ be induced from the invariant quadratic form chosen in \S\ref{sec:linear-algebra} and the diagram; its effect is $e^*_i \mapsto e_{i,*}$ in the standard bases of $X^*(T^\vee)$ and $X_*(T^\vee)$ in \S\ref{sec:linear-algebra}. The coroot of $\alpha \in \Sigma(G^1, T^!)$ can be taken on the level of $\Sigma(G^![s], T^!)$.
		
		Via the descriptions of $B$ and $G^![s]$ (a product of odd $\SO$), the coroot of $\alpha$ is seen to be $B(\alpha)^{-1} \Psi(\alpha)$; i.e.\ short roots become long.
	
		\item The roots of $G^2 \supset T^!$ are characterized by
		\begin{enumerate}
			\item $\xi^*(\beta|_{\mathfrak{a}_{M^!}}) \in R_J$ (so that $\beta \in \Sigma(G_J, T)$, or equivalently $\alpha \in \check{\Sigma}(\widetilde{G_J}^\vee, T^\vee)$),
			\item $\check{\alpha}(t) = 1$ (so that $\alpha \in \Sigma(G_J^![t], T^!)$),
			\item $\alpha(\epsilon[t]) = 1$ (so that $\alpha \in \Sigma(G^2, T^!)$).
		\end{enumerate}
	
		The coroot of $\alpha \in \Sigma(G^2, T^!)$ can be taken in $\Sigma(G_J^![t], T^!)$. Since $G^!_J[t]$ is a product of odd $\SO$ and the effect of $B_J(\cdot)^{-1}$ is to make short roots long, the coroot of $\alpha$ is seen to be $B_J(\alpha)^{-1} \Psi(\alpha)$.

		\item As $s=t$ and $\epsilon[s] = \epsilon[t]$, it follows from the above that $\Sigma(G^1, T^!) = \Sigma(G^2, T^!)$.
	\end{itemize}
	
	We contend that the vertical arrows below of rescaling
	\[\begin{tikzcd}
		\Sigma^J\left( G^![s]_{\epsilon[s]}, T^! \right) \arrow[d, "{\cdot B(\cdot)^{-1}}"'] \arrow[equal, r, "\eqref{eqn:J-G1}"] & \Sigma(G^1, T^!) \arrow[equal, r] & \Sigma(G^2, T^!) \arrow[d, "{\cdot B_J(\cdot)^{-1}}"] \\
		\Sigma^J\left( G^![s]_{\epsilon[s]}, T^!, B \right) & & \Sigma(G^2, T^!, B_J)
	\end{tikzcd}\]
	coincide. Indeed, this follows from the description of $G_J$ (Propositions \ref{prop:GJ-Sp}, \ref{prop:GJ-SO}) and the effects of $B(\cdot)^{-1}$ and $B_J(\cdot)^{-1}$ (see the earlier discussion on coroots). Lemma \ref{prop:J-roots} follows.
\end{proof}

Here is a by-product of the proof: we have actually shown that
\[ \Sigma(G^1, T^!) = \Sigma(G^2, T^!), \quad \check{\Sigma}(G^1, T^!) = \check{\Sigma}(G^2, T^!). \]

\begin{proof}[Proof of Lemma \ref{prop:sigma-sigma}]
	In order to verify \eqref{eqn:sigma-sigma-final}, it suffices to prove:
	\begin{enumerate}[(i)]
		\item there exists an isomorphism $G^1 \simeq G^2$, compatible with the embeddings of $M^!_{\epsilon[s]} = M^!_{\epsilon[t]}$ in both sides;
		\item we have
		\[ \left( Z_{\check{G}^1}^{\Gamma_F} : Z_{{G^![s]_{\epsilon[s]}}^\vee}^{\Gamma_F} \right) e^{G_J^![t]}_{M^!}(\epsilon[t]) e^{G^![s]}_{M^!}(\epsilon[s])^{-1} \left( Z_{{G^!_J[t]}^\vee} : Z_{{G^![s]}^\vee} \right)^{-1} = 1. \]
	\end{enumerate}

	To prove (i), observe that $G^1$ and $G^2$ share the maximal $F$-torus $T^!$, and have the same roots and coroots; this is the by-product of the proof of Lemma \ref{prop:J-roots}.
	
	As for the compatibility of the embeddings of $M^!_{\epsilon[s]} = M^!_{\epsilon[t]}$, note that the images of $\Sigma(M^!, T^!)$ in $\Sigma(G^i, T^!)$ have the same characterizations ($i = 1, 2$), namely $\beta|_{A_{M^!}} = 1$, $\check{\alpha}(s) = \check{\alpha}(t) = 1$, and $\alpha(\epsilon[s]) = \alpha(\epsilon[t]) = 1$ in the notation of the proof of Lemma \ref{prop:J-roots}.
	
	To prove (ii), we may assume $\tilde{G} = \Mp(W)$. Recall that $Z_{G^![s]^\vee} \subset Z_{G^!_J[t]^\vee}$ as subgroups of $Z_{{M^!}^\vee}$ by the proof of Lemma \ref{prop:J-roots}. Consider the diagram
	\[\begin{tikzcd}[column sep=small]
		Z_{{M^!}^\vee} \bigg/ Z_{{G^![s]}^\vee} \arrow[r, "\varphi"] & Z_{{M^!_{\epsilon[s]}}^\vee}^{\Gamma_F} \bigg/ Z_{{G^![s]_{\epsilon[s]}}^\vee}^{\Gamma_F} \arrow[r, "\text{quot}"] & Z_{{M^!_{\epsilon[s]}}^\vee}^{\Gamma_F} \bigg/ Z_{{G^![s]_{\epsilon[s], J_s}}^\vee}^{\Gamma_F} \arrow[equal, r] & Z_{{M^!_{\epsilon[s]}}^\vee}^{\Gamma_F} \bigg/ Z_{\check{G}^1}^{\Gamma_F} \arrow[equal, d, "\because \;\text{(i)}" inner sep=0.6em] \\
		Z_{{M^!}^\vee} \bigg/ Z_{\tilde{G}^\vee}^\circ \arrow[u, "\text{quot}"] \arrow[r, "\text{quot}"'] & Z_{{M^!}^\vee} \bigg/ Z_{{G_J^![t]}^\vee} \arrow[r, "\varphi"'] & Z_{{M^!_{\epsilon[t]}}^\vee}^{\Gamma_F} \bigg/ Z_{{G_J^![t]_{\epsilon[t], J_t}}^\vee}^{\Gamma_F} \arrow[equal, r] & Z_{{M^!_{\epsilon[t]}}^\vee}^{\Gamma_F} \bigg/ Z_{\check{G}^2}^{\Gamma_F}
	\end{tikzcd}\]
	where $\varphi$ stands for the homomorphisms in Definition \ref{def:eLRmu} and $\text{quot}$ denotes the evident quotient homomorphism. All arrows are surjective with finite kernels.
	
	The diagram commutes. Computing the cardinality of the kernel in two ways, we obtain
	\[ \left( Z_{{G^![s]}^\vee} : Z_{\tilde{G}^\vee}^\circ \right) e^{G^![s]}_{M^!}(\epsilon[s])^{-1} \left( Z_{\check{G}^1}^{\Gamma_F} : Z_{{G^![s]_{\epsilon[s]}}^\vee}^{\Gamma_F} \right)
	= \left( Z_{{G^!_J[t]}^\vee} : Z_{\tilde{G}^\vee}^\circ \right) e^{G_J^![t]}_{M^!}(\epsilon[t])^{-1} , \]
	hence
	\[ e^{G^![s]}_{M^!}(\epsilon[s])^{-1} \left( Z_{\check{G}^1}^{\Gamma_F} : Z_{{G^![s]_{\epsilon[s]}}^\vee}^{\Gamma_F} \right) = \left( Z_{{G^!_J[t]}^\vee} : Z_{{G^![s]}^\vee} \right) e^{G_J^![t]}_{M^!}(\epsilon[t])^{-1} . \]
	This is exactly (ii).
\end{proof}

\begin{remark}
	In the counterpart of (ii) in \cite[V.5.6]{MW16-1}, one has to treat a case of type $\mathrm{A}_{2n}$ in which $G^1$ is related to $G^2$ by nonstandard endoscopy. Here we are in the ``typical case'', for which the nonstandard endoscopic datum in question is tautological. 
\end{remark}

\section{Summary for general local fields}\label{sec:rho-main-summary}
As before, let $\rev: \tilde{G} \to G(F)$ be of metaplectic type and let $M \subset G$ be a Levi subgroup. Here $F$ can be any local field of characteristic zero.

The key observation is that the methods of \S\ref{sec:ext-Arch} also work for non-Archimedean $F$. Given such a local field $F$, the method in \S\S\ref{sec:ext-Arch}--\ref{sec:ext-Arch-proof} can be simplified in the following manner.
\begin{enumerate}
	\item One replaces all $D_{\text{tr-orb}, -}$, $SD_{\text{tr-orb}}$ by $D_{\mathrm{geom}, -}$, $SD_{\mathrm{geom}}$, respectively. The objects (O1) --- (O4) are thus all unconditionally defined.
	\item Among the conditions (C1) --- (C8), only (C2) (i.e.\ stabilization of $I_{\tilde{M}}(\tilde{\gamma}, \cdot)$) and (C3) (i.e.\ stabilization of $\rho^{\tilde{G}}_J$) need proof. Again, we are reduced to proving that
	\[ \rho^{\widetilde{G_J}, \Endo}_J(\mathbf{M}^!, \delta) \mapsto \rho^{\tilde{G}, \Endo}_J(\mathbf{M}^!, \delta) \]
	for all $J \in \mathcal{J}^G_M$, $\mathbf{M}^! \in \Endo_{\elli}(\tilde{M})$, elliptic $\mathcal{O}^! \subset M^!(F)$ and $\delta \in SD_{\mathrm{geom}}(M^!, \mathcal{O}^!) \otimes \mes(M^!)^\vee$. In fact, the arguments simplify in the non-Archimedean case.
	\item To be more specific, the results on harmonic analysis are already available in the non-Archimedean case. For example the reduction of local geometric Theorem \ref{prop:local-geometric} from $G$-equisingular to $G$-regular case is in Proposition \ref{prop:local-geometric-regular-reduction}; alternatively, one can appeal to the upcoming Corollary \ref{prop:local-geometric-regular-reduction-nonArch} via Shalika germs. In the proof of Lemma \ref{prop:ext-Arch-red-ell}, the usage of \cite[V.4.4 Lemme]{MW16-1} can also be replaced by the easier non-Archimedean version established in \cite[I.5.10]{MW16-1}, based on the property that $R^!_\epsilon = M^!_\epsilon$.
	
	The real hardcore of the arguments in \S\ref{sec:ext-Arch-proof} is combinatorial, and the proofs carry over verbatim.
\end{enumerate}

Of course, the $G$-regular case of local geometric theorem (Hypothesis \ref{hyp:local-geometric-nonArch}) is still needed as input. The inductive assumptions also remain in force, namely: all the results about $I^{\tilde{G}}_{\tilde{M}}$, $\rho^{\tilde{G}}_J$, etc.\ are established in either of the following cases:
\begin{itemize}
	\item $\tilde{G}$ is replaced by a group of metaplectic type of smaller dimension;
	\item $M$ is replaced by some strictly larger Levi subgroup of $G$ (the case $M = G$ is already known: $I^{\tilde{G}}_{\tilde{G}}$ are the invariant orbital integrals, and $\rho^{\tilde{G}}_{\emptyset} = \identity$).
\end{itemize}

\begin{theorem}\label{prop:local-geometric-reduction}
	Let $F$ be a local field of characteristic zero. Make the inductive assumptions above, and assume that for all $\tilde{\gamma} \in D_{\mathrm{orb}, -}(\tilde{M}) \otimes \mes(M)^\vee$ such that $\Supp(\tilde{\gamma}) \subset \tilde{M}_{G\text{-reg}}$, we have
	\[ I^{\tilde{G}, \Endo}_{\tilde{M}}(\tilde{\gamma}, \cdot) = I^{\tilde{G}}_{\tilde{M}}(\tilde{\gamma}, \cdot). \]
	Then:
	\begin{enumerate}
		\item if $F$ is non-Archimedean, the displayed equality extends to all $\tilde{\gamma} \in D_{\mathrm{geom}, -}(\tilde{M}) \otimes \mes(M)^\vee$;
		\item if $F$ is Archimedean, the displayed equality extends to all
		\[ \tilde{\gamma} \in \left( D_{\text{tr-orb}, -}(\tilde{M}) + D_{\text{geom, $G$-equi}, -}(\tilde{M})\right) \otimes \mes(M)^\vee . \]
	\end{enumerate}
	In particular, the statement of Theorem \ref{prop:local-geometric} holds true under these assumptions.
\end{theorem}
\begin{proof}
	This is merely a summary of the discussions above. In the Archimedean case, the fact that $I^{\tilde{G}}_{\tilde{M}}$ can be well-defined on $D_{\text{tr-orb}, -}(\tilde{M}) + D_{\text{geom, $G$-equi}, -}(\tilde{M})$ is the content of (C8).
\end{proof}

\section{Shalika germs}\label{sec:Shalika-germs}
Fix a non-Archimedean local field $F$ of characteristic zero.

\subsection{Basics}
The Shalika germs for orbital integrals on coverings are discussed in \cite[\S 4]{Li12b}. Here we follow the paradigm of \cite[II.2.1]{MW16-1} to rephrase the theory.

Let $\rev: \tilde{G} \twoheadrightarrow G(F)$ be a group of metaplectic type, and let $\mathcal{O}$ be a semisimple conjugacy class in $G(F)$. The Shalika germ in question is a germ of linear maps
\[ g_{\mathcal{O}}: D_{\mathrm{geom}, -}(\tilde{G}) \otimes \mes(G)^\vee \to D_{\mathrm{geom}, -}(\tilde{G}, \mathcal{O}) \otimes \mes(G)^\vee \]
characterized by
\[ I^{\tilde{G}}(\tilde{\gamma}, f) = I^{\tilde{G}}(g_{\mathcal{O}}(\tilde{\gamma}), f) \]
for all $f \in \orbI_{\asp}(\tilde{G}) \otimes \mes(G)$ and $\tilde{\gamma} \in D_{\mathrm{geom},-}(\tilde{G}) \otimes \mes(G)^\vee$ such that $\Supp(\tilde{\gamma})$ is sufficiently close to $\rev^{-1}(\mathcal{O})$. Being a germ means: two such linear maps are identified if they coincide when $\Supp(\tilde{\gamma}) \subset \rev^{-1}(\mathcal{V})$, where $\mathcal{V}$ is an invariant neighborhood of $\mathcal{O}$.

The \emph{separation property} of Shalika germs is rephrased as follows. For every $\tau \in D_{\mathrm{geom}, -}(\tilde{G}, \mathcal{O}) \otimes \mes(G)^\vee$, there exists $\tilde{\gamma} \in D_{\mathrm{geom}, -}(\tilde{G}) \otimes \mes(G)^\vee$, with $\Supp(\tilde{\gamma}) \subset \tilde{G}_{\mathrm{reg}}$ and close to $\rev^{-1}(\mathcal{O})$, such that $g_{\mathcal{O}}(\tilde{\gamma}) = \tau$. In fact, $\tilde{G}_{\mathrm{reg}}$ may be replaced by $\rev^{-1} (\mathcal{W})$ where $\mathcal{W}$ is any open dense invariant subset of $G_{\mathrm{reg}}(F)$.

By linearity, the formalism extends to the case that $\mathcal{O}$ is a finite union of semisimple conjugacy classes.

The stable case is analogous. Let $G^!$ be a quasisplit $F$-group and let $\mathcal{O}^!$ be a stable semisimple conjugacy class in $G^!(F)$. By replacing $I^{\tilde{G}}$ by $S^{G^!}$, we have the germ of linear maps
\[ Sg_{\mathcal{O}}: SD_{\mathrm{geom}}(G^!) \otimes \mes(G^!)^\vee \to SD_{\mathrm{geom}}(G^!, \mathcal{O}^!) \otimes \mes(G^!)^\vee \]
with a similar characterization and separation property.

We also need the Shalika germs of invariant weighted orbital integrals. The uncovered case is due to Arthur \cite[\S 9]{Ar88LB}. Since the combinatorial ingredients in the weighted orbital integrals on $\tilde{G}$ are the same (see \S\ref{sec:orbint-weighted-nonArch}), whilst the remaining parts are based on harmonic analysis, Arthur's results on germs carry over to $\tilde{G}$.

Specifically, let $M \subset G$ be a Levi subgroup, $\mathcal{O}$ a semisimple conjugacy class in $M(F)$, and denote by $\mathcal{O}^G$ the semisimple conjugacy class in $G(F)$ containing $\mathcal{O}$. We follow \cite[II.2.3]{MW16-1} to define the Shalika germ as a germ of linear maps
\[ g^{\tilde{G}}_{\tilde{M}, \mathcal{O}}: D_{\mathrm{geom}, G\text{-equi}, -}(\tilde{M}) \otimes \mes(M)^\vee \to D_{\mathrm{geom}, -}(\tilde{G}, \mathcal{O}^G) \otimes \mes(G)^\vee \]
characterized recursively by
\[ I^{\tilde{G}}_{\tilde{M}}(\tilde{\gamma}, f) = \sum_{L \in \mathcal{L}(M)} I^{\tilde{G}}_{\tilde{L}}\left( g^{\tilde{L}}_{\tilde{M}, \mathcal{O}}(\tilde{\gamma}), f \right) \]
for all $f \in \orbI_{\asp}(\tilde{G}) \otimes \mes(G)$ and $\tilde{\gamma}$ with $\Supp(\tilde{\gamma})$ close to $\rev^{-1}(\tilde{O})$. When $M=G$ we recover the previous version.
\index{gGMO@$g^{\tilde{G}}_{\tilde{M}, \mathcal{O}}$}

Again, the definition extends to the case when $\mathcal{O}$ is a finite union of semisimple conjugacy classes.

As in Arthur's setting \cite[p.270, Remark]{Ar88LB}, we have the following crucial property:
\[ \mathcal{O}\; \text{is $G$-equisingular and $M \neq G$} \implies g^{\tilde{G}}_{\tilde{M}, \mathcal{O}} = 0 . \]

The weighted germs also satisfy the following descent formula.
\begin{proposition}\label{prop:Shalika-descent}
	Let $L \in \mathcal{L}^G(M)$. Write $\tilde{\gamma}^{\tilde{L}} := \Ind^{\tilde{L}}_{\tilde{M}}(\tilde{\gamma})$ (see \eqref{eqn:parabolic-ind-dist}) and define $\mathcal{O}^L$ to be the semisimple conjugacy class in $L(F)$ containing $\mathcal{O}$. Then
	\[ g^{\tilde{G}}_{\tilde{L}, \mathcal{O}^L}\left( \tilde{\gamma}^{\tilde{L}} \right) = \sum_{L^\dagger \in \mathcal{L}(M)} d^G_M(L, L^\dagger) g^{\tilde{L}^\dagger}_{\tilde{M}, \mathcal{O}}\left( \tilde{\gamma} \right)^{\tilde{G}} . \]
	whenever $\Supp(\tilde{\gamma})$ is sufficiently close to $\rev^{-1}(\mathcal{O})$.
\end{proposition}
\begin{proof}
	Same as \cite[II.2.11]{MW16-1}: this is based on the descent formula of Proposition \ref{prop:orbint-weighted-descent-nonArch}.
\end{proof}

\begin{remark}\label{rem:Shalika-rho}
	Shalika germs are related to the germs $\rho^{\tilde{G}}_J$ in \S\ref{sec:rho-sigma} by the equivalence
	\[ g^{\tilde{G}}_{\tilde{M}, \mathcal{O}}(a\tilde{\gamma}) \sim \sum_{J \in \mathcal{J}^G_M} \rho^{\tilde{G}}_J(\tilde{\gamma}, a)^{\tilde{G}}, \quad \tilde{\gamma} \in D_{\mathrm{geom}, -}(\tilde{M}, \mathcal{O}) \otimes \mes(M)^\vee. \]
	between germs of functions in $a \in A_M(F)$ in general position, $a \to 1$. This is immediate from Definition--Proposition \ref{def:rho-germ} (ii).
\end{remark}

Consider now the stable setting, namely $G^!$ is a quasisplit $F$-group with Levi subgroup $M^!$. We follow \cite[II.2.4]{MW16-1} to define
\[ SD_{\mathrm{geom}, G^!\text{-equi}}(M^!) := SD_{\mathrm{geom}}(M^!) \cap D_{\mathrm{geom}, G^!\text{-equi}}(M^!). \]
Suppose that a system of $B$-functions have been prescribed on $G^!$. Let $\mathcal{O}^!$ be a stable semisimple conjugacy class in $M^!(F)$, and let $\mathcal{O}^{G^!}$ be the stable semisimple conjugacy class in $G^!(F)$ containing $\mathcal{O}^!$. In view of \cite[II.2.4, Théorème; II.2.5 Proposition]{MW16-1}, we have the germ of linear maps
\[ Sg^{G^!}_{M^!, \mathcal{O}^!}(B): SD_{\mathrm{geom}, G^!\text{-equi}}(M^!) \otimes \mes(M^!)^\vee \to SD_{\mathrm{geom}}(G^!, \mathcal{O}^{G^!}) \otimes \mes(G^!)^\vee \]
satisfying
\[ S^{G^!}_{M^!}(\delta, B, f) = \sum_{L^! \in \mathcal{L}(M^!)} I^{G^!}_{L^!}\left( Sg^{L^!}_{M^!, \mathcal{O}^!}(\delta, B|_{L^!}), B, f \right) \]
for all $f$ and $\delta$ with $\Supp(\delta)$ close to $\mathcal{O}^!$.
\index{SgOB@$Sg^{G^{"!}}_{M^{"!}, \mathcal{O}^{"!}}(B)$}

Again, when $M^! = G^!$ we recover $Sg_{\mathcal{O}^!}$. Also,
\[ \mathcal{O}^!\; \text{is $G^!$-equisingular and $M^! \neq G^!$} \implies Sg^{G^!}_{M^!, \mathcal{O}^!} = 0 . \]

\subsection{Statement of the matching}
Let $\tilde{M} \subset \tilde{G}$ as before. Fix a stable semisimple conjugacy class $\mathcal{O}$ in $M(F)$. Denote by $\mathcal{O}^G$ the stable semisimple conjugacy class in $G(F)$ containing $\mathcal{O}$.

\begin{definition}\label{def:gEndo-1}
	\index{gEndoMOshrek@$g^{\tilde{G}, \Endo}_{\tilde{M}, \mathcal{O}}(\mathbf{M}^{"!}, \delta)$}
	Let $\mathbf{M}^! \in \Endo_{\elli}(\tilde{M})$ and let $\mathcal{O}^!$ be the finite union of stable semisimple conjugacy classes in $M^!(F)$ that map to $\mathcal{O}$. Define the germ of linear maps
	\[ g^{\tilde{G}, \Endo}_{\tilde{M}, \mathcal{O}}(\mathbf{M}^!): SD_{\mathrm{geom}, \tilde{G}\text{-equi}}(M^!) \otimes \mes(M^!)^\vee \to D_{\mathrm{geom}, -}(\tilde{G}, \mathcal{O}^G) \otimes \mes(G)^\vee \]
	by the formula
	\[ g^{\tilde{G}, \Endo}_{\tilde{M}, \mathcal{O}}(\mathbf{M}^!, \delta) := \sum_{s \in \Endo_{\mathbf{M}^!}(\tilde{G})} i_{M^!}(\tilde{G}, G^![s]) \trans_{\mathbf{G}^![s], \tilde{G}}\left( Sg^{G^![s]}_{M^!, \mathcal{O}^![s]}(\delta[s], B^{\tilde{G}}) \right) \]
	where $\Supp(\delta)$ is assumed to be close to $\mathcal{O}^!$.
\end{definition}

\begin{lemma}\label{prop:gEndo-aux}
	In the circumstance above, we have
	\[ I^{\tilde{G}, \Endo}_{\tilde{M}}(\mathbf{M}^!, \delta, f) = \sum_{L \in \mathcal{L}(M)} I^{\tilde{G}, \Endo}_{\tilde{L}} \left( g^{\tilde{L}, \Endo}_{\tilde{M}, \mathcal{O}}(\mathbf{M}^!, \delta), f \right) \]
	for all $f \in \orbI_{\asp}(\tilde{G}) \otimes \mes(G)$ and all $\delta$ with $\Supp(\delta)$ close to $\mathcal{O}^!$, where the left hand side is as Definition \ref{def:IEndo-1}.
\end{lemma}
\begin{proof}
	By the properties of stable Shalika germs, $I^{\tilde{G}, \Endo}_{\tilde{M}}(\mathbf{M}^!, \delta, f)$ equals
	\[ \sum_{s, L^!} i_{M^!}(\tilde{G}, G^![s]) S^{G^![s]}_{L^!}\left( Sg^{L^!}_{M^!, \mathcal{O}^![s]}(\delta[s], B^{\tilde{G}}|_{L^!}), B^{\tilde{G}}, f^{G^![s]} \right) \]
	where $s \in \Endo_{\mathbf{M}^!}(\tilde{G})$, $L^! \in \mathcal{L}^{G^![s]}(M^!)$ and $f^{G^![s]} = \Trans_{\mathbf{G}^![s], \tilde{G}}(f)$. Apply Lemma \ref{prop:sL-Ls} with $R = M$ to transform this into a sum over $(L, s)$ where $L \in \mathcal{L}^G(M)$ and $s \in \Endo_{\mathbf{M}^!}(\tilde{G})$ factors uniquely through a pair $(s^L, s_L)$ with
	\[ s^L \in \Endo_{\mathbf{M}^!}(\tilde{L}), \quad s_L \in \Endo_{\mathbf{L}^![s^L]}(\tilde{G}), \]
	so that $L^!$ (resp.\ $\mathbf{G}^![s]$) corresponds to $L^![s^L]$ (resp.\ $\mathbf{G}^![s_L]$).
	
	Lemma \ref{prop:i-transitivity} asserts $i_{M^!}(\tilde{G}, G^![s]) = i_{M^!}(\tilde{L}, L^!) i_{L^!}(\tilde{G}, G^![s_L])$. Lemma \ref{prop:B-Levi} implies $B^{\tilde{G}}|_{L^!} = B^{\tilde{L}}$. We also have $\delta[s] = \delta[s^L][s_L]$; as $z[s_L]$ is central in $L^!(F)$, translation by $z[s_L]$ commutes with $Sg^{L^!}_{M^!, \mathcal{O}^![\star]}(B^{\tilde{L}})$. All in all, $I^{\tilde{G}, \Endo}_{\tilde{M}}(\mathbf{M}^!, \delta, f)$ becomes
	\begin{multline*}
		\sum_{\substack{L \in \mathcal{L}^G(M) \\ s^L \in \Endo_{\mathbf{M}^!}(\tilde{L}) \\ L^! := L^![s^L] }} i_{M^!}(\tilde{L}, L^!) \sum_{\substack{s_L \in \Endo_{\mathbf{L}^!}(\tilde{G}) \\ G^! := G^![s_L] }} i_{L^!}(\tilde{G}, G^!)
		S^{G^!}_{L^!} \left( Sg^{L^!}_{M^!, \mathcal{O}^![s^L]}(\delta[s^L], B^{\tilde{L}})[s_L], B^{\tilde{G}}, f^{G^!} \right) \\
		= \sum_{L \in \mathcal{L}^G(M)} \sum_{\substack{s^L \in \Endo_{\mathbf{M^!}(\tilde{L})} \\ L^! := L^![s^L] }} i_{M^!}(\tilde{L}, L^!) I^{\tilde{G}, \Endo}_{\tilde{L}} \left( \trans_{\mathbf{L}^!, \tilde{L}}\left( Sg^{L^!}_{M^!, \mathcal{O}^![s^L]}(\delta[s^L], B^{\tilde{L}}) \right), f \right) \\
		= \sum_{L \in \mathcal{L}^G(M)} I^{\tilde{G}, \Endo}_{\tilde{L}} \left( g^{\tilde{L}, \Endo}_{\tilde{M}, \mathcal{O}}(\mathbf{M}^!, \delta), f \right).
	\end{multline*}
	This is what we seek.
\end{proof}

\begin{definition-proposition}\label{def:gEndo}
	\index{gEndoMO@$g^{\tilde{G}, \Endo}_{\tilde{M}, \mathcal{O}}(\tilde{\gamma})$}
	Express $\tilde{\gamma} \in D_{\mathrm{geom}, G\text{-equi}, -}(\tilde{M}) \otimes \mes(M)^\vee$ with $\Supp(\tilde{\gamma})$ sufficiently close to $\rev^{-1}(\mathcal{O})$ as
	\[ \tilde{\gamma} = \sum_{\mathbf{M}^!} \trans_{\mathbf{M}^!, \tilde{M}}(\delta_{\mathbf{M}^!}) \]
	where $(\mathbf{M}^!, \delta_{\mathbf{M}^!})$ is as Definition \ref{def:gEndo-1} for each $\mathbf{M}^! \in \Endo_{\elli}(\tilde{M})$. Put
	\[ g^{\tilde{G}, \Endo}_{\tilde{M}, \mathcal{O}}(\tilde{\gamma}) := \sum_{\mathbf{M}^!} g^{\tilde{G}, \Endo}_{\tilde{M}, \mathcal{O}}(\mathbf{M}^!, \delta_{\mathbf{M}^!}).  \]
	Then $g^{\tilde{G}, \Endo}_{\tilde{M}, \mathcal{O}}(\tilde{\gamma})$ depends only on $\tilde{\gamma}$, and it satisfies
	\[ I^{\tilde{G}, \Endo}_{\tilde{M}}(\tilde{\gamma}, f) = \sum_{L \in \mathcal{L}(M)} I^{\tilde{G}, \Endo}_{\tilde{L}}\left( g^{\tilde{L}, \Endo}_{\tilde{M}, \mathcal{O}}(\tilde{\gamma}), f \right), \quad f \in \orbI_{\asp}(\tilde{G}) \otimes \mes(G). \]
\end{definition-proposition}
\begin{proof}
	By Lemma \ref{prop:gEndo-aux}, we have
	\begin{align*}
		I^{\tilde{G}, \Endo}_{\tilde{M}}(\tilde{\gamma}, f) & = \sum_{\mathbf{M}^!} I^{\tilde{G}, \Endo}_{\tilde{M}}(\mathbf{M}^!, \delta_{\mathbf{M}^!}, f) \\
		& = \sum_{L \in \mathcal{L}(M)} I^{\tilde{G}, \Endo}_{\tilde{L}}\left( \sum_{\mathbf{M}^!} g^{\tilde{L}, \Endo}_{\tilde{M}, \mathcal{O}}(\mathbf{M}^!, \delta_{\mathbf{M}^!}), f \right),
	\end{align*}
	the first equality being Definition--Proposition \ref{def:I-geom-Endo}. The first term depends only on $\tilde{\gamma}$, and so do those $\sum_{\mathbf{M}^!} g^{\tilde{L}, \Endo}_{\tilde{M}, \mathcal{O}}(\mathbf{M}^!, \delta_{\mathbf{M}^!})$ with $L \neq G$ by induction. Hence the required case $L=G$ follows by varying $f$.
\end{proof}

\begin{remark}\label{rem:Shalika-rho-Endo}
	In view of Proposition \ref{prop:rhoEndo-expansion}, we deduce inductively an obvious endoscopic analogue of Remark \ref{rem:Shalika-rho}, namely the equivalence
	\[ g^{\tilde{G}, \Endo}_{\tilde{M}, \mathcal{O}}(\mathbf{M}^!, \xi(a) \delta) \sim \sum_{J \in \mathcal{J}^G_M} \rho^{\tilde{G}, \Endo}_J(\mathbf{M}^!, \delta, a)^{\tilde{G}}, \quad \delta \in SD_{\mathrm{geom}, -}(\tilde{M}, \mathcal{O}) \otimes \mes(M)^\vee. \]
	between germs of functions in $a \in A_M(F)$ in general position, $a \to 1$; here $\xi: A_M \rightiso A_{M^!}$ is the natural identification from endoscopy.
\end{remark}

We can now state the matching of weighted Shalika germs.

\begin{theorem}[Cf.\ {\cite[III.8.5 Proposition]{MW16-1}}]\label{prop:germ-matching}
	We have the equality $g^{\tilde{G}, \Endo}_{\tilde{M}, \mathcal{O}} = g^{\tilde{G}}_{\tilde{M}, \mathcal{O}}$ between germs of linear maps from $D_{\mathrm{geom}, G\text{-equi}, -}(\tilde{M}) \otimes \mes(M)^\vee$ to $D_{\mathrm{geom}, -}(\tilde{G}, \mathcal{O}^G) \otimes \mes(G)^\vee$.
\end{theorem}

The simplest case $M=G$ of Theorem \ref{prop:germ-matching} will be established in Corollary \ref{prop:germ-matching-G}. The proof of the general case will be completed in \S\ref{sec:Shalika-matching}. As the first step, let us reduce it to the case when $\mathcal{O}$ is elliptic, using Proposition \ref{prop:Shalika-descent} and its endoscopic counterpart below.

\begin{proposition}\label{prop:gEndo-descent}
	Let $L \in \mathcal{L}^G(M)$ and let $\mathcal{O}^L$ be the semisimple conjugacy class in $L(F)$ containing $\mathcal{O}$. Then
	\[ g^{\tilde{G}, \Endo}_{\tilde{L}, \mathcal{O}^L}\left( \tilde{\gamma}^{\tilde{L}} \right) = \sum_{L^\dagger \in \mathcal{L}(M)} d^G_M(L, L^\dagger) g^{\tilde{L}^\dagger, \Endo}_{\tilde{M}, \mathcal{O}}\left( \tilde{\gamma} \right)^{\tilde{G}} \]
	whenever $\Supp(\tilde{\gamma})$ is sufficiently close to $\rev^{-1}(\mathcal{O})$.
\end{proposition}
\begin{proof}
	Same as \cite[II.2.12]{MW16-1}. Compared to the proof of Proposition \ref{prop:Shalika-descent}, the inputs here are the descent formula of Corollary \ref{prop:descent-orbint-Endo-2} and Definition--Proposition \ref{def:gEndo}.
\end{proof}

Similarly, the following property of $g^{\tilde{G}, \Endo}_{\tilde{M}, \mathcal{O}}$ mirrors that of $g^{\tilde{G}}_{\tilde{M}, \mathcal{O}}$.

\begin{lemma}\label{prop:gEndo-equisingular}
	Suppose that $\mathcal{O}$ is $G$-equisingular and $M \neq G$. Then $g^{\tilde{G}, \Endo}_{\tilde{M}, \mathcal{O}} = 0$.
\end{lemma}
\begin{proof}
	For every $\mathbf{M}^! \in \Endo_{\elli}(\tilde{M})$, let $\mathcal{O}^!$ be a stable semisimple conjugacy class in $M^!(F)$ with $\mathcal{O}^! \mapsto \mathcal{O}$. Then $\mathcal{O}^!$ is $\tilde{G}$-equisingular. Let $s \in \Endo_{\mathbf{M}^!}(\tilde{G})$. By \cite[II.2.4 (2)]{MW16-1}, we have
	\[ Sg^{G^![s]}_{M^!, \mathcal{O}^![s]}(B^{\tilde{G}}) = 0 \]
	since $\mathcal{O}^![s]$ is $G^![s]$-equisingular by Lemma \ref{prop:equisingular-endo}. Hence $g^{\tilde{G}, \Endo}_{\tilde{M}, \mathcal{O}} = 0$.
\end{proof}

\subsection{Applications to the local geometric theorem}
Conserve the notations and assumptions from \S\ref{sec:Shalika-germs}. Fix a semisimple conjugacy class $\mathcal{O}$ in $M(F)$.

The local geometric Theorem \ref{prop:local-geometric} asserts that $I^{\tilde{G}}_{\tilde{M}}(\tilde{\gamma}, f) = I^{\tilde{G}, \Endo}_{\tilde{M}}(\tilde{\gamma}, f)$ for all $\tilde{\gamma}$ and $f$. We assume inductively that
\begin{itemize}
	\item the statement of Theorem \ref{prop:local-geometric} holds true when $G$ (resp.\ $M$) is replaced by a strictly smaller (resp.\ larger) Levi subgroup;
	\item ditto for Theorem \ref{prop:germ-matching}.
\end{itemize}

\begin{lemma}\label{prop:gEndo-vs-local-geom}
	Let $\mathcal{D}$ be a subset of $D_{\mathrm{geom}, G\text{-equi}, -}(\tilde{M}) \otimes \mes(M)^\vee$ with the following separation property: for all $\tau \in D_{\mathrm{geom}, -}(\tilde{M}, \mathcal{O}) \otimes \mes(M)^\vee$, there exists $\tilde{\gamma} \in \mathcal{D}$ with $\Supp(\tilde{\gamma})$ arbitrarily close to $\rev^{-1}(\mathcal{O})$, such that $g^{\tilde{M}}_{\tilde{M}, \mathcal{O}}(\tilde{\gamma}) = \tau$.
	
	Suppose that $I^{\tilde{G}}_{\tilde{M}}(\tilde{\gamma}, \cdot) = I^{\tilde{G}, \Endo}_{\tilde{M}}(\tilde{\gamma}, \cdot)$ for all $\tilde{\gamma} \in \mathcal{D}$. Then the following are equivalent:
	\begin{enumerate}[(i)]
		\item $I^{\tilde{G}}_{\tilde{M}}(\tilde{\gamma}, \cdot) = I^{\tilde{G}, \Endo}_{\tilde{M}}(\tilde{\gamma}, \cdot)$ for all $\tilde{\gamma} \in D_{\mathrm{geom}, -}(\tilde{M}, \mathcal{O}) \otimes \mes(M)^\vee$;
		\item $g^{\tilde{G}}_{\tilde{M}, \mathcal{O}}(\tilde{\gamma}) = g^{\tilde{G}, \Endo}_{\tilde{M}, \mathcal{O}}(\tilde{\gamma})$ for all $\tilde{\gamma} \in \mathcal{D}$ with $\Supp(\tilde{\gamma})$ sufficiently close to $\rev^{-1}(\mathcal{O})$.
	\end{enumerate}
\end{lemma}
\begin{proof}
	In view of various inductive assumptions, for $\tilde{\gamma} \in \mathcal{D}$ with $\Supp(\tilde{\gamma})$ sufficiently close to $\rev^{-1}(\mathcal{O})$ and all $f$, the germ expansions for $I^{\tilde{G}}_{\tilde{M}}$ and $I^{\tilde{G}, \Endo}_{\tilde{M}}$ imply
	\begin{align*}
		0 & = I^{\tilde{G}}_{\tilde{M}}(\tilde{\gamma}, f) - I^{\tilde{G}, \Endo}_{\tilde{M}}(\tilde{\gamma}, f) \\
		& = I^{\tilde{G}}_{\tilde{M}}\left( g^{\tilde{M}}_{\tilde{M}, \mathcal{O}}(\tilde{\gamma}), f \right) - I^{\tilde{G}, \Endo}_{\tilde{M}}\left( g^{\tilde{M}}_{\tilde{M}, \mathcal{O}}(\tilde{\gamma}), f \right) \\
		& \quad + I^{\tilde{G}}_{\tilde{G}}\left( g^{\tilde{G}}_{\tilde{M}, \mathcal{O}}(\tilde{\gamma}) - g^{\tilde{G}, \Endo}_{\tilde{M}, \mathcal{O}}(\tilde{\gamma}), f \right).
	\end{align*}

	Assuming (i), we deduce that $g^{\tilde{G}}_{\tilde{M}, \mathcal{O}}(\tilde{\gamma}) - g^{\tilde{G}, \Endo}_{\tilde{M}, \mathcal{O}}(\tilde{\gamma})$ has vanishing orbital integrals for these $\tilde{\gamma}$, hence (ii) holds.
	
	Assuming (ii), we deduce that $I^{\tilde{G}}_{\tilde{M}}(g^{\tilde{M}}_{\tilde{M}, \mathcal{O}}(\tilde{\gamma}), \cdot) = I^{\tilde{G}, \Endo}_{\tilde{M}}(g^{\tilde{M}}_{\tilde{M}, \mathcal{O}}(\tilde{\gamma}), \cdot)$ for these $\tilde{\gamma}$. By the separation property, this implies (i).
\end{proof}

\begin{corollary}\label{prop:germ-matching-G}
	We have $g^{\tilde{G}}_{\tilde{G}, \mathcal{O}} = g^{\tilde{G}, \Endo}_{\tilde{G}, \mathcal{O}}$ for all stable semisimple conjugacy class $\mathcal{O}$ in $G(F)$.
\end{corollary}
\begin{proof}
	Take $M=G$ and $\mathcal{D} = D_{\mathrm{geom}, -}(\tilde{G}) \otimes \mes(G)^\vee$ in Lemma \ref{prop:gEndo-vs-local-geom}. The identity $I^{\tilde{G}}_{\tilde{G}}(\tilde{\gamma}, \cdot) = I^{\tilde{G}, \Endo}_{\tilde{G}}(\tilde{\gamma}, \cdot)$ is trivially true for all $\tilde{\gamma} \in \mathcal{D}$, and the matching of Shalika germs follows.
\end{proof}

We may also take $\mathcal{D}$ spanned by those $\tilde{\gamma}$ with supports in $\tilde{M}_{G\text{-reg}}$ and sufficiently close to $\rev^{-1}(\mathcal{O})$. This yields the following consequence.

\begin{corollary}\label{prop:local-geometric-regular-reduction-nonArch}
	Suppose that $\mathcal{O}$ is $G$-equisingular and fix an invariant neighborhood $\mathcal{V}$ of $\mathcal{O}$ in $M(F)$. If $I^{\tilde{G}}_{\tilde{M}}(\tilde{\gamma}, \cdot) = I^{\tilde{G}, \Endo}_{\tilde{M}}(\tilde{\gamma}, \cdot)$ for all $\tilde{\gamma} \in D_{\mathrm{geom}, -}(\tilde{M}) \otimes \mes(M)^\vee$ such that $\Supp(\tilde{\gamma}) \subset \rev^{-1}(\mathcal{V} \cap M_{G\text{-reg}}(F))$, then $I^{\tilde{G}}_{\tilde{M}}(\tilde{\gamma}, \cdot) = I^{\tilde{G}, \Endo}_{\tilde{M}}(\tilde{\gamma}, \cdot)$ for all $\tilde{\gamma} \in D_{\mathrm{geom}, -}(\tilde{M}, \mathcal{O}) \otimes \mes(M)^\vee$.
\end{corollary}
\begin{proof}
	The $G$-equisingularity implies $g^{\tilde{G}}_{\tilde{M}, \mathcal{O}} = 0 = g^{\tilde{G}, \Endo}_{\tilde{M}, \mathcal{O}}$ when $M \neq G$, by Lemma \ref{prop:gEndo-equisingular}. On the other hand, we have seen that $g^{\tilde{G}}_{\tilde{G}, \mathcal{O}} = g^{\tilde{G}, \Endo}_{\tilde{G}, \mathcal{O}}$. Lemma \ref{prop:gEndo-vs-local-geom} can thus be applied.
\end{proof}

This gives another proof of the non-Archimedean case of Proposition \ref{prop:local-geometric-regular-reduction}.

\section{Proof of the matching of Shalika germs}\label{sec:Shalika-matching}
Fix $\mathbf{M}^! \in \Endo_{\elli}(\tilde{M})$, a stable semisimple conjugacy class $\mathcal{O}^! \subset M^!(F)$ and take $\epsilon \in \mathcal{O}^!$ such that $M^!_\epsilon$ is quasisplit. Denote by $\mathcal{O}$ the image of $\mathcal{O}^!$, which is a stable semisimple conjugacy class in $M(F)$. For all $s \in \Endo_{\mathbf{M}^!}(\tilde{G})$, we define $\overline{M^!_\epsilon} \subset \overline{G^![s]_{\epsilon[s]}}$ by the convention in \S\ref{sec:nonstandard-endoscopy}.

The first and easiest step is to reduce Theorem \ref{prop:germ-matching} to the elliptic case.

\begin{proof}[Proof of Theorem \ref{prop:germ-matching} in the non-elliptic case]
	We may assume $\tilde{\gamma} = \trans_{\mathbf{M}^!, \tilde{M}}(\delta)$ for some
	\[ \delta \in SD_{\mathrm{geom}, \tilde{G}\text{-equi}}(M^!) \otimes \mes(M^!)^\vee \]
	with $\Supp(\delta)$ close to $\mathcal{O}^!$. Put
	\[ R^! := Z_{M^!}\left( A_{M^!_\epsilon} \right), \quad \mathcal{O}_{R^!} := \mathcal{O}^! \cap R^!(F). \]
	By construction, $R^! \hookrightarrow M^!$ is a Levi subgroup containing $\epsilon$, and $\mathcal{O}^! = (\mathcal{O}_{R^!})^{M^!}$. By Proposition \ref{prop:situation-easier}, these data extend uniquely (up to conjugacy) to
	\[\begin{tikzcd}[column sep=large]
		M^! \arrow[dashed, dash, r, "\text{ell.\ endo.}"] & \tilde{M} \\
		R^! \arrow[hookrightarrow, u, "\text{Levi}"] \arrow[dashed, dash, r, "\text{ell.\ endo.}"'] & \tilde{R} \arrow[hookrightarrow, u, "\text{Levi}"']
	\end{tikzcd} \qquad \mathbf{M}^! = \mathbf{M}^![t], \; t \in \Endo_{\mathbf{R}^!}(\tilde{M}). \]

	Observe that $R^!_\epsilon = M^!_\epsilon$. Indeed, every element of $M^!_\epsilon$ centralizes $A_{M^!_\epsilon}$, hence lies in $Z_{R^!}(\epsilon)$ and we conclude that $M^!_\epsilon \subset R^!_\epsilon$. Also note that $\delta \in \Image \left(S\desc^{M^!, *}_\epsilon\right)$. These observations entail $\delta = \nu^{M^!}$ for some $\nu \in SD_{\mathrm{geom}}(R^!) \otimes \mes(R^!)^\vee$ with $\Supp(\nu)$ close to $\mathcal{O}_{R^!}$. Cf.\ the first step of the proof in \cite[III.7.4 Proposition]{MW16-1}.

	Put $\nu[t^{-1}] := \nu \cdot z[t]^{-1}$ as usual, and suppose that $\mathcal{O}_{R^!}[t^{-1}] \mapsto \mathcal{O}_R$, where $\mathcal{O}_R$ is a stable semisimple conjugacy class in $R(F)$. Then $\mathcal{O} = (\mathcal{O}_R)^M$ and
	\[ \tilde{\gamma} = \trans_{\mathbf{M}^!, \tilde{M}}\left( \nu^{M^!} \right) = \left(\trans_{\mathbf{R}^!, \tilde{R}}( \nu[t^{-1}] )\right)^{\tilde{M}}. \]
	
	Suppose that $\mathcal{O}^!$ is not elliptic in $M^!$, then $R^!$ is a proper Levi subgroup of $M^!$. The descent formulas from Propositions \ref{prop:Shalika-descent}, \ref{prop:gEndo-descent} and inductive assumptions then imply $g^{\tilde{G}}_{\tilde{M}, \mathcal{O}}(\tilde{\gamma}) = g^{\tilde{G}, \Endo}_{\tilde{M}, \mathcal{O}}(\tilde{\gamma})$.
\end{proof}

In what follows, we assume $\epsilon$ is elliptic in $M^!$ and take $\eta \in \mathcal{O}$ such that $M_\eta$ is quasisplit. Note that $\eta$ is elliptic in $M$ by Lemma \ref{prop:ellipticity-transfer}. Define
\begin{equation}\label{eqn:Z-descent}
	\mathcal{Z} := \Ker\left[ \frac{Z_{\tilde{M}^\vee}^\circ}{Z_{\tilde{G}^\vee}^\circ} \xrightarrow{\text{descent}} \frac{Z_{M_\eta^\vee}^{\Gamma_F}}{Z_{G_\eta^\vee}^{\Gamma_F}} \right],
\end{equation}
where ``descent'' is the map in Remark \ref{rem:descent-homomorphism}. It is finite when $\eta$ is elliptic in $G$.

\begin{lemma}\label{prop:Shalika-matching-aux}
	Let $s \in \Endo_{\mathbf{M}^!}(\tilde{G})$. Suppose that $\eta$ is elliptic in $G^![s]$, then $\overline{G^![s]_{\epsilon[s]}}$ is an elliptic endoscopic group of $G_\eta$ if and only if $\epsilon[s]$ is elliptic in $G^![s]$.
\end{lemma}
\begin{proof}
	The following inequalities hold
	\[ \dim A_{\overline{G^!_{\epsilon[s]}}} = \dim A_{G^![s]_{\epsilon[s]}} \geq \dim A_{G^![s]} = \dim A_G = \dim A_{G_\eta}, \]
	and the assertion follows at once.
\end{proof}

\begin{lemma}\label{prop:Shalika-matching-aux0}
	Let $s \in \Endo_{\mathbf{M}^!}(\tilde{G})$. Suppose that $\epsilon[s]$ is elliptic in $G^![s]$. Then $\eta$ is elliptic in $G$ and $\overline{G^![s]_{\epsilon[s]}}$ is an elliptic endoscopic group of $G_\eta$. Let
	\[ c(s) := c^{\overline{G^![s]_{\epsilon[s]}}, G^![s]_{\epsilon[s]}}_{\overline{M^!_\epsilon}, M^!_\epsilon} \in 2^{\Z_{\leq 0}} \]
	be the factor in \cite[Proposition 5.3.3]{Li12a}. Then
	\begin{equation*}
		e^{G^![s]}_{M^!}(\epsilon[s]) i_{M^!}\left( \tilde{G}, G^![s]\right) = |\mathcal{Z}|^{-1} c(s) i_{\overline{M^!_\epsilon}}\left( G_\eta, \overline{G^![s]_{\epsilon[s]}} \right)
	\end{equation*}
	where $e^{G^![s]}_{M^!}(\epsilon[s])$ is as in Definition \ref{def:eLRmu}.
\end{lemma}
\begin{proof}
	Since $\epsilon[s] \leftrightarrow \eta$ with respect to $\mathbf{G}^![s]$, Lemma \ref{prop:ellipticity-transfer} implies $\eta$ is elliptic in $G$, and Lemma \ref{prop:Shalika-matching-aux} implies that $\overline{G^![s]_{\epsilon[s]}}$ is an elliptic endoscopic group of $G_\eta$. Consider the diagram:
	\[\begin{tikzcd}[row sep=tiny, column sep=large]
		& Z_{(M^!)^\vee} \big/ Z_{(G^![s])^\vee} \arrow[r, "\text{descent}"] & Z_{(M^!_\epsilon)^\vee}^{\Gamma_F} \big/ Z_{(G^![s]_{\epsilon[s]})^\vee}^{\Gamma_F} \\
		Z_{\tilde{M}^\vee}^\circ \big/ Z_{\tilde{G}^\vee}^\circ \arrow[ru] \arrow[rd, "\text{descent}"' sloped] & & \\
		& Z_{M_\eta^\vee}^{\Gamma_F} \big/ Z_{G_\eta^\vee}^{\Gamma_F} \arrow[r, "\text{std.\ endo.}"'] & Z_{\overline{M^!_\epsilon}^\vee}^{\Gamma_F} \big/ Z_{\overline{G^![s]_{\epsilon[s]}}^\vee}^{\Gamma_F} \arrow[uu]
	\end{tikzcd}\]
	Every quotient group is connected by \cite[Lemma 1.1]{Ar99}, and every arrow is surjective by comparing dimensions using ellipticity. The rightmost arrow comes from switching factors of type $\mathrm{B}$ and $\mathrm{C}$, namely
	\begin{gather*}
		L = \SO(2n+1), \quad \overline{L} = \Sp(2n), \\
		Z_{\overline{L}^\vee}^{\Gamma_F} = \{1\} \to \{\pm 1\} = Z_{L^\vee}^{\Gamma_F}.
	\end{gather*}

	Count the cardinality of fibers in two ways. The composite along
	\begin{tikzpicture}[baseline=(O), xscale=0.3, yscale=0.2]
		\coordinate (O) at (0, 0.5);
		\draw[->] (0, 0) -- (1, 1) -- (2, 1);
	\end{tikzpicture}
	yields
	\[ i_{M^!}\left( \tilde{G}, G^![s]\right)^{-1} e^{G^![s]}_{M^!}(\epsilon[s])^{-1}, \]
	whilst
	\begin{tikzpicture}[baseline=(O), xscale=0.3, yscale=0.2]
		\coordinate (O) at (0, -0.7);
		\draw[->] (0, 0) -- (1, -1) -- (2, -1) -- (2, 0);
	\end{tikzpicture}
	yields $|\mathcal{Z}| \cdot i_{\overline{M^!_\epsilon}}\left( G_\eta, \overline{G^![s]_{\epsilon[s]}}\right)^{-1} c(s)^{-1}$, the $c(s)^{-1}$ coming from the rightmost arrow. This proves the desired equality.
\end{proof}

It remains to establish the following
\begin{proposition}[Cf.\ {\cite[III.8.1 Proposition]{MW16-1}}]\label{prop:Shalika-matching-elliptic}
	Assume that $\mathcal{O}^!$ is elliptic in $M^!$. For all $\delta \in SD_{\mathrm{geom}, \tilde{G}\text{-equi}}(M^!) \otimes \mes(M^!)^\vee$, with $\Supp(\delta)$ sufficiently close to $\mathcal{O}^!$, we have
	\[ g^{\tilde{G}, \Endo}_{\tilde{M}, \mathcal{O}}(\mathbf{M}^!, \delta) = g^{\tilde{G}}_{\tilde{M}, \mathcal{O}}\left( \trans_{\mathbf{M}^!, \tilde{M}}(\delta) \right). \]
\end{proposition}
\begin{proof}
	Choose $\epsilon$ and $\eta$ as before. Recall that
	\[ g^{\tilde{G}, \Endo}_{\tilde{M}, \mathcal{O}}(\mathbf{M}^!, \delta) = \sum_{s \in \Endo_{\mathbf{M}^!}(\tilde{G})} i_{M^!}(\tilde{G}, G^![s]) \trans_{\mathbf{G}^![s], \tilde{G}} Sg^{G^![s]}_{M^!, \mathcal{O}[s]}\left( \delta[s], B^{\tilde{G}} \right). \]
	Take $\delta_\epsilon \in SD_{\mathrm{geom}}(M^!_\epsilon) \otimes \mes(M^!_\epsilon)^\vee$ with $\Supp(\delta_\epsilon)$ close to $1$ and
	\[ \delta = S\desc^{M^!, *}_\epsilon(\delta_\epsilon), \quad \delta[s] = S\desc^{M^!, *}_{\epsilon[s]}(\delta_\epsilon). \]
	
	For every $s \in \Endo_{\mathbf{M}^!}(\tilde{G})$, we define
	\[ \tau(s) := \begin{cases}
		Sg^{G^![s]_{\epsilon[s]}}_{M^!_\epsilon, \mathrm{unip}} \left( \delta_\epsilon, B^{\tilde{G}}_{\mathcal{O}^![s]} \right), & \text{if} \; \epsilon[s] \;\text{is elliptic in}\; G^![s] \\
		0, & \text{otherwise.}
	\end{cases}\]
	It is of unipotent support, so we can set $\tau(s)^{G^![s]} := S\desc^{G^![s], *}_{\epsilon[s]}\left( \tau(s) \right)$. By \cite[III.4.4 Proposition]{MW16-1},
	\[ Sg^{G^![s]}_{M^!, \mathcal{O}^![s]}\left( \delta[s], B^{\tilde{G}} \right) =
	e^{G^![s]}_{M^!}(\epsilon[s]) \tau(s)^{G^![s]} \]

	Hereafter, we adopt the conventions from \S\ref{sec:descent-endoscopy-Levi} on descent of endoscopic data and transfer factors. For every $y \in \mathcal{Y}$, use the chosen transfer factors to define
	\begin{align*}
		\tau[s, y] & := \underbracket{\trans_{\mathbf{G}^!_{\eta[y]}, G_{\eta[y]}}}_{\text{std.}} \underbracket{\trans_{G^![s]_{\epsilon[s]}, \overline{G^![s]_{\epsilon[s]}}}}_{\text{non-std.}} (\tau(s)) \; \in D_{\mathrm{unip}}(G_{\eta[y]}) \otimes \mes(G_{\eta[y]})^\vee , \\
		\tau[s, y]^{\tilde{G}} & := \desc^{\tilde{G}, *}_{\eta[y]}(\tau[s, y]) \; \in D_{\mathrm{geom}, -}(\tilde{G}) \otimes \mes(G)^\vee .
	\end{align*}
	Applying Lemma \ref{prop:transfer-descent-y} relative to $\mathbf{G}^![s] \in \Endo_{\elli}(\tilde{G})$, with the $d(s, y)$ from \eqref{eqn:dsy} replacing $d(y)$, we see
	\[ \trans_{\mathbf{G}^![s], \tilde{G}} Sg^{G^![s]}_{M^!, \mathcal{O}^![s]}\left( \delta[s], B^{\tilde{G}} \right)
	= e^{G^![s]}_{M^!}(\epsilon[s]) \displaystyle\sum_{y \in \dot{\mathcal{Y}}} d(s, y) \tau[s, y]^{\tilde{G}}. \]
	Note that $\tau[s, y]$ is nonzero only when $\epsilon[s]$ is elliptic in $G^![s]$.
	
	Next, suppose that $\mathbf{M}^! \in \Endo_{\elli}(\tilde{M})$ arises from an elliptic semisimple element $s^\flat \in \tilde{M}^\vee$. Choose $\overline{s^\flat} \in M_\eta^\vee$ that is part of the endoscopic datum $\mathbf{M}^!_\eta$ of $M_\eta$ obtained by descent at $(\epsilon, \eta)$. By Remark \ref{rem:descent-homomorphism}, the $\mathbf{G}_\eta^! \in \Endo(G_\eta)$ obtained by descent at ($\epsilon[s], \eta$) is equivalent to the one associated with $\overline{s} \in \Endo_{\mathbf{M}^!_\eta}(G_\eta)$, where $\overline{s}$ is determined by $s \in \Endo_{\mathbf{M}^!}(\tilde{G})$ as in \eqref{eqn:descent-s}. Also recall that the endoscopic group underlying $\mathbf{G}^!_\eta$ is identified with $\overline{G^![s]_{\epsilon[s]}}$.

	Now let $s$ vary. Put
	\begin{align*}
		\mathcal{S} & := \left\{ s \in \Endo_{\mathbf{M}^!}(\tilde{G}) : \epsilon[s] \;\text{is elliptic in}\; G^![s] \right\}, \\
		\overline{\delta_\epsilon} & := \trans_{M^!_\epsilon , \overline{M^!_\epsilon}}(\delta_\epsilon).
	\end{align*}
	To each $s \in \mathcal{S}$ we attach the constant $c(s)$ of Lemma \ref{prop:Shalika-matching-aux0}. By \cite[III.6.6 Proposition]{MW16-1},
	\[ \trans_{G^![s]_{\epsilon[s]}, \overline{G^![s]_{\epsilon[s]}}} \left(\tau(s) \right) =
	\begin{cases}
		c(s)^{-1} Sg^{\overline{G^![s]_{\epsilon[s]}}}_{\overline{M^!_\epsilon}, \mathrm{unip}}( \overline{\delta_\epsilon} ), & s \in \mathcal{S} \\
		0, & \text{otherwise,}
	\end{cases}\]
	Indeed, the cited result is applicable by our choice of $B$-functions in Definition \ref{def:B-metaplectic}. Therefore
	\[ \tau[s, y] = \begin{cases}
		c(s)^{-1} \trans_{\mathbf{G}^!_{\eta[y]}, G_{\eta[y]}} Sg^{\overline{G^![s]_{\epsilon[s]}}}_{\overline{M^!_\epsilon}, \mathrm{unip}}\left( \overline{\delta_\epsilon} \right) , & s \in \mathcal{S} \\
		0, & s \notin \mathcal{S}.
	\end{cases}\]

	Fix representatives in $\mathcal{Y}$ for $\dot{\mathcal{Y}}$. We infer that $g^{\tilde{G}, \Endo}_{\tilde{M}, \mathcal{O}}(\mathbf{M}^!, \delta)$ equals
	\begin{equation}\label{eqn:Shalika-matching-aux1}
		\sum_{y \in \dot{\mathcal{Y}}} \desc^{\tilde{G}, *}_{\eta[y]} \sum_{s \in \mathcal{S}} i_{M^!}(\tilde{G}, G^![s]) e^{G^![s]}_{M^!}(\epsilon[s]) d(s, y) c(s)^{-1} \trans_{\mathbf{G}^!_{\eta[y]}, G_{\eta[y]}} \left( Sg^{\overline{G^![s]_{\epsilon[s]}}}_{\overline{M^!_\epsilon}, \mathrm{unip}} (\overline{\delta_\epsilon}) \right).
	\end{equation}

	Hereafter, suppose $\eta$ is elliptic in $G$. Interpret $\Endo_{\mathbf{M}^!}(\tilde{G})$ (thus $\mathcal{S}$) as a subset of $s^\flat Z_{\tilde{M}^\vee}^\circ \big/ Z_{\tilde{G}^\vee}^\circ$. Consider the map $s \mapsto \overline{s}$ from $\Endo_{\mathbf{M}^!}(\tilde{G})$ to $\overline{s^\flat} Z_{M_\eta^\vee}^{\Gamma_F} \big/ Z_{G_\eta^\vee}^{\Gamma_F}$ in \eqref{eqn:descent-s}. Lemma \ref{prop:Shalika-matching-aux} implies
	\begin{equation}\label{eqn:Shalika-matching-aux2}
		\mathcal{S} = \left\{ s \in \Endo_{\mathbf{M}^!}(\tilde{G}): \overline{s} \;\text{determines an elliptic endoscopic datum of}\; G_\eta \right\}.
	\end{equation}
	Combining Lemma \ref{prop:Shalika-matching-aux0} and \eqref{eqn:Shalika-matching-aux1}, we arrive at
	\begin{gather*}
		g^{\tilde{G}, \Endo}_{\tilde{M}, \mathcal{O}}(\mathbf{M}^!, \delta) = \sum_{y \in \dot{\mathcal{Y}}} \desc^{\tilde{G}, *}_{\eta[y]}  \left( \sum_{\overline{s} \in \overline{s^\flat} Z_{M_\eta^\vee}^{\Gamma_F} / Z_{G_\eta^\vee}^{\Gamma_F} } x(\overline{s}, y) \trans_{\mathbf{G}^!(\overline{s}), G_{\eta[y]}} Sg^{G^!(\overline{s})}_{\overline{M^!_\epsilon}, \mathrm{unip}} (\overline{\delta_\epsilon}) \right),
	\end{gather*}
	where
	\begin{compactitem}
		\item $\mathbf{G}^!(\overline{s})$ denotes the endoscopic datum of $G_\eta$ (equivalently, for $G_{\eta[y]}$) determined by $\overline{s}$ and the endoscopic datum of $M_\eta$ obtained by descent, taken up to equivalence;
		\item for each $y$ and $\overline{s}$, we set
		\begin{equation}\label{eqn:xsy}
			x(\overline{s}, y) := |\mathcal{Z}|^{-1} \displaystyle\sum_{\substack{ s \in \mathcal{S} \\ s \mapsto \overline{s} }} d(s, y) i_{\overline{M^!_\epsilon}}(G_\eta, G^!(\overline{s})) .
		\end{equation}
	\end{compactitem}

	By definition, $i_{\overline{M^!_\epsilon}}(G_\eta, G^!(\overline{s})) = 0$ when $\mathbf{G}^!(\overline{s})$ is non-elliptic. We will prove in Lemma \ref{prop:xsy} that
	\[ x(\overline{s}, y) = \begin{cases}
		i_{\overline{M^!_\epsilon}}(G_\eta, G^!(\overline{s})) d^M(y), & y \in \dot{\mathcal{Y}}^M \\
		0, & y \in \dot{\mathcal{Y}} \smallsetminus \dot{\mathcal{Y}}^M.
	\end{cases}\]

	Still assuming $\eta$ (hence each $\eta[y]$) is elliptic in $G$, we set
	\[ \delta[y] := \trans_{\overline{M^!_\epsilon}, M_{\eta[y]}}(\overline{\delta_\epsilon}) \]
	and deduce
	\begin{gather*}
		g^{\tilde{G}, \Endo}_{\tilde{M}, \mathcal{O}}(\mathbf{M}^!, \delta) = \sum_{y \in \dot{\mathcal{Y}}^M} d^M(y) \desc^{\tilde{G}, *}_{\eta[y]} \left(g^{G_{\eta[y]}, \Endo}_{M_{\eta[y]}, \mathrm{unip}}(\delta[y])\right).
	\end{gather*}

	By Arthur's result \cite[III.1.5 Corollaire]{MW16-1} in standard endoscopy, we obtain
	\begin{align*}
		g^{\tilde{G}, \Endo}_{\tilde{M}, \mathcal{O}}(\mathbf{M}^!, \delta) & = \sum_{y \in \dot{\mathcal{Y}}^M} d^M(y) \desc^{\tilde{G}, *}_{\eta[y]} \left( g^{G_{\eta[y]}}_{M_{\eta[y]}, \mathrm{unip}}(\delta[y]) \right) \\
		& = \sum_{y \in \dot{\mathcal{Y}}^M} d^M(y) g^{\tilde{G}}_{\tilde{M}, \mathcal{O}}\left( \delta[y]^{\tilde{M}} \right) = g^{\tilde{G}}_{\tilde{M}, \mathcal{O}}\left( \sum_{y \in \dot{\mathcal{Y}}^M} d^M(y) \delta[y]^{\tilde{M}} \right), \\
		\text{with}\; \delta[y]^{\tilde{M}} & := \desc^{\tilde{M}, *}_{\eta[y]}(\delta[y]) ,
	\end{align*}
	the second equality being the descent of ordinary germs \cite[III.4.2 Proposition]{MW16-1} at elliptic elements. The last expression equals $g^{\tilde{G}}_{\tilde{M}, \mathcal{O}}\left( \trans_{\mathbf{M}^!, \tilde{M}} (\delta) \right)$ by Lemma \ref{prop:transfer-descent-y}. This settles the case when $\eta$ is elliptic in $G$.
	
	Finally, assume $\eta$ (hence each $\eta[y]$) is non-elliptic in $G$, then $\mathcal{S} = \emptyset$ by Lemma \ref{prop:Shalika-matching-aux0}, hence \eqref{eqn:Shalika-matching-aux1} gives $g^{\tilde{G}, \Endo}_{\tilde{M}, \mathcal{O}}(\mathbf{M}^!, \delta) = 0$. On the other hand, $g^{\tilde{G}}_{\tilde{M}, \mathcal{O}} \circ \desc^{\tilde{M}, *}_{\eta[y]} = 0$ in this case by \textit{loc.\ cit.}, therefore
	\[ g^{\tilde{G}}_{\tilde{M}, \mathcal{O}}\left( \trans_{\mathbf{M}^!, \tilde{M}} (\delta) \right) = g^{\tilde{G}}_{\tilde{M}, \mathcal{O}}\left( \sum_{y \in \dot{\mathcal{Y}}^M} \delta[y]^{\tilde{M}} \right) = 0. \]
	This settles the remaining cases.
\end{proof}

\begin{lemma}[Cf.\ {\cite[III.8.3]{MW16-1}}]\label{prop:xsy}
	Keep the assumptions of Proposition \ref{prop:Shalika-matching-elliptic}. In addition, assume that $\eta$ is elliptic in $G$. For each $\overline{s} \in \overline{s^\flat} Z_{M_\eta^\vee}^{\Gamma_F} \big/ Z_{G_\eta^\vee}^{\Gamma_F}$ we have
	\[ x(\overline{s}, y) = \begin{cases}
		i_{\overline{M^!_\epsilon}}(G_\eta, G^!(\overline{s})) d^M(y) , & y \in \dot{\mathcal{Y}}^M \\
		0, & y \in \dot{\mathcal{Y}} \smallsetminus \dot{\mathcal{Y}}^M.
	\end{cases}\]
\end{lemma}
\begin{proof}
	We still interpret $\Endo_{\mathbf{M}^!}(\tilde{G})$ as a subset of $s^\flat Z_{\tilde{M}^\vee}^\circ \big/ Z_{\tilde{G}^\vee}^\circ$. By Lemma \ref{prop:Shalika-matching-aux}, $\mathbf{G}^!(\overline{s})$ is elliptic if $\overline{s}$ comes from $\mathcal{S}$. Conversely, suppose that $\overline{s} \in \overline{s^\flat} Z_{M_\eta^\vee}^{\Gamma_F} \big/ Z_{G_\eta^\vee}^{\Gamma_F}$ and $\mathbf{G}^!(\overline{s})$ is an elliptic endoscopic datum of $G_\eta$, then we have
	\begin{equation}\label{eqn:xsy-aux}
		\left\{ s \in \mathcal{S}: s \mapsto \overline{s} \right\} = \left\{ s^\flat t \in s^\flat Z_{\tilde{M}^\vee}^\circ \big/ Z_{\tilde{G}^\vee}^\circ : \overline{s^\flat} \cdot \mathrm{descent}(t) = \overline{s} \right\}.
	\end{equation}
	Indeed, $\subset$ is trivial. As for $\supset$, recall that one constructs $\mathbf{G}^![s] \in \Endo(\tilde{G})$ from $\mathbf{M}^!$ and any $s \in s^\flat Z_{\tilde{M}^\vee}^\circ \big/ Z_{\tilde{G}^\vee}^\circ$. If $\overline{s^\flat} \cdot \mathrm{descent}(s) = \overline{s}$, then $\mathbf{G}^![s]$ must be elliptic, otherwise $\mathbf{G}^!(\overline{s})$ would factor through an endoscopic datum of $L_\eta$ for some proper $L \in \mathcal{L}^G(M)$, which is absurd for the following reason: $\mathbf{G}^!(\overline{s})$ is an elliptic endoscopic datum of $G_\eta$, thus
	\[ \dim A_{L_\eta} \leq \dim A_{G^!(\overline{s})} = \dim A_{G_\eta}. \]
	As $\eta \in G(F)$ is elliptic, $\dim A_{G_\eta} = \dim A_G$, we get $\dim A_{L_\eta} \leq \dim A_G$, which contradicts $\dim A_G < \dim A_L \leq \dim A_{L_\eta}$. The conclusion is that $s \in \Endo_{\mathbf{M}^!}(\tilde{G})$, and this shows \eqref{eqn:xsy-aux}.
	
	As observed in the proof of Lemma \ref{prop:Shalika-matching-aux0}, the map ``descent'' in Remark \ref{rem:descent-homomorphism} is surjective since $\eta$ is elliptic in $M$ and $G$. Hence both sides of \eqref{eqn:xsy-aux} are nonempty.

	By the prior discussion, the definition \eqref{eqn:xsy} can be written as
	\[ x(\overline{s}, y) = |\mathcal{Z}|^{-1} \underbracket{i_{\overline{M^!_\epsilon}}(G_\eta, G^!(\overline{s}))}_{=0 \;\text{unless}\; \mathbf{G}^!(\overline{s}) \;\text{is ell.} } \displaystyle\sum_{\substack{s \in \mathcal{S} \\ s \mapsto \overline{s}}} d(s, y), \]
	and the sum above is taken over a $\mathcal{Z}$-coset. The case $y \in \dot{\mathcal{Y}}^M$ becomes straightforward since $d(s,y) = d^M(y)$ by Proposition \ref{prop:dsy}.

	Assume from now onward that $y \notin \dot{\mathcal{Y}}^M$. Our goal is to show that when $\mathbf{G}^!(\overline{s})$ is elliptic, we have
	\[ \sum_{z \in \mathcal{Z}} d(sz, y) = 0 \]
	for any $s \in \mathcal{S}$ with $s \mapsto \overline{s}$. In this case $\mathbf{G}^!(\overline{s})$ is relevant for $G_{\eta[y]}$ by Proposition \ref{prop:relevance-elliptic} (true for general $y$). Also note that $d(sz, 1) = d^M(1) \neq 0$ since all endoscopic data of the quasisplit group $M_\eta$ are relevant. We are thus reduced to show that
	\[ s \in \mathcal{S}, \; y \notin \dot{\mathcal{Y}}^M \implies \sum_{z \in \mathcal{Z}} \frac{d(sz, y)}{d(sz, 1)} = 0. \]

	By definition,
	\[ \frac{d(sz, y)}{d(sz, 1)} = \frac{\Delta_{\mathbf{G}^![sz], \tilde{G}}(\exp(\overline{Y}) \epsilon[sz], \exp(X[y])\tilde{\eta}[y] ) }{\Delta_{\mathbf{G}^![sz], \tilde{G}}(\exp(\overline{Y}) \epsilon[sz], \exp(X)\tilde{\eta} )} \cdot
	\frac{\Delta(\overline{s}, 1)(\exp(\overline{Y}), \exp(X))}{\Delta(\overline{s}, y)(\exp(\overline{Y}), \exp(X[y]))}. \]
	On the right hand side, the second quotient is nonzero and independent of $z$; denote by $\chi(z)$ the first quotient, as a function in $z \in \mathcal{Z}$. It suffices to show that
	\begin{equation}\label{eqn:xsy-aux0}
		y \notin \dot{\mathcal{Y}}^M \implies \chi(z) = \text{constant} \times \;\text{a nontrivial character of}\; \mathcal{Z}.
	\end{equation}
	
	Since $\mathbf{G}^!(\overline{s})$ is relevant, Proposition \ref{prop:diagram-relevant} says that there exists a diagram $(\epsilon[s], B^!, T^!, B, T, \eta)$ for $\mathbf{G}^![s] \in \Endo_{\elli}(\tilde{G})$, compatibly with descent. In fact, since $\epsilon[s]$ (resp.\ $\eta$) is elliptic in $G^![s]$ (resp.\ $G$), we may take $T^!$ (resp.\ $T$) elliptic in $G^![s]_{\epsilon[s]}$ (resp.\ $G_\eta$).

	Furthermore, \cite[10.2 Lemma]{Ko86} guarantees that there is a representative of $y$ such that the pure inner twist $\Ad(y): G_{\eta, \overline{F}} \rightiso G_{\eta[y], \overline{F}}$ restricts to an isomorphism $T \rightiso T[y]$ between $F$-tori, which determines $\mathrm{inv}(y) \in \Hm^1(F, T)$. By construction, the image $\mathrm{inv}^G(y)$ in $\Hm^1(F, G_\eta)$ of $\mathrm{inv}(y)$ is the class corresponding to $y \in \dot{\mathcal{Y}}$.
	
	The $Y$, $\overline{Y}$, $X$ and $X[y]$ above may be taken in the Lie algebras of these tori; one can suppose that they match under isomorphisms. The relative position $\mathrm{inv}(\cdot, \cdot)$ satisfies
	\[ \mathrm{inv}\left(\exp(X)\eta, \exp(X[y]) \eta[y] \right) = \mathrm{inv}(y) \in \Hm^1(F, T) . \]
	
	Recall that $\exp(X[y])\tilde{\eta}[y]$ is stably conjugate in the metaplectic sense to $\exp(X)\tilde{\eta}$ in $\tilde{G}$. The cocycle property \cite[Proposition 5.13]{Li11} for metaplectic transfer factors extends immediately to groups of metaplectic type. It implies
	\[ \chi(z) := \frac{\Delta_{\mathbf{G}^![sz], \tilde{G}}(\exp(\overline{Y}) \epsilon[sz], \exp(X[y])\tilde{\eta}[y] ) }{\Delta_{\mathbf{G}^![sz], \tilde{G}}(\exp(\overline{Y}) \epsilon[sz], \exp(X)\tilde{\eta} )} = \lrangle{\mathrm{inv}(y), sz}; \]
	the pairing is formed by
	\begin{itemize}
		\item regarding $sz$ as an element of $\check{T}^{\Gamma_F}$ (with the help of diagrams, see \S\ref{sec:diagram}),
		\item the canonical perfect pairing $\Hm^1(F, T) \times \pi_0(\check{T}^{\Gamma_F}) \to \CC^\times$ due to Kottwitz \cite[1.2 Theorem]{Ko86}.
	\end{itemize}
	In particular, $\chi$ is a nonzero constant times a character of $\mathcal{Z}$.
	
	Now we use the functoriality of the perfect pairings
	\[\begin{tikzcd}[row sep=small, column sep=small]
		\lrangle{\cdot, \cdot}: \Hm^1(F, T) \arrow[phantom, r, "\times" description] \arrow[d] & \pi_0\left(\check{T}^{\Gamma_F}\right) \arrow[rd] & \\
		\lrangle{\cdot, \cdot}_G: \Hm^1(F, G_\eta) \arrow[phantom, r, "\times" description] & \pi_0\left( Z_{G_\eta^\vee}^{\Gamma_F}\right) \arrow[u] \arrow[twoheadrightarrow, d] \arrow[r] & \CC^\times . \\
		\lrangle{\cdot, \cdot}_M: \Hm^1(F, M_\eta) \arrow[hookrightarrow, u] \arrow[phantom, r, "\times" description] & \pi_0\left(Z_{M_\eta^\vee}^{\Gamma_F}\right) \arrow[ru] &
	\end{tikzcd}\]

	As $\eta$ is elliptic in $G$, we have $\pi_0\left( Z_{G_\eta^\vee}^{\Gamma_F}\right) = Z_{G_\eta^\vee}^{\Gamma_F}$. We contend that
	\begin{equation}\label{eqn:xsy-aux1}
		\Ker\left[ \pi_0\left( Z_{G_\eta^\vee}^{\Gamma_F}\right) \twoheadrightarrow \pi_0\left( Z_{M_\eta^\vee}^{\Gamma_F} \right) \right] \subset \Image\left[ \mathcal{Z} \to Z_{G_\eta^\vee}^{\Gamma_F} \right].
	\end{equation}
	To see this, use the commutative diagram
	\[\begin{tikzcd}
		Z_{\tilde{M}^\vee}^\circ \arrow[r, "\sim"] \arrow[hookrightarrow, d] & Z_{M^\vee}^\circ \arrow[twoheadrightarrow, r] & Z_{M_\eta^\vee}^{\Gamma_F, \circ} \arrow[hookrightarrow, d] & \\
		Z_{M^\vee} \arrow[rr] & & Z_{M_\eta^\vee}^{\Gamma_F} & Z_{G_\eta^\vee}^{\Gamma_F} \arrow[hookrightarrow, l] ,
	\end{tikzcd}\]
	the surjectivity being due to the ellipticity of $\eta$ in $M$. Suppose $x \in Z_{G_\eta^\vee}^{\Gamma_F}$ maps to $1 \in \pi_0( Z_{M_\eta^\vee}^{\Gamma_F})$, then $x \in Z_{M_\eta^\vee}^{\Gamma_F, \circ}$, hence there exists $z \in Z_{\tilde{M}^\vee}^\circ$ with $z \mapsto x$. This easily leads to \eqref{eqn:xsy-aux1}.

	Let us prove \eqref{eqn:xsy-aux0} now. If $\chi$ is constant on $\mathcal{Z}$, the functoriality of pairings would imply
	\[ \lrangle{ \cdot, \mathrm{inv}^G(y)}_G \; \text{is trivial on}\; \Ker\left[ \pi_0\left( Z_{G_\eta^\vee}^{\Gamma_F} \right) \twoheadrightarrow \pi_0\left( Z_{M_\eta^\vee}^{\Gamma_F} \right) \right]. \]
	Hence $\mathrm{inv}^G(y)$ is in the image of $\Hm^1(F, M_\eta)$. By \eqref{eqn:Y-inclusion}, this contradicts $y \notin \dot{\mathcal{Y}}^M$.
\end{proof}

\begin{remark}\label{rem:Z-nonarch-arch}
	The last part of the proof concerning $\lrangle{\cdot, \cdot}$ and $\mathcal{Z}$ also works for Archimedean $F$, upon replacing $\Hm^1(F, \cdot)$ by $\Hm^1_{\mathrm{ab}}(F, \cdot)$; of course, $\Hm^1(F, T) \rightiso \Hm^1_{\mathrm{ab}}(F, T)$. There is one caveat: for a connected reductive $F$-group $L$, Kottwitz's natural homomorphism
	\begin{gather*}
		\Hm^1_{\mathrm{ab}}(F, L) \to \pi_0\left( Z_{L^\vee}^{\Gamma_F} \right)^D , \\
		(\cdots)^D := \text{the Pontryagin dual}
	\end{gather*}
	is bijective for non-Archimedean $F$, but only injective when $F = \R$; see also \cite[Proposition 1.7.3]{Lab99}. The arguments before will carry over if the surjectivity holds for both $L = G_\eta$ and $L = M_\eta$. Such an adaptation will be needed in \S\ref{sec:pf-jump-Endo-lemmas}.
\end{remark}

\section{Proof of the weighted fundamental lemma}\label{sec:proof-LFP}
Let $F$, $\tilde{M}$, $\tilde{G}$, $K$ be as in \S\ref{sec:LFP}. In particular, we are in the unramified situation defined in \S\ref{sec:LF}. Our goal is to prove Theorem \ref{prop:LFP-general} by reducing it to the $\tilde{G}$-regular case solved in \cite{Li12a}.

Fix $\mathbf{M}^! \in \Endo_{\elli}(\tilde{M})$. The lines $\mes(G)$, $\mes(M)$, $\mes(M^!)$, etc.\ are trivialized by the Haar measures determined by hyperspecial subgroups; see \S\ref{sec:LF}. The protagonists here are $r^{\tilde{G}}_{\tilde{M}}(\tilde{\gamma}, K)$ and $r^{\tilde{G}, \Endo}_{\tilde{M}}(\mathbf{M}^!, \delta)$, where $\tilde{\gamma} \in D_{\mathrm{geom}, -}(\tilde{M})$ and $\delta \in SD_{\mathrm{geom}}(M^!)$. 

Fix a stable semisimple conjugacy class $\mathcal{O}^!$ (resp.\ $\mathcal{O}$) in $M^!(F)$ (resp.\ $M(F)$), and assume $\mathcal{O}^! \mapsto \mathcal{O}$. We have the following expansions into Shalika germs.
\begin{itemize}
	\item If $\Supp(\tilde{\gamma})$ is sufficiently close to $\rev^{-1}(\mathcal{O})$, we have
	\[ r^{\tilde{G}}_{\tilde{M}}(\tilde{\gamma}, K) = \sum_{L \in \mathcal{L}(M)} r^{\tilde{G}}_{\tilde{L}}\left( g^{\tilde{L}}_{\tilde{M}, \mathcal{O}}(\tilde{\gamma}), K \right). \]
	This is the metaplectic analogue of \cite[II.4.5 (1)]{MW16-1}. It is based on the fact that non-invariant weighted orbital integrals have the same expansion into weighted Shalika germs as the invariant ones, see \cite[\S 9]{Ar88LB}. The proof is based on harmonic analysis and has nothing to do with metaplectic covers, so it carries over verbatim.
	\item When $\mathcal{O}$ is $G$-equisingular, the equality above reduces to $r^{\tilde{G}}_{\tilde{M}}(\tilde{\gamma}, K) = r^{\tilde{G}}_{\tilde{M}}\left( g^{\tilde{M}}_{\tilde{M}, \mathcal{O}}(\tilde{\gamma}), K \right)$.
	\item Let $s \in \Endo_{\mathbf{M}^!}(\tilde{G})$. If $\mathcal{O}^!$ is $\tilde{G}$-equisingular and $\Supp(\delta)$ is sufficiently close to $\mathcal{O}^!$, we have
	\[ s^{G^![s]}_{M^!}(\delta[s]) = s^{G^![s]}_{M^!}\left( g^{M^!}_{M^!, \mathcal{O}^![s]}(\delta[s]) \right); \]
	here $g^{M^!}_{M^!, \mathcal{O}^![s]}(\delta[s]) \in SD_{\mathrm{geom}}(M^!, \mathcal{O}^![s])$ by \cite[II.2.2 Lemme]{MW16-1}, so we may take $s^{G^![s]}_{M^!}$. Indeed, $\mathcal{O}^![s]$ is $G^![s]$-equisingular by Lemma \ref{prop:equisingular-endo}, so this reduces to \cite[II.4.5 (2)]{MW16-1}.
	\item Recall that $r^{\tilde{G}, \Endo}_{\tilde{M}}(\mathbf{M}^!, \delta) = \sum_s i_{M^!}(\tilde{G}, G^![s]) s^{G^![s]}_{M^!}(\delta[s])$. Under the assumptions of the previous item,
	\begin{equation*}
		r^{\tilde{G}, \Endo}_{\tilde{M}}(\mathbf{M}^!, \delta) = r^{\tilde{G}, \Endo}_{\tilde{M}}\left( \mathbf{M}^!, g^{M^!}_{M^!, \mathcal{O}^!}(\delta) \right).
	\end{equation*}
	This is because $z[s] \in Z_{M^!}(F)$ implies $g^{M^!}_{M^!, \mathcal{O}^!}(\delta)[s] = g^{M^!}_{M^!, \mathcal{O}^![s]}(\delta[s])$.
\end{itemize}

We shall prove the equisingular case of Theorem \ref{prop:LFP-general} first.

\begin{lemma}\label{prop:LFP-equisingular}
	The statements of Theorem \ref{prop:LFP-general} hold when $\delta \in SD_{\mathrm{geom}}(M^!, \mathcal{O}^!)$ and $\mathcal{O}^!$ is $\tilde{G}$-equisingular.
\end{lemma}
\begin{proof}
	By \cite[II.2.2 Lemme]{MW16-1}, there exists $\delta_{\mathrm{reg}} \in SD_{\mathrm{geom}, \tilde{G}\text{-reg}}(M^!)$ with $\Supp(\delta_{\mathrm{reg}})$ arbitrarily close to $\mathcal{O}^!$, such that $g^{M^!}_{M^!, \mathcal{O}^!}(\delta_{\mathrm{reg}}) = \delta$. The germs expansions above imply
	\begin{gather*}
		r^{\tilde{G}, \Endo}_{\tilde{M}}\left( \mathbf{M}^!, \delta \right) = r^{\tilde{G}, \Endo}_{\tilde{M}}\left( \mathbf{M}^!, g^{M^!}_{M^!, \mathcal{O}^!}(\delta_{\mathrm{reg}}) \right) = r^{\tilde{G}, \Endo}_{\tilde{M}}\left( \mathbf{M}^!, \delta_{\mathrm{reg}} \right) \\
		r^{\tilde{G}}_{\tilde{M}}\left( \trans_{\mathbf{M}^!, \tilde{M}}(\delta_{\mathrm{reg}}), K \right) = r^{\tilde{G}}_{\tilde{M}}\left( g^{\tilde{M}}_{\tilde{M}, \mathcal{O}}\left(\trans_{\mathbf{M}^!, \tilde{M}}(\delta_{\mathrm{reg}})\right), K \right).
	\end{gather*}

	Corollary \ref{prop:germ-matching-G} asserts
	\[ g^{\tilde{M}}_{\tilde{M}, \mathcal{O}}(\trans_{\mathbf{M}^!, \tilde{M}}(\delta_{\mathrm{reg}})) = \trans_{\mathbf{M}^!, \tilde{M}}\left( g^{M^!}_{M^!, \mathcal{O}^!}(\delta_{\mathrm{reg}}) \right). \]
	By the $\tilde{G}$-regular case \cite[Théorème 4.2.1]{Li12a} of Theorem \ref{prop:LFP-general},
	\[ r^{\tilde{G}, \Endo}_{\tilde{M}}\left( \mathbf{M}^!, \delta_{\mathrm{reg}} \right) = r^{\tilde{G}}_{\tilde{M}}\left( \trans_{\mathbf{M}^!, \tilde{M}}(\delta_{\mathrm{reg}}), K \right). \]
	This implies that $r^{\tilde{G}, \Endo}_{\tilde{M}}\left( \mathbf{M}^!, \delta \right) = r^{\tilde{G}}_{\tilde{M}}\left(\trans_{\mathbf{M}^!, \tilde{M}}(\delta), K\right)$, as required.
\end{proof}

In order to access to the general case, we need the following germ expansions from \S\ref{sec:rho-sigma}. Recall the spaces of germs $U^G_M \supset U^{G, +}_M$, $U^{G^![s]}_{M^!} \supset U^{G^![s], +}_{M^!}$ defined before Lemma \ref{prop:U-germ-5}.

\begin{itemize}
	\item Let $\tilde{\gamma} \in D_{\mathrm{geom},-}(\tilde{M})$, we have the equivalence
	\[ r^{\tilde{G}}_{\tilde{M}}(a\tilde{\gamma}, K) \sim \sum_{L \in \mathcal{L}(M)} \sum_{J \in \mathcal{J}^G_M} r^{\tilde{G}}_{\tilde{L}}\left( \rho^{\tilde{L}}_J(\tilde{\gamma}, a)^{\tilde{L}}, K \right) \]
	between germs of functions of $a \in A_M(F)$ in general position, $a \to 1$. Again, this is a fact from harmonic analysis having nothing to do with metaplectic covers; see \cite[p.307]{MW16-1} for explanations.
	
	The terms $a \mapsto r^{\tilde{G}}_{\tilde{L}} \left( \rho^{\tilde{L}}_J(\tilde{\gamma}, a)^{\tilde{L}}, K \right)$ belongs to $U_J^{L \supset M}$. The right hand side is thus in
	\[ r^{\tilde{G}}_{\tilde{M}}(\tilde{\gamma}, K) + U^{G, +}_M \; \subset U^G_M . \]	

	\item The stable side is similar: for each $s \in \Endo_{\mathbf{M}^!}(\tilde{G})$, we use the germ expansion of $s^{G^![s]}_{M^!}$ from \cite[II.4.7]{MW16-1} to obtain
	\[ s^{G^![s]}_{M^!}(a^! \delta[s]) \sim s^{G^![s]}_{M^!}(\delta[s]) + U^{G^![s], +}_{M^!} \; \subset U^{G^![s]}_{M^!}, \]
	where $a^! \in A_{M^!}(F)$ is in general position, $a^! \to 1$.
	
	\item Consider the endoscopic case. Recall that we have $\xi: A_M \rightiso A_{M^!}$. Given the previous case, the upcoming Lemma \ref{prop:LFP-prep} will guarantee that
	\[ r^{\tilde{G}, \Endo}_{\tilde{M}}(\mathbf{M}^!, \xi(a)\delta) \sim r^{\tilde{G}, \Endo}_{\tilde{M}}(\mathbf{M}^!, \delta) + U^{G, +}_M \; \subset U^G_M . \]
\end{itemize}

\begin{lemma}\label{prop:LFP-prep}
	Let $s \in \Endo_{\mathbf{M^!}}(\tilde{G})$ and consider the spaces $U^G_M$, $U^{G^![s]}_{M^!}$.
	\begin{enumerate}[(i)]
		\item For all $u^! \in U^{G^![s]}_{M^!}$, put $u := u^! \circ \xi$. Then $u \in U^G_M$ and it has the same constant term as $u'$.
		\item If $u^! \sim v^!$ in $U^{G^![s]}_{M^!}$, then $u \sim v$ in $U^G_M$.
	\end{enumerate}
	Consequently, $u^! \mapsto u$ induces a linear map from $U^{G^![s]}_{M^!}/\sim$ to $U^G_M /\sim$.
\end{lemma}
\begin{proof}
	First, we show $u^! \in U_{J^!} \implies u \in U_J$ for all $J^! \in \mathcal{J}^{G^![s]}_{M^!}\left( B^{\tilde{G}}_{\mathcal{O}^![s]} \right)$ and $J \in \mathcal{J}^G_M$. Let $\gamma_1, \ldots, \gamma_m \in \Sigma^{G^![s]}\left( A_{M^!}, B^{\tilde{G}}_{\mathcal{O}^![s]} \right)$ be linearly independent with $m = \dim \mathfrak{a}^{G^![s]}_{M^!}$. They determine $J^!$ with image $J$. Set $u^! := \prod_{i=1}^m \log \left| \gamma_i - \gamma_i^{-1} \right| \in U_{J^!}$. To see $u \in U_J$, it suffices to recall Lemma \ref{prop:J-transfer}.
	
	As each $L^! \in \mathcal{L}^{G^![s]}(M^!)$ comes from some $L \in \mathcal{L}^G(M)$ (see Lemma \ref{prop:sL-Ls} with $M=R$), (i) follows by replacing $G^![s]$ by various $L^!$.
	
	The proof of (ii) is essentially the same as \cite[II.4.7 (6)]{MW16-1}. Let us be more precise. The relation $\sim$ between germs is defined in terms of
	\begin{compactitem}
		\item a function $d$ with $d(\exp H) = \|H\|$, which works simultaneously for $A_M$ and $A_{M^!}$ via $\xi$;
		\item estimates of functions on domains of the form $|\alpha(a) - 1| > c \cdot d(a)$ where $\alpha \in \Sigma^G(A_M)$ or $\Sigma^{G^![s]}(A_{M^!})$.
	\end{compactitem}
	The sets $\Sigma^G(A_M)$ and $\Sigma^{G^![s]}(A_{M^!})$ are not in bijection via $\xi$. However, they match under dilation by $B$-functions. Since $B^{\tilde{G}}$ takes value in $\{1, \frac{1}{2}\}$, $p > 2$ and $a \to 1$, one readily sees that the dilation is harmless.
\end{proof}

We are ready to prove the weighted fundamental lemma in general.

\begin{proof}[Proof of Theorem \ref{prop:LFP-general}]
	Fix $\delta \in SD_{\mathrm{geom}}(M^!)$. For $a \in A_M(F)$ in general position, $a \to 1$, the $\tilde{G}$-equisingular matching (Lemma \ref{prop:LFP-equisingular}) implies
	\begin{align*}
		r^{\tilde{G}, \Endo}_{\tilde{M}}\left( \mathbf{M}^!, \xi(a)\delta \right) & = r^{\tilde{G}}_{\tilde{M}}\left( \trans_{\mathbf{M}^!, \tilde{M}}(\xi(a)\delta), K \right) \\
		& = r^{\tilde{G}}_{\tilde{M}}\left( a \trans_{\mathbf{M}^!, \tilde{M}}(\delta), K \right).
	\end{align*}
	By Lemma \ref{prop:LFP-prep}, both sides are germs in $U^G_M$, and we deduce the equivalence of germs
	\[ r^{\tilde{G}}_{\tilde{M}}(\trans_{\mathbf{M}^!, \tilde{M}}(\delta), K) - r^{\tilde{G}, \Endo}_{\tilde{M}}(\mathbf{M}^!, \delta) + (\text{element of}\; U^{G, +}_M) \sim 0 \quad \text{in}\; U^G_M. \]
	Hence Proposition \ref{prop:U-germ} (ii) and Lemma \ref{prop:U-germ-5} yield $r^{\tilde{G}}_{\tilde{M}}(\trans_{\mathbf{M}^!, \tilde{M}}(\delta), K) = r^{\tilde{G}, \Endo}_{\tilde{M}}(\mathbf{M}^!, \delta)$.
\end{proof}

\chapter{Weighted characters: non-Archimedean case, and construction of \texorpdfstring{$\epsilon_{\tilde{M}}$}{epsilonM}}\label{sec:cpt-nonarch}

Let $F$ be a non-Archimedean local field of characteristic zero. We will work with a covering of metaplectic type $\rev: \tilde{G} \to G(F)$ where $G = \prod_{i \in I} \GL(n_i) \times \Sp(W)$; here we fix the symplectic $F$-vector space $(W, \lrangle{\cdot|\cdot})$ and the additive character $\psi$ of $F$ to construct the eightfold covering $\tilde{G}$. We also fix a minimal Levi subgroup $M_0$ of $G$.

What follows is closely modeled on Arthur's strategy in \cite{Ar88-1, Ar03-3}. First off, we will review the unitary weighted characters $J_{\tilde{M}}(\pi, f)$ and the associated invariant distributions $I_{\tilde{M}}(\pi, \cdot)$ on $\tilde{G}$, where $\pi \in \Pi_{\mathrm{unit}, -}(\tilde{M})$. They were already defined in \cite{Li13, Li14b} via the \emph{canonical normalization}, i.e.\ in terms of Harish-Chandra's $\mu$-functions instead of chosen normalizing factors.

Next, we define the compactly supported weighted characters ${}^c J_{\tilde{M}}(\pi, f)$, the corresponding variant ${}^c \phi_{\tilde{M}}$ of $\phi_{\tilde{M}}$, and the resulting invariant distributions ${}^c I_{\tilde{M}}(\pi, f)$ in \S\S\ref{sec:weighted-characters-nonarch}--\ref{sec:cIGM}. There is a family of linear maps ${}^c \theta_{\tilde{M}}$ linking $\phi_{\tilde{M}}$ and ${}^c \phi_{\tilde{M}}$, where $M$ ranges over $\mathcal{L}(M_0)$.

In \S\ref{sec:cIGMEndo}, we define the endoscopic counterparts ${}^c I^{\Endo}_{\tilde{M}}$ and ${}^c \theta^{\Endo}_{\tilde{M}}$, and formulate the corresponding matching Theorem \ref{prop:cpt-supported-equalities} for them.

These definitions are geared to establish the existence of a linear map
\[ \epsilon_{\tilde{M}}: \orbI_{\asp}(\tilde{G}) \otimes \mes(G) \to \orbI_{\mathrm{ac}, \asp, \cusp}(\tilde{M}) \otimes \mes(M) \]
for a given $M \in \mathcal{L}(M_0)$ under inductive assumptions. It is characterized by
\[ I^{\Endo}_{\tilde{M}}(\tilde{\gamma}, f) - I_{\tilde{M}}(\tilde{\gamma}, f) = I^{\tilde{M}}\left( \tilde{\gamma}, \epsilon_{\tilde{M}}(f) \right) \]
for all $\tilde{\gamma} \in D_{\mathrm{geom}, -}(\tilde{M}) \otimes \mes(M)^\vee$. In other words, $\epsilon_{\tilde{M}}$ is the obstruction to the local geometric Theorem \ref{prop:local-geometric}. The point is to consider $I^{\Endo}_{\tilde{M}}(\tilde{\gamma}, f) - I_{\tilde{M}}(\tilde{\gamma}, f)$ as a function in $\tilde{\gamma} \in D_{\mathrm{geom}, G\text{-equi}, -}(\tilde{M}) \otimes \mes(M)^\vee$, where $f$ is kept fixed, and then the construction divides into two steps.
\begin{itemize}
	\item Locally, the matching of weighted Shalika germs establishes the existence of such a $\epsilon_{\tilde{M}}(f)$ around any given stable semisimple class $\mathcal{O} \subset M(F)$.
	\item The global problem is more delicate, since neither $I_{\tilde{M}}(\tilde{\gamma}, f)$ nor ${}^c I_{\tilde{M}}(\tilde{\gamma}, f)$ are compactly supported in $\tilde{\gamma}$ modulo conjugation. This is where ${}^c I_{\tilde{M}}$, ${}^c \theta_{\tilde{M}}$ and their endoscopic avatars enter.
\end{itemize}

By induction, it is legitimate to assume the matching theorems for $I_{\tilde{M}}(\tilde{\gamma}, f)$, ${}^c I_{\tilde{M}}(\tilde{\gamma}, f)$ and ${}^c \theta_{\tilde{M}}$ hold when $M$ (resp.\ $G$) is replaced a larger (resp.\ proper) Levi subgroup. Then the ingredients above and an argument via partition of unity furnish the required $\epsilon_{\tilde{M}}(f)$. This is the content of Proposition \ref{def:epsilonM-nonarch}.

\section{Review of unitary weighted characters}\label{sec:weighted-characters-nonarch}
For unitary genuine representations, their \emph{weighted characters} on coverings are defined in \cite[\S 5.7]{Li12b} and \cite[\S 4.1]{Li14b}. Specifically, let $M \in \mathcal{L}(M_0)$ and let $\pi \in \Pi_{\mathrm{unit}, -}(\tilde{M})$. One defines
\[ J_{\tilde{M}}(\pi, f) = J^{\tilde{G}}_{\tilde{M}}(\pi, f) = \Tr\left( \mathcal{M}_M(\pi, \tilde{P}) I_{\tilde{P}}(\pi, f) \right) \]
where
\begin{itemize}
	\item $f \in C^\infty_{c,\asp}(\tilde{G})$,
	\item $P \in \mathcal{P}(M)$ is arbitrary,
	\item $\mathcal{M}_M(\pi, \tilde{P})$ is the endomorphism of $I_{\tilde{P}}(\pi)$ arising from the $(G,M)$-family
	\begin{equation}\label{eqn:M-GM-family}\begin{gathered}
		\mu_{\tilde{Q}|\tilde{P}}(\pi_\lambda)^{-1} \mu_{\tilde{Q}|\tilde{P}}\left( \pi_{\lambda + \frac{\Lambda}{2} } \right) 	J_{\tilde{Q}|\tilde{P}}(\pi_\lambda)^{-1} J_{\tilde{Q}|\tilde{P}}(\pi_{\lambda + \Lambda}), \\
		Q \in \mathcal{P}(M), \; \Lambda \in \mathfrak{a}_{M, \CC}^*,
	\end{gathered}\end{equation}
	where $\mu_{\tilde{Q}|\tilde{P}}(\pi_\lambda) := \prod_{\alpha \in \Sigma_P^{\mathrm{red}} \cap \Sigma_{Q^-}^{\mathrm{red}}} \mu_\alpha(\pi_\lambda)$ is made from Harish-Chandra's $\mu$-functions; this is called the \emph{canonical normalization} of weighted characters in \cite{Ar98}.
	\index{muQP@$\mu_{\tilde{Q}{"|}\tilde{P}}(\pi)$}
\end{itemize}

The result turns out to be independent of $P$. It also differs from Arthur's earlier recipe in \cite{Ar89a}, which is based on the choice of normalizing factors. When $M=G$, we revert to the usual characters.

Define
\begin{align*}
	\mathfrak{a}_{M, F} & := H_M(M(F)), \\
	\mathfrak{a}_{M, F}^\vee & := \left\{ H \in \mathfrak{a}_M^*: \forall X \in \mathfrak{a}_{M,F},\; \lrangle{H, X} \in 2\pi\Z \right\}.
\end{align*}
We have thus the compact torus $i\mathfrak{a}^*_{M, F} := i\mathfrak{a}_M^* / i\mathfrak{a}_{M, F}^\vee$ and the complex torus $\mathfrak{a}_{M, \CC}^* / i\mathfrak{a}_{M, F}^\vee$.
\index{aMF@$\mathfrak{a}_{M, F}$, $\mathfrak{a}_{M, F}^\vee$, $i\mathfrak{a}_{M, F}^*$}

The function $\lambda \mapsto J_{\tilde{M}}(\pi_\lambda , f)$ on $i\mathfrak{a}^*_M$ is analytic and periodic relative to the lattice $i\mathfrak{a}_{M, F}^\vee$. One defines
\[ J_{\tilde{M}}(\pi, X, f) := \int_{i\mathfrak{a}^*_{M, F}} J_{\tilde{M}}(\pi_\lambda, f) e^{-\lrangle{\lambda, X}} \dd\lambda, \quad X \in \mathfrak{a}_{M, F}. \]
Then $J_{\tilde{M}}(\pi_\lambda, X, f) = J_{\tilde{M}}(\pi, X, f) e^{\lrangle{\lambda, X}}$, and the resulting distribution on $\tilde{G}$ is concentrated at $H_G = X$. Therefore $J_{\tilde{M}}(\pi, X, \cdot)$ extends to $C^\infty_{\mathrm{ac}, \asp}(\tilde{G})$. See \cite[Corollaire 4.6]{Li14b}.
\index{JMpiX@$J_{\tilde{M}}(\pi, f)$, $J_{\tilde{M}}(\pi, X, f)$}

Using the trace Paley--Wiener theorem, the map
\[ \phi_{\tilde{M}}: C^\infty_{\mathrm{ac}, \asp}(\tilde{G}) \otimes \mes(G) \to \orbI_{\mathrm{ac}, \asp}(\tilde{M}) \otimes \mes(M) \]
is defined by mapping $f$ to the function $(\pi, X) \mapsto J_{\tilde{M}}(\pi, X, f)$. See \cite[Théorème 4.13]{Li14b}.
\index{phiM@$\phi_{\tilde{M}}$}

The non-unitary weighted characters will be reviewed in \S\ref{sec:weighted-characters}.

\section{The maps \texorpdfstring{${}^c \phi_{\tilde{M}}$}{cphi} and \texorpdfstring{${}^c \theta_{\tilde{M}}$}{ctheta}}\label{sec:cphi-nonarch}
Let $M \in \mathcal{L}(M_0)$. It makes sense to say that a rational function on $\mathfrak{a}_{M, \CC}^* / i\mathfrak{a}_{M, F}^\vee$ has singularities along ``hyperplanes'' of the form $e^{\lrangle{\cdot, \check{\alpha}}} = \text{const}$, where $\alpha \in \Sigma(A_M)$. See \cite[VIII.1.2]{MW16-2}.

The basic example is
\[ \mathfrak{a}_{R, \CC}^* \ni \lambda \mapsto J^{\tilde{G}}_{\tilde{M}}\left( \pi_\lambda^{\tilde{M}}, f \right) \]
where
\begin{itemize}
	\item $R \subset \mathcal{L}^M(M_0)$,
	\item $\pi$ is a genuine admissible representation of $\tilde{R}$,
	\item $\pi_\lambda^{\tilde{M}} = (\pi_\lambda)^{\tilde{M}}$ is the normalized parabolic induction,
	\item $f \in \orbI_{\asp}(\tilde{G}) \otimes \mes(G)$.
\end{itemize}
It is rational with a finite number of polar hyperplanes of the form $e^{\lrangle{\cdot, \check{\alpha}}} = \text{const}$, where $\alpha \in \Sigma(A_R)$.

As in \cite[VIII.1.1]{MW16-2} or \cite[\S 4]{Ar88-1}, we fix a family of functions $\omega_S: \mathfrak{a}_R \to [0,1]$ for every Levi subgroup $R$ of $G$ and $S \in \mathcal{P}(R)$, such that
\begin{itemize}
	\item $\Supp(\omega_S)$ lies in some translate of the chamber $\mathfrak{a}_S^+$;
	\item $\sum_{S \in \mathcal{P}(R)} \omega_S = 1$;
	\item the families $(R, (\omega_S)_{S \in \mathcal{P}(R)})$ are invariant under $G(F)$-conjugation.
\end{itemize}
Moreover, we can replace $G$ by any of its Levi subgroup and impose the following compatibility: for every $L \in \mathcal{L}(R)$ and $S \in \mathcal{P}^L(R)$, we have
\begin{equation}\label{eqn:omega-compatibility}
	\omega_S = \sum_{\substack{S' \in \mathcal{P}(R) \\ S' \cap L = S}} \omega_{S'}.
\end{equation}

\begin{definition}
	\index{JMpi-c@${}^c J_{\tilde{M}}(\pi, f)$}
	Let $\pi \in D_{\elli, -}(\tilde{R}) \otimes \mes(R)^\vee$. Set
	\begin{multline*}
		{}^c J_{\tilde{M}}\left(  \pi^{\tilde{M}}, f \right) = {}^c J^{\tilde{G}}_{\tilde{M}}\left(  \pi^{\tilde{M}}, f \right) := \sum_{X \in \mathfrak{a}_{R,F}} \sum_{S \in \mathcal{P}^G(R)} \omega_S(X) \\
		\int_{\nu_S + i\mathfrak{a}^*_{R,F}} J_{\tilde{M}}\left( \pi_\lambda^{\tilde{M}}, f \right) e^{-\lrangle{\lambda, X}} \dd\lambda
	\end{multline*}
	where
	\begin{itemize}
		\item $i\mathfrak{a}^*_{R, F}$ carries the Haar measure dual to the counting measure on $\mathfrak{a}_{R, F}$,
		\item $\nu_S \in \mathfrak{a}_R^*$ satisfies $\lrangle{\nu_S, \check{\alpha}} \gg 0$ for all $\alpha \in \Sigma^S(A_R)$.
	\end{itemize}
	This is independent of the choice of $(\nu_S)_S$. In view of \eqref{eqn:Dspec-Ind}, we can thus define ${}^c J_{\tilde{M}}(\pi, f)$ for all $\pi \in D_{\mathrm{temp}, -}(\tilde{M}) \otimes \mes(M)^\vee$.
\end{definition}

\begin{definition-proposition}
	\index{phiM-c@${}^c \phi_{\tilde{M}}$}
	There is a linear map
	\[ {}^c \phi_{\tilde{M}} = {}^c \phi^{\tilde{G}}_{\tilde{M}} : C^\infty_{c, \asp}(\tilde{G}) \otimes \mes(G) \to \orbI_{\asp}(\tilde{M}) \otimes \mes(M(F)) \]
	characterized by $I^{\tilde{M}}\left( \pi, {}^c \phi_{\tilde{M}}(f) \right) = {}^c J_{\tilde{M}}\left( \pi, f \right)$ for all $\pi \in D_{\mathrm{temp}, -}(\tilde{M}) \otimes \mes(M)^\vee$.
\end{definition-proposition}
\begin{proof}
	This is the analogue of \cite[VIII.1.3 Proposition]{MW16-2}, which carries over to coverings because of its analytic nature.
\end{proof}

Similarly, one has
\begin{enumerate}
	\item the descent formula for ${}^c \phi_{\tilde{M}}$ as in \cite[VIII.1.4 (1)]{MW16-2};
	\item the explicit behavior of ${}^c \phi_{\tilde{M}}$ under conjugation, cf.\ \cite[VIII.1.4 (2)]{MW16-2};
	\item the extension to $C^\infty_{\mathrm{ac}, \asp}(\tilde{G}) \otimes \mes(G) \to \orbI_{\asp, \mathrm{ac}}(\tilde{M}) \otimes \mes(M)$ as in \cite[VIII.1.4 (6)]{MW16-2}.
\end{enumerate}
All these proofs carry over to coverings without change.

\begin{definition}[Cf.\ {\cite[VIII.1.5]{MW16-2}}]
	\index{thetaM-c@${}^c \theta_{\tilde{M}}$}
	We define the linear maps
	\[ {}^c \theta_{\tilde{M}} = {}^c \theta^{\tilde{G}}_{\tilde{M}}: \orbI_{\asp}(\tilde{G}) \otimes \mes(G) \to \orbI_{\mathrm{ac}, \asp}(\tilde{M})^{W^G(M)} \otimes \mes(M) \]
	characterized inductively by
	\[ \sum_{L \in \mathcal{L}(M)} {}^c \theta^{\tilde{L}}_{\tilde{M}} \circ {}^c \phi^{\tilde{G}}_{\tilde{L}} = \phi^{\tilde{G}}_{\tilde{M}}. \]
\end{definition}

These linear maps also satisfy a descent formula relative to parabolic descent. See \cite[VIII.1.6 Lemme]{MW16-2}, whose proof carries over to coverings as usual.

\begin{definition}
	\index{Schwartz function}
	We say that $f \in \orbI_{\mathrm{ac}, \asp}(\tilde{G}) \otimes \mes(G)$ is a \emph{Schwartz function} if $X \mapsto I^{\tilde{G}}(\pi, X, f)$ is rapidly decreasing on $\mathfrak{a}_{G,F}$ for each $\pi \in D_{\mathrm{temp},-}(\tilde{G}) \otimes \mes(G)^\vee$.
\end{definition}

The definition is the same as \cite[VIII.1.7]{MW16-2}. If $f$ is Schwartz, one can define
\[ I^{\tilde{G}}(\pi, \lambda, f) := \sum_{X \in \mathfrak{a}_{G,F}} e^{\lrangle{\lambda, X}} I^{\tilde{G}}(\pi, X, f), \quad \lambda \in i\mathfrak{a}_G^* / i\mathfrak{a}_{G,F}^\vee . \]
We record the following basic properties of $I^{\tilde{G}}(\pi, \lambda, f)$ from \textit{loc. cit.}
\begin{itemize}
	\item It is $C^\infty$ in $\lambda$.
	\item Every element of $\orbI_{\asp}(\tilde{G}) \otimes \mes(G)$ is Schwartz. Fourier inversion implies that
	\[ f \in \orbI_{\asp}(\tilde{G}) \otimes \mes(G) \implies I^{\tilde{G}}(\pi, \lambda, f) = I^{\tilde{G}}(\pi_\lambda, f). \]
	\item Suppose that $I^{\tilde{G}}(\pi, \lambda, f)$ extends to a rational function in $\lambda$ on $\mathfrak{a}_{G, \CC}^* / i\mathfrak{a}_{G,F}^\vee$. By setting
	\[ I^{\tilde{G}}(\pi, \nu, X, f) := \int_{\underbracket{\nu + i\mathfrak{a}_{G,F}^*}_{\text{compact}}} I^{\tilde{G}}(\pi, \lambda, f) e^{-\lrangle{\lambda, X}} \dd\lambda, \quad X \in \mathfrak{a}_{G,F}, \]
	provided that $\nu + i\mathfrak{a}_{G,F}$ contains no poles, we have
	\[ I^{\tilde{G}}(\pi, \lambda, f) = \sum_{X \in \mathfrak{a}_{G,F}} e^{\lrangle{\lambda, X}} I^{\tilde{G}}(\pi, \nu, X, f) \quad \text{when} \quad \Re(\lambda) = \nu. \]
\end{itemize}

In the special case $f \in \orbI_{\asp}(\tilde{G}) \otimes \mes(G)$, we have $I^{\tilde{G}}(\pi, \lambda, f) = I^{\tilde{G}}(\pi_\lambda, f)$ for all $\lambda \in \mathfrak{a}^*_{G, \CC}$, thus $I^{\tilde{G}}(\pi, \nu, X, f)$ is independent of $\nu$ by shifting contours; in fact, it equals $I^{\tilde{G}}(\pi, 0, X, f) = I^{\tilde{G}}(\pi, X, f)$.

\begin{proposition}[Cf.\ {\cite[VIII.1.8]{MW16-2}}]\label{prop:cthetaGM-nonarch}
	Let $f \in \orbI_{\asp}(\tilde{G}) \otimes \mes(G)$ and $M \in \mathcal{L}(M_0)$. Then
	\begin{itemize}
		\item ${}^c \theta^{\tilde{G}}_{\tilde{M}}(f)$ is Schwartz on $\tilde{M}$;
		\item for all $\pi \in D_{\mathrm{temp},-}(\tilde{M}) \otimes \mes(M)^\vee$, the function
		\[ \lambda \mapsto I^{\tilde{M}}(\pi, \lambda, {}^c \theta^{\tilde{G}}_{\tilde{M}}(f)) \]
		extends to a rational function on $\mathfrak{a}^*_{M, \CC} / i\mathfrak{a}_{M,F}^\vee$ with only a finite number of polar hyperplanes, all of which are of the form $e^{\lrangle{\cdot, \check{\alpha}}} = \mathrm{const}$ with $\alpha \in \Sigma^G(A_M)$;
		\item if $M \neq G$ and $\pi \in D_{\mathrm{ell},-}(\tilde{M}) \otimes \mes(M)^\vee$, we have
		\[ \sum_{S \in \mathcal{P}(M)} \omega_S(X) I^{\tilde{M}} \left( \pi, \nu_S, X, {}^c \theta^{\tilde{G}}_{\tilde{M}}(f) \right) = 0, \quad X \in \mathfrak{a}_{M,F} \]
		provided that $\nu_S$ is sufficiently deep relative to $\mathfrak{a}_S^+$, for each $S \in \mathcal{P}(M)$.
	\end{itemize}
\end{proposition}

Because of the analytic nature of the arguments, we refer to \textit{loc.\ cit.} for the proofs.

To conclude this section, we consider the stable analogues
\[ {}^c S\theta^{G^!}_{M^!}: S\orbI(G^!) \otimes \mes(G^!) \to S\orbI_{\mathrm{ac}}(M^!)^{W^{G^!}(M^!)} \otimes \mes(M^!) \]
\index{SthetaM-c@${}^c S\theta^{G^{"!}}_{M^{"!}}$}
in \cite[VIII.2]{MW16-2}, where $G^!$ is quasisplit and $M^!$ is a Levi subgroup. Suppose that $M^!$ arises from $\mathbf{M}^! \in \Endo_{\elli}(\tilde{M})$ and $\mathbf{G}^! = \mathbf{G}^![s] \in \Endo_{\elli}(\tilde{G})$ as in \eqref{eqn:s-situation}. We assume that the data $(\omega_P)_{P \in \mathcal{P}(M)}$ and $(\omega_{P^!})_{P^! \in \mathcal{P}(M^!)}$ are chosen compatibly, namely
\begin{equation}\label{eqn:omega-compatibility-Endo}
	\omega_{P^!} = \sum_{\substack{P \in \mathcal{P}(M) \\ P \mapsto P^!}} \omega_P
\end{equation}
where we used the isomorphism $\mathfrak{a}_M \simeq \mathfrak{a}_{M^!}$ from endoscopy, and $P \mapsto P^!$ signifies that the corresponding chambers satisfy $\overline{\mathfrak{a}_{P^!}^+} \supset \overline{\mathfrak{a}_P^+}$. In fact, $\overline{\mathfrak{a}_{P^!}^+}$ is the union of various $\overline{\mathfrak{a}_P^+}$, as there are ``more walls'' in $G$ than in $G^!$. See \cite[VIII.2.1]{MW16-2}.

There are also stable analogues of the notion of Schwartz functions $f^!$. It suffices to consider $S^{G^!}(\pi, \lambda, f^!)$ where $\pi \in SD_{\mathrm{temp}}(G^!) \otimes \mes(G^!)^\vee$. Also recall that the tempered local Langlands correspondence for $G^!$ is available; see \S\ref{sec:spectral-distributions}.

\section{The distributions \texorpdfstring{${}^c I^{\tilde{G}}_{\tilde{M}}(\tilde{\gamma}, \cdot)$}{cIGM}}\label{sec:cIGM}
\index{IGM-c-gamma@${}^c I_{\tilde{M}}(\tilde{\gamma}, f)$}
Let $M \in \mathcal{L}(M_0)$, $\tilde{\gamma} \in D_{\mathrm{geom}, -}(\tilde{M}) \otimes \mes(M)^\vee$. Define inductively the distribution ${}^c I_{\tilde{M}}(\tilde{\gamma}, \cdot) = {}^c I^{\tilde{G}}_{\tilde{M}}(\tilde{\gamma}, \cdot)$ by
\[ \sum_{L \in \mathcal{L}(M)} {}^c I^{\tilde{L}}_{\tilde{M}}\left( \tilde{\gamma}, {}^c \phi^{\tilde{G}}_{\tilde{L}}(f) \right) = J^{\tilde{G}}_{\tilde{M}}(\tilde{\gamma}, f), \quad f \in C^\infty_{c, \asp}(\tilde{G}) \otimes \mes(G). \]

The following facts follow in exactly the same way as in \cite[VIII.1.9]{MW16-2}.
\begin{itemize}
	\item They factor through $\orbI_{\asp}(\tilde{G}) \otimes \mes(G)$.
	\item For all $f \in \orbI_{\asp}(\tilde{G}) \otimes \mes(G)$, there exists a compact subset $\Omega \subset M(F)$ such that ${}^c I^{\tilde{G}}_{\tilde{M}}(\tilde{\gamma}, f) = 0$ whenever $\Supp(\tilde{\gamma}) \cap \rev^{-1}(\Omega) = \emptyset$.
	\item ${}^c I^{\tilde{G}}_{\tilde{M}}(\tilde{\gamma}, f) = \displaystyle\sum_{L \in \mathcal{L}^G(M)} I^{\tilde{L}}_{\tilde{M}}\left( \tilde{\gamma}, {}^c \theta^{\tilde{G}}_{\tilde{L}}(f)\right)$.
\end{itemize}

On the stable side, we also have the stable distributions ${}^c S^{G^!}_{M^!}\left( \delta, B^{\tilde{G}}, \cdot \right)$ for $\delta \in SD_{\mathrm{geom}}(M^!) \otimes \mes(M)^\vee$. This is the content of \cite[VIII.2]{MW16-2}. Here the auxiliary functions $\omega_{P^!}$ are chosen compatibly with the $\omega_P$ on the metaplectic side, as explained in \eqref{eqn:omega-compatibility-Endo}.

\begin{remark}\label{rem:cS-no-B}
	The system $B^{\tilde{G}}$ of $B$-functions in ${}^c S^{G^!}_{M^!}$ will be dropped, since we will ultimately focus on the case of $\delta \in SD_{\mathrm{geom}, \tilde{G}\text{-reg}}(M^!) \otimes \mes(M^!)^\vee$.
\end{remark}
\index{SGM-delta-c@${}^c S^{G^{"!}}_{M^{"!}}\left( \delta, B^{\tilde{G}}, \cdot \right)$}

\section{The distributions \texorpdfstring{${}^c I^{\tilde{G}, \Endo}_{\tilde{M}}(\tilde{\gamma}, \cdot)$}{cIGMEndo} and statement of matching theorems}\label{sec:cIGMEndo}
Let $M \in \mathcal{L}(M_0)$ and $\mathbf{M}^! \in \Endo_{\elli}(\tilde{M})$. For every $s \in \Endo_{\mathbf{M}^!}(\tilde{G})$, we have the automorphism $g \mapsto g[s] := g(\cdot z[s])$ of $S\orbI(M^!) \otimes \mes(M)$.

\begin{definition}\label{def:thetaEndo-nonarch}
	\index{thetaM-Endo-c@${}^c \theta^{\Endo}_{\tilde{M}}$}
	Let ${}^c \theta^{\Endo}_{\tilde{M}}(\mathbf{M}^!, \cdot) = {}^c \theta^{\tilde{G}, \Endo}_{\tilde{M}}(\mathbf{M}^!, \cdot)$ be the linear map
	\begin{align*}
		\orbI_{\asp}(\tilde{G}) \otimes \mes(G) & \to S\orbI_{\mathrm{ac}}(M^!) \otimes \mes(M^!) \\
		f & \mapsto \sum_{s \in \Endo_{\mathbf{M}^!}(\tilde{G})} i_{M^!}(\tilde{G}, G^![s]) {}^c S\theta^{G^![s]}_{M^!}\left( f^{G^![s]} \right)[s]
	\end{align*}
	where $f^{G^![s]} := \Trans_{\mathbf{G}^![s], \tilde{G}}(f)$ as usual.
\end{definition}

Again, we have the following standard properties.
\begin{itemize}
	\item The image of ${}^c \theta^{\Endo}_{\tilde{M}}(\mathbf{M}^!, \cdot)$ is $W(M)$-invariant. Cf.\ \cite[VIII.3.2]{MW16-2}.
	\item For all $f$, the function ${}^c \theta^{\tilde{G}, \Endo}_{\tilde{M}}(\mathbf{M}^!, f)$ is Schwartz. For every $\pi \in SD_{\mathrm{temp}}(M^!) \otimes \mes(M^!)^\vee$, the function
	\[ \lambda \mapsto S^{M^!}\left(\pi, \lambda, {}^c \theta^{\tilde{G}, \Endo}_{\tilde{M}}(\mathbf{M}^!, f) \right) \]
	extends to a rational function on $i\mathfrak{a}_{M^!}^* / i\mathfrak{a}_{M^!, F}^\vee$, with only finitely many polar hyperplanes of the form $e^{\lrangle{\cdot, \check{\alpha}}} = \mathrm{const}$ where $\alpha \in \bigcup_s \Sigma^{G^![s]}(A_{M^!})$. Cf.\ \cite[VIII.3.1 (4)]{MW16-2}.
\end{itemize}

Due to the presence of $z[s]$-twists, the descent formula requires more care.

\begin{proposition}[Cf.\ {\cite[VIII.3.3 Proposition (i)]{MW16-2}}]\label{prop:ctheta-descent}
	Suppose that $R \in \mathcal{L}^M(M_0)$, $\mathbf{R}^! \in \Endo_{\elli}(\tilde{R})$ and $t \in \Endo_{\mathbf{R}^!}(\tilde{M})$ so that $\mathbf{M}^! = \mathbf{M}^![t]$. In other words, we are in the situation
	\[\begin{tikzcd}
		M^! \arrow[dash, r, "\text{endo.}", "\text{ell.}"'] & \tilde{M} \\
		R^! \arrow[dash, r, "\text{endo.}", "\text{ell.}"'] \arrow[hookrightarrow, u, "\text{Levi}"] & \tilde{R} \arrow[hookrightarrow, u, "\text{Levi}"'] .
	\end{tikzcd}\]
	For every $f \in \orbI_{\asp}(\tilde{G}) \otimes \mes(G)$, we have
	\[ {}^c \theta^{\tilde{G}, \Endo}_{\tilde{M}}(\mathbf{M}^!, f)_{R^!}[t] = \sum_{L \in \mathcal{L}^G(R)} d^G_R(M, L) {}^c \theta^{\tilde{L}, \Endo}_{\tilde{R}}\left( \mathbf{R}^!, f_{\tilde{L}} \right). \]
	Here $f \mapsto f_{\tilde{L}}$ denotes the parabolic descent and $(\cdots)_{R^!}$ is its stable avatar.
\end{proposition}
\begin{proof}
	Follow the paradigm of the proof of Proposition \ref{prop:descent-orbint-Endo}. Apart from the combinatorics from the cited result, the only input is the descent formula \cite[VIII.2.3 Lemme]{MW16-2} for ${}^c S\theta^{G^![s]}_{M^!}$.
\end{proof}

\begin{proposition}[Cf.\ {\cite[VIII.3.4]{MW16-2}}]\label{prop:theta-Endo-glue-nonarch}
	Given $M \in \mathcal{L}(M_0)$, there exists a unique linear map
	\[ {}^c \theta^{\Endo}_{\tilde{M}} = {}^c \theta^{\tilde{G}, \Endo}_{\tilde{M}} : \orbI_{\asp}(\tilde{G}) \otimes \mes(G) \to \orbI_{\mathrm{ac}, \asp}(\tilde{M}) \otimes \mes(M) \]
	characterized by
	\[ \Trans_{\mathbf{M}^!, \tilde{M}} \left( {}^c \theta^{\Endo}_{\tilde{M}}(f) \right) = {}^c \theta^{\Endo}_{\tilde{M}}(\mathbf{M}^!, f) \]
	for all $\mathbf{M}^! \in \Endo_{\elli}(\tilde{M})$.
\end{proposition}
\begin{proof}
	It amounts to a straightforward gluing of various ${}^c \theta^{\Endo}_{\tilde{M}}(\mathbf{M}^!, \cdot)$. See the arguments for Definition--Proposition \ref{def:I-geom-Endo}.
\end{proof}

\begin{proposition}\label{prop:thetaEndo-descent}
	Given $M \in \mathcal{L}(M_0)$ and $R \in \mathcal{L}^M(M_0)$, we have the descent formula
	\[ {}^c \theta^{\tilde{G}, \Endo}_{\tilde{M}}(f)_{\tilde{R}} = \sum_{L \in \mathcal{L}^G(R)} d^G_R(M, L) {}^c \theta^{\tilde{L}, \Endo}_{\tilde{R}}(f_{\tilde{L}}). \]
\end{proposition}
\begin{proof}
	It suffices to show that for each $\mathbf{R}^! \in \Endo_{\elli}(\tilde{R})$, both sides in the assertion have the same transfer to $R^!$. Take $t \in \Endo_{\mathbf{R}^!}(\tilde{M})$ and set $\mathbf{M}^! := \mathbf{M}^![t]$. The left hand side becomes ${}^c \theta^{\Endo}_{\tilde{M}}(\mathbf{M}^!, f)_{R^!} [t]$ whilst the right hand side becomes $\sum_{L \in \mathcal{L}^G(R)} d^G_R(M, L) {}^c \theta^{\tilde{L}, \Endo}_{\tilde{R}}(\mathbf{R}^!, f_{\tilde{L}})$. The required equality follows at once from Proposition \ref{prop:ctheta-descent}.
\end{proof}

Moreover, the following properties hold for the endoscopic counterparts.
\begin{itemize}
	\item The function ${}^c \theta^{\Endo}_{\tilde{M}}(f)$ is Schwartz, and $\lambda \mapsto I^{\tilde{G}}(\pi, \lambda, {}^c \theta^{\Endo}_{\tilde{M}}(f))$ extends to a rational function with only finitely many polar hyperplanes of the form $e^{\lrangle{\cdot, \check{\alpha}}} = \mathrm{const}$, where $\alpha \in \Sigma^G(A_M)$.
	
	Indeed, this is the same as (and easier than) \cite[VIII.3.6]{MW16-2}, the main ingredients being the Proposition \ref{prop:thetaEndo-descent} and the case for ${}^c S\theta$. Note that the coroots from $\Sigma^{G^!}(A_{M^!})$ and $\Sigma^G(A_M)$ may differ by a scaling by $2$, but this does not affect the required properties.
	
	\item If $M \neq G$ and $\pi \in D_{\mathrm{ell},-}(\tilde{M}) \otimes \mes(M)^\vee$, then
	\[ \sum_{S \in \mathcal{P}(M)} \omega_S(X) I^{\tilde{M}} \left( \pi, \nu_S, X, {}^c \theta^{\tilde{G}, \Endo}_{\tilde{M}}(f) \right) = 0 \]
	for all $X \in \mathfrak{a}_{M,F}$, provided that $\nu_S$ is sufficiently deep for each $S \in \mathcal{P}(M)$. Again, this reduces to the stable counterparts for various $G^! \supset M^!$: see \cite[VIII.3.7 Proposition]{MW16-2}. Our case is even easier since $\mathfrak{a}_{M, F} = \mathfrak{a}_{M^!, F}$ for all $\mathbf{M}^! \in \Endo_{\elli}(\tilde{M})$.
\end{itemize}

\begin{definition}
	\index{IGM-c-Endo-shrek@${}^c I^{\Endo}_{\tilde{M}}(\mathbf{M}^{"!}, \delta, f)$}
	Let $M \in \mathcal{L}(M_0)$, $\mathbf{M}^! \in \Endo_{\elli}(\tilde{M})$ and $\delta \in SD_{\mathrm{geom}, \tilde{G}\text{-reg}}(M^!) \otimes \mes(M^!)^\vee$, set
	\begin{align*}
		{}^c I^{\Endo}_{\tilde{M}}(\mathbf{M}^!, \delta, f) & = {}^c I^{\tilde{G}, \Endo}_{\tilde{M}}(\mathbf{M}^!, \delta, f) \\
		& := \sum_{s \in \Endo_{\mathbf{M}^!}(\tilde{G})} i_{M^!}(\tilde{G}, G^![s]) {}^c S^{G^![s]}_{M^!}\left( \delta[s], f^{G^![s]} \right)
	\end{align*}
	for all $f \in \orbI_{\asp}(\tilde{G}) \otimes \mes(G^!)$, where $f^{G^![s]} := \Trans_{\mathbf{G}^![s], \tilde{G}}(f)$. See also Remark \ref{rem:cS-no-B}.
\end{definition}

As in Definition--Proposition \ref{def:I-geom-Endo}, it is now routine to obtain the following

\begin{definition-proposition}
	\index{IGM-c-Endo@${}^c I^{\Endo}_{\tilde{M}}(\tilde{\gamma}, f)$}
	Let $M \in \mathcal{L}(M_0)$. There are invariant distributions
	\[ {}^c I^{\Endo}_{\tilde{M}}(\tilde{\gamma}, \cdot) = {}^c I^{\tilde{G}, \Endo}_{\tilde{M}}(\tilde{\gamma}, \cdot): \orbI_{\asp}(\tilde{G}) \otimes \mes(G) \to \CC , \]
	where $\tilde{\gamma} \in D_{\mathrm{geom}, G\text{-reg}, -}(\tilde{M}) \otimes \mes(M)^\vee$, which are characterized by
	\[ {}^c I^{\Endo}_{\tilde{M}}\left( \trans_{\mathbf{M}^!, \tilde{M}}(\delta), f \right) = {}^c I^{\Endo}_{\tilde{M}}(\mathbf{M}^!, \delta, f) \]
	for all $\mathbf{M}^! \in \Endo_{\elli}(\tilde{M})$ and $\delta \in SD_{\mathrm{geom}, \tilde{G}\text{-reg}}(M^!) \otimes \mes(M^!)\vee$.
\end{definition-proposition}

\begin{proposition}\label{prop:cpt-Endo-nonarch}
	For $f$ fixed, there exists a compact subset $\Omega \subset M(F)$ such that ${}^c I^{\Endo}_{\tilde{M}}(\tilde{\gamma}, f) = 0$ whenever $\Supp(\tilde{\gamma}) \cap \rev^{-1}(\Omega) = \emptyset$.
\end{proposition}
\begin{proof}
	This reduces to the similar property enjoyed by ${}^c S^{G^!}_{M^!}$ for various $G^! \supset M^!$, which is in turn deduced from the same property for ${}^c I^{G^!}_{M^!}$.
\end{proof}

We now state the matching between the distributions ${}^c \theta^{\tilde{G}}_{\tilde{M}}$, ${}^c I^{\tilde{G}}_{\tilde{M}}$ and their endoscopic counterparts.

\begin{theorem}\label{prop:cpt-supported-equalities}
	For all $M \in \mathcal{L}(M_0)$, we have
	\begin{enumerate}[(i)]
		\item ${}^c \theta^{\tilde{G}}_{\tilde{M}} = {}^c \theta^{\tilde{G}, \Endo}_{\tilde{M}}$.
		\item ${}^c I^{\tilde{G}}_{\tilde{M}}(\tilde{\gamma}, \cdot) = {}^c I^{\tilde{G}, \Endo}_{\tilde{M}}(\tilde{\gamma}, \cdot)$, for all $\tilde{\gamma} \in D_{\mathrm{geom}, G\text{-reg}, -}(\tilde{M}) \otimes \mes(M)^\vee$.
	\end{enumerate}
\end{theorem}

The case $M=G$ of Theorem \ref{prop:cpt-supported-equalities} is already available. In Corollary \ref{prop:cpt-supported-equalities-aux}, the general case will be reduced to the $G$-regular case of the local geometric Theorem \ref{prop:local-geometric}, which will be established in \S\ref{sec:end-stabilization}.

We proceed to study the comparison between compact and non-compact endoscopic distributions. Fix $M \in \mathcal{L}(M_0)$. The following is the endoscopic counterpart of an identity in \S\ref{sec:cIGM}.

\begin{proposition}[Cf.\ {\cite[VIII.4.1]{MW16-2}}]\label{prop:cpt-noncpt}
	For all $\tilde{\gamma} \in D_{\mathrm{geom}, \tilde{G}\text{-reg}, -}(\tilde{M}) \otimes \mes(M)^\vee$, we have
	\[ {}^c I^{\tilde{G}, \Endo}_{\tilde{M}}(\tilde{\gamma}, \cdot) = \sum_{L \in \mathcal{L}(M)} I^{\tilde{L}, \Endo}_{\tilde{M}}\left( \tilde{\gamma}, {}^c \theta^{\tilde{G}, \Endo}_{\tilde{L}}(\cdot) \right). \]
\end{proposition}
\begin{proof}
	By linearity in $\tilde{\gamma}$, it suffices to fix $\mathbf{M}^! \in \Endo_{\elli}(\tilde{M})$, $\delta \in SD_{\mathrm{geom}, \tilde{G}\text{-reg}}(M^!) \otimes \mes(M^!)^\vee$ such that $\tilde{\gamma} = \trans_{\mathbf{M}^!, \tilde{M}}(\delta)$ and show that
	\begin{equation}\label{eqn:cpt-noncpt-aux}
		\sum_{s \in \Endo_{\mathbf{M}^!}(\tilde{G})} i_{M^!}(\tilde{G}, G^![s]) {}^c S^{G^![s]}_{M^!}\left( \delta[s], f^{G^![s]} \right) = \sum_{L \in \mathcal{L}(M)} I^{\tilde{L}, \Endo}_{\tilde{M}}\left( \mathbf{M}^!, \delta, {}^c \theta^{\tilde{G}, \Endo}_{\tilde{L}}(f) \right).
	\end{equation}

	The stable counterpart \cite[VIII.4.1]{MW16-2} implies that the left-hand side of \eqref{eqn:cpt-noncpt-aux} is
	\[ \sum_{s \in \Endo_{\mathbf{M}^!}(\tilde{G})} i_{M^!}(\tilde{G}, G^![s]) \sum_{L^! \in \mathcal{L}^{G^![s]}(M^!)} S^{L^!}_{M^!}\left( \delta[s], {}^c S\theta^{G^![s]}_{L^!}\left( f^{G^![s]}\right) \right). \]
	We apply Lemma \ref{prop:sL-Ls} (with $R=M$) and Lemma \ref{prop:i-transitivity} to transform the above into
	\begin{multline*}
		\sum_{L \in \mathcal{L}^G(M)} \sum_{s^L \in \Endo_{\mathbf{M}^!}(\tilde{L})} \sum_{s_L \in \Endo_{\mathbf{L}^![s^L]}(\tilde{G})} i_{M^!}(\tilde{L}, L^![s^L]) i_{L^![s^L]}(\tilde{G}, G^![s_L]) \\
		S^{L^![s^L]}_{M^!}\left( \delta[s^L][s_L], {}^c S\theta^{G^![s_L]}_{L^![s^L]}\left( f^{G^![s_L]} \right) \right).
	\end{multline*}
	Since $z[s_L]$ is central in $L^![s^L]$, its twist on $\delta[s^L]$ may be moved onto ${}^c S\theta^{G^![s_L]}_{L^![s^L]}\left( f^{G^![s_L]} \right)$, so we get
	\[ \sum_{(L, s^L)} i_{M^!}(\tilde{L}, L^![s^L]) S^{L^![s^L]}_{M^!}\left( \delta[s^L], A(L, s^L) \right) \]
	where
	\begin{align*}
		A(L, s^L) & := \sum_{s_L} i_{L^![s^L]}(\tilde{G}, G^![s_L]) {}^c S\theta^{G^![s_L]}_{L^![s^L]}\left( f^{G^![s_L]} \right)[s_L] \\
		& = {}^c \theta^{\tilde{G}, \Endo}_{\tilde{L}}\left( \mathbf{L}^![s^L], f \right) = \Trans_{\mathbf{L}^![s^L], \tilde{L}} \left( {}^c \theta^{\tilde{G}, \Endo}_{\tilde{L}}(f) \right).
	\end{align*}

	The left-hand side of \eqref{eqn:cpt-noncpt-aux} thus equals
	\[ \sum_L \sum_{s^L} i_{M^!}(\tilde{L}, L^![s^L]) S^{L^![s^L]}_{M^!}\left( \delta[s^L], \Trans_{\mathbf{L}^![s^L], \tilde{L}} \left( {}^c \theta^{\tilde{G}, \Endo}_{\tilde{L}}(f) \right) \right). \]
	By definition, this equals the right-hand side of \eqref{eqn:cpt-noncpt-aux}.
\end{proof}

\begin{corollary}[Cf.\ {\cite[VIII.4.3]{MW16-2}}]\label{prop:cpt-supported-equalities-aux}
	Suppose that
	\begin{itemize}
		\item $I^{\tilde{G}}_{\tilde{M}}(\tilde{\gamma}, \cdot) = I^{\tilde{G}, \Endo}_{\tilde{M}}(\tilde{\gamma}, \cdot)$ for all $M \in \mathcal{L}(M_0)$ and all $\tilde{\gamma} \in D_{\mathrm{geom}, G\text{-reg}, -}(\tilde{M}) \otimes \mes(M)^\vee$.
		\item the assertions in Theorem \ref{prop:cpt-supported-equalities} are known when $M$ is replaced by any $L \in \mathcal{L}(M)$ with $L \neq M$, or when $G$ is replaced by any $L \in \mathcal{L}(M)$ with $L \neq G$.
	\end{itemize}
	Then the assertions (i) and (ii) in Theorem \ref{prop:cpt-supported-equalities} hold for $M$.
\end{corollary}
\begin{proof}
	The cited proof is reproduced below for the reader's convenience. Expanding ${}^c I^{\tilde{G}}_{\tilde{M}}$ and ${}^c I^{\tilde{G}, \Endo}_{\tilde{M}}$ in terms of non-compact analogues (by Proposition \ref{prop:cpt-noncpt}), the assumptions imply
	\[ {}^c I^{\tilde{G}}_{\tilde{M}}(\tilde{\gamma}, f) - {}^c I^{\tilde{G}, \Endo}_{\tilde{M}}(\tilde{\gamma}, f) = I^{\tilde{M}}\left( \tilde{\gamma}, {}^c \theta^{\tilde{G}}_{\tilde{M}}(f) - {}^c \theta^{\tilde{G}, \Endo}_{\tilde{M}}(f) \right). \]
	Denote the difference inside $I^{\tilde{M}}(\tilde{\gamma}, \cdot)$ as $\varphi$. The equality combined with the compact-support property of ${}^c I^{\tilde{G}}_{\tilde{M}}$ and ${}^c I^{\tilde{G}, \Endo}_{\tilde{M}}$ implies that $\varphi \in \orbI_{\asp}(\tilde{M}) \otimes \mes(M)$. On the other hand, a comparison of the descent formulas for ${}^c \theta^{\tilde{G}}_{\tilde{M}}$ (see \cite[VIII.1.6 Lemme]{MW16-2}) and ${}^c \theta^{\tilde{G}, \Endo}_{\tilde{M}}$ (Proposition \ref{prop:thetaEndo-descent}) shows that $\varphi$ is cuspidal.
	
	If $M = G$, it is clear that $\varphi = 0$. If $M \neq G$, we apply Proposition \ref{prop:cthetaGM-nonarch} and its endoscopic analogue (which follows from the stable case in \cite[VIII.3.7 Proposition]{MW16-2}) to see
	\[ \sum_{S \in \mathcal{P}(M)} \omega_S(X) I^{\tilde{M}} \left( \pi, \nu_S, X, \varphi \right) = 0 \]
	for all $X \in \mathfrak{a}_{M,F}$ and $\pi \in D_{\elli, -}(\tilde{M}) \otimes \mes(M)^\vee$, provided that $\nu_S$ is sufficiently deep for each $S \in \mathcal{P}(M)$. Since $\varphi \in \orbI_{\asp}(\tilde{M}) \otimes \mes(M)$, we can shift $\nu_S$ to $0$ in the displayed equality to obtain $I^{\tilde{M}}(\pi, X, \varphi) = I^{\tilde{M}}(\pi, X, 0, \varphi) = 0$.
	
	By varying $X$ and $\pi$, we deduce $\varphi = 0$ in the case $M \neq G$, since $\varphi \in \orbI_{\asp, \cusp}(\tilde{M}) \otimes \mes(M)$. The assertion follows.
\end{proof}

\section{Construction of \texorpdfstring{$\epsilon_{\tilde{M}}$}{epsilonM}}\label{sec:epsilon-nonarch}
This section is parallel to \cite[VIII.4.4]{MW16-2}. Our situation is simpler since the matching of Shalika germs (Theorem \ref{prop:germ-matching}) is already available to us.

As before, fix $M \in \mathcal{L}(M_0)$ and assume inductively that whenever $\tilde{\gamma}$ is $\tilde{G}$-regular,
\begin{gather*}
	I^{\tilde{L}}_{\tilde{M}}(\tilde{\gamma}, \cdot) = I^{\tilde{L}, \Endo}_{\tilde{M}}(\tilde{\gamma}, \cdot), \quad L \in \mathcal{L}(M), \; L \neq G, \\
	I^{\tilde{G}}_{\tilde{L}}(\tilde{\gamma}, \cdot) = I^{\tilde{G}, \Endo}_{\tilde{L}}(\tilde{\gamma}, \cdot), \quad L \in \mathcal{L}(M), \; L \neq M.
\end{gather*}
Similarly, we assume that the assertions of Theorem \ref{prop:cpt-supported-equalities} hold when $G$ (resp.\ $M$) becomes smaller (resp.\ larger).

\begin{proposition}\label{def:epsilonM-nonarch}
	\index{epsilonM@$\epsilon_{\tilde{M}}$}
	There exists a linear map
	\[ \epsilon_{\tilde{M}}: \orbI_{\asp}(\tilde{G}) \otimes \mes(G) \to \orbI_{\mathrm{ac}, \asp, \cusp}(\tilde{M}) \otimes \mes(M) \]
	characterized uniquely by
	\[ I^{\Endo}_{\tilde{M}}(\tilde{\gamma}, f) - I_{\tilde{M}}(\tilde{\gamma}, f) = I^{\tilde{M}}\left( \tilde{\gamma}, \epsilon_{\tilde{M}}(f) \right) \]
	for all $f \in \orbI_{\asp}(\tilde{G}) \otimes \mes(G)$ and $\tilde{\gamma} \in D_{\mathrm{geom}, -}(\tilde{M}) \otimes \mes(M)^\vee$.
\end{proposition}
\begin{proof}
	We may assume $M \neq G$. The uniqueness of $f \mapsto \epsilon_{\tilde{M}}(f)$ is clear. Let us fix $f$ and prove that for each stable semisimple conjugacy class $\mathcal{O} \subset M(F)$, there exists $\varphi_1 \in \orbI_{\asp}(\tilde{M}) \otimes \mes(M)$ such that
	\begin{equation}\label{eqn:epsilonM-nonarch-aux0}
		I^{\Endo}_{\tilde{M}}(\tilde{\gamma}, f) - I_{\tilde{M}}(\tilde{\gamma}, f) = I^{\tilde{M}}(\tilde{\gamma}, \varphi_1)
	\end{equation}
	for all $\tilde{\gamma} \in D_{\mathrm{geom}, G\text{-equi}, -}(\tilde{M}) \otimes \mes(M)^\vee$ such that $\Supp(\tilde{\gamma})$ is close to $\rev^{-1}(\mathcal{O})$.
	
	The induction hypotheses and the matching of weighted Shalika germs (Theorem \ref{prop:germ-matching}) imply that
	\begin{align*}
		I^{\Endo}_{\tilde{M}}(\tilde{\gamma}, f) - I_{\tilde{M}}(\tilde{\gamma}, f) & =
		\sum_{L \in \mathcal{L}(M)} \left( I^{\Endo}_{\tilde{L}}\left( g^{\tilde{L}, \Endo}_{\tilde{M}, \mathcal{O}} (\tilde{\gamma}), f \right) - I_{\tilde{L}}\left( g^{\tilde{L}}_{\tilde{M}, \mathcal{O}}(\tilde{\gamma}) , f \right) \right) \\
		& = I^{\Endo}_{\tilde{M}}\left( g^{\tilde{M}}_{\tilde{M}, \mathcal{O}}(\tilde{\gamma}), f \right) - I_{\tilde{M}}\left( g^{\tilde{M}}_{\tilde{M}, \mathcal{O}}(\tilde{\gamma}), f \right).
	\end{align*}
	
	Since $D_{\mathrm{geom}, -}(\tilde{M}, \mathcal{O})$ is finite-dimensional, we may choose $\varphi_1$ such that
	\[ I^{\tilde{M}}(\tilde{\eta}, \varphi_1) = I^{\Endo}_{\tilde{M}}(\tilde{\eta}, f) - I_{\tilde{M}}(\tilde{\eta}, f) \]
	for all $\tilde{\eta} \in D_{\mathrm{geom}, -}(\tilde{M}, \mathcal{O}) \otimes \mes(M)^\vee$. Hence
	\[ I^{\Endo}_{\tilde{M}}\left( g^{\tilde{M}}_{\tilde{M}, \mathcal{O}}(\tilde{\gamma}), f \right) - I_{\tilde{M}}\left( g^{\tilde{M}}_{\tilde{M}, \mathcal{O}}(\tilde{\gamma}), f \right)
	= I^{\tilde{M}}\left( g^{\tilde{M}}_{\tilde{M}, \mathcal{O}}(\tilde{\gamma}), \varphi_1 \right) = I^{\tilde{M}}(\tilde{\gamma}, \varphi_1)
	\]
	by the definition of Shalika germs. This proves \eqref{eqn:epsilonM-nonarch-aux0}.

	Next, we claim that there exist a compact subset $\Omega \subset M(F)$ and $\varphi_2 \in \orbI_{\mathrm{ac}, \asp}(\tilde{M}) \otimes \mes(M)$ such that
	\begin{equation}\label{eqn:epsilonM-nonarch-aux1}
		I^{\Endo}_{\tilde{M}}(\tilde{\gamma}, f) - I_{\tilde{M}}(\tilde{\gamma}, f) = I^{\tilde{M}}(\tilde{\gamma}, \varphi_2)
	\end{equation}
	for all $f$ and $\tilde{\gamma} \in D_{\mathrm{geom}, G\text{-reg}, -}(\tilde{M}) \otimes \mes(M)^\vee$ such that $\Supp(\tilde{\gamma}) \cap \rev^{-1}(\Omega) = \emptyset$.
	
	To establish \eqref{eqn:epsilonM-nonarch-aux1}, apply Proposition \ref{prop:cpt-noncpt} and its counterpart for ${}^c I_{\tilde{M}}$ to express the left-hand side as
	\[ {}^c I^{\Endo}_{\tilde{M}}(\tilde{\gamma}, f) - {}^c I_{\tilde{M}}(\tilde{\gamma}, f) -
	\sum_{\substack{L \in \mathcal{L}(M) \\ L \neq G}} \left( I^{\tilde{L}, \Endo}_{\tilde{M}}\left( \tilde{\gamma}, {}^c \theta^{\tilde{G}, \Endo}_{\tilde{L}}(f) \right) - I^{\tilde{L}}_{\tilde{M}}\left(\tilde{\gamma}, {}^c \theta^{\tilde{G}}_{\tilde{L}}(f) \right) \right). \]
	
	By induction hypotheses, $L \neq G$ (resp.\ $L \neq M$) implies $I^{\tilde{L}, \Endo}_{\tilde{M}} = I^{\tilde{L}}_{\tilde{M}}$ (resp.\ ${}^c \theta^{\tilde{G}, \Endo}_{\tilde{L}} = {}^c \theta^{\tilde{G}}_{\tilde{L}}$). Hence the previous expression becomes
	\[ {}^c I^{\Endo}_{\tilde{M}}(\tilde{\gamma}, f) - {}^c I_{\tilde{M}}(\tilde{\gamma}, f) - I^{\tilde{M}}\left( \tilde{\gamma}, {}^c \theta^{\tilde{G}, \Endo}_{\tilde{M}}(f) - {}^c \theta^{\tilde{G}}_{\tilde{M}}(f) \right). \]

	Combining \S\ref{sec:cIGM} and Proposition \ref{prop:cpt-Endo-nonarch}, the term ${}^c I^{\Endo}_{\tilde{M}}(\cdot, f) - {}^c I_{\tilde{M}}(\cdot, f)$ is seen to vanish off some compact $\Omega \subset M(F)$. One takes $\varphi_2 := {}^c \theta^{\tilde{G}, \Endo}_{\tilde{M}}(f) - {}^c \theta^{\tilde{G}}_{\tilde{M}}(f)$ to obtain \eqref{eqn:epsilonM-nonarch-aux1}.
	
	Now we can use an argument of partition of unity to glue $\varphi_1$ (for various $\mathcal{O}$) and $\varphi_2$ from \eqref{eqn:epsilonM-nonarch-aux0} and \eqref{eqn:epsilonM-nonarch-aux1}, and obtain $\varphi \in \orbI_{\mathrm{ac}, \asp}(\tilde{M}) \otimes \mes(M)$ such that
	\begin{equation}\label{eqn:epsilonM-nonarch-aux2}
		I^{\Endo}_{\tilde{M}}(\tilde{\gamma}, f) - I_{\tilde{M}}(\tilde{\gamma}, f) = I^{\tilde{M}}(\tilde{\gamma}, \varphi), \quad \tilde{\gamma} \in D_{\mathrm{geom}, G\text{-reg}, -}(\tilde{M}) \otimes \mes(M)^\vee .
	\end{equation}
	
	We have to extend \eqref{eqn:epsilonM-nonarch-aux2} to general $\tilde{\gamma}$. It suffices to fix $\mathcal{O} \subset M(F)$ and assume that $\Supp(\tilde{\gamma}) = \rev^{-1}(\mathcal{O})$. The separation property of Shalika germs (see \S\ref{sec:Shalika-germs}) furnishes $\tilde{\gamma}'$ with $G$-regular support close to $\rev^{-1}(\mathcal{O})$, such that $g^{\tilde{M}}_{\tilde{M}, \mathcal{O}}(\tilde{\gamma}') = \tilde{\gamma}$. Then
	\begin{align*}
		I^{\tilde{M}}(\tilde{\gamma}, \varphi) = I^{\tilde{M}}(\tilde{\gamma}', \varphi) & = I^{\Endo}_{\tilde{M}}(\tilde{\gamma}', f) - I_{\tilde{M}}(\tilde{\gamma}', f) \\
		& = \sum_{L \in \mathcal{L}(M)} \left( I^{\Endo}_{\tilde{L}}\left( g^{\tilde{L}, \Endo}_{\tilde{M}}(\tilde{\gamma}'), f \right) - I_{\tilde{L}}\left( g^{\tilde{L}}_{\tilde{M}}(\tilde{\gamma}'), f \right) \right).
	\end{align*}
	In the last expression, all the differences with $L \neq M$ vanish by the induction assumption $I^{\Endo}_{\tilde{L}} = I_{\tilde{L}}$ and the matching of Shalika germs. What remains is $I^{\Endo}_{\tilde{M}}(\tilde{\gamma}, f) - I_{\tilde{M}}(\tilde{\gamma}, f)$. This extends \eqref{eqn:epsilonM-nonarch-aux2} to all $\tilde{\gamma}$.
		
	All in all, we take $\epsilon_{\tilde{M}}(f) := \varphi$, which is uniquely characterized and linear in $f$. It remains to show $\epsilon_{\tilde{M}}(f) \in \orbI_{\mathrm{ac}, \asp, \cusp}(\tilde{M}) \otimes \mes(M)$. Indeed, for all $R \in \mathcal{L}^M(M_0)$, $\tilde{\eta} \in D_{\mathrm{geom}, -}(\tilde{R}) \otimes \mes(R)^\vee$ with $R \neq M$, a comparison between descent formulas (Proposition \ref{prop:orbint-weighted-descent-nonArch}, Corollary \ref{prop:descent-orbint-Endo-2}) using induction hypotheses entails
	\[ I^{\tilde{R}}\left(\tilde{\eta}, \epsilon_{\tilde{M}}(f)_{\tilde{R}} \right) = I^{\Endo}_{\tilde{M}}\left( \tilde{\eta}^{\tilde{M}}, f \right) - I_{\tilde{M}}\left( \tilde{\eta}^{\tilde{M}}, f \right) = 0. \]
	This completes the proof.
\end{proof}

\begin{remark}\label{rem:epsilonM-nonarch}
	As a by-product of the proof, we see that $\epsilon_{\tilde{M}}(f)$ can be written as
	\[ {}^c \theta^{\Endo}_{\tilde{M}}(f) - {}^c \theta_{\tilde{M}}(f) + \;\text{an element of}\; \orbI_{\asp}(\tilde{M}) \otimes \mes(M). \]
\end{remark}

\begin{corollary}\label{prop:epsilonM-nonarch-Schwartz}
	Let $f \in \orbI_{\asp}(\tilde{G}) \otimes \mes(G)$. Then $\epsilon_{\tilde{M}}(f)$ is Schwartz, and for every $\pi \in D_{\mathrm{temp}, -}(\tilde{M}) \otimes \mes(M)^\vee$, the function
	\[ \lambda \mapsto I^{\tilde{M}}\left(\pi, \lambda, \epsilon_{\tilde{M}}(f) \right), \quad \lambda \in i\mathfrak{a}_M^* / i\mathfrak{a}_{M, F}^\vee \]
	extends to a rational function on $\mathfrak{a}_{M, \CC}^* / i\mathfrak{a}_{M, F}^\vee$ with only finitely many polar hyperplanes of the form $e^{\lrangle{\cdot, \check{\alpha}}} = \mathrm{const}$ where $\alpha \in \Sigma^G(A_M)$.
\end{corollary}
\begin{proof}
	In view of Remark \ref{rem:epsilonM-nonarch}, this reduces to the corresponding properties of ${}^c \theta_{\tilde{M}}(f)$ (Proposition \ref{prop:cthetaGM-nonarch}), ${}^c \theta^{\Endo}_{\tilde{M}}(f)$ (see \S\ref{sec:cIGMEndo}), and of $\orbI_{\asp}(\tilde{M}) \otimes \mes(M)$ (see \S\ref{sec:cphi-nonarch}).
\end{proof}

\chapter{Weighted characters: Archimedean case}\label{sec:weighted-char-arch}
Let $F$ be a an Archimedean local field. Consider a covering of metaplectic type $\rev: \tilde{G} \to G(F)$ where $G = \prod_{i \in I} \GL(n_i) \times \Sp(W)$; here we fix the auxiliary data $(W, \lrangle{\cdot|\cdot})$, $\psi$ as in \S\ref{sec:cpt-nonarch}. We also fix a maximal compact subgroup $K \subset G(F)$ arising from a Cartan involution $\theta$, in good position relative to $M_0$.

Let $M \in \mathcal{L}(M_0)$. Our goal is to define the compactly supported variant ${}^c I_{\tilde{M}}(\tilde{\gamma}, f)$ of $I_{\tilde{M}}(\tilde{\gamma}, f)$ for $\tilde{\gamma} \in D_{\mathrm{orb}, -}(\tilde{M}) \otimes \mes(M)^\vee$ and $f \in \orbI_{\asp}(\tilde{G}, \tilde{K}) \otimes \mes(G)$, and formulate the matching theorems. In contrast with the non-Archimedean case, a new problem here is that the canonically normalized weighted character $J_{\tilde{M}}(\pi_\lambda, f)$, as a meromorphic family in $\lambda$, may have infinitely many polar hyperplanes. For this reason, we revert to Arthur's earlier recipe of using normalizing factors. This is inspired by \cite[IX.5]{MW16-2}, where a more complicated ``rational normalization'' is used in the twisted setting.

Without loss of generality, we may assume $\tilde{G} = \Mp(W)$ with $\dim W = 2n$. Based on the results of \cite{Li12b} (resp.\ \cite{Ar89a}), the normalizing factors for $\Mp(W)$ (resp.\ $\SO(2n+1)$) are explicitly and canonically given. Furthermore, by matching the infinitesimal characters of discrete series between $\Mp(2n)$ and $\SO(2n+1)$, their normalizing factors differ by a precise power of $\sqrt{2}$. This is the content of \S\ref{sec:normalized-R-arch}, which is based on explicit computations. In this way, we obtain the normalized weighted characters $J^r_{\tilde{M}}(\pi, f)$ for $\pi \in \Pi_{\mathrm{unit}, -}(\tilde{M})$, as well as the variant ${}^r \phi_{\tilde{M}}$ of $\phi_{\tilde{M}}$.

In \S\ref{sec:cphi-arch}, we proceed to define the compactly supported weighted characters ${}^{c, r} J_{\tilde{M}}(\pi, f)$ out from $J^r_{\tilde{M}}(\pi, f)$. We also obtain the variant ${}^c \phi_{\tilde{M}}$ of $\phi_{\tilde{M}}$, from which ${}^c I_{\tilde{M}}(\tilde{\gamma}, f)$ is built. The same construction also applies to the stable side, although this was not documented before.

The triad of maps ${}^c \phi_{\tilde{M}}$, $\phi_{\tilde{M}}$ and ${}^r \phi_{\tilde{M}}$ are related by families of maps ${}^c \theta_{\tilde{M}}$, ${}^{c, r} \theta_{\tilde{M}}$ and ${}^r \theta_{\tilde{M}}$, as $M \in \mathcal{L}(M_0)$ varies. Together with the distributions ${}^c I_{\tilde{M}}(\tilde{\gamma}, \cdot)$, each of them comes with an endoscopic counterpart; see \S\ref{sec:cI-Endo-arch}. The formulation of the respective matching theorems is the content of \S\ref{sec:matching-theta-arch}.

Among these statements, the case of ${}^r \theta_{\tilde{M}}$ (Theorem \ref{prop:rtheta-matching-arch}) is immediately within reach. Roughly speaking, it reduces to the stabilization of some factors $r^{\tilde{G}}_{\tilde{M}}(\pi)$ made from normalizing factors; we seek to express them via similar factors $s^{G^![s]}_{M^!}(\phi^!)$ on the stable side, where $\mathbf{M}^! \in \Endo_{\elli}(\tilde{M})$ and $s \in \Endo_{\mathbf{M}^!}(\tilde{G})$. This is achieved in \S\ref{sec:easy-stabilization} by reducing to the case of $\SO(2n+1)$ settled in \cite{Ar99b}.

\section{Normalizing factors}\label{sec:normalized-R-arch}
\subsection{The setting}
Consider a general covering $\tilde{G} \twoheadrightarrow G(F)$ over a local field $F$ of characteristic zero and a Levi subgroup $M$ of $G$. For $\pi \in \Pi_-(\tilde{M})$, the normalized intertwining operators take the form
\[ R_{\tilde{Q}|\tilde{P}}(\pi_\lambda) = r_{\tilde{Q}|\tilde{P}}(\pi_\lambda)^{-1} J_{\tilde{Q}|\tilde{P}}(\pi_\lambda) : I_{\tilde{P}}(\pi_\lambda) \to I_{\tilde{Q}}(\pi_\lambda), \quad P, Q \in \mathcal{P}(M), \]
where $J_{\tilde{Q}|\tilde{P}}(\pi_\lambda)$ is the standard intertwining operator, and $r_{\tilde{Q}|\tilde{P}}(\pi_\lambda)$ is a meromorphic function in $\lambda \in \mathfrak{a}_{M, \CC}^*$. The \emph{normalizing factors} $r_{\tilde{Q}|\tilde{P}}(\pi_\lambda)$ are subject to the conditions in \cite[\S 3.1]{Li12b}, which are in turn based on \cite[\S 2]{Ar89a}. Their definition is eventually reduced to the case when $M$ is maximal proper and $\pi \in \Pi_{2, -}(\tilde{M})$ (i.e.\ essentially square-integrable, or discrete series).
\index{rQP@$r_{\tilde{Q}{"|}\tilde{P}}(\pi)$, $R_{\tilde{Q}{"|}\tilde{P}}(\pi)$}

Specifically, fix a maximal $F$-torus $T \subset G$ and denote the elements of $\Sigma(A_M)$ as $[\alpha]$, viewed as equivalence classes of $\alpha \in \Sigma(G, T)$ restricting to $[\alpha]$, we have the property
\begin{equation}
	r_{\tilde{Q}|\tilde{P}}(\pi_\lambda) = \prod_{[\alpha] \in \Sigma_Q^{\mathrm{red}} \cap \Sigma_P^{\mathrm{red}}} r_{[\alpha]}(\pi_\lambda),
\end{equation}
where $r_{[\alpha]}$ is the normalizing factor relative to the Levi subgroup $\widetilde{M_{[\alpha]}} \subset \tilde{G}$ containing $\tilde{M}$ minimally with $\pm\alpha$ in $\Sigma(M_{[\alpha]}, T)$.

Now revert to the case where $\tilde{G}$ is of metaplectic type.
\begin{itemize}
	\item For non-Archimedean $F$, the existence of normalizing factors for general coverings is settled in \cite[\S 3.1]{Li12b}, although they are not canonically given except in the unramified case (see \S\ref{sec:spec-nr}). This suffices for the purposes of this work.
	
	\item For Archimedean $F$, we will have to compare them with their counterparts for $\SO(2n+1)$, and this requires explicit formulas.
\end{itemize}

We remark that in the uncovered case, Langlands proposes an explicit choice of normalizing factors in terms of local factors. Although the metaplectic counterpart is accessible through $\Theta$-correspondence, we do not need this level of precision in this work.

Henceforth, $F \in \{\R, \CC\}$ and $G = \Sp(W)$, $\rev: \tilde{G} \to G(F)$. Our strategy is to work out the formula \cite[\S 3.1, (10)]{Li12b}. In view of the reduction procedure for normalizing factors, we fix the following data
\begin{itemize}
	\item a maximal proper Levi subgroup $M$ of $G$ and $P, Q \in \mathcal{P}(M)$ which are opposite to each other;
	\item $\pi \in \Pi_{2, -}(\tilde{M})$.
\end{itemize}
Then there is only one $[\alpha]$, we have $r_{\tilde{Q}|\tilde{P}}(\pi_\lambda) = r_{[\alpha]}(\pi_\lambda)$, and $M_{[\alpha]} = G$. The case when $\tilde{G}$ is of metaplectic type will then follow immediately. 

\begin{notation}
	\index{GammaRC@$\Gamma_{\R}$, $\Gamma_{\CC}$}
	We write $\Gamma_{\R}(z) := \pi^{-z/2} \Gamma\left(\frac{z}{2}\right)$ and $\Gamma_{\CC}(z) := 2(2\pi)^{-z} \Gamma(z)$, for $z \in \CC$.
\end{notation}

\subsection{The real case}
Let $F = \R$. Take a maximal $\R$-torus $T \subset M$ such that $T/A_M$ is anisotropic. Define $\rho^M$ to be the half-sum of positive roots of $(\mathfrak{m}_{\CC}, \mathfrak{t}_{\CC})$ relative to some choice of Borel subgroup over $M_{\CC}$. Denote by $\sigma$ the complex conjugation in $\Gal(\CC|\R)$.

Fix a Cartan involution $\theta$ of $\mathfrak{g}$ so that the associated maximal compact subgroup is in good position relative to $M$. Fix an invariant bilinear form $\mathsf{B}$ on $\mathfrak{g}$ such that $X \mapsto -\mathsf{B}(X, \theta X)$ is positive definite.

\begin{notation}
	Unless otherwise stated, we will use the \emph{twofold} version of the metaplectic covering $\rev: \tilde{G} \to G(\R)$ in what follows.
\end{notation}

Take $\dd\chi \in \mathfrak{t}^*_{\CC}$ such that $\dd \chi + \rho^M$ represents the infinitesimal character of $\pi$. By Harish-Chandra's theory of discrete series, $\chi$ lifts to a genuine character of $\tilde{T}$. Take $\lambda_0 \in \mathfrak{t}_{\CC}^*$ and $\mu \in i\mathfrak{t}^*$ such that
\[ \lrangle{\dd\chi, H} = \lrangle{\lambda_0, H - \sigma H} + \frac{1}{2} \lrangle{\mu - \rho^M, H + \sigma H}, \quad H \in \mathfrak{t}_{\CC}. \]
If $\pi$ is replaced by $\pi_\lambda$ with $\lambda \in i\mathfrak{a}^*_M$, then $(\lambda_0, \mu)$ is replaced by $(\lambda_0 + \lambda, \mu + \lambda)$.

\index{SigmaPGT@$\Sigma_P(G, T)$}
Let $\Sigma_P(G, T) \subset \Sigma(G,T)$ be the subset of roots restricting to $\Sigma_P$. An element of $\Sigma(G, T)$ is said to be \emph{real} if it is $\sigma$-invariant. Since $M$ is maximal and $T/A_M$ is anisotropic, there is at most one real root in $\Sigma_P(G, T)$. Let $\alpha$ be such a real root, and let $H_{\alpha} \in \mathfrak{t}$ be characterized by $\mathsf{B}(H, H_{\alpha}) = \lrangle{\alpha, H}$ for all $H \in \mathfrak{t}$. Take a triplet $(H',X', Y')$ such that $H' = 2\lrangle{\alpha, H_{\alpha}}^{-1} H_{\alpha}$, $X' \in \mathfrak{g}_{\alpha}$, $X' \in \mathfrak{g}_{-\alpha}$, and
\[ [H', X'] = 2X', \quad [H', Y'] = -2Y', \quad \quad [X', Y'] = H' . \]
From the triplet we obtain a homomorphism $\SL(2) \to G$. Denote the pullback along $\rev$ of this homomorphism as $\Phi_\alpha: \widetilde{\SL}(2) \to \tilde{G}$. Calculating with
\[
	H' = \begin{pmatrix} 1 & \\ & -1 \end{pmatrix}, \quad
	X' = \begin{pmatrix} & 1 \\ & \end{pmatrix}, \quad
	Y' = \begin{pmatrix} & \\ 1 & \end{pmatrix}
\]
in $\mathfrak{sl}_2$, where the blank entries are zero, we see that
$\pi (X' - Y') = \bigl(\begin{smallmatrix} 0 & \pi \\ -\pi & 0 \end{smallmatrix}\bigr)$. Put
\[ \gamma := \exp\left( \dd\Phi_\alpha (\pi (X' - Y')) \right) \; \in \tilde{T}. \]

Choose a Borel subgroup $B$ of $G_{\CC}$ such that $T_{\CC} \subset B \subset P_{\CC}$, and that $B \cap M_{\CC}$ is the Borel subgroup of $M_{\CC}$ chosen before. Define $\rho_B$ in this way. Since we work with twofold coverings here ($m=2$), the discussions above imply $\gamma^4 = 1$ and there exists $k \in \{0, 1, 2, 3\}$ such that
\[ (-1)^{\lrangle{\rho_P, \check{\alpha}}} \chi(\gamma) = \exp\left(\frac{2\pi k}{4i}\right). \]

When $\frac{k}{4} = 0$ (resp.\ $\frac{k}{4} = \frac{1}{2}$), put $N_0 := 1$ (resp.\ $N_0 := 0$). Furthermore, we follow \cite[Remarques 3.2.2, 3.2.3]{Li12b} to define
\begin{equation}\label{eqn:normalizing-G}
	\mathsf{G}(z) := \begin{cases}
		\displaystyle\frac{\Gamma_{\R}(z + N_0)}{\Gamma_{\R}(z + N_0 + 1)}, & \frac{k}{4} \in \{0, \frac{1}{2}\} \\
		\sqrt{2} \cdot \displaystyle\frac{\Gamma_{\R}(2z)}{\Gamma_{\R}(2z + 1)}, & \frac{k}{4} \in \{ \frac{1}{4}, \frac{3}{4} \}.
\end{cases}\end{equation}
\index{Gz@$\mathsf{G}(z)$}

Recall that $\rho_P$ is defined as the half-sum of roots appearing in $U_P$, and $\lrangle{\rho_P, \check{\alpha}}$ is interpreted accordingly. This quantity appeared in Harish-Chandra's work \cite[p.178, p.189]{HC76} on $\mu$-functions, and he denotes $\lrangle{\rho_P, \check{\alpha}}$ by $\rho_\alpha$. One can work with $\rho_B$ instead, by virtue of the following result.

\begin{lemma}\label{prop:rho-B-rho-P}
	In the circumstance above, we have $\lrangle{\rho_B, \check{\alpha}} = \lrangle{\rho_P, \check{\alpha}}$.
\end{lemma}
\begin{proof}
	By definition, we have $\rho_B = \rho_P + \rho^M_{B \cap M_{\CC}}$, where $\rho^M_{B \cap M_{\CC}}$ is the half-sum of $B \cap {M_{\CC}}$-positive roots in $M_{\CC}$. Since $T/A_M$ is anisotropic, $\sigma \rho^M_{B \cap M_{\CC}} = -\rho^M_{B \cap M_{\CC}}$. On the other hand, $\sigma\check{\alpha} = \check{\alpha}$ since $\alpha$ is real. It follows that $\lrangle{\rho^M_{B \cap M_{\CC}}, \check{\alpha}} = 0$ since $\lrangle{\cdot, \cdot}$ is $\sigma$-invariant.
\end{proof}

\begin{notation}
	Let us say that we are in ``case I'' if there exists a real root in $\Sigma_P(G, T)$, denoted as $\alpha_0$; otherwise we are in ``case II''. In case I, we pick for each $\Gal(\CC|\R)$-orbit in $\Sigma_P(G, T) \smallsetminus \{\alpha_0\}$ a representative $\alpha$ such that
	\begin{equation}\label{eqn:alpha-normalizing-choice}
		\lrangle{\sigma\mu - \mu, \check{\alpha}} \in \Z_{\leq 0}.
	\end{equation}
\end{notation}

\begin{definition-proposition}\label{prop:normalizing-real-0}
	In the circumstance above, the normalizing factors can be taken as
	\[ r_{\tilde{Q}|\tilde{P}}(\pi_\lambda) = \begin{cases}
		\displaystyle\prod_{\substack{\alpha \in \Sigma_P(G, T) \\ \bmod\; \Gal(\CC|\R) \\ \alpha \neq \alpha_0 }} \frac{\Gamma_{\CC}(\lrangle{\mu + \lambda, \check{\alpha}})}{\Gamma_{\CC}(\lrangle{\mu + \lambda, \;\check{\alpha}} + 1)} \cdot \mathsf{G}(\lrangle{\mu + \lambda, \check{\alpha}_0}), & \text{case I}, \\
		\displaystyle\prod_{\substack{\alpha \in \Sigma_P(G, T) \\ \bmod\; \Gal(\CC|\R) }} \frac{\Gamma_{\CC}(\lrangle{\mu + \lambda, \; \check{\alpha}})}{\Gamma_{\CC}(\lrangle{\mu + \lambda, \;\check{\alpha}} + 1)}, & \text{case II}, \\
	\end{cases}\]
	for all $\lambda \in \mathfrak{a}_{M, \CC}^*$. The representatives $\alpha$ in each orbit are chosen as in \eqref{eqn:alpha-normalizing-choice}.
\end{definition-proposition}
\begin{proof}
	This is the formula \cite[(10)]{Li12b} which intervenes in the Théorème 3.2.1 thereof.
\end{proof}

We have to be more precise. Since $M$ is maximal proper in $G$ and admits a maximal torus $T$ with $T/A_M$ anisotropic, there exist $a \in \{1, 2\}$ and a symplectic subspace $W^\flat \subset W$ of codimension $2a$ such that
\[ M = \GL(a) \times \Sp(W^\flat) \hookrightarrow \Sp(W) = G. \]
Decompose $T$ into $T_{\GL} \times T^\flat$ accordingly, $T^\flat \subset \Sp(W^\flat)$ being an anisotropic torus. Let $n := \frac{1}{2} \dim W$.

\begin{itemize}
	\item As before, choose a Borel subgroup $B$ of $G_{\CC}$ appropriately.
	\item Take a symplectic basis for $W^\flat_{\CC}$, enlarge it to a symplectic basis of $W_{\CC}$, so that $X^*(T_{\CC})$ has the standard basis $e^*_1, \ldots, e^*_n$ satisfying
	\[ X^*(T_{\GL, \CC}) = \bigoplus_{k=1}^a \Z e^*_k, \quad X^*(T^\flat_{\CC}) = \bigoplus_{k=a+1}^n \Z e^*_k , \]
	and the simple roots in $\Sigma(G, T)$ relative to $B$ are $2e^*_n$ and $e^*_i - e^*_{i+1}$, for $1 \leq i < n$.
	\item In particular, $X_*(A_M) = \Z\left(e_{1,*} + \cdots + e_{a,*} \right)$ and $T_{\GL} \simeq \Gm$ (resp.\ $\Res_{\CC|\R} \Gm$) when $a=1$ (resp.\ $a=2$).
	\item In case I, the real root $\alpha_0$ takes the form
	\begin{compactitem}
		\item $2e^*_1$ when $a=1$,
		\item $e^*_1 + e^*_2$ when $a=2$, in which case $\sigma$ swaps $e^*_1$ and $e^*_2$.
	\end{compactitem}
	Indeed, $\sigma e^*_i = -e^*_i$ for all $i > a$ since $T^\flat$ is anisotropic; they cannot appear in $\alpha_0$.
\end{itemize}

Regard $G' := \SO(2n+1)$ as the ``principal'' endoscopic group associated with $(n,0) \in \Endo_{\elli}(\tilde{G})$. With the choices above, one can take in $G'$ a Levi subgroup $M'$, opposite $P', Q' \in \mathcal{P}(M')$, a maximal torus $T' \subset M'$ such that $T'/A_{M'}$ is anisotropic, together with a Borel subgroup $B' \subset G'_{\CC}$ such that $T'_{\CC} \subset B' \subset P'_{\CC}$. Furthermore, there is an isomorphism $T' \rightiso T$ of tori, restricting to $A_M \rightiso A_{M'}$.
\begin{itemize}
	\item In this manner, one can match $\Gal(\CC|\R)$-equivariantly $\Sigma_P(G, T)$ and its counterpart $\Sigma_{P'}(G', T')$. One can take a standard basis $(e^*_1)', \ldots, (e^*_n)'$ for $X^*(T'_{\CC})$ in a parallel manner, so that
	\[ T \rightiso T' \quad\text{induces}\quad e^*_i \mapsto (e^*_i)', \quad \forall 1 \leq i \leq n. \]
	
	\item Moreover, the infinitesimal characters of genuine discrete series of $\tilde{M}$ are in bijection with those of discrete series of $M'(\R)$. The matching upgrades to the level of Harish-Chandra parameters for discrete series: see \cite[\S 11]{Ad98}.
\end{itemize}

This recipe is also explained in \cite[Lemma 7.3.3]{Li19}. Another way is to use the ``diagrams'' in \S\ref{sec:diagram} to produce such correspondences of Borel pairs between $G$ and $\SO(2n+1)$.

The recipe of Definition--Proposition \ref{prop:normalizing-real-0} applies to the normalizing factors for $G'$ as well, by setting $m=1$ everywhere. This is done explicitly in \cite{Ar89a}, which we followed. The resulting normalizing factors are then expressed in terms of local $L$-functions, namely
\[ r_{Q'|P'}(\sigma_\lambda) = L(0, \sigma_\lambda, \rho^\vee_{Q'|P'}) L(1, \sigma_\lambda, \rho^\vee_{Q|P})^{-1} \]
as in \cite[(3.2)]{Ar89a}, where $\sigma$ is a discrete series representation of $M'(\R)$ and $\lambda \in \mathfrak{a}_{M',\CC}^* \simeq \mathfrak{a}_{M, \CC}^*$. In particular, $r_{Q'|P'}(\sigma_\lambda)$ depend only on the infinitesimal character of $\sigma_\lambda$.

All in all, $r_{Q|P}(\pi_\lambda)$ can be compared with its $G'$-counterpart $r_{Q'|P'}(\sigma_\lambda)$, where we take $\sigma$ with infinitesimal character matching that of $\pi$.

\begin{theorem}\label{prop:normalizing-real}
	In the circumstance above, setting $\Sigma_{P, \mathrm{long}}(G, T) := \Sigma_P(G,T) \cap \left\{ \text{long roots}\right\}$, we have
	\[ \dfrac{r_{\tilde{Q}|\tilde{P}}(\pi_\lambda)}{r_{Q'|P'}(\sigma_\lambda)} = \sqrt{2}^{|\Sigma_{P, \mathrm{long}}(G, T)|}. \]
\end{theorem}
\begin{proof}
	This is a comparison between Definition--Proposition \ref{prop:normalizing-real-0} and its counterpart for $G' := \SO(2n+1)$. Denote by $\Sigma^c_P(G, T)$ the subset of non-real roots in $\Sigma_P(G, T)$, and $\Sigma^c_{P, \mathrm{long}}(G, T) := \Sigma^c_P(G, T) \cap \Sigma_{P, \mathrm{long}}(G, T)$. Also observe that if $\alpha \in \Sigma(G,T)$ corresponds to $\beta \in \Sigma(G', T')$, then
	\[ \check{\alpha} = \begin{cases}
		\check{\beta}, & \text{if}\; \alpha \;\text{is short} \\
		\frac{1}{2} \cdot \check{\beta}, & \text{if}\; \alpha \;\text{is long}
	\end{cases}\]
	under the identification $T' \simeq T$.
	
	\textbf{Step 1: Non-real roots.}
	We begin with the terms indexed by $\Sigma^c_P(G, T) / \Gal(\CC|\R)$. Since
	\[ \Gamma_{\CC}(z)/\Gamma_{\CC}(z+1) = \frac{2\pi}{z}, \]
	we infer that
	\[ \frac{\Gamma_{\CC}(\lrangle{\mu + \lambda, \check{\alpha}})}{\Gamma_{\CC}(\lrangle{\mu + \lambda, \;\check{\alpha}} + 1)} = \begin{cases}
		\dfrac{\Gamma_{\CC}(\lrangle{\mu + \lambda, \check{\beta}})}{\Gamma_{\CC}(\lrangle{\mu + \lambda, \;\check{\beta}} + 1)}, & \text{if}\; \alpha \;\text{is short} \\
		2 \cdot \dfrac{\Gamma_{\CC}(\lrangle{\mu + \lambda, \check{\beta}})}{\Gamma_{\CC}(\lrangle{\mu + \lambda, \;\check{\beta}} + 1)}, & \text{if}\; \alpha \;\text{is long}.
	\end{cases}\]
	Taking the product over $\Gal(\CC|\R)$-orbits, they contribute $\sqrt{2}^{|\Sigma^c_{P, \mathrm{long}}(G, T)|}$ to the ratio between normalizing factors.
	
	Note that we are in case II if and only if there is no real root in $\Sigma_{P'}(G', T')$, that is, if and only if we are in ``case II'' on the $G'$-side. We are done when either condition holds, in which case $\Sigma^c_{P, \mathrm{long}}(G, T) = \Sigma_{P, \mathrm{long}}(G, T)$.
	
	Assume hereafter that we are in case I. Let $\alpha_0$ be the unique real root in $\Sigma_P(G, T)$. The goal is to reach the ratio $\sqrt{2}^{|\Sigma_{P, \mathrm{long}}(G, T)|}$.
	
	\textbf{Step 2: Long real root.}
	Suppose $\alpha_0$ is long. The homomorphism $\Phi_{\alpha_0}: \widetilde{\SL}(2) \hookrightarrow \tilde{G}$ constructed earlier is then the embedding of metaplectic coverings associated with $W = W_0 \oplus W'_0$ where $W_0$, $W'_0$ are symplectic vector subspaces with $\dim W_0 = 2$. The covering $\Mp(2) = \widetilde{\SL}(2)$ being non-split, it restricts to the non-trivial twofold covering of the maximal compact subgroup $\SO(2, \R)$. Thus
	\[ \gamma^2 \in \bmu_2 = \Ker(\rev), \quad \gamma^2 \neq 1. \]
	As $\chi$ is genuine, $\chi(\gamma)$ is of order exactly $4$, hence so is $(-1)^{\lrangle{\rho_B, \check{\alpha}_0}} \chi(\gamma)$, and we see $\frac{k}{4} \in \{\frac{1}{4}, \frac{3}{4}\}$. Hence we are in the second case in \eqref{eqn:normalizing-G}.

	Denote by $\gamma'$ the $G'$-counterpart of $\gamma$. The $G'$-counterpart of $\mathsf{G}(z)$ must fall into the first case of \eqref{eqn:normalizing-G} since $(\gamma')^2 = 1$, and we claim that $N_0 = 0$ for $G'$. Indeed, $\beta_0$ is short and $\Phi_{\beta_0}: \SL(2) \to G'$ factors through $\SO(3)$, thus $\gamma' = 1$ here. It remains to observe that $\lrangle{\rho_{B'}, \check{\beta}_0} \equiv 1 \pmod{2}$ by \eqref{eqn:rho-comparison}.
	
	Since $\check{\beta}_0 = 2\check{\alpha}_0$, we conclude from \eqref{eqn:normalizing-G} that the $\mathsf{G}(\lrangle{\mu + \lambda, \check{\alpha}_0})$ for $\tilde{G}$ equals $\sqrt{2}$ times its $G'$-counterpart.

	\textbf{Step 3: Short real root.}
	Next, suppose $\alpha_0$ is short. The only possibility is $\alpha_0 = e^*_1 + e^*_2$ and $a=2$. By \eqref{eqn:rho-comparison},
	\begin{equation}\label{eqn:normalizing-real-aux0}
		\lrangle{\rho_B, \check{\alpha}_0} \equiv \lrangle{\rho_{B'}, \check{\beta}_0} + 1 \pmod{2} .
	\end{equation}

	The Galois-invariant part of $X_*(T_{\CC})$ is generated by $e_{1,*} + e_{2,*}$. Together with $\alpha_0 = e_1^* + e_2^*$ and $\check{\alpha}_0$, this determines the image of the map $\SL(2) \to G$, the uncovered version of $\Phi_{\alpha_0}$, by describing its root datum; in particular, this is an embedding. We refer to \eqref{eqn:SL2-Sp4-embedding} for a down-to-earth description. In this way we construct $\chi: \tilde{T} \to \CC^\times$ and $\gamma \in \tilde{T}$. In fact, $\gamma$ lies over $-1 \in \CC^\times \simeq T_{\GL}(\R)$.

	The recipe above makes no use of long roots, hence it carries over to $T' \subset G'$ verbatim. In particular, we have $\SL(2) \hookrightarrow G'$ and $\chi'$, $\gamma'$. Identify both $T_{\GL}(\R)$ and its $G'$-counterpart $T'_{\GL}(\R)$ with $\CC^\times$. As before, $\gamma' = -1 \in \CC^\times = T'_{\GL}(\R)$. Our goal is to show that $\mathsf{G}(\lrangle{\mu + \lambda, \check{\alpha}_0})$ equals its $G'$-counterpart. In view of \eqref{eqn:normalizing-real-aux0}, this is tantamount to
	\begin{equation}\label{eqn:normalizing-real-aux1}\begin{gathered}
		\gamma^2 = 1 \quad \text{(so we are in the first case of \eqref{eqn:normalizing-G})}, \\
		\text{and} \quad \chi(\gamma) = -\chi'(\gamma').
	\end{gathered}\end{equation}

	The correspondence of infinitesimal characters (and Harish-Chandra parameters) between $\tilde{G}$ and $G'$ is compatible with $T \simeq T'$. Moreover, $\rho^M$ and $\rho^{M'}$ have the same projection to $X^*(T_{\GL, \CC}) \otimes \Q$ and $X^*(T'_{\GL, \CC}) \otimes \Q$, namely $\frac{1}{2}(e^*_1 - e^*_2)$. It follows that
	\[ \dd\chi|_{\mathfrak{t}_{\GL, \CC}} = \dd\chi'|_{\mathfrak{t}'_{\GL, \CC}} . \]
	
	For any subgroup $S \subset G$, denote the twofold (resp.\ eightfold) covering as $\tilde{S} \subset \tilde{G}$ (resp.\ $\tilde{S}^{(8)} \subset \tilde{G}^{(8)}$), then
	\index{SpWtilde-8@$\Mp(W)^{(8)}$}
	\begin{align}
		\label{eqn:normalizing-real-e1}
		\tilde{M}^{(8)} & = \GL(2, \R) \times \Mp(W^\flat)^{(8)} \subset \tilde{G}^{(8)}, \\
		\label{eqn:normalizing-real-e2}
		\tilde{T}^{(8)} & = T_{\GL}(\R) \times \tilde{T}^{\flat, (8)}.
	\end{align}
	Extend $\chi$ to a genuine character $\chi^{(8)}$ of $\tilde{T}^{(8)}$. From \eqref{eqn:normalizing-real-e2} we see $\chi^{(8)}|_{\CC^\times} = \chi'|_{\CC^\times}$.
	
	The image of $\GL(2, \R)$ under \eqref{eqn:normalizing-real-e1} lies in a copy of $\Mp(4)^{(8)}$ inside $\tilde{G}^{(8)}$. Consider the image of $-1 \in \GL(2, \R)$; it equals the canonical preimage of $-1 \in \Sp(4, \R)$ in \cite[Définition 2.8]{Li11}. The discussions above lead to $\chi'(\gamma') = \chi^{(8)}(-1)$. On the other hand, $\gamma$ is also a preimage of $-1$ in $\Mp(4)^{(8)}$. Denote by $z$ the nontrivial element in $\bmu_2$. As $\chi(\gamma) = \chi^{(8)}(\gamma)$ and $\chi^{(8)}$ is genuine, the two equalities in \eqref{eqn:normalizing-real-aux1} reduce to the group-theoretic assertion that
	\begin{equation}\label{eqn:normalizing-real-aux2}
		\gamma = z \cdot (-1) \quad \text{inside}\quad \tilde{T}_{\GL}^{(8)} \subset \Mp(4)^{(8)}.
	\end{equation}
	
	Note that \eqref{eqn:normalizing-real-aux2} is really a statement about $\Mp(4)$. Indeed, $\alpha_0 = e_1^* + e_2^*$ is a short real root for $T_{\GL} \subset \Sp(4)$, hence $\Phi_{\alpha_0}$ and $\gamma$ are constructed using an $\SL(2)$-triple inside $\Sp(4)$. The verification of \eqref{eqn:normalizing-real-aux2} is deferred to Proposition \ref{prop:normalizing-Sp4}.
\end{proof}

We remark that $|\Sigma_{P, \mathrm{long}}(G, T)|$ can be computed using any maximal torus $T$ of $M$: there is no need to assume $T/A_M$ anisotropic, and one can compute it over $\CC$.

\subsection{The complex case}
When $F = \CC$, we have $\tilde{G} \simeq G(\CC) \times \bmu_8$ canonically. Genuine representations of $\tilde{G}$ are nothing but representations of $G(\CC)$. Hence the normalization reduces to the uncovered case of complex groups. The recipe in the real case still applies by regarding $G(\CC)$ as a real group.

One can still take the principal endoscopic group $G' = \SO(2n+1)$ of $\tilde{G}$, and compare $r_{\tilde{Q}|\tilde{P}}(\pi_\lambda)$ with its $G'$-counterpart $r_{Q'|P'}(\sigma_\lambda)$, as we did in the discussions preceding Theorem \ref{prop:normalizing-real}. The situation here simplifies in two aspects:
\begin{enumerate}
	\item there are no real roots anymore,
	\item $M=T$ since $M$ is assumed to have discrete series.
\end{enumerate}

In the statement below, the sets $\Sigma(G, T)$ and $\Sigma_P(G, T)$ are taken in the complex sense: no Weil restrictions $\CC|\R$ are performed. Since the argument is not new, we shall make it brief.

\begin{theorem}\label{prop:normalizing-cplx}
	In the circumstance above, we have
	\[ \dfrac{r_{\tilde{Q}|\tilde{P}}(\pi_\lambda)}{r_{Q'|P'}(\sigma_\lambda)} = 2^{|\Sigma_{P, \mathrm{long}}(G, T)|}, \]
	here $\Sigma_{P, \mathrm{long}}(G, T)$ is the subset of $\Sigma_P(G, T)$ consisting of long roots.
\end{theorem}
\begin{proof}
	This is identical to the step 1 of the proof of Theorem \ref{prop:normalizing-real}. By viewing $G$, $M$, etc.\ as real groups in those arguments, each $\alpha \in \Sigma_{P, \mathrm{long}}(G, T)$ is counted twice, whence the factor $2$ instead of $\sqrt{2}$.
\end{proof}

\subsection{Computations in \texorpdfstring{$\Sp(4)$}{Sp(4)}}
The goal is to complete the proof of Theorem \ref{prop:normalizing-real}. We keep the notations there.

Consider a symplectic $\R$-vector space of dimension $4$, equipped with a symplectic basis $u_{\pm 1}, u_{\pm 2}$. This defines the group $\Sp(4)$ as well as its split maximal torus $T_s$. Denote by $\epsilon^*_i$ (resp.\ $\epsilon_{i,*}$) the corresponding elements in $X^*(T_s)$ (resp.\ $X_*(T_s)$) for $i = 1,2$.

Similarly, we have the group $\Sp(2)$ associated with a symplectic vector space $\R u_1 \oplus \R u_{-1}$, together with $\epsilon^*_1$ and $\epsilon_{1, *}$ as before.

Let us realize $\Sp(4)$ by matrices as
\begin{align*}
	\Sp(4) & = \left\{ g \in \GL(4): {}^t g J g = J \right\}, \\
	J & := \begin{pmatrix} & 1_{2 \times 2} \\ -1_{2 \times 2} & \end{pmatrix},
\end{align*}
that is, with respect to the ordered basis $u_1, u_2, u_{-1}, u_{-2}$. Using the polarization given by $\R u_1 \oplus \R u_2$ and $\R u_{-1} \oplus \R u_{-2}$, we have the embedding
\[ \GL(2) \hookrightarrow \Sp(4), \quad g \mapsto \begin{pmatrix} g & \\ & {}^t g^{-1} \end{pmatrix}. \]

View $\mathbb{G}_{m, \CC}$ as an $\R$-torus by Weil restriction. It embeds into $\GL(2)$: on the level of $\R$-points, this is simply
\[ z = x+iy \mapsto \begin{pmatrix} x & -y \\ y & x \end{pmatrix} \in \GL(2, \R). \]
Accordingly, $\mathbb{G}_{m, \CC}$ also embeds into $\Sp(4)$ as a maximal torus. Denote its image by $T$.

Such a non-split maximal torus has appeared in the proof of Theorem \ref{prop:normalizing-real}. The real root $\alpha_0 = e_1^* + e_2^* \in X^*(T_{\CC})$ is $z \mapsto z\bar{z} = |z|^2$ on the level of $\R$-points, and $\check{\alpha}_0 = e_{1, *} + e_{2, *}$ is simply
\begin{equation}\label{eqn:TGL-Sp4-embedding}
	\Gm \to T \subset \Sp(4): \quad t \mapsto \begin{pmatrix} t \cdot 1_{2 \times 2} & \\ & t^{-1} \cdot 1_{2 \times 2} \end{pmatrix}.
\end{equation}
In particular, the embedding above given by $\check{\alpha}_0 \in X^*(T_{\CC})$ coincides with that given by the coroot $\epsilon_{1,*} + \epsilon_{2, *} \in X^*(T_s)$.

It is readily seen that the $\SL(2)$ triple $(H', X', Y')$ associated with $T \subset \Sp(4)$ and $\alpha_0$ can be taken as
\begin{equation}\label{eqn:SL2-triple-Sp4}\begin{aligned}
	X' = \begin{pmatrix}
		& 1_{2 \times 2} \\
		&
	\end{pmatrix}, \quad
	Y' = \begin{pmatrix}
		& \\
		1_{2 \times 2} &
	\end{pmatrix}, \\
	H' = [X', Y'] = \begin{pmatrix}
		1_{2 \times 2} & \\
		& -1_{2 \times 2}
	\end{pmatrix}.
\end{aligned}\end{equation}
The corresponding homomorphism $\SL(2) \to \Sp(4)$ is therefore
\begin{equation}\label{eqn:SL2-Sp4-embedding}
	\begin{pmatrix} a & b \\ c & d \end{pmatrix} \mapsto
	\begin{pmatrix}
		a & & b & \\
		& a & & b \\
		c & & d & \\
		& c & & d
	\end{pmatrix}.
\end{equation}
Note that \eqref{eqn:SL2-Sp4-embedding} coincides with the embedding of $\SL(2)$ associated with the short root $\epsilon_1^* + \epsilon_2^* \in X^*(T_s)$. We will forget about $T$ and focus on the split maximal torus $T_s$ in what follows.

Now consider the \emph{twofold} (resp.\ \emph{eightfold}) metaplectic covering $\Mp(4)$ (resp.\ $\Mp(4)^{(8)}$). All are defined relative to the same additive character $\psi$ of $\R$.

\begin{lemma}\label{prop:SL2-Sp4-split}
	Let $\widetilde{\SL}(2) \subset \Mp(4)$ be the preimage of $\SL(2, \R)$ under \eqref{eqn:SL2-Sp4-embedding}. Then the covering $\widetilde{\SL}(2) \to \SL(2, \R)$ splits.
\end{lemma}
\begin{proof}
	This can be seen from the Brylinski--Deligne theory of coverings\footnote{See \cite{GG18} for an introduction.}. Instead, we reason directly as follows. Since $\SL(2, \R)$ equals its derived subgroup, it suffices to show that the eightfold version $\widetilde{\SL}(2)^{(8)}$ splits.
	
	Since \eqref{eqn:SL2-Sp4-embedding} is associated with $T_s \subset \Sp(4)$ and the short root $\epsilon^*_1 + \epsilon^*_2$, its image is contained in a $\GL(2)$ associated with some polarization (not the one given by $u_1, u_2$ and $u_{-1}, u_{-2}$). The eightfold covering splits over this copy of $\GL(2, \R)$. Hence $\widetilde{\SL}(2)^{(8)}$ splits as well. 
\end{proof}

Note that the inclusion $\widetilde{\SL}(2) \hookrightarrow \Mp(4)$ is nothing but the map $\Phi_{\alpha_0}$ in the proof of Theorem \ref{prop:normalizing-real}.

The following construction is standard in the general theory of covering groups. We consider only the case over $\R$.

\begin{notation}
	Suppose we are given a connected reductive $\R$-group $G$ with split maximal torus $T_s$. To each root $\alpha \in X^*(T_s)$, suppose that an isomorphism $x_\alpha: \Ga \rightiso U_\alpha$ is chosen, where $U_\alpha \subset G$ is the root subgroup associated with $\alpha$, such that $x_\alpha$ and $x_{-\alpha}$ are opposite in the following sense: there is a surjection $\zeta_{\pm\alpha}$ from $\SL(2)$ onto the derived subgroup of $\lrangle{T_s, U_\alpha, U_{-\alpha}}$, and we require that
	\[ x_\alpha(t) = \zeta_{\pm\alpha} \begin{pmatrix} 1 & t \\ & 1 \end{pmatrix}, \quad
	x_{-\alpha}(t) = \zeta_{\pm\alpha} \begin{pmatrix} 1 & \\ t & 1 \end{pmatrix}. \]

	Consider any covering $\tilde{G} \to G(\R)$ and any root $\alpha \in X^*(T_s)$. Then $x_\alpha$ lifts canonically to a homomorphism of Lie groups $\R \to \tilde{G}$, which we still denote as $x_\alpha$. Define the elements
	\begin{align*}
		w_\alpha(t) & := x_\alpha(t) x_{-\alpha}(-t^{-1}) x_\alpha(t), \\
		h_\alpha(t) & := w_\alpha(t) w_{-\alpha}(-1)
	\end{align*}
	in $\tilde{G}$, for all $t \in \R^\times$. We have $h_\alpha(t) \mapsto \check{\alpha}(t) \in T_s(\R)$.
\end{notation}

Let us apply this formalism to $\Mp(2)$ and $\Mp(4)$.
\begin{itemize}
	\item For $\Mp(2)$ and $\Mp(4)$, the root $\alpha$ can be $2\epsilon^*_i$ and $x_{\pm\alpha}$ is the standard one. We take $i=1$ for $\Mp(2)$.
	\item For $\Mp(4)$, we also allow $\alpha = \epsilon_1^* + \epsilon_2^*$, in which case $x_{\pm\alpha}$ is the one determined by \eqref{eqn:SL2-triple-Sp4} after taking $\dd x_{\pm\alpha}|_{t=0}$.
\end{itemize}

\begin{lemma}\label{prop:Sp4-bisector}
	We have $h_{\epsilon^*_1 + \epsilon^*_2}(t) = h_{2\epsilon^*_1}(t) h_{2\epsilon^*_2}(t)$ for all $t \in \R^\times$.
\end{lemma}
\begin{proof}
	We have
	\[ (\epsilon_1^* + \epsilon^*_2)^\vee = \epsilon_{1, *} + \epsilon_{2, *} = (2\epsilon^*_1)^\vee + (2\epsilon^*_2)^\vee . \]
	
	All coverings in view are Brylinski--Deligne coverings. Consider temporarily the case of $\Mp(2n)$ for general $n \geq 1$, endowed with the chosen split maximal torus $T_s \subset \Sp(2n)$. In order to perform computations, we follow \cite[\S 2.6]{GG18} to ``incarnate'' the covering by choosing
	\begin{itemize}
		\item a ``bisector'' $D \in X^*(T_s) \otimes X^*(T_s)$, also viewed as an integral bilinear form on $X_*(T_s)$,
		\item an arbitrary group homomorphism $\eta: X_*(T_s) \to \R^\times$.
	\end{itemize}
	Specifically, one can take $\eta = \mathbf{1}$ and
	\[ D := \epsilon_1^* \otimes \epsilon_1^* + \cdots + \epsilon_n^* \otimes \epsilon_n^* , \]
	so that $D(y, y) = Q(y) := \sum_{i=1}^n y_i^2$ for all $y = \sum_{i=1}^n y_i \epsilon_{i, *} \in X_*(T_s)$, where $Q: X_*(T_s) \to \Z$ is the quadratic form parameterizing the twofold metaplectic covering.
	
	Using these data with $n = 2$, one may identify $\Mp(4)$ with $\Sp(4, \R) \times \bmu_2$ as sets, with $\bmu_2$ embedded as $\{1\} \times \bmu_2$, and make calculations with explicit cocycles over $T_s(\R)$.
	
	To be more precise, denote by $(\cdot, \cdot)_{\R, 2}$ the quadratic Hilbert symbol over $\R$. The discussion preceding \cite[Proposition 2.5]{GG18} yields
	\begin{align*}
		(\epsilon_{1, *}(t), 1) \cdot (\epsilon_{2,*}(t), 1) & = \left( \epsilon_{1, *}(t) \epsilon_{2,*}(t), \left(t, 1 \right)_{\R, 2} \left(1, t \right)_{\R, 2} \right) \\
		& = \left( \epsilon_{1, *}(t) \epsilon_{2,*}(t), 1 \right),
	\end{align*}
	where we used
	\[ \epsilon^*_1 \epsilon_{1, *}(t) = t = \epsilon^*_2 \epsilon_{2, *}(t), \quad \epsilon^*_1 \epsilon_{2, *}(t) = 1 = \epsilon_2^* \epsilon_{1, *}(t). \]
	
	From \cite[Proposition 2.5]{GG18} and $\eta = 1$, we infer that
	\begin{align*}
		h_{\epsilon^*_1 + \epsilon^*_2}(t) & = \left( (\epsilon_{1,*} + \epsilon_{2,*})(t), \; (\eta(\epsilon_{1,*} + \epsilon_{2,*}), t )_{\R, 2} \right) \\
		& = \prod_{i=1}^2 \left( \epsilon_{i,*}(t), 1 \right).
	\end{align*}
	But the final expression is also $\prod_{i=1}^2 h_{2\epsilon^*_i}(t)$ by the same result applied to $\Mp(2)$.
\end{proof}

\begin{lemma}\label{prop:Sp4-gamma}
	Define $\gamma := \exp\left(\pi(X' - Y')\right) \in \Mp(4)$ using the data in \eqref{eqn:SL2-triple-Sp4}. Then $\gamma$ equals $h_{\epsilon^*_1 + \epsilon^*_2}(-1)$.
\end{lemma}
\begin{proof}
	Their images in $\Sp(4, \R)$ are both $-1$. Consider the sub-covering $\widetilde{\SL}(2) \subset \Mp(4)$ in Lemma \ref{prop:SL2-Sp4-split}. There is an isomorphism between coverings $\iota: \SL(2, \R) \times \bmu_2 \rightiso \widetilde{\SL}(2)$, and $\exp\left(\theta(X' - Y')\right) \in \iota(\SL(2, \R) \times \{1\})$ for all $\theta \in \R$ by continuity. On the other hand, $x_{\pm(\epsilon^*_1 + \epsilon^*_2)}(\R)$ also lands in $\iota(\SL(2, \R) \times \{1\})$. Hence $h_{\epsilon^*_1 + \epsilon^*_2}(-1) \in \iota(\SL(2, \R) \times \{1\})$. The equality follows at once.
\end{proof}

The next result concludes the proof of Theorem \ref{prop:normalizing-real} by establishing \eqref{eqn:normalizing-real-aux2}.

\begin{proposition}\label{prop:normalizing-Sp4}
	Let $\gamma \in \Mp(4)$ be as in Lemma \ref{prop:Sp4-gamma}. Let $-1$ denote the canonical preimage of $-1 \in \Sp(4, \R)$ in $\Mp(4)^{(8)}$. Then $\gamma = z \cdot (-1)$ where $z$ denotes the nontrivial element in $\bmu_2$.
\end{proposition}
\begin{proof}
	There is a nature homomorphism $j: \Mp(2)^{(8)} \times \Mp(2)^{(8)} \to \Mp(4)^{(8)}$, where the first (resp.\ second) $\Mp(2)$ is constructed from the symplectic vector subspace spanned by $u_1, u_{-1}$ (resp.\ $u_2, u_{-2}$). It covers the evident embedding $\Sp(2) \times \Sp(2) \hookrightarrow \Sp(4)$. There is a canonical $-1$ in $\Mp(2)^{(8)}$ lying over $-1 \in \Sp(2, \R)$. We claim that
	\[ \left( \text{the}\; -1 \in \Mp(4)^{(8)} \right) = j(-1, -1). \]
	Indeed, this follows from the characterization \cite[Corollaire 4.6]{Li11} of the canonical $-1$ in terms of the character $\Theta_\psi$ of the Weil representation, together with the decomposition formula \cite[Proposition 4.25]{Li11} for $\Theta_\psi$.

	In parallel, note that the $h_{2\epsilon^*_i}(t)$ in Lemma \ref{prop:Sp4-bisector} can be constructed in the $i$-th copy of $\Mp(2)$ inside $\Mp(4)$ via $j$. Putting $t = -1$ yields
	\[ h_{\epsilon^*_1 + \epsilon^*_2}(-1) = j(h_{2\epsilon^*_1}(-1), h_{2\epsilon^*_2}(-1)). \]
	We are reduced to proving the following assertion: there exists $w \in \bmu_8$ of order $4$ such that for all $i \in \{1, 2\}$, by embedding $\Mp(2)^{(8)}$ into $\Mp(4)^{(8)}$ as the $i$-th copy under $j$, we have
	\[ h_{2\epsilon_i^*}(-1) = w \cdot \left( \text{the}\; -1 \in \Mp(2)^{(8)} \right). \]
	
	This is probably well-known to experts. Let us extract it from the proof of \cite[Proposition 9.1.1]{Li20}.
	\begin{itemize}
		\item In the notation of \textit{loc.\ cit.}, the $-1 \in \Mp(2)^{(8)}$ is denoted as $\sigma_\ell[-1]$: this is the lifting via the Schrödinger model associated with the Lagrangian subspace $\ell := \R u_i$.
		\item In \textit{loc.\ cit.}, one also considers the Lion--Perrin lifting of $-1$ in $\Mp(2)$, denoted by $\sigma_{\mathrm{LP}}[-1]$.
	\end{itemize}
	As elements of $\Mp(2)^{(8)}$, they are shown to satisfy
	\begin{align*}
		\sigma_\ell[-1] \cdot \gamma_\psi(1)^{-2} & = \sigma_\ell[-1] \cdot \frac{\gamma_\psi(-1)}{\gamma_\psi(1)} \\
		& = \sigma_{\mathrm{LP}}[-1] = (-1, -1)_{\R, 2} \cdot h_{2\epsilon^*_i}(-1).
	\end{align*}
	Here we used the property $\gamma_\psi(-a) \gamma_\psi(a) = 1$. Take $w := -\gamma_\psi(1)^{-2}$. It remains to notice that $\gamma_\psi(1) \in \bmu_8$ is of order $8$; for example, just take $a = b = -1$ in $(a, b)_{\R, 2} = \frac{\gamma_\psi(ab) \gamma_\psi(1)}{\gamma_\psi(a) \gamma_\psi(b)}$.
\end{proof}

\section{Review of unitary weighted characters}\label{sec:weighted-characters-arch}
Let $M \in \mathcal{L}(M_0)$. For every $\pi \in \Pi_-(\tilde{M})$ and $f \in C^\infty_{c, \asp}(\tilde{G}, \tilde{K}) \otimes \mes(G)$ (see Definition \ref{def:K-finite}), we have the \emph{normalized weighted character}
\[ J^r_{\tilde{M}}(\pi_\lambda, f) = J^{\tilde{G}, r}_{\tilde{M}}(\pi_\lambda, f) := \Tr\left( \mathcal{R}_{\tilde{M}}(\pi_\lambda, \tilde{P}) I^{\tilde{G}}_{\tilde{P}}(\pi_\lambda, f) \right), \quad P \in \mathcal{P}(M). \]
See \cite[\S 5.7]{Li12b}, which is in turn based on \cite[\S 6]{Ar89a}. Explications:
\begin{itemize}
	\item $\mathcal{R}_{\tilde{M}}(\pi_\lambda, \tilde{P})$ arises from the $(G,M)$-family
	\[ R_{\tilde{Q}|\tilde{P}}(\pi_\lambda)^{-1} R_{\tilde{Q}|\tilde{P}}(\pi_{\lambda + \Lambda}), \quad Q \in \mathcal{P}(M), \; \Lambda \in i\mathfrak{a}_M^* \]
	made from the normalized intertwining operators $R_{\tilde{Q}|\tilde{P}}(\cdots)$;
	\item $J^r_{\tilde{M}}(\pi_\lambda, f)$ is a meromorphic family in $\lambda \in \mathfrak{a}_{M, \CC}^*$. It is independent of $P \in \mathcal{P}(M)$ and analytic at $\lambda$ if $\pi_\lambda$ is unitary.
\end{itemize}

\index{JrMpiX@$J^r_{\tilde{M}}(\pi, f)$, $J^r_{\tilde{M}}(\pi, X, f)$}
For all $\pi \in \Pi_{\mathrm{unit}, -}(\tilde{M})$ and $X \in \mathfrak{a}_M$, we put
\[ J^r_{\tilde{M}}(\pi, X, f) := \int_{i\mathfrak{a}_M^*} J^r_{\tilde{M}}(\pi_\lambda, f) e^{-\lrangle{\lambda, X}} \dd\lambda . \]
\begin{itemize}
	\item It is well-defined since $\lambda \mapsto J^r_{\tilde{M}}(\pi_\lambda, f)$ is Schwartz. This follows from a straightforward generalization of the estimate \cite[Lemma 2.1]{Ar94}, since our normalization factors have the same properties as in the linear case.
	
	Consequently, $X \mapsto J^r_{\tilde{M}}(\pi, X, f)$ is a Schwartz function.
	\item The distribution $J^r_{\tilde{M}}(\pi, X, \cdot)$ is concentrated at the projection of $X$ to $\mathfrak{a}_G$ (Definition \ref{def:concentrated}). Therefore, it extends to $C^\infty_{\mathrm{ac}, \asp}(\tilde{G}, \tilde{K})$.
\end{itemize}

On the other hand, the \emph{canonically normalized weighted characters} are defined in \cite[\S 5.7]{Li12b} and \cite[Définition 4.1]{Li14b}. They are given by
\[ J_{\tilde{M}}(\pi_\lambda, f) = J^{\tilde{G}}_{\tilde{M}}\left( \pi_\lambda, f \right) = \Tr\left( \mathcal{M}_{\tilde{M}}(\pi_\lambda, \tilde{P}) I^{\tilde{G}}_{\tilde{P}}(\pi_\lambda, f) \right), \quad P \in \mathcal{P}(M) \]
\index{JMpiX}
and the corresponding $J_{\tilde{M}}(\pi_\lambda, X, f)$ for all $X \in \mathfrak{a}_M$, where
\begin{itemize}
	\item as in the non-Archimedean case, $\mathcal{M}_{\tilde{M}}(\pi_\lambda, \tilde{P})$ arises from the $(G,M)$-family described in \eqref{eqn:M-GM-family};
	\item $J_{\tilde{M}}(\pi_\lambda, f)$ is a meromorphic family in $\lambda \in \mathfrak{a}_{M, \CC}^*$. It is independent of $P \in \mathcal{P}(M)$ and is analytic at $\lambda$ if $\pi_\lambda$ is unitary; moreover, $\lambda \mapsto J_{\tilde{M}}(\pi_\lambda, f)$ is rapidly decreasing over $i\mathfrak{a}_M^*$ by \cite[Théorème 4.2]{Li14b}. Therefore it makes sense to define $J_{\tilde{M}}(\pi, X, f)$ by Fourier transform, as has been done for $J^r_{\tilde{M}}(\pi, X, f)$; see \cite[Définition 4.4]{Li14b}.
\end{itemize}
The non-unitary situation will be reviewed in \S\ref{sec:weighted-characters}.

To simplify matters, we assume $\pi \in \Pi_{\mathrm{unit}, -}(\tilde{M})$. By \cite[\S 5.7]{Li12b},
\begin{equation}\label{eqn:Jmu-Jr}
	J_{\tilde{M}}(\pi, f) = \sum_{L \in \mathcal{L}(M)} r^{\tilde{L}}_{\tilde{M}}(\pi) J^r_{\tilde{L}}\left( \pi^{\tilde{L}}, f \right)
\end{equation}
where $\pi^{\tilde{L}}$ is the normalized parabolic induction, and $r^{\tilde{G}}_{\tilde{M}}(\pi)$ arises from the $(G,M)$-family
\[ r_{\tilde{Q}}(\Lambda, \pi) := r_{\tilde{Q}^- | \tilde{Q}}(\pi)^{-1} r_{\tilde{Q}^- | \tilde{Q}}(\pi_{\Lambda / 2}), \quad Q \in \mathcal{P}(M), \; \Lambda \in \mathfrak{a}_{M, \CC}^* . \]
The descent of $(G, M)$-families yields $r^{\tilde{R}}_{\tilde{M}}(\pi)$ for all $L \in \mathcal{L}(M)$ and $R \in \mathcal{P}(L)$; see \cite[p.97]{ArIntro}. It turns out that $r^{\tilde{L}}_{\tilde{M}}(\pi) := r^{\tilde{R}}_{\tilde{M}}(\pi)$ depends only on $L$, whence the notation.
\index{rGM@$r^{\tilde{G}}_{\tilde{M}}(\pi)$}

\begin{definition}\label{def:Schwartz-arch}
	\index{Schwartz function}
	We say that $f \in \orbI_{\mathrm{ac}, \asp}(\tilde{M}) \otimes \mes(M)$ is a \emph{Schwartz function} if $X \mapsto I^{\tilde{M}}(\pi, X, f)$ is a Schwartz function on $\mathfrak{a}_M$ for each $\pi \in D_{\mathrm{temp},-}(\tilde{M}) \otimes \mes(M)^\vee$.
\end{definition}

\begin{definition}\label{def:rphi-arch}
	\index{phiM}
	\index{phiM-r@${}^r \phi_{\tilde{M}}$}
	There are continuous linear maps
	\[\begin{tikzcd}[column sep=large]
		C^\infty_{\mathrm{ac}, \asp}(\tilde{G}, \tilde{K}) \otimes \mes(G) \otimes \mes(G) \arrow[r, shift left, "{{}^r \phi_{\tilde{M}} = {}^r \phi^{\tilde{G}}_{\tilde{M}}}"] \arrow[r, shift right, "{\phi_{\tilde{M}} = \phi^{\tilde{G}}_{\tilde{M}}}"'] & \orbI_{\mathrm{ac}, \asp}(\tilde{M}, \tilde{K} \cap \tilde{M}) \otimes \mes(M)
	\end{tikzcd}\]
	characterized by
	\[ I^{\tilde{M}}(\pi, X, {}^r \phi_{\tilde{M}}(\cdot)) = J^r_{\tilde{M}}(\pi, X, \cdot), \quad I^{\tilde{M}}(\pi, X, \phi_{\tilde{M}}(\cdot)) = J_{\tilde{M}}(\pi, X, \cdot) \]
	for all $\pi \in D_{\mathrm{temp}, -}(\tilde{M}) \otimes \mes(M)^\vee$ and $X \in \mathfrak{a}_M$. Their images consist of Schwartz functions defined above.
\end{definition}

\begin{remark}[Extension to Schwartz--Harish-Chandra spaces]\label{rem:SHC}
	\index{CG@$\mathscr{C}_{\asp}(\tilde{G})$}
	It is shown in \cite[Théorème 5.7.2]{Li12b} that $\phi_{\tilde{M}}$ and ${}^r \phi_{\tilde{M}}$ both extend by continuity to linear maps
	\[ \mathscr{C}_{\asp}(\tilde{G}) \otimes \mes(G) \to I\mathscr{C}_{\asp}(\tilde{M}) \otimes \mes(M), \]
	where $\mathscr{C}_{\asp}$ stands for the anti-genuine Schwartz--Harish-Chandra space, and $I\mathscr{C}_{\asp}$ is the corresponding avatar of $\orbI_{\asp}$.
\end{remark}

\begin{proposition}\label{prop:phimu-phir}
	For all $\pi \in \Pi_{\mathrm{temp},-}(\tilde{M})$ and $f \in \mathscr{C}_{\asp}(\tilde{G}) \otimes \mes(G)$, we have
	\[ I^{\tilde{M}}\left(\pi, \phi_{\tilde{M}}(f)\right) = \sum_{L \in \mathcal{L}(M)} r^{\tilde{L}}_{\tilde{M}}(\pi) I^{\tilde{L}}\left(\pi^{\tilde{L}}, {}^r \phi_{\tilde{L}}(f)\right). \]	
\end{proposition}
\begin{proof}
	Immediate from \eqref{eqn:Jmu-Jr}.
\end{proof}

When $\tilde{G} = \Mp(W)$ with $\dim W = 2n$ and $\pi \in \Pi_{\mathrm{temp},-}(\tilde{M})$, it makes sense to compare $r_{\tilde{Q}|\tilde{P}}(\pi)$ with its $\SO(2n+1)$-counterpart $r_{Q'|P'}(\sigma)$ as in Theorems \ref{prop:normalizing-real} and \ref{prop:normalizing-cplx}. Hence one can also compare $r^{\tilde{L}}_{\tilde{M}}(\pi_\lambda)$ with $r^{L'}_{M'}(\sigma_\lambda)$.

\begin{theorem}
	Suppose $\tilde{G} = \Mp(W)$ with $\dim W = 2n$. For all $\pi \in \Pi_{\mathrm{temp},-}(\tilde{M})$ and $L \in \mathcal{L}(M)$, we have
	\[ r^{\tilde{L}}_{\tilde{M}}(\pi) = \text{its $\SO(2n+1)$-counterpart}\; r^{L^!}_{M^!}(\sigma). \]
	The same comparison generalizes to groups of metaplectic type.
\end{theorem}
\begin{proof}
	Theorems \ref{prop:normalizing-real} and \ref{prop:normalizing-cplx} imply that the $(G, M)$-family $\left(r_{\tilde{Q}}(\Lambda, \pi)\right)_{Q \in \mathcal{P}(M)}$ matches its $\SO(2n+1)$-counterpart. Ditto for the $(L, M)$-families $\left(r^{\tilde{R}}_{\tilde{S}}(\Lambda, \pi)\right)_{S \in \mathcal{P}^L(M)}$ so obtained by descent, where $L \in \mathcal{L}(M)$ and $R \in \mathcal{P}(L)$.
	
	By the general formalism of $(L, M)$-families, see \eqref{eqn:GMtheta} and \eqref{eqn:GMlim}, the $r^{\tilde{R}}_{\tilde{M}}(\pi)$ is then
	\[ \lim_{\substack{\Lambda \to 0 \\ \Lambda \in i\mathfrak{a}^*_M}} \sum_{S \in \mathcal{P}^L(M)} r^{\tilde{R}}_{\tilde{S}}(\Lambda, \pi) \theta^L_S(\Lambda)^{-1} . \]
	The same construction applies to $\SO(2n+1)$-counterparts. Since $\mathfrak{a}^{L'}_{S'} \simeq \mathfrak{a}^L_S$ preserves the Haar measures determined by the invariant quadratic forms chosen in Definition \ref{def:invariant-quadratic-form}, it is readily seen that the effects of rescaling coroots cancel with each other in \eqref{eqn:GMtheta}, and $\theta^L_S = \theta^{L'}_{S'}$.
	
	All in all, $r^{\tilde{L}}_{\tilde{M}}(\pi)$ equals its $\SO(2n+1)$-counterparts. The case of groups of metaplectic types follows by extending the matching tautologically on the $\GL$-factors.
\end{proof}

Finally, define the $\CC$-vector space
\begin{equation}\label{eqn:UGM-spec-arch}
	U^G_M := \left\{\begin{array}{l}
	\text{meromorphic functions}\; u: \mathfrak{a}^{G,*}_{M, \CC} \to \CC, \\
	\text{of moderate growth on}\; i\mathfrak{a}^{G,*}_M
\end{array}\right\}. \end{equation}
Moderate growth means that there exists $N \geq 0$ such that $\int_{i\mathfrak{a}^{G,*}_M} |u(\lambda)| (1 + \|\lambda\|)^{-N} \dd\lambda$ is finite, where $\|\cdot\|$ is any Euclidean norm.
\index{UGM}

\begin{theorem}\label{prop:rGM-growth}
	The function $\lambda \mapsto r^{\tilde{G}}_{\tilde{M}}(\pi_\lambda)$ belongs to $U^G_M$ for every $\pi \in \Pi_{\mathrm{unit},-}(\tilde{M})$.
\end{theorem}
\begin{proof}
	We know that $r^{\tilde{G}}_{\tilde{M}}(\pi_\lambda)$ is analytic whenever $\lambda \in i\mathfrak{a}^*_M$. The moderate growth is shown in the last step of the proof of \cite[Lemme 5.17]{Li14b}.
\end{proof}

Note that $r^{\tilde{G}}_{\tilde{M}}(\pi) = r^{\tilde{G}}_{\tilde{M}}(\pi_\lambda)$ whenever $\lambda \in \mathfrak{a}^*_{G, \CC}$.

\section{The maps \texorpdfstring{${}^c \theta^{\tilde{G}}_{\tilde{M}}$}{cphi}, \texorpdfstring{${}^r \theta^{\tilde{G}}_{\tilde{M}}$}{rtheta}, \texorpdfstring{${}^{c,r} \theta^{\tilde{G}}_{\tilde{M}}$}{crtheta}, and the distributions \texorpdfstring{${}^c I^{\tilde{G}}_{\tilde{M}}(\tilde{\gamma}, \cdot)$}{cIGM}}\label{sec:cphi-arch}
Fix the data $\omega_S: \mathfrak{a}_R \to [0,1]$ for each $R \in \mathcal{L}(M_0)$ and $S \in \mathcal{P}(R)$, and similarly for each Levi subgroup of $G$, as in \S\ref{sec:cphi-nonarch}. Also fix $M \in \mathcal{L}(M_0)$.

\begin{definition}\label{def:crJ-arch}
	\index{JMpi-cr@${}^{c, r} J_{\tilde{M}}(\pi, f)$}
	Let $\pi \in D_{\elli, -}(\tilde{R}) \otimes \mes(R)^\vee$ where $R \in \mathcal{L}^M(M_0)$, and $f \in C^\infty_{c, \asp}(\tilde{G}, \tilde{K})$. Set
	\begin{multline*}
		{}^{c,r} J_{\tilde{M}}\left(  \pi^{\tilde{M}}, f \right) = {}^{c,r} J^{\tilde{G}}_{\tilde{M}}\left( \pi^{\tilde{M}}, f \right) := \int_{\mathfrak{a}_R} \dd X \sum_{S \in \mathcal{P}^G(R)} \omega_S(X) \\
		\int_{\nu_S + i\mathfrak{a}^*_R} \dd\lambda \; J^r_{\tilde{M}}\left( \pi_\lambda^{\tilde{M}}, f \right) e^{-\lrangle{\lambda, X}}
	\end{multline*}
	where the Haar measures on $\mathfrak{a}_R$ and $i\mathfrak{a}^*_R$ are dual to each other, and $\nu_S \in \mathfrak{a}_R^*$ satisfies $\lrangle{\nu_S, \check{\alpha}} \gg 0$ for all $\alpha \in \Sigma^S(A_R)$. This is independent of the choice of $(\nu_S)_S$.
\end{definition}

To justify this definition, the only required analytic input is the fact that for all $R \in \mathcal{L}^M(M_0)$ and $\pi \in D_{\elli, -}(\tilde{R}) \otimes \mes(R)^\vee$, the function $\lambda \mapsto J^r_{\tilde{M}}\left( \pi_\lambda^{\tilde{M}}, f \right)$ is
\begin{itemize}
	\item meromorphic with only a \emph{finite} number of polar hyperplanes,
	\item the polar hyperplanes take the form $\lrangle{\cdot, \check{\alpha}} = \text{const}$ where $\alpha \in \Sigma(A_M)$,
	\item rapidly decreasing on vertical strips.
\end{itemize}
These are ultimately based on our conditions \cite[\S 3]{Li12b} on normalizing factors, which are also the same as Arthur's conditions.

Note that when $R$ varies, distributions of the form $\pi^{\tilde{M}}$ generate $D_{\mathrm{temp}, -}(\tilde{M}) \otimes \mes(M)^\vee$ by \eqref{eqn:Dspec-Ind}. Therefore, one obtains ${}^{c,r} J_{\tilde{M}}\left(\pi, f \right)$ for all $\pi \in D_{\mathrm{temp}, -}(\tilde{G}) \otimes \mes(M)^\vee$. They extend to $D_{\mathrm{spec}, -}(\tilde{G}) \otimes \mes(M)^\vee$ by meromorphic continuation.

\begin{definition-proposition}\label{def:cphi-arch}
	\index{phiM-c}
	In the definition above, the function in the integrand is $C^\infty_c$ in $X$. There exists a continuous linear map
	\[ {}^c \phi_{\tilde{M}} = {}^c \phi^{\tilde{G}}_{\tilde{M}}: C^\infty_{c, \asp}(\tilde{G}, \tilde{K}) \otimes \mes(G) \to \orbI_{\asp}(\tilde{M}, \tilde{K} \cap \tilde{M}) \otimes \mes(M) \]
	characterized uniquely by
	\[ I^{\tilde{M}}\left( \pi, {}^c \phi_{\tilde{M}}(\cdot) \right) = {}^{c,r} J_{\tilde{M}}\left( \pi, \cdot\right) \]
	for all $\pi \in D_{\mathrm{temp},-}(\tilde{M}) \otimes \mes(M)^\vee$.
\end{definition-proposition}
\begin{proof}
	The arguments for \cite[IX.5.8 Proposition]{MW16-2} carry over verbatim. It is formally the same as the non-Archimedean case. The required analytic inputs are the same as in Definition \ref{def:crJ-arch}.
\end{proof}

Thus far, we have defined the linear maps:
\begin{align*}
	\phi^{\tilde{G}}_{\tilde{M}} & \quad \text{canonical normalization} \\
	{}^r \phi^{\tilde{G}}_{\tilde{M}} & \quad \text{via normalizing factors} \\
	{}^c \phi^{\tilde{G}}_{\tilde{M}} & \quad \text{the compactly-supported version}
\end{align*}
They are all continuous. The next goal is to explicate their relations. The situation can be summarized by the triangle
\[\begin{tikzcd}
	{}^c \phi \arrow[rr, "{{}^c \theta}"] \arrow[rd, "{{}^{c,r} \theta}"'] & & \phi \\
	& {}^r \phi \arrow[ru, "{{}^r \theta}"'] &
\end{tikzcd}\]
We have to define the edges.

\begin{definition}[Cf.\ {\cite[IX.5.5]{MW16-2}}]\label{def:rtheta-arch}
	\index{thetaM-r@${}^r \theta_{\tilde{M}}$, ${}^{c, r} \theta_{\tilde{M}}$}
	There is a continuous linear map
	\[ {}^r \theta_{\tilde{M}} = {}^r \theta^{\tilde{G}}_{\tilde{M}} : C^\infty_{c, \asp}(\tilde{G}, \tilde{K}) \otimes \mes(G) \to \orbI_{\mathrm{ac}, \asp}(\tilde{M}, \tilde{M} \cap \tilde{K}) \otimes \mes(M) \]
	characterized inductively by $\sum_{L \in \mathcal{L}^G(M)} {}^r \theta^{\tilde{L}}_{\tilde{M}} \circ {}^r \phi^{\tilde{G}}_{\tilde{L}} = \phi^{\tilde{G}}_{\tilde{M}}$. It extends by continuity to $\mathscr{C}_{\asp}(\tilde{G}) \otimes \mes(G) \to I\mathscr{C}_{\asp}(\tilde{M}) \otimes \mes(M)$.
\end{definition}

To complete the inductive definition, we must show that
\begin{itemize}
	\item ${}^r \theta_{\tilde{M}}$ extends to $C^\infty_{\mathrm{ac}, \asp}(\tilde{G}, \tilde{K}) \otimes \mes(G)$ and factors through $\orbI_{\asp, \mathrm{ac}}(\tilde{G}, \tilde{K}) \otimes \mes(G)$;
	\item in the Schwartz--Harish-Chandra setting, it factors through $I\mathscr{C}_{\asp}(\tilde{G}) \otimes \mes(G)$.
\end{itemize}
This will be achieved in \S\ref{sec:rhoGM-spec-arch}.

\begin{definition}[Cf.\ {\cite[IX.5.9]{MW16-2}}]\label{def:crtheta-arch}
	There is a continuous linear map
	\[ {}^{c,r} \theta_{\tilde{M}} = {}^{c,r} \theta^{\tilde{G}}_{\tilde{M}} : C^\infty_{c, \asp}(\tilde{G}, \tilde{K}) \otimes \mes(G) \to \orbI_{\mathrm{ac}, \asp}(\tilde{M}, \tilde{K} \cap \tilde{M}) \otimes \mes(M) \]
	characterized inductively by $\sum_{L \in \mathcal{L}^G(M)} {}^{c,r} \theta^{\tilde{L}}_{\tilde{M}} \circ {}^c \phi^{\tilde{G}}_{\tilde{L}} = {}^r \phi^{\tilde{G}}_{\tilde{M}}$.
\end{definition}

\begin{definition}\label{def:ctheta-arch}
	\index{thetaM-c}
	There is a continuous linear map
	\[ {}^c \theta_{\tilde{M}} = {}^c \theta^{\tilde{G}}_{\tilde{M}} : C^\infty_{c, \asp}(\tilde{G}, \tilde{K}) \otimes \mes(G) \to \orbI_{\mathrm{ac}, \asp}(\tilde{M}, \tilde{K} \cap \tilde{M}) \otimes \mes(M) \]
	characterized inductively by $\sum_{L \in \mathcal{L}^G(M)} {}^c \theta^{\tilde{L}}_{\tilde{M}} \circ {}^c \phi^{\tilde{G}}_{\tilde{L}} = \phi^{\tilde{G}}_{\tilde{M}}$.
\end{definition}

To finish the inductive definitions, we will show in \S\ref{sec:phi-inductive} that both ${}^{c,r} \theta_{\tilde{M}}$ and ${}^c \theta_{\tilde{M}}$ factor through $\orbI_{\asp}(\tilde{G}, \tilde{K}) \otimes \mes(G)$. Nevertheless, the following key result is independent of this assumption.

\begin{lemma}[Cf.\ {\cite[IX.5.12 Lemme]{MW16-2}}]\label{prop:c-rcr}
	Let $f \in C^\infty_{c, \asp}(\tilde{G}, \tilde{K}) \otimes \mes(G)$, then
	\[ {}^c \theta^{\tilde{G}}_{\tilde{M}} = \sum_{L \in \mathcal{L}^G(M)} {}^r \theta^{\tilde{L}}_{\tilde{M}} \circ {}^{c,r} \theta^{\tilde{G}}_{\tilde{L}}(f). \]
\end{lemma}
\begin{proof}
	We omit the proof since it is identical to the one in \textit{loc.\ cit.}, which does not involve any special analytic property.
\end{proof}

\begin{definition}[Cf.\ {\cite[IX.5.13]{MW16-2}}]\label{def:cIGM-arch}
	\index{IGM-c-gamma}
	For all $\tilde{\gamma} \in D_{\mathrm{orb},-}(\tilde{M}) \otimes \mes(M)^\vee$ and $f \in C^\infty_{c, \asp}(\tilde{G}, \tilde{K}) \otimes \mes(G)$, define inductively
	\begin{equation*}
		{}^c I_{\tilde{M}}(\tilde{\gamma}, f) = {}^c I^{\tilde{G}}_{\tilde{M}}(\tilde{\gamma}, f)
		= J_{\tilde{M}}(\tilde{\gamma}, f) - \sum_{\substack{L \in \mathcal{L}(M) \\ L \neq G}} {}^c I^{\tilde{L}}_{\tilde{M}}\left(\tilde{\gamma}, {}^c \phi_{\tilde{L}}(f) \right).
	\end{equation*}
\end{definition}

Again, in order to finish the inductive definition, we have to show that ${}^c I_{\tilde{M}}(\tilde{\gamma}, \cdot)$ factorizes through $\orbI_{\asp}(\tilde{G}, \tilde{K}) \otimes \mes(G)$. This will also be deferred to \S\ref{sec:phi-inductive}.

Granting the inductive definition, the following properties can be verified as in the non-Archimedean case.
\begin{itemize}
	\item For $f$ fixed, there exists a compact subset $\Omega \subset M(F)$ such that ${}^c I_{\tilde{M}}(\tilde{\gamma}, f) = 0$ whenever $\Supp(\tilde{\gamma}) \cap \rev^{-1}(\Omega) = \emptyset$.
	\item For all $R \in \mathcal{L}^M(M_0)$ and $\tilde{\gamma} \in D_{\mathrm{orb}, -}(\tilde{R}) \otimes \mes(R)^\vee$, we have
	\[ {}^c I_{\tilde{M}}\left( \tilde{\gamma}^{\tilde{M}}, f \right) = \sum_{L \in \mathcal{L}^G(R)} d^G_R(M, L) {}^c I^{\tilde{L}}_{\tilde{R}}(\tilde{\gamma}, f_{\tilde{L}}) \]
	where $\tilde{\gamma}^{\tilde{M}}$ is the parabolic induction of $\tilde{\gamma}$ and $f \mapsto f_{\tilde{L}}$ is the parabolic descent.
	\item Let $\tilde{\gamma} \in D_{\mathrm{orb}, -}(\tilde{M}) \otimes \mes(M)^\vee$, we have
	\begin{equation}\label{eqn:cpt-noncpt-arch}
		{}^c I_{\tilde{M}}(\tilde{\gamma}, \cdot) = \sum_{L \in \mathcal{L}(M)} I^{\tilde{L}}_{\tilde{M}}\left(\tilde{\gamma}, {}^c \theta_{\tilde{L}}(\cdot) \right).
	\end{equation}
	Indeed, $\sum_{L \in \mathcal{L}(M)} {}^c \theta^{\tilde{L}}_{\tilde{M}} \circ {}^c \phi_{\tilde{L}} = \phi_{\tilde{M}}$ implies
	\begin{multline*}
		\sum_{L \in \mathcal{L}^G(M)} \sum_{L' \in \mathcal{L}^L(M)} I^{\tilde{L}'}_{\tilde{M}}\left( \tilde{\gamma}, {}^c \theta^{\tilde{L}}_{\tilde{L}'} \circ {}^c \phi_{\tilde{L}}(\cdot) \right) = \sum_{L' \in \mathcal{L}(M)} I^{\tilde{L}'}_{\tilde{M}}\left( \tilde{\gamma}, \sum_{L \in \mathcal{L}^G(L')} {}^c \theta^{\tilde{L}}_{\tilde{L}'} \circ {}^c \phi_{\tilde{L}}(\cdot) \right) \\
		= \sum_{L' \in \mathcal{L}(M)} I^{\tilde{L}'}_{\tilde{M}}\left(\tilde{\gamma}, \phi_{\tilde{L}'}(\cdot) \right) = J_{\tilde{M}}(\tilde{\gamma}, \cdot).
	\end{multline*}
	Compare this with $\sum_{L \in \mathcal{L}(M)} {}^c I^{\tilde{L}}_{\tilde{M}}\left(\tilde{\gamma}, {}^c \phi_{\tilde{L}}(\cdot) \right) = J_{\tilde{M}}(\tilde{\gamma}, \cdot)$, and cancel the terms with $L \neq G$ by induction to obtain the result.
\end{itemize}

\section{The maps \texorpdfstring{$\rho^{\tilde{G}}_{\tilde{M}}$}{rhoGM}}\label{sec:rhoGM-spec-arch}
Fix $M \in \mathcal{L}(M_0)$ and define $U^G_M$ by \eqref{eqn:UGM-spec-arch}.

\begin{lemma}\label{prop:rtheta-rGM}
	For all $\pi \in \Pi_{\mathrm{temp},-}(\tilde{M})$ and $f \in \mathscr{C}_{\asp}(\tilde{G}) \otimes \mes(G)$, we have
	\begin{align*}
		I^{\tilde{M}}\left(\pi, {}^r \theta_{\tilde{M}}(f) \right) & = r^{\tilde{G}}_{\tilde{M}}(\pi) I^{\tilde{G}}(\pi^{\tilde{G}}, f) \\
		& = r^{\tilde{G}}_{\tilde{M}}(\pi) I^{\tilde{M}}(\pi, f_{\tilde{M}}).
	\end{align*}
	Consequently, ${}^r \theta_{\tilde{M}}$ factors through $I\mathscr{C}_{\asp}(\tilde{G}) \otimes \mes(G)$ in the Schwartz--Harish-Chandra setting.
\end{lemma}
\begin{proof}
	Argue by induction. The case $M=G$ is trivial. For $M \neq G$, Proposition \ref{prop:phimu-phir} implies
	\begin{multline*}
		I^{\tilde{M}}\left(\pi, {}^r \theta_{\tilde{M}}(f) \right) = I^{\tilde{M}}\left(\pi, \phi_{\tilde{M}}(f) \right) - \sum_{\substack{L \in \mathcal{L}(M) \\ L \neq G}} I^{\tilde{M}}\left( \pi, {}^r \theta^{\tilde{L}}_{\tilde{M}}\left( {}^r \phi_{\tilde{L}}(f) \right) \right) \\
		= I^{\tilde{M}}\left(\pi, \phi_{\tilde{M}}(f) \right) - \sum_{\substack{L \in \mathcal{L}(M) \\ L \neq G}} r^{\tilde{L}}_{\tilde{M}}(\pi) I^{\tilde{L}}\left( \pi^{\tilde{L}}, {}^r \phi_{\tilde{L}}(f) \right) = r^{\tilde{G}}_{\tilde{M}}(\pi) I^{\tilde{G}}\left( \pi^{\tilde{G}}, f \right).
	\end{multline*}
	The right-hand side depends only on the image of $f$ in $I\mathscr{C}_{\asp}(\tilde{G}) \otimes \mes(G)$.
\end{proof}

\begin{corollary}\label{prop:rtheta-rGM-factorization}
	The map ${}^r \theta_{\tilde{M}}$ extends canonically to $C^\infty_{\mathrm{ac}, \asp}(\tilde{G}) \otimes \mes(G)$ (resp.\ $C^\infty_{\mathrm{ac}, \asp}(\tilde{G}, \tilde{K}) \otimes \mes(G)$), and factors through $\orbI_{\mathrm{ac}, \asp}(\tilde{G}) \otimes \mes(G)$ (resp.\ $\orbI_{\mathrm{ac}, \asp}(\tilde{G}, \tilde{K}) \otimes \mes(G)$).
\end{corollary}
\begin{proof}
	Use the decompositions $\mathfrak{a}_M = \mathfrak{a}_G \oplus \mathfrak{a}^G_M$ and $i\mathfrak{a}_M^* = i\mathfrak{a}_G^* \oplus i\mathfrak{a}^{G,*}_M$. Let $X = X_G + X^G_M \in \mathfrak{a}_M$. Since $r^{\tilde{G}}_{\tilde{M}}(\pi)$ is invariant under $i\mathfrak{a}_G^*$-twists, Lemma \ref{prop:rtheta-rGM} yields
	\begin{multline*}
		I^{\tilde{M}}\left(\pi, X, {}^r \theta_{\tilde{M}}(f)\right) =
		\int_{i\mathfrak{a}^{G,*}_{M}} r^{\tilde{G}}_{\tilde{M}}(\pi_\mu) e^{-\lrangle{\mu, X^G_M}} \int_{i\mathfrak{a}_G^*} e^{-\lrangle{\lambda, X_G}} I^{\tilde{G}}\left( (\pi_\mu^{\tilde{G}})_\lambda, f \right) \dd\lambda \dd\mu \\
		= \int_{i\mathfrak{a}^{G,*}_{M}} r^{\tilde{G}}_{\tilde{M}}(\pi_\mu) e^{-\lrangle{\mu, X^G_M}} I^{\tilde{G}}\left( \pi_\mu^{\tilde{G}}, X_G, f \right) \dd\mu ,
	\end{multline*}
	which depends only on the orbital integrals of $f$, and is concentrated at $H_{\tilde{G}}(\cdot) = X_G$.
\end{proof}

Recall the vector spaces $D_{\mathrm{temp}, -}(\tilde{M})_0$ from Definition \ref{def:Dspec}.

\begin{lemma}[Cf.\ {\cite[IX.5.7 Lemme]{MW16-2}}]\label{prop:rhoGM-arch}
	\index{rhoGM@$\rho^{\tilde{G}}_{\tilde{M}}$}
	There exists a unique linear map
	\[ \rho^{\tilde{G}}_{\tilde{M}}: D_{\mathrm{temp},-}(\tilde{M})_0 \otimes \mes(M)^\vee \to U^G_M \otimes D_{\mathrm{temp},-}(\tilde{M})_0 \otimes \mes(M)^\vee \]
	satisfying
	\[ I^{\tilde{M}}\left( \pi, X, {}^r \theta^{\tilde{G}}_{\tilde{M}}(f) \right) = \int_{i\mathfrak{a}_M^*} I^{\tilde{M}}\left( \rho^{\tilde{G}}_{\tilde{M}}(\pi, \lambda)_\lambda , f_{\tilde{M}} \right) e^{-\lrangle{\lambda, X}} \dd \lambda \]
	for all $X \in \mathfrak{a}_M$, $f \in \mathscr{C}_{\asp}(\tilde{G}) \otimes \mes(G)$ and $\pi \in D_{\mathrm{temp},-}(\tilde{M})_0$.
\end{lemma}
\begin{proof}
	Consider the existence first. We are reduced to the case $\pi \in \Pi_{\mathrm{temp}, -}(\tilde{M})$ such that $\omega_\pi$ is trivial on $A_M(F)^\circ$. In this setting, put
	\[ \rho^{\tilde{G}}_{\tilde{M}}(\pi) := \left[ \lambda \mapsto r^{\tilde{G}}_{\tilde{M}}(\pi_\lambda) \right] \otimes \pi \]
	where $[\cdots]$ belongs to $U^G_M$ by Theorem \ref{prop:rGM-growth}. The general case is defined by linearity.
	
	The required formula for $I^{\tilde{M}}\left( \pi, X, {}^r \theta_{\tilde{M}}(f) \right)$ amounts to
	\[ I^{\tilde{M}}\left( \pi, X, {}^r \theta_{\tilde{M}}(f) \right) = \int_{i\mathfrak{a}_M^*} r^{\tilde{G}}_{\tilde{M}}(\pi_\lambda) I^{\tilde{G}}\left( \pi_\lambda^{\tilde{G}}, f \right) e^{-\lrangle{\lambda, X}} \dd \lambda. \]
	It boils down to the equality of Lemma \ref{prop:rtheta-rGM}.
	
	The uniqueness of $\rho^{\tilde{G}}_{\tilde{M}}$ is shown as in \textit{loc.\ cit.} In view of Fourier inversion, we must show that for all $v \in U^G_M \otimes D_{\mathrm{temp},-}(\tilde{M}) \otimes \mes(M)^\vee$, if
	\[ I^{\tilde{M}}(v(\lambda), f_{\tilde{M}}) = 0 \quad \text{for all}\quad f \in \mathscr{C}_{\asp}(\tilde{G}) \otimes \mes(G), \; \lambda \in i\mathfrak{a}_M^*, \]
	then $v = 0$. To deduce a contradiction, suppose $v = \sum_{i=1}^k u_i \otimes \pi_i$ such that $\pi_1, \ldots, \pi_k$ are linearly independent, $k \geq 1$ and $u_i \neq 0$ for all $i$. By the trace Paley--Wiener theorem. there exists $g_1 \in \orbI_{\asp}(\tilde{M}) \otimes \mes(M)$ such that $I^{\tilde{M}}(\pi_i, g_1) = \delta_{i, 1}$ (Kronecker's $\delta$).
	
	Now apply \cite[IX.5.6 Lemme]{MW16-2} to obtain $\mathcal{H} \subset \mathfrak{a}_{M, \CC}^*$, a finite union of $\mathfrak{a}_{G,\CC}^*$-invariant proper affine subspaces (the ``locus with symmetries''), such that for all $\lambda \in i\mathfrak{a}_M^* \smallsetminus \mathcal{H}$ and all $g \in \orbI_{\asp}(\tilde{M}) \otimes \mes(M)$, there exists $f \in \orbI_{\asp}(\tilde{G}) \otimes \mes(G)$ such that
	\[ I^{\tilde{M}}\left( \pi_{i, \lambda}, f_{\tilde{M}} \right) = I^{\tilde{M}}\left( \pi_{i, \lambda}, g \right), \quad i = 1, \ldots, k. \]
	
	Take $\lambda \in i\mathfrak{a}_M^* \smallsetminus \mathcal{H}$ such that $u_1(\lambda) \neq 0$ and take $g := e^{-\lrangle{\lambda, H_{\tilde{M}}(\cdot)}} g_1$ so that
	\[ I^{\tilde{M}}\left( \pi_{i, \lambda}, f_{\tilde{M}} \right) = I^{\tilde{M}}(\pi_{i, \lambda}, g) = I^{\tilde{M}}(\pi_i, g_1) = \delta_{i,1}, \quad i = 1, \ldots, k. \]
	Hence $I^{\tilde{M}}(v(\lambda), f_{\tilde{M}}) = u_1(\lambda) \neq 0$. Contradiction.
\end{proof}

We shall study the stable version of this map in \S\ref{sec:rho-stability-arch}.

\section{Completion of the inductive definitions}\label{sec:phi-inductive}
Conserve the conventions in \S\ref{sec:cphi-arch}. Suppose that $f \in C^\infty_{c, \asp}(\tilde{G}, \tilde{K}) \otimes \mes(G)$ maps to $0 \in \orbI_{\asp}(\tilde{G}, \tilde{K}) \otimes \mes(G)$. We will show that
\begin{equation}\label{eqn:theta-inductive}\begin{gathered}
	{}^{c,r} \theta_{\tilde{M}}(f) = {}^c \theta_{\tilde{M}}(f) = 0, \\
	{}^c I_{\tilde{M}}(\tilde{\gamma}, f) = 0, \quad \tilde{\gamma} \in D_{\mathrm{orb}, -}(\tilde{M}) \otimes \mes(M)^\vee ,
\end{gathered}\end{equation}
thereby finish the inductive definitions in \S\ref{sec:cphi-arch}. Although the arguments are similar to \cite[IX]{MW16-2}, we give the proofs since our version involves the normalizing factors.

\begin{lemma}\label{prop:cr-descent}
	For all $f \in C^\infty_{c, \asp}(\tilde{G}) \otimes \mes(G)$ and $R \in \mathcal{L}^M(M_0)$, we have
	\[ {}^{c, r} \theta_{\tilde{M}}(f)_{\tilde{R}} = \sum_{L \in \mathcal{L}(R)} d^G_R(M, L) \, {}^{c,r} \theta^{\tilde{L}}_{\tilde{R}}(f_{\tilde{Q}}) \]
	where $Q \in \mathcal{P}(L)$, and $f \mapsto f_{\tilde{Q}}$ denotes the parabolic descent in \eqref{eqn:f_P}.

	The same descent formula holds for ${}^r \theta_{\tilde{M}}$ and ${}^c \theta_{\tilde{M}}$.
\end{lemma}
\begin{proof}
	Observe that ${}^r \phi_{\tilde{M}}$ satisfies a similar descent formula, which boils down to the descent formula for $J^r_{\tilde{M}}$. Let us deduce the descent formula for ${}^c \phi_{\tilde{M}}$ next. For all $\pi \in D_{\mathrm{temp},-}(\tilde{R}) \otimes \mes(R)^\vee$,
	\begin{multline*}
		I^{\tilde{R}}\left( \pi, {}^c \phi_{\tilde{M}}(f)_{\tilde{R}} \right) = I^{\tilde{M}}\left( \pi^{\tilde{M}}, {}^c \phi_{\tilde{M}}(f) \right) = {}^{c, r} J_{\tilde{M}}\left(\pi^{\tilde{M}}, f \right) \\
		= \sum_{L \in \mathcal{L}(R)} d^G_R(M, L) \, {}^{c,r} J^{\tilde{L}}_{\tilde{R}}(\pi, f_{\tilde{Q}}) = I^{\tilde{R}}\left( \pi, \sum_{L \in \mathcal{L}(R)} d^G_R(M, L) \, {}^c \phi^{\tilde{L}}_{\tilde{R}}(f_{\tilde{Q}}) \right)
	\end{multline*}
	where we used the descent formula for ${}^{c,r} J_{\tilde{M}}$ (see \cite[p.940]{MW16-2}). Vary $\pi$ to reach the desired descent formula.
	
	To obtain the descent formula for ${}^{c,r} \theta_{\tilde{M}} = {}^r \phi_{\tilde{M}} - \sum_{L \neq G} {}^{c,r} \theta^{\tilde{L}}_{\tilde{M}} \circ {}^c \phi_{\tilde{L}}$ (Definition \ref{def:crtheta-arch}), one can iterate the arguments of \cite[VIII.1.6 Lemme]{MW16-2}, whose only inputs are the descent formulas for ${}^c \phi_{\tilde{M}}$ and ${}^r \phi_{\tilde{M}}$.
	
	The descent formulas for ${}^r \theta_{\tilde{M}}$ and ${}^c \theta_{\tilde{M}}$ follow in the same way, by using
	\[ {}^r \theta_{\tilde{M}} = \phi_{\tilde{M}} - \sum_{L \neq G} {}^r \theta^{\tilde{L}}_{\tilde{M}} \circ {}^r \phi_{\tilde{L}}, \quad {}^c \theta_{\tilde{M}} = \phi_{\tilde{M}} - \sum_{L \neq G} {}^c \theta^{\tilde{L}}_{\tilde{M}} \circ {}^c \phi_{\tilde{L}} \]
	(Definitions \ref{def:rtheta-arch}, \ref{def:ctheta-arch}) and the available descent formulas for ${}^r \phi_{\tilde{M}}$ and ${}^c \phi_{\tilde{M}}$.
\end{proof}

\begin{lemma}\label{prop:crtheta-Schwartz}
	For all $f \in C^\infty_{c, \asp}(\tilde{G}, \tilde{K}) \otimes \mes(G)$, its image ${}^{c,r} \theta_{\tilde{M}}(f)$ in $\orbI_{\mathrm{ac}, \asp}(\tilde{M}, \tilde{M} \cap \tilde{K}) \otimes \mes(M)$ is Schwartz.
\end{lemma}
\begin{proof}
	This follows inductively from Definition \ref{def:crtheta-arch}, since the image of ${}^r \phi^{\tilde{G}}_{\tilde{L}}$ (resp.\ ${}^c \phi^{\tilde{G}}_{\tilde{M}}$) consists of Schwartz (resp.\ compactly supported) functions, for various $L \in \mathcal{L}^G(M)$.
\end{proof}

Let $f \in C^\infty_{c, \asp}(\tilde{G}, \tilde{K}) \otimes \mes(G)$ and $\pi \in D_{\mathrm{temp},-}(\tilde{M}) \otimes \mes(M)^\vee$. As in the non-Archimedean case \S\ref{sec:cphi-nonarch}, one defines
\[ I^{\tilde{M}} \left(\pi, \lambda, {}^{c,r} \theta_{\tilde{M}}(f) \right) := \int_{\mathfrak{a}_M} e^{\lrangle{\lambda, X}} I^{\tilde{M}}\left( \pi, X, {}^{c,r} \theta_{\tilde{M}}(f) \right) \dd X. \]
for all $\lambda \in i\mathfrak{a}^*_M$. It equals $I^{\tilde{M}}(\pi_\lambda, f)$ if ${}^{c,r} \theta_{\tilde{M}}(f) \in \orbI_{\asp}(\tilde{M}, \tilde{M} \cap \tilde{K}) \otimes \mes(M)$.

In fact, the above extends to a meromorphic function in $\lambda \in \mathfrak{a}^*_{M, \CC}$, with only finitely many singular hyperplanes of the form $\lrangle{\cdot, \check{\alpha}} = \mathrm{const}$. The arguments are the same as in \cite[IX.5.10]{MW16-2}: by reasoning by induction via Definition \ref{def:crtheta-arch}, it reduces ultimately to similar properties for $J^r_{\tilde{M}}$ seen in \S\ref{sec:cphi-arch}.

Therefore, for general $\nu$ off the finitely many singular hyperplanes, it makes sense to set
\[ I^{\tilde{M}}\left( \pi, \nu, X, {}^{c,r} \theta_{\tilde{M}}(f) \right) := \int_{\nu + i\mathfrak{a}_M^*} I^{\tilde{M}}\left( \pi, \lambda, {}^{c,r} \theta_{\tilde{M}}(f) \right) e^{-\lrangle{\lambda, X}} \dd \lambda. \]

\begin{lemma}\label{prop:cr-vanishing}
	If $M \neq G$ and $\pi \in D_{\mathrm{ell}, -}(\tilde{M}) \otimes \mes(M)^\vee$, then
	\[ \sum_{S \in \mathcal{P}(M)} \omega_S(X) I^{\tilde{M}}\left( \pi, \nu_S, X, {}^{c,r} \theta_{\tilde{M}}(f) \right) = 0 \]
	for all $f \in C^\infty_{c, \asp}(\tilde{G}, \tilde{K}) \otimes \mes(G)$ and $X \in \mathfrak{a}_M$, provided that $\nu_S$ is sufficiently positive relative to $S$, for all $S \in \mathcal{P}(M)$.
\end{lemma}
\begin{proof}
	This corresponds to \cite[IX.5.10 Proposition]{MW16-2}. The crucial input is that $\lambda \mapsto J^r_{\tilde{M}}(\pi_\lambda, f)$ has only a finite number of polar hyperplanes.
\end{proof}

\begin{proof}[Completion of the inductive definitions]
	Suppose that $f \in C^\infty_{c,\asp}(\tilde{G}, \tilde{K}) \otimes \mes(G)$ maps to zero in $\orbI_{\asp}(\tilde{G}, \tilde{K}) \otimes \mes(G)$. The goal is to prove \eqref{eqn:theta-inductive}.
	
	First, we may assume $M \neq G$ and that the vanishing results are known when $M$ is replaced by any $L \in \mathcal{L}(M)$ with $M \neq L$. For each $\tilde{\gamma} \in D_{\mathrm{orb},-}(\tilde{M}) \otimes \mes(M)^\vee$ we have by \eqref{eqn:cpt-noncpt-arch}
	\begin{equation}\label{eqn:phi-inductive-aux1}
		{}^c I_{\tilde{M}}(\tilde{\gamma}, f) = \sum_{L \in \mathcal{L}(M)} I^{\tilde{L}}_{\tilde{M}}\left(\tilde{\gamma}, {}^c \theta_{\tilde{L}}(f)\right) = I^{\tilde{M}}\left( \tilde{\gamma}, {}^c \theta_{\tilde{M}}(f) \right);
	\end{equation}
	it also equals $I^{\tilde{M}}\left( \tilde{\gamma}, {}^{c,r} \theta_{\tilde{M}}(f) \right)$ since
	\begin{equation}\label{eqn:phi-inductive-aux2}
		{}^c \theta_{\tilde{M}}(f) = \sum_{L \in \mathcal{L}(M)} {}^r \theta^{\tilde{L}}_{\tilde{M}} \circ {}^{c,r} \theta_{\tilde{L}}(f) = {}^{c,r} \theta_{\tilde{M}}(f)
	\end{equation}
	by Lemma \ref{prop:c-rcr} and induction hypotheses. By the ``compactness'' of ${}^c I_{\tilde{M}}(\cdot, f)$, the expression \eqref{eqn:phi-inductive-aux1} vanishes if $\Supp(\tilde{\gamma}) \cap \rev^{-1}(\Omega) = \emptyset$, where $\Omega \subset M(F)$ is compact and depends only on $f$. It follows that ${}^{c,r} \theta_{\tilde{M}}(f) \in \orbI_{\asp}(\tilde{M}, \tilde{M} \cap \tilde{K}) \otimes \mes(M)$.

	Let $\pi \in D_{\elli, -}(\tilde{M}) \otimes \mes(M)^\vee$, we now know that $\lambda \mapsto I^{\tilde{M}}\left( \pi, \lambda, {}^{c,r} \theta_{\tilde{M}}(f) \right)$ is holomorphic on $\mathfrak{a}_{M, \CC}^*$. Hence the $I^{\tilde{M}}\left( \pi, \nu, X, {}^{c,r} \theta_{\tilde{M}}(f) \right)$, defined a priori for generic $\nu \in \mathfrak{a}_M^*$, is actually independent of $\nu$ by shifting contours. It is thus legitimate to replace all $\nu_S$ in Lemma \ref{prop:cr-vanishing} by $0$. This yields $I^{\tilde{M}}\left( \pi, 0, X, {}^{c,r} \theta_{\tilde{M}}(f) \right) = 0$ for all $X$. By Fourier inversion, $I^{\tilde{M}}\left( \pi, {}^{c,r} \theta_{\tilde{M}}(f) \right) = 0$ for all $\pi \in D_{\elli, -}(\tilde{M}) \otimes \mes(M)^\vee$.

	Lemma \ref{prop:cr-descent} and the induction hypotheses imply the cuspidality of ${}^{c,r} \theta_{\tilde{M}}(f)$. We conclude that ${}^{c,r} \theta_{\tilde{M}}(f) = 0$.

	It follows from \eqref{eqn:phi-inductive-aux2} that ${}^c \theta_{\tilde{M}}(f) = 0$. Hence ${}^c I_{\tilde{M}}(\tilde{\gamma}, f) = 0$ for all $\tilde{\gamma} \in D_{\mathrm{orb}, -}(\tilde{M}) \otimes \mes(M)^\vee$ by \eqref{eqn:phi-inductive-aux1}.
\end{proof}

We also record the following Archimedean counterpart for parts of Proposition \ref{prop:cthetaGM-nonarch}.

\begin{proposition}\label{prop:cthetaGM-arch}
	Let $f \in \orbI_{\asp}(\tilde{G}, \tilde{K}) \otimes \mes(G(F))$ and $M \in \mathcal{L}(M_0)$. Then
	\begin{itemize}
		\item ${}^c \theta_{\tilde{M}}(f)$ is Schwartz on $\tilde{M}$;
		\item for all $\pi \in D_{\mathrm{temp}, -}(\tilde{M}) \otimes \mes(M(F))^\vee$, the function $\lambda \mapsto I^{\tilde{M}}(\pi, \lambda, {}^c \theta_{\tilde{M}}(f))$ on $i\mathfrak{a}_M^*$ extends to a meromorphic function on $\mathfrak{a}^*_{M, \CC}$, with poles on hyperplanes of the form $\lrangle{\cdot, \check{\alpha}} = \text{const}$ for $\alpha \in \Sigma^G(A_M)$, and rapidly decreasing over vertical strips.
	\end{itemize}
\end{proposition}
\begin{proof}
	Same as \cite[IX.5.11]{MW16-2}. Although our ${}^{c, r} J_{\tilde{M}}$ are defined in terms of normalizing factors à la Arthur, whilst a ``rational normalization'' is employed in \textit{loc.\ cit.}, they satisfy the same analytic properties, and this is what one needs.
\end{proof}

\section{The stable versions}\label{sec:Sphi-arch}
In this section, $G^!$ is a direct product of quasisplit groups of the form $\GL(n)$ of $\SO(2n+1)$, for various $n$. This class of groups is closed under passing to Levi subgroups.

We fix a minimal Levi subgroup $M_0^!$ and a maximal compact subgroup $K^! \subset G^!(F)$ in good position relative to $M_0^!$, say arising from a Cartan involution. For such groups, we have
\begin{itemize}
	\item local Langlands correspondence for tempered representations;
	\item normalizing factors for intertwining operators, expressed in terms of $L$-functions;
	\item canonically normalized weighted characters $J_{M^!}(\pi^!, \cdot) = J^{G^!}_{M^!}(\pi^!, \cdot)$, for all $M^! \in \mathcal{L}(M_0^!)$ (see \S\ref{sec:weighted-characters-arch} or \cite{Ar98});
	\item normalized weighted characters $J^r_{M^!}(\pi^!, \cdot) = J^{G^!, r}_{M^!}(\pi^!, \cdot)$ using the normalizing factors above (see \S\ref{sec:weighted-characters-arch}, \cite[\S 7]{Ar89a} or \cite{Ar88-1});
	\item the compactly-supported version ${}^{c,r} J_{M^!}(\pi^!, \cdot) = {}^{c,r} J^{G^!}_{M^!}(\pi^!, \cdot)$, defined in exactly the way of Definition \ref{def:crJ-arch}. 
\end{itemize}

In this way, one defines $\phi^{G^!}_{M^!}$, ${}^r \phi^{G^!}_{M^!}$ and ${}^{c,r} \phi^{G^!}_{M^!}$. One also defines
\begin{gather*}
	{}^r \theta_{M^!} = {}^r \theta^{G^!}_{M^!}, \quad {}^{c,r} \theta_{M^!} = {}^{c,r} \theta^{G^!}_{M^!}, \quad {}^c \theta_{M^!} = {}^c \theta^{G^!}_{M^!}.
\end{gather*}
They are continuous linear maps $\orbI(G^!, K^!) \otimes \mes(G^!) \to \orbI_{\mathrm{ac}}(M^!, K^! \cap M^!) \otimes \mes(M^!)^\vee$. By the arguments of \S\ref{sec:rhoGM-spec-arch}, the map ${}^r \theta_{M^!}$ can even be extended in two ways, namely:
\begin{itemize}
	\item to Schwartz--Harish-Chandra spaces, namely $I\mathscr{C}(G^!) \otimes \mes(G^!) \to I\mathscr{C}(M^!) \otimes \mes(M^!)$;
	\item to $\orbI_{\mathrm{ac}}(G^!, K^!) \otimes \mes(G^!) \to \orbI_{\mathrm{ac}}(M^!, K^! \cap M^!) \otimes \mes(M^!)^\vee$.
\end{itemize}

Similarly, we build ${}^c I^{G^!}_{M^!}\left( \delta, f \right)$ as in Definition \ref{def:cIGM-arch}.

These constructions have already appeared in \cite{Ar88-1}, except that Arthur did not use canonically normalized characters at that time. In contrast, normalizing factors are not used in \cite{MW16-2}; a more complicated ``rational normalization'' is employed instead.

Denote the endoscopic data of $G^!$ as $\mathbf{G}^{!!} \in \Endo(G^!)$. Note that the class of groups $G^!$ is closed under passing to endoscopic groups, and the endoscopy for such $G^!$ involves neither $K$-groups nor $z$-extensions.

\begin{definition}
	\index{SthetaM-c}
	\index{SthetaM-r@${}^{c,r} S\theta^{G^{"!}}_{M^{"!}}$, ${}^r S\theta^{G^{"!}}_{M^{"!}}$}
	For every $M^! \in \mathcal{L}(M_0^!)$ and every placeholder of superscript $\star \in \left\{ r, (c,r), c\right\}$, define inductively
	\[ {}^\star S\theta^{G^!}_{M^!}: \orbI(G^!, K^!) \otimes \mes(G^!) \to S\orbI_{\mathrm{ac}}(M^!, K^! \cap M^!) \otimes \mes(M^!) \]
	by the requirement
	\[ \sum_{s \in \Endo_{\mathbf{M^!}}(G^!)} \iota_{M^!}(G^!, G^{!!}[s]) {}^\star S\theta^{G^{!!}[s]}_{M^!}\left( f^{G^{!!}[s]} \right) = {}^\star \theta^{G^!}_{M^!}(f) \quad \text{projected to}\; S\orbI_{\mathrm{ac}} \]
	where
	\begin{compactitem}
		\item $\mathbf{M}^!$ stands for the tautological endoscopic datum of $M^!$;
		\item to each $s \in \Endo_{\mathbf{M^!}}(G^!)$ we attach the endoscopic endoscopic datum $\mathbf{G}^{!!}[s]$ of $G^!$ by Arthur's recipe;
		\item $f \in \orbI(G^!, K^!) \otimes \mes(G^!)$ and $f^{G^{!!}[s]}$ is its transfer to $G^{!!}[s]$.
	\end{compactitem}
\end{definition}

As usual, for ${}^r S\theta^{G^!}_{M^!}$ can be extended to $I\mathscr{C}(G^!) \otimes \mes(G^!)$ or $\orbI_{\mathrm{ac}}(G^!, K^!) \otimes \mes(G)$. The following stability is required to complete the inductive definition.

\begin{proposition}\label{prop:Stheta-stability-arch}
	For each $\star \in \left\{ r, (c,r), c \right\}$, the map ${}^\star S\theta^{G^!}_{M^!}$ factors through $S\orbI(G^!, K^!) \otimes \mes(G^!)$. Moreover, ${}^r S\theta^{G^!}_{M^!}$ extends canonically to $S\orbI_{\mathrm{ac}}(G^!, K^!) \otimes \mes(G^!)$.
\end{proposition}

For ${}^r S\theta^{G^!}_{M^!}$, the required properties in Proposition \ref{prop:Stheta-stability-arch} follow readily from the upcoming Lemma \ref{prop:sigmaGM-arch}; see the discussions in \S\ref{sec:rhoGM-spec-arch}. For the remaining two cases, the proof is deferred to \S\ref{sec:rho-stability-arch}.

The following result does not presuppose stability. First off, for all $M^! \in \mathcal{L}(M_0^!)$ and $\pi^! \in SD_{\mathrm{temp}}(M^!) \otimes \mes(M^!)^\vee$ one defines $S^{M^!}\left( \pi^!, \nu, X, {}^{c,r} S\theta_{M^!}(\cdot) \right)$ as in \S\ref{sec:phi-inductive}. Here $X \in \mathfrak{a}_{M^!}$, and $\nu \in \mathfrak{a}_{M^!, \CC}^*$ is in general position.

\begin{lemma}\label{prop:Scr-vanishing}
	If $M^! \neq G^!$ and $\pi^! \in SD_{\mathrm{ell}}(M^!) \otimes \mes(M^!)^\vee$, then
	\[ \sum_{S^! \in \mathcal{P}(M^!)} \omega_{S^!}(X) S^{M^!}\left( \pi^!, \nu_S, X, {}^{c,r} S\theta_{M^!}(f) \right) = 0 \]
	for all $f \in C^\infty_c(G^!, K^!) \otimes \mes(G^!)$ and $X \in \mathfrak{a}_{M^!}$, provided that $\nu_{S^!}$ is sufficiently positive relative to $S^!$, for all $S^!$.
\end{lemma}
\begin{proof}
	This corresponds to \cite[IX.6.2 Proposition]{MW16-2}. Again, the crucial input is that $\lambda \mapsto J^r_{M^!}(\pi^!_\lambda, f)$ has only a finite number of polar hyperplanes.
\end{proof}

Define $U^{G^!}_{M^!}$ as in \eqref{eqn:UGM-spec-arch}. The following is the counterpart of Lemma \ref{prop:rhoGM-arch}. In its proof, the stability of ${}^r S\theta^{G^!}_{M^!}$ is not presupposed.

\begin{lemma}\label{prop:sigmaGM-arch}
	\index{sigmaGM@$\sigma^{G^{"!}}_{M^{"!}}$}
	There exists a unique linear map
	\[ \sigma^{G^!}_{M^!}: SD_{\mathrm{temp}}(M)_0 \otimes \mes(M^!)^\vee \to U^{G^!}_{M^!} \otimes D_{\mathrm{temp}}(M^!)_0 \otimes \mes(M^!)^\vee \]
	satisfying
	\[ S^{M^!}\left( \pi^!, X, {}^r S\theta_{M^!}(f) \right) = \int_{i\mathfrak{a}_{M^!}^*} I^{M^!}\left( \sigma^{G^!}_{M^!}(\pi^!, \lambda)_\lambda, f_{M^!} \right) e^{-\lrangle{\lambda, X}} \dd \lambda \]
	for all $X \in \mathfrak{a}_{M^!}$, $f \in \orbI(G^!) \otimes \mes(G^!)$ and $\pi^! \in SD_{\mathrm{temp}}(M^!)_0$.
\end{lemma}
\begin{proof}
	Same as \cite[IX.6.3 Lemme]{MW16-2}. To summarize, the recipe is to set
	\begin{equation}\label{eqn:sigmaGM-arch}
		\sigma^{G^!}_{M^!}(\pi^!) := \rho^{G^!}_{M^!}(\pi^!) - \sum_{\substack{s \in \Endo_{\mathbf{M}^!}(G^!) \\ s \neq 1}} i_{M^!}\left( G^!, G^{!!}[s] \right) \sigma^{G^{!!}[s]}_{M^!}(\pi^!)
	\end{equation}
	inductively, where $\mathbf{M}^!$ stands for the tautological endoscopic datum of $M^!$, and we use the identification $U^{G^{!!}[s]}_{M^!} \simeq U^{G^!}_{M^!}$ via $\mathfrak{a}^{G^{!!}[s]}_{M^!} \simeq \mathfrak{a}^{G^!}_{M^!}$. The arguments for uniqueness are the same as Lemma \ref{prop:rhoGM-arch}.
\end{proof}

The following will be established in \S\ref{sec:rho-stability-arch}.

\begin{lemma}[Cf.\ {\cite[IX.6.4]{MW16-2}}]\label{prop:sigmaGM-stability-arch}
	The image of $\sigma^{G^!}_{M^!}$ is contained in $U^{G^!}_{M^!} \otimes SD_{\mathrm{temp}}(M^!)_0 \otimes \mes(M^!)^\vee$.
\end{lemma}

Now we can state the compactly supported variant of the stable distributions $S^{G^!}_{M^!}(\delta, \cdot)$ on the geometric side.

\begin{definition}
	For all $\delta \in SD_{\mathrm{orb}}(M^!) \otimes \mes(M^!)^\vee$ and $f \in \orbI(G^!, K^!) \otimes \mes(G^!)$, define
	\[ {}^c S^{G^!}_{M^!}\left( \delta, f \right) := {}^c I^{G^!}_{M^!}\left(\delta, f \right) - \sum_{\substack{s \in \Endo_{\mathbf{M}^!}(G^!) \\ s \neq 1}} i_{M^!}\left( G^!, G^{!!}[s] \right) {}^c S^{G^{!!}[s]}_{M^!}\left(\delta, f^{G^{!!}[s]} \right) \]
	inductively, where $f \mapsto f^{G^{!!}[s]}$ is the transfer.
\end{definition}

In order to accomplish the inductive definition, we need the following result on stability, to be proved in \S\ref{sec:rho-stability-arch}.

\begin{lemma}\label{prop:cSGM-stability-arch}
	The continuous linear form ${}^c S^{G^!}_{M^!}\left( \delta, \cdot \right)$ factors through $S\orbI(G^!, K^!) \otimes \mes(G^!)$, for all $\delta$.
\end{lemma}

Assuming the validity of the inductive definition, the compact and non-compact versions are related as below.

\begin{lemma}\label{prop:cpt-noncpt-S-arch}
	For all $\delta \in SD_{\mathrm{orb}}(M^!) \otimes \mes(M^!)^\vee$, we have
	\[ {}^c S^{G^!}_{M^!}\left(\delta, \cdot \right) = \sum_{L^! \in \mathcal{L}(M^!)} S^{L^!}_{M^!}\left(\delta, {}^c S\theta_{L^!}(\cdot) \right). \]
\end{lemma}
\begin{proof}
	Same as the non-Archimedean case \cite[VIII.4.1]{MW16-2}.
\end{proof}

We will only use the case when $\delta$ is of $G^!$-regular support.

\begin{lemma}\label{prop:cpt-noncpt-S}
	We have
	\begin{equation*}
		{}^c S\theta^{G^!}_{M^!} = \sum_{L^! \in \mathcal{L}(M^!)} {}^r S\theta^{L^!}_{M^!} \circ {}^{c,r} S\theta^{G^!}_{L^!}.
	\end{equation*}
\end{lemma}
\begin{proof}
	The arguments are formal; see the reference cited in Lemma \ref{prop:c-rcr}.
\end{proof}

\section{Proof of stability}\label{sec:rho-stability-arch}
Conserve the conventions of \S\ref{sec:Sphi-arch}. In particular, we fix $M^! \in \mathcal{L}^{G^!}(M_0^!)$.

Following \cite[Theorem 5]{Ar99b}, adapted to the local Archimedean setting, one defines a function
\begin{equation*}
	s^{G^!}_{M^!}: \Phi_{\mathrm{bdd}}(M^!) \to \CC.
\end{equation*}
\index{sGM-phi@$s^{G^{"!}}_{M^{"!}}(\phi^{"!})$}

It is the stable version of the factors $r^{G^!}_{M^!}(\pi^!)$ made from normalizing factors where $\pi^! \in \Pi_{\mathrm{temp}}(M^!)$; see \S\ref{sec:weighted-characters-arch}. Specifically, let $\mathbf{M}^{!!} \in \Endo_{\elli}(M^!)$ and $\phi^{!!} \in \Phi_{\mathrm{bdd}}(M^{!!})$ be any $L$-parameter corresponding to $\pi^! \in \Pi_{\mathrm{temp}}(M^!)$ via $\Lgrp{M^{!!}} \to \Lgrp{M^!}$, the following equality is imposed
\begin{equation}
	r^{G^!}_{M^!}(\pi^!) = \sum_{s \in \Endo_{\mathbf{M}^{!!}(G^!)}} i_{M^{!!}}(G^!, G^{!!}[s]) s^{G^{!!}[s]}_{M^{!!}}(\phi^{!!}).
\end{equation}

The inductive definition is thus given by $s^{G^!}_{M^!}(\phi^!) = r^{G^!}_{M^!}(\pi^!) - \sum_{s \neq 1} i_{M^!}(G^!, G^{!!}[s]) s^{G^{!!}[s]}_{M^!}(\phi^!)$, by setting $\mathbf{M}^{!!} = \mathbf{M}^!$. Consequently, $\lambda \mapsto s^{G^!}_{M^!}(\phi_\lambda)$ extends to a meromorphic function on $\mathfrak{a}_{M^! ,\CC}^*$.

\begin{lemma}
	For all $\phi^! \in \Phi_{\mathrm{bdd}}(M^!)$, the function $\lambda \mapsto s^{G^!}_{M^!}(\phi^!_\lambda)$ belongs to $U^{G^!}_{M^!}$. 
\end{lemma}
\begin{proof}
	Since the twists on $\phi^!$ by $\mathfrak{a}_{M^!, \CC}^*$ is compatible with the twists on representations under local Langlands correspondence, this reduces to the moderate growth of $\lambda \mapsto r^{G^!}_{M^!}(\pi^!_\lambda)$; see the references in Theorem \ref{prop:rGM-growth}.
\end{proof}

\begin{proof}[Proof of Lemma \ref{prop:sigmaGM-stability-arch}]
	Fix $\phi^! \in \Phi_{\mathrm{bdd}}(M^!)$, denote the tempered $L$-packet as $\Pi_{\phi^!}$ and consider the corresponding stable character
	\begin{equation}\label{eqn:sigmaGM-stability-arch-aux0}
		\pi^! = \sum_{\sigma \in \Pi_{\phi^!}} m_\sigma \sigma \; \in SD_{\mathrm{temp}}(M^!)_0 \otimes \mes(M)^\vee.
	\end{equation}
	For our groups, we shall actually take $m_\sigma = 1$. The stability of $\sigma^{G^!}_{M^!}(\pi^!)$ will follow from the claim that
	\begin{equation}\label{eqn:sigmaGM-stability-arch-aux}
		\sigma^{G^!}_{M^!}(\pi^!, \lambda) = s^{G^!}_{M^!}(\phi_\lambda) \otimes \pi^!.
	\end{equation}

	Indeed, we may write $r^{G^!}_{M^!}(\sigma_\lambda) = r^{G^!}_{M^!}(\phi^!_\lambda)$ since the normalizing factors depend only on tempered $L$-parameters by the construction in \cite[\S 3]{Ar89a}. By \eqref{eqn:sigmaGM-arch}, the linearity of $\rho^{G^!}_{M^!}$ and induction on $\dim G^!$, we have
	\begin{align*}
		\sigma^{G^!}_{M^!}(\pi^!, \lambda) & = \sum_{\sigma \in \Pi_{\phi^!}} m_\sigma r^{G^!}_{M^!}(\sigma_\lambda) \otimes \sigma - \sum_{\substack{s \in \Endo_{\mathbf{M}^!}(G^!) \\ s \neq 1}} i_{M^!}\left(G^!, G^{!!}[s] \right) s^{G^{!!}[s]}_{M^!}(\phi^!_\lambda) \otimes \pi^! \\
		& = r^{G^!}_{M^!}\left( \phi^!_\lambda \right) \otimes \pi^! - \sum_{\substack{s \in \Endo_{\mathbf{M}^!}(G^!) \\ s \neq 1}} i_{M^!}\left(G^!, G^{!!}[s] \right) s^{G^{!!}[s]}_{M^!}(\phi^!_\lambda) \otimes \pi^! \\
		& = s^{G^!}_{M^!}(\phi^!_\lambda) \otimes \pi^!
	\end{align*}
	by the definition of $s^{G^!}_{M^!}$. This proves \eqref{eqn:sigmaGM-stability-arch-aux}.
\end{proof}

\begin{lemma}\label{prop:crS-descent}
	Let $\star \in \{r, (c, r), c\}$. For all $f \in \orbI(G^!, K^!) \otimes \mes(G^!)$ and $R^! \in \mathcal{L}^{M^!}(M_0^!)$, we have
	\[ {}^{\star} S\theta_{M^!}(f)_{R^!} = \sum_{L^! \in \mathcal{L}^{G^!}(R^!)} e^{G^!}_{R^!}(M^!, L^!) \, {}^{\star} S\theta^{L^!}_{R^!}(f_{L^!}). \]
\end{lemma}
\begin{proof}
	This follows from the descent formulas à la Lemma \ref{prop:cr-descent}, together with the same combinatorial arguments from \cite[VIII.2.3]{MW16-2}. Cf.\ the discussion after \cite[IX.6.1, Proposition]{MW16-2}.
\end{proof}

The notion of Schwartz functions also apply to $S\orbI_{\mathrm{ac}}(M^!) \otimes \mes(M^!)$. Cf.\ Definition \ref{def:Schwartz-arch}.

\begin{lemma}\label{prop:crStheta-Schwartz}
	For all $f \in \orbI(G^!, K^!) \otimes \mes(G^!)$, its image ${}^{c,r} S\theta_{M^!}(f)$ is Schwartz.
\end{lemma}
\begin{proof}
	Same as the proof of Lemma \ref{prop:crtheta-Schwartz}.
\end{proof}

\begin{proof}[Proof of the stablity of ${}^{c,r} S\theta^{G^!}_{M^!}$, ${}^c S\theta^{G^!}_{M^!}$ and ${}^c S^{G^!}_{M^!}(\delta, \cdot)$]
	Given Lemmas \ref{prop:Scr-vanishing}, \ref{prop:crS-descent} and \ref{prop:crStheta-Schwartz}, the proof is formally the same as that in \S\ref{sec:phi-inductive}.
\end{proof}

\section{Intermezzo: an easy stabilization}\label{sec:easy-stabilization}
Take $\tilde{G}$, $\tilde{K}$ and $M \in \mathcal{L}(M_0)$ as in \S\ref{sec:weighted-characters-arch}. On the other hand, let $G^!$, $K^!$ and $M^! \in \mathcal{L}(M_0^!)$ be as in \S\ref{sec:Sphi-arch}.

Recall that we have defined the factors
\begin{itemize}
	\item $r^{\tilde{L}}_{\tilde{M}}(\pi_\lambda)$ as a meromorphic family in $\lambda \in \mathfrak{a}_{M, \CC}^*$, analytic on $i\mathfrak{a}_M^*$, for each $\pi \in \Pi_{\mathrm{temp},-}(\tilde{M})$;
	\item $s^{G^!}_{M^!}(\phi^!_\lambda)$ as a meromorphic family in $\lambda \in \mathfrak{a}_{M^!, \CC}^*$, analytic on $i\mathfrak{a}_{M^!}^*$, for each $\phi^! \in \Phi_{\mathrm{bdd}}(M^!)$.
\end{itemize}

Assume for a moment that $\tilde{G} = \Mp(W)$, with $2n = \dim W$. In the discussions prior to Theorem \ref{prop:normalizing-real}, we saw that one can match various objects defined for $\tilde{G}$ and $G' := \SO(2n+1)$. In particular, the Levi subgroups $M \subset G$ are in bijection with their counterparts $M' \subset G'$ up to conjugacy. One can also match infinitesimal characters between genuine representation of $\tilde{G}$ (resp.\ $\tilde{M}$) and those of representations of $G'$ (resp.\ $M'$). Ditto for the infinitesimal characters of discrete series.

\begin{lemma}
	Suppose that $\pi \in \Pi_{\mathrm{temp}, -}(\tilde{M})$. The $r^{\tilde{G}}_{\tilde{M}}(\pi)$ equals its counterpart of $G' := \SO(2n+1)$ by matching infinitesimal characters.
\end{lemma}
\begin{proof}
	Immediate consequence of Theorems \ref{prop:normalizing-real}, \ref{prop:normalizing-cplx}, and the definition of $r^{\tilde{G}}_{\tilde{M}}(\pi)$.
\end{proof}

Consider now the situation \eqref{eqn:s-situation}, with $\mathbf{G}^! := \mathbf{G}^![s]$ for $s \in \Endo_{\mathbf{M}^!}(\tilde{G})$. The superscript $'$ will stand for the $\SO(2n+1)$-counterparts of various objects.

\begin{theorem}\label{prop:easy-stabilization}
	Suppose that $\pi \in \Pi_{\mathrm{temp}, -}(\tilde{M})$ corresponds to $\phi^! \in \Phi_{2, \mathrm{bdd}}(M^!)$ by matching infinitesimal characters. Then
	\[ r^{\tilde{G}}_{\tilde{M}}(\pi) = \sum_{s \in \Endo_{\mathbf{M}^!}(\tilde{G})} i_{M^!}(\tilde{G}, G^![s]) s^{G^![s]}_{M^!}(\phi^!). \]
\end{theorem}
\begin{proof}
	Our arguments are based on the $G' := \SO(2n+1)$ version
	\begin{equation}\label{eqn:easy-stabilization-SO}
		r^{G'}_{M'}(\pi') = \sum_{s' \in \Endo_{\mathbf{M}'^!}(G')} i_{M^{'!}}\left( G', G'^![s'] \right) s^{G'^![s']}_{M'^!}(\phi'^!)
	\end{equation}
	where $\pi' \in \Pi_2(M')$ has infinitesimal character matching $\pi$. In \cite{Ar99b}, a close analogue of \eqref{eqn:easy-stabilization-SO} is proven for normalizing factors made from Satake parameters. Since both sorts of normalizing factors are by prescribed by $L$-functions, the same combinatorial techniques thereof, especially \cite[Lemma 4, Theorem 5]{Ar99b}, also prove \eqref{eqn:easy-stabilization-SO}. Our case is even easier since their are fewer analytic issues\footnote{More precisely, the $L$-functions here are local.}.

	Since $Z_{\tilde{M}^\vee}^\circ = Z_{(M')^\vee}^\circ$, by \cite[Lemma 1.1]{Ar99} there is a surjection $\Endo_{\mathbf{M}^!}(\tilde{G}) \twoheadrightarrow \Endo_{\mathbf{M}'^!}(G')$. If $s \mapsto s'$, then $G^![s] = G'^![s']$ and one can match $s^{G^![s]}_{M^!}(\phi^!)$ with its $G'^![s']$-counterpart. Also,
	\[ i_{M^!}(\tilde{G}, G^![s]) = 2^{ \# \;\text{of $\SO$ in $M^!$} - \# \;\text{of $\SO$ in $G^![s]$}}. \]
	
	Suppose that there are no $\SO$-factors in $M^!$, so that $M = M'$ and $M^! = M'^!$. One readily checks $2 i_{M^!}(\tilde{G}, G^![s]) = i_{M^{'!}}(G', G'^![s'])$. On the other hand, there is an involution without fixed points on $\Endo_{\mathbf{M}^!}(\tilde{G})$ corresponding to multiplication by $-1 \in \tilde{G}^\vee$, making it a $\{\pm 1\}$-torsor over $\Endo_{\mathbf{M}'^!}(G')$. These two differences compensate each other, and the required identity reduces to the case of $G'$.
	
	If there are $\SO$-factors in $M^!$, then $i_{M^!}(\tilde{G}, G^![s])$ equals its avatar for $\SO(2n+1)$. In this case $\Endo_{\mathbf{M}^!}(\tilde{G})$ is seen to equal its avatar for $G'$. The same argument reduces the equality to the case of $G'$.
\end{proof}

We stated Theorem \ref{prop:easy-stabilization} under the tacit assumption that $\tilde{G} = \Mp(W)$, but it is routine to extend it to any $\tilde{G}$ of metaplectic type.

\section{The endoscopic versions}\label{sec:cI-Endo-arch}
Now take $\tilde{G}$ and $M \in \mathcal{L}(M_0)$ as in \S\ref{sec:weighted-characters-arch}.

In order to define ${}^{c,r} J^{\tilde{L}}_{\tilde{M}}$ and ${}^{c,r} J^{L^!}_{M^!}$ in a compatible way, for all $M \subset L \subset G$ and $M^! \subset L^! \subset G^!$ on the stable side, we choose the auxiliary functions $\omega_S$ and $\omega_{S^!}$ subject to the conditions \eqref{eqn:omega-compatibility}, \eqref{eqn:omega-compatibility-Endo}, etc.

\begin{definition}
	\index{thetaM-Endo-c}
	\index{thetaM-Endo-r@${}^r \theta^{\Endo}_{\tilde{M}}$, ${}^{c, r} \theta^{\Endo}_{\tilde{M}}$}
	Fix $\mathbf{M}^! \in \Endo_{\elli}(\tilde{M})$. Let $\star \in \left\{ r, (c,r), c \right\}$. Define
	\begin{align*}
		{}^\star \theta^{\Endo}_{\tilde{M}}(\mathbf{M}^!, f) & = {}^\star \theta^{\tilde{G}, \Endo}_{\tilde{M}}(\mathbf{M}^!, f) \\
		& = \sum_{s \in \Endo_{\mathbf{M}^!}(\tilde{G})} i_{M^!}(\tilde{G}, G^![s]) {}^\star S\theta^{G^![s]}_{M^!}\left(f^{G^![s]} \right)[s]
	\end{align*}
	where $f \in \orbI_{\asp}(\tilde{G}, \tilde{K}) \otimes \mes(G)$, $f^{G^![s]} := \Trans_{\mathbf{M}^!, \tilde{M}}(f)$ and for every $s$ we define the automorphism $g \mapsto g(\cdot)[s] := g(\cdot z[s])$ of $S\orbI_{\mathrm{ac}}(G^![s], K^!) \otimes \mes(G^![s])$.
\end{definition}

\begin{proposition}[Cf.\ {\cite[VIII.3.4]{MW16-2}}]\label{prop:theta-Endo-glue-arch}
	Let $\star \in \left\{ r, (c,r), c \right\}$. There exists a unique linear map
	\[ {}^\star \theta^{\Endo}_{\tilde{M}} = {}^\star \theta^{\tilde{G}, \Endo}_{\tilde{M}} : \orbI_{\asp}(\tilde{G}, \tilde{K}) \otimes \mes(G) \to \orbI_{\mathrm{ac}, \asp}(\tilde{M}) \otimes \mes(M) \]
	which is continuous and characterized by
	\[ \Trans_{\mathbf{M}^!, \tilde{M}} \circ {}^\star \theta^{\Endo}_{\tilde{M}} = {}^\star \theta^{\Endo}_{\tilde{M}}(\mathbf{M}^!, \cdot) \]
	for all $\mathbf{M}^! \in \Endo_{\elli}(\tilde{M})$.
\end{proposition}
\begin{proof}
	Identical to the non-Archimedean version (Proposition \ref{prop:theta-Endo-glue-nonarch}).
\end{proof}

\begin{remark}
	One can also define ${}^r \theta^{\Endo}_{\tilde{M}}(\mathbf{M}^!, f)$ for $f \in I\mathscr{C}_{\asp}(\tilde{G}) \otimes \mes(G)$.
\end{remark}

Below is the endoscopic version of Lemma \ref{prop:c-rcr}.

\begin{lemma}[Cf.\ {\cite[IX.7.2 Lemme (i)]{MW16-2}}]\label{prop:c-rcr-Endo}
	We have ${}^c \theta^{\tilde{G}, \Endo}_{\tilde{M}} = \sum_{L \in \mathcal{L}(M)} {}^r \theta^{\tilde{L}, \Endo}_{\tilde{M}} \circ {}^{c,r} \theta^{\tilde{G}, \Endo}_{\tilde{L}}$.
\end{lemma}
\begin{proof}
	It suffices to fix $f \in \orbI_{\asp}(\tilde{G}, \tilde{K}) \otimes \mes(G)$, $\mathbf{M}^! \in \Endo_{\elli}(\tilde{M})$ and prove the equality after applying $\Trans_{\mathbf{M}^!, \tilde{M}}$ to both sides. By Lemma \ref{prop:cpt-noncpt-S}, the left-hand side is
	\begin{multline*}
		\sum_{s \in \Endo_{\mathbf{M}^!}(\tilde{G})} i_{M^!}(\tilde{G}, G^![s]) {}^c S\theta^{G^![s]}_{M^!}\left( f^{G^![s]} \right)[s] \\
		= \sum_{s \in \Endo_{\mathbf{M}^!}(\tilde{G})} \sum_{L^! \in \mathcal{L}^{G^![s]}(M^!)} i_{M^!}(\tilde{G}, G^![s]) \left( {}^r S\theta^{L^!}_{M^!} \circ {}^{c,r} S\theta^{G^![s]}_{L^!}\left( f^{G^![s]} \right) \right) [s]
	\end{multline*}
	
	Next, apply Lemma \ref{prop:sL-Ls} (with $R = M$) to transform the sum over $(s, L^!)$ into $(L, r^L, r_L)$ where $L \in \mathcal{L}^G(M)$, $r^L \in \Endo_{\mathbf{M}^!}(\tilde{L})$, $r_L \in \Endo_{\mathbf{L}^![r^L]}(\tilde{G})$, so that
	\[ G^![s] = G^![r_L], \quad L^! = L^![r^L], \quad z[s] = z[r^L] z[r_L], \]
	and by Lemma \ref{prop:i-transitivity}, the sum becomes
	\begin{multline*}
		\sum_L \sum_{r^L} i_{M^!}(\tilde{L}, L^![r^L]) {}^r S\theta^{L^![r^L]}_{M^!}
		\left( \sum_{r_L} i_{L^![r^L]}(\tilde{G}, G^![r_L]) {}^{c,r} S\theta^{G^![r_L]}_{L^![r^L]}\left( f^{G^![r_L]} \right) [r_L]  \right) [r^L] \\
		= \sum_L \sum_{r^L} i_{M^!}(\tilde{L}, L^![r^L]) {}^r S\theta^{L^![r^L]}_{M^!} \left( {}^{c,r} \theta^{\tilde{G}, \Endo}_{\tilde{L}}\left( \mathbf{L}^![r^L], f \right)\right) [r^L] \\
		= \sum_L \sum_{r^L} i_{M^!}(\tilde{L}, L^![r^L]) {}^r S\theta^{L^![r^L]}_{M^!} \left( {}^{c,r} \theta^{\tilde{G}, \Endo}_{\tilde{L}}(f)^{L^![r^L]} \right) [r^L] \\
		= \sum_{L \in \mathcal{L}^G(M)} {}^r \theta^{\tilde{L}, \Endo}_{\tilde{M}}\left(\mathbf{M}^!, {}^{c,r} \theta^{\tilde{G}, \Endo}_{\tilde{L}}(f) \right).
	\end{multline*}
	We have used the fact that $z[r_L]$ belongs to the connected center of $L^![r^L]$, hence the translation $[z_L]$ commutes with ${}^r S\theta^{L^![r^L]}_{M^!}$.
	
	The outcome equals the transfer of $\sum_{L \in \mathcal{L}^G(M)} {}^r \theta^{\tilde{L}, \Endo}_{\tilde{M}} \circ {}^{c,r} \theta^{\tilde{G}, \Endo}_{\tilde{L}}(f)$ to $M^!$, as desired.
\end{proof}

The following is the endoscopic version of Lemma \ref{prop:cr-descent}.

\begin{lemma}\label{prop:cr-descent-Endo}
	Let $\star \in \{r, (c, r), c\}$. For all $f \in \orbI_{\asp}(\tilde{G}, \tilde{K}) \otimes \mes(G)$ and $R \in \mathcal{L}^M(M_0)$, we have
	\[ {}^{\star} \theta^{\tilde{G}, \Endo}_{\tilde{M}}(f)_{\tilde{R}} = \sum_{L \in \mathcal{L}(R)} d^G_R(M, L) \, {}^{\star} \theta^{\tilde{L}, \Endo}_{\tilde{R}}(f_{\tilde{L}}). \]
\end{lemma}
\begin{proof}
	Given Lemma \ref{prop:cr-descent} and its stable counterpart in Lemma \ref{prop:crS-descent}, the proof follows the familiar combinatorial apparatus: see Proposition \ref{prop:descent-orbint-Endo}. 
\end{proof}

\begin{proposition}[Cf.\ {\cite[IX.6.6 (1)]{MW16-2}}]\label{prop:cr-vanishing-Endo}
	Suppose that $M \neq G$. For all $\pi \in D_{\elli, -}(\tilde{M}) \otimes \mes(M)^\vee$ we have
	\[ \sum_{S \in \mathcal{P}(M)} \omega_S(X) I^{\tilde{M}}\left(\pi, \nu_S, X, {}^{c,r} \theta^{\Endo}_{\tilde{M}}(f) \right) = 0 \]
	for all $f$ and all $X$, provided that $\nu_S \gg 0$ relative to $S$.
\end{proposition}
\begin{proof}
	Every such $\pi$ can be written as $\sum_{\mathbf{M}^!, \pi^!} \trans_{\mathbf{M}^!, \tilde{M}}(\pi^!)$, where $\mathbf{M}^! \in \Endo_{\elli}(\tilde{M})$ and $\pi^!$ comes from some $\phi^! \in \Phi_{2, \mathrm{bdd}}(M^!)$. Indeed, this fact is already contained in the discussion in the proof of Theorem \ref{prop:local-character-relation}, which hinges ultimately upon the spectral inversion from \cite[\S 7]{Li19}.
	
	Hence we may assume $\pi = \trans_{\mathbf{M}^!, \tilde{M}}(\pi^!)$. Now it suffices to use the definition of ${}^{c,r} \theta^{\Endo}_{\tilde{M}}$ to reduce to the stable side, i.e.\ to Lemma \ref{prop:Scr-vanishing}.
\end{proof}

Concerning ${}^c \theta^{\Endo}_{\tilde{M}}(f)$, here is an endoscopic counterpart of Proposition \ref{prop:cthetaGM-arch}.

\begin{proposition}\label{prop:cthetaGM-Endo-arch}
	Let $f \in \orbI_{\asp}(\tilde{G}, \tilde{K}) \otimes \mes(G(F))$ and $M \in \mathcal{L}(M_0)$. Then
	\begin{itemize}
		\item ${}^c \theta^{\Endo}_{\tilde{M}}(f)$ is Schwartz on $\tilde{M}$;
		\item for all $\pi \in D_{\mathrm{temp}, -}(\tilde{M}) \otimes \mes(M(F))^\vee$, the function $\lambda \mapsto I^{\tilde{M}}(\pi, \lambda, {}^c \theta^{\Endo}_{\tilde{M}}(f))$ on $i\mathfrak{a}_M^*$ extends to a meromorphic function on $\mathfrak{a}^*_{M, \CC}$, with poles on hyperplanes of the form $\lrangle{\cdot, \check{\alpha}} = \text{const}$ for $\alpha \in \Sigma^G(A_M)$, and rapidly decreasing over vertical strips.
	\end{itemize}
\end{proposition}
\begin{proof}
	It is routine to reduce to the case of ${}^c \theta^{\Endo}_{\tilde{M}}(\mathbf{M}^!, f)$ for various $\mathbf{M}^! \in \Endo_{\elli}(\tilde{M})$, and use the available analogues on the stable side. Besides the formal manipulations, the key inputs in the Archimedean case are
	\begin{itemize}
		\item the spectral inversion, see \S\ref{sec:spectral-transfer} and \cite[\S 7]{Li19}, which expresses $\pi$ as a transfer from various $\mathbf{M}^!$;
		\item the $\tilde{K}$-finite transfer, see Theorem \ref{prop:image-transfer}.
	\end{itemize}

	See also the sketch in \cite[p.1095]{MW16-2}; a detailed discussion in the non-Archimedean context can be found in \cite[VIII.3.6]{MW16-2}.
\end{proof}

Without further ado, we identify $\mathfrak{a}^*_{M, \CC}$ and $\mathfrak{a}^*_{M^!, \CC}$ via $\xi: A_M \rightiso A_{M^!}$.

\begin{definition}\label{def:rhoGM-Endo-arch}
	\index{rhoGM-Endo@$\rho^{\tilde{G}, \Endo}_{\tilde{M}}$}
	Fix $\mathbf{M}^! \in \Endo_{\elli}(\tilde{M})$ and let $\pi^! \in SD_{\mathrm{temp}}(M^!)_0 \otimes \mes(M^!)^\vee$. Define
	\[ \rho^{\tilde{G}, \Endo}_{\tilde{M}}(\mathbf{M}^!, \pi^!, \lambda) := \sum_{s \in \Endo_{\mathbf{M}^!}(\tilde{G})} i_{M^!}(\tilde{G}, G^![s]) \trans_{\mathbf{M}^!, \tilde{M}}\left( \sigma^{G^![s]}_{M^!}(\pi^!, \lambda) \right). \]
	It belongs to $U^G_M \otimes D_{\mathrm{temp},-}(\tilde{M})_0 \otimes \mes(M)^\vee$.
\end{definition}

\begin{lemma}[Cf.\ {\cite[p.1098]{MW16-2}}]\label{prop:rhoGMEndo-arch}
	There exists a unique linear map $\rho^{\tilde{G}, \Endo}_{\tilde{M}}: D_{\mathrm{temp},-}(\tilde{M})_0 \otimes \mes(M)^\vee \to U^G_M \otimes D_{\mathrm{temp},-}(\tilde{M})_0 \otimes \mes(M)^\vee$ such that
	\[ I^{\tilde{M}}\left( \pi, X, {}^r \theta^{\tilde{G}, \Endo}_{\tilde{M}}(f) \right) = \int_{i\mathfrak{a}^*_M} I^{\tilde{M}}\left( \rho^{\tilde{G}, \Endo}_{\tilde{M}}(\pi, \lambda)_\lambda , f_{\tilde{M}} \right) e^{-\lrangle{\lambda, X}} \dd\lambda \]
	for all $\pi \in D_{\mathrm{temp},-}(\tilde{M})_0 \otimes \mes(M)^\vee$ and $f$. When $\pi = \sum_{\mathbf{M}^!, \pi^!} \trans_{\mathbf{M}^!, \tilde{M}}(\pi^!)$ where $(\mathbf{M}^!, \pi^!)$ is as in Definition \ref{def:rhoGM-Endo-arch}, we have
	\begin{equation}\label{eqn:rhoGM-Endo-arch}
		\rho^{\tilde{G}, \Endo}_{\tilde{M}}(\pi) = \sum_{\mathbf{M}^!, \pi^!} \rho^{\tilde{G}, \Endo}_{\tilde{M}}(\mathbf{M}^!, \pi^!).
	\end{equation}
\end{lemma}
\begin{proof}
	The uniqueness is proved in the same way as Lemma \ref{prop:rhoGM-arch}. The uniqueness holds in the following sense: when $\pi$ is given, the exists at most one $\rho^{\tilde{G}, \Endo}_{\tilde{M}}(\pi)$ such that the equality holds for all $f$. Therefore, in order to show existence, it suffices to verify the equality for $\rho^{\tilde{G},\Endo}_{\tilde{M}}(\pi) := \rho^{\tilde{G}, \Endo}_{\tilde{M}}(\mathbf{M}^!, \pi^!)$ when $\pi = \trans_{\mathbf{M}^!, \tilde{M}}(\pi^!)$ for some given $(\mathbf{M}^!, \pi^!)$. It will follow that the expression \eqref{eqn:rhoGM-Endo-arch} is independent of the decomposition of $\pi$.

	Identifying $\mathfrak{a}_M$ with $\mathfrak{a}_{M^!}$, Lemma \ref{prop:sigmaGM-arch} implies
	\begin{multline*}
		I^{\tilde{M}}\left( \pi, X, {}^r \theta^{\tilde{G}, \Endo}_{\tilde{M}}(f) \right) = \sum_{s \in \Endo_{\mathbf{M}^!}(\tilde{G})} i_{M^!}(\tilde{G}, G^![s]) S^{M^!}\left(\pi^!, X, {}^r S\theta^{G^![s]}_{M^!}\left( f^{G^![s]} \right)[s] \right) \\
		= \sum_s i_{M^!}(\tilde{G}, G^![s]) \int_{i\mathfrak{a}_{M^!}^*} S^{M^!}\left( \sigma^{G^![s]}_{M^!}(\pi^!, \lambda)_\lambda, (f^{G^![s]})_{M^!}[s] \right) e^{-\lrangle{\lambda, X}} \dd\lambda \\
		= \sum_s i_{M^!}(\tilde{G}, G^![s]) \int_{i\mathfrak{a}_{M^!}^*} S^{M^!}\left( \sigma^{G^![s]}_{M^!}(\pi^!, \lambda)_\lambda, (f_{\tilde{M}})^{M^!} \right) e^{-\lrangle{\lambda, X}} \dd\lambda \\
		= \sum_s i_{M^!}(\tilde{G}, G^![s]) \int_{i\mathfrak{a}_M^*} I^{\tilde{M}}\left( \trans_{\mathbf{M}^!, \tilde{M}} \left(\sigma^{G^![s]}_{M^!}(\pi^!, \lambda)\right)_\lambda , f_{\tilde{M}} \right) e^{-\lrangle{\lambda, X}} \dd\lambda .
	\end{multline*}
	The last expression equals $\int_{i\mathfrak{a}_M^*} I^{\tilde{M}}\left( \rho^{\tilde{G}, \Endo}_{\tilde{M}}(\mathbf{M}^!, \pi^!, \lambda)_\lambda , f_{\tilde{M}} \right) e^{-\lrangle{\lambda, X}} \dd\lambda$.
\end{proof}

Take $\phi^! \in \Phi_{\mathrm{bdd}}(M^!)_0$ (i.e.\ with central character trivial over $A_{M^!}(F)^\circ$) and the corresponding stable character $\pi^! = \sum_{\sigma \in \Pi_{\phi^!}} m_\sigma \sigma$ as in \eqref{eqn:sigmaGM-stability-arch-aux0}. By \eqref{eqn:sigmaGM-stability-arch-aux} we have
\begin{equation}\label{eqn:rhoGM-spec-matching-arch}
	\rho^{\tilde{G}, \Endo}_{\tilde{M}}(\mathbf{M}^!, \pi^!, \lambda) = \sum_{s \in \Endo_{\mathbf{M}^!}(\tilde{G})} i_{M^!}(\tilde{G}, G^![s]) s^{G^![s]}_{M^!}(\phi^!_\lambda) \otimes \trans_{\mathbf{M}^!, \tilde{M}}(\pi^!).
\end{equation}

\begin{theorem}\label{prop:rhoGM-spec-matching-arch}
	We have $\rho^{\tilde{G}, \Endo}_{\tilde{M}}(\pi) = \rho^{\tilde{G}}_{\tilde{M}}(\pi)$ for all $\pi \in D_{\mathrm{temp},-}(\tilde{M})_0 \otimes \mes(M)^\vee$.
\end{theorem}
\begin{proof}
	Both sides are linear in $\pi$. Let us fix $\mathbf{M}^! \in \Endo_{\elli}(\tilde{M})$ and assume $\pi = \trans_{\mathbf{M}^!, \tilde{M}}(\pi^!)$ where $\pi^! \in SD_{\mathrm{temp},-}(M^!)_0 \otimes \mes(M^!)^\vee$ is the stable character attached to $\phi^! \in \Phi_{\mathrm{bdd}}(M^!)_0$. Write $\pi$ as a linear combination of irreducibles, and take any $\pi_0 \in \Pi_{\mathrm{temp},-}(\tilde{M})$ appearing with nonzero coefficient. Then the infinitesimal character of $\pi$ matches that of $\phi^!$. Recall that $r^{\tilde{G}}_{\tilde{M}}(\pi_0)$ is independent of the choice of $\pi_0$.

	Using \eqref{eqn:rhoGM-spec-matching-arch} together with Theorem \ref{prop:easy-stabilization}, we obtain
	\begin{align*}
		\rho^{\tilde{G}, \Endo}_{\tilde{M}}(\pi, \lambda) & = \rho^{\tilde{G}, \Endo}_{\tilde{M}}(\mathbf{M}^!, \pi^!, \lambda) \\
		& = \sum_{s \in \Endo_{\mathbf{M}^!}(\tilde{G})} i_{M^!}(\tilde{G}, G^![s]) s^{G^![s]}_{M^!}(\phi^!_\lambda) \otimes \trans_{\mathbf{M}^!, \tilde{M}}(\pi^!) \\
		& = r^{\tilde{G}}_{\tilde{M}}(\pi_{0, \lambda}) \otimes \trans_{\mathbf{M}^!, \tilde{M}}(\pi^!) = r^{\tilde{G}}_{\tilde{M}}(\pi_{0, \lambda}) \otimes \pi .
	\end{align*}
	By adjoint spectral transfer \cite[Corollary 7.5.4]{Li19}, it follows that $\rho^{\tilde{G}, \Endo}_{\tilde{M}}(\pi, \lambda) = r^{\tilde{G}}_{\tilde{M}}(\pi_\lambda) \otimes \pi$ for all $\pi \in \Pi_{\mathrm{temp},-}(\tilde{M})_0$. This equals $\rho^{\tilde{G}}_{\tilde{M}}(\pi, \lambda)$ by its construction in Lemma \ref{prop:rhoGM-arch}.
\end{proof}

Theorem \ref{prop:rhoGM-spec-matching-arch} is the analogue of \cite[IX.6.9 Lemme]{MW16-2}, which is proved through a long local-global loop. The proof above is much simpler since we have direct control of the normalizing factors in question.

\begin{definition}
	\index{IGM-c-Endo-shrek}
	For $\mathbf{M}^! \in \Endo_{\elli}(\tilde{M})$, $\delta \in SD_{\mathrm{orb}, \tilde{G}\text{-reg}}(M^!) \otimes \mes(M^!)^\vee$ and $f \in \orbI_{\asp}(\tilde{G}, \tilde{K}) \otimes \mes(G)$, set
	\begin{align*}
		{}^c I^{\Endo}_{\tilde{M}}(\mathbf{M}^!, \delta, f) & = {}^c I^{\tilde{G}, \Endo}_{\tilde{M}}(\mathbf{M}^!, \delta, f) \\
		& := \sum_{s \in \Endo_{\mathbf{M}^!}(\tilde{G})} i_{M^!}(\tilde{G}, G^![s]) {}^c S^{G^![s]}_{M^!}\left(\delta[s], f^{G^![s]} \right).
	\end{align*}
\end{definition}

\begin{definition-proposition}
	\index{IGM-c-Endo}
	Let $\tilde{\gamma} = \sum_{\mathbf{M}^!, \delta} \trans_{\mathbf{M}^!, \tilde{M}}(\delta)$ where the $(\mathbf{M}^!, \delta)$ are as above. Then
	\[ {}^c I^{\Endo}_{\tilde{M}}(\tilde{\gamma}, f) = {}^c I^{\tilde{G}, \Endo}_{\tilde{M}}(\tilde{\gamma}, f) := \sum_{\mathbf{M}^!, \delta} {}^c I^{\tilde{G}, \Endo}_{\tilde{M}}(\mathbf{M}^!, \delta, f) \]
	depends only on $\tilde{\gamma}$ and $f$. This defines ${}^c I^{\Endo}_{\tilde{M}}(\tilde{\gamma}, f)$ for all $\tilde{\gamma} \in D_{\mathrm{orb}, G\text{-reg}, -}(\tilde{M}) \otimes \mes(M)^\vee$ and $f \in \orbI_{\asp}(\tilde{G}, \tilde{K}) \otimes \mes(G)$.
\end{definition-proposition}
\begin{proof}
	Same as Definition--Proposition \ref{def:I-geom-Endo}: it is based on the description of the kernel of the collective transfer
	\[ \bigoplus_{\mathbf{M}^! \in \Endo_{\elli}(\tilde{M})} SD_{\mathrm{geom}}(M^!, \mathcal{O}_{M^!}) \otimes \mes(M^!)^\vee \to D_{\mathrm{geom}, -}(\tilde{M}, \mathcal{O}) \otimes \mes(M)^\vee \]
	in Theorem \ref{prop:Dgeom-preservation}, where $\mathcal{O}$ is a chosen stable semisimple class in $M(F)$, and $\mathcal{O}_{M^!}$ denotes its preimage.
\end{proof}

\begin{remark}
	Following the familiar paradigm, cf.\ Proposition \ref{prop:descent-orbint-Endo} and its proof, one establishes a descent formula for ${}^c I^{\Endo}_{\tilde{M}}\left(\tilde{\gamma}^{\tilde{M}}, \cdot \right)$ that is formally the same as that of ${}^c I_{\tilde{M}}\left( \tilde{\gamma}^{\tilde{M}}, \cdot\right)$.
\end{remark}

The following result is the Archimedean analogue of Proposition \ref{prop:cpt-noncpt}.
\begin{proposition}\label{prop:cpt-noncpt-Endo-arch}
	We have ${}^c I^{\Endo}_{\tilde{M}}\left(\tilde{\gamma}, \cdot \right) = \displaystyle\sum_{L \in \mathcal{L}(M)} I^{\tilde{L}, \Endo}_{\tilde{M}}\left(\tilde{\gamma}, {}^c \theta^{\Endo}_{\tilde{L}}(\cdot) \right)$ for all $\tilde{\gamma} \in D_{\mathrm{orb}, G\text{-reg}, -}(\tilde{M}) \otimes \mes(M)^\vee$.
\end{proposition}
\begin{proof}
	In view of its stable counterpart in Lemma \ref{prop:cpt-noncpt-S-arch}, the arguments for Proposition \ref{prop:cpt-noncpt} carry over verbatim.
\end{proof}

\section{Statement of matching theorems}\label{sec:matching-theta-arch}
Let $\tilde{G}$ and $M \in \mathcal{L}(M_0)$ be as in \S\ref{sec:cI-Endo-arch}. The matching between ${}^r \theta^{\tilde{G}, \Endo}_{\tilde{M}}$ and ${}^r \theta^{\tilde{G}}_{\tilde{M}}$ is now within reach.

\begin{theorem}\label{prop:rtheta-matching-arch}
	We have ${}^r \theta^{\tilde{G}, \Endo}_{\tilde{M}} = {}^r \theta^{\tilde{G}}_{\tilde{M}}$.
\end{theorem}
\begin{proof}
	In view of Lemmas \ref{prop:rhoGM-arch} and \ref{prop:rhoGMEndo-arch}, the equality follows at once from $\rho^{\tilde{G}, \Endo}_{\tilde{M}} = \rho^{\tilde{G}}_{\tilde{M}}$ (Theorem \ref{prop:rhoGM-spec-matching-arch}).
\end{proof}

The proof of the results below will be completed in \S\ref{sec:end-stabilization}, together with all the other main theorems. Their non-Archimedean counterpart is Theorem \ref{prop:cpt-supported-equalities}.

\begin{theorem}\label{prop:theta-matching-arch}
	We have ${}^c \theta^{\tilde{G}, \Endo}_{\tilde{M}} = {}^c \theta^{\tilde{G}}_{\tilde{M}}$ and ${}^{c,r} \theta^{\tilde{G}, \Endo}_{\tilde{M}} = {}^{c,r} \theta^{\tilde{G}}_{\tilde{M}}$. Furthermore, for all $\tilde{\gamma} \in D_{\mathrm{orb}, G\text{-reg}, -}(\tilde{M}) \otimes \mes(M)^\vee$ and $f \in \orbI_{\asp}(\tilde{G}, \tilde{K}) \otimes \mes(G)$, we have
	\[ {}^c I^{\Endo}_{\tilde{M}}(\tilde{\gamma}, f) = {}^c I_{\tilde{M}}(\tilde{\gamma}, f). \]
\end{theorem}

We will give a conditional proof of these assertions under the Hypothesis \ref{hyp:ext-Arch}, that is, under the assumption that $I^{\tilde{G}, \Endo}_{\tilde{M}}(\tilde{\gamma}, f) = I^{\tilde{G}}_{\tilde{M}}(\tilde{\gamma}, f)$ whenever $\tilde{\gamma} \in D_{\mathrm{orb}, G\text{-reg}, -}(\tilde{M}) \otimes \mes(M)^\vee$ and $f \in \orbI_{\asp}(\tilde{G}, \tilde{K}) \otimes \mes(G)$.

For each $\star \in \left\{ r, (c,r), c \right\}$, define
\[ {}^\star \varphi(f) := {}^\star \theta^{\tilde{G}, \Endo}_{\tilde{M}}(f) - {}^\star \theta^{\tilde{G}}_{\tilde{M}}(f) \]
for each $f$. Note that ${}^r \varphi = 0$ is already established in Theorem \ref{prop:rtheta-matching-arch}.

\begin{lemma}\label{prop:theta-matching-arch-prep}
	Assume that Hypothesis \ref{hyp:ext-Arch} holds for $\tilde{M} \subset \tilde{G}$, and that all the main theorems hold when $M$ (resp.\ $G$) is replaced by any $L \in \mathcal{L}(M)$ with $L \neq M$ (resp.\ $L \neq G$). Then
	\[ {}^c \varphi = {}^{c,r} \varphi = 0, \quad {}^c I^{\Endo}_{\tilde{M}}(\tilde{\gamma}, \cdot) = {}^c I_{\tilde{M}}(\tilde{\gamma}, \cdot) \]
	for all $\tilde{\gamma} \in D_{\mathrm{orb}, G\text{-reg}, -}(\tilde{M}) \otimes \mes(M)^\vee$.
\end{lemma}
\begin{proof}
	Let $\tilde{\gamma} \in D_{\mathrm{orb}, G\text{-reg}, -}(\tilde{M}) \otimes \mes(M)^\vee$ and $f \in \orbI_{\asp}(\tilde{G}, \tilde{K}) \otimes \mes(G)$. Apply \eqref{eqn:cpt-noncpt-arch}, Proposition \ref{prop:cpt-noncpt-Endo-arch}, Hypothesis \ref{hyp:ext-Arch} and the induction hypotheses to obtain
	\begin{align*}
		{}^c I^{\Endo}_{\tilde{M}}(\tilde{\gamma}, f) - {}^c I_{\tilde{M}}(\tilde{\gamma}, f) & = \sum_{L \in \mathcal{L}(M)} \left( I^{\tilde{L}, \Endo}_{\tilde{M}}\left(\tilde{\gamma}, {}^c \theta^{\tilde{G}, \Endo}_{\tilde{L}}(f) \right) - I^{\tilde{L}}_{\tilde{M}}\left(\tilde{\gamma}, {}^c \theta^{\tilde{G}}_{\tilde{L}}(f) \right) \right) \\
		& = I^{\Endo}_{\tilde{M}}(\tilde{\gamma}, f) - I_{\tilde{M}}(\tilde{\gamma}, f) + I^{\tilde{M}}\left(\tilde{\gamma}, {}^c \theta^{\tilde{G}, \Endo}_{\tilde{M}}(f) - {}^c \theta^{\tilde{G}}_{\tilde{L}}(f) \right) \\
		& = I^{\tilde{M}}\left(\tilde{\gamma}, {}^c \varphi(f) \right).
	\end{align*}

	As ${}^c \varphi(f) \in \orbI_{\asp}(\tilde{M}, \tilde{K} \cap \tilde{M}) \otimes \mes(M)$, Lemmas \ref{prop:c-rcr}, \ref{prop:c-rcr-Endo} and the inductive hypotheses imply
	\begin{multline*}
		{}^c \varphi = \sum_{L \in \mathcal{L}(M)} \left( {}^r \theta^{\tilde{L}, \Endo}_{\tilde{M}} \circ {}^{c,r} \theta^{\tilde{G}, \Endo}_{\tilde{L}} - {}^r \theta^{\tilde{L}}_{\tilde{M}} \circ {}^{c,r} \theta^{\tilde{G}}_{\tilde{L}} \right) \\
		= \left( {}^r \theta^{\tilde{G}, \Endo}_{\tilde{M}} - {}^r \theta^{\tilde{G}}_{\tilde{M}} \right) + \left( {}^{c,r} \theta^{\tilde{G}, \Endo}_{\tilde{M}} - {}^{c,r} \theta^{\tilde{G}}_{\tilde{M}} \right)
		= {}^r \varphi + {}^{c,r} \varphi = {}^{c,r} \varphi.
	\end{multline*}
	Hence ${}^{c,r} \varphi(f) \in \orbI_{\asp}(\tilde{M}, \tilde{K} \cap \tilde{M}) \otimes \mes(M)$ as well.
	
	For all $\pi \in D_{\mathrm{ell},-}(\tilde{M})_0 \otimes \mes(M)^\vee$ and $\lambda \in i\mathfrak{a}_M^*$, we infer that
	\[ I^{\tilde{M}}\left(\pi_\lambda, {}^c \varphi(f) \right) = I^{\tilde{M}}\left( \pi, \lambda, {}^{c,r} \varphi(f) \right) \]
	and it has holomorphic continuation to all $\lambda \in \mathfrak{a}_{M, \CC}^*$. Its Fourier transform
	\[ I^{\tilde{M}}\left( \pi, \nu, X, {}^{c,r} \varphi(f) \right) :=  \int_{\nu + i\mathfrak{a}_M^*} I^{\tilde{M}}\left( \pi, \lambda, {}^{c,r} \varphi(f) \right) e^{-\lrangle{\lambda, X}} \dd \lambda, \quad X \in \mathfrak{a}_M \]
	is thus independent of the shifting parameter $\nu$; taking $\nu = 0$ yields $I^{\tilde{M}}(\pi, X, {}^{c,r} \varphi(f))$ by Fourier inversion.
	
	Combining Lemmas \ref{prop:cr-vanishing} and \ref{prop:cr-vanishing-Endo}, we have
	\begin{align*}
		0 & = \sum_{S \in \mathcal{P}(M)} \omega_S(X) I^{\tilde{M}}\left( \pi, \nu_S, X, {}^{c,r} \varphi(f) \right) \\
		& = I^{\tilde{M}}\left(\pi, X, {}^{c,r} \varphi(f) \right) \sum_{S \in \mathcal{P}(M)} \omega_S(X) = I^{\tilde{M}}\left(\pi, X, {}^{c,r} \varphi(f) \right)
	\end{align*}
	provided that $\nu_S$ is sufficiently positive relative to $S$, for each $S \in \mathcal{P}(M)$. We conclude that
	\begin{equation}\label{eqn:theta-matching-cond-arch-aux}
		I^{\tilde{M}}\left( \pi, X, {}^{c,r} \varphi(f) \right) = 0, \quad \pi \in D_{\mathrm{ell},-}(\tilde{M})_0 \otimes \mes(M)^\vee, \; X \in \mathfrak{a}_M.
	\end{equation}
	
	Since ${}^{c,r} \theta^{\tilde{G}, \Endo}_{\tilde{M}}$ and ${}^{c,r} \theta^{\tilde{G}}_{\tilde{M}}$ satisfy the same descent formula (Lemmas \ref{prop:cr-descent}, \ref{prop:cr-descent-Endo}), the inductive hypotheses imply that ${}^{c,r} \varphi(f)$ is cuspidal. Combined with \eqref{eqn:theta-matching-cond-arch-aux} for various $\pi$ and $X$, we conclude ${}^{c,r} \varphi(f) = 0$.
	
	Consequently, ${}^c \varphi(f) = 0$. By the first step of the proof, we also have ${}^c I^{\Endo}_{\tilde{M}}(\tilde{\gamma}, f) - {}^c I_{\tilde{M}}(\tilde{\gamma}, f) = 0$. This concludes the proof.
\end{proof}

The non-Archimedean counterpart of the above is Corollary \ref{prop:cpt-supported-equalities-aux}. The observation below will be used in \S\ref{sec:construction-epsilonM-arch}.

\begin{lemma}\label{prop:cphi-cuspidal}
	Let $\star \in \{c, (c, r)\}$. For all $f \in \orbI_{\asp}(\tilde{G}, \tilde{K}) \otimes \mes(G)$, the function ${}^\star \varphi(f)$ is cuspidal.
\end{lemma}
\begin{proof}
	By induction, we may assume ${}^\star \theta^{\tilde{L}}_{\tilde{R}} = {}^\star \theta^{\tilde{L}, \Endo}_{\tilde{R}}$ whenever $R \in \mathcal{L}^M(M_0)$ and $L \in \mathcal{L}^G(R) \smallsetminus \{G\}$. It remains to compare the descent formulas for ${}^\star \theta^{\tilde{G}}_{\tilde{M}}$ (Lemma \ref{prop:cr-descent}, noting that $f_{\tilde{Q}}$ can now be replaced by $f_{\tilde{L}}$) and ${}^\star \theta^{\tilde{G}, \Endo}_{\tilde{M}}$ (Lemma \ref{prop:cr-descent-Endo}): indeed, if $R \subsetneq M$ and $d^G_R(M, L) \neq 0$, then $L \subsetneq G$.
\end{proof}
\chapter{Construction of \texorpdfstring{$\epsilon_{\tilde{M}}$}{epsilonM} in the Archimedean case}\label{sec:epsilon-arch}
The convention here is the same as \S\ref{sec:weighted-char-arch}; in particular, $F$ is Archimedean. Given $M \in \mathcal{L}(M_0)$ and letting $K^M := K \cap M(F)$, our goal is to construct the map
\[ \epsilon_{\tilde{M}}: \orbI_{\asp}(\tilde{G}, \tilde{K}) \otimes \mes(G) \to \orbI_{\mathrm{ac}, \asp, \cusp}(\tilde{M}, \widetilde{K^M}) \otimes \mes(M) \]
based on the results from \S\ref{sec:weighted-char-arch} and inductive assumptions, such that
\[ I^{\Endo}_{\tilde{M}}(\tilde{\gamma}, f) - I_{\tilde{M}}(\tilde{\gamma}, f) = I^{\tilde{M}}\left( \tilde{\gamma}, \epsilon_{\tilde{M}}(f) \right) \]
for all $\tilde{\gamma} \in D_{\mathrm{orb}, G\text{-reg}, -}(\tilde{M}) \otimes \mes(M)^\vee$.

The basic idea is similar to the non-Archimedean case in \S\ref{sec:epsilon-nonarch}. The construction here also breaks into local and global parts. Nonetheless, there are substantial differences for the Archimedean case. First, the use of Shalika germs in the local part is replaced by
\begin{itemize}
	\item stabilization of the differential equations satisfied by weighted orbital integrals (see \S\ref{sec:stable-DE}),
	\item stabilization of the jump relations of weighted orbital integrals (see \S\ref{sec:jump-relations}),
\end{itemize}
since these properties characterize invariant orbital integrals on $\tilde{M}$.

The stabilization of differential equations will be stated for $F \in \{\R, \CC\}$. For the stabilization of jump relations, we only consider $F = \R$ and elliptic maximal tori in $M$; specifically, the statement will involve the notion of semi-regular quadruplets (Definition \ref{def:semi-regular}).

In the global part of the construction, we need the distributions ${}^c I_{\tilde{M}}$ defined through normalized intertwining operators in \S\ref{sec:cphi-arch}, as well as their endoscopic analogues.

The overall strategy follows \cite{Ar00, Ar08} and \cite[IX]{MW16-2}. After constructing the map $\epsilon_{\tilde{M}}$, which a priori takes values in $\orbI_{\mathrm{ac}, \asp, \cusp}(\tilde{M}) \otimes \mes(M)$, we will show in \S\ref{sec:epsilon-KM-finite} that its image lies in the $\widetilde{K^M} \times \widetilde{K^M}$-finite part; the basic tool is the Fourier transform of weighted orbital integrals, cf.\ \cite{Ar94}.

\section{Derivatives of invariant weighted orbital integrals}\label{sec:derivative-orbint}
For each affine $F$-variety $V$, let $\mathcal{O}_V$ be the algebra of regular functions on $V_{\CC}$, and let $\mathrm{Diff}(V)$ be the algebra of algebraic differential operators on $V_{\CC}$. For every vector space $W$, denote by $\Sym V$ its symmetric algebra.

Fix $M \in \mathcal{L}(M_0)$ and a maximal $F$-torus $T \subset M$. There is an embedding $\Sym(\mathfrak{t}_{\CC}) \hookrightarrow \mathrm{Diff}(T)$ of algebras, denoted by $U \mapsto \partial_U$. Hence $\mathcal{O}_{T_{G\text{-reg}}} \cdot \Sym\mathfrak{t}$ embeds as a subalgebra of $\mathrm{Diff}(T_{G\text{-reg}})$. We have the Weyl groups
\[ W(M(F), T(F)) := N_{M(F)}(T(F)) / T(F), \quad W(M,T)(F) := (N_M(T)/T)(F). \]
\index{partialU@$\partial_U$}
\index{WMT@$W(M(F), T(F))$, $W(M, T)(F)$}

For convenience, fix a Haar measure on $T(F)$. For each $f \in \orbI_{\asp}(\tilde{G}) \otimes \mes(G)$, we have the anti-genuine $C^\infty$ function
\[ I_{\tilde{M}}(\cdot, f) = I^{\tilde{G}}_{\tilde{M}}(\cdot, f) : \widetilde{T_{G\text{-reg}}} \to \CC. \]
For each $D \in \mathrm{Diff}(T_{G\text{-reg}})$, it makes sense to study $D I_{\tilde{M}}(\cdot, f)$.

Given $\eta \in T(F)$ and $\tilde{\eta} \in \rev^{-1}(\eta)$, we fix a sufficiently small open neighborhood $\mathcal{U}^\flat \subset \mathfrak{g}_\eta(F)$ of $0$ as in the formalism of descent; see \S\ref{sec:descent-HA}. Take $\mathfrak{u} := \mathcal{U}^\flat \cap \mathfrak{t}(F)$. In particular, we may assume that $\mathfrak{u}$ has compact closure in $\mathfrak{t}(F)$, the exponential is defined on $\mathfrak{u}$, and $\exp(X)\eta \in G_{\mathrm{reg}}(F)$ if and only if $X \in \mathfrak{g}_{\mathrm{reg}}(F)$ for all $X \in \mathfrak{u}$.

\begin{proposition}[Cf.\ {\cite[Lemma 13.2]{Ar88LB} or \cite[IX.3.2]{MW16-2}}]\label{prop:estimate-derivatives-orbint}
	Let $\eta \in T(F)$, $\tilde{\eta} \in \rev^{-1}(\eta)$ and take $\mathfrak{u} \subset \mathfrak{t}(F)$ as above. For each $U \in \Sym(\mathfrak{t}_{\CC})$, there exists $N \in \Z_{\geq 1}$ such that for every $f \in \orbI_{\asp}(\tilde{G}) \otimes \mes(G)$ there is a constant $c = c(N,f) > 0$ satisfying
	\[ \left|\partial_U I_{\tilde{M}}\left( \exp(X)\tilde{\eta}, f \right)\right| \leq c \left| D^{G_\eta}(X) \right|^{-N} \]
	for all $X \in \mathfrak{u} \cap \mathfrak{t}_{G\text{-reg}}(F)$.
\end{proposition}
\begin{proof}
	The case of reductive groups is done in \cite{Ar88LB}. Since the proof is of analytic nature, involving no metaplectic features, we refer to \textit{loc.\ cit.} for details.
\end{proof}

Next, consider the endoscopic version from Definition--Proposition \ref{def:I-geom-Endo}:
\[ I^{\Endo}_{\tilde{M}}(\cdot, f) = I^{\tilde{G}, \Endo}_{\tilde{M}}(\cdot, f) : \widetilde{T_{G\text{-reg}}} \to \CC. \]
It is also anti-genuine, $C^\infty$, and one can consider its derivatives. Some preparations are in order.

Fix a stable conjugacy class $\mathcal{T}$ of embeddings of maximal tori $T \hookrightarrow M$. Then $M(F)$ acts on $\mathcal{T}$ by conjugation. Abusing notations, we will view each $T \in \mathcal{T}$ as a subgroup of $M$.

Given $T$, consider the abelian $2$-group
\index{DT@$\mathfrak{D}(T)$}
\begin{align*}
	\mathfrak{D}(T) & := \Hm^1(F, T) \\
	& = \Ker[\Hm^1(F, T) \to \Hm^1(F, M)] = \mathfrak{D}(T, M; F).
\end{align*}
Each $1$-cocycle can be expressed as $\sigma \mapsto y\sigma(y)^{-1}$ for some $y \in M(\CC)$. The corresponding class will be denoted by $[y] \in \mathfrak{D}(T)$. This induces a bijection
\begin{align*}
	\mathfrak{D}(T) & \xrightarrow{1:1} \mathcal{T}/M(F) =: \underline{\mathcal{T}} \\
	[y] & \mapsto y^{-1}Ty \;\mod \text{conj}.
\end{align*}

According to \cite[Définition 4.1.3]{Li15}, we have the canonical isomorphism
\begin{equation}\label{eqn:stable-conjugacy}
	\Ad(y^{-1}): \widetilde{T}_{M\text{-reg}} \rightiso (\widetilde{y^{-1}Ty})_{M\text{-reg}}, \quad \tilde{\gamma} \mapsto y^{-1}\tilde{\gamma} y
\end{equation}
which covers the stable conjugacy $\gamma \mapsto y^{-1}\gamma y$ on $M(F)$. It depends only on $[y]$.

\index{RT@$\mathfrak{R}(T)$}
Denote by $\mathfrak{R}(T)$ the Pontryagin dual of $\mathfrak{D}(T)$. When $T$ is elliptic in $M$, there is a bijection for each $\gamma \in T_{G\text{-reg}}(F)$:
\begin{equation*}\begin{aligned}
	\mathfrak{R}(T) & \xrightarrow{1:1} \left\{ (\mathbf{M}^!, \delta): \mathbf{M}^! \in \Endo_{\elli}(\tilde{M}), \; \delta \in \Sigma_{G\text{-reg}}(M^!), \; \gamma \leftrightarrow \delta \right\} \\
	& \quad \text{assuming}\; T \subset M: \text{elliptic}.
\end{aligned}\end{equation*}
Indeed, when $\tilde{M} = \tilde{G} = \Mp(W)$ this is \cite[Lemme 5.2.1]{Li15}, and the general case follows easily. More precisely, $\mathbf{M}^!$ arises from an element in $\check{T}$ associated with $\kappa \in \mathfrak{R}(T)$ by Tate--Nakayama duality. Once $\mathbf{M}^!$ is given, the stable class $\delta$ corresponding to $\gamma$ is uniquely determined.

The assignment $\kappa \mapsto (\mathbf{M}^!, \delta)$ can be extended to non-elliptic $T \subset M$ as follows. Take a Levi subgroup $R \subset M$ such that $T$ is elliptic in $R$, and produce the pair $(\mathbf{R}^!, \delta)$ from $\kappa$. What one needs is to realize $\mathbf{R}^!$ as a Levi of some $\mathbf{M}^! \in \Endo_{\elli}(\tilde{M})$; this we can do canonically by taking the unique $u \in \Endo_{\mathbf{R}^!}(\tilde{M})$ with $\delta[u] = \delta$, and set $\mathbf{M}^! := \mathbf{M}^![u]$. In this way, we obtain a map
\begin{equation}\label{eqn:R-vs-M-delta}\begin{aligned}
	\mathfrak{R}(T) & \to \left\{ (\mathbf{M}^!, \delta): \mathbf{M}^! \in \Endo_{\elli}(\tilde{M}), \; \delta \in \Sigma_{G\text{-reg}}(M^!), \; \gamma \leftrightarrow \delta \right\} \\
	\kappa & \mapsto (\mathbf{M}^!_\kappa, \delta_\kappa).
\end{aligned}\end{equation}
It is not bijective in general.

The main purpose of \eqref{eqn:R-vs-M-delta} is to deduce the inversion formula \eqref{eqn:IEndo-T-inversion}. In order to state it, we fix
\begin{itemize}
	\item a system of representatives $T \in \mathcal{T}$ for $\underline{\mathcal{T}}$;
	\item for each $T \in \underline{\mathcal{T}}$ and each class in $\mathfrak{D}(T)$, we fix $y \in M(\CC)$ such that $\sigma \mapsto y\sigma(y)^{-1}$ represents the cohomology class, and $y^{-1}Ty$ is another representative of $\underline{\mathcal{T}}$.
\end{itemize}
Different representatives $T$, $T'$ of $\underline{\mathcal{T}}$ are thus correlated by an isomorphism of tori $\iota_{T, T'}: T' \rightiso T$, arising from $\Ad(y^{-1})$ for the chosen $y$. We shall write
\[ \mathrm{inv}(T, T') := [y] \; \in \mathfrak{D}(T), \quad \text{if}\quad T' = y^{-1}Ty. \]

\begin{lemma}\label{prop:stable-DE-iota}
	The data above can be chosen so that $\iota_{T,T} = \identity$ and $\iota_{T, T'} \iota_{T', T''} = \iota_{T, T''}$ for all $T, T', T'' \in \underline{\mathcal{T}}$.
\end{lemma}
\begin{proof}
	Fix $T$ and take a Levi subgroup $R \subset M$ such that $T$ is elliptic in $R$. There is a group action by $W(R, T)(F)$ on $\mathcal{T}$. In fact, by \cite[Lemma 7.2.2]{Li19}, one can choose a subgroup $W'$ of $W(R, T)(F)$ that realizes all the stable conjugacy, and $|W'| = |\mathfrak{D}(T)|$. Choose representatives for elements of $W'$.
\end{proof}

For all $T \in \underline{\mathcal{T}}$ and $\kappa \in \mathfrak{R}(T)$, the map \eqref{eqn:R-vs-M-delta} yields $(\mathbf{M}^!, \delta) := (\mathbf{M}^!_\kappa , \delta_\kappa)$ and $T^! := M^!_\delta \subset M^!$, where we take a representative of $\delta$ in $M_{G\text{-reg}}(F)$. Given $\kappa$, we may and do fix a ``diagram'' joining $\delta \in T^! \subset M^!$ and $\gamma \in T \subset M$. By Definition \ref{def:diagram}, it affords
\begin{itemize}
	\item the isomorphism $\xi_{T^!, T}: T^! \rightiso T$ between $F$-tori, restricting to the canonical isomorphism $\xi: A_{M^!} \rightiso A_M$;
	\item the isomorphism $\tilde{\xi}_{T^!, T}: T^! \rightiso T$ between $F$-varieties, which is a twist of $\xi_{T^!, T}$ such that $\delta \leftrightarrow \gamma = \tilde{\xi}_{T^!, T}(\delta)$ with respect to $\mathbf{M}^! \in \Endo_{\elli}(\tilde{M})$.
\end{itemize}
Notice that $\xi_{T^!, T}$ and $\tilde{\xi}_{T^!, T}$ have the same differential at $1$.

Furthermore, the diagrams can be chosen to ensure the compatibilities
\[ \iota_{y^{-1}Ty, T} \xi_{T^!, T} = \xi_{T^!, y^{-1}Ty}, \quad \iota_{y^{-1}Ty, T} \tilde{\xi}_{T^!, T} = \tilde{\xi}_{T^!, y^{-1}Ty}, \]
where $\iota_{T,T'}$ are as in Lemma \ref{prop:stable-DE-iota}.

Let us revert to the study of $I^{\Endo}_{\tilde{M}}(\tilde{\gamma}, f)$. Fix the following data
\begin{itemize}
	\item $T \in \underline{\mathcal{T}}$;
	\item $\gamma \in T_{G\text{-reg}}(F)$, $\tilde{\gamma} \in \rev^{-1}(\gamma)$;
	\item for each $\kappa \in \mathfrak{R}(T)$, take the corresponding $(\mathbf{M}^!_\kappa, \delta_\kappa)$ and a ``diagram'' so that
	\[ \delta_\kappa = \tilde{\xi}_{T^!_\kappa, T}^{-1}(\gamma), \quad T^!_\kappa := (M^!_\kappa)_{\delta_\kappa}. \]
	\item a Haar measure of $T(F)$, and the corresponding Haar measures of $(y^{-1}Ty)(F)$ for each $[y]$, and of $T^!_\kappa(F)$ for each $\kappa$; they are compatible via $\iota_{T, y^{-1}Ty}$ and $\xi_{T^!_\kappa, T}$.
\end{itemize}

Given $\kappa$ and $[y]$, the cocycle condition \cite[p.51]{Li19} says
\[ \lrangle{\kappa, [y]} = \Delta^\kappa(\delta_\kappa, y^{-1}\tilde{\gamma} y) \Delta^\kappa(\delta_\kappa, \tilde{\gamma})^{-1} , \quad \Delta^\kappa := \Delta_{\mathbf{M}^!_\kappa, \tilde{M}}. \]

Identify $\tilde{\gamma}$ with an element of $D_{\mathrm{orb}, G\text{-reg}}(\tilde{M}) \otimes \mes(M)^\vee$, and similarly for various $\delta_\kappa$. Making the same identification for $y^{-1}\tilde{\gamma}y$ for every $[y] \in \mathfrak{D}(T)$, we obtain
\[ \trans^\kappa(\delta_\kappa) := \trans_{\mathbf{M}^!_\kappa, \tilde{M}}(\delta_\kappa) = \sum_{[y] \in \mathfrak{D}(T)} \Delta^\kappa(\delta_\kappa, y^{-1}\tilde{\gamma}y) y^{-1}\tilde{\gamma}y , \]
thus Fourier inversion on $\mathfrak{D}(T)$ gives
\begin{align*}
	\tilde{\gamma} & = |\mathfrak{D}(T)|^{-1} \sum_{\kappa \in \mathfrak{R}(T)} \sum_{[y] \in \mathfrak{D}(T)} \lrangle{\kappa, [y]} y^{-1}\tilde{\gamma} y \\
	& = \sum_{\kappa \in \mathfrak{R}(T)} |\mathfrak{D}(T)|^{-1} \Delta^\kappa(\delta_\kappa, \tilde{\gamma})^{-1} \trans^\kappa\left( \delta_\kappa \right) \\
	& = \sum_{\kappa \in \mathfrak{R}(T)} \Delta^\kappa(\tilde{\gamma}, \delta_\kappa) \trans^\kappa \left( \delta_\kappa \right),
\end{align*}
where $\Delta^\kappa(\tilde{\gamma}, \delta_\kappa)$ is the adjoint transfer factor in Definition \ref{def:adjoint-transfer-factor}.

In conclusion, we obtain
\begin{equation}\label{eqn:IEndo-T-inversion}
	I^{\Endo}_{\tilde{M}}(\tilde{\gamma}, f) = \sum_{\kappa \in \mathfrak{R}(T)} \Delta^\kappa(\tilde{\gamma}, \delta_\kappa) \sum_{s \in \Endo_{\mathbf{M}^!_\kappa}(\tilde{G})} i_{M^!_\kappa}(\tilde{G}, G^![s]) S^{G^![s]}_{M^!_\kappa}\left( \delta_\kappa[s], f^{G^![s]} \right).
\end{equation}

\begin{proposition}[Cf.\ {\cite[IX.3.4]{MW16-2}}]\label{prop:estimate-derivatives-orbint-Endo}
	Let $\eta$, $\tilde{\eta}$, $\mathfrak{u}$ be as in Proposition \ref{prop:estimate-derivatives-orbint}. For each $U \in \Sym(\mathfrak{t}_{\CC})$, there exists $N \in \Z_{\geq 1}$ such that for every $f \in \orbI_{\asp}(\tilde{G}) \otimes \mes(G)$ there is a constant $c = c(N,f) > 0$ satisfying
	\[ \left| \partial_U I^{\Endo}_{\tilde{M}}\left( \exp(X)\tilde{\eta}, f \right)\right| \leq c \left| D^{G_\eta}(X) \right|^{-N} \]
	for all $X \in \mathfrak{u} \cap \mathfrak{t}_{G\text{-reg}}(F)$.
\end{proposition}
\begin{proof}
	This is based on \eqref{eqn:IEndo-T-inversion}. For every $\kappa \in \mathfrak{R}(T)$, take
	\[ \epsilon_\kappa := \tilde{\xi}_{T^!_\kappa, T}^{-1}(\eta). \quad Y_\kappa := (\dd\xi_{T^!_\kappa, T})^{-1}(X) \in \mathfrak{t}^!_{\kappa, \tilde{G}\text{-reg}}(F). \]
	Denote by $\partial_{U^!}$ the transportation of $\partial_U$ to $T^!_{\kappa, \tilde{G}\text{-reg}}$ via $\xi_{T^!_\kappa, T}$, for each $\kappa$. Also note that $\Delta^\kappa(\delta_\kappa, \tilde{\gamma})$ is locally constant, hence commutes with these differential operators.
	
	Let $s \in \Endo_{\mathbf{M}^!_\kappa}(\tilde{G})$. By the stable version \cite[IX.3.3]{MW16-2} of Proposition \ref{prop:estimate-derivatives-orbint}, for each $\kappa$ we have
	\[ \partial_{U^!} S^{G^![s]}_{M^!_\kappa}\left( \exp(Y_\kappa)\epsilon_\kappa[s], f^{G^![s]} \right) \leq c^! \left| D^{G^![s]_{\epsilon_\kappa[s]}}(Y_\kappa) \right|^{-N} \]
	for some constant $c^! = c^!\left(\kappa, N, f^{G^![s]}\right) > 0$. By the descent of endoscopic data (see \S\ref{sec:descent-endoscopy}), we are in the situation
	\[\begin{tikzcd}[column sep=huge]
		G^![s]_{\epsilon_\kappa[s]} \arrow[dashed, dash, "\text{nonstd.\ endo.}", r] & \overline{G^![s]_{\epsilon_\kappa[s]}} \arrow[dashed, dash, "\text{endo.}", r] & G_\eta .
	\end{tikzcd}\]

	Non-standard endoscopy rescales roots, whilst the set of roots of $\overline{G^![s]_{\epsilon_\kappa[s]}}$ embeds into that of $G_\eta$, via the transfer of maximal tori. Since $X$ (hence $Y_\kappa$) is close to zero, we have $\left| D^{G_\eta}(X) \right| \leq b \left| D^{G^![s]_{\epsilon_\kappa[s]}}(Y_\kappa) \right|$ for some $b > 0$. Assembling these estimates, we obtain the estimate for $\partial_U I^{\Endo}_{\tilde{M}}(\exp(X)\tilde{\eta}, f)$.
\end{proof}

\section{Differential equations}
\subsection{The operators \texorpdfstring{$\delta^G_M$}{deltaGM}}
For each connected reductive $F$-group $L$, let $\mathfrak{Z}(\mathfrak{l})$ be the center of the universal enveloping algebra $\mathcal{U}(\mathfrak{l}_{\CC})$. The Harish-Chandra isomorphism gives a homomorphism of algebras
\begin{align*}
	\mathfrak{Z}(\mathfrak{l}) & \to \mathfrak{Z}(\mathfrak{t}) = \Sym(\mathfrak{t}_{\CC}) \\
	z & \mapsto z_T
\end{align*}
for each maximal torus $T \subset L$. For any Levi subgroup $M \subset L$, we also have a natural homomorphism
\begin{align*}
	\mathfrak{Z}(\mathfrak{l}) & \to \mathfrak{Z}(\mathfrak{m}) \\
	z & \mapsto z_M,
\end{align*}
with the compatibility $(z_M)_T = z_T$ provided that $T \subset M$. When $M=T$, it coincides with Harish-Chandra homomorphism, thus the notation is consistent.
\index{Z(L)@$\mathfrak{Z}(\mathfrak{l})$}

The $\mathfrak{Z}(\mathfrak{g})$-action on $C^\infty_{c,\asp}(\tilde{G})$ induces actions on $\orbI_{\asp}(\tilde{G})$ and $\orbI_{\asp}(\tilde{G}, \tilde{K})$, etc. Denote these actions by $(z,f) \mapsto zf$.

\begin{proposition}\label{prop:delta-invariant}
	\index{deltaGM@$\delta^G_M$}
	Let $T$ be a maximal torus in $M \in \mathcal{L}(M_0)$. There is a unique family of linear maps
	\[ \delta^L_M = \delta^L_{M, T}: \mathfrak{Z}(\mathfrak{l}) \to \mathrm{Diff}(T_{L\text{-reg}}), \]
	where $L \in \mathcal{L}(M)$, such that
	\begin{itemize}
		\item $\delta^L_M$ has image in $\mathcal{O}_{T_{L\text{-reg}}} \cdot \Sym(\mathfrak{t})$;
		\item $\delta^L_L(z) = z_T$;
		\item for all $L \in \mathcal{L}(M)$, by regarding $I_{\tilde{L}}(\cdot, f)$ as smooth functions over $\tilde{T}_{G\text{-reg}}$, we have
		\[ I_{\tilde{M}}(\cdot, zf) = \sum_{L \in \mathcal{L}(M)} \delta^L_M(z_L) I_{\tilde{L}}(\cdot, f) \]
		for all $z \in \mathfrak{Z}(\mathfrak{g})$ and $f \in \orbI_{\asp}(\tilde{G}) \otimes \mes(G)$.
	\end{itemize}
\end{proposition}
\begin{proof}
	The uniqueness of this family of maps follows inductively from
	\[ \delta^G_M(z) I^{\tilde{G}}(\cdot, f) = I_{\tilde{M}}(\cdot, zf) - \sum_{L \neq G} \delta^L_M(z_L) I_{\tilde{L}}(\cdot, f). \]
	
	As for the existence, one begins with the case where $I_{\tilde{L}}$ is replaced by $J_{\tilde{L}}$. In the case of reductive groups, the corresponding assertion is in \cite[\S 12]{Ar88LB} which is in turn based on \cite{Ar76}. Moreover, the results of \cite{Ar76} apply to coverings as well. The property that $\delta^G_M(z)$ lands in $\mathcal{O}_{T_{G\text{-reg}}} \cdot \Sym(\mathfrak{t})$ is explained in \cite[IX.1.5--1.6]{MW16-2}; the arguments concern only the adjoint action of $T$ on $\mathfrak{g}$. Thus the maps $\delta$ defined for $(M, G)$ and $(\tilde{M}, \tilde{G})$ are equal. This is why we write $\delta^L_M$ instead of $\delta^{\tilde{L}}_{\tilde{M}}$.
	
	The transition from $J_{\tilde{L}}$ to $I_{\tilde{L}}$ is explained in \cite{Ar00}. It boils down to the assertion that
	\[ \phi_{\tilde{L}}(zf) = z_L \phi_{\tilde{L}}(f). \]
	In turn, this reduces to standard properties of intertwining operators, and the arguments apply to coverings. Note that the property holds for both the canonically normalized $\phi_{\tilde{L}}$ and the $\phi^r_{\tilde{L}}$ made from normalized intertwining operators.
\end{proof}

We record some further properties of these maps.

\begin{itemize}
	\item $\delta^G_M(z)$ is $W(M(F), T(F))$-invariant. Indeed, its characterization is invariant.
	\item For all $z, z' \in \mathfrak{Z}(\mathfrak{g})$, by \cite[IX.1.2 (3)]{MW16-2} we have
	\begin{equation}\label{eqn:delta-zz}
		\delta^G_M(zz') = \sum_{L \in \mathcal{L}(M)} \delta^L_M(z_L) \delta^G_L(z').
	\end{equation}
	\item Fix a Borel subgroup $T_{\CC} \subset B \subset G_{\CC}$, and denote by $\Delta$ the resulting set of simple roots. By writing $\delta^G_M(z) = \sum_U q^G_M(U, z) \partial_U$, where $U$ ranges over a basis of $\Sym(\mathfrak{t})$ made up of homogeneous elements and $\partial_U$ stands for the corresponding differential operator, the coefficients $q^G_M(U,z) \in \mathcal{O}_{T_{G\text{-reg}}}$ satisfy
	\begin{equation}\label{eqn:delta-vanishing}
		q^G_M(U,z, \gamma) \to 0 \quad\text{as}\; \gamma \underset{\Delta}{\to} \infty
	\end{equation}
	whenever $M \neq G$, where $\gamma \underset{\Delta}{\to} \infty$ signifies that $|\alpha(\gamma)| \to \infty$ for all $\alpha \in \Delta$.
	\item According to \cite[p.286]{Ar88LB}, $\delta^G_M$ is the restriction of its avatar over $\CC$: it is determined by the triplet $(\mathfrak{g}_{\CC}, T_{\CC}, A_M)$. Therefore, $\delta^G_M$ is of an algebraic nature, and it is essentially sufficient to regard the case of split $T$. 
\end{itemize}

Let $\Omega \in \mathfrak{Z}(\mathfrak{g})$ be the Casimir operator with respect to a chosen invariant bilinear form. By \cite[IX.1.7, p.995]{MW16-2} we know that $C^G_M := \delta^G_M(\Omega)$ belongs to $\mathcal{O}_{T_{G\text{-reg}}}$.
\index{OmegaG@$\Omega_G$}

\subsection{The stable version}
Consider now a quasisplit $F$-group $G^!$ with chosen minimal Levi subgroup $M_0^!$ and maximal compact subgroup $K^!$ in good relative position, we have $\mathfrak{Z}(\mathfrak{g}^!)$ acting on $S\orbI(G^!)$, $S\orbI(G^!, K^!)$, etc. For all $M^! \in \mathcal{L}(M_0^!)$ and maximal tori $T^! \subset M^!$, we have a canonical family of linear maps
\[ S\delta^{L^!}_{M^!} = S\delta^{L^!}_{M^!, T^!}: \mathfrak{Z}(\mathfrak{l}^!) \to \mathrm{Diff}\left(T^!_{L^!\text{-reg}}\right) \]
where $L^! \in \mathcal{L}(M^!)$, such that
\begin{itemize}
	\item $S\delta^{L^!}_{M^!}$ has image in $\mathcal{O}_{T^!_{L^!\text{-reg}}} \cdot \Sym(\mathfrak{t}^!)$, and $S\delta^{L^!}_{M^!}(z)$ is $W(M^!, T^!)(F)$-invariant for all $z$;
	\item $S\delta^{L^!}_{L^!}(z) = z_{T^!}$;
	\item $S_{M^!}(\cdot, zf^!) = \sum_{L^! \in \mathcal{L}(M^!)} S\delta^{L^!}_{M^!}(z_{L^!}) S_{L^!}\left(\cdot, f^! \right)$ for all $f^! \in S\orbI(G^!) \otimes \mes(M^!)$, where both sides are regarded as smooth functions on $T^!_{G^!\text{-reg}}(F)$.
\end{itemize}

The properties of these maps are documented in \cite[IX.2.1--2.2]{MW16-2}, including the effect of the Casimir operator $\Omega^! \in \mathfrak{Z}(\mathfrak{g}^!)$. Namely, $S\delta^{G^!}_{M^!}(\Omega^!)$ is an explicit element $SC^{G^!}_{M^!} \in \mathcal{O}_{T^!_{G^!\text{-reg}}}$, as in \cite[p.1005]{MW16-2}.
\index{SdeltaGM@$S\delta^{G^{"!}}_{M^{"!}}$}

\subsection{The endoscopic version}
For all $\mathbf{G}^! \in \Endo_{\elli}(\tilde{G})$, there is a canonical homomorphism of algebras $\mathfrak{Z}(\mathfrak{g}) \to \mathfrak{Z}(\mathfrak{g}^!)$, denoted by $z \mapsto z^{G^!}$; see \cite[Proposition 7.3.1]{Li19}. It commutes with transfer $f \mapsto f^{G^!} := \Trans_{\mathbf{G}^!, \tilde{G}}(f)$, in the sense that
\begin{equation}\label{eqn:z-commute-transfer}
	(zf)^{G^!} = z^{G^!} f^{G^!}, \quad f \in \orbI_{\asp}(\tilde{G}) \otimes \mes(G), \quad z \in \mathfrak{Z}(\mathfrak{g}).
\end{equation}

The construction goes as follows. Take maximal tori $T^! \subset G^!$ and $T \subset G$ joined by a ``diagram'' as in \S\ref{sec:diagram}; this affords a Weyl-equivariant isomorphism $\xi_{T^!, T}: T^! \rightiso T$ between $F$-tori, and the resulting morphism between varieties
\[ \mathfrak{t^!}^*_{\CC} / W(G^!_{\CC}, T^!_{\CC}) \to \mathfrak{t}^*_{\CC} / W(G_{\CC}, T_{\CC}) \]
induces $z \mapsto z^{G^!}$ via Harish-Chandra's isomorphism. In \textit{loc.\ cit.} one presumes $T^!$ and $T$ are elliptic, but $z \mapsto z^{G^!}$ is readily seen to be independent of all choices. 

The following result is an immediate consequence of this construction.
\begin{lemma}\label{prop:z-transfer-descent}
	Let $M \subset G$ be a Levi subgroup, $\mathbf{M}^! \in \Endo_{\elli}(\tilde{M})$, $s \in \Endo_{\mathbf{M}^!}(\tilde{G})$ and $G^! := G^![s]$. We have
	\[ \left( z^{G^!} \right)_{M^!} = (z_M)^{M^!}, \quad z \in \mathfrak{Z}(\mathfrak{g}). \]
	It is thus legitimate the denote both sides by $z^{M^!}$.
\end{lemma}

Fix $M \in \mathcal{L}(M_0)$. Our goal is an endoscopic version of $\delta^G_M(z)$. In contrast with the case of reductive groups, the metaplectic setting requires a twist on the stable side.

\begin{definition}\label{def:diff-op-s-twist}
	Let $\mathbf{M}^! \in \Endo_{\elli}(\tilde{M})$ and let $T^! \subset M^!$ be a maximal $F$-torus. For each $s \in \Endo_{\mathbf{M}^!}(\tilde{G})$ and $D \in \mathrm{Diff}\left(T^!_{\tilde{G}\text{-reg}}\right)$, let $D[s] \in \mathrm{Diff}(T^!_{\tilde{G}\text{-reg}})$ be the transport of structure of $D$ via translation by $z[s]$; in particular, $(D_1 D_2)[s] = D_1[s] D_2[s]$. Two basic examples:
	\begin{itemize}
		\item if $D$ is a regular function, then $D[s]$ is the regular function $t \mapsto D(t[s])$;
		\item if $D$ arises from $\Sym \mathfrak{t}^!$, then $D[s] = D$.
	\end{itemize}
\end{definition}

Henceforth, we adopt the formalism of \S\ref{sec:derivative-orbint}.  In particular, we fix a stable class $\mathcal{T}$ of embeddings of maximal tori $T \hookrightarrow M$, and isomorphisms $\iota_{T, T'}: T' \rightiso T$ for all $T, T' \in \mathcal{T}$ which arise from stable conjugacy and are subject to the requirements of Lemma \ref{prop:stable-DE-iota}. We also fix Haar measures on all these tori compatibly.

We will consider matrices $A$ of algebraic differential operators $A_{T, T'}$ indexed by $T, T' \in \underline{\mathcal{T}}$. The precise meaning is: $A_{T, T'}$ maps functions over $T'_{G\text{-reg}}$ to functions over $T_{G\text{-reg}}$, and $A_{T, T'} \circ \iota_{T, T'}^*$ belongs to $\mathrm{Diff}(T_{G\text{-reg}})$.

\begin{definition-proposition}[Cf.\ \cite{Ar00} and {\cite[IX.2.3]{MW16-2}}]\label{def:deltaEndo}
	\index{deltaGM-Endo@$\delta^{\tilde{G}, \Endo}_{\tilde{M}}$}
	Given the choices above, for each $z \in \mathfrak{Z}(\mathfrak{g})$ there exists a unique matrix $\delta^{\tilde{G}, \Endo}_{\tilde{M}}(z)$ of differential operators, such that $\delta^{\tilde{G}, \Endo}_{\tilde{M}}$ is linear in $z$ and
	\[ I_{\tilde{M}}^{\Endo}(\tilde{\gamma}, zf) = \sum_{[y] \in \mathfrak{D}(T)} \sum_{L \in \mathcal{L}(M)} \delta^{\tilde{L}, \Endo}_{\tilde{M}}(z_L)_{T, y^{-1}Ty} I^{\Endo}_{\tilde{L}}\left( y^{-1}\tilde{\gamma}y, f \right) \]
	for all $T \in \underline{\mathcal{T}}, \tilde{\gamma} \in \tilde{T}_{G\text{-reg}}$ and $f \in \orbI_{\asp}(\tilde{G}) \otimes \mes(G)$.
\end{definition-proposition}
\begin{proof}
	The uniqueness follows by induction using the postulated identity. Indeed,
	\[ \sum_{[y] \in \mathfrak{D}(T)} \delta^{\tilde{G}, \Endo}_{\tilde{M}}(z)_{T, y^{-1}Ty} I^{\tilde{G}} \left( y^{-1}\tilde{\gamma} y, f \right) = I_{\tilde{M}}^{\Endo}(\tilde{\gamma}, zf) - \left(\text{terms with}\; L \neq G\right). \]
	Locally around a given $\tilde{\gamma}$, one can choose $f$ to isolate terms contributed by different $[y]$. This determines $\delta^{\tilde{G}, \Endo}_{\tilde{M}}(z)_{T, y^{-1}Ty}$.

	Now comes the existence. Given $T, T' \in \underline{\mathcal{T}}$, use the notation from \S\ref{sec:derivative-orbint} to define the operator
	\begin{multline}\label{eqn:deltaEndo}
		\delta^{\tilde{G}, \Endo}_{\tilde{M}}(z)_{T, T'} := |\mathfrak{D}(T)|^{-1} \sum_{\kappa \in \mathfrak{R}(T)} \lrangle{\kappa, \mathrm{inv}(T, T')} \sum_{s \in \Endo_{\mathbf{M}^!_\kappa}(\tilde{G})} i_{M^!_\kappa}(\tilde{G}, G^![s]) \\
		\times \; \tilde{\xi}_{T^!_\kappa, T}^{-1, *} \circ \underbracket{S\delta^{G^![s]}_{M^!_\kappa}\left( z^{G^![s]} \right)[s]}_{\text{Definition \ref{def:diff-op-s-twist}}} \circ \tilde{\xi}_{T^!_\kappa, T'}^* .
	\end{multline}
	In what follows, we verify that \eqref{eqn:deltaEndo} does the job.

	Fix $T \in \underline{\mathcal{T}}$. Let $\gamma \in T_{G\text{-reg}}(F)$ and $\tilde{\gamma} \in \rev^{-1}(\gamma)$. Apply the inversion formula \eqref{eqn:IEndo-T-inversion} to $zf$ and use \eqref{eqn:z-commute-transfer} to deduce that
	\begin{align*}
		I^{\Endo}_{\tilde{M}}(\tilde{\gamma}, zf) & = \sum_{\kappa \in \mathfrak{R}(T)} \Delta^\kappa(\tilde{\gamma}, \delta_\kappa) \sum_{s \in \Endo_{\mathbf{M}^!_\kappa}(\tilde{G})} i_{M^!_\kappa}(\tilde{G}, G^![s]) S^{G^![s]}_{M^!_\kappa}\left( \delta_\kappa[s], (zf)^{G^![s]} \right) \\
		& = \sum_{\kappa \in \mathfrak{R}(T)} \Delta^\kappa(\tilde{\gamma}, \delta_\kappa) \sum_{s \in \Endo_{\mathbf{M}^!_\kappa}(\tilde{G})} i_{M^!_\kappa}(\tilde{G}, G^![s]) \\
		& \quad \sum_{L^! \in \mathcal{L}^{G^![s]}(M^!_\kappa)} \left( S\delta^{L^!}_{M^!_\kappa}(z^{L^!}) S^{L^!}_{M^!_\kappa}\right) \left( \delta_\kappa[s], f^{G^![s]} \right).
	\end{align*}
	
	The next step is familiar. For each $\kappa$, apply Lemma \ref{prop:sL-Ls} with $R = M$ to transform the sum over $(s, L^!)$ into $(L, r^L, r_L)$ where $L \in \mathcal{L}^G(M)$, $r^L \in \Endo_{\mathbf{M}^!_\kappa}(\tilde{L})$, $r_L \in \Endo_{\mathbf{L}^![r^L]}(\tilde{G})$, with
	\[ G^![s] = G^![r_L], \quad L^! = L^![r^L], \quad z[s] = z[r^L] z[r_L]. \]
	We also apply Lemma \ref{prop:i-transitivity} to break the coefficients $i_{M^!_\kappa}(\tilde{G}, G^![s])$. The previous sum becomes the sum of the expression below over $\kappa \in \mathfrak{R}(T)$ (with $(\mathbf{M}^!, \delta) := (\mathbf{M}^!_\kappa, \delta_\kappa)$), which are functions in $\bm{x} \in T^!_{\tilde{G}\text{-reg}}(F)$ evaluated at $\bm{x} = \delta[r^L]$:
	\begin{multline*}
		\Delta_{\tilde{M}, \mathbf{M}^!}(\tilde{\gamma}, \delta) \sum_{L \in \mathcal{L}^G(M)} \sum_{r^L} i_{M^!}(\tilde{L}, L^![r^L]) \\
		\cdot S\delta^{L^![r^L]}_{M^!}(z^{L^![r^L]}) \left( \sum_{r_L} i_{L^![r^L]}(\tilde{G}, G^![r_L]) S^{G^![r_L]}_{L^![r^L]}\left( \bm{x}[r_L], f^{G^![r_L]} \right) \right) \\
		= \Delta_{\tilde{M}, \mathbf{M}^!}(\tilde{\gamma}, \delta) \sum_{L \in \mathcal{L}^G(M)} \sum_{r^L} i_{M^!}(\tilde{L}, L^![r^L])
		S\delta^{L^![r^L]}_{M^!}(z^{L^![r^L]}) I^{\Endo}_{\tilde{L}}\left( \mathbf{L}^![r^L], \bm{x}, f \right).
	\end{multline*}
	Here entered the obvious fact that $S\delta^{L^![r^L]}_{M^!}(z^{L^!})$ is invariant under translation by $z[r_L]$, as the latter is central in $L^![r^L]$.

	Consider $\kappa \in \mathfrak{R}(T)$, $\mathbf{M}^! = \mathbf{M}^!_\kappa$, $L \in \mathcal{L}^G(M)$ and $r^L \in \Endo_{\mathbf{M}^!}(\tilde{L})$. For any maximal torus $T^! \subset M^!$ joined to $T \subset M$ by a diagram, and any $\delta' \in T^!_{\tilde{G}\text{-reg}}(F)$, set $\delta := \delta' \cdot z[r^L]^{-1}$ and $\gamma := \tilde{\xi}_{T^!, T}(\delta)$ For all $[y] \in \mathfrak{D}(T)$, we have $\delta' \leftrightarrow y^{-1}\gamma y$ with respect to $\mathbf{L}^![r^L]$. Let $\tilde{\gamma} \in \rev^{-1}(\gamma)$. The parabolic descent of transfer factors \cite[p.52]{Li19} implies $\Delta_{\mathbf{L}^![r^L], \tilde{L}}(\delta', y^{-1}\tilde{\gamma}y) = \Delta_{\mathbf{M}^!, \tilde{M}}(\delta, y^{-1}\tilde{\gamma}y)$. Hence
	\begin{align*}
		I^{\Endo}_{\tilde{L}}\left( \mathbf{L}^![r^L], \delta', f \right) & = \sum_{[y]} \Delta_{\mathbf{L}^![r^L], \tilde{L}}(\delta', y^{-1}\tilde{\gamma}y) I^{\Endo}_{\tilde{L}}\left( y^{-1}\tilde{\gamma} y, f \right) \\
		& = \sum_{[y]} \Delta_{\mathbf{M}^!, \tilde{M}}(\delta, y^{-1}\tilde{\gamma}y) I^{\Endo}_{\tilde{L}}\left( y^{-1} \tilde{\gamma} y , f \right).
	\end{align*}

	The transfer factor is locally constant, thus commutes with differential operators. For all $L$, $y$ and $\kappa$, define the operator
	\[ D_{L,y,\kappa} := \sum_{t \in \Endo_{\mathbf{M}^!_\kappa}(\tilde{L})} i_{M^!_\kappa}(\tilde{L}, L^![t]) \tilde{\xi}_{T^!_\kappa, T}^{-1, *} \circ \underbracket{S\delta^{L^![t]}_{M^!_\kappa}\left(z^{L^!}\right)[t]}_{\text{Definition \ref{def:diff-op-s-twist}}} \circ \tilde{\xi}_{T^!_\kappa, y^{-1}Ty}^*. \]
	The reasoning above and Definition \ref{def:adjoint-transfer-factor} lead to
	\begin{multline*}
		\sum_L \sum_{\kappa \in \mathfrak{R}(T)} \sum_{[y] \in \mathfrak{D}(T)} \Delta^\kappa(\tilde{\gamma}, \delta) \Delta^\kappa(\delta, y^{-1}\tilde{\gamma} y) D_{L,y,\kappa} I^{\Endo}_{\tilde{L}} \left( y^{-1}\tilde{\gamma} y, f \right) \\
		= \sum_{[y] \in \mathfrak{D}(T)} \sum_L |\mathfrak{D}(T)|^{-1} \sum_{\kappa \in \mathfrak{R}(T)} \lrangle{\kappa, [y]} D_{L,y,\kappa} I^{\Endo}_{\tilde{L}}\left( y^{-1}\tilde{\gamma} y, f \right).
	\end{multline*}
	The left hand side equals $I^{\Endo}_{\tilde{M}}(\tilde{\gamma}, zf)$. In view of \eqref{eqn:deltaEndo} with $G$ replaced by various $L \in \mathcal{L}(M)$, this verifies the required identity.
\end{proof}

In view of Lemma \ref{prop:stable-DE-iota}, the multiplication of matrices of differential operators can be performed by the usual rule.

\begin{lemma}\label{prop:deltaEndo-zz}
	For all $z, z' \in \mathfrak{Z}(\mathfrak{g})$, we have
	\[ \delta^{\tilde{G}, \Endo}_{\tilde{M}}(zz') = \sum_{L \in \mathcal{L}(M)} \delta^{\tilde{L}, \Endo}_{\tilde{M}}(z_L) \delta^{\tilde{G}, \Endo}_{\tilde{L}}(z'). \]
\end{lemma}
\begin{proof}
	Same as the proof for $\delta^G_M$ cited earlier. It uses only Definition--Proposition \ref{def:deltaEndo} and induction.
\end{proof}

\section{Stabilization of differential equations}\label{sec:stable-DE}
Fix a stable class $\mathcal{T}$ of embeddings of maximal tori $T \hookrightarrow M$, as well as the auxiliary choices as in \S\ref{sec:derivative-orbint}. We view $\delta^G_M(z)$ as a scalar matrix of differential operators, namely $\delta^G_M(z)_{T, T'} := 0$ when $T \neq T'$, and $\delta^G_M(z)_{T, T} = \delta^G_{M,T}(z)$.

Our goal is the following statement.

\begin{theorem}\label{prop:stablization-DE}
	For all $M \in \mathcal{L}(M_0)$, we have $\delta^{\tilde{G}, \Endo}_{\tilde{M}} = \delta^G_M$.
\end{theorem}

We need several reduction steps.
\begin{lemma}\label{prop:stabilization-DE-G}
	We have $\delta^{\tilde{G}, \Endo}_{\tilde{G}} = \delta^G_G$.
\end{lemma}
\begin{proof}
	Let $T, T' \in \underline{\mathcal{T}}$. For all $z \in \mathfrak{Z}(\mathfrak{g})$, we have $\delta^G_G(z) = z_T \in (\Sym\mathfrak{t}_{\CC})^{W(G_{\CC}, T_{\CC})}$. On the other hand, \eqref{eqn:deltaEndo} with $M=G$ gives
	\[ \delta^{\tilde{G}, \Endo}_{\tilde{G}}(z)_{T, T'} = |\mathfrak{D}(T)|^{-1} \sum_{\substack{\kappa \in \mathfrak{R}(T) \\ \mathbf{G}^! := \mathbf{G}^!_\kappa}} \lrangle{\kappa, \mathrm{inv}(T, T')}
	\tilde{\xi}_{T^!, T}^{-1, *} \circ S\delta^{G^!}_{G^!}(z^{G^!}) \circ \tilde{\xi}_{T^!, T'}^* . \]
	
	Since $z \mapsto z^{G^!}$ is defined exactly by Harish-Chandra isomorphisms and the isomorphism between tori $\xi_{T^!, T}: T^! \rightiso T$, the above equals $z_T$ times the factor
	\[ |\mathfrak{D}(T)|^{-1} \sum_{\kappa \in \mathfrak{R}(T)} \lrangle{\kappa, \mathrm{inv}(T, T')}. \]
	The assertion follows at once.
\end{proof}

\begin{lemma}\label{prop:stabilization-DE-c}
	Let $\mathfrak{z}(\mathfrak{g})$ denote the center of $\mathfrak{g}_{\CC}$. For $z \in \Sym(\mathfrak{z}(\mathfrak{g}))$, we have
	\[ \delta^{\tilde{G}, \Endo}_{\tilde{M}}(z) = \delta^G_M(z) = \begin{cases}
		0, & M \neq G \\
		z, & M = G
	\end{cases} \]
	where $\mathfrak{z}(\mathfrak{g})$ embeds into $\mathfrak{z}(\mathfrak{t}) = \mathfrak{t}_{\CC}$ in the evident way, for each $T \in \underline{\mathcal{T}}$. Similarly, we have $S\delta^{G^!}_{M^!} = 0$ over $\Sym(\mathfrak{z}(\mathfrak{g}^!))$ when $M^! \subsetneq G^!$ on the endoscopic side.
\end{lemma}
\begin{proof}
	The descriptions of $\delta^{\tilde{G}, \Endo}_{\tilde{M}}(z)$ and $\delta^G_M(z)$ follow inductively from their characterizations.
\end{proof}

Fix $T \in \mathcal{T}$. Choose $m \in M(\CC)$ that conjugates $T_{\CC}$ to the split maximal torus $T_0 = M_0$ prescribed by the chosen symplectic basis. Transport the corresponding Weyl-invariant quadratic form on $X_*(T_0)$ (Definition \ref{def:invariant-quadratic-form}) to $X_*(T_{\CC})$. These quadratic forms are positive definite, and determine an invariant non-degenerate bilinear form on $\mathfrak{g}_{\CC}$. This well-defines the Casimir operator $\Omega_G \in \mathfrak{Z}(\mathfrak{g})$.

The same recipe applies to the elliptic endoscopic groups if the invariant quadratic forms prescribed by Definition \ref{def:invariant-quadratic-form} is used. This determines the corresponding Casimir operators $\Omega_{G^!}$ on the endoscopic side.

By this construction, if $\mathbf{M}^! \in \Endo_{\elli}(\tilde{M})$ and $T \subset M$ is joined to $T^! \subset M^!$ by a ``diagram'' in \S\ref{sec:diagram}, the resulting isomorphism $\xi_{T^!, T}: T^! \rightiso T$ matches their invariant quadratic forms. From the standard description of $\Omega_G$ and $\Omega_{G^!}$ in terms of quadratic forms, one infers that
\begin{equation}\label{eqn:Casimir-transfer}
	\left(\Omega_G\right)^{G^!} \in \Omega_{G^!} + \Sym\left(\mathfrak{z}(\mathfrak{g}^!)\right), \quad \mathbf{G}^! \in \Endo_{\elli}(\tilde{G}).
\end{equation}
In fact, they differ only by an explicit scalar in $\CC$.

\begin{lemma}\label{prop:stabilization-Casimir}
	We have $\delta^{\tilde{G}, \Endo}_{\tilde{M}}(z) = \delta^G_M(z)$ when $z \in \Omega_G + \Sym(\mathfrak{z}(\mathfrak{g}))$.
\end{lemma}
\begin{proof}
	The case $z \in \Sym(\mathfrak{z}(\mathfrak{g}))$ is Lemma \ref{prop:stabilization-DE-c}. The case $M=G$ is Lemma \ref{prop:stabilization-DE-G}. Hence it suffices to treat the case with
	\[ d := \dim \mathfrak{a}^G_M \geq 1, \quad z = \Omega_G . \]
	Note that $\dim \mathfrak{a}^{G^![s]}_{M^!} = d$ for all $\mathbf{M}^! \in \Endo_{\elli}(\tilde{M})$ and $s \in \Endo_{\mathbf{M}^!}(\tilde{G})$. Fix a diagram joining $T \subset M$ to $T^! \subset M^!$, whence the isomorphism $\xi_{T^!, T}$ of tori and its ``twisted'' version $\tilde{\xi}_{T^!, T}$.
	
	By \cite[IX.1.7, IX.2.2]{MW16-2}, there are explicit $C^{G}_M \in \mathcal{O}_{T_{G\text{-reg}}}$, $SC^{G^![s]}_{M^!} \in \mathcal{O}_{T^!_{G^![s]\text{-reg}}}$ such that
	\[ \delta^G_M(\Omega_G) = C^G_M, \quad S\delta^{G^![s]}_{M^!}(\Omega_{G^![s]}) = SC^{G^![s]}_{M^!}. \]
	Moreover, they both vanish when $d \geq 2$.

	Since $S\delta^{G^![s]}_{M^!}$ is zero on $\Sym(\mathfrak{z}(\mathfrak{g}^![s]))$ by Lemma \ref{prop:stabilization-DE-c}, in view of \eqref{eqn:deltaEndo} and \eqref{eqn:Casimir-transfer}, we are reduced to proving the equation of scalars
	\begin{equation*}
		|\mathfrak{D}(T)|^{-1} \sum_{\kappa \in \mathfrak{R}(T)} \lrangle{\kappa, \mathrm{inv}(T, T')} \sum_{s \in \Endo_{\mathbf{M}^!_\kappa}(\tilde{G})} i_{M^!_\kappa}(\tilde{G}, G^![s]) SC^{G^![s]}_{M^!_\kappa}(\delta_\kappa[s])
		= \begin{cases}
			C^G_M(\gamma), & T = T' \\
			0, & T \neq T'
		\end{cases}
	\end{equation*}
	when $d = 1$, for all $T, T' \in \underline{\mathcal{T}}$, $\gamma \in T_{G\text{-reg}}(F)$ and $\delta_\kappa := \tilde{\xi}_{T^!_\kappa, T}^{-1}(\gamma) \in T_{G\text{-reg}}(F)$.

	In turn, it suffices to prove that when $d=1$, for each $\kappa \in \mathfrak{R}(T)$ we have
	\begin{equation}\label{eqn:stable-DE-aux1}
		\sum_{s \in \Endo_{\mathbf{M}^!_\kappa}(\tilde{G})} i_{M^!_\kappa}(\tilde{G}, G^![s]) SC^{G^![s]}_{M^!_\kappa}(\delta_\kappa[s]) = C^G_M(\gamma).
	\end{equation}

	Since $\kappa$ is kept fixed, it will be omitted from the indices hereafter. Without loss of generality, we may assume $\tilde{G} = \Mp(W)$ and $d=1$, so that
	\begin{gather*}
		M = \GL(m) \times \Sp(W^\flat), \quad 2m + \dim W^\flat = \dim W, \quad m \in \Z_{\geq 1}, \\
		\mathbf{M}^! \;\text{corresponds to}\; (m', m''), \quad m' + m'' = \frac{1}{2}\dim W^\flat, \\
		M^! = \GL(m) \times \SO(2m' + 1) \times \SO(2m'' + 1).
	\end{gather*}
	For all $s \in \Endo_{\mathbf{M}^!}(\tilde{G})$, let $(n', n'')$ be the pair corresponding to $\mathbf{G}^![s] \in \Endo_{\elli}(\tilde{G})$, so that
	\[ G^![s] = \SO(2n'+1) \times \SO(2n''+1), \quad n' + n'' = \frac{1}{2} \dim W. \]
	
	The following combinatorial arguments echo the proof of Lemma \ref{prop:J-roots}.
	\begin{itemize}
		\item We may write $\Endo_{\mathbf{M}^!}(\tilde{G}) = \{+, -\}$, according to whether $\GL(m)$ goes into the factor $\SO(2n'+1)$ or $\SO(2n''+1)$ of $G^![s]$.

		\item Fix $s^\flat \in T^\vee$ with $(s^\flat)^2 = 1$ giving rise to $\mathbf{M}^!$. Given $s \in \Endo_{\mathbf{M}^!}(\tilde{G})$, denote by $\hat{s}$ the element of $T^\vee \subset \tilde{G}^\vee$ giving rise to $\mathbf{G}^![s]$, satisfying $\hat{s}^2 = 1$.

		\item For each $\alpha \in \Sigma(G,T)$, there is a corresponding $\check{\beta} \in \Sigma(\tilde{G}^\vee, T^\vee)$; it corresponds to a root $\beta$ in $G^![s]$ if and only if $\check{\beta}(\hat{s}) = 1$. In any case, $\beta|_{\mathfrak{a}_{M}^!}$ is independent of $s$.
		
		\item There exists $\check{\beta}_1 \in \Sigma(\check{G}^!, T^\vee) \smallsetminus \Sigma(\check{M}^!, T^\vee)$ such that for every $\check{\beta} \in \Sigma(\check{G}^!, T^\vee) \smallsetminus \Sigma(\check{M}^!, T^\vee)$, we have
		\[ \exists! \; k = k^{\widehat{G^!}}(\check{\beta}) \in \Z_{\geq 1}, \; \check{\beta}|_{A_{\check{M}^!}} = \pm k \check{\beta}_1 |_{A_{\check{M}^!}}. \]
		
		\item The choice of invariant positive-definite quadratic forms yields Euclidean norms $\|\cdot\|$ on $\mathfrak{a}^*_M$ and $\mathfrak{a}^*_{M^!}$, compatibly with the natural map $\xi: \mathfrak{a}_{M^!} \rightiso \mathfrak{a}_M$.
	\end{itemize}

	Comparing the explicit formulas for $C^G_M(\gamma)$ and $SC^{G^!}_{M^!}(\delta[s])$ in \cite[p.995, p.1005]{MW16-2}, the desideratum \eqref{eqn:stable-DE-aux1} is reduced to the following: for each $\alpha \in \Sigma(G,T) \smallsetminus \Sigma(M,T)$ and the corresponding $\check{\beta}$, we have
	\begin{multline}\label{eqn:stable-DE-aux2}
		\sum_{\substack{s \in \Endo_{\mathbf{M}^!}(\tilde{G}) \\ \check{\beta}(\hat{s}) = 1}} \frac{i_{M^!}(\tilde{G}, G^![s])}{k^{\widehat{G^![s]}}(\check{\beta})} \left\| \beta|_{\mathfrak{a}_{M^!}} \right\| (1 - \beta(\delta[s]))^{-1} (1 - \beta(\delta[s])^{-1})^{-1} \\
		= \left\| \alpha|_{\mathfrak{a}_M} \right\| (1 - \alpha(\gamma))^{-1}(1 - \alpha(\gamma)^{-1})^{-1} .
	\end{multline}
	for all $\gamma \in T_{G\text{-reg}}(F)$ and $\delta \in T^!_{\tilde{G}\text{-reg}}(F)$ as in \eqref{eqn:stable-DE-aux1} (where $\delta$ was denoted as $\delta_\kappa$).

	We need the following quick observations.
	\begin{itemize}
		\item For $\alpha$, $\beta$ as in \eqref{eqn:stable-DE-aux2},
		\begin{equation*}
			\dfrac{\| \alpha|_{\mathfrak{a}_M} \|}{\| \beta|_{\mathfrak{a}_{M^!}} \|} = \begin{cases}
				1, & \alpha: \text{short}, \\
				2, & \alpha: \text{long}.
		\end{cases}\end{equation*}
		\item We have $k^{\widehat{G^![s]}}(\check{\beta}) \in \{1,2\}$. Moreover, $k^{\widehat{G^![s]}}(\check{\beta})=2$ if and only if $\alpha = \pm(\epsilon_i + \epsilon_j)$ in the usual basis (including the case $\alpha = \pm 2\epsilon_i$), such that $i$ and $j$ both fall under the $\GL(m)$-factor in $M^!$, and
		\begin{compactitem}
			\item either $\GL(m) \times \SO(2m'+1) \hookrightarrow \SO(2n'+1)$ when $s = +$, with $m' \geq 1$,
			\item or $\GL(m) \times \SO(2m''+1) \hookrightarrow \SO(2n''+1)$ when $s = -$, with $m'' \geq 1$.
		\end{compactitem}
		\item For all $s \in \Endo_{\mathbf{M}^!}(\tilde{G})$,
		\[ i_{M^!}(\tilde{G}, G^![s]) = 2^{\#\; \text{of $\SO$-factors in}\; M^! - \#\; \text{in $G^![s]$} }. \]
		\item If $\alpha$ is long, we can write $\alpha = \pm 2\epsilon_i$ and $\beta = \pm\epsilon_i$ in the usual bases; in this case
		\[ \beta(\delta[s])^2 = \beta(\delta)^2 = \alpha(\gamma), \quad \beta(\delta[+]) = -\beta(\delta[-]). \]
		If $\alpha$ is short, say $\alpha = \epsilon_i \pm \epsilon_j$ with $i \neq j$, then $\beta$ takes the same form. Claim:
		\begin{equation}\label{eqn:beta-alpha-short}
			\beta(\delta[s]) = \alpha(\gamma).
		\end{equation}
	
		Indeed, this is an assertion over $\CC$ and we may assume $T$ and $T^!$ are split. Express $\delta[s]$ in the standard coordinates for $T^! \subset G^![s]$ as
		\[ \delta[s] = (\underbracket{u_1, \ldots, u_{n'}}_{\SO(2n'+1)} , \underbracket{u_{n'+1}, \ldots, v_{n' + n''}}_{\SO(2n''+1)} ). \]
		Since $\delta[s] \leftrightarrow \gamma$ via $\mathbf{G}^![s] \in \Endo_{\elli}(\tilde{G})$, without loss of generality we can express
		\[ \gamma = \left( u_1, \ldots, u_{n'}, -u_{n'+1}, \ldots, -u_{n'+n''} \right). \]
		As $\beta \in \Sigma(G^![s], T^!)$, we have either $1 \leq i, j \leq n'$ or $n' < i, j \leq n'+n''$. Hence the minus signs in $\gamma$ have no effect: we have $\beta(\delta[s]) = u_i u_j^{\pm 1} = \alpha(\gamma)$, whence \eqref{eqn:beta-alpha-short}. 
	\end{itemize}

	Let us verify \eqref{eqn:stable-DE-aux2} when $\alpha$ is long. In this case, $\check{\beta}(\hat{s}) = 1$ for both $s$. Put
	\[ u := \beta(\delta[+]), \quad i(\pm) := i_{M^!}(\tilde{G}, G^![\pm]), \quad k(\pm) := k^{\widehat{G^![\pm]}}(\check{\beta}). \]
	As already observed, $\beta(\delta[-]) = -u$ and $\alpha(\gamma) = u^2$. The left hand side is $\|\beta|_{\mathfrak{a}_{M^!}}\|$ times
	\begin{equation}\label{eqn:stable-DE-aux4}
		\frac{i(+)}{k(+)} (1-u)^{-1}(1-u^{-1})^{-1} + \frac{i(-)}{k(-)} (1+u)^{-1}(1+u^{-1})^{-1}.
	\end{equation}

	In the table below, we assume $m', m'' \geq 1$ in each row. By the discussions above, we have
	\[\begin{array}{|c|c|c|c|} \hline
		\mathbf{M}^! & M^! & (i(+), k(+)) & (i(-), k(-)) \\ \hline
		(0, 0) & \GL(m) & \left(\frac{1}{2}, 1\right) & \left(\frac{1}{2}, 1\right) \\
		(m', 0) & \GL(m) \times \SO(2m'+1) & (1, 2) & \left(\frac{1}{2}, 1\right) \\
		(0, m'') & \GL(m) \times \SO(2m''+1) & \left(\frac{1}{2}, 1\right) & (1, 2)  \\
		(m', m'') & \GL(m) \times \SO(2m'+1) \times \SO(2m''+1) & (1, 2) & (1, 2) \\ \hline
	\end{array}\]
	In each case, $\frac{i(\pm)}{k(\pm)} = \frac{1}{2}$ and \eqref{eqn:stable-DE-aux4} equals $2 (1 - u^2)^{-1}(1 - u^{-2})^{-1}$, which equals the right hand side of \eqref{eqn:stable-DE-aux2} divided by $\|\beta|_{\mathfrak{a}_{M^!}}\|$.
	
	Consider the case of a short root $\alpha$, say $\alpha = \epsilon_i \pm \epsilon_j$ in the usual basis, where $i \neq j$; one can write $\beta$ in the same form. Define $i(\pm)$ and $k(\pm)$ as before. Since $\alpha \notin \Sigma(M, T)$, at least one of $i$, $j$ must fall under the $\GL(m)$ factor. In the table below, let us say $\check{\beta}$ is of $\GL$-type if both indices $i$, $j$ fall under $\GL(m)$; of $\SO^+$ type (resp.\ $\SO^-$ type) if one indexes falls under $\GL(m)$ and the other under $\SO(2m'+1)$ (resp.\ $\SO(2m''+1)$). Again, we assume $m', m'' \geq 1$ in each row of the table below.
	\[\begin{array}{|c|c|c|c|c|} \hline
		\mathbf{M}^! & \text{type of root} & (i(+), k(+)) & (i(-), k(-)) \\ \hline
		(0, 0) & \GL & \left(\frac{1}{2}, 1\right) & \left(\frac{1}{2}, 1\right) \\
		(m', 0) & \SO^+ & (1, 1) & \text{none} \\
		(m', 0) & \GL & (1, 2) & \left(\frac{1}{2}, 1\right) \\
		(0, m'') & \SO^- & \text{none} & (1, 1)  \\
		(0, m'') & \GL & \left(\frac{1}{2}, 1\right) & (1, 2) \\
		(m', m'') & \SO^+ & (1, 1) & \text{none} \\
		(m', m'') & \SO^- & \text{none} & (1, 1) \\
		(m', m'') & \GL & (1, 2) & (1, 2) \\ \hline
	\end{array}\]

	By \eqref{eqn:beta-alpha-short}, we have $\beta(\delta[s]) = \alpha(\gamma)$ for each $s$ with $\check{\beta}(\hat{s}) = 1$. This observation reduces \eqref{eqn:stable-DE-aux2} into to the assertion
	\[ \sum_{\substack{s \in \{+, -\} \\ \check{\beta}(\hat{s})=1 }} i(s) k(s)^{-1} = 1, \]
	which is apparently true by the table. The proof is by now complete.
\end{proof}

\begin{proof}[Proof of Theorem \ref{prop:stablization-DE}]
	We have to show the vanishing of the matrix $A(z) := \delta^{\tilde{G}, \Endo}_{\tilde{M}}(z) - \delta^G_M(z)$, for all $z \in \mathfrak{Z}(\mathfrak{g})$. This is based on the following inputs.
	\begin{itemize}
		\item The case of Casimir operator in Lemma \ref{prop:stabilization-Casimir}.
		\item The formula \eqref{eqn:delta-zz} for $\delta^G_M(zz')$ and the endoscopic version in Lemma \ref{prop:deltaEndo-zz}.
		\item The property \eqref{eqn:delta-vanishing} for the coefficients of $\delta^G_M$, and its counterparts for the operators $S\delta^{G^![s]}_{M^!}$.
	\end{itemize}
	The remaining arguments, albeit non-trivial, are the same as \cite[IX.2.5]{MW16-2}.
\end{proof}

\begin{corollary}[Cf.\ {\cite[Corollary 4.1]{Ar00}}]\label{prop:stable-DE}
	For all $f \in \orbI_{\asp}(\tilde{G}) \otimes \mes(G)$ and $z \in \mathfrak{Z}(\mathfrak{g})$, we have $I^{\Endo}_{\tilde{M}}(\cdot, zf) = \sum_{L \in \mathcal{L}(M)} \delta^L_M(z_L) I^{\Endo}_{\tilde{L}}(\cdot, f)$ as $C^\infty$ functions on $\widetilde{T_{G\text{-reg}}}$.
\end{corollary}
\begin{proof}
	This follows from $\delta^{\tilde{G}, \Endo}_{\tilde{M}} = \delta^G_M$ and the characterization of $\delta^{\tilde{G}, \Endo}_{\tilde{M}}$.
\end{proof}

\section{Jump relations in the real case}\label{sec:jump-relations}
We take $F = \R$ throughout.

\subsection{The case of invariant weighted orbital integrals: semi-regularity}\label{sec:jump-invariant}
The exposition here combines \cite{Ar08} and \cite[IX.4.1]{MW16-2}. Take $\tilde{G}$ and $M \in \mathcal{L}(M_0)$ as before. Let $T \subset M$ be a maximal torus and assume $\eta \in T(\R)$. In the study of jump relations, we make the following assumptions:
\begin{enumerate}[(i)]
	\item $T$ is elliptic in $M$;
	\item $G_{\eta, \mathrm{SC}} \simeq \SL(2)$;
	\item the preimage $T_d$ of $T$ in $G_{\eta, \mathrm{SC}}$ is isomorphic to $\Gm$.
\end{enumerate}
By (ii), there exists $\alpha \in \Sigma(G, T)$ such that $\alpha(\eta) = 1$, unique up to sign. Denote by $M_\alpha$ the Levi subgroup of $G$ obtained by adjoining the roots $\pm \alpha$ to $M$.

\begin{definition}[Semi-regularity]\label{def:semi-regular}
	\index{semi-regular}
	For simplicity's sake, the quadruplet $(\eta, T, M, G)$ is said to be \emph{semi-regular} when the conditions (i) --- (iii) above hold.
\end{definition}

The next result will imply, among others, that $M_\alpha$ contains $M$ as a maximal proper Levi subgroup, since it says $M_\eta \subsetneq (M_\alpha)_\eta$. In particular, $\alpha \notin \Sigma(M, T)$.

\begin{lemma}\label{prop:semi-regular}
	When $(\eta, T, M, G)$ is semi-regular, we have
	\begin{gather*}
		M_\eta = T, \quad (M_\alpha)_\eta = G_\eta, \quad A_{M_\eta} = A_M, \quad A_{M_\alpha} = A_{(M_\alpha)_\eta} = A_{G_\eta},
	\end{gather*}
	In particular, $\eta$ is $G$-equisingular (Definition \ref{def:equisingular}) and elliptic as an element of $M_\alpha$ or $M$; also, $(\eta, T, M, M_\alpha)$ is a semi-regular quadruplet.
\end{lemma}
\begin{proof}
	We begin by showing $M_\eta = T$. The inclusion $T \subset M_\eta$ is clear. Conversely, $M_\eta$ is a Levi subgroup of $G_\eta$. Observe that $G_\eta$ has semisimple rank $1$ by (ii) and has $T$ as a proper Levi by (iii). We cannot have $M_\eta = G_\eta$, otherwise $T$ would be non-elliptic in $M_\eta$ by considering $T_d \subset M_{\eta, \mathrm{SC}} = G_{\eta, \mathrm{SC}}$, thus non-elliptic in $M$; this contradicts (i).

	We have the inclusion of Levi subgroups $M_\eta \subset (M_\alpha)_\eta$ of $G_\eta$, which is proper since $\alpha$ goes into a root of $G_{\eta, \mathrm{SC}} \simeq \SL(2)$. Hence $(M_\alpha)_\eta = G_\eta$ by (ii).
	
	Now we have $A_{M_\eta} = A_T = A_M \subset A_{M_\eta}$, hence $A_{M_\eta} = A_M$. Moreover, we have just seen that $(M_\alpha)_\eta \supsetneq M_\eta$; this implies
	\[ \dim A_{M_\alpha} \leq \dim A_{(M_\alpha)_\eta} < \dim A_{M_\eta} = \dim A_M = \dim A_{M_\alpha} + 1. \]
	It follows that $A_{M_\alpha} = A_{(M_\alpha)_\eta} = A_{G_\eta}$.
\end{proof}

In what follows, we fix
\begin{itemize}
	\item a root $\alpha \in \Sigma(G, T)$ as above,
	\item $T_c \subset G_{\eta, \mathrm{SC}}$: anisotropic maximal torus, unique up to conjugacy over $\R$,
	\item $u \in G_{\eta, \mathrm{SC}}(\CC)$ such that $u T_{d, \CC} u^{-1} = T_{c, \CC}$,
	\item $H_d \in \mathfrak{t}_d(\R) \smallsetminus \{0\}$ and $H_c \in \mathfrak{t}_c(\R) \smallsetminus \{0\}$ such that $u H_d u^{-1} = iH_c$.
\end{itemize}

Extend the image of $T_c$ under $G_{\eta, \mathrm{SC}} \to G_\eta$ to a maximal torus $\underline{T}$ of $G_\eta$. It is also a maximal torus of $G$.

\begin{lemma}
	The torus $\underline{T}$ is an elliptic maximal torus of $M_\alpha$.
\end{lemma}
\begin{proof}
	By construction, $\dim A_{\underline{T}} = \dim A_T - 1$. Hence it suffices to notice that $\underline{T} \subset M_\alpha$, which follows from $(M_\alpha)_\eta = G_\eta$.
\end{proof}

The results above are essentially in \cite[IX.4.1 (4)]{MW16-2} (see also \cite[pp.196--197]{Ar08}), whose $\underline{\tilde{M}}$ corresponds to our $M_\alpha$.

Furthermore, denote by $w_d \in W(G_\eta, T)(\R)$ and $w_c \in W(G_\eta, \underline{T})(\R)$ the unique non-trivial elements in these Weyl groups. Choose Haar measures on $T(\R)$ and $\underline{T}(\R)$ by taking differential forms of top degrees, such that $\Ad(u): \mathfrak{t}^*_{\CC} \rightiso \underline{\mathfrak{t}}^*_{\CC}$ preserves these forms up to $\pm i$.

The situation is clearly modeled on the baby case $M = \Gm \hookrightarrow \SL(2) = G$, $\eta = \pm 1$, so that $C := \Ad(u)$ realizes a \emph{Cayley transform}.

Now fix $\tilde{\eta} \in \rev^{-1}(\eta)$. In Definition \ref{def:invariant-quadratic-form} we have chosen an explicit invariant positive-definite quadratic form on $X^*(T_0)$. It determines quadratic forms on $X^*(T_{\CC})$ and $X^*(\underline{T}_{\CC})$ compatibly with $C$, as in \S\ref{sec:stable-DE}; the corresponding Euclidean norms are both denoted as $\|\cdot\|$.

In what follows, the weighted orbital integrals will be defined using the Haar measures on $T(\R)$, $\underline{T}(\R)$ and $\mathfrak{a}^G_M$, $\mathfrak{a}^G_{\underline{M}}$ determined by the invariant quadratic forms above.

\begin{definition}
	\index{IM-mod@$I^{\mathrm{mod}}_{\tilde{M}}$}
	\index{M-underline@$\underline{M}$}
	For all $f \in \orbI_{\asp}(\tilde{G}) \otimes \mes(G)$ and $X \in \mathfrak{t}_d(\R)$ or $\mathfrak{t}_c(\R)$ which is $G$-regular and sufficiently close to $0$, set $\underline{M} := M_\alpha$ and
	\begin{align*}
		I^{\mathrm{mod}}_{\tilde{M}}(\exp(X)\tilde{\eta}, f) & = I^{\tilde{G}, \mathrm{mod}}_{\tilde{M}}(\exp(X)\tilde{\eta}, f) \\
		& := I_{\tilde{M}}\left( \exp(X)\tilde{\eta}, f \right)
		+ \left\| \check{\alpha} \right\| \log|\alpha(X)| \cdot I_{\underline{\tilde{M}}}(\exp(X)\tilde{\eta}, f) .
	\end{align*}
	
	Note that if the quadratic form on $X^*(T_0)$ is scaled by $\lambda > 0$, then both $I_{\tilde{M}}$ and $I_{\tilde{M}}^{\mathrm{mod}}$ are scaled by $\lambda^{\dim \mathfrak{a}^G_M}$.
\end{definition}
%%%%%%%%%%%%%%%%%%%%%%%%%%%%%%%%%%%%%%

It is a basic fact in the endoscopy for real groups that the stable conjugacy class of $H_c \in \mathfrak{g}_\eta(\R)$ contains either one or two $G_\eta(\R)$-conjugacy classes. Define
\[ c(\eta) = c^G(\eta) := \begin{cases}
	1, & \text{if the stable class of $H_c$ contains two conjugacy classes}, \\
	0, & \text{otherwise.}
\end{cases}\]
In terms of Galois cohomology, we have $2^{c(\eta)} = |\mathfrak{D}(\underline{T}, G_\eta; \R)| = |\mathfrak{D}(\underline{T}, \underline{M}_\eta; \R)$.
\index{c(eta)@$c(\eta)$}

\begin{proposition}[Jump relations, cf.\ {\cite[IX.4.1 Proposition]{MW16-2}}]\label{prop:jump-relation-I}
	Let $f \in \orbI_{\asp}(\tilde{G}) \otimes \mes(G)$ and $U \in \Sym(\mathfrak{t}_{\CC})$, $\underline{U} \in \Sym(\underline{\mathfrak{t}}_{\CC})$. The corresponding differential operators on the maximal tori are denoted by $\partial_U$ and $\partial_{\underline{U}}$, respectively.
	\begin{enumerate}[(i)]
		\item The limits of $\partial_U I^{\mathrm{mod}}_{\tilde{M}}\left( \exp(rH_d)\tilde{\eta}, f \right)$ as $r \to 0\pm$ both exist.
		\item If $w_d U = U$, the two limits above coincide.
		\item The limits of $\partial_{\underline{U}} I_{\underline{\tilde{M}}}\left( \exp(rH_c)\tilde{\eta}, f \right)$ as $r \to 0\pm$ both exist.
		\item If $w_d U = -U$, put $\underline{U} := uUu^{-1} \in \Sym(\underline{\mathfrak{t}}_{\CC})$. Then
		\begin{equation*}\begin{split}
			\lim_{r \to 0+} \partial_U I^{\mathrm{mod}}_{\tilde{M}}\left( \exp(rH_d)\tilde{\eta}, f \right) & -
			\lim_{r \to 0-} \partial_U I^{\mathrm{mod}}_{\tilde{M}}\left( \exp(rH_d)\tilde{\eta}, f \right) \\
			& = 2^{1 + c(\eta)} \pi i \|\check{\alpha}\| \lim_{r \to 0+} \partial_{\underline{U}} I_{\underline{\tilde{M}}}\left( \exp(rH_c)\tilde{\eta}, f \right) \\
			& = -2^{1 + c(\eta)} \pi i \|\check{\alpha}\| \lim_{r \to 0-} \partial_{\underline{U}} I_{\underline{\tilde{M}}}\left( \exp(rH_c)\tilde{\eta}, f \right).
		\end{split}\end{equation*}
	\end{enumerate}
\end{proposition}
\begin{proof}
	When $G=M$, all these assertions are due to Harish-Chandra. For general $M$, the arguments from \textit{loc.\ cit.}\ or \cite[\S 13]{Ar88LB} carry over verbatim. It does not involve any special property of coverings.
\end{proof}

\subsection{The stable version}\label{sec:jump-stable}
Now we consider a quasisplit $\R$-group $G^!$ with Levi subgroup $M^!$. Let $T^! \subset M^!$ be a maximal torus and $\epsilon \in T^!(\R)$. Assume that $(\epsilon, T^!, M^!, G^!)$ is semi-regular in the sense of Definition \ref{def:semi-regular}. There exists $\beta \in \Sigma(G^!, T^!) \smallsetminus \Sigma(M^!, T^!)$ such that $\beta(\epsilon) = 1$, unique up to sign. Denote by $M^!_\beta$ the Levi subgroup obtained by adjoining the roots $\pm \beta$ to $M^!$.

From these data, we define $T^!_d \subset G^!_{\epsilon, \mathrm{SC}}$, $T^!_c \subset G^!_{\epsilon, \mathrm{SC}}$, $H^!_d$, $H^!_c$ and $u \in G^!_{\epsilon, \mathrm{SC}}(\CC)$ such that $uT_{d, \CC} u^{-1} = T_{c, \CC}$ as in \S\ref{sec:jump-invariant}. We also have the Weyl group elements $w_d$ and $w_c$ for $T^!_d$ and $T^!_c$ respectively.

Following \cite[IX.4.2]{MW16-2}, define
\begin{equation}\label{eqn:iGM-jump}\begin{aligned}
	Z^{\check{\beta}}_{(M^!)^\vee} & := \left\{ z \in Z_{(M^!)^\vee}^{\Gamma_{\R}} : \check{\beta}(z) = 1 \right\} \supset Z_{(M^!_\beta)^\vee}^{\Gamma_{\R}}, \\
	i^{G^!}_{M^!}(\epsilon) & := \left( Z_{(M^!)^\vee}^{\check{\beta}} : Z_{(M^!_\beta)^\vee}^{\Gamma_{\R}} \right)^{-1};
\end{aligned}\end{equation}
taking $\Gamma_{\R}$-invariants is actually superfluous in our case.

As before, we have the invariant positive-definite quadratic form on $X^*(T^!_{\CC})$ induced by the choice made in Definition \ref{def:invariant-quadratic-form}, so that $\|\check{\beta}\|$ makes sense. Same for $X*(\underline{T}^!_{\CC})$.

\begin{definition}\label{def:Smod}
	\index{SM-mod@$S^{\mathrm{mod}}_{M^{"!}}$}
	For all $f^! \in S\orbI(G^!) \otimes \mes(G^!)$ and $Y \in \mathfrak{t}^!_d(\R)$ or $\mathfrak{t}^!_c(\R)$ which is $G^!$-regular and sufficiently close to $0$, set $\underline{M}^! := M^!_\beta$ and use the invariant quadratic forms above to define
	\begin{align*}
		S^{\mathrm{mod}}_{M^!}\left(\exp(Y)\epsilon, f^!\right) & = S^{G^!, \mathrm{mod}}_{M^!}\left(\exp(Y)\epsilon, f^!\right) \\
		& := S_{M^!}\left( \exp(Y)\epsilon, f^! \right)
		+ i^{G^!}_{M^!}(\epsilon) \left\| \check{\beta} \right\| \log|\beta(Y)| \cdot S_{\underline{M}^!}\left(\exp(Y)\epsilon, f^!\right) .
	\end{align*}
\end{definition}

The following is a restatement of \cite[IX.4.2 Proposition]{MW16-2}.

\begin{proposition}\label{prop:jump-stable}
	For every $U^! \in \Sym(\mathfrak{t}^!_{\CC})$, the following variants of jump relations in Proposition \ref{prop:jump-relation-I} hold.
	\begin{enumerate}[(i)]
		\item The limits of $\partial_{U^!} S^{\mathrm{mod}}_{M^!}\left( \exp(rH^!_d)\epsilon, f^! \right)$ as $r \to 0\pm$ both exist.
		\item If $w_d U^! = U^!$, the two limits above coincide.
		\item The limits of $\partial_{\underline{U}^!} S_{\underline{M}^!}\left( \exp(rH^!_c)\epsilon, f^! \right)$ as $r \to 0\pm$ both exist.
		\item If $w_d U^! = -U^!$, put $\underline{U}^! := uU^! u^{-1} \in \Sym(\underline{\mathfrak{t}^!}_{\CC})$. Then
		\begin{equation*}\begin{split}
			\lim_{r \to 0+} \partial_{U^!} S^{\mathrm{mod}}_{M^!}\left( \exp(rH^!_d)\epsilon, f^! \right) & -
			\lim_{r \to 0-} \partial_{U^!} S^{\mathrm{mod}}_{M^!}\left( \exp(rH^!_d)\epsilon, f^! \right) \\
			& = 2 \pi i \|\check{\beta}\| i^{G^!}_{M^!}(\epsilon) \lim_{r \to 0+} \partial_{\underline{U}^!} S_{\underline{M}^!}\left( \exp(rH^!_c)\epsilon, f^! \right) \\
			& = -2 \pi i \|\check{\beta}\| i^{G^!}_{M^!}(\epsilon) \lim_{r \to 0-} \partial_{\underline{U}^!} S_{\underline{M}^!}\left( \exp(rH^!_c)\epsilon, f^! \right).
		\end{split}\end{equation*}
	\end{enumerate}
\end{proposition}

\subsection{The endoscopic version: desiderata}
Consider $\rev: \tilde{G} \to G(\R)$, $M \subset G$, $\eta$, $\tilde{\eta}$, $\alpha$, $T$, $\underline{T}$ as in the context of Proposition \ref{prop:jump-relation-I}.

Fix Haar measures on $T(\R)$, $\underline{T}(\R)$, in a compatible way. For every $\mathbf{M}^! \in \Endo_{\elli}(\tilde{M})$, the Haar measures of the corresponding maximal tori $T^!$ are chosen compatibly with the isomorphisms $T^! \simeq T$ arising from any diagram.

The invariant positive-definite quadratic forms on $T^!$ and $T$ are also compatibly chosen; see the explanations in \S\ref{sec:stable-DE}.

\begin{definition}\label{def:IMEndo-mod}
	\index{IMEndo-mod@$I^{\Endo, \mathrm{mod}}_{\tilde{M}}$}
	For all $f \in \orbI_{\asp}(\tilde{G}) \otimes \mes(G)$ and $X \in \mathfrak{t}_{G\text{-reg}}(\R)$ which is sufficiently close to $0$, set $\underline{M} := M_\alpha$ and
	\begin{align*}
		I^{\Endo, \mathrm{mod}}_{\tilde{M}}(\exp(X)\tilde{\eta}, f) & = I^{\tilde{G}, \Endo, \mathrm{mod}}_{\tilde{M}}(\exp(X)\tilde{\eta}, f) \\
		& := I^{\Endo}_{\tilde{M}}\left( \exp(X)\tilde{\eta}, f \right)
		+ \left\| \check{\alpha} \right\| \log|\alpha(X)| \cdot I^{\Endo}_{\underline{\tilde{M}}}(\exp(X)\tilde{\eta}, f) .
	\end{align*}
	It equals $I^{\mathrm{mod}}_{\tilde{M}}(\exp(X)\tilde{\eta}, f)$ if the local geometric Theorem \ref{prop:local-geometric} holds.
\end{definition}

\begin{theorem}\label{prop:jump-Endo}
	The jump relations in Proposition \ref{prop:jump-relation-I} continue to hold in the endoscopic setting, in the sense below.
		\begin{enumerate}[(i)]
		\item The limits of $\partial_U I^{\Endo, \mathrm{mod}}_{\tilde{M}}\left( \exp(rH_d)\tilde{\eta}, f \right)$ as $r \to 0\pm$ both exist.
		\item If $w_d U = U$, the two limits above coincide.
		\item The limits of $\partial_{\underline{U}} I^{\Endo}_{\underline{\tilde{M}}}\left( \exp(rH_c)\tilde{\eta}, f \right)$ as $r \to 0\pm$ both exist.
		\item If $w_d U = -U$, put $\underline{U} := uUu^{-1} \in \Sym(\underline{\mathfrak{t}}_{\CC})$. Then
		\begin{equation*}\begin{split}
			\lim_{r \to 0+} \partial_U I^{\Endo, \mathrm{mod}}_{\tilde{M}}\left( \exp(rH_d)\tilde{\eta}, f \right) & -
			\lim_{r \to 0-} \partial_U I^{\Endo, \mathrm{mod}}_{\tilde{M}}\left( \exp(rH_d)\tilde{\eta}, f \right) \\
			& = 2^{1 + c(\eta)} \pi i \|\check{\alpha}\| \lim_{r \to 0+} \partial_{\underline{U}} I^{\Endo}_{\underline{\tilde{M}}}\left( \exp(rH_c)\tilde{\eta}, f \right) \\
			& = -2^{1 + c(\eta)} \pi i \|\check{\alpha}\| \lim_{r \to 0-} \partial_{\underline{U}} I^{\Endo}_{\underline{\tilde{M}}}\left( \exp(rH_c)\tilde{\eta}, f \right).
		\end{split}\end{equation*}
	\end{enumerate}
\end{theorem}

The proof will be achieved in \S\S\ref{sec:pf-jump-Endo}--\ref{sec:pf-jump-Endo-lemmas}.

\section{Combinatorial preparations}
Retain the conventions from \S\ref{sec:jump-relations} about $\tilde{G}$, $\tilde{M}$, etc. Fix the following data;
\begin{itemize}
	\item $\mathbf{M}^! \in \Endo_{\elli}(\tilde{M})$,
	\item an elliptic maximal torus $T^! \subset M^!$, isomorphic to $T \subset M$ via a chosen diagram joining $\epsilon \in T^!(\R)$ and $\eta \in T(\R)$; in particular $\epsilon \leftrightarrow \eta$ with respect to $\mathbf{M}^!$.
\end{itemize}

Furthermore, we take $\underline{M} \in \mathcal{L}^G(M)$, $t \in \Endo_{\mathbf{M}^!}(\underline{\tilde{M}})$ and $u \in \Endo_{\mathbf{\underline{M}}^!}(\underline{\tilde{G}})$ where $\mathbf{\underline{M}}^! := \mathbf{\underline{M}}^![t]$; recall that $t$ and $u$ determine $s \in \Endo_{\mathbf{M}^!}(\tilde{G})$, so that we are in the following situation (cf.\ \eqref{eqn:s-situation}):
\begin{equation}\label{eqn:semireg-situation}\begin{tikzcd}
	\mathbf{G}^![s] \arrow[dash, dashed, r, "\text{endo.}", "\text{ell.}"'] & \tilde{G} \\
	\mathbf{\underline{M}}^! \arrow[dash, dashed, r, "\text{endo.}", "\text{ell.}"'] \arrow[hookrightarrow, u, "\text{Levi}"] & \underline{\tilde{M}} \arrow[hookrightarrow, u, "\text{Levi}"'] \\
	\mathbf{M}^! \arrow[dash, dashed, r, "\text{endo.}", "\text{ell.}"'] \arrow[hookrightarrow, u, "\text{Levi}"] & \tilde{M} \arrow[hookrightarrow, u, "\text{Levi}"']
\end{tikzcd} \quad \mathbf{G}^![s] = \mathbf{G}^![u]. \end{equation}

They are subject to the following conditions:
\begin{itemize}
	\item $(\eta, T, M, G)$ (resp.\ $(\epsilon[s], T^!, M^!, G^![s])$) is semi-regular (Definition \ref{def:semi-regular}), so that we obtain an $\alpha \in \Sigma(G, T) \smallsetminus \Sigma(M, T)$ (resp.\ $\beta \in \Sigma(G^![s], T^!) \smallsetminus \Sigma(M^!, T^!)$), canonical up to sign;
	\item $\underline{M} = M_\alpha$, $\underline{M}^! = M^!_\beta$;
	\item $\alpha$ and $\beta$ are chosen so that they correspond under $\xi_{T^!, T}: T^! \rightiso T$, up to a rescaling factor when $\alpha$ is long.
\end{itemize}

Using the diagram, we may identify $X_*(T_{\CC})$ with $X^*(T^\vee)$ where $T^\vee \subset \tilde{G}^\vee$ is part of the root datum of $\tilde{G}^\vee$; see Remark \ref{rem:diagram-dual}.

\begin{definition}[Cf.\ \eqref{eqn:iGM-jump}]\label{def:ZMvee-alpha}
	\index{ZMvee-alpha@$Z_{\tilde{M}^\vee}^{\hat{\alpha}^*}$}
	In this situation, let $\hat{\alpha}^*: Z_{\tilde{M}^\vee}^\circ \to \CC^\times$ be the restriction of $\check{\alpha} \in X^*(T^\vee)$. This is independent of the choice of diagrams. Set
	\[ Z_{\tilde{M}^\vee}^{\hat{\alpha}^*} := \Ker(\hat{\alpha}^*). \]
\end{definition}

\begin{lemma}\label{prop:alpha-Z-Mbar}
	The homomorphism $\hat{\alpha}^*$ is trivial on $Z_{\underline{\tilde{M}}^\vee}^\circ$, and $\left( Z_{\tilde{M}^\vee}^{\hat{\alpha}^*} : Z_{\underline{\tilde{M}}^\vee}^\circ \right)$ is finite.
\end{lemma}
\begin{proof}
	If we work in the non-metaplectic dual $G^\vee$ of $G$, then $\hat{\alpha}^*$ can be defined from $\check{\alpha} \in X_*(T_{\CC})$ as usual, and the corresponding statements are in \cite[p.194]{Ar08}. In fact, in this case one readily shows that $\check{\alpha}$ is trivial on $Z_{\underline{M}^\vee}$.
	
	In the metaplectic case, the $\SO$-factors in $G^\vee$ and its Levi subgroups are flipped into $\Sp$-factors, but they all share the same maximal torus $T^\vee$. Hence $Z_{\tilde{L}^\vee} \supset Z_{L^\vee}$ as subgroups of $T^\vee$ for all standard Levi subgroups $L \subset G$, and $Z_{\tilde{L}^\vee}^\circ = Z_{L^\vee}^\circ$. Therefore, the aforementioned results carry over.
\end{proof}

As in \S\ref{sec:jump-relations}, $X_*(T_{\CC})$ and $X_*(T^!_{\CC})$ carry compatible invariant Euclidean norms, both denoted by $\|\cdot\|$.

\begin{lemma}[Cf.\ {\cite[IX.4.5 (5)]{MW16-2}}]\label{prop:jump-m}
	Let $i^{G^![s]}_{M^!}(\epsilon[s])$ be as in  \eqref{eqn:iGM-jump}. We have
	\[ i_{M^!}\left( \tilde{G}, G^![s] \right) i^{G^![s]}_{M^!}(\epsilon[s]) \cdot \left\|\check{\beta}\right\| =
	\left( Z_{\tilde{M}^\vee}^{\hat{\alpha}^*} : Z_{\underline{\tilde{M}}^\vee}^\circ \right)^{-1}
	i_{\underline{M}^!}\left( \tilde{G}, G^![s] \right) \cdot \left\|\check{\alpha}\right\|. \]
\end{lemma}
\begin{proof}
	Use \cite[Lemma 1.1]{Ar99} to obtain the natural surjection $Z_{\tilde{M}^\vee}^\circ / Z_{\tilde{G}^\vee}^\circ \twoheadrightarrow Z_{(M^!)^\vee} / Z_{G^![s]^\vee}$. Denote its kernel as $K_1$. One readily checks that $|K_1|^{-1} = i_{M^!}\left( \tilde{G}, G^![s] \right)$ as in \cite[p.230 (2)]{MW16-1}.

	We contend that the image of $Z_{\tilde{M}^\vee}^{\hat{\alpha}^*} / Z_{\tilde{G}^\vee}^\circ$ is contained in $Z_{(M^!)^\vee}^{\check{\beta}} / Z_{G^![s]^\vee}$. When $\alpha$ is short, $\check{\alpha}$ transports to $\check{\beta}$ under $T^! \simeq T$; in this case $Z_{\tilde{M}^\vee}^{\hat{\alpha}^*} / Z_{\tilde{G}^\vee}^\circ$ is actually the preimage of $Z_{(M^!)^\vee}^{\check{\beta}} / Z_{G^![s]^\vee}$. When $\alpha$ is long, $\check{\alpha}$ transports to $\frac{1}{2} \check{\beta}$, and the containment is clear.

	Set $K'_1 := K_1 \cap (Z_{\tilde{M}^\vee}^{\hat{\alpha}^*} / Z_{\tilde{G}^\vee}^\circ)$; we have just seen that $K_1 = K'_1$ if $\alpha$ is short. From these and Lemma \ref{prop:alpha-Z-Mbar}, we obtain the following commutative diagram of abelian groups, with exact rows:
	\[\begin{tikzcd}
        & 1 \arrow[d] & 1 \arrow[d] & 1 \arrow[d] & & \\
		1 \arrow[r] & K_3 \arrow[d] \arrow[r] & Z_{\underline{\tilde{M}}^\vee}^\circ / Z_{\tilde{G}^\vee}^\circ \arrow[d] \arrow[r] & Z_{(\underline{M}^!)^\vee} / Z_{G^![s]^\vee} \arrow[r] \arrow[d] & 1 \arrow[r] \arrow[d] & 1 \\
		1 \arrow[r] & K'_1 \arrow[r] \arrow[d] & Z_{\tilde{M}^\vee}^{\hat{\alpha}^*} / Z_{\tilde{G}^\vee}^\circ \arrow[r] \arrow[d] & Z_{(M^!)^\vee}^{\check{\beta}} / Z_{G^![s]^\vee} \arrow[r] \arrow[d] & C_1 \arrow[r] \arrow[d] & 1 \\
		1 \arrow[r] & K_2 \arrow[r] \arrow[d] & Z_{\tilde{M}^\vee}^{\hat{\alpha}^*} / Z_{\underline{\tilde{M}}^\vee}^\circ \arrow[r] \arrow[d] & Z_{(M^!)^\vee}^{\check{\beta}} / Z_{(\underline{M}^!)^\vee} \arrow[r] \arrow[d] & C_2 \arrow[r] & 1 \\
        & 1 & 1 & 1 & &
	\end{tikzcd}\]
	where $K_2, K_3$ (resp.\ $C_1$, $C_2$) are defined to be the kernels (resp.\ cokernels); they are all finite. The second and the third columns are readily seen to be exact, hence so is the first column by the Snake Lemma.
	
	Next, observe that $|K_3|^{-1} = i_{\underline{M}^!}(\tilde{G}, G^![s])$ as in the case of $|K_1|^{-1}$. Hence
	\begin{align*}
		i_{M^!}(\tilde{G}, G^![s]) & = |K_1|^{-1} = |K'_1|^{-1} (K_1 : K'_1)^{-1} \\
		& = |K_2|^{-1} |K_3|^{-1} (K_1 : K'_1)^{-1} \\
		& = |K_2|^{-1} i_{\underline{M}^!}(\tilde{G}, G^![s]) (K_1 : K'_1)^{-1}.
	\end{align*}
	It remains to prove that
	\begin{equation*}
		|K_2| (K_1 : K'_1) = \left( Z_{\tilde{M}^\vee}^{\hat{\alpha}^*} : Z_{\underline{\tilde{M}}^\vee}^\circ \right)  i^{G^![s]}_{M^!}(\epsilon[s]) \cdot \frac{\|\check{\beta}\|}{\|\check{\alpha}\|} .
	\end{equation*}

	Using the third row of the diagram and \eqref{eqn:iGM-jump}, we see $|K_2| = \left( Z_{\tilde{M}^\vee}^{\hat{\alpha}^*} : Z_{\underline{\tilde{M}}^\vee}^\circ \right)  i^{G^![s]}_{M^!}(\epsilon[s]) |C_2|$, thus we are reduced to proving
	\begin{equation*}
		|C_2| (K_1 : K'_1) = \frac{\|\check{\beta}\|}{\|\check{\alpha}\|}.
	\end{equation*}

	When $\alpha$ is short, we have seen that $K_1 = K'_1$, $\|\check{\alpha}\| = \|\check{\beta}\|$ whilst $C_2 = \{1\}$ (upon replacing $\tilde{G}$ by $\underline{\tilde{M}}$). The required equality follows at once.
	
	Hereafter, suppose that $\alpha$ is long. Write $\tilde{G} = \prod_{i \in I} \GL(n_i) \times \Mp(2n)$. Without loss of generality, we may express $\alpha = 2\epsilon_i$ under the usual basis for $\Sp(2n)$. The index $i$ must fall under a $\GL$-factor of $M$ that embeds into $\Sp(2n)$. Moreover, $\check{\alpha} = \check{\epsilon}_i$ and $\check{\beta} = 2\check{\epsilon}_i$. It is clear that
	\begin{gather*}
		Z_{\tilde{M}^\vee}^\circ = Z_{(M^!)^\vee}^\circ , \quad Z_{\tilde{M}^\vee}^{\hat{\alpha}^*} = Z_{\underline{\tilde{M}}^\vee}^\circ = Z_{(\underline{M}^!)^\vee}^\circ, \\
		(K_1 : K'_1 ) = \left( Z_{(M^!)^\vee}^\circ \cap Z_{G^![s]^\vee} : Z_{(\underline{M}^!)^\vee}^\circ \cap Z_{G^![s]^\vee} \right).
	\end{gather*}
	Thus it remains to verify in this case that
	\begin{equation}\label{eqn:jump-m-aux}
		\left( Z_{(M^!)^\vee}^{\check{\beta}} : Z_{(\underline{M}^!)^\vee} \right) \left( Z_{(M^!)^\vee}^\circ \cap Z_{G^![s]^\vee} : Z_{(\underline{M}^!)^\vee}^\circ \cap Z_{G^![s]^\vee} \right) = 2.
	\end{equation}

	Observe that \eqref{eqn:jump-m-aux} involves only the objects on the endoscopic side. We may write
	\[ G^![s] = \prod_{i \in I} \GL(n_i) \times \SO(2n'+1) \times \SO(2n''+1), \quad n' + n'' = n. \]
	The first step is to reduce to the case $I = \emptyset$ and $n'' = 0$. Indeed, $M^!$ and $\underline{M}^!$ decompose accordingly, and the construction of $\underline{M}^!$ takes place inside either $\SO(2n'+1)$ or $\SO(2n''+1)$, on which $\beta$ lives. Hence we may rename $G^![s]$ to $G^!$ and assume $G^! = \SO(2n+1)$.
	
	Accordingly, we can write $M^! = \prod_{j=1}^k \GL(m_j) \times \SO(2m+1)$ where $m \in \Z_{\geq 0}$, such that $\check{\beta}$ factors through the dual of $\GL(m_1)$, so that $\underline{M}^!$ is obtained by merging $\GL(m_1)$ with $\SO(2m'+1)$ to form a larger Levi subgroup of $G^!$. As a referee pointed out, we actually have $m_1 = 1$: by semi-regularity, $\epsilon_i$ is Galois-invariant and $i$ is in the part of the basis landing in $\GL(m_1)$, thus $m_1 = 1$ since $T^! \cap \GL(m_1)$ must be an elliptic maximal torus of $\GL(m_1)$.
	\begin{itemize}
		\item Suppose $m=0$, then the first index in \eqref{eqn:jump-m-aux} is $1$ since
		\[ Z_{(M^!)^\vee}^{\check{\beta}} = \{\pm 1\} \times \prod_{j \geq 2} \CC^\times = Z_{(\underline{M}^!)^\vee}. \]
		On the other hand, $Z_{(M^!)^\vee}^\circ \cap Z_{(G^!)^\vee} = Z_{(G^!)^\vee} = \{\pm 1\}$ and $Z_{(\underline{M}^!)^\vee}^\circ \cap Z_{(G^!)^\vee} = \{1\}$, so the second index is $2$. Hence \eqref{eqn:jump-m-aux} is verified.
		\item Suppose $m \geq 1$, then
		\[ Z_{(M^!)^\vee}^{\check{\beta}} = \{\pm 1\} \times \prod_{j \geq 2} \CC^\times \times \{\pm 1\} ,\]
		whilst $Z_{(\underline{M}^!)^\vee}$ has only one $\{\pm 1\}$-factor diagonally embedded, so the first index is $2$. On the other hand,
		\[ Z_{(M^!)^\vee}^\circ = \prod_{j \geq 1} \CC^\times \times \{1\} \]
		intersects trivially with $Z_{(G^!)^\vee} \simeq \{\pm 1\}$, so the second index is $1$. Again, \eqref{eqn:jump-m-aux} is verified.
	\end{itemize}

	Summing up, the case of long roots is completed.
\end{proof}

A crucial question remains: how to pass from the second to the first column of \eqref{eqn:semireg-situation}? This is addressed by the following result.

\begin{lemma}\label{prop:jump-t}
	Given a semi-regular quadruplet $(\eta, T, M, G)$ and $\alpha$, define $\underline{M} := M_\alpha$. Let $\mathbf{M}^! \in \Endo_{\elli}(\tilde{M})$ and take a diagram joining $\epsilon \in T^! \subset M^!$ to $\eta \in T \subset M$. Suppose that $M^!_\epsilon$ is quasisplit.
	\begin{enumerate}[(i)]
		\item Take $s^\flat \in T^\vee$ that gives rise to $\mathbf{M}^!$. Interpret $\Endo_{\mathbf{M}^!}(\underline{\tilde{M}})$ as a subset of $s^\flat Z_{\tilde{M}^\vee}^\circ / Z_{\underline{\tilde{M}}^\vee}^\circ$ (Proposition \ref{prop:Endo-s-interpretation}). Then those $t \in \Endo_{\mathbf{M}^!}(\underline{\tilde{M}})$ such that $(\epsilon[t], T^!, M^!, \underline{M}^![t])$ is semi-regular form a coset under $Z_{\tilde{M}^\vee}^{\hat{\alpha}^*} / Z_{\underline{\tilde{M}}^\vee}^\circ$.
		\item As a consequence,
		\[ \left( Z_{\tilde{M}^\vee}^{\hat{\alpha}^*} : Z_{\underline{\tilde{M}}^\vee}^\circ \right) = \left|\left\{ t \in \Endo_{\mathbf{M}^!}(\underline{\tilde{M}}) : (\epsilon[t], T^!, M^!, \underline{M}^![t])\; \text{is semi-regular} \right\}\right|. \]
	\end{enumerate}
\end{lemma}
\begin{proof}
	Recall from \S\ref{sec:dual} that for each $t \in s^\flat Z_{\tilde{M}^\vee}^\circ / Z_{\underline{\tilde{M}}^\vee}^\circ$, we can construct $\underline{\mathbf{M}}^![t] \in \Endo(\underline{\tilde{M}})$. In Definition \ref{def:ZMvee-alpha} we considered the coroot $\check{\alpha} \in X^*(T^\vee)$ attached to $\alpha$. Recall that we set $\check{\beta} = \check{\alpha}$ (resp.\ $\check{\beta} = 2\check{\alpha}$) when $\alpha$ is short (resp.\ long), so that $\check{\beta} \in \Sigma(\tilde{G}^\vee, T^\vee)$. We claim that in order to have $t \in \Endo_{\mathbf{M}^!}(\underline{\tilde{M}})$ and $(\epsilon[t], T^!, M^!, \underline{M}^!)$ be semi-regular, it is necessary and sufficient that
	\begin{itemize}
		\item $\check{\beta}(t) = 1$ where $t$ is identified as an element of $T^\vee$, and
		\item the resulting $\beta \in \Sigma(\underline{M}^!, T^!) \smallsetminus \Sigma(M^!, T^!)$ satisfies $\beta(\epsilon[t]) = 1$.
	\end{itemize}

	The necessity of these conditions is clear. Consider now their sufficiency. The ellipticity of $T^!$ in $M^!$ is automatic. The first condition above implies $\beta \in \Sigma(\underline{M}^![t], T^!)$, whilst $\beta \notin \Sigma(M^!, T^!)$ since $\alpha \notin \Sigma(M, T)$. By comparing semisimple ranks, we see that $\underline{\mathbf{M}}^![t]$ must be an elliptic endoscopic datum of $\underline{\tilde{M}}$. Next, recall from \S\ref{sec:descent-endoscopy} that
	\[\begin{tikzcd}[column sep=large]
		\underline{M}^!{[t]}_{\epsilon[t]} \arrow[dashed, dash, r, "\text{non.\ std.}"] \arrow[dashed, dash, r, "\text{endo.}"'] &
		\overline{\underline{M}^!{[t]}_{\epsilon[t]}} \arrow[dashed, dash, r, "\text{std.}"] \arrow[dashed, dash, r, "\text{endo.}"'] &
		\underline{M}_\eta ,
	\end{tikzcd}\]
	the rightmost group has semisimple rank $1$; ditto for the leftmost one by the second condition (since it has the roots $\pm\beta$). It is routine to deduce that $\overline{\underline{M}^![t]_{\epsilon[t]}}$ is a quasisplit inner twist\footnote{It is actually an isomorphism, since $\underline{M}_\eta$ is quasisplit.} of $\underline{M}_\eta$. From this and the semi-regularity of $(\eta, T, M, \underline{M})$, it is immediate that
	\begin{itemize}
		\item the preimage $T^!_d$ of $T^!$ in $\underline{M}^![t]_{\epsilon[t], \mathrm{SC}}$ is isomorphic to $\Gm$;
		\item $\underline{M}^![t]_{\epsilon[t], \mathrm{SC}} \simeq \SL(2)$.
	\end{itemize}

	To complete the proof, we apply the claim in the following way. Those $t \in s^\flat Z_{\tilde{M}^\vee}^\circ / Z_{\underline{\tilde{M}}^\vee}^\circ$ such that $\check{\beta}(t) = 1$ form a coset under $(Z_{\tilde{M}^\vee}^\circ)^{\check{\beta} = 1} / Z_{\underline{\tilde{M}}^\vee}^\circ$. If $\alpha$ is short, by \eqref{eqn:beta-alpha-short} we have $\beta(\epsilon[t]) = \alpha(\eta) = 1$ and $(Z_{\tilde{M}^\vee}^\circ)^{\check{\beta} = 1} = Z_{\tilde{M}^\vee}^{\hat{\alpha}^*}$. If $\alpha$ is long, then $\beta$ is short in $\underline{M}^![t]$, and it is routine to check that the condition $\beta(\epsilon[t]) = 1$ cuts out a $Z_{\tilde{M}^\vee}^{\hat{\alpha}^*} / Z_{\underline{\tilde{M}}^\vee}^\circ$-coset.
\end{proof}

\section{Proof of Theorem \ref{prop:jump-Endo}}\label{sec:pf-jump-Endo}
Let $(\eta, T, M, G)$ be a semi-regular quadruplet. Take a corresponding root $\alpha \in \Sigma(G, T) \smallsetminus \Sigma(M, T)$ and put $\underline{M} := M_\alpha$.

\begin{definition}\label{def:semi-reg-s}
	Given $\kappa \in \mathfrak{R}(T)$ and the corresponding $\mathbf{M}^! \in \Endo_{\elli}(\tilde{M})$ given by Tate--Nakayama duality, fix once and for all a diagram joining $\eta \in T \subset M$ and $\epsilon \in T^! \subset M^!$. Suppose that $s \in \Endo_{\mathbf{M}^!}(\tilde{G})$ and $(\epsilon[s], T^!, M^!, G^![s])$ is semi-regular, with the root $\beta \in \Sigma(G^![s], T^!) \smallsetminus \Sigma(M^!, T^!)$ corresponding to $\alpha$ in the sense that for $Y = \dd\xi_{T^!, T}(X) \in \mathfrak{t}^!(\R)$,
	\[ \beta(Y) = \begin{cases}
		\alpha(X), & \alpha: \text{short}, \\
		\frac{1}{2} \alpha(X), & \alpha: \text{long.}
	\end{cases}\]
	In this case we simply say that $s \in \Endo_{\mathbf{M}^!}(\tilde{G})$ is \emph{semi-regular}, and write $\underline{M}^! := M^!_\beta$.
\end{definition}

\begin{remark}\label{rem:semi-reg-equisingular}
	Suppose that $s$ is semi-regular. By Lemmas \ref{prop:equisingular-endo} and \ref{prop:semi-regular}, $\epsilon[s]$ is $G^![s]$-equisingular as an element of $\underline{M}^!$. The semi-regularity of $(\epsilon[s], T^!, M^!, G^![s])$ also implies that $M^!_\epsilon = M^!_{\epsilon[s]} = T^!$; in particular, $M^!_\epsilon = M^!_{\epsilon[s]}$ is quasisplit.
\end{remark}

Fix $\tilde{\eta} \in \rev^{-1}(\eta)$ and $f \in \orbI_{\asp}(\tilde{G}) \otimes \mes(G)$. For $X \in \mathfrak{t}_{G\text{-reg}}(\R)$ which is close to $0$, we can view $\exp(X)\tilde{\eta}$ as an element of $D_{\mathrm{orb}, \mathrm{G}\text{-reg}, -}(\tilde{M}) \otimes \mes(M)^\vee$ by using the prescribed Haar measure on $T(\R)$. Denote this element as $\tilde{\gamma}(X)$.

We shall adopt the theory of \S\ref{sec:derivative-orbint} in its Lie algebra version. For each $\kappa \in \mathfrak{R}(T)$, take the $(\mathbf{M}^!, T^!, Y) = (\mathbf{M}^!_\kappa, T^!_\kappa, Y_\kappa)$ (the pair $(Y, T^!)$ is taken up to stable conjugacy), such that $Y \in \mathfrak{t}^!(\R)$ and $X = \dd\xi_{T^!, T}(Y)$. Thanks to \cite{Ko82}, $(Y, T^!)$ can be chosen within the stable class so that $\epsilon = \epsilon_\kappa := \tilde{\xi}_{T^!, T}(\delta)$ makes $M^!_\epsilon$ quasisplit.

Transport the Haar measure to $T^!(\R)$ and define
\[ \delta_\kappa(X) := \Delta(\exp(X)\tilde{\eta}, \exp(Y)\epsilon) \exp(Y)\epsilon \; \in SD_{\mathrm{orb}}(M^!) \otimes \mes(M)^\vee \]
where $\Delta(\exp(X)\tilde{\eta}, \exp(Y)\epsilon)$ is the adjoint transfer factor (Definition \ref{def:adjoint-transfer-factor}). This leads to
\begin{equation}\label{eqn:gammaX-deltaX}
	\tilde{\gamma}(X) = \sum_{\kappa \in \mathfrak{R}(T)} \trans_{\mathbf{M}^!, \tilde{M}}\left( \delta_\kappa(X) \right).
\end{equation}

By Remark \ref{rem:semi-reg-equisingular} and \cite[V.1.4 (2)]{MW16-1}, there exists $\varphi \in S\orbI(\underline{\mathfrak{m}}^!) \otimes \mes(\underline{M}^!)$ such that
\[ S^{G^![s]}_{\underline{M}^!}\left(\exp(Y) \epsilon[s], f^{G^![s]}\right) = S^{\underline{M}^!}\left(\exp(Y) \epsilon[s], \varphi \right), \]
for $Y$ close to $0$. Since $T^!_d \subset \underline{M}^!_{\epsilon[s], \mathrm{SC}} \simeq \SL(2)$ is split, the condition $\mathrm{I}_2(\underline{\mathfrak{m}}^!_{\epsilon[s]})$ in \cite{Bo94a} extends this stable orbital integral to a $C^\infty$-function around $Y=0$. Hence one can replace $\beta(Y)$ by $\alpha(X)$ in Definition \ref{def:Smod} to obtain
\begin{multline*}
	S^{G^![s], \mathrm{mod}}_{M^!}\left(\exp(Y)\epsilon[s], f^{G^![s] }\right) =
	S^{G^![s]}_{M^!}\left( \exp(Y)\epsilon[s], f^{G^![s]} \right) \\
	+ i^{G^![s]}_{M^!}(\epsilon[s]) \left\| \check{\beta} \right\| \log|\alpha(X)| \cdot S^{G^![s]}_{\underline{M}^!}\left(\exp(Y)\epsilon[s], f^{G^![s]}\right) + C^\infty
\end{multline*}
as germs of functions in $Y$ (or equivalently, in $X$), where $C^\infty$ stands for a $C^\infty$ function around $X=0$.

Lemma \ref{prop:semi-regular} also implies $M_\eta = T$, hence $\Delta(\exp(Y)\epsilon, \exp(X)\tilde{\eta})$ and its adjoint extend to locally constant functions around $X=Y=0$, by the results in \S\ref{sec:descent-endoscopy}.

On the other hand, since $G_\eta$ has semisimple rank $1$ and we are in the situation
\[\begin{tikzcd}[column sep=large]
	{G^![s]}_{\epsilon[s]} \arrow[dashed, dash, r, "\text{non.\ std.}"] \arrow[dashed, dash, r, "\text{endo.}"'] &
	\overline{{G^![s]}_{\epsilon[s]}} \arrow[dashed, dash, r, "\text{std.}"] \arrow[dashed, dash, r, "\text{endo.}"'] &
	G_\eta ,
\end{tikzcd}\]
if $s \in \Endo_{\mathbf{M}^!}(\tilde{G})$ is not semi-regular, then the argument in the first part of the proof of Lemma \ref{prop:jump-t} shows that $G^![s]_{\epsilon[s]}$ is a torus, i.e.\ $\epsilon[s]$ is $G^![s]$-regular. In this case, a familiar reasoning shows that $S^{G^![s]}_{M^!}\left( \exp(Y)\epsilon[s], f^{G^![s]} \right)$ extends to a $C^\infty$ function around $Y=0$.

Fix $\kappa \in \mathfrak{R}(T)$ together with the resulting data $(\mathbf{M}^!, \epsilon, Y)$ and $\delta_\kappa(X)$. All in all, by varying $s \in \Endo_{\mathbf{M}^!}(\tilde{G})$, we have
\begin{multline}\label{eqn:Dkappa-0}
	I^{\tilde{G}, \Endo}_{\tilde{M}}\left( \mathbf{M}^!, \delta_\kappa(X), f \right) = C^\infty + \sum_{s: \text{semi-regular}} i_{M^!}(\tilde{G}, G^![s]) S^{G^![s], \mathrm{mod}}_{M^!}\left( \delta_\kappa(X)[s], f^{G^![s]} \right) \\
	- \log|\alpha(X)| \sum_{s: \text{semi-regular}} i_{M^!}(\tilde{G}, G^![s]) \left\|\check{\beta}\right\| i^{G^![s]}_{M^!}(\epsilon[s]) S^{G^![s]}_{\underline{M}^!}\left( \delta_\kappa(X)[s]^{\underline{M}^!} , f^{G^![s]} \right) \\
	= C^\infty + \underbracket{B_\kappa(X)}_{\text{the first sum}} - \log|\alpha(X)| \cdot \underbracket{D_\kappa(X)}_{\text{the second sum}}.
\end{multline}

\begin{remark}\label{rem:s-t-semi-regular}
	The sums over semi-regular $s$ are to be unfolded in the following manner. We claim that $s$ factorizes uniquely into $t \in \Endo_{\mathbf{M}^!}(\underline{\tilde{M}})$ and $u \in \Endo_{\underline{\mathbf{M}}^!}(\tilde{G})$ where $\underline{\mathbf{M}}^! := \underline{\mathbf{M}}^![t]$, so that $\mathbf{G}^![s] = \mathbf{G}^![u]$; see the discussions in \S\ref{sec:descent-orbint-Endo}. Indeed, this is because $\underline{M} = M_\alpha$ whilst the $\beta \in X^*(T^!_{\CC})$ corresponding to $\alpha$ is already a root of $G^![s]$.
	\begin{itemize}
		\item Let us indicate such a factorization by $s \mid t$. When $s \mid t$, it has been remarked that $\epsilon[s] = \epsilon[t][u]$ is $G^![s]$-equisingular as an element of $\underline{M}^!$, hence
		\[ \underline{M}^!_{\epsilon[t]} = \underline{M}^!_{\epsilon[s]} = G^![s]_{\epsilon[s]}, \quad \underline{M}^!_{\epsilon[t], \mathrm{SC}} = G^![s]_{\epsilon[s], \mathrm{SC}} \simeq \SL(2). \]
		This implies that $(\epsilon[t], T^!, M^!, \underline{M}^!)$ is a semi-regular quadruplet. Therefore $t \in \Endo_{\mathbf{M}^!}(\underline{\tilde{M}})$ is semi-regular.
		\item Conversely, given any semi-regular $t \in \Endo_{\mathbf{M}^!}(\underline{\tilde{M}})$, every $u \in \Endo_{\underline{\mathbf{M}}^!}(\tilde{G})$ gives rise to the situation \eqref{eqn:semireg-situation}. Indeed, let $s \in \Endo_{\mathbf{M}^!}(\tilde{G})$ be obtained from $t$ and $u$. Since $\epsilon[t] \in \underline{M}^!(\R)$ is $\tilde{G}$-equisingular (it corresponds to $\eta \in \underline{M}$), Lemma \ref{prop:equisingular-endo} implies $\epsilon[s] = \epsilon[t][u]$ is $G^![s]$-equisingular, i.e.\ $G^![s]_{\epsilon[s]} = \underline{M}^!_{\epsilon[s]} = \underline{M}^!_{\epsilon[t]}$. Hence $(\epsilon[s], T^!, M^!, G^![s])$ is semi-regular as well, and $s \in \Endo_{\mathbf{M}^!}(\tilde{G})$ is semi-regular.
	\end{itemize}
\end{remark}

Going back to \eqref{eqn:Dkappa-0}, for all $\kappa \in \mathfrak{R}(T)$ and semi-regular $t \in \Endo_{\mathbf{M}^!}(\underline{\tilde{M}})$, we define
\[ D_\kappa(X, t) := \sum_{s \mid t} i_{M^!}(\tilde{G}, G^![s]) \left\| \check{\beta} \right\| i^{G^![s]}_{M^!}(\epsilon[s]) S^{G^![s]}_{\underline{M}^!}\left( \delta_\kappa(X)[s]^{\underline{M}^![t]}, f^{G^![s]} \right) \]
so that $D_\kappa(X) = \sum_t D_\kappa(X, t)$. By Lemma \ref{prop:jump-m},
\begin{align*}
	D_\kappa(X, t) & = \left( Z_{\tilde{M}^\vee}^{\hat{\alpha}^*} : Z_{\underline{\tilde{M}}^\vee}^\circ \right)^{-1} \left\|\check{\alpha}\right\| \sum_{s \mid t} i_{\underline{M}^![t]}\left( \tilde{G}, G^![s]\right) S^{G^![s]}_{\underline{M}^![t]}\left(\delta_\kappa(X)[s]^{\underline{M}^![t]}, f^{G^![s]} \right) \\
	& = \left( Z_{\tilde{M}^\vee}^{\hat{\alpha}^*} : Z_{\underline{\tilde{M}}^\vee}^\circ \right)^{-1} \left\|\check{\alpha}\right\| I^{\tilde{G}, \Endo}_{\underline{\tilde{M}}}\left( \underline{\mathbf{M}}^![t], \delta_\kappa(X)[t]^{\underline{M}^![t]}, f \right) \\
	& = \left( Z_{\tilde{M}^\vee}^{\hat{\alpha}^*} : Z_{\underline{\tilde{M}}^\vee}^\circ \right)^{-1} \left\|\check{\alpha}\right\| I^{\tilde{G}, \Endo}_{\underline{\tilde{M}}}\left( \trans_{\underline{\mathbf{M}}^![t], \underline{\tilde{M}}} \left( \delta_\kappa(X)[t]^{\underline{M}^![t]} \right), f \right).
\end{align*}

By commuting parabolic induction and transfer up to $z[t]$-twist, then summing over $t$ using Lemma \ref{prop:jump-t}, we deduce
\[ D_\kappa(X) = \sum_{t: \text{semi-regular}} D_\kappa(X, t) = \left\| \check{\alpha} \right\| I^{\tilde{G}, \Endo}_{\underline{\tilde{M}}}\left( \trans_{\mathbf{M}^!, \tilde{M}}(\delta_\kappa(X))^{\underline{\tilde{M}}}, f \right). \]

Substituting into \eqref{eqn:Dkappa-0} and \eqref{eqn:gammaX-deltaX}, we arrive at
\begin{equation}\label{eqn:IGEndo-kappa}
	I^{\tilde{G}, \Endo}_{\tilde{M}}\left( \tilde{\gamma}(X), f \right) = C^\infty + \sum_{\kappa \in \mathfrak{R}(T)} B_\kappa(X) - \left\| \check{\alpha} \right\| \log|\alpha(X)| I^{\tilde{G}, \Endo}_{\underline{\tilde{M}}}\left( \tilde{\gamma}(X), f \right).
\end{equation}

Plugging \eqref{eqn:IGEndo-kappa} into the Definition \ref{def:IMEndo-mod} of $I^{\tilde{G}, \Endo, \mathrm{mod}}_{\tilde{M}}$, we obtain
\begin{equation}\label{eqn:IGEndomod-kappa}
	I^{\tilde{G}, \Endo, \mathrm{mod}}_{\tilde{M}}\left( \tilde{\gamma}(X), f \right) = C^\infty + \sum_{\kappa \in \mathfrak{R}(T)} B_\kappa(X), \quad \text{with the same $C^\infty$ term.}
\end{equation}

Now take $U \in \Sym(\mathfrak{t}_{\CC})$. For each $\kappa$, denote by $U^!$ its transportation via $\xi_{T^!, T}: T^! \rightiso T$; also set $\mathbf{M}^! = \mathbf{M}^!_\kappa$ as usual and put
\begin{align*}
	B_\kappa(X, U) & := \sum_{s: \text{semi-regular}} i_{M^!}\left( \tilde{G}, G^![s] \right) \partial_{U^!} S^{G^![s], \mathrm{mod}}_{M^!}\left( \delta_\kappa(X)[s], f^{G^![s]} \right), \\
	B(X, U) & := \sum_{\kappa \in \mathfrak{R}(T)} B_\kappa(X, U).
\end{align*}

Since the transfer factor $\Delta$ is locally constant, it commutes with $\partial_{U^!}$, and \eqref{eqn:IGEndomod-kappa} yields
\[ \partial_U I^{\tilde{G}, \Endo, \mathrm{mod}}_{\tilde{M}}\left( \tilde{\gamma}(rH_d), f \right) = C^\infty + B(rH_d, U), \quad r \to 0. \]

\begin{proof}[Proof of Theorem \ref{prop:jump-Endo} except the first equality of (iv)]
	The assertion (i) is about
	\[ \lim_{r \to 0\pm} \partial_U I^{\tilde{G}, \Endo, \mathrm{mod}}_{\tilde{M}}\left( \tilde{\gamma}(rH_d), f \right). \]
	In view of the foregoing discussions and the expression for $B_\kappa(X, U)$, it reduces to the stable counterparts (Proposition \ref{prop:jump-stable}), for each $\kappa \in \mathfrak{R}(T)$ and semi-regular $s \in \Endo_{\mathbf{M}^!}(\tilde{G})$. Recall that the $s$ which are not semi-regular give rise to $C^\infty$-functions around $r=0$, thus are immaterial for the jump relations.
	
	The assertion (ii) also reduces to the stable side, since the condition $w_d U = U$ transcribes to the stable side for each $\kappa$ and semi-regular $s$.
	
	The assertion (iii) asserts that the limits of $\partial_{\underline{U}} I^{\tilde{G}, \Endo}_{\underline{\tilde{M}}}\left( \exp(rH_c)\tilde{\eta}, f \right)$  as $r \to 0\pm$ both exist, where $\underline{U} \in \Sym(\underline{\mathfrak{t}}_{\CC})$. Let us express $I^{\tilde{G}, \Endo}_{\underline{\tilde{M}}}\left( \exp(rH_c)\tilde{\eta}, f \right)$ by stable weighted orbital integrals on various $\underline{M}^!$. As in the discussions prior to \eqref{eqn:Dkappa-0}, the $G$-equisingularity of $\eta \in \underline{M}$ and the existence of such limits for invariant stable orbital integrals for various $\underline{M}^!$ will imply the existence of both limits. Specifically, we rely on the fact that the maximal torus $\underline{T}^! \subset \underline{M}^!_{\epsilon[t]}$ which transfers to $\underline{T} \subset \underline{M}_\eta$ has only imaginary compact roots, and apply the condition $\mathrm{I}_2(\underline{\mathfrak{m}}^!_{\epsilon[t]})$ in \cite{Bo94a}.
	
	Now turn to the second equality of (iv). It is tantamount to
	\[ \lim_{r \to 0+} \partial_{\underline{U}} I^{\tilde{G}, \Endo}_{\underline{\tilde{M}}}\left( \exp(rH_c)\tilde{\eta}, f \right) = - \lim_{r \to 0-} \partial_{\underline{U}} I^{\tilde{G}, \Endo}_{\underline{\tilde{M}}}\left( \exp(rH_c)\tilde{\eta}, f \right) \]
	when $w_d U = -U$ and $\underline{U} := uUu^{-1}$. As before, this is readily reduced to the stable side for each $\kappa$ and semi-regular $s$.
\end{proof}

We proceed to prove the first equality of Theorem \ref{prop:jump-Endo} (iv). Some extra preparations are in order.

Consider $\kappa \in \mathfrak{R}(T)$ and the data $(\mathbf{M}^!, T^!, \epsilon)$ which depends on $\kappa$. For a semi-regular $s \in \Endo_{\mathbf{M}^!}(\tilde{G})$, we denote by $\beta(s)$ (resp.\ $\check{\beta}(s)$) the root (resp.\ coroot) of $\underline{M}^!$ corresponding to $\alpha$. Take the unique $t \in \Endo_{\mathbf{M}^!}(\underline{\tilde{M}})$ such that $s \mid t$, and form $\underline{\mathbf{M}}^! := \underline{\mathbf{M}}^![t]$.

The maximal torus $T_c \subset G_{\eta, \mathrm{SC}} = \underline{M}_{\eta, \mathrm{SC}}$ extends to $\underline{T} \subset \underline{M}_\eta$. By Remark \ref{rem:s-t-semi-regular}, $t$ is automatically semi-regular. Claim:
\begin{itemize}
	\item $\underline{T}$ is the transfer of a maximal torus $\underline{T}^! \subset \underline{M}^!$, and there is a diagram joining $\epsilon[t] \in \underline{T}^!$ to $\eta \in T$ with respect to $\underline{\mathbf{M}}^!$;
	\item the Cayley transforms $C$ and $C^!$ can be chosen to render
	\[\begin{tikzcd}
		T^!_{\CC} \arrow[r, "\sim"] \arrow[d, "{C^!}"'] & T_{\CC} \arrow[d, "C"] \\
		\underline{T}^!_{\CC} \arrow[r, "\sim"'] & \underline{T}_{\CC}
	\end{tikzcd} \quad \begin{tikzcd}
		\Sym(\mathfrak{t}^!_{\CC}) \arrow[r, "\sim"] \arrow[d, "{C^!}"'] & \Sym(\mathfrak{t}_{\CC}) \arrow[d, "C"] \\
		\Sym(\underline{\mathfrak{t}}^!_{\CC}) \arrow[r, "\sim"'] & \Sym(\underline{\mathfrak{t}}_{\CC})
	\end{tikzcd}\]
	commutative, where the horizontal arrows are induced by $\xi_{\underline{T}^!, \underline{T}}$ and $\xi_{T^!, T}$, and defined over $\R$.
\end{itemize}

Indeed, we have seen that $\underline{M}_\eta$ identifies with $\overline{\underline{M}^!_{\epsilon[t]}}$; see also the proof of Lemma \ref{prop:jump-t}. Thus the problem reduces to an easy comparison between $\SL(2)$ and $\mathrm{PGL}(2)$; the existence of suitable diagrams is ensured by Proposition \ref{prop:diagram-relevant}.

Choose the Haar measures on $\underline{T}^!$ and $\underline{T}$ compatibly. In all what follows, we assume
\[ U \in \Sym(\mathfrak{t}_{\CC}), \quad w_d U = -U, \quad \underline{U} := C(U), \quad \underline{U}^! = C^!(U^!). \]

For all $\underline{X} \in \underline{\mathfrak{t}}_{G\text{-reg}}(\R)$ and semi-regular $t \in \Endo_{\mathbf{M}^!}(\underline{\tilde{M}})$, we also obtain $\underline{Y} \in \underline{\mathfrak{t}}^!_{G\text{-reg}}(\R)$ such that $\underline{X} = \dd\xi_{\underline{T}^!, \underline{T}}(\underline{Y})$. This is clearly analogous to the case of $X = \dd\xi_{T^!, T}(Y)$ considered earlier. Let us summarize this situation by $Y \leftrightarrow X$ and $\underline{X} \leftrightarrow \underline{Y}$.

\begin{lemma}\label{prop:equisingular-two}
	Take $\kappa \in \mathfrak{R}(T)$, the corresponding $\mathbf{M}^! = \mathbf{M}^!_\kappa \in \Endo_{\elli}(\tilde{M})$ and semi-regular $t \in \Endo_{\mathbf{M}^!}(\underline{\tilde{M}})$ (Definition \ref{def:semi-reg-s}).
	\begin{enumerate}[(i)]
		\item The resulting $\epsilon[t]$ and $\eta$ are in equisingular correspondence with respect to $\underline{\mathbf{M}}^! := \underline{\mathbf{M}^!}[t] \in \Endo_{\elli}(\underline{\tilde{M}})$ in the sense of \cite[Définition 4.2.1]{Li15}, extended in the obvious way to groups of metaplectic type.
		\item For $Y \leftrightarrow X$ and $\underline{X} \leftrightarrow \underline{Y}$ as above, the following limits exists
		\begin{align*}
			\Delta_{\mathbf{M}^!, \tilde{M}}(\epsilon, \tilde{\eta}) & = \lim_{Y \to 0} \Delta_{\mathbf{M}^!, \tilde{M}}\left( \exp(Y)\epsilon, \exp(X)\tilde{\eta} \right), \\
			\Delta_{\underline{\mathbf{M}}^!, \underline{\tilde{M}}}(\epsilon[t], \tilde{\eta}) & = \lim_{\underline{Y} \to 0} 	\Delta_{\underline{\mathbf{M}}^!, \underline{\tilde{M}}}\left( \exp(\underline{Y})\epsilon[t], \exp(\underline{X})\tilde{\eta} \right).
		\end{align*}
		\item Moreover, $\Delta_{\mathbf{M}^!, \tilde{M}}(\epsilon, \tilde{\eta}) = \Delta_{\underline{\mathbf{M}}^!, \underline{\tilde{M}}}(\epsilon[t], \tilde{\eta})$.
	\end{enumerate}
\end{lemma}
\begin{proof}
	Since $(\epsilon[t], T^!, M^!, \underline{M}^!)$ and $(\eta, T, M, \underline{M})$ are both semi-regular, the comparison between $\underline{M}_\eta$ and $\overline{\underline{M}^!_{\epsilon[t]}}$ made above implies that $\epsilon[t]$ and $\eta$ are in equisingular correspondence; see the cited definition. This proves (i).

	Consider (ii). Lemma \ref{prop:semi-regular} implies $\eta$ is regular in $M$, hence $\epsilon$ is regular in $M^!$. The existence of the first limit follows at once. In view of \cite[Définition 4.2.5]{Li15}, the second limit also exists: it defines the equisingular transfer factor in \textit{loc.\ cit.}
	
	By \textit{loc.\ cit.}, the second limit is also attained by $\Delta_{\underline{\mathbf{M}}^!, \underline{\tilde{M}}}\left( \exp(Y)\epsilon[t], \exp(X)\tilde{\eta} \right)$ as $Y \to 0$. By parabolic descent for transfer factors, it equals the limit of $\Delta_{\mathbf{M}^!, \tilde{M}}\left( \exp(Y)\epsilon, \exp(X)\tilde{\eta} \right)$. Therefore $\Delta_{\mathbf{M}^!, \tilde{M}}(\epsilon, \tilde{\eta}) = \Delta_{\underline{\mathbf{M}}^!, \underline{\tilde{M}}}(\epsilon[t], \tilde{\eta})$ and (iii) follows.
\end{proof}

Consider $\kappa \in \mathfrak{R}(T)$, $\mathbf{M}^! := \mathbf{M}^!_\kappa$ and a semi-regular $t \in \Endo_{\mathbf{M}^!}(\underline{\tilde{M}})$ and $\underline{Y} \leftrightarrow \underline{X}$. Using Lemma \ref{prop:equisingular-two}, define
\begin{align*}
	\delta_\kappa(\underline{X}, t) & := |\mathfrak{D}(T)|^{-1} \Delta_{\mathbf{M}^!, \tilde{M}}(\epsilon, \tilde{\eta})^{-1} \cdot \exp(\underline{Y})\epsilon[t] \; \in SD_{\mathrm{orb}}(\underline{M}^!) \otimes \mes(\underline{M}^\vee), \\
	\delta_{\underline{X}, U, \kappa, t} & := \partial_{\underline{U}^!} S^{\underline{M}^![t]}\left( \delta_\kappa(\underline{X}, t), \cdot \right) \; \in SD_{\mathrm{geom}}(\underline{M}^!) \otimes \mes(\underline{M}^!)^\vee.
\end{align*}

When $s \mid t$, we obtain $u \in \Endo_{\underline{\mathbf{M}}^!}(\tilde{G})$ such that $\mathbf{G}^![s] = \mathbf{G}^![u]$ (see Remark \ref{rem:s-t-semi-regular}). Set
\begin{align*}
	\underline{B}_\kappa(\underline{X}, U, t) & := \sum_{s \mid t} i_{M^!}(\tilde{G}, G^![s]) i^{G^![s]}_{M^!}(\epsilon[s]) \left\| \check{\beta}(s) \right\| \cdot \left\| \check{\alpha} \right\|^{-1} \\
	& \cdot \partial_{\underline{U}^!} S^{G^![s], \mathrm{mod}}_{\underline{M}^!} \left( \delta_\kappa(\underline{X}, t)[u], f^{G^![s]} \right), \\
	\underline{B}_\kappa(\underline{X}, U) & := \sum_{t: \text{semi-regular}} \underline{B}_\kappa(\underline{X}, U, t) \\
	& = \sum_{s: \text{semi-regular}}  i_{M^!}(\tilde{G}, G^![s]) i^{G^![s]}_{M^!}(\epsilon[s]) \left\| \check{\beta}(s) \right\| \cdot \left\| \check{\alpha} \right\|^{-1} \\
	& \cdot |\mathfrak{D}(T)|^{-1} \Delta_{\mathbf{M}^!, \tilde{M}}(\epsilon, \tilde{\eta})^{-1} \\
	& \cdot \partial_{\underline{U}^!} S^{G^![s], \mathrm{mod}}_{\underline{M}^!} \left( \exp(\underline{Y}) \epsilon[s], f^{G^![s]} \right).
\end{align*}

From Lemma \ref{prop:jump-m} we deduce
\begin{equation}\label{eqn:stable-jump-mod4}
	\underline{B}_\kappa(\underline{X}, U, t) = \left( Z_{\tilde{M}^\vee}^{\hat{\alpha}^*} : Z_{\underline{\tilde{M}}^\vee}^\circ \right)^{-1} I^{\tilde{G}, \Endo}_{\underline{\tilde{M}}}\left( \trans_{\underline{\mathbf{M}}^![t], \underline{\tilde{M}}} \left( \delta_{\underline{X}, U, \kappa, t}\right), f \right).
\end{equation}

Now plug $X = rH_d \in \mathfrak{t}_d(\R) \subset \mathfrak{t}(\R)$ into \eqref{eqn:IGEndomod-kappa} to obtain
\begin{multline}\label{eqn:stable-jump-mod1}
	\lim_{r \to 0+} \partial_U I^{\tilde{G}, \Endo, \mathrm{mod}}_{\tilde{M}}\left( \tilde{\gamma}(rH_d), f \right) - \lim_{r \to 0-} \partial_U I^{\tilde{G}, \Endo, \mathrm{mod}}_{\tilde{M}}\left( \tilde{\gamma}(rH_d), f \right) \\
	= \lim_{r \to 0+} \sum_{\kappa \in \mathfrak{R}(T)} \left( B_\kappa(rH_d, U) - B_\kappa(-rH_d, U) \right)
\end{multline}
provided that the limit on the right hand side exists, which we will show later on.

In parallel, take $\underline{X} = rH_c \in \mathfrak{t}_c(\R) \subset \underline{\mathfrak{t}}(\R)$. On the stable side, $Y = rH^!_d$ and $\underline{Y} = rH^!_c$ where $\dd\xi_{\underline{T}^!, \underline{T}}(H^!_c) = H_c$, etc. Our choice of Cayley transforms in \S\ref{sec:jump-relations} ensures that $C^!(H^!_d) = iH^!_c$.

\begin{lemma}\label{prop:stable-jump-mod2}
	For all $\kappa \in \mathfrak{R}(T)$, we have $B_\kappa(rH_d, U) - B_\kappa(-rH_d, U) - 2\pi i\left\|\check{\alpha}\right\| \underline{B}_\kappa(rH_c, U) \to 0$ as $r \to 0+$.
\end{lemma}
\begin{proof}
	Set $\mathbf{M}^! = \mathbf{M}^!_\kappa$. Using $\Delta_{\mathbf{M}^!, \tilde{M}}(\epsilon, \tilde{\eta}) = \displaystyle\lim_{Y \to 0} \Delta_{\mathbf{M}^!, \tilde{M}}\left( \exp(Y)\epsilon, \exp(X)\tilde{\eta} \right)$ (Lemma \ref{prop:equisingular-two}) and the local constancy of transfer factors, we have
	\[ \delta_\kappa(\pm rH_d) = |\mathfrak{D}(T)|^{-1} \Delta_{\mathbf{M}^!, \tilde{M}}(\epsilon, \tilde{\eta})^{-1} \exp(\pm rH_d)\epsilon \]
	as $r \to 0+$. In turn, this implies
	\begin{multline*}
		B_\kappa(rH_d, U) - B_\kappa(-rH_d, U) - 2\pi i\left\|\check{\alpha}\right\| \underline{B}_\kappa(rH_c, U) \\
		= \sum_{s: \text{semi-regular}} i_{M^!}(\tilde{G}, G^![s]) |\mathfrak{D}(T)|^{-1} \Delta_{\mathbf{M}^!, \tilde{M}}(\epsilon, \tilde{\eta})^{-1} X(s, r),
	\end{multline*}
	where
	\begin{multline*}
		X(s, r) := \partial_{U^!} S^{G^![s], \mathrm{mod}}_{M^!}\left( \exp(rH_d)\epsilon[s], f^{G^![s]} \right) \\
		- \partial_{U^!} S^{G^![s], \mathrm{mod}}_{M^!}\left( \exp(-rH^!_d)\epsilon[s], f^{G^![s]} \right) \\
		- 2\pi i\left\| \check{\beta}(s) \right\| i^{G^![s]}_{M^!}(\epsilon[s]) \partial_{\underline{U}^!} S^{G^![s]}_{\underline{M}^![t]} \left( \exp(rH^!_c)\epsilon[s], f^{G^![s]} \right)
	\end{multline*}
	with $s \mid t$ as usual. The stable jump relation in Proposition \ref{prop:jump-stable} (iv) implies $\lim_{r \to 0+} X(s, r) = 0$ for all semi-regular $s$.
\end{proof}

Define the distributions
\begin{equation}\label{eqn:stable-jump-mod3}\begin{aligned}
	\tilde{\gamma}_{\underline{X}, U, \kappa} & := \left( Z_{\tilde{M}^\vee}^{\hat{\alpha}^*} : Z_{\underline{\tilde{M}}^\vee}^\circ \right)^{-1} \sum_{\substack{t \in \Endo_{\mathbf{M}^!}(\underline{\tilde{M}}) \\ \text{semi-regular}}} \trans_{\underline{\mathbf{M}}^![t], \underline{\tilde{M}}}\left( \delta_{\underline{X}, U, \kappa, t} \right) \; \in D_{\mathrm{geom}, -}(\underline{\tilde{M}}) \otimes \mes(\underline{M})^\vee , \\
	\tilde{\gamma}_{\underline{X}, U} & := \sum_{\kappa \in \mathfrak{R}(T)} \gamma_{\underline{X}, U, \kappa}, \\
	\underline{B}(\underline{X}, U) & := \sum_{\kappa \in \mathfrak{R}(T)} \underline{B}_\kappa(\underline{X}, U) = I^{\tilde{G}, \Endo}_{\underline{\tilde{M}}}\left( \tilde{\gamma}_{\underline{X}, U}, f \right) \quad \text{by \eqref{eqn:stable-jump-mod4}.}
\end{aligned}\end{equation}

Our goal is to understand $\tilde{\gamma}_{\underline{X}, U}$ and $\underline{B}(\underline{X}, U)$.

\begin{lemma}\label{prop:equisingular-eta}
	Fix $\kappa \in \mathfrak{R}(T)$ and a semi-regular $t \in \Endo_{\mathbf{M}^!}(\underline{\tilde{M}})$ where $\mathbf{M}^! = \mathbf{M}^!_\kappa$; let $\underline{\mathbf{M}}^! := \underline{\mathbf{M}}^![t]$ be the resulting elliptic endoscopic datum of $\underline{\tilde{M}}$. Then
	\begin{enumerate}[(i)]
		\item the abelianization map $\mathrm{ab}^1_{\underline{M}_\eta}: \Hm^1(\R, \underline{M}_\eta) \to \Hm^1_{\mathrm{ab}}(\R, \underline{M}_\eta)$ is bijective;
		\item the natural homomorphism
		\[ q: \mathfrak{D}(\underline{T}) := \Hm^1(\R, \underline{T}) \to \Hm^1_{\mathrm{ab}}(\R, \underline{M}_\eta) \]
		is surjective, and its fibers have cardinality $2^{c(\eta)}$.
	\end{enumerate}
\end{lemma}
\begin{proof}
	Since $\underline{M}_{\eta, \mathrm{SC}} \simeq \SL(2)$ has trivial $\Hm^1$, by \cite[Proposition 1.6.7]{Lab99} we infer that $\mathrm{ab}^1_{\underline{M}_\eta}$ is bijective.
	
	The surjectivity of $q$ is then a general fact since $\underline{T}$ is a fundamental torus in $\underline{M}_\eta$; see \cite[10.2 Lemma]{Ko86}. Since $c(\eta) = |\mathfrak{D}(\underline{T}, G_\eta; \R)| = |\mathfrak{D}(\underline{T}, \underline{M}_\eta; \R)|$, it gives the cardinality of the fibers of $q$.
\end{proof}

Recall that $\Hm^1(\R, \underline{M})$ is trivial. For every $[y] \in \mathfrak{D}(\underline{T})$, represent it as the coboundary of some $y \in \underline{M}(\CC)$ and let $\underline{U}[y]$ be the transportation of $\underline{U}$ to $y^{-1} \underline{T} y$ via $\Ad(y)$. Now take $\underline{X} = rH_c$ and $\underline{Y} = rH^!_c$, where $r \in \R \smallsetminus \{0\}$ is close to $0$.

The stable conjugacy \cite[Définition 4.1.3]{Li15} in $\underline{\tilde{M}}$ affords $y^{-1}\tilde{\eta}y \in \rev^{-1}(y^{-1}\eta y)$.

\index{RL@$\mathfrak{R}(L)$}
For each connected reductive $\R$-group $L$, denote by $\mathfrak{R}(L)$ the Pontryagin dual of $\Hm^1_{\mathrm{ab}}(\R, L)$. The datum $(\epsilon[t], \eta, \underline{\mathbf{M}}^!)$ determines, through \cite[Lemme 5.2.1]{Li15}, the \emph{endoscopic character}
\[ \kappa_t \in \mathfrak{R}(\underline{M}_\eta) := \text{the Pontryagin dual of}\; \Hm^1_{\mathrm{ab}}(\R, \underline{M}_\eta). \]
If we replace $\eta$ by $y^{-1}\eta y$ where $[y] \in \mathfrak{D}(\underline{T})$, then $\underline{M}_\eta$ and $\underline{M}_{y^{-1}\eta y}$ are related by pure inner twists, and the characters $\kappa_t$ are compatible with the resulting isomorphism $\mathfrak{R}(\underline{M}_\eta) \simeq \mathfrak{R}(\underline{M}_{y^{-1}\eta y})$.

\begin{lemma}\label{prop:equisingular-Delta}
	For $\kappa$, $t$, $[y]$, $\kappa_t$ as above, the limits in Lemma \ref{prop:equisingular-two} satisfy
	\[ \Delta_{\underline{\mathbf{M}}^![t], \underline{\tilde{M}}}\left(\epsilon[t], y^{-1}\tilde{\eta} y\right) \cdot \Delta_{\underline{\mathbf{M}}^![t], \underline{\tilde{M}}}\left( \epsilon[t], \tilde{\eta} \right)^{-1} = \lrangle{ \kappa_t, q[y] }. \]
	This property characterizes $\kappa_t$ by varying $[y] \in \mathfrak{D}(\underline{T})$.
\end{lemma}
\begin{proof}
	The first part is the cocycle condition for equisingular transfer factors in \cite[Théorème 4.1.2 (iii)]{Li15}. The second part follows from Lemma \ref{prop:equisingular-eta}.
\end{proof}

Lemmas \ref{prop:equisingular-two}, \ref{prop:equisingular-Delta} imply
\begin{multline*}
	\trans_{\underline{\mathbf{M}}^![t], \underline{\tilde{M}}}\left( \delta_{rH_c, U, \kappa, t} \right) = \\
	|\mathfrak{D}(T)|^{-1} \sum_{[y] \in \mathfrak{D}(\underline{T})}
	\frac{\Delta_{\underline{\mathbf{M}}^![t], \underline{\tilde{M}}}\left(\epsilon[t], y^{-1}\tilde{\eta} y\right)}{\Delta_{\underline{\mathbf{M}}^![t], \underline{\tilde{M}}}\left( \epsilon[t], \tilde{\eta} \right)} \cdot
	\partial_{\underline{U}[y]} I^{\underline{\tilde{M}}}\left( \exp(ry^{-1}H_c y) (y^{-1}\tilde{\eta} y), \cdot \right)
	\\
	= |\mathfrak{D}(T)|^{-1} \sum_{[y] \in \mathfrak{D}(\underline{T})}
	\lrangle{\kappa_t, q[y]} \underbracket{\partial_{\underline{U}[y]} I^{\underline{\tilde{M}}}\left( \exp(ry^{-1}H_c y) (y^{-1}\tilde{\eta} y), \cdot \right)}_{\text{independent of}\; t, \kappa}.
\end{multline*}
Summing over $t$ and then over $\kappa$, we obtain
\begin{multline*}
	\tilde{\gamma}_{rH_c, U} = \sum_{\kappa} \sum_t \left( Z_{\tilde{M}^\vee}^{\hat{\alpha}^*} : Z_{\underline{\tilde{M}}^\vee}^\circ \right)^{-1} |\mathfrak{D}(T)|^{-1} \\
	\sum_{[y] \in \mathfrak{D}(\underline{T})}
	\lrangle{\kappa_t, q[y]} \partial_{\underline{U}[y]} I^{\underline{\tilde{M}}}\left( \exp(ry^{-1}H_c y) (y^{-1}\tilde{\eta} y), \cdot \right) \\
	= \sum_{[y] \in \mathfrak{D}(\underline{T})} d(q[y]) \partial_{\underline{U}[y]} I^{\underline{\tilde{M}}}\left( \exp(ry^{-1}H_c y) (y^{-1}\tilde{\eta} y), \cdot \right)
\end{multline*}
where
\begin{equation}\label{eqn:dqy}
	d(q[y]) := \left( Z_{\tilde{M}^\vee}^{\hat{\alpha}^*} : Z_{\underline{\tilde{M}}^\vee}^\circ \right)^{-1} |\mathfrak{D}(T)|^{-1} \sum_{\kappa \in \mathfrak{R}(T)} \sum_{\substack{t \in \Endo_{\mathbf{M}^!}(\underline{\tilde{M}}) \\ \text{semi-regular} }} \lrangle{\kappa_t, q[y]}.
\end{equation}

When $q[y] = 0$, we may and do choose the representative $y \in \underline{M}(\CC)$ of $[y] \in \mathfrak{D}(\underline{T})$ to ensure $y^{-1}\tilde{\eta} y = \tilde{\eta}$. The upcoming Lemma \ref{prop:dqy-vanishing} leads to
\begin{equation}\label{eqn:stable-jump-mod5}
	\tilde{\gamma}_{rH_c, U} = \sum_{[y] \in q^{-1}(0)} \partial_{\underline{U}[y]} I^{\underline{\tilde{M}}}\left( \exp(ry^{-1}H_c y) \tilde{\eta}, \cdot \right).
\end{equation}

\begin{proof}[Proof of the first equality of Theorem \ref{prop:jump-Endo} (iv)]
	Given \eqref{eqn:stable-jump-mod1} and Lemma \ref{prop:stable-jump-mod2}, it suffices to show that
	\[ \lim_{r \to 0+} \underline{B}(rH_c, U) = 2^{c(\eta)} \lim_{r \to 0+} \partial_{\underline{U}} I^{\tilde{G}, \Endo}_{\underline{\tilde{M}}}\left( \exp(rH_c)\tilde{\eta}, f \right). \]
	In view of \eqref{eqn:stable-jump-mod3} and \eqref{eqn:stable-jump-mod5}, we are reduced to proving
	\[ \sum_{[y] \in q^{-1}(0)} \lim_{r \to 0+} \partial_{\underline{U}[y]} I^{\tilde{G}, \Endo}_{\underline{\tilde{M}}} \left( \exp(ry^{-1}H_c y) \tilde{\eta}, f \right) = 2^{c(\eta)} \lim_{r \to 0+} \partial_{\underline{U}} I^{\tilde{G}, \Endo}_{\underline{\tilde{M}}}\left( \exp(rH_c)\tilde{\eta}, f \right). \]
	
	For each $[y] \in q^{-1}(0)$, take its representative $y$ from $W(\underline{M}_\eta, \underline{T})(\R)$ so that $y^{-1} \underline{T}y = \underline{T}$ and $y^{-1}\tilde{\eta} y = \tilde{\eta}$. To see this is possible, we may work inside $\underline{M}_{\eta, \mathrm{SC}} \simeq \SL(2)$ and recall \cite[Définition 4.1.3]{Li15}. By working in $\mathfrak{sl}(2)$, this also shows that only two cases can occur:
	\begin{itemize}
		\item $y^{-1}H_c y = H_c$, in which case $\underline{U}[y] = \underline{U}$,
		\item $y^{-1}H_c y = -H_c$, in which case $\underline{U}[y] = w_c \underline{U}$.
	\end{itemize}
	In the first case, $\partial_{\underline{U}[y]} I^{\tilde{G}, \Endo}_{\underline{\tilde{M}}} \left( \exp(ry^{-1}H_c y) \tilde{\eta}, f \right) = \partial_{\underline{U}} I^{\tilde{G}, \Endo}_{\underline{\tilde{M}}}\left( \exp(rH_c)\tilde{\eta}, f \right)$. In the second case, we infer from $w_d U = -U$ and by applying $C$ that $\underline{U}[y] = w_c \underline{U} = -\underline{U}$. Therefore
	\begin{align*}
		\lim_{r \to 0+} \partial_{\underline{U}[y]} I^{\tilde{G}, \Endo}_{\underline{\tilde{M}}} \left( \exp(ry^{-1}H_c y) \tilde{\eta}, f \right)
		& = -\lim_{r \to 0-} \partial_{\underline{U}} I^{\tilde{G}, \Endo}_{\underline{\tilde{M}}} \left( \exp(rH_c) \tilde{\eta}, f \right) \\
		& = \lim_{r \to 0+} \partial_{\underline{U}} I^{\tilde{G}, \Endo}_{\underline{\tilde{M}}} \left( \exp(rH_c) \tilde{\eta}, f \right)
	\end{align*}
	the last equality being a consequence of the second equality of Theorem \ref{prop:jump-Endo} (iv). Since $|q^{-1}(0)| = 2^{c(\eta)}$ by Lemma \ref{prop:equisingular-eta}, the proof is complete.
\end{proof}

\section{A lemma of vanishing}\label{sec:pf-jump-Endo-lemmas}
We proceed to complete the proof of Theorem \ref{prop:jump-Endo} (iv). Retain the conventions before.

Remember that our aim is to prove \eqref{eqn:stable-jump-mod5}. The first goal is to analyze the endoscopic character $\kappa_t \in \mathfrak{R}(\underline{M}_\eta)$ attached to $\kappa \in \mathfrak{R}(T)$ and a semi-regular $t \in \Endo_{\mathbf{M}^!}(\underline{\tilde{M}})$, where $\mathbf{M}^! := \mathbf{M}^!_\kappa$ as usual. Fix $s^\flat \in T^\vee$ with $(s^\flat)^2 = 1$ that gives rise to $\mathbf{M}^!$. Embed $\Endo_{\mathbf{M}^!}(\underline{\tilde{M}})$ as a subset of $s^\flat Z_{\tilde{M}^\vee}^\circ \big/ Z_{\underline{\tilde{M}}^\vee}^\circ$.

The allowed semi-regular $t$ form a coset under $Z_{\tilde{M}^\vee}^{\hat{\alpha}^*} / Z_{\underline{\tilde{M}}^\vee}^\circ$ by Lemma \ref{prop:jump-t}. On the other hand, $\mathfrak{R}(\underline{M}_\eta)$ can be identified with $\pi_0( Z_{\underline{M}_\eta^\vee}^{\Gamma_{\R}})$.

\begin{proposition}\label{prop:Z-descent-arch}
	Consider the homomorphism of ``descent'' $Z_{\tilde{M}^\vee}^\circ / Z_{\underline{\tilde{M}}^\vee}^\circ \to Z_{M_\eta^\vee}^{\Gamma_{\R}} / Z_{\underline{M}_\eta^\vee}^{\Gamma_{\R}}$, cf.\ \eqref{eqn:Z-descent} or \cite[p.362]{MW16-1}. Its kernel equals $Z_{\tilde{M}^\vee}^{\hat{\alpha}^*} / Z_{\underline{M}^\vee}^\circ$.
\end{proposition}
\begin{proof}
	Define $Z_{M^\vee}^{\check{\alpha}}$ as in \eqref{eqn:iGM-jump}. As an easy fact explained in \cite[p.1049]{MW16-2}, the ``descent'' homomorphism in the uncovered case
	\[ Z_{M^\vee} / Z_{\underline{M}^\vee} \to Z_{M_\eta^\vee}^{\Gamma_{\R}} / Z_{\underline{M}_\eta^\vee}^{\Gamma_{\R}} \]
	has kernel equal to $Z_{M^\vee}^{\check{\alpha}} / Z_{\underline{M}^\vee}$. The metaplectic version of the descent homomorphism is obtained by restriction via $Z_{\tilde{M}^\vee}^\circ = Z_{M^\vee}^\circ \subset Z_{M^\vee}$. It remains to recall Definition \ref{def:ZMvee-alpha}.
\end{proof}

\begin{lemma}\label{prop:kappa-t-variance}
	Given $\kappa \in \mathfrak{D}(T)$ and a semi-regular $t \in \Endo_{\mathbf{M}^!}(\underline{\tilde{M}})$, if $t$ is multiplied by $z \in Z_{\tilde{M}^\vee}^{\hat{\alpha}^*} / Z_{\underline{\tilde{M}}^\vee}^\circ$, then $\kappa_{zt}$ is $\kappa_t$ multiplied by the image of $z$ under the ``descent'' homomorphism (see Lemma \ref{prop:Z-descent-arch})
	\[ Z_{\tilde{M}^\vee}^{\hat{\alpha}^*} \to Z_{\underline{M}_\eta^\vee}^{\Gamma_{\R}} . \]
\end{lemma}
\begin{proof}
	Dually to $\mathfrak{D}(\underline{T}) \twoheadrightarrow \Hm^1_{\mathrm{ab}}(\R, \underline{M}_\eta)$, we have the inclusion $\mathfrak{R}(\underline{M}_\eta) \hookrightarrow \mathfrak{R}(\underline{T}) \simeq \pi_0(\underline{T}^{\vee, \Gamma_{\R}})$.
	
	To the elliptic maximal torus $\underline{T} \subset \underline{M}$, and the pair $(\mathbf{M}^![t], \underline{T}^!)$ determined by $(\kappa, t)$, we may attach the endoscopic character denoted by $\underline{\kappa}(t) \in \mathfrak{R}(\underline{T})$; see \cite[p.556]{Li11}. It can be compared with $\kappa_t$ as follows. Use Proposition \ref{prop:ss-parameters} to write
	\[ \underline{M}_\eta = U \times \Sp(W_+) \times \Sp(W_-) \]
	where $U$ is a direct product of unitary groups or $\GL$ over $\R$ or $\CC$, and $W_\pm$ accounts for the eigenvalues $\pm 1$ of $\eta$ in the symplectic factors of $\underline{M}$. Accordingly,
	\[ \underline{T} = \underline{T}_U \times \underline{T}_+ \times \underline{T}_- , \quad \underline{\kappa}(t) = \left( \underline{\kappa}(t)_U, \underline{\kappa}(t)_+, \underline{\kappa}(t)_- \right). \]
	
	Since $\Hm^1(\R, \Sp(W_\pm))$ is trivial, $\mathfrak{R}(\underline{M}_\eta) = \mathfrak{R}(U) \hookrightarrow \mathfrak{R}(\underline{T}_U)$. According to \cite[Lemme 4.2.10]{Li15}, the image of $\kappa_t$ in $\mathfrak{R}(\underline{T})$ equals $\underline{\kappa}(t)_U$.
	
	By describing $\underline{\kappa}(t)$ via $\underline{T}^\vee$ and Tate--Nakayama duality, one sees that $\underline{\kappa}(tz)$ is simply $\underline{\kappa}(t)$ multiplied by the image of $z$ in $\pi_0(Z_{\underline{M}_\eta}^{\Gamma_{\R}})$. As $\pi_0(Z_{\underline{M}_\eta^\vee}^{\Gamma_{\R}})$ = $\pi_0(Z_{U^\vee}^{\Gamma_{\R}})$, multiplication by $z$ only affects $\underline{\kappa}(t)_U$. This is the desired description of $\kappa_{zt}$.
\end{proof}

We finalize the proof of Theorem \ref{prop:jump-Endo} with the following result.

\begin{lemma}\label{prop:dqy-vanishing}
	For all $y \in \mathfrak{D}(\underline{T})$, the factor $d(q[y])$ in \eqref{eqn:dqy} satisfies
	\[ d(q[y]) = \begin{cases}
		1, & q[y] = 0 \\
		0, & q[y] \neq 0.
	\end{cases}\]
\end{lemma}
\begin{proof}
	Some preparations are in order. Thanks to \cite[1.2 Theorem]{Ko86} and \cite[Proposition 1.7.3]{Lab99}, for each connected reductive $\R$-group $L$ there are natural maps
	\[ \Hm^1(\R, L) \xrightarrow{\mathrm{ab}^1_L} \Hm^1_{\mathrm{ab}}(\R, L) \hookrightarrow \pi_0\left( Z_{L^\vee}^{\Gamma_{\R}} \right)^D \]
	where $(\cdots)^D$ means Pontryagin dual. The pairing so obtained between $\Hm^1_{\mathrm{ab}}(\R, L)$ and $\pi_0( Z_{L^\vee}^{\Gamma_{\R}})$ will still be denoted by $\lrangle{\cdot, \cdot}$. When $L_{\mathrm{der}}$ is simply connected, by \cite[1.2 Theorem]{Ko86} the composite $\Hm^1(\R, L) \to \pi_0( Z_{L^\vee}^{\Gamma_{\R}} )^D$ is surjective, thus $\Hm^1_{\mathrm{ab}}(\R, L) \to \pi_0( Z_{L^\vee}^{\Gamma_{\R}} )^D$ is bijective in this case.
	
	Now resume the proof. Let us embed $\Hm^1_{\mathrm{ab}}(\R, M_\eta)$ as a subset of $\Hm^1_{\mathrm{ab}}(\R, \underline{M}_\eta)$, and divide the proof into three cases.
	
	First off, suppose that $q[y] \in \Hm^1_{\mathrm{ab}}(\R, \underline{M}_\eta)$ but $q[y] \notin \Hm^1_{\mathrm{ab}}(\R, M_\eta)$. By Lemma \ref{prop:jump-t}, when $\kappa$ is fixed, the sum over $t$ in \eqref{eqn:dqy} runs over a coset under
	\[ \mathcal{Z} := Z_{\tilde{M}^\vee}^{\hat{\alpha}^*} / Z_{\underline{\tilde{M}}^\vee}^\circ . \]

	We claim that the $t$-sum vanishes. Using Lemma \ref{prop:kappa-t-variance} and Kottwitz's pairing $\lrangle{\cdot, \cdot}$, this is tantamount to
	\[ \sum_{z \in \mathcal{Z}} \lrangle{q[y], z} = 0 \]
	where we abuse notation to identify $z$ with its image in $\pi_0(Z_{\underline{M}_\eta^\vee}^{\Gamma_{\R}})$.

	Recall that $\eta$ is elliptic in $\underline{M}$. In view of the description of $\mathcal{Z}$ in Proposition \ref{prop:Z-descent-arch}, which takes the same form as \eqref{eqn:Z-descent}, the last part of the proof of Lemma \ref{prop:xsy} can be repeated to show that $\lrangle{q[y], \cdot}|_{\mathcal{Z}}$ is nontrivial, by replacing $G$ by $\underline{M}$ and $\Hm^1(F, \cdot)$ by $\Hm^1_{\mathrm{ab}}(\R, \cdot)$. This adaptation to $F = \R$ is given in Remark \ref{rem:Z-nonarch-arch} by assuming the surjectivity of Kottwitz's maps. The surjectivity does hold for $M_\eta$ and $\underline{M}_\eta$, since the description of centralizers in $\Sp$ (Proposition \ref{prop:ss-parameters}) and $\GL$ entail that their derived subgroups are simply connected.

	Next, suppose that $q[y] \in \Hm^1_{\mathrm{ab}}(\R, M_\eta) \smallsetminus \{0\}$. We claim that
	\[ \lrangle{\kappa_t, q[y]} = \lrangle{\kappa, q[y]} \]
	for all semi-regular $t \in \Endo_{\mathbf{M}^!}(\underline{\tilde{M}})$, where the right hand side is the perfect pairing between $\mathfrak{R}(T)$ and $\mathfrak{D}(T)$, by noting that $M_\eta = T$ (Lemma \ref{prop:semi-regular}). To see this, we take $y \in M_\eta(\CC)$ to represent $[y]$ and invoke the characterization of $\kappa_t$ in Lemma \ref{prop:equisingular-Delta}. It implies:
	\begin{align*}
		\lrangle{\kappa_t, q[y]} & = \Delta_{\underline{\mathbf{M}}^![t], \underline{\tilde{M}}}\left(\epsilon[t], y^{-1}\tilde{\eta} y\right) \cdot \Delta_{\underline{\mathbf{M}}^![t], \underline{\tilde{M}}}\left( \epsilon[t], \tilde{\eta} \right)^{-1} \\
		& = \Delta_{\mathbf{M}^!, \tilde{M}}\left(\epsilon, y^{-1}\tilde{\eta} y\right) \cdot \Delta_{\mathbf{M}^!, \tilde{M}}\left( \epsilon, \tilde{\eta} \right)^{-1} \\
		& = \lrangle{\kappa, q[y]};
	\end{align*}
	the second (resp.\ third) equality follows from the parabolic descent (resp.\ cocycle property) of metaplectic transfer factors, as $\eta \in M_{\mathrm{reg}}$.
	
	By the claim above, the terms in the $t$-sum in \eqref{eqn:dqy} does not depend on $t$; furthermore, the $\kappa$-sum vanishes since $q[y] \neq 0$.
	
	Finally, suppose that $q[y] = 0$. In this case, we have $d(q[y]) = 1$ from Lemma \ref{prop:jump-t}.
\end{proof}

\section{Construction of \texorpdfstring{$\epsilon_{\tilde{M}}$}{epsilonM}}\label{sec:construction-epsilonM-arch}
We revert to the setting of $F \in \{\R, \CC\}$. The following arguments are similar to \cite[IX.8.2, 8.3]{MW16-2}.

We fix Haar measures to get rid of the lines $\mes(\cdots)$ and their duals. We also fix Haar measures on the maximal tori of $M$, in a manner that is compatible with stable conjugacy.

Let $f \in \orbI_{\asp}(\tilde{G})$. For every maximal torus $T \subset M$, define the genuine function
\[ \varphi_{f, \tilde{T}}(\tilde{\gamma}) = \varphi_{f, \tilde{T}, \tilde{G} \supset \tilde{M}}(\tilde{\gamma}) := I^{\Endo}_{\tilde{M}}(\tilde{\gamma}, f) - I_{\tilde{M}}(\tilde{\gamma}, f), \quad \tilde{\gamma} \in \tilde{T}_{G\text{-reg}}. \]
Then $\varphi_{f, m\tilde{T}m^{-1}}(m\tilde{\gamma}m^{-1}) = \varphi_{f, \tilde{T}}(\tilde{\gamma})$ whenever $m \in M(F)$. If $T$ is not elliptic in $M$, then $\varphi_{f, \tilde{T}} = 0$ by induction, since $I^{\Endo}_{\tilde{M}}(\tilde{\gamma}, f)$ and $I_{\tilde{M}}(\tilde{\gamma}, f)$ satisfy the same descent formulas (Propositions \ref{prop:orbint-weighted-descent-Arch}, \ref{prop:descent-orbint-Endo}).

\begin{lemma}\label{prop:phiT-DE}
	For all $z \in \mathfrak{Z}(\mathfrak{g})$, we have $\varphi_{zf, \tilde{T}} = z_T \varphi_{f, \tilde{T}}$.
\end{lemma}
\begin{proof}
	Proposition \ref{prop:delta-invariant} and Corollary \ref{prop:stable-DE} imply that
	\[ \varphi_{zf, \tilde{T}} = \sum_{L \in \mathcal{L}(M)} \delta^L_M(z_L) \varphi_{f, \tilde{T},  \tilde{G} \supset \tilde{L}}. \]
	By induction, we may assume $\varphi_{f, \tilde{T}, \tilde{G} \supset \tilde{L}} = 0$ when $L \supsetneq M$. On the other hand, $\delta^M_M(z_M) = z_T$.
\end{proof}

Suppose $\eta \in T(F)$ and $\tilde{\eta} \in \rev^{-1}(\eta)$. Assume $T$ to be elliptic in $M$. The roots in $\Sigma(M_\eta, T)$ are all imaginary, and those in $\Sigma(G_\eta, T) \smallsetminus \Sigma(M_\eta, T)$ are non-imaginary since they have non-trivial restriction to $A_{M_\eta} \subset T$. We refer to \cite[p.249]{KV95} for the definitions pertaining to imaginary roots. Choose a Borel subgroup of $M_{\eta, \CC}$ containing $T_{\CC}$ and define the sign
\[ \Delta_\eta(X) := \prod_{\substack{\alpha \in \Sigma(M_\eta, T) \\ \alpha > 0}} \sgn\left( i\alpha(X) \right), \quad X \in \mathfrak{t}_{M_\eta\text{-reg}}(F). \]
If $F=\CC$, the ellipticity of $T$ forces $M = T$, in which case $\Delta_\eta = 1$.

Given the data above, define
\[ \mathfrak{t}_1 := \left\{ X \in \mathfrak{t}: \forall \alpha \in \Sigma(G_\eta, T) \smallsetminus \Sigma(M_\eta, T), \; \alpha(X) \neq 0 \right\}. \]

\begin{lemma}\label{prop:epsilonM-arch-bdd}
	Let $f$, $T$ and $\tilde{\eta}$ be as above, where $T \subset M$ is an elliptic maximal torus.
	\begin{enumerate}[(i)]
		\item There exists an open neighborhood $\mathcal{V}$ of $0$ in $\mathfrak{t}(F)$ such that
		\[ X \mapsto \Delta_\eta(X) \varphi_{f, \tilde{T}}(\exp(X) \tilde{\eta}), \quad X \in \mathfrak{t}_{G_\eta\text{-reg}}(F) \]
		extends to a $C^\infty$ function on $\mathfrak{t}_1(\R) \cap \mathcal{V}$.

		\item Furthermore, for all connected component $\Omega$ of $\mathfrak{t}_1(F)$ and all $U \in \Sym(\mathfrak{t}_{\CC})$, the function $\partial_U \Delta_\eta(\cdot) \varphi_{f, \tilde{T}}(\exp(\cdot) \tilde{\eta})$ extends to a bounded continuous function over $\overline{\Omega} \cap \mathcal{V}$.
	\end{enumerate}
\end{lemma}
\begin{proof}
	Assume $F=\R$ until the end of the proof. Set
	\begin{align*}
		\mathfrak{t}_2 & := \left\{ X \in \mathfrak{t}: \forall \alpha \in \Sigma(G_\eta, T), \; \alpha(X) \neq 0 \right\}, \\
		\mathfrak{t}_3 & := \left\{ X \in \mathfrak{t}_2: \forall w \in W(G_{\CC}, T_{\CC}) \smallsetminus \{1\}, \; w\eta = \eta \implies wX \neq X \right\}.
	\end{align*}
	Therefore $\mathfrak{t}_1 \supset \mathfrak{t}_2 \supset \mathfrak{t}_3 = \mathfrak{t}_{G_\eta\text{-reg}}$. They are complements of finite unions of proper subspaces.
	
	Let $X_0 \in \mathfrak{t}_1(\R)$ be close to $0$. Then $\eta_0 := \exp(X_0)\eta$ is $G$-equisingular: see the discussions on equisingularity in \cite[II.1.2]{MW16-1} or \cite{Ar88LB}.
	
	Near the equisingular element $\tilde{\eta}_0 := \exp(X_0)\tilde{\eta}$, we know that $I_{\tilde{M}}(\cdot, f)$ and $I_{\tilde{M}}^{\Endo}(\cdot, f)$ are both invariant orbital integrals over $\tilde{M}$. The characterization of (cuspidal) orbital integrals from Harish-Chandra or \cite{Bo94a} affords a $C^\infty$ extension of $\Delta_{\eta_0}(\cdot) \varphi_{f, \tilde{T}}(\exp(\cdot)\tilde{\eta}_0)$ around $0$. On the other hand, we claim that when $X$ is sufficiently close to $0$ relative to $X_0$,
	\begin{equation}\label{eqn:epsilonM-arch-bdd-aux}
		\Delta_{\eta_0}(X) = (\text{const}) \cdot \Delta_\eta(X + X_0).
	\end{equation}
	Indeed, as $X_0$ is close to $0$, we have $\Sigma(M_{\eta_0}, T) \subset \Sigma(M, T)$, and $\alpha(X_0) = 0$ for all $\alpha \in \Sigma(M_{\eta_0}, T)$. Thus
	\begin{align*}
		\Delta_\eta(X + X_0) & = \prod_{\substack{\alpha \in \Sigma(M_{\eta_0}, T) \\ \alpha > 0}} \sgn(i\alpha(X + X_0)) \prod_{\substack{\alpha \in \Sigma(M_\eta, T) \smallsetminus \Sigma(M_{\eta_0}, T) \\ \alpha > 0}} \underbracket{\sgn(i\alpha(X + X_0))}_{\text{const}} \\
		& = (\text{const}) \cdot \prod_{\substack{\alpha \in \Sigma(M_{\eta_0}, T) \\ \alpha > 0}} \sgn(i\alpha(X)) \\
		& = (\text{const}) \cdot \Delta_{\eta_0}(X).
	\end{align*}
	
	In view of \eqref{eqn:epsilonM-arch-bdd-aux}, the function $\Delta_\eta(\cdot) \varphi_{f, \tilde{T}}(\exp(\cdot)\tilde{\eta})$ has a $C^\infty$ extension around $X_0$. Vary $X_0$ to conclude that
	\begin{equation*}
		\Delta_\eta(\cdot) \varphi_{f, \tilde{T}}(\exp(\cdot)\tilde{\eta}) \quad \text{has $C^\infty$ extension over}\; \mathfrak{t}_1(\R) \cap \mathcal{V} ,
	\end{equation*}
	where $\mathcal{V}$ is an open neighborhood of $0$ in $\mathfrak{t}(\R)$. This settles (i).
	
	Consider any $X \in \mathfrak{t}_3(\R)$ close to $0$. Then $\exp(X)\eta$ is $G$-regular. By the Propositions \ref{prop:estimate-derivatives-orbint}, \ref{prop:estimate-derivatives-orbint-Endo}, upon shrinking $\mathcal{V}$ we have the following estimate: for every $U \in \Sym(\mathfrak{t}_{\CC})$, there exists $N_U \in \Z_{\geq 1}$ such that for all $f$, there exists $c = c(f,U) > 0$ with
	\[ \left| \partial_U \varphi_{f, \tilde{T}}(\exp(X)\tilde{\eta}) \right| \leq c\left| D^{G_\eta}(X)\right|^{-N_U} \quad X \in \mathfrak{t}_3(\R) \cap \mathcal{V}. \]

	Take generators $U_1, \ldots, U_k$ of $\Sym(\mathfrak{t}_{\CC})$ as a $\mathfrak{Z}(\mathfrak{g})$-module, $k \in \Z_{\geq 1}$. Let $N := \max_{1 \leq i \leq k} N_{U_i}$. In view of Lemma \ref{prop:phiT-DE}, we conclude that for all $U \in \Sym(\mathfrak{t}_{\CC})$ and $f$, there exists $c = c(f, U) > 0$ such that
	\[ \left| \partial_U \varphi_{f, \tilde{T}}(\exp(X)\tilde{\eta}) \right| \leq c\left| D^{G_\eta}(X)\right|^{-N} \quad X \in \mathfrak{t}_3(\R) \cap \mathcal{V}. \]
	
	We are now in a position to apply a result of Langlands \cite[Lemma 6.2]{Ar08} to conclude that the estimate above actually holds with $N=0$ and some other constant $c > 0$.
	
	Since $\Delta_\eta$ is locally constant and commutes with $\partial_U$ over $\mathfrak{t}_3(\R)$, we obtain a similar estimate
	\[ \left| \partial_U \Delta_\eta(\cdot) \varphi_{f, \tilde{T}}(\exp(X)\tilde{\eta}) \right| \leq c = c(f, U) \]
	over $\mathfrak{t}_3(\R) \cap \mathcal{V}$. By (i), the same estimate extends over $\mathfrak{t}_1(\R) \cap \mathcal{V}$.
	
	Let $\Omega$ be a connected component of $\mathfrak{t}_1(\R)$. It then follows from a classical technique of Harish-Chandra (see \cite[Lemma I-3-21]{Va77}) that $\partial_U \Delta_\eta(\cdot) \varphi_{f, \tilde{T}}(\exp(\cdot)\tilde{\eta})$ extends to a bounded continuous function on $\overline{\Omega} \cap \mathcal{V}$. This settles (ii).
	
	Now assume $F = \CC$. This forces $M=T$ and $\Delta_\eta = 1$. The same arguments carry over, and the aforementioned technique of Harish-Chandra simplifies since all roots are complex, thus defining ``walls'' in $\mathfrak{t}$ of real codimension $2$; see \textit{loc.\ cit.}
\end{proof}

\begin{proposition}\label{prop:epsilonM-arch-local}
	Let $f \in \orbI_{\asp}(\tilde{G})$, $\eta \in M(F)_{\mathrm{ss}}$ and $\tilde{\eta} \in \rev^{-1}(\eta)$. There exist an open neighborhood $\mathcal{U}$ of $\eta$ in $M(F)$ and $\varphi \in \orbI_{\mathrm{cusp}, \asp}(\tilde{M})$, such that $I^{\tilde{M}}(\tilde{\gamma}, \varphi) = I^{\Endo}_{\tilde{M}}(\tilde{\gamma}, f) - I_{\tilde{M}}(\tilde{\gamma}, f)$ for all $\tilde{\gamma} \in D_{\text{orb}, G\text{-reg},-}(\tilde{M})$ such that $\Supp(\tilde{\gamma})$ intersects $\rev^{-1}\mathcal{U}$.
\end{proposition}
\begin{proof}
	The goal is to show that over some invariant neighborhood of $\tilde{\eta}$, there exists $\varphi \in \orbI_{\asp}(\tilde{M})$ such that $(\varphi_{f, \tilde{T}})_T$ arises from $I^{\tilde{M}}(\cdot, \varphi)$, where $T$ ranges over conjugacy classes of maximal tori of $M$. By descent (see \S\ref{sec:descent-HA}), it suffices to consider those conjugacy classes with a representative $T \subset M_\eta$. If such a $\varphi$ exists, it is unique as a germ around $\rev^{-1}(\eta)$, and cuspidal since $\varphi_{f, \tilde{T}} = 0$ for non-elliptic $T$.
	
	Consider the statement that for all $T \subset M_\eta$,
	\begin{equation}\label{eqn:epsilonM-arch-aux0}
		\Delta_\eta(\cdot) \varphi_{f, \tilde{T}}(\exp(\cdot)\tilde{\eta}) \; \text{extends to a $C^\infty$ function around $0 \in \mathfrak{t}(F)$}.
	\end{equation}
	Note that the case of non-elliptic $T$ is trivial since $\varphi_{f, \tilde{T}} = 0$ in this case.

	Let us show how \eqref{eqn:epsilonM-arch-aux0} implies our goal. This is based on the characterization of orbital integrals \cite{Bo94a, Bo94b}. Indeed, after descent to Lie algebras,
	\begin{itemize}
		\item the boundedness condition $\mathrm{I}_1(\mathfrak{m}_\eta)$ in \textit{loc.\ cit.} is the content of Lemma \ref{prop:epsilonM-arch-bdd} (ii), as we work over an arbitrarily small neighborhood of $\eta$;
		\item the smoothness condition $\mathrm{I}_2(\mathfrak{m}_\eta)$ is immediate from \eqref{eqn:epsilonM-arch-aux0}, since $\Sigma(M_\eta, T)$ consists of imaginary roots, thus $\Delta_\eta$ is exactly the sign appearing in \textit{loc.\ cit.};
		\item the jump relations $\mathrm{I}_3(\mathfrak{m}_\eta)$ simplify since $\varphi_{f, \tilde{T}} = 0$ for non-elliptic $T$;
		\item the support condition $\mathrm{I}_4(\mathfrak{m}_\eta)$ can be enforced upon a cut-off around $0$;
	\end{itemize}
	
	Assume $F = \R$ until the end of the proof.
	
	In order to verify \eqref{eqn:epsilonM-arch-aux0}, according to the technique of Harish-Chandra \cite[Lemma I-3-21]{Va77} and Lemma \ref{prop:epsilonM-arch-bdd}, it suffices to extend $\partial_U \left( \Delta_\eta(\cdot) \varphi_{f, \tilde{T}}(\exp(\cdot)\tilde{\eta})\right)$ across the ``walls'' of real codimension $1$ defined by various $\alpha \in \Sigma(G_\eta, T) \smallsetminus \Sigma(M_\eta, T)$, where $U \in \Sym(\mathfrak{t}_{\CC})$. As seen before Lemma \ref{prop:epsilonM-arch-bdd}, $\alpha$ is not imaginary. If $\alpha$ is complex, the corresponding ``wall'' has codimension two, which can be excluded. Suppose $\alpha$ is real. Take $X_0 \in \mathcal{V}$ such that
	\[ \left\{ \beta \in \Sigma(G_\eta, T): \beta(X_0) = 0 \right\} = \{\alpha, -\alpha\}. \]
	Since $\Delta_\eta$ is constant near $X_0$, it suffices to extend $\partial_U \varphi_{f, \tilde{T}}(\exp(\cdot)\tilde{\eta})$ by continuity across $\alpha=0$ around any such point $X_0$.
	
	Consider any $X_1$ from the hyperplane $\alpha(\cdot) = 0$, sufficiently close to $X_0$. Define $\eta_1 := \exp(X_1)\eta$ and $\tilde{\eta}_1 := \exp(X_1)\tilde{\eta}$. These choices imply that $(\eta_1, T, M, G)$ is a semi-regular quadruplet (Definition \ref{def:semi-regular}), with $\pm\alpha$ being the corresponding roots. Define $\underline{M} := M_\alpha \supset \underline{T}$ and $H_c, H_d$ as in the discussion in \S\ref{sec:jump-relations} of jump relations. We are led to study the jumps of $\partial_U \varphi_{f, \tilde{T}}(\exp(\cdot)\tilde{\eta}_1)$.
	
	Note that $I^{\tilde{G}, \Endo}_{\underline{\tilde{M}}}(\exp(rH_d)\tilde{\eta}_1, f) = I^{\tilde{G}}_{\underline{\tilde{M}}}(\exp(rH_d)\tilde{\eta}_1, f)$ by induction, hence
	\begin{multline*}
		\left( \lim_{r \to 0+} - \lim_{r \to 0+}\right) \partial_U \varphi_{f, \tilde{T}}\left(\exp(rH_d) \tilde{\eta}_1\right) \\
		= \left( \lim_{r \to 0+} - \lim_{r \to 0+}\right) \left(\partial_U I^{\tilde{G}, \Endo}_{\tilde{M}}\left( \exp(rH_d)\tilde{\eta}_1, f \right) - \partial_U I^{\tilde{G}}_{\tilde{M}}\left( \exp(rH_d)\tilde{\eta}_1, f \right)\right) \\
		= \left( \lim_{r \to 0+} - \lim_{r \to 0+}\right) \left(\partial_U I^{\tilde{G}, \Endo, \mathrm{mod}}_{\tilde{M}}\left( \exp(rH_d)\tilde{\eta}_1, f \right) - \partial_U I^{\tilde{G}, \mathrm{mod}}_{\tilde{M}}\left( \exp(rH_d)\tilde{\eta}_1, f \right)\right).
	\end{multline*}

	We are now in a position to apply the Propositions \ref{prop:jump-relation-I} and \ref{prop:jump-Endo}: they imply that the final expression vanishes. Indeed, since $\dim A_{\underline{M}} < \dim A_M$, induction assumptions imply
	\[ \partial_{\underline{U}} I^{\tilde{G}, \Endo}_{\underline{\tilde{M}}}(\exp(rH_c)\tilde{\eta}_1, f) = \partial_{\underline{U}} I^{\tilde{G}}_{\underline{\tilde{M}}}(\exp(rH_c)\tilde{\eta}_1, f). \]
	
	Finally, consider the case $F = \CC$. This forces $M=T$ and $\Delta_\eta = 1$. As before, we are led to extend $\partial_U \varphi_{f, \tilde{T}}(\exp(\cdot)\tilde{\eta})$ across walls of real codimension $1$. However, all roots of are complex, defining walls of real codimension $2$ in $\mathfrak{t}$, hence there is no need to bother with jump relations.
\end{proof}

Next, suppose that $f \in \orbI_{\asp}(\tilde{G}, \tilde{K})$. In \S\ref{sec:matching-theta-arch} we defined
\[ {}^c \varphi(f) := {}^c \theta^{\tilde{G}, \Endo}_{\tilde{M}}(f) - {}^c \theta^{\tilde{G}}_{\tilde{M}}(f) \; \in \orbI_{\mathrm{ac}, \asp}(\tilde{M}, \tilde{K} \cap \tilde{M}). \]
It is cuspidal by Lemma \ref{prop:cphi-cuspidal}.

As in the non-Archimedean case, we assume the validity of Theorem \ref{prop:theta-matching-arch} when $M$ (resp.\ $G$) is replaced by a bigger (resp.\ proper) Levi subgroup of $G$.

\begin{proposition}\label{prop:epsilon-M-arch-prep}
	For all $f \in \orbI_{\asp}(\tilde{G}, \tilde{K})$, there exists a unique $\phi \in \orbI_{\cusp, \asp}(\tilde{M})$ such that $I^{\tilde{M}}\left( \tilde{\gamma}, \phi - {}^c \varphi(f)\right) = I^{\Endo}_{\tilde{M}}(\tilde{\gamma}, f) - I_{\tilde{M}}(\tilde{\gamma}, f)$ for all $\tilde{\gamma} \in D_{\mathrm{orb}, G\text{-reg}, -}(\tilde{M})$.
\end{proposition}
\begin{proof}
	Recall from the proof of Lemma \ref{prop:theta-matching-arch-prep} that the inductive assumptions imply
	\begin{align*}
		{}^c I^{\Endo}_{\tilde{M}}(\tilde{\gamma}, f) - {}^c I_{\tilde{M}}(\tilde{\gamma}, f) & = \sum_{L \in \mathcal{L}(M)} \left( I^{\tilde{L}, \Endo}_{\tilde{M}}\left(\tilde{\gamma}, {}^c \theta^{\tilde{G}, \Endo}_{\tilde{L}}(f) \right) - I^{\tilde{L}}_{\tilde{M}}\left(\tilde{\gamma}, {}^c \theta^{\tilde{G}}_{\tilde{L}}(f) \right) \right) \\
		& = \underbracket{I^{\tilde{M}}\left(\tilde{\gamma}, {}^c \varphi(f) \right)}_{L=M} +
		\underbracket{ I^{\Endo}_{\tilde{M}}(\tilde{\gamma}, f) - I_{\tilde{M}}(\tilde{\gamma}, f) }_{L=G} .
	\end{align*}

	By Lemma \ref{prop:epsilonM-arch-local}, the term labeled by $L=G$ is locally an invariant orbital integral over $\tilde{M}$, hence so is ${}^c I^{\Endo}_{\tilde{M}}(\cdot, f) - {}^c I_{\tilde{M}}(\cdot, f)$. However, the latter function is also of compact support modulo conjugacy. Using a partition of unity, we see that ${}^c I^{\Endo}_{\tilde{M}}(\cdot, f) - {}^c I_{\tilde{M}}(\cdot, f) = I^{\tilde{M}}(\cdot, \phi)$ for some $\phi \in \orbI_{\asp}(\tilde{M})$. Cuspidality follows from a comparison of descent formulas for ${}^c I^{\Endo}_{\tilde{M}}$ and ${}^c I_{\tilde{M}}$: see the proof of Lemma \ref{prop:cphi-cuspidal}.
	
	Substituting this back into the displayed equations, we obtain $I^{\tilde{M}}\left( \cdot, \phi - {}^c \varphi(f)\right) = I^{\Endo}_{\tilde{M}}(\cdot , f) - I_{\tilde{M}}(\tilde{\gamma}, f)$. The uniqueness is clear.
\end{proof}

We now allow the Haar measures on $G(F)$ and $M(F)$ to vary.

\begin{definition-proposition}\label{def:epsilonM-arch}
	\index{epsilonM}
	Define the linear map
	\begin{align*}
		\epsilon_{\tilde{M}}: \orbI_{\asp}(\tilde{G}, \tilde{K}) \otimes \mes(G) & \to \orbI_{\mathrm{ac}, \asp, \cusp}(\tilde{M}) \otimes \mes(M) \\
		f & \mapsto \phi - {}^c \varphi(f)
	\end{align*}
	where ${}^c \varphi(f)$ and $\phi$ are as in Proposition \ref{prop:epsilon-M-arch-prep}. It satisfies
	\[ \epsilon_{\tilde{M}}(f) \in \orbI_{\cusp, \asp}(\tilde{M}) \otimes \mes(M) + \orbI_{\mathrm{ac}, \cusp, \asp}(\tilde{M}, \tilde{K} \cap \tilde{M}) \otimes \mes(M) , \]
	and is uniquely characterized by the identity
	\[ I^{\tilde{M}}(\tilde{\gamma}, \epsilon_{\tilde{M}}(f)) = I^{\Endo}_{\tilde{M}}(\tilde{\gamma}, f) - I_{\tilde{M}}(\tilde{\gamma}, f), \quad \tilde{\gamma} \in D_{\mathrm{orb}, G\text{-reg}, -}(\tilde{M}) \otimes \mes(M)^\vee . \]
	Moreover, it fulfills the following properties:
	\begin{enumerate}[(i)]
		\item $\epsilon_{\tilde{M}}(zf) = z_M \epsilon_{\tilde{M}}(f)$ for all $z \in \mathfrak{Z}(\mathfrak{g})$;
		\item $\epsilon_{\tilde{M}}(f)$ is Schwartz (Definition \ref{def:Schwartz-arch}), and for all $\pi \in D_{\elli, -}(\tilde{M}) \otimes \mes(M)^\vee$, the Fourier transform of $\mathfrak{a}_M \ni X \mapsto I^{\tilde{M}}(\pi, X, \epsilon_{\tilde{M}}(f))$, denoted by $\mathfrak{a}_M^* \ni \lambda \mapsto I^{\tilde{M}}(\pi, \lambda, \epsilon_{\tilde{M}}(f))$, has a meromorphic continuation to $\mathfrak{a}_{M, \CC}^*$ such that
		\begin{itemize}
			\item it is rapidly decreasing on vertical strips,
			\item its poles are supported on hyperplanes of the form $\lrangle{\cdot, \check{\alpha}} = \mathrm{const}$ where $\alpha \in \Sigma(A_M)$, not necessarily finite in number.
		\end{itemize}
	\end{enumerate}
\end{definition-proposition}
\begin{proof}
	It remains to prove (i) and (ii). The assertion (i) concerns the local behavior of $\tilde{\gamma} \mapsto I^{\Endo}_{\tilde{M}}(\tilde{\gamma}, f) - I_{\tilde{M}}(\tilde{\gamma}, f)$, say around a semisimple element $\tilde{\eta}$ of $\tilde{M}$. Fix a maximal torus $T \subset M$ with $\eta \in T(F)$; Lemma \ref{prop:phiT-DE} implies that
	\[ I^{\tilde{M}}(\cdot, \epsilon_{\tilde{M}}(zf)) = z_T I^{\tilde{M}}\left( \cdot, \epsilon_{\tilde{M}}(f) \right) = I^{\tilde{M}}\left( \cdot, z_M \epsilon_{\tilde{M}}(f) \right) \]
	as germs of functions on $\tilde{T}_{G\text{-reg}}$ around $\tilde{\eta}$. Assertion (i) follows at once.
	
	Consider (ii). We have $\epsilon_{\tilde{M}}(f) = \phi -  {}^c \theta^{\tilde{G}, \Endo}_{\tilde{M}}(f) + {}^c \theta^{\tilde{G}}_{\tilde{M}}(f)$. Since $\phi \in \orbI_{\cusp, \asp}(\tilde{M}) \otimes \mes(M)$, the Fourier transform of $X \mapsto I^{\tilde{M}}(\pi, X, \phi)$ has the desired properties, as given by the classical Paley--Wiener theorem. As for ${}^c \theta^{\tilde{G}}_{\tilde{M}}(f)$ and ${}^c \theta^{\tilde{G}, \Endo}_{\tilde{M}}(f)$, it suffices to combine Propositions \ref{prop:cthetaGM-arch} and \ref{prop:cthetaGM-Endo-arch}.
\end{proof}

\section{\texorpdfstring{$\tilde{K} \cap \tilde{M}$}{K cap M}-finiteness}\label{sec:epsilon-KM-finite}
Retain the conventions from \S\ref{sec:construction-epsilonM-arch}. Let
\[ K^M := K \cap M(F), \quad \widetilde{K^M} := \rev^{-1}(K^M) \subset \tilde{M}. \]
Note that $K^M$ is a maximal compact subgroup of $M(F)$, in good position relative to $M_0(F)$.

The goal is to prove the following sequel to Definition--Proposition \ref{def:epsilonM-arch}.

\begin{theorem}\label{prop:epsilonM-KM-finite}
	The map $\epsilon_{\tilde{M}}$ takes values in $\orbI_{\mathrm{ac}, \asp, \cusp}(\tilde{M}, \widetilde{K^M}) \otimes \mes(M)$.
\end{theorem}

For the ease of exposition, we assume $F = \R$ until the end of this section. The same arguments work over $F = \CC$ after tiny modifications, and become simpler; this will be explained in Remark \ref{rem:epsilon-KM-finite-C}.

The proof below is essentially the same as \cite[IX.8.7--8.8]{MW16-2}. Since the arguments are involved, we give an overview below. The first ingredient concerns the Fourier transform of the invariant orbital integrals of $\epsilon_{\tilde{M}}(f)$, which is based on \cite[Theorem 4.1]{Ar94}.

To state the next few results, recall the sets $T_{\mathrm{disc}, -}(\tilde{L})_0 \subset T_{\mathrm{disc}, -}(\tilde{L})$ defined in \S\ref{sec:disc-parameter}, together with the measure on $\mathbb{S}^1 \backslash T_{\mathrm{disc}, -}(\tilde{L})$, where $L \in \mathcal{L}(M_0)$.

\begin{proposition}[Cf.\ {\cite[IX.8.7 Proposition]{MW16-2}}]\label{prop:Fourier-epsilonM-arch}
	Let $T$ be an elliptic maximal torus of $M$. To each $R \in \mathcal{L}(M_0)$ and $\sigma \in T_{\elli, -}(\tilde{R})_0$, one may attach a finite set $\mathcal{H}^{\tilde{G}, \Endo}_{\tilde{R}}(\sigma)$ of affine subspaces of $i\mathfrak{a}_R^*$, subject to the symmetry condition
	\[ \forall w \in W^G_0, \quad \mathcal{H}^{\tilde{G}, \Endo}_{w\tilde{R}}(w\sigma) = \mathcal{H}^{\tilde{G}, \Endo}_{\tilde{R}}(\sigma). \]
	Moreover, for all $H \in \mathcal{H}^{\tilde{G}, \Endo}_{\tilde{R}}(\sigma)$ whose vectorial part we denote by $i\vec{H}$, there exists a unique anti-genuine function $\xi_{\tilde{R}, \sigma, H}$ on $\widetilde{T}_{G\text{-reg}} \times H$ such that
	\begin{enumerate}[(i)]
		\item $\xi_{\tilde{R}, \sigma, H}$ is $C^\infty$;
		\item for all $D \in \Sym(\mathfrak{t}_{\CC})$ and $\Delta \in \Sym(\vec{H}_{\CC})$, viewed as differential operators, and for all compact subset $\Gamma \subset \tilde{T}$, there exists $d, N \in \Z_{\geq 1}$ and $c > 0$ such that
		\[ \left| D\Delta \xi_{\tilde{R}, \sigma, H}(\tilde{\gamma}, \lambda) \right| \leq \mathbf{d}_{\mathrm{reg}}(\gamma)^{-d} \left( 1 + \|\lambda\| \right)^N, \quad \tilde{\gamma} \in \Gamma \cap \tilde{G}_{\mathrm{reg}}, \; \lambda \in H, \]
		where $\|\cdot\|$ is a chosen Euclidean norm on $H$ and $\mathbf{d}_{\mathrm{reg}}$ is defined in \cite[p.1130]{MW16-2} --- it measures the distance to the walls;
		\item they satisfy the $W^G_0$-symmetry;
		\item for all $f \in \orbI_{\asp}(\tilde{G}, \tilde{K}) \otimes \mes(G)$ and $\tilde{\gamma} \in \widetilde{T}_{G\text{-reg}}$, we have
		\[ I^{\tilde{M}}\left( \tilde{\gamma}, \epsilon_{\tilde{M}}(f)\right) = \sum_{\substack{R, \sigma \\ H \in \mathcal{H}^{\tilde{G}, \Endo}_{\tilde{R}}(\sigma) }} \int_H \xi_{\tilde{R}, \sigma, H}(\tilde{\gamma}, \lambda) f_{\tilde{R}}(\sigma_\lambda) \dd \lambda. \]
	\end{enumerate}
\end{proposition}

More precisely, $\mathcal{H}^{\tilde{G}, \Endo}_{\tilde{R}}(\sigma)$ contains all $\nu + i\mathfrak{a}_L^* \subset i\mathfrak{a}_R^*$ such that $L \in \mathcal{L}(R)$, $\nu \in i\mathfrak{a}_R^*$, and $\sigma_\nu$ induces to an element of $T_{\mathrm{disc}, -}(\tilde{L})_0$; see Lemma \ref{prop:elli-disc-finiteness}.

\begin{proof}
	We give a quick sketch: the details are in \cite[IX.8.4--8.7]{MW16-2}. The key assertion is surely (iv), which expresses $I^{\Endo}_{\tilde{M}}(\tilde{\gamma}, f) - I_{\tilde{M}}(\tilde{\gamma}, f)$. Consider the case of $I_{\tilde{M}}(\tilde{\gamma}, f)$ first: we seek its \emph{Fourier transform} in the sense of Arthur, worked out in \cite[Theorem 4.1]{Ar94} for reductive groups. The proof is of an analytic nature. The required ingredients for coverings such as normalizing factors and the invariant local trace formula are supplied by \cite{Li12b}; see also \S\ref{sec:normalized-R-arch} for a review. The result à la Arthur takes the abstract form
	\[ I_{\tilde{M}}\left( \tilde{\gamma}, f \right) = \sum_{L \in \mathcal{L}(M_0)} \frac{|W^L_0|}{|W^G_0|} \int_{\mathbb{S}^1 \backslash T_{\mathrm{disc},-}(\tilde{L})} I_{\tilde{M}}(\tilde{\gamma}, \tau) f_{\tilde{L}}(\tau) \dd\tau \]
	where $I_{\tilde{M}}(\tilde{\gamma}, \tau)$ is a function with the properties of Weyl-symmetry, smoothness, etc., and $I_{\tilde{M}}(\tilde{\gamma}, z\tau) = z I_{\tilde{M}}(\tilde{\gamma}, \tau)$ for all $z \in \mathbb{S}^1$ and $\tau \in T_{\mathrm{disc}, -}(\tilde{L})$, so that the integral makes sense. By the definition in \S\ref{sec:disc-parameter}, the right hand side is a discrete sum of integrals over $i\mathfrak{a}_L^*$.

	To obtain the relatively concrete version à la \cite[IX.8.4 Proposition]{MW16-2}, take the $\lambda \in i\mathfrak{a}_L^*$ such that $\tau_{-\lambda} \in T_{\mathrm{disc}, -}(\tilde{L})_0$, and use Lemma \ref{prop:elli-disc-finiteness} to express $\tau_{-\lambda}$ as the parabolic induction of $\sigma_\nu$ where $\sigma \in T_{\elli, -}(\tilde{R})_0$, $R \in \mathcal{L}^L(M_0)$ and $\nu \in i\mathfrak{a}_R^*$. The resulting integrals take place over affine subspaces of $i\mathfrak{a}_R^*$ from sets $\mathcal{H}^{\tilde{G}}_{\tilde{R}}(\sigma)$, for various $(R, \sigma)$; these indexing sets are finite by Lemma \ref{prop:elli-disc-finiteness}.

	We also need the Fourier transform for $I^{\Endo}_{\tilde{M}}(\tilde{\gamma}, f)$. Given the corresponding results on the stable side \cite[IX.8.5]{MW16-2}, it suffices to recall the spectral transfer and its adjoint established in \cite[\S 7]{Li19} and reviewed in \S\ref{sec:spectral-transfer}. The remaining arguments are similar to \cite[IX.8.6]{MW16-2}, but simpler since they involve neither twists, $z$-extensions nor K-groups. Moreover, the geometric transfer factors here are locally constant, thus commutes with $D$.
	
	Accordingly, $\mathcal{H}^{\tilde{G}}_{\tilde{R}}(\sigma)$ is enlarged to another finite set $\mathcal{H}^{\tilde{G}, \Endo}_{\tilde{R}}(\sigma)$, on which the integrals take place.
\end{proof}

Consider an elliptic maximal torus $T \subset M$. We identify various maximal tori in $G_{\CC}$ by conjugation, up to Weyl-group actions. For every $\sigma \in T_{\mathrm{ell}, -}(\tilde{R})$, its infinitesimal character $\mu(\sigma)$ is an element of $\mathfrak{t}^*_{\CC} / W^G$ where $W^G := W(G_{\CC}, T_{\CC})$. In the same vein, every $H \in \mathcal{H}^{\tilde{G}, \Endo}_{\tilde{R}}(\sigma)$ can be identified as a real affine subspace of $\mathfrak{t}^*_{\CC}$.

\begin{proposition}[Cf.\ {\cite[IX.8.8 Proposition]{MW16-2}} (ii)]\label{prop:Fourier-exp-polynomial}
	Let $R$, $\sigma$, and $H \in \mathcal{H}^{\tilde{G}, \Endo}_{\tilde{R}}(\sigma)$ be given.
	If $\dim H > \dim A_M$, then $\xi_{\tilde{R}, \sigma, H} = 0$.
\end{proposition}
\begin{proof}
	It suffices to redo the cite proof and we only sketch the ideas. Given Proposition \ref{prop:Fourier-epsilonM-arch}, the arguments are based on
	\begin{itemize}
		\item the cuspidality of $\epsilon_{\tilde{M}}(f)$ and the fact that $\epsilon_{\tilde{M}}(zf) = z_M \epsilon_{\tilde{M}}(f)$, for all $f \in \orbI_{\asp}(\tilde{G}, \tilde{K}) \otimes \mes(G)$ and $z \in \mathfrak{Z}(\mathfrak{g})$ (Definition--Proposition \ref{def:epsilonM-arch} (i));
		\item Harish-Chandra's techniques, including the use of differential equations and jump relations --- see \cite[Part I]{Va77};
		\item elementary manipulation of exponential polynomials.
	\end{itemize}
	Among these ingredients, only the first item requires metaplectic input.
\end{proof}

For the next result, recall the Definition \ref{def:PW-spaces} of various Paley--Wiener spaces.

\begin{lemma}[Cf.\ {\cite[IX.8.9 Proposition]{MW16-2}}]\label{prop:epsilonM-KM-finite-prep}
	Fix $R \in \mathcal{L}(M_0)$ and $\sigma \in T_{\elli, -}(\tilde{R})_0$. Let $\mathcal{F}_R$ be the classical Paley--Wiener space of functions on $\mathfrak{a}_{R, \CC}^*$, embedded as a subspace of $\mathrm{PW}_{\asp}(\tilde{R})$ via the choice of $\sigma$. Use the symmetrization map to obtain the maps
	\[\begin{tikzcd}
		\mathcal{F}_R \arrow[twoheadrightarrow, r, "\mathrm{sym}"] & \mathcal{F}_R^{W^G(R)} \underset{\text{via}\; \sigma}{\subset} \mathrm{PW}_{\asp}(\tilde{R})^{W^G(R)} \arrow[hookrightarrow, r] & \mathrm{PW}_{\asp}(\tilde{G}) \\
		& & \orbI_{\asp}(\tilde{G}, \tilde{K}) \otimes \mes(G) \arrow[u, "\sim" sloped, "\mathrm{pw}"'] . 
	\end{tikzcd}\]

	Given $\varphi \in \mathcal{F}_R$, let $f_{\varphi} := \mathrm{pw}^{-1}\left(\mathrm{sym}(\varphi)\right) \in \orbI_{\asp}(\tilde{G}, \tilde{K}) \otimes \mes(G)$. Then we have
	\begin{enumerate}[(i)]
		\item if $\dim A_R < \dim A_M$, then $\epsilon_{\tilde{M}}(f_\varphi) = 0$;
		\item if $\dim A_R \geq \dim A_M$, then $\epsilon_{\tilde{M}}(f_\varphi)$ is $\widetilde{K^M}$-finite.
	\end{enumerate}
\end{lemma}

Before proving it, let us show how it implies our desideratum.

\begin{proof}[Proof of Theorem \ref{prop:epsilonM-KM-finite}]
	It suffices to observe that by the definition of $\mathrm{PW}_{\asp}(\tilde{G})$, every $f \in \orbI_{\asp}(\tilde{G}, \tilde{K}) \otimes \mes(G)$ is a finite sum of functions $f_\varphi$ for various $(R, \sigma)$ and $\varphi \in \mathcal{F}_R$.
\end{proof}

The Haar measure on $\mathfrak{a}_M$ has been chosen. Denote the classical Fourier transform of Schwartz functions by
\begin{align*}
	\mathcal{S}(i\mathfrak{a}_M^*) & \rightiso \mathcal{S}(\mathfrak{a}_M) \\
	\alpha & \mapsto \hat{\alpha}.
\end{align*}

Fix an elliptic maximal torus $T \subset M$ as before. Write $W^G := W(G_{\CC}, T_{\CC})$, and similarly for $W^M$ and $W^R$. Given $\pi \in T_{\elli, -}(\tilde{M})$, fix any representative of the $W^M$-orbit $\mu(\pi)$ in $\mathfrak{t}^*_{\CC}$, and similarly for $\mu(\sigma)$. Let $\mathfrak{Z}(\mathfrak{m})$ act on $\mathcal{S}(i\mathfrak{a}_M^*)$ by multipliers, namely
\[ (z\alpha)(\lambda) = z\left( \mu(\pi) + \lambda \right) \alpha(\lambda), \quad \lambda \in i\mathfrak{a}^*_M \]
where we regard $z$ as a $W^M$-invariant polynomial on $\mathfrak{t}^*_{\CC}$ via the Harish-Chandra isomorphism; this is independent of representatives. We deduce a $\mathfrak{Z}(\mathfrak{m})$-action on $\mathcal{S}(\mathfrak{a}_M)$, denoted by $\rho$, determined by
\[ \widehat{z\alpha} = \rho(z) \hat{\alpha}, \quad \alpha \in \mathcal{S}(i\mathfrak{a}_M^*), \; z \in \mathfrak{Z}(\mathfrak{m}). \]

With due care on the subscript ``$\mathrm{ac}$'', a routine verification as in \cite[p.1141]{MW16-2} gives
\begin{equation}\label{eqn:epsilonM-KM-finite-aux0}
	I^{\tilde{M}}(\pi, \cdot, zh) = \rho(z) I^{\tilde{M}}(\pi, \cdot, h), \quad h \in \orbI_{\mathrm{ac}, \asp}(\tilde{M}) \otimes \mes(M), \; z \in \mathfrak{Z}(\mathfrak{m}),
\end{equation}
both sides viewed as functions on $\mathfrak{a}_M$; moreover, when $h \in \orbI_{\asp}(\tilde{M}) \otimes \mes(M)$ we have
\[ I^{\tilde{M}}(\pi_\lambda, zh) = z(\mu(\pi) + \lambda) I^{\tilde{M}}(\pi_\lambda, h) . \]

Similarly, $\mathfrak{Z}(\mathfrak{r})$ acts on $\mathcal{S}(i\mathfrak{a}_R^*)$ by $(w\beta)(\lambda) = w(\mu(\sigma) + \lambda)\beta(\lambda)$, for all $w \in \mathfrak{Z}(\mathfrak{r})$ and $\beta \in \mathcal{S}(i\mathfrak{a}_R^*)$. This action leaves $\mathcal{F}_R$ invariant. It is routine to check that
\begin{equation}\label{eqn:epsilonM-KM-finite-aux1}
	f_{z_R \varphi} = z f_{\varphi}, \quad z \in \mathfrak{Z}(\mathfrak{g}), \; \varphi \in \mathcal{F}_R.
\end{equation}

\begin{proof}[Proof of Lemma \ref{prop:epsilonM-KM-finite-prep}]
	For the reader's convenience, the cited proof is reproduced and explicated below. Fix $\pi \in T_{\elli, -}(\tilde{M})_0$. Recall that $I^{\tilde{M}}(\pi, \nu, h)$ is defined as the Fourier transform of $I^{\tilde{M}}(\pi, X, h)$, for all $h$ in the image of $\epsilon_{\tilde{M}}$. It is meromorphic in $\nu \in \mathfrak{a}_{M,\CC}^*$. Define the linear form
	\[ \ell_\nu: \mathcal{F}_R \to \CC, \quad \varphi \mapsto I^{\tilde{M}}\left( \pi, \nu, \epsilon_{\tilde{M}}(f_\varphi) \right). \]

	Applying Fourier transform to \eqref{eqn:epsilonM-KM-finite-aux0}, we deduce that for all $z \in \mathfrak{Z}(\mathfrak{g})$,
	\[ I^{\tilde{M}}\left( \pi, \nu, \epsilon_{\tilde{M}}(zf_\varphi) \right) = I^{\tilde{M}}(\pi, \nu, z_M \epsilon_{\tilde{M}}(f_\varphi)) = z(\mu(\pi) + \nu) I^{\tilde{M}}(\pi, \nu, \epsilon_{\tilde{M}}(f_\varphi)). \]

	With these definitions and \eqref{eqn:epsilonM-KM-finite-aux1}, we obtain
	\begin{equation*}
		\ell_\nu(z_R \varphi) = z(\mu(\pi) + \nu) \ell_\nu(\varphi).
	\end{equation*}

	Define the ideal $J_\nu := \left\{ z \in \mathfrak{Z}(\mathfrak{g}) : z(\mu(\pi) + \nu)=0 \right\}$, so that $\ell_\nu(J_\nu \mathcal{F}_R) = 0$. Also, define the following finite subset of $\mathfrak{t}_{\CC}^*$:
	\[ \mathrm{int}_\nu(\pi, \sigma) := W^G \cdot \left( W^M \mu(\pi) + \nu \right) \cap \left( W^R \mu(\sigma) + \mathfrak{a}_{R, \CC}^* \right). \]
	By the general result \cite[IV.2.5 Lemme]{MW16-1} on Paley--Wiener functions, we see that
	\begin{itemize}
		\item if $\mathrm{int}_\nu(\pi, \sigma) = \emptyset$, then $J_\nu \mathcal{F}_R = \mathcal{F}_R$, and consequently $\ell_\nu = 0$;
		\item if $\mathrm{int}_\nu(\pi, \sigma) \neq \emptyset$, then there exists $N \in \Z_{\geq 1}$ such that
		\[ J_\nu \mathcal{F}_R \supset \left\{\begin{array}{r|l}
			& \forall \lambda \in \mathfrak{a}_{R, \CC}^* , \\
			\varphi \in \mathcal{F}_R & (W^R\mu(\sigma) + \lambda) \cap W^G \cdot \left( W^M \mu(\pi) + \nu \right) \neq \emptyset \\
			& \quad \implies \mathrm{ord}_\lambda(\varphi) \geq N
		\end{array}\right\}. \]
	\end{itemize}

	For each $H \in \mathcal{H}^{\tilde{G}, \Endo}_{\tilde{R}}(\sigma)$, we have its complexification $H_{\CC} \subset \mathfrak{a}_{R, \CC}^*$ and vectorial part $\vec{H}$. Set
	\[ \mathcal{H}^{\tilde{G}, \Endo}_{\tilde{R}}(\sigma, \leq \mathfrak{a}_M) := \left\{ H \in \mathcal{H}^{\tilde{G}, \Endo}_{\tilde{R}}(\sigma): \dim H \leq \dim \mathfrak{a}_M \right\}. \]
	We claim that given any $\nu \in i\mathfrak{a}_M^*$,
	\begin{equation}\label{eqn:epsilonM-KM-finite-prep-claim}\begin{split}
		\text{If} \quad & \forall H \in \mathcal{H}^{\tilde{G}, \Endo}_{\tilde{R}}(\sigma, \leq \mathfrak{a}_M), \;
		\left\{\begin{array}{r|l}
			\lambda \in H_{\CC} & W^R \mu(\sigma) + \lambda \;\text{intersects} \\
			& W^G \cdot (W^M \mu(\pi) + \nu)
		\end{array}\right\} = \emptyset, \\
		\text{then} & \quad \ell_\nu = 0.
	\end{split}\end{equation}
	
	Let us admit the validity of \eqref{eqn:epsilonM-KM-finite-prep-claim}. Let $E$ be the set consisting of $\nu \in \mathfrak{a}^*_{M, \CC}$ satisfying
	\begin{equation}\label{eqn:epsilonM-KM-finite-aux2}
		\exists H \in \mathcal{H}^{\tilde{G}, \Endo}_{\tilde{R}}(\sigma, \leq \mathfrak{a}_M), \quad W^G \cdot\left( W^M \mu(\pi) + \nu \right) \cap \left( W^R \mu(\sigma) + H_{\CC} \right) \neq \emptyset.
	\end{equation}

	Suppose $\nu \in i\mathfrak{a}^*_M$ and $\ell_\nu \neq 0$. Since $\nu \mapsto \ell_\nu(f)$ is Schwartz over $i\mathfrak{a}^*_M$, it follows that $\ell_\nu \neq 0$ for $\nu$ in a nonempty open subset of $i\mathfrak{a}^*_M$, which must be contained in $E$ by \eqref{eqn:epsilonM-KM-finite-prep-claim}. Therefore $E$ cannot be contained in any finite union of proper affine $\CC$-subspaces in $\mathfrak{a}_{M,\CC}^*$.
	
	By construction, $E$ is the finite union over $H \in \mathcal{H}^{\tilde{G}, \Endo}_{\tilde{R}}(\sigma, \leq \mathfrak{a}_M)$, $w \in W^G$, $w' \in W^M$, $w'' \in W^R$ of the affine $\CC$-subspaces
	\[ E_{H, w, w', w''} := \left( w^{-1}\left( w'' \mu(\sigma) + H_{\CC} \right) - w' \mu(\pi) \right) \cap \mathfrak{a}^*_{M, \CC} , \]
	hence there exists $(H, w, w', w'')$ such that $E_{H, w, w', w''} = \mathfrak{a}_{M,\CC}^*$, and so
	\[ \mathfrak{a}_{M,\CC}^* \subset w^{-1} \left( w''\mu(\sigma) + H_{\CC} \right) - w'\mu(\pi). \]
	By comparing the vectorial parts, we deduce $\mathfrak{a}_{M, \CC}^* \subset w^{-1} \vec{H}_{\CC}$. Since $\dim H \leq \dim \mathfrak{a}_M$, we infer that $\mathfrak{a}_{M,\CC}^* = w^{-1}\vec{H}_{\CC}$. Therefore there exists $i\mu_H \in H$ such that
	\[ w' \mu(\pi) \equiv w^{-1}(w''\mu(\sigma) + i\mu_H) \pmod{\mathfrak{a}^*_{M, \CC}}. \]
	Since $\pi \in T_{\elli, -}(\tilde{M})_0$, the congruence actually determines $w' \mu(\pi)$ from $w^{-1}(w'' \mu(\sigma) + i\mu_H)$. Summing up, when $\ell_\nu \neq 0$,
	\begin{itemize}
		\item there exists $H \in \mathcal{H}^{\tilde{G}}_{\tilde{R}}(\sigma)$ such that $\dim H = \dim \mathfrak{a}_M$,
		\item $W^M \mu(\pi)$ is constrained in a finite set determined by $\sigma$.
	\end{itemize}

	The first property forces $\dim A_R \geq \dim A_M$, and then the second property constrains $\pi \in T_{\elli, -}(\tilde{M})_0$ in a finite set determined by $\sigma$. The $\widetilde{K^M}$-finiteness of $\epsilon_{\tilde{M}}(f_\varphi)$ follows from the $\widetilde{K^M}$-finite trace Paley--Wiener theorem.
	
	It remains to justify \eqref{eqn:epsilonM-KM-finite-prep-claim}. Under its hypotheses, there exists a polynomial function $Q$ on $\mathfrak{a}_{R, \CC}^*$ such that
	\begin{itemize}
		\item $\mathrm{ord}_\lambda(Q-1) \geq N$ for all $\lambda$ such that $W^R \mu(\sigma) + \lambda$ intersects $W^G \cdot (W^M \mu(\pi) + \nu)$,
		\item $Q$ vanishes over $H_{\CC}$, for each $H \in \mathcal{H}^{\tilde{G}, \Endo}_{\tilde{R}}(\sigma, \leq \mathfrak{a}_M)$,
		\item $Q$ is $W^G(R)$-invariant.
	\end{itemize}
	Note that the first two conditions are symmetric under $W^G(R)$.
	
	Let $\varphi \in \mathcal{F}_R$. We have $(1-Q)\varphi \in J_\nu \mathcal{F}_R$, so one deduces from $\ell_\nu(J_\nu \mathcal{F}_R) = 0$ that
	\[ \epsilon_{\tilde{M}}(f_\varphi) = \epsilon_{\tilde{M}}\left( f_{(1-Q)\varphi} \right) + \epsilon_{\tilde{M}}\left( f_{Q\varphi} \right) = \epsilon_{\tilde{M}}\left( f_{Q\varphi}\right). \]
	To show $\ell_\nu(\varphi) = 0$, it suffices to show $I^{\tilde{M}}\left(\tilde{\gamma}, \epsilon_{\tilde{M}}(f_{Q\varphi})\right) = 0$ for all $\tilde{\gamma} \in \tilde{M}_{G\text{-reg}}$. Since $\epsilon_{\tilde{M}}(f_{Q\varphi})$ is cuspidal, and all elliptic maximal tori in $M$ are conjugate to $T$ over $\R$, it suffices to deal with $\tilde{\gamma} \in \widetilde{T}_{G\text{-reg}}$. Now apply Proposition \ref{prop:Fourier-epsilonM-arch} to $f_{Q\varphi}$: the resulting formula becomes
	\begin{align*}
		I^{\tilde{M}}\left(\tilde{\gamma}, \epsilon_{\tilde{M}}(f_{Q\varphi})\right) & = \sum_{H \in \mathcal{H}^{\tilde{G}, \Endo}_{\tilde{R}}(\sigma)} \int_H \xi_{\tilde{R}, \sigma, H}(\tilde{\gamma}, \lambda) (f_{Q\varphi})_{\tilde{R}}(\sigma_\lambda) \dd\lambda \\
		& = \sum_{H \in \mathcal{H}^{\tilde{G}, \Endo}_{\tilde{R}}(\sigma)} |W^G(R)|^{-1} \sum_{w \in W^G(R)} \int_H \xi_{\tilde{R}, \sigma, H}(\tilde{\gamma}, \lambda) Q(w\lambda) \varphi(w\lambda) \dd\lambda .
	\end{align*}
	Since $Q$ is $W^G(R)$-invariant, we are reduced to showing that $\xi_{\tilde{R}, \sigma, H}(\tilde{\gamma}, \lambda) Q(\lambda) = 0$ for all $H$ and $\lambda$ as above.
	\begin{itemize}
		\item If $\dim H \leq \dim A_M$, we have $Q(\lambda) = 0$ by the choice of $Q$.
		\item If $\dim H > \dim A_M$, the vanishing of $\xi_{\tilde{R}, \sigma, H}(\tilde{\gamma}, \lambda)$ is ensured by Proposition \ref{prop:Fourier-exp-polynomial}.
	\end{itemize}
	This completes the verification of \eqref{eqn:epsilonM-KM-finite-prep-claim}.
\end{proof}

\begin{remark}\label{rem:epsilon-KM-finite-C}
	The proof above of Theorem \ref{prop:epsilonM-KM-finite} also applies to the case $F = \CC$, and actually becomes simpler. First off, the critical inputs from \S\ref{sec:construction-epsilonM-arch} work for $F = \CC$. Secondly, when $F = \CC$ we may assume $M$ is a maximal torus by the cuspidality of $\epsilon_{\tilde{M}}(f)$, so that $T = M$.
	
	The proof of Proposition \ref{prop:Fourier-epsilonM-arch} carries over verbatim. In fact, the Fourier transform of $I_{\tilde{M}}(\tilde{\gamma}, f)$ à la Arthur works over all local fields $F$ of characteristic zero; the reference \cite{Li12b} is written in this generality. The spectral transfer for $\tilde{G}$ and its adjoint reviewed in \S\ref{sec:spectral-transfer} also apply to $F = \CC$.
	
	Since $M=T$, Proposition \ref{prop:Fourier-exp-polynomial} is not needed here.	In the discussions preceding Proposition \ref{prop:Fourier-exp-polynomial}, we shall replace $(G_{\CC}, T_{\CC})$ by $(G, T)$. However, $\mathfrak{t}_{\CC}$ must be defined by regarding $\mathfrak{t}$ as an $\R$-vector space and then complexify, so that $\mathfrak{t}^*_{\CC}/W(G, T)$ accommodate all infinitesimal characters.
	
	The proof of Lemma \ref{prop:epsilonM-KM-finite-prep} also carries over verbatim. The notation there can be simplified since $M = T$, but the core arguments are the same.
\end{remark}
\chapter{Global trace formula: the geometric side}\label{sec:TF-geom}
Let $F$ be a number field. We will first review the adélic coverings of metaplectic type
\[ 1 \to \bmu_8 \to \tilde{G} \to G(\A_F) \to 1 \]
in \S\S\ref{sec:adelic-coverings}--\ref{sec:metaplectic-type-adelic}. There is a natural decomposition $\tilde{G} = \tilde{G}^1 \times A_{G, \infty}$; see \eqref{eqn:AGinfty}. Take a finite set $V$ of places containing all the ramified places for $\tilde{G}$; define $\tilde{G}_V = \tilde{G}_V^1 \times A_{G, \infty}$ accordingly. Following the paradigm of \cite{MW16-2}, the geometric side of the trace formula for $\tilde{G}$ in \cite{Li14b} will be rephrased as
\[ I_{\mathrm{geom}}(f) = \sum_{M \in \mathcal{L}(M_0)} \frac{|W^M_0|}{|W^G_0|} \sum_{\mathcal{O}} I_{\tilde{M}_V}\left(A^{\tilde{M}}(V, \mathcal{O}_V), f \right), \]
where
\begin{itemize}
	\item $f \in \orbI_{\asp}(\tilde{G}_V, \tilde{K}_V) \otimes \mes(G(F_V))$,
	\item $I_{\tilde{M}_V}(\cdot, f)$ denotes the ``semi-local'' weighted orbital integral of $f$, which decomposes into local ones through a splitting formula (see \S\ref{sec:semi-local-geom-distributions}),
	\item $\mathcal{O}_V$ ranges over semisimple conjugacy classes in $M(F_V)$,
	\item $A^{\tilde{M}}(V, \mathcal{O}_V) \in D_{\mathrm{orb}, -}(\tilde{M}_V^1, \mathcal{O}_V) \otimes \mes(M(F_V))$, viewed as the ``coefficients'' in the expansion of $I_{\mathrm{geom}}(f)$.
\end{itemize}

In the parlance of \cite{Li14b}, the above is the compressed form of $I_{\mathrm{geom}}$. The coefficients can be expressed in terms of
\begin{itemize}
	\item the pre-compressed coefficients $A^{\tilde{M}}(S, \mathcal{O})_{\elli}$ where $\mathcal{O}$ is a semisimple class in $M(F)$, and $S$ is a finite set of places, sufficiently large with respect to $\mathcal{O}$ and containing $V$;
	\item semi-local unramified weighted orbital integrals $r^{\tilde{G}}_{\tilde{M}}(\tilde{\gamma}, K_{S \smallsetminus V})$.
\end{itemize}
The precise conditions on $S$ involve the notion of $S$-admissibility. In turn, $A^{\tilde{M}}(S, \mathcal{O})_{\elli}$ can be expressed in terms of the unipotent coefficient in the trace formula of $M_\sigma$, for appropriate $\sigma \in \mathcal{O}$. These reductions are explained in \S\ref{sec:semi-local-geom-distributions} in detail.

All of $I_{\mathrm{geom}}(f)$, $I_{\tilde{M}_V}(\cdot, f)$ and $A^{\tilde{G}}(V, \mathcal{O}_V)$ have endoscopic counterparts. The matching theorems are stated in \S\ref{sec:Igeom-global-matching}. Grosso modo, the stabilization of $I_{\mathrm{geom}}(f)$ is divided into two parts. The semi-local part is the stabilization of $I_{\tilde{M}_V}(\cdot, f)$, which will be reduced to the local geometric Theorem \ref{prop:local-geometric} in \S\ref{sec:semi-local-geom-Endo}. The global part is the stabilization of the coefficients $A^{\tilde{G}}(V, \mathcal{O}_V)$, deferred to \S\ref{sec:gdesc}. The stabilization of $A^{\tilde{G}}(V, \mathcal{O}_V)$ will also play a role in the study of the spectral side; see \S\ref{sec:spec-stabilization-special}.

We remark that the system of $B$-functions $B^{\tilde{G}}$ in \S\ref{sec:local-geometric} fades away in the global context; see the proof of Proposition \ref{prop:Igeom-Endo-reindexed}.

Finally, we will extend the weighted fundamental lemma to the semi-local setting in \S\ref{sec:semi-local-WFL}. This result will be needed in \S\ref{sec:gdesc}.

\section{Adélic coverings in general}\label{sec:adelic-coverings}
Below is a summary of the abstract formalism of adélic coverings presented in \cite[\S 3.3]{Li14a}. This is a supplement to the local case in \S\ref{sec:covering}.

Let $F$ be a number field and let $G$ be a connected reductive $F$-group. Consider a central extension of locally compact groups
\[ 1 \to \bmu_m \to \tilde{G} \xrightarrow{\rev} G(\A_F) \to 1 \]
where $m \in \Z_{\geq 1}$. For each set $S$ of places of $F$, pull-back to $G(F_S)$ yields a central extension
\[ 1 \to \bmu_m \to \tilde{G}_S \xrightarrow{\rev_S} G(F_S) \to 1. \]
\index{GtildeS@$\tilde{G}_S$}
Abusing notations, we often write $\mathbf{p}_S$ instead of $\mathbf{p}$. In the special case $S = \{v\}$, we obtain a local covering $\rev_v: \tilde{G}_v \to G(F_v)$.

The following auxiliary data are fixed:
\begin{itemize}
	\item a splitting $i$ of $\rev$ over the rational points, namely a commutative diagram
	\[\begin{tikzcd}
		1 \arrow[r] & \bmu_m \arrow[r] & \tilde{G} \arrow[r, "\rev"] & G(\A_F) \arrow[r] & 1 \\
		& & & G(F) \arrow[hookrightarrow, u] \arrow[lu, "i"] &
	\end{tikzcd}\]
	realizing $G(F)$ as a discrete subgroup of $\tilde{G}$;
	\item a finite subset $V_{\mathrm{ram}}$ of places of $F$ containing $\{v: v \mid \infty \}$;
	\index{Vram@$V_{\mathrm{ram}}$}
	\item a smooth connected $\mathfrak{o}_{V_{\mathrm{ram}}}$-model of $G$, so that $K_v := G(\mathfrak{o}_v)$ is hyperspecial for each $v \notin V_{\mathrm{ram}}$
\end{itemize}
The adélic coverings are subject to the following conditions.
\begin{enumerate}[(G1)]
	\item For all $v \notin V_{\mathrm{ram}}$, there is a section $s_v: K_v \to \tilde{G}_v$ of $\rev_v$ over $K_v$, which we fix.
	\item For all $v \notin V_{\mathrm{ram}}$, the ``unramified conditions'' of \cite[Définition 3.1.1]{Li14a} hold for $\tilde{G}_v$ and $s_v$. Roughly speaking, they are group-theoretical conditions that guarantee the Satake isomorphism for $\tilde{G}_v \supset s_v(K_v)$ (see below). They also include the assumption that $m$ is coprime to the residual characteristic of $F_v$.
	\item For all neighborhood $\widetilde{\mathcal{V}} \subset \tilde{G}$ of $1$, there exists a finite set of places $S \supset V_{\mathrm{ram}}$ such that $\widetilde{\mathcal{V}} \supset s_v(K_v)$ for all $v \notin S$.
\end{enumerate}

We shall view $K_v$ as an open compact subgroup of $\tilde{G}_v$ and omit $s_v$ in what follows, for each $v \notin V_{\mathrm{ram}}$. As shown in \textit{loc.\ cit.}, for every set $S$ of places, multiplication induces
\[ \Resprod_{v \in S} \tilde{G}_v \bigg/ \mathbf{N}_S \rightiso \tilde{G}_S, \]
where the restricted product is taken with respect to $K_v$ for $v \in S \smallsetminus V_{\mathrm{ram}}$, and
\[ \mathbf{N}_S := \left\{ (z_v)_{v \in S} \in \bigoplus_{v \in S} \bmu_m : \prod_v z_v = 1 \right\}. \]
When $S$ is the set of all places, we write $\mathbf{N} = \mathbf{N}_S$ so that $\Resprod_v \tilde{G}_v / \mathbf{N} \rightiso \tilde{G}$.

We record the following implications.
\begin{itemize}
	\item These conditions are preserved under passage to Levi subgroups.
	\item For each parabolic subgroup $P = MU$ of $G$, the splitting of $\rev$ over $U(\A_F)$ is compatible with $i$ over $U(F)$; see \cite[Appendice I]{MW94}.
	\item For all $v \notin V_{\mathrm{ram}}$, define the anti-genuine spherical Hecke algebra at $v$ as
	\[ \mathcal{H}_{\asp}(K_v \backslash \tilde{G}_v / K_v) := \left\{ f \in C^\infty_{c, \asp}(\tilde{G}_v): \text{bi-invariant under}\; K_v \right\}. \]
	\index{HKGK}
	The convolution product is given by the Haar measure on $\tilde{G}_v$ specified by \eqref{eqn:measure-cover} and $\mes(K_v) = 1$. The unit is given by the function
	\[ f_{K_v} \in \mathcal{H}_{\asp}(K_v \backslash \tilde{G}_v / K_v), \quad \Supp(f_{K_v}) = \widetilde{K_v}, \; f|_{K_v} = 1. \]
	A substitute of Satake isomorphism for $\mathcal{H}_{\asp}(K_v \backslash \tilde{G}_v / K_v)$ is established in \cite[Proposition 3.2.5]{Li14a}. In particular, it is known to be a finitely generated commutative $\CC$-algebra.
	
	Of course, there are genuine versions $\mathcal{H}_{-}(K_v \backslash \tilde{G}_v / K_v)$ with similar properties.
	\item The genuine automorphic forms, automorphic representations and their factorizations into local representations of $\tilde{G}_v$ can be developed in this context. See \cite{Li14a} or \cite{MW94}.
\end{itemize}

The axioms above are verified for all adélic coverings arising from the $\mathrm{K}_2$-torsors of Brylinski--Deligne \cite{BD01}. For the adélic coverings of metaplectic type (Definition \ref{def:metaplectic-type}) encountered in this work, they are well-known.

Let $F_\infty := \prod_{v \mid \infty} F_v$. Since $A_G$ is split, we may write $A_G = \left( A_{G, \Q} \right)_F$ for a split $\Q$-torus $A_{G, \Q}$ of rank $\dim \mathfrak{a}_G$. The canonical embedding $\R = \Q_\infty \to F_\infty$ induces
\[ A_{G, \Q}(\R) \hookrightarrow A_{G, \Q}(F_\infty) = A_G(F_\infty). \]
Denote by $A_{G, \infty}$ the identity component of $A_{G, \Q}(\R)$, viewed as a Lie subgroup of $A_G(F_\infty)$. Then multiplication induces an isomorphism of topological groups
\begin{equation}\label{eqn:AGinfty}
	G(\A_F)^1 \times A_{G, \infty} \rightiso G(\A_F).
\end{equation}
\index{AGinfty@$A_{G, \infty}$}

On the other hand, let $H_{\tilde{G}}$ be the composition $\tilde{G} \xrightarrow{\rev} G(\A_F) \xrightarrow{H_G} \mathfrak{a}_G$, and set $\tilde{G}^1 := \Ker H_{\tilde{G}}$. The adélic covering splits over $A_{G, \infty}$ since $A_{G, \infty} \simeq \R^{\dim \mathfrak{a}_G}$ is simply connected. Hence multiplication also induces
\begin{equation}\label{eqn:HC-decomp-global}
	\tilde{G}^1 \times A_{G, \infty} \rightiso \tilde{G}.
\end{equation}

More generally, when $V$ is a finite set of places of $F$, containing all Archimedean ones, we may set $\tilde{G}_V^1 := \Ker H_{\tilde{G}}|_{\tilde{G}_V}$ and obtain
\begin{equation}\label{eqn:HC-decomp-semilocal}
	\tilde{G}_V^1 \times A_{G, \infty} \simeq \tilde{G}.
\end{equation}

For the given $V$, choose maximal compact subgroups $K_v \subset G(F_v)$ for every $v \in V$, in good positions relative to the chosen minimal Levi $M_0$, and let $K_V := \prod_{v \in V} K_v$. As usual, $\tilde{K}_V := \rev^{-1}(K_V)$ and $\tilde{K}_v := \rev^{-1}(K_v)$. Define the spaces
\[ C^\infty_{c, \asp}(\tilde{G}_V) \twoheadrightarrow \orbI_{\asp}(\tilde{G}_V), \quad C^\infty_{c, \asp}(\tilde{G}_V, \tilde{K}_V) \twoheadrightarrow \orbI_{\asp}(\tilde{G}_V, \tilde{K}_V) \]
as we did in Definitions \ref{def:orbI} and \ref{def:K-finite}.

For each function $f$ on $\tilde{G}$ (resp.\ $\tilde{G}_V$) and $b \in C^\infty_c(\mathfrak{a}_G)$, put
\[ f^b := b(H_{\tilde{G}}(\cdot)) f . \]
\index{fb}

\begin{definition}\label{def:concentrated-global}
	\index{concentration}
	As in the local case (Definition \ref{def:concentrated}), we say a linear functional $T: C^\infty_{c, \asp}(\tilde{G}) \to \CC$ (resp.\ $T: C^\infty_{c, \asp}(\tilde{G}_V) \to \CC$) is \emph{concentrated} at $Z \in \mathfrak{a}_G$ if $T(f) = T(f^b)$ for all $f$ and $b \in C^\infty_c(\mathfrak{a}_G)$ such that $b(Z) = 1$.
\end{definition}

\begin{definition}
	Following \cite[\S 3.5]{Li14b}, define the subspace $C^\infty_{\mathrm{ac}, \asp}(\tilde{G}_V)$ of $C^\infty(\tilde{G}_V)$ so that $f$ belongs to $C^\infty_{\mathrm{ac}, \asp}(\tilde{G}_V)$ if and only if $f^b \in C^\infty_{c, \asp}(\tilde{G}_V)$ for all $b \in C^\infty_c(\mathfrak{a}_G)$. The invariant orbital integrals $I^{\tilde{G}_V}(\tilde{\gamma}, \cdot)$ are well-defined on these spaces. Similarly, we define $C^\infty_{\mathrm{ac}, \asp}(\tilde{G}_V, \tilde{K}_V)$. Define $\orbI_{\asp, \mathrm{ac}}(\tilde{G}_V)$ and $\orbI_{\asp, \mathrm{ac}}(\tilde{G}_V, \tilde{K}_V)$ accordingly.
\end{definition}

\index{semi-local}
We will be interested in genuine (resp.\ invariant genuine) distributions, i.e.\ linear functionals on the spaces $C^\infty_{c, \asp}(\tilde{G}_V)$ (resp.\ $\orbI_{\asp}(\tilde{G}_V)$), or more commonly, their $\tilde{K}_V \times \tilde{K}_V$-finite variants. By stipulation, these linear functionals must be continuous on the Archimedean $\otimes$-slot; they will be called \emph{semi-local}. Define $D_-(\tilde{G}_V)$ to be the space of invariant semi-local genuine distributions on $\tilde{G}_V$.

To keep things canonical, the spaces of test functions (resp.\ distributions) will often get twisted by the line $\mes(G(F_V))$ (resp.\ $\mes(G(F_V))$) in the semi-local setting.

\section{Groups of metaplectic type: adélic case}\label{sec:metaplectic-type-adelic}
In the rest of this section, we take $F$ to be a number field and fix an additive character $\psi = \prod_v \psi_v$ of $F \backslash \A_F$. Consider a symplectic $F$-vector space $(W, \lrangle{\cdot|\cdot})$. For each place $v$ of $F$, we put $W_v := W \dotimes{F} F_v$ and obtain the local metaplectic covering $\rev_v: \Mp(W_v) \to \Sp(W_v)$ with respect to $\psi_v$. The aim is to construct an adélic metaplectic covering
\[ 1 \to \bmu_8 \to \Mp(W, \A_F) \xrightarrow{\rev} \Sp(W, \A_F) \to 1 \]
in the sense of \S\ref{sec:covering}, whose local components are the $\Mp(W_v)$.

Fix an $\mathfrak{o}_F$-lattice $L \subset W$. When $v \nmid \infty$, put $L_v := L \dotimes{\mathfrak{o}} \mathfrak{o}_v$. We can take
\[ V_{\mathrm{ram}} := \left\{\begin{array}{r|l}
	v: \text{place} & v \mid 2, \;\text{or}\; v \mid \infty, \\
	& \text{or}\; L_v \neq L_v^*
\end{array}\right\}, \]
where the dual lattice $L_v^*$ is defined in \S\ref{sec:metaplectic-type}. This is a finite set.
\index{Vram}

When $v \notin V_{\mathrm{ram}}$, by \S\ref{sec:metaplectic-type} we know that $\Sp(W_v, \mathfrak{o}_v)$ is hyperspecial and there is a prescribed section $s_v$ of $\rev_v$ over $\Sp(W_v, \mathfrak{o}_v)$. The adélic metaplectic covering is constructed as follows. We define
\begin{align*}
	\mathbf{N} & := \left\{ (z_v)_v \in \bigoplus_v \bmu_8 : \prod_v z_v = 1 \right\}, \\
	\Mp(W, \A_F) & := \Resprod_v \Mp(W_v) \big/ \mathbf{N},
\end{align*}
where $\Resprod_v$ is the restricted product with respect to $\Sp(W_v, \mathfrak{o}_v)$ for $v \notin V_{\mathrm{ram}}$. Then $\Resprod_v \rev_v$ factors through a surjection $\rev: \Mp(W, \A_F) \to \Sp(W, \A_F)$ with $\Ker(\rev) \simeq \bmu_8$ canonically.
\index{SpWtildeA@$\Mp(W, \A_F)$}

\begin{itemize}
	\item One can verify that $\rev: \Mp(W, \A_F) \to \Sp(W, \A_F)$ is independent of the choice of $L$.
	\item There exists a canonical splitting $i: \Sp(W, F) \to \Mp(W, \A_F)$ of $\rev$ constructed in terms of the Schrödinger model. Moreover, $i(-1)$ equals the product of those $-1 \in \Mp(W_v)$ in \eqref{eqn:minus-1-lifting}; see \cite[Corollaire 2.17]{Li11}.
	\item For each Levi subgroup $M = \prod_{i \in I} \GL(n_i) \times \Sp(W^\flat)$ of $G$, let $\tilde{M} := \rev^{-1}(M(\A_F))$ and denote by $\rev^\flat: \Mp(W^\flat, \A_F) \to \Sp(W^\flat, \A_F)$ the metaplectic covering associated to $W^\flat$. Then $\rev|_{\tilde{M}}$ is canonically isomorphic to
	\[ (\identity, \rev^\flat ): \prod_{i \in I} \GL(n_i, \A_F) \times \Mp(W^\flat, \A_F) \to \prod_{i \in I} \GL(n_i, \A_F) \times \Sp(W^\flat, \A_F). \]
	The splittings over $\GL(n_i, F)$ or over the adélic points of unipotent radicals are also compatible with the above decomposition.
\end{itemize}

The maps $\rev$ and $(\identity, \rev^\flat)$ satisfy the axioms for adélic coverings in \S\ref{sec:covering}. This motivates the following

\begin{definition}
	\index{group of metaplectic type}
	The adélic groups of metaplectic type are defined as coverings of the form
	\[ \prod_{i \in I} \GL(n_i, \A_F) \times \Mp(W, \A_F), \]
	equipped with the splitting $(\identity, i)$ over its $F$-points.
\end{definition}

As in the local case, $\rev: \Mp(W, \A_F) \to \Sp(W, \A_F)$ can be reduced to a twofold covering, which are more common in the literature.

\begin{remark}
	We will often fix a symplectic basis for $(W, \lrangle{\cdot|\cdot})$, which gives rise to a minimal Levi subgroup $M_0$. In this case, it is customary to take the $\mathfrak{o}_F$-lattice $L \subset W$ to be generated by that basis. The resulting hyperspecial subgroups $K_v$ are then in good position relative to $M_0$ for all $v \notin V_{\mathrm{ram}}$.
\end{remark}

\section{Semi-local geometric distributions}\label{sec:semi-local-geom-distributions}
Following \S\ref{sec:metaplectic-type-adelic}, we consider a covering of metaplectic type $\rev: \tilde{G} \to G(\A_F)$, including the data $(W, \lrangle{\cdot|\cdot})$, $\psi$, etc. We also fix the relevant data giving rise to a minimal Levi subgroup $M_0$, for example a symplectic basis for $W$.

Take a finite set $V \supset V_{\mathrm{ram}} := V_{\mathrm{ram}}(\tilde{G})$ of places of $F$. Note that $V_{\mathrm{ram}}$ may be chosen so that $V_{\mathrm{ram}}(\tilde{G}) \supset V_{\mathrm{ram}}(\tilde{M})$ for all $M \in \mathcal{L}(M_0)$. For each $v \in V$, choose a maximal compact subgroup $K_v \subset G(F_v)$ in good position relative to $M_0$.

We make use of the formalism of semi-local genuine distributions introduced in \S\ref{sec:adelic-coverings}, with $\tilde{K}_V \times \tilde{K}_V$-finite test functions. The subspaces $D_{\mathrm{orb}, -}(\tilde{G}_V)$, etc.\ of $D_-(\tilde{G}_V)$ (resp.\ the sets $\Pi_{\mathrm{unit}, -}(\tilde{M}_V)$, etc.) are interpreted in the obvious way.

Let $M \in \mathcal{L}(M_0)$ and $f \in \orbI_{\asp}(\tilde{G}_V, \tilde{K}_V) \otimes \mes(G(F_V))$. Following \cite[\S 4.5]{Li14b}, the invariant distributions on the geometric side of the invariant trace formula are
\[ I_{\tilde{M}_V}(\tilde{\gamma}, f) = I^{\tilde{G}_V}_{\tilde{M}_V}(\tilde{\gamma}, f), \]
where $\tilde{\gamma} \in D_{\text{orb}, -}(\tilde{M}_V) \otimes \mes(M(F_V))^\vee$. These are semi-local genuine invariant distributions, and they are concentrated at $0 \in \mathfrak{a}_G$, thus extend to $\orbI_{\asp, \mathrm{ac}}(\tilde{G}_V, \tilde{K}_V) \otimes \mes(G(F_V))$.
\index{IMVgamma@$I_{\tilde{M}_V}(\tilde{\gamma}, f)$}

As in the local case, $I_{\tilde{M}_V}(\tilde{\gamma}, f)$ also has the support property: given $f \in \orbI_{\asp}(\tilde{G}_V, \tilde{K}_V)$ and a compact subset $\Xi$ of $\mathfrak{a}_M$, there exists a compact subset $C_V \subset M(F_V)$ such that $I_{\tilde{M}_V}(\tilde{\gamma}, f) = 0$ and $H_M(\gamma) \in \Xi$ together imply $\Supp(\tilde{\gamma}) \cap \rev^{-1}(C_V) \neq \emptyset$. Cf.\ \cite[VI.1.12]{MW16-2}.

Let $V \neq \emptyset$ be a finite set of places of $F$. Consider the data
\[ R^V = (R^v)_{v \in V}, \quad R^v \subset G: \text{Levi subgroup over}\; F_v. \]
We write $L^V \in \mathcal{L}(R^V)$ to mean $L^v \supset R^v$ for all $v$. For $M \in \mathcal{L}(M_0)$, we write $M \in \mathcal{L}(R^V)$ to mean $M \supset R^v$ for all $v$. Set
\begin{align*}
	\mathfrak{a}_{L^V} & := \bigoplus_{v \in V} \mathfrak{a}_{L^v}, \\
	 \mathfrak{a}^G_{L^v} & := \mathfrak{a}_G^\perp \subset \mathfrak{a}_{L^v}, \\
	\mathfrak{a}^G_{R^V} & := \bigoplus_{v \in V} \mathfrak{a}^G_{R^v}.
\end{align*}
Therefore, $L^V \in \mathcal{L}(R^V)$ implies $\mathfrak{a}^G_{R^V} \supset \mathfrak{a}^G_{L^V}$.

For all $M, L^V \in \mathcal{L}(R^V)$, we define $d^G_{R^V}(M, L^V) \in \R_{\geq 0}$ by the following variant of \eqref{eqn:d}: let $\mathrm{diag}$ denote the diagonal embedding, and
\begin{equation}\label{eqn:d-V}
	d^G_{R^V}(M, L^V) := \begin{cases}
		\dfrac{\text{the Haar measure on}\; \mathrm{diag}(\mathfrak{a}^G_M) \oplus \mathfrak{a}^G_{L^V}}{\text{that on}\; \mathfrak{a}^G_{R^V}}, & \mathfrak{a}^G_{R^V} = \mathrm{diag}(\mathfrak{a}^G_M) \oplus \mathfrak{a}^G_{L^V,} \\
		0, & \text{otherwise.}
	\end{cases}
\end{equation}
These definitions follow Arthur's conventions and \cite[VI.1.4]{MW16-2}.

\begin{proposition}\label{prop:semi-local-orbint-splitting}
	Let $M \in \mathcal{L}(M_0)$, $f = \prod_{v \in V} f_v \in \orbI_{\asp}(\tilde{G}_V, \tilde{K}_V) \otimes \mes(G(F_V))$ and
	\[ \tilde{\gamma} = \prod_{v \in V} \tilde{\gamma}_v \in D_{\mathrm{orb}, -}(\tilde{M}_V) \otimes \mes(M(F_V))^\vee. \]
	
	Define $M^V$ by $M^v := M$ for all $v \in V$. Then we have
	\[ I_{\tilde{M}_V}(\tilde{\gamma}, f) = \sum_{L^V \in \mathcal{L}(M^V)} d^G_{M^V}\left( M, L^V \right) \prod_{v \in V} I^{\tilde{L}^v}_{\tilde{M}^v}(\tilde{\gamma}_v, f_{v, \tilde{L}^v}). \]
	For any decomposition $V = V_1 \sqcup V_2$ and the corresponding $\tilde{\gamma} = \tilde{\gamma}_1 \tilde{\gamma}_2$, $f = f_1 f_2$, we have
	\[ I_{\tilde{M}_V}(\tilde{\gamma}, f) = \sum_{L_1, L_2 \in \mathcal{L}(M)} d^G_M(L_1, L_2) I^{\tilde{L}_1}_{\tilde{M}}\left(\tilde{\gamma}_1, f_{1, \tilde{L}_1}\right) I^{\tilde{L}_2}_{\tilde{M}}\left(\tilde{\gamma}_2, f_{2, \tilde{L}_2}\right). \]
\end{proposition}
\begin{proof}
	Same as \cite[VI.1.11]{MW16-2}, based on standard properties of product $(G, M)$-families.
\end{proof}

\begin{remark}\label{rem:semi-local-orbint-ext}
	Consider the space
	\[ D_{\mathrm{geom}, G_\infty\text{-equi}, -}(\tilde{M}_V) := \bigotimes_{\substack{v \in V \\ v \mid \infty}} D_{\mathrm{geom}, G\text{-equi}, -}(\tilde{M}_v) \otimes \bigotimes_{\substack{v \in V \\ v \nmid \infty}} D_{\mathrm{geom}, -}(\tilde{M}_v). \]
	Define $I_{\tilde{M}_V}(\tilde{\gamma}, f)$ for all $\tilde{\gamma} \in D_{\mathrm{geom}, G_\infty\text{-equi}, -}(\tilde{M}_V) \otimes \mes(M(F_V))^\vee$ by the first formula of Proposition \ref{prop:semi-local-orbint-splitting}, as the definition at each $v \mid \infty$ is available by Definition \ref{def:weighted-I-integral-arch}. This is compatible with the original definition.
	
	In a similar vein, we extend $I_{\tilde{M}_V}(\tilde{\gamma}, f)$ to $\tilde{\gamma} \in D_{\text{tr-orb}, -}(\tilde{M}_V) \otimes \mes(M(F_V))^\vee$, where
	\[ D_{\text{tr-orb}, -}(\tilde{M}_V) := \bigotimes_{\substack{v \in V \\ v \mid \infty}} \underbracket{D_{\text{tr-orb}, -}(\tilde{M}_v)}_{\text{Definition \ref{def:tr-orb}}} \otimes \bigotimes_{\substack{v \in V \\ v \nmid \infty}} D_{\text{geom}, -}(\tilde{M}_v), \]
	by turning the first formula of Proposition \ref{prop:semi-local-orbint-splitting} into a definition. To ensure that $I_{\tilde{M}^v}^{\tilde{L}^v}(\tilde{\gamma}_v, \cdot)$ is defined at $v \mid \infty$ and $\tilde{\gamma}_v \in D_{\text{tr-orb}, -}(\tilde{M}_v) \otimes \mes(M(F_v))^\vee$, it suffices to make the Hypothesis \ref{hyp:ext-Arch}: see Theorem \ref{prop:ext-Arch-prog-realize}.
	\index{Dtr-orb}
	
	Under the same hypothesis, we can actually extend $I_{\tilde{M}_V}(\tilde{\gamma}, f)$ to all $\tilde{\gamma}$ in the following space
	\begin{align*}
		D_{\diamondsuit, -}(\tilde{M}_V) & := D_{\text{geom}, G_\infty\text{-equi}, -}(\tilde{M}_V) + D_{\text{tr-orb}, -}(\tilde{M}_V) \\
		& \subset D_{\mathrm{geom}, -}(\tilde{M}_V),
	\end{align*}
	which is also a $\otimes$-product over $V$. This is guaranteed by (C8) in \S\ref{sec:ext-Arch-prog}.
	\index{DdiamondMV@$D_{\diamondsuit, -}$}
\end{remark}

Next, suppose that $G^!$ is a quasisplit $F$-group, equipped with a Borel pair $(T^!, B^!)$ over $F$ and an invariant quadratic form on $X_*(T^!_{\overline{F}}) \otimes \R$. Note that $M^!_0 := T^!$ is a minimal Levi subgroup. We also choose maximal compact subgroups $K^!_v \subset G^!(F_v)$ in good positions relative to $M^!_0$, for each $v \in V$, and set $K^!_V := \prod_{v \in V} K^!_v$.

We assume that $G^!$ comes with an integral model that is unramified outside a finite set $V_{\mathrm{ram}}(G^!)$ of places, containing all Archimedean ones. Given $V \supset V_{\mathrm{ram}}(G^!)$, $M^! \in \mathcal{L}(M_0^!)$, the semi-local distributions above have stable avatars
\[ S^{G^!_V}_{M^!_V}(\delta, f^!) = S_{M^!_V}(\delta, f^!), \]
\index{SGMV-delta-weighted@$S^{G^{"!}_V}_{M^{"!}_V}(\delta, f^{"!})$}
where
\[ \delta \in SD_{\text{orb}, -}(M^!_V) \otimes \mes(M^!(F_V))^\vee, \quad f^! \in S\orbI_{\mathrm{ac}}(G^!(F_V), K^!_V) \otimes \mes(G(F_V)). \]
We also extend the definition of $S_{M^!_V}(\delta, f^!)$ to $\delta \in SD_{\diamondsuit}(M^!_V) \otimes \mes(M^!(F_V))^\vee$ where
\begin{align*}
	SD_{\diamondsuit}(M_V^!) & := SD_{\text{geom}, G^!_\infty\text{-equi}, -}(M_V^!) + SD_{\text{tr-orb}}(M_V^!), \\
	SD_{\text{tr-orb}}(M_V^!) & := \bigotimes_{\substack{v \in V \\ v \mid \infty}} SD_{\text{tr-orb}}(M^!_v) \otimes \bigotimes_{\substack{v \in V \\ v \nmid \infty}} SD_{\text{geom}}(M^!_v).
\end{align*}
as in Remark \ref{rem:semi-local-orbint-ext}, as subspaces of $SD_{\mathrm{geom}}(M_V^!)$. Thanks to \cite{MW16-1, MW16-2}, this extension is unconditional.
\index{SDtr-orb}
\index{SD-diamondsuit@$SD_{\diamondsuit}$}

The analogue of Proposition \ref{prop:semi-local-orbint-splitting} holds in the stable setting, provided that $d^G_{R^V}(M, L^V)$ is replaced by the analogue of \eqref{eqn:e}:
\[ e^{G^!}_{R^{V!}}(M^!, L^{V!}) := \begin{cases}
	k^{G^!}_{R^{V!}}(M^!, L^{V!})^{-1} d^{G^!}_{R^{V!}}(M^!, L^{V!}), & \text{if}\; d^{G^!}_{R^{V!}}(M^!, L^{V!}) \neq 0, \\
	0, & \text{otherwise}
\end{cases}\]
where following \cite[VI.4.2]{MW16-2}, we put
\[ k^{G^!}_{R^{V!}}(M^!, L^{V!}) := \# \Ker\left[ \frac{Z_{(M^!)^\vee}^{\Gamma_F}}{Z_{(G^!)^\vee}^{\Gamma_F}} \to
\prod_{v \in V} \frac{Z_{R_v^\vee}^{\Gamma_{F_v}}}{Z_{L_v^\vee}^{\Gamma_{F_v}}} \right]. \]

We can also take a system of $B$-functions for $G^!$ (Definition \ref{def:system-B-function}). The global case makes no difference, except that $V$ is so large that $v \notin V$ implies the residual characteristic of $F_v$ is coprime to all scaling factors $B_\epsilon(\beta)$. This is surely possible.

The semi-local stable distributions have the variants $S_{M^!_V}(\delta, B, f^!)$. However, we will make them disappear in the stable trace formula: see the proof of Proposition \ref{prop:Igeom-Endo-reindexed}.

\section{Geometric side of the trace formula}\label{sec:geom-side}
Keep the conventions from \S\ref{sec:semi-local-geom-distributions}. Let $\Gamma(\tilde{G}_V)$ (resp.\ $\Gamma(\tilde{G}_V^1)$) be the set of conjugacy classes in $\tilde{G}_V$ (resp.\ in $\tilde{G}_V^1$); see \eqref{eqn:HC-decomp-semilocal}. If $S \supset V$ are finite sets of places, then $F^V_S := \prod_{v \in S \smallsetminus V} F_v$ and $\tilde{G}^V_S := \rev^{-1}\left( G(F^V_S) \right)$. Define $K^V_S$ similarly, which embeds as a compact open subgroup of $\tilde{G}^V_S$.
\index{GtildeVS@$\tilde{G}^V_S$}

\subsection{Distributions and coefficients}
Fix an invariant quadratic form on $X^*(M_0)$ as in Definition \ref{def:invariant-quadratic-form}, noting that $M_0$ is a split maximal $F$-torus. Using $H_G: A_{G, \infty} \rightiso \mathfrak{a}_G$, we equip $A_{G, \infty}$ with the Haar measure induced by the invariant positive-definite quadratic form on $\mathfrak{a}_0$. Therefore, the lines $\mes(G(F_V))$, $\mes(G(F_V)/A_{G, \infty})$ and $\mes(G(F_V)^1)$ of measures are identified; similarly for $G(F_V)$ replaced by $G(\A_F)$.

In \cite[Théorème 6.1]{Li14b}, the geometric side of the invariant trace formula is phrased as the invariant genuine distribution $I_{\mathrm{geom}} = I^{\tilde{G}}_{\mathrm{geom}}$ coming with the expansion
\begin{equation}\label{eqn:Igeom-orig}\begin{aligned}
	I_{\mathrm{geom}}(f^1) & := \sum_{M \in \mathcal{L}(M_0)} \frac{|W^M_0|}{|W^G_0|} \sum_{\tilde{\gamma} \in \Gamma(\tilde{M}^1, V)} a^{\tilde{M}}(\tilde{\gamma}) I_{\tilde{M}_V}(\tilde{\gamma}, f^1), \\
	f^1 & \in \orbI_{\asp}(\tilde{G}_V/A_{G, \infty} ,\tilde{K}_V) \otimes \mes(G(F_V)) \\
	& \subset \orbI_{\asp, \mathrm{ac}}(\tilde{G}_V,\tilde{K}_V) \otimes \mes(G(F_V)).
\end{aligned}\end{equation}
The sum is finite when $f^1$ is given; see the support property mentioned in the beginning of \S\ref{sec:semi-local-geom-distributions}..
\index{Igeom@$I_{\mathrm{geom}}$}

The detailed definitions of the subset $\Gamma(\tilde{M}^1, V)$ of $\Gamma(\tilde{M}_V^1)$ and the coefficients $a^{\tilde{M}}(\tilde{\gamma})$ can be found in \cite[\S 5.2]{Li14b}. Roughly speaking, $\Gamma(\tilde{M}^1, V)$ consists of conjugacy classes of global origin. On the other hand, $a^{\tilde{M}}(\tilde{\gamma})$ are certain ``compressed'' forms of the geometric coefficients from \cite[\S 6.5]{Li14a}. The relevant Haar measures were chosen in \textit{loc.\ cit.} To reconcile with the conventions here, it is more appropriate to view $\Gamma(\tilde{M}^1, V)$ as a subset of $D_{\text{orb}, -}(\tilde{M}_V) \otimes \mes(M(F_V))^\vee$.

Before explicating the coefficients which are of global nature, we give an easy rephrasing of \eqref{eqn:Igeom-orig} in a way that is closer to \cite{Ar88-2, MW16-2}.

\begin{enumerate}
	\item First, $I_{\tilde{M}_V}(\tilde{\gamma}, \cdot)$ is concentrated at $0 \in \mathfrak{a}_G$ for all $\tilde{\gamma} \in \Gamma(\tilde{M}^1, V)$, thus \eqref{eqn:Igeom-orig} extends automatically from $f^1$ to all $f^\dagger \in \orbI_{\asp, \mathrm{ac}}(\tilde{G}_V, \tilde{K}_V) \otimes \mes(G(F_V))$, namely $I_{\mathrm{geom}}(f^\dagger) = I_{\mathrm{geom}}(f^{\dagger, b})$ for any $b \in C^\infty_c(\mathfrak{a}_G)$ with $b(0)=1$. This step is in \cite[Lemme 6.3]{Li14b}.
	
	In particular, the test functions can be drawn from $\orbI_{\asp}(\tilde{G}_V, \tilde{K}_V) \otimes \mes(G(F_V))$.
	
	\item The $I_{\tilde{M}_V}(\tilde{\gamma}, f^\dagger)$ so obtained depends only on $f^\dagger|_{\tilde{G}_V^1}$. To see this, note that $\Supp I_{\tilde{M}_V}(\tilde{\gamma}, \cdot) \subset \tilde{G}_V^1$, and apply the description of such distributions in \cite[Theorem 2.3.5]{Ho03} locally, by recalling Definition \ref{def:concentrated-global} and that $\tilde{G}_V = \tilde{G}^1_V \times A_{G, \infty}$ to rule out derivatives in the normal direction.
	
	\item Summing up, the study of $I_{\mathrm{geom}}(f^1)$ is equivalent to that of
	\begin{equation}\label{eqn:Igeom-0}\begin{aligned}
		I_{\mathrm{geom}}(f) & := \sum_{M \in \mathcal{L}(M_0)} \frac{|W^M_0|}{|W^G_0|} \sum_{\tilde{\gamma} \in \Gamma(\tilde{M}^1, V)} a^{\tilde{M}}(\tilde{\gamma}) I_{\tilde{M}_V}(\tilde{\gamma}, f), \\
		f & \in \orbI_{\asp}(\tilde{G}_V ,\tilde{K}_V) \otimes \mes(G(F_V));
	\end{aligned}\end{equation}
	more precisely, $I_{\mathrm{geom}}(f) = I_{\mathrm{geom}}(f^1)$ holds whenever $f|_{\tilde{G}_V^1} = f^1|_{\tilde{G}_V^1}$.
	\index{Igeom}
	
	This step is essentially done in \cite[Théorème 6.4]{Li14b}, except that only those
	\[ f \in \orbI_{\asp}(\tilde{G}_V^1, \tilde{K_V}) \otimes C^\infty_c(A_{G, \infty}) \otimes \mes(G(F_V)) \]
	were considered in \textit{loc.\ cit.}; that constraint is not necessary.
\end{enumerate}

\begin{remark}\label{rem:cpt-supp-test-fcn-geom}
	In the rest of this work, we will get rid of the $A_{G, \infty}$-invariant test functions $f^1$, and work with test functions $f \in \orbI_{\asp}(\tilde{G}_V, \tilde{K}_V) \otimes \mes(G(F_V))$ instead. In other words, we use the invariant trace formula à la \cite[Théorème 6.4]{Li14b}. We emphasize again that $I_{\mathrm{geom}}(f)$ depends only on $f|_{\tilde{G}_V^1}$.
\end{remark}

We continue to recast \eqref{eqn:Igeom-0} into the form of \cite[VI.2]{MW16-2}, in order to use the method therein. We will also need the supplements from \S\ref{sec:Aell}, giving the necessary reviews and improvements on the construction of these coefficients.

At this stage, let us remark that every $\tilde{\gamma} \in \Gamma(\tilde{M}^1, V)$ is extracted from an element of $M(F)$. That element is not unique, but the $M(F)$-conjugacy class of its semisimple part is: it boils down to the notion of $(M,S)$-equivalence in \cite[Définition 6.5.1]{Li14a}. By collecting terms according to semisimple conjugacy classes $\mathcal{O}_V$ in $M(F_V)$ (Definition \ref{def:geom-O-dist}), we arrive at
\begin{equation}\label{eqn:Igeom-1}
	I_{\mathrm{geom}}(f) = \sum_{M \in \mathcal{L}(M_0)} \frac{|W^M_0|}{|W^G_0|} \sum_{\substack{\mathcal{O}_V \in M(F_V) /\text{conj} \\ \text{semisimple}}} I_{\tilde{M}_V}\left(A^{\tilde{M}}(V, \mathcal{O}_V), f\right),
\end{equation}
where $A^{\tilde{M}}(V, \mathcal{O}_V) \in D_{\mathrm{orb}, -}(\tilde{M}_V^1, \mathcal{O}_V) \otimes \mes(M(F_V))^\vee$ are uniquely determined. As in \eqref{eqn:Igeom-orig}, the sum is finite when $f$ is given.

It will follow from Remark \ref{rem:Acoeff-indep-K} that $A^{\tilde{G}}(V, \mathcal{O}_V)$ is independent of the choice of $K_V$. On the other hand, it does depend on $M_0$ and $K^V$.
\index{AGVO@$A^{\tilde{G}}(V, \mathcal{O}_V)$, $A^{\tilde{G}}(V, \mathcal{O})$}

\begin{remark}\label{rem:Igeom-A}
	More generally, $A^{\tilde{M}}(V, \mathcal{O}_V)$ is defined by linearity when $\mathcal{O}_V$ is a finite union of semisimple conjugacy classes. Hence
	\[ I_{\mathrm{geom}}(f) = \sum_{M \in \mathcal{L}(M_0)} \frac{|W^M_0|}{|W^G_0|} \sum_{\substack{\mathcal{O}_V \in M(F_V)/\mathrm{st. conj} \\ \text{semisimple}}} I_{\tilde{M}_V}\left( A^{\tilde{M}}(V, \mathcal{O}_V), f \right) \]
	for all $f \in \orbI_{\asp}(\tilde{G}_V,\tilde{K}_V) \otimes \mes(G(F_V))$.
\end{remark}

In order to explicate the coefficients $A^{\tilde{G}}(V, \mathcal{O}_V)$, we proceed in steps.

\subsection{Admissibility conditions}\label{sec:S-admissibility}
Let $F$ be a number field and $H$ be a connected reductive $F$-group, endowed with an integral structure outside $V_{\mathrm{ram}} = V_{\mathrm{ram}}(H)$, so that we have hyperspecial subgroups $K_v \subset H(F_v)$ for each $v \notin V_{\mathrm{ram}}$.

Let $S$ be a finite set of places of $F$ such that $S \supset V_{\mathrm{ram}}$. We give a review of the notion of $S$-admissibility following \cite[Définition 5.6.1]{Li14a} or \cite[\S 1]{Ar02}. First, define an invariant morphism
\begin{gather*}
	\mathcal{D} = (D_0, \ldots, D_d): H \to \Ga^{d+1}, \quad d := \dim H, \\
	\det\left( 1 + t - \Ad(x) \middle| \mathfrak{h} \right) = D_0(x) + D_1(x) t + \cdots + D_d(x) t^d \; \in F[t].
\end{gather*}

By invariance, $\mathcal{D}(x)$ depends only on the semisimple part of $x$. Note that $D_d(x)=1$. For all $X \in F^{d+1} \smallsetminus \{0\}$, let $X_{\mathrm{min}}$ denote its first nonzero coordinate. Taking a maximal $F$-torus $T \subset H$ containing $x$, one readily sees that
\begin{equation}\label{eqn:Dmin-root}
	\mathcal{D}(x)_{\mathrm{min}} = \prod_{\substack{\alpha \in \Sigma(H, T) \\ \alpha(x) \neq 1}} \left( 1  - \alpha(x) \right).
\end{equation}
Hence $\mathcal{D}(x)_{\mathrm{min}}$ is the Weyl discriminant of $x$.

\begin{definition}\label{def:admissible-subset}
	\index{admissible subset}
	Let $C_S$ be subset of $F_S^{d+1} \smallsetminus \{0\}$. Denote $\mathfrak{o}^S := \prod_{v \notin S} \mathfrak{o}_v$. We say $C_S$ is $S$-admissible if for all $X \in F^{d+1} \cap \left(C_S \times (\mathfrak{o}^S)^{d+1}\right)$, we have $|X_{\mathrm{min}}|_v = 1$ for all $v \notin S$.
	
	We say a subset $E_S \subset H(F_S)$ is admissible if $\mathcal{D}(E_S)$ is. We say a subset $E \subset H(\A_F)$ is $S$-admissible if $\mathcal{D}(E) \subset C_S \times (\mathfrak{o}^S)^{d+1}$ for some admissible $C_S$.
\end{definition}

It is routine to see that if $E$ is $S$-admissible, then $E$ is $S'$-admissible for all $S' \supset S$.

\begin{definition}\label{def:SO}
	For every semisimple conjugacy class $\mathcal{O} \subset H(F)$, we take the minimal finite set $S(\mathcal{O})$ of places of $F$ satisfying
	\begin{itemize}
		\item $S(\mathcal{O})$ contains all Archimedean places,
		\item for all $v \notin S(\mathcal{O})$, the image in $H(F_v)$ of some (equivalently, every) $\sigma \in \mathcal{O}$ is a compact element,
		\item $\mathcal{O}$ is $S(\mathcal{O})$-admissible.
	\end{itemize}
	For the existence of such a minimal choice $S(\mathcal{O})$, see the Remark below.
\end{definition}

\begin{remark}\label{rem:admissible-root}
	Suppose that a finite set of places $S$ verifies the first two conditions above. The final admissibility condition on $S(\mathcal{O})$ is then equivalent to that for some (equivalently, all) $\sigma \in \mathcal{O}$, a maximal $F$-torus $T \ni \sigma$ and all $v \notin S$, we have
	\[ \forall \alpha \in \Sigma(H, T), \quad \text{either}\; \alpha(\sigma) = 1 \;\text{or}\; 1 - \alpha(\sigma) \in \overline{\mathfrak{o}_v}^\times. \]
	This follows from \eqref{eqn:Dmin-root} since $\alpha(\sigma) \in \overline{\mathfrak{o}_v}^\times$ by the compactness of $\sigma$.
\end{remark}

If $M$ is a Levi subgroup of $H$ and $\mathcal{O}_M \subset \mathcal{O}$, then $S(\mathcal{O}) \supset S(\mathcal{O}_M)$.

\begin{definition}\label{def:SOK}
	\index{SOK@$S(\mathcal{O})$, $S(\mathcal{O}, K)$, $S(\mathcal{O}, K)'$}
	Let $\mathcal{O}$ be a semisimple conjugacy class in $G(F)$. Denote by $\mathcal{O}_v$ (resp.\ $\mathcal{O}_v^{\mathrm{st}}$) the conjugacy class (resp.\ stable conjugacy class) in $H(F_v)$ determined by $\mathcal{O}$, for every place $v$. We define the finite sets
	\begin{align*}
		S(\mathcal{O}, K)' & := S(\mathcal{O}) \cup V_{\mathrm{ram}} \cup \left\{ v \notin V_{\mathrm{ram}} : \mathcal{O}_v \cap K_v = \emptyset \right\}, \\
		S(\mathcal{O}, K) & := S(\mathcal{O}) \cup V_{\mathrm{ram}} \cup \left\{ v \notin V_{\mathrm{ram}} : \mathcal{O}_v^{\mathrm{st}} \cap K_v = \emptyset \right\}.
	\end{align*}
\end{definition}

It follows that
\[ S(\mathcal{O}) \subset S(\mathcal{O}, K) \subset S(\mathcal{O}, K)' . \]

Note that $S(\mathcal{O})$ and $S(\mathcal{O}, K)$ depend only on the stable class of $\mathcal{O}$, thus are more manageable for our purposes.

\begin{lemma}\label{prop:admissible-Endo}
	Let $\rev: \tilde{G} \to G(\A_F)$ be of metaplectic type. Let $\mathbf{G}^! \in \Endo_{\elli}(\tilde{G})$ and let $\mathcal{O}^! \subset G^!(F)$ be a stable semisimple conjugacy class, with image $\mathcal{O} \subset G(F)$. We have $S(\mathcal{O}) \cup \{v: v \mid 2\} \supset S(\mathcal{O}^!)$.
\end{lemma}
\begin{proof}
	Let $\sigma \in \mathcal{O}^!$ and $v \notin S(\mathcal{O})$, assuming that $v \nmid 2$. By looking at eigenvalues, $\sigma$ is seen to be a compact element in $G^!(F_v)$. It remains to show $\sigma$ is $S(\mathcal{O})$-admissible. Fix a maximal torus $T^!$ containing $\sigma$. Let $v \notin S(\mathcal{O})$. In view of Remark \ref{rem:admissible-root}, this is tantamount to $\beta(\sigma) - 1 \in \mathfrak{o}_{\overline{F_v}}^\times$ for all root $\beta \in \Sigma(G^!, T^!)$ with $\beta(\sigma) \neq 1$.
	
	Without loss of generality, let us assume $G = \Sp(W)$. If $\beta$ is a long root, the desideratum about $\beta(\sigma) - 1$ is the same as its counterpart for $\mathcal{O} \subset G(F)$. If $\beta$ is short, the counterpart for $\mathcal{O} \subset G(F)$ says that either $\beta(\sigma) = \pm 1$ or $\beta(\sigma)^2 - 1 \in \mathfrak{o}_{\overline{F_v}}^\times$.
	
	In the former case, $\beta(\sigma) = \pm 1$ meets our desideratum since $v \nmid 2$. In the latter case, since $\beta(\sigma) \pm 1 \in \mathfrak{o}_{\overline{F_v}}$, this still implies our desideratum.
\end{proof}

\subsection{Reduction to elliptic and admissible classes}\label{sec:Aell}
Consider the geometric side $J_{\mathrm{geom}}(\dot{f}) = \sum_{\mathcal{O}} J_{\mathcal{O}}(\dot{f})$ of the coarse trace formula \cite[\S 6.1]{Li14a}, where $\mathcal{O}$ ranges over semisimple conjugacy classes of $G(F)$, and $\dot{f} = \prod_v f_v \in C^\infty_{c, \asp}(\tilde{G}) \otimes \mes(G(\A_F))$. The distributions $J_{\mathcal{O}}(\dot{f})$ depend only on $\dot{f}|_{\tilde{G}^1}$.
\index{Jgeom@$J_{\mathrm{geom}}$}

The \emph{refined geometric expansion} \cite[Théorème 6.5.8]{Li14a} of $J_{\mathrm{geom}}(\dot{f})$ amounts to an expression in terms of semi-local weighted orbital integrals $J_{\tilde{M}_S}(\cdot, f_S)$ with global coefficients $a^{\tilde{M}}(S, \cdot)$, summed over $M \in \mathcal{L}(M_0)$ as usual.

\index{AGSOell@$A^{\tilde{G}}(S, \mathcal{O})_{\elli}$}
In what follows, we will focus on the \emph{elliptic coefficients} in the refined geometric expansion. Given $\mathcal{O}$ and
\[ S \supset V_{\mathrm{ram}} \cup S(\mathcal{O}), \]
there is a coefficient $A^{\tilde{G}}(S, \mathcal{O})_{\elli}$, to be explicated below. Given $f_S \in C^\infty_{c, \asp}(\tilde{G}_S) \otimes \mes(G(F_S))$, set $\dot{f} := f_S \prod_{v \notin S} f_{K_v}$. Following \cite[Lemme 5.4]{Li14b} and the improvements to be given in Remark \ref{rem:condition-S}, such genuine distributions $A^{\tilde{G}}(S, \mathcal{O})_{\elli}$ are characterized by
\begin{itemize}
	\item $A^{\tilde{G}}(S, \mathcal{O})_{\elli} \neq 0$ only when $\mathcal{O}$ is elliptic and $S \supset S(\mathcal{O}, K)'$;
	\item the semisimple part of $\Supp\left( A^{\tilde{G}}(S, \mathcal{O})_{\elli} \right)$ comes from $\mathcal{O}$;
	\item for all $f_S$ such that $\Supp(f_S) \subset G(F_S)$ is admissible, we have
	\begin{equation}\label{eqn:Aell-characterization}
		\sum_{\mathcal{O}} I^{\tilde{G}_S}\left(A(S, \mathcal{O})_{\elli}, f_S \right) = I_{\elli}(\dot{f})
	\end{equation}
	where $I_{\elli}$ denotes the elliptic part of the refined geometric expansion, i.e.\ the terms with $M=G$.
\end{itemize}

Given a compact subset $\Delta \subset G(\A_F)$ of the form $\Delta_S \times K^S$, there are only finitely many $\mathcal{O}$ with $A^{\tilde{G}}(S, \mathcal{O})_{\elli} \neq 0$ and $\mathcal{O} \cap \Delta \neq \emptyset$; see \cite[Théorème 6.5.9]{Li14a}. In \textit{loc.\ cit.}, these coefficients are indexed by $\tilde{\gamma} \in \Gamma(\tilde{G}_S^1)$ and denoted as $a^{\tilde{G}}_{\elli}(\tilde{\gamma})$. Here we collect terms according to $\mathcal{O}$.
\index{aGelli@$a^{\tilde{G}}_{\elli}$}

Hereafter, assume $\mathcal{O}$ to be elliptic and $S \supset S(\mathcal{O}, K)'$. Let $\sigma \in \mathcal{O}$. In this setting, the $G(F_v)$-class of $\sigma$ cuts $K_v$ for all $v \notin S$.

\begin{remark}\label{rem:condition-S}
	The prior conditions follow \cite[\S 3]{Ar88-2} and \cite[\S 2]{Ar02}; they depend only on $\mathcal{O}$ and not on $\sigma$. In \cite{Li14a, Li14b}, one followed Arthur's older work \cite{Ar86} and make stronger conditions on $S$: instead of just assuming $\sigma$ cuts $K_v$ up to $G(F_v)$-conjugacy, one assumes that there exist $\sigma \in \mathcal{O}$ and $M_1 \in \mathcal{L}(M_0)$ such that
	\begin{itemize}
		\item $\sigma^S \in K^S$;
		\item $M_1$ contains $\sigma$, but no proper parabolic subset of $M_1$ contains $\sigma$.
	\end{itemize}
	Putting $M_{1, \sigma} := M_1 \cap G_\sigma$, these requirements also imply that $K_v \cap G_\sigma(F_v)$ and $M_{1, \sigma}$ are in good relative position for each $v \notin S$. The role of $M_1$ will be explained in \S\ref{sec:Aunip}.
	
	When writing \cite{Li14b}, the author was somewhat cavalier in stating these conditions on $\sigma$, $S$ and the space $\mathcal{H}_{\mathrm{adm}, \asp}(\tilde{G}_S, A_{G, \infty})$ of admissible test functions. This is often immaterial since given $\mathcal{O}$, the requirements are met by first taking $\sigma$ and $M_1$, then enlarging $S$. But the study of stable conjugacy seems to require the more flexible version $S \supset S(\mathcal{O}, K)'$.
	
	In the generality of \cite{Li14a} for covering groups, the assumption $\sigma^S \in K^S$ was used in extracting the $S$-component
	\begin{equation}\label{eqn:extraction-S}
		G(F) \ni \sigma = \tilde{\sigma}_S \sigma^S \in \tilde{G}_S \times K^S \subset \tilde{G},
	\end{equation}
	where $K^S \hookrightarrow \tilde{G}^S$ via the chosen sections. The choice of $M_1$ and \eqref{eqn:extraction-S} enter the refined expansion for $J_{\mathcal{O}}$, specifically in the descent to $A^{G_\sigma}_{\mathrm{unip}}(S)$; we will return to the details later on.
	
	If $\sigma^S$ is lifted to $\tilde{G}^S$ via some conjugate of $K^S$, the result might differ by $\Ker(\rev)$, a priori. Therefore, if $\sigma^S$ is only conjugate to $K^S$, one has to fix not only $\sigma \in \mathcal{O}$ and $M_1$, but also $x^S \in G(F^S)$ such that $\sigma^S \in x^S K^S (x^S)^{-1}$, in order to adopt Definition \ref{def:SOK}. The same ambiguity also enters the twisted trace formula with a character \cite[VI.2.2--2.3]{MW16-2}.
	
	Fortunately, the choice of $x^S$ is not a big issue for our $\tilde{G}$: assume $\tilde{G} = \Mp(W, \A_F)$ without loss of generality, then the liftings from $x^S K^S (x^S)^{-1}$ and $K^S$ coincide on their intersection; this is the content of \cite[Lemma A.3.1]{LMS16}.
\end{remark}

Now we revert to the given $V \supset V_{\mathrm{ram}}$. Take $S \supset V \cup S(\mathcal{O})$, then $A^{\tilde{G}}(S, \mathcal{O})_{\elli}$ decomposes into a finite sum
\begin{gather*}
	A^{\tilde{G}}(S, \mathcal{O})_{\elli} = \sum_i k_i \otimes A_i, \\
	k_i \in D_{\mathrm{geom}, -}(\tilde{G}^V_S) \otimes \mes(G(F^V_S))^\vee, \quad A_i \in D_{\mathrm{orb}, -}(\tilde{G}_V) \otimes \mes(G(F_V))^\vee,
\end{gather*}
Of course, this applies to Levi subgroups of $G$ as well. Set
\begin{equation}\label{eqn:Aell-induction}\begin{aligned}
	A^{\tilde{G}}(V, \mathcal{O}) & := \sum_{M \in \mathcal{L}(M_0)} \frac{|W^M_0|}{|W^G_0|} \sum_{\mathcal{O}_M \mapsto \mathcal{O}} \sum_i r^{\tilde{G}}_{\tilde{M}}(k_i, K_S^V) A_i^{\tilde{G}_V}, \\
	A^{\tilde{G}}(V, \mathcal{O}_V) & := \sum_{\mathcal{O} \mapsto \mathcal{O}_V} A^{\tilde{G}}(V, \mathcal{O}),
\end{aligned}\end{equation}
\index{AGVO}
where $\mathcal{O}_V$ is any semisimple conjugacy class in $G(F_V)$, and
\begin{itemize}
	\item $\mathcal{O}_M$ ranges over semisimple conjugacy classes in $M(F)$ contained in $\mathcal{O}$;
	\item $\sum_i k_i \otimes A_i$ is a decomposition of $A^{\tilde{M}}(S, \mathcal{O}_M)_{\elli}$, for each $\mathcal{O}_M$;
	\item $r^{\tilde{G}}_{\tilde{M}}\left(k_i, K_S^V \right)$ are the semi-local unramified weighted orbital integrals along $k_i$, to be reviewed in \S\ref{sec:semi-local-WFL};
	\item $A_i^{\tilde{G}_V}$ is the parabolic induction of $A_i$ from $\tilde{M}_V$ to $\tilde{G}_V$.
\end{itemize}

Indeed, this is simply \cite[Définition 5.6]{Li14b}, rephrased à la \cite[VI.2.3 (8)]{MW16-2}. By \cite[Proposition 5.10]{Li14b}, $A^{\tilde{G}}(V, \mathcal{O}_V)$ depends solely on $K^V$ (which is fixed throughout) and $\mathcal{O}_V$. Again, the coefficients $a^{\tilde{G}}(\tilde{\gamma})$ therein are indexed by $\Gamma(\tilde{G}_V^1)$, and here we collect terms according to $\mathcal{O}$.

\begin{remark}\label{rem:Acoeff-indep-K}
	By \cite[Remarque 6.5.6]{Li14a}, $A^{\tilde{G}}(S, \mathcal{O})_{\elli}$ does not depend on the maximal compact $K_S \subset G(F_S)$. Therefore, $A^{\tilde{G}}(V, \mathcal{O})$ and $A^{\tilde{G}}(V, \mathcal{O}_V)$ do not depend on $K_V$.
\end{remark}

\begin{remark}\label{rem:r-support}
	The set $\mathcal{K}(\tilde{M}^V_S)$ in \cite[Définition 5.6]{Li14b} was erroneously defined. It should be the set of conjugacy classes $k$ in $\tilde{M}^V_S$ whose induction $k^G$ intersects $K^V_S$. By \cite[Corollary 6.2]{Ar88LB}, the classes $k$ with $r^{\tilde{G}}_{\tilde{M}}(k, K^V_S) \neq 0$ must satisfy this property.
	
	Since the induction is induced by $M \hookrightarrow G$ on the semisimple parts\footnote{In contrast, it involves an induction of Lusztig--Spaltenstein type on the unipotent part.}, we infer that $A^{\tilde{G}}(V, \mathcal{O})$ is zero unless the $G(F^S)$-conjugacy class determined by $\mathcal{O}$ intersects $K^V$. This easily implies that the $\mathcal{O}$-sum in \eqref{eqn:Aell-induction} is finite.
\end{remark}

\subsection{Reduction to the unipotent case}\label{sec:Aunip}
Let $\underline{G}$ be a connected reductive $F$-group, equipped with a smooth model over the ring of $S'$-integers, $S' \supset S \supset V_{\mathrm{ram}}(\underline{G})$. Suppose we are given
\begin{itemize}
	\item a hyperspecial subgroup $\underline{K}_v \subset \underline{G}(F_v)$ for each $s \notin S$, equal to $\underline{G}(\mathfrak{o}_v)$ for almost all $v \notin S'$;
	\item a minimal Levi subgroup $\underline{M}_1 \subset \underline{G}$, in good position relative to $\underline{K}_v$ for all $v \notin S$.
\end{itemize}
The unipotent coefficient
\[ A^{\underline{G}}_{\mathrm{unip}}(S) = A^{\underline{G}}(S, \{1\})_{\elli} \; \in D_{\mathrm{unip}}(\underline{G}(F_S)) \otimes \mes(\underline{G}(F_S))^\vee \]
in the trace formula for $\underline{G}$ is studied by Arthur. We also write $A^{\underline{G}}_{\mathrm{unip}}(S, \underline{M}_1, \underline{K}^S)$ to emphasize the auxiliary data.
\index{AGunipS@$A^{\underline{G}}_{\mathrm{unip}}(S)$}

More flexibility can be gained by removing the condition on good relative positions. Choose $\underline{M}_1$; for each $v \notin S$, choose a conjugate $\underline{K}'_v$ of $\underline{K}_v$ that is in good position relative to $\underline{M}_1$, and $\underline{K}'_v = G(\mathfrak{o}_v)$ for almost all $v \notin S'$. Following \cite[VI.2.2 (5)]{MW16-2}, one defines
\[ A^{\underline{G}}_{\mathrm{unip}}\left( S, \underline{K}^S \right) := A^{\underline{G}}_{\mathrm{unip}}\left( S, \underline{M}_1, (\underline{K}')^S \right). \]
By \textit{loc.\ cit.}, this turns out to be independent of $\underline{M}_1$.

Suppose that $\mathcal{O} \subset G(F)$ is elliptic and $S \supset S(\mathcal{O}, K)'$ (Definition \ref{def:SOK}). Take a decomposition $\tilde{\sigma}_S \in \tilde{G}_S$ as in \eqref{eqn:extraction-S}. We shall also choose a standard $M_1 \in \mathcal{L}(M_0)$ intersecting $\mathcal{O}$ minimally, and assume $\sigma \in \mathcal{O} \cap M_1(F)$.

Taking $\underline{G} := G_\sigma$, the unipotent coefficient
\[ A^{G_\sigma}_{\mathrm{unip}}(S) \in D_{\mathrm{unip}}(G_\sigma(F_S)) \otimes \mes(G_\sigma(F_S))^\vee \]
is defined by taking $\underline{K}^S := G_\sigma(F_S) \cap x^{-1} K^S x$ for a suitable $x \in G(F^S)$.

\begin{proposition}[Cf.\ {\cite[VI.2.3 (7)]{MW16-2}}]\label{prop:Aell-descent}
	Assume $\mathcal{O} \subset G(F)$ is elliptic and $S \supset S(\mathcal{O}, K)'$. For all $f_S \in C^\infty_{c, \asp}(G(F_S)) \otimes \mes(G(F_S))$, we have
	\begin{equation}\label{eqn:Aell-descent}
		I^{\tilde{G}_S}\left( A^{\tilde{G}}(S, \mathcal{O})_{\elli}, f_S \right) = \int_{G_\sigma(F_S) \backslash G(F_S)} I^{G_\sigma}\left( A^{G_\sigma}_{\mathrm{unip}}(S), {}^y f|_{\widetilde{G_\sigma}_S}\right) \dd y ,
	\end{equation}
	where ${}^y f_{\widetilde{G_\sigma}_S}$ is the $\mes(G(F_S))$-valued function $\tilde{\gamma} \mapsto f_S \left( y^{-1} \tilde{\sigma}_S \tilde{\gamma} y\right)$ on $\widetilde{G_\sigma}_S := \rev^{-1}(G_\sigma(\A_S))$. The integration makes sense as the covering $\widetilde{G_\sigma}$ splits canonically over $G_{\sigma, \mathrm{unip}}(\A_F)$ by \cite[Proposition 2.2.1, Lemme 3.1.2]{Li14a}.
	
	When $f_S$ is supported on a sufficiently small invariant neighborhood of $\sigma$, we also have
	\[ I^{\tilde{G}_S}\left( A^{\tilde{G}}(S, \mathcal{O})_{\elli}, f_S \right) = I^{G_{\sigma, S}} \left( A^{G_\sigma}_{\mathrm{unip}}(S), \desc^{\tilde{G}}_{\tilde{\sigma}_S}(f_S) \right) \]
	where $\desc^{\tilde{G}}_{\tilde{\sigma}_S}$ is the semi-local version of the map defined \S\ref{sec:descent-HA}; we switch freely between $G_\sigma$ and $\mathfrak{g}_\sigma$ via the exponential map.
\end{proposition}
\begin{proof}
	It suffices to consider $f_S$ with admissible support, since $S \supset S(\mathcal{O})$. Next, recall from \S\ref{sec:Aell} that the left hand side equals the part with $M=G$ in the refined expansion of $J_{\mathcal{O}}$. We go back into the proof of this expansion \cite[Lemme 6.5.4]{Li14a}.
	
	To simplify matters, first make the stronger assumptions as in Remark \ref{rem:condition-S}; in particular, $\sigma_v \in K_v$ for all $v \notin S$. The commutator pairings $[\cdot, \cdot]$ in \textit{loc.\ cit.} are trivial for metaplectic groups. Since $Z_G(\sigma) = G_\sigma$, we have $\iota^G(\sigma) = 1$ and the penultimate displayed formula therein reads
	\begin{multline*}
		J_{\mathcal{O}}(\dot{f}) = \sum_{M \in \mathcal{L}_\sigma^0(M_1)} \sum_{u \in \Gamma_{\mathrm{unip}}(M_\sigma(F), S)} \frac{|W^{M_\sigma}|}{|W^{G_\sigma}_0|} a^{M_\sigma}(S, \dot{u}) \\
		\cdot \int_{G_\sigma(F_S) \backslash G(F_S)} \left( \sum_{R \in \mathcal{F}^{G_\sigma}(M_\sigma)} 	J^{M_R}_{M_\sigma}\left(\dot{u}, \Phi_{R, y, T_1} \right)\right) \dd y.
	\end{multline*}
	Then \cite[Proposition 6.3.4]{Li14a} is applied to transform the $R$-sum into a weighted orbital integral $J_{\tilde{M}_S}(\cdots, f_S)$. Therefore, the required terms correspond exactly to $M = G$ in the sum over $\mathcal{L}_\sigma^0(M_1)$. The resulting $\Phi_{G_\sigma, y, T_1}$ (see also \cite[p.100]{Li14a}) can be replaced by ${}^y f_{\widetilde{G_\sigma}_S}$ since we only care about its orbital integrals. This yields the right hand side of \eqref{eqn:Aell-descent}.
	
	For the second assertion, recall the construction of $\desc^{\tilde{G}}_{\tilde{\sigma}_S}$, eg.\ \cite[Proposition 4.2.1]{Li12b}. This yields the integral over $G_\sigma(F_S) \backslash G(F_S)$ since $Z_G(\sigma) = G_\sigma$.
	
	Consider the general case. The assumption now is that $(g^S)^{-1} \sigma^S g^S \in K^S$ for some $g^S \in G(F^S)$. We need slight modifications on some formulas in Arthur's proof \cite[pp.203--204]{Ar86} which is followed by \cite[Lemme 6.5.4]{Li14a}. Details can be found in the first part proof of \cite[VI.2.3 Proposition]{MW16-2}.
\end{proof}

\subsection{The stable version}
Suppose $G^!$ is a quasisplit $F$-group, endowed with auxiliary data such as a minimal Levi subgroup $M^!_0$, etc., as in the last part of \S\ref{sec:semi-local-geom-distributions}. Take $V \supset V_{\mathrm{ram}}(G^!)$.

Consider a stable semisimple conjugacy class $\mathcal{O}^!$ (resp.\ $\mathcal{O}^!_V$) in $G^!$ (resp.\ in $G^!(F_V)$); see Definition \ref{def:stable-conjugacy}. By \cite[VI.5.2, VII.1.12]{MW16-2}, which is in turn based on Arthur's works on stabilization, we have
\begin{equation}\label{eqn:SAell-sum}\begin{gathered}
	SA^{G^!}(V, \mathcal{O}^!) \in SD_{\text{tr-orb}}(G^!(F_V)) \otimes \mes(G^!(F_V))^\vee , \\
	SA^{G^!}(V, \mathcal{O}^!_V) := \sum_{\mathcal{O}^! \mapsto \mathcal{O}^!_V} SA^{G^!}(V, \mathcal{O}^!).
\end{gathered}\end{equation}
\index{SAVO@$SA^{G^{"!}}(V, \mathcal{O}^{"!})$, $SA^{G^{"!}}(V, \mathcal{O}^{"!}_V)$, $SA^{G^{"!}}(S, \mathcal{O}^{"!})_{\elli}$}

In parallel with \eqref{eqn:Aell-induction}, $SA^{G^!}(V, \mathcal{O}^!)$ can also be reduced to its elliptic version
\[ SA^{G^!}(S, \mathcal{O}^!)_{\elli}, \quad S \supset V_{\mathrm{ram}}(G^!) \cup S(\mathcal{O}^!), \]
and its avatars for Levi subgroups. This is done in \cite[p.780]{MW16-2}. As before, $SA^{G^!}(S, \mathcal{O}^!)_{\elli}$ is nonzero only when $\mathcal{O}^!$ is elliptic and $S \supset S(\mathcal{O}^!, K^!)$ (Definition \ref{def:SOK}), where $K^!_v$ are hyperspecial subgroups. The reduction formula is akin to \eqref{eqn:Aell-induction}. It reads
\begin{equation}\label{eqn:SAell-induction}
	SA^{G^!}(V, \mathcal{O}^!) = \sum_{M^! \in \mathcal{L}(M^!_0)} \frac{|W^{M^!}_0|}{|W^{G^!}_0|} \sum_{\mathcal{O}_M^! \mapsto \mathcal{O}^!} \sum_i s^{G^!}_{M^!}(Sk_i) SA_i^{G^!},
\end{equation}
where $S \supset V \cup S(\mathcal{O}^!)$ and
\begin{itemize}
	\item $\sum_i Sk_i \otimes SA_i$ is a decomposition of $A^{M^!}(S, \mathcal{O}_M^!)_{\elli}$ for each $\mathcal{O}_M^!$,
	\item $SA_i^{G^!}$ is the parabolic induction of $SA_i$ to $G^!(F_V)$,
	\item $s^{G^!}_{M^!}(Sk_i)$ is the stable version of the unramified weighted orbital integral along $Sk_i$, whose definition is deferred to \S\ref{sec:semi-local-WFL}.
\end{itemize}

\begin{remark}\label{rem:s-support}
	In parallel with Remark \ref{rem:r-support}, $s^{G^!}_{M^!}(Sk_i) \neq 0$ only when the induction of $Sk_i$ to $G^!$ intersects $(K^!)^V_S$; see \cite[II.4.2 (2)]{MW16-1}. Again, this implies that $SA^{G^!}(V, \mathcal{O}^!) \neq 0$ only when the stable conjugacy class of $\mathcal{O}^!$ in $G^!(F_v)$ intersects $K^!_v$ for each $v \notin V$. This also guarantees that the $\mathcal{O}^!$-sum in \eqref{eqn:SAell-sum} is finite.
\end{remark}

\index{Sgeom@$S_{\mathrm{geom}}$}
Let $S^{G^!}_{\mathrm{geom}} = S_{\mathrm{geom}}$ be the stable distribution in the geometric side of the stable trace formula \cite{Ar02}. Following \cite[VI.5.9--5.10]{MW16-2}, we have the expansion
\begin{equation}\label{eqn:Sgeom-0}\begin{aligned}
	S_{\mathrm{geom}} & = \sum_{M^! \in \mathcal{L}(M^!_0)} \frac{|W^{M^!}_0|}{|W^{G^!}_0|} \sum_{\mathcal{O}^!} S^{G^!_V}_{M^!_V}\left( SA^{M^!}(V, \mathcal{O}^!), \cdot \right) \\
	& = \sum_{M^! \in \mathcal{L}(M^!_0)} \frac{|W^{M^!}_0|}{|W^{G^!}_0|} \sum_{\mathcal{O}^!_V} S^{G^!_V}_{M^!_V}\left( SA^{M^!}(V, \mathcal{O}^!_V), \cdot \right)
\end{aligned}\end{equation}
where
\begin{itemize}
	\item $\mathcal{O}^!$ (resp.\ $\mathcal{O}^!_V$) ranges over stable semisimple conjugacy classes in $M^!(F)$ (resp.\ in $M^!(F_V)$);
	\item $S^{G^!_V}_{M^!_V}$ stands for the semi-local stable distributions reviewed in \S\ref{sec:semi-local-geom-distributions}.
\end{itemize}

As in $I_{\mathrm{geom}}$, all the semi-local stable distributions in \eqref{eqn:Sgeom-0} are concentrated at $0 \in \mathfrak{a}_{G^!}$, and the test functions will be drawn from $S\orbI(G^!(F_V)) \otimes \mes(G^!(F_V))$. Furthermore, \eqref{eqn:Sgeom-0} becomes a finite sum when a test function $f^!$ is plugged in.

\begin{lemma}\label{prop:SA-central-twist}
	Suppose that $G^!$ is a direct product of groups of the form $\GL(n)$ and $\SO(2n+1)$. Let $z \in A_{G^!}(F)$ be such that $z^2 = 1$. Then
	\begin{gather*}
		S(z\mathcal{O}^!) = S(\mathcal{O}^!), \quad S(z\mathcal{O}^!, K^!) = S(\mathcal{O}^!, K^!), \\
		SA^{G^!}(V, \mathcal{O}^!)[z] = SA^{G^!}(V, z\mathcal{O}^!)
	\end{gather*}
	for all $\mathcal{O}^!$, where $(\cdots)[z]$ means the distribution is twisted by multiplication by $z$. Similarly, $SA^{G^!}(V, \mathcal{O}^!_V)[z] = SA^{G^!}(V, z\mathcal{O}^!_V)$ for all $\mathcal{O}_V^!$.
\end{lemma}
\begin{proof}
	We are immediately reduced to the case $G^! = \GL(n)$, so that $SA^{G^!}(V, \cdot) = A^{G^!}(V, \cdot)$. It is clear that $z$ belongs to all hyperspecial subgroups of $G^!(F_v)$. One can now trace through the definitions of $S(\cdot)$ and $A^{G^!}(V, \cdot)$ to obtain the required identities.
\end{proof}

Finally, the elliptic coefficients $SA^{G^!}(S, \mathcal{O}^!)_{\elli}$ descends to the unipotent case in a way similar to Proposition \ref{prop:Aell-descent}. Let $SA^{G^!}_{\mathrm{unip}}(S) := SA^{G^!}(S, \{1\})_{\elli}$. By \cite[VI.2.2, 5.2]{MW16-2}, $SA^{G^!}_{\mathrm{unip}}(S)$ is independent of all choices of minimal Levi or hyperspecial subgroups.

\begin{theorem}\label{prop:SAell-descent}
	Suppose that $\mathcal{O}^!$ is elliptic and $\epsilon \in \mathcal{O}^!$, such that $G^!_\epsilon$ is quasisplit. Denote by $\Xi_\epsilon$ the $F$-group scheme $Z_{G^!}(\epsilon)/G^!_\epsilon$. Suppose $S \supset S(\mathcal{O^!}, K^!)$, then
	\begin{multline}\label{eqn:SAell-descent}
		S^{G^!_S}\left( SA^{G^!}(S, \mathcal{O}^!)_{\elli}, f_S \right) \\
		= \left|\Xi_\epsilon(F)\right|^{-1} \tau(G^!) \tau(G^!_\epsilon)^{-1} S^{G^!_\epsilon}\left( SA^{G^!_\epsilon}_{\mathrm{unip}}(S), S\desc^{G^!}_{\epsilon_S}(f_S) \right).
	\end{multline}
	for all $f_S \in S\orbI(G^!(F_S)) \otimes \mes(G^!(F_S))$ supported in a sufficiently small stably invariant neighborhood of $\epsilon$. Here
	\begin{compactitem}
		\item $\tau(G^!)$ and $\tau(G^!_\epsilon)$ denotes the normalized Tamagawa numbers, defined using compatible Haar measures on $A_{G^!, \infty} \simeq A_{G^!_\epsilon, \infty}$.
		\item $S\desc^{G^!}_{\epsilon_S}$ is the semi-local version of the map defined in \S\ref{sec:descent-HA}; we switch freely between $G^!_\epsilon$ and $\mathfrak{g}^!_\epsilon$ via the exponential map.
	\end{compactitem}
\end{theorem}
\begin{proof}
	It boils down to \cite[VII.3.2 Théorème]{MW16-2}. The $|\Xi_\epsilon(F)|$ there should be corrected to $|\Xi_\epsilon(F)|^{-1}$; cf.\ \cite[p.626 (7)]{MW16-2}.
\end{proof}

\begin{remark}
	As in Remark \ref{rem:Acoeff-indep-K}, $SA^{G^!}(S, \mathcal{O}^!)_{\elli}$ (resp.\ $SA^{G^!}(V, \mathcal{O}^!)$) does not depend on $K_S$ (resp.\ $K_V$).
\end{remark}

\section{Statement of the global matching}\label{sec:Igeom-global-matching}
Consider a covering of metaplectic type $\rev: \tilde{G} \to G(\A_F)$, together with all the auxiliary data involved in the semi-local distributions in \S\ref{sec:semi-local-geom-distributions}.

For any connected reductive group $H$ over $F$, let $\tau(H)$ denote its normalized Tamagawa number. Recall that $\tau(\GL(m)) = 1$ for all $m$. The following repeats Definition \ref{def:iotaGG-general}.

\begin{definition}[Cf.\ {\cite[p.747]{Li15}}]\label{def:iota-const}
	\index{iotaGG}
	Let $\mathbf{G}^! \in \Endo_{\elli}(\tilde{G})$. Set $\iota(\tilde{G}, G^!) := \tau(G) \tau(G^!)^{-1}$. More concretely, if $G \simeq \prod_{i \in J} \GL(n_j) \times \Sp(2n^\flat)$ and $\mathbf{G}^!$ is given by the pair $(n', n'') \in \Z_{\geq 0}^2$ with $n' + n'' = n^\flat$, then
	\[ \iota(\tilde{G}, G^!) = \begin{cases}
		\frac{1}{4}, & n', n'' \geq 1 \\
		\frac{1}{2}, & n' n'' = 0, \; n^{\flat} \geq 1 \\
		1, & n^{\flat} = 0.
	\end{cases}\]
	If $\mathbf{G}^! \in \Endo(\tilde{G}) \smallsetminus \Endo_{\elli}(\tilde{G})$, we set $\iota(\tilde{G}, G^!) = 0$.
\end{definition}

For all $\mathbf{G}^! \in \Endo(\tilde{G})$, note that $V_{\mathrm{ram}} \supset V_{\mathrm{ram}}(G^!) \cup  \{v: v \mid 2\}$ for all $\mathbf{G}^!$.
\begin{itemize}
	\item The transfer $\Trans_{\mathbf{G}^!, \tilde{G}}$ and its dual $\trans_{\mathbf{G}^!, \tilde{G}}$ extend to the semi-local setting, by taking the $\otimes$-product over $v \in V$.
	\item Furthermore, by Theorem \ref{prop:image-transfer}, $\Trans_{\mathbf{G}^!, \tilde{G}}$ restricts to
	\[ \Trans_{\mathbf{G}^!, \tilde{G}}: \orbI_{\asp}\left(\tilde{G}_V, \tilde{K}_V\right) \otimes \mes(G(F_V)) \to S\orbI\left(G^!(F_V), K^!_V\right) \otimes \mes(G^!(F_V)). \]
\end{itemize}

Consider $I^{\tilde{G}}_{\mathrm{geom}}$ (resp.\ $S^{G^!}_{\mathrm{geom}}$), the distribution of the geometric side of the trace formula for $\tilde{G}$ (resp.\ the stable trace formula for $G^!$), reviewed in \S\ref{sec:geom-side}. Let $V \supset V_{\mathrm{ram}}$ be a finite set of places of $F$, and $f \in \orbI_{\asp}\left(\tilde{G}_V, \tilde{K}_V\right) \otimes \mes(G(F_V))$. Define
\[ I^{\Endo}_{\mathrm{geom}}(f) = I^{\tilde{G}, \Endo}_{\mathrm{geom}}(f) := \sum_{\mathbf{G}^! \in \Endo_{\elli}(\tilde{G})} \iota(\tilde{G}, G^!) S^{G^!}_{\mathrm{geom}}\left( \Trans_{\mathbf{G}^!, \tilde{G}}(f)\right). \]
\index{IEndogeom@$I^{\Endo}_{\mathrm{geom}}$}

Our long-term goal is the
\begin{theorem}[Cf.\ {\cite[VI.5.10]{MW16-2}}]\label{prop:geom-stabilization}
	For $V \supset V_{\mathrm{ram}}$ as above, we have
	\[ I^{\tilde{G}}_{\mathrm{geom}} = I^{\tilde{G}, \Endo}_{\mathrm{geom}}. \]
\end{theorem}

Its proof will be completed at the end of this work; see \S\ref{sec:end-stabilization}.

The global information in $I_{\mathrm{geom}}$ is encoded in the coefficients reviewed in \S\ref{sec:geom-side}. They have the endoscopic counterparts below.

\begin{definition}\label{def:AEndo}
	\index{AEndoVO@$A^{\tilde{G}, \Endo}(V, \mathcal{O})$, $A^{\tilde{G}, \Endo}(V, \mathcal{O}_V)$}
	Suppose that $V \supset V_{\mathrm{ram}}$ and let $\mathcal{O}$ (resp.\ $\mathcal{O}_V$) be a stable semisimple conjugacy class in $G(F)$ (resp.\ $G(F_V)$). Set
	\begin{align*}
		A^{\tilde{G}, \Endo}(V, \mathcal{O}) & := \sum_{\mathbf{G}^! \in \Endo_{\elli}(\tilde{G})} \iota(\tilde{G}, G^!) \sum_{\mathcal{O}^! \mapsto \mathcal{O}} \trans_{\mathbf{G}^!, \tilde{G}} \left( SA^{G^!}(V, \mathcal{O}^!) \right), \\
		A^{\tilde{G}, \Endo}(V, \mathcal{O}_V) & := \sum_{\mathbf{G}^! \in \Endo_{\elli}(\tilde{G})} \iota(\tilde{G}, G^!) \sum_{\mathcal{O}^!_V \mapsto \mathcal{O}_V} \trans_{\mathbf{G}^!, \tilde{G}} \left( SA^{G^!}(V, \mathcal{O}^!_V) \right) \\
		& \in D_{\text{tr-orb}, -}(\tilde{G}_V) \otimes \mes(G(F_V))^\vee .
	\end{align*}
	Stable conjugacy is needed for comparing classes in $G$ and in its endoscopic groups. For the finiteness of the $\mathcal{O}^!$-sum, see Remark \ref{rem:s-support}.
	
	In order to make comparison with the non-endoscopic avatars, we also define $A^{\tilde{G}}(V, \mathcal{O})$ (resp.\ $A^{\tilde{G}}(V, \mathcal{O}_V)$) for stable classes $\mathcal{O}$ in $G(F)$ (resp.\ $\mathcal{O}_V$ in $G(F_V)$) by collecting terms.
\end{definition}

By recalling the effect of $\trans_{\mathbf{G}^!, \tilde{G}}$ on supports (Theorem \ref{prop:Dgeom-preservation}), $A^{\tilde{G}, \Endo}(V, \mathcal{O}_V)$ is seen to equal the sum of the $D_{\mathrm{orb}, -}(\tilde{G}_V, \mathcal{O}_V)$-components of $A^{\tilde{G}, \Endo}(V, \mathcal{O})$, for various $\mathcal{O}$; cf.\ the non-endoscopic version \eqref{eqn:Aell-induction}.

The distribution $A^{\tilde{G}, \Endo}(V, \mathcal{O}_V)$ also depends on $K^V$ in the sense that $\Trans_{\mathbf{G}^!, \tilde{G}}$ (or $\Delta$) is normalized with respect to $K^V$; see \cite[Propositions 5.15, 5.17]{Li11}.

\begin{theorem}[Cf.\ {\cite[VI.5.4 Théorème]{MW16-2}}]\label{prop:matching-coeff-A}
	For any stable semisimple conjugacy class $\mathcal{O}$ in $G(F)$, we have
	\[ A^{\tilde{G}, \Endo}(V, \mathcal{O}) = A^{\tilde{G}}(V, \mathcal{O}) . \]
	Consequently, for any stable semisimple conjugacy class $\mathcal{O}_V$ in $G(F_V)$, we have
	\[ A^{\tilde{G}, \Endo}(V, \mathcal{O}_V) = A^{\tilde{G}}(V, \mathcal{O}_V). \]
\end{theorem}

The proof of Theorem \ref{prop:matching-coeff-A} is deferred to \S\ref{sec:gdesc-ascent}.

\section{Semi-local endoscopic constructions}\label{sec:semi-local-geom-Endo}
Fix $V \supset V_{\mathrm{ram}}$. Let $M \in \mathcal{L}(M_0)$ and $\mathbf{M}^! \in \Endo_{\elli}(\tilde{M})$. From $s \in \Endo_{\mathbf{M}^!}(\tilde{G})$ we construct
\begin{itemize}
	\item the constant $i_{M^!}(\tilde{G}, G^![s])$ from Definition \ref{def:i-const},
	\item the system of $B$-functions $B^{\tilde{G}}$ on $G^![s]$ from Definition \ref{def:B-metaplectic},
\end{itemize}
both performed in the global setting. For $\delta \in SD_{\diamondsuit}(M^!(F_V)) \otimes \mes(M(F_V))^\vee$, we can also define $\delta[s] := \delta \cdot z[s]$ by performing the twist by $z[s]$ at each place in $V$ simultaneously.

\begin{definition}\label{def:IEndoMV}
	\index{IMVEndo-delta-weighted@$I^{\Endo}_{\tilde{M}_V}(\mathbf{M}^{"!}, \delta, f)$}
	For all $\delta \in SD_{\diamondsuit}(M^!(F_V)) \otimes \mes(M(F_V))^\vee$, set
	\begin{align*}
		I^{\Endo}_{\tilde{M}_V}(\mathbf{M}^!, \delta, f) & = I^{\tilde{G}_V, \Endo}_{\tilde{M}_V}(\mathbf{M}^!, \delta, f) \\
		& := \sum_{s \in \Endo_{\elli}(\tilde{G})} i_{M^!}(\tilde{G}, G^![s]) S^{G^![s]_V}_{M^!_V}\left( \delta[s], B^{\tilde{G}}, f^{G^![s]} \right),
	\end{align*}
	where $f^{G^![s]}$ is the transfer of $f \in \orbI_{\asp}(\tilde{G}_V, \tilde{K}_V) \otimes \mes(G(F_V))$ to $G^![s]$.
\end{definition}

\begin{remark}\label{rem:IEndo-local-nonell}
	For all $v \in V$, we deduce from $\mathbf{M}^!$ the $\mathbf{M}^!_v \in \Endo_{\elli}(\tilde{M}_v)$. In contrast with the endoscopy for reductive groups, $\mathbf{M}^!_v$ is still elliptic. This is basically because $\tilde{G}^\vee$ and $\tilde{M}^\vee$ have trivial Galois action. Consequently, $I^{\Endo}_{\tilde{M}_v}\left( \mathbf{M}^!_v, \delta_v, f_v \right)$ is well-defined. We will often abbreviate $\mathbf{M}^!_v$ as $\mathbf{M}^!$.
\end{remark}

\begin{definition}\label{def:IMVEndo-gamma-weighted}
	\index{IMVEndo-gamma-weighted@$I^{\Endo}_{\tilde{M}_V}(\tilde{\gamma}, f)$}
	In view of Proposition \ref{prop:semi-local-orbint-splitting}, for all
	\begin{align*}
		\tilde{\gamma} & = \prod_v \tilde{\gamma}_v \in D_{\mathrm{geom}, G_\infty\text{-equi}, -}(\tilde{M}_V) \otimes \mes(M(F_V))^\vee, \\
		f & = \prod_v f_v \in \orbI_{\asp}(\tilde{G}_V, \tilde{K}_V) \otimes \mes(G(F_V))
	\end{align*}
	we define
	\begin{align*}
		I^{\Endo}_{\tilde{M}_V}(\tilde{\gamma}, f) & = I^{\tilde{G}_V, \Endo}_{\tilde{M}_V}(\tilde{\gamma}, f) \\
		& := \sum_{\substack{L^V \in \mathcal{L}(M^V) \\ \forall v, \; M^v := M}} d^G_{M^V}\left( M, L^V \right) \prod_{v \in V} 	I^{\tilde{L}^v, \Endo}_{\tilde{M}^v}\left( \tilde{\gamma}_v, f_{v, \tilde{L}^v}\right).
	\end{align*}
\end{definition}

Two quick remarks are in order.
\begin{itemize}
	\item One cannot define $I^{\Endo}_{\tilde{M}_V}\left(\tilde{\gamma}, \cdot\right)$ by expressing $\tilde{\gamma}$ as transfer, since $\tilde{\gamma}$ is of a local nature, and the global endoscopic data do not suffice for this purpose.
	\item If Hypothesis \ref{hyp:ext-Arch} holds at all Archimedean places, the same definition can be made for $\tilde{\gamma} \in D_{\diamondsuit, -}(\tilde{M}_V) \otimes \mes(M(F_V))^\vee$.
\end{itemize}

The following results guarantee that the definitions are coherent. We begin with some preparations. Choose $s^\flat \in \tilde{M}^\vee$ with $(s^\flat)^2 = 1$ giving rise to $\mathbf{M}^! \in \Endo_{\elli}(\tilde{M})$. We express every $s \in \Endo_{\mathbf{M}^!}(\tilde{G})$ as $s = ts^\flat$ following the identification \eqref{eqn:Endo-s-dual}, where $t \in Z_{\tilde{M}^\vee}^\circ / Z_{\tilde{G}^\vee}^\circ$.

For every $L_1, L_2 \in \mathcal{L}(M)$, set
\begin{gather*}
	S(\tilde{L}_1, \tilde{L}_2) := \left\{\begin{array}{r|l}
		s \in \Endo_{\mathbf{M}^!}(\tilde{G}) & \text{decomposable into}\; (s^{L_1}, s_{L_1}), (s^{L_2}, s_{L_2}) \\
		& \text{with}\; \mathbf{L}_j^![s^{L_j}] \in \Endo_{\elli}(\tilde{L}_j), \; j=1,2
	\end{array}\right\},
\end{gather*}
the decomposition of $s$ being understood as in \S\ref{sec:descent-orbint-Endo}, with $s^{L_j} \in \Endo_{\mathbf{M}^!}(\tilde{L}_j)$ and $s_{L_j} \in \Endo_{\mathbf{L}_j^![s^{L_j}]}(\tilde{G})$. 

When $d^G_M(L_1, L_2) \neq 0$, the natural homomorphism
\begin{equation}\label{eqn:q-L1L2}
	q: Z_{\tilde{M}^{\vee}}^\circ / Z_{\tilde{G}^{\vee}}^\circ \to \left( Z_{\tilde{M}^{\vee}}^\circ / Z_{\tilde{L}_1^{\vee}}^\circ \right) \oplus \left(Z_{\tilde{M}^{\vee}}^\circ / Z_{\tilde{L}_2^{\vee}}^\circ\right)
\end{equation}
is surjective with finite kernel.

\begin{lemma}\label{prop:semilocal-matching-L1L2}
	Suppose that $d^G_M(L_1, L_2) \neq 0$ and $s = s^\flat t \in S(\tilde{L}_1, \tilde{L}_2)$. Set
	\[ (t_1, t_2) := q(t) , \quad s_j := t_j s^\flat, \quad j = 1,2 \]
	Then $s_j \in \Endo_{\mathbf{M}^!}(\tilde{L}_j)$ and $\mathbf{L}_j^![s^{L_j}] = \mathbf{L}_j^![s_j]$. Moreover, $z\left[ s_1 \right] z\left[ s_2 \right] = z[s]$ in this case.
\end{lemma}
\begin{proof}
	Without loss of generality, we may assume that
	\[ \tilde{G} = \Mp(W), \quad M = \prod_{i \in I} \GL(n_i) \times \Mp(W^\flat). \]
	The datum $t$ (or $s$) amounts to prescribe embeddings of $\GL(n_i)$ into the two factors of $G^![s] = \SO(2n'+1) \times \SO(2n''+1)$. For $j \in \{1, 2\}$, taking the projection $t_j$ means that we regard only those $\GL(n_i)$ embedding into the $\Sp$-factor of $L_j$; from this we build $\mathbf{L}_j^![s_j] \in \Endo_{\elli}(\tilde{L}_j)$.
	
	This is also the recipe in \S\ref{sec:descent-orbint-Endo} for constructing $s^{L_j}$ from $s$. All in all, $\mathbf{L}_j^![s^{L_j}] = \mathbf{L}_j^![s_j]$. The first assertion follows.
	
	The condition that $\mathfrak{a}^{L_1}_M \oplus \mathfrak{a}^{L_2}_M = \mathfrak{a}^G_M$ implies that via $s \in S(\tilde{L}_1, \tilde{L}_2)$, every $\GL(n_i)$ embeds into the $\Sp$-factor of $L_j$ for a unique $j \in \{1,2\}$. This implies the second assertion.
\end{proof}

\begin{lemma}\label{prop:semilocal-matching-q}
	Suppose that $d^G_M(L_1, L_2) \neq 0$ and $t \in Z_{\tilde{M}^\vee}^\circ / Z_{\tilde{G}^\vee}^\circ$. Let $(t_1, t_2) := q(t)$ and suppose that $s_j := s^\flat t_j \in \Endo_{\mathbf{M}^!}(\tilde{L}_j)$. Then $s := s^\flat t$ belongs to $S(\tilde{L}_1, \tilde{L}_2)$.
\end{lemma}
\begin{proof}
	According to Remark \ref{rem:Gs-nonelliptic}, we may always decompose $s$ into $(s^{L_j}, s_{L_j})$ for $j=1,2$; the point to verify is ellipticity. Since $s_j = s^{L_j}$ by construction, $\mathbf{L}^!_j[s^{L_j}] \in \Endo_{\elli}(\tilde{L}_j)$ follows from assumption. It remains to show $s \in \Endo_{\mathbf{M}^!}(\tilde{G})$. Since $G^![s] \supset L_j^![s_j]$ for $j=1,2$, we have $\mathfrak{a}_{G^![s]} \subset \mathfrak{a}_{L_1^![s_1]} \cap \mathfrak{a}_{L_2^![s_2]}$, and the latter space corresponds to $\mathfrak{a}_{L_1} \cap \mathfrak{a}_{L_2} = \mathfrak{a}_G$ under $\mathfrak{a}_{M^!} \simeq \mathfrak{a}_M$.
\end{proof}

\begin{proposition}[Cf.\ {\cite[VI.4.5 Proposition]{MW16-2}}]\label{prop:semilocal-matching-coherence}
	Let $\mathbf{M}^! \in \Endo_{\elli}(\tilde{M})$ and $f = \prod_{v \in V} f_v$.
	\begin{enumerate}[(i)]
		\item For all $\delta = \prod_{v \in V} \delta_v \in SD_{\diamondsuit}(M^!(F_V)) \otimes \mes(M(F_V))^\vee$, we have
		\[ I^{\Endo}_{\tilde{M}_V}(\mathbf{M}^!, \delta, f) = \sum_{\substack{L^V \in \mathcal{L}(M^V) \\ \forall v, \; M^v := M}} d^G_{M^V}\left( M, L^V \right) \prod_{v \in V} I^{\tilde{L}^v, \Endo}_{\tilde{M}^v}\left( \mathbf{M}^!, \delta_v, f_{v, \tilde{L}^v} \right). \]
		\item For all $\delta \in SD_{\mathrm{geom}, \tilde{G}_\infty\text{-equi}}(M^!_V) \otimes \mes(M(F_V))^\vee$, we have
		\[ I^{\Endo}_{\tilde{M}_V}\left( \trans_{\mathbf{M}^!, \tilde{M}}(\delta), f \right) = I^{\Endo}_{\tilde{M}_V}(\mathbf{M}^!, \delta, f). \]
		If Hypothesis \ref{hyp:ext-Arch} holds at all $v \mid \infty$, then $\delta \in SD_{\diamondsuit}(M^!_V) \otimes \mes(M^!(F_V))^\vee$ can be allowed.
	\end{enumerate}
\end{proposition}
\begin{proof}
	Let us begin with (ii). We may assume $\delta = \prod_v \delta_v$. By definition and the local results,
	\begin{multline*}
		I^{\Endo}_{\tilde{M}_V}\left( \trans_{\mathbf{M}^!, \tilde{M}}(\delta), f \right) = \sum_{L^V} d^G_{M^V}(M, L^V) \prod_{v \in V} I^{\tilde{L}_v, \Endo}_{\tilde{M}_v}\left( \trans_{\mathbf{M}^!, \tilde{M}}(\delta_v), f_{v, \tilde{L}^v} \right) \\
		= \sum_{L^V} d^G_{M^V}(M, L^V) \prod_{v \in V} I^{\tilde{L}^v, \Endo}_{\tilde{M}}\left( \mathbf{M}^!, \delta_v, f_{v, \tilde{L}^v} \right).
	\end{multline*}
	Assuming (i), the last expression equals $I^{\Endo}_{\tilde{M}_V}(\mathbf{M}^!, \delta, f)$, as desired.
	
	We break (i) into two assertions.
	\begin{enumerate}
		\item Given a decomposition $V = V_1 \sqcup V_2$ and $f = f_1 f_2$ accordingly, we have
		\[ I^{\Endo}_{\tilde{M}_V}\left( \mathbf{M}^!, \delta, f \right) = \sum_{L_1, L_2 \in \mathcal{L}(M)} d^G_M(L_1, L_2) I^{\tilde{L}_{1, V_1}, \Endo}_{\tilde{M}_{V_1}}\left( \mathbf{M}^!, \delta_1, f_{1, \tilde{L}_1} \right) I^{\tilde{L}_{2, V_2}, \Endo}_{\tilde{M}_{V_2}}\left( \mathbf{M}^!, \delta_2, f_{2, \tilde{L}_2} \right). \]
		\item When $V = \{v\}$, we have
		\[ I^{\Endo}_{\tilde{M}_V}\left(\mathbf{M}^!, \delta, f \right) =  \sum_{L^v \in \mathcal{L}(M_v)} d^G_M\left( M, L^v \right) I^{\tilde{L}^v, \Endo}_{\tilde{M}_v}\left( \mathbf{M}^!_v, \delta_v, f_{v, \tilde{L}_v}\right). \]
	\end{enumerate}
	To see that their conjunction implies (i), it suffices to repeat the first step of the proof of \cite[VI.4.5 Proposition]{MW16-2}, which is of a combinatorial nature.
	
	For the second assertion, recall that $\mathbf{M}^!_v$ is elliptic in our setting. The remaining arguments are thus identical to Proposition \ref{prop:descent-orbint-Endo} with $M=R$. We focus on the first assertion in what follows.

	The systems of $B$-functions in question will all be labeled as $B$, as justifiable by Lemma \ref{prop:B-Levi}. Using the factors from \eqref{eqn:e} and \cite[VI.4.2 Proposition (ii)]{MW16-2}, we obtain
	\begin{multline*}
		I^{\Endo}_{\tilde{M}}\left(\mathbf{M}^!, \delta, f\right) = \sum_{s \in \Endo_{\mathbf{M}}^!(\tilde{G})} i_{M^!}(\tilde{G}, G^![s]) \sum_{L_1^!, L_2^! \in \mathcal{L}^{G^![s]}(M)} e^{G^![s]}_{M^!}(L_1^!, L_2^!) \\
		S^{L^!_{1, V_1}}_{M^!_{V_1}}\left( \delta[s], B, (f_1^{G^![s]})_{L^!_1} \right) S^{L^!_{2, V_2}}_{M^!_{V_2}}\left( \delta[s], B, (f_2^{G^![s]})_{L^!_2} \right).
	\end{multline*}
	By the paradigm of \S\ref{sec:descent-orbint-Endo}, each $(s, L_i^!)$ determines $L_i \in \mathcal{L}(M)$ and $s$ decomposes into $s^{L_i} \in \Endo_{\mathbf{M}^!}(\tilde{L}_i)$ and $s_{L_i} \in \Endo_{\mathbf{L}_i^![s^{L_i}]}(\tilde{G})$, such that
	\begin{gather*}
		\mathfrak{a}_{L_i^!} \simeq \mathfrak{a}_{L_i} \;\text{under}\; \mathfrak{a}_{M^!} \simeq \mathfrak{a}_M, \\
		L_i^! = L_i^![s^{L_i}], \quad (f^{G^![s]})_{L_i^!} [s_{L_i}] = (f_{\tilde{L}_i})^{L_i^!}, \quad i=1,2.
	\end{gather*}
	Moreover, one can sum over $L_i$ and $(s^{L_i}, s_{L_i})$. For every $L_1, L_2 \in \mathcal{L}(M)$, set
	\begin{align*}
		X(\tilde{L}_1, \tilde{L}_2) & := \sum_{s \in S(\tilde{L}_1, \tilde{L}_2)} i_{M^!}(\tilde{G}, G^![s]) e^{G^![s]}_{M^!}\left(L_1^![s^{L_1}], L_2^![s^{L_2}] \right) \\
		& \cdot S^{L^!_{1, V_1}}_{M^!_{V_1}}\left( \delta[s], B, (f_1^{G^![s]})_{L^!_1} \right) S^{L^!_{2, V_2}}_{M^!_{V_2}}\left( \delta[s], B, (f_2^{G^![s]})_{L^!_2} \right).
	\end{align*}
	
	For all $s \in S(\tilde{L}_1, \tilde{L}_2)$, the endoscopic data $\mathbf{L}_i^![s^{L_i}]$ and $\mathbf{G}^![s]$ are all elliptic. In this case, \eqref{eqn:e} says that $d^G_M(L_1, L_2) = 0 \implies X(\tilde{L}_1, \tilde{L}_2) = 0$.
		
	Henceforth, assume $d^G_M(L_1, L_2) \neq 0$. Then
	\begin{multline*}
		X(\tilde{L}_1, \tilde{L}_2) = d^G_M(L_1, L_2) \sum_{s \in S(\tilde{L}_1, \tilde{L}_2)} i_{M^!}(\tilde{G}, G^![s])
		\underbracket{k^{G^![s]}_{M^!}\left( L_1^![s^{L_1}], L_2^![s^{L_2}] \right)^{-1}}_{\text{see \eqref{eqn:e}}} \\
		\cdot S^{L^!_{1, V_1}}_{M^!_{V_1}}\left( \delta[s], B, (f_1^{G^![s]})_{L^!_1} \right) S^{L^!_{2, V_2}}_{M^!_{V_2}}\left( \delta[s], B, (f_2^{G^![s]})_{L^!_2} \right).
	\end{multline*}

	Fix $s$. Let $q$ be as in \eqref{eqn:q-L1L2} and $(s_1, s_2)$ be as in Lemma \ref{prop:semilocal-matching-L1L2}. We claim that
	\begin{multline}\label{eqn:semilocal-matching-coherence}
		i_{M^!}(\tilde{G}, G^![s]) k^{G^![s]}_{M^!}\left( L_1^![s^{L_1}], L_2^![s^{L_2}] \right)^{-1} \\
		= \left|\Ker(q)\right|^{-1} i_{M^!}(\tilde{L}_1, L_1^![s_1]) i_{M^!}(\tilde{L}_2, L_2^![s_2]).
	\end{multline}

	Assume the validity of \eqref{eqn:semilocal-matching-coherence}. By Lemmas \ref{prop:semilocal-matching-L1L2}, \ref{prop:semilocal-matching-q} and $|\Ker(q)| = k^G_M(L_1, L_2)$, the sum in $X(\tilde{L}_1, \tilde{L}_2)$ can be transformed into $d^G_M(L_1, L_2)$ times the product over $i=1,2$ of
	\[ \sum_{s_i \in \Endo_{\mathbf{M}^!}(\tilde{L}_i)} i_{M^!}(\tilde{L}_i, L_i^![s_i]) S^{L_i^![s_i]_{V_i}}_{M^!_{V_i}}\left( \delta_i[s_i], B, (f_{\tilde{L}_i})^{L_i^![s_i]} \right) = I^{\tilde{L}_{i, V_i}, \Endo}_{\tilde{M}_{V_i}}\left( \mathbf{M}^!, \delta_i, f_{i, \tilde{L}_i} \right). \]
	Summing over $(L_1, L_2)$ entails the required assertion.

	We turn to the proof of \eqref{eqn:semilocal-matching-coherence}. Consider the diagram
	\[\begin{tikzcd}
		\dfrac{Z_{\tilde{M}^\vee}^\circ}{Z_{\tilde{G}^\vee}^\circ} \arrow[r, "q"] \arrow[d, "r"'] & \dfrac{Z_{\tilde{M}^\vee}^\circ}{Z_{\tilde{L}_1^\vee}^\circ} \oplus \dfrac{Z_{\tilde{M}^\vee}^\circ}{Z_{\tilde{L}_2^\vee}^\circ} \arrow[d, "r_{12}"] \\
		\dfrac{Z_{\check{M}^!}}{Z_{G^![s]^\vee}} \arrow[r, "{q_{12}}"'] & \dfrac{Z_{\check{M}^!}}{Z_{L_1^![s_1]^\vee}} \oplus \dfrac{Z_{\check{M}^!}}{Z_{L_2^![s_2]^\vee}}
	\end{tikzcd}\]
	the $q$ and $q_{12}$ being quotient homomorphisms, and $r$, $r_{12}$ from endoscopy. It commutes, and all arrows are surjective with finite kernels. Specifically,
	\begin{align*}
		\# \Ker(q_{12}) & = k^{G^![s]}_{M^!}\left( L_1^![s_1], L_2^![s_2] \right), \\
		\# \Ker(r) & = i_{M^!}\left(\tilde{G}, G^![s]\right)^{-1}, \\
		\# \Ker(r_{12}) &= i_{M^!}\left(\tilde{L}_1, L_1^![s_1]\right)^{-1} i_{M^!}\left(\tilde{L}_2, L_2^![s_2] \right)^{-1}.
	\end{align*}
	Counting fibers in two ways yields \eqref{eqn:semilocal-matching-coherence}.
\end{proof}

We are ready to state the semi-local geometric matching, and reduce it to the local case. Recall that $I_{\tilde{M}}(\tilde{\gamma}, \cdot)$ and $I^{\Endo}_{\tilde{M}}(\tilde{\gamma}, \cdot)$ are defined for all $\tilde{\gamma} \in D_{\diamondsuit, -}(\tilde{M}_V) \otimes \mes(M(F_V))^\vee$ if the Hypothesis \ref{hyp:ext-Arch} is verified at all Archimedean places in $V$.

\begin{theorem}\label{prop:semilocal-geometric}
	Assume that the local geometric Theorem \ref{prop:local-geometric} holds for all $v \in V$.
	\begin{enumerate}[(i)]
		\item For all $\tilde{\gamma} \in D_{\diamondsuit, -}(\tilde{M}_V) \otimes \mes(M(F_V))^\vee$, we have
		\[ I^{\Endo}_{\tilde{M}_V}\left(\tilde{\gamma}, \cdot \right) = I_{\tilde{M}_V}\left(\tilde{\gamma}, \cdot \right). \]
		\item For all $\mathbf{M}^! \in \Endo_{\elli}(\tilde{M})$ and $\delta \in SD_{\text{tr-orb}}(M^!(F_V)) \otimes \mes(M^!(F_V))^\vee$, we have
		\[ I_{\tilde{M}_V}\left( \trans_{\mathbf{M}^!, \tilde{M}}(\delta), \cdot \right) = I^{\Endo}_{\tilde{M}_V}\left( \mathbf{M}^!, \delta, \cdot \right). \]
	\end{enumerate}
\end{theorem}
\begin{proof}
	Consider (i). By linearity, we may assume $\tilde{\gamma} = \prod_{v \in V} \tilde{\gamma}_v$. By Remark \ref{rem:semi-local-orbint-ext}, the distribution $I_{\tilde{M}_V}\left(\tilde{\gamma}, \cdot \right)$ splits into local ones. On the other hand, $I^{\Endo}_{\tilde{M}_V}\left(\tilde{\gamma}, \cdot \right)$ splits in the same way by its very definition.
	
	As for (ii), we may assume $\delta = \prod_{v \in V} \delta_v$ so that $\trans_{\mathbf{M}^!, \tilde{M}}(\delta) \in D_{\diamondsuit, -}(\tilde{M}_V) \otimes \mes(M(F_V))^\vee$ decomposes accordingly. Hence $I_{\tilde{M}_V}\left(\trans_{\mathbf{M}^!, \tilde{M}}(\delta), \cdot \right)$ splits into local avatars. One the other hand, Proposition \ref{prop:semilocal-matching-coherence} (i) gives a similar splitting formula for $I^{\Endo}_{\tilde{M}_V}\left( \mathbf{M}^!, \delta, \cdot \right)$. It remains to apply the local results at each $v \in V$.
\end{proof}

The following result is unconditional.

\begin{proposition}\label{prop:Igeom-Endo-reindexed}
	For all $f \in \orbI_{\asp}\left(\tilde{G}_V, \tilde{K}_V\right) \otimes \mes(G(F_V))$, we have
	\[ I^{\Endo}_{\mathrm{geom}}(f) = \sum_{M \in \mathcal{L}(M_0)} \frac{|W^M_0|}{|W^G_0|} \sum_{\mathbf{M}^! \in \Endo_{\elli}(\tilde{M})} \iota(\tilde{M}, M^!) \sum_{\mathcal{O}^!} I^{\Endo}_{\tilde{M}_V}\left(\mathbf{M}^!, SA^{M^!}(V, \mathcal{O}^!), f^{G^!} \right) \]
	where $\mathcal{O}^!$ ranges over stable semisimple conjugacy classes in $M^!(F)$, and $f^{G^!} := \Trans_{\mathbf{G}^!, \tilde{G}}(f)$.
\end{proposition}
\begin{proof}
	Given $f$, plug
	\[ S(\mathbf{G}^!, M^!) := \sum_{\mathcal{O}^!} S^{G^!_V}_{M^!_V}\left( SA^{M^!}(V, \mathcal{O}^!) , f^{G^!} \right). \]
	into Proposition \ref{prop:combinatorial-summation} to obtain
	\begin{multline*}
		I^{\Endo}_{\mathrm{geom}}(f) = \sum_{\mathbf{G}^! \in \Endo_{\elli}(\tilde{G})} \iota(\tilde{G}, G^!) \sum_{M^! \in \mathcal{L}^{G^!}(M^!_0)} \frac{|W^{M^!}_0|}{|W^{G^!}_0|} \sum_{\mathcal{O}^!} S^{G^!_V}_{M^!_V}\left( SA^{M^!}(V, \mathcal{O}^!) , f^{G^!} \right) \\
		= \sum_{M \in \mathcal{L}^G(M_0)} \frac{|W^M_0|}{|W^G_0|} \sum_{\mathbf{M}^! \in \Endo_{\elli}(\tilde{M})} \iota(\tilde{M}, M^!) \\
		\sum_{\mathcal{O}^!} \sum_{s \in \Endo_{\mathbf{M}^!}(\tilde{G})} i_{M^!}(\tilde{G}, G^![s]) S^{G^![s]_V}_{M^!_V}\left( SA^{M^!}(V, \mathcal{O}^!) , f^{G^![s]} \right).
	\end{multline*}
	
	Since $V_{\mathrm{ram}}$ contains all dyadic and all Archimedean places, by \cite[VI.5.3 Lemme]{MW16-2}
	\[ S^{G^!_V}_{M^!_V}\left(SA^{M^!}(V, \mathcal{O}^!), f^{G^!}\right) = S^{G^!_V}_{M^!_V}\left(SA^{M^!}(V, \mathcal{O}^!), B^{\tilde{G}},  f^{G^!}\right). \]
	
	Furthermore, $SA^{M^!}(V, \mathcal{O}^!)$ can be replaced by $SA^{M^!}(V, \mathcal{O}^!)[s]$ since $\sum_{\mathcal{O}^!} = \sum_{\mathcal{O}^![s]}$ and Lemma \ref{prop:SA-central-twist} ensures $SA^{M^!}(V, \mathcal{O}^!)[s] = SA^{M^!}(V, \mathcal{O}^![s])$.
	
	All in all, by Definition \ref{def:IEndoMV} we see
	\begin{equation*}
		I^{\Endo}_{\mathrm{geom}}(f) = \sum_{M \in \mathcal{L}(M_0)} \frac{|W^M_0|}{|W^G_0|} \sum_{\mathbf{M}^!} \iota(\tilde{M}, M^!) \sum_{\mathcal{O}^!} I^{\Endo}_{\tilde{M}_V}\left(\mathbf{M}^!, SA^{M^!}(V, \mathcal{O}^!), f^{G^!} \right)
	\end{equation*}
	as desired.
\end{proof}

\begin{corollary}\label{prop:global-matching-reduction}
	Assume the validity of Theorem \ref{prop:matching-coeff-A}, i.e.\ the stabilization of coefficients, and the local geometric Theorem \ref{prop:local-geometric}. Then the statement in Theorem \ref{prop:geom-stabilization} holds, namely $I_{\mathrm{geom}} = I^{\Endo}_{\mathrm{geom}}$.
\end{corollary}
\begin{proof}
	In the formula of Proposition \ref{prop:Igeom-Endo-reindexed}, the $\mathcal{O}^!$-sum may be fibered over stable semisimple conjugacy classes $\mathcal{O}$ in $M(F)$. Apply Theorem \ref{prop:semilocal-geometric} and Definition \ref{def:AEndo} to transform it into
	\begin{equation*}
		I^{\Endo}_{\mathrm{geom}}(f) = \sum_{M \in \mathcal{L}(M_0)} \frac{|W^M_0|}{|W^G_0|} \sum_{\mathcal{O}} I_{\tilde{M}_V}\left( A^{\tilde{M}, \Endo}(V, \mathcal{O}), f \right).
	\end{equation*}

	As we are assuming Theorem \ref{prop:matching-coeff-A}, $A^{\tilde{M}, \Endo}(V, \mathcal{O}_V) = A^{\tilde{M}}(V, \mathcal{O}_V)$. This equals the expansion of $I_{\mathrm{geom}}(f)$ from Remark \ref{rem:Igeom-A}.
\end{proof}

\section{Semi-local weighted fundamental lemma}\label{sec:semi-local-WFL}
The following results are needed for \S\ref{sec:gdesc}.

Let $\tilde{G} \supset \tilde{M}$ be as before, with $M \in \mathcal{L}(M_0)$. Let $U$ be a finite set of places of $F$ such that $U \cap V_{\mathrm{ram}} = \emptyset$. Regard $K_U = \prod_{u \in U} K_u$ as a subgroup of $\tilde{G}_U$. We assume that $K_u$ is in good position relative to $M$ for each $u$.

Use the Haar measure on $G(F_U)$ such that $\mes(K_U) = 1$, and trivialize the lines $\mes(\cdots)$. Same for the stable side. Let $f_{K_U}$ be the unit of $\mathcal{H}_{\asp}( K_U \backslash \tilde{G}_U / K_U)$.

\begin{definition}
	\index{rGMgammaK}
	\index{rGEndoMdelta}
	As in the local theory (Definitions \ref{def:r-unramified}, \ref{def:rEndo-unramified}), we set
	\begin{align*}
		r^{\tilde{G}}_{\tilde{M}}(\tilde{\gamma}, K_U) & := J^{\tilde{G}}_{\tilde{M}}\left(\tilde{\gamma}, f_{K_U}\right), \\
		r^{\tilde{G}, \Endo}_{\tilde{M}}(\mathbf{M}^!, \delta) & := \sum_{s \in \Endo_{\mathbf{M}^!}(\tilde{G})} i_{M^!}(\tilde{G}, G^![s]) s^{G^![s]}_{M^!}(\delta[s]).
	\end{align*}
	for all $\mathbf{M}^! \in \Endo_{\elli}(\tilde{M})$ and all
	\[ \tilde{\gamma} \in D_{\mathrm{geom}, -}(\tilde{M}_U), \quad \delta \in SD_{\mathrm{geom}}(M^!(F_U)) \]
	where $s^{G^![s]}_{M^!}$ is the stable avatar of $r^{G^![s]}_{M^!}$ defined in \cite[VII.2.2]{MW16-2}; they do not depend on any auxiliary choice. Despite their semi-local nature, the subscripts $U$ are suppressed in $r^{\tilde{G}}_{\tilde{M}}$ and $s^{G^![s]}_{M^!}$ for typographical reasons.
\end{definition}

The non-endoscopic version $r^{\tilde{G}}_{\tilde{M}}(\tilde{\gamma}, K_U)$ already appeared in the invariant trace formula \cite{Li14b} for $\tilde{G}$. What remains is its stabilization, namely a \emph{semi-local weighted fundamental lemma}.

\begin{theorem}\label{prop:semi-local-WFL}
	\index{weighted fundamental lemma}
	Fix $U$ and $M \subset G$ as above. For any given datum $(\mathbf{M}^!, \delta)$, we have $r^{\tilde{G}, \Endo}_{\tilde{M}}\left( \mathbf{M}^!, \delta \right) = r^{\tilde{G}}_{\tilde{M}}\left(\trans_{\mathbf{M}^!, \tilde{M}}(\delta), K_U\right)$.
\end{theorem}
\begin{proof}
	Our strategy is to reduce this to the case $U = \{u\}$, which is simply Theorem \ref{prop:LFP-general}. This is achieved by the splitting formulas below. Consider $\tilde{\gamma}$ and $(\mathbf{M}^!, \delta)$ with $\tilde{\gamma} = \prod_u \tilde{\gamma}_u$ and $\delta = \prod_u \delta_u$. Define $M^U$ by $M^u := M$ for all $u \in U$. We contend that
	\begin{align*}
		r^{\tilde{G}}_{\tilde{M}}\left( \tilde{\gamma}, K_U \right) & = \sum_{L^U \in \mathcal{L}(M^U)} d^G_{M^U}\left( M, L^U \right) \prod_{u \in U} r^{\tilde{L}^u}_{\tilde{M}}\left( \tilde{\gamma}_u, K_u \cap L^u \right), \\
		r^{\tilde{G}, \Endo}_{\tilde{M}}\left(\mathbf{M}^!, \delta\right) & = \sum_{L^U \in \mathcal{L}(M^U)} d^G_{M^U}\left( M, L^U \right) \prod_{u \in U} r^{\tilde{L}^u, \Endo}_{\tilde{M}}\left(\mathbf{M}^!, \delta_u\right).
	\end{align*}
	
	The first identity reduces to the splitting formula for semi-local weighted orbital integrals, see \cite[VI.1.9]{MW16-2}.
	
	The second identity is based on the splitting formula for $s^{G^![s]}_{M^!}$ given in \cite[p.779]{MW16-2}. To deduce the splitting formula for $r^{\tilde{G}, \Endo}_{\tilde{M}}$, it remains to reiterate the proof of Proposition \ref{prop:semilocal-matching-coherence} (i).
\end{proof}
\chapter{Global descent}\label{sec:gdesc}
The goal is to prove the stabilization of coefficients in $I_{\mathrm{geom}} = I^{\tilde{G}}_{\mathrm{geom}}$, namely the identity $A^{\tilde{G}}(V, \mathcal{O}) = A^{\tilde{G}, \Endo}(V, \mathcal{O})$ (Theorem \ref{prop:matching-coeff-A}). Since the matching is tautological on $\GL$-factors, we may and do assume $\tilde{G} = \Mp(W)$.

Overall, we shall follow the strategy termed \emph{global descent} in \cite{Ar01} and \cite[VII]{MW16-2}. The first step is to reduce the problem to its elliptic version $A^{\tilde{G}}(S, \mathcal{O})_{\elli} = A^{\tilde{G}, \Endo}(S, \mathcal{O})_{\elli}$, where $\mathcal{O}$ is a stable semisimple conjugacy class in $G(F)$ and $S$ is sufficiently large relative to $\mathcal{O}$. This reduction is the content of \S\ref{sec:reduction-A-ell}, and relies on the weighted fundamental lemma.

The next step is to perform the Harish-Chandra descent globally. For this purpose, a finer description of the stable semisimple classes and endoscopic data under descent is necessary. It will permit us to break and descend $A^{\tilde{G}, \Endo}(S, \mathcal{O})_{\elli}$ to semisimple centralizers on the stable side. This is the content of \S\S\ref{sec:ss-reparameterization}--\ref{sec:descent-coeff}.

In order to match the expression so obtained with $A^{\tilde{G}}(S, \mathcal{O})_{\elli}$, we need the matching of unipotent coefficients in standard and nonstandard endoscopy. Thanks to the works of Arthur \cite{Ar03-3} and Mo\!eglin--Waldspurger \cite{MW16-2}, these identities are already available to us. In addition, we also need more sophisticated properties of the metaplectic transfer factors in \cite{Li11, Li15}, as well as the behavior of transfer factors in standard endoscopy under pure inner twists (see \S\ref{sec:endoscopy-linear}). These technical details occupy \S\S\ref{sec:gdesc-coh}--\ref{sec:gdesc-ascent}.

\section{Reduction to elliptic case}\label{sec:reduction-A-ell}
Fix a finite set of places $V$ of $F$ such that $V \supset V_{\mathrm{ram}} := V_{\mathrm{ram}}(\tilde{G})$. Concerning Theorem \ref{prop:matching-coeff-A}, it suffices to prove its global form, namely
\[ A^{\tilde{G}, \Endo}(V, \mathcal{O}) = A^{\tilde{G}}(V, \mathcal{O}) \]
for all stable semisimple conjugacy class $\mathcal{O}$ in $G(F)$.

Let $\mathbf{G}^! \in \Endo_{\elli}(\tilde{G})$ and let $\mathcal{O}^!$ be a stable semisimple conjugacy class in $G^!(F)$. For $S \supset V \cup S(\mathcal{O}^!)$, take a decomposition as in \eqref{eqn:SAell-induction}:
\begin{gather*}
	SA^{G^!}(S, \mathcal{O}^!)_{\elli} = \sum_{i=1}^{n(\mathcal{O}^!)} Sk_i^{G^!}(\mathcal{O}^!) \otimes SA_i^{G^!}(\mathcal{O}^!).
\end{gather*}

Let $\mathcal{O}$ be a stable semi-simple conjugacy class in $G(F)$ such that $S \supset V \cup S(\mathcal{O})$, and suppose that $\mathcal{O}^! \mapsto \mathcal{O}$. Lemma \ref{prop:admissible-Endo} implies $S \supset V \cup S(\mathcal{O}^!)$ since $V$ contains all dyadic places. Define
\begin{equation}\label{eqn:AEndo-ell-decomp}\begin{aligned}
	A^{\tilde{G}, \Endo}(S, \mathcal{O})_{\elli} & := \sum_{\substack{\mathbf{G}^! \in \Endo_{\elli}(\tilde{G}) \\ \mathcal{O}^! \mapsto \mathcal{O}}} \iota(\tilde{G}, G^!) \trans_{\mathbf{G}^!, \tilde{G}}\left( SA^{G^!}(S, \mathcal{O}^!)_{\elli} \right) \\
	& = \sum_{i=1}^{n(\mathcal{O})} k_i^{\tilde{G}, \Endo}(\mathcal{O}) \otimes A^{\tilde{G}, \Endo}_i(\mathcal{O}),
\end{aligned}\end{equation}
\index{AEndoSOell@$A^{\tilde{G}, \Endo}(S, \mathcal{O})_{\elli}$}
where the $\mathcal{O}^!$-sum is finite, and
\begin{align*}
	k_i^{\tilde{G}, \Endo}(\mathcal{O}) & \in D_{\mathrm{geom}, -}(\tilde{G}^V_S) \otimes \mes(G(F^V_S))^\vee , \\
	A^{\tilde{G}, \Endo}_i(\mathcal{O}) & \in D_{\mathrm{geom}, -}(\tilde{G}_V) \otimes \mes(G(F_V))^\vee
\end{align*}
are expressed in terms of the transfers of various $Sk_i^{G^!}(\mathcal{O}^!)$, $SA^{G^!}_i(\mathcal{O}^!)$ and $\iota(\tilde{G}, G^!)$. Note that $A^{\tilde{G}, \Endo}(S, \mathcal{O})_{\elli} \neq 0$ only when $\mathcal{O}$ is elliptic.

These recipes apply to Levi subgroups as well, with superscript $\tilde{G}$ modified accordingly. We are going to compare this with the $A^{\tilde{G}, \Endo}(V, \mathcal{O})$ from Definition \ref{def:AEndo}. 

\begin{proposition}[Cf.\ {\cite[VII.2.3 Proposition (i)]{MW16-2}}]\label{prop:matching-coeff-reduction-ell}
	For $\mathcal{O}$ and $S$ as above, we have
	\begin{equation*}
		A^{\tilde{G}, \Endo}(V, \mathcal{O}) = \sum_{M \in \mathcal{L}(M_0)} \frac{|W^M_0|}{|W^G_0|} \sum_{\substack{\mathcal{O}_M \\ \mathcal{O}_M \mapsto \mathcal{O}}} \sum_{i=1}^{n(\mathcal{O}_M)} r^{\tilde{G}}_{\tilde{M}}\left( k^{\tilde{M}, \Endo}_i(\mathcal{O}_M), K_S^V \right) A^{\tilde{M}, \Endo}_i(\mathcal{O}_M)^{\tilde{G}}
	\end{equation*}
	where $(\cdots)^{\tilde{G}} := (\cdots)^{\tilde{G}_V}$ denotes the parabolic induction of distributions.
\end{proposition}
\begin{proof}
	For each $\mathbf{G}^! \in \Endo_{\elli}(\tilde{G})$ and $\mathcal{O}^! \mapsto \mathcal{O}$, by \cite[VII.2.3 Proposition (ii)]{MW16-2} applied to $G^!$, we see $SA^{G^!}(V, \mathcal{O}^!)$ is equal to
	\[ \sum_{M^! \in \mathcal{L}^{G^!}(M_0^!)} \frac{|W^{M^!}_0|}{|W^{G^!}_0|} \sum_{\mathcal{O}^!_M \mapsto \mathcal{O}^!} \sum_{i=1}^{n(\mathcal{O}^!_M)} s^{G^!}_{M^!}\left( Sk_i^{M^!}(\mathcal{O}_M^!) \right) SA^{M^!}_i(\mathcal{O}_M^!)^{G^!} \]
	using chosen decompositions of $SA^{M^!}(S, \mathcal{O}_M^!)_{\elli}$. Now apply transfer to obtain
	\[ A^{\tilde{G}, \Endo}(V, \mathcal{O}) = \sum_{\mathbf{G}^!} \iota(\tilde{G}, G^!) \sum_{M^! \in \mathcal{L}(M^!_0)} \frac{|W^{M^!}_0|}{|W^{G^!}_0|} S(\mathbf{G}^!, M^!) \]
	where for the $\mathcal{O}$ chosen,
	\[ S(\mathbf{G}^!, M^!) := \sum_{\mathcal{O}^! \mapsto \mathcal{O}} \sum_{\mathcal{O}_M^! \mapsto \mathcal{O}^!} \sum_{i=1}^{n(\mathcal{O}_M^!)} s^{G^!}_{M^!}\left( Sk_i^{M^!}(\mathcal{O}_M^!) \right) \trans_{\mathbf{G}^!, \tilde{G}} \left( SA^{M^!}_i(\mathcal{O}_M^!)^{G^!} \right). \]
	The operations are justified by the fact that $S \supset S(\mathcal{O}) \supset S(\mathcal{O}^!) \supset S(\mathcal{O}_M^!)$.
	
	Proposition \ref{prop:combinatorial-summation} is applicable here, and yields
	\[ A^{\tilde{G}, \Endo}(V, \mathcal{O}) = \sum_{M \in \mathcal{L}(M_0)} \frac{|W^M_0|}{|W^G_0|} \sum_{\mathbf{M}^! \in \Endo_{\elli}(\tilde{M})} \iota(\tilde{M}, M^!) \sum_{s \in \Endo_{\mathbf{M}^!}(\tilde{G})} i_{M^!}(\tilde{G}, G^![s]) S(\mathbf{G}^![s], M^!). \]

	Using Proposition \ref{prop:central-twist-corr} (resp.\ Proposition \ref{prop:Levi-central-twist}) in the global (resp.\ semi-local) setting to rearrange the sum,
	\[ S(\mathbf{G}^![s], M^!) = \sum_{\mathcal{O}_M \mapsto \mathcal{O}} \sum_{\mathcal{O}_M^! \mapsto \mathcal{O}_M} \sum_{i=1}^{n(\mathcal{O}_M^!)} s^{G^!}_{M^!}\left( Sk_i^{M^!}(\mathcal{O}_M^!)[s] \right) \trans_{\mathbf{M}^!, \tilde{M}} \left( SA^{M^!}_i(\mathcal{O}_M^!) \right)^{\tilde{G}}. \]

	Collect terms to obtain $A^{\tilde{G}, \Endo}(V, \mathcal{O}) = \sum_M \frac{|W^M_0|}{|W^G_0|} \sum_{\mathcal{O}_M \mapsto \mathcal{O}} B(\tilde{M}, \mathcal{O}_M)$, where
	\begin{align*}
		B(\tilde{M}, \mathcal{O}_M) & := \sum_{\mathbf{M}^!} \iota(\tilde{M}, M^!) \sum_{\mathcal{O}_M^! \mapsto \mathcal{O}_M} b(M^!, \mathcal{O}_M^!), \\
		b(M^!, \mathcal{O}_M^!) & := \sum_{i=1}^{n(\mathcal{O}_M^!)} r^{\tilde{G}, \Endo}_{\tilde{M}}\left( \mathbf{M}^!, Sk_i^{M^!}(\mathcal{O}_M^!) \right) \trans_{\mathbf{M}^!, \tilde{M}}\left( SA^{M^!}_i(\mathcal{O}_{M^!}) \right)^{\tilde{G}} \\
		& = \sum_{i=1}^{n(\mathcal{O}_M^!)} r^{\tilde{G}}_{\tilde{M}}\left( \trans_{\mathbf{M}^!, \tilde{M}} (Sk_i^{M^!}(\mathcal{O}_M^!)), K_S^V \right) \trans_{\mathbf{M}^!, \tilde{M}}\left( SA^{M^!}_i(\mathcal{O}_{M^!}) \right)^{\tilde{G}} ,
	\end{align*}
	the last equality coming from Theorem \ref{prop:semi-local-WFL}. By recalling the formation of \eqref{eqn:AEndo-ell-decomp}, we see
	\[ B(\tilde{M}, \mathcal{O}_M) = \sum_{i=1}^{n(\mathcal{O}_M)} r^{\tilde{G}}_{\tilde{M}}\left( k^{\tilde{M}, \Endo}_i(\mathcal{O}_M), K_S^V \right) A_i^{\tilde{M}, \Endo}(\mathcal{O}_M)^{\tilde{G}}. \]
	This suffices to conclude the proof.
\end{proof}

\begin{corollary}
	Assume that
	\[ A^{\tilde{M}, \Endo}(S, \mathcal{O}_M)_{\elli} = A^{\tilde{M}}(S, \mathcal{O}_M)_{\elli} \]
	for all $M \in \mathcal{L}(M_0)$, $\mathcal{O}_M \mapsto \mathcal{O}$ and sufficiently large $S \supset S(\mathcal{O}) \cup V$. Then $A^{\tilde{G}, \Endo}(V, \mathcal{O}) = A^{\tilde{G}}(V, \mathcal{O})$.
\end{corollary}
\begin{proof}
	For all $M$ and $\mathcal{O}_M$, the assumption allows us to take $k_i = k_i^{\tilde{M}, \Endo}(\mathcal{O}_M)$ and $A_i = A_i^{\tilde{M}, \Endo}(\mathcal{O}_M)$ in the decomposition $A^{\tilde{M}}(S, \mathcal{O}_M)_{\elli} = \sum_i k_i \otimes A_i$ intervening in \eqref{eqn:Aell-induction}; the distributions with superscript $\Endo$ being the ones in \eqref{eqn:AEndo-ell-decomp}. The assertion follows immediately by comparing Proposition \ref{prop:matching-coeff-reduction-ell} and \eqref{eqn:Aell-induction}.
\end{proof}

Our goal is to establish Theorem \ref{prop:matching-coeff-A}. Hereafter, we fix the following data:
\begin{itemize}
	\item $\mathcal{O}$: an elliptic stable semisimple conjugacy class in $G(F)$;
	\item $\eta \in \mathcal{O}$ such that $G_\eta$ is quasisplit;
	\item $\mathcal{E}$: an $F$-pinning of $G_\eta$;
	\item since it is legitimate to enlarge $S$, we may assume $S \supset S(\mathcal{C}, K)'$ (Definition \ref{def:SOK}), where $\mathcal{C}$ is the $G(F)$-conjugacy class containing $\eta$, and that $\eta \in K_v$ and $K_v \cap G_\eta(F_v)$ of $G_\eta(F_v)$ coincides with the hyperspecial subgroup determined by $\mathcal{E}$ via Bruhat--Tits theory, for each $v \notin S$ (see \S\ref{sec:Aell}).
\end{itemize}

We are thus reduced to showing that
\begin{equation}\label{eqn:matching-coeff-Aell}
	A^{\tilde{G}, \Endo}(S, \mathcal{O})_{\elli} = A^{\tilde{G}}(S, \mathcal{O})_{\elli}
\end{equation}

\section{Re-parameterization}\label{sec:ss-reparameterization}
Let $H$ be any quasisplit connected reductive $F$-group, with a chosen Borel pair $(B,T)$ over $F$. Given $\mu \in T(\overline{F})$, we define:
\begin{center}\begin{tabular}{|c|c|} \hline
	$W^H$ & absolute Weyl group \\
	$\Sigma = \Sigma(H, T)$ & absolute roots \\
	$Z^1(\Gamma_F, W^H)$ & continuous $1$-cocycles \\
	$\Sigma(\mu) = \Sigma(H_\mu, T)$ & $\left\{ \alpha \in \Sigma: \alpha(\mu)=1 \right\}$ \\
	$\Sigma_+(\mu)$ & positive roots relative to $B \cap H_\mu$ \\
	$W^H(\mu) = W^{H_\mu}$ & absolute Weyl group for $T \subset H_\mu$ \\
	& (generated by reflections from $\Sigma_+(\mu)$) \\ \hline
\end{tabular}\end{center}

The shorthand $W(\mu) := W^H(\mu)$ will also be used. Define
\begin{align*}
	\underline{\mathrm{Stab}}(H(F)) & := \left\{\begin{array}{r|l}
		(\mu, \omega) & \mu \in T(\overline{F}), \; \omega \in Z^1(\Gamma_F, W^H) \;\text{s.t.} \\
		& \forall \sigma \in \Gamma_F, \; \omega(\sigma) \sigma \;\text{fixes}\; \mu
	\end{array}\right\}, \\
	\mathrm{Stab}(H(F)) & := \left\{\begin{array}{r|l}
		(\mu, \omega) & \mu \in T(\overline{F}), \; \omega \in Z^1(\Gamma_F, W^H) \;\text{s.t.} \\
		& \forall \sigma \in \Gamma_F, \; \omega(\sigma) \sigma \;\text{fixes}\; \mu, \\
		& \quad \omega(\sigma)\sigma \;\text{leaves}\; \Sigma_+(\mu) \;\text{invariant} 
	\end{array}\right\}.
\end{align*}
\index{Stab@$\underline{\mathrm{Stab}}$, $\mathrm{Stab}$}
\begin{itemize}
	\item We say $(\mu_1, \omega_1), (\mu_2, \omega_2) \in \underline{\mathrm{Stab}}(H(F))$ are \emph{equivalent} if $\mu_1 = \mu_2$ and $\omega_1(\sigma) \in W(\mu_1) \omega_2(\sigma)$ for all $\sigma$.
	
	\item Given $(\mu, \omega) \in \underline{\mathrm{Stab}}(H(F))$, we may adjust each $\omega(\sigma)$ by a unique element $\theta(\sigma) \in W(\mu)$ to make $\theta(\sigma)\omega(\sigma)\sigma$ fix $\Sigma_+(\mu)$. The result $\theta\omega$ turns out to be a cocycle, and it follows that $\mathrm{Stab}(H(F))$ affords a system of representatives for equivalence classes in $\underline{\mathrm{Stab}}(H(F))$.
	
	\item We say $(\mu, \omega), (\mu', \omega') \in \mathrm{Stab}(H(F))$ are \emph{conjugate} if there exists $w \in W^H$ such that $\mu' = w\mu w^{-1}$, $\omega'(\sigma) = w \omega(\sigma) \sigma(w^{-1})$ for all $\sigma$ and $\Sigma_+(\mu') = w \Sigma_+(\mu)$.
\end{itemize}

Given $\eta \in H(F)_{\mathrm{ss}}$, one can construct a pair $(\mu, \omega) \in \mathrm{Stab}(H(F))$ in the following way. Pick a Borel pair $(B', T')$ over $\overline{F}$ for $H$ such that $\eta \in T'(\overline{F})$, and complete it into a pinning $\mathcal{E}' = (B', T', (E'_\alpha)_\alpha)$ for $H_{\overline{F}}$. Choose a cochain $u_{\mathcal{E}'}: \Gamma_F \to H_{\mathrm{SC}}(\overline{F})$ (resp.\ $u_\eta: \Gamma_F \to H_{\mathrm{SC}, \eta}(\overline{F})$) such that $\Ad(u_{\mathcal{E}'}(\sigma)) \circ \sigma$ (resp.\ $\Ad(u_\eta(\sigma)) \circ \sigma$) preserves $\mathcal{E}'$ (resp.\ the Borel pair $(B' \cap H_\eta, T')$ of $H_\eta$) for all $\sigma$. Let $v_\eta(\sigma) := u_\eta(\sigma) \circ u_{\mathcal{E}}(\sigma)^{-1}$, which normalizes $T'$. On the other hand, take $h \in H(\overline{F})$ such that $\Ad(h)(B', T') = (B, T)$ and set $\mu := \Ad(h)(\eta) \in T(\overline{F})$. It is verified in \cite[VII.1.2 Proposition]{MW16-2} that
\begin{itemize}
	\item the image of $v_\eta(\sigma) T'$ under $\Ad(h)$ is an element of $W^H$, denoted as $\omega(\sigma)$, and $\omega \in Z^1(\Gamma_F, W^H)$;
	\item the pair $(\mu, \omega)$ is independent of the choices of $(E'_\alpha)_\alpha$, $u_{\mathcal{E}'}$, $u_\eta$, $h$ and it belongs to $\mathrm{Stab}(H(F))$;
	\item the conjugacy class of $(\mu, \omega)$ is independent of the choice of $(B', T')$;
	\item any two $\eta_1, \eta_2 \in H(F)_{\mathrm{ss}}$ are stably conjugate if and only if the conjugacy classes in $\mathrm{Stab}(H(F))$ so obtained are equal.
\end{itemize}

\begin{proposition}\label{prop:Stab-bijection}
	The recipe above yields
	\[ \mathrm{Stab}(H(F)) \twoheadrightarrow \left(\mathrm{Stab}(H(F))/\text{conj.} \right) \xleftarrow{1:1} H(F)_{\mathrm{ss}} \big/ \text{st.\ conj.} \]
\end{proposition}
\begin{proof}
	This is a paraphrase of \cite[VII.1.1--1.3]{MW16-2} in the untwisted quasisplit case.
\end{proof}

Given $(\mu, \omega) \in \mathrm{Stab}(H(F))$, by \cite[Lemma 3.3]{Ko82} we may choose $\epsilon \in H(F)_{\mathrm{ss}}$ in the corresponding stable conjugacy class such that $H_\epsilon$ is quasisplit.

\begin{lemma}\label{prop:mu-omega-fiber}
	Given $(\mu, \omega) \in \mathrm{Stab}(H(F))$ and a corresponding $\epsilon \in H(F)_{\mathrm{ss}}$ such that $H_\epsilon$ is quasisplit, set
	\begin{align*}
		\mathrm{Fix}(\mu, \omega) & := \left\{\begin{array}{r|l}
			w \in W^H & w\mu = \mu, \; \omega(\sigma) = w \omega(\sigma) \sigma(w^{-1}) \\
			& w \Sigma_+(\mu) = \Sigma_+(\mu)
		\end{array}\right\}, \\
		\Xi^H_\epsilon & := Z_H(\epsilon) / H_\epsilon \quad \text{(a finite $F$-group scheme)}.
	\end{align*}
	\index{Xi-epsilon}
	Then $\left| \mathrm{Fix}(\mu, \omega) \right| = |\Xi^H_\epsilon(F)|$. Moreover, the fiber of the surjection in Proposition \ref{prop:Stab-bijection} containing $(\mu, \omega)$ has cardinality equal to
	\[ |W^H| \cdot |W^H(\mu)|^{-1} \cdot |\Xi^H_\epsilon(F)|^{-1}. \]
\end{lemma}
\begin{proof}
	According to \cite[E-35, 4.1]{SS70}, $\Xi^H_\epsilon(\overline{F})$ is in bijection with those $w \in W^H$ satisfying $w\mu = \mu$ and $w \Sigma_+(\mu) = \Sigma_+(\mu)$. Unraveling the definitions, one can verify that $\Gamma_F$ acts on these $w$ via $w \mapsto \omega(\sigma) \sigma(w) \omega(\sigma)^{-1}$. This shows $|\mathrm{Fix}(\mu, \omega)| = \left| \Xi^H_\epsilon(F) \right|$.
	
	The second assertion follows easily: see \cite[p.813]{MW16-2}.
\end{proof}

The results above apply to local fields as well, and there is an obvious localization map $\mathrm{Stab}(H(F)) \to \mathrm{Stab}(H(F_v))$ for each place $v$ of $F$.

\begin{proposition}\label{prop:nr-Stab}
	Suppose that $H$ is an unramified group over a non-Archimedean local field $E$ with $\mathrm{char}(E)=0$. Thus the torus $T$ in the chosen Borel pair has a canonical $\mathfrak{o}_E$-model. Enlarge $(B, T)$ to an $E$-pinning to define a hyperspecial subgroup $K$ of $H(E)$ via Bruhat--Tits theory. Let the stable semisimple conjugacy class $\mathcal{O} \subset H(E)$ be parameterized by $(\mu, \omega) \in \mathrm{Stab}(H(E))$. Consider the conditions on $(\mu, \omega)$:
	\begin{enumerate}[(i)]
		\item $\mu \in T(\mathfrak{o}_{\overline{E}})$;
		\item for every $\alpha \in \Sigma$, either $\alpha(\mu) = 1$ or $\alpha(\mu) - 1 \in \mathfrak{o}_{\overline{E}}^\times$;
		\item the cocycle $\omega$ is trivial on the inertia $I_E \lhd \Gamma_E$.
	\end{enumerate}
	Assuming (i) and (ii), then $\mathcal{O}$ intersects $K$ if and only if (iii) holds.
\end{proposition}
\begin{proof}
	This rephrases \cite[VII.1.6 Lemme]{MW16-2}.
\end{proof}

\begin{corollary}\label{prop:S-nr}
	Let $F$ be a number field and $H$ be a connected reductive $F$-group, with chosen hyperspecial subgroups $K_v \subset H(F_v)$ for each $v \notin V_{\mathrm{ram}}$. Let $\mathcal{O} \subset H(F)$ be the stable class parameterized by $(\mu, \omega) \in \mathrm{Stab}(H(F))$. Then $S(\mathcal{O}, K)$ (Definition \ref{def:SOK}) equals the union of $V_{\mathrm{ram}}$ and the set of $v \notin V_{\mathrm{ram}}$ such that one of the conditions (i) --- (iii) in Proposition \ref{prop:nr-Stab} fails at $v$.
\end{corollary}
\begin{proof}
	Compare (i) --- (iii) with the Definitions \ref{def:SO} and \ref{def:SOK}. The condition (i) amounts to compactness in $G(F_v)$, and (ii) amounts to admissibility.
\end{proof}

Let us go back to the endoscopy for $\tilde{G} = \Mp(W)$. Note that $\Xi^G_\eta$ is trivial for all semisimple $\eta \in G(F)$.

Also note that since $G$ is split, the $\omega: \Gamma_F \to W^G$ in $(\mu, \omega) \in \mathrm{Stab}(G(F))$ are actually homomorphisms. Ditto for $(\mu^!, \omega^!) \in \mathrm{Stab}(G^!(F))$ where $\mathbf{G}^! \in \Endo_{\elli}(\tilde{G})$.

\begin{definition}
	\index{iotaGGO@$\iota(\tilde{G}, G^{"!}, \mathcal{O}^{"!})$}
	Let $\mathbf{G}^! \in \Endo_{\elli}(\tilde{G})$. Consider a stable semisimple class $\mathcal{O}^!$ in $G^!(F)$ with preimage $(\mu^!, \omega^!) \in \mathrm{Stab}(G^!(F))$. Set
	\begin{align*}
		\iota(\tilde{G}, G^!, \mathcal{O}^!) & = \iota(\tilde{G}, G^!, \mu^!, \omega^!) := \iota(\tilde{G}, G^!) \cdot |W^{G^!}(\mu^!)| \cdot |\mathrm{Fix}(\mu^!, \omega^!)| \\
		& = \iota(\tilde{G}, G^!) \cdot |W^{G^!_\epsilon}| \cdot |\Xi^{G^!}_\epsilon(F)|,
	\end{align*}
	where $\epsilon \in \mathcal{O}^!$ is chosen so that $G^!_\epsilon$ is quasisplit.
\end{definition}

Take $S \supset S(\mathcal{O}^!) \cup V_{\mathrm{ram}}$ and $f^!_S = \prod_{v \in S} f^!_v \in S\orbI(G^!(F_S)) \otimes \mes(G^!(F_S))$, supported in a small stably invariant neighborhood of $\epsilon$. Plugging these into \eqref{eqn:SAell-descent} and using Definition \ref{def:iota-const}, we obtain the
\begin{lemma}\label{prop:SAell-descent2}
	For $f^!_S$ as above. Then
	\begin{multline*}
		\iota(\tilde{G}, G^!, \mathcal{O}^!) S^{G^!_S}\left( SA^{G^!}(S, \mathcal{O}^!)_{\elli}, f^!_S \right) \\
		= \iota(\tilde{G}, G^!) |W^{G^!_\epsilon}| \cdot \frac{\tau(G^!)}{\tau(G^!_\epsilon)} \cdot S^{G^!_{\epsilon, S}}\left( SA^{G^!_\epsilon}_{\mathrm{unip}}(S), S\desc_\epsilon(f^!_S) \right) \\
		= |W^{G^!_\epsilon}| \cdot \frac{\tau(G)}{\tau(G^!_\epsilon)} \cdot S^{G^!_{\epsilon, S}}\left( SA^{G^!_{\epsilon, S}}_{\mathrm{unip}}(S), S\desc_\epsilon(f^!_S) \right).
	\end{multline*}
\end{lemma}

These constructions are somewhat parallel to \cite[VII.5.8]{MW16-2}.

Consider now the dual $\tilde{G}^\vee$ and its standard maximal torus $T^\vee$. Using \S\ref{sec:dual}, we have
\begin{equation}\label{eqn:EndoTG}
	\Endo_{\elli, T^\vee}(\tilde{G}) := \left\{ s \in T^\vee: s^2 = 1 \right\} \twoheadrightarrow \Endo_{\elli}(\tilde{G}).
\end{equation}
Abusing notations, we write $\mathbf{G}^! \in \Endo_{\elli, T^\vee}(\tilde{G})$ to mean an element $s$ of $\Endo_{\elli, T^\vee}(\tilde{G})$, giving rise to $\mathbf{G}^! \in \Endo_{\elli}(\tilde{G})$.

Recall that the elements in $\Endo_{\elli}(\tilde{G})$ are in bijection with pairs $(n', n'') \in \Z_{\geq 0}^2$ with $n' + n'' = n$. The fiber of $\Endo_{\elli, T^\vee}(\tilde{G}) \to \Endo_{\elli}(\tilde{G})$ over $(n' ,n'')$ or $\mathbf{G}^!$ has cardinality $\frac{n!}{(n')! (n'')!} = |W^G| / |W^{G^!}|$. Thus
\begin{equation*}
	A^{\tilde{G}, \Endo}(S, \mathcal{O})_{\elli} = |W^G|^{-1} \sum_{\mathbf{G}^! \in \Endo_{\elli, T^\vee}(\tilde{G})} |W^{G^!}| \iota(\tilde{G}, G^!) \sum_{\mathcal{O}^! \mapsto \mathcal{O}} \trans_{\mathbf{G}^!, \tilde{G}}\left( SA^{G^!}(S, \mathcal{O}^!)_{\elli} \right).
\end{equation*}

\begin{definition}[Cf.\ {\cite[p.764]{MW16-2}}]\label{def:Stab-trans}
	The datum $\mathbf{G}^! \in \Endo_{\elli, T^\vee}(\tilde{G})$ lifts the map of stable semisimple classes $\mathcal{O}^! \mapsto \mathcal{O}$ to a map $\mathrm{Stab}(G^!(F)) \to \mathrm{Stab}(G(F))$
	as follows.
	
	Let $(\mu^!, \omega^!) \in \mathrm{Stab}(G^!(F))$. Consider the map $\Psi$ in Definition \ref{def:corr-orbits} which is equivariant with respect to $W^{G^!} \hookrightarrow W^G$ and $\Gamma_F$. Set $\mu := \Psi(\mu^!)$. Then $\omega^!(\sigma)\sigma$ still fixes $\mu$.

	Since $\omega^!(\sigma)\sigma$ does not necessarily fix $\Sigma_+(\mu)$, this furnishes only an element in $\underline{\mathrm{Stab}}(G(F))$. Denote the corresponding representative in $\mathrm{Stab}(G(F))$ as $(\mu, \omega) \in \mathrm{Stab}(G(F))$; it is the required image of $(\mu^!, \omega^!)$.
\end{definition}

Therefore Lemma \ref{prop:mu-omega-fiber} implies that when $\mathcal{O}$ and $\mathbf{G}^! \in \Endo_{\elli, T^\vee}(\tilde{G})$ are chosen, 
\begin{multline*}
	\sum_{\mathcal{O}^! \mapsto \mathcal{O}} \trans_{\mathbf{G}^!, \tilde{G}}\left( SA^{G^!}(S, \mathcal{O}^!)_{\elli} \right) \\
	= \sum_{\substack{(\mu, \omega) \in \mathrm{Stab}(G(F)) \\ (\mu, \omega) \mapsto \mathcal{O} }} \sum_{\substack{ (\mu^!, \omega^!) \in \mathrm{Stab}(G^!(F)) \\ (\mu^!, \omega^!) \mapsto (\mu, \omega) }} \frac{|W^{G^!}(\mu^!)|}{|W^{G^!}|} \cdot |\Xi^{G^!}_\epsilon(F)| \cdot \trans_{\mathbf{G}^!, \tilde{G}}\left( SA^{G^!}(S, \mathcal{O}^!)_{\elli} \right).
\end{multline*}
where $\mathcal{O}^!$ denotes the stable class determined by $(\mu^!, \omega^!)$ and $\epsilon \in \mathcal{O}^!$ makes $G^!_\epsilon$ quasisplit. Hence
\begin{equation*}
	A^{\tilde{G}, \Endo}(S, \mathcal{O})_{\elli} = |W^G|^{-1} \sum_{\mathbf{G}^! \in \Endo_{\elli, T^\vee}(\tilde{G})} \sum_{(\mu, \omega) \mapsto \mathcal{O}} \sum_{(\mu^!, \omega^!) \mapsto (\mu, \omega)} \iota(\tilde{G}, G^!, \mathcal{O}^!) \trans_{\mathbf{G}^!, \tilde{G}}\left( SA^{G^!}(S, \mathcal{O}^!)_{\elli} \right).
\end{equation*}

Denote by $\mathrm{Fiber}(\mathcal{O})$ the fiber over $\mathcal{O}$ inside $\mathrm{Stab}(G(F))$. Set
\begin{equation}\label{eqn:AEndo-muomega}
	A^{\tilde{G}, \Endo}(S, \mu, \omega)_{\elli} := |W^G(\mu)|^{-1} \sum_{\mathbf{G}^! \in \Endo_{\elli, T^\vee}(\tilde{G})} \sum_{(\mu^!, \omega^!) \mapsto (\mu, \omega)} \iota(\tilde{G}, G^!, \mathcal{O}^!) \trans_{\mathbf{G}^!, \tilde{G}}\left( SA^{G^!}(S, \mathcal{O}^!)_{\elli} \right).
\end{equation}
\index{AGEndoS-muomega@$A^{\tilde{G}, \Endo}(S, \mu, \omega)_{\elli}$}

The earlier discussions, $\Xi^G_\eta = \{1\}$ for all $\eta \in \mathcal{O}$ and $|\mathrm{Fiber}(\mathcal{O})| = |W^G| \cdot |W^G(\mu)|^{-1}$ (Lemma \ref{prop:mu-omega-fiber}) lead to
\[ A^{\tilde{G}, \Endo}(S, \mathcal{O})_{\elli} = |\mathrm{Fiber}(\mathcal{O})|^{-1} \sum_{(\mu, \omega) \in \mathrm{Fiber}(\mathcal{O})} A^{\tilde{G}, \Endo}(S, \mu, \omega)_{\elli}. \]

All in all, our aim \eqref{eqn:matching-coeff-Aell} is reduced to the assertion that
\begin{equation}\label{eqn:AEndo-repara}
	A^{\tilde{G}, \Endo}(S, \mu, \omega)_{\elli} = A^{\tilde{G}}(S, \mathcal{O})_{\elli}.
\end{equation}
for all $(\mu, \omega) \in \mathrm{Fiber}(\mathcal{O})$.

\section{Global descent of endoscopic data}\label{sec:gdesc-endoscopy}
Fix $(\mu, \omega) \in \mathrm{Fiber}(\mathcal{O})$. Let $\eta \in \mathcal{O}$ and $S$ be chosen as at the end of \S\ref{sec:reduction-A-ell}; in particular, $G_\eta$ is quasisplit.

We identify $(G_\eta)_{\overline{F}}$ with the $\overline{F}$-group $\overline{G} := G_\mu$ equipped with the Borel pair $(\overline{B}, \overline{T})$ given by intersecting $(B, T)$, and the Galois action on $\overline{T}$ is twisted by $\omega \in Z^1(\Gamma_F, W^G)$. The resulting quasisplit $F$-group is isomorphic to $G_\eta$. Moreover, $T^\vee \subset \tilde{G}^\vee$, which also identifies with $X^*(T_{\overline{F}}) \otimes \CC^\times$, is the standard maximal torus in $\overline{G}^\vee \simeq G_\eta^\vee$.

\begin{definition}\label{def:ET}
	\index{ETell@$E_{T^\vee, \elli}(G_\eta, S)$}
	Following \cite[p.812]{MW16-2}, we define $E_{T^\vee, \elli}(G_\eta, S)$ to be the set of endoscopic data $\mathbf{H} = (H, \mathcal{H}, s)$ of $G_\eta$ such that
	\begin{itemize}
		\item $\mathbf{H}$ is unramified outside $S$,
		\item $\mathbf{H}$ is elliptic,
		\item $s \in T^\vee$.
	\end{itemize}
	For $\mathbf{H}_1, \mathbf{H}_2 \in E_{T^\vee, \elli}(G_\eta, S)$, we write $\mathbf{H}_1 \sim \mathbf{H}_2$ to indicate that $s_1, s_2 \in T^\vee$ are in the same $Z_{G_\eta^\vee}^{\Gamma_F}$-coset. Note that $\sim$ does not affect the underlying endoscopic group $H$.
\end{definition}

\begin{remark}\label{rem:H-hat}
	For every $\mathbf{H} \in E_{T^\vee, \elli}(G_\eta, S)$, the group $H$ carries a decomposition into a direct product of classical groups. This one can see by first decomposing $G_\eta$ by Proposition \ref{prop:ss-parameters}, then applying the explicit description \cite{Wal10} of elliptic endoscopic data to each direct factor. One can realize $H$ as $\overline{\hat{H}}$ for some $F$-group $\hat{H}$, in the notation of \S\ref{sec:nonstandard-endoscopy}.
\end{remark}

Let $\mathbf{G}^! \in \Endo_{\elli}(\tilde{G})$. We say $(\mu^!, \omega^!) \in \mathrm{Stab}(G^!(F))$ is elliptic if the corresponding stable conjugacy class $\mathcal{O}^!$ is. The elliptic pairs form a subset $\mathrm{Stab}(G^!(F))_{\elli}$. Given $\mathcal{O}^! \mapsto \mathcal{O}$, take $\epsilon \in \mathcal{O}^!$ with $G^!_\epsilon$ quasisplit. We have seen in \S\ref{sec:descent-endoscopy} how $\mathbf{G}^!$ can be descended at $(\eta, \epsilon)$ up to equivalence: the situation after descent is
\[\begin{tikzcd}[column sep=huge]
	G^!_\epsilon \arrow[dashed, dash, "\text{nonstd.\ endo.}", r] & \overline{G^!_\epsilon} \arrow[dashed, dash, "\text{endo.}", r] & G_\eta
\end{tikzcd}\]
where $\overline{G^!_\epsilon}$ is as in Notation \ref{nota:Lbar}. The underlying endoscopic datum $\mathbf{G}^!_\eta$ of $G_\eta$ is well-defined only up to equivalence.

Below is an refinement that yields an endoscopic datum $\mathbf{H}$ of $G_\eta$, not just an equivalence class thereof. The point is that one must start with $\mathbf{G}^! \in \Endo_{\elli, T^\vee}(\tilde{G})$ together with $(\mu^!, \omega^!)$ that maps to $(\mu, \omega) \in \mathrm{Fiber}(\mathcal{O})$.

As a preparation, let us explain the descent of the datum $s \in T^\vee$ in $\Endo_{\elli, T^\vee}(\tilde{G})$ which satisfies $s^2 = 1$. Decompose $G_\eta$ into
\[ G_\eta = G_\eta^+ \times G_\eta^- \times G_\eta^u \]
corresponding to the eigenvalues $+1$, $-1$ and $\neq \pm 1$ of $\eta$, respectively. Decompose $T^\vee \subset \overline{G}^\vee$ accordingly and write $s = (s^+, s^-, s^u)$. Set
\begin{equation}\label{eqn:sbars}
	\overline{s} := (s^+, -s^-, s^u) \in T^\vee.
\end{equation}

Note that $G_\eta^-$ (resp.\ $(G_\eta^-)^\vee$) is of the form $\Sp(2n_-)$ (resp.\ $\SO(2n_- + 1, \CC)$), where $2n_-$ is the multiplicity of $-1$ as an eigenvalue of $\eta$. The corresponding factor of $T^\vee$ is identified with the standard maximal torus in $\SO(2n_- + 1, \CC)$; this is how we interpret $-s^-$.

The $\overline{s}$ so obtained will be part of the descended endoscopic datum $\mathbf{H}$, to be explained below.

\begin{proposition}\label{prop:desc-endoscopy-refined}
	The descent of endoscopic data can be upgraded to a canonical map
	\begin{equation}\label{eqn:JH-map}
		\left\{\begin{array}{r|l}
			(\mathbf{G}^!, \mu^!, \omega^!) & \mathbf{G}^! \in \Endo_{\elli, T^\vee}(\tilde{G}), \; (\mu^!, \omega^!) \in \mathrm{Stab}(G^!(F))_{\elli} \\
			& (\mu^!, \omega^!) \mapsto (\mu, \omega) \\
			& S \supset S(\mathcal{O}^!, K^!) \;\text{where}\; \mathrm{Fiber}(\mathcal{O}^!) \ni (\mu^!, \omega^!)
		\end{array}\right\} \to  E_{T^\vee, \elli}(G_\eta, S).
	\end{equation}
\end{proposition}
\begin{proof}
	Consider $s \in T^\vee$ with $s^2 = 1$ and $(\mu^!, \omega^!) \mapsto (\mu, \omega)$. Define $\overline{s} \in T^\vee \subset \overline{G}^\vee$ by \eqref{eqn:sbars}.

	Take $H^\vee := Z_{\overline{G}^\vee}(\overline{s})^\circ$ equipped with the Borel pair by intersecting $(\overline{B}, \overline{T})$. By the construction of $(\mu^!, \omega^!) \mapsto (\mu, \omega)$, there exists a unique continuous map $\omega_H: \Gamma_F \to W^G(\mu)$ such that
	\[ \omega^!(\sigma) = \omega_H(\sigma) \omega(\sigma), \quad \sigma \in \Gamma_F.  \]
	One readily verifies that $\omega_H$ is a $1$-cocycle with respect to the $\Gamma_F$-action on $W^G(\mu)$ given by $\sigma \mapsto \Ad(\omega(\sigma)) \circ \sigma$. One can pull $\omega_H$ back to $\Weil{F}$.

	Take a representative $g(w) \in \overline{G}^\vee$ for $\omega_H(w)$ for all $w \in \Weil{F}$, continuous in $w$. Set
	\begin{equation}\label{eqn:H-cocycle}
		\mathcal{H} := \lrangle{H^\vee, g(w) \rtimes w : w \in \Weil{F}} \; \subset \Lgrp{G_\eta}.
	\end{equation}
	It sits in an extension $1 \to H^\vee \to \mathcal{H} \to \Weil{F} \to 1$, thus $\Weil{F}$ maps continuously into $\mathrm{Out}(H^\vee)$, and so does $\Gamma_F$. Equip $H^\vee$ with the corresponding pinned $\Gamma_F$-action. The restriction of that action on $T^\vee$ is
	\[ \sigma \mapsto \xi_H(\sigma) \omega_H(\sigma), \quad \sigma \in \Gamma_F \]
	where $\xi_H(\sigma) \in W^{H^\vee}$ is the unique element such that $\xi_H(\sigma) \omega_H(\sigma)$ fixes the Borel pair for $H^\vee$. Construct the corresponding quasisplit $F$-group $H$. The image of $(\mathbf{G}^!, \mu^!, \omega^!)$ equals
	\[ \mathbf{H} := (H, \mathcal{H}, \overline{s}). \]
	
	It is routine to check that $\mathbf{H}$ is equivalent to the descended endoscopic datum $\mathbf{G}^!_\eta$ of $G_\eta$ in \cite[\S 7]{Li11}.
	
	The ellipticity of $(\mu^!, \omega^!)$ implies that of $(\mu, \omega)$ by Lemma \ref{prop:ellipticity-transfer}. Since $\mathbf{G}^!$ is elliptic and $(\mu^!, \omega^!)$, $(\mu, \omega)$ both parameterize elliptic classes, $\mathbf{H}$ is elliptic as well.

	It remains to show that $\mathbf{H}$ is unramified outside $S$. Suppose $v \notin S$. Using Corollary \ref{prop:S-nr} and the assumptions on $S$, we see that $\omega$ and $\omega^!$ are both trivial on $I_{F_v}$, and $G_\eta$, $G^!_\epsilon$ (thus $\overline{G^!_\epsilon}$) are unramified by \cite[7.1 Proposition]{Ko86}. In particular, $\omega_H|_{I_{F_v}} = 1$ and $1 \rtimes I_{F_v} \subset \mathcal{H}$.
\end{proof}

\begin{remark}\label{rem:H-fiber-determination}
	The construction above follows \cite[VII.1.7 and 5.2]{MW16-2}. It also implies $(\mu^!, \omega^!)$ is determined by $\mathbf{H}$ and $(\mu, \omega)$. Indeed, $\mu^! = \Psi^{-1}(\mu)$ is known. On the other hand, $W^H \simeq W^{H^\vee}$, and the identification between $H$ and $\overline{G^!_\epsilon}$ gives an identification $\Sigma_+(\mu^!) \simeq \Sigma^H_+$ that involves only some ``nonstandard'' rescaling. The construction of $H$ implies
	\[\omega^!(\sigma) \in W^H \cdot \text{(the pinned $\Gamma_F$-action on $H^\vee$)} \cdot \omega(\sigma) \]
	for all $\sigma \in \Gamma_F$. Since $\omega^!(\sigma) \sigma$ fixes $\Sigma_+(\mu^!)$, and when $w \in W^H$ we have $w\Sigma_+(\mu^!) = \Sigma_+(\mu^!)$ if and only if $w=1$, the above determines $\omega^!(\sigma)$.
\end{remark}

\begin{lemma}\label{prop:J-mult-one}
	The map \eqref{eqn:JH-map} is bijective.
\end{lemma}
\begin{proof}
	Identify $G_{\eta, \overline{F}}$ with $\overline{G}$, and similarly for their dual groups. Fix $\mathbf{H} \in E_{T^\vee, \elli}(G_\eta, S)$.
	
	We begin by proving that the fiber over $\mathbf{H}$ has at most one element. The determination of $s$ is clear from \eqref{eqn:sbars} since $T^\vee$ is shared by $\tilde{G}^\vee$ and $\overline{G}^\vee$. The determination of $(\mu^! , \omega^!)$ has been given in Remark \ref{rem:H-fiber-determination}.

	It remains to show that the fiber over every $\mathbf{H} = (H, \mathcal{H}, \overline{s})$ is non-empty. For each direct factor of $G_\eta$, say a symplectic or unitary group, $\overline{s}^2 = 1$ always holds; cf.\ the explicit description in \cite[\S 1.8]{Wal10} for the local case. As $T^\vee$ is shared by $\tilde{G}^\vee$ and $\overline{G}^\vee$, we define $s$ from \eqref{eqn:sbars} and obtain $\mathbf{G}^! \in \Endo_{\elli, T^\vee}(\tilde{G})$.
	
	By the general theory of endoscopy, using $\mathcal{H} \hookrightarrow \Lgrp{G_\eta}$ one can obtain a continuous map $g: \Weil{F} \to \overline{G}^\vee$ is given such that
	\begin{compactitem}
		\item \eqref{eqn:H-cocycle} holds;
		\item $\Ad(g(\cdot))$ stabilizes the Borel pair of $H^\vee$;
		\item it induces a continuous $1$-cocycle $\omega_H: \Gamma_F \to W^{\overline{G}}$, where the $\Gamma_F$-action on $W^{\overline{G}}$ is induced by its action on $\overline{T} \subset \overline{G}$, involving a twist by $\omega$.
	\end{compactitem}
	Explicit choices of $g$ can be found in \cite{Wal10}. If $\mathbf{H}$ is unramified at a place $v$, then $\omega_H|_{I_{F_v}} = 1$.
	
	Next, return to the description of $\mathrm{Stab}(G^!(F)) \to \mathrm{Stab}(G(F))$ in Definition \ref{def:Stab-trans}. Define
	\[ \mu^! := \Psi^{-1}(\mu), \quad \underline{\omega^!} := \omega_H \omega : \Gamma_F \to W^G. \]
	We claim that
	\[ (\mu^!, \underline{\omega^!}) \in \underline{\mathrm{Stab}}(G^!(F)). \]
	Once this is verified, its representative $(\mu^!, \omega^!) \in \mathrm{Stab}(G^!(F))$ will map to $(\mu, \omega)$ under Definition \ref{def:Stab-trans}, and $(\mathbf{G}^!, \mu^!, \omega^!) \mapsto \mathbf{H}$ by construction. Moreover, the ellipticity of $(\mu^!, \omega^!)$ will follow easily from that of $\mathbf{H}$ and $(\eta, \omega)$.

	\begin{enumerate}
		\item To prove the claim, observe first that by \eqref{eqn:H-cocycle}, the $\Gamma_F$-action on $T^\vee \subset H^\vee$ determined by $1 \to H^\vee \to \mathcal{H} \to \Weil{F} \to 1$ comes from the one on $T^\vee \subset \overline{G}^\vee$ (namely the homomorphism $\omega: \Gamma_F \to W^G$) followed by $\omega_H$; this is just $\underline{\omega}^!$. Hence $\underline{\omega}^!$ defines a homomorphism $\Gamma_F \to W^{H^\vee}$.
	
		As a consequence, we also deduce $\underline{\omega}^!(\sigma) \in W^{\overline{G}^\vee} = W^{\overline{G}}$ for all $\sigma \in \Gamma_F$.

		\item Let us show $\underline{\omega}^!(\sigma)$ stabilizes $s \in T^\vee$ for all $\sigma \in \Gamma_F$. As $Z_{\tilde{G}^\vee}(s)$ is connected, it will then follow that $\underline{\omega^!}$ takes value in $W^{G^!}$, by \cite[E-35, 4.1]{SS70}.
		
		Inspect the relation \eqref{eqn:sbars} between $s$ and $\overline{s}$: they only differ by $-1$ on the factor $(G_\eta^-)^\vee$ of $\overline{G}^\vee \simeq G_\eta^\vee$. The meaning of $-1$ is explicated in the discussions after \eqref{eqn:sbars}; as such, it is evidently fixed by the Weyl group of $(G^-_\eta)^\vee$. Hence it suffices to show $\underline{\omega}^!(\sigma)\overline{s} = \overline{s}$.
		
		The explicit description in \cite{Wal10} decomposes $H = H_+ \times H_-$ according to the signs in $\overline{s}$, and one checks case-by-case that $\overline{s} \in (H^\vee)^{\Gamma_F}$. As noted in the first step, this entails $\underline{\omega}^!(\sigma)\overline{s} = \overline{s}$.

		\item Also, $\underline{\omega^!}(\sigma)\sigma$ stabilizes $\mu \in T(\overline{F})$ since $\omega_H(\sigma) \in W^{\overline{G}}$ and $\omega_H(\sigma) \omega(\sigma) \sigma\mu = \omega_H(\sigma)\mu = \mu$. As $\Psi$ is equivariant with respect to $\Gamma_F$ and $W^{G^!} \hookrightarrow W^G$, it stabilizes $\mu^!$ as well. The claim follows.
	\end{enumerate}
	
	Finally, let us verify $S \supset S(\mathcal{O}^!, K^!)$ using Corollary \ref{prop:S-nr}. For $v \notin S$, we have $\mu \in T(\mathfrak{o}_{\overline{F_v}})$ and $\omega|_{I_{F_v}} = 1$. Thus $\mu^! = \Psi^{-1}(\mu) \in T^!(\mathfrak{o}_{\overline{F_v}})$; this verifies the condition (i). Consider (ii): if $\beta$ is a long (resp.\ short) root of $G^!$, corresponding to a short (resp.\ long) root $\alpha$ of $G$, then $\beta(\mu^!) - 1 = \alpha(\mu) - 1$ (resp.\ $(\beta(\mu^!) - 1)(\beta(\mu^!) + 1) = \alpha(\mu) - 1$), so $\beta(\mu^!) - 1$ is either zero or invertible in $\mathfrak{o}_{F_v}$ since so is $\alpha(\mu) - 1$. The minus signs in $\Psi$ do not bother here.
	
	Consider (iii). Since $\mathbf{H}_v$ is unramified, $\omega_H|_{I_{F_v}} = 1$, hence $\underline{\omega^!}|_{I_{F_v}} = 1$. The same property passes easily from $(\mu^!, \underline{\omega^!})$ to its representative $(\mu^!, \omega^!) \in \mathrm{Stab}(G^!(F))$.
\end{proof}

\begin{definition}
	\index{JH@$\mathcal{J}(\mathbf{H})$, $\mathcal{J}(\mathbf{H}/\sim)$}
	For all $\mathbf{H} \in E_{T^\vee, \elli}(G_\eta, S)$, let
	\begin{align*}
		\mathcal{J}(\mathbf{H}) & := \text{the fiber of \eqref{eqn:JH-map} over}\; \mathbf{H}, \\
		\mathcal{J}(\mathbf{H}/\sim) & := \bigsqcup_{\mathbf{H}' \sim \mathbf{H}} \mathcal{J}(\mathbf{H}').
	\end{align*}
\end{definition}

It follows from Lemma \ref{prop:J-mult-one} that $\mathcal{J}(\mathbf{H})$ is a singleton for every $\mathbf{H} \in E_{T^\vee, \elli}(G_\eta, S)$.

Continue the discussions on \eqref{eqn:AEndo-repara}. In order apply the method of descent, recall that $SA^{G^!}(S, \mathcal{O}^!) = 0$ unless $\mathcal{O}^!$ is elliptic and $S \supset S(\mathcal{O}^!, K^!)$. It is thus justified to fiber the sum in \eqref{eqn:AEndo-muomega} over $E_{T^\vee, \elli}(G_\eta, S)$ to obtain
\begin{equation*}
	A^{\tilde{G}, \Endo}(S, \mu, \omega)_{\elli} = \sum_{\mathbf{H} \in E_{T^\vee, \elli}(G_\eta, S)} A^{\tilde{G}, \Endo}(S, \mathbf{H})_{\elli},
\end{equation*}
where
\begin{equation}\label{eqn:AGEndoH}
	A^{\tilde{G}, \Endo}(S, \mathbf{H})_{\elli} :=
	|W^G(\mu)|^{-1} \sum_{(\mathbf{G}^!, \mu^!, \omega^!) \in \mathcal{J}(\mathbf{H})} \iota(\tilde{G}, G^!, \mathcal{O}^!) \trans_{\mathbf{G}^!, \tilde{G}}\left( SA^{G^!}(S, \mathcal{O}^!)_{\elli} \right)
\end{equation}
with $\mathcal{O}^! \subset G^!(F)$ satisfies $(\mu^!, \omega^!) \in \mathrm{Fiber}(\mathcal{O}^!)$.
\index{AGEndoH@$A^{\tilde{G}, \Endo}(S, \mathbf{H})_{\elli}$, $A^{\tilde{G}, \Endo}(S, \mathbf{H}/\sim)_{\elli}$}

Hence, for all $(\mu, \omega) \in \mathrm{Fiber}(\mathcal{O})$, we get
\begin{equation}\label{eqn:AGEndoHsim}\begin{aligned}
	A^{\tilde{G}, \Endo}(S, \mu, \omega)_{\elli} & = \sum_{\mathbf{H}/\sim} A^{\tilde{G}, \Endo}(S, \mathbf{H}/\sim)_{\elli}, \\
	A^{\tilde{G}, \Endo}(S, \mathbf{H}/\sim)_{\elli} & := \sum_{\substack{\mathbf{H}' \sim \mathbf{H}}} A^{\tilde{G}, \Endo}(S, \mathbf{H}')_{\elli}.
\end{aligned}\end{equation}

\section{Descent of coefficients}\label{sec:descent-coeff}
The scenario is the same as \S\ref{sec:gdesc-endoscopy}. In order to localize the problem, for every place $v$ of $F$ we use Definition \ref{def:Y-set} to define the pointed set
\[ \mathcal{Y}_v := \mathcal{Y}^{G_{F_v}}(\eta_v) \; \subset G(\overline{F_v}). \]
It is acted by $G_\eta(\overline{F_v}) \times G(F_v)$. Each $G_\eta(\overline{F_v})$-orbit $G_\eta(\overline{F}_v) y_v$ in $\mathcal{Y}_v$ gives rise to the stable conjugate $\eta[y_v] \in G(F_v)$, whereas $\eta[y_v g_v] = g_v^{-1} \eta[y_v] g_v$ for all $g_v \in G(F_v)$.

Each $y_v \in \mathcal{Y}_v$ affords a pure inner twist
\[ \Ad(y_v): G_{\eta[y_v], \overline{F_v}} \rightiso G_{\eta, \overline{F_v}}, \]
making every $\mathbf{H} \in E_{T^\vee, \elli}(G_\eta, S)$ into an elliptic endoscopic datum of $G_{\eta[y_v]}$. Those $y_v$ for which $\mathbf{H}$ is relevant form a subset $\mathcal{Y}_v^{\mathrm{rel}}(\mathbf{H}) \subset \mathcal{Y}_v$. Define
\[ \mathcal{Y}_S := \prod_{v \in S} \mathcal{Y}_v, \quad \mathcal{Y}_S^{\mathrm{rel}}(\mathbf{H}) := \prod_{v \in S} \mathcal{Y}_v^{\mathrm{rel}}(\mathbf{H}). \]
\index{YS@$\mathcal{Y}_S$, $\mathcal{Y}_S^{\mathrm{rel}}(\mathbf{H})$}

For all $y_S \in \mathcal{Y}_S$, define $\eta[y_S] := \prod_{v \in S} \eta[y_v] \in G(F_S)$.

\begin{lemma}
	We have $\mathcal{Y}_S^{\mathrm{rel}}(\mathbf{H}) \neq \emptyset$ for all $\mathbf{H}$; in fact, it contains the base-point $(1_v)_{v \in S}$.
\end{lemma}
\begin{proof}
	The simple reason is that $G_\eta$ is assumed to be quasisplit.
\end{proof}

Choose a representative in $E_{T^\vee, \elli}(G_\eta, S)$ for each class in $E_{T^\vee, \elli}(G_\eta, S) / \sim$. Given such a representative $\mathbf{H} = (H, \mathcal{H}, \overline{s})$ and $y_S \in \mathcal{Y}_S^{\mathrm{rel}}(\mathbf{H})$, take the semi-local transfer factors
\begin{align*}
	\Delta[y_S] & = \Delta_{\mathbf{H}}[y_S] := \prod_{v \in S} \Delta_v[y_v], \\
	\Delta_v [y_v] & = \Delta_v[y_v]: H_{\mathrm{reg}}(F_v) \times G_{\eta[y_v], \mathrm{reg}}(F_v) \to \CC
\end{align*}
determined by the quasisplit pure inner twists $\Ad(y_v)$ and the chosen $F$-pinning $\mathcal{E}$ of $G_\eta$. See \S\ref{sec:endoscopy-linear} for a review.
\index{Delta-yS@$\Delta[y_S]$}

\begin{notation}\index{Hstd@$\mathbf{H}_{\mathrm{std}}$}
	In what follows, we denote such ``standard'' representatives in a $\sim$ class by $\mathbf{H}_{\mathrm{std}}$.
\end{notation}

Based on the assumptions imposed on $\eta \in G(F)$ at the end of \S\ref{sec:reduction-A-ell}, we define $\eta^S \in K^S$ and $\tilde{\eta}_S = \prod_{v \in S} \tilde{\eta}_v$ via \eqref{eqn:extraction-S}.

Let $f_S = \prod_{v \in S} f_v \in \orbI_{\asp}(\tilde{G}_S) \otimes \mes(G(F_S))$, supported in a $\bmu_8$-invariant, small and stably invariant neighborhood of $\tilde{\eta}_S$. For all $\mathbf{H} = (H, \mathcal{H}, s) \in E_{T^\vee, \elli}(G_\eta, S)$, $y_S \in \mathcal{Y}_S^{\mathrm{rel}}(\mathbf{H})$ and $v \in S$, take $\hat{H}$ with $H \simeq \overline{\hat{H}}$ (Remark \ref{rem:H-hat}), so that we have the maps
\begin{equation}\label{eqn:AEndo-f-transfer} \begin{tikzcd}[column sep=small]
	\orbI_{\asp}(\tilde{G}_v) \otimes \mes(G(F_v)) \arrow[d, "{\desc^{\tilde{G}}_{\tilde{\eta}[y_v]}}"'] & f_v \arrow[mapsto, d] \\
	\orbI(G_{\eta[y_v]}(F_v)) \otimes \mes(G_{\eta[y_v]}(F_v)) \arrow[d, "\text{transfer}"'] & {f_v[y_v]} \arrow[mapsto, d] \\
	S\orbI(H(F_v)) \otimes \mes(H(F_v)) \arrow[d, "\sim" sloped, "\text{nonstandard transfer}"'] & {\overline{f}[y_v]} \arrow[mapsto, d] \\
	S\orbI(\hat{H}(F_v)) \otimes \mes(\hat{H}(F_v)) & {\hat{f}[y_v]}
\end{tikzcd}\end{equation}
\index{fyv@$f_v[y_v]$, $\overline{f}[y_v]$, $\hat{f}[y_v]$}
and taking products over $v$ yields $f_S \mapsto f_S[y_S] \mapsto \overline{f}[y_S] \mapsto \hat{f}[y_S]$. Here,
\begin{itemize}
	\item the metaplectic stable conjugate $\tilde{\eta}[y_v] \in \tilde{G}_v$ of $\tilde{\eta}_v$ is defined in \cite[Définition 4.1.3]{Li15}, which depends only on $G_\eta(\overline{F_v}) y_v$ and satisfies $\tilde{\eta}[y_v g_v] = g_v^{-1} \tilde{\eta}[y_v] g_v$ for $g_v \in G(F_v)$;
	\item the transfer from $G_{\eta[y_v]}(F_v)$ to $H(F_v)$ uses the representative $\mathbf{H}_{\mathrm{std}}$ in the $\sim$ class of $\mathbf{H}$ and the $\Delta_{\mathbf{H}_{\mathrm{std}}}[y_S]$ defined earlier, noting that $\sim$ has no effect on the endoscopic groups and the correspondence of orbits;
	\item the nonstandard transfer between $H(F_v)$ and $\hat{H}(F_v)$ is described in \S\ref{sec:nonstandard-endoscopy}; it depends only on the identification $H \simeq \overline{\hat{H}}$.
\end{itemize}

\begin{definition}\label{def:djyS}
	\index{djyS@$d_j(y_S)$}
	To $j = (\mathbf{G}^!, \mu^!, \omega^!) \in \mathcal{J}(\mathbf{H}/\sim)$, we have the transfer factor
	\[ \Delta_j = \prod_{v \in S} \Delta_{\mathbf{G}^!_v, \tilde{G}_v} \]
	for $\tilde{G}_S$. Let $y_S \in \mathcal{Y}_S^{\mathrm{rel}}(\mathbf{H}) = \mathcal{Y}_S^{\mathrm{rel}}(\mathbf{H}_{\mathrm{std}})$, define
	\[ d_j(y_S) := \prod_{v \in S} d_j(y_v) \in \CC^\times \]
	by \eqref{eqn:dy}. In other words, $d_j(y_v)$ is the ratio of the descent at $(\epsilon, \tilde{\eta}[y_v])$ of $\Delta_j$ over $\Delta_{\mathbf{H}_{\mathrm{std}}}[y_v]$.
\end{definition}

Given $\mathbf{H}$, choose a system $\dot{\mathcal{Y}}_S^{\mathrm{rel}}(\mathbf{H})$ of representatives of $G_\eta(\overline{F_S}) \backslash \mathcal{Y}_S^{\mathrm{rel}}(\mathbf{H}) / G(F_S)$ in $\mathcal{Y}_S^{\mathrm{rel}}(\mathbf{H})$. Note that $\mathcal{Y}_S^{\mathrm{rel}}(\mathbf{H})$ and $\dot{\mathcal{Y}}_S^{\mathrm{rel}}(\mathbf{H})$ depend only on the class $\mathbf{H}/\sim$, since the notion of relevance is not affected by automorphisms of endoscopic data. Thus we can write $\dot{\mathcal{Y}}_S^{\mathrm{rel}}(\mathbf{H} / \sim)$, etc.

\begin{lemma}\label{prop:AEndoH-descent-0}
	For all class $\mathbf{H}/\sim$ and $f_S$ as above, we have
	\begin{multline*}
		I^{\tilde{G}_S}\left( A^{\tilde{G}, \Endo}(S, \mathbf{H}/\sim)_{\elli}, f_S \right) \\
		= |W^G(\mu)|^{-1} |W^H| \cdot \frac{\tau(G)}{\tau(H)} \cdot \sum_{y_S \in \dot{\mathcal{Y}}_S^{\mathrm{rel}}(\mathbf{H}/\sim)} S^{H_S}\left( SA^H_{\mathrm{unip}}(S), \overline{f}[y_S] \right) \cdot \sum_{j \in \mathcal{J}(\mathbf{H}/\sim)} d_j(y_S).
	\end{multline*}
\end{lemma}
\begin{proof}
	For $j = (\mathbf{G}^!, \mu^!, \omega^!) \in \mathcal{J}(\mathbf{H}/\sim)$, write $f_j^! := \Trans_{\mathbf{G}^!, \tilde{G}}(f_S)$. Consider any $\mathbf{H}$ in the $\sim$ class. Plug Lemma \ref{prop:SAell-descent2} into \eqref{eqn:AGEndoH} to obtain
	\begin{multline*}
		I^{\tilde{G}_S}\left( A^{\tilde{G}, \Endo}(S, \mathbf{H})_{\elli}, f_S \right) \\
		= |W^G(\mu)|^{-1} \sum_{j \in \mathcal{J}(\mathbf{H})} |W^{G^!_\epsilon}| \cdot \frac{\tau(G)}{\tau(G^!_\epsilon)} \cdot S^{G^!_{\epsilon, S}}\left( SA^{G^!_\epsilon}_{\mathrm{unip}}(S), S\desc^{G^!}_\epsilon (f_j^!) \right).
	\end{multline*}

	Next, take $j \in \mathcal{J}(\mathbf{H})$. The descent of endoscopy in \S\ref{sec:descent-endoscopy} furnishes an isomorphism $\theta_j: \hat{H} \rightiso G^!_\epsilon$. By the results cited there --- specifically, the transposed form of Lemma \ref{prop:transfer-descent-y} --- we have
	\[ \theta_j^* \left( S\desc^{G^!}_\epsilon \left(f_j^! \right) \right) = \sum_{y_S \in \dot{\mathcal{Y}}_S^{\mathrm{rel}}(\mathbf{H} / \sim)} d_j(y_S) \hat{f}[y_S] . \]
	In view of the canonicity of $SA_{\mathrm{unip}}$ \cite[VI.2.2, 5.2]{MW16-2}, we may pass to $\hat{H}$ and deduce
	\[ S^{G^!_{\epsilon, S}}\left( SA^{G^!_\epsilon}_{\mathrm{unip}}(S), S\desc^{G^!}_\epsilon (f_j^!) \right) = \sum_{y_S \in \dot{\mathcal{Y}}_S^{\mathrm{rel}}(\mathbf{H} / \sim)} d_j(y_S) S^{\hat{H}_S}\left( SA^{\hat{H}}_{\mathrm{unip}}(S), \hat{f}[y_S] \right). \]
	By construction, $\overline{f}[y_S]$ and $\hat{f}[y_S]$ depend only on $\mathbf{H}/\sim$ and $y_S$, whereas $\tau(G^!_\epsilon) = \tau(\hat{H})$. The only remnant of $j$ now resides in $d_j(y_S)$.
	
	On the other hand, by \cite[VI.5.6 Théorème]{MW16-2} (nonstandard transfer of $SA_{\mathrm{unip}}$ for simply connected groups) and \cite[VII.4.6 Proposition]{MW16-2} (passage from $\SO(2m+1)$ to $\Spin(2m+1)$ up to Tamagawa numbers), for all $y_S$ we have
	\[ S^{\hat{H}_S}\left( SA^{\hat{H}}_{\mathrm{unip}}(S), \hat{f}[y_S] \right) = \frac{\tau(\hat{H})}{\tau(H)} \cdot S^{H_S}\left( SA^H_{\mathrm{unip}}(S), \overline{f}[y_S] \right). \]
	
	We also have $W^{\hat{H}} \simeq W^H$. Finally, summing over $\mathbf{H}$ in the $\sim$ class yields the desired equality.
\end{proof}

\section{Cohomological arguments}\label{sec:gdesc-coh}
Retain the conventions from \S\ref{sec:gdesc-endoscopy}. For all $v$, denote by $[y_v] \in \Hm^1_{\mathrm{ab}}(F_v, G_\eta)$ the class determined by $y_v \in \mathcal{Y}_v$ via $G_\eta(\overline{F_v}) \backslash \mathcal{Y}_v / G(F_v) \to \Hm^1(F_v, G_\eta)$ composed with the abelianization map $\mathrm{ab}^1$.

For all $v \notin S$, define
\[ \mathcal{Y}_v^{\mathrm{nr}} := \left\{ y_v \in \mathcal{Y}_v: \eta[y_v] \;\text{is $G(F_v)$-conjugate to}\; K_v \right\}. \]

\begin{lemma}\label{prop:Ynr}
	The set $\mathcal{Y}_v^{\mathrm{nr}}$ equals the $G_\eta(\overline{F_v}) \times G(F_v)$-orbit containing $1$. Moreover, $\mathcal{Y}_v^{\mathrm{nr}} \subset \mathcal{Y}_v^{\mathrm{rel}}(\mathbf{H})$ for any $\mathbf{H}$.
\end{lemma}
\begin{proof}
	Since $S \supset S(\mathcal{O})$ and $\eta \in K_v$, the first assertion is immediate from a result of Kottwitz \cite[7.1 Proposition]{Ko86}. That result also implies $G_{\eta[y_v]}$ is unramified, thus quasisplit, and the second assertion follows.
\end{proof}

\begin{definition}
	\index{YF@$\mathcal{Y}_F$, $\mathcal{Y}_F^{S\text{-nr}}$, $\mathcal{Y}_F[y_S]$}
	Let $\mathcal{Y}_F$ be the set of $y \in G(\overline{F})$ such that $\eta[y] := y\eta y^{-1} \in G(F)$, and $\Ad(y): G_{\eta[y], \overline{F}} \to G_{\eta, \overline{F}}$ is an inner twist; in fact, it is a pure inner twist since $Z_G(\eta) = G_\eta$ in this case. Denote the natural map $\mathcal{Y}_F \to \mathcal{Y}_S$ as $x \mapsto x_S$. Define
	\[ \mathcal{Y}_F^{S\text{-nr}} := \left\{ x \in \mathcal{Y}_F : v \notin S \implies x_v \in \mathcal{Y}_v^{\mathrm{nr}} \right\} . \]
	For all $y_S \in \mathcal{Y}_S$, define 
	\[ \mathcal{Y}_F[y_S] := \left\{ x \in \mathcal{Y}_F^{S\text{-nr}} : x_S \in G_\eta(\overline{F_S}) y_S G(F_S) \right\}. \]
	As before, we fix systems of representatives $\dot{\mathcal{Y}}_F^{S\text{-nr}} \subset \dot{\mathcal{Y}}_F$ of their $G_\eta(\overline{F}) \times G(F)$-orbits.
\end{definition}

\begin{lemma}\label{prop:YF-singleton}
	For all $y_S \in \mathcal{Y}_S$, the quotient set $G_\eta(\overline{F}) \backslash \mathcal{Y}_F[y_S] / G(F)$ is either empty or a singleton.
\end{lemma}
\begin{proof}
	Identify $G_\eta(\overline{F}) \backslash \mathcal{Y}_F / G(F)$ with $\Hm^1(F, G_\eta)$ as pointed sets. Lemma \ref{prop:Ynr} identifies $G_\eta(\overline{F}) \backslash \mathcal{Y}_F[y_S] / G(F)$ as a fiber of
	\[ \Hm^1(F, G_\eta) \rightiso \Hm_{\mathrm{ab}}^1(F, G_\eta) \dtimes{\Hm_{\mathrm{ab}}^1(\A_F, G_\eta)} \Hm^1(\A_F, G_\eta) \xrightarrow{\text{pr}_2} \Hm^1(\A_F, G_\eta) ; \]
	the first arrow being an isomorphism by \cite[Théorème 1.6.10]{Lab99}. Such a fiber, if non-empty, is identifiable with a fiber of
	\[ \Hm_{\mathrm{ab}}^1(F, G_\eta) \to \Hm_{\mathrm{ab}}^1(\A_F, G_\eta). \]
	
	Recall from \S\ref{sec:ss-classes} that $G_\eta$ is a product of symplectic and unitary factors up to restriction of scalars. It suffices to show $\Ker_{\mathrm{ab}}^1(F, I) = 0$ for any unitary factor $I$ in $G_\eta$; the other cases are trivially true.
	
	Without loss of generality, assume that $I$ is a unitary group relative to a quadratic extension $E|F$ of fields. By \cite[Exemple 3.1.2]{Li15}, the complex $I_{\mathrm{ab}}$ computing $\Hm^1_{\mathrm{ab}}(\cdot, I)$ is quasi-isomorphic to the norm-one torus relative to $E|F$. It remains to apply the local-global principle in such a situation, see \cite[Theorem 6.11]{PR94}.
\end{proof}

Consider any $\mathbf{H}$ and take $j \in \mathcal{J}(\mathbf{H})$. In Definition \ref{def:djyS}, the denominator of $d_j(y_S)$ was $\Delta_{\mathbf{H}_{\mathrm{std}}}[y_S]$, where $\mathbf{H}_{\mathrm{std}}$ is the standard representative in $\mathbf{H}/\sim$. We begin by inspecting the variant $\hat{d}_j(y_S)$ defined using the ``direct descendant'' $\mathbf{H}$ of $j$ via Proposition \ref{prop:desc-endoscopy-refined}. Results from \S\ref{sec:endoscopy-rigid} will be needed.

\begin{lemma}\label{prop:direct-d-trivial}
	For all $\mathbf{H}$, $j \in \mathcal{J}(\mathbf{H})$ and $y_S \in \mathcal{Y}_S^{\mathrm{rel}}(\mathbf{H})$, define $\hat{d}_j(y_S) = \prod_v \hat{d}_j(y_v) \in \CC^\times$, which is the ratio of the descent of $\Delta_j$ over the transfer factor $\Delta_{\mathbf{H}}[y_S]$. Then $\hat{d}_j(y_S) = 1$.
\end{lemma}
\begin{proof}
	For each place $v$ of $F$ and $y_v \in \mathcal{Y}_v^{\mathrm{rel}}(\mathbf{H})$, the ratio $\hat{d}_j(y_v)$ of $\Delta_j$ over $\Delta_{\mathbf{H}}[y_v]$ is determined in \cite[Théorème 7.10]{Li11}. Indeed, the endoscopic datum of $G_{\eta[y_v]}$ used there is readily seen to be $\mathbf{H}$ localized at $v$; the point is that in \textit{loc.\ cit.} one picked out the $\overline{s}$ in \eqref{eqn:sbars}.
	
	The ratio $\hat{d}_j(y_v)$ is a \emph{good constant} in the parlance of \textit{loc.\ cit.} Let us be more precise.
	\begin{itemize}
		\item Assume $y_v \in G_\eta(\overline{F_v}) G(F_v)$. Then $\hat{d}_j(y_v) = 1$ in the unramified setting, and this is indeed the case for all $v \notin S$ since the assumption on the $F$-pinning $\mathcal{E}$ of $G_\eta$ implies that $\Delta_{\mathbf{H}}[y_v]$ is correctly normalized in the unramified setting.
		\item The constant $\hat{d}_j(y_v)$ is independent of the pure inner twist $\Ad(y_v)$. To see this, note that in the arguments leading up to \cite[Théorème 7.10]{Li11}, the implied constants are made from
		\begin{itemize}
			\item constants appearing in the parameter of the \emph{stable} classes of $\eta$ and $\epsilon$ (see \S\ref{sec:ss-classes});
			\item Weil's constants $\gamma_{\psi_v}$, where $\psi_v$ is the $v$-component of $\psi: \A_F/F \to \mathbb{S}^1$;
			\item Hilbert symbols;
			\item values at the ``unitary part'' of the data of the virtual character $\omega_{\psi_v}^+ - \omega_{\psi_v}^-$ where $\omega_{\psi_v}^\pm$ are Weil's representations (see the first step in the proof of \cite[Proposition 7.20]{Li11}), and such a value is stably invariant in the metaplectic sense by \cite[Définition 5.6]{Li11}.
			\footnote{
				There is a typo in the proof of \cite[Proposition 7.21]{Li11}: in the case $\bullet = u$, one should replace $N_{L/F}(u+1)$ by its original version $\det_F(u+1 | W'_u)$. This has no consequences elsewhere.
			}
		\end{itemize}
	\end{itemize}
	
	Consequently, $\hat{d}_j(y_v) = \hat{d}_j(1_v)$ for all $y_v$. We claim that $\prod_v \hat{d}_j(1_v) = 1$. Indeed,
	\begin{compactitem}
		\item the $F$-pinning $\mathcal{E}$ of $G_\eta$ is assumed to be ``unramified'' outside $S$,
		\item Weil's constants, Hilbert symbols, and the virtual character $\omega_{\psi_v}^+ - \omega_{\psi_v}^-$ (see \cite[\S 4.7]{Li11}) all have product formula;
		\item the other data appearing in $\hat{d}_j(1_v)$ are all $F$-rational.
	\end{compactitem}
	
	All in all, $\hat{d}_j(y_S) = \prod_{v \in S} \hat{d}_j(1_v) = \prod_v \hat{d}_j(1_v) = 1$.
\end{proof}

\begin{remark}
	The phenomenon above is coherent with the cocycle property \cite[Théorème 4.2.12]{Li11} for $\Delta_j$, when $\eta$ and $\epsilon$ are in equisingular correspondence. The reason is that \cite[(5.1)]{Kal16} implies
	\[ \Delta_{\mathbf{H}}[y_v] = \lrangle{\overline{s}, [y_v]}_{\mathrm{Kott}} \Delta_{\mathbf{H}}[1_v] = \lrangle{s, [y_v]}_{\mathrm{Kott}} \Delta_{\mathbf{H}}[1_v] \]
	although $\mathbf{H}$ is now a principal endoscopic datum; indeed, the difference between $s$ and $\overline{s}$ in \eqref{eqn:sbars} does not affect $\lrangle{\cdot, \cdot}_{\mathrm{Kott}}$. On the other hand, $\Delta_{\mathbf{H}}[y_v]$ is replaced by the constant $1$ in \textit{loc.\ cit.}, thus independent of $y_v$. The product formula was also stated in \textit{loc.\ cit.} in the equisingular case.
\end{remark}

In view of Lemma \ref{prop:J-mult-one}, for $j \in \mathcal{J}(\mathbf{H})$ and $z \in Z_{\overline{G}^\vee}^{\Gamma_F}$, it makes sense to define
\[ zj := \text{the unique element in}\; \mathcal{J}(z\mathbf{H}). \]
This yields a free action of $Z_{G_\eta^\vee}^{\Gamma_F}$ on $\mathcal{J}(\mathbf{H}/\sim)$.

\begin{lemma}\label{prop:dj-Kott}
	Let $y_S \in \mathcal{Y}_S^{\mathrm{rel}}(\mathbf{H})$.
	\begin{enumerate}[(i)]
		\item For all $j \in \mathcal{J}(\mathbf{H}/\sim)$ and $z \in Z_{\overline{G}^\vee}^{\Gamma_F}$, we have
		\[ d_{zj}(y_S) = d_j(y_S) \prod_{v \in S} \lrangle{[y_v] , z}_{\mathrm{Kott}}. \]
		\item If $\mathcal{Y}_F[y_S] \neq \emptyset$, then $d_j(y_S) = 1$.
	\end{enumerate}
\end{lemma}
\begin{proof}
	Consider (i). In view of Lemma \ref{prop:direct-d-trivial}, it suffices to show that
	\[ \frac{d_{zj}(y_S)}{\hat{d}_{zj}(y_S)} = \frac{d_j(y_S)}{\hat{d}_j(y_S)} \prod_{v \in S} \lrangle{[y_v] , z}_{\mathrm{Kott}}. \]

	Denote by $d_j(y_v)$ the component of $d_j(y_S)$ at $v \in S$, and so forth. We are reduced to showing
	\begin{equation}\label{eqn:dj-Kott}
		\frac{d_{zj}(y_v)}{\hat{d}_{zj}(y_v)} = \frac{d_j(y_v)}{\hat{d}_j(y_v)} \cdot \lrangle{[y_v] , z}_{\mathrm{Kott}}, \quad v \in S.
	\end{equation}
	
	Take $w \in Z_{\overline{G}^\vee}^{\Gamma_F}$ such that $\mathbf{H} = w\mathbf{H}_{\mathrm{std}}$. Fix $v \in S$. Theorem \ref{prop:Delta-under-translation} implies
	\begin{equation}\label{eqn:j-hat}
		\frac{d_j(y_v)}{\hat{d}_j(y_v)} = \frac{\Delta_{\mathbf{H}}[y_v]}{\Delta_{\mathbf{H}_{\mathrm{std}}}[y_v]} =  \lrangle{[y_v], w }_{\mathrm{Kott}}.
	\end{equation}
	
	Consider $zj$, which maps to $z\mathbf{H} = zw\mathbf{H}_{\mathrm{std}}$; the same reasoning leads to
	\[ \frac{d_{zj}(y_v)}{\hat{d}_{zj}(y_v)} = \frac{\Delta_{z\mathbf{H}}[y_v]}{\Delta_{\mathbf{H}_{\mathrm{std}}}[y_v]} = \lrangle{[y_v], zw }_{\mathrm{Kott}}. \]
	The identity \eqref{eqn:dj-Kott} follows at once.
	
	Consider (ii). Both $[y_S] \in \Hm_{\mathrm{ab}}^1(\A_F, G_\eta)$ and $w$ are of global origin, and the local pairings satisfy product formula. Hence \eqref{eqn:j-hat} entails
	\[ \frac{d_j(y_S)}{\hat{d}_j(y_S)} = \prod_{v \in S} \lrangle{[y_v], w}_{\mathrm{Kott}} = \prod_v \lrangle{[y_v], w}_{\mathrm{Kott}} = 1. \]
	It remains to recall from Lemma \ref{prop:direct-d-trivial} that $\hat{d}_j(y_S) = 1$.
\end{proof}

\begin{corollary}\label{prop:dj-inversion}
	For all $\mathbf{H}/\sim$ and $y_S \in \mathcal{Y}_S^{\mathrm{rel}}(\mathbf{H} / \sim)$, we have
	\[ \sum_{j \in \mathcal{J}(\mathbf{H}/\sim)} d_j(y_S) = \begin{cases}
		0, & \text{if} \quad \mathcal{Y}_F[y_S] = \emptyset, \\
		\left| Z_{\overline{G}^\vee}^{\Gamma_F}\right|, & \text{if} \quad \mathcal{Y}_F[y_S] \neq \emptyset.
	\end{cases}\]
\end{corollary}
\begin{proof}
	Take any $j_0 \in \mathcal{J}(\mathbf{H} / \sim)$. The sum is indexed by $z \in Z_{\overline{G}^\vee}^{\Gamma_F}$ where $z$ corresponds to $d_{zj_0}(y_S)$. By Lemma \ref{prop:dj-Kott}, this gives
	\[ \sum_{j \in \mathcal{J}(\mathbf{H}/\sim)} d_j(y_S) = d_{j_0}(y_S) \sum_{z \in Z_{\overline{G}^\vee}^{\Gamma_F}} \prod_v \lrangle{[y_v], z}_{\mathrm{Kott}} . \]
	Moreover, that Lemma asserts $d_j(y_S) = 1$ for all $j \in \mathcal{J}(\mathbf{H}/\sim)$ when $\mathcal{Y}_F[y_S] \neq \emptyset$. It remains to address the case $\mathcal{Y}_F[y_S] = \emptyset$.
	
	Let us adopt the formalism of abelianized Galois cohomologies from \S\ref{sec:Galois-cohomology}. In our case, we will take $I = G_\eta$ and $H = G = \Sp(W)$, so that the complex $G_{\mathrm{ab}}$ is trivial. Thus
	\begin{gather*}
		\mathfrak{D}(G_\eta, G; \cdot) = \Hm^1(\cdot, G_\eta), \quad \mathfrak{E}(G_\eta, G; \cdot) = \Hm_{\mathrm{ab}}^1(\cdot, G_\eta).
	\end{gather*}

	On the other hand, to $y_S$ is attached a class in $\mathfrak{D}(G_\eta, G; \A_F)$, whose components outside $S$ are the base-points. The condition $\mathcal{Y}_F[y_S] \neq \emptyset$ is equivalent to that the class of $y_S$ arises from $\mathfrak{D}(G_\eta, G; F)$.
	
	Denote by $(\cdots)^D$ the Pontryagin dual. By \cite[Proposition 1.7.3]{Lab99} we also have
	\begin{align*}
		\mathfrak{E}(G_\eta, G; \A_F/F) & = \Hm_{\mathrm{ab}}^1(\A_F/F, G_\eta) \simeq \pi_0\left( Z_{\overline{G}^\vee}^{\Gamma_F}\right)^D, \\
		\mathfrak{R}(G_\eta, G; F)_1 & = \mathfrak{R}(G_\eta, G; F) = \Hm_{\mathrm{ab}}^1(\A_F/F, G_\eta)^D \\
		& \simeq \pi_0\left( Z_{\overline{G}^\vee}^{\Gamma_F} \right) .
 	\end{align*}
	canonically. Via the map $\mathfrak{E}(G_\eta, G; \A_F) \to \mathfrak{E}(G_\eta, G; \A_F/F)$ in \eqref{eqn:E-local-global}, the duality pairing between $\mathfrak{R}(G_\eta, G; F)_1$ and $\mathfrak{E}(G_\eta, G; \A_F/F)$ is induced by the products of $\lrangle{\cdot, \cdot}_{\mathrm{Kott}}$ over places $v$.

	Now Corollary \ref{prop:Langlands-obs} entails that
	\[ \mathcal{Y}_F[y_S] = \emptyset \implies \sum_{z \in Z_{\overline{G}^\vee}^{\Gamma_F}} \prod_v \lrangle{[y_v], z}_{\mathrm{Kott}} = 0 \]
	since $\prod_v \lrangle{[y_v], \cdot}$ is then a non-trivial character. This concludes the proof.
\end{proof}

The argument above should be compared with its local version in \S\ref{sec:Shalika-matching}.

Note that $\left| Z_{\overline{G}^\vee}^{\Gamma_F} \right| = \tau(G_\eta)$. This is a consequence of Kottwitz's formula for Tamagawa numbers \cite[Corollaire 1.7.4]{Lab99} and the fact $\Ker_{\mathrm{ab}}^1(F, G_\eta) = 0$, which has been observed in the proof of Lemma \ref{prop:YF-singleton}. As a by-product, $\tau(G_\eta) = \tau(G_{\eta[y]})$ for all $y \in \mathcal{Y}_F$.

Next, when a class $\mathbf{H}$ is given, we impose relevance conditions on $\mathcal{Y}_F^{S\text{-nr}}$ to define
\[ \mathcal{Y}_F^{S\text{-nr}}(\mathbf{H}) := \bigcup_{y_S \in \mathcal{Y}_S^{\mathrm{rel}}(\mathbf{H})} \mathcal{Y}_F[y_S] \subset \mathcal{Y}_F^{S\text{-nr}}, \]
and let $\dot{\mathcal{Y}}_F^{S\text{-nr}}(\mathbf{H}) \subset \dot{\mathcal{Y}}_F^{S\text{-nr}}$ be the system of representatives of $G_\eta(\overline{F}) \times G(F)$-orbits. Once again, they depend only on the class $\mathbf{H}/\sim$, and we shall accordingly write $\dot{\mathcal{Y}}_F^{S\text{-nr}}(\mathbf{H}/\sim)$, etc.
\index{YFSnrH@$\mathcal{Y}_F^{S\text{-nr}}(\mathbf{H}/\sim)$}

Recall \eqref{eqn:AEndo-f-transfer}. Given $\mathbf{H}/\sim$ in $E_{T^\vee, \elli}(G_\eta, S) / \sim$ and $y_S \in \mathcal{Y}_S$, let us denote
\[ \overline{f}^{\mathbf{H}/\sim}[y_S] := \begin{cases}
	\overline{f}[y_S], & \text{if} \quad y_S \in \mathcal{Y}_S^{\mathrm{rel}}(\mathbf{H} / \sim) \\
	0, & \text{otherwise}.
\end{cases}\]
Combining Lemmas \ref{prop:AEndoH-descent-0}, \ref{prop:YF-singleton} and Corollary \ref{prop:dj-inversion}, we obtain
\begin{multline}\label{eqn:AGEndo-rat}
	I^{\tilde{G}_S}\left( A^{\tilde{G}, \Endo}(S, \mathbf{H}/\sim)_{\elli}, f_S \right) \\
	= |W^G(\mu)|^{-1} |W^H| \cdot \frac{\tau(G)}{\tau(H)} \cdot \tau(G_{\eta[y]}) \sum_{y \in \dot{\mathcal{Y}}_F^{S\text{-nr}}} S^{H_S}\left( SA^H_{\mathrm{unip}}(S), \overline{f}^{\mathbf{H}/\sim}[y_S] \right);
\end{multline}
as usual, $\mathbf{H} = (H, \mathcal{H}, \overline{s})$.

\section{Global ascent}\label{sec:gdesc-ascent}
The starting point is \eqref{eqn:AGEndoHsim} and \eqref{eqn:AGEndo-rat}. Summing the latter formula over all classes $\mathbf{H}/\sim$ and noting $\tau(G)=1$, we obtain
\begin{equation}\label{eqn:AGEndo-X}\begin{gathered}
		I^{\tilde{G}_S}\left( A^{\tilde{G}, \Endo}(S, \mu, \omega)_{\elli}, f_S \right) = |W^G(\mu)|^{-1} \sum_{y \in \dot{\mathcal{Y}}_F^{S\text{-nr}}} X[y], \\
		X[y] := \sum_{\substack{\mathbf{H}/\sim \\ \in E_{T^\vee, \elli}(G_\eta, S) / \sim}} \frac{\tau(G_{\eta[y]})}{\tau(H)} \cdot |W^H| \cdot S^{H_S}\left( SA^H_{\mathrm{unip}}(S), \overline{f}^{\mathbf{H}/\sim}[y_S] \right).
\end{gathered}\end{equation}

\index{EellS@$\mathcal{E}_{\elli}(G_{\eta[y]}, S)$}
For each $y \in \dot{\mathcal{Y}}_F^{S\text{-nr}}$, denote by $\mathcal{E}_{\elli}(G_{\eta[y]}, S)$ the set of equivalence classes of endoscopic data of $G_{\eta[y]}$ that are unramified outside $S$ and elliptic. There is an evident map
\begin{align*}
	E_{T^\vee, \elli}(G_\eta, S) & \to \mathcal{E}_{\elli}(G_{\eta[y]}, S) \\
	\mathbf{H} & \mapsto [\mathbf{H}]
\end{align*}
that factors through $\sim$ classes.

\begin{lemma}\label{prop:E-surj}
	The map $\left(E_{T^\vee, \elli}(G_\eta, S) / \sim\right) \to \mathcal{E}_{\elli}(G_{\eta[y]}, S)$ induced by $\mathbf{H} \mapsto [\mathbf{H}]$ is surjective, and each fiber has cardinality equal to
	\[ |W^{G_\eta}| \cdot |W^H|^{-1} |\mathrm{Out}(\mathbf{H})|^{-1} . \]
\end{lemma}
\begin{proof}
	This assertion has nothing to do with coverings. An easy proof can be found in \cite[p.812]{MW16-2}.
\end{proof}

On the other hand, the coefficients in the stabilization of the trace formula for $G_{\eta[y]}$ read
\begin{equation}\label{eqn:iota-descent}
	\iota(G_{\eta[y]}, H) = \frac{\tau(G_{\eta[y]})}{\tau(H)} \cdot |\mathrm{Out}(\mathbf{H})|^{-1}
\end{equation}
for every $[\mathbf{H}] \in \mathcal{E}_{\elli}(G_{\eta[y]}, S)$.

Before proceeding to the next step, note that for all $\mathbf{H} \in E_{T^\vee, \elli}(G_\eta, S)$ and $y \in \dot{\mathcal{Y}}_F^{S\text{-nr}}(\mathbf{H})$, say $y \in \mathcal{Y}_F[y_S]$ for some $y_S \in \mathcal{Y}_S^{\mathrm{rel}}(\mathbf{H})$, the data defining $\Delta_{\mathbf{H}}[y_S]$ are all global, and everything is unramified outside $S$ relative to $K_v \cap G_\eta(F_v)$. This implies that $\Delta_{\mathbf{H}}[y_S]$ depends only on $[\mathbf{H}]$ and $y$; in fact, it is the global transfer factor mentioned in \S\ref{sec:endoscopy-rigid}. Therefore we can write $\trans_{\mathbf{H}, G_{\eta[y]}}$ unambiguously, for each $y$ and $[\mathbf{H}]$; it is the $\otimes$-product over $v \in S$ of local transfers.

On the other hand, the same $\otimes$-product over $v \in S$ yields $\trans_{\mathbf{H}, G_{\eta[y]}}(\cdot) = 0$ when $y_S \in \mathcal{Y}_S \smallsetminus \mathcal{Y}_S^{\mathrm{rel}}(\mathbf{H})$.

\begin{proof}[Proof of Theorem \ref{prop:matching-coeff-A}]
	Continue the foregoing discussions. Substituting Lemma \ref{prop:E-surj} and \eqref{eqn:iota-descent} into \eqref{eqn:AGEndo-X}, we see
	\begin{align*}
		X[y] & = |W^{G_\eta}| \sum_{[\mathbf{H}] \in \Endo_{\elli}(G_{\eta[y]}, S)} \iota(G_{\eta[y]}, H) I^{G_{\eta[y], S}}\left( \trans_{\mathbf{H}, G_{\eta[y]}} SA^H_{\mathrm{unip}}(S), f[y_S] \right) \\
		& = |W^G(\mu)| \cdot I^{G_{\eta[y], S}}\left( \sum_{[\mathbf{H}] \in \Endo_{\elli}(G_{\eta[y]}, S)} \iota(G_{\eta[y]}, H) \trans_{\mathbf{H}, G_{\eta[y]}} SA^H_{\mathrm{unip}}(S), \; f[y_S] \right),
	\end{align*}
	where $f[y_S] := \prod_{v \in S} f_v[y_v]$ as in \eqref{eqn:AEndo-f-transfer}, with $f_v[y_v] := \desc^{\tilde{G}}_{\tilde{\eta}[y_v]}(f_v)$. Also recall that
	\[ \trans_{\mathbf{H}, G_{\eta[y]}} SA^H_{\mathrm{unip}}(S) \neq 0 \implies y_S \in \mathcal{Y}_S^{\mathrm{rel}}(\mathbf{H}). \]

	Our normalization of transfer factors $\Delta_{\mathbf{H}}[y_S]$ is compatible with \cite[VI.3]{MW16-2}. Now apply Arthur's theorem on the matching of coefficients in \cite[p.239, Global Theorem 1]{Ar02} to arrive at
	\begin{equation*}
		X[y] = |W^G(\mu)| \cdot I^{G_{\eta[y], S}}\left( A^{G_{\eta[y]}}_{\mathrm{unip}}(S), f[y_S] \right).
	\end{equation*}
	
	We are going to ascend to $\tilde{G}$. Note that for $y \in \dot{\mathcal{Y}}_F^{S\text{-nr}}$, according to \cite[Lemme 4.1.7]{Li15}, the product over local stable conjugates $\tilde{\eta}_v[y_v] \in \tilde{G}_v$ coincides with $\widetilde{\eta[y]}_S$, the component in the decomposition $\eta[y] = \widetilde{\eta[y]}_S \eta[y]^S$ in \eqref{eqn:extraction-S}. For discussions on how $\eta[y]^S$, an ordinary conjugate of $\eta^S \in K^S$, is lifted to $\tilde{G}^S$, see the last part of Remark \ref{rem:condition-S}.

	For each $y \in \mathcal{Y}_F$, let $\mathcal{O}_y \subset G(F)$ denote the semisimple $G(F)$-conjugacy class containing $\eta[y]$. Taking account of Proposition \ref{prop:Aell-descent}, now \eqref{eqn:AGEndo-X} becomes
	\[ I^{\tilde{G}_S}\left( A^{\tilde{G}, \Endo}(S, \mu, \omega)_{\elli}, f_S \right) = \sum_{y \in \dot{\mathcal{Y}}_F^{S\text{-nr}}} I^{\tilde{G}_S}\left( A^{\tilde{G}}(S, \mathcal{O}_y)_{\elli}, f_S \right). \]
	
	We claim that the sum above equals the same sum over the larger set $\dot{\mathcal{Y}}_F$. Indeed, by \S\ref{sec:Aell}, $A^{\tilde{G}}(S, \mathcal{O}_y) \neq 0$ only when $\mathcal{O}_y$ cuts $K^S$, so
	\[ y \in \dot{\mathcal{Y}}_F \smallsetminus \dot{\mathcal{Y}}_F^{S\text{-nr}} \implies A^{\tilde{G}}(S, \mathcal{O}_y)_{\elli} = 0 \]
	by the definition of $\mathcal{Y}_v^{\mathrm{nr}}$.

	The enlarged sum is thus over all $G(F)$-conjugacy classes inside the stable class $\mathcal{O}$. By varying $f_S$, we obtain $A^{\tilde{G}, \Endo}(S, \mu, \omega)_{\elli} = A^{\tilde{G}}(S, \mathcal{O})_{\elli}$. This is exactly \eqref{eqn:AEndo-repara}.
\end{proof}

\chapter{Stabilization of local trace formula}\label{sec:LTF}
Let $F$ be a local field of characteristic zero. Consider a covering of metaplectic type $\rev: \tilde{G} \to G(F)$ where $G = \prod_{i \in I} \GL(n_i) \times \Sp(W)$, with chosen $(W, \lrangle{\cdot|\cdot})$, symplectic basis and non-trivial additive character $\psi$ of $F$.

The invariant local trace formula for $\tilde{G}$ established in \cite{Li12b} is cast into the form
\[ I_{\mathrm{geom}}(\overline{f_1}, f_2) = I_{\mathrm{disc}}(\overline{f_1}, f_2), \quad f_1, f_2 \in \orbI_{\asp}(\tilde{G}) \otimes \mes(G). \]

Philosophically, $\overline{f_1}$ should be viewed as an element of $\orbI_{\asp}(\tilde{G}^\dagger) \otimes \mes(G)$, where $\tilde{G}^\dagger$ is the \emph{antipodal covering} of $\tilde{G}$. This means that $\tilde{G}^\dagger = \tilde{G}$ as topological groups, but the embeddings of $\bmu_8$ differ by the automorphism $z \mapsto z^{-1}$.

In \S\ref{sec:antipodal}, we will show that when $G = \Sp(W)$, there is a canonical isomorphism between $\tilde{G}^\dagger$ and the covering group associated with $(W, -\lrangle{\cdot|\cdot}, \psi)$, as topological central extensions of $G(F)$ by $\bmu_8$. Using this property, we can relate the metaplectic transfer factors for $\tilde{G}$ and $\tilde{G}^\dagger$, then show that
\[ I_{\tilde{M}^\dagger}(\tilde{\gamma}, \overline{f}) = \overline{I_{\tilde{M}}(\tilde{\gamma}, f) }, \quad
I^{\Endo}_{\tilde{M}^\dagger}(\tilde{\gamma}, \overline{f}) = \overline{I^{\Endo}_{\tilde{M}}(\tilde{\gamma}, f)} \]
for all $f \in \orbI_{\asp}(\tilde{G}) \otimes \mes(G)$.

The structure of the invariant local trace formula is reviewed in \S\ref{sec:review-LTF}. In \S\S\ref{sec:LTF-geom-Endo}--\ref{sec:LTF-spec-Endo}, we define endoscopic counterparts for the distributions in the invariant local trace formula. The stabilization of local trace formula will be a consequence of the $\tilde{G}$-regular case of the local geometric Theorem \ref{prop:local-geometric}.

On the other hand, given $M \in \mathcal{L}(M_0)$, assume inductively that the $\tilde{G}$-regular case of Theorem \ref{prop:local-geometric} is known when $M$ (resp.\ $G$) is replaced by a larger (resp.\ proper) Levi subgroup, we will formulate a key geometric Hypothesis \ref{hyp:key-geometric}, which produces a smooth function $\epsilon(\mathbf{M}^!, \cdot)$ on $M^!_{G\text{-reg}}(F)$ for all $\mathbf{M}^! \in \Endo_{\elli}(\tilde{M})$, such that $\epsilon_{\tilde{M}} = 0$ if and only if $\epsilon(\mathbf{M}^!, \cdot) = 0$ for all $\mathbf{M}^!$. Then we will prove that
\begin{itemize}
	\item $\epsilon(\mathbf{M}^!, \cdot) + \overline{\epsilon(\mathbf{M}^!, \cdot)} = 0$ by using the local trace formula (Lemma \ref{prop:epsilon-vanishing-1});
	\item $\epsilon(\mathbf{M}^!, \cdot) = \overline{\epsilon(\mathbf{M}^!, \cdot)}$ by showing that ``the endoscopic transfer is isomorphic to its complex conjugate'', via the MVW-involution on metaplectic groups (Lemma \ref{prop:epsilon-vanishing-2}).
\end{itemize}
The combination of these facts implies $\epsilon_{\tilde{M}} = 0$.

This reduces both the stabilization of local trace formula and the local geometric Theorem \ref{prop:local-geometric} to the key geometric Hypothesis \ref{hyp:key-geometric}. That hypothesis will ultimately be settled by a local--global argument in \S\ref{sec:proof-key-geometric}, after proving a special case of global spectral stabilization; the notion of $\tilde{M}$-cuspidality introduced in Definition \ref{def:M-cuspidal} will continue to play a key role there.

\section{The antipodal covering}\label{sec:antipodal}
To begin with, consider a general covering
\[ 1 \to \bmu_m \to \tilde{G} \xrightarrow{\rev} G(F) \to 1 \]
as in \S\ref{sec:covering}, where $G$ is any connected reductive $F$-group and $m \in \Z_{\geq 1}$.

\begin{definition}
	\index{Gtilde-dagger@$\tilde{G}^\dagger$}
	The antipodal covering $\tilde{G}^\dagger$ of $\tilde{G}$ is defined as its push-out via the map $\bmu_m \rightiso \bmu_m$ given by $z \mapsto z^{-1}$. In other words,
	\[ \tilde{G}^\dagger = \dfrac{\bmu_m \times \tilde{G}}{(zw, \tilde{x}) \sim (z, w^{-1} \tilde{x})}, \quad z, w \in \bmu_m ,\; \tilde{x} \in \tilde{G}, \]
	and $\bmu_m \hookrightarrow \tilde{G}^\dagger$ (resp.\ $\tilde{G}^\dagger \to G(F)$) is given by $z \mapsto [z, 1]$ (resp.\ $[z, \tilde{x}] \mapsto \rev(\tilde{x})$).
	
	Alternatively, it is often convenient to identify $\tilde{G}^\dagger$ with $\tilde{G}$, with the same $\rev$, but the embedding $\bmu_m \hookrightarrow \tilde{G}$ is modified by $z \mapsto z^{-1}$.
\end{definition}

Recalling the definition in \S\ref{sec:orbital-integrals} of genuine functions (resp.\ anti-genuine representations) over $\tilde{G}$, we have
\[ C^\infty_-(\tilde{G}) = C^\infty_{\asp}(\tilde{G}^\dagger), \quad \Pi_{\asp}(\tilde{G}) = \Pi_-(\tilde{G}^\dagger). \]
Also, by writing $f \mapsto \overline{f}$ for the complex conjugation,
\[ C^\infty_-(\tilde{G}) = \left\{ \overline{f}: f \in C^\infty_{\asp}(\tilde{G}) \right\}. \]
Ditto for $\orbI_-(\tilde{G})$ and so forth. The contragredient (resp.\ complex conjugate) of a genuine representation $\pi$ of $\tilde{G}$ is a genuine representation of $\tilde{G}^\dagger$.

These operations are obviously compatible with inclusion of Levi subgroups $\tilde{M} \hookrightarrow \tilde{G}$. The maps $C^\infty_{c, \asp}(\tilde{G}) \twoheadrightarrow \orbI_{\asp}(\tilde{G})$ and $\orbI_{\asp}(\tilde{G}) \otimes \mes(G) \to \orbI_{\asp}(\tilde{M}) \otimes \mes(M)$ are readily seen to be ``defined over $\R$'', hence commutes with complex conjugation.

Hereafter, we return to the case where $\rev: \tilde{G} \to G(F)$ is a coverings of metaplectic type; in particular, we take $m=8$. What is special about coverings of metaplectic type is the relation between antipodals and endoscopy. Without loss of generality, let us assume $\tilde{G} = \Mp(W)$. Recall that the covering is determined by $\psi \circ \lrangle{\cdot|\cdot}: W \times W \to \mathbb{S}^1$.

\begin{definition-proposition}\label{def:G-minus}
	\index{Gtilde-minus@$\tilde{G}_-$}
	\index{Xi@$\Xi$}
	Let $\tilde{G}_-$ denote the covering of $\Sp(W)$ obtained from the datum $(W, \lrangle{\cdot|\cdot}, \psi^{-1})$, or equivalently from $(W, -\lrangle{\cdot|\cdot}, \psi)$. There is then a unique isomorphism $\Xi: \tilde{G}^\dagger \rightiso \tilde{G}_-$ of topological groups, such that the following diagram commutes:
	\[\begin{tikzcd}
		1 \arrow[r] & \bmu_8 \arrow[r] \arrow[d, "\identity"'] & \tilde{G}^\dagger \arrow[d, "\Xi"] \arrow[r] & G(F) \arrow[d, "\identity"] \arrow[r] & 1 \\
		1 \arrow[r] & \bmu_8 \arrow[r] & \tilde{G}_- \arrow[r] & G(F) \arrow[r] & 1 .
	\end{tikzcd}\]
	Equivalently, by identifying $\tilde{G}^\dagger$ with $\tilde{G}$ as groups, the diagram can also be rendered as
	\[\begin{tikzcd}
		1 \arrow[r] & \bmu_8 \arrow[r] \arrow[d, "{z \mapsto z^{-1}}"'] & \tilde{G} \arrow[d, "\Xi"] \arrow[r] & G(F) \arrow[d, "\identity"] \arrow[r] & 1 \\
		1 \arrow[r] & \bmu_8 \arrow[r] & \tilde{G}_- \arrow[r] & G(F) \arrow[r] & 1 .
	\end{tikzcd}\]
\end{definition-proposition}
\begin{proof}
	Uniqueness is straightforward: any two such isomorphisms of coverings differ by a homomorphism $G(F) \to \bmu_8$, which must be trivial.
	
	As for the existence, it suffices to deal with the second diagram. Fix a polarization $W = \ell \oplus \ell'$ where $\ell$ and $\ell'$ are Lagrangian subspaces. Recall the description of $\tilde{G}$ in terms of the Maslov cocycle: according to \cite{Th06} and \cite[Théorème 2.6]{Li11}, it arises from the section
	\[ G(F) \to \tilde{G}, \quad x \mapsto (x, M_\ell[x]) \]
	in the notation of \textit{loc.\ cit.}, where $M_\ell[x]$ is some endomorphism of the Schwartz space $\mathcal{S}(\ell')$, and $\bmu_8$ acts by dilation on the second component. Multiplication is described by
	\[ (x, M_\ell[x]) (y, M_\ell[y]) = \gamma_\psi(\underbracket{\tau(\ell, y\ell, xy\ell)}_{\text{Maslov index}}) \left(xy, M_\ell[xy] \right). \]
	
	Recall from \textit{loc.\ cit.} that the Maslov index here is understood as a class in the Witt group of quadratic $F$-vector spaces. Since
	\[ \gamma_{\psi^{-1}}(q) = \gamma_\psi(-q) = \gamma_\psi(q)^{-1}, \quad (V, q): \text{any quadratic $F$-vector space}, \]
	we see that replacing $\psi$ (resp.\ $\lrangle{\cdot|\cdot}$) by $\psi^{-1}$ (resp.\ $-\lrangle{\cdot|\cdot}$) has the effect of flipping the Maslov cocycle, which yields $\tilde{G}^\dagger$. This concludes the proof.
\end{proof}

Note that $\tilde{G}$ and $\tilde{G}_-$ have the same endoscopic data, and the correspondence of semisimple classes is also the same. What differs is the transfer factor.

\begin{lemma}\label{prop:Delta-minus}
	For every $f \in \orbI_{\asp}(\tilde{G}) \otimes \mes(G)$ and $\mathbf{G}^!$, we have $\overline{f}^{G^!} = \overline{f^{G^!}}$ where $f^{G^!} := \Trans_{\mathbf{G}^!, \tilde{G}}(f)$ and $\overline{f}^{G^!} := \Trans_{\mathbf{G}^!, \tilde{G}^\dagger}(\overline{f})$.
\end{lemma}
\begin{proof}
	The transfer is characterized by
	\[ f^{G^!}(\delta) = \sum_{\gamma \in \Gamma_{\text{reg}}(G)} \Delta_{\mathbf{G}^!, \tilde{G}}(\delta, \tilde{\gamma}) f(\tilde{\gamma}) \]
	for all $\delta \in \Sigma_{G\text{-reg}}(G^!)$. Using $\Xi: \tilde{G} \rightiso \tilde{G}_-$ and taking complex conjugates, the assertion reduces to the claim that
	\begin{equation}\label{eqn:Delta-minus}
		\Delta_{\mathbf{G}^!, \tilde{G}_-}(\delta, \Xi(\tilde{\gamma})) = \overline{\Delta_{\mathbf{G}^!, \tilde{G}}(\delta, \tilde{\gamma})}.
	\end{equation}
	
	Let us write $\tilde{\gamma} = z(\gamma, M_\ell[\gamma])$ as in the proof of Definition--Proposition \ref{def:G-minus}, where $\ell \subset W$ is Lagrangian and $z \in \bmu_8$. By construction, $\Xi(\tilde{\gamma}) = z^{-1} (\gamma, M_\ell[\gamma])$, but everything is now taken with respect to $\psi$ and $(W, -\lrangle{\cdot|\cdot})$. To prove \eqref{eqn:Delta-minus}, we may assume $z=1$.
	
	By \cite[\S 5.3]{Li11}, the ingredients in $\Delta$ not affected by flipping $\lrangle{\cdot|\cdot}$ are all real; the affected one is just the character of Weil representation $\omega_\psi^+ \oplus \omega_\psi^-$, or more precisely its value at regular semisimple elements divided by $|\cdots|$. This value lies in $\bmu_8$. In Maktouf's character formula \cite[Corollaire 4.4]{Li11}, this is given through Weil's constant $\gamma_\psi(\cdot)$ attached to a Maslov index $\tau(\Gamma_\gamma, \Gamma_1, \ell \oplus \ell)$, which gets multiplied by $-1$ if $\lrangle{\cdot|\cdot}$ is replaced by $-\lrangle{\cdot|\cdot}$. As $\gamma_\psi(-q) = \gamma_\psi(q)^{-1}$ for any quadratic $F$-vector space $(V,q)$, we obtain \eqref{eqn:Delta-minus}.
\end{proof}

\begin{proposition}\label{prop:IM-cplx-conj}
	For all $M \in \mathcal{L}(M_0)$, $\tilde{\gamma} \in \tilde{M}$ and $f \in \orbI_{\asp}(\tilde{G}) \otimes \mes(G)$. Fix Haar measures needed to define weighted orbital integrals, and identify $\tilde{\gamma}$ with an element of $\tilde{M}^\dagger$. Then
	\begin{equation*}
		I_{\tilde{M}^\dagger}(\tilde{\gamma}, \overline{f}) = \overline{I_{\tilde{M}}(\tilde{\gamma}, f) }, \quad
		I^{\Endo}_{\tilde{M}^\dagger}(\tilde{\gamma}, \overline{f}) = \overline{I^{\Endo}_{\tilde{M}}(\tilde{\gamma}, f) }.
	\end{equation*}
	In fact, the first equation holds for coverings in general.
\end{proposition}
\begin{proof}
	Consider the first equality. It is clear that $J_{\tilde{M}^\dagger}(\tilde{\gamma}, \overline{f}) = \overline{J_{\tilde{M}}(\tilde{\gamma}, f)}$, as the weighted orbital integrals are defined over $\R$. To show the equality, it suffices to show that $\phi_{\tilde{L}^\dagger}(\overline{f}) = \overline{\phi_{\tilde{L}}(f)}$ for all Levi subgroup $L$, or equivalently,
	\[ J_{\tilde{L}^\dagger}(\overline{\pi}, \overline{f}) = \overline{J_{\tilde{L}}(\pi, f)} \]
	for all $\pi \in \Pi_{\mathrm{temp}, -}(\tilde{L})$, with complex conjugate $\overline{\pi} \simeq \check{\pi}$.
	
	Recall the definition of canonically normalized tempered weighted character in \cite[\S 5.7]{Li12b}: fix $P \in \mathcal{P}(L)$,
	\begin{align*}
		J_{\tilde{L}}(\pi, f) & := \Tr\left( \mathcal{M}_{\tilde{L}}(\pi, \tilde{P}) I_{\tilde{P}}(\pi, f) \right), \\
		\mathcal{M}_{\tilde{L}}(\pi, \tilde{P}) & := \lim_{\Lambda \to 0} \sum_{Q \in \mathcal{P}(L)} \mathcal{M}_{\tilde{Q}}(\Lambda, \pi, \tilde{P}) \theta_Q(\Lambda)^{-1}, \quad \Lambda \in \mathfrak{a}^*_{L, \CC}, \\
		\mathcal{M}_{\tilde{Q}}(\Lambda, \pi, \tilde{P}) & := \mu_{\tilde{Q}|\tilde{P}}(\pi)^{-1} \mu_{\tilde{Q}|\tilde{P}}(\pi_{\Lambda/2}) J_{\tilde{Q}|\tilde{P}}(\pi)^{-1} J_{\tilde{Q}|\tilde{P}}(\pi_\Lambda).
	\end{align*}
	All these objects are well-defined upon perturbing $\pi$ by $i\mathfrak{a}^*_{L, F}$.
	\begin{itemize}
		\item The standard intertwining operators are defined over $\R$ in the sense that $J_{\tilde{Q}^\dagger|\tilde{P}^\dagger}(\overline{\pi_\Lambda}) = \overline{J_{\tilde{Q}|\tilde{P}}(\pi_\Lambda)}$ for all $\Pi_-(\tilde{L})$, which can be checked by looking at the integrals in the range of convergence, as both sides are anti-holomorphic in $\Lambda$.
		\item So are the Harish-Chandra $\mu$-functions: see \cite[Proposition 2.4.3]{Li12b}.
		\item Finally, $\theta_Q(\Lambda) := \mathrm{vol}(\mathfrak{a}^G_L / \Z \Delta_Q^\vee)^{-1} \prod_{\alpha \in \Delta_Q} \lrangle{\Lambda, \check{\alpha}}$ is defined over $\R$.
	\end{itemize}
	The required identity for $\phi_{\tilde{L}^\dagger}$ follows at once. These arguments apply to coverings in general.
	
	There are also stable counterparts of the first equality for $S^{G^!}_{M^!}(\delta, f^!)$, where $G^!$ (resp.\ $M^!$) is an elliptic endoscopic group of $\tilde{G}$ (resp.\ $\tilde{M}$). This is implicit in \cite[\S 9]{Ar99}, and follows by a standard inductive argument.
	
	To deduce the second equality, it suffices to combine the previous observation, Lemma \ref{prop:Delta-minus} and the first equality.
\end{proof}

The coverings $\tilde{G}$ and $\tilde{G}_-$ or $\tilde{G}^\dagger$ are related by an explicitly constructed isomorphism. Denote by $\GSp(W)$ the group of symplectic similitudes on $W$, and by $\nu: \GSp(W) \to \Gm$ the similitude character.

\begin{proposition}[{\cite[p.36]{MVW87}}]\label{prop:MVW}
	Take any $g \in \GSp(W)$ with $\nu(g) = -1$. Then there exists an isomorphism $\Theta_g: \tilde{G} \rightiso \tilde{G}_-$ of topological groups making the diagram below commutative:
	\[\begin{tikzcd}
		1 \arrow[r] & \bmu_8 \arrow[r] \arrow[d, "\identity"'] & \tilde{G} \arrow[d, "{\Theta_g}"] \arrow[r] & G(F) \arrow[d, "{\Ad(g)}"] \arrow[r] & 1 \\
		1 \arrow[r] & \bmu_8 \arrow[r] & \tilde{G}_- \arrow[r] & G(F) \arrow[r] & 1 .
	\end{tikzcd}\]
\end{proposition}
\begin{proof}
	In the cited reference, this is stated for the coverings obtained by pushing-out along $\bmu_8 \hookrightarrow \CC^\times$. Choose any Lagrangian subspace $\ell$ of $W$ with respect to $\lrangle{\cdot|\cdot}$, or equivalently for $-\lrangle{\cdot|\cdot}$. Using the construction of $\tilde{G}$ (or $\tilde{G}_-$) via Schrödinger models (see the proof of Definition--Proposition \ref{def:G-minus}), the map in \textit{loc.\ cit.} is simply $(x, M_\ell[x]) \mapsto (\Ad(g)x, M_\ell[x])$. Thus it restricts to $\tilde{G} \rightiso \tilde{G}_-$.
\end{proof}

The isomorphism above depends on the choice of $g$, but only up to conjugacy: $\Theta_{h'gh} = \Ad(h') \Theta_g \Ad(h)$ for all $h, h' \in G(F)$. It is customary to fix $g$ and use the shorthand $\Theta = \Theta_g$. The following result will also be needed.

\begin{lemma}\label{prop:MVW-adapted}
	Given $M \in \mathcal{L}(M_0)$, written as $M = \prod_{j \in J} \GL(m_j) \times \Sp(W^\flat)$, one can choose $g \in \GSp(W)$ with $\nu(g) = -1$ such that $\Ad(g)(M) = M$, and that $\Ad(g)$ acts as identity (resp.\ as in Proposition \ref{prop:MVW}) on each $\GL(m_j)$ (resp.\ on $\Sp(W^\flat)$).
	
	Such $\Theta_g$ is said to be \emph{adapted to $\tilde{M}$}.
\end{lemma}
\begin{proof}
	Describe $M$ by a decomposition $W = W^\flat \oplus \bigoplus_{j=1}^r (\ell_j \oplus \ell'_j)$ where $W^\flat$ and $\ell_j \oplus \ell'_j$ are symplectic vector subspaces of $W$, and each $\ell_j$ is totally isotropic of dimension $m_j$. Then $M \simeq \prod_{j \in J} \GL(\ell_j) \times \Sp(W^\flat)$. Take $g = \left((g_j)_{j \in J}, g^\flat\right) \in \GL(W)$ such that $g^\flat \in \GSp(W^\flat)$ satisfies $\nu(g^\flat) = -1$, and
	\[ g_j = \begin{pmatrix} \identity_{\ell_j} & \\ & -\identity_{\ell_j} \end{pmatrix} \in \GL(\ell_j \oplus \ell'_j) \]
	for all $j \in J$. Then $\nu(g) = -1$ with respect to $(W, \lrangle{\cdot|\cdot})$, and $g$ has the desired behavior.
\end{proof}

By choosing $g$ in this way, $\Theta_g$ can be restricted to $\tilde{M}$. The proof above actually implies that $g$ can be chosen so that $\Theta_g^2 = \identity$.

\begin{remark}\label{rem:MVW}
	In view of the previous result, one can define the automorphisms $\Theta$ for any group of metaplectic type, by requiring them to be inner on each $\GL$-factor. Nonetheless, what matters in \S\ref{sec:key-geometric} is just the $W^G(M)$-coset of $\Theta|_{\tilde{M}}$ modulo inner automorphisms. In view of the more general scenario, in which $\Ad(g)$ is to be replaced by a Chevalley involution of $G$, it is probably better to require that $\Theta$ is the transpose-inverse on each $\GL$-factor. In the case of $G = \Sp(W)$ and its Levi subgroup $M$, transpose-inverses are realizable by $W^G(M)$, so this does not affect the upcoming applications.
\end{remark}

\section{Local trace formulas}\label{sec:review-LTF}

We are going to state the invariant local trace formula for $\tilde{G}$ in the form of \cite[\S 5.8]{Li12b}. Subsequently, we will review the stable version for the endoscopic groups of $\tilde{G}$.

\subsection{Invariant local trace formula}\label{sec:invariant-local}
When talking about normalized intertwining operators and the resulting weighted characters, etc., the normalizing factors for genuine representations of $\tilde{G}$ and $\tilde{G}^\dagger$ are taken to be \emph{complementary}; see \cite[p.85]{Ar93} and \cite[p.849]{Li12b}. Specifically, we impose
\[ r_{\tilde{Q}|\tilde{P}}(\pi^\vee) = r_{\tilde{P}|\tilde{Q}}(\pi), \quad \pi \in \Pi_{\mathrm{unit}, -}(\tilde{M}), \quad P, Q \in \mathcal{P}(M), \]
where $\pi^\vee$ is viewed as an element of $\Pi_{\mathrm{unit}, -}(\tilde{M}^\dagger)$. With this choice, the $R$-groups and their central extensions in \S\ref{sec:spectral-distributions} satisfy
\begin{equation}\label{eqn:R-group-contragredient}
	\tilde{R}_{\pi^\vee} \simeq \tilde{R}_\pi, \quad R_{\pi^\vee} \simeq R_\pi
\end{equation}
canonically; see \cite[p.92, p.95]{Ar93}.

Such normalizing factors for $\tilde{G}^\dagger$ were called \emph{weak normalizing factors} in \S 3.1 of \textit{loc.\ cit.}, as the condition (R7) is not assumed.

For each $M \in \mathcal{L}(M_0)$, we define the spaces $T_-(\tilde{M})$, $T_{\elli, -}(\tilde{M})$, etc.\ and the distributions $\Theta_\tau$ using the normalizing factors; see \S\ref{sec:spectral-distributions}. The anti-genuine versions $T_{\asp}(\tilde{M})$, etc.\ are defined using the complementary factors, and will be identified with $T_-(\tilde{M}^\dagger)$, etc.

\begin{definition}
	\index{tau-vee@$\tau^\vee$}
	The map
	\begin{align*}
		\tilde{T}_-(\tilde{G}) & \to \tilde{T}_{\asp}(\tilde{G}) = \tilde{T}_-(\tilde{G}^\dagger) \\
		(M, \pi, r) & \mapsto (M, \pi^\vee, r),
	\end{align*}
	is well-defined in view of \eqref{eqn:R-group-contragredient}. It factors as $T_-(\tilde{G}) \to \tilde{T}_-(\tilde{G}^\dagger)$ as well as their quotients by $\mathbb{S}^1$; the resulting maps will all be denoted as $\tau \mapsto \tau^\vee$.
\end{definition}

In this manner, one gets $\Theta_{\tau^\vee} \in D_{\mathrm{spec}, -}(\tilde{G}^\dagger) \otimes \mes(G)^\vee = D_{\mathrm{spec}, \asp}(\tilde{G}) \otimes \mes(G)^\vee$ from each $\tau \in T_-(\tilde{G})$.

Denote by $\mathscr{C}(\tilde{G})$ the Schwartz--Harish--Chandra space of $\tilde{G}$, and denote by $\mathscr{C}_-(\tilde{G})$ (resp.\ $\mathscr{C}_{\asp}(\tilde{G})$) its genuine (resp.\ anti-genuine) part. Their quotient spaces $I\mathscr{C}_-$ and $I\mathscr{C}_{\asp}$ are defined by the familiar recipe. The invariant local trace formula for $\tilde{G}$ is the following identity.\footnote{In \cite[p.851]{Li12b}, one takes test functions from $\mathscr{C}_{\asp}(\tilde{G} \times \tilde{G}) \subset \mathscr{C}(\tilde{G} \times \tilde{G})$, and $f_1 \otimes f_2$ belongs to this space. An obvious typo appeared there: the subscript $\asp$ should mean the space of functions that are anti-genuine under $z \mapsto (z^{-1}, 1)$ and $z \mapsto (1, z)$.}

\begin{theorem}[{\cite[Théorème 5.8.6]{Li12b}}]
	For all $h_1 \in I\mathscr{C}_-(\tilde{G}) \otimes \mes(G)$ and $h_2 \in I\mathscr{C}_{\asp}(\tilde{G}) \otimes \mes(G)$, we have
	\[ I_{\mathrm{geom}}(h_1, h_2) = I_{\mathrm{disc}}(h_1, h_2). \]
\end{theorem}

Both sides are separately continuous in $(h_1, h_2)$ with respect to the topology on $\mathscr{C}(\tilde{G})$; this is proved in \cite[Proposition 5.7.6]{Li12b}. The equality will also be written as $I^{\tilde{G}}_{\mathrm{geom}} = I^{\tilde{G}}_{\mathrm{disc}}$ in case of possible confusions. The distributions are explicated below.
\begin{description}
	\index{Igeom-local@$I_{\mathrm{geom}}(\cdot, \cdot)$, $I_{\tilde{M}, \mathrm{geom}}(\cdot, \cdot)$}
	\index{Idisc-local@$I_{\mathrm{disc}}(\cdot, \cdot)$}
	\item[Geometric side] For each $M \in \mathcal{L}(M_0)$, equip $A_M(F)$ with the Haar measure determined by the invariant positive-definite quadratic form on $\mathfrak{a}_M \subset \mathfrak{a}_0$. For each elliptic maximal torus $T \subset M$, use the Haar measure such that $\mes(T(F)/A_M(F)) = 1$. We set
	\begin{align*}
		I_{\tilde{M}, \mathrm{geom}}(h_1, h_2) & =
		I^{\tilde{G}}_{\tilde{M}, \mathrm{geom}}(h_1, h_2) \\
		& := \int_{\Gamma_{G\text{-reg, ell}}(M)} I_{\tilde{M}}\left(\gamma, h_1, h_2 \right) \dd\gamma
	\end{align*}
	where $\Gamma_{G\text{-reg, ell}}(M)$ is equipped with the measure à la Weyl integration formula, and $I_{\tilde{M}}\left(\gamma, h_1, h_2 \right) = I_{\tilde{M}}\left(\gamma, h_1 \otimes h_2 \right)$ is defined in \cite[p.852]{Li12b} using the canonically normalized invariant weighted characters. Here we follow the paradigm of the semi-local distributions $I_{\tilde{M}_V}(\tilde{\gamma}, \cdots)$ in \S\ref{sec:semi-local-geom-distributions}, by imaging that $V = \{v_1, v_2\}$ with the same local fields $F$. The new feature here is that it depends only on $\gamma$ and not on $\tilde{\gamma} \in \rev^{-1}(\gamma)$.
	
	Now set
	\begin{equation*}
		I_{\mathrm{geom}}(h_1, h_2) := \sum_{M \in \mathcal{L}(M_0)} \frac{|W^M_0|}{|W^G_0|} (-1)^{\dim \mathfrak{a}^G_M} I_{\tilde{M}, \mathrm{geom}}(h_1, h_2) .
	\end{equation*}
	
	Compared with the general setting in p.823 of \textit{loc.\ cit.}, here we have $A_M(F)^\dagger := A_M(F)$ and $\iota(\gamma) = 1$.
	
	\item[Spectral side] We have
	\begin{equation*}
		I_{\mathrm{disc}}(h_1, h_2) := \int_{\mathbb{S}^1 \backslash T_{\mathrm{disc}, -}(\tilde{G})} i(\tau) \Theta_{\tau^\vee}(h_1) \Theta_\tau(h_2) \dd\tau ,
	\end{equation*}
	where we choose representatives $\tau \in T_{\mathrm{disc}, -}(\tilde{G})$ for elements of $\mathbb{S}^1 \backslash T_{\mathrm{disc}, -}(\tilde{G})$. We refer to \S\ref{sec:disc-parameter} for the definition of $T_{\mathrm{disc}, -}(\tilde{G}) \subset T_-(\tilde{G})$ and the Radon measure on $\mathbb{S}^1 \backslash T_{\mathrm{disc}, -}(\tilde{G})$.
	
	\index{itau@$i(\tau)$}
	For $\tau = (L, \pi, r)$ and $P \in \mathcal{P}(L)$, the function $i(\tau)$ appearing above is given by
	\[ i(\tau) = i^{\tilde{G}}(\tau) := \sum_{t \in W_\pi^0 r \cap W^G(L)_{\mathrm{reg}}} \epsilon_\sigma(t) \left|\det\left( 1 - t \middle| \mathfrak{a}^G_L \right)\right|^{-1} , \]
	where $\epsilon_\sigma(t) = (-1)^{\ell(w_0)}$ if $t = w_0 r$ with $w_0 \in W_\pi^0$ and $r \in R_\pi$; see \S\ref{sec:disc-parameter}. The sign $\ell(w_0)$ is defined since by the theory of $R$-groups, $W_\pi^0$ is known to be the Weyl group of the root system generated $\left\{ \beta \in \Sigma_P^{\mathrm{red}}: \mu_\beta(\pi)=0 \right\}$.

	It is convenient to set $i(\tau) = 0$ if $\tau \notin \mathbb{S}^1 \backslash T_{\mathrm{disc}, -}(\tilde{G})$. Note that the coefficients $i(\tau)$ are invariant under $\mathbb{S}^1$ and $i\mathfrak{a}^*_{G, F}$-twists.
\end{description}

\begin{remark}
	This is the formula considered in \cite{Ar03-3}. It differs from the ``symmetric'' version utilized in \S\ref{sec:epsilon-KM-finite}.
\end{remark}

What we will use is a slightly different formulation. The distributions above can also be written as
\begin{equation}\label{eqn:invariant-local-f}
	I_{\mathrm{geom}}(\overline{f_1}, f_2), \quad I_{\mathrm{disc}}(\overline{f_1}, f_2), \quad f_1, f_2 \in I\mathscr{C}_{\asp}(\tilde{G}).
\end{equation}

Let us relate them with the genuine distributions on $\tilde{G}$ already defined. Consider the geometric side. Take any $\tilde{\gamma} \in \rev^{-1}(\gamma)$ for a representative of $\gamma \in \Gamma_{G\text{-reg, ell}}(M)$ in $M(F)$. There is a splitting formula
\[ I_{\tilde{M}}(\gamma, h_1, h_2) = \sum_{L_1, L_2 \in \mathcal{L}(M)} d^G_M(L_1, L_2) I^{\tilde{L}_1}_{\tilde{M}}(\tilde{\gamma}, h_{1, \tilde{L}_1}) I^{\tilde{L}_2}_{\tilde{M}}(\tilde{\gamma}, h_{2, \tilde{L}_2}); \]
see the discussions in \S\ref{sec:semi-local-geom-distributions} or the proof of \cite[Corollaire 5.8.7]{Li12b}. Note that the right hand side does not depend on the choice of $\tilde{\gamma}$. Since $d^G_M(L_1, L_2) \in \R_{\geq 0}$, whilst $\mathscr{C}(\tilde{G}) \to I\mathscr{C}(\tilde{G})$ and $h_i \mapsto h_{i, \tilde{L}_i}$ are both ``defined over $\R$'', the first equality of Proposition \ref{prop:IM-cplx-conj} (extended to Schwartz--Harish-Chandra spaces by continuity) yields
\[ I_{\tilde{M}}(\gamma, \overline{f_1}, f_2) = \sum_{L_1, L_2 \in \mathcal{L}(M)} d^G_M(L_1, L_2) \overline{I^{\tilde{L}_1}_{\tilde{M}}(\tilde{\gamma}, f_{1, \tilde{L}_1})} I^{\tilde{L}_2}_{\tilde{M}}(\tilde{\gamma}, f_{2, \tilde{L}_2}) \]
for all $f_1, f_2 \in I\mathscr{C}_{\asp}(\tilde{G})$.

The situation on the spectral side is similar. Consider $\tau = (M, \pi, r) \in \tilde{T}_-(\tilde{G})$ and fix $P \in \mathcal{P}(M)$. For every $\tilde{x} \in \tilde{G}$, the endomorphism $I_{\tilde{P}}(\pi^\vee, \tilde{x})$ equals the transpose-inverse of $I_{\tilde{P}}(\pi, \tilde{x})$. On the other hand, the formula (2.1) in \cite[p.86]{Ar93} says that the normalized intertwining operator $R_{\tilde{P}}(r, \pi^\vee)$ equals the transpose-inverse of $R_{\tilde{P}}(r, \pi)$. All these operators are unitary. Hence the definition of $\Theta_{\tau^\vee}$ and $\Theta_\tau$ entails that
\[ I_{\mathrm{disc}}(\overline{f_1}, f_2) = \int_{\mathbb{S}^1 \backslash T_{\mathrm{disc}, -}(\tilde{G})} i(\tau) \overline{\Theta_\tau(f_1)} \Theta_\tau(f_2) \dd\tau . \]
This is the form of the invariant local trace formula employed in \cite[\S 6]{Ar03-3} and \cite[X.3.3]{MW16-2}.

\subsection{Stable side}
Consider the stable side $G^! \supset M^!$, in the usual notations. In particular, $G^!$ is a direct product of factors of the form $\GL(a)$ or $\SO(2b+1)$. Again, we fix an invariant positive-definite quadratic form on $\mathfrak{a}_0^!$ in order to normalize Haar measures. In the upcoming applications, they will come from the prescribed invariant quadratic form on $\mathfrak{a}_0$ via endoscopy.

Denote by $S\mathscr{C}(G^!)$ the stable analogue of $I\mathscr{C}(G^!)$. The stable version of the invariant local trace formula is established in \cite[Corollary 6.4]{Ar03-3}, taking the form
\[ S^{G^!}_{\mathrm{geom}}(\overline{f^!_1}, f^!_2) = S^{G^!}_{\mathrm{disc}}(\overline{f^!_1}, f^!_2) \]
for all $f^!_1, f^!_2 \in S\mathscr{C}(G^!) \otimes \mes(G^!)$. See also \cite[X.3.2]{MW16-2}.

Both sides are separately continuous in $(f^!_1, f^!_2)$ with respect to the topologies on $\mathscr{C}(G^!)$. Note that we are stating only the analogue of \eqref{eqn:invariant-local-f}, i.e.\ the version involving complex conjugates. The precise construction is given inductively from the non-stable case settled in \cite[p.189]{Ar94}. Take the geometric side for example,
\[ S^{G^!}_{\mathrm{geom}}\left( \overline{f^!_1}, f^!_2 \right) = I^{G^!}_{\mathrm{geom}}(\overline{f^!_1}, f^!_2) - \sum_{\substack{\mathbf{G}^{!!} \in \Endo_{\elli}(G^!) \\ G^{!!} \neq G^!}} \iota(G^!, G^{!!}) S^{G^{!!}}_{\mathrm{geom}}\left( \overline{f_1^{!!}}, f_2^{!!} \right) \]
where $f_i^{!!}$ denotes the transfer of $f_i^!$ to the endoscopic group $G^{!!}$, for $i=1,2$. We refer to \cite[p.832]{Ar03-3} for details.

\begin{description}
	\index{Sgeom-local@$S^{G^{"!}}_{\mathrm{geom}}(\cdot, \cdot)$, $S^{G^{"!}}_{M^{"!}, \mathrm{geom}}(\cdot, \cdot)$}
	\index{Sdisc-local@$S^{G^{"!}}_{\mathrm{disc}}(\cdot, \cdot)$}
	\item[Geometric side] Following \cite[(10.11)]{Ar03-3} or \cite[p.1150]{MW16-2}, for each $M^! \in \mathcal{L}(M^!_0)$ we set
	\begin{align*}
		S_{M^!, \mathrm{geom}}(\overline{f^!_1}, f^!_2) & = S^{G^!}_{M^!, \mathrm{geom}}(\overline{f^!_1}, f^!_2) \\
		& := \int_{\Sigma_{G^!\text{-reg, ell}}(M^!)} \left| \mathfrak{E}(M^!_\delta, M^!; F) \right|^{-1} S^{G^!}_{M^!}(\delta, \overline{f^!_1}, f^!_2) \dd\delta ,
	\end{align*}
	where $\mathfrak{E}(M^!_\delta, M^!; F)$ is as in \S\ref{sec:Galois-cohomology}, and $\Sigma_{G^!\text{-reg, ell}}(M^!)$ is endowed with the measure à la stable Weyl integration formula. Again, $S^{G^!}_{M^!}(\delta, \cdots)$ follows the paradigm of semi-local versions in \S\ref{sec:semi-local-geom-distributions} with $V = \{v_1, v_2\}$, with the same local field $F$. In particular, there is a stable splitting formula
	\[ S^{G^!}_{M^!}(\delta, \overline{f^!_1}, f^!_2) = \sum_{L^!_1, L^!_2 \in \mathcal{L}(M^!)} e^{G^!}_{M^!}(L^!_1, L^!_2) \overline{S^{L^!_1}_{M^!}(\delta, f_{1, L^!_1})} S^{L^!_2}_{M^!}(\delta, f_{2, L^!_2}). \]
	
	The geometric side is then expressed as
	\[ S^{G^!}_{\mathrm{geom}}(\overline{f^!_1}, f^!_2) = \sum_{M^! \in \mathcal{L}(M^!_0)} \frac{|W^{M^!}_0|}{|W^{G^!}_0|} (-1)^{\dim \mathfrak{a}^{G^!}_{M^!}} S^{G^!}_{M^!, \mathrm{geom}}(\overline{f^!_1}, f^!_2) . \]
	\item[Spectral side] By \cite[Proposition 6.3 (b), Corollary 6.4 (b)]{Ar03-3}, there are uniquely determined coefficients $s^{G^!}(\phi)$ such that
	\[ S^{G^!}_{\mathrm{disc}}(\overline{f^!_1}, f^!_2) = \int_{\Phi_{\mathrm{disc}}(G^!)} s^{G^!}(\phi) \overline{S\Theta_\phi(f^!_1)} S\Theta_\phi(f^!_2) \dd\phi , \]
	where $\Phi_{\mathrm{disc}}(G^!)$ is a subspace of $\Phi_{\mathrm{bdd}}(G^!)$, discrete modulo $i\mathfrak{a}^*_{G^!, F}$, that supports $s^{G^!}$. The measure on $\Phi_{\mathrm{disc}}(G^!)$ is defined by the same recipe as for $T_{\mathrm{disc}, -}(\tilde{G})/\mathbb{S}^1$; see \eqref{eqn:Tdisc-measure}.
	\index{sGphi@$s^{G^{"!}}(\phi)$}
\end{description}

The coefficients $s^{G^!}(\phi)$ are invariant under $i\mathfrak{a}^*_{G^!, F}$-twists. This follows inductively from the corresponding property of $i^{G^!}$.

\begin{remark}\label{rem:simple-stable-LTF}
	When $f_1^!$ and $f_2^!$ are both cuspidal, the identity $S^{G^!}_{\mathrm{geom}}(\overline{f^!_1}, f^!_2) = S^{G^!}_{\mathrm{disc}}(\overline{f^!_1}, f^!_2)$ simplifies into
	\begin{multline*}
		\int_{\Sigma_{\elli}(G^!)} \left| \mathfrak{E}(M^!_\delta, M^!; F) \right|^{-1} \overline{f_1^!(\delta)} f_2^!(\delta) \dd\delta \\
		=  \int_{\Phi_{\mathrm{bdd}, 2}(G^!)} s^{G^!}(\phi) \overline{S\Theta_\phi(f^!_1)} S\Theta_\phi(f^!_2) \dd\phi .
	\end{multline*}
	Recall that $f_i^!(\delta) = S^{G^!}(\delta, f_i^!)$ by convention, and the integral over $\Phi_{\mathrm{bdd}, 2}(G^!)$ is discrete modulo $i\mathfrak{a}^*_{G, F}$. The identity can be seen by using the stable splitting formula for $S^{G^!}_{M^!}(\delta, \overline{f^!_1}, f^!_2)$ to show that it vanishes when $M^! \neq G^!$.
\end{remark}

\section{Endoscopic constructions: geometric side}\label{sec:LTF-geom-Endo}
We revert to the study of a covering of metaplectic type $\rev: \tilde{G} \to G(F)$.

\index{iotaGG}
Recall from Definition \ref{def:iotaGG-general} that for every $\mathbf{G}^! \in \Endo(\tilde{G})$, we have the constant
\[ \iota(\tilde{G}, G^!) := \left( Z_{(G^!)^\vee} : Z_{\tilde{G}^\vee}^\circ \right)^{-1}
	= \begin{cases}
		2^{- \#(\SO\text{-factors in}\; G^!)}, & \text{if}\; \mathbf{G}^! \in \Endo_{\elli}(\tilde{G}), \\
		0, & \text{otherwise}.
\end{cases}\]
The explicit formula for $\iota(\tilde{G}, G^!)$ is the same as the global case (Definition \ref{def:iota-const}).

For each $M \in \mathcal{L}(M_0)$ and $\mathbf{M}^! \in \Endo_{\elli}(\tilde{M})$, equip $\mathfrak{a}_{M^!}$ with the invariant quadratic form induced from the one on $\mathfrak{a}_M$, which is in turn restricted from $\mathfrak{a}_0$. The isomorphism $A_{M^!}(F) \rightiso A_M(F)$ preserves the corresponding Haar measures. This determines the Haar measure on $T^!(F)$ for each elliptic maximal torus $T^! \subset M^!$, such that $\mes(T^!(F)/A_{M^!}(F)) = 1$.

If one starts with $\mathbf{G}^! \in \Endo_{\elli}(\tilde{G})$ and take the induced invariant quadratic form on $\mathfrak{a}^!_0$, then restriction to $\mathfrak{a}_{M^!}$ for various $M^! \subset G^!$ gives the same Haar measures.

Recall that $G = \prod_{i \in I} \GL(n_i) \times \Sp(W)$. The definition of endoscopic datum is not affected by change of $\psi$ or dilation of $(W, \lrangle{\cdot|\cdot})$, hence $\Endo_{\elli}(\tilde{G}) = \Endo_{\elli}(\tilde{G}^\dagger)$, and so forth. On the other hand, the transfer factors are affected by such changes.

\begin{definition}
	For all $f_1, f_2 \in \orbI_{\asp}(\tilde{G}) \otimes \mes(G)$ and all $\mathbf{G}^! \in \Endo_{\elli}(\tilde{G})$, let
	\begin{itemize}
		\item $\overline{f_1}^{G^!}$ be the transfer of $\overline{f_1} \in \orbI_{\asp}(\tilde{G}^\dagger) \otimes \mes(G)$ via $\Trans_{\mathbf{G}^!, \tilde{G}^\dagger}$, using the datum $(W, \lrangle{\cdot|\cdot}, \psi^{-1})$ (or equivalently $(W, -\lrangle{\cdot|\cdot}, \psi)$) via $\Xi: \tilde{G}^\dagger \rightiso \tilde{G}^-$ (Definition--Proposition \ref{def:G-minus});
		\item $f_2^{G^!}$ be the transfer of $f_2$ via $\Trans_{\mathbf{G}^!, \tilde{G}}$.
	\end{itemize}
	
	We define
	\begin{align*}
		I^{\Endo}_{\mathrm{geom}}(\overline{f_1}, f_2) & = I^{\tilde{G}, \Endo}_{\mathrm{geom}}(\overline{f_1}, f_2) \\
		& := \sum_{\mathbf{G}^! \in \Endo_{\elli}(\tilde{G})} \iota(\tilde{G}, G^!) S^{G^!}_{\mathrm{geom}}\left( \overline{f_1}^{G^!}, f_2^{G^!} \right).
	\end{align*}
	\index{IEndogeom-local@$I^{\Endo}_{\mathrm{geom}}(\cdot, \cdot)$, $I^{\Endo}_{\tilde{M}, \mathrm{geom}}(\cdot, \cdot)$}
	
	For each $M \in \mathcal{L}(M_0)$ and $\mathbf{M}^! \in \Endo_{\elli}(\tilde{M})$, we define
	\begin{align*}
		I^{\Endo}_{\tilde{M}, \mathrm{geom}}(\mathbf{M}^!, \overline{f_1}, f_2) & = I^{\tilde{G}, \Endo}_{\tilde{M}, \mathrm{geom}}(\mathbf{M}^!, \overline{f_1}, f_2) \\
		& := \sum_{s \in \Endo_{\mathbf{M}^!}(\tilde{G})} i_{M^!}(\tilde{G}, G^![s]) S^{G^![s]}_{M^!, \mathrm{geom}}\left( \overline{f_1}^{G^![s]}, f_2^{G^![s]} \right).
	\end{align*}
	Note that the $z[s]$-twists do not intervene here.
\end{definition}

\begin{definition}
	Let
	\begin{align*}
		I^{\Endo}_{\tilde{M}}\left(\mathbf{M}^!, \delta, \overline{f_1}, f_2 \right) & := \sum_{s \in \Endo_{\mathbf{M}^!}(\tilde{G})} \iota_{M^!}(\tilde{G}, G^![s]) S^{G^![s]}_{M^!}\left(\delta[s], \overline{f_1}^{G^![s]}, f_2^{G^![s]}\right) \\
		I^{\Endo}_{\tilde{M}}\left(\gamma, \overline{f_1}, f_2 \right) & := \sum_{L_1, L_2 \in \mathcal{L}(M)} d^G_M(L_1, L_2) \overline{I^{\tilde{L}_1, \Endo}_{\tilde{M}}(\tilde{\gamma}, f_{1, \tilde{L}_1})} I^{\tilde{L}_2, \Endo}_{\tilde{M}}(\tilde{\gamma}, f_{2, \tilde{L}_2})
	\end{align*}
	where $\mathbf{M}^! \in \Endo_{\elli}(\tilde{M})$ and $\delta \in \Sigma_{G\text{-reg, ell}}(M^!)$ (in the first case), or $\gamma \in \Gamma_{\text{reg}}(M)$ and $\tilde{\gamma} \in \rev^{-1}(\gamma)$ is arbitrary (in the second case). The paradigm follows the semi-local constructions in \S\ref{sec:semi-local-geom-Endo}, say with $V = \{v_1, v_2\}$, and the details will not be repeated here.
\end{definition}

In particular, the proof of Proposition \ref{prop:semilocal-matching-coherence} (i) adapts to this context, and yields the splitting formula
\begin{equation}\label{eqn:IMEndo-double-splitting}
	I^{\Endo}_{\tilde{M}}\left(\mathbf{M}^!, \delta, \overline{f_1}, f_2 \right) := \sum_{L_1, L_2 \in \mathcal{L}(M)} d^G_M(L_1, L_2) \overline{I^{\tilde{L}_1, \Endo}_{\tilde{M}}(\mathbf{M}^!, \delta, f_{1, \tilde{L}_1})} I^{\tilde{L}_2, \Endo}_{\tilde{M}}(\mathbf{M}^!, \delta, f_{2, \tilde{L}_2}).
\end{equation}

We are ready to state the stabilization of the geometric side, as follows.

\begin{theorem}\label{prop:stabilization-geom-LTF}
	Let $f_1, f_2 \in \orbI_{\asp}(\tilde{G}) \otimes \mes(G)$.
	\begin{enumerate}[(i)]
		\item For all $M \in \mathcal{L}(M_0)$, we have
		\[ \sum_{\mathbf{M}^! \in \Endo_{\elli}(\tilde{M})} \iota(\tilde{M}, M^!) I^{\Endo}_{\tilde{M}, \mathrm{geom}}(\mathbf{M}^!, \overline{f_1}, f_2) = I_{\tilde{M}, \mathrm{geom}}(\overline{f_1}, f_2). \]
		\item We have
		\[ I^{\Endo}_{\mathrm{geom}}(\overline{f_1}, f_2) = I_{\mathrm{geom}}(\overline{f_1}, f_2). \]
	\end{enumerate}
\end{theorem}

A conditional proof based on the key geometric Hypothesis \ref{hyp:key-geometric} will be given in \S\ref{sec:key-geometric}.

\begin{remark}\label{rem:SHC-LTF}
	Using the trace Paley--Wiener theorems for $I\mathscr{C}_{\asp}(\tilde{G})$ and $S\mathscr{C}(G^!)$, together with the spectral transfer reviewed in \S\ref{sec:spectral-transfer}, one can show that $\Trans_{\mathbf{G}^!, \tilde{G}}$ and $\Trans_{\mathbf{G}^!, \tilde{G}^\dagger}$ extend by continuity to Schwartz--Harish-Chandra spaces. The stabilization above extends automatically to $f_1, f_2 \in I\mathscr{C}_{\asp}(\tilde{G})$ by continuity.
	
	Also, by the continuity of these distributions, it suffices to prove the matching in the case of $\tilde{G} = \Mp(W)$, since the counterpart for the $\GL$-factors is tautological.
\end{remark}

Before embarking on the proof of Theorem \ref{prop:stabilization-geom-LTF}, observe that (i) $\implies$ (ii) in Theorem \ref{prop:stabilization-geom-LTF}. To see this, apply Proposition \ref{prop:combinatorial-summation} to obtain
\begin{multline}\label{eqn:IEndo-geom-sum}
	I^{\Endo}_{\mathrm{geom}}(\overline{f_1}, f_2) = \sum_{\mathbf{G}^!} \iota(\tilde{G}, G^!) \sum_{M^!} \frac{|W^{G^!}_0|}{|W^{M^!}_0|}
	\underbracket{(-1)^{\dim \mathfrak{a}^{G^!}_{M^!}} S^{G^!}_{M^!, \mathrm{geom}}\left(\overline{f}_1^{G^![s]}, f_2^{G^![s]}\right)}_{=: S(\mathbf{G}^!, M^!)} \\
	= \sum_M \frac{|W^G_0|}{|W^M_0|} (-1)^{\dim \mathfrak{a}^G_M} \sum_{\mathbf{M}^!} \iota(\tilde{M}, M^!) \sum_s i_{M^!}(\tilde{G}, G^![s]) S^{G^![s]}_{M^!, \mathrm{geom}}\left(\overline{f}_1^{G^![s]}, f_2^{G^![s]}\right) \\
	= \sum_M \frac{|W^G_0|}{|W^M_0|} (-1)^{\dim \mathfrak{a}^G_M} \sum_{\mathbf{M}^!} \iota(\tilde{M}, M^!) I^{\Endo}_{\tilde{M}, \mathrm{geom}}(\mathbf{M}^!, \overline{f_1}, f_2).
\end{multline}
This equals $I_{\mathrm{geom}}(\overline{f_1}, f_2)$ if (i) holds.

\begin{lemma}[{Cf.\ \cite[pp.276--277]{Ar99}}]\label{prop:LTF-Endo-IP}
	For all $M \in \mathcal{L}(M_0)$ and $f_1, f_2 \in \orbI_{\asp}(\tilde{G}) \otimes \mes(G)$,
	\[ \sum_{\mathbf{M}^! \in \Endo_{\elli}(\tilde{M})} \iota(\tilde{M}, M^!) I^{\Endo}_{\tilde{M}, \mathrm{geom}}(\mathbf{M}^!, \overline{f_1}, f_2)
	= \int_{\Gamma_{G\text{-reg, ell}}(M)} I^{\Endo}_{\tilde{M}}(\gamma, \overline{f_1}, f_2) \dd\gamma. \]
\end{lemma}
\begin{proof}
	Put $n(\delta) := \left|\mathfrak{E}(M^!_\delta, M^!; F)\right|^{-1}$ for each $(\mathbf{M}^!, \delta)$. Since $\delta \mapsto \delta[s]$ does not change the measure on $\Sigma_{\mathrm{reg}}(M^!)$ and $n(\delta) = n(\delta[s])$, our assertion is tantamount to
	\begin{multline*}
		\sum_{\mathbf{M}^!} \iota(\tilde{M}, M^!) \int_{\Sigma_{G\text{-reg, ell}}(M^!)}  n(\delta) I^{\Endo}_{\tilde{M}}\left(\mathbf{M}^!, \delta, \overline{f_1}, f_2 \right) \dd\delta \\
		= \int_{\Gamma_{G\text{-reg, ell}}(M)} I^{\Endo}_{\tilde{M}}(\gamma, \overline{f_1}, f_2) \dd\gamma .
	\end{multline*}
	\index{n-delta@$n(\delta)$}
	
	Applying the definition of $I^{\Endo}_{\tilde{M}}(\gamma, \overline{f_1}, f_2)$ and \eqref{eqn:IMEndo-double-splitting}, we are reduced to proving that
	\begin{multline*}
		\sum_{\mathbf{M}^!} \iota(\tilde{M}, M^!) \int_{\Sigma_{G\text{-reg, ell}}(M^!)} n(\delta) \overline{I^{\tilde{L}_1, \Endo}_{\tilde{M}}(\mathbf{M}^!, \delta, f_{1, \tilde{L}_1})} I^{\tilde{L}_2, \Endo}_{\tilde{M}}(\mathbf{M}^!, \delta, f_{2, \tilde{L}_2}) \dd\delta \\
		= \int_{\Gamma_{G\text{-reg, ell}}(M)} \overline{I^{\tilde{L}_1, \Endo}_{\tilde{M}}(\tilde{\gamma}, f_{1, \tilde{L}_1})} I^{\tilde{L}_2, \Endo}_{\tilde{M}}(\tilde{\gamma}, f_{2, \tilde{L}_2}) \dd\gamma
	\end{multline*}
	where $\tilde{\gamma} \in \rev^{-1}(\gamma)$ is arbitrary, and $L_1, L_2 \in \mathcal{L}(M_0)$ satisfy $d^G_M(L_1, L_2) \neq 0$.
	
	The above follows from the endoscopic inner product formula in \cite[Corollary 5.2.5]{Li19} for $\tilde{M}$. Indeed, it suffices to check that all the functions to be paired are $L^2$. By ellipticity, it suffices to estimate them in the direction of $A_M(F)$, or of $A_{M^!}(F)$ by expanding them in terms of stable analogues. Ultimately, this reduces to the estimate for weighted orbital integrals in \cite[(5.7)]{Ar94}. For non-Archimedean case, the proof of that estimate requires Howe's conjecture for coverings, which is established in \cite{Luo17}.
\end{proof}

\begin{proposition}[{Cf.\ \cite[X.3.2 Proposition]{MW16-2}}]\label{prop:LTF-geom-reduction}
	Assume inductively that the local geometric Theorem \ref{prop:local-geometric} for $G$-regular elements holds if $\tilde{G}$ is replaced by $\tilde{L}$ for every $L \in \mathcal{L}(M) \smallsetminus \{G\}$. Let $f_1, f_2 \in \orbI_{\asp}(\tilde{G}) \otimes \mes(G)$. For every $M \in \mathcal{L}(M_0)$ and $\gamma \in \Gamma_{G\text{-reg, ell}}(M)$, define $\mathfrak{d}_M(\gamma, \overline{f_1}, f_2)$ to be
	\[ \left( \overline{I^{\Endo}_{\tilde{M}}(\tilde{\gamma}, f_1) - I_{\tilde{M}}(\tilde{\gamma}, f_1) }\right) I^{\tilde{M}}(\tilde{\gamma}, f_{2, \tilde{M}}) + \overline{I^{\tilde{M}}(\tilde{\gamma}, f_{1, \tilde{M}})} \left( I^{\Endo}_{\tilde{M}}(\tilde{\gamma}, f_2) - I_{\tilde{M}}(\tilde{\gamma}, f_2) \right) \]
	where $\tilde{\gamma} \in \rev^{-1}(\gamma)$ is arbitrary.
	\begin{enumerate}[(i)]
		\item For all $M \in \mathcal{L}(M_0)$, we have
		\begin{multline*}
			\sum_{\mathbf{M}^! \in \Endo_{\elli}(\tilde{M})} \iota(\tilde{M}, M^!) I^{\Endo}_{\tilde{M}, \mathrm{geom}}(\mathbf{M}^!, \overline{f_1}, f_2) - I_{\tilde{M}, \mathrm{geom}}(\overline{f_1}, f_2) \\
			= \int_{\Gamma_{G\text{-reg, ell}}(M)} \mathfrak{d}_M(\gamma, \overline{f_1}, f_2) \dd\gamma .
		\end{multline*}
		\item We have
		\begin{multline*}
			I^{\Endo}_{\mathrm{geom}}(\overline{f_1}, f_2) - I_{\mathrm{geom}}(\overline{f_1}, f_2) \\
			= \sum_{M \in \mathcal{L}(M_0)} \frac{|W^M_0|}{|W^G_0|} (-1)^{\dim \mathfrak{a}^G_M} \int_{\Gamma_{G\text{-reg, ell}}(M)} \mathfrak{d}_M(\gamma, \overline{f_1}, f_2) \dd\gamma .
		\end{multline*}
		\item We have
		\begin{multline*}
			I^{\Endo}_{\mathrm{geom}}(\overline{f_1}, f_2) - I_{\mathrm{geom}}(\overline{f_1}, f_2) \\
			= \sum_M |W^G(M)|^{-1} (-1)^{\dim \mathfrak{a}^G_M} \int_{\Gamma_{G\text{-reg, ell}}(M)} \mathfrak{d}_M(\gamma, \overline{f_1}, f_2) \dd\gamma,
		\end{multline*}
		where $M$ ranges over $W^G_0$-orbits in $\mathcal{L}(M_0)$, i.e.\ over conjugacy classes of Levi subgroups of $G$.
	\end{enumerate}
\end{proposition}
\begin{proof}
	We only have to show (i); one can then deduce (ii) from \eqref{eqn:IEndo-geom-sum}, and (iii) will follow from invariance and the familiar combinatorics. By Lemma \ref{prop:LTF-Endo-IP}, the left hand side of (i) equals
	\[ \int_{\Gamma_{G\text{-reg, ell}}(M)} \left( I^{\Endo}_{\tilde{M}}(\gamma, \overline{f_1}, f_2) - I_{\tilde{M}}(\gamma, \overline{f_1}, f_2) \right) \dd\gamma. \]
	
	The integrand may be expanded by splitting $I_{\tilde{M}}(\gamma, \overline{f_1}, f_2)$ and $I^{\Endo}_{\tilde{M}}(\gamma, \overline{f_1}, f_2)$: both are expressed as a sum over $L_1, L_2 \in \mathcal{L}(M)$ weighted by $d^G_M(L_1, L_2) \neq 0$. The summands take the form $\overline{I^{\tilde{L}_1, \Endo}_{\tilde{M}}} I^{\tilde{L}_2, \Endo}_{\tilde{M}} - \overline{I^{\tilde{L}_1}_{\tilde{M}}} I^{\tilde{L}_2}_{\tilde{M}}$. By assumption, only the terms with either $L_1 = G$ or $L_2 = G$ can survive. The only possibilities are $(L_1, L_2) = (G, M)$ or $(M, G)$. This yields $\mathfrak{d}_M(\gamma, \overline{f_1}, f_2)$.
\end{proof}

As a consequence, the $G$-regular case of the local geometric Theorem \ref{prop:local-geometric} will imply Theorem \ref{prop:stabilization-geom-LTF} (i), which in turn implies (ii) as already seen. We also obtain the following consequence.

\begin{corollary}\label{prop:LTF-stabilization-cusp}
	Retain the assumptions in Proposition \ref{prop:LTF-geom-reduction}. If $f_1, f_2 \in \orbI_{\asp, \cusp}(\tilde{G}) \otimes \mes(G)$, then $I_{\mathrm{geom}}(\overline{f_1}, f_2) = I^{\Endo}_{\mathrm{geom}}(\overline{f_1}, f_2)$.
\end{corollary}
\begin{proof}
	In this case, $\mathfrak{d}_M = 0$ unless $M = G$, but $I_{\tilde{G}} = I^{\Endo}_{\tilde{G}}$ is already known.
\end{proof}

\section{Endoscopic constructions: spectral side}\label{sec:LTF-spec-Endo}
Now comes the endoscopic construction on the spectral side. Consider
\[ I^{\Endo}_{\mathrm{disc}}(\overline{f_1}, f_2) := \sum_{\mathbf{G}^! \in \Endo_{\elli}(\tilde{G})} \iota(\tilde{G}, G^!) S^{G^!}_{\mathrm{disc}}(\overline{f_1}^{G^!}, f_2^{G^!}) \]
where $\overline{f_1}^{G^!} := \Trans_{\mathbf{G}^!, \tilde{G}^\dagger}(\overline{f_1})$ and $f_2^{G^!} := \Trans_{\mathbf{G}^!, \tilde{G}}(f_2)$.
\index{IEndo-disc-local@$I^{\Endo}_{\mathrm{disc}}(\cdot, \cdot)$}

Recall that $S^{G^!}_{\mathrm{disc}}$ is an integral over $\Phi_{\mathrm{disc}}(G^!)$, a space discrete modulo $i\mathfrak{a}^*_{G^!, F}$. Therefore, $I^{\Endo}_{\mathrm{disc}}(\overline{f_1}, f_2)$ is an integral taken over the space
\[ T^{\Endo}_{\mathrm{disc}}(\tilde{G}) := \bigsqcup_{\mathbf{G}^! \in \Endo_{\elli}(\tilde{G})} \Phi_{\mathrm{disc}}(G^!), \]
which is equipped with the Radon measure such that
\[ \int_{T^{\Endo}_{\mathrm{disc}}(\tilde{G})} \beta = \sum_{\mathbf{G}^! \in \Endo_{\elli}(\tilde{G})} \int_{\Phi_{\mathrm{disc}}(G^!)} \beta(\phi) \dd\phi \]
for all $\beta \in C_c(T^{\Endo}_{\mathrm{disc}}(\tilde{G}))$.

On the other hand, the spectral transfer factor $\Delta(\phi, \tau)$ makes sense for each $\phi = (\mathbf{G}^!, \phi) \in T^{\Endo}_{\mathrm{disc}}(\tilde{G})$ and $\tau \in T_-(\tilde{G})$. Define the $\mathbb{S}^1$-stable subset
\[ T_{\Endo\text{-disc}, -}(\tilde{G}) := \left\{\begin{array}{r|l}
	\tau \in T_-(\tilde{G}) & \exists \; (\mathbf{G}^!, \phi) \in \Phi^{\Endo}_{\mathrm{disc}}(\tilde{G}) \\
	& \Delta(\phi, \tau) \neq 0
\end{array}\right\}
\cup T_{\mathrm{disc}, -}(\tilde{G})
\]
of $T_-(\tilde{G})$, make it into a space that is discrete modulo $\mathbb{S}^1 \times i\mathfrak{a}^*_{G, F}$. As $i\mathfrak{a}^*_{G^!, F} \simeq i\mathfrak{a}^*_{G, F}$, this space can be equipped with the Radon measure\footnote{For metaplectic groups $\tilde{G} = \Mp(W)$, there are no outer automorphisms of endoscopic origin; see Remark \ref{rem:centerless}.} such that each $i\mathfrak{a}^*_{G, F}$-orbit containing $(L, \pi, r)$ carries the quotient measure divided by $|Z_{R_\pi}(r)|$.

\begin{lemma}\label{prop:spectral-square-change-of-variables}
	With the measures defined above, the formula of ``change of variables''
	\begin{multline}\
		\int_{\substack{\tau \in \mathbb{S}^1 \backslash T_{\Endo\text{-disc}, -}(\tilde{G}) \\ \tau = (L, \pi, r) }} \sum_{\phi \in T^{\Endo}_{\mathrm{disc}}(\tilde{G})} \beta(\phi) |\Delta(\phi, \tau)|^2 \alpha(\tau) |Z_{R_\pi}(r)| \dd\tau \\
		= \int_{T^{\Endo}_{\mathrm{disc}}(\tilde{G})} \sum_{\tau \in \mathbb{S}^1 \backslash T_{\Endo\text{-disc}, -}(\tilde{G})} \beta(\phi) |\Delta(\phi, \tau)|^2 \alpha(\tau) \dd\phi
	\end{multline}
	holds for all $\alpha \in C_c(\mathbb{S}^1 \backslash T_{\Endo\text{-disc}, -}(\tilde{G}))$ and $\beta \in C_c(T^{\Endo}_{\mathrm{disc}}(\tilde{G}))$.
\end{lemma}
\begin{proof}
	This is formally similar to the elliptic version \cite[Lemma 5.3]{Ar96}, with the same proof. Compared with \eqref{eqn:Tdisc-measure}, the measure $\dd\tau$ in \textit{loc.\ cit.} does not involve $|Z_{R_\pi}(r)|^{-1}$, whence the factor $|Z_{R_\pi}(r)| \dd\tau$ in our version.
	\footnote{Besides, $|Z_{R_\pi}(r)|$ is also missing in the elliptic version stated in \cite[Lemma 6.2.7]{Li19}. That result is only used in the proof of \cite[Proposition 8.3.5]{Li19}, and does not affect the proof since $Z_{R_\pi}(r)$ is not altered under $i\mathfrak{a}^*_{G, F}$-twists.}
\end{proof}

\begin{definition}
	\index{iEndo-tau@$i^{\Endo}(\tau)$}
	For all $\tau = (L, \pi, r) \in \tilde{T}_-(\tilde{G})$, set
	\[ i^{\Endo}(\tau) := i^{\tilde{G}, \Endo}(\tau) := \sum_{\phi \in \Phi^{\Endo}_{\mathrm{disc}}(\tilde{G})} \iota(\tilde{G}, G^!) |\Delta(\phi, \tau)|^2 s^{G^!}(\phi) |Z_{R_\pi}(r)| \]
	where $\Delta(\phi, \tau)$ is the spectral transfer factor in Definition \S\ref{def:spectral-transfer-factor}. This depends only on the $\mathbb{S}^1$-orbit of $\tau$, and vanishes outside $T_{\Endo\text{-disc}, -}(\tilde{G})$.
\end{definition}

The following is analogous to \cite[Proposition 6.3 (a)]{Ar03-3}.

\begin{lemma}\label{prop:IEndo-disc-LTF}
	For all $f_1, f_2 \in \orbI_{\asp}(\tilde{G}) \otimes \mes(G)$. we have
	\[ I^{\Endo}_{\mathrm{disc}}(\overline{f_1}, f_2) = \int_{\mathbb{S}^1 \backslash T_{\Endo\text{-disc}, -}(\tilde{G})} i^{\Endo}(\tau) \overline{\Theta_\tau(f_1)} \Theta_\tau(f_2) \dd\tau . \]
\end{lemma}
\begin{proof}
	Lemma \ref{prop:Delta-minus} implies $\overline{f_1}^{G^!} = \overline{f_1^{G^!}}$. The left hand side is hence the integral over $(\mathbf{G}^!, \phi) \in T^{\Endo}_{\mathrm{disc}}(\tilde{G})$ of
	\[ \iota(\tilde{G}, G^!) s^{G^!}(\phi) \overline{S\Theta_\phi(f_1^!)} S\Theta_\phi(f_2^!). \]
	
	Via spectral transfer, the expression above equals
	\begin{equation*}
		\iota(\tilde{G}, G^!) s^{G^!}(\phi) \sum_{\tau_1, \tau_2 \in T_-(\tilde{G})/\mathbb{S}^1} \overline{\Delta(\phi, \tau_1)} \Delta(\phi, \tau_2) \overline{\Theta_{\tau_1}(f_1)} \Theta_{\tau_2}(f_2).
	\end{equation*}
	
	In the integration over $\phi$, only the terms with $\tau_1 = \tau_2 \in T_{\Endo\text{-disc}, -}(\tilde{G})$ contribute. This is due to the orthogonality relation in \cite[Definition 6.2.4, Lemma 6.2.6]{Li19}; for the Archimedean case, see also \S 7.5 of \textit{loc.\ cit.} It remains to apply Lemma \ref{prop:spectral-square-change-of-variables}.
\end{proof}

The assertion of Theorem \ref{prop:stabilization-geom-LTF} (ii) is equivalent to that $I_{\mathrm{disc}}(\overline{f_1}, f_2) = I^{\Endo}_{\mathrm{disc}}(\overline{f_1}, f_2)$ for all $f_1, f_2$. It also entails the stabilization of the spectral coefficients $i(\tau)$. We record this formally as follows.

\begin{theorem}[{Cf.\ \cite[Corollary 6.4]{Ar03-3}}]\label{prop:stabilization-disc-LTF}
	Assuming the validity of Theorem \ref{prop:stabilization-geom-LTF} (ii), for all $f_1$, $f_2$ we have
	\[ I_{\mathrm{disc}}(\overline{f_1}, f_2) = I^{\Endo}_{\mathrm{disc}}(\overline{f_1}, f_2), \]
	and $i^{\Endo}(\tau) = i(\tau)$ for all $\tau$. In particular, under this assumption, $i^{\Endo}(\tau)$ is nonzero only when $\tau \in T_{\mathrm{disc}, -}(\tilde{G})$.
\end{theorem}
\begin{proof}
	The first equality is already seen. Now compare $I_{\mathrm{disc}}(\overline{f_1}, f_2)$ with the integral in Lemma \ref{prop:IEndo-disc-LTF}, noting that $\mathbb{S}^1 \backslash T_{\mathrm{disc}, -}(\tilde{G})$ is a union of connected components in $\mathbb{S}^1 \backslash T_{\Endo\text{-disc}, -}(\tilde{G})$, in a way compatible with the Radon measures. Since $f_1$ and $f_2$ are arbitrary, it follows that $i^{\Endo}(\tau) = i(\tau)$.
\end{proof}

\begin{proposition}\label{prop:LTF-spec-reduction}
	Under the assumptions of Proposition \ref{prop:LTF-geom-reduction}, we have
	\begin{equation*}
		I^{\Endo}_{\mathrm{disc}}(\overline{f_1}, f_2) - I_{\mathrm{disc}}(\overline{f_1}, f_2)
		= \int_{i\mathfrak{a}^*_{G, F}} \sum_{M, \chi_M, \tau} a(\tilde{M}, \chi_M, \tau) \overline{\Theta_{\tau_\lambda}(f_{1, \tilde{M}})} \Theta_{\tau_\lambda}(f_{2, \tilde{M}}) \dd\lambda
	\end{equation*}
	for all $f_1, f_2$, where
	\begin{itemize}
		\item $M$ ranges over representatives of conjugacy classes of Levi subgroups of $G$,
		\item $\chi_M$ ranges over $i\mathfrak{a}^*_{G, F}$-orbits of genuine characters $\widetilde{A_M} \to \mathbb{S}^1$,
		\item $\tau$ ranges over $\mathbb{S}^1 \backslash T_{\elli, -}(\tilde{M}) / i\mathfrak{a}^*_{G, F}$, for which we choose representatives in $T_{\elli, -}(\tilde{M})$, and we require that $\tau$ has central character equal to $\chi_M$ on $\widetilde{A_M}$, modulo $i\mathfrak{a}^*_{G, F}$;
		\item $a(\tilde{M}, \chi_M, \cdot)$ are identically zero for all but finitely many pairs $(M, \chi_M)$.
	\end{itemize}

	For all $M$ as above, we may choose $(\chi_M, \tau) \mapsto a(\tilde{M}, \chi_M, \tau)$ to be $W^G(M)$-invariant.
\end{proposition}
\begin{proof}
	Such an expression follows by sorting out the integrals in $I_{\mathrm{disc}}$ and $I^{\Endo}_{\mathrm{disc}}$, using Lemma \ref{prop:IEndo-disc-LTF}. To verify the vanishing of $a(\tilde{M}, \chi_M, \cdot)$ for all but finitely many pairs $(M, \chi_M)$, we use the stable counterparts in \cite[p.1160]{MW16-2} for the contribution from $I^{\Endo}_{\mathrm{disc}}$. For the contribution from $I_{\mathrm{disc}}$, simply repeat the argument from \cite[X.3.3 Lemme, p.1157]{MW16-2}; the only needed input is the first part of Lemma \ref{prop:elli-disc-finiteness}.
	
	The $W^G(M)$-invariance of $a(\tilde{M}, \chi_M, \tau)$ can be ensured by averaging, since $\Theta_{\tau_\lambda}(f_{i, \tilde{M}})$ is $W^G(M)$-invariant for $i=1,2$.
	
	Finally, to show that only the terms with $M \neq G$ contribute, we use the vanishing of the left hand side when $f_1$, $f_2$ are both cuspidal (Corollary \ref{prop:LTF-stabilization-cusp}), together with the trace Paley--Wiener theorem which implies that $\orbI_{\asp, \cusp}(\tilde{G})$ separates points in $\mathbb{S}^1 \backslash T_{\elli, -}(\tilde{G})$.
\end{proof}

\begin{definition}\label{def:M-cuspidal}
	Let $M \in \mathcal{L}(M_0)$. We say that $f \in \orbI_{\asp}(\tilde{G}) \otimes \mes(G)$ is $\tilde{M}$-cuspidal if we have $f_{\tilde{L}} = 0$ for all Levi subgroups $L$ of $G$ that does not contain any conjugate of $M$.
\end{definition}

For example, cuspidality is equivalent to $\tilde{G}$-cuspidality, and every function is $M_0$-cuspidal. The next two results are taken from \cite[X.3.4]{MW16-2}.

\begin{corollary}\label{prop:M-cuspidal-two}
	Fix $M \in \mathcal{L}(M_0)$ and an $\tilde{M}$-cuspidal $f_2 \in \orbI_{\asp}(\tilde{G}) \otimes \mes(G)$. If $I^{\Endo}_{\mathrm{disc}}(\overline{f_1}, f_2) = I_{\mathrm{disc}}(\overline{f_1}, f_2)$ holds for all $f_1 \in \orbI_{\asp}(\tilde{G}) \otimes \mes(G)$ which is $\tilde{M}$-cuspidal, then it holds for all $f_1 \in \orbI_{\asp}(\tilde{G}) \otimes \mes(G)$.
\end{corollary}
\begin{proof}
	Apply the expansion furnished by Proposition \ref{prop:LTF-spec-reduction} and the trace Paley--Wiener theorem.
\end{proof}

\begin{corollary}
	Suppose that $I^{\Endo}_{\tilde{L}}(\tilde{\gamma}, \cdot) = I_{\tilde{L}}(\tilde{\gamma}, \cdot)$ holds for all $L \supset M$ and all $\tilde{\gamma} \in \Gamma_{G\text{-reg, ell}}(\tilde{L})$. Then $I^{\Endo}_{\tilde{L}}(\tilde{\gamma}, f) = I_{\tilde{L}}(\tilde{\gamma}, f)$ holds for all $L \in \mathcal{L}(M_0)$, $\tilde{\gamma} \in \Gamma_{G\text{-reg, ell}}(\tilde{L})$ and $\tilde{M}$-cuspidal $f$.
\end{corollary}
\begin{proof}
	Suppose $f_1, f_2$ are both $\tilde{M}$-cuspidal. The expansion of Proposition \ref{prop:LTF-geom-reduction} involves only terms indexed by $L \in \mathcal{L}(M_0)$ that contains a conjugate of $M$. In particular, $I^{\Endo}_{\mathrm{geom}}(\overline{f_1}, f_2) = I_{\mathrm{geom}}(\overline{f_1}, f_2)$ holds. Therefore, Corollary \ref{prop:M-cuspidal-two} implies that $I^{\Endo}_{\mathrm{geom}}(\overline{f_1}, f_2) = I_{\mathrm{geom}}(\overline{f_1}, f_2)$ if only $f_2$ is assumed to be $\tilde{M}$-cuspidal.
	
	Looking back at the expansion in Proposition \ref{prop:LTF-geom-reduction}, with $f_2$ being $\tilde{M}$-cuspidal, it now becomes
	\[ \sum_{L \in \mathcal{L}(M_0)} \frac{|W^L_0|}{|W^G_0|} (-1)^{\dim \mathfrak{a}^G_L} \int_{\Gamma_{G\text{-reg, ell}}(L)} \overline{I^{\tilde{L}}(\tilde{\gamma}, f_{1, \tilde{L}})} \left( I^{\Endo}_{\tilde{L}}(\tilde{\gamma}, f_2) - I_{\tilde{L}}(\tilde{\gamma}, f_2) \right) \dd\gamma. \]
	
	Observe that $I^{\Endo}_{\tilde{L}}(\tilde{\gamma}, f_2)$ and $I_{\tilde{L}}(\tilde{\gamma}, f_2)$ are both $W^G(L)$-invariant in $\tilde{\gamma}$. As $f_1$ is arbitrary, we deduce $I^{\Endo}_{\tilde{L}}(\tilde{\gamma}, f_2) = I_{\tilde{L}}(\tilde{\gamma}, f_2)$ for all $L$ and $\tilde{\gamma} \in \Gamma_{G\text{-reg, ell}}(\tilde{L})$.
\end{proof}

\section{Key geometric hypothesis}\label{sec:key-geometric}
Fix $M \in \mathcal{L}(M_0)$. We assume inductively that
\begin{equation}\label{eqn:key-geometric-inductive}
	I^{\tilde{S}, \Endo}_{\tilde{L}}(\tilde{\gamma}, \cdot) = I^{\tilde{S}}_{\tilde{L}}(\tilde{\gamma}, \cdot)
\end{equation}
for all $\tilde{\gamma} \in \Gamma_{S\text{-reg}}(\tilde{L})$, where $L \subset S$ are Levi subgroups of $G$ containing $M$ such that
\[ \text{either}\quad S \subsetneq G \quad \text{or}\quad M \subsetneq L. \]

Let $f \in \orbI_{\asp}(\tilde{G}) \otimes \mes(G)$. Recall the $\epsilon_{\tilde{M}}(f) \in \orbI_{\mathrm{ac}, \asp, \cusp}(\tilde{M}) \otimes \mes(M)$ constructed in Definition--Propositions \ref{def:epsilonM-nonarch} and \ref{def:epsilonM-arch}. For all $\mathbf{M}^! \in \Endo_{\elli}(\tilde{M})$, one has the function
\[ \delta \mapsto \Trans_{\mathbf{M}^!, \tilde{M}}\left( \epsilon_{\tilde{M}}(f)\right)(\delta), \quad \delta \in \Sigma_{\mathrm{reg}}(M^!) \]
given by stable orbital integrals; see also \eqref{eqn:collective-transfer-fcn}. On the other hand, it also makes sense to consider
\[ \delta \mapsto \Trans_{\mathbf{M}^!, \tilde{M}}\left(f_{\tilde{M}}\right)(\delta). \]

We now state the \emph{key geometric hypothesis}.

\begin{hypothesis}\label{hyp:key-geometric}
	\index{key geometric hypothesis}
	\index{epsilonMshrek@$\epsilon(\mathbf{M}^{"!}, \delta)$}
	For all $\mathbf{M}^! \in \Endo_{\elli}(\tilde{M})$, there exists a smooth function
	\[ \epsilon(\mathbf{M}^!, \cdot): M^!_{G\text{-reg}}(F) \to \CC \]
	characterized uniquely by
	\[ \Trans_{\mathbf{M}^!, \tilde{M}}\left(\epsilon_{\tilde{M}}(f)\right)(\delta) = \epsilon(\mathbf{M}^!, \delta) \Trans_{\mathbf{M}^!, \tilde{M}}\left(f_{\tilde{M}}\right)(\delta) \]
	for all $\delta \in \Sigma_{G\text{-reg}}(M^!)$ and all $f \in \orbI_{\asp}(\tilde{G}) \otimes \mes(G)$.
\end{hypothesis}

\begin{remark}\label{rem:key-geometric-invariance}
	The group $W^G(M)$ acts on the $\GL$-factors of $M^!$, and $\epsilon(\mathbf{M}^!, \cdot)$ is $W^G(M)$-invariant. To see this, we need the $W^G(M)$-invariance of $f_{\tilde{M}}$ and $\epsilon_{\tilde{M}}(f)$; for the latter, we use the invariance of $I_{\tilde{M}}$ under automorphisms \cite[Lemma 3.3]{Ar98}, as well as its endoscopic version.
\end{remark}

\begin{remark}
	Since $\epsilon_{\tilde{M}}(f)$ is cuspidal for all $f$, we have $\epsilon(\mathbf{M}^!, \delta) = 0$ unless $\delta$ is elliptic.
\end{remark}

Our aim is to show that $I^{\Endo}_{\tilde{M}}(\tilde{\gamma}, f) = I_{\tilde{M}}(\tilde{\gamma}, f)$ for all $\tilde{\gamma}$ and $f$, or equivalently that $\epsilon_{\tilde{M}}(f) = 0$. By adjoint geometric transfer (Proposition \ref{prop:inverse-transfer-0}), this amounts to showing $\Trans_{\mathbf{M}^!, \tilde{M}}(\epsilon_{\tilde{M}}(f)) = 0$ for all $\mathbf{M}^! \in \Endo_{\elli}(\tilde{M})$. This will be achieved, under Hypothesis \ref{hyp:key-geometric}, by showing that $\epsilon(\mathbf{M}^!, \delta) = 0$ for all $(\mathbf{M}^!, \delta)$. We shall proceed in steps.

\begin{lemma}[Cf.\ {\cite[Lemma 6.5]{Ar03-3}}]\label{prop:epsilon-vanishing-1}
	Under the inductive assumptions \eqref{eqn:key-geometric-inductive} and the key geometric Hypothesis \ref{hyp:key-geometric}, we have $\epsilon(\mathbf{M}^!, \delta) + \overline{\epsilon(\mathbf{M}^!, \delta)} = 0$ for all $(\mathbf{M}^!, \delta)$.
\end{lemma}
\begin{proof}
	We may and do assume that $M$ is a proper standard Levi subgroup, and $\delta$ is elliptic in $M^!$. Take $f_1, f_2 \in \orbI_{\asp}(\tilde{G}) \otimes \mes(G)$ and expand $I^{\Endo}_{\mathrm{geom}}(\overline{f_1}, f_2) - I_{\mathrm{geom}}(\overline{f_1}, f_2) = I^{\Endo}_{\mathrm{disc}}(\overline{f_1}, f_2) - I_{\mathrm{disc}}(\overline{f_1}, f_2)$ in two ways, namely via Propositions \ref{prop:LTF-geom-reduction} (iii) and \ref{prop:LTF-spec-reduction}, to obtain
	\begin{multline*}
		\sum_{L:\; \text{Levi}/\text{conj}} |W^G(L)|^{-1} (-1)^{\dim \mathfrak{a}^G_L} \int_{\Gamma_{G\text{-reg, ell}}(L)} \dd\gamma \cdot \\
		\overline{\left(I^{\Endo}_{\tilde{L}}(\tilde{\gamma}, f_1) - I_{\tilde{L}}(\tilde{\gamma}, f_1)\right)} I^{\tilde{L}}(\tilde{\gamma}, f_{2, \tilde{L}})
		+ \overline{I^{\tilde{L}}(\tilde{\gamma}, f_{1, \tilde{L}})} \left(I^{\Endo}_{\tilde{L}}(\tilde{\gamma}, f_2) - I_{\tilde{L}}(\tilde{\gamma}, f_2)\right) \\
		= \int_{i\mathfrak{a}^*_{G, F}} \sum_{\substack{L, \chi_L, \tau \\ L/\text{conj} \\ L \neq G}} a(\tilde{L}, \chi_L, \tau) \overline{\Theta_{\tau_\lambda}(f_{1, \tilde{L}})} \Theta_{\tau_\lambda}(f_{2, \tilde{L}}) \dd\lambda .
	\end{multline*}

	Fix representatives for $L$'s in the sum above. Consider any $L' \supsetneq M$. We contend that $a(\tilde{L}', \chi_{L'}, \tau) = 0$ for all $\chi_{L'}$ and $\tau$. Take $f_1, f_2$ to be $\tilde{L}'$-cuspidal. The geometric (i.e.\ left hand) side of the displayed equation vanishes by inductive assumptions. By choosing $f_1 = f_2$ suitably and using the trace Paley--Wiener theorem plus the properties stated in Proposition \ref{prop:LTF-spec-reduction}, it follows that $a(\tilde{L}', \chi_{L'}, \tau) = 0$ as desired.
	
	Next, take $f_1, f_2$ to be $\tilde{M}$-cuspidal. On the spectral side, only the terms with $L=M$ survive. On the geometric side, only the terms with $L \supset M$ survive; furthermore, inductive assumptions kill the terms with $L \neq M$. Hence the geometric side can be simplified into a sum of two inner products taken over $\Gamma_{G\text{-reg, ell}}(\tilde{M})$, and we also have
	\[ I^{\Endo}_{\tilde{M}}(\tilde{\gamma}, \cdot) - I_{\tilde{M}}(\tilde{\gamma}, \cdot) = I^{\tilde{M}}(\tilde{\gamma}, \epsilon_{\tilde{M}}(\cdot)). \]

	In turn, by using
	\begin{itemize}
		\item the endoscopic inner product formula for $\tilde{M}$ in \cite[Corollary 5.2.5]{Li19},
		\item the key geometric Hypothesis \ref{hyp:key-geometric},
	\end{itemize}
	the geometric side is seen to equal
	\begin{multline}\label{eqn:key-geometric-aux-0}
		(\text{nonzero constant}) \cdot \sum_{\mathbf{M}^! \in \Endo_{\elli}(\tilde{M})} \iota(\tilde{M}, M^!) \\
		\int_{\Sigma_{\tilde{M}\text{-reg, ell}}(M^!)} n(\delta) \left( \epsilon(\mathbf{M}^!, \delta) + \overline{\epsilon(\mathbf{M}^!, \delta)}\right)
		\overline{\Trans_{\mathbf{M}^!, \tilde{M}} (f_{1, \tilde{M}})(\delta)} \cdot \Trans_{\mathbf{M}^!, \tilde{M}}(f_{2, \tilde{M}})(\delta) \dd\delta .
	\end{multline}

	Fix an $\mathbf{M}^! \in \Endo_{\elli}(\tilde{M})$. In what follows, we take $f_1$ such that
	\begin{itemize}
		\item the $W^G(M)$-invariant function $\delta' \mapsto \Trans_{\mathbf{M}^!, \tilde{M}} (f_{1, \tilde{M}})(\delta')$ on $\Sigma_{\mathrm{reg}}(M^!)$ has compact supported in $\Sigma_{\elli, G\text{-reg}}(M^!)$;
		\item for any $\mathbf{M}^{!!} \in \Endo_{\elli}(\tilde{M})$ with $\mathbf{M}^{!!} \neq \mathbf{M}^!$, we have $\Trans_{\mathbf{M}^{!!}, \tilde{M}} (f_{1, \tilde{M}}) = 0$.
	\end{itemize}

	In fact, we will take a sequence of $f_1$ such that $\Trans_{\mathbf{M}^!, \tilde{M}} (f_{1, \tilde{M}})$ approximates the Dirac measure attached to the $W^G(M)$-orbit of any given $\delta \in \Sigma_{\elli, G\text{-reg}}(M^!)$. The existence such a sequence follows from Theorem \ref{prop:image-transfer} and the characterization of orbital integrals in \cite{Bo94a, Vi81} (cf.\ the proof of Proposition \ref{prop:epsilonM-arch-local}).
	
	With these assumptions, we infer that
	\[ c_1 := \left( \overline{\epsilon(\mathbf{M}^!, \delta)} + \epsilon(\mathbf{M}^!, \delta) \right) \Trans_{\mathbf{M}^!, \tilde{M}} (f_{1, \tilde{M}})(\delta) \in S\orbI_{\cusp}(M^!) \otimes \mes(M^!), \]
	whilst by $\tilde{M}$-cuspidality of $f_2$, we also have
	\[ c_2 := \Trans_{\mathbf{M}^!, \tilde{M}}(f_{2, \tilde{M}}) \in S\orbI_{\cusp}(M^!) \otimes \mes(M^!). \]

	Hence the integral in \eqref{eqn:key-geometric-aux-0} can be expressed using $S^{M^!}_{\mathrm{geom}}(\overline{c_1}, c_2) = S^{M^!}_{\mathrm{disc}}(\overline{c_1}, c_2)$, simplified according to Remark \ref{rem:simple-stable-LTF}, as an inner product of stable characters indexed by $\Phi_{\mathrm{bdd}, 2}(M^!)$, defined using linear combinations of integrals over $i\mathfrak{a}^*_{M^!, F} \simeq i\mathfrak{a}^*_{M, F}$.

	On the other hand, the expression
	\[ I^{\Endo}_{\mathrm{disc}}(\overline{f_1}, f_2) - I_{\mathrm{disc}}(\overline{f_1}, f_2) = \int_{i\mathfrak{a}^*_{G, F}} \sum_{\chi_M, \tau} a(\tilde{M}, \chi_M, \tau) \overline{\Theta_{\tau_\lambda}(f_{1, \tilde{M}})} \Theta_{\tau_\lambda}(f_{2, \tilde{M}}) \]
	can be translated to $\mathbf{M}^!$ on the endoscopic side via spectral transfer. Keep $f_1$ fixed and vary the $\tilde{M}$-cuspidal test function $f_2$. Since $M \neq G$, a familiar comparison with the previous integrals over $i\mathfrak{a}^*_{M, F}$ shows that both expansions are zero. In particular, the integral over $\Sigma_{G\text{-reg}, \elli}(M^!)$ in \eqref{eqn:key-geometric-aux-0} vanishes for such $(f_1, f_2)$.
	
	Notice that the three terms and $n(\delta)$ in the integrand of \eqref{eqn:key-geometric-aux-0} are all $W^G(M)$-invariant by Remark \ref{rem:key-geometric-invariance}. Given $(\mathbf{M}^!, \delta)$, let $f_{1, \tilde{M}}$ approach the Dirac measure attached to $W^G(M) \delta$ and vary $f_2$. This implies $\epsilon(\mathbf{M}^!, \delta) + \overline{\epsilon(\mathbf{M}^!, \delta)} = 0$, as asserted.
\end{proof}

\begin{lemma}\label{prop:epsilon-vanishing-2}
	Under the inductive assumptions \eqref{eqn:key-geometric-inductive} and the key geometric Hypothesis \ref{hyp:key-geometric}, we have $\epsilon(\mathbf{M}^!, \delta) = \overline{\epsilon(\mathbf{M}^!, \delta)}$ for all $(\mathbf{M}^!, \delta)$.
\end{lemma}
\begin{proof}
	By taking complex conjugates in the characterization of $\epsilon(\mathbf{M}^!, \cdot)$ and using Lemma \ref{prop:Delta-minus}, we obtain
	\[ \Trans_{\mathbf{M}^!, \tilde{M}^\dagger}(\overline{\epsilon_{\tilde{M}}(f)})(\delta) = \overline{\epsilon(\mathbf{M}^!, \delta)} \Trans_{\mathbf{M}^!, \tilde{M}^\dagger}(\overline{f}_{\tilde{M}^\dagger})(\delta). \]
	
	Proposition \ref{prop:IM-cplx-conj} implies $\overline{\epsilon_{\tilde{M}}(f)} = \epsilon_{\tilde{M}^\dagger}(\overline{f})$ for all $f$. Therefore, we see that Hypothesis \ref{hyp:key-geometric} also holds for $\tilde{M}^\dagger \subset \tilde{G}^\dagger$. Moreover, by denoting the antipodal version of $\epsilon(\mathbf{M}^!, \cdot)$ as $\epsilon(\mathbf{M}^!, \cdot)^\dagger$, we have $\epsilon(\mathbf{M}^!, \cdot)^\dagger = \overline{\epsilon(\mathbf{M}^!, \cdot)}$.
	
	Next, describe $M$ by a decomposition $W = \bigoplus_{j \in J} (\ell_j \oplus \ell'_j) \oplus W^\flat$ into symplectic vector subspaces. Consider the isomorphism $\Theta_g$ in Proposition \ref{prop:MVW} that is adapted to $\tilde{M}$ (Lemma \ref{prop:MVW-adapted}). Note that $g$ induces isomorphisms of symplectic $F$-vector spaces
	\[ (W, \lrangle{\cdot|\cdot}) \rightiso (W, -\lrangle{\cdot|\cdot}), \quad (W^\flat, \lrangle{\cdot|\cdot}) \rightiso (W^\flat, -\lrangle{\cdot|\cdot}) \]
	and is inner on each $\GL$-factor of $\tilde{M}$. The endoscopy for metaplectic groups is defined in terms of linear algebraic data, when $\psi$ is kept fixed; consequently, $\Theta_g$ transports the whole formalism of endoscopy for $\tilde{M} \subset \tilde{G}$ to $\tilde{M}_- \subset \tilde{G}_-$.
	
	Identifying $\tilde{M}^\dagger \subset \tilde{G}^\dagger$ with $\tilde{M}_- \subset \tilde{G}_-$ by Proposition \ref{prop:MVW}, we conclude that $\epsilon(\mathbf{M}^!, \cdot) = \epsilon(\mathbf{M}^!, \cdot)^\dagger$. The assertion follows immediately.
\end{proof}

\begin{corollary}\label{prop:local-geometric-G-reg}
	Under the inductive assumptions \eqref{eqn:key-geometric-inductive} and the key geometric Hypothesis \ref{hyp:key-geometric}, we have $I_{\tilde{M}}(\tilde{\gamma}, \cdot) = I^{\Endo}_{\tilde{M}}(\tilde{\gamma}, \cdot)$ for all $\tilde{\gamma} \in \Gamma_{G\text{-reg}}(\tilde{M})$; equivalently, $\epsilon_{\tilde{M}}(f) = 0$.
	
	In particular, the assertions in Theorem \ref{prop:stabilization-geom-LTF} hold under these assumptions.
\end{corollary}
\begin{proof}
	Combining Lemmas \ref{prop:epsilon-vanishing-1} and \ref{prop:epsilon-vanishing-1}, we obtain $\epsilon(\mathbf{M}^!, \cdot) = 0$ for all $\mathbf{M}^! \in \Endo_{\elli}(\tilde{M})$. This implies $\epsilon_{\tilde{M}}(f) = 0$, and the first part follows. As for the second part, see the discussions after Proposition \ref{prop:LTF-geom-reduction}.
\end{proof}

In the local setting, we have reduced everything to the Hypothesis \ref{hyp:key-geometric}. Its verification is the topic of \S\ref{sec:proof-key-geometric}.

\begin{remark}
	The proof of $\epsilon_{\tilde{M}}(f) = 0$ given above can be compared with the simple case treated in \cite[p.1169]{MW16-2} for twisted trace formulas. This should not be over-generalized: the case of general coverings will require more sophisticated techniques, provided that a theory of endoscopy is available.
\end{remark}
\chapter{Global trace formula: the spectral side}\label{sec:spectral-side}
This part marks the end of stabilization. Let $F$ be a number field and consider a covering $\rev: \tilde{G} \to G(\A_F)$ of metaplectic type, together with auxiliary data chosen as in \S\ref{sec:TF-geom}. Take a finite set $V$ of places of $F$ containing all ramified places.

Let us disrupt the ordering of sections and begin with \S\ref{sec:Ispec}. The spectral side of the invariant trace formula \cite{Li14b} is phrased as
\[ I_{\mathrm{spec}}(f) = \sum_{t \geq 0} I_t(f), \quad f \in \orbI_{\asp}(\tilde{G}_V, \tilde{K}_V) \otimes \mes(G(F_V)), \]
where $t$ specifies the height of the imaginary part of infinitesimal characters in the spectral expansion. Each $I_t$ can be decomposed into $\sum_{\substack{\nu \in \mathfrak{h}^*_{\CC}/W_\infty \\ \|\Im(\nu)\| = t}} I_\nu$, according to infinitesimal characters.

Similarly to the geometric side $I_{\mathrm{geom}}$, these distributions are ``compressed'' and it is also necessary to deal with the pre-compressed form, built upon the discrete parts
\[ I_{\mathrm{disc}, t}(\dot{f}) = \sum_{\pi \in \Pi_{t, -}(\tilde{G}^1)} a^{\tilde{G}}_{\mathrm{disc}}(\pi) \Tr \pi(\dot{f}|_{\tilde{G}^1}) \]
where $\dot{f} \in \orbI_{\asp}(\tilde{G}, \tilde{K}) \otimes \mes(G(\A_F))$, as well as their counterparts $I^{\tilde{M}}_{\mathrm{disc}, t}$ for various $M \in \mathcal{L}(M_0)$. We shall view $I_{\mathrm{disc}, t}$ as a formal linear combination of genuine irreducible characters of $\tilde{G}$. Its evaluation at $\dot{f}$ makes sense; it is convenient to justify such operations by employing the notion of \emph{semi-finiteness} (Definition \ref{def:semi-finite-distribution}) from \cite[X.5.3]{MW16-2}, which we adopt here. Let $\rho$ be such a formal linear combination, then $\rho$ can be decomposed according to
\begin{itemize}
	\item the infinitesimal characters $\nu$, written as $\rho = \sum_\nu \rho_\nu$;
	\item the Satake parameters $c^V$ off $V$, if we only look at the part $\rho^{V\text{-nr}}$ unramified outside $V$, so that $\rho^{V\text{-nr}} = \sum_{\nu, c^V} \rho_{\nu, c^V}$.
\end{itemize}
Applied to $\rho := I_{\mathrm{disc}, t}$, this yields
\[ I_{\mathrm{disc}, \nu}(\dot{f}), \quad I_{\mathrm{disc}, \nu, c^V}(\dot{f}). \]
We also obtain distributions $I_{\mathrm{disc}, \nu, c^V}$, etc.\ on $\tilde{G}_V$, by evaluating the earlier distributions at $\dot{f} := f \prod_{v \notin V} f_{K_v}$ for $f \in \orbI_{\asp}(\tilde{G}_V, \tilde{K}_V) \otimes \mes(G(F_V))$.

Given $\nu$, we can express $I_\nu(f)$ as
\[ I_\nu(f) = \sum_{M \in \mathcal{L}(M_0)} \frac{|W^M_0|}{|W^G_0|} I_{\tilde{M}}^{\mathrm{glob}}\left(I^{\tilde{M}}_{\mathrm{disc}, \nu}, f\right), \]
where $I_{\tilde{M}}^{\mathrm{glob}}$ is the (invariant) \emph{global weighted character} to be defined in \S\ref{sec:spec-nr}. The global weighted character involves the following two ingredients introduced in \S\ref{sec:weighted-characters}:
\begin{itemize}
	\item the semi-local weighted characters $I_{\tilde{M}_V}(\pi_V, \lambda, X, f)$ for not necessarily unitary representations $\pi_V \in \Pi_-(\tilde{M}_V)$;
	\item the factors $r^{\tilde{G}}_{\tilde{M}}(c^V)$ related to Satake parameters $c^V$ of $\tilde{M}$, assuming $c^V$ is automorphic or merely quasi-automorphic.
\end{itemize}

All these objects admit endoscopic counterparts. This step requires the transfer of centers of universal enveloping algebras, obtained in \cite{Li19}, as well as the spherical fundamental lemma \cite{Luo18}. The stabilization of the spectral side of the trace formula is stated in \S\ref{sec:spectral-stabilization}. It boils down to the stabilization of $I_{\tilde{M}}^{\mathrm{glob}}$ and $I_{\mathrm{disc}, \nu, c^V}$.

The case of $I_{\tilde{M}}^{\mathrm{glob}}$ is relatively easy. One ingredient is a special case of the stabilization of local weighted characters, settled in \S\ref{sec:local-spec}; we remark that the general case will be deduced from the local geometric Theorem \ref{prop:local-geometric}. The other ingredient is the stabilization of $r^{\tilde{G}}_{\tilde{M}}(c^V)$: this is some sort of spectral weighted fundamental lemma, and can be proved directly by the method in \S\ref{sec:easy-stabilization}.

The case of $I_{\mathrm{disc}, \nu, c^V}$ is the hardcore of this work. In Lemma \ref{prop:geom-diff-h}, we will use the available stabilization for $A^{\tilde{G}}(V, \mathcal{O})$ to reduce it to the stabilization of $I_{\mathrm{geom}}$. Recall that the stabilization of $I_{\mathrm{geom}}$ is a consequence of the key geometric Hypothesis \ref{hyp:key-geometric}, as already seen in \S\ref{sec:key-geometric}.

In \S\ref{sec:spec-stabilization-special}, we pick $M \in \mathcal{L}(M_0)$ and impose local and global inductive assumptions on various matching theorems. We then prove the stabilization of $I_{\mathrm{disc}, \nu}$ in Proposition \ref{prop:Idisc-stabilization-special} for $f = \prod_v f_v$ which is $\tilde{M}$-cuspidal (Definition \ref{def:M-cuspidal}) at two distinct places, and has vanishing singular orbital integrals at some non-Archimedean place.

In \S\ref{sec:proof-key-geometric} a local-global argument reduces the key geometric Hypothesis \ref{hyp:key-geometric} to the aforementioned special case of spectral stabilization. This will enable us to prove all the remaining matching theorems in this work, thereby completing the long induction process. The remaining steps are summarized in \S\ref{sec:end-stabilization}.

The strategy sketched above roughly follows \cite{Ar03-3}, but also incorporates some viewpoints from \cite{MW16-2}. However, the result of Finis--Lapid--Müller \cite{FLM11} on the convergence of $I_{\mathrm{spec}}(f)$ has not been written down for coverings. Therefore, instead of following the paradigm of \cite{MW16-2}, we revert to Arthur's original approach via the multiplier convergence estimates, which suffices for our purposes here.

\section{Invariant weighted characters}\label{sec:weighted-characters}
We will extend the definition of weighted characters in \S\ref{sec:weighted-characters-nonarch}, \S\ref{sec:weighted-characters-arch} to non-unitary case. Next, we will introduce the invariant weighted characters and their semi-local counterparts, for not necessarily unitary representations. The semi-local version will intervene in the spectral side of the invariant trace formula.

\subsection{Local version}
To begin with, consider a local field $F$ of characteristic zero. We work with a covering of metaplectic type
\[ \rev: \tilde{G} \to G(F), \quad G = \prod_{i \in I} \GL(n_i) \times \Sp(W). \]
We also fix a symplectic basis for $(W, \lrangle{\cdot|\cdot})$, the resulting minimal Levi subgroup $M_0$ and positive-definite invariant quadratic form on $\mathfrak{a}_0$, and a maximal compact subgroups $K \subset G(F)$ in good position relative to $M_0$.

Below is a review of \cite{Ar88-2} or \cite[X.4.2]{MW16-2} in the metaplectic setting. Let $M \in \mathcal{L}(M_0)$.
\begin{description}
	\item[Non-Archimedean case] Define $\mathfrak{a}_{M, F}$ and its Pontryagin dual $i\mathfrak{a}^*_{M, F} = i\mathfrak{a}_M^* / i\mathfrak{a}_{M, F}^\vee$ as in \S\ref{sec:weighted-characters-nonarch}. Set $\mathfrak{a}^*_{M, F, \CC}$ to be the $\CC$-torus $\mathfrak{a}_{M, \CC}^* / i\mathfrak{a}_{M, F}^\vee$.
	\item[Archimedean case] Define $\mathfrak{a}_{M, F} := \mathfrak{a}_M$ and its Pontryagin dual $i\mathfrak{a}^*_{M, F} := i\mathfrak{a}^*_M$. Set $\mathfrak{a}^*_{M, F, \CC}$ to be the $\CC$-vector space $\mathfrak{a}^*_{M, \CC}$.
\end{description}
In each case, we have $\mathfrak{a}_{M, F} = H_M(M(F))$. Equip $\mathfrak{a}_{M ,F}$ and $i\mathfrak{a}^*_{M,F}$ with dual Haar measures; in the non-Archimedean case, we use the counting measure on $\mathfrak{a}_{M, F}$

Let $\pi \in \Pi_-(\tilde{M})$. The associated canonically normalized weighted character is
\[ J_{\tilde{M}}(\pi_\lambda, f) = J_{\tilde{M}}^{\tilde{G}}(\pi_\lambda, f), \quad f \in C^\infty_{c, \asp}(\tilde{G}, \tilde{K}) \otimes \mes(G(F)), \]
viewed as a meromorphic (resp.\ rational) family as $\lambda$ varies over $\mathfrak{a}^*_{M, F, \CC}$. For generic $\lambda$ such that $\lambda + i\mathfrak{a}^*_{M, F}$ is disjoint from the singularities, one defines
\index{JMpiX-lambda@$J_{\tilde{M}}(\pi, \lambda, X, f)$}
\begin{align*}
	J_{\tilde{M}}(\pi, \lambda, X, f) & = J_{\tilde{M}}^{\tilde{G}}(\pi, \lambda, X, f) \\
	& := \int_{\lambda + i\mathfrak{a}^*_{M, F}} J_{\tilde{M}}(\pi_\mu, f) e^{-\lrangle{\mu, X}} \dd\mu
\end{align*}
where $X \in \mathfrak{a}_M$. By shifting contours, $\lambda \mapsto J_{\tilde{M}}(\pi, \lambda, X, f)$ is seen to be locally constant around each point where it is defined. Also, $J_{\tilde{M}}(\pi, \lambda, X, \cdot)$ is concentrated at $H_{\tilde{G}}(\cdot) = X_G$, where $X_G$ is the projection of $X$ along the natural map $\mathfrak{a}_{M, F} \to \mathfrak{a}_{G, F}$. Thus $J_{\tilde{M}}(\pi, \lambda, X, \cdot)$ extends to $C^\infty_{\mathrm{ac}, \asp}(\tilde{G}) \otimes \mes(G(F))$.

The invariant version $I_{\tilde{M}}(\pi, \lambda, X, f) = I^{\tilde{G}}_{\tilde{M}}(\pi, \lambda, X, f)$ is constructed inductively by
\[ \sum_{L \in \mathcal{L}(M_0)} I^{\tilde{L}}_{\tilde{M}}\left( \pi, \lambda, X, \phi_{\tilde{L}}(f) \right) = J_{\tilde{M}}(\pi, \lambda, X, f) \]
for $\lambda$ off the singular hyperplanes (infinitely many for Archimedean $F$). They are concentrated at $H_{\tilde{G}}(\cdot) = X_G$, thus defined for $f \in \orbI_{\asp, \mathrm{ac}}(\tilde{G}, \tilde{K}) \otimes \mes(G(F))$. As usual, $\phi_{\tilde{L}}$ is defined in terms of canonically normalized tempered weighted characters.
\index{IMpiX-lambda@$I_{\tilde{M}}(\pi, \lambda, X, f)$}

In the case $\pi \in \Pi_{\mathrm{unit}, -}(\tilde{M})$, the weighted character $J_{\tilde{M}}(\pi_\lambda, f)$ is analytic for $\lambda \in i\mathfrak{a}^*_{M, F}$, and the distributions
\begin{equation}\label{eqn:I-0-0}
	J_{\tilde{M}}(\pi, X, f) := J_{\tilde{M}}(\pi, 0, X, f) , \quad I_{\tilde{M}}(\pi, X, f) := I_{\tilde{M}}(\pi, 0, X, f)
\end{equation}
have been introduced in \S\ref{sec:weighted-characters-nonarch} and \S\ref{sec:weighted-characters-arch}. These definitions also extend by linearity to the case when $\pi$ is merely a virtual genuine representation. Their behavior under parabolic induction is described as follows.

\begin{proposition}\label{prop:inv-weighted-descent}
	Suppose that $R \in \mathcal{L}^M(M_0)$ and $\sigma$ is a virtual genuine representation of $\tilde{R}$. Then its normalized parabolic induction $\sigma^{\tilde{M}}$ to $\tilde{M}$ satisfies
	\[ I_{\tilde{M}}\left(\sigma^{\tilde{M}}, \lambda, X, f\right) = \sum_{L \in \mathcal{L}(R)} d^G_R(L, M) \int_{\substack{Y \in \mathfrak{a}_{R, F} \\ Y_M = X}} I^{\tilde{L}}_{\tilde{R}}(\sigma, \lambda, Y, f_{\tilde{L}}) \dd Y \]
	for all $\lambda$ off the singularities, $X$ and $f$.
\end{proposition}
\begin{proof}
	Same as \cite[p.1184]{MW16-2}. It is of an analytic-combinatorial nature.
\end{proof}

\subsection{Semi-local version}
Now let $F$ be a number field and consider
\[ \rev: \tilde{G} \to G(\A_F), \quad G = \prod_{i \in I} \GL(n_i) \times \Sp(W) \]
as in \S\ref{sec:metaplectic-type-adelic}. Denote by $\tilde{G}_\infty$ the preimage of $\prod_{v \mid \infty} G(F_v)$ under $\rev$. We also fix a symplectic basis for $(W, \lrangle{\cdot|\cdot})$, the resulting $M_0$, etc.\ as well as maximal compact subgroups $K_v \subset G(F_v)$ in good position relative to $M_0$ for each place $v$, such that $K_v$ is hyperspecial whenever $v \notin V_{\mathrm{ram}}$.

Take a finite set $V$ of places of $F$, such that $V$ contains the Archimedean places. In particular, $V$ satisfies the \emph{closure property} \cite[Définition 3.1]{Li14b}.

In view of that property, one can still define $\mathfrak{a}_{M, V}$, its Pontryagin dual $i\mathfrak{a}^*_{M, V}$, and $\mathfrak{a}^*_{M, V, \CC}$ as in the local case, by setting $\mathfrak{a}_{M, V} := H_M(M(F_V))$. In any case, for $v \in V$ we have natural maps
\[ \mathfrak{a}_{M, F_v} \hookrightarrow \mathfrak{a}_{M, V} \hookrightarrow \mathfrak{a}_M, \quad \mathfrak{a}^*_{M, \CC} \twoheadrightarrow \mathfrak{a}^*_{M, V, \CC} \twoheadrightarrow \mathfrak{a}^*_{M, F_v, \CC}. \]

Let $M \in \mathcal{L}(M_0)$. Consider
\begin{itemize}
	\item $\pi_V \in \Pi_-(\tilde{M}_V)$,
	\item $X \in \mathfrak{a}_{M,V}$ and $f \in \orbI_{\asp}(\tilde{G}_V, \tilde{K}_V) \otimes \mes(G(F_V))$,
	\item $\lambda \in \mathfrak{a}^*_{M, V, \CC}$ in general position, that is, lying off the singular hyperplanes.
\end{itemize}

Iterating the constructions in the local case, we get the semi-local weighted characters
\[ J_{\tilde{M}_V}(\pi_{V, \lambda}, f) = J_{\tilde{M}_V}^{\tilde{G}_V}(\pi_{V, \lambda}, f), \quad f \in C^\infty_{c, \asp}(\tilde{G}_V, \tilde{K}_V) \otimes \mes(G(F_V)), \]
meromorphic in $\lambda$ in general position, from which we construct
\begin{gather*}
	J_{\tilde{M}_V}\left( \pi_V, \lambda, X, f \right) = J^{\tilde{G}_V}_{\tilde{M}_V}\left( \pi_V, \lambda, X, f \right), \\
	I_{\tilde{M}_V}\left( \pi_V, \lambda, X, f \right) = I^{\tilde{G}_V}_{\tilde{M}_V}\left( \pi_V, \lambda, X, f \right).
\end{gather*}
\index{IMVpiX-lambda@$I_{\tilde{M}_V}(\pi_V, \lambda, X, f)$}

As before, they are locally constant in $\lambda$, concentrated at $H_{\tilde{G}}(\cdot) = X_G \in \mathfrak{a}_{G,V}$, and extend to $f$ in $C^\infty_{\asp, \mathrm{ac}}(\tilde{G}_V, \tilde{K}_V) \otimes \mes(G(F_V))$ and $\orbI_{\asp, \mathrm{ac}}(\tilde{G}_V, \tilde{K}_V) \otimes \mes(G(F_V))$, respectively.

The following splitting formula is the spectral analogue of Proposition \ref{prop:semi-local-orbint-splitting}. The only difference comes from the parameter $X$.

\begin{proposition}\label{prop:IM-spec-splitting}
	Let $\pi_V = \bigotimes_v \pi_v \in \Pi_-(\tilde{M}_V)$ and $f = \prod_v f_v \in \orbI_{\asp}(\tilde{G}_V, \tilde{K}_V) \otimes \mes(G(F_V))$. For all $X \in \mathfrak{a}_{M, V}$ and $\lambda \in \mathfrak{a}^*_{M, V, \CC}$ off singularities, let
	\[ \mathcal{A}_X := \left\{ (X_v)_v \in \prod_{v \in V} \mathfrak{a}_{M, F_v} : \sum_{v \in V} X_{v, M} = X \right\}. \]
	Define $M^V$ so that $M^v = M$ for all $v \in V$. Then we have
	\begin{equation*}
		I_{\tilde{M}_V}(\pi_V, \lambda, X, f) = \sum_{L^V \in \mathcal{L}(M^V)} d^G_{M^V}\left( M, L^V \right)
		\int_{\mathcal{A}_X} \prod_{v \in V} I^{\tilde{L}^v}_{\tilde{M}^v}(\pi_v, \lambda, X_v, f_{v, \tilde{L}^v}) \dd (X_v)_v
	\end{equation*}
	for an appropriate measure on $\mathcal{A}_X$, where we used the maps $\mathfrak{a}^*_{M, V, \CC} \to \mathfrak{a}^*_{M, F_v, \CC}$.
\end{proposition}
\begin{proof}
	Same as \cite[p.1198]{MW16-2}.
\end{proof}

For coverings, the following result can be found in \cite[\S 4.1]{Li14b}.

\begin{definition-proposition}
	Let $\pi_V \in \Pi_{\mathrm{unit}, -}(\tilde{M}_V)$. Then
	\[ J_{\tilde{M}_V}(\pi_V, \lambda, \cdot, \cdot), \quad I_{\tilde{M}_V}(\pi_V, \lambda, \cdot, \cdot) \]
	are both defined for $\lambda$ in some neighborhood of $i\mathfrak{a}^*_{M, V}$. In this case, we define
	\begin{equation}\label{eqn:I-0-0-V}
		J_{\tilde{M}_V}(\pi_V, X, f) := J_{\tilde{M}_V}(\pi_V, 0, X, f) , \quad I_{\tilde{M}_V}(\pi_V, X, f) := I_{\tilde{M}_V}(\pi_V, 0, X, f).
	\end{equation}
\end{definition-proposition}

The above is just the semi-local version of \eqref{eqn:I-0-0}. Also note that for unitary $\pi_V$ and $f \in \orbI_{\asp}(\tilde{G}_V, \tilde{K}_V) \otimes \mes(G(F_V))$, one can define
\[ I_{\tilde{M}_V}(\pi_V, f) := \int_{\mathfrak{a}_{M, V}} I_{\tilde{M}_V}(\pi_V, X, f) \dd X. \]

Finally, for a quasisplit group $G^!$, its Levi subgroup $M^!$, and a stable virtual character $\pi_V^!$ of $M^!(F_V)$, we have the stable counterpart
\[ S^{G^!_V}_{M^!_V}(\pi^!_V, \lambda, X, \cdot), \]
We refer to \cite[X.4.3]{MW16-2} for the precise construction and properties.
\index{SpiVX-lambda@$S^{G^{"!}_V}_{M^{"!}_V}(\pi_V^{"!}, \lambda, X, \cdot)$}

\begin{remark}\label{rem:weighted-character-A}
	The invariant distributions $I_{\tilde{M}_V}(\pi_V, \lambda, X, f)$ also appeared in \cite[Hypothèse 4.17]{Li14b}, albeit with a minor difference. In that context, one also fixes a connected Lie subgroup $A \subset A_{G, \infty}$, isomorphic to $\mathfrak{a} \subset A_{G, \infty}$ under $H_G$. Endow $\mathfrak{a}$ (thus $A$) with the Haar measure determined by the chosen positive-definite invariant quadratic form. In defining $J_{\tilde{M}_V}(\pi_V, \lambda, X, \cdot)$ and $I_{\tilde{M}_V}(\pi_V, \lambda, X, \cdot)$, the assumptions become:
	\begin{itemize}
		\item $V$ contains all Archimedean places (so that $\mathfrak{a}_{M, V} = \mathfrak{a}_M$);
		\item $\pi_V \in \Pi_-(\tilde{M}_V /A)$;
		\item $X \in \mathfrak{a}_M/\mathfrak{a}$ and $f^A \in C^\infty_{c, \asp}(\tilde{G}_V /A, \tilde{K}_V) \otimes \mes(G(F_V))$ or with $C^\infty_{c, \asp}$ replaced by $\orbI_{\asp}$;
		\item $\lambda \in \mathfrak{a}^*_{M, \CC}$ is orthogonal to $\mathfrak{a}$ and lies in general position.
	\end{itemize}
	For example, in this setting we put
	\begin{align*}
		J_{\tilde{M}_V}(\pi_V, \lambda, X, f^A) & = J^{\tilde{G}_V}_{\tilde{M}_V}(\pi, \lambda, X, f^A) \\
		& := \int_{\lambda + i\mathfrak{a}^*_M \cap \mathfrak{a}^\perp} J_{\tilde{M}_V}(\pi_\mu, f^A) e^{-\lrangle{\mu, X}} \dd\mu.
	\end{align*}
	In particular, for $M=G$ and $A = A_{G, \infty}$, the integral degenerates to $I^{\tilde{G}_V}(\pi, f^A)$.

	When $\pi_V$ is unitary, $J_{\tilde{M}_V}(\pi_V, \lambda, X, f^A)$ and $I_{\tilde{M}}(\pi_V, \lambda, X, f^A)$ are defined for $\lambda$ in some neighborhood of $i\mathfrak{a}^*_M$. The same constructions applies in the stable context, for all $M^! \subset G^!$ and $A^! \subset A_{G^!, \infty}$.
	
	The study of these distributions can be reduced to the case $A = \{1\}$, and then to the local case. We will return to this issue in \S\ref{sec:Ispec}.
\end{remark}

\begin{remark}\label{rem:weighted-characters-general-G}
	The foregoing constructions of $J^{\tilde{G}}_{\tilde{M}}$, $I^{\tilde{G}}_{\tilde{M}}$, etc.\ apply to all adélic coverings subject to the requirements in \cite{Li14b} (in particular, the Hypothèse 3.12 therein).
\end{remark}

\section{Local spectral stabilization}\label{sec:local-spec}
We continue to work in the local setting of \S\ref{sec:weighted-characters}. We will also use the maps ${}^c \phi_{\tilde{M}}$, ${}^c \theta_{\tilde{M}}$ defined in \S\ref{sec:cphi-nonarch} and \S\ref{sec:cphi-arch}; they intervened in studying the compactly supported version of invariant weighted orbital integrals ${}^c I_{\tilde{M}}(\tilde{\gamma}, \cdot)$. The focus here is on the spectral aspect.

Let $\pi$ be a genuine irreducible representation of $\tilde{M}$, or more generally a genuine virtual character. The genuine invariant distributions ${}^c I_{\tilde{M}}(\pi, \lambda, \cdot) = {}^c I^{\tilde{G}}_{\tilde{M}}(\pi, \lambda, \cdot)$ are characterized by
\[ \sum_{L \in \mathcal{L}(M)} {}^c I^{\tilde{L}}_{\tilde{M}}(\pi, \lambda, {}^c \phi_{\tilde{L}}(f)) = J_{\tilde{M}}(\pi_\lambda, f); \]
they are meromorphic in $\lambda$. We also have
\begin{align*}
	{}^c I_{\tilde{M}}(\pi, \lambda, X, f) & = {}^c I^{\tilde{G}}_{\tilde{M}}(\pi, \lambda, X, f) \\
	& = \int_{\lambda + i\mathfrak{a}^*_{M, F}} e^{-\lrangle{\mu, X}} \cdot {}^c I_{\tilde{M}}(\pi, \mu, f) \dd\mu ,
\end{align*}
which is locally constant in $\lambda$ off the singularities.
\index{IMpiX-lambda-c@${}^c I_{\tilde{M}}(\pi, \lambda, X, f)$}

Recall that ${}^c \theta_{\tilde{M}} = {}^c \theta^{\tilde{G}}_{\tilde{M}}: \orbI_{\asp}(\tilde{G}) \otimes \mes(G(F)) \to \orbI_{\asp, \mathrm{ac}}(\tilde{M}) \otimes \mes(M(F))$ satisfies
\[ \phi_{\tilde{M}} = \sum_{L \in \mathcal{L}(M)} {}^c \theta^{\tilde{L}}_{\tilde{M}} \circ {}^c \phi_{\tilde{L}}. \]
In Propositions \ref{prop:cthetaGM-nonarch} and \ref{prop:cthetaGM-arch}, it is also shown that when $\pi$ is tempered, $X \mapsto I^{\tilde{M}}(\pi, X, {}^c \theta_{\tilde{M}}(f))$ is the Fourier transform of a function on $i\mathfrak{a}^*_{M, F}$, denoted reasonably by $\lambda \mapsto I^{\tilde{M}}\left( \pi_\lambda, {}^c \theta_{\tilde{M}}(f)\right)$ hereafter, which extends meromorphically to $\mathfrak{a}^*_{M, F, \CC}$.

\begin{lemma}\label{prop:cI-tempered-theta}
	For all $\pi$, $X$ and $f$, we have
	\[ {}^c I_{\tilde{M}}(\pi, \lambda, X, f) = \sum_{L \in \mathcal{L}(M)} I^{\tilde{L}}_{\tilde{M}}\left(\pi, \lambda, X, {}^c \theta_{\tilde{L}}(f) \right) \]
	for $\lambda$ in general position. If $\pi$ is tempered, then we have
	\[ {}^c I_{\tilde{M}}(\pi, \lambda, f) = I^{\tilde{M}}\left( \pi_\lambda, {}^c \theta_{\tilde{M}}(f)\right). \]
\end{lemma}
\begin{proof}
	Same as \cite[pp.1186--1187]{MW16-2}. It is based on standard properties of the maps ${}^c \theta_{\tilde{M}}$, etc., having nothing to do with special properties of the covering.
\end{proof}

Consider now a quasisplit group $G^!$ with Levi subgroup $M^!$. In \cite[X.4.3]{MW16-2}, one defined the stable counterparts
\[ {}^c S^{G^!}_{M^!}(\pi^!, \lambda, X, \cdot) \]
of the aforementioned distributions, where $\pi^!$ is now a stable virtual representation of $M^!(F)$, $X \in \mathfrak{a}_{M, F}$ and $\lambda \in \mathfrak{a}^*_{M, F, \CC}$ (off the singularities). Their stability is established in \textit{loc.\ cit.} In the next definition, we take $M^!$ (resp.\ $G^! = G^![s]$) arising from $\mathbf{M}^! \in \Endo_{\elli}(\tilde{M})$ (resp.\ $s \in \Endo_{\mathbf{M}^!}(\tilde{G})$).
\index{SpiX-lambda-c@${}^c S^{G^{"!}}_{M^{"!}}(\pi^{"!}, \lambda, X, \cdot)$}

\begin{definition}\label{def:cIEndo-4-var}
	\index{IEndoMpiX-lambda@$I^{\Endo}_{\tilde{M}}(\pi, \lambda, X, f)$}
	\index{IEndoMpiX-lambda-c@${}^c I^{\Endo}_{\tilde{M}}(\pi, \lambda, X, f)$}
	Let $\mathbf{M}^! \in \Endo_{\elli}(\tilde{M})$ and $\pi^!, \lambda, X$ be as above. Identify $\mathfrak{a}_{M, F}$ and $\mathfrak{a}_{M^!, F}$ (resp.\ $\mathfrak{a}^*_{M, F, \CC}$ and $\mathfrak{a}^*_{M^!, F, \CC}$) by endoscopy. Consider a stable virtual character $\pi^!$ of $M^!(F)$. For all $f \in \orbI_{\asp}(\tilde{G}, \tilde{K}) \otimes \mes(G(F))$, put $f^{G^![s]} := \Trans_{\mathbf{G}^![s], \tilde{G}}(f)$ and define
	\begin{align*}
		I^{\Endo}_{\tilde{M}}(\pi^!, \lambda, X, f) & = I^{\tilde{G}, \Endo}_{\tilde{M}}(\pi^!, \lambda, X, f) \\
		& := \sum_{s \in \Endo_{\mathbf{M}^!}(\tilde{G})} i_{M^!}(\tilde{G}, G^![s]) S^{G^![s]}_{M^!}(\pi^![s], \lambda, X, f^{G^![s]}), \\
		{}^c I^{\Endo}_{\tilde{M}}(\pi^!, \lambda, X, f) & = {}^c I^{\tilde{G}, \Endo}_{\tilde{M}}(\pi^!, \lambda, X, f) \\
		& := \sum_{s \in \Endo_{\mathbf{M}^!}(\tilde{G})} i_{M^!}(\tilde{G}, G^![s]) {}^c S^{G^![s]}_{M^!}(\pi^![s], \lambda, X, f^{G^![s]}).
	\end{align*}
\end{definition}

A more consistent notation would be $I^{\Endo}_{\tilde{M}}(\mathbf{M}^!, \pi^!, \lambda, X, f)$ and so on. The only reason for not doing so is typographic.

Recall Definition \ref{def:thetaEndo-nonarch} for ${}^c \theta^{\Endo}_{\tilde{M}}(\cdot)$ and ${}^c \theta^{\Endo}_{\tilde{M}}(\mathbf{M}^!, \cdot)$. The equality ${}^c \theta^{\Endo}_{\tilde{M}} = {}^c \theta_{\tilde{M}}$ is the content of Theorems \ref{prop:cpt-supported-equalities} and \ref{prop:theta-matching-arch}, to be proved in \S\ref{sec:end-stabilization}.

\begin{lemma}\label{prop:cIEndo-theta}
	Assume that the local geometric Theorem \ref{prop:local-geometric} holds when $M$ (resp.\ $G$) is replaced by any $L \in \mathcal{L}(M)$ with $L \neq M$ (resp.\ $L \neq G$). Then
	\begin{multline*}
		{}^c I^{\Endo}_{\tilde{M}}\left(\pi^!, \lambda, X, f \right) = \sum_{L \in \mathcal{L}(M)} I^{\tilde{L}, \Endo}_{\tilde{M}}\left(\pi^!, \lambda, X, {}^c \theta_{\tilde{L}}(f) \right) \\
		+ I^{M^!}\left(\pi^!, \lambda, X, {}^c \theta^{\Endo}_{\tilde{M}}(\mathbf{M}^!, f) - {}^c \theta_{\tilde{M}}(f)^{M^!} \right) ,
	\end{multline*}
	where ${}^c \theta_{\tilde{M}}(f)^{M^!} := \Trans_{\mathbf{M}^!, \tilde{M}} \left( {}^c \theta_{\tilde{M}}(f)\right)$.
\end{lemma}
\begin{proof}
	By the same combinatorial argument as in the first step in proving Proposition \ref{prop:cpt-noncpt}, the left hand side gets transformed into
	\begin{multline*}
		\sum_{L \in \mathcal{L}^G(M)} \sum_{s^L \in \Endo_{\mathbf{M}^!}(\tilde{L})} \sum_{s_L \in \Endo_{\mathbf{L}^![s^L]}(\tilde{G})} i_{M^!}(\tilde{L}, L^![s^L]) i_{L^![s^L]}(\tilde{G}, G^![s_L]) \\
		S^{L^![s^L]}_{M^!}\left( \pi^! [s^L], \lambda, X, {}^c S\theta^{G^![s_L]}_{L^![s^L]}\left( f^{G^![s_L]} \right)[s_L] \right)
	\end{multline*}
	as $z[s_L]$ is central in $L^![s^L]$; see the explanations in the cited proof. This equals
	\begin{multline*}
		\sum_{L \in \mathcal{L}^G(M)} \sum_{s^L \in \Endo_{\mathbf{M}^!}(\tilde{L})} i_{M^!}(\tilde{L}, L^![s^L]) S^{L^![s^L]}_{M^!}\left( \pi^![s^L], \lambda, X, {}^c \theta^{\Endo}_{\tilde{L}}(\mathbf{L}^![s^L], f) \right) \\
		= \sum_{L \in \mathcal{L}^G(M)} I^{\tilde{L}, \Endo}_{\tilde{M}}\left( \pi^!, \lambda, X, {}^c \theta^{\Endo}_{\tilde{L}}(f) \right).
	\end{multline*}
	
	According to Corollary \ref{prop:cpt-supported-equalities-aux} (for non-Archimedean $F$) or Lemma \ref{prop:theta-matching-arch-prep} (for Archimedean $F$), our assumptions imply ${}^c \theta^{\Endo}_{\tilde{L}} = {}^c \theta_{\tilde{L}}$ for all $L \supsetneq M$. Hence the above equals
	\[ \sum_{L \supsetneq M} I^{\tilde{L}, \Endo}_{\tilde{M}}\left(\pi^! , \lambda, X, {}^c \theta_{\tilde{L}}(f) \right) + I^{\tilde{M}, \Endo}_{\tilde{M}}\left(\pi^!, \lambda, X, {}^c \theta^{\Endo}_{\tilde{M}}(f) \right). \]
	The last term is also $I^{M^!}\left(\pi^!, \lambda, X, {}^c \theta^{\Endo}_{\tilde{M}}(\mathbf{M}^!, f)\right)$.
\end{proof}

The behavior of $I^{\Endo}_{\tilde{M}}(\pi^!, \lambda, X, f)$ (resp.\ ${}^c I^{\Endo}_{\tilde{M}}(\pi^!, \lambda, X, f)$) under parabolic induction is described by a descent formula akin to that of $I_{\tilde{M}}(\pi^!, \lambda, X, f)$ (resp.\ ${}^c I_{\tilde{M}}(\pi^!, \lambda, X, f)$); see Proposition \ref{prop:inv-weighted-descent}. As usual, the proofs follow the combinatorial paradigm of Proposition \ref{prop:descent-orbint-Endo}.

\begin{definition}\label{def:virtual-central-character}
	Let $\pi = c_1 \pi_1 + \cdots + c_n \pi_n$ be a virtual representation of $\tilde{M}$, where $\pi_1, \ldots, \pi_n$ are distinct genuine irreducible representations. We say $\pi$ is unitary (resp.\ tempered) if each constituent $\pi_i$ is. We say $\pi$ has central character $\omega$ if each $\pi_i$ has the same central character $\omega$, and say $\pi$ has unitary central character if $\omega$ is unitary.
	
	The same definition applies to stable virtual representations for the endoscopic groups of $\tilde{M}$, as well as to the semi-local setting.
\end{definition}

Consider $\mathbf{M}^! \in \Endo_{\elli}(\tilde{G})$ and a stable virtual representation $\pi^!$ of $M^!(F)$.

\begin{lemma}\label{prop:spectral-trick-aux}
	Inside $SD_{\mathrm{spec}}(M^!) \otimes \mes(M^!(F))^\vee$, one has a decomposition of the form
	\[ \pi^! = \sum_{\substack{R^! \subset M^! \\ \text{Levi} /\text{conj.} }} \left( \sigma_{R^!} \right)^{M^!} \]
	where $\sigma_{R^!} \in SD_{\mathrm{spec}}(R^!) \otimes \mes(R^!(F))^\vee$ are stable virtual representations of $R^!(F)$, such that
	\begin{enumerate}[(i)]
		\item $\sigma_{R^!} = \sigma'_{R^!, \mu_{R^!}}$ for some $\sigma'_{R^!} \in SD_{\mathrm{ell}}(R^!) \otimes \mes(R^!(F))$ and $\mu_{R^!} \in \mathfrak{a}^*_{R^!, F, \CC}$;
		\item $\sigma_{M^!}$ is tempered provided that $\pi^!$ has unitary central character.
	\end{enumerate}
\end{lemma}
\begin{proof}
	Such a decomposition satisfying (i) is furnished by \eqref{eqn:SDspec-Ind} in general. Furthermore, we may decompose each $\sigma'_{R^!}$ using the parameters from $\Phi_{2, \mathrm{bdd}}(R^!)$. Now assume $\pi^!$ has unitary central character $\omega^!$. By regrouping the sum according to central characters, we may assume that the central character of each constituent in $\sigma_{R^!} = \sigma'_{R^!, \mu_{R^!}}$ equals $\omega^!$ after restriction to $Z_{M^!}(F)$. In particular, this implies $\mu_{M^!}$ is imaginary and $\sigma_{M^!} \in SD_{\mathrm{temp}}(M^!(F)) \otimes \mes(M^!(F))$.
\end{proof}

\begin{corollary}\label{prop:spectral-trick}
	Assume inductively that the assertion $I^{\tilde{S}, \Endo}_{\tilde{M}}(\tilde{\gamma}, \cdot) = I^{\tilde{S}}_{\tilde{M}}(\tilde{\gamma}, \cdot)$ in the local geometric Theorem \ref{prop:local-geometric} holds for all $S \in \mathcal{L}(M)$ with $S \neq G$. Let $\pi^!$ be a stable virtual representation of $M^!(F)$, and put $\pi := \trans_{\mathbf{M}^!, \tilde{M}}(\pi^!) \in D_{\mathrm{spec}, -}(\tilde{M}) \otimes \mes(M(F))^\vee$.
	\begin{enumerate}[(i)]
		\item Suppose that we also have $I^{\tilde{G}, \Endo}_{\tilde{L}}(\tilde{\gamma}, \cdot) = I^{\tilde{G}}_{\tilde{L}}(\tilde{\gamma}, \cdot)$ for all $L \in \mathcal{L}(M)$, including $L=M$. Then we have
		\[ I^{\Endo}_{\tilde{M}}\left(\pi^!, \lambda, X, \cdot \right) = I_{\tilde{M}}\left(\pi, \lambda, X, \cdot \right) \]
		for all $X \in \mathfrak{a}_{M,F}$ and $\lambda$ is in general position.
		\item Suppose that the equalities of the form
		\[ I^{\tilde{S}, \Endo}_{\tilde{R}}\left(\sigma^!, \lambda, X, \cdot \right) = I^{\tilde{S}}_{\tilde{R}}\left(\sigma, \lambda, X, \cdot \right) \]
		as in (i) holds whenever $S \in \mathcal{L}^G(M_0)$ is proper, $R \in \mathcal{L}^S(M_0)$, $\mathbf{R}^! \in \Endo_{\elli}(\tilde{R})$ and $\sigma^!$ replaces the role of $\pi^!$, which we may assume by induction. Assume furthermore that $\pi^!$ has unitary central character. Then $I^{\Endo}_{\tilde{M}}\left(\pi^!, \lambda, X, \cdot \right) = I_{\tilde{M}}\left(\pi, \lambda, X, \cdot \right)$ still holds for $\lambda$ in general position, but constrained in a small neighborhood of $ i\mathfrak{a}^*_{M, F}$.
	\end{enumerate}
\end{corollary}
\begin{proof}
	Consider (i). By Corollary \ref{prop:cpt-supported-equalities-aux} (for non-Archimedean $F$) or Lemma \ref{prop:theta-matching-arch-prep} (for Archimedean $F$), the assumptions imply that ${}^c \theta^{\Endo}_{\tilde{M}}(\mathbf{M}^!, f) = {}^c \theta_{\tilde{M}}(f)^{M^!}$. Thus Lemma \ref{prop:cIEndo-theta} implies
	\[ {}^c I^{\Endo}_{\tilde{M}}\left(\pi^!, \lambda, X, f \right) = \sum_{L \in \mathcal{L}(M)} I^{\tilde{L}, \Endo}_{\tilde{M}}\left(\pi^!, \lambda, X, {}^c \theta_{\tilde{L}}(f) \right). \]
	Compare this with Lemma \ref{prop:cI-tempered-theta}:
	\[ {}^c I_{\tilde{M}}\left(\pi, \lambda, X, f \right) = \sum_{L \in \mathcal{L}(M)} I^{\tilde{L}}_{\tilde{M}}\left(\pi, \lambda, X, {}^c \theta_{\tilde{L}}(f) \right). \]
	By induction, we are reduced to proving that ${}^c I^{\Endo}_{\tilde{M}}\left(\pi^!, \lambda, X, f \right) = {}^c I_{\tilde{M}}\left(\pi, \lambda, X, f \right)$. In turn, this reduces to proving that
	\[ \sum_{s \in \Endo_{\elli}(\tilde{M})} i_{M^!}(\tilde{G}, G^![s]) {}^c S^{G^![s]}_{M^!}\left( \pi^![s], \mu, f^{G^![s]} \right) = {}^c I_{\tilde{M}}\left( \pi, \mu, f \right) \]
	as meromorphic families in $\mu \in \mathfrak{a}^*_{M, \CC} \simeq \mathfrak{a}^*_{M^!, \CC}$.
	
	Using Lemma \ref{prop:spectral-trick-aux}, descent formulas and meromorphic continuation, we are reduced to the case that $\pi^!$ (thus $\pi$) is tempered. It is therefore justified to take $\mu$ imaginary.

	Using the second part of Lemma \ref{prop:cI-tempered-theta} as well as its stable analogue (see the last part of the proof in \cite[p.1191]{MW16-2}), this is in turn reduced to the equality
	\[ S^{M^!}\left(\pi^!_\mu, \sum_{s \in \Endo_{\mathbf{M}^!}(\tilde{G})} i_{M^!}(\tilde{G}, G^![s]) {}^c S\theta^{G^![s]}_{M^!}(f^{G^![s]}) [s] \right) = I^{\tilde{M}}\left(\pi_\mu, {}^c \theta_{\tilde{M}}(f) \right); \]
	recall that the left hand side is $S^{M^!}\left(\pi_\mu^!, {}^c \theta^{\Endo}_{\tilde{M}}(\mathbf{M}^!, f) \right)$,
	
	Note that ${}^c \theta_{\tilde{M}}$ and ${}^c S\theta^{G^![s]}_{M^!}$ take value in Schwartz functions; the meaning of $I^{\tilde{M}}(\pi_\mu, \cdot)$ and $S^{M^!}(\pi^!_\mu, \cdot)$ applied to them is explicated before Lemma \ref{prop:cI-tempered-theta}. It remains to use ${}^c \theta^{\Endo}_{\tilde{M}}(\mathbf{M}^!, f) = {}^c \theta_{\tilde{M}}(f)^{M^!}$.

	Consider (ii). Decompose $\pi^!$ by Lemma \ref{prop:spectral-trick-aux} into $\sum_{\substack{R^! \subset M^! \\ /\text{conj.}}} \left( \sigma_{R^!}\right)^{M^!}$. Let us compare $I^{\Endo}_{\tilde{M}}(\cdots)$ and $I_{\tilde{M}}(\cdots)$ for each $R^!$ separately.
	\begin{itemize}
		\item For $R^! = M^!$, we know $\sigma_{M^!}$ is tempered. Thus $I^{\Endo}_{\tilde{M}}(\cdot)$ and $I_{\tilde{M}}(\cdots)$ associated with $\sigma_{M^!}$ and $\lambda$ in a small neighborhood $\mathcal{U}$ of $i\mathfrak{a}^*_{M, F}$ are both zero, unless $M=G$ (see \cite[Lemma 3.1]{Ar88-1}, which is used in \cite[\S 5.8]{Li12b}). In this case the matching becomes trivial for $\lambda$ in $\mathcal{U}$.

		\item For $R^! \neq M^!$, we complete $R^!$ into a data $R \in \mathcal{L}^M(M_0) \smallsetminus \{M\}$ and $\mathbf{R}^! \in \Endo_{\elli}(\tilde{R})$; see Remark \ref{rem:s-situation-variant}. We can then use descent formulas to expand both sides in terms of $I^{\tilde{S}, \Endo}_{\tilde{R}}(\cdots)$ and $I^{\tilde{S}}_{\tilde{R}}(\cdots)$ where $S \subsetneq G$, so that the equalities in the assumptions can be applied.
	\end{itemize}

	This yields the required equality, with the caveat that the terms with $R^! \neq M^!$ might create singularities in $\mathcal{U}$.
\end{proof}

Note that the assumption $I^{\tilde{S},\ \Endo}_{\tilde{R}}(\cdots) = I^{\tilde{S}}_{\tilde{R}}(\cdots)$ in (ii) is independent of $M$. This is crucial for later applications.

We state the unconditional version of (i) as follows. It is a spectral analogue of the local geometric Theorem \ref{prop:local-geometric}, and will be proved in \S\ref{sec:end-stabilization}.

\begin{theorem}\label{prop:local-spectral}
	Let $\mathbf{M}^! \in \Endo_{\elli}(\tilde{M})$ and let $\pi^!$ be a stable virtual representation of $M^!(F)$. Set $\pi := \trans_{\mathbf{M}^!, \tilde{M}}(\pi^!) \in D_{\mathrm{spec}, -}(\tilde{M}) \otimes \mes(M(F))^\vee$. Then
	\[ I^{\Endo}_{\tilde{M}}\left(\pi^!, \lambda, X, \cdot \right) = I_{\tilde{M}}\left(\pi, \lambda, X, \cdot \right) \]
	for all $X \in \mathfrak{a}_{M, F}$ and $\lambda$ in general position.
\end{theorem}

We conclude by considering the semi-local variant. Let $\tilde{G}$, $M \in \mathcal{L}^G(M_0)$ and $\mathbf{M}^! \in \Endo_{\elli}(\tilde{M})$ be as before, but they are now defined over a number field $F$.

\begin{definition}
	Take $V \supset V_{\mathrm{ram}}$. Let $\pi^!_V$ be a stable virtual representation of $M^!(F_V)$. As in Definition \ref{def:cIEndo-4-var}, set
	\begin{equation}\label{eqn:IEndo-semilocal-spectral}\begin{aligned}
		I^{\Endo}_{\tilde{M}_V}(\pi^!_V, \lambda, X, f) & = I^{\tilde{G}_V, \Endo}_{\tilde{M}_V}(\pi^!_V, \lambda, X, f) \\
		& := \sum_{s \in \Endo_{\mathbf{M}^!}(\tilde{G})} i_{M^!}(\tilde{G}, G^![s]) S^{G^![s]}_{M^!}(\pi^!_V[s], \lambda, X, f^{G^![s]}),
	\end{aligned}\end{equation}
	for all $f \in \orbI_{\asp}(\tilde{G}_V, \tilde{K}_V) \otimes \mes(G(F_V))$, $X \in \mathfrak{a}_M$ and $\lambda \in \mathfrak{a}^*_{M, \CC}$ off the singularities. It is locally constant in $\lambda$. The twist $\pi^!_V \mapsto \pi^!_V[s]$ by $z[s] \in M^!(F)$ on stable virtual characters is performed on each place of $v \in V$ simultaneously.
	\index{IEndoMVpiX-lambda@$I^{\Endo}_{\tilde{M}_V}(\pi^{"!}_V, \lambda, X, f)$}
\end{definition}

Now comes the semi-local version of Theorem \ref{prop:local-spectral}. Its verification is also deferred to \S\ref{sec:end-stabilization}.

\begin{theorem}\label{prop:semilocal-spectral}
	For $\mathbf{M}^!$ and $\pi^!_V$ as above, define $\pi_V := \trans_{\mathbf{M}^!, \tilde{M}}(\pi^!_V)$. Identifying $\mathfrak{a}_M$ and $\mathfrak{a}_{M^!}$ via endoscopy, we have
	\[ I_{\tilde{M}_V}(\pi^!_V, \lambda, X, f) = I^{\Endo}_{\tilde{M}_V}(\pi^!_V, \lambda, X, f), \]
	for all $f$, $X$ and $\lambda$ off the singularities so that both sides are well-defined.
\end{theorem}

\section{Unramified parameters and global weighted characters}\label{sec:spec-nr}
Consider a global covering $\rev: \tilde{G} \to G(\A_F)$ of metaplectic type. Take $V \supset V_{\mathrm{ram}}$.

Since $W_v$ is equipped with a self-dual lattice with respect to $\psi_v \circ \lrangle{\cdot|\cdot}$ for each $v \notin V$, it makes sense to talk about Satake parameters $c_v$ (resp.\ $c^V = (c_v)_{v \notin V}$) of genuine spherical irreducible representations of $\tilde{G}_v$ (resp.\ those of $\tilde{G}$ which are unramified outside $V$), with respect to the hyperspecial subgroup $K_v$ (resp.\ $K^V$). One identifies $c_v$ with a semisimple conjugacy class in $\tilde{G}^\vee$. Satake parameters of $\tilde{G}$ outside $V$ can be twisted by $\mathfrak{a}_{G, \CC}^*$, denoted as $c_v \mapsto c_{v, \lambda}$ and $c^V \mapsto c^V_\lambda$.
\index{cV@$c^V$}

We say $c^V$ is \emph{automorphic} if it occurs in the discrete part of the trace formula for $\tilde{G}$; see \eqref{eqn:Ispec-1}. This implies that $c^V$ is induced from the Satake parameter of some constituent of $L^2_{\mathrm{disc}, -}(M(F) \backslash \tilde{M} /A_{M, \infty})$, for some Levi subgroup $M$; in particular, automorphic Satake parameters are unitarizable.

Let $\mathbf{G}^! \in \Endo_{\elli}(\tilde{G})$. For each $v \notin V_{\mathrm{ram}}$, take any hyperspecial subgroup $K^!_v \subset G^!(F_v)$. Recall from \S\ref{sec:LF} that when restricted to spherical functions, $\Trans_{\mathbf{G}^!_v, \tilde{G}_v}$ is realized by the homomorphism
\[ b_{\mathbf{G}^!_v, \tilde{G}_v}: \mathcal{H}_{\asp}(K_v \backslash \tilde{G}_v / K_v ) \to \mathcal{H}(K^!_v \backslash G^!(F_v) / K^!_v). \]
This induces a transfer $c_v \mapsto c_v^!$ of Satake parameters from $G^!(F_v)$ to $\tilde{G}_v$.

For the endoscopic groups of $G$, there is broader notion of \emph{quasi-automorphic} Satake parameters introduced in \cite[p.1176]{MW16-2}; quasi-automorphy is preserved under transfer from an elliptic endoscopic datum. This can also be adapted to $\tilde{G}$ as follows.

\begin{definition}\label{def:quasi-automorphic-Satake}
	\index{quasi-automorphic}
	A Satake parameter $c^V$ of $\tilde{G}^V$ is called \emph{quasi-automorphic} if
	\begin{itemize}
		\item either $c^V$ is automorphic,
		\item or $c^V$ is the transfer of a quasi-automorphic Satake parameter of $G^!$, for some $\mathbf{G}^! \in \Endo_{\elli}(\tilde{G})$.
	\end{itemize}
\end{definition}

Now enter partial $L$-functions. Let $M \in \mathcal{L}(M_0)$; without loss of generality, assume $M$ is standard. Consider the adjoint $Z_{\tilde{M}^\vee}^\circ$-action on $\tilde{\mathfrak{g}}^\vee := \Lie \tilde{G}^\vee$. For all nontrivial character $\alpha$ of $Z_{\tilde{M}^\vee}^\circ$, denote by $\tilde{\mathfrak{g}}^\vee[\alpha]$ the corresponding eigenspace. Suppose that $\tilde{\mathfrak{g}}^\vee[\alpha] \neq \{0\}$. For each place $v \notin V$ and a Satake parameter $c_v$ for $\tilde{G}_v$, we have the local $L$-factor $L(\alpha, c_v, s)$ associated with the adjoint representation of $\tilde{M}^\vee$ on $\tilde{\mathfrak{g}}^\vee[\alpha]$.

\begin{lemma}\label{prop:qauto-r-conv}
	Let $c^V$ be a quasi-automorphic Satake parameter of $\tilde{M}$. For each character $\alpha$ of $Z_{\tilde{M}^\vee}^\circ$ with $\tilde{\mathfrak{g}}^\vee[\alpha] \neq \{0\}$, the partial $L$-function
	\[ L(\alpha, c^V, s) := \prod_{v \notin V} L(\alpha, c_v, s) \]
	is absolutely convergent whenever $\Re(s) \gg 0$. It extends to a meromorphic function on $\CC$.
\end{lemma}
\begin{proof}
	The recipe is exactly the same as \cite{Ar99b}. We may and do assume $\tilde{G} = \Mp(W, \A_F)$ with $\dim W = 2n$. Suppose that $c^V$ is transferred from some $\mathbf{M}^! \in \Endo_{\elli}(\tilde{M})$. Consider the principal endoscopic datum $\mathbf{M}'$ of $\tilde{M}$, associated with a pair of the form $(n, 0)$. For all $v \notin V$, take hyperspecial subgroups $K^{M'}_v \subset M'(F_v)$, $K^{M^!}_v \subset M^!(F_v)$ and let $K^M_v := K_v \cap M(F_v)$, then $b_{\mathbf{M}^!_v, \tilde{M}_v}$ factors into
	\[ \mathcal{H}_{\asp}(K^M_v \backslash \tilde{M}_v / K^M_v ) \rightiso \mathcal{H}(K^{M'}_v \backslash M'_v / K^{M'}_v) \to \mathcal{H}(K^{M^!}_v \backslash M^!(F_v) / K^{M^!}_v) \]
	where the first arrow is $b_{\mathbf{M}^!_v, \tilde{M}_v}$, and the second is the one for the corresponding endoscopic datum of $M'$.
	
	It follows that $c^V$ is also quasi-automorphic for $M'$, and $L(\alpha, c_v, s)$ equals its counterpart for $M'$ for all $v \notin V$. The convergence and meromorphic continuation are thus reduced to the known case of $M'$.
	
	Assume that $c^V$ is automorphic. In this case, the convergence is proved in the same way as in \cite[13.2 Theorem]{Bo79}, namely through unitarity. As for the meromorphic continuation, the result follows from part of the Langlands--Shahidi method \cite[Theorem 8.4]{Gao18} for $\tilde{M} \subset \tilde{G}$.
\end{proof}

\index{rGM-cV@$r^{\tilde{G}}_{\tilde{M}}(c^V)$}
When $c^V$ is an automorphic Satake parameter of $\tilde{M}$, one has defined in \cite[\S 5.4]{Li14b} the meromorphic family
\[ r^{\tilde{G}}_{\tilde{M}}(c^V_\lambda), \quad \lambda \in \mathfrak{a}^*_{M, \CC}. \]
It has been shown in \cite[Lemme 5.17]{Li14b} to be holomorphic when $\Re\lrangle{\lambda, \check{\alpha}} \gg 0$ for all $\alpha \in \Sigma_P^{\mathrm{red}}$, of moderate growth and analytic over $i\mathfrak{a}^*_M$, hence analytic over an open neighborhood of $i\mathfrak{a}^*_M$. Recall from \cite[\S 5.4]{Li14b} that $r^{\tilde{G}}_{\tilde{M}}(c^V_\lambda)$ is defined through unramified normalizing factors for $c^V$. By the metaplectic Gindikin--Karpelevich formula \cite[Theorem 7.10]{Gao18}, these normalizing factors are expressed in terms of $L$-factors à la Langlands--Shahidi. Specifically, $r^{\tilde{G}}_{\tilde{M}}(c^V_\lambda)$ arises from the $(G, M)$-family
\[ r_{\tilde{Q}}(\Lambda, c^V) := r_{\tilde{Q}^- | \tilde{Q}}(c^V)^{-1} r_{\tilde{Q}^- | \tilde{Q}}(c^V_{\Lambda / 2}), \quad Q \in \mathcal{P}(M), \; \Lambda \in \mathfrak{a}_{M, \CC}^* , \]
with $Q^-$ opposite to $Q$ and
\[ r_{\tilde{Q}|\tilde{P}}(c^V_\lambda) = \prod_{\alpha \in \Sigma_{P^{\vee, -}}^{\mathrm{red}} \cap \Sigma_{Q^{\vee}}^{\mathrm{red}}} \frac{L(\alpha, c^V_\lambda , 0)}{L(\alpha, c^V_\lambda , 1)}, \quad P, Q \in \mathcal{P}(M), \]
where $P \leftrightarrow P^\vee$ is the bijection between $\mathcal{P}(M)$ and $\mathcal{P}(\tilde{M}^\vee)$ explained in \cite[p.216]{Ar99}.

The definition of $r^{\tilde{G}}_{\tilde{M}}(c^V_\lambda)$ via $L$-functions extends to the case of quasi-automorphic $c^V$, say transferred from some $\mathbf{M}^! \in \Endo_{\elli}(\tilde{M})$. As in the proof of Lemma \ref{prop:qauto-r-conv}, we take the principal endoscopic group $M'$ (resp.\ $G'$) of $\tilde{M}$ (resp.\ $\tilde{G}$). Using the fact that $\tilde{G}^\vee = (G')^\vee$ and $\tilde{M}^\vee = (M')^\vee$, it follows that $r^{\tilde{G}}_{\tilde{M}}(c^V_\lambda)$ still have the same analytic properties alluded to above. Indeed, this reduces to the case for $M' \subset G'$ addressed in \cite[Lemma 3.2]{Ar02} or \cite[p.1182]{MW16-2}.

Consider now the stable versions. Let $G^!$ be any elliptic endoscopic group of $\tilde{G}$, with Levi subgroup $M^!$. For every quasi-automorphic Satake parameter $c^{V, !}$ of $M^!$, we have the meromorphic family
\[ s^{G^!}_{M^!}(c^{V, !}_\lambda), \quad \lambda \in \mathfrak{a}^*_{M^!, \CC} \]
with similar analytic properties. See \cite[X.4.1]{MW16-2} or \cite[Theorem 5]{Ar99b}.
\index{sGM-cV@$s^{G^{"!}}_{M^{"!}}(c^{V, "!})$}

\begin{definition}
	\index{rGMEndo-cV@$r^{\tilde{G}, \Endo}_{\tilde{M}}(\mathbf{M}^{"!}, c^{V, "!})$}
	Let $\mathbf{M}^! \in \Endo_{\elli}(\tilde{M})$. Let $c^{V, !}$ be a quasi-automorphic Satake parameter of $M^!$. Define
	\[ r^{\tilde{G}, \Endo}_{\tilde{M}}(\mathbf{M}^!, c^{V, !}_\lambda) := \sum_{s \in \Endo_{\mathbf{M}^!}(\tilde{G})} i_{M^!}(\tilde{G}, G^![s]) s^{G^![s]}_{M^!}(c^{V, !}_\lambda). \]
\end{definition}

The stabilization of these meromorphic families is within reach.

\begin{theorem}\label{prop:spectral-fundamental-lemma}
	Let $\tilde{M} \subset \tilde{G}$, $\mathbf{M}^!$ and $c^{V, !}$ be as above. Denote by $c^V$ the resulting quasi-automorphic Satake parameter of $\tilde{M}$. Then we have the equality of meromorphic families
	\[ r^{\tilde{G}, \Endo}_{\tilde{M}}(\mathbf{M}^!, c^{V, !}_\lambda) = r^{\tilde{G}}_{\tilde{M}}(c^V_\lambda) \]
	in $\lambda$, where we identify $\mathfrak{a}^*_{M, \CC}$ and $\mathfrak{a}^*_{M^!, \CC}$ via endoscopy.
\end{theorem}
\begin{proof}
	Without loss of generality, one may assume $\tilde{G} = \Sp(W, \A_F)$ with $\dim W = 2n$. Following the paradigm for its local Archimedean version in Theorem \ref{prop:easy-stabilization} and the discussions above, this reduces readily to the case for $\SO(2n+1)$, namely \cite[Theorem 5]{Ar99b}. There is an extra ingredient in the present global setting, namely the meromorphic continuation of partial $L$-functions involved in $r^{\tilde{G}}_{\tilde{M}}$ and $r^{\tilde{G}, \Endo}_{\tilde{M}}$; this is furnished by Lemma \ref{prop:qauto-r-conv}.
\end{proof}

We proceed to study global weighted characters. Fix $V \supset V_{\mathrm{ram}}$ and $M \in \mathcal{L}(M_0)$. The notion of unitarity and central characters of virtual characters will be as in Definition \ref{def:virtual-central-character}, except that we now work with representations of $\tilde{M}_V$. Consider a pair $(\pi_V, c^V)$ subject to one of the following sets of conditions:
\begin{enumerate}[(a)]
	\item The $\pi_V$ is a virtual genuine representation of $\tilde{M}_V$, assumed to be unitary with central character (Definition \ref{def:virtual-central-character}); $c^V$ is an automorphic Satake parameter of $\tilde{M}$ off $V$.
	\item There exists $\mathbf{M}^! \in \Endo_{\elli}(\tilde{M})$ such that $(\pi_V, c^V)$ is the transfer of a similar pair $(\pi^!_V, c^{V, !})$ from $M^!$, assuming that $\pi^!_V$ is stable with unitary central character, but the Satake parameter $c^{V, !}$ is only assumed to be quasi-automorphic.
\end{enumerate}

For every $L \in \mathcal{L}(M)$, denote by $\pi_V^{\tilde{L}}$ the normalized parabolic induction of $\pi_V$ to $\tilde{L}_V$.

\begin{definition}\label{def:global-weighted-character}
	\index{IMglob@$I^{\mathrm{glob}}_{\tilde{M}}$, $I^{\mathrm{glob}, \Endo}_{\tilde{M}}$}
	Given $(\pi_V, c^V)$ and $f \in \orbI_{\asp}(\tilde{G}_V, \tilde{K}_V) \otimes \mes(G(F_V))$, we impose the following conditions are $\pi_V$:
	\begin{itemize}
		\item $I_{\tilde{L}_V}\left( (\pi_{V, \lambda})^{\tilde{L}} , \nu, X, f \right)$ is defined for $\nu$ in a neighborhood of $i\mathfrak{a}_{L, V}^*$, whenever $\lambda \in i\mathfrak{a}^{L, *}_M$ and $X \in \mathfrak{a}_{L, V}$;
		\item $\lambda \mapsto I_{\tilde{L}_V}\left( (\pi_{V, \lambda})^{\tilde{L}} , 0, 0, f \right)$ is of rapid decay on $i\mathfrak{a}^{L, *}_M$.
	\end{itemize}
	For instance, they are satisfied in the case (a) above. The (invariant) \emph{global weighted character} is then defined as
	\begin{align*}
		I^{\mathrm{glob}}_{\tilde{M}}\left( \pi_V \otimes c^V, f \right) & = I^{\tilde{G}, \mathrm{glob}}_{\tilde{M}}\left( \pi_V \otimes c^V, f \right) \\
		& := \sum_{L \in \mathcal{L}(M)} \int_{i\mathfrak{a}^{L, *}_M} \dd \lambda \cdot r^{\tilde{L}}_{\tilde{M}}(c^V_\lambda) I_{\tilde{L}_V}\left( (\pi_{V, \lambda})^{\tilde{L}} , 0, 0, f \right).
	\end{align*}
\end{definition}

The definition above follows \cite[X.4.9]{MW16-2}. We have seen that $\lambda \mapsto r^{\tilde{G}}_{\tilde{M}}(c^V_\lambda)$ is analytic and of moderate growth on $i\mathfrak{a}^*_M$, hence the integral converges.

The stable avatars of the $I^{\mathrm{glob}}_{\tilde{M}}$ have been defined in \cite[p.1201]{MW16-2}; here we denote them as
\[ S^{G^!, \mathrm{glob}}_{M^!}\left(\pi_V^! \otimes c^{V, !}, f^!\right) := \sum_{L^! \in \mathcal{L}^{G^!}(M^!)} \int_{i\mathfrak{a}^{L^!, *}_{M^!}} \dd \lambda \cdot s^{L^!}_{M^!}(c^{V, !}_\lambda) S^{G^!}_{L^!_V}\left( (\pi^!_{V, \lambda})^{L^!} , 0, 0, f^! \right), \]
under similar assumptions on $S^{G^!_V}_{L^!_V}\left( \left(\pi^!_{V, \lambda}\right)^{L^!}, \nu, X, f^! \right)$. This leads to the following endoscopic version.
\index{SMglob@$S^{G^{"!}, \mathrm{glob}}_{M^{"!}}$}

\begin{definition}\label{def:global-weighted-character-Endo}
	Assume that the assertion $I^{\tilde{S}, \Endo}_{\tilde{L}}(\tilde{\gamma}, \cdot) = I^{\tilde{S}}_{\tilde{L}}(\tilde{\gamma}, \cdot)$ of the local geometric Theorem \ref{prop:local-geometric} holds whenever $S \subsetneq G$, which is covered by induction. Suppose that we are in the case (b) before, i.e.\ $(\pi_V, c^V)$ is the transfer of $(\pi^!_V, c^{V, !})$ from some $\mathbf{M}^! \in \Endo_{\elli}(\tilde{M})$, which we fix. We also assume that
	\[ S^{G^![s], \mathrm{glob}}_{M^!}\left( \pi_V^![s] \otimes c^{V, !}, f^{G^![s]} \right) \]
	is defined for all $s \in \Endo_{\mathbf{M}^!}(\tilde{G})$. Set
	\begin{align*}
		I^{\mathrm{glob}, \Endo}_{\tilde{M}}\left(\mathbf{M}^!, \pi_V^! \otimes c^{V, !}, f \right) & = I^{\tilde{G}, \mathrm{glob}, \Endo}_{\tilde{M}}\left(\mathbf{M}^!, \pi_V^! \otimes c^{V, !}, f \right) \\
		& :=\sum_{s \in \Endo_{\mathbf{M}^!}(\tilde{G})} i_{M^!}(\tilde{G}, G^![s]) S^{G^![s], \mathrm{glob}}_{M^!}\left( \pi_V^![s] \otimes c^{V, !}, f^{G^![s]} \right),
	\end{align*}
	where $f \in \orbI_{\asp}(\tilde{G}_V, \tilde{K}_V) \otimes \mes(G(F_V))$ and $f^{G^![s]} := \Trans_{\mathbf{G}^![s], \tilde{G}}(f)$.
\end{definition}

We are ready to stabilize the global weighted characters.

\begin{proposition}[Cf.\ {\cite[X.4.9 Proposition]{MW16-2}}]\label{prop:stabilization-global-weighted-character}
	Keep the assumptions of Definition \ref{def:global-weighted-character-Endo}, and assume that all matching of distribution holds when $G$ is replaced by a proper Levi subgroup, as covered by induction. Then $I^{\mathrm{glob}}_{\tilde{M}}(\pi_V \otimes c^V, \cdot)$ is also defined, and
	\[ I^{\mathrm{glob}, \Endo}_{\tilde{M}}(\mathbf{M}^!, \pi_V^! \otimes c^{V, !}, \cdot) = I^{\mathrm{glob}}_{\tilde{M}}(\pi_V \otimes c^V, \cdot). \]
\end{proposition}
\begin{proof}
	Fix $f \in \orbI_{\asp}(\tilde{G}_V, \tilde{K}_V) \otimes \mes(G(F_V))$. Apply the definitions to obtain
	\begin{multline*}
		I^{\mathrm{glob}, \Endo}_{\tilde{M}}\left(\mathbf{M}^!, \pi_V^! \otimes c^{V, !}, f \right) = \sum_{s \in \Endo_{\mathbf{M}^!}(\tilde{G}) } i_{M^!}(\tilde{G}, G^![s]) \sum_{L^! \in \mathcal{L}^{G^![s]}(M^!)} \\
		\int_{i\mathfrak{a}^{L^!, *}_{M^!}} \dd\lambda \cdot s^{L^!}_{M^!}(c^{V, !}_\lambda) S^{G^![s]_V}_{L^!_V}\left( (\pi^!_{V, \lambda})[s]^{L^!}, 0, 0, f^{G^![s]} \right).
	\end{multline*}

	Apply the combinatorial Lemma \ref{prop:sL-Ls} with $R=M$ to transform the sum over $(s, L^!)$ to a sum over $(L, s)$, and decompose $s$ into $(s^L, s_L)$ where $s^L \in \Endo_{\mathbf{M}^!}(\tilde{L})$ and $s_L \in \Endo_{\mathbf{L}^![s^L]}(\tilde{G})$. Also decompose $i_{M^!}(\tilde{G}, G^![s])$ by Lemma \ref{prop:i-transitivity}, from which we obtain
	\begin{equation*}
		I^{\mathrm{glob}, \Endo}_{\tilde{M}}\left(\mathbf{M}^!, \pi_V^! \otimes c^{V, !}, f \right) = \sum_{L \in \mathcal{L}(M)} \int_{i\mathfrak{a}^{L, *}_M} \sum_{s^L} i_{M^!}(\tilde{L}, L^![s^L]) s^{L^![s^L]}_{M^!}(c^{V, !}_\lambda) E(\lambda, s^L) \dd \lambda ,
	\end{equation*}
	with
	\begin{align*}
		E(\lambda, s^L) & := \sum_{s_L} i_{L^![s^L]}(\tilde{G}, G^![s_L]) S^{G^![s_L]_V}_{L^![s^L]_V}\left( (\pi^!_{V, \lambda}[s^L])^{L^![s^L]}[s_L], 0, 0, f^{G^![s_L]} \right) \\
		& = I^{\Endo}_{\tilde{L}_V}\left(\mathbf{L}^![s^L] (\pi^!_{V, \lambda}[s^L])^{L^![s^L]}, 0, 0, f \right) \quad \text{(see \eqref{eqn:IEndo-semilocal-spectral})}
	\end{align*}
	where we used the observation that $z[s_L]$-translation commutes with $(\cdots)^{L^![s^L]}$, by centrality.

	Recall that the transfer of $(\pi^!_{V, \lambda}[s^L])^{L^![s^L]}$ to $\tilde{L}_V$ equals $(\pi_{V, \lambda})^{\tilde{L}}$. We claim that for $\nu \in \mathfrak{a}^*_{L, \CC}$ close to $i\mathfrak{a}^*_{L, V}$ and in general position, we have
	\begin{equation}\label{eqn:stabilization-semilocal-special}
		I^{\Endo}_{\tilde{L}_V}\left(\mathbf{L}^![s^L], (\pi^!_{V, \lambda}[s^L])^{L^![s^L]}, \nu, X, f \right) = I_{\tilde{L}_V}\left( \pi_{V, \lambda}^{\tilde{L}}, \nu, X, f \right).
	\end{equation}

	To prove this special instance of Theorem \ref{prop:semilocal-spectral}, we reduce to the local setting by using
	\begin{itemize}
		\item the splitting formula for $I_{\tilde{L}_V}(\cdots)$ (Proposition \ref{prop:IM-spec-splitting});
		\item the splitting formula for $I^{\Endo}_{\tilde{L}_V}$ --- the required combinatorial arguments are the same as in the proof of Proposition \ref{prop:semilocal-matching-coherence} (i).
	\end{itemize}
	Both splitting formulas take the same form. Once reduced to the local setting, Corollary \ref{prop:spectral-trick} (ii) is applicable since by assumption,
	\begin{itemize}
		\item the equalities in Corollary \ref{prop:spectral-trick} (i) hold if $G$ is replaced by a proper Levi subgroup, and
		\item $\pi^!_V$ has unitary central character.
	\end{itemize}
	
	Consequently, for $\nu$ as above, $I_{\tilde{L}_V}\left( \pi_{V, \lambda}^{\tilde{L}}, \nu, X, f \right)$ inherits the required properties (non-singular for imaginary $\nu$, rapid decay in $\lambda$) from the stable side. In particular, \eqref{eqn:stabilization-semilocal-special} extends to $(\nu, X) = (0, 0)$. Hence $E(\lambda, s^L) = I_{\tilde{L}_V}\left( \pi_{V, \lambda}^{\tilde{L}}, 0, 0, f\right)$.

	On the other hand, the sum of $i_{M^!}(\tilde{L}, L^![s^L]) s^{L^![s^L]}_{M^!}(c^{V, !}_\lambda)$ over $s^L$ equals $r^{\tilde{L}}_{\tilde{M}}(c^V_\lambda)$, by Theorem \ref{prop:spectral-fundamental-lemma}. All in all, we conclude that $I^{\mathrm{glob}}_{\tilde{M}}(\pi_V \otimes c^V, \cdot)$ is defined, and $I^{\mathrm{glob}, \Endo}_{\tilde{M}}\left(\pi_V^! \otimes c^{V, !}, f \right) = I^{\mathrm{glob}}_{\tilde{M}}(\pi_V \otimes c^V, f)$.
\end{proof}

\section{Digression on endoscopic coverings}\label{sec:endo-covering}
The following digressions are inspired by a comment from Fan Gao. It is not used elsewhere in this work.

Take $G = \Sp(W)$ with $\dim W = 2n$. Recall that $\tilde{G}^{(2)}$ denotes the twofold metaplectic covering, a Brylinski--Deligne covering of $G(\A_F)$. For any Levi subgroup $M$ of $G$, define $\tilde{M}^{(2)}$ to be the preimage of $M(\A_F)$ in $\tilde{G}^{(2)}$.

The proof of Lemma \ref{prop:qauto-r-conv} relies critically on the meromorphy of $L$-functions of Langlands--Shahidi type for $\tilde{G}$, or equivalently for $\tilde{G}^{(2)}$, in the case when $c^V$ is automorphic. In order to adapt Shahidi's proof to this case, Gao constructed certain \emph{endoscopic coverings} of $\tilde{G}^{(2)}$ in \cite[\S 8.3]{Gao18}, which differ from our definitions in \S\ref{sec:endoscopic-data}. In fact, he treated Brylinski--Deligne coverings of split groups in general. In order to reconcile these points of view, let us recall Gao's construction briefly.

Consider a standard parabolic subgroup $P = MU \subset G$ and an automorphic Satake parameter $c^V$ of $\tilde{M}$, or equivalently of $\tilde{M}^{(2)}$. It is routine to reduce the proof of meromorphy of $L(\alpha, c^V, s)$ to the case of maximal $P$. Let $\tilde{P}^\vee = \tilde{M}^\vee U^\vee$ be the corresponding standard parabolic of $\tilde{G}^\vee$. There is a unique reduced positive root $\alpha_1$ of $(Z_{\tilde{M}^\vee}^\circ, \tilde{G}^\vee)$, and all the positive roots are
\[ k\alpha_1, \quad k = 1, \ldots, r. \]
Let $\Ad_k$ be the adjoint representation of $\tilde{M}^\vee$ on $\tilde{\mathfrak{g}}^\vee[k\alpha_1]$. Then the adjoint action of $\tilde{M}^\vee$ on $\mathfrak{u}^\vee$ decomposes into $\bigoplus_{k=1}^r \Ad_k$.

Fix $1 \leq k \leq r$. In \textit{loc.\ cit.} is constructed a subgroup $\tilde{G}_k^\vee \subset \tilde{G}^\vee$ together with a maximal standard parabolic $\tilde{P}_k^\vee = \tilde{M}^\vee U_k^\vee$. It has the key property that $\Ad_k$ equals the representation $\Ad_1$ associated with $\tilde{P}_k^\vee \subset \tilde{G}^\vee_k$. 
	
Specifically, $\tilde{G}_k^\vee$ is the centralizer of some $s_k \in Z_{\tilde{M}^\vee}^\circ$. This element represents a generator of $\Ker(k\alpha_1) / \Ker(\alpha_1) \simeq \bmu_k$. We claim that $s_k$ is elliptic in $\tilde{G}^\vee$. Indeed, we have
\[ \dim Z_{\tilde{M}^\vee}^\circ - 1 = \dim Z_{\tilde{G}^\vee}^\circ \leq \dim Z_{\tilde{G}_k^\vee}^\circ \leq \dim Z_{\tilde{M}^\vee}^\circ . \]
If $s_k$ is not elliptic, the first inequality must be strict, and this would force $\tilde{M}^\vee = \tilde{G}_k^\vee$; on the other hand, $k\alpha_1$ is a root of $(\tilde{G}_k^\vee, Z_{\tilde{M}^\vee}^\circ)$ (see \cite[p.1141]{Ar99}). Contradiction.

It follows from the claim that $s_k^2 = 1$, and
\[ \tilde{G}^\vee_k = \Sp(2n', \CC) \times \Sp(2n'', \CC), \quad n' + n'' = n. \]

Let us describe $\tilde{G}_k^\vee$ in terms of a sub-root datum of that of $\tilde{G}^\vee$, say with the same lattices of characters and cocharacters, but fewer roots and coroots. By rescaling the (co-)roots and passing to the dual, one obtains $M \subset G_k \subset G$. Now the twofold Brylinski--Deligne covering $\tilde{G}_k$ is built in two steps.
\begin{itemize}
	\item First, take an ``incarnation'' of $\tilde{G}$ by a pair $(D, \eta)$, where $D$ is an integral bilinear form on $X_*$ and $\eta$ is homomorphism from the coroot lattice of $G$ to $\R^\times$. We refer to \cite[\S 2.6]{GG18} for an overview of incarnations.
	\item Secondly, take $\tilde{G}_k$ to be the Brylinski--Deligne covering of $G_k(\A_F)$ incarnated by $(D, \eta_k)$, where $\eta_k$ is the restriction of $\eta$ to the coroot lattice of $G_k$. It is shown in \cite[\S 8.3]{Gao18} that $\tilde{G}_k$ restricts to $\tilde{M}^{(2)}$.
\end{itemize}

The standard choice of $(D, \eta)$ is given in the proof of Lemma \ref{prop:Sp4-bisector}: take $\eta = \mathbf{1}$ and $D = \sum_{i=1}^n \epsilon^*_i \otimes \epsilon^*_i$ in the usual basis. We may also assume that $\epsilon^*_1, \ldots, \epsilon^*_{n'}$ (resp.\ $\epsilon^*_{n'+1}, \ldots, \epsilon^*_n$) lie in the $\Sp(2n')$ (resp.\ $\Sp(2n'')$) factor of $\tilde{G}^\vee$. Thus $(D, \eta_k)$ breaks into two factors.

With the shorthand $\Mp(2n)^{(2)}$ for $\Mp(2n, \A_F)^{(2)}$, it is now clear that $\tilde{G}_k$ is nothing but
\[ \Mp(2n')^{(2)} \utimes{\bmu_2} \Mp(2n'')^{(2)} := \Mp(2n')^{(2)} \times \Mp(2n'')^{(2)} \big/ \{ (z, z^{-1}) : z \in \bmu_2 \}. \]
Clearly, $\tilde{G}_1 = \tilde{G}^{(2)}$, and $\dim G_k < \dim G$ if $k > 1$. If we replace $\Mp(2n')^{(2)}$ (resp.\ $\Mp(2n'')^{(2)}$) by $\SO(2n'+1)$ (resp.\ $\SO(2n''+1)$) by the familiar philosophy, we revert to the elliptic endoscopic groups defined in \S\ref{sec:endoscopic-data}.

\begin{remark}
	Let $c^V$ be an automorphic Satake parameter of $\tilde{M}$. As in Shahidi's case, the meromorphy of $L(\alpha, c^V, s)$ follows easily from the construction of $\tilde{M}^{(2)} \subset \tilde{G}_k \subset \tilde{G}^{(2)}$ for various $1 \leq k \leq r$: the argument is based on the meromorphy of constant terms of Eisenstein series and induction on dimensions. Although the endoscopic coverings appeared, the proof does not require any advanced property of endoscopy, such as the correspondence of semisimple classes or transfer.
\end{remark}

Note that the construction of $\tilde{G}_k$ works over local fields as well.

\begin{remark}
	The coverings $\Mp(2n')^{(2)} \utimes{\bmu_2} \Mp(2n'')^{(2)}$ are also the elliptic endoscopic groups of $\Mp(2n)^{(2)}$ advocated by D.\ Renard in \cite{Re99}.
\end{remark}

\section{Spectral side of the trace formula}\label{sec:Ispec}
We continue to work with a covering of metaplectic type $\rev: \tilde{G} \to G(\A_F)$. Denote by $\tilde{G}_\infty$ the preimage of $\prod_{v \mid \infty} G(F_v)$ under $\rev$. It also decomposes into $\tilde{G}_\infty^1 \times A_{G, \infty}$.

\subsection{The invariant side}
As in \cite[p.1040]{Li13} or \cite[Sect 4]{Ar88-2}, fix a Cartan subalgebra $\mathfrak{h} = i\mathfrak{h}_K \oplus \mathfrak{h}_0$ of the split real form of $\mathfrak{g}_{\infty, \CC}$; define its complexification and absolute Weyl group
\index{hC@$\mathfrak{h}_{\CC}$}
\index{Winfty@$W_\infty$}
\[ \mathfrak{h}_{\CC} \subset \mathfrak{g}_{\infty, \CC}, \quad W_\infty := W(\mathfrak{g}_{\infty, \CC}, \mathfrak{h}_{\CC}). \]
Note that $W_\infty$ acts on $\mathfrak{h}$. Fix a $W_\infty$-invariant Euclidean norm $\|\cdot\|$ on $\mathfrak{h}$, and denote the dual Hermitian norm on $\mathfrak{h}^*_{\CC}$ by $\|\cdot\|$. Using these data, one can talk about the height $\|\Im\nu\|$ for any $\nu \in \mathfrak{h}^*_{\CC} / W_\infty$. The same $\|\cdot\|$ will be used for all $\tilde{M}$, where $M \in \mathcal{L}(M_0)$.

Let $V \supset V_{\mathrm{ram}}$. In \cite[Théorème 6.1]{Li14b}, the spectral side of the invariant trace formula is the invariant genuine distribution $I^{\tilde{G}}_{\mathrm{spec}} = I_{\mathrm{spec}}$ coming with the expansion
\begin{equation}\label{eqn:Ispec-0}\begin{aligned}
	I_{\mathrm{spec}} & = \sum_{t \geq 0} I_t , \\
	I_t & = \sum_{M \in \mathcal{L}(M_0)} \frac{|W^M_0|}{|W^G_0|} \int_{\Pi_{t, -}(\tilde{M}^1, V)} a^{\tilde{M}}(\pi) I_{\tilde{M}_V}(\pi, 0, \cdot) \dd\pi \\
	& = \sum_{\nu: \|\Im(\nu)\| = t} I_\nu .
\end{aligned}\end{equation}
\index{Ispec@$I_{\mathrm{spec}}$, $I_t$, $I_\nu$}
Here
\begin{itemize}
	\item $\Pi_{t, -}(\tilde{M}^1, V)$ is a space of irreducible unitary genuine representations $\pi$ of $\tilde{M}_V^1$, or equivalently those of $\tilde{M}_V/ A_{G, \infty}$;
	\item $\nu$ specifies the infinitesimal characters of such representations $\pi$;
	\item $a^{\tilde{M}}(\pi)$ are certain coefficients and $\dd \pi$ is a suitable Radon measure.
\end{itemize}

In \textit{loc.\ cit.}, these distributions are evaluated on test functions
\[ f^1 \in \orbI_{\asp}(\tilde{G}_V/A_{G, \infty} ,\tilde{K}_V) \otimes \mes(G(F_V)). \]
For details about these objects, see \cite[\S 5.5]{Li14b}. We remind the reader that the $I_{\tilde{M}_V}(\pi, 0, \cdot)$ here is the $I_{\tilde{M}_V}(\pi, 0, 0, \cdot)$ in \S\ref{sec:weighted-characters}, with $A := A_{G, \infty}$; see \eqref{eqn:I-0-0-V} and Remark \ref{rem:weighted-character-A}.

\begin{remark}
	The convergence of $I_{\mathrm{spec}} = \sum_{t \geq 0} \sum_{M \in \mathcal{L}(M_0)} \int_{\Pi_{t, -}(\tilde{M}^1, V)}$ is understood as an iterated integral. We do not claim the convergence as a multiple integral, since the results of Finis--Lapid--Müller \cite{FLM11} have not been adapted to metaplectic coverings yet, although such a result is expected.
\end{remark}

In parallel with the geometric side, we shall rephrase $I_{\mathrm{spec}}$, $I_t$ and $I_\nu$ as linear functionals on $\orbI_{\asp}(\tilde{G}_V, \tilde{K}_V) \otimes \mes(G(F_V))$. As explained in \S\ref{sec:geom-side}, this is essentially given by \cite[Théorème 6.4]{Li14b}: it suffices to notice that $I_{\tilde{M}_V}(\pi, 0, \cdot)$ is concentrated at $0 \in \mathfrak{a}_G$, as seen in \cite[\S 4.5]{Li14b}.

\begin{remark}\label{rem:cpt-supp-test-fcn-spec}
	In what follows, we will only evaluate $I_{\mathrm{spec}}$, $I_t$ and $I_\nu$ on
	\[ f \in \orbI_{\asp}(\tilde{G}_V, \tilde{K}_V) \otimes \mes(G(F_V)); \]
	they depend only on $f|_{\tilde{G}_V^1}$. The $A_{G, \infty}$-invariant test functions $f^1$ will not be used anymore. Cf.\ Remark \ref{rem:cpt-supp-test-fcn-geom}.
	
	The invariant trace formula $I_{\mathrm{geom}} = I_{\mathrm{spec}}$ and their simple forms in \cite[\S 6.3]{Li14b} hold for such test functions $f$. They remain valid for $f \in \orbI_{\mathrm{ac}, \asp}(\tilde{G}_V, \tilde{K}_V) \otimes \mes(G(F_V))$, since what matters is only $f|_{\tilde{G}^1_V}$.
\end{remark}

\subsection{Discrete part}
The genuine invariant distributions $I_t$ are \emph{compressed} --- they live on $\tilde{G}_V$. The \emph{pre-compressed} distribution in \cite[pp.1060--1062]{Li13} or \cite[(5.10)]{Li14b} is
\begin{equation}\label{eqn:Ispec-1}\begin{aligned}
	I_{\mathrm{disc}, t}(\dot{f}) = I^{\tilde{G}}_{\mathrm{disc}, t}(\dot{f}) & := \sum_{M \in \mathcal{L}(M_0)} \sum_{s \in W^G(M)_{\mathrm{reg}}} \frac{|W^M_0|}{|W^G_0|} \\
	& \quad \left| \det(1-s | \mathfrak{a}^G_M) \right|^{-1} \Tr\left( M_{P|P}(s, 0) \mathcal{I}_{\tilde{P}, \mathrm{disc}, t}(0, \dot{f}) \right) \\
	& = \sum_{\pi \in \Pi_{-, t}(\tilde{G}^1)} a^{\tilde{G}}_{\mathrm{disc}}(\pi) \Tr \pi(\dot{f}|_{\tilde{G}^1}) \\
	& = \sum_{\pi \in \Pi_{-, t}(\tilde{G}^1)} a^{\tilde{G}}_{\mathrm{disc}}(\pi) \int_{i\mathfrak{a}^*_G} \Tr \pi_\mu(\dot{f}) \dd\mu
\end{aligned}\end{equation}
where $t \geq 0$, $P \in \mathcal{P}(M)$ and $\dot{f} \in \orbI_{\asp}(\tilde{G}, \tilde{K}) \otimes \mes(G(\A_F))$. In \textit{loc.\ cit.}, a variant $J_{\mathrm{disc}, t}(\dot{f}^1)$ was considered, where $A_{G, \infty}$-invariance was imposed on $\dot{f}^1$.
\index{Idisc@$I_{\mathrm{disc}, t}$}
\index{adisc@$a^{\tilde{G}}_{\mathrm{disc}}$}

The distribution $I_{\mathrm{disc}, t}$ is of an adélic nature. It is expressed as an infinite linear combination of irreducible genuine characters $\Theta_{\pi}$ of $\tilde{G}^1$. However, the evaluation on each $\dot{f}$ is well-defined: the results in \cite[\S 7]{Li13} guarantee that $I_{\mathrm{disc}, t}$ is \emph{semi-finite} in the sense of \cite[X.5.3]{MW16-2}, recapitulated as follows.

\begin{definition}\label{def:semi-finite-distribution}
	\index{semi-finite}
	\index{nu-pi@$\nu_{\pi}$}
	Consider a formal linear combination $\rho = \sum_{\pi} c_{\pi} \pi$ of irreducible admissible genuine representations of $\tilde{G}$; call $\pi$ a constituent of $\rho$ if $c_{\pi} \neq 0$, and denote by $\nu_{\pi}$ its infinitesimal character at Archimedean places. We say $\rho$ is \emph{semi-finite} if the following are satisfied. 
	\begin{enumerate}[(i)]
		\item Let $\pi$ be any constituent, then its complex-conjugate $\overline{\pi}$ and contragredient $\check{\pi}$ have the same infinitesimal character. This is automatically the case if $\pi$ is unitary.
	
		\item Given a finite set $\Gamma$ of $\tilde{K}_\infty$-types, there exists $r = r(\Gamma)$ such that for all constituent $\pi$ with some $\tilde{K}_\infty$-type in $\Gamma$, we have $\|\Re(\nu_{\pi})\| \leq r$.
	
		\item Given $R \geq 0$, a finite set $V \supset V_{\mathrm{ram}}$ of places, and a compact open subgroup $\mathcal{K}_V$ of $\prod_{\substack{v \in V \\ v \nmid \infty}} \tilde{G}_v$, consider the constituents $\pi$ that are unramified outside $V$, with a specified $\tilde{K}_\infty \times \mathcal{K}_V$-type, and $\|\Im \nu_{\pi}\| \leq R$.	Then there are only finitely many infinitesimal characters $\nu$ and Satake parameters $c^V$ that can intervene in these constituents.
	
		\item For each compact open subgroup $\mathcal{K}$ of $\prod'_{v \nmid \infty} \tilde{G}_v$, there are at most finitely many constituents $\pi$ with specified infinitesimal characters and $\mathcal{K}$-types.
	
		\item Given a finite set $\Upsilon$ of $\tilde{K}_\infty \times \mathcal{K}$-types and an infinitesimal character $\nu$, let $\rho^\Upsilon_\nu$ be the formal sum of constituents of $\rho$ that contains some $\tilde{K}_\infty \times \mathcal{K}$-type from $\Upsilon$ and has infinitesimal character $\nu$; it is actually a virtual character by (iv). Then $\sum_\nu \rho^\Upsilon_\nu(\dot{f})$ converges absolutely for every $\dot{f}$.
	\end{enumerate}
\end{definition}

In fact, the properties satisfied by $I_{\mathrm{disc}, t}$ are much stronger than those considered in \cite[X.5.3]{MW16-2}, since the height $t$ is specified here; see \cite[\S 7]{Li13}.

\begin{notation}\label{nota:semi-finite-decomposition}
	Given a semi-finite $\rho$, we denote by $\rho^{V\text{-nr}}$ the part of $\rho$ that is unramified outside $V$. The general operations in \cite[X.5.3]{MW16-2} on semi-finite distributions furnish the canonical decompositions
	\[ \rho = \sum_\nu \rho_\nu, \quad \rho^{V\text{-nr}} = \sum_\nu \sum_{c^V} \rho_{\nu, c^V} \]
	where $\nu$ ranges over $\mathfrak{h}_{\CC}^* / W_\infty$ and $c^V$ ranges over Satake parameters off $V$; when applied to $\dot{f}$, they yield absolutely convergent sums.
	
	Let $\rho$ be semi-finite. The same token $\rho$ will also be used to denote the following distribution
	\begin{equation}\label{eqn:rho-decomp}
		\rho(f) := \rho(f f_{K^V}), \quad f \in \orbI_{\asp}(\tilde{G}_V, \tilde{K}_V) \otimes \mes(G(F_V))
	\end{equation}
	on $\tilde{G}_V$; it depends only on $\rho^{V\text{-nr}}$.
\end{notation}

As a special case, one arrives at the genuine invariant distributions
\begin{equation*}
	\index{Idisc-nu-cV@$I_{\mathrm{disc}, \nu, c^V}$, $I_{\mathrm{disc}, \nu}$}
	I_{\mathrm{disc}, \nu, c^V} = I^{\tilde{G}}_{\mathrm{disc}, \nu, c^V}, \quad
	I_{\mathrm{disc}, \nu} = I^{\tilde{G}}_{\mathrm{disc}, \nu},
\end{equation*}
on $\tilde{G}$, as well as on $\tilde{G}_V$ in the notation above\footnote{The precise meaning will be clear according to the context.}, where $\nu$ ranges over $\mathfrak{h}_{\CC}^* / W_\infty$ with $\|\Im \nu\| = t$, and $c^V$ ranges over Satake parameters off $V$.

Note that $I_{\mathrm{disc}, \nu, c^V}$ is nonzero only when $c^V$ is automorphic.

It is convenient to express $I_t$ (resp.\ $I_\nu$) in terms of $I^{\tilde{M}}_{\mathrm{disc}, t}$ (resp.\ $I^{\tilde{M}}_{\mathrm{disc}, \nu}$), for various $M \in \mathcal{L}(M_0)$. This is phrased in terms of global weighted characters (Definition \ref{def:global-weighted-character}).

\begin{proposition}[Cf.\ {\cite[p.1206]{MW16-2}}]\label{prop:It-Levi-decomp}
	Let $\nu \in \mathfrak{h}^*_{\CC} /W_\infty$. We have
	\begin{equation*}
		I_\nu \left(f \right) = \sum_{M \in \mathcal{L}(M_0)} \frac{|W^M_0|}{|W^G_0|} I^{\mathrm{glob}}_{\tilde{M}}\left( I^{\tilde{M}}_{\mathrm{disc}, \nu}, f \right)
	\end{equation*}
	for all $f \in \orbI_{\asp}(\tilde{G}_V, \tilde{K}_V) \otimes \mes(G(F_V))$, where $I^{\tilde{M}}_{\mathrm{disc}, \nu}$ is viewed as a distribution on $\tilde{G}_V$ and $I^{\mathrm{glob}}_{\tilde{M}}\left( I^{\tilde{M}}_{\mathrm{disc}, \nu}, f \right)$ is understood as the convergent sum of terms $I^{\mathrm{glob}}_{\tilde{M}}(\pi_V \otimes c^V, f)$.
	
	Surely, there are similar results with the subscripts $\nu$ replaced by any $t \geq 0$.
\end{proposition}
\begin{proof}
	Let $t := \|\Im \nu\|$. Recalling \cite[Définition 5.18]{Li14b}, the expression of $I_\nu(f)$ in \eqref{eqn:Ispec-0} (with $M$ replaced by $L$) unfolds into
	\begin{multline*}
		\sum_{L \in \mathcal{L}(M_0)} \frac{|W^L_0|}{|W^G_0|} \sum_{M \in \mathcal{L}^L(M_0)} \frac{|W^M_0|}{|W^L_0|} \sum_{\pi_V, c^V} a^{\tilde{M}}_{\mathrm{disc}}(\pi_V \otimes c^V) \int_{i\mathfrak{a}^{L,*}_M} r^{\tilde{L}}_{\tilde{M}}(c^V_\lambda) I_{\tilde{L}_V}\left( \pi_V^{\tilde{L}}, 0, 0, f \right) \dd\lambda \\
		= \sum_{M \in \mathcal{L}(M_0)} \frac{|W^M_0|}{|W^G_0|} \sum_{\pi_V, c^V} a^{\tilde{M}}_{\mathrm{disc}}(\pi_V \otimes c^V) \sum_{L \in \mathcal{L}(M)} \int_{i\mathfrak{a}^{L,*}_M} r^{\tilde{L}}_{\tilde{M}}(c^V_\lambda) I_{\tilde{L}_V}\left( \pi_V^{\tilde{L}}, 0, 0, f \right) \dd\lambda
	\end{multline*}
	where $\pi_V$ (resp.\ $c^V$) ranges over $\Pi_{\mathrm{disc}, t, -}(\tilde{M}^1, V)$ (resp.\ Satake parameters of $\tilde{M}$ off $V$) such that the Archimedean component of $\pi_V$ has infinitesimal character $\nu$. Notice that
	\begin{itemize}
		\item in \textit{loc.\ cit.}, one puts $a^{\tilde{M}}_{\mathrm{disc}}(\pi_{V, \lambda} \otimes c^V_\lambda)$ inside the integral, which is redundant since the only reasonable definition of $a^{\tilde{M}}_{\mathrm{disc}}(\pi_V \otimes c^V)$ is invariant under $i\mathfrak{a}^{L,*}_M$-twists --- it depends only on $\pi_V \otimes c^V$ restricted to $\tilde{M}^1$;
		\item one wrote $I_{\tilde{L}_V}(\pi_V^{\tilde{L}}, 0, \cdot)$ instead of $I_{\tilde{L}_V}(\pi_V^{\tilde{L}}, 0, 0, \cdot)$ in \textit{loc.\ cit.}; see \eqref{eqn:I-0-0-V}.
	\end{itemize}

	By the fact that $r^{\tilde{G}}_{\tilde{G}}(c^V) = 1$ and \eqref{eqn:Ispec-1}, $I^{\tilde{G}}_{\mathrm{disc}, \nu}(f)$ equals the part with $M=G$ in the expansion above. This is true if $\tilde{G}$ is replaced by a Levi subgroup. Hence the expression for $I_\nu(f)$ follows.
\end{proof}

\begin{remark}
	The foregoing results about the spectral side of the trace formula apply to all adélic coverings subject to the requirements in \cite{Li14b}. Cf.\ Remark \ref{rem:weighted-characters-general-G}.
\end{remark}

\subsection{The stable side}
Consider $\mathbf{G}^! \in \Endo_{\elli}(\tilde{G})$. We choose a maximal compact subgroup $K^!_v$ of $G^!(F_v)$ for each place $v$, and assume $K^!_v$ to be hyperspecial when $v \notin V_{\mathrm{ram}}$. Observe that $A_{G^!, \infty} \simeq A_{G, \infty}$ via endoscopy.

To define the height of infinitesimal characters, we need to take $\mathfrak{h}^!$ together with an Hermitian norm on $\mathfrak{h}^{!,*}_{\CC}$; we transfer $\mathfrak{h}^!$ to the $\mathfrak{h}$ for $\tilde{G}$, then take the corresponding norm. With these data, we obtain the stable distribution $S^{G^!}_t$ on $G^!(F_V)$ and
\[ S^{G^!}_{\mathrm{spec}} = \sum_{t \geq 0} S^{G^!}_t, \quad S^{G^!}_t = \sum_{\nu^!: \|\Im\nu^!\| = t} S^{G^!}_{\nu^!} . \]
As in the invariant side, $A_{G^!, \infty}$-invariance is not imposed on test functions here.
\index{Sspec@$S^{G^{"!}}_{\mathrm{spec}}$, $S^{G^{"!}}_t$, $S^{G^{"!}}_{\nu^{"!}}$}

The discrete version $I_{\mathrm{disc}, t}$ also has a stable counterpart $S^{G^!}_{\mathrm{disc}, t}$. This is the achievement of Arthur's work \cite{Ar02}, especially the Global Theorem 2 therein. Moreover, $S^{G^!}_{\mathrm{disc}, t}$ is semi-finite. The decompositions in Notation \ref{nota:semi-finite-decomposition} also apply to the stable context; see \cite[X.5.5, 5.8]{MW16-2}. Therefore, we obtain the stable distributions
\begin{equation*}
	\index{Sdisc@$S^{G^{"!}}_{\mathrm{disc}, t}$, $S^{G^{"!}}_{\mathrm{disc}, \nu^{"!}}$, $S^{G^{"!}}_{\mathrm{disc}, \nu^{"!}, c^{V, "!}}$}
	S^{G^!}_{\mathrm{disc}, \nu^!}, \quad S^{G^!}_{\mathrm{disc}, \nu^!, c^{V, !}}
\end{equation*}
on $G^!(\A_F)$, where $V \supset V_{\mathrm{ram}}(G^!)$ is finite, $\nu^!$ (resp.\ $c^{V, !}$) stand for the infinitesimal characters (resp.\ Satake parameters off $V$) for $G^!$. They are also viewed as distributions on $G^!(F_V)$ through the map
\begin{align*}
	S\orbI(G^!(F_V), K^!_V) \otimes \mes(G^!(F_V)) & \to S\orbI(G^!(\A_F), K^!) \otimes \mes(G^!(\A_F)) \\
	f & \mapsto \dot{f}^! := f^! \mathbf{1}_{K^{!,V}}.
\end{align*}

Endoscopy affords a map $\mathfrak{h}^{!, *}_{\CC}/W^!_\infty \to \mathfrak{h}^*_{\CC} / W_\infty$, denoted by $\nu^! \mapsto \nu$, that preserves heights $\|\cdot\|$; see \cite[Proposition 7.3.1]{Li19}. Endoscopy also affords a map of Satake parameters $c^{V, !} \mapsto c^V$. From this, we deduce the decompositions
\begin{equation}\label{eqn:S-spec-decomp}\begin{aligned}
	S^{G^!}_\nu & := \sum_{\nu^! \mapsto \nu} S^{G^!}_{\nu^!}, \\
	S^{G^!}_{\mathrm{disc}, \nu, c^V} & := \sum_{\nu^! \mapsto \nu} \sum_{c^{V, !} \mapsto c^V} S^{G^!}_{\mathrm{disc}, \nu^!, c^{V, !}}. \\
	S^{G^!}_{\mathrm{disc}, t}\ & = \sum_{\nu, c^V} S^{G^!}_{\mathrm{disc}, \nu, c^V}
\end{aligned}\end{equation}
where we assume $\|\Im(\nu)\| = \|\Im(\nu^!)\| = t$.

There is a stable counterpart of Proposition \ref{prop:It-Levi-decomp}.

\begin{proposition}[Cf.\ {\cite[X.5.8]{MW16-2}}]\label{prop:St-Levi-decomp}
	Let $\nu \in \mathfrak{h}^*_{\CC}/W_\infty$. We have
	\[ S_\nu \left(f^! \right) = \sum_{M^! \in \mathcal{L}(M^!_0)} \frac{|W^{M^!}_0|}{|W^{G^!}_0|} S^{G^!, \mathrm{glob}}_{M^!}\left( S^{M^!}_{\mathrm{disc}, \nu}, f^! \right) \]
	for all $f^! \in S\orbI(G^!(F_V), K^!_V) \otimes \mes(G^!(F_V))$, where $S^{M^!}_{\mathrm{disc}, \nu}$ is viewed as a stable distribution on $M^!(F_V)$ and $S^{G^!, \mathrm{glob}}_{M^!}\left( S^{M^!}_{\mathrm{disc}, \nu}, f^! \right)$ is understood as a convergent sum.
\end{proposition}

\section{Statement of spectral stabilization}\label{sec:spectral-stabilization}
Let $\rev: \tilde{G} \to G(\A_F)$ be as in \S\ref{sec:Ispec}.

\begin{proposition}
	Let $\nu \in \mathfrak{h}_{\CC}^* / W_\infty$ and $t \geq 0$. For all $\mathbf{G}^! \in \Endo_{\elli}(\tilde{G})$, the transfers of $S^{G^!}_{\mathrm{disc}, \nu}$ and $S^{G^!}_{\mathrm{disc}, t}$ to $\tilde{G}$, defined as formal sums, are semi-finite. More precisely, $\trans_{\mathbf{G}^!, \tilde{G}} \left( S^{G^!}_{\mathrm{disc}, t} \right)$ is the formal sum of $\trans_{\mathbf{G}^!, \tilde{G}} \left( S^{G^!}_{\mathrm{disc}, \nu} \right)$ over $\nu$ with $\|\Im(\nu)\| = t$, and by fixing a finite set of places $V \supset V_{\mathrm{ram}}$, we have
	\[ \trans_{\mathbf{G}^!, \tilde{G}} \left( S^{G^!}_{\mathrm{disc}, \nu} \right)^{V\text{-nr}} = \sum_{c^V} \trans_{\mathbf{G}^!, \tilde{G}} \left( S^{G^!}_{\mathrm{disc}, \nu, c^V} \right). \]
\end{proposition}
\begin{proof}
	This is the metaplectic avatar of \cite[X.5.7]{MW16-2}. The only non-formal inputs for the proof in \textit{loc.\ cit.} are the spherical fundamental lemma \cite{Luo18} for places outside $V$, together with the spectral transfer at each $v \in V$, which is already covered by \S\ref{sec:spectral-transfer}.
\end{proof}

\begin{definition}[Cf.\ {\cite[X.5.9]{MW16-2}}]\label{def:IdiscEndo}
	\index{IEndo-disc@$I^{\Endo}_{\mathrm{disc}, \nu}$, $I^{\Endo}_{\mathrm{disc}, \nu, c^V}$}
	Fix a finite set of places $V \supset V_{\mathrm{ram}}$. Let $\nu \in \mathfrak{h}_{\CC}^* / W_\infty$ and let $c^V$ be a Satake parameter off $V$. Put
	\[ I^{\Endo}_{\mathrm{disc}, \nu, c^V} = I^{\tilde{G}, \Endo}_{\mathrm{disc}, \nu, c^V} := \sum_{\mathbf{G}^! \in \Endo_{\elli}(\tilde{G})} \iota(\tilde{G}, G^!) \trans_{\mathbf{G}^!, \tilde{G}} \left( S^{G^!}_{\mathrm{disc}, \nu, c^V} \right), \]
	and define $I^{\Endo}_{\mathrm{disc}, \nu} = I^{\tilde{G}, \Endo}_{\mathrm{disc}, \nu}$, $I^{\Endo}_{\mathrm{disc}, t} = I^{\tilde{G}, \Endo}_{\mathrm{disc}, t}$ similarly, where $t \geq 0$. They also give rise to distributions on $\tilde{G}_V$, denoted by the same tokens, by evaluation at test functions $\dot{f} = f f_{K^V}$.
\end{definition}

We will mainly work with the corresponding distributions on $\tilde{G}_V$. Note that by the fundamental lemma for units, for all $f \in \orbI_{\asp}(\tilde{G}_V, \tilde{K}_V) \otimes \mes(G(F_V))$ we still have formulas like
\begin{align*}
	I^{\Endo}_{\mathrm{disc}, \nu, c^V}(f) & = \sum_{\mathbf{G}^! \in \Endo_{\elli}(\tilde{G})} \iota(\tilde{G}, G^!) \trans_{\mathbf{G}^!, \tilde{G}} \left( S^{G^!}_{\mathrm{disc}, \nu, c^V} \right)(f) \\
	& = \sum_{\mathbf{G}^! \in \Endo_{\elli}(\tilde{G})} \iota(\tilde{G}, G^!) S^{G^!}_{\mathrm{disc}, \nu, c^V} \left(\Trans_{\mathbf{G}^!, \tilde{G}}(f) \right),
\end{align*}
and so forth. The goal of spectral stabilization is the following statement.

\begin{theorem}\label{prop:spectral-stabilization-disc}
	For all $V \supset V_{\mathrm{ram}}$, $\nu$ and $c^V$ as in Definition \ref{def:IdiscEndo}, we have $I_{\mathrm{disc}, \nu, c^V} = I^{\Endo}_{\mathrm{disc}, \nu, c^V}$.
	
	Consequently, we also have $I_{\mathrm{disc}, \nu} = I^{\Endo}_{\mathrm{disc}, \nu}$, even as distributions on $\tilde{G}$ by enlarging $V$.
\end{theorem}

The proof of this Theorem is deferred to \S\ref{sec:end-stabilization}.

Before proving Theorem \ref{prop:spectral-stabilization-disc}, some reductions are in order. Let us assume that the equality $I^{\tilde{S}, \Endo}_{\tilde{L}}(\tilde{\gamma}, \cdot) = I^{\tilde{S}}_{\tilde{L}}(\tilde{\gamma}, \cdot)$ of the local geometric Theorem \ref{prop:local-geometric} holds whenever $S \subsetneq G$. This one can assume by induction. Likewise, we assume the validity of Theorem \ref{prop:spectral-stabilization-disc} when $\tilde{G}$ is replaced by any proper Levi subgroup.

For each $f \in \orbI_{\asp}(\tilde{G}_V, \tilde{K}_V) \otimes \mes(G(F_V))$, let $\dot{f} := f f_{K^V}$.

\begin{lemma}\label{prop:diff-I-vs-Idisc}
	Under the assumptions above, we have
	\[ I_\nu(f) - \sum_{\mathbf{G}^!} \iota(\tilde{G}, G^!) S^{G^!}_\nu(f^{G^!}) = I^{\tilde{G}}_{\mathrm{disc}, \nu}(f) - \sum_{\mathbf{G}^!} \iota(\tilde{G}, G^!) S^{G^!}_{\mathrm{disc}, \nu}(f^{G^!}) \]
	for all $\nu \in \mathfrak{h}^*_{\CC} / W_\infty$ and $f$, where $f^{G^!} := \Trans_{\mathbf{G}^!, \tilde{G}}(f)$.
\end{lemma}
\begin{proof}	
	Expand $I_\nu$ by Proposition \ref{prop:It-Levi-decomp}, which yields a sum over $M \in \mathcal{L}^G(M_0)$. In parallel, expand $S^{G^!}_\nu$ by Proposition \ref{prop:St-Levi-decomp} for each $\mathbf{G}^!$. To compare them, note that Proposition \ref{prop:combinatorial-summation} transforms the left hand side into
	\begin{multline*}
		\sum_M \frac{|W^M_0|}{|W^G_0|} I^{\mathrm{glob}}_{\tilde{M}}\left( I^{\tilde{M}}_{\mathrm{disc}, \nu}, f \right)
		- \sum_{\mathbf{G}^!} \iota(\tilde{G}, G^!) \sum_{M^!} \frac{|W^{M^!}_0|}{|W^{G^!}_0|} S^{G^!, \mathrm{glob}}_{M^!}\left( S^{M^!}_{\mathrm{disc}, \nu}, f^! \right) \\
		= \sum_{M \in \mathcal{L}(M_0)} \frac{|W^M_0|}{|W^G_0|} \left( I^{\mathrm{glob}}_{\tilde{M}}\left( I^{\tilde{M}}_{\mathrm{disc}, \nu}, f \right) - \sum_{\mathbf{M}^! \in \Endo_{\elli}(\tilde{M})} \iota(\tilde{M}, M^!) I^{\mathrm{glob}, \Endo}_{\tilde{M}}\left( S^{M^!}_{\mathrm{disc}, \nu}, f \right) \right).
	\end{multline*}
	The fact that $S^{G^!, \mathrm{glob}}_{M^!}\left( S^{M^!}_{\mathrm{disc}, \nu}, f^! \right)$ is defined is part in Arthur's works on stabilization; see for example \cite[\S 8, \S 10]{Ar02}. Unsurprisingly, it follows inductively from the fact that $I^{G^!, \mathrm{glob}}_{M^!}\left( I^{M^!}_{\mathrm{disc}, \nu}, f^! \right)$ is defined (based on unitarity of $I^{M^!}_{\mathrm{disc}, \nu}$), and the stability of $S^{G^!, \mathrm{glob}}_{M^!}$.

	Furthermore, by
	\begin{itemize}
		\item the stabilization of global weighted characters (Proposition \ref{prop:stabilization-global-weighted-character}) which does not require any assumption about $M$ now,
		\item the assumption that $I^{\tilde{M}}_{\mathrm{disc}, \nu, c^V} = I^{\tilde{M}, \Endo}_{\mathrm{disc}, \nu, c^V}$ for all $M \subsetneq G$ and all $c^V$, which stabilizes $I^{\tilde{M}}_{\mathrm{disc}, \nu}$,
	\end{itemize}
	we see that only the terms with $M = G$ survive. This yields the right hand side.
\end{proof}

The spectral stabilization extends to all terms in the expansion of $I_{\mathrm{spec}}$. Its proof will be given in \S\ref{sec:end-stabilization}.

\begin{theorem}\label{prop:spectral-stabilization}
	\index{IEndo-spec@$I^{\Endo}_{\mathrm{spec}}$, $I^{\Endo}_t$, $I^{\Endo}_\nu$}
	For all $\nu \in \mathfrak{h}^*_{\CC} / W_\infty$ and $t \geq 0$, we have
	\begin{align*}
		I_\nu(f) & = \left( I^{\Endo}_\nu(f) := \sum_{\mathbf{G}^!} \iota(\tilde{G}, G^!) S^{G^!}_\nu(f^{G^!}) \right), \\
		I_t(f) & = \left( I^{\Endo}_t(f) := \sum_{\mathbf{G}^!} \iota(\tilde{G}, G^!) S^{G^!}_t(f^{G^!}) \right), \\
		I_{\mathrm{spec}}(f) & = \left( I^{\Endo}_{\mathrm{spec}}(f) := \sum_{\mathbf{G}^!} \iota(\tilde{G}, G^!) S^{G^!}_{\mathrm{spec}}(f^{G^!}) \right),
	\end{align*}
	where $\mathbf{G}^!$ ranges over $\Endo_{\elli}(\tilde{G})$ in each case.
\end{theorem}

Among the three equalities above, it suffices to prove the one for $I_\nu$. In view of Lemma \ref{prop:diff-I-vs-Idisc}, that assertion reduces to the stabilization of $I_{\mathrm{disc}, \nu}$, i.e.\ reduces to Theorem \ref{prop:spectral-stabilization-disc}.

\section{Spectral stabilization in a special case}\label{sec:spec-stabilization-special}
Fix a finite set $V$ of places such that $V \supset V_{\mathrm{ram}}$, and fix $M \in \mathcal{L}(M_0)$. We write $\widetilde{K^M} := \tilde{K} \cap \tilde{M}$. Without loss of generality, we also assume $M$ to be proper.

Denote by $\mathcal{E}(\mathfrak{h})$ the space of distributions of compact support on $\mathfrak{h}$. For every $\alpha \in \mathcal{E}(\mathfrak{h})$, its Fourier transform $\hat{\alpha} \in C^\infty(i\mathfrak{h}^*)$ extends meromorphically to $\mathfrak{h}^*_{\CC}$. Hereafter, take $\alpha \in C^\infty_c(\mathfrak{h})^{W_\infty}$. Such distributions serve as multipliers at Archimedean places: one can define the endomorphism $f \mapsto f_\alpha$ of $\orbI_{\asp}(\tilde{G}_V, \tilde{K}_V) \otimes \mes(G(F_V))$ such that
\index{f-alpha@$f_\alpha$}
\[ I^{\tilde{G}_V}(\pi_V, f_\alpha) = \hat{\alpha}(\nu) I^{\tilde{G}_V}(\pi_V, f) \]
for all irreducible unitarizable representation $\pi_V$ of $\tilde{G}_V$ whose Archimedean component $\pi_\infty \in \Pi_{\mathrm{unit}, -}(\tilde{G}_\infty)$ has infinitesimal character $\nu$. This is the content of \cite[Théorème 4.4]{Li13}, which just reformulates Arthur's results \cite{Ar83}.

Consider $f \in \orbI_{\asp, \mathrm{ac}}(\tilde{G}_V, \tilde{K}_V) \otimes \mes(G(F_V))$ and suppose that for each $\pi_V \in \Pi_{\mathrm{temp}, -}(\tilde{G}_V)$, there exists a function $\lambda \mapsto I^{\tilde{G}}(\pi_V, \lambda, f)$ on $i\mathfrak{a}^*_G$, smooth of rapid decay, such that
\[ I^{\tilde{G}_V}\left( \pi_V, X, f \right) = \int_{i\mathfrak{a}^*_G} e^{-\lrangle{\lambda, X}} I^{\tilde{G}_V}\left(\pi_V, \lambda, f \right) \dd\lambda, \quad X \in \mathfrak{a}_G. \]
For example, this is the case when $f \in \orbI_{\asp}(\tilde{G}_V, \tilde{K}_V) \otimes \mes(G(F_V))$. Then by following \cite[(6.1)]{Ar88-2}, one can still construct $f_\alpha$ with the same properties, characterized by
\[ I^{\tilde{G}_V}\left(\pi_V, X, f_\alpha \right) = \int_{i\mathfrak{a}^*_G} e^{-\lrangle{\lambda, X}} \hat{\alpha}(\nu + \lambda) I^{\tilde{G}_V}\left(\pi_V, \lambda, f\right) \dd\lambda. \]

It follows that $\hat{\alpha}(\nu + \lambda) I^{\tilde{G}_V}\left(\pi_V, \lambda, f\right) = I^{\tilde{G}_V}(\pi_V, \lambda, f_\alpha)$.

\begin{definition}\label{def:hV}
	\index{hfV@$h^f_{v_0}$, $h^f_V$}
	For all $f = \prod_v f_v \in \orbI_{\asp}(\tilde{G}_V, \tilde{K}_V) \otimes \mes(G(F_V))$ and $v_0 \in V$, set
	\begin{align*}
		h^f_{v_0} & := \epsilon_{\tilde{M}_{v_0}}(f_{v_0}) \prod_{v \in V \smallsetminus \{v_0\}} f_{v, \tilde{M}_v} \\
		& \in \orbI_{\asp, \mathrm{ac}}(\tilde{M}_V, (\widetilde{K^M})_V) \otimes \mes(M(F_V)), \\
		h^f_V & := \sum_{v_0 \in V} h^f_{v_0}.
	\end{align*}
\end{definition}

The functions $h^f_{v_0}$ and $h^f_V$ satisfy the properties mentioned above; in particular, one can consider $(h^f_V)_\alpha$, etc. Indeed, this follows from the corresponding properties of $\epsilon_{\tilde{M}_{v_0}}(f_{v_0})$ stated in Corollary \ref{prop:epsilonM-nonarch-Schwartz} (for non-Archimedean $v_0$) and in the discussions after Definition--Proposition \ref{def:epsilonM-arch} (for Archimedean $v_0$).

\begin{lemma}[Cf.\ {\cite[X.7.1 Lemme + Corollaire]{MW16-2}}]\label{prop:hV-Fourier}
	Fix $v_0 \in V$. Let $f \in \orbI_{\asp}(\tilde{G}_V, \tilde{K}_V) \otimes \mes(G(F_V))$.
	Let $\pi_V$ be a genuine admissible representation of $\tilde{M}_V / A_{M, \infty}$ of finite length, unitarizable with infinitesimal character $\nu$. There exists a meromorphic function
	\[ \lambda \mapsto I^{\tilde{M}}\left( \pi_V, \lambda, h^f_{v_0} \right), \quad \lambda \in \mathfrak{a}^*_{M, \CC}, \]
	which is analytic and of rapid decay on $i\mathfrak{a}^*_M$, such that
	\begin{equation*}
		I^{\tilde{M}_V}\left( \pi_V, X, h^f_{v_0} \right) = \int_{i\mathfrak{a}^*_M} e^{-\lrangle{\lambda, X}} I^{\tilde{M}}\left( \pi_V, \lambda, h^f_{v_0} \right) \dd\lambda
	\end{equation*}
	for all $X \in \mathfrak{a}_M$, and
	\begin{align*}
		I^{\tilde{M}}\left(\pi_V, \lambda, h^{zf}_{v_0} \right) & = z\left(\nu + \lambda\right) I^{\tilde{M}}\left(\pi_V, \lambda, h^f_{v_0} \right), \\
		I^{\tilde{M}}\left(\pi_V, X, h^{f_\alpha}_{v_0} \right) & = \int_{i\mathfrak{a}^*_M} e^{-\lrangle{\lambda, X}} \hat{\alpha}\left(\nu + \lambda\right) I^{\tilde{M}}\left(\pi_V, \lambda, h^f_{v_0} \right) \dd\lambda.
	\end{align*}
	for all $z \in \mathfrak{Z}(\mathfrak{g}_\infty)$ and $\alpha \in C^\infty_c(\mathfrak{h})^{W_\infty}$. Consequently, $h^{f_\alpha}_{v_0} = \left( h^f_{v_0} \right)_\alpha$.
\end{lemma}
\begin{proof}
	These assertions are of a local, analytic nature, based on Fourier transform and the trace Paley--Wiener theorems at each place. The proof is ultimately based on the properties of $\epsilon_{\tilde{M}_{v_0}}$ mentioned just before this Lemma, for each $v_0 \in V$. Therefore, the arguments in the cited reference carry over.
\end{proof}

Summing over $v_0$, we obtain $h^{f_\alpha}_V = \left( h^f_V \right)_\alpha$ for all $\alpha \in C^\infty_c(\mathfrak{h})^{W_\infty}$.

\begin{remark}\label{rem:moderate-function}
	A priori, $h^f_V$ does not necessarily belong to $\orbI_{\asp}(\tilde{M}_V, (\widetilde{K^M})_V) \otimes \mes(M(F_V))$ --- it only lies in $\orbI_{\asp, \cusp, \mathrm{ac}}(\cdots)$. Nonetheless, $h^f_V$ can still be plugged into the distribution $I^{\tilde{M}}_{\mathrm{geom}}$, and the simple trace formulas apply to $h^f_V$ as well; they all depend on $h^f_V |_{\tilde{G}_V^1}$.
	
	On the other hand, we will need to know that $h^f_V$ is \emph{moderate} in the sense of \cite[p.545]{Ar88-2} for certain estimates in the proof of Proposition \ref{prop:Idisc-stabilization-special} involving multipliers. To show $h^f_V$ is moderated, simply observe that the uniform growth condition (i) in \textit{loc.\ cit.} follows from Lemma \ref{prop:hV-Fourier}, whilst the support condition (ii) is trivial by cuspidality. The same arguments have been used in \cite[p.804]{Ar03-3}.
\end{remark}

We now move to the trace formula. It is legitimate to make the following assumptions.
\begin{description}
	\item[Local inductive assumption] Let $L \subset S$ be Levi subgroups of $G$ containing $M$ such that
	\[ \text{either}\quad S \subsetneq G \quad \text{or}\quad M \subsetneq L. \]
	We assume that for all $\tilde{\gamma} \in \tilde{L}_{S\text{-reg}}$, we have
	\[ I^{\tilde{S}, \Endo}_{\tilde{L}}(\tilde{\gamma}, \cdot) = I^{\tilde{S}}_{\tilde{L}}(\tilde{\gamma}, \cdot). \]
	\item[Global inductive assumption] For all proper Levi subgroup $L \in \mathcal{L}(M)$, $\nu \in \mathfrak{h}^*_{\CC} / W_\infty$ and Satake parameters $c^V$ off $V$, we have
	\[ I^{\tilde{L}}_{\mathrm{disc}, \nu, c^V} = I^{\tilde{L}, \Endo}_{\mathrm{disc}, \nu, c^V}. \]
	As seen in Theorem \ref{prop:spectral-stabilization-disc}, this implies $I^{\tilde{L}}_{\mathrm{disc}, \nu} = I^{\tilde{L}, \Endo}_{\mathrm{disc}, \nu}$ by allowing $V$ to vary.
\end{description}

Recall that all these distributions are evaluated on various $f \in \orbI_{\asp}(\tilde{G}_V, \tilde{K}_V) \otimes \mes(G(F_V))$, and depend only on $f|_{\tilde{G}_V^1}$.

\begin{lemma}\label{prop:disc-diff-geom}
	Let $f = \prod_{v \in V} f_v \in \orbI_{\asp}(\tilde{G}_V, \tilde{K}_V) \otimes \mes(G(F_V))$. Then
	\[ \sum_{t \geq 0} \sum_{\substack{\mu \in \mathfrak{h}^*_{\CC}/W_\infty \\ \|\Im(\mu)\|=t}} \left( I^{\Endo}_{\mathrm{disc}, \mu}(f) - I_{\mathrm{disc}, \mu}(f) \right) = I^{\Endo}_{\mathrm{geom}}(f) - I_{\mathrm{geom}}(f). \]
\end{lemma}
\begin{proof}
	By Lemma \ref{prop:diff-I-vs-Idisc}, the left hand side equals
	\begin{equation*}
		\sum_{t \geq 0} \sum_{\mu: \|\Im(\mu)\|=t} \left( I^{\Endo}_\mu(f) - I_\mu(f) \right) = I^{\Endo}_{\mathrm{spec}}(f) - I_{\mathrm{spec}}(f).
	\end{equation*}
	The invariant (resp.\ stable) trace formula for $\tilde{G}$ (resp.\ $G^!$ with $\mathbf{G}^! \in \Endo_{\elli}(\tilde{G})$) asserts $I_{\mathrm{geom}} = I_{\mathrm{spec}}$ (resp.\ $S^{G^!}_{\mathrm{geom}} = S^{G^!}_{\mathrm{spec}}$); see \cite{Li14b, Ar03-3} respectively. The right hand side of the displayed equation is therefore equal to $I^{\Endo}_{\mathrm{geom}}(f) - I_{\mathrm{geom}}(f)$.
\end{proof}

In favorable situations, the difference $I^{\Endo}_{\mathrm{geom}}(f) - I_{\mathrm{geom}}(f)$ in Lemma \ref{prop:disc-diff-geom} can be evaluated as follows.

\begin{lemma}\label{prop:geom-diff-h}
	Let $f = \prod_{v \in V} f_v \in \orbI_{\asp}(\tilde{G}_V, \tilde{K}_V) \otimes \mes(G(F_V))$ be subject to the conditions that
	\begin{itemize}
		\item $f$ is $\tilde{M}$-cuspidal (Definition \ref{def:M-cuspidal}) at two distinct places;
		\item there is a non-Archimedean place $v \in V$ such that for all $\tilde{\gamma} \in D_{\mathrm{geom}, -}(\tilde{G}_v) \otimes \mes(G(F_v))$, we have $I^{\tilde{G}_v}(\tilde{\gamma}_v, f_v) \neq 0$ only if $\Supp(\tilde{\gamma}) \subset \tilde{G}_{\mathrm{reg}}$.
	\end{itemize}
	Then
	\[ I^{\Endo}_{\mathrm{geom}}(f) - I_{\mathrm{geom}}(f) = |W^G(M)|^{-1} I^{\tilde{M}}_{\mathrm{geom}}(h^f_V). \]
\end{lemma}
\begin{proof}
	We assume $M$ is a standard Levi subgroup, and $f$ will be kept fixed throughout the proof. First off, the formula \eqref{eqn:Igeom-1} can be expressed as
	\[ I_{\mathrm{geom}}(f) = \sum_{L: \text{Levi} \; / \text{conj.}} |W^G(L)|^{-1} \sum_{\mathcal{O}^L} I_{\tilde{L}_V}\left(A^{\tilde{L}}(V, \mathcal{O}^L), f\right) \]
	where $\mathcal{O}^L$ ranges over stable semisimple classes in $L(F)$, by using the invariance of $I_{\tilde{L}_V}$ under $W^G_0$. Likewise, the formula for $I^{\Endo}_{\mathrm{geom}}$ in Proposition \ref{prop:Igeom-Endo-reindexed} is also a sum indexed by standard Levi subgroups. Hence $I^{\Endo}_{\mathrm{geom}}(f) - I_{\mathrm{geom}}(f)$ is the sum over standard Levi subgroups $L \subset G$ of $|W^G(L)|^{-1}$ times
	\begin{equation}\label{eqn:geom-diff-h-aux1}
		\sum_{\mathcal{O}^L} \left( \sum_{\mathbf{L}^! \in \Endo_{\elli}(\tilde{L})} \iota(\tilde{L}, L^!) \sum_{\mathcal{O}^{L, !} \mapsto \mathcal{O}^L} I^{\Endo}_{\tilde{L}_V}\left(\mathbf{L}^!, SA^{L^!}(\mathcal{O}^{L, !}, V), f \right) - I_{\tilde{L}_V}\left(A^{\tilde{L}}(\mathcal{O}^L, V), f\right) \right).
	\end{equation}
		
	We begin with the part $L=G$. By Theorem \ref{prop:matching-coeff-A} for $\tilde{G}$,
	\[ A^{\tilde{G}}(V, \mathcal{O}^G) = \sum_{\mathbf{G}^! \in \Endo_{\elli}(\tilde{G})} \iota(\tilde{G}, G^!) \sum_{\mathcal{O}^{G, !} \mapsto \mathcal{O}^G} \trans_{\mathbf{G}^!, \tilde{G}} \left( SA^{G^!}(V, \mathcal{O}^{G, !}) \right). \]
	Applying $I_{\tilde{G}_V}(\cdot, f) = I^{\Endo}_{\tilde{G}_V}(\cdot, f)$, the expression \eqref{eqn:geom-diff-h-aux1} is seen to vanish in this case.
	
	Now consider a proper standard Levi subgroup $L$. Let $D_{\diamondsuit, -}(\tilde{L}_V, \mathcal{O}^L)$ be $D_{\diamondsuit, -}(\tilde{L}_V) \cap D_{\mathrm{geom}, -}(\tilde{L}_V, \mathcal{O}^L)$, and so forth. Given $\mathcal{O}^L$ and
	\[ \tilde{\gamma} = \prod_{v \in V} \tilde{\gamma}_v \in D_{\mathrm{geom}, G_\infty\text{-equi}, -}(\tilde{L}_V, \mathcal{O}^L) \otimes \mes(L(F_V))^\vee , \]
	we define $L^V$ by $L^v := L$ for all $v \in V$ and infer from Remark \ref{rem:semi-local-orbint-ext} that
	\[ I_{\tilde{L}_V}(\tilde{\gamma}, f) = \sum_{S^V \in \mathcal{L}(L^V)} d^G_{L^V}\left( M, S^V \right) \prod_{v \in V} I^{\tilde{S}^v}_{\tilde{L}^v}(\tilde{\gamma}_v, f_{v, \tilde{S}^v}). \]
	Claim: when $d^G_{L^V}(L, S^V) \neq 0$, we are in either one of the following cases.
	\begin{description}
		\item[Old case] $S^v \neq G$ for all $v \in V$.
		
		\item[New case] $S^{v_0} = G$ for some $v_0 \in V$, whilst $S^v = L$ for all $v \in V \smallsetminus \{v_0\}$. In this case $d^G_{L^V}(L, S^V) = 1$.
	\end{description}
	In order to verify the claim, recall \eqref{eqn:d-V}: if $d^G_{L^V}(L, S^V) \neq 0$ then
	\[ \bigoplus_{v \in V} \mathfrak{a}^G_L = \mathrm{diag}(\mathfrak{a}^G_L) + \bigoplus_{v \in V} \mathfrak{a}^G_{S^v}. \]
	If there exists some $v_0$ with $S^{v_0} = G$, then
	\[ \sum_{v \neq v_0} \dim \mathfrak{a}^G_{S^v} + \dim\mathfrak{a}^G_L = |V| \dim\mathfrak{a}^G_L ; \]
	together with $\dim \mathfrak{a}^G_{S^v} \leq \dim \mathfrak{a}^G_L$, this implies $S^v = L$ for all $v \neq v_0$. The computation of $d^G_{L^V}(L, S^V)$ in the new case is straightforward.

	Therefore, $I_{\tilde{L}_V}(\tilde{\gamma}, f) = I_{\mathrm{new}}(\tilde{\gamma}) + I_{\mathrm{old}}(\tilde{\gamma})$ where
	\begin{align*}
		I_{\mathrm{new}}(\tilde{\gamma}) & := \sum_{v_0 \in V} I_{\tilde{L}_{v_0}}\left( \tilde{\gamma}_{v_0}, f_{v_0} \right) \prod_{v \neq v_0} I^{\tilde{L}_v}\left( \tilde{\gamma}_v, f_{v, \tilde{L}_v} \right), \\
		I_{\mathrm{old}}(\tilde{\gamma}) & := \sum_{S^V: \text{old}} d^G_{L^V}(L, S^V) \prod_{v \in V} I^{\tilde{S}^v}_{\tilde{L}^v}\left( \tilde{\gamma}_v, f_{v, \tilde{S^v}} \right).
	\end{align*}

	The terms in $I_{\mathrm{old}}(\tilde{\gamma})$ extend to all factorizable $\tilde{\gamma} \in D_{\diamondsuit, -}(\tilde{L}_V, \mathcal{O}^L) \otimes \mes(L(F_V))^\vee$ by induction. The crucial observation here is that $I_{\mathrm{new}}(\tilde{\gamma})$ also extends. Indeed, if $v_0$ is non-Archimedean, $I_{\tilde{L}_{v_0}}(\tilde{\gamma}_{v_0}, f_{v_0})$ is always defined; if $v_0$ is Archimedean, the assumption at some other non-Archimedean place implies that $\mathcal{O}^L$ (thus all $\gamma_v$) is strongly $G$-regular, otherwise the product over $V \smallsetminus \{v_0\}$ would be zero, so in this case $I_{\tilde{L}_{v_0}}(\tilde{\gamma}_{v_0}, f_{v_0})$ is also defined.

	Since $I_{\mathrm{new}}(\tilde{\gamma})$ and $I_{\mathrm{old}}(\tilde{\gamma})$ are multi-linear in $(\tilde{\gamma}_v)_v$, we conclude that $I_{\tilde{L}_V}(\tilde{\gamma}, f)$ extends canonically to all $\tilde{\gamma}$ in the space
	\[ D_{\diamondsuit, -}(\tilde{L}_V, \mathcal{O}^L) \otimes \mes(L(F_V))^\vee \]
	since the space is a $\otimes$-product over $V$.
	
	By Definition \ref{def:IMVEndo-gamma-weighted} and exactly the same method, we extend $I^{\Endo}_{\tilde{L}}(\tilde{\gamma}, f)$ to all $\tilde{\gamma}$ in the space displayed above. Specifically, $I^{\Endo}_{\tilde{L}}(\tilde{\gamma}, f)$ is the sum of
	\begin{align*}
		I^{\Endo}_{\mathrm{new}}(\tilde{\gamma}) & := \sum_{v_0 \in V} I^{\Endo}_{\tilde{L}}\left( \tilde{\gamma}_{v_0}, f_{v_0} \right) \prod_{v \neq v_0} I^{\tilde{L}_v}\left( \tilde{\gamma}_v, f_{v, \tilde{L}} \right) \\
		& \quad (\text{nonzero only when}\; \mathcal{O}_L \;\text{is strongly $G$-regular}), \\
		I^{\Endo}_{\mathrm{old}}(\tilde{\gamma}) & := \sum_{S^V: \text{old}} d^G_{L^V}(L, S^V) \prod_{v \in V} I^{\tilde{S}^v, \Endo}_{\tilde{L}^v}\left( \tilde{\gamma}_v, f_{v, \tilde{S^v}} \right) \\
		& = \sum_{S^V: \text{old}} d^G_{L^V}(L, S^V) \prod_{v \in V} I^{\tilde{S}^v}_{\tilde{L}^v}\left( \tilde{\gamma}_v, f_{v, \tilde{S^v}} \right).
	\end{align*}

	Furthermore, if $\mathbf{L}^! \in \Endo_{\elli}(\tilde{L})$ and $\mathcal{O}^{L, !}$ is a stable semisimple class in $L^!(F)$ mapping to $\mathcal{O}^L$, then
	\begin{equation}\label{eqn:geom-diff-h-aux0}
		I^{\Endo}(\mathbf{L}^!, \delta, f) = I^{\Endo}(\trans_{\mathbf{L}^!, \tilde{L}}(\delta), f)
	\end{equation}
	for all $\delta \in \left( SD_{\diamondsuit}(L^!_V) \cap SD_{\mathrm{geom}}(L^!_V, \mathcal{O}^{L, !}) \right) \otimes \mes(L^!(F_V))^\vee$. Indeed, this follows by comparing $I^{\Endo}_{\mathrm{new}}(\tilde{\gamma}) + I^{\Endo}_{\mathrm{old}}(\tilde{\gamma})$ with the splitting formula in Proposition \ref{prop:semilocal-matching-coherence} (i).

	Given a proper Levi $L$ and stable class $\mathcal{O}_L$, their contribution to \eqref{eqn:geom-diff-h-aux1} can be simplified by Theorem \ref{prop:matching-coeff-A} and \eqref{eqn:geom-diff-h-aux0} into
	\[ \left( I^{\Endo}_{\tilde{L}_V} - I_{\tilde{L}_V} \right)\left( A^{\tilde{L}}(\mathcal{O}^L, V), f \right), \]
	where $I^{\Endo}_{\tilde{L}}$ and $I_{\tilde{L}}$ are defined by the recipe alluded to above. In order to analyze this expression, we express $A^{\tilde{L}}(\mathcal{O}^L, V)$ as a linear combination of factorizable elements $\tilde{\gamma}$, and compare $I^{\Endo}_{\tilde{L}}(\tilde{\gamma}, f)$ with $I_{\tilde{L}}(\tilde{\gamma}, f)$: the ``old'' parts are equal, whist
	\[ I^{\Endo}_{\mathrm{new}}(\tilde{\gamma}) - I_{\mathrm{new}}(\tilde{\gamma}) = \sum_{v_0 \in V} \left( I^{\Endo}_{\tilde{L}_{v_0}}( \tilde{\gamma}_{v_0}, f_{v_0}) - I_{\tilde{L}_{v_0}}( \tilde{\gamma}_{v_0}, f_{v_0}) \right) \prod_{v \neq v_0} I^{\tilde{L}_v}\left( \tilde{\gamma}_v, f_{v, \tilde{L}_v} \right). \]
	If $L$ does not contain $M$, the above vanishes by $\tilde{M}$-cuspidality at two places. If $L \supsetneq M$, we also get zero by induction. If $L = M$, there are two cases for the summand indexed by $v_0 \in V$:
	\begin{enumerate}
		\item either $\mathcal{O}^L$ is strongly $G$-regular or $v_0$ is non-Archimedean, then summand is
		\[ I^{\tilde{M}_{v_0}}\left(\tilde{\gamma}_{v_0}, \epsilon_{\tilde{M}_{v_0}}(f_{v_0})\right) \prod_{v \neq v_0} I^{\tilde{M}_v}\left( \tilde{\gamma}_v, f_{v, \tilde{M}} \right) = I^{\tilde{M}_V}\left( \tilde{\gamma}, h^f_{v_0} \right); \]
		\item $v_0$ is Archimedean and $\mathcal{O}^L$ is not strongly $G$-regular, then the summand is zero by the assumption at some non-Archimedean place and our conventions in the ``new case'', whereas $I^{\tilde{M}_V}\left( \tilde{\gamma}, h^f_{v_0} \right)$ is also zero.
	\end{enumerate}
	All in all, we obtain $I^{\tilde{M}_V}\left( \tilde{\gamma}, h^f_V \right)$ by summing over $v_0 \in V$.

	By summing over proper $L$, $\mathcal{O}^L$ and $\tilde{\gamma}$ with suitable coefficients, the conclusion is that $I^{\Endo}_{\mathrm{geom}}(f) - I_{\mathrm{geom}}(f)$ equals $|W^G(M)|^{-1}$ times the expansion of $I^{\tilde{M}}_{\mathrm{geom}}\left(h^f_V \right)$.
\end{proof}

Let $\theta_\infty$ be a Cartan involution of $\mathfrak{g}_\infty$ adapted to the choices of $\mathfrak{h}$ and $K_\infty$. Following \cite[p.558]{Ar96}, we define
\[ \mathfrak{h}^*_u := \left\{ \nu \in \mathfrak{h}^*_{\CC}: \exists \omega \in W_\infty \theta_\infty, \; \omega^2 = \identity, \; \overline{\nu} = \omega \nu  \right\}. \]
It is a $W_\infty$-invariant subset of $\mathfrak{h}^*_{\CC}$ that contains the infinitesimal characters of all unitarizable genuine irreducible representations of $\tilde{G}_\infty$. It also contains the transfers of similar subsets from various $\mathbf{G}^! \in \Endo_{\elli}(\tilde{G})$.

For any $r, T \in \R_{\geq 0}$, define the $W_\infty$-invariant subsets
\begin{align*}
	\index{hu-rT@$\mathfrak{h}^*_u$, $\mathfrak{h}^*_u(r)$, $\mathfrak{h}^*_u(r, T)$}
	\mathfrak{h}^*_u(r, T) & := \left\{ \nu \in \mathfrak{h}^*_u : \|\Re(\nu)\| \leq r, \; \|\Im(\nu)\| \geq T \right\}, \\
	\mathfrak{h}^*_u(r) & := \mathfrak{h}^*_u(r, 0).
\end{align*}

\begin{lemma}\label{prop:hur}
	Given $f \in \orbI_{\asp}(\tilde{G}_V, \tilde{K}_V) \otimes \mes(G(F_V))$, there exists $r \geq 0$, depending only on the $\tilde{K}_\infty$-types of $f$, such that $I_{\mathrm{disc}, \nu}(f) = 0 = I^{\Endo}_{\mathrm{disc}, \nu}(f)$ unless $\nu \in \mathfrak{h}^*_u(r) / W_\infty$.
	
	Ditto for $S^{G^!}_{\mathrm{disc}, \nu}$ on the stable side, for each $\mathbf{G}^! \in \Endo_{\elli}(\tilde{G})$.
\end{lemma}
\begin{proof}
	This is a variant of \cite[Lemma 4.1]{Ar03-3} for $\tilde{G}$. For the case of $I_{\mathrm{disc}, \nu}$, it suffices to repeat the arguments in \textit{loc.\ cit.} which are only based on unitarity.

	The case of $S^{G^!}_{\mathrm{disc}, \nu}$ is already done in \textit{loc.\ cit.} It is routine to deduce the case of $I^{\Endo}_{\mathrm{disc}, \nu}$, provided that the $K^!_\infty$-types of $\Trans_{\mathbf{G}^!, \tilde{G}}(f)$ can be controlled in terms of $\tilde{K}_\infty$-types of $f$, for various $\mathbf{G}^!$. This follows from the trace Paley--Wiener theorems and standard properties of Archimedean spectral transfer; see \S\ref{sec:spectral-transfer}, especially the finite-support property for spectral transfer factors.
\end{proof}

\begin{proposition}[Cf.\ {\cite[Proposition 4.2 (a)]{Ar03-3}}]\label{prop:Idisc-stabilization-special}
	Let $f = \prod_v f_v \in \orbI_{\asp}(\tilde{G}_V, \tilde{K}_V) \otimes \mes(G(F_V))$ be subject to the conditions that
	\begin{itemize}
		\item $f$ is $\tilde{M}$-cuspidal at two distinct places;
		\item there exists a non-Archimedean $v \in V$ such that for all $\tilde{\gamma} \in D_{\mathrm{geom}, -}(\tilde{G}_v) \otimes \mes(G(F_v))$, we have $I^{\tilde{G}_v}(\tilde{\gamma}, f_v) \neq 0$ only if $\Supp(\tilde{\gamma}) \subset \tilde{G}_{v, \mathrm{reg}}$.
	\end{itemize}
	Then
	\[ I_{\mathrm{disc}, \nu}(f) = I^{\Endo}_{\mathrm{disc}, \nu}(f) \]
	for all $\nu \in \mathfrak{h}^*_{\CC} / W_\infty$.
\end{proposition}
\begin{proof}
	In view of the invariant (resp.\ stable) trace formula for $\tilde{G}$ (resp.\ various $G^!$), Lemmas \ref{prop:disc-diff-geom} and \ref{prop:geom-diff-h} imply
	\begin{equation}\label{eqn:Idisc-stabilization-special-aux}
		\sum_{t \geq 0} \sum_{\substack{\mu \in \mathfrak{h}^*_{\CC}/W_\infty \\ \|\Im(\mu)\| = t}} \left( I^{\Endo}_{\mathrm{disc}, \mu}(f) - I_{\mathrm{disc}, \mu}(f) \right)
		= |W^G(M)|^{-1} I^{\tilde{M}}_{\mathrm{geom}} \left( h^f_V \right).
	\end{equation}
	
	Observe that $h^f_V$ is cuspidal at two distinct places, since $\epsilon_{\tilde{M}_v}(f_v)$ is cuspidal for all $v$. The version 1 of the simple trace formula \cite[Théorème 6.7]{Li14b} implies $I^{\tilde{M}}_{\mathrm{geom}} \left( h^f_V \right) = \sum_{t \geq 0} I^{\tilde{M}}_{\mathrm{disc}, t}(\dot{h}^f_V)$.
	
	In brief, the idea is that both sides of
	\[ \sum_{t \geq 0} \sum_{\mu: \|\Im(\mu)\|=t } \left( I^{\Endo}_{\mathrm{disc}, \mu}(f) - I_{\mathrm{disc}, \mu}(f) \right) = |W^G(M)|^{-1} \sum_{t \geq 0} I^{\tilde{M}}_{\mathrm{disc}, t} \left( h^f_V \right) \]
	must vanish since the left (resp.\ right) hand side is a sum of integrals over affine subspaces of $\mathfrak{h}^*$ under $i\mathfrak{a}_G^*$ (resp.\ $i\mathfrak{a}_M^*$), whilst $M \neq G$. The required equality then follows by isolating infinitesimal characters. To make this rigorous, we appeal to multipliers.
	
	In order to simplify the usage of multipliers, let us we assume $\tilde{G} = \Mp(W)$ hereafter, so that $\mathfrak{a}^*_G  = \{0\}$. This is legitimate since the matching on $\GL$-factors is tautological.
	\footnote{For the general case, one has to work with multipliers in $C^\infty_c(\mathfrak{h}^1)^{W_\infty}$ where $\mathfrak{h}^1$ is the direct summand of $\mathfrak{h}$ cut out by $\tilde{G}^1_\infty$. The results in the beginning of this section should also be adapted to $\tilde{G}_V^1 \simeq \tilde{G}_V / A_{G, \infty}$. This is comparable with the framework in the proof of \cite[Proposition 4.2 (a)]{Ar03-3}.}
	
	Given $T \geq 0$ and $\alpha \in C^\infty_c(\mathfrak{h})^{W_\infty}$, we replace $f$ by $f_\alpha$ to get
	\begin{align*}
		I^{\tilde{M}}_{\mathrm{geom}} \left( h^{f_\alpha}_V \right) & = \sum_{0 \leq t \leq T} I^{\tilde{M}}_{\mathrm{disc}, t} \left( h^{f_\alpha}_V \right) + \sum_{t > T} I^{\tilde{M}}_{\mathrm{disc}, t} \left( h^{f_\alpha}_V \right) \\
		& = \sum_{0 \leq t \leq T} I^{\tilde{M}}_{\mathrm{disc}, t} \left( (\dot{h}^f_V)_\alpha \right) + \sum_{t > T} I^{\tilde{M}}_{\mathrm{disc}, t} \left( (h^f_V)_\alpha \right).
	\end{align*}
	
	Take $r \geq 0$ as in Lemma \ref{prop:hur}, large enough with respect to both $(\tilde{G}, f)$ and $(\tilde{M}, h_V^f)$; only the $\tilde{K}_V$-types of these functions matter here. This choice is not affected by the replacement of $f$ by $f_\alpha$. We claim that there exists $C, k > 0$ such that for all $N, T > 0$ and $\alpha \in C^\infty_c(\mathfrak{h})^{W_\infty}$ supported in $\|\cdot\| \leq N$, we have
	\begin{equation}\label{eqn:Idisc-stabilization-special-aux0}
		\sum_{t > T} I^{\tilde{M}}_{\mathrm{disc}, t} \left( (h^f_V)_\alpha \right) \leq Ce^{kN} \sup_{\mu \in \mathfrak{h}^*_u(r, T)} \left|\hat{\alpha}(\mu)\right|.
	\end{equation}
	
	The claim takes the same form as Arthur's \emph{multiplier convergence estimates} in \cite[(3.3)]{Ar02}, see also \cite[p.2432]{Li14b}. Nonetheless, there are two subtle differences.
	\begin{itemize}
		\item We are using multipliers from $C^\infty_c(\mathfrak{h})^{W_\infty}$ rather than the $C^\infty_c(\mathfrak{h}^{M, 1})$ used in \cite{Ar02}, where $\mathfrak{h}^{M, 1}$ is the direct summand of $\mathfrak{h}$ cut out by $\tilde{M}^1_\infty$. Otherwise it would be difficult treat multipliers of $\tilde{G}_V$ and $\tilde{M}_V$ on an equal footing.
		\item In parallel, we use test functions on $\tilde{M}_V$ and not $\tilde{M}_V^1$; however, the multiplier convergence estimate only hold for moderate functions on $\tilde{M}_V$, which have appeared in Remark \ref{rem:moderate-function} and form a subspace of $\orbI_{\asp, \mathrm{ac}}$ containing $\orbI_{\asp}$.
	\end{itemize}
	The relevant version of the estimates is stated in \cite[\S 6]{Ar88-2}; see Lemma 6.3 and the remark after Corollary 6.5 therein. The proofs in \textit{loc.\ cit.} are of an analytic nature, and the ingredients for coverings can be found in \cite{Li13}. Surely, the same applies when $\tilde{M}$ is replaced by $\tilde{G}$.
	
	Thanks to Remark \ref{rem:moderate-function}, we know that $(h^f_V)_\alpha = h^{f_\alpha}_V$ is moderate. Hence \eqref{eqn:Idisc-stabilization-special-aux0} holds.
	
	Denote $\alpha_m := \underbracket{\alpha \star \cdots \star \alpha}_{m \;\text{times}}$ for all $m \geq 1$. Note that $\alpha_m$ is supported in $\|\cdot\| \leq mN$. Given $N, T > 0$, let us take $\alpha$ such that
	\begin{itemize}
		\item $\alpha$ is supported in $\|\cdot\| \leq N$,
		\item $\hat{\alpha}$ restricts to $\mathfrak{h}^*_u(r) \to [0, 1]$,
		\item $\hat{\alpha}(\nu) = 1$ if and only if $\mu \in W_\infty \cdot \nu$.
	\end{itemize}

	By the classical Paley--Wiener theorem, $\hat{\alpha}$ is rapidly decreasing over vertical strips. Hence we may take $T$ sufficiently large to ensure $|\hat{\alpha}| \leq e^{-2kN}$ over $\mathfrak{h}^*_u(r, T)$. Hence
	\begin{equation}\label{eqn:Idisc-stabilization-special-aux1}\begin{aligned}
		\sum_{t > T} I^{\tilde{M}}_{\mathrm{disc}, t} \left( (h^f_V)_{\alpha_m} \right) & \leq Ce^{kmN} \sup_{\mu \in \mathfrak{h}^*_u(r, T)} \left|\hat{\alpha}(\mu)\right|^m \\
		& \leq Ce^{-kmN} \to 0, \quad \text{as} \quad m \to +\infty.
	\end{aligned}\end{equation}
	
	On the other hand, we claim that
	\begin{equation}\label{eqn:Idisc-stabilization-special-aux2}
		\lim_{m \to +\infty} \sum_{t \geq 0} \left( I^{\Endo}_{\mathrm{disc}, t}(f_{\alpha_m}) - I_{\mathrm{disc}, t}(f_{\alpha_m}) \right) = I^{\Endo}_{\mathrm{disc}, \nu}(f) - I_{\mathrm{disc}, \nu}(f).
	\end{equation}
	To see this, note that when $t$ is given, it follows from \eqref{eqn:Ispec-1} and dominated convergence that
	\[ \lim_{m \to +\infty} \left( I^{\Endo}_{\mathrm{disc}, t}(f_{\alpha_m}) - I_{\mathrm{disc}, t}(f_{\alpha_m}) \right) = \begin{cases}
		I^{\Endo}_{\mathrm{disc}, \nu}(f) - I_{\mathrm{disc}, \nu}(f), & \text{if}\; \|\Im(\nu)\| = t \\
		0, & \text{otherwise.}
	\end{cases}\]
	We take a sufficiently large $T'$, break the $t$-sum into $\sum_{0 \leq t \leq T'}$ and $\sum_{t > T'}$, then show as in \eqref{eqn:Idisc-stabilization-special-aux1} that $\sum_{t > T'} \to 0$ by using
	\begin{itemize}
		\item the multiplier convergence estimate \eqref{eqn:Idisc-stabilization-special-aux0} with $\tilde{M}$ (resp.\ $h_V^f$) replaced by $\tilde{G}$ (resp.\ $f$);
		\item its stable counterparts \cite[Proposition 10.5 (b)]{Ar02} for all $G^!$;
	\end{itemize}
	According to \cite[Lemme 7.2]{Li13}, combined with the construction of discrete spectrum and Lemma \ref{prop:hur}, the $t$-sum over $[0, T']$ is finite and involves only finitely many infinitesimal characters. This proves \eqref{eqn:Idisc-stabilization-special-aux2}.

	The combination of \eqref{eqn:Idisc-stabilization-special-aux}, \eqref{eqn:Idisc-stabilization-special-aux1} and \eqref{eqn:Idisc-stabilization-special-aux2} leads to
	\[ I^{\Endo}_{\mathrm{disc}, \nu}(f) - I_{\mathrm{disc}, \nu}(f) = |W^G(M)|^{-1} \lim_{m \to +\infty} \sum_{0 \leq t \leq T} I^{\tilde{M}}_{\mathrm{disc}, t} \left( (h^f_V)_{\alpha_m} \right). \]
	
	It remains to show that the limit above is zero. Observe that the $\tilde{K}_V$-type of $(h^f_V)_{\alpha_m}$ does not depend on $m$. As noticed earlier (in the case of $\tilde{G}$), the $t$-sum over $[0, T]$ is finite and involves only finitely many infinitesimal characters. In view of the semi-finiteness of $I^{\tilde{M}}_{\mathrm{disc}, t}$, the last displayed equality in Lemma \ref{prop:hV-Fourier} expresses $\sum_{0 \leq t \leq T} I^{\tilde{M}}_{\mathrm{disc}, t} \left( (h^f_V)_{\alpha_m}\right)$ as a sum (uniformly finite in $m$) of terms
	\[ \text{const} \cdot \int_{i\mathfrak{a}^*_M} e^{-\lrangle{\lambda, X}} \hat{\alpha}^m \left(\nu + \lambda\right) I^{\tilde{M}}\left(\pi_V, \lambda, h^f_{v_0} \right) \dd\lambda  \]
	where $v_0 \in V$, $\pi_V \in \Pi_{\mathrm{unit}, -}(\tilde{M}_V / A_{M, \infty})$, and $\nu \in \mathfrak{h}^*_u(r)$ is a representative of the infinitesimal character of $\pi_V$. By the choice of $\alpha$ and the dominated convergence theorem, we see
	\[ \lim_{m \to +\infty} \sum_{0 \leq t \leq T} I^{\tilde{M}}_{\mathrm{disc}, t} \left( (h^f_V)_{\alpha_m}\right) = 0 , \]
	as desired.
\end{proof}

\section{Proof of the key geometric hypothesis}\label{sec:proof-key-geometric}
We set out to prove Hypothesis \ref{hyp:key-geometric}. This is a local problem, so we take $F$ to be a local field of characteristic zero, $\rev: \tilde{G} \to G(F)$ be a covering of metaplectic type. Let $\tilde{M} \subset \tilde{G}$ and $\mathbf{M}^! \in \Endo_{\elli}(\tilde{M})$ be as in that statement.

For simplicity, let us assume $G = \Sp(W)$ and $M = \prod_{i \in I} \GL(\ell_i) \times \Sp(W^\flat)$ with respect to an orthogonal decomposition of symplectic $F$-vector spaces
\[ W = W^\flat \oplus \bigoplus_{i \in I} (\ell_i \oplus \ell'_i), \]
where $\ell_i$ and $\ell'_i$ are all totally isotropic. Identify $\ell'_i$ with the dual of $\ell_i$, for all $i$. Recall from the description of $\tilde{G}$ via Maslov cocycles (the proof of Definition--Proposition \ref{def:G-minus}) that $\tilde{G}$ is unaltered if $\psi$ (or $\lrangle{\cdot|\cdot}$) is dilated by $F^{\times 2}$.

\begin{lemma}\label{prop:globalization}
	Let $W = W^\flat \oplus \bigoplus_{i \in I} (\ell_i \oplus \ell'_i)$ and $M \subset G$ be as above. Let $T$ be an elliptic maximal $F$-torus in $M$. For all $r \geq 1$, there exist
	\begin{itemize}
		\item a number field $\mathring{F}$ with distinct places $u_0, \ldots, u_r$;
		\item a nontrivial additive character $\mathring{\psi}: \mathring{F} \backslash \A_{\mathring{F}} \to \mathbb{S}^1$;
		\item a symplectic $\mathring{F}$-vector space $(\mathring{W}, \lrangle{\cdot|\cdot})$ together with an orthogonal decomposition into symplectic subspaces $\mathring{W}^\flat \oplus \bigoplus_{i \in I} (\mathring{\ell}_i \oplus \mathring{\ell}'_i)$, where $\mathring{\ell}_i$ and $\mathring{\ell}'_i$ are all totally isotropic;
		\item a maximal $\mathring{F}$-torus $\mathring{T} \subset \mathring{M} := \prod_{i \in I} \GL(\mathring{\ell}_i) \times \Sp(\mathring{W}^\flat) \subset \mathring{G} := \Sp(\mathring{W})$.
	\end{itemize}

	Denote by $\widetilde{\mathring{M}} \subset \widetilde{\mathring{G}}$ the corresponding coverings. The local and global data are compatible for every $0 \leq j \leq r$ in the following sense.
	\begin{itemize}
		\item There exists an isomorphism $\mathring{F}_{u_j} \simeq F$, which we fix; the global data localize at $u_j$ to $T \subset M \subset G$ and $W$, $(\ell_i, \ell'_i)_{i \in I}$, up to isomorphisms.
		\item On the other hand, there exist $a_j \in F^\times$ such that
		\[ \mathring{\psi}(t) = \psi(a_j^2 t), \quad t \in F^\times \simeq F_{u_j}^\times . \]
		\item In particular, $\tilde{G} \simeq \widetilde{\mathring{G}}_{u_j}$, $\tilde{M} \simeq \widetilde{\mathring{M}}_{u_j}$, and $\mathring{T}$ is an elliptic maximal $\mathring{F}$-torus in $\mathring{M}$.
		\item The map $\mathring{T}(\mathring{F}) \to \mathring{T}(\mathring{F}_{u_j}) \simeq T(F)$ has dense image.
		\item The natural map $\Hm^1(\mathring{F}_{u_j}, \mathring{T}_{u_j}) \to \Hm^1(\mathring{F} \backslash \A_{\mathring{F}}, \mathring{T})$ is bijective.
	\end{itemize}
\end{lemma}
\begin{proof}
	This is essentially \cite[Proposition 8.4.1]{Li19}. It is comprised of two parts:
	\begin{inparaenum}[(i)]
		\item globalization of fields, and
		\item globalization of linear-algebraic data such as $(W, \lrangle{\cdot|\cdot})$, etc.
	\end{inparaenum}
	Compared with \textit{loc.\ cit.}, the only new ingredient here is the Levi subgroup $M \subset G$ and $(\ell_i, \ell'_i)_{i \in I}$; the required arguments for them are similar, and actually simpler since they only involve the part (i).
	
	Although one assumes $F$ is non-Archimedean in \textit{loc.\ cit.}, the arguments carry over verbatim.
\end{proof}

\begin{proof}[Proof of Hypothesis \ref{hyp:key-geometric}]
	In order to well-define orbital integrals, we choose Haar measures on $M^!_\delta(F)$ for all $\delta \in M^!_{\mathrm{reg}}(F)$ in a coherent way, and similarly for other groups.

	First off, it suffices to prove that for any given $\delta \in M^!_{G\text{-reg}}(F)$, there exists $\epsilon(\mathbf{M}^!, \delta) \in \CC$ such that
	\begin{equation}\label{eqn:key-geometric-pf-aux-0}
		\Trans_{\mathbf{M}^!, \tilde{M}}\left(\epsilon_{\tilde{M}}(f)\right)(\delta) = \epsilon(\mathbf{M}^!, \delta) \Trans_{\mathbf{M}^!, \tilde{M}}\left(f_{\tilde{M}}\right)(\delta)
	\end{equation}
	for all $f$. By varying $f$, it will follow that $\epsilon(\mathbf{M}^!, \delta)$ is a uniquely defined smooth function in $\delta$.

	We may take $\delta$ from a given maximal $F$-torus $T^! \subset M^!$. As $\epsilon_{\tilde{M}}(\cdot)$ is cuspidal, so is $\Trans_{\mathbf{M}^!, \tilde{M}}(\epsilon_{\tilde{M}}(\cdot))$ and we may assume $T^!$ is elliptic. We claim that it suffices to verify \eqref{eqn:key-geometric-pf-aux-0} for $\delta$ in a dense subset of $T^!_{G\text{-reg}}(F)$: introduce $g \in \orbI_{\asp}(\tilde{G}, \tilde{K}) \otimes \mes(G(F))$ and reformulate \eqref{eqn:key-geometric-pf-aux-0} as an equality between functions on $T^!_{G\text{-reg}}(F)$:
	\begin{equation}\label{eqn:key-geometric-pf-aux-1}
		\Trans_{\mathbf{M}^!, \tilde{M}}\left(\epsilon_{\tilde{M}}(f)\right)(\delta) \cdot
		\Trans_{\mathbf{M}^!, \tilde{M}}\left(g_{\tilde{M}}\right)(\delta)
		= \Trans_{\mathbf{M}^!, \tilde{M}}\left(\epsilon_{\tilde{M}}(g)\right)(\delta) \cdot
		\Trans_{\mathbf{M}^!, \tilde{M}}\left(f_{\tilde{M}}\right)(\delta)
	\end{equation}
	for all $f, g$. Both sides being continuous in $\delta$, the claim follows.

	Suppose that $T^! \subset M^!$ transfers to $T \subset M$. Then $T$ is elliptic as well. Apply Lemma \ref{prop:globalization} with $r \geq 2$ to globalize these data. We also need the globalization of endoscopic data, which is easier since $\Endo_{\elli}(\widetilde{\mathring{M}}) = \Endo_{\elli}(\widetilde{\mathring{M}}_u)$ at each place $u$. We can be more precise: $\mathbf{M}^!$ and the stable class of $T^! \subset M^!$ are determined by $T$ together with a character $\kappa$ of $\Hm^1(F, T)$; the global case is similar by considering characters $\mathring{\kappa}$ of $\Hm^1(\mathring{F} \backslash \A_{\mathring{F}}, \mathring{T})$ instead. By the property concerning $\Hm^1$ in Lemma \ref{prop:globalization}, we can globalize $(T, \kappa)$ to $(\mathring{T}, \mathring{\kappa})$, hence globalize $(\mathbf{M}^!, T^!)$ into $(\mathbf{\mathring{M}}^! , \mathring{T}^!)$.

	Consider a $\widetilde{\mathring{G}}$-regular $\mathring{\delta} \in \mathring{T}^!(\mathring{F})$, and take $\mathring{\gamma} \in \mathring{T}(\mathring{F})$ with $\mathring{\delta} \leftrightarrow \mathring{\gamma}$. Define
	\begin{equation*}\begin{aligned}
		\mathcal{O}_1: & \text{the conjugacy class in}\; \mathring{G}(\mathring{F}) \;\text{containing}\; \mathring{\gamma}, \\
		\mathcal{O}: & \text{the conjugacy class in}\; \mathring{M}(\mathring{F}) \;\text{containing}\; \mathring{\gamma}, \\
		\mathcal{O}^!_1: & \text{the stable conjugacy class in}\; \mathring{G}^!(\mathring{F}) \;\text{containing}\; \mathring{\delta}, \\
		\mathcal{O}^!: & \text{the stable conjugacy class in}\; \mathring{M}^!(\mathring{F}) \;\text{containing}\; \mathring{\delta}.
	\end{aligned}\end{equation*}
	
	Take a finite set of places $V \supset V_{\mathrm{ram}} \cup \{u_0, \ldots, u_r\}$ of $\mathring{F}$, sufficiently large in the precise sense that $V \supset S(\mathcal{O}_1, K)'$ (Definition \ref{def:SOK}), with hyperspecial subgroups $(K_u)_{u \notin V_{\mathrm{ram}}}$ chosen as in \S\ref{sec:geom-side}. We also choose maximal compact subgroups $K_u \subset \mathring{G}(\mathring{F}_u)$ for each $u \in V$, in good positions relative to $\mathring{M}$. By the description of the image of transfer (Theorem \ref{prop:image-transfer}) and trace Paley--Wiener theorems, we may choose\footnote{The symbol $f$ has been used, thus we write $\varphi$.}
	\[ \varphi = f_{u_0} \prod_{u \neq u_0} f_u \; \in \orbI_{\asp}(\widetilde{\mathring{G}}_V, \tilde{K}_V) \otimes \mes(\mathring{G}(\mathring{F}_V)) \]
	with the following properties. First, write $\tilde{G}_v := \widetilde{\mathring{G}}_v$, $\tilde{M}_v := \widetilde{\mathring{M}}_v$ for all places $v$, and identify $\tilde{M}$ with $\tilde{M}_{u_j}$ for all $0 \leq j \leq r$.
	\begin{itemize}
		\item There exists a non-Archimedean place $v \in V \smallsetminus \{u_0\}$ such that $f_v$ comes from a $C^\infty_{c, \asp}$-function on the $\tilde{G}_{v, \mathrm{reg}}$.
		\item The functions $f_{u_1}$ and $f_{u_2}$ are $\tilde{M}$-cuspidal (Definition \ref{def:M-cuspidal}).
		\item For all $0 \leq j \leq r$ and all $\mathbf{M}^{!!} \in \Endo_{\elli}(\tilde{M})$ not equivalent to $\mathbf{M}^!$, we have $\Trans_{\mathbf{M}^{!!}, \tilde{M}}\left( f_{u_j, \tilde{M}} \right) = 0$.
		\item Let $j \in \{1, 2\}$. For every $\delta_j^* \in \Sigma_{\mathrm{reg}}(M^!)$ such that $\delta_j^* \leftrightarrow \mathring{\gamma}_{u_j}$ relative to $\mathbf{M}^!$ but $\delta_j^* \neq \mathring{\delta}_{u_j}$ modulo $W^G(M)$, we require that $\Trans_{\mathbf{M}^!, \tilde{M}}(f_{u_j, \tilde{M}})(\delta_j^*) = 0$.
		\item For all $v \in V \smallsetminus \{u_0\}$, we require that $\Trans_{\mathring{\mathbf{M}}^!_v, \tilde{M}_v}(f_{v, \tilde{M}_v})(\mathring{\delta}_v) = 1$.
	\end{itemize}
	
	There are no conditions on $f_{u_0}$.
	
	We contend that upon imposing further constraints on $\varphi$ at the places in $V \smallsetminus \{u_0\}$, which may depend on the support of $f_{u_0}$, we may ensure that there is a constant $c \neq 0$ such that
	\begin{equation}\label{eqn:key-geometric-aux-1}
		I^{\widetilde{\mathring{M}}}_{\mathrm{geom}}(h_V^\varphi) = c \cdot  S^{\mathring{M}^!_V}\left(\mathring{\delta}_V , \Trans_{\mathring{\mathbf{M}}^!, \mathring{\tilde{M}}}(h_V^\varphi) \right)
	\end{equation}
	where we fix a Haar measure on $\mathring{M}^!_{\mathring{\delta}}(F_V)$ to define the stable orbital integral.

	Granting \eqref{eqn:key-geometric-aux-1}, the definition of $h^\varphi_V$ entails
	\begin{multline*}
		I^{\widetilde{\mathring{M}}}_{\mathrm{geom}}(h_V^\varphi) = x \cdot \Trans_{\mathbf{M}^!, \tilde{M}}\left( \epsilon_{\tilde{M}}(f_{u_0})\right)(\mathring{\delta}_{u_0}) \underbracket{\prod_{v \in V \smallsetminus \{ u_0\}} \Trans_{\mathring{\mathbf{M}}^!_v, \tilde{M}_v}(f_{v, \tilde{M}}) \left( f_{v, \tilde{M}}\right)(\mathring{\delta}_v)}_{= 1} \\
		+ y \cdot \Trans_{\mathbf{M}^!, \tilde{M}}\left( f_{u_0, \tilde{M}}\right)(\mathring{\delta}_{u_0}),
	\end{multline*}
	where $x, y$ are constants, $x \neq 0$. On the other hand, as $\varphi$ is $\tilde{M}$-cuspidal at $u_1$ and $u_2$, and coming from $C^\infty_{c, \asp}(\tilde{G}_{v, \mathrm{reg}})$ for some non-Archimedean place $v \in V$,	Proposition \ref{prop:Idisc-stabilization-special} together with Lemmas \ref{prop:disc-diff-geom} and \ref{prop:geom-diff-h} imply
	\begin{align*}
		I^{\widetilde{\mathring{M}}}_{\mathrm{geom}}(h_V^\varphi) & = |W^{\mathring{G}}(\mathring{M})| \left( I^{\Endo}_{\mathrm{geom}}(\varphi) - I_{\mathrm{geom}}(\varphi)\right) \\
		& = |W^{\mathring{G}}(\mathring{M})| \sum_{t \geq 0} \sum_{\substack{\mu \in \mathfrak{h}^*_{\CC}/W_\infty \\ \|\Im(\mu)\|=t }} \left( I^{\Endo}_{\mathrm{disc}, \mu}(\varphi) - I_{\mathrm{disc}, \mu}(\varphi) \right) = 0.
	\end{align*}
	
	The above implies that given $\delta := \mathring{\delta}_{u_0} \in \mathring{T}^!(F_{u_0}) \simeq T^!(F)$, for all $f$ we have
	\[ \Trans_{\mathbf{M}^!, \tilde{M}}(f_{\tilde{M}})(\delta) = 0 \implies \Trans_{\mathbf{M}^!, \tilde{M}}(\epsilon_{\tilde{M}}(f))(\delta) = 0. \]
	Hence \eqref{eqn:key-geometric-pf-aux-0} holds for these $\delta$. The case of general $\delta \in T^!_{\tilde{M}\text{-reg}}(F)$ follows by density.

	Proceed to the verification of \eqref{eqn:key-geometric-aux-1}. Enlarging $V$ if necessary, we may take two non-Archimedean places $v \neq v'$ in $V \smallsetminus \{u_0, \ldots, u_r\}$ and assume that $f_v$ (resp.\ $f_{v'}$) comes from $C^\infty_{c, \asp}(\tilde{G}_{v, \mathrm{reg}})$ (resp.\ $C^\infty_{c, \asp}(\tilde{G}_{v', \mathrm{reg}})$).
	
	Recall that $h^\varphi_V = \sum_{v_0 \in V} h^\varphi_{v_0}$. Since $\epsilon_{\tilde{M}_{v_0}}(\cdot)$ is cuspidal, every $h^\varphi_{v_0}$ is by construction cuspidal at two places, and supported on regular semisimple locus at another place. All in all, the version 3 of the simple trace formula \cite[Théorème 6.7]{Li14b} for $\widetilde{\mathring{M}}$ yields, in the formalism of \eqref{eqn:Igeom-1}, the equality
	\begin{equation}\label{eqn:h-TF-simple}
		I^{\widetilde{\mathring{M}}}_{\mathrm{geom}}\left( h_V^\varphi \right) = \sum_{\substack{\mathcal{O}^* \subset \mathring{M}_{\mathrm{reg}}(\mathring{F}) \\ / \text{st.\ conj.}} } I^{\widetilde{\mathring{M}}_V}\left( A^{\widetilde{\mathring{M}}}(V, \mathcal{O}^*), h_V^\varphi \right).
	\end{equation}
	
	Next, the stabilization of the coefficients $A^{\widetilde{\mathring{M}}}(V, \mathcal{O}^*)$ (Theorem \ref{prop:matching-coeff-A}) entails
	\begin{equation}\label{eqn:h-STF-simple}
		I^{\widetilde{\mathring{M}}}_{\mathrm{geom}}\left( h_V^\varphi \right) =
		\sum_{\substack{\mathring{\mathbf{M}}^{!!} \in \Endo_{\elli}(\widetilde{\mathring{M}}) \\ \mathcal{O}^{!!, *}: \mathring{G}\text{-reg}}} \iota(\widetilde{\mathring{M}}, \mathring{M}^{!!})
		S^{\mathring{M}^{!!}_V}\left( SA^{\mathring{M}^{!!}}(V, \mathcal{O}^{!!, *}), \Trans_{\mathring{\mathbf{M}}^{!!}, \widetilde{\mathring{M}}} (h^\varphi_V) \right).
	\end{equation}
	
	Claim: those $\mathcal{O}^* \subset \mathring{M}_{\mathrm{reg}}(\mathring{F})$ with
	$I^{\widetilde{\mathring{M}}_V}( A^{\widetilde{\mathring{M}}}(V, \mathcal{O}^*), h_V^\varphi ) \neq 0$
	can be constrained in a finite subset $\Gamma$ depending only on the support of $\varphi$ (viewed as a function $\tilde{G}_V \to \CC$). One can resort to Lemma \ref{prop:geom-diff-h}: it suffices to control those $\mathcal{O}^*$ contributing to the $\mathring{M}$-part in the expansions of $I^{\Endo}_{\mathrm{geom}}(\varphi)$ and $I_{\mathrm{geom}}(\varphi)$ in terms of $\Supp(\varphi)$; note that such an $\mathcal{O}^*$ must cut $K^V$. This can be done by using a support property \cite[Theorem 3.3]{Ar88-2} for $I_{\tilde{M}}(\cdot, \varphi)$, as well as its stable analogues proven in a similar way.
	
	The above is precisely the argument in \cite[p.849]{Ar03-3}. Alternatively, one can argue by the cuspidality of $\epsilon_{\tilde{M}_v}$ as in \cite[p.1239, Remarque]{MW16-2}.

	By the claim, we can further shrink $\Supp(f_v)$ or $\Supp(f_{v'})$ without changing $\Gamma$. In this way, we can arrange that all stable classes $\mathcal{O}^*$ contributing to \eqref{eqn:h-TF-simple} intersect the $W^{\mathring{G}}(\mathring{M})$-orbit of $\mathcal{O}$. Those $\mathcal{O}^{!!, *}$ contributing to \eqref{eqn:h-STF-simple} must map to such a class $\mathcal{O}^*$.

	At each $u_j$, the localization map $\Endo_{\elli}(\widetilde{\mathring{M}}) \rightiso \Endo_{\elli}(\tilde{M})$ is bijective. Our choice of test functions pins down $(\mathring{\mathbf{M}}^{!!}, \mathcal{O}^{!!, *})$: modulo $W^{\mathring{G}}(\mathring{M})$, it must be $(\mathring{\mathbf{M}}^!, \mathcal{O}^!)$. Since $\mathcal{O}$ is $\mathring{G}$-regular and $h_V^\varphi$ is $W^{\mathring{G}}(\mathring{M})$-invariant, this leads to
	\begin{equation*}
		I^{\widetilde{\mathring{M}}}_{\mathrm{geom}}(h_V^\varphi) = |W^{\mathring{G}}(\mathring{M})| \iota(\widetilde{\mathring{M}}, \mathring{M}^! ) S^{\mathring{M}^!_V}\left( SA^{\mathring{M}^!}(V, \mathcal{O}^!), \Trans_{\mathring{\mathbf{M}}^!, \mathring{\tilde{M}}}(h_V^\varphi) \right).
	\end{equation*}
	
	By construction, $V \supset S(\mathcal{O}_1, K)' \supset S(\mathcal{O}_1)$ (see \S\ref{sec:geom-side}), whilst $S(\mathcal{O}_1) \cup V_{\mathrm{ram}} \supset S(\mathcal{O}^!_1)$ by Lemma \ref{prop:admissible-Endo} and $S(\mathcal{O}^!_1) \supset S(\mathcal{O}^!)$. Hence $SA^{\mathring{M}^!}(V, \mathcal{O}^!) = SA^{\mathring{M}^!}(V, \mathcal{O}^!)_{\elli}$. In order to reach \eqref{eqn:key-geometric-aux-1}, it remains to show that
	\begin{equation}\label{eqn:h-I-S}
		S^{\mathring{M}^!_V}\left( SA^{\mathring{M}^!}(V, \mathcal{O}^!)_{\elli}, \cdot\right) = (\text{nonzero constant}) \cdot S^{\mathring{M}^!_V}\left( \mathring{\delta}_V, \cdot \right).
	\end{equation}

	Since $\mathcal{O}^! \subset \mathring{M}^!_{\mathrm{reg}}(\mathring{F})$, one reduces \eqref{eqn:h-I-S} to a standard fact about the elliptic regular part of the stable trace formula for $\mathring{M}^!$; see eg.\ \cite[9.6 Theorem]{Ko86}. One can also deduce \eqref{eqn:h-I-S} easily from the descent formula in Theorem \ref{prop:SAell-descent}.
\end{proof}

The arguments above are similar to \cite[\S 7]{Ar03-3} and \cite[X.7.7]{MW16-2}, but much simpler.

\section{End of stabilization}\label{sec:end-stabilization}
Now that the key geometric Hypothesis \ref{hyp:key-geometric} is established, we can prove all the main theorems in this work, under inductive assumptions.

\begin{enumerate}
	\item \textbf{Stabilization of local trace formula (Theorems \ref{prop:stabilization-geom-LTF}, \ref{prop:stabilization-disc-LTF}).}
	
	This is the content of Corollary \ref{prop:local-geometric-G-reg}.
	
	\item \textbf{The local geometric Theorem \ref{prop:local-geometric} and its semi-local version, Theorem \ref{prop:semilocal-geometric}.}
	
	Theorem \ref{prop:local-geometric-reduction} reduces the local theorem to the $\tilde{G}$-regular case; in turn, the latter is a consequence of the key geometric Hypothesis. The semi-local theorem has been established in \S\ref{sec:semi-local-geom-Endo} by assuming the local one.
	
	\item \textbf{Local non-Archimedean Theorem \ref{prop:rho-main-nonArch} for stabilizing $\rho_J$.}
	
	This is also reduced to the $\tilde{G}$-regular case of the local geometric Theorem: see \S\ref{sec:rho-main-summary}.
	
	\item \textbf{Realization of the local program in \S\ref{sec:ext-Arch-prog} and \S\ref{sec:rho-main-summary}, including the stabilization of $\rho_J$ in the Archimedean case.}
	
	Proofs have been given there, by assuming the $\tilde{G}$-regular case of the local geometric Theorem.
	
	\item \textbf{Global Theorem \ref{prop:geom-stabilization} for stabilizing $I_{\mathrm{geom}}$.}
	
	Apply Corollary \ref{prop:global-matching-reduction} and the available stabilization of the coefficients in $I_{\mathrm{geom}}$ (Theorem \ref{prop:matching-coeff-A}).
	
	\item \textbf{Global Theorem \ref{prop:spectral-stabilization-disc} for stabilizing $I_{\mathrm{disc}, \nu, c^V}$ and $I_{\mathrm{disc}, \nu}$.}
	
	It suffices to show $I_{\mathrm{disc}, \nu} = I^{\Endo}_{\mathrm{disc}, \nu}$ as distributions on $\tilde{G}_V$, since the $c^V$-parts can then be isolated. By Lemma \ref{prop:disc-diff-geom} and Theorem \ref{prop:geom-stabilization} (the previous item),
	\[ \sum_{t \geq 0} \sum_{\substack{\mu \in \mathfrak{h}^*_{\CC}/W_\infty \\ \|\Im(\mu)\| = t}} \left( I^{\Endo}_{\mathrm{disc}, \mu}(f) - I_{\mathrm{disc}, \mu}(f) \right) = I^{\Endo}_{\mathrm{geom}}(f) - I_{\mathrm{geom}}(f) = 0 \]
	for all $f$. Replacing $f$ by $f_\alpha$ with a suitable $\alpha \in C^\infty_c(\mathfrak{h})^{W_\infty}$ and its powers under convolution, we can isolate the $\nu$-part and deduce $I^{\Endo}_{\mathrm{disc}, \nu}(f) = I_{\mathrm{disc}, \nu}(f)$ as in the proof of Proposition \ref{prop:Idisc-stabilization-special}.
	
	\item \textbf{Global Theorem \ref{prop:spectral-stabilization} for stabilizing $I_\nu$ and $I_t$}
	
	It suffices to apply Lemma \ref{prop:diff-I-vs-Idisc} and Theorem \ref{prop:spectral-stabilization-disc} (the previous item).
	
	\item \textbf{The local spectral Theorem \ref{prop:local-spectral} and its semi-local version, Theorem \ref{prop:semilocal-spectral}.}
	
	For the local theorem, note that the premise for Corollary \ref{prop:spectral-trick} is now verified. The semi-local theorem reduces to the local one via splitting formulas as in the geometric setting; see the proof of Proposition \ref{prop:stabilization-global-weighted-character}.
	
	\item \textbf{Local non-Archimedean Theorem \ref{prop:cpt-supported-equalities} for stabilizing ${}^c I_{\tilde{M}}$ and ${}^c \theta_{\tilde{M}}$.}
	
	This follows inductively from the $G$-regular case of the local geometric Theorem \ref{prop:local-geometric}, by virtue of Corollary \ref{prop:cpt-supported-equalities-aux}.
	
	\item \textbf{Local Archimedean Theorem \ref{prop:theta-matching-arch} for stabilizing ${}^c I_{\tilde{M}}$, ${}^c \theta_{\tilde{M}}$ and ${}^{c,r}  \theta_{\tilde{M}}$.}
	
	Similar to the non-Archimedean counterpart, by appealing to Lemma \ref{prop:theta-matching-arch-prep}.
\end{enumerate}

\printindex

% Below: for Biblatex...

\printbibliography[heading=bibintoc]

% Below: for bibtex - preferred on arXiv
%\bibliographystyle{abbrv}
%\bibliography{STF-Mp}

\end{document}